\documentclass[letterpaper,11pt,oneside]{book}
\usepackage[T1]{fontenc}
\usepackage[utf8]{inputenc}
\usepackage{lmodern}
\usepackage{graphicx}
\usepackage[english]{babel}
\usepackage{amsfonts,amsrefs,latexsym,amsmath,amssymb,amsthm,mathrsfs,verbatim}
\usepackage{fancyhdr}
\usepackage[plainpages=false]{hyperref}
\usepackage{url,color}
\usepackage{pifont}
\usepackage{upgreek}
\usepackage{times}
\usepackage{calligra}
\usepackage{marvosym}
\usepackage[percent]{overpic}
\usepackage[hypcap=false]{caption}
\usepackage{xcolor}
\usepackage{wasysym}

\newcommand{\changefont}{%
    \fontsize{12}{11}\selectfont
}
\numberwithin{equation}{section}

\newtheorem{theorem}{Theorem}[section]
\newtheorem{proposition}{Proposition}[section]
\newtheorem{lemma}[proposition]{Lemma}
\newtheorem{corollary}[proposition]{Corollary}

\theoremstyle{definition}
\newtheorem{definition}{Definition}[section]
\newtheorem{remark}{Remark}[section]

\DeclareRobustCommand{\Annoy}{{\breve{\underline{L}}}}

\newcommand{\rgeo}{\varrho}

\newcommand{\myexp}{e}

\newcommand{\Tboot}{T_{(Bootstrap)}}
\newcommand{\Tlocal}{T_{(Local)}}

\newcommand{\mytr}{{\mbox{\upshape{tr}}_{\mkern-2mu \gsphere}}}
\newcommand{\myspacetimetr}{{\mbox{\upshape{tr}}_{\mkern-2mu g}}}

\newcommand{\gsphere}{g \mkern-8.5mu / }
\newcommand{\rescaledmetricsphere}{h \mkern-10mu / }
\newcommand{\gt}{\underline{g}}
\newcommand{\gtinverse}{\underline{g}^{-1}}
\newcommand{\ginversesphere}{\gsphere^{-1}}

\newcommand{\hsphere}{h \mkern-8.5mu / }
\newcommand{\hinversesphere}{\hsphere^{-1}}

\newcommand{\Euct}{e}

\newcommand{\Eucsphereunit}{{e \mkern-8.5mu / }}

\newcommand{\msphere}{m \mkern-10.5mu / }
\newcommand{\minversesphere}{\msphere^{-1}}

\newcommand{\sphereproject}{{\Pi \mkern-12mu / } \, }
\newcommand{\Sigmatproject}{\underline{\Pi}}

\newcommand{\Euctvolform}{\upsilon_{\Euct}}
\newcommand{\gtvolform}{\upsilon_{\gt}}

\newcommand{\vol}{\varpi}
\newcommand{\tvol}{\underline{\varpi}}
\newcommand{\conevol}{\overline{\varpi}}
\newcommand{\spherevol}{\upsilon_{g \mkern-8.5mu /}}
\newcommand{\rescaledspherevol}{\upsilon_{h \mkern-9mu /}}
\newcommand{\Eucspherevol}{\upsilon_{e \mkern-8.5mu /}}
\newcommand{\argspherevol}[1]{\upsilon_{{g \mkern-8.5mu /}#1}}
\newcommand{\rescaledargspherevol}[1]{\upsilon_{{h \mkern-9mu /}#1}}
\newcommand{\argEucspherevol}[1]{\upsilon_{{e \mkern-8.5mu /}#1}}

\newcommand{\sphereGamma}{\Gamma \mkern-10mu / }

\newcommand{\Fried}{\digamma}

\newcommand{\D}{\mathscr{D}}

\newcommand{\angD}{ {\nabla \mkern-14mu / \,} }
\newcommand{\angDarg}[1]{{\angD_{\mkern-3mu #1}}}

\newcommand{\angDsquaredarg}[2]{ {\angD_{\mkern-3mu #1 #2}^2} }
\newcommand{\angfreeDsquared}{ { / \mkern-10mu \hat{\nabla}^2 }}
\newcommand{\angfreeDsquaredarg}[2]{\angfreeDsquared_{\mkern-14mu #1 #2}}

\newcommand{\angcheckD}{ {{\check{\nabla} \mkern-14mu /} \,}}
\newcommand{\angcheckDarg}[1]{ {{\check{\nabla} \mkern-14mu /}_{#1} \,}}

\newcommand{\angdiv}{\mbox{\upshape{div} $\mkern-17mu /$\,}}
\newcommand{\angLap}{ {\Delta \mkern-12mu / \, } }
\newcommand{\angFlatLap}{ {\Delta \mkern-12mu / \, }_{\msphere} }

\newcommand{\angdiff}{ {{d \mkern-9mu /} }}
\newcommand{\angdiffarg}[1]{ {d \mkern-9mu /}_{#1} }
\newcommand{\angdiffuparg}[1]{ {d \mkern-9mu /}^{#1} }

\newcommand{\angLie}{ { \mathcal{L} \mkern-10mu / } }

\newcommand{\angfreeLie}{\hat{\angLie}}
\newcommand{\angfreeLiearg}[1]{{/ \mkern-11.5mu {\hat{\mathcal{L}}_{#1}}} }
\newcommand{\angfreeLietwoarg}[2]{{/ \mkern-11.5mu {\hat{\mathcal{L}}_{#1}^{#2}} }}
\newcommand{\SigmatLie}{\underline{\mathcal{L}}}
\newcommand{\Lie}{\mathcal{L}}

\newcommand{\angJ}{ {\mathscr{J} \mkern-14mu / \, } }
\newcommand{\angnormalJ}{ {J \mkern-10mu / \,} }

\newcommand{\angpi}{ { \pi \mkern-10mu / }}

\newcommand{\angk}{ { {k \mkern-10mu /} \, } }
\newcommand{\angkarg}[1]{ {{k \mkern-10mu /}_{#1} \, } }
\newcommand{\angkdoublearg}[2]{ {{k \mkern-10mu /}_{#1 #2} \, } }
\newcommand{\angkmixedarg}[2]{ {{k \mkern-10mu /}_{#1}^{\ #2} \, } }
\newcommand{\angktriplearg}[3]{ {{k \mkern-10mu /}_{#1 #2}^{#3} \, } }
\newcommand{\angkuparg}[1]{ { {k \mkern-10mu /}^{#1} \, } }
 
\newcommand{\angkfreetriplearg}[3]{ {{\hat{k} \mkern-10mu /}_{#1 #2}^{#3} \, } }

\newcommand{\angG}{ {{G \mkern-12mu /} \, }}
\newcommand{\angGarg}[1]{ {{G \mkern-12mu /}_{\mkern 1mu #1} \, }}

\newcommand{\angGdoublearg}[2]{ {{G \mkern-12mu /}_{#1 #2} \, }}
\newcommand{\angGprime}{ {{ {G'} \mkern-16mu /} \, \, }}
\newcommand{\angGprimearg}[1]{ {{ {G'} \mkern-16mu /}_{\mkern 1mu #1} \, }}
\newcommand{\angGprimemixedarg}[2]{ {{ {G'} \mkern-16mu /}_{\mkern-1mu #1}^{\ \mkern 4mu #2} \, }}
\newcommand{\angGprimedoublearg}[2]{ {{ {G'} \mkern-16mu /}_{\mkern 1mu #1 #2} \, }}
\newcommand{\angGmixedarg}[2]{ {{G \mkern-12mu /}_{#1}^{\ #2} \, }}

\newcommand{\angH}{ {{H \mkern-13mu /} \, }}
\newcommand{\angHarg}[1]{ {{H \mkern-13mu /}_{\mkern 1mu #1} \, }}
\newcommand{\angHuparg}[1]{ {{H \mkern-13mu /}^{#1} \, }}
\newcommand{\angHdoublearg}[2]{ {{H \mkern-13mu /}_{#1 #2} \, }}

\newcommand{\angHmixedarg}[2]{ {{H \mkern-12mu /}_{#1}^{\ #2} \, }}

\newcommand{\angxi}{ { {\xi \mkern-9mu /}  \, }}
\newcommand{\angxiprime}{ { {  {\xi'} \mkern-13mu /}  \, }}
\newcommand{\angxiarg}[1]{ {{\xi \mkern-9mu /}_{#1}  \, }}
\newcommand{\angxiprimearg}[1]{ { {  {\xi'} \mkern-13mu /}_{#1}  \, }}

\newcommand{\deform}[1]{{^{(#1)} \mkern-1mu \pi}}
\newcommand{\deformarg}[3]{{^{(#1)} \mkern-1mu \pi_{#2 #3}}}
\newcommand{\deformuparg}[3]{{^{(#1)} \mkern-1mu \pi^{#2 #3}}}
\newcommand{\deformmixedarg}[3]{{^{(#1)} \mkern-1mu \pi_{#2}^{\ #3}}}

\newcommand{\angdeform}[1]{{^{(#1)} \mkern-2mu \angpi}}
\newcommand{\angdeformoneformarg}[2]{{^{(#1)} \mkern-2mu \angpi_{#2}}}
\newcommand{\angdeformoneformupsharparg}[2]{{^{(#1)} \mkern-2mu {\pi \mkern-10mu /}_{#2}^{\#}}}
\newcommand{\angdeformarg}[3]{{^{(#1)} \mkern-2mu \angpi_{#2 #3}}}
\newcommand{\angdeformuparg}[3]{{^{(#1)} \mkern-2mu{\pi\mkern-10mu /}^{#2 #3}}}
\newcommand{\angdeformmixedarg}[3]{{^{(#1)} \mkern-2mu {\pi \mkern-10mu /}_{#2}^{\ #3}}}
\newcommand{\angdeformfree}[1]{{^{(#1)} \mkern-2mu \hat{\angpi}}}
\newcommand{\angdeformupsharparg}[1]{{^{(#1)} \mkern-2mu {\pi \mkern-10mu /}^{\#}}}
\newcommand{\angdeformupdoublesharparg}[1]{{^{(#1)} \mkern-2mu {\pi \mkern-10mu /}^{\# \#}}}
\newcommand{\angdeformfreearg}[3]{{^{(#1)} \mkern-2mu {\hat{\pi}\mkern-10mu /}_{#2 #3}}}

\newcommand{\angdeformfreeupdoublesharparg}[1]{{^{(#1)} \mkern-2mu {\hat{\pi} \mkern-10mu /}^{\# \#}}}
\newcommand{\angdeformfreeuparg}[3]{{^{(#1)} \mkern-2mu {\hat{\pi}\mkern-10mu /}^{#2 #3}}}
\newcommand{\angdeformfreemixedarg}[3]{{^{(#1)} \mkern-2mu {\hat{\pi}\mkern-10mu /}_{#2}^{\ #3}}}

\newcommand{\Lgeo}{L_{(Geo)}}
\newcommand{\Lunit}{L}
\newcommand{\uLgood}{\breve{\underline{L}}}
\newcommand{\uLunit}{\underline{L}}
\newcommand{\Rad}{\breve{R}}
\newcommand{\Radunit}{R}

\newcommand{\Timenormal}{N}

\newcommand{\Mult}{T}
\newcommand{\Mor}{\widetilde{K}}

\newcommand{\Roteuc}{\Rot_{(Flat)}}
\newcommand{\Roteucarg}[1]{\Rot_{(Flat;#1)}}
\newcommand{\Rot}{O}
\newcommand{\RotRadcomponent}[1]{\uprho_{(#1)}}

\newcommand{\FutFailFac}{ {^{(+)} \mkern-1mu \aleph} }
\newcommand{\PastFailFac}{ {^{(-)} \mkern-1mu \aleph} }
\newcommand{\InitialFutFailFac}{ {^{(+)} \mkern-1mu \mathring{\aleph}} }

\newcommand{\Jcurrent}[1]{^{(#1)} \mkern-10mu \mathscr{J}}
\newcommand{\Jenergycurrent}[1]{^{(#1)} \mkern-3mu J}

\newcommand{\Cur}{\mathscr{R}}

\newcommand{\Rsphere}{\mathfrak{R}}
\newcommand{\Riemannsphere}{\mathfrak{R}}
\newcommand{\Gauss}{\mathfrak{K}}

\newcommand{\enmomtensor}{Q}

\newcommand{\enzero}{\mathbb{E}}
\newcommand{\enone}{\widetilde{\mathbb{E}}}
\newcommand{\flzero}{\mathbb{F}}
\newcommand{\flone}{\widetilde{\mathbb{F}}}

\newcommand{\Morint}{\widetilde{\mathbb{K}}}

\newcommand{\totzeromax}[1]{\mathbb{Q}_{(#1)}}
\newcommand{\totonemax}[1]{\widetilde{\mathbb{Q}}_{(#1)}}
\newcommand{\totMormax}[1]{\widetilde{\mathbb{K}}_{(#1)}}

\newcommand{\chifullmod}{\mathscr{X}}
\newcommand{\chifullmodarg}[1]{{^{(#1)} \mkern-4mu \mathscr{X}}}
\newcommand{\chipartialmod}{\widetilde{\mathscr{X}}}
\newcommand{\chipartialmodarg}[1]{{^{(#1)} \mkern-4mu \widetilde{\mathscr{X}}}}
\newcommand{\chipartialmodinhom}{\widetilde{\mathfrak{X}}}
\newcommand{\chipartialmodinhomarg}[1]{{^{(#1)} \mkern-2mu \widetilde{\mathfrak{X}}}}
\newcommand{\mupartialmod}{\widetilde{\mathscr{M}}}
\newcommand{\mupartialmodarg}[1]{{^{(#1)} \mkern-4mu \widetilde{\mathscr{M}} \, } }
\newcommand{\mupartialmodargsharp}[1]{{^{(#1)} \mkern-4mu \widetilde{\mathscr{M}}^{\#} \, }}
\newcommand{\mupartialmodinhom}{\widetilde{\mathfrak{M}}}
\newcommand{\mupartialmodinhomarg}[1]{{^{(#1)} \mkern-1mu \widetilde{\mathfrak{M}} \, }}
\newcommand{\mupartialmodinhomargsharp}[1]{{^{(#1)} \mkern-1mu \widetilde{\mathfrak{M}}^{\#} \, }}

\newcommand{\waveinhom}{\mathfrak{F}}

\newcommand{\inhomarg}[1]{{^{(#1)} \mkern-2mu \waveinhom}}
\newcommand{\inhomleftexparg}[2]{{^{(#1)} \mkern-.5mu #2 }}

\newcommand{\chifullmodsourcearg}[1]{\inhomleftexparg{#1}{\mathfrak{I}}}

\newcommand{\chipartialmodsourcearg}[1]{\inhomleftexparg{#1}{\mathfrak{B}}}
\newcommand{\mupartialmodsourcearg}[1]{\inhomleftexparg{#1}{\mathfrak{J}}}

\newcommand{\basicenergyerror}[1]{{^{(#1)} \mkern-.5mu \mathfrak{P}}}
\newcommand{\basicenergyerrorarg}[2]{{^{(#1)} \mkern-.5mu \mathfrak{P}}_{(#2)}}

\newcommand{\rotationcommutationerror}[3]{{^{(#2,#3)} \mkern-.5mu #1}}
\newcommand{\rotationcommutationerrorarg}[4]{{^{(#3,#4)} \mkern-.5mu #1_{#2}}}
\newcommand{\rotationcommutationerrorupsharp}[3]{{^{(#2,#3)} \mkern-.5mu #1^{\#}}}

\newcommand{\IBPerror}[1]{\mbox{\upshape{Error}}_{(#1)}}
\newcommand{\IBPseconderror}[1]{\widetilde{\mbox{\upshape{Error}}}_{(#1)}}

\newcommand{\Vplus}[2]{{^{(+)} \mkern-.5mu  \mathcal{V}_{#1}^{#2}}}
\newcommand{\Vminus}[2]{{^{(-)} \mkern-.5mu  \mathcal{V}_{#1}^{#2}}}
\newcommand{\Sigmaplus}[3]{{^{(+)} \mkern-.5mu  \Sigma_{#1;#2}^{#3}}}
\newcommand{\Sigmaminus}[3]{{^{(-)} \mkern-.5mu  \Sigma_{#1;#2}^{#3}}}

\newcommand{\Conone}{A}
\newcommand{\Cononestar}{A_*}
\newcommand{\Littleconone}{a}
\newcommand{\Contwo}{B}

\newcommand{\Littlecontwo}{b}
\newcommand{\Conthree}{P}

\newcommand{\smoothfunction}{f}
\newcommand{\newsmoothfunction}{\eta}

\newcommand{\angDleftexp}[1]{{^{(#1)} \mkern-1mu \angD}}
\newcommand{\angDleftexparg}[2]{{^{(#1)} \mkern-1mu \angD_{#2}}}
\newcommand{\xileftexp}[1]{{{^{(#1)} \mkern-1mu \xi}}}
\newcommand{\xileftexparg}[3]{{^{(#1)} \mkern-1mu \xi_{#2 #3}}}
\newcommand{\gsphereleftexp}[1]{{{^{(#1)} \mkern-1mu \gsphere}}}
\newcommand{\Scalarsphereleftexp}[1]{{{^{(#1)} \Rsphere}}}
\newcommand{\gsphereleftexparg}[3]{{^{(#1)} \mkern-1mu \gsphere_{#2 #3}}}
\newcommand{\ginversesphereleftexp}[1]{{^{(#1)} \mkern-1mu \ginversesphere}}
\newcommand{\ginversesphereleftexparg}[3]{{(^{(#1)} \mkern-1mu \ginversesphere)^{#2 #3}}}
\newcommand{\sphereGammaleftexp}[1]{{^{(#1)} {\Gamma \mkern-10mu / } }}
\newcommand{\sphereGammaleftexparg}[4]{{^{(#1)} {\Gamma \mkern-10mu / }_{#2 \ #3}^{\ #4}  }}

\newcommand{\topboxedmulterrorone}[1]{\boxed{\mathbf{I)}}_{(#1)}}
\newcommand{\topmulterrortwo}[1]{\mathbf{II)}_{(#1)}}
\newcommand{\topmulterrorthree}[1]{\mathbf{III)}_{(#1)}}
\newcommand{\topmulterrorfour}[1]{\mathbf{IV)}_{(#1)}}
\newcommand{\topboxedmulterrorfive}[1]{\boxed{\mathbf{V)}}_{(#1)}}
\newcommand{\topmulterrorsix}[1]{\mathbf{VI)}_{(#1)}}
\newcommand{\topmulterrorseven}[1]{\mathbf{VII)}_{(#1)}}
\newcommand{\topmulterroreight}[1]{\mathbf{VIII)}_{(#1)}}

\newcommand{\topboxedmorerrorone}[1]{\boxed{\mathbf{\widetilde{I})}}_{(#1)}}
\newcommand{\topmorerrortwo}[1]{\mathbf{\widetilde{II})}_{(#1)}}
\newcommand{\topmorerrorthree}[1]{\mathbf{\widetilde{III})}_{(#1)}}
\newcommand{\topmorerrorfour}[1]{\mathbf{\widetilde{IV})}_{(#1)}}
\newcommand{\topboxedmorerrorfive}[1]{\boxed{\mathbf{\widetilde{V})}}_{(#1)}}
\newcommand{\topmorerrorsix}[1]{\mathbf{\widetilde{VI})}_{(#1)}}
\newcommand{\topmorerrorseven}[1]{\mathbf{\widetilde{VII})}_{(#1)}}
\newcommand{\topmorerroreight}[1]{\mathbf{\widetilde{VIII})}_{(#1)}}

\newcommand{\multerrorzero}[1]{\mathbf{0)}_{(#1)}}
\newcommand{\multerrorone}[1]{\mathbf{i)}_{(#1)}}
\newcommand{\multerrortwo}[1]{\mathbf{ii)}_{(#1)}}
\newcommand{\multerrorthree}[1]{\mathbf{iii)}_{(#1)}}
\newcommand{\multerrorfour}[1]{\mathbf{iv)}_{(#1)}}
\newcommand{\multerrorfive}[1]{\mathbf{v)}_{(#1)}}
\newcommand{\multerrorsix}[1]{\mathbf{vi)}_{(#1)}}

\newcommand{\morerrorzero}[1]{\mathbf{\widetilde{0})}_{(#1)}}
\newcommand{\morerrorone}[1]{\mathbf{\widetilde{i})}_{(#1)}}
\newcommand{\morerrortwo}[1]{\mathbf{\widetilde{ii})}_{(#1)}}
\newcommand{\morerrorthree}[1]{\mathbf{\widetilde{iii})}_{(#1)}}
\newcommand{\morerrorfour}[1]{\mathbf{\widetilde{iv})}_{(#1)}}
\newcommand{\morerrorfive}[1]{\mathbf{\widetilde{v})}_{(#1)}}

\newcommand{\shockfunction}{\mathbb{S}}
\newcommand{\timeminushalfshockfunction}{\check{\mathbb{S}}}

\newcommand{\halfdata}{\check{\Psi}}
\newcommand{\halftimedata}{\check{\Psi}_0}

\newcommand{\shockparameter}{\upalpha_*}

\newcommand{\scaledPsi}[1]{{^{(#1)} \mkern-2mu \Psi}}
\newcommand{\scaleddata}[1]{{^{(#1)} \mkern-2mu \mathring{\Psi}}}
\newcommand{\scaledtimedata}[1]{{^{(#1)} \mkern-2mu \mathring{\Psi}_0}}

\newcommand{\shockset}{\mathfrak{S}}

\newcommand{\myarray}[2][]{\left(
		\begin{array}{lr}
    	 #1 \\
    	 #2 
     \end{array} \right)}

\newcommand{\threemyarray}[3][]{\left(
		\begin{array}{lr}
    	 #1 \\
    	 #2 \\
    	 #3
     \end{array} \right)}

\newcommand{\fourmyarray}[4][]{\left(
		\begin{array}{lr}
    	 #1 \\
    	 #2 \\
    	 #3 \\
    	 #4
     \end{array} \right)}

\newcommand{\fivemyarray}[5][]{\left(
		\begin{array}{lr}
    	 #1 \\
    	 #2 \\
    	 #3 \\
    	 #4 \\
    	 #5
     \end{array} \right)}

\newcommand{\sevenmyarray}[7][]{\left(
		\begin{array}{lr}
    	 #1 \\
    	 #2 \\
    	 #3 \\
    	 #4 \\
    	 #5 \\
    	 #6 \\
    	 #7
     \end{array} \right)}

\newcommand{\ninemyarray}[9][]{\left(
		\begin{array}{lr}
    	 #1 \\
    	 #2 \\
    	 #3 \\
    	 #4 \\
    	 #5 \\
    	 #6 \\
    	 #7 \\
    	 #8 \\
    	 #9
     \end{array} \right)}

\newenvironment{dedication}
{
   \cleardoublepage
   \thispagestyle{empty}
   \vspace*{\stretch{1}}
   \hfill\begin{minipage}[t]{0.66\textwidth}
   \raggedright
}
{
   \end{minipage}
   \vspace*{\stretch{3}}
\texttt{}   \clearpage
}

\makeatletter
\renewcommand{\@chapapp}{}

\makeatother

\title{{\fontsize{.5in}{.6in}\selectfont  
\textbf{Shock Formation in Small-Data Solutions to $3D$ Quasilinear Wave Equations}}  
 \\ 
}
\author{\Large{\textsc{Jared Speck}}\thanks{$^{*}$Massachusetts Institute of Technology, Department of Mathematics. 
\texttt{jspeck@math.mit.edu}}}

\begin{document}

\frontmatter
\pagenumbering{alph}
\maketitle
\pagenumbering{roman}

\pagestyle{fancy}


\begin{dedication}
To Genevieve, for sleeping (just) long enough to allow me to complete this project.
\end{dedication}


\thispagestyle{fancy}
\tableofcontents
\thispagestyle{fancy}

\listoffigures
\thispagestyle{fancy}

\mainmatter

\chapter*{Preface}
\thispagestyle{fancy}
In his 2007 monograph, D. Christodoulou proved a breakthrough result giving a 
detailed description of the formation of shocks in solutions to the relativistic Euler equations
in three spatial dimensions. He assumed that
the data have small $H^N$ norm, where $N$ is a sufficiently large integer.
To deduce the shock formation, he also assumed that the data
verify a signed integral inequality. 
In the present monograph, we extend Christodoulou's framework and use it to prove 
that shock singularities often develop 
in initially small, regular
solutions to two important classes of quasilinear wave equations in three spatial dimensions.
Our work also generalizes and unifies earlier work on singularity formation initiated by F. John
in the 1970's and continued by
L. H\"{o}rmander,
S. Alinhac, 
and many others.
Specifically, we study \textbf{i)} covariant scalar wave equations of the form 
$\square_{g(\Psi)} \Psi = 0$ and \textbf{ii)} non-covariant scalar wave equations
of the form $(h^{-1})^{\alpha \beta}(\partial \Phi) \partial_{\alpha} \partial_{\beta} \Phi = 0.$
Our main result shows that whenever the nonlinear terms fail Klainerman's classic null condition,
shocks develop in solutions arising from an open set of small data.
Hence, within the classes \textbf{i)} and \textbf{ii)},
our work can be viewed as a sharp converse to the well-known result of
Christodoulou and Klainerman, which showed
that when the classic null condition is verified, 
small-data global existence holds.

\bigskip

\noindent \textbf{Mathematics Subject Classification (2010)} Primary: 35L67; Secondary: 35L05, 35L10, 35L72

\section*{Acknowledgements}
I am grateful for the support offered by NSF grant \# DMS-1162211 
and by a Solomon Buchsbaum grant administered by the Massachusetts Institute of Technology.
I would like to thank the American Institute of Mathematics for funding three SQuaRE workshops
on the formation of shocks, which greatly aided the development of many of the ideas
in this monograph. I would like to thank Gustav Holzegel, Sergiu Klainerman, Jonathan Luk, 
Willie Wong, and Shiwu Yang for participating in the workshops and for their helpful 
contributions. I offer special thanks to Willie Wong for creating Figure~\ref{F:MAXDEVBOUND}.
I am grateful for the support of my PhD advisors, Michael Kiessling and Shadi
Tahvildar-Zadeh, who encouraged me to read Christodoulou's monograph.
I would also like to thank Hans Lindblad for sharing his insight on Alinhac's work.

\mbox{}\\
\noindent Jared Speck \\
\noindent \url{http://math.mit.edu/~jspeck/}



\chapter{Introduction and Overview of the Two Main Theorems}
\label{C:INTRO}
\thispagestyle{fancy}
In \cite{dC2007}, D. Christodoulou proved a breakthrough result giving a detailed
description of the formation of shocks in small-data solutions to the relativistic Euler equations
in three spatial dimensions.
In this monograph, we develop an extended version of his framework 
that enables us to derive sharp information about the global behavior of solutions
to two important classes of quasilinear wave equations (described just below)
in three spatial dimensions. Our results apply to data belonging to $H^N$ for a sufficiently large integer $N$
(see \eqref{E:INTROSMALLDATA})
and they show that for equations belonging to the two classes,
initially small, regular solutions often develop 
(see Remark~\ref{R:SHOCKFORMATIONFORALARGESETOFSMALLDATA})
shock singularities\footnote{We must make suitable assumptions on the nonlinear
terms to ensure that shocks form.} 
in finite time. 
We assume that the equations are quadratic perturbations of the standard linear
wave equation $\square_m \Psi = 0$ on Minkowski spacetime $(\mathbb{R}^{1+3},m),$
where $m$ is the Minkowski metric and $\square_m = - \partial_t^2 + \Delta$
relative to standard rectangular coordinates.
The two classes of wave equations that we treat are 
\textbf{i)} covariant wave equations of the form
\[
	\square_{g(\Psi)} \Psi = 0,
\]
where $\square_{g(\Psi)}$ is the covariant wave operator of the Lorentzian metric $g(\Psi),$
and \textbf{ii)}
non-covariant wave equations\footnote{Throughout, we use Einstein's summation convention.}
of the form
\[
	(h^{-1})^{\alpha \beta}(\partial \Phi) \partial_{\alpha} \partial_{\beta} \Phi = 0,
\]
where $h(\partial \Phi)$ is a Lorentzian metric and $\partial_{\alpha}$
denotes a Minkowski rectangular coordinate derivative.
We show below that these two classes of equations are intimately related.
We restrict our attention to initial data that are 
compactly supported in the Euclidean unit ball, 
and we study the solution only in the region
determined by the portion of the data lying in the exterior of a Euclidean sphere 
of radius $1 - U_0$ centered at the origin in $\mathbb{R}^3$ 
(see Figure~\ref{F:REGION} on pg. \pageref{F:REGION}),
where $0 < U_0 < 1$ is a parameter.
We expect that our work could be extended to apply to non-compactly supported data
and to a larger solution region. 
An extensive account of related work as well as an overview of our proofs are  
provided in the companion survey article \cite{gHsKjSwW2014}. 
Hence, in the present introduction, we provide a relatively sparse
account of the relevant literature and a briefer overview our proofs.

Fritz John gave a non-constructive proof \cite{fJ1981}  
showing that a large class of equations 
that are quadratic perturbations of the standard linear
wave equation on $\mathbb{R}^{1+3}$
must fail to have global $C^3$ solutions, 
as long as the data are smooth, compactly supported, and nontrivial
(without any smallness assumption).
His proof applied to some semilinear equations including
\[
	\square_m \Phi = - (\partial_t \Phi)^2,
\]
to some quasilinear equations including
\[
	\square_m \Phi = - (\partial_t \Phi) \partial_t^2 \Phi,
\]
and to related equations with nonlinearities verifying a signed technical condition.
Later work by Alinhac\footnote{Despite the title, Alinhac's article addresses both the cases of two and three spatial dimensions} 
\cite{sA2001b} and Christodoulou \cite{dC2007} 
(see also the follow-up work Christodoulou-Miao \cite{dCsM2012})
revealed that for some related classes of 
quasilinear wave equations having a nontrivial intersection with John's class, 
for solutions launched by a set of small data, the singularity 
is caused by the intersection of characteristic hypersurfaces. That is,
the singularity is a shock.
We discuss the precise classes of equations and data covered by Alinhac and Christodoulou and
the differences in their approaches in Sect.~\ref{S:COMPARISONWITHRELATEDWORK}.
Although the sign of the nonlinearity played an important role in John's proof, 
it did not in the proofs of Alinhac or Christodoulou, nor does it in the present work.
It is interesting to note that as of the present, 
the nature of the breakdown 
of solutions in the semilinear equations covered by John's work \cite{fJ1981} 
remains poorly understood. 
Our main objectives in the monograph are the following.
 \renewcommand{\theenumi}{\textbf{\Roman{enumi}}}
\begin{enumerate}
	\item We identify criteria for the structure of the equations and 
		\emph{stable} criteria on small initial data
		that guarantee that the solution develops a shock singularity in finite time.
		In particular, our results recover, as special cases, 
		the most important aspects
		of the shock-formation results of Christodoulou and Alinhac.
		Specifically, for the equations we study,
		the global behavior of small-data solutions
		is determined by the presence or absence of 
		certain tensorial components in the wave equations.
		If the tensorial components are absent,
		then the nonlinearities verify Klainerman's classic null condition \cite{sk1984}, and the methods of
		\cite{sK1986} and \cite{dC1986a} yield small-data global existence.
		If they are present, then the classic null condition \cite{sk1984} fails, and our 
		work yields shock formation for solutions launched by an open set of small data
		(see Remark~\ref{R:SHOCKFORMATIONFORALARGESETOFSMALLDATA}).
		Moreover, we show that for general small data,  
		shocks are the only kind of singularities that can develop
		in the region of interest.
	\item For the shock-forming solutions, we provide a detailed description of the dynamics from 
		$t=0$ all the way up to and including the constant-time
		hypersurface subset $\Sigma_{T_{(Lifespan);U_0}}^{U_0}$ 
		(see definition \eqref{E:SIGMATU} and \eqref{E:TLIFESPAN})
		where the first shock singularity point lies.
		To provide this description, we extend the groundbreaking framework
		developed by Christodoulou \cite{dC2007},
		which he used to prove sharp shock formation results for irrotational regions of relativistic fluids.
		We have also significantly simplified some
		aspects of his approach, in part by being more selective in 
		our use of geometry. That is, we perform fully geometric 
		decompositions only when deriving the most delicate estimates. 
		In deriving many of the less delicate estimates, we use
		a simpler, less geometric approach, 
		which has spared us a great deal of analysis.
		
		As does \cite{dC2007}, our results provide,
		in particular, a sharp description 
		of the quantities that remain regular up to and including the 
		time of first shock formation.
		An important virtue of our estimates is that
		as in \cite{dC2007}*{Ch.15},
		they could be extended to give a detailed description of the maximal future development of the 
		portion of the data in the exterior of $S_{0,U_0} \subset \Sigma_0;$
		see Figure~\ref{F:MAXDEVBOUND} on pg. \pageref{F:MAXDEVBOUND}.
		This is one key advantage of Christodoulou's approach
		over Alinhac's, which is tailored only to see the first singularity point.\footnote{Hence, Alinhac's results only apply
		to data for which there is a unique first singularity point.}
		One would need a precise description of the solution up to the boundary of the maximal development,
		not just information near the first singularity point,
		if one wanted to try to weakly continue the solution under suitable extension criteria.
		
\end{enumerate}

\section{First description of the two theorems} 
\label{S:FIRSTDESCRIPTIONOFTHETHEOREMS}
Our work is divided into two main theorems.

\renewcommand{\theenumi}{\textbf{\arabic{enumi}}}
\begin{enumerate}
	\item \noindent{\textbf{Main Theorem of the Monograph: A Sharp Classical Lifespan Theorem.}}
	In Theorem~\ref{T:LONGTIMEPLUSESTIMATES}, 
	we show that
	for small-data solutions to covariant wave equations of the form $\square_{g(\Psi)} \Psi = 0,$
	there is a scalar function $\upmu$ that dictates the global behavior of the solution in the region of interest.
	The function $\upmu,$ 
	which we discuss in great detail below,
	has a simple geometric interpretation:
	$1/\upmu$ measures the foliation density of a family of outgoing null cones 
	(which are characteristic hypersurfaces of the dynamic metric).
	In short, $\upmu$ starts out near $1$ and either 
	it goes to $0$ 
	in finite time and a shock forms due to the intersection of characteristic hypersurfaces
	(that is, due to the blow-up of the foliation density),
	or it remains positive and, because of
	dispersive effects, no singularity forms in the region.
	Furthermore, the theorem provides many quantitative
	estimates that are verified up to and including
	the constant-time hypersurface subset $\Sigma_{T_{(Lifespan);U_0}}^{U_0}$ 
	(see definition \eqref{E:SIGMATU} and \eqref{E:TLIFESPAN})
	where the first shock singularity point lies.
	The proof of Theorem~\ref{T:LONGTIMEPLUSESTIMATES}
	is based on a long bootstrap argument that occupies the majority of the monograph.
	Our proof is based on Christodoulou's framework \cite{dC2007}, but we
	develop some alternate strategies that significantly shorten and simplify some
	aspects of the analysis and allow us to address a larger class of equations.
	In Appendix \eqref{A:EQUIVALENTPROBLEM}, we outline how to extend
	Theorem~\ref{T:LONGTIMEPLUSESTIMATES} so that it applies to non-covariant wave
	equations of the form $(h^{-1})^{\alpha \beta}(\partial \Phi) \partial_{\alpha} \partial_{\beta} \Phi = 0.$
	
	\begin{remark}[$\upmu^{-1}$ \textbf{degeneracy is responsible for the length}] \label{R:WHYSOLONG}
	The main reasons that the present work and the work
	\cite{dC2007} are so long
	are \textbf{i)} the high-order $L^2$ estimates
	of Theorem~\ref{T:LONGTIMEPLUSESTIMATES}
	are rather degenerate
	with respect to powers of $\upmu^{-1}$
	(see Sect.~\ref{SS:L2HIERARCHY} for an overview)
	and \textbf{ii)} an intimately related difficulty: the top-order $L^2$ estimates
  are very difficult to derive without ``losing derivatives''
	(see Sect.~\ref{SS:EIKONALTOPORDER} for an overview).
	The degeneracy of the high-order $L^2$ estimates with respect to powers of $\upmu^{-1}$
	\textbf{is the main new feature} that 
	distinguishes the problem of small-data shock formation 
	from other global results for nonlinear wave equations
	that are found in the literature. 
	
\end{remark}

\begin{remark}[\textbf{Key innovations of Christodoulou's framework}]
	The most important innovations of the framework of \cite{dC2007}
	are that it provides a means to establish the degenerate high-order $L^2$ estimates
	as well as a means to show that the degeneracy does not
	propagate down to the lower levels.
\end{remark}
	
		\begin{remark}[\textbf{Strong null condition}]
		Our work can easily be extended to allow for the presence of additional
		``harmless'' quadratic semilinear terms in the wave equations.
		For example, a sufficient condition in the region $\lbrace t \geq 0 \rbrace$
		is that the terms verify the 
		future strong null condition of Def.~\ref{D:STRONGNULL}.
		Such terms remain small throughout the future evolution
		and do not interfere with the shock-formation processes.
		In contrast, the results of our theorems are not stable under the addition of
		general cubic semilinear terms. Roughly, the reason is that
		our framework is tailored
		to precisely control the dangerous quadratic interactions,
		whereas for solutions that blow-up, general cubic terms would dominate
		the dynamics near the singularity.
	\end{remark}
\item \noindent{\textbf{Second Theorem: Shock-Formation for Nearly Spherically Symmetric Small Data.}}
		We state this theorem as Theorem~\ref{T:STABLESHOCKFORMATION}.
		It is relatively easy to prove
		thanks to the difficult estimates of Theorem~\ref{T:LONGTIMEPLUSESTIMATES}.
		The theorem shows that there exists an open set of small data that launch solutions
		such that $\upmu$ vanishes in finite time, thus yielding the onset of a shock.
		Specifically, Theorem~\ref{T:STABLESHOCKFORMATION} shows that shocks form in
		solutions corresponding to small, nontrivial, compactly supported, 
		nearly spherically symmetric data. For technical reasons, the theorem
		applies to data specified at time $-1/2$
		and supported in a Euclidean ball of radius\footnote{
		By making straightforward
		modifications throughout the monograph, this result can be extended to allow
		for an arbitrary finite radius of support.} 
		$1/2.$ Moreover, with some additional effort, 
		we could extend Theorem~\ref{T:STABLESHOCKFORMATION} 
		to show shock formation in solutions 
		corresponding to a significantly larger class of data;
		see Remark~\ref{R:SHOCKFORMATIONFORALARGESETOFSMALLDATA}.
\end{enumerate}

\begin{remark}[\textbf{Shock formation for small compactly supported data}]
	\label{R:SHOCKFORMATIONFORALARGESETOFSMALLDATA}
		In Sect.~\ref{SS:INTROSHOCKFORMINGDATACOMPARISON}, 
		we outline how to extend Theorem~\ref{T:STABLESHOCKFORMATION} 
		to show finite-time shock formation in solutions corresponding to
		all sufficiently small compactly supported nontrivial data.
		More precisely, the discussion in Sect.~\ref{SS:INTROSHOCKFORMINGDATACOMPARISON}
		leaves open the possibility that the small size of the data that is
		sufficient for proving the sharp classical lifespan theorem (Theorem~\ref{T:LONGTIMEPLUSESTIMATES})
		might not be sufficient for proving shock formation.
		That is, our proof of shock formation depends on the profile of the data,
		which we might need to multiply by a small rescaling factor 
		in order to ensure the blow-up.
\end{remark}

\section{The basic structure of the equations}
\label{S:INTRORECTANGULARCOORDINATES}
As we mentioned at the beginning, the first class of problems that we study 
is Cauchy problems for covariant scalar wave equations of the form
\begin{align} \label{E:WAVEGEO}
	\square_g \Psi & = 0,
		\\
	(\Psi|_{t=0}, \partial_t \Psi|_{t=0}) & = (\mathring{\Psi},\mathring{\Psi}_0).
		\label{E:INTRODATA}
\end{align}
Above, $g = g(\Psi)$ is a Lorentzian metric, 
\begin{align} \label{E:INTROCOVARIANTWAVEOPERATOR}
	\square_g := (g^{-1})^{\alpha \beta} \D_{\alpha \beta}^2
\end{align}
is the covariant wave operator\footnote{Relative to an arbitrary coordinate system,
$\square_g \Psi = \frac{1}{\sqrt{|\mbox{\upshape{det}} \mbox{$g$}}|} \partial_{\alpha}(\sqrt{|\mbox{\upshape{det}} g|}(g^{-1})^{\alpha \beta} \partial_{\beta} \Psi).$} corresponding to $g,$ and 
$\D$ is the Levi-Civita connection corresponding to $g.$ 
In order to make more precise statements, we now formally introduce
one of the two coordinate systems that we use in our analysis. 
Specifically, throughout this monograph, 
$\lbrace x^{\mu} \rbrace_{\mu = 0,1,2,3}$ denotes a fixed rectangular coordinate system
on Minkowski spacetime $(\mathbb{R}^{1 + 3},m)$ relative to which 
$m_{\mu \nu} = \mbox{diag}(-1,1,1,1).$ $x^0$ is the time coordinate and $(x^1,x^2,x^3)$ are the spatial coordinates.
Greek ``spacetime'' indices vary from $0$ to $3$ and lowercase Latin ``spatial'' indices vary from $1$ to $3.$
Repeated indices are summed over their respective ranges. 
The symbol $\partial_{\nu}$ denotes the rectangular coordinate vectorfield $\frac{\partial}{\partial x^{\nu}}.$
We often use the alternate notation $t := x^0,$ $\partial_t := \frac{\partial}{\partial x^0}.$
Upon rescaling the metric,\footnote{Actually, rescaling the metric leads to the presence of an additional
semilinear term on the right-hand side of equation \eqref{E:WAVEGEO}. However, this additional term 
has a very good null structure, in the sense of Lemma~\ref{L:SPECIALNULLSTRUCTUREINHOMOGENEOUS}.
In the region $\lbrace t \geq 0 \rbrace,$ 
we call this good structure the ``future strong null condition''
(see Def.~\ref{D:STRONGNULL}).
As will become clear, such a good term makes only a tiny contribution
to the future dynamics of small-data solutions (even those solutions that form shocks) and hence we ignore it.} 
we may assume that
\begin{align} \label{E:GINVERSE00ISMINUSONE}
	(g^{-1})^{00}(\Psi) & = - 1,
\end{align}
which simplifies some of our calculations. 

We assume that relative to standard Minkowski-rectangular coordinates,
$g$ is a quadratic perturbation of the Minkowski metric $m,$ 
$(\mu,\nu = 0,1,2,3):$
\begin{align} \label{E:LITTLEGDECOMPOSED}
	g_{\mu \nu} 
	= g_{\mu \nu}(\Psi)
	& := m_{\mu \nu} 
		+ g_{\mu \nu}^{(Small)}(\Psi),
\end{align}
where 
\begin{align} \label{E:MINKOWSKIRELATIVETORECTANGULAR}
	m_{\mu \nu} := \mbox{diag}(-1,1,1,1),
\end{align}
and the $g_{\mu \nu}^{(Small)}$ are smooth functions of $\Psi$ 
(at least when $\Psi$ is sufficiently small)
verifying 
\begin{align} \label{E:METRICNONLINEARITYVANISHESWHENPSIIS0}
	g_{\mu \nu}^{(Small)}(0) = 0.
\end{align}
We prove our main theorems under the assumption that the size 
$\mathring{\upepsilon}$ 
of the data
is sufficiently small, where
\begin{align} \label{E:INTROSMALLDATA}
	\mathring{\upepsilon} 
	= \mathring{\upepsilon}[(\mathring{\Psi},\mathring{\Psi}_0)]
	:= \| \mathring{\Psi} \|_{H_{\Euct}^{25}(\Sigma_0^1)}
		+ \| \mathring{\Psi}_0 \|_{H_{\Euct}^{24}(\Sigma_0^1)}.
\end{align}
In \eqref{E:INTROSMALLDATA},
$\| \cdot \|_{H_{\Euct}^M(\Sigma_0^1)}$
denotes the standard Euclidean Sobolev norm 
corresponding to order $\leq M$ rectangular spatial derivatives
on the Euclidean unit ball $\Sigma_0^1 \subset \mathbb{R}^3.$
The number of derivatives in \eqref{E:INTROSMALLDATA}, $25,$
could be reduced with some additional effort.
We were motivated to explicitly keep track of the number of derivatives
that we use in proving our results
because, as we will see, the number is connected to certain crucially 
important structural features
of our equations and our estimates.

\begin{remark}[$N_{(Christodoulou)}$]
\label{R:CHRISTODOULOUSN}
In \cite{dC2007}, Christodoulou did not 
provide an explicit estimate for the number of derivatives $N_{(Christodoulou)}$ on the data needed 
for his main results. In his bootstrap argument, $N_{(Christodoulou)}^{-1}$ was a parameter that was
chosen to be sufficiently small to counter other large, non-explicit constants.
\end{remark}

For simplicity, we do not study the solution along complete 
constant-time hypersurfaces. Rather, we fix a constant 
$U_0 \in (0,1)$ and study the solution in the 
future region determined by the portion of the nontrivial data belonging to the 
following annular region in $\mathbb{R}^3:$
$\Sigma_0^{U_0}:=\lbrace (x^1,x^2,x^3) \ | \ 1 - U_0 \leq r \leq 1 \rbrace,$
where $r = \sqrt{\sum_{a=1}^3 (x^a)^2}$ is the standard Euclidean
radial coordinate. In doing so, we avoid the origin, where
$r$ vanishes and some of our estimates degenerate.
The region of interest, which is evolutionarily determined by
the nontrivial data in $\Sigma_0^{U_0},$
is trapped between the two outgoing
null cones $\mathcal{C}_{U_0}$ and $\mathcal{C}_0,$ where the latter cone is flat
because $\Psi$ completely vanishes in its exterior (see Figure~\ref{F:REGION});
we explain this region in much greater detail below.

\begin{center}
\begin{overpic}[scale=.333]{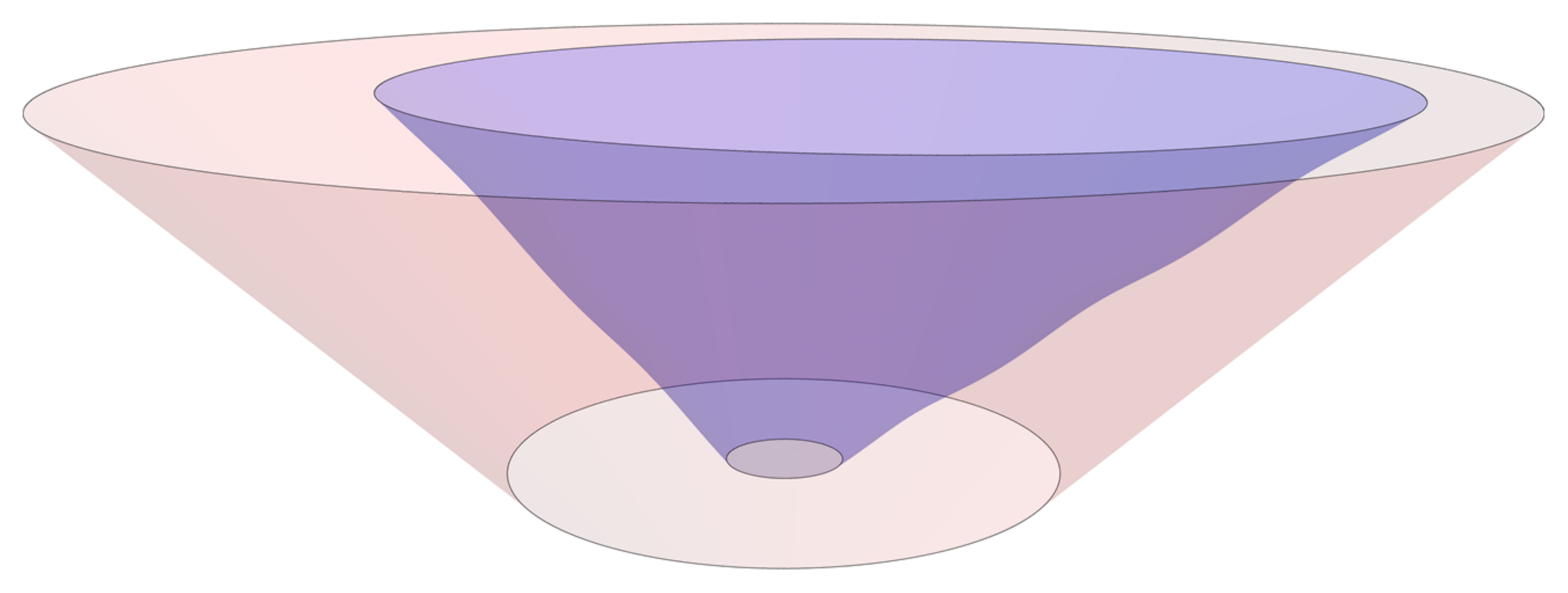} 
\put (66,16) {\large$\displaystyle \mathcal{C}_{U_0}$}
 \put (74,11.5) {\large$\displaystyle \mathcal{C}_0$}
\put (42,3.5) {nontrivial data}
\put (70,3.5) {trivial data}
\put (90,15) {\large$\displaystyle \Psi \equiv 0$}
\end{overpic}
\captionof{figure}{The region of interest}
\label{F:REGION}
\end{center}

In brief, our main goal in this monograph is to show that if the nonlinearities
in \eqref{E:WAVEGEO} fail Klainerman's classic null condition \cite{sk1984}
(see Sect.~\ref{S:STRUCTUREOFNONLINEARITIES}), 
then nontrivial data of arbitrarily small size 
verifying an \emph{open} condition launch a solution 
that forms a shock in finite time. 
Furthermore, we give a detailed description of the various phases of the dynamics of all small-data solutions, 
up to and including the time of first shock formation for those that form shocks. 
We stress that the covariant wave equation \eqref{E:WAVEGEO}
features both quasilinear and semilinear terms when it
is written relative to local coordinates; see equation \eqref{E:NEWFORMWAVEEQUATIONINCOORDINATES}
However, for equation \eqref{E:WAVEGEO},
failure of the classic null condition for both the quasilinear terms and the semilinear terms
is \emph{determined by a single factor}, a function $\FutFailFac$ that we define below in \eqref{E:INTROFAILFACT}.

We now briefly discuss the second class of equations that we study,
that is, non-covariant equations of the form
\begin{align} \label{E:ONEDERIVATIVEQUASILINEARWAVE}
	(h^{-1})^{\alpha \beta}(\partial \Phi) \partial_{\alpha} \partial_{\beta} \Phi & = 0.
\end{align}
In analogy with the case of equation \eqref{E:WAVEGEO}, we assume that
\begin{align} \label{E:LITTLEHDECOMPOSED}
	h_{\mu \nu} 
	& := m_{\mu \nu} 
		+ h_{\mu \nu}^{(Small)}(\partial \Phi),
\end{align}
where the $h_{\mu \nu}^{(Small)}$ are smooth functions of\footnote{Throughout, 
$\partial \Phi$ denotes the gradient of $\Phi$ with respect to the rectangular spacetime coordinates.} 
$\partial \Phi$ 
(at least when $\partial \Phi$ is sufficiently small)
verifying 
\begin{align} \label{E:LITTLEHVANISHESATZERO}
	h_{\mu \nu}^{(Small)}(0) = 0.
\end{align}
It turns out that the study of solutions to \eqref{E:LITTLEHDECOMPOSED}
can be effectively reduced to the study of equation \eqref{E:WAVEGEO}. 
We provide an outline of the reduction
in Appendix~\ref{A:EQUIVALENTPROBLEM}.
In short, by differentiating \eqref{E:ONEDERIVATIVEQUASILINEARWAVE} with rectangular 
coordinate derivatives $\partial_{\nu},$
the question of the long-time behavior of solutions to \eqref{E:ONEDERIVATIVEQUASILINEARWAVE}
can be transformed into an equivalent question of the long-time behavior
of solutions to a coupled system of equations in the unknowns 
$\Psi_{\nu} := \partial_{\nu} \Phi.$
The main point is that the system comprises scalar equations 
that are closely related to equation \eqref{E:WAVEGEO}.
Furthermore, in the small-data regime, 
the system turns out to be rather weakly coupled except in a few key aspects.
This structure was first observed by Christodoulou in \cite{dC2007} for a class of 
wave equations that derive from a Lagrangian. 
In this monograph, we show that the structure survives 
for general wave equations of the form \eqref{E:ONEDERIVATIVEQUASILINEARWAVE}-\eqref{E:LITTLEHDECOMPOSED}.
This fact is based on the availability of some
good (null) structure present in certain semilinear terms
(see Lemma~\ref{L:SPECIALNULLSTRUCTUREINHOMOGENEOUS}), which 
completely vanished in Christodoulou's work \cite{dC2007}.
In total, except for a handful of key aspects, the difference between the coupled system
and the scalar equation \eqref{E:WAVEGEO} is small. Hence, in this
monograph, we provide detailed proofs only in the case of the scalar equation 
\eqref{E:WAVEGEO}. We remark that one important additional ingredient 
is needed to understand the long-time behavior of solutions
to \eqref{E:ONEDERIVATIVEQUASILINEARWAVE}. We explain this
remark in Sect.~\ref{S:MAINNEWESTIMATE} and provide the
ingredient in Prop.~\ref{P:EACHRADPSINUFACTORISESSENTIALLYTWICELUPMU}.

\section{The structure of the equation relative to rectangular coordinates}
	\label{S:STRUCTUREOFNONLINEARITIES}
Relative to the rectangular coordinates introduced in Sect.~\ref{S:INTRORECTANGULARCOORDINATES}, 
equation \eqref{E:WAVEGEO} can be expressed as
\begin{align} \label{E:WAVEEQUATIONINCOORDINATES}
	(g^{-1})^{\alpha \beta} \partial_{\alpha} \partial_{\beta} \Psi
	+ (g^{-1})^{\alpha \beta} (g^{-1})^{\kappa \lambda} \Gamma_{\alpha \kappa \beta} \partial_{\lambda} \Psi
	& = 0,
\end{align}
where (recall the decomposition \eqref{E:LITTLEGDECOMPOSED})
\begin{align}
	\Gamma_{\alpha \kappa \beta}
	= \Gamma_{\alpha \kappa \beta}(\Psi,\partial \Psi)
	& := 
		\frac{1}{2}
		\left\lbrace
			\partial_{\alpha} g_{\kappa \beta}^{(Small)}
			+ \partial_{\beta} g_{\alpha \kappa}^{(Small)}
			- \partial_{\kappa} g_{\alpha \beta}^{(Small)}
		\right\rbrace
\end{align}
are the lowered Christoffel symbols of $g$ relative to the rectangular coordinates. 
In order to fully clarify the structure
of the nonlinearities, we introduce the following smooth functions of $\Psi,$ which
play a fundamental role in our analysis.

\begin{definition}[\textbf{Derivatives of the metric component functions with respect to} $\Psi$]
\label{D:BIGGANDBIGGPRIME}
Relative to the rectangular coordinates, we define
the smooth functions $G_{\mu \nu}(\cdot)$ as follows:
\begin{align} \label{E:BIGGDEF}
	G_{\mu \nu}
	= G_{\mu \nu}(\Psi)
	:= \frac{d}{d \Psi} g_{\mu \nu}^{(Small)}(\Psi).
\end{align}
 
We also define the 
following smooth functions 
$G_{\mu \nu}'(\cdot),$ which play 
a supporting role in our analysis:
\begin{align} \label{E:BIGGPRIMEDEF}
	G_{\mu \nu}'
	= G_{\mu \nu}'(\Psi)
	:= \frac{d}{d \Psi} G_{\mu \nu}(\Psi).
\end{align}
\end{definition}

Using \eqref{E:BIGGDEF}, we can express
\begin{align}
	2 \Gamma_{\alpha \kappa \beta}
	& = G_{\kappa \beta} \partial_{\alpha} \Psi
		+ G_{\alpha \kappa} \partial_{\beta} \Psi
		- G_{\alpha \beta} \partial_{\kappa} \Psi
\end{align}
and hence equation \eqref{E:WAVEEQUATIONINCOORDINATES} can be rewritten as
\begin{align} \label{E:NEWFORMWAVEEQUATIONINCOORDINATES}
	(g^{-1})^{\alpha \beta} \partial_{\alpha} \partial_{\beta} \Psi
	+ \frac{1}{2}
		(g^{-1})^{\alpha \beta} (g^{-1})^{\kappa \lambda} 
		\left\lbrace
				2 G_{\alpha \kappa} \partial_{\beta} \Psi
			-  G_{\alpha \beta} \partial_{\kappa} \Psi
		\right\rbrace
	\partial_{\lambda} \Psi
	& = 0.
\end{align}

\section{The classic null condition}
\label{S:WHENNULLCONDITIONFAILS}
Klainerman's well-known result \cite{sK1986},
which he derived using the vectorfield method,
shows that the equation
\eqref{E:NEWFORMWAVEEQUATIONINCOORDINATES} has global small-data solutions
if the nonlinearities verify his classic \emph{null condition}
\cite{sk1984};
see also \cite{dC1986a} for Christodoulou's alternate proof, which is based on the conformal method. 
There are several equivalent ways of formulating the classic null condition.
The standard way is to use Klainerman's original definition \cite{sk1984}, 
which involves Taylor expanding the nonlinearities
around $(\Psi, \partial \Psi, \partial^2 \Psi)=\mathbf{0}$
and keeping only the quadratic part.
In the case of equation \eqref{E:NEWFORMWAVEEQUATIONINCOORDINATES},
we compute that the quadratic terms are, up to constant factors, 
as follows:
\begin{subequations}
\begin{align}
	& G_{\kappa \lambda}(\Psi=0) 
	(m^{-1})^{\alpha \kappa} 
	(m^{-1})^{\beta \lambda}
	\Psi 
	\partial_{\alpha} \partial_{\beta} \Psi,
		\label{E:QUASILINEARFAILSNULL} \\
	& G_{\kappa \lambda}(\Psi=0)
		(m^{-1})^{\kappa \lambda}
		(m^{-1})^{\alpha \beta}
		\partial_{\alpha} \Psi \partial_{\beta} \Psi,
		\label{E:VERIFIESNULL} \\
	& G_{\kappa \lambda}(\Psi=0)
		(m^{-1})^{\alpha \kappa}
		(m^{-1})^{\beta \lambda}
		\partial_{\alpha} \Psi \partial_{\beta} \Psi.
		\label{E:SEMILINEARFAILSNULL}
\end{align}
\end{subequations}
By definition, the quasilinear terms
\eqref{E:QUASILINEARFAILSNULL} are said to verify the classic null
condition if and only if for every Minkowski-null covector\footnote{Recall that 
$\xi$ is Minkowski-null if and only if $(m^{-1})^{\alpha \beta} \xi_{\alpha} \xi_{\beta} = 0.$} 
$\xi,$
we have $G_{\kappa \lambda}(\Psi=0)
	(m^{-1})^{\alpha \kappa} 
	(m^{-1})^{\beta \lambda} 
	(m^{-1})^{\beta \lambda} \xi_{\alpha} \xi_{\beta} = 0.$
For the semilinear terms \eqref{E:VERIFIESNULL}, the
defining condition is
$G_{\kappa \lambda}(\Psi=0)
		(m^{-1})^{\kappa \lambda}
		(m^{-1})^{\alpha \beta} \xi_{\alpha} \xi_{\beta} = 0,$
while for \eqref{E:SEMILINEARFAILSNULL} the defining condition is
$G_{\kappa \lambda}(\Psi=0)
		(m^{-1})^{\alpha \kappa}
		(m^{-1})^{\beta \lambda} \xi_{\alpha} \xi_{\beta} = 0.$
It follows that \eqref{E:VERIFIESNULL} always verifies
the classic null condition, while
\eqref{E:QUASILINEARFAILSNULL} and \eqref{E:SEMILINEARFAILSNULL}
verify it if and only if 
$G_{\alpha \beta}(\Psi=0)\ell^{\alpha} \ell^{\alpha} = 0$
for every Minkowski-null vector $\ell.$

We now examine the classic null condition from a different perspective,
one more closely connected to our analysis of shock-forming solutions
in the region $\lbrace t \geq 0 \rbrace.$
To proceed, we consider, for the purpose of illustration, the vectorfield frame
\begin{align} \label{E:MINKOWSKIFRAME}
	\lbrace 
		\Lunit_{(Flat)} := \partial_t + \partial_r, 
		\Radunit_{(Flat)} := - \partial_r,
		X_{(Flat);1}, 
		X_{(Flat);2} 
	\rbrace,
\end{align}
where $X_{(Flat);1}$ and $X_{(Flat);2}$
are a local frame on the Euclidean spheres of constant $t$ and $r$
and $\partial_r$ is the standard Euclidean radial derivative.
Note that relative to the frame
\eqref{E:MINKOWSKIFRAME}, we can decompose the inverse Minkowski metric as follows:
\begin{align} \label{E:MINKINFLATFRAME}
	(m^{-1})^{\alpha \beta}
	& = - \Lunit_{(Flat)}^{\alpha} \Lunit_{(Flat)}^{\beta}
		- (\Lunit_{(Flat)}^{\alpha} \Radunit_{(Flat)}^{\beta} 
				+ \Radunit_{(Flat)}^{\alpha} \Lunit_{(Flat)}^{\beta})
		+ (\minversesphere)^{AB} X_{(Flat);A}^{\alpha} X_{(Flat);B}^{\beta},
\end{align}
where the $2 \times 2$ matrix $(\minversesphere)^{AB}$ is 
the inverse of the $2 \times 2$ matrix $\msphere_{AB} := m_{\alpha \beta} X_{(Flat);A}^{\alpha} X_{(Flat);B}^{\beta}.$
With the help of \eqref{E:MINKINFLATFRAME} 
and \eqref{E:QUASILINEARFAILSNULL}-\eqref{E:SEMILINEARFAILSNULL},
it is straightforward to see that
for equation \eqref{E:NEWFORMWAVEEQUATIONINCOORDINATES},
Klainerman's classic null condition is verified if and only the quadratic part 
of the nonlinearities, 
when expressed in terms of the frame derivatives,
does not contain any quasilinear terms proportional to\footnote{Throughout, if $X$ is a vectorfield and $f$ is a function,
then $X f := X^{\alpha} \partial_{\alpha} f$ denotes the derivative of $f$ in the direction $X.$} 
$\Psi \Radunit_{(Flat)} (\Radunit_{(Flat)} \Psi)$ or semilinear terms
proportional to $(\Radunit_{(Flat)} \Psi)^2.$
Simple calculations yield that up to constant factors,
both of these terms have the same
coefficient, namely the
\emph{future null condition failure factor}
$\FutFailFac$ defined by
\begin{align} \label{E:INTROFAILFACT}
	\FutFailFac 
	:= \underbrace{G_{\alpha \beta} (\Psi=0)}_{\mbox{constants}}
	\Lunit_{(Flat)}^{\alpha} \Lunit_{(Flat)}^{\beta}.
\end{align}
Note that $\FutFailFac$ can be viewed as a function 
depending only on $\theta = (\theta^1,\theta^2),$ 
where $\theta^1$ and $\theta^2$ are local angular 
coordinates corresponding to standard spherical coordinates 
on Minkowski spacetime.
The following important result
follows from the methods of \cite{sK1986} and \cite{dC1986a}:
if the terms $\Psi \Radunit_{(Flat)} (\Radunit_{(Flat)} \Psi)$ 
and $(\Radunit_{(Flat)} \Psi)^2$ are absent
from equation \eqref{E:NEWFORMWAVEEQUATIONINCOORDINATES}
when it is expanded relative to the frame \eqref{E:MINKOWSKIFRAME},
then small-data \emph{future-global} existence holds.
Hence, up to constant factors,
$\FutFailFac$ is the \underline{coefficient} of the terms
in equation \eqref{E:NEWFORMWAVEEQUATIONINCOORDINATES}
that are possible obstructions to future-global existence.
The relevant point is that in the region $\lbrace t \geq 0 \rbrace,$
$\Psi \Radunit_{(Flat)} (\Radunit_{(Flat)} \Psi)$
and $(\Radunit_{(Flat)} \Psi)^2$
decay more slowly than the other quadratic nonlinear terms
and hence have the potential to eventually 
cause singularities to form
(see however, Remark~\ref{R:QUASILINEARTERMSAREHARMLESSBYTHEMSELVES} concerning the product $\Psi \Radunit_{(Flat)} (\Radunit_{(Flat)} \Psi)$). All other quadratic terms decay sufficiently fast as $t \to \infty$
to allow for small-data future-global solutions.
Moreover, from the discussion in the previous paragraph,
we conclude that Klainerman's classic null condition 
is verified by the nonlinear terms in equation \eqref{E:NEWFORMWAVEEQUATIONINCOORDINATES} 
if and only if\footnote{\label{FN:ALLPOSSIBLENULLVECTORS}
The key point is that all possible future-directed Minkowski-null vectors 
$\ell$ from the previous paragraph are achieved by the vectorfield
$\Lunit_{(Flat)}$ in \eqref{E:INTROFAILFACT} as it varies over the region $\lbrace t \geq 0 \rbrace.$}
$\FutFailFac \equiv 0.$
Our main shock-formation theorem, 
Theorem~\ref{T:STABLESHOCKFORMATION}, 
provides a converse to the results of \cite{sK1986} and \cite{dC1986a} 
for the equations under consideration:
if $\FutFailFac \not \equiv 0,$ 
then equation \eqref{E:NEWFORMWAVEEQUATIONINCOORDINATES} exhibits 
future shock formation,
caused by the presence of the quadratic terms
$\Psi \Radunit_{(Flat)} (\Radunit_{(Flat)} \Psi)$
and $(\Radunit_{(Flat)} \Psi)^2,$
in solutions launched by an open set of small data.


\begin{remark}[\textbf{Most equations of type \eqref{E:NEWFORMWAVEEQUATIONINCOORDINATES} have} $\FutFailFac \not \equiv 0$]
	\label{R:MOSTEQNSFAILNULLCONDITION}
	It follows from the above discussion that for equation \eqref{E:NEWFORMWAVEEQUATIONINCOORDINATES},
	$\FutFailFac$ completely vanishes if and only if, relative to the rectangular coordinates,
	$g_{\alpha \beta} = (1 + f(\Psi))m_{\alpha \beta} + \mathcal{O}(\Psi^2)$
	for some smooth function $f(\Psi)$ verifying $f(0) = 0.$
	Hence, for equations of type \eqref{E:NEWFORMWAVEEQUATIONINCOORDINATES}, 
	the classic null condition is restrictive and holds only in very special cases.
\end{remark}

\begin{remark}[\textbf{The quasilinear terms by themselves do not cause small-data singularities}]
	\label{R:QUASILINEARTERMSAREHARMLESSBYTHEMSELVES}
	The results of Alinhac \cite{sA2003} and Lindblad \cite{hL2008} show
	that if we modify equation \eqref{E:WAVEEQUATIONINCOORDINATES}
	by deleting the semilinear terms
	that fail the classic null condition, then the modified equation admits small-data global solutions.\footnote{
	More precisely, Lindblad's work \cite{hL2008} treats the general case,
	while Alinhac's work \cite{sA2003} addresses the specific equation $-\partial_t^2 \Psi + (1 + \Psi)^2 \Delta \Psi = 0.$}
	Their result holds even though the modified equation
	can contain a quasilinear term proportional to $\Psi \Radunit_{(Flat)} (\Radunit_{(Flat)} \Psi),$
	which fails the classic null condition.
	In the presence of such a term, the global solutions can have distorted asymptotics.
	In particular, the true characteristics can significantly deviate from the 
	Minkowskian characteristics as $t \to \infty.$ Roughly, the analog of $\upmu$
	in \cite{sA2003} and \cite{hL2008} can become arbitrarily small as $t \to \infty,$
	but it never vanishes in finite time.
	On the other hand, if we delete the quasilinear terms from 
	equation \eqref{E:WAVEEQUATIONINCOORDINATES} 
	but retain the semilinear terms,
	then the results of John \cite{fJ1981} imply that under certain structural assumptions 
	on the semilinear terms,
	all nontrivial solutions must break down in finite time.
	However, the singularity formation mechanism is not known.
	
	We stress that \emph{the results of Alinhac \cite{sA2003} and Lindblad \cite{hL2008} do not extend to equation}
	\eqref{E:ONEDERIVATIVEQUASILINEARWAVE}. More precisely, even though there are no semilinear
	terms present in this equation, it exhibits small-data shock formation
	when its nonlinearities fail the classic null condition; see Remark~\ref{R:BIGDIFFERENCE}.
\end{remark}

\begin{remark}[\textbf{Shock formation for $\lbrace t \leq 0 \rbrace$}]
We could also study past shock formation. In this case, 
the function $\FutFailFac$ from \eqref{E:INTROFAILFACT}
needs to be modified to account for the fact that
$\Lunit_{(Flat)}$ is not outward pointing as we head towards the past;
see Remark~\ref{R:PASTFAILUREFACTOR}.
\end{remark}

\section{Basic geometric constructions}
\label{S:INTROBASICGEO}
In deriving our main results, we adopt the framework of \cite{dC2007}
and perform the vast majority of our
analysis relative to a new system of \emph{geometric coordinates}: 
\begin{align} \label{E:INTROGEOMETRICCOORDINATES}
	(t,u,\vartheta^1,\vartheta^2).
\end{align}
The coordinate $t$ is the Minkowski time coordinate from Sect.~\ref{S:INTRORECTANGULARCOORDINATES}, 
while the ``dynamic coordinate'' $u$ is an \emph{outgoing eikonal function}.
Specifically, $u$ a solution to the eikonal equation
\begin{align} \label{E:INTROEIKONAL}
	(g^{-1})^{\alpha \beta}(\Psi) \partial_{\alpha} u \partial_{\beta} u = 0
\end{align}
with level sets that are outgoing to the future (see Figure~\ref{F:SOLIDREGION}).
We define the initial value of $u$ by 
\begin{align}
	u|_{t=0} = 1 - r,
\end{align}
where $r$ is the standard radial coordinate on $\mathbb{R}^3.$
The solution to \eqref{E:INTROEIKONAL} is a perturbation of the flat
eikonal function $u_{(Flat)} = 1 + t - r.$
The level sets of $u,$ which we denote by $\mathcal{C}_u,$
are null (characteristic) hypersurfaces of the Lorentzian metric $g.$
The $\mathcal{C}_u$ 
intersect the constant Minkowski-time
hypersurfaces $\Sigma_t$ in spheres $S_{t,u}.$
We denote the Riemannian metric that $g$ induces on $S_{t,u}$
by $\gsphere.$
We denote the annular region in $\Sigma_t$ that is trapped between
the inner sphere $S_{t,u}$ and the outer sphere $S_{t,0}$ by $\Sigma_t^u.$
We denote the portion of $\mathcal{C}_u$ in between
$\Sigma_0$ and $\Sigma_t$ by $\mathcal{C}_u^t.$
We denote the solid spacetime region in between 
$\Sigma_0,$ $\Sigma_t,$ $\mathcal{C}_u^t,$ 
and $\mathcal{C}_0^t$ by $\mathcal{M}_{t,u}$
(see Figure~\ref{F:SOLIDREGION}). More precisely,
we define $\mathcal{M}_{t,u}$ to be 
``open at the top.'' That is, $\Sigma_t^u$ is not part of $\mathcal{M}_{t,u}.$
In contrast, we define $\mathcal{C}_u^t$ to be ``closed at the top.''
The functions $\vartheta^1,\vartheta^2$ from \eqref{E:INTROGEOMETRICCOORDINATES}
are local coordinates on the $S_{t,u}$ that, 
as we explain below, are easy to construct.

\begin{center}
\begin{overpic}[scale=.3]{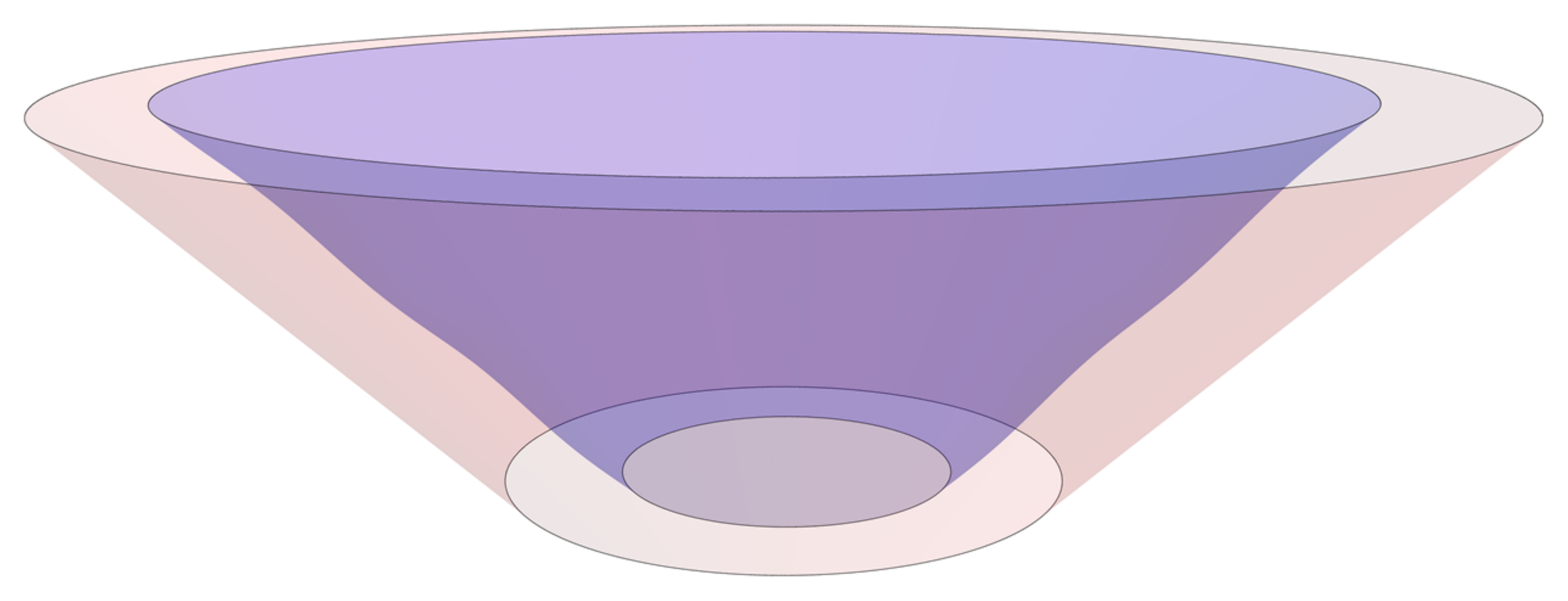} 
	\put (67,15) {\large $\displaystyle \mathcal{C}_u^t$}
  \put (71,10) {\large $\displaystyle \mathcal{C}_0^t$}
	\put (76.5,21) {\large $\displaystyle \mathcal{M}_{t,u}$}
	\put (4,29.5) {\large $\displaystyle \Sigma_t^u$}
	\put (33.2,5) {\large $\displaystyle \Sigma_0^u$}
	\put (82,30) {\large $\displaystyle S_{t,u}$}
	\put (93.5,28) {\large $\displaystyle S_{t,0}$}
	\put (56,5.7) {\large $\displaystyle S_{0,u}$}
	\put (61.5,2.5) {\large $\displaystyle S_{0,0}$}
\end{overpic}
\captionof{figure}{Surfaces and regions}
\label{F:SOLIDREGION}
\end{center}

We stress the following important point: 
the eikonal equation \eqref{E:INTROEIKONAL} is an evolution equation for $u$ that 
must be solved \emph{in conjunction} with the wave equation \eqref{E:WAVEGEO}. That is, 
we are studying the wave equation \eqref{E:WAVEGEO}
coupled to \eqref{E:INTROEIKONAL}. In rectangular coordinates, the 
wave equation completely decouples. However, as we will see, it
is extremely useful to study the wave equation relative to 
a dynamic frame constructed with the help of $u.$
Relative to this frame, \eqref{E:WAVEGEO} and \eqref{E:INTROEIKONAL} 
and their higher-order versions
become highly coupled.

Intuitively, one might expect that the intersection of the level sets of
$u$ should correspond to a shock singularity. 
This is precisely what happens
in the case of Burgers' equation in $1+1$ dimensions:
\begin{align} \label{E:BURGERS}
	\partial_t \Psi + \Psi \partial_x \Psi = 0,
\end{align}
where the role of the eikonal function is played by $\Psi$ itself.
For illustration, we now recall the basic argument leading to shock formation for
solutions to \eqref{E:BURGERS}.
The characteristics corresponding to \eqref{E:BURGERS} 
are the integral curves of the following ODE system in $\gamma(t) := (\gamma^0(t),\gamma^1(t)):$
\begin{align} \label{E:BURGERODE}
	\frac{d}{dt} \gamma^0(t) & = 1, \qquad \frac{d}{dt} \gamma^1(t) = \Psi \circ \gamma(t).
\end{align}
Equation \eqref{E:BURGERS} implies that $\Psi$ is constant along each solution curve $\gamma,$
and we therefore conclude from \eqref{E:BURGERODE} that each such curve is in fact a straight 
line in the $(t,x)$ plane with constant slope $\Psi(0,\gamma^1(0)).$
It follows that $\Psi(t, x_0 + t \mathring{\Psi}(x_0)) = \mathring{\Psi}(x_0),$
where $\mathring{\Psi}(x_0) := \Psi(0,x_0).$ From the chain rule, we see that
$\partial_x \Psi(t, x_0 + t \mathring{\Psi}(x_0)) = (1 + t \Psi'(x_0))^{-1} \mathring{\Psi}'(x_0).$
Hence, if $\mathring{\Psi}'(x_0) < 0$ for some point $x_0,$ 
then $\partial_x \Psi$ blows up 
by the time $t = - [\Psi'(x_0)]^{-1}$ due to the intersection of the characteristics emanating 
near $x_0.$ In contrast, $\Psi$ itself remains bounded.
Another way to think about the dynamics is that $\Psi$ is well-behaved along the 
directions tangent to the characteristics but becomes singular in the transversal
directions. As we will see, these basic features are also partially present 
in the shock-forming wave equation solutions that we study in this monograph.

In more than one spatial dimension, wave equations do not reduce to 
transport equations along the characteristics.
Nonetheless, it is well-established that 
in some regimes,
eikonal functions can be used
to derive sharp results for quasilinear wave-like equations in more than
one spatial dimension. For example, 
for small-data solutions to wave equations in three spatial dimensions with nonlinearities
that verify the classic null condition,
the solution's lower-order derivatives
verify a transport equation with sources that decay at an integrable-in-time rate
along the null generators of the characteristic hypersurfaces,
and this structure allows one to derive sharp estimates.
The first instance of the use of an eikonal function to solve a global nonlinear problem
is found in the celebrated proof 
of Christodoulou-Klainerman \cite{dCsK1993} of the global stability of 
Minkowski spacetime as a solution to Einstein's equations. 
Eikonal functions also played a central role in
the Klainerman-Rodnianski \cite{sKiR2003} proof of low regularity local 
well-posedness for a class of quasilinear wave equations.
They also played a fundamental role in
both Alinhac's and Christodoulou's proofs of small-data shock 
formation. Moreover, in the shock-formation proofs, 
the eikonal function played a more essential role than
it did in the proof of the stability of Minkowski spacetime.
Specifically, Lindblad and Rodnianski gave a second proof 
\cite{hLiR2010} of the stability of Minkowski spacetime,
relative to wave coordinates,
that did not rely upon a true eikonal function. Instead, 
they
closed their small-data global existence proof by deriving 
estimates relative to the 
background Minkowskian geometry\footnote{More precisely, some of the transport equation estimates in \cite{hLiR2010}
were derived along first-order corrections to the outgoing Minkowskian characteristics.} with the help of a Minkowski eikonal function
$u_{(Flat)} = t - r.$ In contrast, 
as we will see, 
for the solutions studied in this monograph,
small-data shock formation \emph{exactly corresponds} to the intersection 
of the level sets of a true eikonal function verifying \eqref{E:INTROEIKONAL}.
It is difficult to imagine an alternative sharp proof 
of small-data shock formation
that references only the background Minkowskian geometry; 
it seems that one would need at least a very close approximation of a true eikonal function.

Intimately connected to the eikonal function is the following important vectorfield, 
which can be expressed relative to rectangular coordinates as follows: 
\begin{align} \label{E:INTROLGEO}
	\Lgeo^{\nu} := - (g^{-1})^{\nu \alpha} \partial_{\alpha} u.
\end{align}
It is easy to verify that $\Lgeo$ is null and geodesic, that is, that
$g(\Lgeo,\Lgeo) = 0$ and $\D_{\Lgeo} \Lgeo = 0$
(recall that $\D$ is the Levi-Civita connection of $g$ and thus $\D g = 0$).
Furthermore, a related quantity that we mentioned above, 
the \emph{inverse foliation density} $\upmu,$
is the most important object in this monograph
from the point of view of small-data shock formation:
\begin{align} \label{E:INTROUPMU}
	\upmu := - \frac{1}{(g^{-1})^{\alpha \beta} \partial_{\alpha} u \partial_{\beta} t}
	= \frac{1}{\Lgeo^0}.
\end{align}
The function $1/\upmu$ is a measure of the density of the 
stacking of the level sets of $u$ relative to $\Sigma_t.$ 
In the case of the background solution $\Psi \equiv 0,$ 
we have $\upmu \equiv 1.$ For perturbed solutions, 
the stacking density becomes infinite 
(that is, the level sets of $u$ intersect and a shock forms)
when $\upmu = 0.$ In Figure~\ref{F:RESCALEDFRAME} on pg. \pageref{F:RESCALEDFRAME}, 
we exhibit a solution in which
$\upmu$ has become very small in a certain region because the 
density of the level sets of $u$ has become large
(that is, a shock has nearly formed in this figure).
One of our primary goals 
in this monograph is to prove that 
for a class of small-data solutions,
$\upmu$ becomes $0$ in finite time before any other kind of singularity occurs.
Actually, in addition to exhibiting a class of small data for which
$\upmu$ becomes $0$ in finite time, 
we prove in Theorem~\ref{T:LONGTIMEPLUSESTIMATES}
an analog of the main result of \cite{dC2007}, 
which is stronger: we show that for small data,
\textbf{in the region of interest, the only possible way a singularity can
form is for $\upmu$ to become $0$ in finite time.} 

We now introduce the following rescaled version of $\Lgeo,$
which plays a fundamental role in our analysis:
\begin{align} \label{E:INTROL}
	\Lunit := \upmu \Lgeo.
\end{align}
Note that $\Lunit t = \Lunit^0 = 1.$ 
Our interest in $\Lunit$ lies in the fact 
that our proof shows that when a shock forms,
the rectangular components $\Lgeo^{\nu}$ blow-up, while 
the rectangular components $\Lunit^{\nu}$ remain near 
those of $\Lunit_{(Flat)} := \partial_t + \partial_r.$
For this reason and many others that will become apparent,
$\Lunit$ is useful for analyzing solutions.

To complete our geometric
coordinate system, we complement $t,u$ with local coordinates $(\vartheta^1,\vartheta^2)$ 
on the initial Euclidean sphere $S_{0,0}$ and propagate them inward to the $S_{0,u}$ by
solving $- \partial_r \vartheta^A = 0$ and then 
to the future by solving $\Lunit \vartheta^A = 0,$
$(A=1,2).$
We denote the coordinate vectorfield corresponding to 
$\vartheta^1$ by $X_1 := \frac{\partial}{\partial \vartheta^1}|_{t,u,\vartheta^2}$
and similarly for $X_2.$ It follows that the 
$X_A$ are tangent to the spheres $S_{t,u}.$ 

\begin{remark}[$\upmu$ \textbf{is connected to the Jacobian determinant of the change of variables map}]
\label{R:UPMUCONNECTIONTODETERMINANT}
For the solutions under consideration, the Jacobian determinant of the change of variables map $\Upsilon$ from
geometric to rectangular coordinates vanishes precisely when
$\upmu$ vanishes
(see Lemma~\ref{L:JACOBIAN}).
Hence, small-data shock formation can alternatively be viewed as a breakdown
in the map $\Upsilon^{-1}.$
\end{remark}


\section{The rescaled frame and dispersive \texorpdfstring{$C^0$}{sup-norm} estimates}
\label{S:RESCALEDFRAME}
The following basic principle, first exhibited by
Christodoulou \cite{dC2007}, 
is a key ingredient in our analysis: 
the lower-order derivatives of the solutions
in the directions $\Lunit, X_1, X_2,$ which are tangent to the 
outgoing null cones $\mathcal{C}_u,$ 
exhibit dispersive behavior,
even as a shock forms!
Furthermore, when expressed in terms of the Minkowskian time coordinate $t,$
the decay rates of these derivatives are the same
as those of solutions to the linear wave equation.
This is somewhat analogous to the regular behavior observed in 
solutions to
Burgers' equation
along the characteristics
(see Sect.~\ref{S:INTROBASICGEO}).
Furthermore, 
in directions transversal to the 
$\mathcal{C}_u,$ the solutions also behave in a linearly
dispersive fashion, but 
\emph{only after one rescales the transversal direction by $\upmu.$} 
The main point is that \emph{the rescaling by $\upmu$ eliminates the semilinear
term in the wave equation that is most dangerous in the region $\lbrace t \geq 0 \rbrace$
because of its slow decay}; 
see Remark~\ref{R:RESCALINGGIVESTHENULLCONDITION}.
 
Hence, as will become apparent, 
a good choice for a rescaled transversal vectorfield is 
the \emph{radially inward} vectorfield $\Rad,$ which verifies
\begin{align}
	\Rad u &= 1, && \Rad t = 0, && g(\Rad,X_1) = g(\Rad,X_2) = 0.
\end{align} 
It is not difficult to verify using \eqref{E:GINVERSE00ISMINUSONE} that (see Lemma~\ref{L:NULLPAIRBASICPROPERTIES})
\begin{align}
	g(\Rad, \Rad) 
	& = \upmu^2, 
	&& g(\Lunit, \Rad) = - \upmu.
\end{align}
An important aspect of our proof is showing that rectangular components of the vectorfield
\begin{align}
	\Radunit := \upmu^{-1} \Rad
\end{align}
remain near those of
$\Radunit_{(Flat)} := - \partial_r$
throughout the entire
course of the evolution.
Thus, when $\upmu$ vanishes, the vectorfield $\Rad,$ when viewed
as a vectorfield on Minkowski spacetime expressed relative to rectangular
coordinates, also vanishes.

In total, the above vectorfields form a useful \emph{rescaled frame} with span
equal to $\mbox{span}\lbrace \partial_{\alpha} \rbrace_{\alpha=0,1,2,3}$
at each point where $\upmu > 0$ (see Figure~\ref{F:RESCALEDFRAME}):
\begin{align} \label{E:RESCALEDFRAME}
	\lbrace \Lunit, \Rad, X_1, X_2 \rbrace.
\end{align}	
We carry out most of our analysis relative to the frame \eqref{E:RESCALEDFRAME}.
Occasionally, when performing calculations, 
we find it convenient to replace $\Rad$ with the following vectorfield:
\begin{align}	
	\uLgood := \upmu \Lunit + 2 \Rad,
\end{align}
which is future-directed, $g-$orthogonal to the $S_{t,u},$ ingoing, null, and normalized by $g(\Lunit, \uLgood) = - 2 \upmu.$

\begin{center}
\begin{overpic}[scale=.3]{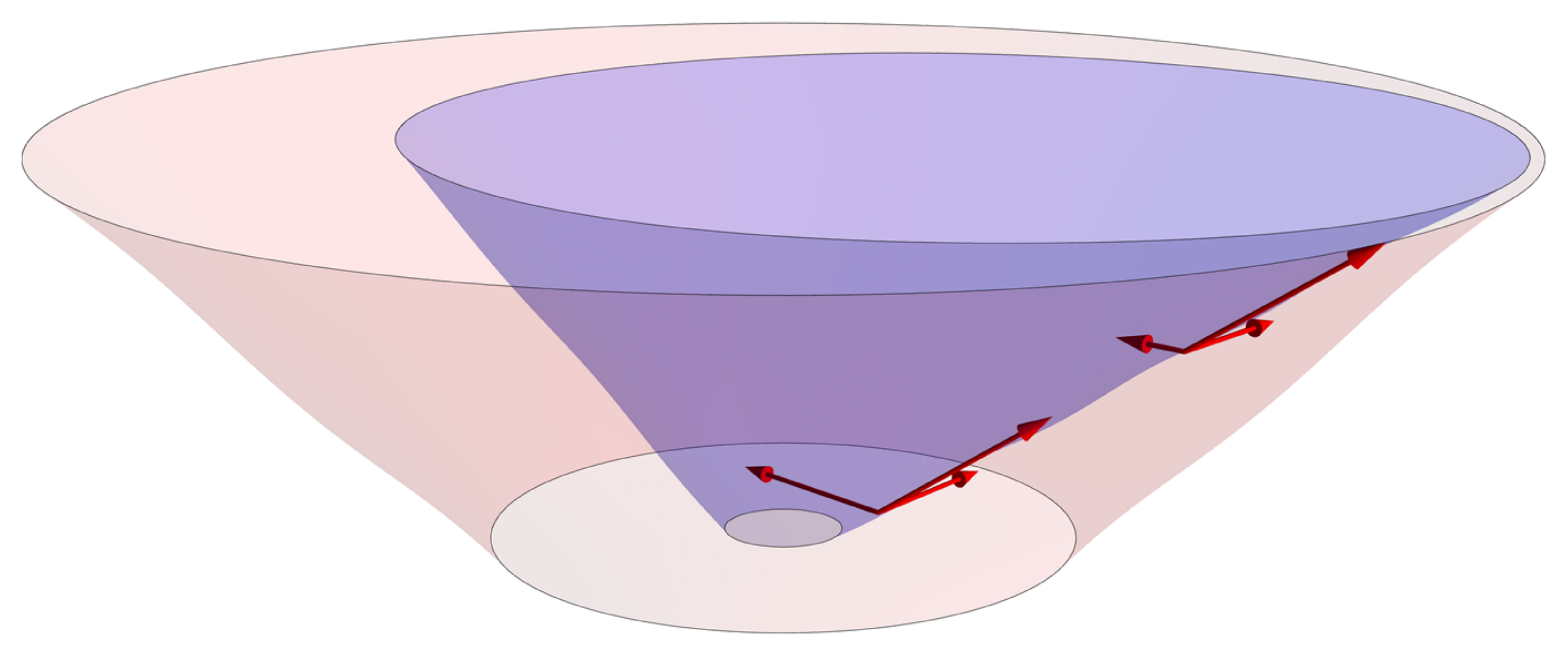}
	\put (49.5,11.5) {\large$\displaystyle \Rad$}
	\put (60,14) {\large$\displaystyle \Lunit$}
	\put (60.5,8) {\large$\displaystyle X_1, X_2$}
	\put (71,21) {\large$\displaystyle \Rad$}
	\put (81.1,24.5) {\large$\displaystyle \Lunit$}
	\put (79,17.5) {\large$\displaystyle X_1, X_2$}
	\put (90,29) {\large$\displaystyle \upmu$ is small}
	\put (52,4.5) {\large$\displaystyle \upmu \approx 1$}
\end{overpic}
\captionof{figure}{The rescaled frame at two distinct points, embedded in $(m,\mathbb{R}^{1+3})$}
\label{F:RESCALEDFRAME}
\end{center}

In the proof of our sharp classical lifespan theorem (Theorem~\ref{T:LONGTIMEPLUSESTIMATES}), 
to capture the dispersive behavior relative to the rescaled frame \eqref{E:RESCALEDFRAME},
we make the following bootstrap assumptions
(see Ch.~\ref{C:C0BOUNDBOOTSTRAP} for a more precise statement)
on a spacetime region of the form $\mathcal{M}_{\Tboot,U_0}$ 
(see Figure~\ref{F:SOLIDREGION} on pg. \pageref{F:SOLIDREGION}):
\begin{subequations}
\label{E:HP}
\begin{align}
	& \| \Lunit \Psi \|_{C^0(\Sigma_t^{U_0})}, \qquad \| \angdiff \Psi \|_{C^0(\Sigma_t^{U_0})} \leq \frac{\varepsilon}{(1 + t)^2},
		\label{E:INTROTANGENTIALARERAPIDLYDECAYING} \\
	& \| \Psi \|_{C^0(\Sigma_t^{U_0})}, \qquad \| \Rad \Psi \|_{C^0(\Sigma_t^{U_0})} \leq \frac{\varepsilon}{1 + t}.
		\label{E:INTROTRANSVERSALC0BOUND}
\end{align}
\end{subequations}
In \eqref{E:INTROTANGENTIALARERAPIDLYDECAYING} and throughout, 
$\angdiff f$ denotes the angular differential of the scalar-valued function $f,$
viewed as a function of the geometric angular coordinates $(\vartheta^1, \vartheta^2).$
Furthermore, $\varepsilon$ is a small number whose smallness is shown, 
near the end of the proof of the sharp classical lifespan theorem,
to be controlled by the size of the data (see \eqref{E:INTROSMALLDATA}).
We show that the assumptions \eqref{E:INTROTANGENTIALARERAPIDLYDECAYING}-\eqref{E:INTROTRANSVERSALC0BOUND},
which reflect the decay rates of solutions to the \emph{linear wave equation},
hold all the way up to and including the time of first shock formation.
We make similar assumptions for some of the higher derivatives of $\Psi,$
but not all the way up to top order. More precisely, 
our bootstrap assumptions imply that
each additional $\Lunit$ or angular differentiation results in a decay rate gain of
$(1 + t)^{-1}.$
These bootstrap assumptions
help us control the inessential error terms
that arise in our analysis. At the end of the proof of Theorem~\ref{T:LONGTIMEPLUSESTIMATES},
we derive improvements of
these bootstrap assumptions
from our $L^2$ estimates with the help of Sobolev embedding.

It is important to understand how the covariant wave operator looks relative to the rescaled frame.
In Prop.~\ref{P:GEOMETRICWAVEOPERATORFRAMEDECOMPOSED}, we show that
\begin{align} \label{E:INTROWAVEEQNRESCALEDFRAME}
	\upmu \square_{g(\Psi)} \Psi
	& = - \Lunit(\upmu \Lunit \Psi + 2 \Rad \Psi)
			+ \upmu \angLap \Psi
			- \mytr \upchi \Rad \Psi
			+ \mbox{Error},
\end{align}
where $\angLap$ is the covariant Laplacian 
corresponding to the metric $\gsphere$ on the $S_{t,u}$ induced by $g,$
and $\mbox{Error}$ denotes quadratically small terms that 
decay at an integrable-in-time rate. The error terms
depend on the up-to-first-order derivatives $\Psi$ and 
the up-to-second-order derivatives of $u.$ 
In \eqref{E:INTROWAVEEQNRESCALEDFRAME}, $\upchi$ is a crucially important
symmetric type $\binom{0}{2}$ tensorfield on the spheres that verifies
\begin{align} \label{E:INTROCHIDEF}
	\upchi_{AB} 
	& := \upchi(X_A,X_B)
		= g(\D_A \Lunit, X_B),
\end{align}
and $\mytr \upchi = (\ginversesphere)^{AB} \upchi_{AB}.$
Roughly, $\upchi$ can be viewed as a second derivative of $u$
that involves an angular derivative.
In proving our sharp classical lifespan theorem,
we encounter the following serious difficulty:
it takes a huge amount of effort to derive top-order $L^2$ 
estimates for $\upchi.$ One issue is that a naive approach results in 
a loss of derivatives. Another issue is that 
even though we are able to find way to avoid derivative loss,
our estimates for the top derivatives
are very degenerate and can blow-up like a power of $1/\upmu;$  
see Sects.~\ref{SS:EIKONALTOPORDER} and \ref{SS:L2HIERARCHY} for additional details.
In the case $\Psi \equiv 0$ we have $\mytr \upchi_{(Flat)} = 2 r^{-1},$
whereas for perturbed solutions, we have
$\mytr \upchi = 2 \rgeo^{-1} + \mbox{Error},$
where $\mbox{Error}$ is of size $\varepsilon \ln(\myexp + t) (1 + t)^{-2}$ and
\begin{align} \label{E:INTRORGEO}
	\rgeo(t,u) := 1 - u + t
\end{align}
is a geometric radial variable, tailored to the eikonal function.
On the bootstrap region $\mathcal{M}_{\Tboot,U_0},$ 
we have the following
estimate, which we often use:
$\rgeo \approx 1 + t.$
In view of the fact that $\Lunit \rgeo = 1$ and $\Rad \rgeo = - 1,$
we deduce from \eqref{E:INTROWAVEEQNRESCALEDFRAME} that
\begin{align} \label{E:INTROALTERNATEWAVEEQNRESCALEDFRAME}
	\upmu \square_{g(\Psi)} \Psi
	& = - \Lunit 
				\left\lbrace
					\upmu \Lunit (\rgeo \Psi) + 2 \Rad (\rgeo \Psi)
				\right\rbrace
			+ \rgeo \upmu \angLap \Psi
			+ \mbox{Error},
\end{align}
where $\mbox{Error}$ denotes quadratically small terms that 
decay at an integrable-in-time rate.
Equation \eqref{E:INTROALTERNATEWAVEEQNRESCALEDFRAME} plays an important role in our proof of small-data
shock formation; see Sect.~\ref{S:LOWERBOUNDLEADINGTOBLOWUP}
and Lemma~\ref{L:ALMOSTSPHERICALLYSYMMETRICWAVEEQUATION}.

\begin{remark}[\textbf{``Eliminating'' the dangerous semilinear term by rescaling}]
	\label{R:RESCALINGGIVESTHENULLCONDITION}
	A crucially important structural feature of
	equations \eqref{E:INTROWAVEEQNRESCALEDFRAME} and \eqref{E:INTROALTERNATEWAVEEQNRESCALEDFRAME}
	is that they do not contain any quadratic term proportional to $(\Rad \Psi)^2.$
	This is quite remarkable since the wave equation \eqref{E:NEWFORMWAVEEQUATIONINCOORDINATES}, 
 	when expressed relative to the Minkowskian frame \eqref{E:MINKOWSKIFRAME},
 	can contain a term proportional to $(\Radunit_{(Flat)} \Psi)^2.$
	The point is that by using the rescaled frame 
	\eqref{E:RESCALEDFRAME} 
	and bringing the $\upmu$ factor ``inside'' the outer
	$\Lunit$ differentiation 
	in \eqref{E:INTROWAVEEQNRESCALEDFRAME} and \eqref{E:INTROALTERNATEWAVEEQNRESCALEDFRAME},
	we have generated cancellation that ``eliminates'' the term $(\Rad \Psi)^2.$
	That is, we have been able to eliminate the 
	quadratic semilinear term with the slowest decay rate as $t \to \infty;$ 
	see the discussion in Sect.~\ref{S:WHENNULLCONDITIONFAILS}.
	The price we pay is that
	\textbf{i)} our rescaled frame
		degenerates relative to the rectangular coordinate frame 
		$\lbrace \partial_{\alpha} \rbrace_{\alpha=0,1,2,3}$
		as $\upmu$ shrinks
	and \textbf{ii)}
	the terms $\mbox{Error}$
	on the right-hand side of \eqref{E:INTROALTERNATEWAVEEQNRESCALEDFRAME}
	have a somewhat complicated structure because they depend not only
	on $\Psi$ and its first derivatives, but also 
	on the up-to-second-order derivatives of $u;$
	see Prop.~\ref{P:GEOMETRICWAVEOPERATORFRAMEDECOMPOSED}.
	Hence, in order to verify that the terms $\mbox{Error}$ are in fact
	error terms, we must obtain a good understanding of the asymptotic
	behavior of $u$ and its derivatives.
\end{remark}

\section[The coupled system and \texorpdfstring{$C^0$}{sup-norm} estimates for
\texorpdfstring{$u$}{the eikonal function}]
{The basic structure of the coupled system and \texorpdfstring{$C^0$}{sup-norm} estimates for the eikonal function quantities}
\label{S:C0FOREIKONAL}
As we stressed above, we are studying the wave equation
\eqref{E:WAVEGEO} coupled to the eikonal equation \eqref{E:INTROEIKONAL}.
In practice, rather than directly studying the eikonal equation,
we study evolution equations verified by its first derivatives, 
which are represented by the rectangular components 
$\Lgeo^{\nu}$ (see \eqref{E:INTROLGEO}). 
More precisely, as we will see, it is convenient
to instead study the equations verified by 
the quantities $\upmu$ and $\Lunit^i,$ $(i = 1,2,3),$
defined by \eqref{E:INTROUPMU} and \eqref{E:INTROL}.
We can derive evolution equations for these quantities
by writing down the geodesic equation $\D_{\Lgeo} \Lgeo = 0$
relative to the rectangular coordinates and then translating it
into equations for $\upmu$ and $\Lunit^i.$ We carry this out in
full detail in Chapter~\ref{C:TRANSPORTEQUATIONSFORTHEEIKONALFUNCTIONQUANTITIES}.
We find that the following transport equations hold:
\begin{subequations} \label{E:INTROTRANSPORT}
\begin{align}
	\Lunit \upmu
	& = \frac{1}{2} G_{\Lunit \Lunit} \Rad \Psi + \mbox{Error}(\upmu, \Lunit \Psi, \Lunit^1,\Lunit^2,\Lunit^3),
		\label{E:INTROLUNITUPMUHEURISTIC} \\
	\Lunit \Lunit^i 
	& = \mbox{Error}^i(\Lunit \Psi, \angdiff \Psi, \Lunit^1,\Lunit^2,\Lunit^3),
		\label{E:INTROLUNITLUNITIHEURISTIC}
\end{align}
\end{subequations}
where $G_{\Lunit \Lunit} := G_{\alpha \beta}(\Psi)\Lunit^{\alpha} \Lunit^{\beta},$
$G_{\alpha \beta}$ is as defined in \eqref{E:BIGGDEF},
and $\mbox{Error}$ denotes small terms that decay at an integrable-in-time rate.

In total, we can use \eqref{E:HP},
\eqref{E:INTROTRANSPORT},
and a small-data assumption to derive 
$C^0$ estimates for
$\upmu,$ $\Lunit^i,$ and their lower-order derivatives.
Like the bootstrap assumptions \eqref{E:HP}, these estimates complement our $L^2$ estimates 
(which we discuss below)
and allow us to control various error terms. An important aspect of these $C^0$
estimates is that
$\upmu - 1$ and its lower-order derivatives 
can shrink or grow like $\mathcal{O}(\varepsilon \ln(\myexp + t)).$
The possible logarithmic-in-time shrinking behavior of $\upmu$ 
is in fact the source of small-data shock-formation;
we illustrate this in more detail in Sect.~\ref{S:UPMUHEURISTICS}.
The same statement holds for the difference between 
$\rgeo \Lunit^i$ and its Minkowskian value $x^i,$ that is, for
$\rgeo \Lunit^i - x^i.$ Similar statements hold for many other 
scalar and tensorial quantities that play a role in this monograph.
We derive most of these $C^0$ estimates in Chapter~\ref{C:C0BOUNDBOOTSTRAP}.

\section{Lower bounds for \texorpdfstring{$\Rad \Psi$}{the rescaled radial derivative of the solution} in the case of shock formation}
\label{S:LOWERBOUNDLEADINGTOBLOWUP}
We now heuristically explain the behavior of $\Psi$
for the shock-forming solutions of interest. 
The behavior that we describe in this section is one of the key ingredients 
in the proof of our small-data shock-formation theorem,
Theorem~\ref{T:STABLESHOCKFORMATION}.
The main idea is that there exists an open set of data 
whose corresponding solutions
verify the following \emph{signed lower bound}, 
which is valid for sufficiently large times
and for some choices of $u$ and $\vartheta$
(that is, along some integral curves of $\Lunit$):
\begin{align} \label{E:BLOWUPLOWERBOUND}
	\Rad \Psi(t,u,\vartheta) 
	& \gtrsim \frac{\varepsilon}{1 + t}.
\end{align}
Furthermore, for solutions corresponding to a different open set of data, 
the opposite signed bound holds:
$\Rad \Psi(t,u,\vartheta) \lesssim - \varepsilon (1 + t)^{-1}.$
Since $\Rad = \upmu \Radunit,$ \eqref{E:BLOWUPLOWERBOUND} is equivalent to
\begin{align} \label{E:ALTBLOWUPLOWERBOUND}
	\Radunit \Psi(t,u,\vartheta)
	& \gtrsim \frac{1}{\upmu} \frac{\varepsilon}{(1 + t)}.
\end{align}
Clearly, the bound \eqref{E:ALTBLOWUPLOWERBOUND} implies that
the near-Euclidean-unit-length derivative $\Radunit \Psi$
blows up when $\upmu \to 0.$
Moreover, the bound \eqref{E:ALTBLOWUPLOWERBOUND} 
(with the correct sign, depending on the nonlinearities)
is in fact the main ingredient needed to show that $\upmu$ in fact 
goes to $0.$ We elaborate upon this in Sect.~\ref{S:UPMUHEURISTICS}.
All of our shock-forming solutions feature blow-up in $\Radunit \Psi$ 
for exactly these reasons.

To explain how we derive \eqref{E:ALTBLOWUPLOWERBOUND}, we impose the wave equation by
setting the right-hand side of \eqref{E:INTROALTERNATEWAVEEQNRESCALEDFRAME}
equal to $0,$ which yields
\begin{align} \label{E:FRAMEWAVEINTRO}
	\Lunit 
	\left\lbrace
		\upmu \Lunit (\rgeo \Psi) 
		+ 
		2 \Rad (\rgeo \Psi)
	\right\rbrace
	& = \rgeo \upmu \angLap \Psi
		+ \mbox{Error},
\end{align}
where the terms $\mbox{Error}$ are quadratically small and decay at an integrable-in-time rate. 
For convenience, in this monograph,
we prove the bound \eqref{E:ALTBLOWUPLOWERBOUND}
in complete detail only for an open set of nearly spherically symmetric small data
(see the proof of Theorem~\ref{T:STABLESHOCKFORMATION}).
However, in Sect.~\ref{SS:INTROSHOCKFORMINGDATACOMPARISON}, we provide an outline of
how to enlarge the class of small data to which our shock-formation arguments apply.
By ``nearly spherically symmetric,'' we mean that we assume 
that the lower-order angular derivatives of the data
are even smaller than the small radial derivative, and
we prove that this condition is propagated by the solution
(see Cor.~\ref{C:ANGULARDDERIVATIVEESTIMATES}).
Furthermore, from the higher-order analog of \eqref{E:HP}, 
we find that even though $\rgeo \upmu \angLap \Psi$ 
is only linearly small,
it decays at an integrable-in-time rate.
The net effect is that for such solutions, 
the right-hand side of \eqref{E:FRAMEWAVEINTRO}
has a tiny amplitude and is integrable in time,
and \eqref{E:FRAMEWAVEINTRO} can effectively be treated as a transport equation along the integral curves of 
$\Lunit = \frac{\partial}{\partial t}.$ 
Hence, integrating \eqref{E:FRAMEWAVEINTRO} 
with respect to time at fixed $u$ and $\vartheta$
and ignoring
the small right-hand side, we obtain, relative to the geometric coordinates:
\begin{align} \label{E:INTROWAVEEQNTRANSPORTINTEGRATEDESTIMATE}
		\left\lbrace
			\upmu \Lunit (\rgeo \Psi)
			+ 2 \Rad (\rgeo \Psi)
		\right\rbrace
		(t,u,\vartheta)
	& \approx f[\mathring{\Psi}, \mathring{\Psi}_0](u,\vartheta),
\end{align}
where $f[\mathring{\Psi}, \mathring{\Psi}_0](u,\vartheta)$
is equal to the term in braces on the left-hand side of
\eqref{E:FRAMEWAVEINTRO} evaluated at $t=0,$
and $(\mathring{\Psi}, \mathring{\Psi}_0)$ are the data \eqref{E:INTRODATA}.
Furthermore, using \eqref{E:INTROTANGENTIALARERAPIDLYDECAYING}
and \eqref{E:INTROWAVEEQNTRANSPORTINTEGRATEDESTIMATE}
and the identities $\Lunit \rgeo = 1$
and $\Rad \rgeo = - 1,$
we find that for sufficiently large times, we have
\begin{align}  \label{E:INTROWAVEEQNTRANSPORTINTEGRATEDESTIMATEMIDTIMEANDLATERESTIMATE}
	\Rad \Psi(t,u,\vartheta) 
	& \approx
	\frac{1}{2} \frac{1}{\rgeo(t,u)} f[\mathring{\Psi}, \mathring{\Psi}_0](u,\vartheta).
\end{align}
For suitable data,
\eqref{E:INTROWAVEEQNTRANSPORTINTEGRATEDESTIMATEMIDTIMEANDLATERESTIMATE}
implies the desired bound \eqref{E:ALTBLOWUPLOWERBOUND}.

\section{The main ideas behind the vanishing of \texorpdfstring{$\upmu$}{the inverse foliation density}}
\label{S:UPMUHEURISTICS}
We now explain why the estimates of Sect.~\ref{S:LOWERBOUNDLEADINGTOBLOWUP}
imply that $\upmu$ vanishes in finite time.
We first address the important factor $G_{\Lunit \Lunit}$ in equation
\eqref{E:INTROLUNITUPMUHEURISTIC}. In Lemma~\ref{L:NULLCONDFACT}, 
we show that
$G_{\Lunit \Lunit}(t,u,\vartheta)$ is well-approximated by 
$\FutFailFac(t=0,u=0,\vartheta),$ where $\FutFailFac,$ the future null condition failure factor from \eqref{E:INTROFAILFACT},
is actually independent\footnote{When $t=0,$
the geometric angular coordinate $\vartheta$ coincides with the standard Euclidean spherical coordinate $\theta$
and hence, as we explained just below equation \eqref{E:INTROFAILFACT}, 
$\FutFailFac(0,u,\vartheta)$ can be viewed as a function of $\vartheta$ alone.} of $u$ at $t=0.$
We remark that one significant difference between our work and that of \cite{dC2007} is that
in \cite{dC2007}, the analog of $\FutFailFac$ is a (non-zero) constant. 
The presence of an angularly dependent $\FutFailFac$ 
in our work here creates additional complications in the analysis.
In view of the above discussion 
and the transport equation \eqref{E:INTROLUNITUPMUHEURISTIC},
we see that for the class of data that lead to the bound \eqref{E:BLOWUPLOWERBOUND}, 
and for angles $\vartheta$ with $\FutFailFac(t=0,u=0,\vartheta) < 0,$
we have
\begin{align} \label{E:LUNITUPMUBLOWUPEQUATIONHEURISTIC}
	\Lunit \upmu(t,u,\vartheta)
	& \lesssim - \varepsilon \left|\FutFailFac(t=0,u=0,\vartheta) \right| \frac{1}{1 + t} 
		+ \mbox{Error}, 
\end{align}
where the error terms are small in magnitude compared to the first term.
Integrating \eqref{E:LUNITUPMUBLOWUPEQUATIONHEURISTIC} 
along the integral curves
of $\Lunit,$ 
ignoring the error terms, 
and using the initial condition
$\upmu = 1 + \mathcal{O}(\varepsilon),$ we obtain that
\begin{align} \label{E:UPMUBLOWUPEQUATIONHEURISTIC}
	\upmu(t,u,\vartheta) \lesssim 1 - \varepsilon \left|\FutFailFac(t=0,u=0,\vartheta) \right| \ln(\myexp + t).
\end{align}
Hence, for such data, we conclude from inequality \eqref{E:UPMUBLOWUPEQUATIONHEURISTIC}
that $\upmu$ vanishes at a time $T_{(Lifespan)} \sim \exp (c \varepsilon^{-1}).$
Furthermore, for angles $\vartheta$ with $\FutFailFac(t=0,u=0,\vartheta) > 0,$
we can also derive the bound \eqref{E:UPMUBLOWUPEQUATIONHEURISTIC}
and the vanishing of $\upmu$ in finite time
by using the remark made in the first sentence after equation \eqref{E:BLOWUPLOWERBOUND}.

\section{The role of Theorem~\ref{T:LONGTIMEPLUSESTIMATES} in justifying the heuristics}
\label{S:JUSTIFYINGTHEHEURISTICS}
In Sects.~\ref{S:LOWERBOUNDLEADINGTOBLOWUP} and \ref{S:UPMUHEURISTICS},
we sketched a proof of our small-data shock-formation theorem,
Theorem~\ref{T:STABLESHOCKFORMATION}.
However, in order to justify the heuristic arguments given above,
we must 
\textbf{i)} show that coupled system 
\eqref{E:WAVEGEO}, \eqref{E:INTROLUNITUPMUHEURISTIC}, \eqref{E:INTROLUNITLUNITIHEURISTIC}
for 
$\Psi,$
$\upmu,$ 
$\Lunit^i$
is well-posed in the sense that we can control
suitable numbers of derivatives of all quantities 
without losing derivatives
and
\textbf{ii)} show that all of the ``error terms'' mentioned above
contribute only small corrections to the heuristic analysis. In particular,
we have to make sure that the blow-up behavior \eqref{E:ALTBLOWUPLOWERBOUND} 
does not couple back into the error terms and cause them to become
large near the shock, when $\upmu$ is small.
The task \textbf{ii)} is highly coupled to \textbf{i)} and occupies the bulk of our work.
As we have mentioned,
our main results in this direction are provided by
Theorem~\ref{T:LONGTIMEPLUSESTIMATES}, which is the main
theorem in the monograph. Roughly, the theorem states that 
the solution persists unless $\upmu$ vanishes, at which point a shock has formed.
Furthermore, the theorem provides a slew of quantitative estimates
that must hold until $\upmu$ vanishes. 
In fact, many quantities 
can be extended \underline{relative to the geometric coordinates}
as $C^k$ functions for some $k \geq 0$
all the way to the constant-time
hypersurface subset $\Sigma_{T_{(Lifespan);U_0}}^{U_0}$ 
where the first shock singularity point lies.

\subsection{\texorpdfstring{$L^2$}{Square integral} estimates for \texorpdfstring{$\Psi$}{the wave variable}}
The only known methods for controlling solutions 
to quasilinear wave equations are heavily based on $L^2$ estimates.
To derive such estimates for $\Psi,$ we use 
both the vectorfield multiplier method and the vectorfield 
commutator methods. 
The multiplier method roughly corresponds 
to applying a well-chosen vectorfield differential operator to $\Psi,$
multiplying both sides of the wave equation by 
this quantity, and integrating by parts over 
the spacetime region $\mathcal{M}_{t,u}$
with the help of the divergence theorem.
We present the details in Chapter~\ref{C:DIVERGENCETHEOREM}.
In this monograph, we use two distinct multipliers:
\begin{subequations}
\begin{align}
	\Mult 
	& := (1 + 2 \upmu) \Lunit + 2 \Rad,
		\label{E:MULTINTRO} \\
	\Mor 
		&:= \rgeo^2 \Lunit.	
		\label{E:MORINTRO}
\end{align}
\end{subequations}
The timelike multiplier \eqref{E:MULTINTRO} was used in \cite{dC2007} while the 
Morawetz-type multiplier \eqref{E:MORINTRO} 
is a slightly modified version of a multiplier 
used in \cite{dC2007}.

We first discuss the estimates generated by \eqref{E:MULTINTRO}. 
In this case, energies $\enzero[\Psi](t,u)$ 
and cone fluxes $\flzero[\Psi](t,u)$ 
naturally arise from applying the divergence theorem
to solutions of $\square_{g(\Psi)} \Psi = 0$ on the
spacetime region $\mathcal{M}_{t,u}$ bounded by
$\Sigma_0^u,$ $\Sigma_t^u,$ $\mathcal{C}_u^t,$ and $\mathcal{C}_0^t$
(see Figure~\ref{F:SOLIDREGION} on pg. \pageref{F:SOLIDREGION}).
These quantities 
have the following $L^2$ coerciveness properties\footnote{It is important that the coefficient
of $\| \Rad \Psi \|_{L^2(\Sigma_t^u)}^2$ in \eqref{E:INTROMULTENERGYCOERCIVITY} is equal to $1.$
This constant affects the number of derivatives that we need to close our estimates;
see the sentence just above inequality \eqref{E:GRONWALLREADYINTROTOPORDERENERGYIDCARICATURE}.}
(see Lemma~\ref{L:COERCIVEENERGIESANDFLUXES}):
\begin{subequations}
		\begin{align} \label{E:INTROMULTENERGYCOERCIVITY}
				\enzero[\Psi](t,u)
				& \geq
							C^{-1} \| \Psi \|_{L^2(S_{t,u})}^2
								+
							C^{-1} \| \Psi \|_{L^2(\Sigma_t^u)}^2
								+
							C^{-1} \| \sqrt{\upmu} \Lunit \Psi \|_{L^2(\Sigma_t^u)}^2
								+
							C^{-1} \| \upmu \Lunit \Psi \|_{L^2(\Sigma_t^u)}^2
							\\
					& \ \ 
								+ 
							\| \Rad \Psi \|_{L^2(\Sigma_t^u)}^2
								+ 
							C^{-1} \| \sqrt{\upmu} \angdiff \Psi \|_{L^2(\Sigma_t^u)}^2
								+ 
							C^{-1} \| \upmu \angdiff \Psi \|_{L^2(\Sigma_t^u)}^2,
						\notag \\
				\flzero[\Psi](t,u)
					& \geq
						C^{-1} \| \Lunit \Psi \|_{L^2(\mathcal{C}_u^t)}^2
						+ C^{-1} \| \sqrt{\upmu} \Lunit \Psi \|_{L^2(\mathcal{C}_u^t)}^2
						+ C^{-1} \| \sqrt{\upmu} \angdiff \Psi \|_{L^2(\mathcal{C}_u^t)}^2.
						\label{E:INTROMULTCONEFLUXCOERCIVITY} 
		\end{align}	
		\end{subequations}
In \eqref{E:INTROMULTENERGYCOERCIVITY}-\eqref{E:INTROMULTCONEFLUXCOERCIVITY},
quantities inside the norms are considered to be functions of $(t',u',\vartheta^1,\vartheta^2)$ and
\begin{align} \label{E:INTROVOLFORMS}
	d \spherevol := \sqrt{\mbox{\upshape{det}} \gsphere} \, d \vartheta^1 d \vartheta^2,
		\qquad
	d \tvol := d \spherevol \, du',
		\qquad
	d \conevol := d \spherevol \, dt',
		\qquad
	d \vol := d \spherevol du' \, dt'
\end{align}
are volume forms on
$S_{t,u},$
$\Sigma_t^u,$
$\mathcal{C}_u^t,$
and $\mathcal{M}_{t,u}$
respectively. 
The second and fourth forms 
are ``rescaled'' in the sense that
the canonical forms 
induced by $g$ are respectively
$\upmu \, d \tvol$
and $\upmu \, d \vol.$
The virtue of the rescaled forms 
is that they remain regular 
(that is, finite and non-zero)
up to and including 
the shock singularity. Hence, any degeneration of the energies and fluxes
at the shock is caused by an explicitly displayed $\upmu-$dependent factor going to $0.$
		
In the case of the multiplier $\Mult,$ the divergence theorem yields 
(see Prop.~\ref{P:DIVTHMWITHCANCELLATIONS})
the following energy-flux identity
for functions $\Psi$ vanishing along the flat outer cone $\mathcal{C}_0:$
\begin{align} \label{E:INTROMTUDIVERGENCETHM}
	\enzero[\Psi](t,u)
	+ 
	\flzero[\Psi](t,u)
	& = \enzero[\Psi](0,u)
			- \int_{\mathcal{M}_{t,u}}
					(\Mult \Psi) (\upmu \square_{g(\Psi)} \Psi)
				\, d \vol
			- 
			\frac{1}{2} 
			\int_{\mathcal{M}_{t,u}}
				\upmu \enmomtensor^{\alpha \beta}[\Psi] \deformarg{\Mult}{\alpha}{\beta}
			\, d \vol.
\end{align}
In \eqref{E:INTROMTUDIVERGENCETHM}, 
\begin{align}
	\enmomtensor_{\mu \nu}[\Psi]
	:= \D_{\mu} \Psi \D_{\nu} \Psi
	- \frac{1}{2} g_{\mu \nu} (g^{-1})^{\alpha \beta} \D_{\alpha} \Psi \D_{\beta} \Psi
\end{align}
denotes the energy-momentum tensorfield corresponding to $\Psi$ and
\begin{align} \label{E:INTRODEFORMATIONTENSORDEF}
	\deform{V}_{\mu \nu} := \D_{\mu} V_{\nu} + \D_{\nu} V_{\mu}
\end{align} 
denotes the deformation tensor of a vectorfield $V.$
The last term on the right-hand side of \eqref{E:INTROMTUDIVERGENCETHM} can be decomposed,
relative to the rescaled frame $\lbrace \Lunit, \Rad, X_1, X_2 \rbrace,$ 
into a collection of error terms that have to be bounded
back in terms of the left-hand side. 
Then from Gronwall's inequality, 
we can derive an a priori estimate for 
$	\enzero[\Psi](t,u) + \flzero[\Psi](t,u).$
This strategy, which at this level of generality is standard, 
is the main idea of the proof of Theorem~\ref{T:LONGTIMEPLUSESTIMATES}.

To close our estimates, we also use the aforementioned vectorfield commutator method. 
More precisely, we commute the wave equation 
\eqref{E:WAVEGEO} with a
family of vectorfields and derive the identity 
\eqref{E:INTROMTUDIVERGENCETHM} for the commuted quantities.
Our commutation vectorfields are the set
\begin{align} \label{E:COMMUTATIONSETINTRO}
	\mathscr{Z} 
	& := \lbrace \rgeo \Lunit, \Rad, \Rot_{(1)}, \Rot_{(2)}, \Rot_{(3)} \rbrace,
\end{align}
where the $\Rot_{(l)}$ are rotation vectorfields tangent to the spheres $S_{t,u}.$
The set \eqref{E:COMMUTATIONSETINTRO} has span equal to $\mbox{span}\lbrace \partial_{\alpha} \rbrace_{\alpha=0,1,2,3}$
at each point where $\upmu > 0.$
Christodoulou used a similar commutation set in \cite{dC2007}, 
the difference being that we use $\rgeo \Lunit$ in place of Christodoulou's commutation vectorfield
$(1 + t) \Lunit$ because the former vectorfield
exhibits slightly better commutation properties with various operators.
Specifically, to prove our sharp classical lifespan theorem, we commute the wave equation up to 
$24$ times with all possible combinations of vectorfields in $\mathscr{Z}.$
The set $\mathscr{Z}$ is carefully constructed so that the error terms on the right-hand
side of the identity \eqref{E:INTROMTUDIVERGENCETHM} for the differentiated quantities 
$\mathscr{Z}^M \Psi,$ $M \leq 24,$
are ``controllable.''
It is not possible for us, in this introduction, to fully describe 
the properties that
controllable error terms should have, but we note that in
Remark~\ref{R:DEFTENSORCALCULATIONSIMPORTANTASPECTS},
we highlight some of the important ones.
We construct the rotation vectorfields $\Rot_{(l)}$
(see Ch.~\ref{C:ROTATIONS})
by projecting the Euclidean rotations
$\Roteucarg{l} := \epsilon_{lab} x^a \partial_b$
onto the spheres $S_{t,u},$ that is, by 
writing them as linear combinations of 
the non-rescaled frame vectors $\lbrace \Radunit, X_1, X_2 \rbrace$
and then discarding the $\Radunit$ component.
Because the $\Roteucarg{l}$ generally contain
a (small) non-zero component of $\Radunit,$ they are unsuitable
for studying solutions near points where $\upmu$ vanishes;
by \eqref{E:ALTBLOWUPLOWERBOUND}, $\Roteucarg{l} \Psi$ can blow-up at such points.
It is important to note that all commutation vectorfields $Z \in \mathscr{Z}$ depend on the 
first derivatives of the eikonal function $u.$ 
Hence, the $Z$ are \emph{dynamically constructed}
and depend on the solution itself, through $u$
and its first derivatives.

When we use the energy-flux identity \eqref{E:INTROMTUDIVERGENCETHM} to derive a priori 
$L^2$ estimates,
we encounter the following four main difficulties. 
\begin{enumerate}
	\item The energy-flux quantities 
	\eqref{E:INTROMULTENERGYCOERCIVITY}-\eqref{E:INTROMULTCONEFLUXCOERCIVITY}
		are not sufficient for controlling some of the error integrals
		appearing on the right-hand side of \eqref{E:INTROMTUDIVERGENCETHM}.
		One reason is that some of the error terms 
		lack good $t-$weights. The same remarks apply to the error integrals
		corresponding to the higher-order quantities $\mathscr{Z}^M \Psi.$
	\item The energy defined by \eqref{E:INTROMULTENERGYCOERCIVITY} 
		\textbf{controls only $\upmu-$weighted versions
		of $\Lunit \Psi$ and $\angdiff \Psi$} 
		and hence becomes very weak if $\upmu$ decays to $0$
		(and similarly for the higher-order analogs of \eqref{E:INTROMULTENERGYCOERCIVITY}).
		However, 
		\textbf{some of the error terms lack $\upmu$ weights}
		and hence appear to be uncontrollable 
		when $\upmu$ is near $0.$
		The cone flux \eqref{E:INTROMULTCONEFLUXCOERCIVITY} fixes this problem for the $\Lunit$ derivatives,
		but not for the angular derivatives.
	\item To control the error terms, we need good $C^0$ estimates for the lower-order 
		derivatives of $\Psi$ 
		and good $C^0$ estimates for the lower-order derivatives of quantities
		constructed out of the eikonal function,
		such as $\upmu,$ $\Lunit^i,$ and $\upchi.$
	\item We need to derive suitable $L^2$ estimates to control the higher-order derivatives 
		of quantities constructed out of the eikonal function.
\end{enumerate}
The first difficulty above is standard and often appears in
the global analysis of solutions to nonlinear wave equations.
However, in the problem of small-data shock formation, 
the last three difficulties are
quite challenging, especially the last one.
In his work \cite{dC2007}, 
Christodoulou developed new strategies 
to overcome them. We have extended them so that
they apply to the equations of interest in this monograph.

To handle difficulty (1), 
the lack of good $t-$weights,
we also derive an analog 
of the energy-flux identity \eqref{E:INTROMTUDIVERGENCETHM}
in the case of the Morawetz-type multiplier \eqref{E:MORINTRO}.
Actually, the correct analog of the divergence theorem identity \eqref{E:INTROMTUDIVERGENCETHM} 
is somewhat more involved in this case because we need some lower-order correction terms 
to obtain useful estimates; see Prop.~\ref{P:DIVTHMWITHCANCELLATIONS}.
The corresponding energies $\enone$ and fluxes $\flone$ are coercive in the following sense
(see Lemma~\ref{L:COERCIVEENERGIESANDFLUXES}):
		\begin{subequations}
		\begin{align}
				\enone[\Psi](t,u)
				& \geq C^{-1} (1 + t)^2 
									\left\| \sqrt{\upmu} 
										\left(
											\Lunit \Psi 
											+ \frac{1}{2} \mytr \upchi \Psi 
										\right)
									\right\|_{L^2(\Sigma_t^u)}^2
							+ 
								\|\rgeo^2 
									\sqrt{\upmu} 
									\angdiff \Psi	
								\|_{L^2(\Sigma_t^u)}^2,
					 \label{E:INTROENONECOERCIVENESS} \\
				\flone[\Psi](t,u)
				& \geq 
					C^{-1}
					\left\| 
						(1 + t') 
						\left(
							\Lunit \Psi 
							+ \frac{1}{2} \mytr \upchi \Psi 
						\right)
					\right\|_{L^2(\mathcal{C}_u^t)}^2.
					\label{E:INTROFLUXONECOERCIVENESS}
		\end{align}	
		\end{subequations}
The good $t-$weights in \eqref{E:INTROENONECOERCIVENESS}-\eqref{E:INTROFLUXONECOERCIVENESS}
are suitable for handling most error terms.

To handle difficulty (2), 
the unfortunate presence of $\upmu-$weights in front of $|\angdiff \Psi|^2$ in the energies and cone fluxes,
we use the remarkable observation made in
\cite{dC2007}: in the case of the Morawetz multiplier \eqref{E:MORINTRO}, 
the divergence theorem yields an \emph{additional positive spacetime integral} that is activated when
$\upmu$ becomes small and that provides control over the angular derivatives, 
\emph{without the problematic $\upmu$ weights}.
Specifically, the spacetime integral $\Morint[\Psi]$ defined by\footnote{
In \eqref{E:INTROCOERCIVEMORDEF}, 
$[\Lunit \upmu]_{-} := |\Lunit \upmu|$ when $\Lunit \upmu < 0$ and $[\Lunit \upmu]_{-} := 0$ when $\Lunit \upmu \geq 0.$}
\begin{align} \label{E:INTROCOERCIVEMORDEF} 
		\Morint[\Psi](t,u)
		& :=
	 	\frac{1}{2}
	 	\int_{\mathcal{M}_{t,u}}
			\rgeo^2
			[\Lunit \upmu]_{-}
			|\angdiff \Psi|^2
		\, d \vol,
\end{align}		
appears on the right-hand side of the analog of the identity \eqref{E:INTROMTUDIVERGENCETHM}
with a negative sign. Bringing it over to the left-hand side, we obtain a spacetime Morawetz integral
that is coercive in the following sense
(see Lemma~\ref{L:MORAWETZSPACETIMECOERCIVITY}):
\begin{align} \label{E:INTROMORAWETZSPACETIMECOERCIVITY}
		\Morint[\Psi](t,u) 
		& \gtrsim
		\int_{\mathcal{M}_{t,u}}
			\mathbf{1}_{\lbrace \upmu \leq 1/4 \rbrace}
			\frac{1 + t'}{\ln(\myexp + t')}
			|\angdiff \Psi|^2
		\, d \vol.
\end{align}
Not only does the estimate \eqref{E:INTROMORAWETZSPACETIMECOERCIVITY}
solve the problem of $\upmu$ weights for $|\angdiff \Psi|^2,$ 
it also has favorable $t-$ weights.
We stress that we must use the integrals
$\Morint[\mathscr{Z}^M \Psi]$ to close the order $M$
a priori $L^2$ estimates. That is,
we use the Morawetz integrals 
at all orders, not just at the top order.
The coerciveness estimate \eqref{E:INTROMORAWETZSPACETIMECOERCIVITY} follows
directly from the following critically important\footnote{The value $1/4$ in inequality \eqref{E:INTROSMALLMUIMPLIESLMUISNEGATIVE}
is convenient but not optimal; it can be enlarged to $1$ minus a function that depends on the size of the data.}
estimate for $\Lunit \upmu$ (see \eqref{E:SMALLMUIMPLIESLMUISNEGATIVE}):
\begin{align} \label{E:INTROSMALLMUIMPLIESLMUISNEGATIVE}
	\upmu(t,u,\vartheta) \leq \frac{1}{4}
	\implies
	\Lunit \upmu(t,u,\vartheta) \leq - c \frac{1}{(1 + t) \ln(\myexp + t)}.
\end{align} 

\begin{remark}[\textbf{The ``point of no return'' implication of \eqref{E:INTROSMALLMUIMPLIESLMUISNEGATIVE}}]
	\label{R:POINTOFNORETURN}
	Inequality \eqref{E:INTROSMALLMUIMPLIESLMUISNEGATIVE} has another
	interesting consequence. Since the right-hand side of \eqref{E:INTROSMALLMUIMPLIESLMUISNEGATIVE}
	is not integrable in $t,$ if $\upmu(t,u,\vartheta) \leq 1/4,$
	then $\upmu$ must continue to shrink 
	(along the integral curve of $\Lunit$ corresponding to fixed $u$ and $\vartheta$)
	as $t$ increases until it eventually becomes $0$ and a shock forms. 
\end{remark}

To handle difficulty (3), that of obtaining good $C^0$ estimates
for lower-order derivatives,
we use different strategies for $\Psi$ and the eikonal function
quantities. For $\Psi,$ we simply use a Sobolev embedding estimate
(see Cor.~\ref{C:C0BOUNDSOBOLEVINTERMSOFENERGIES}) 
to improve the $C^0$ bootstrap assumptions \eqref{E:HP}
after we derive suitable $L^2$ estimates.
As we mentioned in Sect.~\ref{S:C0FOREIKONAL}, 
to derive $C^0$ estimates for the eikonal function quantities,
we integrate the transport equations \eqref{E:INTROLUNITUPMUHEURISTIC}-\eqref{E:INTROLUNITLUNITIHEURISTIC}
for $\upmu$ and $\Lunit^i$
and use a small-data assumption together with
the good $C^0$ bootstrap assumptions for $\Psi.$ 
Actually, to control the error terms, 
we have to derive $C^0$ estimates for 
a huge number of quantities related to the eikonal function, 
and the previous sentence, though correct in spirit,
is only a caricature of the estimates needed.
We derive many of these $C^0$ estimates in Chapter~\ref{C:C0BOUNDBOOTSTRAP}.
Though they are tedious to derive, they are not particularly difficult.

The final difficulty (4) is by far the most challenging one.
The main challenge, which was first
identified in 
\cite{dCsK1993}
and \cite{sKiR2003}
and which appears in \cite{dC2007} and the present monograph in a more severe form,\footnote{As we explain in 
Sect.~\ref{SS:INTROTOPORDEERENERGYESTIMATES}, the new difficulty in
the small-data shock formation problem is that the price one pays for avoiding top-order derivative loss 
is a the introduction of a factor of $\upmu^{-1},$ which leads to degenerate top-order $L^2$ estimates.}
is that one has to work very hard to control the top-order
derivatives of $u$ in $L^2;$ a naive approach would result 
in a loss of one derivative of $u$ relative to $\Psi$ and would seem to suggest 
that the estimates will not close. 
To further explain this difficulty, we commute
the wave equation \eqref{E:WAVEGEO} with a commutation vectorfield 
$Z \in \mathscr{Z}$ (see \eqref{E:COMMUTATIONSETINTRO}) to obtain the following schematic expression
(see \eqref{E:BOXZCOM} for a precise expression):
\begin{align} \label{E:INTROWAVECOMMUTEDSCHEMATIC}
	\upmu \square_{g(\Psi)} (Z \Psi)
	& = \upmu \D \deform{Z} \cdot \D \Psi
		+ \upmu \deform{Z} \cdot \D^2 \Psi,
\end{align}
where $\deform{Z}$ is the deformation tensor of $Z$ (see \eqref{E:INTRODEFORMATIONTENSORDEF}).
It is crucially important to note that for $Z \in \mathscr{Z},$
the term $\upmu \D \deform{Z}$ 
from \eqref{E:INTROWAVECOMMUTEDSCHEMATIC}
depends on the \emph{third derivatives} of the eikonal function $u.$
For example, in the case that $Z$ is a rotation vectorfield $\Rot,$ we obtain
(see Prop.~\ref{P:DEFORMATIONTENSORFRAMECOMPONENTS}, 
Lemma~\ref{L:WAVEONCECOMMUTEDBASICSTRUCTURE}, 
and Prop.~\ref{P:COMMUTATIONCURRENTDIVERGENCEFRAMEDECOMP}) 
\begin{align} \label{E:INTROWAVEEQUATIONROTATIONCOMMUTED}
	\upmu \square_{g(\Psi)} (\Rot \Psi)
	& = (\Rad \Psi) \Rot \mytr \upchi
		+ \cdots
		\sim (\Rad \Psi) \partial^3 u
			+ \cdots.
\end{align}
Furthermore, one can show 
(see Cor.~\ref{C:ALPHARENORMALIZED} 
and the proof of Cor.~\ref{C:FIRSTRENORMALIZEDTRANSPORTEQUATION})
that $\Rot \mytr \upchi$ verifies the transport equation
\begin{align} \label{E:TRANSPORTROTTRCHIDERIVATIVELOSS}
	\Lunit \Rot \mytr \upchi
	& = \frac{1}{2} G_{\Lunit \Lunit} \angLap \Rot \Psi
		+ \cdots.
\end{align}
By applying the divergence theorem identity \eqref{E:INTROMTUDIVERGENCETHM} 
with $\Rot \Psi$ in the role of $\Psi$
and using equation \eqref{E:INTROWAVEEQUATIONROTATIONCOMMUTED}, we can obtain control over some
second-order derivatives of $\Psi$ in $L^2,$ 
assuming that we can control the right-hand side of 
\eqref{E:INTROWAVEEQUATIONROTATIONCOMMUTED} in $L^2.$
However, to obtain $L^2$ estimates for
$\Rot \mytr \upchi,$ we must 
integrate equation \eqref{E:TRANSPORTROTTRCHIDERIVATIVELOSS}
along the integral curves of $\Lunit,$
which seems to require that we have control over
\emph{three derivatives} of $\Psi$ in $L^2.$
This derivative-loss problem cannot be overcome by simply further commuting the wave equation.
Hence, our estimates will not close
without an additional ingredient, which we describe
in Sect.~\ref{SS:EIKONALTOPORDER}.

\begin{remark}[\textbf{The importance of error terms that never appear}]
	\label{R:AVOIDDERIVATIVELOSS}
	The following property of our commutation vectorfield set $\mathscr{Z}$
	is crucially important: after commuting the 
	operator $\upmu \square_{g(\Psi)}$ one time with $Z \in \mathscr{Z}$
	(see \eqref{E:INTROWAVECOMMUTEDSCHEMATIC}), 
	we never produce the terms
	$\angdiff \Rad \upmu$ or $\Rad \Rad \upmu;$
	this fact follows from the proof of Prop.~\ref{P:IDOFKEYDIFFICULTENREGYERRORTERMS}.
	Like the term $\Rot \mytr \upchi$
	in equation \eqref{E:INTROWAVEEQUATIONROTATIONCOMMUTED},
	these terms are third-order derivatives of the eikonal function $u$
	and thus could cause derivative loss for the reasons explained just above.
	However, there is a serious discrepancy between 
	$\Rot \mytr \upchi$ and the other two terms.
	As we describe in Sect.~\ref{SS:EIKONALTOPORDER},
	there is a way to avoid the derivative loss
	corresponding to the term $\Rot \mytr \upchi.$
	In fact, the argument we provide can be extended to show that 
	for all $Z \in \mathscr{Z},$ we can avoid the derivative loss corresponding 
	to $Z \mytr \upchi.$ In contrast, we are not aware of any method that would allow us to
	avoid the derivative loss that would occur in the presence of
	$\angdiff \Rad \upmu$ or $\Rad \Rad \upmu.$
\end{remark}

\subsection{Top-order \texorpdfstring{$L^2$}{square integral} estimates for \texorpdfstring{$u$}{the eikonal function}}
\label{SS:EIKONALTOPORDER}
As an example, we now explain how we overcome the derivative loss for the eikonal 
function quantity $\Rot \mytr \upchi$ in \eqref{E:INTROWAVEEQUATIONROTATIONCOMMUTED},
which depends on third derivatives of $u.$
For convenience, we imagine that \eqref{E:INTROWAVEEQUATIONROTATIONCOMMUTED}
is an equation for a top-order quantity.
The main idea, which was first revealed in \cite{sKiR2003}, is that
by virtue of the wave equation $\square_{g(\Psi)} \Psi = 0,$
$\Lunit \mytr \upchi$ has a remarkable structure. 
Specifically, up to lower-order derivative terms, 
$\Lunit \mytr \upchi$
is equal to the curvature component 
$-(\ginversesphere)^{AB} \Cur_{\Lunit A \Lunit B},$
where $\Cur$ is the Riemann curvature tensor of $g.$ 
The special structure revealed in \cite{sKiR2003} is:
by using the wave equation for $\Psi,$ 
one can show that $-(\ginversesphere)^{AB} \Cur_{\Lunit A \Lunit B}$
can be
written as perfect $\Lunit$ derivative 
of the first derivatives of $\Psi$
plus lower-order terms; see Cor.~\ref{C:ALPHARENORMALIZED}.
Putting these $\Lunit$ derivatives back on the left,
we obtain a transport equation along the integral curves of $\Lunit$
for a ``modified''
version of $\mytr \upchi$ that does not lose
derivatives relative to $\Psi$ (see Cor.~\ref{C:FIRSTRENORMALIZEDTRANSPORTEQUATION}).
Furthermore, this special structure essentially survives under commutations, 
and in total, we obtain an equation of the following schematic form
(see Prop.~\ref{P:TOPORDERTRCHIJUNKRENORMALIZEDTRANSPORT} for the precise version):
\begin{align} \label{E:INTROMODIFIEDROTCHITRANSPORTEQUATION}
	\Lunit
		\left\lbrace
			\upmu \Rot \mytr \upchi 
			- \Rot (G_{\Lunit \Lunit} \Rad \Psi)
			+ \cdots
		\right\rbrace
	& \sim \upmu \partial^2 \Psi + \cdots.
\end{align}

Since the right-hand side of \eqref{E:INTROMODIFIEDROTCHITRANSPORTEQUATION} depends on
only two derivatives of $\Psi,$
it thus appears that we have solved the problem of losing derivatives.
However, on the right-hand side of \eqref{E:INTROMODIFIEDROTCHITRANSPORTEQUATION}, 
present in the $\cdots,$
lies another difficult term, roughly of the form
$\upmu \angfreeLiearg{\Rot} \hat{\upchi} \cdot \hat{\upchi},$
where $\hat{\upchi}$ is the trace-free part of $\upchi$
and $\angfreeLie$ is the projection of the trace-free Lie derivative operator
onto the $S_{t,u}$
(see Def.~\ref{D:ANGFREELIEDEF}).
The difficulty is that like $\Rot \mytr \upchi,$
$\angfreeLiearg{\Rot} \hat{\upchi}$ is 
(schematically) of the form $\partial^3 u$
and hence, as we described above, has the potential to cause derivative loss. 
The saving grace is that, as was first shown in \cite{sKiR2003}
based in part on the ideas of \cite{dCsK1993}, 
one can derive an elliptic PDE on the spheres $S_{t,u}$ 
of the form $\upmu \angdiv \hat{\upchi} = \cdots,$
where $\cdots$ denotes controllable source terms that do not create derivative loss.
The elliptic estimates imply that\footnote{Here and throughout, $\angD$ denotes the
Levi-Civita connection corresponding to the metric $\gsphere$ on $S_{t,u}.$} 
$\upmu \angD \hat{\upchi}$ can be 
suitably bounded in $L^2$ (see Lemma~\ref{L:TYPE02ELLIPTICESTIMATESPHERES}),
and from this estimate, one can easily 
recover the desired $L^2$ estimates for $\angfreeLiearg{\Rot} \hat{\upchi}.$

\subsection{A hierarchy of \texorpdfstring{$L^2$}{square integral} estimates}
\label{SS:L2HIERARCHY}
Now that we have sketched how we avoid losing derivatives in both 
$\Psi$ and $u,$
we now provide a quantitative overview of the $L^2$
estimates that we derive in our sharp classical lifespan theorem
(see Theorem~\ref{T:LONGTIMEPLUSESTIMATES} for the complete details). 
Our estimates involve the quantity
\begin{align}	\label{E:INTROUPMUSTAR}
	\upmu_{\star}(t,u) := \min \lbrace 1, \min_{\Sigma_t^u} \upmu \rbrace,
\end{align}
which captures the ``worst-case'' behavior of $\upmu$ in the region $\Sigma_t^u.$
Much like in \cite{dC2007}, we derive the following hierarchy of $L^2$ estimates 
for $\Psi,$ where $\enzero[\cdot]$ and $\flzero[\cdot]$ are $L^2-$controlling in the sense of
\eqref{E:INTROMULTENERGYCOERCIVITY}
and
\eqref{E:INTROMULTCONEFLUXCOERCIVITY},
$\mathring{\upepsilon}$ is the size of the data 
(which we assume to be small)
as defined in \eqref{E:INTROSMALLDATA},
$u \in [0,U_0],$
and \emph{the estimates hold up to and including the first time that $\upmu_{\star}(\cdot,u)$ vanishes,\footnote{Hence, the estimates hold for all time if no shock forms.}}
which corresponds to the time of first shock formation in a strip of eikonal function width $u:$
\begin{subequations}
\begin{align}
	\enzero^{1/2}[\mathscr{Z}^{24} \Psi](t,u)
	+ \flzero^{1/2}[\mathscr{Z}^{24} \Psi](t,u)
	& \lesssim \mathring{\upepsilon} \ln^{\Cononestar}(\myexp + t) \upmu_{\star}^{-8.75}(t,u),
		\label{E:INTROTOPORDERENERGYESTIMATE} \\
	\enzero^{1/2}[\mathscr{Z}^{23} \Psi](t,u)
	+ \flzero^{1/2}[\mathscr{Z}^{23} \Psi](t,u)
	& \lesssim \mathring{\upepsilon} \upmu_{\star}^{-7.75}(t,u),
		\label{E:INTROJUSTBELOWTOPORDERENERGYESTIMATE} \\
	& \cdots,
		\notag \\
	\enzero^{1/2}[\mathscr{Z}^{16} \Psi](t,u)
	+ \flzero^{1/2}[\mathscr{Z}^{16} \Psi](t,u)
	& \lesssim \mathring{\upepsilon} \upmu_{\star}^{-.75}(t,u),	
		\\
	\enzero^{1/2}[\mathscr{Z}^{15} \Psi](t,u)
	+ \flzero^{1/2}[\mathscr{Z}^{15} \Psi](t,u)
	& \lesssim \mathring{\upepsilon},	
			\label{E:INTROORDERNONDEGENERATEMIDORDERENERGYESTIMATE} \\
	& \cdots, 
		\notag \\
	\enzero^{1/2}[\Psi](t,u)
	+ \flzero^{1/2}[\Psi](t,u)
	& \lesssim \mathring{\upepsilon}.
	\label{E:INTROORDER0ORDERENERGYESTIMATE}
\end{align}	
\end{subequations}
In the above estimates,
$\Cononestar > 0$ is a large constant,
and
$\mathscr{Z}^M$ denotes an arbitrary $M^{th}$ order operator
constructed out of the vectorfields $Z$ belonging to the commutation set $\mathscr{Z}$ 
from \eqref{E:COMMUTATIONSETINTRO}.
The exponent $8.75$ in \eqref{E:INTROTOPORDERENERGYESTIMATE} is connected to certain
structural constants in the equations
that we describe below, and its size is intimately tied to the number
of derivatives we need to close our estimates.
We also derive a similar hierarchy for the Morawetz multiplier quantities
$\enone,$ $\flone,$ and $\Morint$ 
(see \eqref{E:INTROMORAWETZSPACETIMECOERCIVITY})
and for the norm $\| \cdot \|_{L^2(\Sigma_t^u)}$ 
of quantities such as 
$\upmu,$ $\Lunit^i,$
and $\upchi$ constructed out of the eikonal function $u.$

\begin{remark}[\textbf{The number of derivatives needed}]
	We have used $25$ derivatives to close our estimates partially out of convenience.
	We need many derivatives because
	of the degenerate top-order estimate \eqref{E:INTROTOPORDERENERGYESTIMATE}
	and because our analysis allows us to gain only one power of $\upmu_{\star}$ as we descend
	in our $L^2$ estimates. As we explain below, the gain of one power at each step seems to be 
	a fundamental, unavoidable feature, as does the need to descend to a non-$\upmu_{\star}^{-1}-$degenerate
	level (such as \eqref{E:INTROORDERNONDEGENERATEMIDORDERENERGYESTIMATE}).
	However, we expect that the total number of derivatives
	could be reduced with additional effort.
\end{remark}

\emph{The proof of \eqref{E:INTROTOPORDERENERGYESTIMATE}-\eqref{E:INTROORDER0ORDERENERGYESTIMATE}
differs in some crucially important ways from the familiar bootstrap-type argument for deriving $L^2$ estimates.} 
One new feature, inherited from the framework of \cite{dC2007}, is that
even though $\upmu_{\star}$ is a ``$0^{th}$ order'' quantity,
we cannot ``bootstrap'' its behavior because
we cannot generally say how it will behave at late times.
Hence, we have to find a way to derive the estimates
\eqref{E:INTROTOPORDERENERGYESTIMATE}-\eqref{E:INTROORDER0ORDERENERGYESTIMATE}
independent of whether or not $\upmu_{\star}$ goes to $0.$
As we describe in 
Sect.~\ref{SS:INTROLUPMUOVERUPMU}.
some of this analysis is based on a posteriori estimates.

A second feature of crucial importance is that the lower-order estimates 
\eqref{E:INTROORDERNONDEGENERATEMIDORDERENERGYESTIMATE}-\eqref{E:INTROORDER0ORDERENERGYESTIMATE}
do not degenerate at all, even if $\upmu_{\star}$ goes to $0.$ 
That is, from the vantage point of the \underline{rescaled} frame, 
the lower-order derivatives of $\Psi$ remain regular
\emph{all the way to the shock}. Equivalently,
relative to the geometric coordinates $(t,u,\vartheta^1,\vartheta^2),$ 
$\Psi$ remains regular.
The non-degenerate estimates 
\eqref{E:INTROORDERNONDEGENERATEMIDORDERENERGYESTIMATE}-\eqref{E:INTROORDER0ORDERENERGYESTIMATE}
are especially important because they can be used to derive, 
through a Sobolev embedding result
(see Cor.~\ref{C:C0BOUNDSOBOLEVINTERMSOFENERGIES}),
improvements of the $C^0$ bootstrap assumptions \eqref{E:INTROTANGENTIALARERAPIDLYDECAYING}-\eqref{E:INTROTRANSVERSALC0BOUND}.
As we have noted, in practice,
we use the $C^0$ bootstrap assumptions to show that the terms we have deemed ``error terms'' 
are in fact small.

\subsection{Top-order \texorpdfstring{$L^2$}{square integral} estimates}
\label{SS:INTROTOPORDEERENERGYESTIMATES}
We now explain how we derive the top-order estimate
\eqref{E:INTROTOPORDERENERGYESTIMATE} and 
in particular why it is rather degenerate
relative to $\upmu_{\star}^{-1}.$
We consider the commuted wave equation
\eqref{E:INTROWAVEEQUATIONROTATIONCOMMUTED} and,
for the sake of illustration,
we imagine that 
$\upmu \square_{g(\Psi)} \Rot \Psi =  (\Rad \Psi) \Rot \mytr \upchi + \cdots$ 
is the top-order wave equation. We now sketch a proof of why
the difficult term $(\Rad \Psi) \Rot \mytr \upchi$
leads to the degenerate energy estimate
\eqref{E:INTROTOPORDERENERGYESTIMATE}.
To proceed, we apply the divergence theorem identity
\eqref{E:INTROMTUDIVERGENCETHM} with $\Rot \Psi$ in place of $\Psi$
and retain only the difficult
error integral generated by the term $(\Rad \Psi) \Rot \mytr \upchi$
(which comes from the first integral on the right-hand side of \eqref{E:INTROMTUDIVERGENCETHM}).
Furthermore, we split $\Mult \Rot \Psi = 2 \Rad \Rot \Psi + (1 + 2 \upmu) \Lunit \Rot \Psi$
(see \eqref{E:MULTINTRO})
and retain only the most difficult piece $2 \Rad \Rot \Psi$ involving 
transversal derivatives of $\Psi.$
In total, we arrive at the following integral inequality
(recall that $\mathring{\upepsilon}$ is the size of the data \eqref{E:INTROSMALLDATA}):
\begin{align} \label{E:INTROTOPORDERENERGYIDCARICATURE}
	\enzero[\Rot \Psi](t,u)
	+ \flzero[\Rot \Psi](t,u)
	& \leq C \mathring{\upepsilon}^2
		- \int_{\mathcal{M}_{t,u}}
				(2 \Rad \Psi)(\Rot \mytr \upchi) \Rad \Rot \Psi
			\, d \vol
		+ \cdots.
\end{align}
As we explained in Sect.~\ref{SS:EIKONALTOPORDER},
in order to avoid losing derivatives in \eqref{E:INTROTOPORDERENERGYIDCARICATURE},
we have to replace the $\Rot \mytr \upchi$ term with the modified quantity
in braces on the left-hand side of \eqref{E:INTROMODIFIEDROTCHITRANSPORTEQUATION}.
This replacement leads to the introduction of a factor of $\upmu^{-1}$ whose importance we discuss below.
We then have to estimate the error integrals corresponding to both the modified quantity
and the discrepancy between it and $\Rot \mytr \upchi.$ 
The analysis of the former error integral is extremely technical
(see the proof of Lemma~\ref{L:DANGEROUSTOPORDERMULTERRORINTEGRAL}) and in fact
constitutes the most difficult analysis in the monograph.
In order to avoid hijacking the discussion, we instead focus on the discrepancy 
error integral, which is also difficult and is sufficient to illustrate
the main ideas behind our derivation of \eqref{E:INTROTOPORDERENERGYESTIMATE}.
The discrepancy is equal to 
$\upmu^{-1}$ multiplied by the $\Psi-$dependent quantities in braces on the left-hand
side of \eqref{E:INTROMODIFIEDROTCHITRANSPORTEQUATION}. 
The factor $\upmu^{-1}$ is of supreme importance and is at the heart of difficulty.
In \eqref{E:INTROMODIFIEDROTCHITRANSPORTEQUATION}, the only explicitly written 
discrepancy term 
is $\Rot (G_{\Lunit \Lunit} \Rad \Psi)$ because the remaining
terms in $\cdots$ involve $\mathcal{C}_u-$tangential derivatives of $\Psi$ and are much easier
handle. Furthermore, only the top-order part
$G_{\Lunit \Lunit} \Rad \Rot \Psi$ is difficult to analyze, so
we only discuss this part.
Upon making these substitutions in \eqref{E:INTROTOPORDERENERGYIDCARICATURE}, we obtain
 \begin{align} \label{E:INTROTOPORDERENERGYIDCARICATUREHARDTERM}
	\enzero[\Rot \Psi](t,u)
	+ \flzero[\Rot \Psi](t,u)
	& \leq C \mathring{\upepsilon}^2
		+
		2 \int_{\mathcal{M}_{t,u}}
				\frac{1}{\upmu} G_{\Lunit \Lunit} (\Rad \Psi) (\Rad \Rot \Psi)^2
			\, d \vol
		+ \cdots.
\end{align}
If we were to simply substitute the dispersive bound \eqref{E:INTROTRANSVERSALC0BOUND}
for $|\Rad \Psi|$ into \eqref{E:INTROTOPORDERENERGYIDCARICATUREHARDTERM}, 
then the resulting inequality would lead, by a Gronwall argument 
(based on an inequality in the spirit of \eqref{E:INTROUPMUGAINSTHROUGHTIMEINTEGRATION} below), 
to the a priori estimate 
$\enzero^{1/2}[\Rot \Psi](t,u) + \flzero^{1/2}[\Rot \Psi](t,u) \lesssim \mathring{\upepsilon} \upmu_{\star}^{- c \mathring{\upepsilon} \ln(\myexp + t)},$
which is inadequate for closing our estimates due to the growing exponent.
In particular, this line of reasoning 
would not allow us derive the  
a priori estimate
\eqref{E:INTROTOPORDERENERGYESTIMATE}.
To circumvent this difficulty, we adopt Christodoulou's strategy and
use equation \eqref{E:INTROLUNITUPMUHEURISTIC} to replace
$G_{\Lunit \Lunit} \Rad \Psi$ with $2 \Lunit \upmu$ plus some
additional error integrals that are relatively easy to control
because their integrands are decaying sufficiently fast in time.
Then, thanks to the coerciveness estimate \eqref{E:INTROMULTENERGYCOERCIVITY},
\eqref{E:INTROTOPORDERENERGYIDCARICATUREHARDTERM} becomes
\begin{align} \label{E:GRONWALLREADYINTROTOPORDERENERGYIDCARICATURE}
	\enzero[\Rot \Psi](t,u)
	+ \flzero[\Rot \Psi](t,u)
	& \leq C \mathring{\upepsilon}^2
		- 4 \int_{\mathcal{M}_{t,u}}
					\frac{\Lunit \upmu}{\upmu}(\Rad \Rot \Psi)^2
			\, d \vol
		+ \cdots
		\\
	& \leq C \mathring{\upepsilon}^2
			+ 4 \int_{t'=0}^t
						\sup_{\Sigma_{t'}^u}
							\left|\frac{\Lunit \upmu}{\upmu} \right|
						\enzero[\Rot \Psi](t',u)
					\, dt'
			\notag	
			+ \cdots
\end{align}

We now highlight two crucially important features of the estimate \eqref{E:GRONWALLREADYINTROTOPORDERENERGYIDCARICATURE}:
\textbf{i)} the factor $4$ is a universal ``structural constant'' that does not depend on how many times we commute
the equations and
\textbf{ii)} the ``Gronwall factor''
$\sup_{\Sigma_{t'}^u} \left|\frac{\Lunit \upmu}{\upmu}\right|$
on the right-hand side of \eqref{E:GRONWALLREADYINTROTOPORDERENERGYIDCARICATURE}
has a special structure that allows us to close our
estimates. 
To see the first hint of the structure, we imagine that
$\upmu$ depends only on $t$ and that $\upmu$ is small and decreasing to $0$
along the integral curves of $\Lunit = \frac{\partial}{\partial t}.$
Then $\upmu^{-1} \left|\Lunit \upmu \right| = \frac{\partial}{\partial t} \ln \upmu^{-1}.$
Hence, applying Gronwall's inequality to 
\eqref{E:GRONWALLREADYINTROTOPORDERENERGYIDCARICATURE}
and ignoring the terms $\cdots,$
we would deduce that 
\begin{align} \label{E:INTROENERGYCARICATUREESTIMATE}
	\enzero[\Rot \Psi](t,u)
	+ \flzero[\Rot \Psi](t,u)
	& \lesssim 
		\mathring{\upepsilon}^2
		\upmu^{-4}(t).
\end{align}
Roughly, this drastically simplified argument captures the main reason 
that the top-order estimate \eqref{E:INTROTOPORDERENERGYESTIMATE}
is degenerate. The reason that the power of $\upmu_{\star}^{-1}$ in 
\eqref{E:INTROTOPORDERENERGYESTIMATE} is larger than the power suggested by inequality \eqref{E:INTROENERGYCARICATUREESTIMATE}
is that we have ignored
the presence of some additional difficult error integrals
that arise in the analysis of
\eqref{E:GRONWALLREADYINTROTOPORDERENERGYIDCARICATURE}
and also in the analysis of the Morawetz energy-flux quantities $\enone$ and $\flone.$
We also see why it is important that the structural constant $4$ 
does not change as we repeatedly commute the equations; if the constant rapidly grew,
then the power of $\upmu_{\star}^{-1}$ on the right-hand side
of the $L^2$ estimates \eqref{E:INTROTOPORDERENERGYESTIMATE}-\eqref{E:INTROJUSTBELOWTOPORDERENERGYESTIMATE} etc.
would also grow, and we would have no reason to expect that we could eventually descend
to a non-degenerate level \eqref{E:INTROORDERNONDEGENERATEMIDORDERENERGYESTIMATE}
and close the whole process.

\subsection{The behavior of \texorpdfstring{$\frac{\Lunit \upmu}{\upmu}$}{the critical Gronwall factor}}
\label{SS:INTROLUPMUOVERUPMU}
It is not easy to rigorously derive an estimate of the form
\eqref{E:INTROENERGYCARICATUREESTIMATE} from inequality \eqref{E:GRONWALLREADYINTROTOPORDERENERGYIDCARICATURE}.
Because the estimate is at the heart of 
our sharp classical lifespan theorem
(Theorem~\ref{T:LONGTIMEPLUSESTIMATES}),
we now provide some additional details.
As we will see, the most serious complication is:
$\upmu(t,u,\vartheta)$ is highly $(u,\vartheta)-$dependent,
even though it is approximately monotonic at fixed $(u, \vartheta).$ 
In particular, $\upmu$
can either logarithmically grow or shrink in $t$
depending on $(u,\vartheta).$
Consequently, as we explain below,
in order to derive suitable estimates for $\upmu$ that account for 
its full spectrum of possible behaviors,
we must complement our a priori estimates with \emph{a posteriori estimates}.

The key ingredient driving 
the approximate monotonicity of $\upmu$ is its 
smaller-than-expected weighted acceleration along the integral curves of
$\Lunit.$ More precisely, we have the following better-than-expected estimate:
\begin{align} \label{E:INTROLUNITUPMU}
	\left|
		\Lunit (\rgeo \Lunit \upmu)
	\right|
	& \lesssim 
		\mathring{\upepsilon}
		\frac{\ln(\myexp + t)}{(1 + t)^2}, 
\end{align}
where the right-hand side is integrable in $t.$
A naive approach based only on the dispersive-type bootstrap assumption \eqref{E:HP}
would lead to the non-integrable rate $(1 + t)^{-1}$
on the right-hand side of \eqref{E:INTROLUNITUPMU}, 
which is not sufficient for deriving any of our main results.
The source of the gain is that the
only potentially slowly decaying term on the right-hand side of
the equation $\Lunit (\rgeo \Lunit \upmu) = \cdots$ 
(see equation \eqref{E:INTROLUNITUPMUHEURISTIC})
is $G_{\Lunit \Lunit} \Lunit(\rgeo \Rad \Psi),$
and furthermore, for solutions to $\square_{g(\Psi)} \Psi = 0,$
$\Lunit(\rgeo \Rad \Psi)$ decays at an integrable-in-time rate
(see Lemma~\ref{L:FASTERTHANEXPECTEDPSIDECAY}).
This fact is familiar from the case of the linear wave equation in three spatial dimensions,
whose solutions $\Psi$ verify
$\left|\Lunit_{(Flat)}(r \Radunit_{(Flat)} \Psi)\right| \lesssim (1 + t)^{-2}.$
By twice integrating \eqref{E:INTROLUNITUPMU} along 
the integral curves of $\Lunit$ from
the initial data hypersurface to a late time $t$ 
and using\footnote{Throughout Ch.~\ref{C:INTRO}, we often write ``$A \sim B$'' to imprecisely indicate
that $A$ is well-approximated by $B.$} 
$\upmu(0,\cdot) \sim 1,$
we see that
for times $0 \leq s \leq t,$
$\Lunit \upmu$ and $\upmu$ can be heuristically 
modeled by
\begin{subequations}
\begin{align}
	\Lunit \upmu(s,u,\vartheta) 
	& \sim 
		- \updelta_{t,u,\vartheta} 
		\frac{1}{\rgeo(s,u)},
		\label{E:INTROLUPMUMODEL} \\
	\upmu(s,u,\vartheta)
	& \sim 1 
		- \updelta_{t,u,\vartheta} 
			\ln \left( \frac{\rgeo(s,u)}{\rgeo(0,u)} \right),
			\label{E:INTROUPMUMODEL}
\end{align}
\end{subequations}
where $\updelta_{t,u,\vartheta} := - [\rgeo \Lunit \upmu](t,u,\vartheta).$
We stress that the estimates 
\eqref{E:INTROLUPMUMODEL}-\eqref{E:INTROUPMUMODEL}
are not purely a priori in nature;
they involve the quantity
$\updelta_{t,u,\vartheta},$
which is a posteriori in nature.
For the sake of illustration, we imagine that $\updelta_{t,u,\vartheta} \geq 0$
for all $t,u,\vartheta$ and set
$\updelta_t := \sup_{(u',\vartheta) \in [0,u] \times \mathbb{S}^2} \updelta_{t,u',\vartheta}.$ 
This corresponds to a scenario in which $\upmu$ is shrinking along 
all integral curves of $\Lunit$ at a curve-dependent rate.
We then find that the difficult factor on the right-hand 
side of \eqref{E:GRONWALLREADYINTROTOPORDERENERGYIDCARICATURE}
can be bounded as follows:
\begin{align} \label{E:INTROMODELDIFFICULTINTEGRATINGFACTOR}
\int_{t'=0}^t
	\sup_{\Sigma_{t'}^u}
	\left|
		\frac{\Lunit \upmu}{\upmu}
	\right|						
\, dt'
\leq 
\int_{s=0}^t
	\frac{\updelta_t }
	{\rgeo(s,u) \left\lbrace 1 - \updelta_t \ln \left( \frac{\rgeo(s,u)}{\rgeo(0,u)} \right) \right\rbrace}					
\, ds
= \ln 
	\left\lbrace 
		1 - \updelta_t \ln \left( \frac{\rgeo(t,u)}{\rgeo(0,u)} \right)	
	\right\rbrace
\sim \ln \upmu_{\star}^{-1}(t,u).
\end{align}
We stress that the important feature that allows us to arrive at the right-hand side 
of \eqref{E:INTROMODELDIFFICULTINTEGRATINGFACTOR}
is that the \emph{same constant} $\updelta_t$ appears in the both the
numerator and denominator in the $\, ds$ integral.
Clearly, the bound \eqref{E:INTROMODELDIFFICULTINTEGRATINGFACTOR}, when
combined with inequality \eqref{E:GRONWALLREADYINTROTOPORDERENERGYIDCARICATURE} and Gronwall's inequality,
leads to the desired estimate
$\enzero[\Rot \Psi](t,u)
	+ \flzero[\Rot \Psi](t,u)
	\lesssim 
		\mathring{\upepsilon}
		\upmu_{\star}^{-4}(t,u).$
In summary, the only feature of \eqref{E:INTROTOPORDERENERGYESTIMATE} 
not captured by our heuristic analysis is the 
growth factor $\ln^{\Cononestar}(\myexp + t);$ this factor, 
which is not very important because it is present only in the top-order estimates,
appears because along some integral curves of $\Lunit,$ 
$\upmu$ can grow and 
the ratio $\frac{\Lunit \upmu}{\upmu}$ 
can decay at the slow rate $\lbrace (1 + t) \ln(\myexp + t) \rbrace^{-1}$
(see inequality \eqref{E:POSITIVEPARTOFLMUOVERMUISSMALL}).

\subsection{The descent scheme for the below-top-order \texorpdfstring{$L^2$}{square integral} estimates}
If we used the strategy of Sect.~\ref{SS:INTROTOPORDEERENERGYESTIMATES} to estimate
$\enzero[\cdot]$ and $\flzero[\cdot]$ at all orders, 
then all of our $L^2$ estimates would experience the same $\upmu_{\star}^{-1}$
degeneracy as the top-order estimate \eqref{E:INTROTOPORDERENERGYESTIMATE}.
If this were the best we could do, then there would be no hope of 
improving the bootstrap assumptions \eqref{E:HP} through Sobolev embedding, 
and our proof would fall apart.
Hence, it is extremely important 
that the below-top-order $L^2$ estimates 
\eqref{E:INTROJUSTBELOWTOPORDERENERGYESTIMATE}-\eqref{E:INTROORDER0ORDERENERGYESTIMATE}
become less degenerate in $\upmu_{\star}^{-1}$ as we descend.
The main idea behind deriving the improved estimates is that at the lower levels, we can avoid 
using modified quantities and simply allow the loss of one derivative. We then need
to make sure that this procedure somehow results in a gain of a power of
$\upmu_{\star}.$ The reasons that we can in fact gain are 
\textbf{i)} we have sharp information about how $\upmu_{\star}$ behaves
(in the difficult case where $\upmu$ is shrinking, we roughly have the caricature estimate \eqref{E:INTROUPMUMODEL}) 
and \textbf{ii)} by avoiding the use of modified quantities, we see a large gain in the powers of $t$ available 
in the dangerous error integrand on the right-hand side of \eqref{E:INTROTOPORDERENERGYIDCARICATURE}.
To illustrate the scheme, let us imagine that the third-order derivatives of $\Psi$ 
are top-order and that we are trying to deduce improved $L^2$ estimates for the second-order 
derivatives of $\Psi.$ For simplicity, we imagine that we have proved the 
following top-order $L^2$ estimates for the third-order derivatives, 
which are consistent with \eqref{E:INTROTOPORDERENERGYESTIMATE}
up to the unimportant factor $\ln^{\Cononestar}(\myexp + t):$
\begin{align} \label{E:INTROTOPORDERENERGYCARICATUREBOUND}
	\enzero^{1/2}[\mathscr{Z}^2 \Psi](t,u)
	+ \flzero^{1/2}[\mathscr{Z}^2 \Psi](t,u)
	& \lesssim \mathring{\upepsilon} \upmu_{\star}^{-\Contwo}(t,u),
		\\
	\enone^{1/2}[\mathscr{Z}^2 \Psi](t,u)
	+ \flone^{1/2}[\mathscr{Z}^2 \Psi](t,u)
	& \lesssim \mathring{\upepsilon} \upmu_{\star}^{-\Contwo}(t,u),
		\label{E:INTROTOPORDERMORAWETZCARICATUREBOUND}
\end{align}
where $\Contwo > 1.$ We now sketch a proof of 
how to use 
\eqref{E:INTROTOPORDERENERGYCARICATUREBOUND}-\eqref{E:INTROTOPORDERMORAWETZCARICATUREBOUND}
to derive estimates for the just-below-top-order quantities
$\enzero[\Rot \Psi](t,u) + \flzero[\Rot \Psi](t,u)$
that are similar to 
\eqref{E:INTROTOPORDERENERGYCARICATUREBOUND}
but with the exponent $\Contwo$ \emph{reduced by one}, as in \eqref{E:INTROJUSTBELOWTOPORDERENERGYESTIMATE}.

To derive the desired estimate, we estimate the error integral 
$- \int_{\mathcal{M}_{t,u}}
				(2 \Rad \Psi)(\Rot \mytr \upchi) \Rad \Rot \Psi
			\, d \vol$
on the right-hand
side of \eqref{E:INTROTOPORDERENERGYIDCARICATURE} in a different way.
To proceed, we commute \eqref{E:TRANSPORTROTTRCHIDERIVATIVELOSS}
with $\rgeo^2$ to derive\footnote{Recall that $\rgeo \approx 1 + t$ in the region of interest.} the equation\footnote{The point  
of including the factor $\rgeo^2$ is that $\rgeo^2 \Rot \mytr \upchi$ verifies a good equation because
$\rgeo^2$ leads to the cancellation of a linear term; 
see the proof of Lemma~\ref{L:POINTWISEESTIMATESFORLDERIVATIVESOFCHIJUNKINTERMSOFOTHERVARIABLES}.}
\begin{align} \label{E:WITHRGEOFACTORTRANSPORTROTTRCHIDERIVATIVELOSS}
	\Lunit (\rgeo^2 \Rot \mytr \upchi)
	& = \frac{1}{2} \rgeo^2 G_{\Lunit \Lunit} \angLap \Rot \Psi
		+ \cdots,
\end{align}
use the following property (see \eqref{E:ANGLAPFUNCTIONPOINTWISEINTERMSOFROTATIONS})
of our rotation vectorfields (familiar from the case of Minkowski spacetime):
\begin{align} \label{E:INTROANGLAPPSIINTERMSOFANGDIFFPSI}
	\left|
		\rgeo^2 \angLap \Psi
	\right|
	& \lesssim (1 + t) \sum_{l=1}^3 |\angdiff \Rot_{(l)} \Psi|,
\end{align}
and use the coerciveness property \eqref{E:INTROENONECOERCIVENESS}
to deduce the schematic inequality
\begin{align} \label{E:INTROLRGEOSQUAREDROTTRCHIINL2}
	\left\|
		\Lunit (\rgeo^2 \Rot \mytr \upchi)
	\right\|_{\Sigma_t^u}
		\lesssim 
		(1 + t) 
		\left\| 
			\angdiff \Rot \Psi
		\right\|_{\Sigma_t^u}
		+ \cdots
	& \lesssim 
		\upmu_{\star}^{-1/2}(t,u)
		\enone^{1/2}[\Rot \Psi](t,u)
		+ \cdots
			\\
	& \lesssim 
		\mathring{\upepsilon} \upmu_{\star}^{-\Contwo-1/2}(t,u)
		+ \cdots.
		\notag
\end{align}
Note that we have generated an additional factor $\upmu_{\star}^{-1/2}$
on the right-hand side of \eqref{E:INTROLRGEOSQUAREDROTTRCHIINL2}
because of the factor $\sqrt{\upmu}$ in \eqref{E:INTROENONECOERCIVENESS}.
It is not too difficult 
(see Lemma~\ref{L:L2NORMSOFTIMEINTEGRATEDFUNCTIONS})
to integrate
\eqref{E:INTROLRGEOSQUAREDROTTRCHIINL2} in time 
(recall that $\Lunit = \frac{\partial}{\partial t}$)
and to account for the fact that the $S_{t,u}$ area form
$d \spherevol$ is pointwise of size $\rgeo^2$
to deduce that
\begin{align} \label{E:INTRORGEOSQUAREDROTTRCHIINL2FIRSTESTIMATE}
	\| \rgeo^2 \Rot \mytr \upchi \|_{\Sigma_t^u}
	& \lesssim
	\mathring{\upepsilon} (1 + t) 
	\int_{t'=0}^t  
		\frac{1}{(1 + t') \upmu_{\star}^{\Contwo + 1/2}(t',u)} 
	\, dt'
		+ \cdots.
\end{align}

The important point concerning 
the time integral in \eqref{E:INTRORGEOSQUAREDROTTRCHIINL2FIRSTESTIMATE}
is that we can prove the following estimate, 
which allows us to \underline{gain a power of $\upmu_{\star}$ by integrating in time},
at the expense of introducing a factor $\ln(\myexp + t),$
which is harmless below top order because of the good weight $\rgeo^2$ available
on the left-hand side of \eqref{E:INTRORGEOSQUAREDROTTRCHIINL2FIRSTESTIMATE}:
\begin{align} \label{E:INTROUPMUGAINSTHROUGHTIMEINTEGRATION}
	\int_{t'=0}^t 
		\frac{1}{(1 + t') \upmu_{\star}^{\Contwo + 1/2}(t',u)} 
	\, dt'
	& \lesssim \ln(\myexp + t) \upmu_{\star}^{-\Contwo + 1/2}(t,u).
\end{align}
The estimate \eqref{E:INTROUPMUGAINSTHROUGHTIMEINTEGRATION} 
can be proved by using arguments similar to but much simpler
than the ones we used in sketching the proof of \eqref{E:INTROMODELDIFFICULTINTEGRATINGFACTOR};
see inequality \eqref{E:LOGLOSSKEYMUINTEGRALBOUND} and its proof.
%
Inserting \eqref{E:INTROUPMUGAINSTHROUGHTIMEINTEGRATION} into
\eqref{E:INTRORGEOSQUAREDROTTRCHIINL2FIRSTESTIMATE}
and using $\rgeo \approx 1 + t,$ 
we find that
\begin{align} \label{E:INTRORGEOSQUAREDROTTRCHIINL2SECONDESTIMATE}
	\| \Rot \mytr \upchi \|_{\Sigma_t^u}
	& \lesssim \frac{\ln(\myexp + t)}{1 + t} \upmu_{\star}^{-\Contwo + 1/2}(t,u)
		+ \cdots.
\end{align}

We now insert the estimate \eqref{E:INTRORGEOSQUAREDROTTRCHIINL2SECONDESTIMATE}
into the error integral on the right-hand side of \eqref{E:INTROTOPORDERENERGYIDCARICATURE}
and use Cauchy-Schwarz, 
the bootstrap assumption \eqref{E:INTROTRANSVERSALC0BOUND},
and the coerciveness estimate \eqref{E:INTROMULTENERGYCOERCIVITY}
to deduce
\begin{align} \label{E:INTROKEYDESCENTSPACETIMEINTEGRALESTIMATE}
	\left|
		\int_{\mathcal{M}_{t,u}}
			2 (\Rad \Psi) (\Rot \mytr \upchi) \Rad \Rot \Psi
		\, d \vol
	\right|
	& \lesssim 
		\mathring{\upepsilon}
		\int_{t'=0}^t
			\frac{1}{1 + t'}
			\| \Rot \mytr \upchi \|_{\Sigma_{t'}^u}
			\| \Rad \Rot \Psi \|_{\Sigma_{t'}^u}
		\, dt'
			\\
	& \lesssim 
		\mathring{\upepsilon}^2
		\int_{t'=0}^t
			\frac{\ln(\myexp + t')}{(1 + t')^2 \upmu_{\star}^{\Contwo-1/2}(t',u)}
			\enzero^{1/2}[\Rot \Psi](t',u)
		\, dt'
		+ \cdots.
		\notag
\end{align}
Hence, from \eqref{E:INTROTOPORDERENERGYIDCARICATURE}
and \eqref{E:INTROKEYDESCENTSPACETIMEINTEGRALESTIMATE}, we have
\begin{align} \label{E:GRONWALLREADYBELOWTOPORDERENERGYIDCARICATURE}
	& \sup_{s \in [0,t]}
		\left\lbrace
			\enzero[\Rot \Psi](s,u)
			+ \flzero[\Rot \Psi](s,u)
		\right\rbrace
		\\
	& \leq C \mathring{\upepsilon}^2
		+ 
		\sup_{s \in [0,t]}
		\enzero^{1/2}[\Rot \Psi](s,u)
		\mathring{\upepsilon}^2
		\int_{t'=0}^t
			\frac{\ln(\myexp + t')}{(1 + t')^2 \upmu_{\star}^{\Contwo-1/2}(t',u)}
		\, dt'
		+ \cdots.
		\notag
\end{align}

We now use inequality \eqref{E:GRONWALLREADYBELOWTOPORDERENERGYIDCARICATURE}
and an argument similar to the one used to deduce \eqref{E:INTROUPMUGAINSTHROUGHTIMEINTEGRATION}
in order to gain a power of $\upmu_{\star}$ through time integration,
but we now take advantage of the additional decay available in the
integrand of \eqref{E:GRONWALLREADYBELOWTOPORDERENERGYIDCARICATURE}
to avoid the $\ln(\myexp + t)$ factor in \eqref{E:INTROUPMUGAINSTHROUGHTIMEINTEGRATION}
(see \eqref{E:LESSSINGULARTERMSMUINTEGRALBOUND} for the precise estimate).
We thus conclude the desired reduction in the power of $\upmu_{\star}^{-1}:$
\begin{align} \label{E:INTROWEGAINEDATLOWERORDER}
	\enzero^{1/2}[\Rot \Psi](t,u)
	+ \flzero^{1/2}[\Rot \Psi](t,u)
	& \lesssim 
		\mathring{\upepsilon}
		\upmu_{\star}^{-\Contwo + 3/2}(t,u)
		+ \cdots.
\end{align}
We remark that the estimate \eqref{E:INTROWEGAINEDATLOWERORDER}
is slightly misleading in the sense that it suggests
that we can gain $\upmu_{\star}^{3/2}$ at each stage in the descent.
In reality, there are additional error integrals on the 
right-hand side of \eqref{E:INTROTOPORDERENERGYIDCARICATURE}
that only allow us to gain
$\upmu_{\star},$ 
and hence we have the hierarchy
\eqref{E:INTROTOPORDERENERGYESTIMATE}-\eqref{E:INTROORDER0ORDERENERGYESTIMATE}.

\section{Comparison with related work}
\label{S:COMPARISONWITHRELATEDWORK}
We now compare our work with that of Alinhac and Christodoulou.

\subsection{Alinhac's shock-formation results}
Alinhac obtained some important results that significantly 
advanced our understanding of singularity formation in 
solutions to wave equations. Specifically,
he was the first to understand the nature of the first singularity 
in small-data solutions to a class of quasilinear wave equations 
\cites{sA1999a,sA1999b,sA2001b,sA2002}.
As we describe below, he proved 
finite-time shock formation when the data are
small, compactly supported, and verify 
some non-degeneracy conditions
that generically hold.\footnote{For equations that are invariant under Euclidean rotations,
some small data containing a spherically symmetric angular sector
do not verify Alinhac's non-degeneracy assumptions.}
He studied equations in both two and three spatial dimensions, 
but we discuss here only the case of three spatial dimensions,
which, despite the title, is addressed in \cite{sA2001b}.
We now briefly summarize Alinhac's results. His work applied to 
equation \eqref{E:ONEDERIVATIVEQUASILINEARWAVE} 
under the assumptions \eqref{E:LITTLEHDECOMPOSED}-\eqref{E:LITTLEHVANISHESATZERO}. 
That is, he studied the Cauchy problem\footnote{Following the conventions of \cite{dC2007}, 
we denote the dynamic metric by $h = h(\partial \Phi)$ in this section and the next one.}
\begin{align}
	(h^{-1})^{\alpha \beta}(\partial \Phi) \partial_{\alpha} \partial_{\beta} \Phi 
	& = 0,
		\label{E:ALINHACWAVE} \\
	(\Phi|_{t=0},\partial_t \Phi|_{t=0}) 
	& = (\mathring{\Phi}, \mathring{\Phi}_0).
		\label{E:ALINHACDATA}
\end{align}
More precisely, he studied solutions corresponding to a one-parameter family 
of compactly supported data of the form $(\uplambda \mathring{\Phi}, \uplambda \mathring{\Phi}_0),$
where $\uplambda \in (0,1)$ was chosen to be sufficiently small
(and the amount of smallness needed depends on $(\mathring{\Phi}, \mathring{\Phi}_0)$).
The correct analog of the future null condition failure factor \eqref{E:INTROFAILFACT}
for equation \eqref{E:ALINHACWAVE} is
\begin{align} \label{E:INTRONEWNULLCONDITIONFAILUREFACTOR}
		\FutFailFac
		&:=
		m_{\kappa \lambda} H_{\mu \nu}^{\kappa}(\partial \Phi = 0) 
		\Lunit_{(Flat)}^{\lambda} \Lunit_{(Flat)}^{\mu} \Lunit_{(Flat)}^{\nu},
	\end{align}
where $m_{\mu \nu}$ is the Minkowski metric
and relative to rectangular coordinates, we have
\begin{align} \label{E:INTROBIGH}
	H_{\mu \nu}^{\lambda}(\partial \Phi)
	:= \frac{\partial}{\partial (\partial_{\lambda} \Phi)} h_{\mu \nu}^{(Small)}(\partial \Phi).
\end{align}
In particular, when $\FutFailFac \equiv 0,$ Klainerman's classic null condition \cite{sk1984} is verified and 
his work \cite{sK1986} (or Christodoulou's work \cite{dC1986a}) yields small-data global existence.
As in the case of \eqref{E:INTROFAILFACT},
the quantity $\FutFailFac$ defined in \eqref{E:INTRONEWNULLCONDITIONFAILUREFACTOR} 
can be viewed as a function 
depending only on $\theta = (\theta^1,\theta^2),$ 
where $\theta^1$ and $\theta^2$ are local angular 
coordinates corresponding to standard spherical coordinates 
on Minkowski spacetime.

	Alinhac's results partially confirm a conjecture of Fritz John,
	where by ``partially,'' we mean
	that his results require non-degeneracy assumptions on the data.
	To further explain these assertions, we first state the definition of the Radon
	transform of a function $f$ on $\mathbb{R}^3.$
	Given any such $f$ and points 
	$q \in \mathbb{R},$ $\theta \in \mathbb{S}^2 \subset \mathbb{R}^3,$ 
	we define 
	\begin{align} \label{E:RADONTRANSFORMOFF}
		\mathcal{R}[f](q,\theta)
		& := 
			\int_{P_{q,\theta}} 
				f(y) 
			\, d \sigma_{q,\theta}(y),
	\end{align}
	where $P_{q,\theta} := \lbrace y \in \mathbb{R}^3 \ | \ \Euct(\theta,y) = q \rbrace,$
	$d \sigma_{q,\theta}(y)$ denotes the area form induced on the plane $P_{q,\theta}$ by the 
	Euclidean metric $\Euct$ on $\mathbb{R}^3,$ and $\Euct(\theta,y)$ is the Euclidean inner
	product of $\theta$ and $y$ (where both are viewed as vectors in $\mathbb{R}^3,$ the former being unit length).
	Alinhac's condition for shock formation
	involves the following function of $(q,\theta),$ 
	which also depends on the data $(\mathring{\Phi}, \mathring{\Phi}_0):$
	\begin{align} \label{E:FRIEDNALNDERRADIATIONFIELD}
		\Fried[(\mathring{\Phi}, \mathring{\Phi}_0)]
		(q,\theta)
		& := 
			-
			\frac{1}{4 \pi} 
				\frac{\partial}{\partial q}
				\mathcal{R}[\mathring{\Phi}](q,\theta)
			+
			\frac{1}{4\pi} 
			\mathcal{R}[\mathring{\Phi}_0](q,\theta).
	\end{align}
	The function $\Fried[(\mathring{\Phi}, \mathring{\Phi}_0)]$
	is Friedlander's radiation field for the solution to the 
	\emph{linear} wave equation 
	$\square_m \Phi = 0$
	corresponding to the data
	$(\mathring{\Phi}, \mathring{\Phi}_0).$
	That is, the $r-$weighted linear solution $r \Phi_{(Linear)},$
	relative to standard spherical coordinates $(t,r,\theta)$
	on Minkowski spacetime,
	is asymptotic to $\Fried[(\mathring{\Phi}, \mathring{\Phi}_0)](q=r-t,\theta)$
	in the region of interest.
	The relevance of $\Fried$ lies in the fact that 
	$\Fried[(\mathring{\Phi}, \mathring{\Phi}_0)](q=r-t,\theta)$
	provides a good approximation 
	to the \emph{nonlinear} solution
	at time $\mathring{\upepsilon}^{-1}$
	(where $\mathring{\upepsilon}$ is the small size of the data),
	long before any singularity has formed.
	Roughly speaking, 
	at time $\mathring{\upepsilon}^{-1},$
	for a class of small-data solutions that form shocks,
	the dangerous term that eventually causes the shock formation
	has acquired the ``shock-driving'' sign, at least along some
	integral curves of the outgoing null vectorfield $\Lunit.$ 
	It takes a certain amount of time to
	see the sign of the dangerous term because one must wait for
	certain $\mathcal{C}_u-$tangent derivatives to sufficiently decay
	before its sign becomes visible;
	see Sect.~\ref{SS:INTROSHOCKFORMINGDATACOMPARISON} 
	for additional discussion.
	Hence, for a class of small-data solutions,
	whether or not shock formation occurs
	can be detected from the state of the dangerous term
	relatively early in the evolution, 
	at time $\mathring{\upepsilon}^{-1}.$ Furthermore, 
	the state of the dangerous term 
	at time $\mathring{\upepsilon}^{-1}$
	can in turn be determined from
	the data with the help of
	$\Fried[(\mathring{\Phi}, \mathring{\Phi}_0)](q=r-t,\theta).$
	
	A connection between $\Fried$ and the lifespan of the
	solution to the nonlinear equation \eqref{E:ALINHACWAVE} was first observed by
	John \cite{fJ1987} and
	H\"{o}rmander \cite{lH1987}. In these works, they 
	proved that for initial data of the form
	$(\uplambda \mathring{\Phi}, \uplambda \mathring{\Phi}_0),$
	the classical lifespan $T_{(Lifespan);\uplambda}$ of the solution 
	to equation \eqref{E:ALINHACWAVE} verifies
	\begin{align} \label{E:JOHNHORMANDERLOWERBND}
		\liminf_{\uplambda \downarrow 0}
		\uplambda \ln T_{(Lifespan);\uplambda}
		\geq 
		\frac{1}
		{\sup_{(q,\theta) \in \mathbb{R} \times \mathbb{S}^2}
		\frac{1}{2} 
		\FutFailFac(\theta) 
		\frac{\partial^2}{\partial q^2} \Fried[(\mathring{\Phi}, \mathring{\Phi}_0)]
		(q,\theta)}.
	\end{align}
	John conjectured that inequality \eqref{E:JOHNHORMANDERLOWERBND} is sharp 
	and that in the limit $\uplambda \downarrow 0,$ 
	the lower bound also serves as an upper bound.
	He made significant progress towards proving his conjecture
	\cites{fJ1989,fJ1990} by showing that near the expected blow-up
	time, the second rectangular derivatives of the solution start to grow. 
	
	We now motivate John's conjecture and the result \eqref{E:JOHNHORMANDERLOWERBND}.
	We first note that an important aspect of the works \cite{fJ1987} and \cite{lH1987} 
	is that they show that the dispersive decay rates corresponding to the linear wave equation are verified
	by small-data solutions to the nonlinear equation \eqref{E:ALINHACWAVE} until 
	very close to the conjectured singularity time.
	Their proof was based on Klainerman's Minkowskian vectorfield method, 
	the Minkowskian eikonal function $r - t,$ 
	and Friedlander's radiation field $\Fried$ for the linear wave equation,
	whose role we further explain below.
	Hence, one can obtain a good approximation to equation \eqref{E:ALINHACWAVE},
	valid until near the conjectured first singularity time,
	by expanding it relative to the Minkowskian frame \eqref{E:MINKOWSKIFRAME}
	and keeping only the quadratic term that fails the classic null condition
	(see the discussion in Sect.~\ref{S:WHENNULLCONDITIONFAILS})
	and the related linear term that drives its evolution.
	Also using the approximation $\frac{r}{t} \sim 1,$
	we find that a suitable approximate equation is as follows:
	\begin{align} \label{E:MOTIVATINGJOHNSCONJECTUREEQN}
		\Lunit_{(Flat)} (r \Radunit_{(Flat)} \Phi)
		+ \frac{1}{2} 
			\frac{1}{t} 
			\FutFailFac (r \Radunit_{(Flat)} \Phi) \Radunit_{(Flat)} (r \Radunit_{(Flat)} \Phi)
		& = 0.
	\end{align}
	The approximation $\frac{r}{t} \sim 1$ used 
	in motivating \eqref{E:MOTIVATINGJOHNSCONJECTUREEQN}
	is accurate at large times in the Minkowskian wave zone $|r - t| \approx 1,$ 
	which, in the case of compactly supported data, is the region
	where one expects the first singularity to form.
	Note also that in obtaining \eqref{E:MOTIVATINGJOHNSCONJECTUREEQN}, 
	we have discarded Euclidean angular derivatives 
	(which decay quickly, at least until near the shock).
	Equation \eqref{E:MOTIVATINGJOHNSCONJECTUREEQN} is a Burgers-type
	equation in the unknown $r \Radunit_{(Flat)} \Phi$ 
	whose solutions can blow up along the integral curves
	of $\Lunit_{(Flat)}.$ From the point of view of understanding the blow-up
	of solutions to the nonlinear wave equation \eqref{E:ALINHACWAVE},
	the most relevant data for 
	the approximating equation \eqref{E:MOTIVATINGJOHNSCONJECTUREEQN}
	is not $r \Radunit_{(Flat)}\Phi|_{t=0}.$
	The reason is that equation \eqref{E:MOTIVATINGJOHNSCONJECTUREEQN} 
	does not take into account the influence of
	some of the linear terms in equation \eqref{E:ALINHACWAVE}, such as the
	Euclidean angular derivatives, which can be influential in the early phase of the dynamics.
	As we explain in more in more detail in Sect.~\ref{SS:INTROSHOCKFORMINGDATACOMPARISON},
	this early phase lasts roughly until time $\uplambda^{-1}$ 
	(recall that $\uplambda$ can roughly be viewed as the size of the data).
	Moreover, the ``data'' induced at time $\uplambda^{-1}$ by the approximately 
	linear evolution is effectively determined, up to small errors, 
	by $\Fried$ as follows
	relative to standard spherical coordinates $(t,r,\theta)$
	on Minkowski spacetime:
	\begin{align} \label{E:APPROXIMATEDATA}
	r \Radunit_{(Flat)}\Phi(\uplambda^{-1},r,\theta)
	\sim
	- \FutFailFac(\theta) \frac{\partial}{\partial q} \Fried[(\mathring{\Phi}, \mathring{\Phi}_0)](r - \uplambda^{-1},\theta).
	\end{align}
	By modifying the argument we used to prove blow-up for the standard Burgers equation \eqref{E:BURGERS},
	we compute that the solution to the model equation \eqref{E:MOTIVATINGJOHNSCONJECTUREEQN} corresponding to the data \eqref{E:APPROXIMATEDATA}
	blows up along the integral curve of $\Lunit_{(Flat)}$ emanating from the point with Minkowskian spherical coordinates $(\uplambda^{-1},r,\theta)$
	by the time $\sim \exp \left\lbrace \frac{1}{2} \FutFailFac(\theta) \frac{\partial^2}{\partial q^2} \Fried[(\mathring{\Phi}, \mathring{\Phi}_0)](r - 		
		\uplambda^{-1},\theta)
	\right\rbrace,$ whenever the term inside the exponential is positive.
	Hence, the lifespan of the solution to equation \eqref{E:MOTIVATINGJOHNSCONJECTUREEQN} 
	corresponding to the data \eqref{E:APPROXIMATEDATA} is well-approximated by 
	the right-hand side of \eqref{E:JOHNHORMANDERLOWERBND}.
	
	It is easy to see that for compactly supported data,
	the right-hand of \eqref{E:JOHNHORMANDERLOWERBND} is $\geq 0.$
	It is natural to wonder if there any nontrivial compactly supported 
	data such that the right-hand side of \eqref{E:JOHNHORMANDERLOWERBND} is equal
	to $0.$ The results of John \cite{fJ1987} imply that such data would lead to a global solution.
	However, John showed \cite{fJ1987}*{pg. 98} that when $\FutFailFac \not \equiv 0,$ there are no such data. 
	We state his important observation as a proposition and sketch its simple proof.
	
 \begin{proposition} [\textbf{Only trivial data cause the right-hand side of \eqref{E:JOHNHORMANDERLOWERBND} to vanish}]
		\label{P:JOHNSCRITERIONISALWAYSSATISFIEDFORCOMPACTLYSUPPORTEDDATA}
		Let $\mathring{\Phi}, \mathring{\Phi}_0 \in C_c^{\infty}(\mathbb{R}^3).$ 
		Assume that $\FutFailFac \not \equiv 0$ and that
		\begin{align} \label{E:JOHNSQUANTITYVANISHES}
		\sup_{(q,\theta) \in \mathbb{R} \times \mathbb{S}^2}
		\FutFailFac(\theta) 
		\frac{\partial^2}{\partial q^2} \Fried[(\mathring{\Phi}, \mathring{\Phi}_0)]
		(q,\theta)
		& = 0.
		\end{align}
		Then $(\mathring{\Phi}, \mathring{\Phi}_0) = (0,0).$
	\end{proposition}
	
	\begin{proof}[Sketch of proof]
	We will deduce from the assumptions of the proposition that
	$\Fried[(\mathring{\Phi}, \mathring{\Phi}_0)] \equiv 0.$
	Then using the identity
	$\Fried[(\mathring{\Phi}, \mathring{\Phi}_0)](q,\theta) \pm \Fried[(\mathring{\Phi}, \mathring{\Phi}_0)](-q,-\theta) = 0,$
	where $- \theta$ denotes the polar opposite of $\theta$ on $\mathbb{S}^2,$
	it is straightforward to see that $\mathring{\Phi}$ and $\mathring{\Phi}_0$ have vanishing Radon transforms.
	The proposition then follows from the fact that
	a compactly supported function with a vanishing Radon 
	transform must completely vanish.
	
	To show that $\Fried[(\mathring{\Phi}, \mathring{\Phi}_0)] \equiv 0,$
	we partition $\mathbb{S}^2$ into an open set $\mathcal{O}$ and a closed set $\mathcal{C}$ as follows:
	\begin{align}
		\mathbb{S}^2 & = \mathcal{O} \cup \mathcal{C}, &
			\\
		\mathcal{O} & := \lbrace \theta \in \mathbb{S}^2 \ | \ \FutFailFac(\theta) \neq 0 \rbrace,
		\qquad \mathcal{C} := \lbrace \theta \in \mathbb{S}^2 \ | \ \FutFailFac(\theta) = 0 \rbrace.
	\end{align} 
	For a fixed $\theta \in \mathcal{O},$
	we deduce from \eqref{E:JOHNSQUANTITYVANISHES}
	that $\frac{\partial^2}{\partial q^2} \Fried[(\mathring{\Phi}, \mathring{\Phi}_0)](q,\theta)$
	is either non-positive or non-negative for $q \in \mathbb{R}.$
	Since $\Fried[(\mathring{\Phi}, \mathring{\Phi}_0)](q,\theta)$ vanishes for large $|q|$
	(this fact is essentially the sharp version of Huygens' principle), 
	it is straightforward to see that $\Fried[(\mathring{\Phi}, \mathring{\Phi}_0)](q,\theta) = 0$ for $q \in \mathbb{R}.$
	We next note that since $\FutFailFac$ is a nontrivial analytic function on $\mathbb{S}^2,$ 
	$\mathcal{C}$ must have an empty interior and thus $\mathcal{C} = \partial \mathcal{C}.$
	Hence, for fixed $\theta \in \mathcal{C},$ we conclude from continuity and the property
	$\Fried[(\mathring{\Phi}, \mathring{\Phi}_0)]|_{\mathbb{R} \times \mathcal{O}} \equiv 0$
	that $\Fried[(\mathring{\Phi}, \mathring{\Phi}_0)](q,\theta) = 0$ for $q \in \mathbb{R}.$
	We have thus shown that $\Fried[(\mathring{\Phi}, \mathring{\Phi}_0)] \equiv 0$ as desired.
	\end{proof}
	
	We now summarize the main features of Alinhac's results
	as Theorem~\ref{T:ALINHAC}. The results stated in the theorem
	are a partial summary of Theorems $2$ and $3$ of
	\cite{sA2001b} in the case of three spatial dimensions.
	We stress that \emph{his main contribution was confirming that John's conjecture holds for
	a large set of data verifying uniqueness and non-degeneracy conditions};
 	see just below equation \eqref{E:ALINHACDATABLOWUPFUNCTION}.
	
	\begin{theorem} [\textbf{Alinhac}]
	\label{T:ALINHAC}
	Let $(\uplambda \mathring{\Phi}, \uplambda \mathring{\Phi}_0)
	\in C^{\infty}(\mathbb{R}^3) \times C^{\infty}(\mathbb{R}^3)$ 
	be a one-parameter family of 
	compactly supported data for the quasilinear wave equation \eqref{E:ALINHACWAVE}
	under the assumptions \eqref{E:LITTLEHDECOMPOSED}-\eqref{E:LITTLEHVANISHESATZERO}.
	Let $\Phi_{\uplambda}$ denote the corresponding solution.
	Assume that the (data-dependent) function
	\begin{align} \label{E:ALINHACDATABLOWUPFUNCTION}
		\frac{1}{2} 
		\FutFailFac(\theta)
		\frac{\partial^2}{\partial q^2} \Fried[(\mathring{\Phi}, \mathring{\Phi}_0)]
		(q,\theta)
	\end{align}
	has a unique, strictly positive, non-degenerate maximum at $(q_*,\theta_*),$
	where the future null condition failure factor $\FutFailFac$ is defined in \eqref{E:INTRONEWNULLCONDITIONFAILUREFACTOR}.
	If $\uplambda$ is sufficiently small and positive 
	(where the amount of smallness needed depends on $(\mathring{\Phi}, \mathring{\Phi}_0)$),
	then we have the following conclusions.
	
	\medskip
	
	\noindent \underline{\textbf{Sharp classical lifespan result.}}
	The classical lifespan $T_{(Lifespan);\uplambda}$ of $\Phi_{\uplambda}$
	is finite and verifies
	\begin{align} \label{E:ALINHACLIMITINGLIFESPAN}
		\lim_{\uplambda \downarrow 0}
			\uplambda \ln T_{(Lifespan);\uplambda}
			& =  \frac{1}
					{\frac{1}{2} 
		\FutFailFac(\theta_*) 
		\frac{\partial^2}{\partial q^2} \Fried[(\mathring{\Phi}, \mathring{\Phi}_0)](q_*,\theta_*)}.
	\end{align}
	
	\medskip
	
	\noindent \underline{\textbf{Sharp description near the unique first blow-up point.}}
	Along the constant-time hypersurface $\Sigma_{T_{(Lifespan);\uplambda}}$ 
	of first shock formation, there is a unique point $p_{(Blow-up);\uplambda}$
	where the solution blows up. In addition, with $C = C[(\mathring{\Phi}, \mathring{\Phi}_0)],$
	$\Phi_{\uplambda}$ and its first rectangular derivatives verify the bound
	\begin{align}
			\left|
				\Phi_{\uplambda}
			\right|
			+
			\sum_{\alpha=0}^3
			\left|
				\partial_{\alpha} \Phi_{\uplambda}
			\right| 
			& \leq C
				\uplambda \frac{1}{1 + t},
	\end{align}
	and in particular, they remain finite along $\Sigma_{T_{(Lifespan);\uplambda}}.$ 
	In the \textbf{complement} of a small neighborhood of $p_{(Blow-up);\uplambda}$
	intersected with $\lbrace t \leq T_{(Lifespan);\uplambda} \rbrace,$
	the second-order rectangular derivatives
	of $\Phi_{\uplambda}$ verify a similar bound.
	In contrast, the following blow-up behavior occurs for $t$ close to and $\leq T_{(Lifespan);\uplambda}:$
			\begin{align} \label{E:ALINHACSOLUTIONREMAINSC2BLOWUP}
				C^{-1} \left(t \ln \frac{T_{(Lifespan);\uplambda}}{t} \right)^{-1}
				\leq
				\sum_{\alpha,\beta=0}^3 
				\left\|
					\partial_{\alpha} \partial_{\beta} \Phi_{\uplambda}
				\right\|_{C^0(\Sigma_t)}
				& \leq C \left(t \ln \frac{T_{(Lifespan);\uplambda}}{t} \right)^{-1}.
			\end{align}
	Furthermore, in a small neighborhood of $p_{(Blow-up);\uplambda}$
	intersected with $\lbrace t \leq T_{(Lifespan);\uplambda} \rbrace,$
	there exist an eikonal function $u$ 
	and a related system of geometric coordinates, 
	one of which is $u,$
	such that relative to the geometric coordinates, 
	$\Phi_{\uplambda}$ and its higher derivatives extend smoothly to 
	$\Sigma_{T_{(Lifespan);\uplambda}},$ even at the point
	$p_{(Blow-up);\uplambda}.$
	Within this region,
	the change of variables map from geometric to rectangular coordinates is smooth and invertible
	except at the point $p_{(Blow-up);\uplambda},$ where its Jacobian determinant vanishes.
\end{theorem}

We now highlight two merits of Alinhac's results. First, he was the first to show 
that failure of the classic null condition in equation \eqref{E:ALINHACWAVE} often leads, 
in the small-data regime, 
to the formation of a singularity caused by the intersection of the level sets
of a true eikonal function. 
A second merit is that, for reasons explained below, his proof is relatively short.

The main limitation of Alinhac's results is that his framework 
is fundamentally tailored to the first point of intersection of the 
characteristics and thus it cannot easily be extended 
to provide information about the maximal development of the data.
We further discuss the origin of this limitation two paragraphs below.
In contrast, as we describe in Sects.~\ref{SS:INTROSUMMARYCHRISTODOULOUWORK}
and \ref{SS:INTROSHOCKFORMINGDATACOMPARISON}, 
Christodoulou's approach and the approach of the present monograph
allow one to obtain a detailed description of the maximal future development of the 
data in the exterior of the sphere $S_{0,U_0} \subset \Sigma_0^{U_0}.$
This sharp information is essential for attacking the problem of 
extending the solution, in a generalized sense, beyond the shock.
Furthermore, Alinhac's results are fundamentally based on
his non-degeneracy assumptions on the data,
which are stated just below \eqref{E:ALINHACDATABLOWUPFUNCTION}
and which cause the solution to form a shock.
Thus, his results do not provide a true analog 
of the sharp classical lifespan theorem
proved in \cite{dC2007} and in Theorem~\ref{T:LONGTIMEPLUSESTIMATES} of the present monograph, 
which imply that singularities in small-data solutions can only be caused by the vanishing of $\upmu.$
Moreover, his assumptions are not verified by some data 
containing a spherically symmetric sector
when the wave equation is invariant under Euclidean rotations.

We now sketch some of the main ideas behind Alinhac's approach, which is
quite different than the approach of \cite{dC2007} and the present monograph.
As we will see in the next paragraph, the most important part of his analysis is based on a type of 
``backwards approach,'' as opposed to the ``forwards approach'' associated to solving a pure Cauchy problem,
which was adopted in \cite{dC2007} and the present monograph.
To proceed with our sketch, 
we first recall John's conjecture that the lifespan lower bound proved in
\cite{fJ1987} and \cite{lH1987} is also an upper bound in the small-data limit.
Alinhac's proof is grounded in the belief that 
one should expand the actual lifespan $T_{(Lifespan);\uplambda}$ 
of the solution to \eqref{E:ALINHACWAVE}-\eqref{E:ALINHACDATA}
as the lifespan predicted by John's conjecture plus 
a small error to be solved for. 
Thus, he splits the nonlinear evolution into two stages.
In the first stage, he simply quotes the
results that we described above when motivating
John's conjecture via the approximating equation \eqref{E:MOTIVATINGJOHNSCONJECTUREEQN}.
That is, Alinhac uses the results of \cite{fJ1987} and \cite{lH1987},
which were based on Klainerman's Minkowskian vectorfield method, 
the Minkowskian eikonal function $r - t,$ 
and Friedlander's radiation field for the linear wave equation,
to follow the solution almost all the way
to the lifespan predicted by John's conjecture.

In the second stage, near the very end of the solution's classical lifespan,
Alinhac constructs, via an iteration scheme based on Nash-Moser-type estimates,
a true eikonal function that enables him to follow the solution all the way 
to the first singularity. The ``initial'' data for this stage are
of course inherited from the state of the solution at the end of the first stage.
Because the blow-up time $T_{(Lifespan);\uplambda}$ 
is not precisely known, his proof involves
a free boundary (because the upper boundary of the domain is $T_{(Lifespan);\uplambda}$). 
A related difficulty is that the lifespan
can slightly vary from iterate to iterate, so that there 
is an extra step involved to fix the domain.
The reason that Alinhac's proof is so short is that he avoids
using many of the intricate geometric structures
present in Christodoulou's framework.
For example, his iteration scheme produces a true eikonal function
for the nonlinear wave equation only in the limit;
the iterates themselves involve approximate eikonal functions
whose level sets are null hypersurfaces for 
the metric evaluated along the \emph{previous} iterate.
Furthermore, his commutation vectorfields
also depend on the previous iterate.
The net effect is that his linearized equations and commutation vectorfields
do not have the same good structure enjoyed by the nonlinear equations studied
in \cite{dC2007} and the present monograph, and thus he does not recover
(or need) the same sharp top-order $L^2$ estimates for the eikonal function
required by those works. Consequently, Alinhac's energy estimates
for the iterates lose derivatives relative to the previous iterate.
Nonetheless, he is able to close his proof by deriving tame $L^2-$type estimates
for his linearized equations and using the Nash-Moser framework.
Alinhac's iteration scheme is fundamentally based on his condition ``(H);''
see \cite{sA1999b}*{pg.15}. Roughly speaking, condition (H)
requires that each iterate has characteristics that intersect at exactly one
point belonging to its constant-time hypersurface of first blow-up.
In particular, for spherically symmetric data in the case of equations 
invariant under the Euclidean rotations,
condition (H) already fails for the zeroth iterate
and hence the next iterate cannot be constructed.
For similar reasons, his framework does not allow one to study constant-time 
hypersurfaces $\Sigma_t$ containing a submanifold along which $\upmu$ is $0,$
as can happen when $\Sigma_t$ lies to the future
of the constant-time hypersurface of first blow-up;
see the ``singular part'' of the boundary of the maximal development
described in Theorem~\ref{T:CHRISTODOULOUSHOCKFORMATION}.
In total, \emph{the restrictive nature of condition (H) is the reason 
that Alinhac's framework applies only
to data that verify his uniqueness and non-degeneracy assumptions and that
it does not reveal the complete structure of the maximal development of the data.}

\subsection{Christodoulou's shock formation results}
\label{SS:INTROSUMMARYCHRISTODOULOUWORK}
In \cite{dC2007} Christodoulou proved several landmark results
in relativistic fluid mechanics; see also the work \cite{dCsM2012}
of Christodoulou-Miao for the same results proved in the case of non-relativistic fluid mechanics.
Specifically, Christodoulou proved analogs of our 
Theorem~\ref{T:LONGTIMEPLUSESTIMATES}
and
Theorem~\ref{T:STABLESHOCKFORMATION}
for all of the physically relevant scalar quasilinear wave equations that arise in irrotational 
relativistic fluid mechanics in Minkowski spacetime. 
In addition, his work went somewhat beyond these
two theorems in that he also gave a detailed description of the maximal future development
of the data given in the exterior of the sphere $S_{0,U_{0}} \subset \Sigma_0^{U_0}$
whenever $0 \leq U_0 < 1/2.$
Moreover, he extended his result to cover a class of small fluid 
data that have non-vanishing vorticity. However, the main aspects of his work
addressed only a region in which the fluid is irrotational. 
In the irrotational region, the relativistic Euler (fluid) equations 
reduce to a scalar quasilinear wave equation for a potential
function $\Phi.$ Again, we stress that the difficult part of his argument was
his proof of an analog of Theorem~\ref{T:LONGTIMEPLUSESTIMATES},
and that the remaining aspects his work are easier to derive.
We also stress that although Alinhac had already proved his small-data shock formation results 
(summarized in Theorem~\ref{T:ALINHAC})
for a larger class of equations,\footnote{One has to take into account the trivial differences in normalization, 
noted below, in order
to see that Christodoulou's equations fall under the scope of Alinhac's work.}
there was great novelty in Christodoulou's thoroughness of his description of the
dynamics and in particular, in his description of the solution along
the boundary of the maximal development of the data.
A particularly attractive feature of Christodoulou's detailed description is 
that it is suitable as a starting point for
trying to extend the solution, in a generalized sense, beyond the shock.

We now provide some details. Christodoulou's wave equations 
were the Euler-Lagrange equations for
Lagrangians $\mathcal{L}$ of the form
$\mathcal{L} = \mathcal{L}(\upsigma),$ 
where 
\begin{align} \label{E:ENTHALPHYSQUARED}
	\upsigma := - (m^{-1})^{\alpha \beta} \partial_{\alpha} \Phi \partial_{\beta} \Phi
\end{align}
and $(m^{-1})^{\alpha \beta} = \mbox{diag}(-1,1,1,1)$ is the standard 
reciprocal Minkowski metric.
As is explained in \cite{dC2007}, 
in order to obtain a relativistic fluid interpretation from 
$\mathcal{L}(\upsigma),$ it suffices to make the 
following five positivity assumptions:
\begin{align} \label{E:LAGRANGIANPOSITIVITY}
	\upsigma,
		\,
	\mathscr{L}(\upsigma), 
		\,
	\frac{d \mathscr{L}}{d \upsigma},
		\,
	\frac{d}{d \upsigma}\left(\mathscr{L}/ \sqrt{\upsigma} \right),
		\,
	\frac{d^2\mathscr{L}}{d \upsigma^2} > 0. 
\end{align}
The assumptions \eqref{E:LAGRANGIANPOSITIVITY}
imply that $\Phi$ can be interpreted as
a potential function for an irrotational relativistic fluid with physically desirable
properties such as having a positive pressure,
a speed of sound in between $0$ and $1,$ etc.
The Euler-Lagrange equation corresponding to $\mathcal{L}(\upsigma)$ is
\begin{align} \label{E:CHRISTFLUIDEL}
		\partial_{\alpha} \left( \frac{\partial \mathcal{L}}{\partial (\partial_{\alpha} \Phi)} \right)
		= - 2 \partial_{\alpha} \left( \frac{\partial \mathcal{L}}{\partial \upsigma} (m^{-1})^{\alpha \beta} \partial_{\beta} \Phi \right)
		& = 0,
\end{align}
and the background solutions to \eqref{E:CHRISTFLUIDEL}
of interest are $\Phi = kt,$ where $k$ is a non-zero constant.
These are the solutions that correspond to non-zero constant fluid 
states.\footnote{Perturbations of the non-zero constant states are much easier to study
than compact perturbations of the vacuum state.
The reason is that the Euler equations become very degenerate along vacuum boundaries.}
Relative to Minkowski-rectangular coordinates, equation \eqref{E:CHRISTFLUIDEL} 
can be expressed as
\begin{align}
	(h^{-1})^{\alpha \beta}(\partial \Phi) \partial_{\alpha} \partial_{\beta} \Phi & = 0,
\end{align}
where the \emph{reciprocal acoustical metric} $h^{-1}$ is defined by
\begin{align}
	(h^{-1})^{\alpha \beta}(\partial \Phi)
	& =(m^{-1})^{\alpha \beta}
		- F
			(m^{-1})^{\alpha \kappa} (m^{-1})^{\beta \lambda} 
			\partial_{\kappa} \Phi \partial_{\lambda} \Phi,
				\\
	F = F(\upsigma) 
	& := \frac{2}{G} \frac{d G}{d \upsigma},
		\\
	G = G(\upsigma) 
	& := 2 \frac{d \mathcal{L}}{d \upsigma}.
\end{align}

Because the background solutions are non-zero, there are some
superficial differences in the
way that Christodoulou's solutions look compared to our solutions and those
of Alinhac. The differences essentially correspond to different normalization choices
and are not of fundamental importance.
For example, the propagation speed corresponding to  
the background solution is not $1$ as in our work and that of Alinhac, but
is instead
\begin{align}
	\upeta_0 := \upeta(\upsigma = k^2),
\end{align}
where $\upeta > 0$ is the function of $\upsigma$ defined by
\begin{align}
	\upeta^2 & = \upeta^2 (\upsigma) := 1 - \upsigma H,
		\\
	H & = H(\upsigma) := \frac{F}{1 + \upsigma F}.
\end{align}
Using \eqref{E:LAGRANGIANPOSITIVITY}, 
it is straightforward to show that the
\emph{speed of sound} $\upeta$ verifies
$0 < \upeta < 1.$
Christodoulou did not assume that $(h^{-1})^{00}$ is equal to $-1$ as 
we did in \eqref{E:GINVERSE00ISMINUSONE},
but instead there is a dynamic lapse function $\upalpha > 0$ defined by
\begin{align}
	\upalpha^{-2}
	= \upalpha^{-2}(\partial \Phi)
	& := - (h^{-1})^{00}(\partial \Phi).
\end{align}
From the point of view of the behavior of small perturbations of the background solutions,
the relevant inverse background metric
is not the standard inverse Minkowski metric in the form $(m^{-1})^{\mu \nu} = \mbox{diag}(-1,1,1,1),$ 
but is instead a flat inverse metric that takes the following form relative
to the rectangular coordinates:
$(h^{-1})^{\alpha \beta}(\partial_t \Phi = k, \partial_1 \Phi= \partial_2 \Phi= \partial_3 \Phi = 0).$
More precisely, we have 
\begin{align}
	h(\partial_t \Phi = k, \partial_1 \Phi= \partial_2 \Phi= \partial_3 \Phi = 0)
= - \upeta_0^2 dt^2 + \sum_{a=1}^3 (dx^a)^2.
\end{align}
The eikonal function corresponding to the background solution is
\begin{align}
	u_{(Flat)} & = 1 - r + \upeta_0 t,
\end{align}
where $r = \sqrt{\sum_{a=1}^3 (x^a)^2}.$ 
The inverse foliation density corresponding to the background solution is
\begin{align}
	\upmu_{(Flat)} & = \upeta_0.
\end{align}
The outgoing and ingoing null vectorfields corresponding to the background solution are
\begin{align}
	\Lunit_{(Flat)} 
	& = \partial_t + \upeta_0 \partial_r,
	&&
	\uLgood_{(Flat)} 
	= \upeta_0^{-1} \partial_t - \partial_r.
\end{align}
The analog of the future null condition failure factor \eqref{E:INTRONEWNULLCONDITIONFAILUREFACTOR} 
is
\begin{align} \label{E:CHRISTODOULOUNULLCONDITIONFAILUREFACTOR}
	\frac{d H}{d \upsigma}(\upsigma = k^2).
\end{align}
Unlike in the general case of \eqref{E:INTRONEWNULLCONDITIONFAILUREFACTOR},
the quantity in
\eqref{E:CHRISTODOULOUNULLCONDITIONFAILUREFACTOR} is a constant.
Christodoulou showed that 
$\frac{d H}{d \upsigma}(\upsigma = k^2)$ vanishes for all non-zero $k$ 
if and only if, up to 
trivial normalization constants,
\begin{align} \label{E:EXCEPTIONALLAGRANGIAN}
	\mathcal{L}(\upsigma) = 1 - \sqrt{1 - \upsigma}.
\end{align}
The Lagrangian \eqref{E:EXCEPTIONALLAGRANGIAN} is therefore the only 
member of the above family of Lagrangians such that
the quadratic nonlinearities that arise in 
expanding its Euler-Lagrange equation \eqref{E:CHRISTFLUIDEL} 
around the background solutions $\Phi = kt$ verify Klainerman's classic null condition.

We now summarize Christodoulou's results. 
The results stated below as Theorem~\ref{T:CHRISTODOULOUSHOCKFORMATION} 
are a conglomeration of
\cite{dC2007}*{Theorem 13.1 on pg. 888, Theorem 14.1 on pg. 903, Proposition 15.3 on pg. 974,
and the Epilogue on pg. 977}.
The quantities such as $\Lunit$, $\uLgood,$ etc. that appear in the theorem
are essentially the same as the quantities that we use throughout our monograph, up to 
the differences in normalization pointed out in the previous paragraph.

\begin{theorem}[\textbf{Christodoulou}]
	\label{T:CHRISTODOULOUSHOCKFORMATION}
	Let $\upsigma$ be as defined in \eqref{E:ENTHALPHYSQUARED}.
	Assume that the Lagrangian $\mathcal{L}(\upsigma)$ verifies 
	the positivity conditions
	\eqref{E:LAGRANGIANPOSITIVITY}
	in a neighborhood of $\upsigma = k^2,$ where $k$ is a non-zero constant, but that
	$\mathcal{L}(\upsigma)$ is not the exceptional Lagrangian
	\eqref{E:EXCEPTIONALLAGRANGIAN}.
	Consider the following Cauchy problem for the quasilinear (Euler-Lagrange) wave equation
	corresponding to $\mathcal{L},$ expressed relative to rectangular coordinates:
	\begin{align} \label{E:CHRISTODOULOUWAVE}
		\partial_{\alpha} \left( \frac{\partial \mathcal{L}}{\partial (\partial_{\alpha} \Phi)} \right)
		& = 0, 
			\\
		(\Phi|_{t=0},\partial_t \Phi|_{t=0}) & = (\mathring{\Phi}, \mathring{\Phi}_0).
	\end{align}
	Assume that the data are small perturbations of the data corresponding to 
	the non-zero constant-state solution
	$\Phi = k t$ and that the perturbations are compactly
	supported in the Euclidean unit ball.
	Let $U_0 \in (0,1/2)$ and let
	\begin{align}
		\mathring{\upepsilon}
		= \mathring{\upepsilon}[(\mathring{\Phi}, \mathring{\Phi}_0)]
		:=
		\| \mathring{\Phi}_0 - k \|_{H^N(\Sigma_0^{U_0})}
		+ \sum_{i=1}^3 \| \partial_i \mathring{\Phi} \|_{H^N(\Sigma_0^{U_0})}
	\end{align}
	denote the size of the data,
	where $N$ is a sufficiently large integer.\footnote{A numerical value of $N$ was not provided in \cite{dC2007}.}
	
	\medskip
	
	\noindent \underline{\textbf{Sharp classical lifespan result.}}
	If $\mathring{\upepsilon}$ is sufficiently small, 
	then a sharp classical lifespan theorem in analogy with
	Theorem~\ref{T:LONGTIMEPLUSESTIMATES} holds.
	
	\medskip
	
	\noindent \underline{\textbf{Small-data shock formation.}}
	We define the following data-dependent functions of $u|_{\Sigma_0} = 1 - r:$
	\begin{align}
	& \mathcal{E}[(\mathring{\Phi}, \mathring{\Phi}_0)](u)
		\\
	& \ \ 
		:= \sum_{\Psi \in \lbrace \partial_t \Phi - k, \partial_1 \Phi, \partial_2 \Phi, \partial_3 \Phi \rbrace} 
		\int_{\Sigma_0^u}
					\left\lbrace
						\upalpha^{-2} \upmu
						(\upeta_0^{-1} + \upalpha^{-2} \upmu)
						(\Lunit \Psi)^2
						+ (\uLgood \Psi)^2
						+ (\upeta_0^{-1} + 2 \upalpha^{-2} \upmu) \upmu |\angdiff \Psi|^2
					\right\rbrace
			\, d \tvol,
			\notag
	\end{align}	
	\begin{align} \label{E:CHRISTODOULOUDATAFUNDTION}
		\mathcal{S}[(\mathring{\Phi}, \mathring{\Phi}_0)](u)
		& :=
		\int_{S_{0,u}}
			r 
			\left\lbrace
				(\mathring{\Phi}_0 - k)
				- \upeta_0 \partial_r \mathring{\Phi}
			\right\rbrace
		 \, d \Eucspherevol
			+	
			\int_{\Sigma_0^u}
			 \left\lbrace
				2 (\mathring{\Phi}_0 - k)
				- \upeta_0 \partial_r \mathring{\Phi}
			\right\rbrace
		 \, d^3 x,
	\end{align}
	where 
	$d \tvol$ is defined in \eqref{E:INTROVOLFORMS},
	$d \Eucspherevol$ denotes the Euclidean area form on the sphere $S_{0,u}$ of 
	Euclidean radius $r = 1 - u,$ 
	and $d^3 x$ denotes the standard flat volume form on $\mathbb{R}^3.$
	Assume that 
	\begin{align}
		\ell := \frac{d H}{d \upsigma}(\upsigma = k^2) > 0. 
	\end{align}	
	There exist constants $C > 0$ and $C' > 0,$ 
	independent of $U \in (0,U_0],$
	such that if $\mathring{\upepsilon}$ is sufficiently small and
	if for some $U \in (0,U_0]$ we have
	\begin{align} \label{E:SHOCKFUNCTIONMUSTBESUFFICIENTLYLARGE}
		\mathcal{S}[(\mathring{\Phi}, \mathring{\Phi}_0)](U)
		& \leq - C \mathring{\upepsilon} \mathcal{E}^{1/2}[(\mathring{\Phi}, \mathring{\Phi}_0)](U) < 0,
	\end{align}
	then a shock forms in the solution\footnote{That is, $\Phi$ and its first rectangular derivatives remain bounded, while
	some second-order rectangular derivative blows up due to the vanishing of $\upmu.$} 
	$\Phi$ 
	and the first shock in the maximal development of the portion of the data
	in the exterior of $S_{0,U} \subset \Sigma_0^U$
	originates in the hypersurface region $\Sigma_{T_{(Lifespan)};U}^U$ 
	(see Def.~\ref{D:HYPERSURFACESANDCONICALREGIONS}), where 
	\begin{align}
		T_{(Lifespan);U} < \exp\left( C' \frac{U}{\left|k^3 \ell \mathcal{S}[(\mathring{\Phi}, \mathring{\Phi}_0)](U)\right|} \right).
	\end{align}
	A similar result holds if $\ell < 0;$ in this case, we delete the 
	``$-$'' sign in \eqref{E:SHOCKFUNCTIONMUSTBESUFFICIENTLYLARGE} and change 
	``$\leq$'' and ``$<$'' to ``$\geq$'' and ``$>.$''
	
	\medskip
	
	\noindent \underline{\textbf{Description of the boundary of the maximal development}.}
	For a large set of small-data shock-forming solutions verifying some technical non-degeneracy conditions,
	the boundary $\mathcal{B}$ of the maximal development of the data in the exterior of $S_{0,U} \subset \Sigma_0^U$
	is a union $\mathcal{B} = (\partial_- \mathcal{H} \cup \mathcal{H}) \cup \underline{\mathcal{C}},$
	where $\partial_- \mathcal{H} \cup \mathcal{H}$ is the singular part and 
	$\underline{\mathcal{C}}$ is the regular part; see Figure~\ref{F:MAXDEVBOUND}.
	$\upmu$ vanishes along the singular part and is positive on the regular part, and 
	the solution and its rectangular derivatives up to a certain order 
	extend continuously in rectangular coordinates to the regular part.
	Each component of $\partial_- \mathcal{H}$ is a smooth $2-$dimensional embedded submanifold of Minkowski spacetime,
	spacelike with respect to the dynamic metric $h.$ The corresponding component of 
	$\mathcal{H}$ is a smooth, embedded, $3-$dimensional
	submanifold in Minkowski spacetime ruled by $h-$null curves
	with past endpoints on $\partial_- \mathcal{H}.$ The corresponding component 
	$\underline{\mathcal{C}}$ is the incoming null hypersurface
	corresponding to $\partial_- \mathcal{H},$ and it is ruled by incoming $h-$null geodesics 
	with past endpoints on $\partial_- \mathcal{H}.$
\end{theorem}

\begin{center}
\begin{overpic}[scale=1]{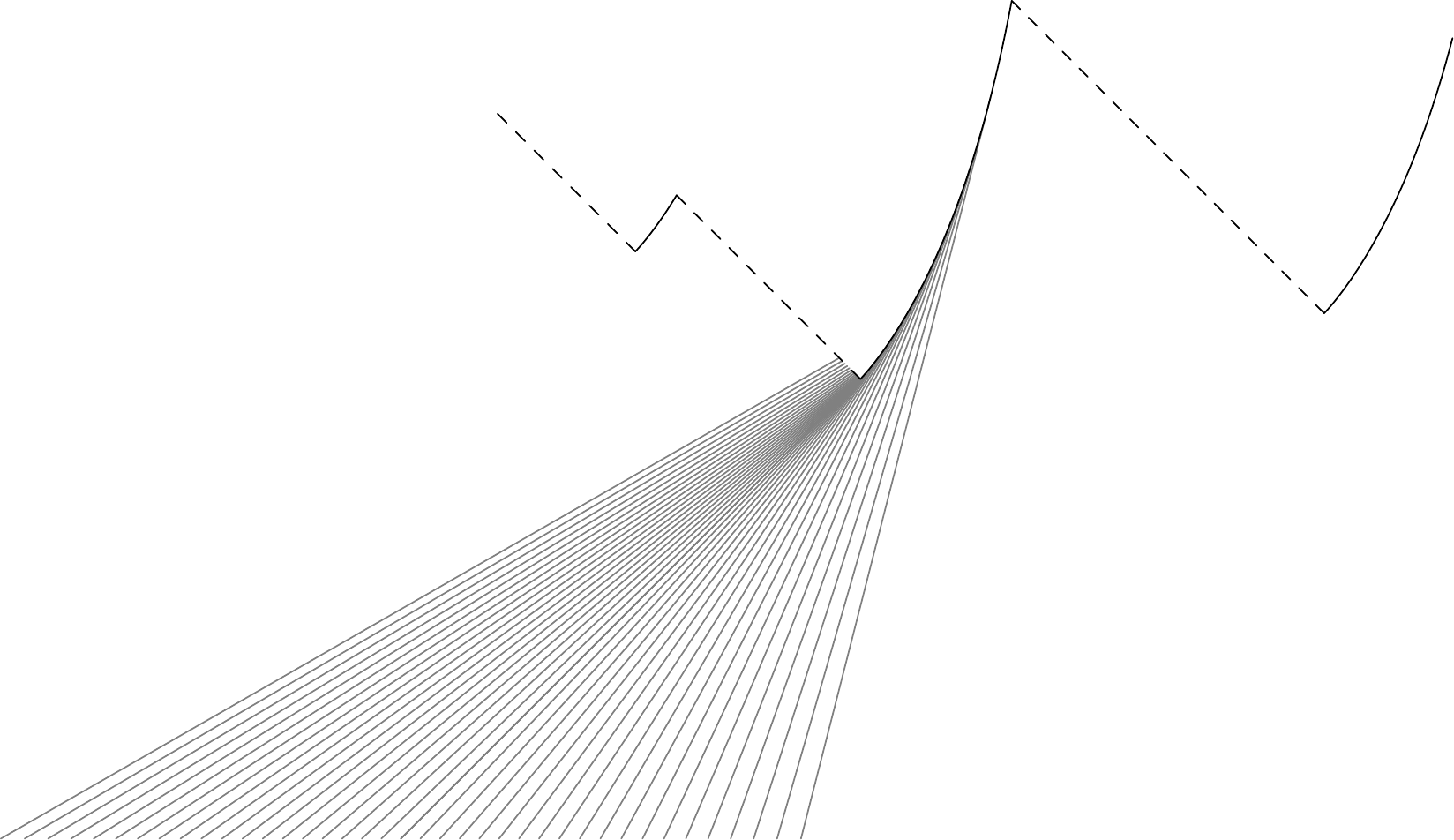}
        \put (62,44) {\large$\mathcal{H}$}
	\put (93,44) {\large$\mathcal{H}$}
	\put (88,33) {\large$\partial_-\mathcal{H}$}
	\put (78,44) {\large$\underline{\mathcal{C}}$}
\end{overpic}
\captionof{figure}{A cross section of the maximal development. The gray lines
are level sets of the eikonal function $u$ near the first
blow-up point. The dotted lines are the regular
boundary \underline{$\mathcal{C}$}. 
The black curves are the boundary $\mathcal{H},$ 
whose past endpoints are $\partial_- \mathcal{H}.$} 
 \label{F:MAXDEVBOUND}
\end{center}

Again, we stress that a key advantage of Christodoulou's approach is that
the estimates of his sharp classical lifespan theorem allow 
him to understand the maximal development of the data in
the exterior of the sphere $S_{0,U_{0}} \subset \Sigma_0^{U_0},$
including the behavior of the solution along the boundary.
We also mention again that, based on the estimates of Theorem~\ref{T:LONGTIMEPLUSESTIMATES}, 
our analysis of shock-forming solutions to equations \eqref{E:WAVEGEO} and \eqref{E:ONEDERIVATIVEQUASILINEARWAVE} 
could be extended, without too much additional effort,
to provide an analogous description of the maximal development.

\subsection{Further discussion on the shock-forming data}
\label{SS:INTROSHOCKFORMINGDATACOMPARISON}
We now compare our conditions on the data for shock formation
with those of Alinhac and Christodoulou.
We first note that there is no obvious connection between any of the
three classes of shock-forming data identified in
Theorems~\ref{T:ALINHAC},
\ref{T:CHRISTODOULOUSHOCKFORMATION},
and \ref{T:STABLESHOCKFORMATION}.
However, the previous discussion
essentially shows that among the three theorems, 
Alinhac's treats the largest set of data:
recall that John conjectured
(see the discussion just below \eqref{E:JOHNHORMANDERLOWERBND}), 
based on the results of his Prop.~\ref{P:JOHNSCRITERIONISALWAYSSATISFIEDFORCOMPACTLYSUPPORTEDDATA}, 
that finite-time blow-up should occur 
in all solutions to equation \eqref{E:ALINHACWAVE}
corresponding to nontrivial small compactly supported data, 
and Alinhac was able to partially confirm the conjecture
by proving it under some non-degeneracy assumptions on the data
(see just below \eqref{E:ALINHACDATABLOWUPFUNCTION}).
As we describe below, the results of the present monograph can be
extended to allow us to complete the proof
of John's conjecture by eliminating Alinhac's non-degeneracy assumptions. 
Using the same ideas, we can also prove the analogous
result for Christodoulou's wave equations \eqref{E:CHRISTODOULOUWAVE}
and our equations $\square_{g(\Psi)} \Psi = 0.$ 
However, there is a sense in which our results are not complete;
as was the case for Alinhac, in order for our proof of shock formation
to go through, we may need to shrink the amplitude of the data
to an extent that depends on their profile. Hence,
our analysis leaves open the possibility
that there exist nontrivial data small enough such that the sharp classical
lifespan (Theorem~\ref{T:LONGTIMEPLUSESTIMATES}) applies, 
but such that a shock does not form in the solution.

We can also derive a ``localized very long time existence result'' 
for data such that the quantity \eqref{E:ALINHACDATABLOWUPFUNCTION}  
(or its analogs for the other equations)
is non-positive in a suitable annular region. 
For such data, there is a region of spacetime in which 
the term that typically dominates the 
behavior of $\Lunit \upmu$ has a non-negative sign
and thus cannot cause $\upmu$ to decay.
Later in this section, we outline
how to use this fact to prove that in a certain region of spacetime, 
the solution exists beyond the lifespan $\exp\left(\frac{1}{C \mathring{\upepsilon}} \right)$
that holds for general small data
(see inequality \eqref{E:TLIFESPANLOWERBOUND}).

Before addressing the issues noted above, 
we first note that under some structural assumptions on the nonlinearities, 
one can derive an analog of Christodoulou's
shock-formation criterion \eqref{E:SHOCKFUNCTIONMUSTBESUFFICIENTLYLARGE}
for equations \eqref{E:WAVEGEO} and \eqref{E:ONEDERIVATIVEQUASILINEARWAVE}.
Specifically, a sufficient assumption on the nonlinearities would be 
that the future null condition failure factor $\FutFailFac = \FutFailFac(\theta)$ from
\eqref{E:INTROFAILFACT}/\eqref{E:INTRONEWNULLCONDITIONFAILUREFACTOR}
takes on a strictly positive or negative sign for all $\theta \in \mathbb{S}^2.$
Then, thanks to the sharp estimates of Theorem~\ref{T:LONGTIMEPLUSESTIMATES},
an analog of inequality \eqref{E:SHOCKFUNCTIONMUSTBESUFFICIENTLYLARGE}
should be sufficient to allow Christodoulou's arguments leading to his shock-formation theorem \cite{dC2007}*{Theorem 14.1 on pg. 903}
to go through nearly verbatim. His arguments are based on using the estimates of his version of Theorem~\ref{T:LONGTIMEPLUSESTIMATES}
to derive estimates for quantities that are averaged over the spheres $S_{t,u}.$ For example, 
using the estimates of Theorem~\ref{T:LONGTIMEPLUSESTIMATES} and
applying Christodoulou's arguments to solutions of \eqref{E:WAVEGEO} under the assumption of 
an analog of inequality \eqref{E:SHOCKFUNCTIONMUSTBESUFFICIENTLYLARGE},
we could deduce the existence of some \emph{unknown} integral curve of $\Lunit$ along which
we have a large-time lower bound of the form $\pm \Rad \Psi \gtrsim \mathring{\upepsilon} (1 + t)^{-1}$
(we can achieve either sign for $\Rad \Psi,$ depending on our choice of data).
Inserting this estimate into equation \eqref{E:INTROLUNITUPMUHEURISTIC}
and using the fact that $G_{\Lunit \Lunit}(t,u,\vartheta)$ is well-approximated
by $\FutFailFac(t=0,u=0,\vartheta)$
(see Lemma~\ref{L:NULLCONDFACT}),
we could deduce the following large-time estimate:
\begin{align} \label{E:APPLYCHARGUMENTTOOUREQUATIONS}
	\pm \Lunit \upmu(t,u_*,\vartheta_*) 
	\geq c \FutFailFac(t=0,u=0,\vartheta_*) \mathring{\upepsilon} \frac{1}{1 + t} + \mbox{Error},
\end{align}
where $(u_*,\vartheta_*)$ are the coordinates
corresponding to the unknown integral curve,
and for $\mathring{\upepsilon}$ sufficiently small, 
the term $\mbox{Error}$ is dominated by the first
term on the right-hand side of \eqref{E:APPLYCHARGUMENTTOOUREQUATIONS}.
The sign $\pm$ in \eqref{E:APPLYCHARGUMENTTOOUREQUATIONS}
depends on the sign of the analog of the left-hand side of \eqref{E:SHOCKFUNCTIONMUSTBESUFFICIENTLYLARGE}.
Hence, when $\FutFailFac(t=0,u=0,\vartheta)$ takes on a definite sign\footnote{Recall that at $t=0,$
the geometric angular coordinate $\vartheta$ coincides with the standard Euclidean spherical coordinate $\theta.$}  
for all $\vartheta \in \mathbb{S}^2,$ 
the data are sufficiently small,
and the data are chosen to generate the sign ``$-$''
on the left-hand side of \eqref{E:APPLYCHARGUMENTTOOUREQUATIONS},
we can integrate inequality \eqref{E:APPLYCHARGUMENTTOOUREQUATIONS} in time
to deduce that $\upmu$ must vanish in finite time.

We now return to the issue noted in the first paragraph of Sect.~\ref{SS:INTROSHOCKFORMINGDATACOMPARISON}:
sketching a proof that we can use the analog of Theorem~\ref{T:LONGTIMEPLUSESTIMATES}
for equation \eqref{E:ALINHACWAVE}
to eliminate Alinhac's uniqueness and non-degeneracy assumptions on the data
(see just below \eqref{E:ALINHACDATABLOWUPFUNCTION}).
Specifically, we can prove small-data shock formation 
for equation \eqref{E:ALINHACWAVE} by assuming only John's criterion
\begin{align}  \label{E:ALINHACRELAXED}
	\mbox{the function from \eqref{E:ALINHACDATABLOWUPFUNCTION}
	is positive at one point \ } 
	(q_*,\theta_*)
\end{align}
and by perhaps shrinking the amplitude of the data by a small constant factor if necessary,
as Alinhac did in his proof Theorem~\ref{T:ALINHAC}.\footnote{
As will become clear, shrinking the amplitude of the data can help the shock-driving term
dominate error terms.}
Furthermore, by Prop.~\ref{P:JOHNSCRITERIONISALWAYSSATISFIEDFORCOMPACTLYSUPPORTEDDATA}, 
the condition \eqref{E:ALINHACRELAXED} is 
\emph{always satisfied when the data are nontrivial and compactly supported.}
John's criterion \eqref{E:ALINHACRELAXED}
can also be modified in a straightforward fashion
to apply to Christodoulou's equations \eqref{E:CHRISTODOULOUWAVE}
and our equations $\square_{g(\Psi)} \Psi = 0.$
Moreover, the criterion is sharp in a way that we make precise in the next paragraph.
We now sketch proofs of these claims. 
For definiteness, we focus only on the equation $\square_{g(\Psi)} \Psi = 0,$
and we consider data $(\mathring{\Psi}, \mathring{\Psi}_0)$ 
that are compactly supported in the Euclidean unit ball $\Sigma_0^1.$
Throughout this discussion, 
$\mathring{\upepsilon}$ denotes the (small) size of the data $(\mathring{\Psi}, \mathring{\Psi}_0).$
The main idea of the proof is that under
the correct version of John's criterion for the equation $\square_{g(\Psi)} \Psi = 0,$
we can find a point $p$ belonging to a set
of the form $\Sigma_{\mathring{\upepsilon}^{-1}}^{U_0}$ such that 
\emph{along the integral curve} of $\Lunit$ passing through $p,$ 
we have the following estimate for times beyond $\mathring{\upepsilon}^{-1}:$
\begin{align} \label{E:INTROREFINEDLUPMURELATION}
	\Lunit \upmu(t,u,\vartheta) 
	& = 
	-\frac{1}{2} \rgeo^{-1}(t,u) 
	\left\lbrace 
		\FutFailFac 
		\frac{\partial}{\partial q} \Fried[(\mathring{\Psi}, \mathring{\Psi}_0)]
	\right\rbrace|_p
	+ \mbox{Error},
\end{align}
where the constant-valued term 
$
\left\lbrace 
		\FutFailFac 
		\frac{\partial}{\partial q} \Fried[(\mathring{\Psi}, \mathring{\Psi}_0)]
	\right\rbrace|_p
$
in \eqref{E:INTROREFINEDLUPMURELATION} is \emph{positive}
and $\mbox{Error} = o\left(\mathring{\upepsilon} \right) (1 + t)^{-1}.$
Hence, integrating \eqref{E:INTROREFINEDLUPMURELATION}
with respect to $t$ starting from time $\mathring{\upepsilon}^{-1},$
we infer that $\upmu$ will vanish in finite time,
thus yielding shock formation.\footnote{This is the point in the argument 
where we might need to rescale the amplitude of the data in order
to ensure that the first term on the right-hand side of 
\eqref{E:INTROREFINEDLUPMURELATION} dominates the second;
the first term shrinks linearly in the scaling factor, 
while the second shrinks at a superlinear rate.}
To find a viable point $p,$ we first note that,
as will become clear, the quantity
$- \frac{1}{2} \FutFailFac(\theta) 
\frac{\partial}{\partial q} \Fried[(\mathring{\Psi}, \mathring{\Psi}_0)]
	(q = r - \mathring{\upepsilon}^{-1},\theta)$
is the correct analog of the quantity \eqref{E:ALINHACDATABLOWUPFUNCTION} in the present context,
where $\Fried[(\mathring{\Psi}, \mathring{\Psi}_0)]$ denotes Friedlander's radiation field (see \eqref{E:FRIEDNALNDERRADIATIONFIELD}).
Hence, if the data verify John's criterion, that is, they are such that the product
$\FutFailFac(\theta) \frac{\partial}{\partial q} \Fried[(\mathring{\Psi}, \mathring{\Psi}_0)](q,\theta)$	
is positive at some point $(q_*,\theta_*),$ 
then there exists a point $p$ belonging to a set 
of the form\footnote{One minor difficulty 
that would have to be addressed in a complete proof
is that $\Sigma_{\mathring{\upepsilon}^{-1}}^{U_0}$ might correspond to a
value $U_0 > 1,$ which is a set that we have not even defined. 
However, our results can be extended to apply to such sets.
The main task would be to modify the construction of the eikonal
function $u$ so that on $\Sigma_{\mathring{\upepsilon}^{-1}},$
it is allowed to take on some values larger than $1.$
Alternatively, we could avoid extending the construction of the eikonal function
by instead starting with data that are given along $\Sigma_{-1/2}$
and compactly supported in the Euclidean ball of radius $1/2$
centered at the origin.
We adopt this latter approach in Theorem \ref{T:STABLESHOCKFORMATION}.
}  
$\Sigma_{\mathring{\upepsilon}^{-1}}^{U_0}$ such that 
$-\frac{1}{2} \left\lbrace \FutFailFac \frac{\partial}{\partial q} \Fried[(\mathring{\Psi}, \mathring{\Psi}_0)] \right\rbrace|_p$
is negative; this point $p$ will suffice for our purposes.
In view of the evolution equation \eqref{E:INTROLUNITUPMUHEURISTIC}
for $\upmu,$ to conclude \eqref{E:INTROREFINEDLUPMURELATION},
the main task is to show that for times beyond $\mathring{\upepsilon}^{-1},$
the important term $\frac{1}{2} [G_{\Lunit \Lunit} \rgeo \Rad \Psi](t,r,\theta)$
is well-approximated by 
$-\frac{1}{2} \rgeo^{-1}(t,u) \left\lbrace \FutFailFac \frac{\partial}{\partial q} \Fried[(\mathring{\Psi}, \mathring{\Psi}_0)] \right\rbrace|_p$
along the integral curve of $\Lunit$ passing through $p.$
There are two main steps in the proof of the approximation.
The first step is to show that
on the time interval $[0,\mathring{\upepsilon}^{-1}],$
the linear dynamics dominate the behavior of the nonlinear solution
and as a consequence, 
the product
$\frac{1}{2} [G_{\Lunit \Lunit} \rgeo \Rad \Psi](t=\mathring{\upepsilon}^{-1},r,\theta)$
is well-approximated by $- \frac{1}{2} \FutFailFac(\theta) 
\frac{\partial}{\partial q} \Fried[(\mathring{\Psi}, \mathring{\Psi}_0)]
	(q = r - \mathring{\upepsilon}^{-1},\theta).$ 
Here, $(t,r,\theta)$ are standard spherical coordinates on Minkowski spacetime.
We describe this step in more detail two paragraphs below.
To explain the second step, we let $(\mathring{\upepsilon}^{-1},u,\vartheta)$ be the geometric coordinates of $p.$
We now integrate equation \eqref{E:FRAMEWAVEINTRO} along 
the integral curve of $\Lunit$ emanating from $p$
and use the estimates of Lemma~\ref{L:NULLCONDFACT}
to derive refined version of 
\eqref{E:INTROWAVEEQNTRANSPORTINTEGRATEDESTIMATEMIDTIMEANDLATERESTIMATE}
showing that for all times beyond $\mathring{\upepsilon}^{-1},$ 
the product 
$\frac{1}{2} G_{\Lunit \Lunit}(t,u,\vartheta) \Rad \Psi(t,u,\vartheta)$
is well-approximated by 
$-\frac{1}{2} \rgeo^{-1}(t,u) G_{\Lunit \Lunit}(\mathring{\upepsilon}^{-1},u,\vartheta) \Rad \Psi(\mathring{\upepsilon}^{-1},u,\vartheta).$
Note that by \eqref{E:INTROTANGENTIALARERAPIDLYDECAYING} and its higher-order analogs,
for times beyond $\mathring{\upepsilon}^{-1},$
the term $\rgeo \upmu \angLap \Psi$
from equation \eqref{E:FRAMEWAVEINTRO}
has had time to sufficiently decay
and thus it only makes a small
contribution when we integrate 
from $p.$ Combining these two steps, we deduce the 
desired approximation, which completes our sketch of a proof of \eqref{E:INTROREFINEDLUPMURELATION} 


We now sketch the proof of the
``localized very long time existence result'' 
mentioned above, which holds
in regions $\FutFailFac \frac{\partial}{\partial q} \Fried[(\mathring{\Psi}, \mathring{\Psi}_0)]$
is non-positive. To illustrate what we mean, 
we assume that the data are such that 
there exists a $q_1 \in (0,1)$
and an annular region
$\mathcal{N} := \lbrace (q,\theta) \ | \ q \in [q_1,1] \times \mathbb{S}^2 \rbrace$ 
with $- \FutFailFac \frac{\partial}{\partial q} \Fried[(\mathring{\Psi}, \mathring{\Psi}_0)]|_{\mathcal{N}} \geq 0.$ 
We consider the family of integral curves of $\Lunit$ that
start at time $\mathring{\upepsilon}^{-1}$ and 
emanate from the subset of points in $\Sigma_{\mathring{\upepsilon}^{-1}}$
corresponding to\footnote{
We note that $q$ should be thought of as the Minkowski null coordinate $q = r - t.$ 
Thus, the Minkowskian spacetime region corresponding to $\mathcal{N}$ 
is the region between the flat inner cone $\lbrace r - t = q_1 \rbrace$
and the flat outer cone $\mathcal{C}_0 = \lbrace r - t = 1 \rbrace.$}  
$\mathcal{N}.$ Fixing any point $p$ belonging to the subset 
and integrating inequality \eqref{E:INTROREFINEDLUPMURELATION}
(which also holds for this point $p$)
along the corresponding integral curve (starting from time $\mathring{\upepsilon}^{-1}$),
we deduce that any shrinking of $\upmu$ 
could only be caused by the term $\mbox{Error} = o\left(\mathring{\upepsilon} \right) (1 + t)^{-1}.$
In particular, the amplitude of $\mbox{Error}$ goes to $0$ strictly faster than $\mathring{\upepsilon}.$
Hence, it takes at least $\exp\left\lbrace \frac{1}{o\left(\mathring{\upepsilon} \right)} \right\rbrace$ amount of time
for $\upmu$ to vanish in this region,
and by Theorem~\ref{T:LONGTIMEPLUSESTIMATES}, the sharp classical lifespan theorem,
no singularity of any kind will form while $\upmu$ is positive.
As we mentioned above, this existence time is longer than the time 
$\exp\left(\frac{1}{C \mathring{\upepsilon}} \right)$
that holds for general small data
(see inequality \eqref{E:TLIFESPANLOWERBOUND}).

We now give a brief overview explaining the so-called ``first step'' stated two paragraphs above, 
namely that relative to standard spherical coordinates $(t,r,\theta)$ on Minkowski spacetime,
$[G_{\Lunit \Lunit} \rgeo \Rad \Psi](t,r,\theta)$ is well-approximated by
$- \FutFailFac \frac{\partial}{\partial q} \Fried[(\mathring{\Psi}, \mathring{\Psi}_0)](q = r - t,\theta)$ at time $t = \mathring{\upepsilon}^{-1}.$
Since Lemma~\ref{L:NULLCONDFACT} shows that
$G_{\Lunit \Lunit}$ is well-approximated by $\FutFailFac,$
the main issue that we must explain is the following estimate, 
valid for small data:
\begin{align} \label{E:FRIEDLANDERAPPROXIMATION}
	\left\|
		\rgeo \Rad \Psi
		+ 
		\frac{\partial}{\partial q}
		\Fried[(\mathring{\Psi}, \mathring{\Psi}_0)]
	\right\|_{C^0(\Sigma_{\mathring{\upepsilon}^{-1}}^{U_0})}
	& \leq C \mathring{\upepsilon}^2 \ln\left( \frac{1}{\mathring{\upepsilon}} \right).
\end{align}
To derive \eqref{E:FRIEDLANDERAPPROXIMATION}, 
one can rely on two standard kinds of estimates.
First, one needs an estimate showing that at time $\mathring{\upepsilon}^{-1},$ the solution 
$\Psi_{(Linear)}$ to the \emph{linear} wave equation $\square_m \Psi = 0$
with initial data $(\mathring{\Psi}, \mathring{\Psi}_0)$
is well-approximated by $r^{-1} \Fried[(\mathring{\Psi}, \mathring{\Psi}_0)]$
and that similar results hold for the higher-order $(t,r,\theta)$ coordinate derivatives of
$\Psi_{(Linear)}.$
Such estimates can be derived with the help of the fundamental solution to the linear 
wave equation; see, for example, \cite{lH1997}.
The second kind are estimates showing that at time $\mathring{\upepsilon}^{-1},$ 
\emph{long before any singularity can form},
the nonlinear solution $\Psi$ is well-approximated by $\Psi_{(Linear)},$
that $\rgeo$ is well-approximated by $r,$
that $\upmu$ is well-approximated by $1,$
that $\Lunit$ is well-approximated by $\partial_t + \partial_r,$
that $\Radunit$ is well-approximated by $-\partial_r,$
etc. The estimates for $\Psi$ can be derived with the help of Klainerman's 
Minkowskian vectorfield method approach, as described in \cite{sK1986}.
The estimates for $\upmu$ and $\Lunit^i$ can then be derived
with the help of the transport equations
\eqref{E:INTROLUNITUPMUHEURISTIC}-\eqref{E:INTROLUNITLUNITIHEURISTIC}
and the dispersive $C^0$ estimates for $\Psi$
(see Sect.~\ref{S:RESCALEDFRAME}).
The remaining estimates then follow without much difficulty;
all of these estimates can be derived 
by using arguments that are explained 
in detail in Chapter~\ref{C:C0BOUNDBOOTSTRAP}.
Combining such estimates, one can deduce \eqref{E:FRIEDLANDERAPPROXIMATION}
without much difficulty.

It would be interesting to better understand the connection between 
our nearly spherically symmetric assumptions on the data stated in
Theorem~\ref{T:STABLESHOCKFORMATION}, 
Christodoulou's shock-formation condition \eqref{E:SHOCKFUNCTIONMUSTBESUFFICIENTLYLARGE},
and John's criterion \eqref{E:ALINHACRELAXED}. 
We speculate that either our condition or Christodoulou's 
can be used to directly prove that the 
corresponding analog of \eqref{E:ALINHACRELAXED} must hold
in the relevant region, but we do not investigate this possibility here.
We conclude by mentioning that it would be interesting
to investigate whether or not \emph{all} (nontrivial) compactly supported data 
that are small enough such that Theorem~\ref{T:LONGTIMEPLUSESTIMATES} applies
must necessarily lead to shock formation;
the limitation of the argument outlined in this section
is that we might need to shrink the amplitude of the data
in order to know that a shock will form
(see just below \eqref{E:ALINHACRELAXED}).

\section{Outline of the monograph}
In Appendix \ref{A:EQUIVALENTPROBLEM}, we outline how to extend our results
for covariant equations $\square_{g(\Psi)} \Psi = 0$
to non-covariant equations 
$(h^{-1})^{\alpha \beta}(\partial \Phi) \partial_{\alpha} \partial_{\beta} \Phi = 0.$
In Appendix \ref{A:NOTATION}, we collect some of the important notation 
and conventions that we use throughout the monograph.

As we have emphasized, our main goal in this monograph is to prove
Theorem~\ref{T:LONGTIMEPLUSESTIMATES}, the sharp classical lifespan theorem.
Our proof of shock formation for nearly spherically symmetric small data, 
Theorem~\ref{T:STABLESHOCKFORMATION}, then follows without much additional effort.
We prove Theorem~\ref{T:STABLESHOCKFORMATION} 
in Chapter~\ref{C:PROOFOFSHOCKFORMATION}, the final section of the monograph.
The monograph and indeed, the main bootstrap argument culminating in the proof 
of Theorem~\ref{T:LONGTIMEPLUSESTIMATES}, are organized in an essentially linear fashion.
We now outline the main steps in its proof.

\noindent \textbf{Step 1: Deriving equations.}
Chapters \ref{C:BASICGEOMETRICCONSTRUCTIONS}-\ref{C:RENORMALIZEDEIKONALFUNCTIONQUANTITIES}
are dedicated to geometry and algebra. In short, in these sections, we define
all of the quantities that play a role in our analysis and we derive the equations that they
satisfy. In particular, in Prop.~\ref{P:DIVTHMWITHCANCELLATIONS},
we derive energy-flux identities, which are the starting
point for our derivation of a priori $L^2-$type estimates for $\Psi$ and its derivatives.
The most difficult analysis occurs in Chapter~\ref{C:RENORMALIZEDEIKONALFUNCTIONQUANTITIES},
where we take special care to derive evolution equations for suitably modified 
eikonal function quantities.
As we described in Sect.~\ref{SS:EIKONALTOPORDER},
we must use the modified quantities at the top order in order
to avoid derivative loss and other difficulties.

\noindent \textbf{Step 2: Bootstrap assumptions and easy $C^0$ and pointwise estimates.}
In Chapter~\ref{C:C0BOUNDBOOTSTRAP}, we state $C^0$ bootstrap assumptions for
$\Psi$ and its lower order derivatives on a ``bootstrap region'' of the form 
$\mathcal{M}_{\Tboot,U_0}$ (see Figure~\ref{F:SOLIDREGION} on pg. \pageref{F:SOLIDREGION}). 
We also state similar assumptions for 
$\upmu,$ 
$\Lunit^i,$ 
and $\upchi$
(more precisely, for re-centered versions of these variables in 
which we subtract off their background values),
as well as the positivity assumption $\upmu > 0.$
We then use these bootstrap assumptions, the equations from Step 1,
and a small-data assumption to derive $C^0$ and pointwise estimates for 
most of the quantities defined in Step 1. In particular, we derive $C^0$ estimates
for the re-centered versions of
$\upmu,$ 
$\Lunit^i,$ 
$\upchi,$
and their lower-order derivatives 
that are improvements over the bootstrap assumptions,
thus closing this portion of the bootstrap argument.
These pointwise estimates are tedious but relatively easy to derive.
We save the most difficult pointwise estimates for Steps 3 and 5.

\noindent \textbf{Step 3: Sharp estimates for $\upmu.$}
In Chapter~\ref{C:SHARPESTIMATESFORUPMU}, 
we derive sharp estimates, 
far more detailed than those of Step 2, 
for the inverse foliation density $\upmu.$
These estimates play a key role in our Gronwall argument for
the energy-flux quantities. 
The most difficult aspects of this analysis 
involve a posteriori estimates for $\upmu$, as described in
Sect.~\ref{SS:INTROLUPMUOVERUPMU}.

\noindent \textbf{Step 4: Coercive $L^2$ quantities.}
In Chapter~\ref{C:L2COERCIVENESS}, we exhibit the coerciveness
of the $L^2-$type energy-flux quantities defined in Step 1.
We then use the energy-flux quantities as building blocks to
construct our fundamental $L^2-$controlling quantities.
These are the quantities that we use to control $\Psi$ in $L^2.$

\noindent \textbf{Step 5: Pointwise estimates for the error terms.}
In Step 5, we derive pointwise estimates for the 
error integrands appearing in the energy-flux  
identities of Prop.~\ref{P:DIVTHMWITHCANCELLATIONS}.
More precisely, we derive these pointwise estimates in the case
of $\Psi$ and also the cases of the higher-order derivatives
$\mathscr{Z}^N \Psi,$ where $\mathscr{Z}^N$
is an $N^{th}$ order string of commutation vectorfields,
which we constructed in Step 1.
Many of the error integrands that we must bound arise because
when we commute the wave equation $\upmu \square_{g(\Psi)} \Psi = 0$
with $\mathscr{Z}^N,$ we generate a huge number of error terms.
These pointwise estimates are a precursor to Step 8,
in which we derive estimates
for the corresponding error integrals. We divide the pointwise estimates 
into Ch.~\ref{C:POINTWISEBOUNDSFOREASYERRORINTEGRANDS} (easy estimates)
and Ch.~\ref{C:POINTWISEESTIMATESDIFFICULTERRORINTEGRANDS} (difficult estimates).
In Chapter~\ref{C:POINTWISEBOUNDSFOREASYERRORINTEGRANDS},
we identify those integrand factors that are ``harmless.''
The harmless terms have a negligible effect on the dynamics and 
are easy to treat in Step 8.
The remaining terms, which are the ones that 
would cause derivative loss and other problems if handled improperly,
require special care.
In Chapter~\ref{C:POINTWISEESTIMATESDIFFICULTERRORINTEGRANDS}, we derive 
pointwise estimates for these remaining difficult 
error integrand terms. To derive these difficult estimates, we use
the modified eikonal function quantities constructed in Step 1.

\noindent \textbf{Step 6: Elliptic estimates and Sobolev embedding.}
Ch.~\ref{C:ELLIPTIC} is an interlude. 
We derive elliptic estimates for various quantities on the Riemannian manifolds
$(S_{t,u},\gsphere).$
We use the elliptic estimates in Step 8 to control
some of the top-order error terms in $L^2.$ 
We also derive a Sobolev embedding result for functions on $(S_{t,u},\gsphere).$ 
This is the tool that we eventually use, 
after we have derived suitable $L^2$ estimates,
to improve the $C^0$ bootstrap assumptions for 
$\Psi$ made in Step 2.

\noindent \textbf{Step 7: Below-top-order $L^2$ estimates for the eikonal function quantities.}
In Chapter~\ref{C:EIKONALBELOWTOPORDERSOBOLEVESTIMATES}, we use 
the equations from Step 1 and the $C^0$ and pointwise estimates from Step 2
to derive a priori $L^2$ estimates for the below-top-order
derivatives of quantities constructed out of the eikonal function,
such as the re-centered versions of
$\upmu,$ $\Lunit^i,$ etc;
such terms appear in Step 5 when we commute the wave equation.
The right-hand sides of these estimates involve
the fundamental $L^2-$controlling quantities constructed in Step 4. 
These estimates are relatively 
easy to derive without using the modified quantities
constructed in Step 1 because we allow them to lose one derivative
relative to $\Psi.$

\noindent \textbf{Step 8: A priori $L^2$ estimates for $\Psi$ and the eikonal function quantities 
up to top-order.}
This step, which we carry out in Chapter~\ref{C:ERRORTERMSOBOLEV},
is the hardest one. We estimate, in $L^2,$ 
the remaining error terms appearing in the 
commuted wave equation. In particular, we estimate the difficult
top-order eikonal function quantities that we ignored in Step 7.
As in Step 7, the right-hand sides of these estimates involve
the fundamental $L^2-$controlling quantities constructed in Step 4. 
We then use the estimates for the error terms to derive, by
a long Gronwall argument based on integral identities such as
\eqref{E:INTROMTUDIVERGENCETHM}, a priori estimates
for the fundamental $L^2-$controlling quantities
on the bootstrap region $\mathcal{M}_{\Tboot,U_0}.$ 
This section contains the main estimates we need to prove
Theorem~\ref{T:LONGTIMEPLUSESTIMATES}, our
sharp classical lifespan theorem.

\noindent \textbf{Step 9: Local well-posedness.}
Ch.~\ref{C:LOCALWELLPOSEDNESS} is another interlude.
We sketch a proof of local well-posedness for the equations of interest,
and we establish continuation criteria for avoiding blow-up
of the solutions. In particular, the results of Ch.~\ref{C:LOCALWELLPOSEDNESS} imply that 
the bootstrap assumptions from Step 2 are satisfied 
for at least a short time.

\noindent \textbf{Step 10: The sharp classical lifespan theorem.}
In Chapter~\ref{C:LONGTIMEEXISTENCE}, we state and prove
Theorem~\ref{T:LONGTIMEPLUSESTIMATES}, which is the main
theorem of the monograph. The difficult part of the 
argument was carried out in Step 8. Roughly, the
theorem shows that for small data, the \emph{only} way a
singularity could form is for $\upmu$ to go to $0$ in finite time,
which would necessarily signify the onset of shock formation.
The theorem also provides a collection of estimates that are 
verified by the solution, whether or not it forms a shock.
We also prove Cor.~\ref{C:ANGULARDDERIVATIVEESTIMATES},
which shows that 
if the initial data have ``very small'' angular derivatives,
then this condition is propagated by the solution;
we use the corollary in Step 11.

\noindent \textbf{Step 11: Shock formation.}
In Chapter~\ref{C:PROOFOFSHOCKFORMATION}, 
we prove Theorem~\ref{T:STABLESHOCKFORMATION}, which
is our small-data shock-formation theorem.
Roughly, the theorem states that sufficiently small, nearly 
spherically symmetric data lead to 
finite-time shock formation.
The theorem is not difficult to prove, thanks to
the estimates of Theorem~\ref{T:LONGTIMEPLUSESTIMATES}
and Cor.~\ref{C:ANGULARDDERIVATIVEESTIMATES}.
We can also extend Theorem~\ref{T:STABLESHOCKFORMATION} 
to show shock formation in solutions 
corresponding to a significantly larger class of data;
see Sect.~\ref{SS:INTROSHOCKFORMINGDATACOMPARISON}
for an outline of a proof.

\section{Suggestions on how to read the monograph}
This monograph makes extensive use of geometry.
For an introduction to the basic concepts in Lorentzian and Riemannian geometry
that play a role in our analysis, readers may consult \cite{bO1983} and \cite{pP2006}.

To become acquainted with the new difficulties that
set this work apart from more standard global results for nonlinear wave equations,
we suggest
starting with 
Prop.~\ref{P:DEFORMATIONTENSORFRAMECOMPONENTS},
Lemma~\ref{L:WAVEONCECOMMUTEDBASICSTRUCTURE},
and
Prop.~\ref{P:COMMUTATIONCURRENTDIVERGENCEFRAMEDECOMP}.
From these, one can see the kind of error terms that arise when we commute 
the wave equation with our commutation vectorfields.
The error terms lead to the presence of error integrals in the 
$L^2-$type energy-flux identities of
Prop.~\ref{P:DIVTHMWITHCANCELLATIONS}.
We identify the most difficult error terms in
Lemma~\ref{L:WAVENTIMESCOMMUTEDBASICSTRUCTURE},
Prop.~\ref{P:IDOFKEYDIFFICULTENREGYERRORTERMS},
and Cor.~\ref{C:REDUCTIONOFPROOFTOPURESPATIALCOMMUTATORS}.
We bound the corresponding difficult error integrals in 
Lemma~\ref{L:DANGEROUSTOPORDERMULTERRORINTEGRAL}
and Lemma~\ref{L:DANGEROUSTOPORDERMORERRORINTEGRAL}.
Because of the degenerate nature of the
Gronwall-type Lemma~\ref{L:FUNDAMENTALGRONWALL},
the bounds from these two lemmas 
lead to degenerate a priori estimates for the high-order $L^2$
quantities. \textbf{These degenerate a priori estimates
are the primary non-standard aspects of our work and the work \cite{dC2007}.}


\chapter{Initial Data, Basic Geometric Constructions, and the Future Null Condition Failure Factor} 
\label{C:BASICGEOMETRICCONSTRUCTIONS}
\thispagestyle{fancy}
In Chapter~\ref{C:BASICGEOMETRICCONSTRUCTIONS}, we discuss the initial data of interest for 
the covariant wave equation\footnote{We formally define the covariant wave operator $\square_{g(\Psi)}$ in Def.~\ref{D:WAVEOPERATORSANDLAPLACIANS}.} \eqref{E:WAVEGEO}, that is, for the equation $\square_{g(\Psi)} \Psi = 0.$
We also define the eikonal function $u,$
some intimately related families of surfaces, 
a related set of geometric coordinates,
and several vectorfield frames,
all of which play an important role in our analysis.
In particular, we define the inverse foliation density $\upmu,$
which in a certain sense is the main object of study in this monograph.
As we discussed in Sect.~\ref{S:LOWERBOUNDLEADINGTOBLOWUP},
for the solutions of interest,
the formation of a shock and the blow-up of a radial derivative
of $\Psi$ both precisely correspond to the vanishing of $\upmu$.
We also define the ``future null condition failure factor'' $\FutFailFac,$ which, 
roughly speaking, is the coefficient of the terms that drive small-data shock formation
in the region $\lbrace t \geq 0 \rbrace.$
Finally, we introduce some schematic notation that
we often use to capture the essential features of our equations 
and estimates.

\begin{remark}[\textbf{There is some redundancy}]
Some of the discussion and definitions in Chapter~\ref{C:BASICGEOMETRICCONSTRUCTIONS}
also appear in Chapter~\ref{C:INTRO}.
For pedagogical reasons, we have chosen to repeat some material here.
\end{remark}

\section{Initial data} \label{S:DATA}
The Cauchy hypersurface of interest is 
$\Sigma_0 := \lbrace (t,x^1,x^2,x^3) \in \mathbb{R}^4 \ | \ t = 0 \rbrace.$
The data for the wave equation \eqref{E:WAVEGEO} are 
\begin{align}
	\mathring{\Psi} := \Psi|_{\Sigma_0}, \qquad \mathring{\Psi}_{0} := \partial_t \Psi|_{\Sigma_0}.
\end{align}
We assume that $\mathring{\Psi},\mathring{\Psi}_0$ are compactly supported in the Euclidean unit ball.
Let $U_0$ be any real number verifying
\begin{align} \label{E:U0}
	0 < U_0 < 1.
\end{align}	
We view $U_0$ as a parameter that is fixed until Sect.~\ref{S:SHOCKFORMINGDATA}.
Our use of the notation ``$U_0$'' is connected to the eikonal function $u,$ which we define
in Sect.~\ref{S:EIKONAL}.

Let $r = \sqrt{\sum_{a=1}^3(x^a)^2}$ 
denote the standard Euclidean radial coordinate on $\Sigma_0.$ 
We study the future-behavior of the solution in the nontrivial region that corresponds to
the portion of the nontrivial data lying in an annular region, centered at the origin, 
of inner Euclidean radius
$1 - U_0$ and outer Euclidean radius $1$ (that is, the thickness of the region is $U_0$):
\begin{align} \label{E:SIGMA0U0DEF}
	\Sigma_0^{U_0} := \lbrace p \in \Sigma_0 \ | \ 1 - U_0 \leq r(p) \leq 1 \rbrace.
\end{align}
The spacetime region that we study is trapped between the two outgoing
null cones $\mathcal{C}_{U_0}$ and $\mathcal{C}_0,$ where the latter cone is flat
because $\Psi$ completely vanishes in its exterior (see Figure~\ref{F:REGION} on pg. \pageref{F:REGION}).

\begin{remark}[\textbf{The role of $U_0$}] \label{R:CHOOSINGU0}
	In Sect.~\ref{S:SHOCKFORMINGDATA},
	for given initial data,
	we choose the data-dependent parameter $U_0$ 
	so that the shock singularity forms in 
	the region trapped between the two outgoing null cones
	$\mathcal{C}_{U_0}$ and $\mathcal{C}_0.$
\end{remark}

\section{The eikonal function and the geometric radial variable} \label{S:EIKONAL}

As we explained in Sect.~\ref{S:INTROBASICGEO}, the eikonal function plays 
a fundamental role in our analysis of solutions.

\begin{definition} [\textbf{The eikonal function}]
The eikonal function $u$ increases towards the future (that is, $\partial_t u > 0$) and verifies
\begin{align} \label{E:OUTGOINGEIKONAL}
	(g^{-1})^{\alpha \beta} \partial_{\alpha} u \partial_{\beta} u & = 0.
\end{align}

We define the initial condition of $u$ as follows, where $(x^1,x^2,x^3)$ are the rectangular coordinates on $\Sigma_0:$
\begin{align} \label{E:EIKONALFUNCTIONINITIALCONDITIONS}
	u|_{\Sigma_0}(x^1,x^2,x^3) = 1 - r(x^1,x^2,x^3), \qquad r(x^1,x^2,x^3):= \sqrt{\sum_{a=1}^3 (x^a)^2}.
\end{align}

\end{definition}

We often use the following geometric radial variable in our analysis.

\begin{definition} [\textbf{The geometric radial variable}]
We define the geometric radial variable $\rgeo$ by
\begin{align} \label{E:GEOMETRICRADIAL}
	\rgeo
	&:= 1 - u + t.
\end{align}

\end{definition}

The following subsets of spacetime play an important role in our analysis; see Figure~\ref{F:SOLIDREGION}
on pg. \pageref{F:SOLIDREGION}.

\begin{definition} [\textbf{Subsets of spacetime}]
\label{D:HYPERSURFACESANDCONICALREGIONS}
We define the following spacetime subsets:
\begin{subequations}	
\begin{align}
	\Sigma_{t'} & := \lbrace (t,x^1,x^2,x^3) \in \mathbb{R}^4  \ | \ t = t' \rbrace, 
		\label{E:SIGMAT} \\
	\Sigma_{t'}^{u'} & := \lbrace (t,x^1,x^2,x^3) \in \mathbb{R}^4  \ | \ t = t', \ 0 \leq u(t,x^1,x^2,x^3) \leq u' \rbrace, 
		\label{E:SIGMATU} 
		\\
	\mathcal{C}_{u'}^{t'} & := \lbrace (t,x^1,x^2,x^3) \in \mathbb{R}^4 \ | \ u(t,x^1,x^2,x^3) = u' \rbrace \cap \lbrace (t,x^1,x^2,x^3) 
		\in \mathbb{R}^4  \ | \ 0 \leq t \leq t' 		
		\rbrace, 
		\label{E:CUT} \\
	S_{t',u'} 
		&:= \mathcal{C}_{u'}^{t'} \cap \Sigma_{t'}^{u'}
		= \lbrace (t,x^1,x^2,x^3) \in \mathbb{R}^4 \ | \ t = t', \ u(t,x^1,x^2,x^3) = u' \rbrace, 
			\label{E:STU} \\
	\mathcal{M}_{t',u'} & := \cup_{u \in [0,u']} \mathcal{C}_u^{t'} \cap \lbrace (t,x^1,x^2,x^3) \in \mathbb{R}^4  \ | \ t < t' \rbrace.
		\label{E:MTUDEF}
\end{align}
\end{subequations}
We refer to the $\Sigma_t$ and $\Sigma_t^u$ as ``constant time slices,'' the $\mathcal{C}_u^t$ as ``outgoing null cones,''
and the $S_{t,u}$ as ``spheres.'' 
We sometimes use the notation $\mathcal{C}_u$ in place of $\mathcal{C}_u^t$ 
when we are not concerned with the truncation time $t.$
\end{definition}

\begin{remark}
	We make the following remarks.
	\begin{itemize}
		\item It follows from standard domain of dependence considerations that for any $t \geq 0,$ 
			the solution corresponding to the data
			from Sect.~\ref{S:DATA} agrees with the trivial solution $\Psi \equiv 0$
			in the exterior of $\mathcal{C}_0^t.$ We therefore do not comment further on this region.
		\item The $S_{t,u}$ are not generally round Euclidean spheres even though
			the $S_{0,u}$ are.
		\item Note that $\mathcal{M}_{t,u}$ is ``open at the top.''
	\end{itemize}
\end{remark}

We now define some important geometric objects related to the above subsets.

\begin{definition}[\textbf{First fundamental forms}] \label{D:FIRSTFUND}
	We define the following first fundamental forms.
	\begin{itemize}
		\item $\gt$ denotes the Riemannian metric on $\Sigma_t$ induced by $g.$ 
			That is, $\gt(X,Y) = g(X,Y)$ for all $\Sigma_t-$tangent vectors $X$ and $Y.$
			$\gt^{-1}$ denotes the corresponding inverse metric.
		\item $\gsphere$ denotes the Riemannian metric on $S_{t,u}$ induced by $g.$
		  That is, $\gsphere(X,Y) = g(X,Y)$ for all $S_{t,u}-$tangent vectors $X$ and $Y.$
			$\ginversesphere$ denotes the corresponding inverse metric.
	\end{itemize}
\end{definition}

\begin{definition}[\textbf{Levi-Civita connections}] 
\label{D:CONNECTIONS}
	We use the following notation for the Levi-Civita connections associated $g,$ $m,$ and $\gsphere.$
	\begin{itemize}
		\item $\D$ denotes the Levi-Civita connection of the spacetime metric $g.$
		\item $\nabla$ denotes the Levi-Civita connection of the Minkowski metric $m.$
		\item $\angD$ denotes the Levi-Civita connection of $\gsphere.$
	\end{itemize}
\end{definition}

\begin{definition}[\textbf{Covariant wave operators and Laplacians}] 
	\label{D:WAVEOPERATORSANDLAPLACIANS}
	We use the following standard notation.
	\begin{itemize}
		\item $\square_g := (g^{-1})^{\alpha \beta} \D_{\alpha \beta}^2$ denotes the covariant wave operator 
			corresponding to the spacetime metric $g.$
		\item $\angLap := (\ginversesphere)^{AB} \angDsquaredarg{A}{B}$ denotes the covariant Laplacian 
			corresponding to $\gsphere.$ 
	\end{itemize}
\end{definition}

\section{Frame vectorfields and the inverse foliation density}
\label{S:FRAMEVECTORFIELDS}
We now define the gradient vectorfield $\Lgeo$ associated to the eikonal function.

\begin{definition} [\textbf{The outgoing null geodesic vectorfield}]
Let $u$ be the eikonal function \eqref{E:OUTGOINGEIKONAL}. We define the vectorfield $\Lgeo$ by
\begin{align} \label{E:LGEOEQUATION}
	\Lgeo^{\nu} & := - (g^{-1})^{\nu \alpha} \partial_{\alpha} u.
\end{align}
\end{definition}

Note that by \eqref{E:OUTGOINGEIKONAL}, we have
\begin{align} \label{E:LGEOISGEODESIC}
	\D_{\Lgeo} \Lgeo & = 0.
\end{align}
Also, by \eqref{E:OUTGOINGEIKONAL}, $\Lgeo$ is null 
(that is, $g(\Lgeo, \Lgeo)=0$)
and $g-$orthogonal to the level sets $\mathcal{C}_u$ of the eikonal function.

We now introduce the most important quantity in our analysis: 
the inverse foliation density of the level sets of $u$
relative to the hypersurfaces $\Sigma_t.$

\begin{definition}[\textbf{The inverse foliation density of the outgoing cones}] 
\label{D:UPMU}
Let $t$ be the Minkowskian time coordinate and $u$ the eikonal function. 
We define the inverse foliation density 
(of the level sets of $u$ relative to the hypersurfaces $\Sigma_t$)
$\upmu$ by
\begin{align} \label{E:UPMUDEF}
	\upmu & := - \frac{1}{(g^{-1})^{\alpha \beta} \partial_{\alpha} t \partial_{\beta} u} 
	= - \frac{1}{(g^{-1})^{0 \alpha} \partial_{\alpha} u} 
	= \frac{1}{\Lgeo^0},
\end{align}
where the last two equalities hold in the \textbf{rectangular} coordinate system.
\end{definition}

When $\upmu$ vanishes, many quantities, including the rectangular components 
$\Lgeo^{\nu},$ blow-up like $\upmu^{-1}.$
The outgoing null vectorfield $\Lunit,$ defined just below, is a rescaled version of $\Lgeo$ that 
``removes the singular factor $\upmu^{-1}.$'' We prove that $\Lunit$ remains regular and near
$\Lunit_{(Flat)} = \partial_t + \partial_r$ throughout the evolution. 
Similarly, the vectorfield $\uLunit$ defined below is an ingoing null vectorfield
that remains regular and near $\uLunit_{(Flat)} = \partial_t - \partial_r$ throughout the evolution.
In contrast, the rectangular components of the ingoing null vectorfield
$\uLgood = \upmu \uLunit$ vanishes precisely when the shock forms.

\begin{definition} [$\Lunit,$ $\uLgood,$ \textbf{and} $\uLunit$]
	\label{D:NULLVECTORFIELDS}
	We define the rescaled null vectorfield $\Lunit$ as follows:
	\begin{align} \label{E:LUNIT}
		\Lunit & := \upmu \Lgeo. 
	\end{align}
	We then define $\uLgood$ to be the (unique) null vectorfield that is 
	$g-$orthogonal to $S_{t,u}$ and that is normalized by
	\begin{align} \label{E:ULGOODNORMALIZATION}
		g(\uLgood, \Lunit) & = - 2 \upmu.
	\end{align}
	
	Finally, we define
	\begin{align} \label{E:BA}
		\uLunit 
		& := \upmu^{-1} \uLgood.
	\end{align}
	
\end{definition}

In the next lemma, we reveal some basic properties of $\Lunit$ and $\uLgood.$

\begin{lemma}[\textbf{Basic properties of the null pair $\Lunit, \uLgood$}]
\label{L:NULLPAIRBASICPROPERTIES}
The following identities hold, where $t$ is the Minkowski time coordinate:
\begin{subequations}
\begin{align}
	\Lunit u & = 0, 
		& \Lunit t = \Lunit^0 & = 1,
		\label{E:LAPPLIEDTOUANDT}
		 \\
	\uLgood u & = 2,
		&  \uLgood t = \uLgood^0 & = \upmu.
		\label{E:ULAPPLIEDTOUANDT}
\end{align}
\end{subequations}

\end{lemma}

\begin{proof}
	From \eqref{E:OUTGOINGEIKONAL} and \eqref{E:LGEOEQUATION}, we see that
	$\Lgeo u = 0.$ It thus follows from \eqref{E:LUNIT} that $\Lunit u = 0$
	as desired. The identity $\Lunit t = 1$ follows from 
	\eqref{E:LGEOEQUATION}, 
	\eqref{E:UPMUDEF},
	and \eqref{E:LUNIT}. We have thus shown \eqref{E:LAPPLIEDTOUANDT}.
	We next note that $\uLgood u = - g(\uLgood, \Lgeo) = - \upmu^{-1} g(\uLgood, \Lunit).$
	We conclude from \eqref{E:ULGOODNORMALIZATION} that $\uLgood u = 2$ as desired.
	We have thus shown \eqref{E:ULAPPLIEDTOUANDT}.
	
	To show that $\uLgood t = \upmu,$ we define (relative to rectangular coordinates)
	$V^{\nu} := - (g^{-1})^{\nu \alpha} \partial_{\alpha} t.$
	We note that $V$ is future-directed and 
	$g-$orthogonal to $\Sigma_t$ and hence to $S_{t,u}.$ 
	Since $\Lunit$ and $\uLgood$ span the $g-$orthogonal
	complement of $S_{t,u},$ there exist scalars $a$ and $b$ such that
	\begin{align} \label{E:VINTERMSOFNULLPAIR}
		V^{\nu} = a \Lunit^{\nu} + b \uLgood^{\nu}.
	\end{align}
	Contracting \eqref{E:VINTERMSOFNULLPAIR} against $\Lunit_{\nu},$ 
	using the fact that $\Lunit$ 
	is null, 
	and using
	\eqref{E:ULGOODNORMALIZATION}
	and
	\eqref{E:LAPPLIEDTOUANDT},
	we find that
	$1 = \Lunit t = - \Lunit_{\nu} V^{\nu} = 2 \upmu b.$
	Similarly, we find that $\uLgood t = 2 \upmu a.$
	Also using the fact that $\uLgood$ is null, we deduce that
	\begin{align} \label{E:GNORMFIRSTCOMPUTATION}
		g_{\alpha \beta} V^{\alpha} V^{\beta}
		& = - 4 \upmu ab = - 2 a.
	\end{align}
	On the other hand, by \eqref{E:GINVERSE00ISMINUSONE}, we have that
	\begin{align}  \label{E:GNORMSECONDCOMPUTATION}
		g_{\alpha \beta} V^{\alpha} V^{\beta}
		& = (g^{-1})^{\alpha \beta} \partial_{\alpha} t \partial_{\beta} t
		= (g^{-1})^{00} = - 1. 
	\end{align}	
	From \eqref{E:GNORMFIRSTCOMPUTATION} and \eqref{E:GNORMSECONDCOMPUTATION}
	it follows that $a = 1/2$ and hence $\uLgood t = \upmu$ as desired.
\end{proof}

We now define an \emph{inward-pointing} $\Sigma_t-$tangent vectorfield $\Rad$
and an \emph{upward pointing} $\Sigma_t-$normal vectorfield $\Timenormal.$

\begin{definition}[\textbf{The vectorfields $\Rad,$ $\Radunit,$ and $\Timenormal$}]
\label{D:RADIALVECTORFIELDSANDTIMELIKENORMAL}
We define the vectorfields 
$\Rad,$ 
$\Radunit,$
and $\Timenormal$ as follows:
\begin{subequations}
\begin{align}
	\Rad 
	& := \frac{1}{2} \lbrace - \upmu \Lunit + \uLgood \rbrace, 
		\label{E:RADINTERMSOFLANDULGOOD} \\
	\Radunit 
	& := \upmu^{-1} \Rad 
		= \frac{1}{2} \lbrace - \Lunit + \uLunit \rbrace, 
		\label{E:RADUNIT} \\
	\Timenormal 
	& := \frac{1}{2} \lbrace \Lunit + \uLunit \rbrace
		= \Lunit + \Radunit.
		\label{E:UNITSIGMATNORMAL}
\end{align}
\end{subequations}
\end{definition}

The Minkowskian analogs of \eqref{E:RADUNIT} and \eqref{E:UNITSIGMATNORMAL}
are $- \partial_r$ and $\partial_t$ respectively.

In the next lemma, we reveal some basic properties of $\Lunit,$ $\Rad,$ $\Radunit,$ and $\Timenormal.$

\begin{lemma} [\textbf{Basic properties of $\Lunit,$ $\Rad,$ $\Radunit,$ and $\Timenormal$}]
\label{L:BASICPROPERTIESOFLUNITRADANDTIMENORMAL}
$\Rad$ is $g-$orthogonal to the $S_{t,u}$ and tangent to $\Sigma_t.$
$\Timenormal$ is timelike, future-directed, $g-$orthogonal to the $S_{t,u},$ 
and $g-$orthogonal to $\Sigma_t.$

In addition, the following identities hold:
\begin{subequations}
\begin{align}
	g(\Lunit, \Lunit)  
	& = 0,
		\label{E:LISNULL} \\
	g(\Rad,\Rad) 
	& = \upmu^2, 
		\label{E:RADLENGTH} \\
	g(\Radunit,\Radunit) 
	& = 1, 
		\label{E:HATRADLENGTH} 
		\\
	g(\Lunit, \Rad)  
	& = - \upmu,
		\label{E:LRADINNERPRODUCT} 
		\\
	g(\Lunit, \Radunit)  
	& = - 1.
		\label{E:LRADUNITINNERPRODUCT} 
\end{align}
\end{subequations}

Furthermore, we have
\begin{align}
		g(\Timenormal,\Timenormal) & = - 1.
		\label{E:TIMENORMALLENGTH}
\end{align}


In addition, we have
\begin{subequations}
\begin{align}
	\Rad t & = 0, & \Rad u & = 1, 
		\label{E:RADUISONERADTISZERO} \\
	\Timenormal t & = 1, & \Timenormal u & = 1.
		\label{E:TIMENORMALAPPLIEDTOTANDU}
\end{align}
\end{subequations}

Finally, the commutator vectorfield $[\Lunit, \Rad]$ is $S_{t,u}-$tangent.

\end{lemma}

\begin{proof}
\eqref{E:LISNULL} is a restatement of the fact that $\Lunit$ is null.
\eqref{E:RADLENGTH} follows from \eqref{E:RADINTERMSOFLANDULGOOD}, the fact
that $\Lunit$ and $\uLgood$ are null, and \eqref{E:ULGOODNORMALIZATION}.
\eqref{E:HATRADLENGTH} follows from \eqref{E:RADUNIT}	and \eqref{E:RADLENGTH}.
\eqref{E:LRADINNERPRODUCT} follows from \eqref{E:RADINTERMSOFLANDULGOOD}, \eqref{E:LISNULL}, and \eqref{E:ULGOODNORMALIZATION}.
\eqref{E:LRADUNITINNERPRODUCT} follows from \eqref{E:RADUNIT} and \eqref{E:LRADINNERPRODUCT}.

\eqref{E:TIMENORMALLENGTH} follows from \eqref{E:UNITSIGMATNORMAL}, \eqref{E:LISNULL}, \eqref{E:HATRADLENGTH}, and \eqref{E:LRADUNITINNERPRODUCT} .

\eqref{E:RADUISONERADTISZERO} follows from \eqref{E:RADINTERMSOFLANDULGOOD}, \eqref{E:LAPPLIEDTOUANDT}, and \eqref{E:ULAPPLIEDTOUANDT}.
\eqref{E:TIMENORMALAPPLIEDTOTANDU} follows from \eqref{E:UNITSIGMATNORMAL}, \eqref{E:LAPPLIEDTOUANDT}, and \eqref{E:ULAPPLIEDTOUANDT}.

The fact that $\Rad$ is $\Sigma_t-$tangent follows from the identity $\Rad t = 0.$
The fact that $\Rad$ is $g-$orthogonal to $S_{t,u}$ follows from \eqref{E:RADINTERMSOFLANDULGOOD}
and the fact that $\Lunit$ and $\uLgood$ are $g-$orthogonal to $S_{t,u}.$ 

To obtain the desired properties of $\Timenormal,$ 
we note that the proof of Lemma~\ref{L:NULLPAIRBASICPROPERTIES}
reveals that $\Timenormal$ is equal to the vectorfield $V^{\nu} := - (g^{-1})^{\nu \alpha} \partial_{\alpha} t,$
which is timelike, future-directed, $g-$orthogonal to the $S_{t,u},$ 
and $g-$orthogonal to $\Sigma_t.$ 

To obtain the fact that $[\Lunit, \Rad]$ is $S_{t,u}-$tangent, 
we use \eqref{E:LAPPLIEDTOUANDT} and \eqref{E:RADUISONERADTISZERO}
to conclude that $0 = [\Lunit, \Rad] t = [\Lunit, \Rad] u.$

\end{proof}

In the next lemma,
we derive some important identities involving the rectangular spatial derivatives of
the eikonal function.

\begin{lemma}[\textbf{Identities involving the rectangular spatial derivatives of} $u$]
	The eikonal function $u$ verifies the following equation:
	\begin{align} \label{E:ALTERNATEEIKONAL}
		(\gt^{-1})^{ab} \partial_a u \partial_b u 
		 = \upmu^{-2},
	\end{align}
	where $\gt^{-1}$ is as in Def.~\ref{D:FIRSTFUND}.
	
	Furthermore, the rectangular spatial derivatives 
	of $u$ verify the following equation $(i=1,2,3):$
	\begin{align} \label{E:PARTIALIUINTERMSOFRADUNIT}
		\upmu \partial_i u 
		& = \Radunit_i.
	\end{align}	
\end{lemma}

\begin{proof}
To prove \eqref{E:ALTERNATEEIKONAL}, we
use definition \eqref{E:UPMUDEF} and the assumption $(g^{-1})^{00} = - 1$
(that is, \eqref{E:GINVERSE00ISMINUSONE})
to deduce that 
$\upmu^{-2} = (\partial_t u)^2 - 2 (g^{-1})^{0a}\partial_t u \partial_a u + ((g^{-1})^{0a} \partial_a u)^2.$
From this identity and the eikonal equation \eqref{E:OUTGOINGEIKONAL},
we deduce the identity 
\begin{align} \label{E:FIRSTALTERNATEEIKONAL}
	(g^{-1})^{ab} \partial_a u \partial_b u + ((g^{-1})^{0a} \partial_a u)^2 = \upmu^{-2}.
\end{align}
Furthermore, it is straightforward to verify that under 
the assumption $(g^{-1})^{00} = - 1,$
we have the following matrix identity $(i,j = 1,2,3):$
\begin{align} \label{E:USEFULMATRIXIDENTITY}
	(\gt^{-1})^{ij} 
	& = (g^{-1})^{ij} 
		+ (g^{-1})^{0i}
			(g^{-1})^{0j}.
\end{align}
From \eqref{E:FIRSTALTERNATEEIKONAL} and \eqref{E:USEFULMATRIXIDENTITY},
we deduce the desired identity \eqref{E:ALTERNATEEIKONAL}.

Next, we consider the one-forms on $\Sigma_t$ 
with rectangular components 
$(\partial_1 u, \partial_2 u, \partial_3 u)$ and 
$(\Radunit_1, \Radunit_2, \Radunit_3).$
They are both inward pointing 
and $\gt-$orthogonal to the $S_{t,u}.$ Hence, we
must have $\partial_i u = z \Radunit_i$ with $z > 0.$
From the fact that $(\gtinverse)^{ab} \Radunit_a \Radunit_b = \gt(\Radunit,\Radunit) = 1,$
we deduce that $z = \sqrt{(\gtinverse)^{ab} \partial_a u \partial_b u}.$
The desired identity \eqref{E:PARTIALIUINTERMSOFRADUNIT}
now follows from this identity and \eqref{E:ALTERNATEEIKONAL}.

\end{proof}



\section{Geometric coordinates}
\label{S:GEOCOORDS}
In this section, we construct the geometric coordinate system $(t,u, \vartheta^1, \vartheta^2)$
that plays a fundamental role in our analysis.
To begin, we fix an atlas 
$\lbrace (\mathbb{D}_i,\vartheta_i^1, \vartheta_i^2) \rbrace_{i=1,2}$
on the \emph{Euclidean} unit sphere $S_{0,0} \simeq \mathbb{S}^2,$
where the $\mathbb{D}_i$ are open subsets of $\mathbb{S}^2$
with $\mathbb{S}^2 = \mathbb{D}_1 \cup \mathbb{D}_2.$

\begin{definition} [\textbf{Standard atlas on} $\mathbb{S}^2$]
	\label{D:ATLASONS2}
	We refer to the above atlas as the standard atlas on $\mathbb{S}^2.$
\end{definition}	

\begin{remark}[\textbf{Suppression of the precise coordinate chart}]	
	\label{R:NOCOORDINATECHART}
	Throughout the monograph, we typically suppress 
	the set $\mathbb{D}_i$ 
	and the index $i$ 
	and simply refer to the above coordinates 
	as $(\vartheta^1, \vartheta^2).$
\end{remark}

We now (inwardly) extend the local coordinate functions $(\vartheta^1, \vartheta^2)$ 
to the region $\Sigma_0^{U_0} = \cup_{u \in [0,U_0]} S_{0,u} \subset \Sigma_0^{U_0}$ by using the
\emph{Euclidean} radial vectorfield $- \partial_r = - \frac{x^a}{r} \partial_a$ to transport them,
$(A=1,2):$
\begin{align}
	- \partial_r \vartheta^A|_{\Sigma_0^{U_0}} & = 0.
\end{align}
We then extend $(\vartheta^1,\vartheta^2)$ to $S_{t,u}$ for $t > 0$ and $u \in [0,U_0]$
by using the $g-$null vectorfield $\Lunit$ to transport them.
That is, 
\begin{align}
	\Lunit \vartheta^A & = 0.
\end{align}

\begin{definition}[\textbf{Geometric coordinates}]
	\label{D:GEOMETRICCOORDINATES}
	We refer to the functions $(t,u, \vartheta^1, \vartheta^2)$ 
	induced on spacetime regions of the form $\mathcal{M}_{T,U_0}$
	as the \emph{geometric coordinates}.
\end{definition}

\begin{remark}[\textbf{Different ways to think of} $\vartheta$]
	Note that we are slightly abusing notation by using
	the same symbol $(\vartheta^1,\vartheta^2)$ to denote local coordinate functions on
	$\mathbb{S}^2$ and also the corresponding coordinate functions 
	induced on $\mathcal{M}_{T,U_0}.$
	Equivalently, we are identifying the coordinate function pair
	$(\vartheta^1,\vartheta^2)$ on $\mathcal{M}_{T,U_0}$
	with the corresponding point $\vartheta$ belonging to
	$S_{0,0} \simeq \mathbb{S}^2.$
	We often make these identifications throughout the monograph;
	the precise meaning of the symbol ``$\vartheta$'' will
	always be clear from context.
\end{remark}

Since $\Lunit t =1,$
it follows that relative to the geometric coordinates,
we have that
\begin{align} \label{E:LISDDT}
	\Lunit = \frac{\partial}{\partial t}|_{u,\vartheta^1,\vartheta^2}.
\end{align}
We often use the identity \eqref{E:LISDDT} in our analysis.

\begin{definition}[\textbf{$X_1$ and $X_2$}]
	$X_1$ and $X_2$ are the locally defined $S_{t,u}-$tangent vectorfields defined by
\begin{align} \label{E:STUFRAME}
	X_1 & := \frac{\partial}{\partial \vartheta^1}|_{t,u,\vartheta^2},
		&&	
	X_2 := \frac{\partial}{\partial \vartheta^2}|_{t,u,\vartheta^1}.
\end{align}
\end{definition}

\begin{remark}[\textbf{The meaning of} $\xi_A$]
	We often denote the contraction of a one-form
	$\xi$ against the $S_{t,u}$ frame vectors
	$\lbrace X_1, X_2 \rbrace$
	by using the abbreviated notation 
	$\xi_A := \xi_{X_A} = \xi_a X_{A}^a.$
	We use similar abbreviations when contracting other
	types of tensors against vectors;
	see Sect.~\ref{S:CONTRACTIONNOTATION}.
\end{remark}

\begin{remark}[\textbf{$C^k-$equivalent differential structures until shock formation}]
We often identify
spacetime regions of the form
$\mathcal{M}_{t,U_0}$ 
(see definition \eqref{E:MTUDEF})
with the region
$[0,t) \times [0,U_0] \times \mathbb{S}^2$
corresponding to the geometric coordinates.
This identification is 
justified by the fact that 
during the classical lifespan of the solution,
the differential structure on 
$\mathcal{M}_{t,U_0}$ corresponding to the 
geometric coordinates is $C^k-$equivalent,
for some large integer $k,$ to the
differential structure on 
$\mathcal{M}_{t,U_0}$
corresponding to the rectangular coordinates.
The equivalence is captured by the fact that the change
of variables map $\Upsilon$ (see Sect.~\ref{S:CHOV})
from geometric to rectangular coordinates is differentiable 
with a differentiable inverse, until a shock forms;
see Theorem~\ref{T:LONGTIMEPLUSESTIMATES}.
However, at points where $\upmu$ vanishes
and the rectangular derivatives of $\Psi$ blow-up
(see inequality \eqref{E:RADUNITPSIBLOWSUP}),
the injectivity of $\Upsilon$ breaks down
and the equivalence of the differential structures breaks down
as well.
\end{remark}

\section{Frames}
In this section, we define the vectorfield frames that we use in our analysis.

\begin{definition}[\textbf{Frames}]
	We define the following frames,
	where the vectorfields are defined
	in Sects.~\ref{S:FRAMEVECTORFIELDS} and \ref{S:GEOCOORDS}.
	\begin{itemize}
		\item $\lbrace \Lunit, \Rad, X_1, X_2 \rbrace$ denotes the rescaled frame.
		\item $\lbrace \Lunit, \Radunit, X_1, X_2 \rbrace$ denotes the non-rescaled frame.
		\item $\lbrace \Lunit, \uLgood, X_1, X_2 \rbrace$ denotes the rescaled null frame.
		\item $\lbrace \Lunit, \uLunit, X_1, X_2 \rbrace$ denotes the non-rescaled null frame.
	\end{itemize}
\end{definition}

\begin{remark}[\textbf{The span of the frame vectorfields}]
	The analysis of Sect.~\ref{S:CHOV}
	can be used to show that for the small-data solutions 
	that we study,
	$\lbrace \Lunit, \Rad, X_1, X_2 \rbrace$
	and
	$\lbrace \Lunit, \uLgood, X_1, X_2 \rbrace$
	have span equal to $\mbox{span} \lbrace \partial_{\alpha} \rbrace_{\alpha = 0,1,2,3}$
	at each point where $\upmu > 0.$
	However, these two frames degenerate
	relative to $\lbrace \partial_{\alpha} \rbrace_{\alpha = 0,1,2,3}$ when $\upmu$ vanishes;
	see Remark~\ref{R:UPMUCONNECTIONTODETERMINANT}.
	In contrast, the estimates of Theorem~\ref{T:LONGTIMEPLUSESTIMATES}
	can be used to show that for small-data solutions, 
	the frames $\lbrace \Lunit, \Radunit, X_1, X_2 \rbrace$
	and
	$\lbrace \Lunit, \uLunit, X_1, X_2 \rbrace$
	have span equal to $\mbox{span} \lbrace \partial_{\alpha} \rbrace_{\alpha = 0,1,2,3}$
	and do not degenerate when $\upmu$ vanishes.
\end{remark}

\section{The future null condition failure factor}
We now define the function $\FutFailFac,$
which depends on the structure of the nonlinearities and 
is of fundamental importance in determining whether or not
shocks can form to the future.

\begin{definition}[\textbf{The future null condition failure factor}]
\label{D:FAILUREFACTOR}
	We define the scalar-valued \emph{future null condition failure factor} $\FutFailFac$ by
	\begin{align} \label{E:FAILUREFACTOR}
		\FutFailFac 
		&:= \underbrace{G_{\alpha \beta}(\Psi = 0)}_{\mbox{constants}} \Lunit_{(Flat)}^{\alpha}	
		\Lunit_{(Flat)}^{\beta},
	\end{align}
	where $G_{\alpha \beta}$ is defined in \eqref{E:BIGGDEF},
	$\Lunit_{(Flat)} = \partial_t + (x^a/r) \partial_a$ is the standard 
	outgoing Minkowski-null Minkowski-geodesic vectorfield, 
	and $r$ is the standard Euclidean radial coordinate on $\mathbb{R}^3.$
	
	We also define $\InitialFutFailFac,$ relative to the geometric coordinates, by
	\begin{align} \label{E:INITIALFAILUREFACTOR}
		\InitialFutFailFac(t,u,\vartheta)
		& =
		\InitialFutFailFac(\vartheta) 
		:= \FutFailFac(t=0,u=0,\vartheta).
	\end{align} 
	
\end{definition}
As we explained in Sect.~\ref{S:STRUCTUREOFNONLINEARITIES},
$\FutFailFac$ is the coefficient of the dangerous slow decaying quadratic 
terms in the wave equation $\square_{g(\Psi)} \Psi = 0$
in the region $\lbrace t \geq 0 \rbrace.$
In particular, when $\FutFailFac \equiv 0,$ Klainerman's
work \cite{sK1986} and Christodoulou's work \cite{dC1986a}
both yield small-data global existence.

\begin{remark}[\textbf{The dependence of $\FutFailFac$ on various coordinates}]
Note that $\FutFailFac$ can be viewed as a function 
depending only on
$\theta = (\theta^1,\theta^2),$ 
where $\theta^1$ and $\theta^2$ are local angular 
coordinates corresponding to standard spherical coordinates 
on Minkowski spacetime. 
Furthermore, at $t=0,$
our geometric coordinates $(\vartheta^1, \vartheta^2)$
coincide with $(\theta^1, \theta^2),$
and $\Lunit_{(Flat)}^{\alpha}|_{t=0}$ does not depend on $u.$
It follows that relative to the geometric coordinates $(t,u,\vartheta),$
the right-hand side of \eqref{E:INITIALFAILUREFACTOR} is a function of $\vartheta$ alone.
\end{remark}

\begin{remark}[\textbf{The role of} $\InitialFutFailFac$]
	\label{R:MATHRINGALEPH}
	The point of introducing $\InitialFutFailFac$ is that
	it is a good approximation to $\FutFailFac$
	(see Lemma~\ref{L:NULLCONDFACT}) 
	that has the added advantage of being constant along the integral curves of $\Lunit$
	(along which $u$ and $\vartheta$ are fixed).
	This property sometimes makes $\InitialFutFailFac$ slightly easier to work with  
	compared to $\FutFailFac.$
\end{remark}

\begin{remark}[\textbf{Past null condition failure factor}]
\label{R:PASTFAILUREFACTOR}
	We could also study shock formation in the region $\lbrace t \leq 0 \rbrace.$
	In this case, the relevant analog of $\FutFailFac$
	is the function $\PastFailFac,$ which is defined by replacing
	$\Lunit_{(Flat)}$ with $-\partial_t + \partial_r$
	in equation \eqref{E:FAILUREFACTOR}.
	Note that $-\partial_t + \partial_r$ is an outgoing Minkowski-null vectorfield
	in the region $\lbrace t \leq 0 \rbrace.$
	It is straightforward to see that
	$\FutFailFac$ is nontrivial
	if and only if $\PastFailFac$
	is nontrivial. In fact, 
	because the quantities $G_{\alpha \beta}(\Psi = 0)$
	in \eqref{E:FAILUREFACTOR} are constants,
	$\FutFailFac$
	and 
	$- \PastFailFac$
	have the same range.
\end{remark}

\section{Expressions for the spacetime metric}
In this section, we provide some expressions for the metrics and volume forms
that we use in our analysis. We begin by noting that from
\eqref{E:RADUISONERADTISZERO}, it follows that 
there exists an $S_{t,u}-$tangent vectorfield $\Xi$ such that
\begin{align} \label{E:RADINTERMSOFGEOMETRICCOORDINATEPARTIALDERIVATIVES}
	\Rad = \frac{\partial}{\partial u} - \Xi, 
\end{align}
where the $-$ sign is a convention. In the next lemma, we show that $\Xi$ verifies a transport equation 
along the integral curves of $\Lunit.$

\begin{lemma}[\textbf{Transport equation verified by $\Xi$}]
	\label{L:VECTORFIELDXITRANSPORT}
	The $S_{t,u}-$tangent vectorfield $\Xi$ verifies the following transport equation:
\begin{align} \label{E:XIEVOLUTION}
	\Lunit \Xi^A = [\Lunit, \Xi]^A = -[\Lunit, \Rad]^A,
\end{align}
where $[\cdot,\cdot]$ denotes the Lie bracket of two vectorfields.
\end{lemma}

\begin{proof}
	Lemma~\ref{L:VECTORFIELDXITRANSPORT}
	follows easily from equation \eqref{E:RADINTERMSOFGEOMETRICCOORDINATEPARTIALDERIVATIVES}
	and the identities 
	$[\frac{\partial}{\partial u}, \Lunit] = [\Lunit, X_A] = [\frac{\partial}{\partial u}, X_A] = 0.$
\end{proof}

We now provide the form of the spacetime metric relative to the geometric coordinates 
$(t,u,\vartheta^1,\vartheta^2).$

\begin{lemma}[\textbf{Form of the spacetime metric relative to geometric coordinates}]
We can express the spacetime metric 
relative to the geometric coordinates
$(t,u,\vartheta^1,\vartheta^2)$ as
\begin{align} \label{E:METRICRELATIVETOGEOMETRICCOORDINATES}
	g 
	& = - 2 \upmu dt du
		+ \upmu^2 du^2
		+ \gsphere_{AB}(d \vartheta^A + \Xi^A du) (d \vartheta^B + \Xi^B du).
\end{align}	
\end{lemma}
\begin{proof}
	Since $\Lunit = \frac{\partial}{\partial t}$ is null,
	it follows from \eqref{E:LISDDT} that the coefficient of $dt^2$ in \eqref{E:METRICRELATIVETOGEOMETRICCOORDINATES} is $0.$
	Next, since $g(\Lunit, X_A)=0,$ it follows from \eqref{E:LISDDT}
	that the coefficient of $dt d \vartheta^A$ in \eqref{E:METRICRELATIVETOGEOMETRICCOORDINATES} is $0.$
	Next, using 
	\eqref{E:LISDDT},
	\eqref{E:LRADINNERPRODUCT},
	and \eqref{E:RADINTERMSOFGEOMETRICCOORDINATEPARTIALDERIVATIVES}, 
	we compute that $g(\frac{\partial}{\partial t}, \frac{\partial}{\partial u}) 
	= g(\Lunit, \Rad) - \Xi^A g(\Lunit, X_A) = g(\Lunit, \Rad) = - \upmu.$
	This yields the term $- 2 \upmu dt du.$
	Similarly, it follows from \eqref{E:RADLENGTH} 
	and \eqref{E:RADINTERMSOFGEOMETRICCOORDINATEPARTIALDERIVATIVES}
	that
	$g(\frac{\partial}{\partial u}, \frac{\partial}{\partial u}) =  \upmu^2 + \Xi^A \Xi^B g(X_A, X_B) 
	= \upmu^2 + \gsphere_{AB} \Xi^A \Xi^B,$
	which yields the term $(\upmu^2 + \gsphere_{AB} \Xi^A \Xi^B)du^2.$
	Next, since
	$g(X_A,X_B) = \gsphere(X_A, X_B) = \gsphere_{AB},$ it follows 
	that the coefficient of $d \vartheta^A d \vartheta^B$ in \eqref{E:METRICRELATIVETOGEOMETRICCOORDINATES} is $\gsphere_{AB}.$
	Finally, since \eqref{E:RADINTERMSOFGEOMETRICCOORDINATEPARTIALDERIVATIVES} implies that
	$0 = g(\Rad, X_A) 
	= g(\frac{\partial}{\partial u}, X_A) - \Xi^B g(X_B,X_A)
	= g(\frac{\partial}{\partial u}, X_A) - \gsphere_{AB} \Xi^B,$ 
	we conclude that the coefficient 
	of $du d \vartheta^A$ in \eqref{E:METRICRELATIVETOGEOMETRICCOORDINATES} is 
	equal to $2 \gsphere_{AB} \Xi^B.$
	
\end{proof}

We now provide expressions for the geometric volume form factors of $g$ and $\gt.$

\begin{corollary}[\textbf{The geometric volume form factors of} $g$ \textbf{and} $\gt$]
\label{C:SPACETIMEVOLUMEFORMWITHUPMU}
Let $\gt$ and $\gsphere$ be as in Def.~\ref{D:FIRSTFUND}.
The following identity is verified by the spacetime metric $g:$
\begin{align} \label{E:SPACETIMEVOLUMEFORMWITHUPMU}
	|\mbox{\upshape{det}} g| 
	& = \upmu^2 \mbox{\upshape{det}} \gsphere,
\end{align}
where the determinant on the left-hand side is taken
relative to the geometric coordinates
$(t,u,\vartheta^1,\vartheta^2)$
and the determinant on the right-hand side 
is taken
relative to the geometric coordinates $(\vartheta^1,\vartheta^2).$

Furthermore, the following identity is verified by the 
first fundamental form $\gt$ of $\Sigma_t^{U_0}:$
\begin{align} \label{E:SIGMATVOLUMEFORMWITHUPMU}
	\mbox{\upshape{det}} \gt 
	& = \upmu^2 \mbox{\upshape{det}} \gsphere.
\end{align}
where the determinant on the left-hand side is taken
relative to the geometric coordinates
$(u,\vartheta^1,\vartheta^2)$ induced on $\Sigma_t^{U_0}$
and the determinant on the right-hand side 
is taken relative to the geometric coordinates $(\vartheta^1,\vartheta^2).$

\end{corollary}

\begin{proof}
Since both sides of \eqref{E:SPACETIMEVOLUMEFORMWITHUPMU}
have the same transformation properties
with respect to changes of angular coordinates of the form
$\widetilde{\vartheta}^1 = f^1(u,\vartheta^1, \vartheta^2),$
$\widetilde{\vartheta}^2 = f^2(u,\vartheta^1, \vartheta^2)$
(where the $f^A$ are the coordinate transformation functions),
it suffices to show that
\eqref{E:SPACETIMEVOLUMEFORMWITHUPMU}
holds for a well-chosen version 
$(\widetilde{\vartheta}^1, \widetilde{\vartheta}^2)$
of such angular coordinates.
To this end, we fix $t$ and endow
$\Sigma_t^{U^0}$ with local coordinates
$(u, \widetilde{\vartheta}^1, \widetilde{\vartheta}^2)$
(where $u$ is the usual eikonal function)
in such a way that $\Rad = \frac{\partial}{\partial u}|_{\widetilde{\vartheta}^1, \widetilde{\vartheta}^2}.$
To achieve this construction, we set $\widetilde{\vartheta}^A := \vartheta^A$
along $S_{t,0}$ and then (inwardly) transport 
$\widetilde{\vartheta}^A$ with $\Rad,$ that is, $\Rad \widetilde{\vartheta}^A = 0.$
We then extend the coordinates $\widetilde{\vartheta}^A$ off of $\Sigma_t^{U_0}$
by demanding that $[\Lunit, \widetilde{\vartheta}^A] = 0.$
Relative to these new local coordinates $(t,u,\widetilde{\vartheta}^1,\widetilde{\vartheta}^2),$ 
the identity \eqref{E:METRICRELATIVETOGEOMETRICCOORDINATES} holds 
along $\Sigma_t^{U^0},$ but with
$\gsphere_{AB}$ replaced by $g(\frac{\partial}{\partial \widetilde{\vartheta}^A}, \frac{\partial}{\partial \widetilde{\vartheta}^B}),$
$d \vartheta$ replaced by $d \widetilde{\vartheta},$
and $\Xi$ replaced by $0.$ Hence, relative to these new coordinates,
the identity \eqref{E:SPACETIMEVOLUMEFORMWITHUPMU} 
follows from a simple computation based on the block form 
of the metric (along $\Sigma_t^{U_0}$).

The identity \eqref{E:SIGMATVOLUMEFORMWITHUPMU} can be proved in a similar fashion, and we omit the details.

\end{proof}

\begin{lemma}[\textbf{The metrics in terms of the frame vectorfields}]
\label{L:SPACETIMEMETRICFRAMEVECTORFIELDS}
Let $\Lunit$ and $\uLgood$ be the vectorfields 
from Def.~\ref{D:NULLVECTORFIELDS},
and let $\gsphere$ be the metric induced on $S_{t,u}$ by $g,$ 
as in Def.~\ref{D:FIRSTFUND}. Then the following identities hold
(we note that the meaning of $\gsphere_{\mu \nu}$ and $(\ginversesphere)^{\mu \nu}$
is clarified in Remark~\ref{R:STUTENSORSDUALIDNETITIES}):
\begin{subequations}
\begin{align} 
	g_{\mu \nu} 
	& =  - \frac{1}{2} \upmu^{-1} \Lunit_{\mu} \uLgood_{\nu} 
		- \frac{1}{2} \upmu^{-1} \uLgood_{\mu} \Lunit_{\nu} + \gsphere_{\mu \nu},
		\label{E:GNULLFRAME} \\
	(g^{-1})^{\mu \nu} 
	& = - \frac{1}{2} \upmu^{-1} \left(\Lunit^{\mu} \uLgood^{\nu} 
			+ \uLgood^{\mu} \Lunit^{\nu} \right) 
			+ (\ginversesphere)^{\mu \nu}.
			\label{E:GINVERSENULLFRAME}
\end{align}
\end{subequations}

In addition, let
$\Rad$ and $\Radunit$ be the vectorfields
from Def.~\ref{D:RADIALVECTORFIELDSANDTIMELIKENORMAL}.
Then the following identities hold:
\begin{subequations}
\begin{align} \label{E:METRICFRAMEDECOMP}
	g_{\mu \nu} 
	& = - \Lunit_{\mu} \Lunit_{\nu}
			- (\Lunit_{\mu} \Radunit_{\nu} 
			+ \Radunit_{\mu} \Lunit_{\nu})
			+ \gsphere_{\mu \nu} \\
		& = - \Lunit_{\mu} \Lunit_{\nu}
			- \upmu^{-1}
				(
					\Lunit_{\mu} \Rad_{\nu} 
					+ \Rad_{\mu} \Lunit_{\nu}
				)
			+ \gsphere_{\mu \nu},
				\label{E:METRICFRAMEDECOMPLUNITRADUNITFRAME} \\
	(g^{-1})^{\mu \nu}  
		& = 
			- \Lunit^{\mu} \Lunit^{\nu}
			- (
					\Lunit^{\mu} \Radunit^{\nu} 
					+ \Radunit^{\mu} \Lunit^{\nu}
				)
			+ (\ginversesphere)^{\mu \nu}
			\label{E:GINVERSEFRAMEWITHRECTCOORDINATESFORGSPHEREINVERSE}	\\
		& = 
		- \Lunit^{\mu} \Lunit^{\nu}
		- \upmu^{-1} 
			(
				\Lunit^{\mu} \Rad^{\nu} 
				+ \Rad^{\mu} \Lunit^{\nu}
			 )
		+ (\ginversesphere)^{\mu \nu}.
		\label{E:GINVERSELUNITRADUNITFRAME}
\end{align}
\end{subequations}

Finally, we have
\begin{align} \label{E:GINVERSSPHEREINTERMSOFVECTORFIELDS}
	(\ginversesphere)^{\mu \nu} 
	& = (\ginversesphere)^{AB} X_A^{\mu} X_B^{\nu}.
\end{align}

\end{lemma}

\begin{proof}
The identity \eqref{E:GNULLFRAME} can easily be verified by contracting both sides 
against all possible pairs of vectors belonging to the rescaled null frame 
$\lbrace \Lunit, \uLgood, X_1, X_2 \rbrace$
and computing that both sides agree.
The identity \eqref{E:GINVERSENULLFRAME} then follows from raising
the indices of both sides of \eqref{E:GNULLFRAME} with $g^{-1}$
and from the identity $(g^{-1})^{\mu \alpha} (g^{-1})^{\nu \beta} \gsphere_{\alpha \beta} = (\ginversesphere)^{\mu \nu}.$

The identities 
\eqref{E:METRICFRAMEDECOMP}-\eqref{E:GINVERSELUNITRADUNITFRAME}
then follow from
using \eqref{E:RADINTERMSOFLANDULGOOD} to substitute for 
$\uLgood$ in \eqref{E:GNULLFRAME} and \eqref{E:GINVERSENULLFRAME}.

\end{proof}

\section{Contraction and component notation}
	\label{S:CONTRACTIONNOTATION}
In this section, we define some contraction and component notation that we use throughout the monograph.
\begin{definition}[\textbf{Contraction and component notation}]
	\label{D:CONTRACTION}
	If $\xi$ is a type $\binom{0}{2}$ spacetime tensor and 
	$V,$ $W$ are
	vectors, then we define the contraction
	\begin{align} \label{E:DOWNCONTRACTION}
		\xi_{VW} 
		&:= \xi_{\alpha \beta}V^{\alpha} W^{\beta}.
	\end{align}
	Similarly, if $\xi$ is a type $\binom{2}{0}$ spacetime tensor,
	then we define the contraction
	\begin{align}
		\xi_{VW}  \label{E:TYPEUPPERCONTRACTION}
		&:= \xi^{\alpha \beta} V_{\alpha} W_{\beta}.
	\end{align}
	We use similar contraction notation for tensors $\xi$ of any type.
	
	If $V$ is a spacetime vector, then we write
	\begin{align} \label{E:UPCOMPONENTS}
		V = V^{\Lunit} \Lunit + V^{\Rad} \Rad + V^A X_A
	\end{align}
	to denote the decomposition of $V$ relative to the rescaled frame
	$\lbrace \Lunit, \Rad, X_1, X_2 \rbrace.$
\end{definition}

\begin{remark}[$S_{t,u}$ \textbf{contraction abbreviations}]
	As we noted above, we often use abbreviations such as
	$\xi_A := \xi_{X_A}$ when contracting against the $S_{t,u}$ frame vectors
	$X_1$ and $X_2.$
\end{remark}

\begin{remark}[\textbf{Expansion of a vectorfield relative to the rescaled frame}]
	\label{R:UPDOWNFRAMETRANSLATION}
	Note that by Lemma~\ref{L:NULLPAIRBASICPROPERTIES},
	we have
	$V^{\Lunit} = - V_{\Lunit} - \upmu^{-1} V_{\Rad},$
	$V^{\Rad} = - \upmu^{-1} V_{\Lunit},$
	and $V^A = (\ginversesphere)^{AB}V_B.$
\end{remark}

\section{Projection operators and tensors along submanifolds}
In this section, we define two important projection operators.
We then provide closely related definitions of $\Sigma_t^{U_0}$ tensors and $S_{t,u}$ tensors.

\begin{definition}[\textbf{Projection operators}]
\label{D:PROJECTIONS}
We define 
the $\Sigma_t^{U_0}$ projection operator $\Sigmatproject$
and the $S_{t,u}$ projection operator
$\sphereproject$
to be the following type $\binom{1}{1}$ spacetime tensorfields  
$(\mu,\nu = 0,1,2,3):$
\begin{subequations}
\begin{align} 
	\Sigmatproject_{\nu}^{\ \mu} 
	&:=	\delta_{\nu}^{\ \mu}
			- \Timenormal_{\nu} \Timenormal^{\mu} 
		= \delta_{\nu}^{\ \mu}
			+ \delta_{\nu}^{\ 0} \Timenormal^{\mu},
			\label{E:SIGMATPROJECTION} \\
	\sphereproject_{\nu}^{\ \mu} 
	&:=	\delta_{\nu}^{\ \mu}
			+ \Radunit_{\nu} \Lunit^{\mu} 
			+ \Lunit_{\nu} (\Lunit^{\mu} + \Radunit^{\mu}) 
		= \delta_{\nu}^{\ \mu}
			- \delta_{\nu}^{\ 0} \Lunit^{\mu} 
			+  \Lunit_{\nu} \Radunit^{\mu}.
			\label{E:SPHEREPROJECTION}
	\end{align}
	\end{subequations}
	Above, $\Timenormal = \Lunit + \Radunit$ denotes the future-directed unit-normal to $\Sigma_t^{U_0}$
	and $\delta_{\nu}^{\ \mu}$ is the standard Kronecker delta.
	Furthermore, the second equality in \eqref{E:SPHEREPROJECTION} holds by virtue of
	\eqref{E:GINVERSE00ISMINUSONE}.
\end{definition}

It is straightforward to verify that
$\Sigmatproject_{\nu}^{\ \mu} \Timenormal^{\nu} = \Sigmatproject_{\nu}^{\ \mu} \Timenormal_{\mu} = 0$
and that $\Sigmatproject_{\nu}^{\ \mu} V^{\nu} = V^{\mu}$ for $\Sigma_t^{U_0}-$tangent vectors $V.$
Similarly, it is straightforward to verify that
$\sphereproject_{\nu}^{\ \mu} \Lunit^{\nu} 
= \sphereproject_{\nu}^{\ \mu} \Lunit_{\mu} 
= \sphereproject_{\nu}^{\ \mu} \Radunit^{\nu} 
= \sphereproject_{\nu}^{\ \mu} \Radunit_{\mu}
= 0$
and that $\sphereproject_{\nu}^{\ \mu} Y^{\nu} = Y^{\mu}$ for $S_{t,u}-$tangent vectors $Y.$

Furthermore, it is straightforward to 
verify that the following alternate expressions hold relative to the rectangular coordinates,
$(i,j=1,2,3$ and $\mu,\nu = 0,1,2,3):$
\begin{align} \label{E:ALTERNATEPROJECTIONEXPRESSIONS}
	\sphereproject_j^{\ i} 
		&	=	\delta_j^{\ i}
			- \Radunit_j \Radunit^i,
			\qquad \sphereproject_{\mu}^{\ 0} = 0,
			\qquad \sphereproject_0^{\ \nu} = 0,
\end{align}
where $\delta_j^{\ i}$ is the standard Kronecker delta.
In particular, the ``$0$ components'' of $\sphereproject$
since they all vanish. We often use this fact in our analysis.

\begin{definition}[\textbf{Projections of tensors onto} $\Sigma_t^{U_0}$ and $S_{t,u}$]
\label{D:PROJECTIONOFATENSOR}
If $\xi_{\nu_1 \cdots \nu_n}^{\nu_1 \cdots \mu_m}$
is a type $\binom{m}{n}$ spacetime tensor,
then the projections of $\xi$ onto 
$\Sigma_t^{U_0}$
and $S_{t,u}$
are 
respectively
the type $\binom{m}{n}$ tensors 
$\Sigmatproject \xi$
and
$\sphereproject \xi$ 
with the following components:
\begin{subequations}
\begin{align} 
(\Sigmatproject \xi)_{\nu_1 \cdots \nu_n}^{\mu_1 \cdots \mu_m}
& :=\Sigmatproject_{\nu_1}^{\ \widetilde{\nu}_1} \cdots \Sigmatproject_{\nu_n}^{\ \widetilde{\nu}_n} 
	\Sigmatproject_{\widetilde{\mu}_1}^{\ \mu_1} \cdots \Sigmatproject_{\widetilde{\mu}_m}^{\ \mu_m}
	\xi_{\widetilde{\nu}_1 \cdots \widetilde{\nu}_n}^{\widetilde{\mu}_1 \cdots \widetilde{\mu}_m},
		\\
(\sphereproject \xi)_{\nu_1 \cdots \nu_n}^{\mu_1 \cdots \mu_m}
& :=\sphereproject_{\nu_1}^{\ \widetilde{\nu}_1} \cdots \sphereproject_{\nu_n}^{\ \widetilde{\nu}_n} 
	\sphereproject_{\widetilde{\mu}_1}^{\ \mu_1} \cdots \sphereproject_{\widetilde{\mu}_m}^{\ \mu_m}
	\xi_{\widetilde{\nu}_1 \cdots \widetilde{\nu}_n}^{\widetilde{\mu}_1 \cdots \widetilde{\mu}_m}.
	\label{E:STUPROJECTIONOFATENSOR}
\end{align}
\end{subequations}
\end{definition}

\begin{definition}[\textbf{$S_{t,u}$ projection notation}]
	\label{D:STUSLASHPROJECTIONNOTATION}
	If $\xi$ is a spacetime tensor, then we define
	\begin{align} \label{E:TENSORSTUPROJECTED}
		\angxi := \sphereproject \xi.
	\end{align}
	
	If $\xi$ is a symmetric type $\binom{0}{2}$ spacetime tensor and $V$ is a spacetime
	vectorfield, then we define
	\begin{align} \label{E:TENSORVECTORANDSTUPROJECTED}
		\angxiarg{V} 
		& := \sphereproject (\xi_V),
	\end{align}
	where $\xi_V$ is the spacetime one-form with components
	$\xi_{\alpha \nu} V^{\alpha},$ $(\nu = 0,1,2,3).$
\end{definition}

In our analysis, we estimate two kinds of quantities: scalar functions and $S_{t,u}$ tensorfields. 

\begin{definition}[$\Sigma_t^{U_0}$ and $S_{t,u}$ \textbf{tensors}]
\label{D:HYPERSURFACETENSORS}
Let $\xi$ be a spacetime tensor. We say that $\xi$ is 
a $\Sigma_t^{U_0}$ tensor if
\begin{align} \label{E:SIGMATUTENSOR}
\Sigmatproject \xi 
& = \xi.
\end{align}
Similarly, we say that $\xi$ is
an $S_{t,u}$ tensor if 
\begin{align} \label{E:STUTENSOR}
\sphereproject \xi 
& = \xi.
\end{align}

\end{definition}
It is easy to show that $\xi$ is a $\Sigma_t^{U_0}$ tensor if and only if
any contraction of any downstairs index of $\xi$ against $\Timenormal$ necessarily results in $0,$
and any contraction of any upstairs index of $\xi$ against the $g-$dual of $\Timenormal$ 
necessarily results in $0.$
Similarly, it is easy to show that $\xi$ is an $S_{t,u}$ tensor if and only if
any contraction of any downstairs index of $\xi$ against either $\Lunit$ or $\Radunit$ necessarily results in $0,$
and any contraction of any upstairs index of $\xi$ against either the $g-$dual of $\Lunit$ or the $g-$dual of $\Radunit$ 
necessarily results in $0.$

\begin{remark}[\textbf{Inherent $S_{t,u}$ tensors vs. $S_{t,u}$ tensors embedded in spacetime}] \label{R:STUTENSORSDUALIDNETITIES}
	Throughout this monograph, we alternate back and forth between viewing $S_{t,u}$ tensors 
	as tensors that are inherent to the two-dimensional spheres $S_{t,u}$ with corresponding uppercase Latin indices $A,$ $B,$ etc.,
	and viewing $S_{t,u}$ tensors as spacetime tensors with the property \eqref{E:STUTENSOR}.
	That is, roughly speaking, we identify
	$\xi_{\widetilde{B}_1 \cdots \widetilde{B}_n}^{\widetilde{A}_1 \cdots \widetilde{A}_m} 
	\simeq \xi_{\nu_1 \cdots \nu_n}^{\mu_1 \cdots \mu_m}.$
	When it comes to performing calculations, both points of view can be helpful, depending on the situation at hand.
	Similar remarks apply to $\Sigma_t^{U_0}$ tensors.
\end{remark}

\section{The trace and trace-free parts of tensors}
\label{S:TRACEANDTRACEFREEPARTS}
We now provide some standard definitions connected to the trace and trace-free part
of tensors.

\begin{definition}[\textbf{Trace and trace-free parts of tensors}]
	\label{D:TRACEFREE}
	If $\xi$ is a symmetric type $\binom{0}{2}$ spacetime tensor, then
	\begin{align} \label{E:SPACETIMETRACEDEF}
		\myspacetimetr \xi := (g^{-1})^{\alpha \beta} \xi_{\alpha \beta}	
	\end{align}
	denotes its $g-$trace.
	
	If $\xi$ is a symmetric type $\binom{0}{2}$ $S_{t,u}$ tensor, then
	\begin{subequations}
	\begin{align} \label{E:TRACEDEF}
		\mytr \xi := (\ginversesphere)^{AB} \xi_{AB}	
	\end{align}
	denotes its $\gsphere-$trace and
	\begin{align} \label{E:TRACEFREE}
		\hat{\xi}_{AB} 
		& := \xi_{AB} - \frac{1}{2} \mytr \xi \gsphere_{AB}
	\end{align}
	\end{subequations}
	denotes its trace-free part.
	
	If $\xi$ and $\omega$ are $S_{t,u}$ covectors, then
	\begin{align} \label{E:SYMMETRIZEDTRACEFREETENSORPRODUCT}
		(\xi \hat{\otimes} \omega)_{AB} 
		& := \frac{1}{2}(\xi_A \omega_B + \xi_B \omega_A) - \frac{1}{2} \mytr (\xi \otimes \omega) \gsphere_{AB}
	\end{align}
	denotes the symmetrized trace-free part of the type $\binom{0}{2}$ $S_{t,u}$ tensor 
	$(\xi \otimes \omega)_{AB} := \xi_A \omega_B.$
	
\end{definition}

\section{Angular differential}
We now define the angular differential of a scalar-valued function.

\begin{definition}[\textbf{Angular differential}]
\label{D:ANGDIFFDEF}
If $f$ is a function, we define $\angdiff f$ to be the $S_{t,u}$ one-form 
\begin{align} \label{E:ANGDIFFDEF}
	\angdiff f := \sphereproject d f,
\end{align}
where $d f$ is the standard spacetime differential of $f$
and we are using the notation of Def.~\ref{D:PROJECTIONOFATENSOR}.
\end{definition}
Note that $\angdiffarg{A} f = X_A f$ and hence $\angdiff f$ can be viewed
as the standard angular differential of $f$ viewed as a function of the
geometric angular coordinates $(\vartheta^1, \vartheta^2).$

\begin{definition} [\textbf{Clarification of the meaning of} $\angdiffarg{i}$]
	\label{D:ANGDIFFIFUNCTION}
	If $f$ is a function, then for $i = 1,2,3,$
	\begin{align} \label{E:ANGDIFFIFUNCTION}
		\angdiffarg{i} f 
		:= (\angdiff f) \cdot \partial_i
		= \sphereproject_i^{\ a} \partial_a f
	\end{align}
	denotes the $i^{th}$ component of the $S_{t,u}$ one-form $\angdiff f$
	relative to the rectangular coordinate frame.
\end{definition}

In the next lemma, we compute the angular differential of the rectangular coordinate functions.

\begin{lemma}[\textbf{Angular differential of the $x^i$}]
Let $\sphereproject$ be the $S_{t,u}$ projection from Def.~\ref{D:PROJECTIONS}
and let $x^i$ be the rectangular spatial coordinate function.
Then the following identities hold relative to the rectangular coordinates
$(i,j = 1,2,3):$
\begin{align} \label{E:RECTANGULARANGDIFFCOORDINATE}
	\angdiffarg{j} x^i & = \sphereproject_j^{\ i}.
\end{align}

Furthermore, we have (for $A = 1,2$ and $i = 1,2,3$)
\begin{align} \label{E:ANGDIFFAXI}
	\angdiffarg{A} x^i & = X_A^i.
\end{align}

\end{lemma}

\begin{proof}
	To prove \eqref{E:RECTANGULARANGDIFFCOORDINATE},
	we use the definition of $\angdiff x^j$ and the identity $\partial_{\alpha} x^i = \delta_{\alpha}^i$ 
	to compute that
	$\angdiffarg{j} x^i = \sphereproject_j^{\ \alpha} \partial_{\alpha} x^i = \sphereproject_j^{\ i}$ as desired.
	To prove \eqref{E:ANGDIFFAXI}, we contract $X_A^j$ against both sides of \eqref{E:RECTANGULARANGDIFFCOORDINATE}.
	The left-hand side clearly results in $\angdiffarg{A} x^i,$ while
	since $X_A$ is $S_{t,u}-$tangent, the right-hand side results in $X_A^i$ as desired.
\end{proof}

\section{Musical notation}
\label{S:SHARPNOTATION}
In this section, we define our use of the symbols ``$\#$'' and ``$\flat.$''

\begin{definition}[\textbf{Sharp and flat notation}]
	\label{D:SHARPNOTATION}
	If $\xi$ is an $S_{t,u}$ covector, then we define
	$\xi^{\#}$ to be the $\gsphere-$dual of $\xi,$
	which is an $S_{t,u}-$tangent vector. That is,
	$(\xi^{\#})^A := (\ginversesphere)^{AB} \xi_B.$
	We often abbreviate $\xi^A := (\xi^{\#})^A.$
	
	If $Y$ is an $S_{t,u}-$tangent vector, then
	we define $Y_{\flat}$ to be the $\gsphere-$dual of $Y,$
	which is an $S_{t,u}$ covector. That is,
	$(Y_{\flat})_A := (\ginversesphere)_{AB} Y^B.$
	We often abbreviate $Y_A := (Y_{\flat})_A.$
	
	Similarly, if $\xi$ is a symmetric type $\binom{0}{2}$ tensor,
	then we define its $\gsphere-$dual $\xi^{\#}$ to be the type $\binom{1}{1}$ $S_{t,u}$ tensor
	with $A,B$ component equal to $(\ginversesphere)^{AC} \xi_{CB},$
	and its $\gsphere-$double dual $\xi^{\# \#}$ to be the symmetric type $\binom{2}{0}$ $S_{t,u}$ tensor
	with $A,B$ component equal to $(\ginversesphere)^{AC} (\ginversesphere)^{BD} \xi_{CD}.$
	
	We use similar notation to denote the $\gsphere-$duals of general type 
	$\binom{m}{0}$ and type $\binom{0}{n}$ $S_{t,u}$ tensors, and we use
	abbreviations similar to the ones mentioned above for vectors and covectors.
\end{definition}

\section{Pointwise norms}
Unless we explicitly state otherwise, 
we measure the pointwise norms of $S_{t,u}$ tensors 
relative to the metric $\gsphere.$ In the following definition,
we make this precise.

\begin{definition}[\textbf{Norm of} $S_{t,u}$ \textbf{tensorfields relative to} $\gsphere$]
\label{D:POINTWISENORMS}
Let $\xi$ be a type $\binom{m}{n}$ $S_{t,u}$ tensor 
with rectangular components 
$\xi_{\nu_1 \cdots \nu_n}^{\mu_1 \cdots \mu_m}$  
and $S_{t,u}$ components $\xi_{B_1 \cdots B_n}^{A_1 \cdots A_m}.$
We define $|\xi|^2,$ the square of the norm of $\xi,$ as follows:
\begin{align} \label{E:POINTWISENORMS}
	|\xi|^2
	& := 
	(\ginversesphere)^{\nu_1 \widetilde{\nu}_1} \cdots (\ginversesphere)^{\nu_n \widetilde{\nu}_n}
	\gsphere_{\mu_1 \widetilde{\mu}_1} \cdots \gsphere_{\mu_m \widetilde{\mu}_m}
	\xi_{\nu_1 \cdots \nu_n}^{\mu_1 \cdots \mu_m}
	\xi_{\widetilde{\nu}_1 \cdots \widetilde{\nu}_n}^{\widetilde{\mu}_1 \cdots \widetilde{\mu}_m}
		\\
	& = 
		(\ginversesphere)^{B_1 \widetilde{B}_1} \cdots (\ginversesphere)^{B_n \widetilde{B}_n}
		\gsphere_{A_1 \widetilde{A}_1} \cdots \gsphere_{A_m \widetilde{A}_m}
		\xi_{B_1 \cdots B_n}^{A_1 \cdots A_m}
		\xi_{\widetilde{B}_1 \cdots \widetilde{B}_n}^{\widetilde{A}_1 \cdots \widetilde{A}_m}.
		\notag
\end{align}
\end{definition}

\section{Lie derivatives and projected Lie derivatives}
In this section, we provide the standard definition of
the Lie derivative of a tensorfield. We then define
several related Lie derivative operators involving 
projections.

\begin{definition}[\textbf{Lie derivatives}]
\label{D:LIEDERIVATIVE}
If $V^{\mu}$ is a spacetime vectorfield and  
$\xi_{\nu_1 \cdots \nu_n}^{\mu_1 \cdots \mu_m}$
is a type $\binom{m}{n}$ spacetime tensorfield,
then relative the rectangular coordinates,
the Lie derivative of $\xi$ with respect to $V$ is
the type $\binom{m}{n}$ spacetime tensorfield $\Lie_V \xi$ with the following components:
\begin{align} \label{E:LIEDERIVATIVE}
\Lie_V \xi_{\nu_1 \cdots \nu_n}^{\mu_1 \cdots \mu_m}
& := V^{\alpha} \partial_{\alpha} \xi_{\nu_1 \cdots \nu_n}^{\mu_1 \cdots \mu_m}
	- \sum_{a=1}^m \xi_{\nu_1 \cdots \nu_n}^{\mu_1 \cdots \mu_{a-1} \alpha \mu_{a+1} \cdots \mu_m} \partial_{\alpha} V^{\mu_a}
	+ \sum_{b=1}^n \xi_{\nu_1 \cdots \nu_{b-1} \alpha \nu_{b+1} \cdots \nu_n}^{\mu_1 \cdots \mu_m} \partial_{\nu_b} V^{\alpha}.
\end{align}

In addition, when
$V$ and $W$ are both vectorfields,
we often use the standard Lie bracket notation
$[V,W] := \Lie_V W.$
\end{definition}


In our analysis, we often Lie differentiate tensorial products of $S_{t,u}$ tensors
and apply the Leibniz rule for Lie derivatives.
The non-$S_{t,u}$ components that arise upon Lie differentiating cancel 
out of such products.
This fact motivates the following definition.

\begin{definition}[\textbf{$\Sigma_t^{U_0}-$projected and $S_{t,u}-$projected Lie derivatives}]
\label{D:PROJECTEDLIE}
Let $\xi$ be a tensorfield.
We define 
$\SigmatLie_V \xi$ and
$\angLie_V \xi$ to respectively be the following 
$\Sigma_t^{U_0}$ and
$S_{t,u}$ tensorfields:
\begin{subequations}
\begin{align} \label{E:SIGMATTPROJECTEDLIE}
	\SigmatLie_V \xi
	& := \Sigmatproject \Lie_V \xi,
		\\
	\angLie_V \xi
	& := \sphereproject \Lie_V \xi,
	\label{E:STUPROJECTEDLIE}
\end{align}
\end{subequations}
where we are using the notation of Def.~\ref{D:PROJECTIONOFATENSOR}.
\end{definition}

\begin{lemma}[\textbf{Alternate expression for Lie derivatives of $S_{t,u}$ tensorfields}]
\label{L:ALTERNATELIEDERIVATIVEFORSTU}
If $\xi_{\nu_1 \cdots \nu_n}^{\mu_1 \cdots \mu_m}$
is a type $\binom{m}{n}$ $S_{t,u}$ tensorfield
and $X$ is an $S_{t,u}-$tangent vectorfield,
then relative to the rectangular coordinates, we have
\begin{align} \label{E:ALTERNATELIEDERIVATIVEFORSTU}
	\angLie_X \xi_{\nu_1 \cdots \nu_n}^{\mu_1 \cdots \mu_m}
		& = X^{\alpha} \angD_{\alpha} \xi_{\nu_1 \cdots \nu_n}^{\mu_1 \cdots \mu_m}
	- \sum_{a=1}^m \xi_{\nu_1 \cdots \nu_n}^{\mu_1 \cdots \mu_{a-1} \alpha \mu_{a+1} \cdots \mu_m} \angD_{\alpha} X^{\mu_a}
	+ \sum_{b=1}^n \xi_{\nu_1 \cdots \nu_{b-1} \alpha \nu_{b+1} \cdots \nu_n}^{\mu_1 \cdots \mu_m} \angD_{\nu_b} X^{\alpha}.
\end{align}
\end{lemma}

\begin{proof}
	By the torsion-free property of the Levi-Civita connection $\D$ of $g,$
	equation \eqref{E:LIEDERIVATIVE} holds for $V := X$ and with $\partial$ replaced by $\D$ on the 
	right-hand side. We then project both sides of the identity onto $S_{t,u}.$
	For $S_{t,u}$ tensorfields $\xi,$ $\D \xi$ and $\angD \xi$ differ only by terms that are $g-$orthogonal to $S_{t,u}.$
	Hence, the $S_{t,u}$ projection allows us to replace $\D$ with $\angD.$
	We have thus proved \eqref{E:ALTERNATELIEDERIVATIVEFORSTU}.
\end{proof}

\section{Second fundamental forms}
In this section, we provide the standard definition of the second
fundamental forms of $\Sigma_t^{U_0}$ and $S_{t,u}$ relative to 
the metric $g.$ We then provide some useful expressions for these
tensorfields.

\begin{definition}[\textbf{Second fundamental form of} $\Sigma_t^{U_0}$]
	\label{D:SECONDFUNDSIGMATDEF} 
	Let $\Timenormal = \Lunit + \Radunit$ be the future-directed unit normal to $\Sigma_t.$
	We define the second fundamental form $k$ of $\Sigma_t^{U_0}$ relative to $g$ to be
	the following type $\binom{0}{2}$ tensorfield:
	\begin{align}
		k 
		&:= \frac{1}{2} \SigmatLie_{\Timenormal} g.
			\label{E:SECONDFUNDSIGMATDEF} 
	\end{align}

\end{definition}

\begin{definition}[\textbf{Null second fundamental form of} $S_{t,u}$]
\label{D:CHI}
We define the null second fundamental form $\upchi$ of $S_{t,u}$ relative to $g$
to be the following type $\binom{0}{2}$ tensorfield:
\begin{align}
	\upchi 
	&:= \frac{1}{2} \angLie_{\Lunit} g.
	\label{E:CHIDEF}
\end{align}

\end{definition}

It follows in a straightforward fashion from the above two definitions 
and Def.~\ref{D:FIRSTFUND}
that the following alternate expressions hold:
\begin{align}
	k 
	& = \frac{1}{2} \SigmatLie_{\Timenormal} \gt,
		\label{E:KALTDEF} \\
	\upchi 
	& = \frac{1}{2} \angLie_{\Lunit} \gsphere.
		\label{E:CHIALTDEF}
\end{align}
In the next lemma, we provide additional expressions for $\upchi,$ $\angk,$ and $\angkarg{R},$
where $\angk$ and $\angkarg{R}$ are defined in terms of $k$ by Def.~\ref{D:STUSLASHPROJECTIONNOTATION}.

\begin{lemma}[\textbf{Alternate expressions for} $\upchi,$ $\angk,$ and $\angkarg{R}$]
	$\upchi$ and $\angk$ are symmetric type $\binom{0}{2}$ $S_{t,u}$ tensorfields that verify the following identities:
	\begin{align}
		\upchi_{AB}
		& = g(\D_A \Lunit, X_B),
			\label{E:CHIINTERMSOFCOVARIANTDERIVATIVES} \\
		\angkdoublearg{A}{B}
			& = g(\D_A \Timenormal, X_B).
			\label{E:ANGKINTERMSOFCOVARIANTDERIVATIVES}
	\end{align}
	
	Furthermore, $\angkarg{R}$ is an $S_{t,u}$ one-form that verifies the following identity:
	\begin{align} \label{E:ANGKRINTERMSOFCOVARIANTDERIVATIVES}
		\angkdoublearg{R}{A}
		& = g(\D_A \Lunit, \Rad).
	\end{align}
	
\end{lemma}

\begin{proof}
The symmetry properties 
follow trivially from the definitions.
We now prove \eqref{E:ANGKINTERMSOFCOVARIANTDERIVATIVES}
Using definition \eqref{E:SECONDFUNDSIGMATDEF},
the Leibniz rule,
the torsion-free property 
$\D_V W - \D_W V = [V,W],$ 
the identity $[X_A, X_B] = 0,$
and the identity $g(\Timenormal, X_A) = 0,$
we deduce \eqref{E:ANGKINTERMSOFCOVARIANTDERIVATIVES} as follows:
\begin{align}
2 \angkdoublearg{A}{B} 
& = \SigmatLie_{\Timenormal} g(X_A,X_B)
	\\
& = \Timenormal [g(X_A,X_B)] 
		- g([\Timenormal, X_A],X_B) 
		- g(X_A,[\Timenormal, X_B])
		\notag \\	
& = g(\D_{\Timenormal} X_A,X_B)
	+  g(X_A,\D_{\Timenormal} X_B)
	- g([\Timenormal, X_A],X_B) 
	- g(X_A,[\Timenormal, X_B])
		\notag \\
& = g(\D_A \Timenormal, X_B)
	+  g(\D_B \Timenormal, X_A)
	= g(\D_A \Timenormal, X_B)
	- g(\Timenormal, \D_B X_A)
		\notag \\
& = g(\D_A \Timenormal, X_B)
	- g(\Timenormal, \D_A X_B)
		= 2 g(\D_A \Timenormal, X_B).
	\notag
\end{align}

The proofs of 
\eqref{E:CHIINTERMSOFCOVARIANTDERIVATIVES} 
and
\eqref{E:ANGKRINTERMSOFCOVARIANTDERIVATIVES} 
are similar
(the proof of the former is based on definition \eqref{E:CHIDEF}),
and we omit the details.

\end{proof}

\section{Components of \texorpdfstring{$G$}{the first derivative of the metric with respect to the solution} 
relative to the non-rescaled frame \texorpdfstring{$\lbrace \Lunit, \Radunit, X_1, X_2 \rbrace$}{}}
In this section, we define the components of $G$ relative to the 
non-rescaled frame $\lbrace \Lunit, \Radunit, X_1, X_2 \rbrace.$

\begin{definition}[\textbf{Components of} $H$ \textbf{relative to the non-rescaled frame} $\lbrace \Lunit, \Radunit, X_1, X_2 \rbrace$]
\label{D:GFRAMECOMPONENTS}
Given any symmetric type $\binom{0}{2}$ spacetime tensorfield $H$ with rectangular components 
$H_{\mu \nu},$ we define the frame components of 
$H_{\mu \nu}$ (relative to the frame $\lbrace \Lunit, \Radunit, X_1, X_2 \rbrace$)
to be the scalar-valued functions
$H_{\Lunit \Lunit}:= H_{\alpha \beta} \Lunit^{\alpha} \Lunit^{\beta},$
$H_{\Lunit \Radunit} := H_{\alpha b} \Lunit^{\alpha} \Radunit^b,$
and
$H_{\Radunit \Radunit} := H_{ab} \Radunit^a \Radunit^b,$
the $S_{t,u}$ one-forms with frame components
$\angHdoublearg{\Lunit}{A} := H_{\alpha b} \Lunit^{\alpha} X_A^b,$
$\angHdoublearg{\Radunit}{A} := H_{\alpha b} \Radunit^{\alpha} X_A^b,$
and the symmetric type $\binom{0}{2}$ $S_{t,u}$ tensorfield with frame components
$\angHdoublearg{A}{B} := H_{a b} X_A^a X_B^b.$
We often respectively use the following abbreviations for the latter three tensorfields:
\begin{align} \label{E:FRAMECOMPONENTABBREVIATIONS}
	\angHarg{\Lunit},
	\angHarg{\Radunit},
	\angH.
\end{align}

We rarely need to distinguish between the various 
frame components of $H.$ Hence, to simplify the presentation,
it is convenient to 
place all of these components into an array
and to view them as a single entity.
Specifically, we define
$H_{(Frame)}$ and $H_{(Frame)}^{\#}$ 
to be the following arrays: 
\begin{subequations}
\begin{align}
H_{(Frame)} 
	&:= \left(
				H_{\Lunit \Lunit},
				H_{\Lunit \Radunit},
				H_{\Radunit \Radunit},
				\angHarg{\Lunit},
				\angHarg{\Radunit},
				\angH
			\right),
			\\
	H_{(Frame)}^{\#} 
	&:= \left(
				H_{\Lunit \Lunit},
				H_{\Lunit \Radunit},
				H_{\Radunit \Radunit},
				\angHmixedarg{\Lunit}{\#},
				\angHmixedarg{\Radunit}{\#},
				\angHuparg{\, \#} 
			\right).
\end{align}
\end{subequations}
\end{definition}

\begin{remark}[\textbf{$H$ is always $G$ or $G'$}]
	\label{R:ALWAYSGFRAMEORGPRIMEFRAME}
	In this monograph, the tensorfield
	$H_{\mu \nu}$ from Def.~\ref{D:GFRAMECOMPONENTS} will
	always be equal to either
	$G_{\mu \nu}$ or $G_{\mu \nu}',$ which are defined in 
	\eqref{E:BIGGDEF} and \eqref{E:BIGGPRIMEDEF}.
\end{remark}

\begin{definition}[\textbf{Norm of} $H_{(Frame)}$]
We define
\begin{align}
|H_{(Frame)}| 
	&:= 	|H_{\Lunit \Lunit}|
				+ |H_{\Lunit \Radunit}|
				+ |H_{\Radunit \Radunit}|
				+ |\angHarg{\Lunit}|
				+ |\angHarg{\Radunit}|
				+ |\angH|,
\end{align}
and similarly for $|H_{(Frame)}^{\#}|.$
\end{definition}

\begin{definition}[\textbf{Derivatives of} $H_{(Frame)}$]
If $V$ is a vectorfield, then we define
\begin{align}
\angLie_V H_{(Frame)} 
	&:= \left(
			V	H_{\Lunit \Lunit},
			V	H_{\Lunit \Radunit},
			V	H_{\Radunit \Radunit},
			\angLie_V	\angHarg{\Lunit},
			\angLie_V	\angHarg{\Radunit},
			\angLie_V	\angH
			\right),
\end{align}
and similarly for $\angLie_V H_{(Frame)}^{\#}.$
We use similar notation if the Lie derivative operator
is replaced with the $S_{t,u}$ covariant derivative operator 
$\angD.$

\end{definition}

\section{The change of variables map \texorpdfstring{$\Upsilon$}{}}
\label{S:CHOV}
In the next lemma, we compute the Jacobian determinant of the change of variables map
$(t,u,\vartheta^1,\vartheta^2) \overset{\Upsilon}{\rightarrow} (x^0,x^1,x^2,x^3).$
Later, we use the lemma to show that for 
the solutions of interest, $\Upsilon$ is regular as long as $\upmu > 0.$

\begin{lemma}[\textbf{The Jacobian determinant of the map} $(t,u,\vartheta^1,\vartheta^2) \rightarrow (x^0,x^1,x^2,x^3)$]
\label{L:JACOBIAN}
Let $\Upsilon:[0,T) \times [0,U_0] \times \mathbb{S}^2 \rightarrow \mathcal{M}_{T,U_0},$
$\Upsilon(t,u,\vartheta^1,\vartheta^2) = (x^0,x^1,x^2,x^3)$
be the change of variables map from geometric to rectangular coordinates.
Let $\gt$ and $\gsphere$ be as in Def.~\ref{D:FIRSTFUND}.
Then the Jacobian determinant of $\Upsilon$ can be expressed as
\begin{align} \label{E:JACOBIAN}
	\mbox{\upshape{det}}
	\frac{\partial (x^0,x^1,x^2,x^3)}{\partial (t,u,\vartheta^1, \vartheta^2)}
	= \upmu (\mbox{\upshape{det}} \gt)^{-1/2} \sqrt{\mbox{\upshape{det}} \gsphere},
\end{align}
where $(\mbox{\upshape{det}} \gt)^{-1/2}$ is a smooth function of $\Psi$ in a neighborhood of $0$ 
and verifies $(\mbox{\upshape{det}} \gt)^{-1/2}(\Psi=0) = 1.$ In \eqref{E:JACOBIAN},
$\mbox{\upshape{det}} \gt$ is taken relative to the rectangular spatial coordinates 
and	$\mbox{\upshape{det}} \gsphere$ is 
taken relative to the local geometric angular coordinates $(\vartheta^1,\vartheta^2).$
\end{lemma}

\begin{proof}
	The analysis already carried out thus far in Chapter~\ref{C:BASICGEOMETRICCONSTRUCTIONS} 
	implies that for $A = 1,2$ and $i = 1,2,3,$ we have
	\begin{align} \label{E:CHOVMATRIXCOMPUTATIONS}
		& \frac{\partial x^0}{\partial t} = 1, \,
			\frac{\partial x^0}{\partial u} = 0, \,
			\frac{\partial x^0}{\partial \vartheta^A} = 0, 
			\\
		& \frac{\partial x^i}{\partial t} = \Lunit^i, \,
			\frac{\partial x^i}{\partial u} = \Rad^i + \Xi^i, \,
			\frac{\partial x^i}{\partial \vartheta^A} = \angdiffarg{A} x^i = X_A^{\alpha} \partial_{\alpha} x^i = X_A^i, 
			\notag
	\end{align}
	where $\Xi$ is the $S_{t,u}-$tangent vectorfield from \eqref{E:RADINTERMSOFGEOMETRICCOORDINATEPARTIALDERIVATIVES}.
	
	Using \eqref{E:CHOVMATRIXCOMPUTATIONS}, we compute that
	\begin{align} \label{E:CHOVJACOBIAN}
		\mbox{\upshape{det}}
		\frac{\partial (x^0,x^1,x^2,x^3)}{\partial (t,u,\vartheta^1, \vartheta^2)}
		& = 
		\mbox{\upshape{det}}
		\left(
		\begin{array}{cccc}
			1 & 0 & 0 & 0 \\
			\Lunit^1 & \Rad^1 + \Xi^1 & X_1^1 & X_2^1 \\
			\Lunit^2 & \Rad^2 + \Xi^2 & X_1^2 & X_2^2 \\
			\Lunit^3 & \Rad^3 + \Xi^3 & X_1^3 & X_2^3
		\end{array}	
		\right) \\
		& = 
		\mbox{\upshape{det}}
		\left(
		\begin{array}{ccc}
			\Rad^1 & X_1^1 & X_2^1 \\
			\Rad^2 & X_1^2 & X_2^2 \\
			\Rad^3 & X_1^3 & X_2^3
		\end{array}	
		\right) 
		+
		\mbox{\upshape{det}}
		\left(
		\begin{array}{ccc}
			\Xi^1 & X_1^1 & X_2^1 \\
			\Xi^2 & X_1^2 & X_2^2 \\
			\Xi^3 & X_1^3 & X_2^3
		\end{array}	
		\right) 
		 \notag \\
	& = 
		\upmu 
		\mbox{\upshape{det}}
		\left(
		\begin{array}{ccc}
			\Radunit^1 & X_1^1 & X_2^1 \\
			\Radunit^2 & X_1^2 & X_2^2 \\
			\Radunit^3 & X_1^3 & X_2^3
		\end{array}	
		\right),
		\notag
	\end{align}
	where to deduce the last equality,
	we used the fact that $\Xi \in \mbox{span}\lbrace X_1, X_2 \rbrace.$
	The last determinant on the right-hand side of \eqref{E:CHOVJACOBIAN} can be 
	interpreted as 
	$\Euctvolform(\Radunit,X_1,X_2),$ where 
	$\Euctvolform$ is the volume form of the Euclidean metric $\Euct$ on $\Sigma_t$
	(where $\Euct_{ij} = \delta_{ij}$ relative to the rectangular spatial coordinates).
	We next recall the following standard fact: 
	$\gtvolform 
		= \sqrt{\mbox{\upshape{det}} \gt} \Euctvolform,$
	where $\gtvolform$ is the volume form of $\gt$ (see Def.~\ref{D:FIRSTFUND})
	and the determinant of $\gt$ is taken relative
	to rectangular coordinates. Since $\gt_{ij} = \delta_{ij} + g_{ij}^{(Small)}(\Psi),$
	it follows that $(\mbox{\upshape{det}} \gt)^{-1/2} = 1 + f(\Psi),$ where $f$ is smooth in a neighborhood of $0$ 
	and vanishes at $\Psi = 0.$ We furthermore note that since $\Radunit$ is the $\gt-$unit normal to
	$S_{t,u},$ it follows that $\gtvolform(\Radunit,X_1,X_2) = \spherevol(X_1,X_2) = \sqrt{\mbox{\upshape{det}} \gsphere}.$
	Here, $\spherevol$ is the volume form of $\gsphere$ (see Def.~\ref{D:FIRSTFUND})
	and	$\mbox{\upshape{det}} \gsphere$ 
	is taken relative to the geometric angular coordinates $(\vartheta^1,\vartheta^2).$
	Combining these identities, we conclude the desired identity \eqref{E:JACOBIAN}.
	
\end{proof}

\section{Area forms, volume forms, and norms}
	\label{S:VOLFORMSANDSOBOLEVNORMS}
In this section, we define the area forms, volume forms, and norms
that we use during our $L^2$ analysis of solutions.
We begin by noting that the results of Cor.~\ref{C:SPACETIMEVOLUMEFORMWITHUPMU} imply 
that relative to the geometric coordinates,
the geometric area and volume forms induced on 
$S_{t,u}$
$\Sigma_t^{U_0},$
and
$\mathcal{M}_{\Tboot,U_0}$
are respectively 
$d \argspherevol{(t',u',\vartheta)} := \sqrt{\mbox{\upshape{det}} \gsphere(t',u',\vartheta)} d \vartheta,$
$\upmu d \argspherevol{(t',u',\vartheta)} \, du',$
and $\upmu d \argspherevol{(t',u',\vartheta)} \, du' \, dt'.$
However, because the factor $\upmu$ in the latter two forms is extremely important,
we now define rescaled volume forms in which the factor is missing.
Throughout the monograph, 
all integrals are defined relative to the rescaled forms.
In particular, \textbf{we always explicitly indicate the factor in integrands when it is present}. 
We stress that we do \emph{not} rescale the area form on $S_{t,u}.$

\begin{definition}[\textbf{Area forms and rescaled volume forms on} 
	\label{D:AREAANDRESCALEDVOLUMEFORM}
	$S_{t,u},$ $\Sigma_t^u,$ $\mathcal{C}_u^t,$ \textbf{and} $\mathcal{M}_{t,u}$]
	We respectively define the following area and volume forms
	on
	$S_{t,u},$ $\Sigma_t^u,$ $\mathcal{C}_u^t,$ and $\mathcal{M}_{t,u}:$
	\begin{subequations}
	\begin{align}
		d \argspherevol{(t',u',\vartheta)} & = \sqrt{\mbox{\upshape{det}} \gsphere(t',u',\vartheta)} d \vartheta,
			\label{E:SPHEREAREAFORM} \\
		d \tvol & = d \argspherevol{(t',u',\vartheta)} \, du',
			\label{E:RESCALEDSIGMATVOLUMEFORM} \\
		d \conevol & = d \argspherevol{(t',u',\vartheta)} \, dt',
			\label{E:RESCALEDCONEVOLUMEFORM} \\
		d \vol & = d \argspherevol{(t',u',\vartheta)} \, du' \, dt'.
			\label{E:RESCALEDSPACETIMEVOLUMEFORM}
	\end{align}
	\end{subequations}
	
\end{definition}

We now define our Sobolev norms. The important point is that 
all norms are defined relative to the area and volume forms of 
Def.~\ref{D:AREAANDRESCALEDVOLUMEFORM}.
\begin{definition}[\textbf{Sobolev norms relative to the area and rescaled volume forms}]
	We define the norms
	$\| \cdot \|_{L^2(S_{t,u})},$
	$\| \cdot \|_{L^2(\Sigma_t^u)},$
	$\| \cdot \|_{L^2(\mathcal{C}_u^t)},$
	and 
	$\| \cdot \|_{\mathcal{M}_{t,u}}$
	as follows:
	\begin{subequations}
	\begin{align}
		\| f \|_{L^2(S_{t,u})}^2
		& := \int_{\mathbb{S}^2} f^2(t,u,\vartheta) d \spherevol
			= \int_{S_{t,u}} f^2 d \spherevol,
			\label{E:STUL2NORMDEF} \\
		\| f \|_{L^2(\Sigma_t^u)}^2
		& := \int_{u'=0}^u \int_{\mathbb{S}^2} f^2(t,u',\vartheta) d \spherevol \, du'
			= \int_{\Sigma_t^u} f^2 \, d \tvol,
			\label{E:SIGMATL2NORMDEF} \\
		\| f \|_{L^2(\mathcal{C}_u^t)}^2
		& := \int_{t'=0}^t \int_{\mathbb{S}^2} f^2(t',u,\vartheta) d \spherevol \, dt'
			= \int_{\mathcal{C}_u^t} f^2 \, d \conevol,
				\label{E:CONEL2NORMDEF} \\
		\| f \|_{\mathcal{M}_{t,u}}^2
		& := \int_{t'=0}^t \int_{u'=0}^u \int_{\mathbb{S}^2} f^2(t',u',\vartheta) d \spherevol \, du' \, dt'
			= \int_{\mathcal{M}_{t,u}} f^2 \, d \vol.
			\label{E:SPACETIMEL2NORMDEF}
	\end{align}
	\end{subequations}
	
	We similarly define $L^2$ norms over subsets $\Omega$ of the above sets.
	For example, if $\Omega \subset \Sigma_t^u,$ then
	\begin{align} \label{E:SUBSETL2NORM}
	\| f \|_{L^2(\Omega)}^2
		& := \int_{\Omega} f^2 \, d \tvol.
	\end{align}
	The relevant area/volume form 
	corresponding to the norm $\| \cdot \|_{L^2(\Omega)}$
	will always be clear from context.
\end{definition}

We now define our $C^0$ spaces and the corresponding norms.

\begin{definition}[$C^0$ \textbf{spaces and norm}]
	\label{D:C0NORM}
	Given any subset $\Omega$ of spacetime, 
	$C^0(\Omega)$ denotes the set of functions 
	of the geometric coordinates $(t,u,\vartheta)$
	that are continuous on $\Omega.$
	We define the corresponding norm
	$\| \cdot \|_{C^0(\Omega)}$ 
	as follows:
	\begin{align} \label{E:C0NORM}
		\| f \|_{C^0(\Omega)} 
		& := \sup_{p \in \Omega} |f(p)|.
	\end{align}
\end{definition}

\begin{remark}[\textbf{We never explicitly state that $f \in C^0(\Omega)$}] \label{R:C0NORMCONVENTION}
Whenever we state an estimate
that implies that $\| f \|_{C^0(\Omega)} < \infty,$
we adopt the convention that the estimate 
also carries with it the implication $f \in C^0(\Omega).$
That is, by our conventions, we only write estimates 
of the form $\| f \|_{C^0(\Omega)} < \infty$
whenever $f \in C^0(\Omega).$
\end{remark}

%

\section{Schematic notation}
\label{S:SCHEMATIC}
We often use schematic notation in our effort to 
highlight only the important features
of our equations and inequalities. In particular, we often 
encounter products of terms such that the numerical coefficients
and the precise tensorial structure are irrelevant from the point 
of view of proving our main sharp classical lifespan theorem. 
We indicate the basic structure 
of such terms in a very crude fashion by using ``schematic array notation''
\emph{without explicit reference to the precise frame components of the terms}.
For the terms whose precise structure is important,
\emph{we always either explicitly write their frame components or
explicitly mention that the term is ``exact'' including the numerical constants.}
We now give an example. Below, we derive equation \eqref{E:KABGOOD}, 
which precisely reads
\begin{align} \label{E:SCHEMATICKSPHEREGOOD}
	\angktriplearg{A}{B}{(Tan-\Psi)} 
	& = \frac{1}{2} \angGdoublearg{A}{B} \Lunit \Psi
			- \frac{1}{2} \angGdoublearg{\Lunit}{B} \angdiffarg{A} \Psi
			- \frac{1}{2} \angGdoublearg{\Lunit}{A} \angdiffarg{B} \Psi
			- \frac{1}{2} \angGdoublearg{\Radunit}{B} \angdiffarg{A} \Psi
			- \frac{1}{2} \angGdoublearg{\Radunit}{A} \angdiffarg{B} \Psi.
\end{align}
In view of Def.~\ref{D:GFRAMECOMPONENTS},
in order to indicate \eqref{E:SCHEMATICKSPHEREGOOD} schematically, 
we write
\begin{align} \label{E:KSHPEREGOODWRITTENSCHEMATICALLY}
	\angk 
	&= G_{(Frame)} 
		\myarray[\Lunit \Psi]
			{\angdiff \Psi}.
\end{align}
Note that we have suppressed all numerical factors in \eqref{E:KSHPEREGOODWRITTENSCHEMATICALLY}.
Whenever we use such schematic notation, the reader should interpret this as meaning that the numerical factors 
have no substantial effect on the estimates that
we derive for the schematic terms.
Also note that certain combinations suggested by
\eqref{E:KSHPEREGOODWRITTENSCHEMATICALLY} do not actually occur on the right-hand side of 
the precise formula \eqref{E:SCHEMATICKSPHEREGOOD}.
For example, there is no such term of the form
$G_{\Lunit \Radunit} \Lunit \Psi,$ and in fact, such a term does not even make sense 
since the left-hand side of \eqref{E:SCHEMATICKSPHEREGOOD} is a type $\binom{0}{2}$ $S_{t,u}$ tensorfield.
The reader should interpret this as meaning that
even if there were such a term, 
it would not affect our inequalities except for possibly 
altering them by unimportant constant factors.

We now give a few more examples. If we want to indicate the precise structure of 
the left-hand side of \eqref{E:SCHEMATICKSPHEREGOOD} and also the first term on the right, 
we write
\begin{align} \label{E:KSHPEREGOODWRITTENPARTIALLYSCHEMATICALLY}
	\angkdoublearg{A}{B} 
	&= \frac{1}{2} \angGdoublearg{A}{B} \Lunit \Psi
			+ G_{(Frame)} \angdiff \Psi.
\end{align}

To indicate the type $\binom{1}{0}$ $S_{t,u}$ tensorfield 
$\epsilon_{lca} g_{ab}^{(Small)} x^c \angdiffuparg{\#} x^b$ schematically, 
where $\epsilon_{\cdots}$ is the fully antisymmetric symbol normalized by $\epsilon_{123}=1$ 
and $g_{ab}^{(Small)} = g_{ab}^{(Small)}(\Psi)$ is defined in \eqref{E:LITTLEGDECOMPOSED},
we write
\begin{align}
	\smoothfunction(\Psi) \Psi x \angdiffuparg{\#} x,
\end{align}
where $\smoothfunction$ is a smooth function of $\Psi$
(recall that $g_{ab}^{(Small)}(0) = 0$).

To indicate the type $\binom{1}{1}$ $S_{t,u}$ tensorfield
$\frac{1}{2} \angGmixedarg{B}{A} \Lunit \Psi
- \frac{1}{2} \angGmixedarg{\Lunit}{A} \angdiffarg{B} \Psi
+ \frac{1}{2} \angGprimemixedarg{\Lunit}{A} \angdiffarg{B} \Psi$
schematically, we write
\begin{align} \label{E:ANEXAMPLESCHEMATICEXPRESSION}
	\myarray
		[G_{(Frame)}^{\#}]
		{G_{(Frame)}^{'\#}}
	\myarray[\Lunit \Psi]
		{\angdiff \Psi}.
\end{align}
Again, we stress that even though 
the notation \eqref{E:ANEXAMPLESCHEMATICEXPRESSION} suggests the presence
of some terms that are not actually in the precise expression,
the reader should interpret this as meaning that
the presence or absence of such terms will not affect any of the estimates
that we derive except for possibly 
altering them by unimportant constant factors.

The justification of our schematic treatment of many terms is based on the
fact that our estimates are, with only a few exceptions that we clearly point out,
\emph{not sensitive} to the precise tensorial structures present.
This will only become fully clear when we actually carry out the analysis.


\chapter{Transport Equations for the Eikonal Function Quantities}
\label{C:TRANSPORTEQUATIONSFORTHEEIKONALFUNCTIONQUANTITIES}
\thispagestyle{fancy}
In Chapter~\ref{C:TRANSPORTEQUATIONSFORTHEEIKONALFUNCTIONQUANTITIES}, we first
calculate the Christoffel symbols of $g$ relative to the Minkowski rectangular coordinates.
We then derive various transport equations for 
and angular differentials of $\upmu$ and the rectangular components 
$\Lunit^i,$ $
\Radunit^i$ of the frame vectorfields 
$\Lunit$ and $\Radunit.$
These quantities are first derivatives of the eikonal function, so in effect, we are deriving
commuted versions of the eikonal equation.
We also introduce the ``re-centered'' variables 
$\Lunit_{(Small)}^i,$ 
$\Radunit_{(Small)}^i,$ 
and $\upchi^{(Small)},$
which are closely connected to
$\Lunit^i,$ 
$\Radunit^i,$ 
and $\upchi$ but vanish for the background solution $\Psi \equiv 0$
(because to form the re-centered variables, 
we subtract off the background values from
$\Lunit^i,$ 
$\Radunit^i,$ 
and $\upchi$).
We also derive transport equations for and angular differentials of
some of the re-centered variables.
These transport equations are the evolution equations that we use to estimate
the eikonal function quantities, except at the top order.
We then provide an expression for $\upchi^{(Small)}$ in terms of $\Psi$
and $\Lunit_{(Small)}^i.$ In the remainder of the monograph, we use the expression for 
analyzing all derivatives of $\upchi_{(Small)}$ except those of top order.
Finally, we provide some Lie derivative identities that we 
use later in the monograph.

\section{Covariant derivatives and Christoffel symbols relative to the rectangular coordinates}

We often use the following lemma to compute the covariant derivatives of various vectorfields
relative to the rectangular coordinates.

\begin{lemma}[\textbf{Covariant derivatives and Christoffel symbols relative to the rectangular coordinates}]
Let $V$ be a spacetime vectorfield. Then relative to the rectangular coordinates, we have
\begin{subequations}
\begin{align} \label{E:RECTANGULARCOVARIANTDERIVATIVE}
	\D_{\mu} V^{\nu}
	& = \partial_{\mu} V^{\nu}
		+ (g^{-1})^{\nu \kappa} \Gamma_{\mu \kappa \alpha} V^{\alpha},
		\\
	\Gamma_{\alpha \kappa \beta}
	& = \frac{1}{2} 
			\left\lbrace
				\partial_{\alpha} g_{\kappa \beta}
				+ \partial_{\beta} g_{\alpha \kappa}
				- \partial_{\kappa} g_{\alpha \beta}
			\right\rbrace
		= \frac{1}{2} 
			\left\lbrace
				G_{\kappa \beta} \partial_{\alpha} \Psi
				+  G_{\alpha \kappa} \partial_{\beta} \Psi
				- G_{\alpha \beta}\partial_{\kappa} \Psi
			\right\rbrace. 
			\label{E:RECTANGULARCHRISTOFFEL}
\end{align}
\end{subequations}
\end{lemma}

\begin{proof}
	The identity \eqref{E:RECTANGULARCOVARIANTDERIVATIVE}
	and the first equality in \eqref{E:RECTANGULARCHRISTOFFEL}
	are standard identities from differential geometry.
	The second equality in \eqref{E:RECTANGULARCHRISTOFFEL}
	follows from the first equality, the chain rule,
	\eqref{E:LITTLEGDECOMPOSED},
	and \eqref{E:BIGGDEF}.
\end{proof}

\section{Transport equation for the inverse foliation density}
In this section, we derive the transport equation verified by $\upmu.$

\begin{lemma} [\textbf{The transport equation verified by $\upmu$}] \label{L:UPMUFIRSTTRANSPORT}
The inverse foliation density $\upmu$ defined in \eqref{E:UPMUDEF} verifies the following transport equation:
\begin{align} \label{E:UPMUFIRSTTRANSPORT}
	\Lunit \upmu 
	& = \upomega^{(Trans-\Psi)} 
		+ \upmu \upomega^{(Tan-\Psi)}
		:= \upomega,
\end{align}	
where
\begin{subequations}
\begin{align}
	\upomega^{(Trans-\Psi)} 
	& := \frac{1}{2} G_{\Lunit \Lunit} \Rad \Psi,
		\label{E:LITTLEOMEGATRANS} \\
	\upomega^{(Tan-\Psi)} 
	& := - \frac{1}{2} G_{\Lunit \Lunit} \Lunit \Psi
		- G_{\Lunit \Radunit} \Lunit \Psi.
		\label{E:LITTLEOMEGATAN}
\end{align}
\end{subequations}

\end{lemma}

\begin{proof}
Relative to rectangular coordinates, the $0$ component of equation \eqref{E:LGEOISGEODESIC} is
\begin{align} \label{E:LGEOGEODESIC0COMPONENT}
	\Lgeo \Lgeo^0 
	& = - (g^{-1})^{0 \gamma} \Gamma_{\alpha \gamma \beta} \Lgeo^{\alpha} \Lgeo^{\beta}.
\end{align}
Multiplying \eqref{E:LGEOGEODESIC0COMPONENT} by $\upmu^3,$ 
referring to definition \eqref{E:LUNIT},
using the decomposition \eqref{E:GINVERSEFRAMEWITHRECTCOORDINATESFORGSPHEREINVERSE}
and the identities $(\ginversesphere)^{0 \gamma} = 0,$ $\Lunit^0 = 1,$ $\Radunit^0=0,$ $\Rad = \upmu \Radunit,$
and using equation \eqref{E:RECTANGULARCHRISTOFFEL}, we deduce that
\begin{align} \label{E:SECONDVERSIONLGEOGEODESIC0COMPONENT}
\upmu^2 \Lunit \Lgeo^0 
& = \frac{1}{2}
		\left\lbrace
			\upmu \Lunit^{\gamma}
			+ \Rad^{\gamma} 
		\right\rbrace
	  \left\lbrace
			G_{\gamma \beta} \partial_{\alpha} \Psi
			+ G_{\alpha \gamma} \partial_{\beta} \Psi
			- G_{\alpha \beta}\partial_{\gamma} \Psi
		\right\rbrace 
		\Lunit^{\alpha} \Lunit^{\beta}
			\\
& = \frac{1}{2} \upmu G_{\Lunit \Lunit} \Lunit \Psi
		+ \upmu G_{\Lunit \Radunit} \Lunit \Psi
		- \frac{1}{2} G_{\Lunit \Lunit} \Rad \Psi.
		\notag
\end{align}
The desired equation \eqref{E:UPMUFIRSTTRANSPORT} now follows from \eqref{E:SECONDVERSIONLGEOGEODESIC0COMPONENT}
and the identities $\upmu^2 \Lunit \Lgeo^0 = \upmu^2 \Lunit (\frac{1}{\upmu}) = - \Lunit \upmu.$

\end{proof}

\section{Re-centered variables}

We now define re-centered versions 
$\Lunit_{(Small)}^i,$ 
$\Radunit_{(Small)}^i,$ 
and
$\upchi^{(Small)}$
of $\Lunit,$ 
$\Radunit,$ and 
$\upchi.$
The re-centered quantities vanish when $\Psi \equiv 0$
because to form them, we subtract off the background values from
$\Lunit^i,$ 
$\Radunit^i,$ 
and $\upchi.$

\begin{definition}[$\Lunit_{(Small)}^i,$ $\Radunit_{(Small)}^i,$ \textbf{and} $\upchi^{(Small)}$]
\label{D:RENORMALIZEDVARIABLES}
We define the scalar-valued functions 
$\Lunit_{(Small)}^i$ and $\Radunit_{(Small)}^i,$ $(i=1,2,3),$
and the symmetric type $\binom{0}{2}$ $S_{t,u}$ tensorfield $\upchi^{(Small)}$ as follows:
\begin{subequations}
\begin{align}
	\Lunit_{(Small)}^i
	& := \Lunit^i
		- \frac{x^i}{\rgeo},
		\label{E:LUNITJUNK} \\
	\Radunit_{(Small)}^i
	& := \Radunit^i
		+ \frac{x^i}{\rgeo},
		\label{E:RADUNITJUNK} \\
	\upchi^{(Small)}
		& := \upchi
		- \frac{\gsphere}{\rgeo}.
		\label{E:CHIJUNKDEF}
\end{align}	
\end{subequations}
	
\end{definition}

We often use the following lemma to reduce the analysis of
the rectangular components of $\Radunit$ to those of $\Lunit.$

\begin{lemma}[\textbf{Algebraic relationships between the rectangular components of $\Lunit$ and $\Radunit$}]
	\label{L:LUNITRADUNITALGEBRAICRELATION}
	The following identities hold relative to the rectangular coordinates, $(\mu = 0,1,2,3$ and $i=1,2,3):$
	\begin{subequations}
	\begin{align}
		\Radunit_{\mu} 
		& = - \Lunit_{\mu} - \delta_{\mu}^0,
		 \label{E:DOWNSTAIRSRADUNITPLUSLUNITISDELTA} \\
		\Radunit^i
		& = - \Lunit^i
			- (g^{-1})^{0i},
			\label{E:UPSTAIRSSRADUNITPLUSLUNITISAFUNCTIONOFPSI} 
	\end{align}
	\end{subequations}
	where $\delta_{\mu}^0$ is the standard Kronecker delta.
\end{lemma}
\begin{proof}
To prove \eqref{E:DOWNSTAIRSRADUNITPLUSLUNITISDELTA}, 
we use \eqref{E:UNITSIGMATNORMAL} and \eqref{E:TIMENORMALLENGTH} to deduce that
the one-form with rectangular components
$\Radunit_{\mu} + \Lunit_{\mu}$ is future-directed, length $-1,$ and $g-$normal to $\Sigma_t.$
We now note that the one-form $\xi$ with rectangular components $\xi_{\mu} = - \delta_{\mu}^0$ 
also is future-directed, 
$g-$normal to $\Sigma_t,$ 
and has length $(g^{-1})^{00} = -1$ (see \eqref{E:GINVERSE00ISMINUSONE}].
Hence, the identity \eqref{E:DOWNSTAIRSRADUNITPLUSLUNITISDELTA} holds.

To prove \eqref{E:UPSTAIRSSRADUNITPLUSLUNITISAFUNCTIONOFPSI}, we simply raise
the indices of \eqref{E:DOWNSTAIRSRADUNITPLUSLUNITISDELTA} with $g^{-1}$
to deduce that $\Radunit^i = - \Lunit^i - (g^{-1})^{i0}.$
\end{proof}

Note that Lemma~\ref{L:LUNITRADUNITALGEBRAICRELATION} implies that for $i=1,2,3,$ we have
\begin{align} \label{E:RADUNITJUNKLIKELMINUSUNITJUNK}
	\Radunit_{(Small)}^i
	& = - \Lunit_{(Small)}^i - (g^{-1})^{0i}.
\end{align}

\section{Transport equations for the rectangular frame components}

In this section, we derive transport equations for and angular differentials of the rectangular components
of $\Lunit,$ $\Radunit,$ and their re-centered versions.

\begin{proposition}[\textbf{Transport equations for and angular differentials of the rectangular frame components}] \label{P:FRAMECOMPONENTSTRANSPORT}
Let $\Lunit^i,$ $\Radunit^i,$ $\Lunit_{(Small)}^i,$ and $\Radunit_{(Small)}^i$
be the (scalar-valued) rectangular components given by
Defs.~\ref{D:NULLVECTORFIELDS}, \ref{D:RADIALVECTORFIELDSANDTIMELIKENORMAL}, and \ref{D:RENORMALIZEDVARIABLES}.
Then the following transport equations are verified by $\Lunit^i$ and $\Lunit_{(Small)}^i,$
$(i=1,2,3):$
\begin{subequations}
\begin{align}  \label{E:LUNITPROP}
	\Lunit \Lunit^i
	 & = \frac{1}{2} G_{\Lunit \Lunit} (\Lunit \Psi) \Radunit^i
			- \angGdoublearg{\Lunit}{A} (\ginversesphere)^{AB} (\angdiffarg{B} x^i) \Lunit \Psi
			+ \frac{1}{2} G_{\Lunit \Lunit} (\ginversesphere)^{AB} (\angdiffarg{A} x^i) \angdiffarg{B} \Psi,
			\\
	\Lunit (\rgeo \Lunit_{(Small)}^i)
	& = - \frac{1}{2} \rgeo G_{\Lunit \Lunit} (\Lunit \Psi) \Lunit_{(Small)}^i
			\label{E:LUNITJUNKLPROP} \\
	& \ \ - \frac{1}{2} G_{\Lunit \Lunit} (\Lunit \Psi) x^i
			- \frac{1}{2} \rgeo G_{\Lunit \Lunit}  (\Lunit \Psi) (g^{-1})^{0i} 
			\notag \\
	& \ \ + \frac{1}{2} \rgeo G_{\Lunit \Lunit} (\ginversesphere)^{AB} (\angdiffarg{A} \Psi) \angdiffarg{B} x^i
		- \rgeo (\ginversesphere)^{AB} \angGdoublearg{\Lunit}{A} (\angdiffarg{B} x^i) \Lunit \Psi.
		\notag 
\end{align}
\end{subequations}

Furthermore, there exist $S_{t,u}$ one-forms $\lambda^{(Tan-\Psi)}$ and $\uptheta^{(Tan-\Psi)}$
and symmetric type $\binom{0}{2}$ $S_{t,u}$ tensorfields
$\Lambda^{(Tan-\Psi)}$ 
and
$\Theta^{(Tan-\Psi)}$ 
such that the following expressions hold for the
angular differential 
(see Def.~\ref{D:ANGDIFFDEF})
of $\Lunit^i,$ 
$\Radunit^i,$ 
$\Lunit_{(Small)}^i,$ 
and $\Lunit^i$
$(i=1,2,3):$
\begin{subequations}
\begin{align}  
	\angdiffarg{A} \Lunit^i
	& = (\ginversesphere)^{BC} \upchi_{AB} \angdiffarg{C} x^i
	+ (\ginversesphere)^{BC} \Lambda_{AB}^{(Tan-\Psi)} \angdiffarg{C} x^i
	+ \uplambda_A^{(Tan-\Psi)} \Radunit^i,
		\label{E:ANGDIFFLI} \\
	\angdiffarg{A} \Radunit^i
	& = - (\ginversesphere)^{BC} \upchi_{AB} \angdiffarg{C} x^i
		+ (\ginversesphere)^{BC} \Theta_{AB}^{(Tan-\Psi)} \angdiffarg{C} x^i
		+ \theta_A^{(Tan-\Psi)} \Radunit^i,
		\label{E:ANGDIFFRADUNITI} \\
	\angdiffarg{A} (\Lunit_{(Small)}^i)
	& = (\ginversesphere)^{BC} \upchi^{(Small)}_{AB} \angdiffarg{C} x^i
		+ \uplambda_A^{(Tan-\Psi)} \Radunit_{(Small)}^i
		- \frac{x^i}{\rgeo} \uplambda_A^{(Tan-\Psi)}
		+ (\ginversesphere)^{BC} \Lambda_{AB}^{(Tan-\Psi)} \angdiffarg{C} x^i,
			\label{E:ANGDIFFLJUNKI} \\
	\angdiffarg{A} \Radunit_{(Small)}^i
	& = - (\ginversesphere)^{BC} \upchi^{(Small)}_{AB} \angdiffarg{C} x^i
		+ \uptheta_A^{(Tan-\Psi)} \Radunit_{(Small)}^i
		- \frac{x^i}{\rgeo} \uptheta_A^{(Tan-\Psi)}
		+ (\ginversesphere)^{BC} \Theta_{AB}^{(Tan-\Psi)} \angdiffarg{C} x^i,
		\label{E:ANGDIFFRADJUNKI} 
\end{align}
\end{subequations}
where
\begin{subequations}
	\begin{align}
		\uplambda_A^{(Tan-\Psi)} 
		& := - G_{\Lunit \Radunit} \angdiffarg{A} \Psi
			- \frac{1}{2} G_{\Radunit \Radunit} \angdiffarg{A} \Psi,
			\label{E:LITTLELAMBDAGOOD}  \\
		\uptheta_A^{(Tan-\Psi)} 
		& := - \frac{1}{2} G_{\Radunit \Radunit} \angdiffarg{A} \Psi,
			\label{E:LITTLETHETAGOOD} \\
		\Lambda_{AB}^{(Tan-\Psi)} 
		& := - \frac{1}{2} \angGdoublearg{A}{B} \Lunit \Psi
			+ \frac{1}{2} \angGdoublearg{\Lunit}{A} \angdiffarg{B} \Psi
			- \frac{1}{2} \angGdoublearg{\Lunit}{B} \angdiffarg{A} \Psi,
			\label{E:BIGLAMBDAGOOD} \\
		\Theta_{AB}^{(Tan-\Psi)} 
		& := \frac{1}{2} \angGdoublearg{A}{B} \Lunit \Psi
			- \frac{1}{2} \angGdoublearg{\Lunit}{A} \angdiffarg{B} \Psi 
			- \frac{1}{2} \angGdoublearg{\Lunit}{B} \angdiffarg{A} \Psi
			- \angGdoublearg{\Radunit}{B} \angdiffarg{A} \Psi.
			\label{E:BIGTHETAGOOD}
		\end{align}
\end{subequations}

\end{proposition}

\begin{proof}
We first prove \eqref{E:LUNITPROP}.
With $\nabla$ denoting the Levi-Civita connection of the Minkowski metric, we compute that
\begin{align} \label{E:LLUNITRECTANGLUARCHRISTOFFELEXPANSION}
	\Lunit \Lunit^{\nu}
	= \nabla_{\Lunit} \Lunit^{\nu}
	& = \D_{\Lunit} \Lunit^{\nu} 
		- (g^{-1})^{\nu \kappa} \Lunit^{\alpha} \Lunit^{\beta} \Gamma_{\alpha \kappa \beta}
		\\
	& = \D_{\Lunit} \Lunit^{\nu}	
		- \frac{1}{2} 
			(g^{-1})^{\nu \kappa} 
			\Lunit^{\alpha} 
			\Lunit^{\beta}
			\left\lbrace
				G_{\kappa \beta} \partial_{\alpha} \Psi
				+ G_{\alpha \kappa} \partial_{\beta} \Psi
				- G_{\alpha \beta} \partial_{\kappa} \Psi
			\right\rbrace.
			\notag
\end{align}		
Since $\Lunit^0 = 1,$ it follows that $\nabla_{\Lunit} \Lunit^{\nu}$ is $\Sigma_t-$tangent. 
Hence, we can expand $(i=1,2,3)$
\begin{align} \label{E:LLUNITIFRAMEEXPANSION}
	\nabla_{\Lunit} \Lunit^i = z \Radunit^i + Y^i,
\end{align}
where $z$ is a scalar and $Y$ is an $S_{t,u}-$tangent vectorfield.
Next, using
\eqref{E:LGEOISGEODESIC}, 
\eqref{E:LUNIT},
and \eqref{E:UPMUFIRSTTRANSPORT},
we deduce that
\begin{align} \label{E:DLLFIRSTCOMPUTATION}
	\D_{\Lunit} \Lunit^{\nu}
	& = \upmu^{-1} \upomega^{(Trans-\Psi)} \Lunit^{\nu}
		+ \upomega^{(Tan-\Psi)} \Lunit^{\nu}.
\end{align}
Contracting \eqref{E:LLUNITIFRAMEEXPANSION} against $\Rad_i$ 
and using
\eqref{E:LITTLEOMEGATRANS},
\eqref{E:LITTLEOMEGATAN},
\eqref{E:LLUNITRECTANGLUARCHRISTOFFELEXPANSION}, 
\eqref{E:DLLFIRSTCOMPUTATION},
and the identities $\Rad_i \Radunit^i = \upmu$
and $\Rad_i \Lunit^i = - \upmu,$
we deduce
\begin{align}
	\upmu z 
	& =  g(\nabla_{\Lunit} \Lunit,\Rad)
	= g(\D_{\Lunit} \Lunit,\Rad)
		- \frac{1}{2} 
			\Rad^{\kappa} 
			\Lunit^{\alpha} 
			\Lunit^{\beta}
			\left\lbrace
				G_{\kappa \beta} \partial_{\alpha} \Psi
				+ G_{\alpha \kappa} \partial_{\beta} \Psi
				- G_{\alpha \beta} \partial_{\kappa} \Psi
			\right\rbrace
			\label{E:LITTLEZCOMPUTATION} \\
		& = - \upomega^{(Trans-\Psi)} 
				- \upmu \upomega^{(Tan-\Psi)}
				- \upmu G_{\Lunit \Radunit} \Lunit \Psi	
				+ \frac{1}{2} G_{\Lunit \Lunit} \Rad \Psi
			= \frac{1}{2} \upmu G_{\Lunit \Lunit} \Lunit \Psi.
			\notag
\end{align}

Similarly, we compute that $Y^i = (\ginversesphere)^{AB} g(Y,X_A) X_B^i,$ where
\begin{align}
	g(Y,X_A)
	& = g(\nabla_{\Lunit} \Lunit,X_A)
		= g(\D_{\Lunit} \Lunit,X_A)
		- \frac{1}{2} 
			X_A^{\kappa} 
			\Lunit^{\alpha} 
			\Lunit^{\beta}
			\left\lbrace
				G_{\kappa \beta} \partial_{\alpha} \Psi
				+ G_{\alpha \kappa} \partial_{\beta} \Psi
				- G_{\alpha \beta} \partial_{\kappa} \Psi
			\right\rbrace
			\label{E:NABLALLINNERPRODUCTWITHXA} \\
		& = - \angGdoublearg{\Lunit}{A} \Lunit \Psi
			+ \frac{1}{2} G_{\Lunit \Lunit} \angdiffarg{A} \Psi.
			\notag
\end{align}	
From 
\eqref{E:LLUNITIFRAMEEXPANSION},
\eqref{E:LITTLEZCOMPUTATION},
\eqref{E:NABLALLINNERPRODUCTWITHXA},
and the identity $\angdiffarg{A} x^i = X_A^i,$
we conclude that
\begin{align}
	\Lunit \Lunit^i
	& =\nabla_{\Lunit} \Lunit^i 
		= \frac{1}{2} G_{\Lunit \Lunit} (\Lunit \Psi) \Radunit^i
			- \angGmixedarg{\Lunit}{A} (\angdiffarg{A} x^i) \Lunit \Psi
			+ \frac{1}{2} G_{\Lunit \Lunit} \angdiffuparg{A} \Psi \angdiffarg{A} x^i.
				\label{E:LLUNITIFIRSTCOMPUTATION} 
\end{align}	
We have thus proved \eqref{E:LUNITPROP}.

The desired identity \eqref{E:LUNITJUNKLPROP} then follows from \eqref{E:LUNITPROP}, 
Def.~\ref{D:RENORMALIZEDVARIABLES}, 
the identities 
$\Lunit x^i = \Lunit^i,$
$\Lunit \rgeo = 1,$
and
\eqref{E:RADUNITJUNKLIKELMINUSUNITJUNK},
and straightforward computations.

The proofs of \eqref{E:ANGDIFFLI} 
and \eqref{E:ANGDIFFLJUNKI}
are similar and rely on the identity
\eqref{E:CHIINTERMSOFCOVARIANTDERIVATIVES}
in place of \eqref{E:DLLFIRSTCOMPUTATION};
we omit the details.

The identities 
\eqref{E:ANGDIFFRADUNITI}
and \eqref{E:ANGDIFFRADJUNKI}
then follow from
\eqref{E:ANGDIFFLI}, 
\eqref{E:ANGDIFFLJUNKI},
Lemma~\ref{L:LUNITRADUNITALGEBRAICRELATION},
the chain rule identity
$- \angdiffarg{A} (g^{-1})^{0i} = (g^{-1})^{0\alpha} (g^{-1})^{i \beta} G_{\alpha \beta} \angdiffarg{A} \Psi,$
\eqref{E:GINVERSEFRAMEWITHRECTCOORDINATESFORGSPHEREINVERSE},
and straightforward computations.
	
\end{proof}

\section{An expression for \texorpdfstring{$\upchi_{(Small)}$}{the re-centered null second fundamental form} in terms of other quantities}
\label{S:CHISMALLINTERMSOFOTHERQUANTITIES}
In the next lemma, we show that $\upchi_{(Small)}$ is an auxiliary variable in the sense that
it is completely determined in terms of the up-to-second-order derivatives of $\Psi$ 
and the up-to-first-order derivatives of $\Lunit_{(Small)}^i,$ $(i=1,2,3).$

\begin{lemma}[\textbf{Expression for} $\upchi_{(Small)}$ \textbf{in terms of other quantities}]
\label{L:CHIJUNKINTERMSOFOTHERVARIABLES}
The quantities $\upchi^{(Small)}$ and $\mytr \upchi^{(Small)}$ can be expressed as follows
(see Def.~\ref{D:ANGDIFFIFUNCTION}):
\begin{subequations}
\begin{align} \label{E:CHIJUNKINTERMSOFOTHERVARIABLES}
	\upchi^{(Small)}_{AB} 
	& = 
		g_{ab} (\angdiffarg{A} x^a) \angdiffarg{B} \Lunit_{(Small)}^b
	 	- \Lambda_{AB}^{(Tan-\Psi)},
			\\
	\mytr \upchi^{(Small)}
	& = 
	 	\angdiffarg{a} \Lunit_{(Small)}^a
		- (\ginversesphere)^{AB} \Lambda_{AB}^{(Tan-\Psi)},
		\label{E:TRCHIJUNKINTERMSOFOTHERVARIABLES}
\end{align}
\end{subequations}
where $\Lambda_{AB}^{(Tan-\Psi)}$ is the $S_{t,u}$ tensorfield defined in \eqref{E:BIGLAMBDAGOOD}.

\end{lemma}

\begin{proof}
To derive \eqref{E:CHIJUNKINTERMSOFOTHERVARIABLES}, 
we contract \eqref{E:ANGDIFFLJUNKI} against
$g_{ij} \angdiffarg{D} x^j$ 
and use the identities
$g_{ij} (\angdiffarg{C} x^i) \angdiffarg{D} x^j = \gsphere_{CD}$
and
$g_{ij} (\Radunit_{(Small)}^i - \frac{x^i}{\rgeo}) \angdiffarg{D} x^j
= g_{ij} \Radunit^i \angdiffarg{D} x^j
= 0.$

To derive \eqref{E:TRCHIJUNKINTERMSOFOTHERVARIABLES}, we contract
\eqref{E:CHIJUNKINTERMSOFOTHERVARIABLES} against
$(\ginversesphere)^{AB}$ and use the identity
$(\ginversesphere)^{AB} g_{ab} (\angdiffarg{A} x^a) \angdiffarg{B}  = \angdiffarg{b}.$

\end{proof}

\section{Some identities involving deformation tensors and Lie derivatives}
\label{S:LIEDERIVATIVEIDENTITIES}
In this section, we provide some deformation tensor and Lie derivative identities
that play a role in our analysis.
We begin by providing the standard definition of a deformation tensor of a vectorfield.

\begin{definition}[\textbf{Deformation tensor}] \label{D:DEFORMTENSDEFINED}
We associate the following type $\binom{0}{2}$ tensorfield $\deform{V}$ to a vectorfield $V:$
\begin{align} \label{E:DEFORMTENSDEFINED}
	\deformarg{V}{\mu}{\nu} 
	:= \Lie_V g_{\mu \nu} 
	= \D_{\mu} V_{\nu} + \D_{\nu} V_{\mu}.
\end{align}
\end{definition}

\begin{lemma}[\textbf{Vectorfield commutator properties}]
\label{L:VECTORFIELDCOMMUTATORS}
Let $Y$ be any $S_{t,u}-$tangent vectorfield. Then 
$[\Lunit, \Rad],$
$[\Lunit, Y],$ 
and $[\Rad, Y]$ are also $S_{t,u}-$tangent,
and the following identities hold:
\begin{subequations}
\begin{align}
	[\Lunit, \Rad] 
	& = \Lie_{\Lunit} \Rad 
		= \angLie_{\Lunit} \Rad 
		= \angdeformoneformupsharparg{\Rad}{\Lunit}
		= - \angdeformoneformupsharparg{\Lunit}{\Rad},
		\label{E:LCOMMUTERADISSTUTANGENT} \\
	[\Lunit, Y] 
	& = \Lie_{\Lunit} Y 
		= \angLie_{\Lunit} Y 
		= \angdeformoneformupsharparg{Y}{\Lunit},
		\label{E:LCOMMUTETANGENTISTANGENT} \\
	[\Rad, Y] 
	& = \Lie_{\Rad} Y
		= \angLie_{\Rad} Y 
		= \angdeformoneformupsharparg{Y}{\Rad}.
		\label{E:RADCOMMUTETANGENTISTANGENT}
\end{align}
\end{subequations}
Above, the vectorfields $\angdeformoneformupsharparg{V}{W}$ 
are the $\gsphere-$duals of the one-forms
$\angdeformoneformarg{V}{W}$
defined by
\eqref{E:TENSORVECTORANDSTUPROJECTED}.

\end{lemma}

\begin{proof}
We first prove the equalities in \eqref{E:RADCOMMUTETANGENTISTANGENT}.
To prove the first two,
it suffices to show that 
$[\Rad, Y] t = [\Rad, Y] u = 0.$ These identities follow easily 
from the identities $Y t = Y u = \Rad t = 0$ and $\Rad u = 1.$
To deduce the final equality in \eqref{E:RADCOMMUTETANGENTISTANGENT},
it suffices to show that 
$\deformarg{Y}{\Rad}{A} = g([\Rad, Y],X_A)$
for $A=1,2.$ To this end, we
use the torsion-free property
$[\Rad, Y] = \D_{\Rad} Y - \D_Y \Rad,$
the identity $g(\Rad, X_A) = 0,$
and the fact that $[X_A,Y]$
is $S_{t,u}-$tangent
to compute that
\begin{align} 
	\deformarg{Y}{\Rad}{A} 
	: = g(\D_{\Rad} Y,X_A) + g(\D_A Y,\Rad)
	& = g([\Rad, Y],X_A) + g(\D_Y \Rad,X_A) + g(\D_A Y,\Rad)
		\\
	& = g([\Rad, Y],X_A) - g(\D_Y X_A, \Rad) + g(\D_A Y,\Rad)	
		\notag \\
	& = g([\Rad, Y],X_A) + g([X_A,Y],\Rad)
		\notag \\
	& = g([\Rad, Y],X_A)
		\notag
\end{align}
as desired.
The identities in \eqref{E:LCOMMUTETANGENTISTANGENT} follow similarly
with the help of the identities $\Lunit t = 1$ and $\Lunit u = 0.$
The identities in \eqref{E:LCOMMUTERADISSTUTANGENT} then follow similarly
with the help of the fact that 
by 
\eqref{E:LCOMMUTETANGENTISTANGENT} 
and
\eqref{E:RADCOMMUTETANGENTISTANGENT},
$[\Lunit,X_A]$ 
and 
$[\Rad,X_A]$ are $S_{t,u}-$tangent.

\end{proof}

\begin{lemma}[\textbf{Basic properties of the $S_{t,u}$ projection operator}]
Let $\sphereproject_{\nu}^{\ \mu}$ be the $S_{t,u}$ projection operator
from Def.~\ref{D:PROJECTIONS}.
If $Y$ is an $S_{t,u}-$tangent vectorfield,
and if $V \in \lbrace \Lunit, \rgeo \Lunit, \Rad \rbrace$
or $V$ is also an $S_{t,u}-$tangent vectorfield,
then we have $(\mu = 0,1,2,3):$
\begin{align} \label{E:LIEDERIVATIVEOFPROJECTIONSHASNOSTUCOMPONENET}
	Y^{\alpha} \Lie_V \sphereproject_{\alpha}^{\ \mu} & = 0.
\end{align}
Furthermore, we have
\begin{align} \label{E:PROJECTEDLIEDERIVATIVEOFSPHEREPROJECTIONSIS0}
	(\angdiffarg{A} x^a) \Lie_V \sphereproject_a^{\ \mu}
	& = 0.
\end{align}
Finally, we have $(\mu, \nu = 0,1,2,3):$
\begin{align} \label{E:STUPROJECTONCONTRACTEDAGAINSTLIEDERIVATIVESTUPROJECTIONIS0}
	\sphereproject_{\nu}^{\ \alpha} \Lie_V \sphereproject_{\alpha}^{\ \mu} & = 0.
\end{align}

\end{lemma}

\begin{proof}
Since $Y^{\mu} = Y^{\alpha} \sphereproject_{\alpha}^{\ \mu},$ the Leibniz rule
implies that
$\Lie_V Y^{\mu} 
= \Lie_V(Y^{\alpha} \sphereproject_{\alpha}^{\ \mu}) 
= (\Lie_V Y^{\alpha}) \sphereproject_{\alpha}^{\ \mu}
+ Y^{\alpha} \Lie_V \sphereproject_{\alpha}^{\ \mu}.$ 
On the other hand, Lemma~\ref{L:VECTORFIELDCOMMUTATORS} implies that
$\Lie_V Y^{\mu} = (\Lie_V Y^{\alpha}) \sphereproject_{\alpha}^{\ \mu}.$
The desired identity \eqref{E:LIEDERIVATIVEOFPROJECTIONSHASNOSTUCOMPONENET} now
follows from equating these two identities for $\Lie_V Y^{\mu}.$

To prove \eqref{E:PROJECTEDLIEDERIVATIVEOFSPHEREPROJECTIONSIS0}, we note 
that $\angdiffarg{A} x^a = X_A^b \partial_b x^a = X_A^a.$ Since for each fixed $A$ 
the vectorfield with rectangular spatial components $X_A^a$ is $S_{t,u}-$tangent, 
\eqref{E:PROJECTEDLIEDERIVATIVEOFSPHEREPROJECTIONSIS0}
follows from \eqref{E:LIEDERIVATIVEOFPROJECTIONSHASNOSTUCOMPONENET}.

To prove \eqref{E:STUPROJECTONCONTRACTEDAGAINSTLIEDERIVATIVESTUPROJECTIONIS0},
we simply note that for each fixed $\nu,$
$\sphereproject_{\nu}^{\ \alpha}$ can be viewed as an $S_{t,u}-$tangent
vectorfield. Hence,
\eqref{E:STUPROJECTONCONTRACTEDAGAINSTLIEDERIVATIVESTUPROJECTIONIS0}
follows from \eqref{E:LIEDERIVATIVEOFPROJECTIONSHASNOSTUCOMPONENET}.


\end{proof}

\begin{corollary}[\textbf{Sometimes} $S_{t,u}$ \textbf{projection is redundant}]
\label{C:ANGLIEVFORGOODVOFTENSOREQUALSANGLIEOFSTUPROJECTEDTENSOR}
If $\xi$ is a type $\binom{0}{n}$ spacetime tensorfield and
if $V \in \lbrace \Lunit, \rgeo \Lunit, \Rad \rbrace$
or $V$ is an $S_{t,u}-$tangent vectorfield, then 
\begin{align} \label{E:ANGLIEOFTENSOREQUALSANGLIEOFSTUPROJECTEDTENSOR}
	\angLie_V \xi = \angLie_V \angxi.
\end{align}

\end{corollary}

\begin{proof}
	 The identity \eqref{E:ANGLIEOFTENSOREQUALSANGLIEOFSTUPROJECTEDTENSOR} is
	 an easy consequence of
	 \eqref{E:STUPROJECTONCONTRACTEDAGAINSTLIEDERIVATIVESTUPROJECTIONIS0}.
\end{proof}

\begin{corollary}[\textbf{Basic commutation formula for} $S_{t,u}-$\textbf{projected Lie derivatives}]
\label{C:STUPROJECTEDDERIVATIVESCOMMUTATOR}
If $\xi$ is a type $\binom{0}{n}$ $S_{t,u}$ tensorfield,
if $V \in \lbrace \Lunit, \rgeo \Lunit, \Rad \rbrace$
or $V$ is an $S_{t,u}-$tangent vectorfield, 
and if $W \in \lbrace \Lunit, \rgeo \Lunit, \Rad \rbrace$
or $W$ is an $S_{t,u}-$tangent vectorfield, then 
\begin{align} \label{E:STUPROJECTEDDERIVATIVESCOMMUTATOR}
	[\angLie_V, \angLie_W] \xi 
	& = \angLie_{[V,W]} \xi.
\end{align}
\end{corollary}

\begin{proof}
	We give the proof in the case that $\xi$ is an $S_{t,u}$
	one-form; the general case can be handled in the same way.
	The standard commutation property for Lie derivatives is
	$[\Lie_V, \Lie_W] \xi = \Lie_{[V,W]} \xi.$
	Thus, a simple calculation yields that
	\begin{align}  \label{E:ALMOSTPROVEDSTUPROJECTEDDERIVATIVESCOMMUTATOR}
		[\angLie_V, \angLie_W] \xi_{\mu} 
		& = \angLie_{[V,W]} \xi_{\mu}
			+ \sphereproject_{\mu}^{\ \alpha} 
				(\Lie_V \sphereproject_{\alpha}^{\ \nu}) 
				\Lie_W \xi_{\nu}
			- \sphereproject_{\mu}^{\ \alpha} 
				(\Lie_W \sphereproject_{\alpha}^{\ \nu}) 
				\Lie_V \xi_{\nu},
	\end{align}
	where $\sphereproject_{\mu}^{\ \alpha}$ is the $S_{t,u}$ projection tensorfield from Def.~\ref{D:PROJECTIONS}.
	Using \eqref{E:STUPROJECTONCONTRACTEDAGAINSTLIEDERIVATIVESTUPROJECTIONIS0}, we see that
	the last two terms on the right-hand side of
	\eqref{E:ALMOSTPROVEDSTUPROJECTEDDERIVATIVESCOMMUTATOR} vanish.
	We have thus proved \eqref{E:STUPROJECTEDDERIVATIVESCOMMUTATOR}.
\end{proof}

\begin{lemma}[\textbf{Expressions for $S_{t,u}-$projected Lie derivatives}]
	\label{L:PROJECTEDLIEDERIVATIVESINTERMSOFOTHERVARIABLES}
	Let $\xi$ be a spacetime one-form with rectangular components $\xi_{\alpha} = \xi_{\alpha}(\Psi)$ 
	that are functions of $\Psi,$
	and let $\xi'$ be the one-form with rectangular components
	$\xi_{\alpha}' = \xi_{\alpha}'(\Psi) := \frac{d}{d \Psi} \xi_{\alpha}(\Psi).$ Let
	$\angxi,$ $\angxiprime$ denote the projections of $\xi,$ $\xi'$ onto $S_{t,u}$
	(see Def.~\ref{D:STUSLASHPROJECTIONNOTATION}). 
	Let $Z \in \lbrace \Lunit, \rgeo \Lunit, \Rad, \Rot_{(1)}, \Rot_{(2)}, \Rot_{(3)} \rbrace.$  
	Then
	\begin{align}
		\angLie_Z \angxi_A
		& := (\angLie_Z \angxi)_A
			= \angxiprimearg{A} Z \Psi
			+ \xi_a \angdiffarg{A} Z^a.
	\end{align}
	
	Furthermore, if $V \in \lbrace \Lunit, \Radunit, \frac{\partial}{\partial x^j} \rbrace_{j=1,2,3}$
 	then with 
	$\xi_V := \xi_{\alpha} V^{\alpha}$ 
	and $\xi_V' := \xi_{\alpha}' V^{\alpha},$
	we have
	\begin{align}
		Z \xi_V
		& := Z (\xi_V)
			= \xi_V' Z \Psi
			+ \xi_a Z V^a.
	\end{align}
	
	In addition, if $V \in \lbrace \Lunit, \Radunit, \frac{\partial}{\partial x^j} \rbrace_{j=1,2,3},$
	then with $G_{\mu \nu}(\Psi)$ and $G_{\mu \nu}'(\Psi)$ as defined in 
	Def.~\ref{D:BIGGANDBIGGPRIME}, we have
	\begin{subequations}
	\begin{align}
		Z G_{VW}
		& :=
		Z (G_{VW})
		= G_{VW}' Z \Psi
			+ G_{aW} Z V^a
		 	+ G_{Va} Z W^a,
		 \label{E:GFRAMESCALARDERIVATIVE}	\\
		\angLie_Z \angGdoublearg{V}{A}
		& :=
		(\angLie_Z (\angGarg{V}))_A
		= \angGprimedoublearg{V}{A} Z \Psi
			+ \angGdoublearg{A}{a} Z V^a
			+ G_{Va} \angdiffarg{A} Z^a,
				\label{E:GFRAMEONEFORMDERIVATIVE} \\
		\angLie_Z \angGdoublearg{A}{B}
		&:=
		(\angLie_Z \angG)_{AB}
		= \angGprimedoublearg{A}{B} Z \Psi
			+ \angGdoublearg{A}{a} \angdiffarg{B} Z^a
			+ \angGdoublearg{a}{B} \angdiffarg{A} Z^a.
			\label{E:GFRAMETYPE02DERIVATIVE}
	\end{align}
	\end{subequations}
	Above, $\angGarg{V}$ is the $S_{t,u}$ one-form 
	formed by projecting the one-form with rectangular components $G_{V \nu}$ onto $S_{t,u},$
	$\angG$ is the symmetric type $\binom{0}{2}$ $S_{t,u}$ tensorfield
	formed by projecting the symmetric type $\binom{0}{2}$ tensor with rectangular components $G_{\mu \nu}$ onto $S_{t,u},$
	$\angGprime$ is the symmetric type $\binom{0}{2}$ $S_{t,u}$ tensorfield
	formed by projecting the symmetric type $\binom{0}{2}$ tensor with rectangular components $G_{\mu \nu}'$ onto $S_{t,u},$
	and similarly for the other quantities (see Def.~\ref{D:STUSLASHPROJECTIONNOTATION}).
\end{lemma}
\begin{proof}
We prove only \eqref{E:GFRAMEONEFORMDERIVATIVE}; the proofs of the remaining identities are essentially the same.
We first note that the $i^{th}$ rectangular component of the $S_{t,u}$ one-form $\angGarg{V}$
is $G_{\alpha \beta} V^{\alpha} \sphereproject_i^{\beta},$ 
where $\sphereproject$ is the $S_{t,u}$ projection from Def.~\ref{D:PROJECTIONS}.
Hence, using the Leibniz and chain rules, the fact that the $0$ components of $\sphereproject$
are trivial, and the fact that the $0$ component of $V$ is constant,
we have
\begin{align} \label{E:FIRSTRELATIONGFRAMEONEFORMDERIVATIVE}
	\Lie_Z (G_{\alpha \beta} V^{\alpha} \sphereproject_i^{\ \beta})
	& = V^{\alpha} \sphereproject_i^{\ \beta} \Lie_Z G_{\alpha \beta}  
		+ \sphereproject_i^{\ \beta} G_{\alpha \beta} \Lie_Z V^{\alpha}
		+ G_{\alpha \beta} V^{\alpha} \Lie_Z \sphereproject_i^{\ \beta}
			\\
	& = V^{\alpha} \sphereproject_i^{\ \beta} \underbrace{Z G_{\alpha \beta}}_{G_{\alpha \beta}' Z \Psi} 
			+ \sphereproject_i^{\ \beta} V^{\alpha} G_{\alpha \gamma} \partial_{\beta} Z^{\gamma}
			+ \sphereproject_i^{\ \beta} G_{\beta \gamma} V Z^{\gamma}
			+ \sphereproject_i^{\ \beta} G_{\alpha \beta} Z V^{\alpha}
			\notag \\
	& \ \ - \sphereproject_i^{\ \beta} G_{\alpha \beta} V Z^{\alpha}
			+ G_{\alpha \beta} V^{\alpha} \Lie_Z \sphereproject_i^{\ \beta}.
			\notag
\end{align}
Contracting the left-hand side of \eqref{E:FIRSTRELATIONGFRAMEONEFORMDERIVATIVE}
against $\angdiffarg{A} x^i$ yields $(\angLie_Z (\angGarg{V}))_A.$
On the other hand, 
contracting the right-hand side of \eqref{E:FIRSTRELATIONGFRAMEONEFORMDERIVATIVE} 
against $\angdiffarg{A} x^i$ and using the identity $(\angdiffarg{A} x^i) \Lie_Z \sphereproject_i^{\ \beta} = 0$
(that is, \eqref{E:PROJECTEDLIEDERIVATIVEOFSPHEREPROJECTIONSIS0})
yields 
$\angGprimedoublearg{V}{A} Z \Psi
+ G_{V \gamma} \angdiffarg{A} Z^{\gamma} 
+ \angGdoublearg{A}{\gamma} V Z^{\gamma}
+	\angGdoublearg{a}{A} Z V^a
- \angGdoublearg{\alpha}{A} V Z^{\alpha}.$
Noting the cancellation of the third and fifth terms 
and noting that $\angdiff_A Z^0 = 0,$
we see that the desired identity \eqref{E:GFRAMEONEFORMDERIVATIVE} thus follows.

\end{proof}

\begin{corollary}[\textbf{The structure of the} $\Lunit$ \textbf{derivatives of some components of} $G_{(Frame)}$]
\label{C:GFAMELDERIVATIVE}
The following identities hold:
	\begin{subequations}	
\begin{align}
		\Lunit G_{\Lunit \Lunit}
		& =	
			\myarray
				[G_{(Frame)} G_{(Frame)}^{\#}]
				{G_{(Frame)}'}
			\myarray
				[\Lunit \Psi]
				{\angdiff \Psi},
				\label{E:LDERIVATIVEGLL} \\
	\angLie_{\Lunit} \angGdoublearg{\Lunit}{A}
	& = \angGmixedarg{\Lunit}{B} \upchi_{AB} 
	 	+ \myarray
				[G_{(Frame)} G_{(Frame)}^{\#}]
				{G_{(Frame)}'}
			\myarray
				[\Lunit \Psi]
				{\angdiff \Psi},
			 \label{E:LDERIVATIVEGLSPHERE} \\	
	\angLie_{\Lunit} \angGdoublearg{A}{B}
	& = \angGdoublearg{A}{C} \upchi_B^{\ C}
		+ \angGdoublearg{B}{C} \upchi_A^{\ C}
		+ \myarray
				[G_{(Frame)} G_{(Frame)}^{\#}]
				{G_{(Frame)}'}
			\myarray
				[\Lunit \Psi]
				{\angdiff \Psi},
		\label{E:LDERIVATIVEGSPHERE}
	\end{align}
	\end{subequations}
	where $G_{(Frame)},$ $G_{(Frame)}^{\#}$ and $G_{(Frame)}'$ are as in
	Defs.~\ref{D:BIGGANDBIGGPRIME} and \ref{D:GFRAMECOMPONENTS}.
	
\end{corollary}

\begin{proof}
We prove only \eqref{E:LDERIVATIVEGLSPHERE} since the proofs of the other identities
are similar. We first use \eqref{E:GFRAMESCALARDERIVATIVE} with $V = Z := \Lunit$ to deduce that
\begin{align} \label{E:FIRSTRELATIONLDERIVATIVEGLSPHERE}
		\angLie_{\Lunit} \angGdoublearg{\Lunit}{A}
		& = \angGprimedoublearg{\Lunit}{A} \Lunit \Psi
			+ \angGdoublearg{A}{a} \Lunit \Lunit^a
			+ G_{\Lunit a} \angdiffarg{A} \Lunit^a.
\end{align}
Inserting \eqref{E:LUNITPROP} and \eqref{E:ANGDIFFLI} into \eqref{E:FIRSTRELATIONLDERIVATIVEGLSPHERE},
we deduce that
\begin{align} \label{E:SECONDRELATIONLDERIVATIVEGLSPHERE}
\angLie_{\Lunit} \angGdoublearg{\Lunit}{A}
	& = \angGprimedoublearg{\Lunit}{A} \Lunit \Psi
		+ \frac{1}{2} G_{\Lunit \Lunit} \angGdoublearg{A}{\Radunit} \Lunit \Psi
		- \angGmixedarg{A}{B} \angGdoublearg{\Lunit}{B}  \Lunit \Psi
		+ \frac{1}{2} G_{\Lunit \Lunit} \angGmixedarg{A}{B} \angdiffarg{B} \Psi
			\\
	& \ \ + \angGmixedarg{\Lunit}{B} \upchi_{AB} 
		+ \angGmixedarg{\Lunit}{B} \Lambda_{AB}^{(Tan-\Psi)} 
		+ G_{\Lunit \Radunit} \uplambda_A^{(Tan-\Psi)}.
		\notag
\end{align}
The desired identity \eqref{E:LDERIVATIVEGLSPHERE} now follows from
\eqref{E:SECONDRELATIONLDERIVATIVEGLSPHERE},	
\eqref{E:LITTLELAMBDAGOOD},
and \eqref{E:BIGLAMBDAGOOD}.
	
\end{proof}


\chapter[Connection Coefficients and Geometric Decompositions of the Wave Operator]{Connection Coefficients of the Rescaled Frames and Geometric Decompositions of the Wave Operator}
\label{C:CONNECTIONCOEFF}
\thispagestyle{fancy}
In Chapter~\ref{C:CONNECTIONCOEFF}, we compute the connection coefficients of the two 
rescaled frames $\lbrace \Lunit, \Rad, X_1, X_2 \rbrace$ 
and $\lbrace \Lunit, \uLgood, X_1, X_2 \rbrace.$
We refer to these frames as ``rescaled'' because the vectorfields
$\Rad$ and $\uLgood$ are adapted to the behavior of $\upmu$ and hence can
significantly deviate from their Minkowskian counterparts,
which are $- \partial_r$ and $\partial_t - \partial_r.$
We then provide several decompositions of 
the $\upmu-$weighted covariant wave operator $\upmu \square_{g(\Psi)}$
relative to these frames.

\section{Connection coefficients of the rescaled frame \texorpdfstring{$\lbrace \Lunit, \Rad, X_1, X_2 \rbrace$}{}}
	
	In the next lemma, we compute the connection coefficients of the rescaled frame $\lbrace \Lunit, \Rad, X_1, X_2 \rbrace.$
	
	\begin{lemma}[\textbf{Connection coefficients of the rescaled frame $\lbrace \Lunit, \Rad, X_1, X_2 \rbrace$} 
		\textbf{and their decomposition into} $\upmu^{-1}-$\textbf{singular and} $\upmu^{-1}-$\textbf{regular pieces}]
	\label{L:CONNECTIONLRADFRAME}
	Let $\upzeta$ be the $S_{t,u}$ one-form defined by (see the identity \eqref{E:ANGKRINTERMSOFCOVARIANTDERIVATIVES})
	\begin{align} \label{E:ZETADEF}
		\upzeta_A := \angkdoublearg{R}{A}
		= g(\D_A \Lunit, \Radunit) = \upmu^{-1} g(\D_A \Lunit, \Rad).
	\end{align}
	Then the covariant derivatives of the frame vectorfields can be expressed as follows,
	where the tensorfields $k,$ $\upchi,$ and $\upomega$ are defined in 
	\eqref{E:SECONDFUNDSIGMATDEF}, 
	\eqref{E:CHIDEF},
	and \eqref{E:UPMUFIRSTTRANSPORT}:
	\begin{subequations}
	\begin{align}
		\D_{\Lunit} \Lunit 
		& = \upmu^{-1} \upomega \Lunit, 
			\label{E:DLL} \\
		\D_{\Rad} \Lunit 
		& = - \upomega \Lunit 
			+ \upmu \upzeta^A X_A
			+ (\angdiffuparg{A} \upmu) X_A, 
			\label{E:DRADL} \\
		\D_A \Lunit 
		& = - \upzeta_A \Lunit
			+ \upchi_A^{\ B} X_B, 
			\label{E:DAL} \\
		\D_{\Lunit} \Rad 
		& = - \upomega \Lunit 
			- \upmu \upzeta^A X_A, 
			\label{E:DLRAD} \\
		\D_{\Rad} \Rad 
		& = \upmu \upomega \Lunit
			+ \left\lbrace \upmu^{-1} \Rad \upmu  + \upomega \right\rbrace \Rad
			- \upmu (\angdiffuparg{A} \upmu) X_A, 
			\label{E:DRADRAD} \\
		\D_A \Rad 
		& = \upmu \upzeta_A \Lunit
			+ \upzeta_A \Rad
			+ \upmu^{-1} (\angdiffarg{A} \upmu) \Rad
			+ \upmu \angkmixedarg{A}{B} X_B
			- \upmu \upchi_A^{\ B} X_B,
			\label{E:DARAD} \\
		\D_{\Lunit} X_A 
		& = \D_A \Lunit, 
			\label{E:DLA} \\
		\D_A X_B 
		& = \angDarg{A} X_B
			+ \angkdoublearg{A}{B} \Lunit
			+ \upmu^{-1} \upchi_{AB} \Rad.
			\label{E:DAXBINTERMSOFANGDAXB}
	\end{align}
	\end{subequations}
	
	Furthermore, we can decompose the frame components of the tensorfields
	$k$ and $\upzeta$ into 
	$\upmu^{-1}-$singular and $\upmu^{-1}-$regular
	pieces as follows: 
	\begin{subequations}
	\begin{align}
		\upzeta_A & = \upmu^{-1} \upzeta_A^{(Trans-\Psi)}
		+ \upzeta_A^{(Tan-\Psi)},
		 \label{E:ZETADECOMPOSED} \\
	\angkdoublearg{A}{B} 
	& = \upmu^{-1} \angktriplearg{A}{B}{(Trans-\Psi)}
		+ \angktriplearg{A}{B}{(Tan-\Psi)},
		\label{E:ANGKDECOMPOSED}
	\end{align}
\end{subequations}
where
\begin{subequations}
	\begin{align}
		\upzeta_A^{(Trans-\Psi)} 
		& :=
			- \frac{1}{2} \angGdoublearg{\Lunit}{A} \Rad \Psi,
			\label{E:ZETATRANSVERSAL} \\
		\angktriplearg{A}{B}{(Trans-\Psi)} 
		& := \frac{1}{2} \angGdoublearg{A}{B} \Rad \Psi,
			\label{E:KABTRANSVERSAL}
	\end{align}	
\end{subequations}
and
\begin{subequations}
\begin{align}
	\upzeta_A^{(Tan-\Psi)}
	& := \frac{1}{2} \angGdoublearg{\Radunit}{A} \Lunit \Psi
			- \frac{1}{2} G_{\Lunit \Radunit} \angdiffarg{A} \Psi
			- \frac{1}{2} G_{\Radunit \Radunit} \angdiffarg{A} \Psi,
		\label{E:ZETAGOOD} \\
	\angktriplearg{A}{B}{(Tan-\Psi)} 
	& := \frac{1}{2} \angGdoublearg{A}{B} \Lunit \Psi
			- \frac{1}{2} \angGdoublearg{\Lunit}{B} \angdiffarg{A} \Psi
			- \frac{1}{2} \angGdoublearg{\Lunit}{A} \angdiffarg{B} \Psi
			- \frac{1}{2} \angGdoublearg{\Radunit}{B} \angdiffarg{A} \Psi
			- \frac{1}{2} \angGdoublearg{\Radunit}{A} \angdiffarg{B} \Psi.
			\label{E:KABGOOD}
\end{align}
\end{subequations}
\end{lemma}

\begin{proof}
	As an example, we prove \eqref{E:DARAD}. The remaining
	identities in \eqref{E:DLL}-\eqref{E:DAXBINTERMSOFANGDAXB} 
	can be proved using similar arguments.
	To proceed, we expand $\D_A \Rad$ relative to the frame as follows:
	\begin{align} \label{E:DARADFRAMEEXPANDED}
		\D_A \Rad 
		& = a^{\Lunit} \Lunit 
				+ a^{\Rad} \Rad 
				+ a^B X_B,
	\end{align}
	where $a^{\Lunit},$ $a^{\Rad},$ and $a^B$ are scalars to be determined.
	To compute $a^{\Rad},$ we first
	take the inner product of the left-hand side of \eqref{E:DARADFRAMEEXPANDED} with $\Lunit$ and compute that
	$g(\D_A \Rad, \Lunit) = g(\D_A (\upmu \Radunit), \Lunit) = (\angdiffarg{A} \upmu) g(\Radunit, \Lunit) - \upmu g(\Radunit, \D_A \Lunit)
	= - \angdiffarg{A} \upmu - \upmu \upzeta_A.$ On the other hand, the inner product of the right-hand side of \eqref{E:DARADFRAMEEXPANDED} with 	$\Lunit$
	is $- \upmu a^{\Rad}.$ Equating these two quantities, we deduce that $a^{\Rad} = \upmu^{-1} \angdiffarg{A} \upmu + \upzeta_A$ as desired.
	
	To compute $a^{\Lunit},$ we first 
	take the inner product of the left-hand side of \eqref{E:DARADFRAMEEXPANDED} with $\Rad$ and compute that
	$g(\D_A \Rad, \Rad) = \frac{1}{2} \angdiffarg{A} g(\Rad,\Rad) = \frac{1}{2} \angdiffarg{A} (\upmu^2) = \upmu \angdiffarg{A} \upmu.$
	On the other hand, the inner product of the right-hand side of \eqref{E:DARADFRAMEEXPANDED} with $\Rad$
	is $- \upmu a^{\Lunit} + \upmu^2 a^{\Rad}.$ Equating these two quantities and using the above expression for $a^{\Rad},$
	we deduce that $a^{\Lunit} = \upmu \upzeta_A$ as desired.
	
	To compute $a^B,$ 
	we take the inner product of the left-hand side of \eqref{E:DARADFRAMEEXPANDED} with $X_B$ and 
	use \eqref{E:CHIINTERMSOFCOVARIANTDERIVATIVES} and \eqref{E:ANGKINTERMSOFCOVARIANTDERIVATIVES} to
	compute that
	$a^C \gsphere_{BC} = g(\D_A \Rad, X_B) = g(\D_A (\upmu \Radunit), X_B) = \upmu g(\D_A \Radunit, X_B)
	= \upmu g(\D_A \Timenormal, X_B) - \upmu g(\D_A \Lunit, X_B)
	= \upmu \angkdoublearg{A}{B} - \upmu \upchi_{AB}$ as desired.
	
	As a second example, we prove 
	\eqref{E:ANGKDECOMPOSED},
	\eqref{E:KABTRANSVERSAL},
	and \eqref{E:KABGOOD}.
	The identities 
	\eqref{E:ZETADECOMPOSED},
	\eqref{E:ZETATRANSVERSAL},
	and \eqref{E:ZETAGOOD}
	can be proved using similar arguments.
	To proceed, we use definition \eqref{E:SECONDFUNDSIGMATDEF} 
	and the chain rule identity $\Timenormal g_{ij} = G_{ij} \Timenormal \Psi$
	to deduce that relative to rectangular coordinates, we have
	\begin{align} \label{E:SIGMATSECONDFUNDRECT}
		2 k_{ij} 
		&= G_{ij} \Timenormal \Psi
			+ g_{i \alpha} \partial_j \Timenormal^{\alpha}
			+ g_{j \alpha} \partial_i \Timenormal^{\alpha}.
	\end{align}
	Contracting \eqref{E:SIGMATSECONDFUNDRECT} against $X_A^i X_B^j,$
	recalling that $\Timenormal = \Lunit + \Radunit = \Lunit + \upmu^{-1} \Rad,$
	and using the decomposition \eqref{E:METRICFRAMEDECOMP},
	we deduce that
	\begin{align} \label{E:SIGMATSECONDFUNDAG}
		2 \angkdoublearg{A}{B} 
		&= \angGdoublearg{A}{B} 
				\left\lbrace 
					\upmu^{-1} \Rad \Psi
					+ \Lunit \Psi
				\right\rbrace
			+ \gsphere_{Aa} 
				\angdiffarg{B} 
					\left\lbrace 
						\Lunit^a
						+ \Radunit^a
					\right\rbrace
			+ \gsphere_{Ba} 
				\angdiffarg{A}
					\left\lbrace 
						\Lunit^a
						+ \Radunit^a
				\right\rbrace.
	\end{align}
	To conclude the desired identities
	\eqref{E:ANGKDECOMPOSED},
	\eqref{E:KABTRANSVERSAL},
	and \eqref{E:KABGOOD},
	we use equations
	\eqref{E:ANGDIFFLI}	
	and
	\eqref{E:ANGDIFFRADUNITI}
	to substitute for 
	$\angdiffarg{A} \Lunit^a,$
	$\angdiffarg{B} \Lunit^a,$
	$\angdiffarg{A} \Radunit^a,$
	and
	$\angdiffarg{B} \Radunit^a$
	in \eqref{E:SIGMATSECONDFUNDAG},
	use the identities
	$(\ginversesphere)^{BC} \gsphere_{Aa} \angdiffarg{C} x^a = \delta_A^B$
	and $\gsphere_{Aa} \Radunit^a = 0,$
	and perform straightforward calculations.
	
\end{proof}

\begin{corollary}[\textbf{An expression for $[\Lunit, \Rad]$}]
The following identity holds:
\begin{align}
	\Lie_{\Lunit} \Rad^A 
	& = [\Lunit, \Rad]^A
		= - (\ginversesphere)^{AB} \angdiffarg{B} \upmu
			- 2 \upmu (\ginversesphere)^{AB} \upzeta_B.
			\label{E:LRADCOMM}
\end{align}
	\end{corollary}
\begin{proof}
The corollary follows from the torsion-free property
$[\Lunit, \Rad] = \D_{\Lunit}\Rad - \D_{\Rad}\Lunit,$
\eqref{E:DRADL}, and \eqref{E:DLRAD}.
\end{proof}

\section{Connection coefficients of the rescaled null frame \texorpdfstring{$\lbrace \Lunit, \Annoy, X_1, X_2 \rbrace$}{}}

For convenience, we perform some of our computations relative to the rescaled null frame
$\lbrace \Lunit, \uLgood, X_1, X_2 \rbrace.$
In the following lemma, we provide the connection coefficients of this frame.

\begin{lemma}[\textbf{Connection coefficients of the rescaled null frame $\lbrace \Lunit, \uLgood, X_1, X_2 \rbrace$}]
	\label{L:CONNECTIONCOEFFICENTSOFNULLFRAME}
	We define the $S_{t,u}$ tensorfields 
	\begin{subequations}
	\begin{align}
	\upeta &:= \upzeta + \upmu^{-1} \angdiff \upmu,
		\label{E:UPETATENSOR} \\
	\underline{\upchi}
		& := 2 \upmu \angk - \upmu \upchi.
		\label{E:UCHI}
	\end{align}
	\end{subequations}
	
	Then the covariant derivatives of the rescaled null frame vectorfields can be expressed as follow:
	\begin{subequations}
	\begin{align}
		\D_{\Lunit} \Lunit & = \upmu^{-1}(\Lunit \upmu) \Lunit, 
			\\
		\D_{\uLgood} \Lunit & = - (\Lunit \upmu) \Lunit 
			+ 2 \upmu \upeta^A X_A, 
			\\
		\D_A \Lunit & = - \upzeta_A \Lunit
			+ \upchi_A^{\ B} X_B, 
			\\
		\D_{\Lunit} \uLgood & = - 2 \upmu \upzeta^A X_A, 
			\\
		\D_{\uLgood} \uLgood & = \lbrace \upmu^{-1} \uLgood \upmu 
			+ L\upmu \rbrace \uLgood
			- 2 \upmu (\angdiffuparg{A} \upmu) X_A, 
			\\
		\D_A \uLgood & = \upeta_A \uLgood
			+ \underline{\upchi}_A^{\ B} X_B, 
			\\
		\D_{\Lunit} X_A & = \D_A \Lunit, 
			\\
		\D_A X_B & = \angDarg{A} X_B
			+ \frac{1}{2} \upmu^{-1} \underline{\upchi}_{AB} \Lunit
			+ \frac{1}{2} \upmu^{-1} \upchi_{AB} \uLgood.
	\end{align}
	\end{subequations}
	
\end{lemma}

\begin{proof}
Lemma~\ref{L:CONNECTIONCOEFFICENTSOFNULLFRAME} follows from
the identity $\uLgood = \upmu \Lunit + 2 \Rad,$
Lemma~\ref{L:CONNECTIONLRADFRAME},
and straightforward computations.
\end{proof}

\section{Frame decomposition of the \texorpdfstring{$\upmu-$}{inverse-foliation-density-}weighted wave operator}
In the next proposition, we decompose the $\upmu-$weighted
wave operator $\upmu \square_{g(\Psi)}$ relative to the 
rescaled frame $\lbrace \Lunit, \Rad, X_1, X_2 \rbrace.$

\begin{proposition} [\textbf{Frame decomposition of $\upmu \square_{g(\Psi)} f$}]
	\label{P:GEOMETRICWAVEOPERATORFRAMEDECOMPOSED}
	Let $f$ be a function.
	Then relative to the frame $\lbrace \Lunit, \Rad, X_1, X_2 \rbrace,$
	$\upmu \square_{g(\Psi)} f$ can be expressed in either of the following two forms:
	\begin{subequations}
	\begin{align} \label{E:LONOUTSIDEGEOMETRICWAVEOPERATORFRAMEDECOMPOSED}
		\upmu \square_{g(\Psi)} f 
			& = - \Lunit(\upmu \Lunit f + 2 \Rad f)
				+ \upmu \angLap f
				\\
			& \ \ \underbrace{- \mytr \upchi \Rad f}_{\mbox{slowest decay}}
				- \mytr  \angkuparg{(Trans-\Psi)} \Lunit f
				- \upmu \mytr  \angkuparg{(Tan-\Psi)} \Lunit f
				\notag \\
			& \ \ 
				- 2 \upzeta^{(Trans-\Psi) \#} \cdot \angdiff f
				- 2 \upmu \upzeta^{(Tan-\Psi) \#} \cdot \angdiff f, 
				\notag \\
		\upmu \square_{g(\Psi)} f  
			& = - (\upmu \Lunit + 2 \Rad)(\Lunit f)
				+ \upmu \angLap f
				\label{E:LONINSIDEGEOMETRICWAVEOPERATORFRAMEDECOMPOSED} \\
			& \ \ \underbrace{- \mytr \upchi \Rad f}_{\mbox{slowest decay}}
				- \upomega^{(Trans-\Psi)} \Lunit f
				- \upmu \upomega^{(Tan-\Psi)} \Lunit f
					\notag \\
			& \ \ - \mytr  \angkuparg{(Trans-\Psi)} \Lunit f
				- \upmu \mytr  \angkuparg{(Tan-\Psi)} \Lunit f
				\notag \\
			& \ \ 
				+ 2 \upzeta^{(Trans-\Psi) \#} \cdot \angdiff f
				+ 2 \upmu \upzeta^{(Tan-\Psi) \#} \cdot \angdiff f
				+ 2 (\angdiffuparg{\#} \upmu) \cdot \angdiff f.
				\notag
	\end{align}
	\end{subequations}
	
	Furthermore, with $\uLgood = \upmu \Lunit + 2 \Rad,$ 
	we also have the following alternate decompositions:
	\begin{subequations}
	\begin{align} \label{E:WAVEEQUATIONSHOCKFORMATIONVERSION}
			\Lunit 
			\left\lbrace
				\uLgood
					(\rgeo f)
			\right\rbrace
			& = - \rgeo \upmu \square_{g(\Psi)} f 
				\\
			& \ \ + \rgeo \angLap f
				+ \rgeo (\upmu - 1) \angLap f
				\notag \\
			& \ \ + 2 (\upmu - 1) \Lunit f
				+ \upomega^{(Trans-\Psi)} f
				+ \upmu \upomega^{(Tan-\Psi)} f
				\notag \\
			& \ \ - \rgeo \mytr  \angkuparg{(Trans-\Psi)} \Lunit f
				- \rgeo \upmu \mytr  \angkuparg{(Tan-\Psi)} \Lunit f
				\notag \\
			& \ \ - \rgeo \mytr \upchi^{(Small)} \Rad f
				- 2 \rgeo \upzeta^{(Trans-\Psi) \#} \cdot \angdiff f
				- 2 \rgeo \upmu \upzeta^{(Tan-\Psi) \#} \cdot \angdiff f,
				\notag
		\end{align}
		\begin{align}	 \label{E:WAVEOPDECOMPRADRGEOLPSIPLUSHALFTRCHIPSI}
				\Rad 
				\left\lbrace
					\rgeo 
					\left(
						\Lunit f
						+ \frac{1}{2} \mytr \upchi f
					\right)
				\right\rbrace 
			& = - \frac{1}{2} \upmu \Lunit(\rgeo \Lunit f)
				+ \frac{1}{2} \rgeo \upmu \angLap f
				\\
			& \ \ - \frac{1}{2} \rgeo \upmu \square_{g(\Psi)} f  
					\notag \\
			& \ \ - \frac{1}{2} \rgeo \upomega^{(Trans-\Psi)} \Lunit f
				- \frac{1}{2} \rgeo \upmu \upomega^{(Tan-\Psi)} \Lunit f
					\notag \\
			& \ \ - \frac{1}{2} \rgeo \mytr  \angkuparg{(Trans-\Psi)} \Lunit f
				- \frac{1}{2} \rgeo \upmu \mytr  \angkuparg{(Tan-\Psi)} \Lunit f
				\notag \\
			& \ \ 
				+ \rgeo \upzeta^{(Trans-\Psi) \#} \cdot \angdiff f
				+ \rgeo \upmu \upzeta^{(Tan-\Psi) \#} \cdot \angdiff f
				+ \rgeo (\angdiffuparg{\#} \upmu) \cdot \angdiff f
				\notag \\
			& \ \ 
				+ \frac{1}{2} \upmu \Lunit f
				- \Lunit f
				\notag \\
			& \ \ + \frac{1}{2} \rgeo (\Rad \mytr \upchi^{(Small)}) f
				- \frac{1}{2} \mytr \upchi^{(Small)} f.
				\notag	
			\end{align}		
			\end{subequations}				
	In the above expressions, the $S_{t,u}$ tensorfields
	$\upchi,$ 
	$\upchi^{(Small)},$
	$\upomega^{(Trans-\Psi)},$
	$\upomega^{(Tan-\Psi)},$
	$\upzeta^{(Trans-\Psi)},$
	$\angkuparg{(Trans-\Psi)},$ 
	$\upzeta^{(Tan-\Psi)},$
	and
	$\angkuparg{(Tan-\Psi)}$
	are defined by
	\eqref{E:CHIDEF},
	\eqref{E:CHIJUNKDEF},
	\eqref{E:LITTLEOMEGATRANS},
	\eqref{E:LITTLEOMEGATAN},	
	\eqref{E:ZETATRANSVERSAL},
	\eqref{E:KABTRANSVERSAL},
	\eqref{E:ZETAGOOD},
	and \eqref{E:KABGOOD}.
\end{proposition}

\begin{proof}
	Using \eqref{E:GINVERSENULLFRAME}, 
	Lemmas \ref{L:CONNECTIONLRADFRAME}
	and \ref{L:CONNECTIONCOEFFICENTSOFNULLFRAME},
	and the Leibniz rule identities
	$X_A(X_B f) 
	= (\D_A X_B) f 
		+
		\D_{AB}^2 f 
	= (\angDarg{A} X_B) f 
		+ 
		\angDsquaredarg{A}{B} f,$
	we compute that
	\begin{align} 
		\upmu \square_{g(\Psi)} f
		& = \upmu (g^{-1})^{\alpha \beta} \D_{\alpha \beta}^2 f
			= - \D_{\Lunit \uLgood}^2 f
			+ \upmu (\ginversesphere)^{AB} \D_{AB}^2 f
				\label{E:FIRSTCOMPUTATIONOFFRAMEDECOMPOFGEOMETRICWAVEOP} \\
		& = - \Lunit(\uLgood f)
			+ (\D_{\Lunit} \uLgood) f 
			+ \upmu (\ginversesphere)^{AB} \angDsquaredarg{A}{B} f
			- \upmu (\ginversesphere)^{AB} \angkdoublearg{A}{B} \Lunit f
			- (\ginversesphere)^{AB} \upchi_{AB} \Rad f
			\notag
				\\
		& = - \Lunit(\uLgood f)
			- 2 \upmu \upzeta^{\#} \cdot \angdiff f
			+ \upmu \angLap f
			- \upmu \mytr  \angk \Lunit f
			- \mytr \upchi \Rad f
				\notag \\
		& = - \Lunit(\uLgood f)
			+ \upmu \angLap f
			- \mytr \upchi \Rad f
			- \mytr  \angkuparg{(Trans-\Psi)} \Lunit f
			- \upmu \mytr  \angkuparg{(Tan-\Psi)} \Lunit f
			\notag	\\
		& \ \ 
			- 2 \upzeta^{(Trans-\Psi) \#} \cdot \angdiff f
			- 2 \upmu \upzeta^{(Tan-\Psi) \#} \cdot \angdiff f.
			\notag
		\end{align}
		From \eqref{E:FIRSTCOMPUTATIONOFFRAMEDECOMPOFGEOMETRICWAVEOP}
		and the identity $\uLgood = \upmu \Lunit + 2 \Rad,$
		we deduce \eqref{E:LONOUTSIDEGEOMETRICWAVEOPERATORFRAMEDECOMPOSED}.
		The proof of \eqref{E:LONINSIDEGEOMETRICWAVEOPERATORFRAMEDECOMPOSED} is similar;
		we simply interchange the order of $\Lunit$ and $\uLgood$ in
		\eqref{E:FIRSTCOMPUTATIONOFFRAMEDECOMPOFGEOMETRICWAVEOP}.
		
		To prove \eqref{E:WAVEEQUATIONSHOCKFORMATIONVERSION},
		we multiply both sides of \eqref{E:LONOUTSIDEGEOMETRICWAVEOPERATORFRAMEDECOMPOSED}
		by $\rgeo,$ use the identities $\Lunit \rgeo = 1$ and $\Rad \rgeo = - 1,$
		use Lemma~\ref{L:UPMUFIRSTTRANSPORT} to substitute for $\Lunit \upmu,$
		and carry out straightforward computations.
		To prove \eqref{E:WAVEOPDECOMPRADRGEOLPSIPLUSHALFTRCHIPSI},
		we multiply both sides of \eqref{E:LONINSIDEGEOMETRICWAVEOPERATORFRAMEDECOMPOSED}
		by $\rgeo$ and use a similar argument together with 
		the decomposition $\mytr \upchi = 2 \rgeo^{-1} + \mytr \upchi^{(Small)}.$
\end{proof}


\chapter{Construction of the Rotation Vectorfields and Their Basic Properties}
\label{C:ROTATIONS}
\thispagestyle{fancy}
In Chapter~\ref{C:ROTATIONS}, we construct the $S_{t,u}-$tangent rotation vectorfields 
$\Rot$ and derive some of their basic properties. 
Our construction is not difficult: we simply project away the dangerous $\Radunit$ component
present in the Euclidean rotations. If we did not remove this dangerous component,
then the rotational derivatives of $\Psi$ could blow-up like $\upmu^{-1}$ 
at the shock-formation 
points. Hence, they would be useless for detecting dispersive behavior
that persists up until the shock.

\section{Construction of the rotation vectorfields}
\label{S:CONSTRUCTIONOFROTATIONS}
To begin the detailed construction, we
first recall that the type $\binom{1}{1}$ $S_{t,u}$ tensorfield $\sphereproject,$ 
which projects onto the spheres $S_{t,u}$ (see Def.~\ref{D:PROJECTIONS}), 
has the following non-zero rectangular components (see \eqref{E:ALTERNATEPROJECTIONEXPRESSIONS}),
$(i,j=1,2,3):$
\begin{align} \label{E:SPHEREPROJECTSPATIALRECTANGULAR}
	\sphereproject_j^{\ i} := (\gtinverse)^{ia} \gsphere_{aj}
		= \delta_j^{\ i} 
			- g_{ja} \Radunit^a \Radunit^i.
\end{align}

\begin{definition}[\textbf{Rotation vectorfields}]
\label{D:ROTATION}
For $l = 1,2,3,$ we define the \textbf{Euclidean} rotation $\Roteucarg{l}$ to be
the $\Sigma_t-$tangent vectorfield with the following rectangular spatial components,
$(i=1,2,3):$
\begin{align} \label{E:EULCIDEANROTATION}
	\Roteucarg{l}^i := \epsilon_{lai} x^a,
\end{align}
where $\epsilon_{ijk}$ is the fully antisymmetric symbol
normalized by $\epsilon_{123} = 1.$

We define $\Rot_{(l)}$ to be the $S_{t,u}-$tangent vectorfield
with the following rectangular spatial components,
$(i=1,2,3):$
\begin{align} \label{E:GEOMETRICROTATION}
	\Rot_{(l)}^i := \sphereproject_a^{\ i} \Roteucarg{l}^a.
\end{align}
\end{definition}

\section{Basic properties of the rotation vectorfields}
\label{S:BASICPROPSOFROTATIONS}
Some of our most delicate analysis involves the components of the Euclidean
rotations in the direction of $\Radunit.$ This component is captured in the next definition.
\begin{definition}[\textbf{The $\Radunit$ component of $\Roteucarg{l}$}]
	\label{D:RADUNITCOMPONENTOFROT}
		We define the following scalar-valued functions $\RotRadcomponent{l},$ 
		$(l = 1,2,3):$
		\begin{align} \label{E:RADUNITCOMPONENTOFROT}
			\RotRadcomponent{l} := g(\Roteucarg{l}, \Radunit).
		\end{align}
\end{definition}

In the next lemma, we decompose $\Roteucarg{l}$ into $\Rot_{(l)}$
plus an error term in the radial direction $\Radunit.$ We also provide an 
explicit expression for the radial component, which we denote by $\RotRadcomponent{l}.$

\begin{lemma}[\textbf{Decomposition $\Rot_{(l)}$ into $\Roteucarg{l}$ plus an error, and an expression for $\RotRadcomponent{l}$}]
	\label{L:ROTATIONDECOMPOSITIONINTOEUCLIDEANPLUSRADCOMPONENT}
The following identities hold, where
$g_{ij}^{(Small)}$ and $\Radunit_{(Small)}^i$
are defined by \eqref{E:LITTLEGDECOMPOSED} and \eqref{E:RADUNITJUNK}:
\begin{align}
	\Roteucarg{l} 
	& = \Rot_{(l)} + \RotRadcomponent{l} \Radunit,
		\label{E:ROTATIONDECOMPOSITIONINTOEUCLIDEANPLUSRADCOMPONENT} \\
	\RotRadcomponent{l} 
	& = g_{bc} \epsilon_{lab} x^a \Radunit^c
		= \epsilon_{lab} x^a \Radunit_{(Small)}^b
		+ \epsilon_{lab} g_{bc}^{(Small)} x^a \Radunit_{(Small)}^c
		- \frac{1}{\rgeo} g_{bc}^{(Small)} \epsilon_{lab} x^a x^c.
		\label{E:EUCLIDEANROTATIONRADCOMPONENT} 
\end{align}
\end{lemma}
\begin{proof}
The identity \eqref{E:ROTATIONDECOMPOSITIONINTOEUCLIDEANPLUSRADCOMPONENT}
follows easily 
from definition \eqref{E:GEOMETRICROTATION}, 
definition \eqref{E:RADUNITCOMPONENTOFROT}, 
and the fact that $g(\Radunit,\Radunit)= 1.$

The first equality in \eqref{E:EUCLIDEANROTATIONRADCOMPONENT} follows directly from
expressing the right-hand side of \eqref{E:RADUNITCOMPONENTOFROT}
in rectangular spatial coordinates and definition \eqref{E:EULCIDEANROTATION}.
To deduce the second equality, we insert the decompositions
$g_{ab} = \delta_{ab} + g_{ab}^{(Small)}$
and $\Radunit^c = - \rgeo^{-1} x^c + \Radunit_{(Small)}^c$
into the first equality. By the anti-symmetry of $\epsilon_{\cdots},$
the ``large'' term $\rgeo^{-1} \delta_{bc} \epsilon_{lab} x^a x^c$ vanishes
and hence the second equality follows.
\end{proof}

\begin{lemma}[\textbf{The rectangular components of $\Rot_{(l)}$}]
	\label{L:ROTATIONRECATNGULARCOMPONENT}
	The rectangular spatial components $\Rot_{(l)}^i,$ 
	$(i=1,2,3),$
	can be expressed as follows:
	\begin{align} \label{E:ROTATIONRECATNGULARCOMPONENT}
	\Rot_{(l)}^i
	& = \epsilon_{lai} x^a
			- \RotRadcomponent{l} \Radunit^j 
		= \epsilon_{lai} x^a
			+ \RotRadcomponent{l} \frac{x^i}{\rgeo} 
			- \RotRadcomponent{l} \Radunit_{(Small)}^i.
	\end{align}
\end{lemma}
\begin{proof}
	Lemma~\ref{L:ROTATIONRECATNGULARCOMPONENT} follows easily from
	definition \eqref{E:GEOMETRICROTATION}, 
	the identity
	\eqref{E:ROTATIONDECOMPOSITIONINTOEUCLIDEANPLUSRADCOMPONENT},
	and the decomposition \eqref{E:RADUNITJUNK}.
\end{proof}

\begin{lemma}[\textbf{An expression for the} $S_{t,u}$ \textbf{components of} $\Rot_{(l)}$]
	\label{L:ROTATATIONACOMPONENTALTERNATEEXPRESSION}
	The following identities hold,
	where $g_{ij}^{(Small)}$
	is defined by \eqref{E:LITTLEGDECOMPOSED},
	$(A=1,2$ and $\angdiffuparg{A} x^b = (\ginversesphere)^{AB} \angdiffarg{B} x^b):$
	\begin{align} \label{E:ROTATATIONACOMPONENTALTERNATEEXPRESSION}
		\Rot_{(l)}^A
		& = \epsilon_{lca} g_{ab} x^c \angdiffuparg{A} x^b
			=
			\epsilon_{lab} x^a \angdiffuparg{A} x^b	
			+ \epsilon_{lca} g_{ab}^{(Small)} x^c \angdiffuparg{A} x^b.
	\end{align}
\end{lemma}
\begin{proof}
	First, using 
	definition \eqref{E:EULCIDEANROTATION},
	\eqref{E:ROTATIONDECOMPOSITIONINTOEUCLIDEANPLUSRADCOMPONENT},
	and the identities 
	$g(\Radunit, X_B) = 0$
	and
	$X_B^i = \angdiffarg{B} x^i,$
	we compute that
	$\Rot_{(l)B} 
		:= \gsphere(\Rot_{(l)}, X_B) 
		= g(\Roteucarg{l}, X_B) 
		= g_{ab} \epsilon_{lca} x^c X_B^b
		= \epsilon_{lca} g_{ab} x^c \angdiffarg{B} x^b.
	$
	The first equality in \eqref{E:ROTATATIONACOMPONENTALTERNATEEXPRESSION}
	now follows from contracting the previous identity against $(\ginversesphere)^{AB}.$
	The second then follows from the decomposition
	$g_{ab} = \delta_{ab} + g_{ab}^{(Small)}.$
\end{proof}

In the next lemma, we derive an expression 
for the Lie derivatives $\Lie_{\Rot_{(l)}} \Rot_{(m)}.$
The most important aspect of the lemma is
that \emph{the expression does not involve the dangerous transversal derivative} 
$\Radunit \Psi = \upmu^{-1} \Rad \Psi.$

	\begin{lemma}[\textbf{Expression for the commutators $[\Rot_{(m)}, \Rot_{(n)}]$}]
	\label{L:ROTATATIONCOMMUTATORS}
	Let $\Theta^{(Tan-\Psi)}$ be the type $\binom{0}{2}$ $S_{t,u}-$tangent tensorfield 
	\eqref{E:BIGTHETAGOOD}, and let $\RotRadcomponent{l}$ $(l = 1,2,3),$ 
	be the scalar-valued functions defined in \eqref{E:RADUNITCOMPONENTOFROT}. 
	Then the following identity holds:
	\begin{align} \label{E:ROTATATIONCOMMUTATORS}
		\Lie_{\Rot_{(l)}} \Rot_{(m)}
		= [\Rot_{(l)}, \Rot_{(m)}] 
		= - \epsilon_{lmn} \Rot_{(n)}
			+ \rotationcommutationerrorupsharp{W}{l}{m},
	\end{align}
	where $\rotationcommutationerror{W}{l}{m}$ is the $S_{t,u}$ one-form given by
	\begin{align}
	\rotationcommutationerrorarg{W}{A}{l}{m}
	& = - \RotRadcomponent{l} \Rot_{(m)}^B \upchi_{BA}^{(Small)} 
			+ \RotRadcomponent{m} \Rot_{(l)}^B \upchi_{BA}^{(Small)}
			\label{E:ROTATIONCOMMUTATORSERRORONEFORM} \\
		& \ \ 
			+ \RotRadcomponent{l} \Rot_{(m)}^B \Theta_{BA}^{(Tan-\Psi)}
			- \RotRadcomponent{m} \Rot_{(l)}^B \Theta_{BA}^{(Tan-\Psi)}
				\notag \\
		& \ \ - \RotRadcomponent{l} \epsilon_{mab} \Radunit_{(Small)}^a g_{bc} \angdiffarg{A} x^c
			+ \RotRadcomponent{m} \epsilon_{lab} \Radunit_{(Small)}^a g_{bc} \angdiffarg{A} x^c.
			\notag
	\end{align}
\end{lemma}	
\begin{proof}
	Using \eqref{E:ROTATIONDECOMPOSITIONINTOEUCLIDEANPLUSRADCOMPONENT},
	the first equation \eqref{E:ROTATIONRECATNGULARCOMPONENT}, 
	and the decomposition \eqref{E:RADUNITJUNK},
	we compute that
\begin{align} \label{E:ROTLITTLELAPPLIEDTOROTLITTLEMIFIRSTCOMPUTATION}
\Rot_{(l)} \Rot_{(m)}^i 
	& = \Rot_{(l)}^c \epsilon_{mci} 
		- \RotRadcomponent{m} \Rot_{(l)} \Radunit^i
		- (\Rot_{(l)} \RotRadcomponent{m}) \Radunit^i 
			\\
	& = \Roteucarg{l}^c \epsilon_{mci} 
		+ \rgeo^{-1} \RotRadcomponent{l} x^c \epsilon_{mci}
		- \RotRadcomponent{l} \Radunit_{(Small)}^c \epsilon_{mci}
		+ \rgeo^{-1} \RotRadcomponent{m} \epsilon_{lci} x^c 
		+ \rgeo^{-1} \RotRadcomponent{l} \RotRadcomponent{m} \Radunit^i
			\notag \\
	& \ \ - \RotRadcomponent{m} \Rot_{(l)} \Radunit_{(Small)}^i
		- (\Rot_{(l)} \RotRadcomponent{m}) \Radunit^i.
		\notag
\end{align}
Using equation \eqref{E:ANGDIFFRADJUNKI} to substitute 
for $\Rot_{(l)} \Radunit_{(Small)}^i = \Rot_{(l)}^A \angdiffarg{A} \Radunit_{(Small)}^i$
on the right-hand side of \eqref{E:ROTLITTLELAPPLIEDTOROTLITTLEMIFIRSTCOMPUTATION}, we deduce that
\begin{align} \label{E:ROTLITTLELAPPLIEDTOROTLITTLEMITHCOMPONENT}
	\Rot_{(l)} \Rot_{(m)}^i
	& = \Roteucarg{l}^c \epsilon_{mci}
		- \RotRadcomponent{l} \Radunit_{(Small)}^c \epsilon_{mci}
		+ \rgeo^{-1} \RotRadcomponent{l} \epsilon_{mci} x^c 
		+ \rgeo^{-1} \RotRadcomponent{m} \epsilon_{lci} x^c
			\\
	& \ \
		+ \RotRadcomponent{m} (\ginversesphere)^{BC} \upchi_{AB}^{(Small)} \Rot_{(l)}^A \angdiffarg{C} x^i
		- \RotRadcomponent{m} (\ginversesphere)^{BC} \Theta_{AB}^{(Tan-\Psi)} \Rot_{(l)}^A \angdiffarg{C} x^i
		\notag	\\
	& \ \ 
		+ \rgeo^{-1} \RotRadcomponent{l} \RotRadcomponent{m} \Radunit^i
		- \RotRadcomponent{m} \theta_A^{(Tan-\Psi)} \Rot_{(l)}^A \Radunit^i
		- (\Rot_{(l)} \RotRadcomponent{m}) \Radunit^i.
		\notag
\end{align}		
We now compute
$[\Lie_{\Rot_{(l)}} \Rot_{(m)}]^i = \Rot_{(l)} \Rot_{(m)}^i - \Rot_{(m)} \Rot_{(l)}^i$
by interchanging the roles of $l$ and $m$ in \eqref{E:ROTLITTLELAPPLIEDTOROTLITTLEMITHCOMPONENT}
and then subtracting the two identities. It is straightforward to compute that
the first term on the right-hand side of \eqref{E:ROTLITTLELAPPLIEDTOROTLITTLEMITHCOMPONENT}
generates the term $[\Roteucarg{l},\Roteucarg{m}]^i
= - \epsilon_{lmn} \Roteucarg{n}^i = - \epsilon_{lmn} \Rot_{(n)}^i + \uptau_{(n)} \epsilon_{lmn} \Radunit^i.$
Moreover, since $[\Rot_{(l)}, \Rot_{(m)}]$ is $S_{t,u}-$tangent, we can then apply the 
$S_{t,u}$ projection tensorfield $\sphereproject_i^{\ j}$ \eqref{E:SPHEREPROJECTSPATIALRECTANGULAR}
to both sides of the resulting expression, 
which preserves the vector $[\Rot_{(l)}, \Rot_{(m)}]$ on the left-hand side
and annihilates all terms that are proportional to the vector $\Radunit.$
In total, we deduce from this line of reasoning that $(j=1,2,3)$
\begin{align} \label{E:ROTLITTLELLITTLEMCOMMUTATORFIRSTEXPRESSION}
	[\Rot_{(l)}, \Rot_{(m)}]^j
	& = - \epsilon_{lmn} \Rot_{(n)}^j
		- \RotRadcomponent{l} \Radunit_{(Small)}^c \sphereproject_i^{\ j} \epsilon_{mci}
		+ \RotRadcomponent{m} \Radunit_{(Small)}^c \sphereproject_i^{\ j} \epsilon_{lci}
		\\
	& \ \
		+ \RotRadcomponent{m} (\ginversesphere)^{BC} \upchi_{AB}^{(Small)} \Rot_{(l)}^A \angdiffarg{C} x^j
		- \RotRadcomponent{l} (\ginversesphere)^{BC} \upchi_{AB}^{(Small)} \Rot_{(m)}^A \angdiffarg{C} x^j
			\notag \\
	& \ \
		- \RotRadcomponent{m} (\ginversesphere)^{BC} \Theta_{AB}^{(Tan-\Psi)} \Rot_{(l)}^A \angdiffarg{C} x^j
		+ \RotRadcomponent{l} (\ginversesphere)^{BC} \Theta_{AB}^{(Tan-\Psi)} \Rot_{(m)}^A \angdiffarg{C} x^j.
		\notag
\end{align}
Clearly, the first term on the right-hand side of
\eqref{E:ROTLITTLELLITTLEMCOMMUTATORFIRSTEXPRESSION}
is the first term on the right-hand side of \eqref{E:ROTATIONCOMMUTATORSERRORONEFORM} as desired.
To show that the remaining terms on the right-hand side of \eqref{E:ROTLITTLELLITTLEMCOMMUTATORFIRSTEXPRESSION}
are (when viewed as $S_{t,u}-$tangent vectors with rectangular components $(\cdot)^j$)
$\gsphere-$dual to $\rotationcommutationerror{W}{l}{m},$ 
we contract them against
$g_{jk} \angdiffarg{D} x^k$ 
and use the identities
$g_{jk} (\angdiffarg{D} x^k) \sphereproject_i^{\ j} = g_{ic} \angdiffarg{D} x^c$
and
$g_{jk} (\angdiffarg{C} x^j) \angdiffarg{D} x^k = \gsphere_{CD}.$

\end{proof}


\chapter{Definition of the Commutation Vectorfields and Deformation Tensor Calculations}
\label{C:DEFORMATIONTENSORCALCULATIONS}
\thispagestyle{fancy}
In Chapter~\ref{C:DEFORMATIONTENSORCALCULATIONS}, we define the set $\mathscr{Z}$ of commutation vectorfields
that we use to commute the wave equation $\upmu \square_{g(\Psi)} \Psi = 0.$
We also compute the components of the deformation tensors $\deform{Z}$ of 
the vectorfields $Z \in \mathscr{Z}$ relative to the rescaled frame $\lbrace \Lunit, \Rad, X_1, X_2 \rbrace.$ 
In order to prove our sharp classical lifespan theorem, 
we must precisely understand the structure of some of these components. 
The reason is that the derivatives 
of the $\deform{Z}$ appear as inhomogeneous terms in the commuted wave equation
(see Lemma~\ref{L:WAVEONCECOMMUTEDBASICSTRUCTURE} and 
Prop.~\ref{P:COMMUTATIONCURRENTDIVERGENCEFRAMEDECOMP}),
and some of the corresponding frame components are difficult to analyze.
In fact, some of them seem to be on the border of the kinds of terms
that would allow the proof of our sharp classical lifespan theorem to go through.

\section{The commutation vectorfields}
To prove our sharp classical lifespan theorem, we commute the 
covariant wave equation $\upmu \square_{g(\Psi)} \Psi = 0$
with the vectorfields belonging to the following set $\mathscr{Z}.$

\begin{definition}[\textbf{Set of commutation vectorfields}]
\label{D:DEFSETOFCOMMUTATORVECTORFIELDS}
	We define the set $\mathscr{Z}$ of commutation vectorfields as follows:
\begin{align}  \label{E:DEFSETOFCOMMUTATORVECTORFIELDS}
		\mathscr{Z} := \lbrace \rgeo \Lunit, \Rad, \Rot_{(1)}, \Rot_{(2)}, \Rot_{(3)} \rbrace.
\end{align}

\end{definition}

In our analysis, we also use the following subsets of $\mathscr{Z}.$

\begin{definition}[\textbf{Spatial and rotation commutation subsets}]
\label{D:DEFSETOFSPATIALCOMMUTATORVECTORFIELDS}
We define the subsets $\mathscr{S}, \mathscr{O} \subset \mathscr{Z}$ of spatial commutation vectorfields
and rotation commutation vectorfields as follows:
\begin{subequations}
\begin{align}  \label{E:DEFSETOFSPATIALCOMMUTATORVECTORFIELDS}
		\mathscr{S} 
		& := \lbrace \Rad, \Rot_{(1)}, \Rot_{(2)}, \Rot_{(3)} \rbrace,
			\\
		\mathscr{O} 
		& := \lbrace \Rot_{(1)}, \Rot_{(2)}, \Rot_{(3)} \rbrace.
		\label{E:DEFSETOFROTATIONCOMMUTATORVECTORFIELDS}
\end{align}
\end{subequations}
\end{definition}

\section{Deformation tensor calculations}
We recall that the deformation tensor
$\deform{V}$ of a vectorfield $V$ is defined in 
Def.~\ref{D:DEFORMTENSDEFINED}.
In this section, we calculate 
$\deform{V}$
for various vectorfields $V.$
We begin with the next lemma, which shows 
that for some important vectorfields $V,$ 
the $S_{t,u}$ projection of $\deform{V}$
is the same as $\angLie_V \gsphere.$

\begin{lemma}[\textbf{Connection between projected lie derivatives of} $\gsphere$ \textbf{and} $\angpi$]
	\label{L:CONNECTIONBETWEENANGLIEOFGSPHEREANDDEFORMATIONTENSORS}
	If $V \in \lbrace \Lunit, \rgeo \Lunit, \Rad \rbrace$
	or $V$ is an $S_{t,u}-$tangent vectorfield, then 
	\begin{subequations}
	\begin{align} \label{E:CONNECTIONBETWEENANGLIEOFGSPHEREANDDEFORMATIONTENSORS}
		\angLie_V \gsphere_{AB} 
		& = \angdeformarg{Z}{A}{B},
			\\
		(\angLie_V \ginversesphere)^{AB} 
		& = - \angdeformuparg{Z}{A}{B}.
			\label{E:CONNECTIONBETWEENANGLIEOFGINVERSESPHEREANDDEFORMATIONTENSORS}
	\end{align}
	\end{subequations}
\end{lemma}

\begin{proof}
The identity \eqref{E:CONNECTIONBETWEENANGLIEOFGSPHEREANDDEFORMATIONTENSORS} follows from
Cor.~\ref{C:ANGLIEVFORGOODVOFTENSOREQUALSANGLIEOFSTUPROJECTEDTENSOR} and
Def.~\ref{D:DEFORMTENSDEFINED}. The identity \eqref{E:CONNECTIONBETWEENANGLIEOFGINVERSESPHEREANDDEFORMATIONTENSORS}
then follows from \eqref{E:CONNECTIONBETWEENANGLIEOFGSPHEREANDDEFORMATIONTENSORS}
and the identity 
$(\angLie_V \ginversesphere)^{AB} = - (\ginversesphere)^{AC} (\ginversesphere)^{BD} \Lie_V \gsphere_{CD}.$
\end{proof}

We now calculate the frame components of the 
deformation tensors $\deformarg{V}{\mu}{\nu}$
of various vectorfields including the commutation vectorfields $\mathscr{Z}.$
The main result is the following proposition.

\begin{remark}[\textbf{Some important structural features}] \label{R:DEFTENSORCALCULATIONSIMPORTANTASPECTS}
	The following three aspects of the proposition are highly important.
	\begin{itemize}
		\item The dangerous transversal derivative $\Radunit \Psi = \upmu^{-1} \Rad \Psi$
		is completely absent from the right-hand sides of the expressions.
		\item $\deformarg{Z}{\Lunit}{\Lunit} = 0$ for all $Z \in \mathscr{Z}.$  
			This identity in particular implies the absence of the dangerous
			quadratic term $\upmu^{-1} \deformarg{Z}{\Lunit}{\Lunit} \Rad \Psi,$
			whose derivatives would appear
			on the right-hand side of the identity \eqref{E:DIVCOMMUTATIONCURRENTDECOMPOSITION}.
			The presence of such a term could in principle destroy our proof of 
			sharp classical lifespan theorem
			because we would have no obvious way to control the dangerous factor $\upmu^{-1}.$
		\item $\upmu^{-1} 
											\left\lbrace 
												\deformarg{Z}{\Lunit}{\Rad}
												+ Z \upmu 
											\right\rbrace$
			is either $0$ or $-1$ for all $Z \in \mathscr{Z}.$
			We elaborate upon the importance of this fact in
			Remark~\ref{R:POTENTIALLYDANGEROUSTERMISNOTDANGEROUS}.
	\end{itemize}
\end{remark}

\begin{proposition}[\textbf{Expressions for the frame components of various deformation tensors}]
\label{P:DEFORMATIONTENSORFRAMECOMPONENTS}
We have the following identities,
where we are using the notation of 
Sects.~\ref{S:TRACEANDTRACEFREEPARTS}
and \ref{S:SCHEMATIC}
and the quantities $\smoothfunction$ are smooth functions of $\Psi:$
\begin{subequations}
\begin{align}
	\deformarg{\Rad}{\Lunit}{\Lunit} & = 0, 
		\\
	\deformarg{\Rad}{\Rad}{\Radunit} 
	& = 2 \Rad \upmu,
		\label{E:RADDEFORMRADRADUNIT} \\
	\deformarg{\Rad}{\Lunit}{\Rad} 
	& = - \Rad \upmu,
		\label{E:RADDEFORMLRADUNIT} \\
	\angdeformarg{\Rad}{\Lunit}{A}
	& = - \angdiffarg{A} \upmu 
		+ G_{(Frame)}
			\threemyarray
				[\upmu \Lunit \Psi]
				{\Rad \Psi}
				{\upmu \angdiff \Psi},
		\label{E:RADDEFORMLA} \\
	\angdeformarg{\Rad}{\Rad}{A} 
	& = 0,
		\label{E:RADDEFORMRADA} \\
	\angdeformfreearg{\Rad}{A}{B}
	& = - 2 \upmu \hat{\upchi}_{AB}^{(Small)}
		+ G_{(Frame)}
			\hat{\otimes}
			\threemyarray
				[\upmu \Lunit \Psi]
				{\Rad \Psi}
				{\upmu \angdiff \Psi},
		\label{E:RADDEFORMTRFREEANG} \\
	\mytr \angdeform{\Rad} 
	& = - \frac{4}{\rgeo} \upmu 
						- 2 \upmu \mytr \upchi^{(Small)}
						+ 2 \mytr  \angkuparg{(Trans-\Psi)}
						+ 2 \upmu \mytr  \angkuparg{(Tan-\Psi)}
				\label{E:RADDEFORMANG}
				\\
	& \ \ = - \frac{4}{\rgeo} \upmu
			- 2 \upmu \mytr \upchi^{(Small)}
			+ G_{(Frame)}^{\#}
			\threemyarray
				[\upmu \Lunit \Psi]
				{\Rad \Psi}
				{\upmu \angdiff \Psi},
				\notag
\end{align}
\end{subequations}

\begin{subequations}
\begin{align}	
	\deformarg{\Lunit}{\Lunit}{\Lunit}
	& = 0, 
		\\
	\deformarg{\Lunit}{\Rad}{\Radunit} 
	& = 2 \Lunit \upmu
		= G_{(Frame)}
			\myarray
				[\upmu \Lunit \Psi]
				{\Rad \Psi},
		\\
	\deformarg{\Lunit}{\Lunit}{\Rad} 
	& = - \Lunit \upmu 
		= G_{(Frame)}
			\myarray
				[\upmu \Lunit \Psi]
				{\Rad \Psi},
		\\
	\angdeformarg{\Lunit}{\Lunit}{A}
	& = 0, 
		\\
	\angdeformarg{\Lunit}{\Rad}{A} 
	& = \angdiffarg{A} \upmu
	 		+ G_{(Frame)}
			\threemyarray
				[\upmu \Lunit \Psi]
				{\Rad \Psi}
				{\upmu \angdiff \Psi}, \\
	\angdeformfreearg{\Lunit}{A}{B} 
	& = 2 \hat{\upchi}_{AB}^{(Small)},
	 \label{E:LDEFORMTRFREEANG}	\\
	\mytr  \angdeform{\Lunit}
	& = 2 \mytr \upchi^{(Small)}
		+ \frac{4}{\rgeo},
		\label{E:LDEFORMTRANG}
\end{align}	
\end{subequations}

\begin{subequations}
\begin{align}	
	\deformarg{\rgeo \Lunit}{\Lunit}{\Lunit}
	& = 0, 
		\label{E:RGEOLDFORMLL} \\
	\deformarg{\rgeo \Lunit}{\Rad}{\Radunit} 
	& = 2 \rgeo \Lunit \upmu + 2 \upmu
	= 	\rgeo
			G_{(Frame)}
			\myarray
				[\upmu \Lunit \Psi]
				{\Rad \Psi}
			+ 2 \upmu,
		\label{E:RGEOLDEFORMRADRADUNIT} \\
	\deformarg{\rgeo \Lunit}{\Lunit}{\Rad} 
	& = - \rgeo \Lunit \upmu - \upmu
		=  	\rgeo
				G_{(Frame)}
				\myarray
				[\upmu \Lunit \Psi]
				{\Rad \Psi}
			- \upmu,
		\label{E:RGEOLDEFORMLRAD} \\
	\angdeformarg{\rgeo \Lunit}{\Lunit}{A}
	& = 0, 
		\label{E:RGEOLDEFORMLA} \\
	\angdeformarg{\rgeo \Lunit}{\Rad}{A} 
	& = \rgeo \angdiffarg{A} \upmu
	 	+ \rgeo
	 		G_{(Frame)}
			\threemyarray
				[\upmu \Lunit \Psi]
				{\Rad \Psi}
				{\upmu \angdiff \Psi}, 
			\label{E:RGEOLDEFORMRADA} \\
	\angdeformfreearg{\rgeo \Lunit}{A}{B} 
	& = 2 \rgeo \hat{\upchi}_{AB}^{(Small)},
		\label{E:RGEOLDEFORMTRFREESPHERE} \\
	\mytr  \angdeform{\rgeo \Lunit} 
	& = 2 \rgeo \mytr \upchi^{(Small)}
		+ 4,
		\label{E:RGEOLDEFORMTRSPHERE} 
\end{align}	
\end{subequations}

\begin{subequations}
\begin{align}
	\deformarg{\Rot_{(l)}}{\Lunit}{\Lunit} 
	& = 0, 
		\label{E:ROTDEFORMLL} \\
	\deformarg{\Rot_{(l)}}{\Rad}{\Radunit} 
	& = 2 \Rot_{(l)} \upmu, 
		\label{E:ROTDEFORMRADRADUNIT} \\
	\deformarg{\Rot_{(l)}}{\Lunit}{\Rad} 
	& = - \Rot_{(l)} \upmu, 
		\label{E:ROTDEFORMLRAD} \\
	\angdeformarg{\Rot_{(l)}}{\Lunit}{A} 
	& = - \upchi_{AB}^{(Small)} \Rot_{(l)}^B 
		+ \angdeformarg{\Rot_{(l)};Error}{\Lunit}{A},
			\label{E:ROTDEFORMLSPHERE} \\
	\angdeformarg{\Rot_{(l)}}{\Rad}{A} 
	& = \upmu \upchi_{AB}^{(Small)} \Rot_{(l)}^B 
				+ \RotRadcomponent{l} \angdiffarg{A} \upmu
				+ \angdeformarg{\Rot_{(l)};Error}{\Rad}{A},
				\label{E:ROTDEFORMRADSPHERE} \\
	\angdeformfreearg{\Rot_{(l)}}{A}{B}
	& = 2 \RotRadcomponent{l} \hat{\upchi}_{AB}^{(Small)}
		+ \angdeformfreearg{\Rot_{(l)};Error}{A}{B},
		\label{E:ROTDEFORMSPHERETRACEFREE} \\
	\mytr  \angdeform{\Rot_{(l)}}
	& = 2 \RotRadcomponent{l} \mytr \upchi^{(Small)}
		+ \mytr  \angdeform{\Rot_{(l)};Error},
		\label{E:ROTDEFORMSPHERETRACE}
	\end{align}
	\end{subequations}

\begin{subequations}
\begin{align}
	\angdeformoneformarg{\Rot_{(l)};Error}{\Lunit}
	& = G_{(Frame)} 
			\myarray
				[\Rot_{(l)}]
				{\RotRadcomponent{l}}
			\myarray
				[\Lunit \Psi]
				{\angdiff \Psi}
		+ \smoothfunction(\Psi) \Lunit_{(Small)} \angdiff x,
					\label{E:ERRORROTDEFORMLSPHERE} \\
	\angdeformoneformarg{\Rot_{(l)};Error}{\Rad}
	& = G_{(Frame)} 
				\myarray
					[\upmu 	\Rot_{(l)}]
					{\upmu \RotRadcomponent{l}}
			\myarray
				[\Lunit \Psi]
				{\angdiff \Psi}
		+ G_{(Frame)} 
			\RotRadcomponent{l} 
			\Rad \Psi
		+ \upmu
			\smoothfunction(\Psi) 
			\myarray[\Psi]
				{\Lunit_{(Small)}}	 
				\angdiff x,
				\label{E:ERRORROTDEFORMRADSPHERE}	\\				
	\angdeformfree{\Rot_{(l)};Error}
	& =  G_{(Frame)} 
			\hat{\otimes}
			\myarray
				[\Rot_{(l)}]
				{\RotRadcomponent{l}}
			\myarray
				[\Lunit \Psi]
				{\angdiff \Psi}
			+ \smoothfunction(\Psi) 
				\Psi 
				\angdiff x \hat{\otimes} \angdiff x,
		\label{E:ERRORROTDEFORMSPHERETRACEFREE} \\
\mytr  \angdeform{\Rot_{(l)};Error}
	& =  G_{(Frame)}
			\ginversesphere
			\myarray
				[	\Rot_{(l)}]
				{\RotRadcomponent{l}}
			\myarray
				[\Lunit \Psi]
				{\angdiff \Psi}
		+ \frac{\RotRadcomponent{l}}{\rgeo}
			+ \smoothfunction(\Psi) 
			 	\Psi 
			 	(\angdiffuparg{\#} x) \angdiff x.
			 	\label{E:ERRORROTDEFORMSPHERETRACE}
\end{align}
\end{subequations}
\end{proposition}

\begin{remark}[\textbf{Clarification of the meaning of $\Rot_{(l)}$ in the above schematic relations}]
	\label{R:CLARIFICATIONOFROT}
	In \eqref{E:ERRORROTDEFORMLSPHERE}-\eqref{E:ERRORROTDEFORMSPHERETRACE},
	the vectorfield $\Rot_{(l)}$ appearing in the first array is not acting as a differential operator,
	but rather as a tensorfield that is being contracted against other tensorfields
	such as $G_{(Frame)}$ and $\angdiff \Psi.$
\end{remark}

Before proving the proposition, we first compute some covariant derivatives of
the $S_{t,u}$ projection tensorfield $\sphereproject$ relative to frame 
vectorfields $\lbrace \Lunit, \Rad, X_1, X_2 \rbrace.$

\begin{lemma}[\textbf{Frame covariant derivatives of the spherical projection tensorfield} $\sphereproject$]
\label{L:FRAMECOVARIANTDERIVATIVESOFSPHERICALPROJECTION}
Let $\sphereproject$ be the type $\binom{1}{1}$ $S_{t,u}$ projection tensorfield defined in \eqref{E:SPHEREPROJECTION}.
Then we have the following identities:
\begin{subequations}
\begin{align}
	(\D_{\Lunit} \sphereproject) \Radunit
	& = \upmu^{-1} \upzeta^{(Trans-\Psi) A} X_A
		+ \upzeta^{(Tan-\Psi) A} X_A,
		\label{E:DLSPHEREPROJECTAPPLIEDTORADUNIT}
				\\
	(\D_{\Lunit} \sphereproject) X_A
	& = - \upmu^{-1} \upzeta_A^{(Trans-\Psi)} \Lunit
		- \upzeta_A^{(Tan-\Psi)} \Lunit,
		\label{E:DLSPHEREPROJECTAPPLIEDTOXA}
				\\
	(\D_{\Rad} \sphereproject) \Radunit
	& = (\angdiffuparg{A} \upmu) X_A,
		\label{E:DRADSPHEREPROJECTAPPLIEDTORADUNIT} \\
		(\D_{\Rad} \sphereproject) X_A
	& = \upzeta_A^{(Trans-\Psi)} \Lunit
		+ \upmu \upzeta_A^{(Tan-\Psi)} \Lunit
		+ \upzeta_A^{(Trans-\Psi)} \Radunit
		+ \upmu \upzeta_A^{(Tan-\Psi)} \Radunit
		+ (\angdiffarg{A} \upmu) \Radunit,
		\label{E:DRADSPHEREPROJECTAPPLIEDTOA} \\
	(\D_A \sphereproject) \Radunit
	& = \upchi_A^{\ B} X_B
		- \upmu^{-1} \angkmixedarg{A}{(Trans-\Psi) B} X_B
		- \angkmixedarg{A}{(Tan-\Psi) B} X_B,
		\label{E:DASPHEREPROJECTAPPLIEDTORADUNIT} \\
	(\D_A \sphereproject) X_B
	& = \upmu^{-1} \angktriplearg{A}{B}{(Trans-\Psi)} \Lunit
		+ \angktriplearg{A}{B}{(Tan-\Psi)} \Lunit
		+ \upchi_{AB} \Radunit.
		\label{E:DASPHEREPROJECTAPPLIEDTOB}
\end{align}
\end{subequations}
	In the above expressions, the $S_{t,u}$ tensorfields
	$\upchi,$ 
	$\upzeta^{(Trans-\Psi)},$
	$\angkuparg{(Trans-\Psi)},$ 
	$\upzeta^{(Tan-\Psi)},$
	and
	$\angkuparg{(Tan-\Psi)}$
	are defined by
	\eqref{E:CHIDEF},
	\eqref{E:ZETATRANSVERSAL},
	\eqref{E:KABTRANSVERSAL},
	\eqref{E:ZETAGOOD},
	and \eqref{E:KABGOOD}.
\end{lemma}

\begin{proof}
The main idea of the proof is to use the decompositions provided by Lemma~\ref{L:CONNECTIONLRADFRAME}.
As examples, we prove \eqref{E:DLSPHEREPROJECTAPPLIEDTORADUNIT} and \eqref{E:DRADSPHEREPROJECTAPPLIEDTOA}. 
The remaining identities in \eqref{E:DLSPHEREPROJECTAPPLIEDTORADUNIT}-\eqref{E:DASPHEREPROJECTAPPLIEDTOB}
can be proved using similar arguments.
To prove \eqref{E:DLSPHEREPROJECTAPPLIEDTORADUNIT},
we differentiate the identity
$\sphereproject \Rad = 0$
and use the identity $\Rad = \upmu \Radunit$
to deduce that
$(\D_{\Lunit} \sphereproject) \Radunit 
= \upmu^{-1} (\D_{\Lunit} \sphereproject) \Rad
= - \upmu^{-1} \sphereproject \D_{\Lunit} \Rad.$
The desired identity \eqref{E:DLSPHEREPROJECTAPPLIEDTORADUNIT}
now follows easily from the previous identity, \eqref{E:DLRAD}, and \eqref{E:ZETADECOMPOSED}.

To prove \eqref{E:DRADSPHEREPROJECTAPPLIEDTOA},
we differentiate the identity $\sphereproject X_A = X_A$ to deduce 
$(\D_{\Rad} \sphereproject) X_A = \D_{\Rad} X_A - \sphereproject  \D_{\Rad} X_A.$
The desired identity \eqref{E:DASPHEREPROJECTAPPLIEDTORADUNIT} now
follows easily from the previous identity, 
the fact that $\D_{\Rad} X_A - \D_A \Rad = [\Rad,X_A]$ is $S_{t,u}-$tangent
(see \eqref{E:RADCOMMUTETANGENTISTANGENT}),
\eqref{E:DARAD}, 
and
\eqref{E:ZETADECOMPOSED}.

\end{proof}

\begin{proof}[Proof of Prop.~\ref{P:DEFORMATIONTENSORFRAMECOMPONENTS}]
The main idea of the proof is to use the decompositions provided by Lemmas \ref{L:CONNECTIONLRADFRAME}
and \ref{L:FRAMECOVARIANTDERIVATIVESOFSPHERICALPROJECTION}. We give four examples, which
by far involve the most difficult computations.

We begin by proving \eqref{E:ROTDEFORMSPHERETRACEFREE} and \eqref{E:ROTDEFORMSPHERETRACE}.
These identities and \eqref{E:ROTDEFORMLSPHERE} 
are the only ones in the proposition for which 
it requires some effort to see the cancellation of some terms
involving the dangerous transversal derivative $\Radunit \Psi = \upmu^{-1} \Rad \Psi.$
To proceed, we use Defs.~\ref{D:DEFORMTENSDEFINED} and \ref{D:ROTATION}
to deduce that
\begin{align} \label{E:FIRSTRELATIONROTSPHEREAB}
	\angdeformarg{\Rot_{(l)}}{A}{B}
	& = g(\D_A [\sphereproject \Roteucarg{l}], X_B) 
		+ g(\D_B [\sphereproject \Roteucarg{l}], X_A) 
		\\
	& = g([\D_A \sphereproject] \Roteucarg{l}, X_B) 
		+ g([\D_B \sphereproject] \Roteucarg{l}, X_A)
		\notag	\\
	& \ \
		+ g(\D_A \Roteucarg{l}, X_B) 
		+ g(\D_B \Roteucarg{l}, X_A).
			\notag
\end{align}
From \eqref{E:ROTATIONDECOMPOSITIONINTOEUCLIDEANPLUSRADCOMPONENT},
\eqref{E:DASPHEREPROJECTAPPLIEDTORADUNIT},
\eqref{E:DASPHEREPROJECTAPPLIEDTOB},
\eqref{E:ANGKDECOMPOSED},
\eqref{E:KABTRANSVERSAL}, 
and \eqref{E:KABGOOD},
we deduce that
\begin{align} \label{E:DASPHEREPROJECTROTEUCB}
g([\D_A \sphereproject] \Roteucarg{l}, X_B)
+ g([\D_B \sphereproject] \Roteucarg{l}, X_A)
& = 2 \RotRadcomponent{l} \upchi_{AB}
	- 2 \upmu^{-1} \RotRadcomponent{l} \angktriplearg{A}{B}{(Trans-\Psi)}
	- 2 \RotRadcomponent{l} \angktriplearg{A}{B}{(Tan-\Psi)}
		\\
& = 2 \RotRadcomponent{l} \upchi_{AB}
	- \upmu^{-1} \RotRadcomponent{l} \angGdoublearg{A}{B} \Rad \Psi
		+ \RotRadcomponent{l}
			G_{(Frame)} 
			\myarray
				[\Lunit \Psi]
				{\angdiff \Psi}.
	\notag
\end{align}
Next, we use 
\eqref{E:RECTANGULARCOVARIANTDERIVATIVE},
\eqref{E:RECTANGULARCHRISTOFFEL}, 
Def.~\ref{D:ROTATION}, 
\eqref{E:ROTATIONDECOMPOSITIONINTOEUCLIDEANPLUSRADCOMPONENT},
the identity $\Radunit = \upmu^{-1} \Rad,$
and the identity $\angdiffarg{A} x^i = X_A^i$ to compute 
(relative to rectangular coordinates) that
\begin{align} \label{E:DAROTEUCB}
	g(\D_A \Roteucarg{l}, X_B)
	+ g(\D_B \Roteucarg{l}, X_A)
	& = \epsilon_{lca} g_{ab} \angdiffarg{A} x^c \angdiffarg{B} x^b
		+ \epsilon_{lca} g_{ab} \angdiffarg{B} x^c \angdiffarg{A} x^b
		+ \angGdoublearg{A}{B} \Roteucarg{l} \Psi.
\end{align}
We now insert the decompositions $\Roteucarg{l} = \Rot_{(l)} + \RotRadcomponent{l} \Radunit$
and $g_{ab} = \delta_{ab} + g_{ab}^{(Small)}$ 
into \eqref{E:DAROTEUCB}
and observe that 
by the antisymmetry of $\epsilon_{\cdots},$ the $\delta_{ab}$ part cancels
from the first two terms on the right-hand side of \eqref{E:DAROTEUCB}.
Hence, we have
\begin{align} \label{E:SECONDFORMDAROTEUCB}
	g(\D_A \Roteucarg{l}, X_B)
	& + g(\D_B \Roteucarg{l}, X_A)	
		\\
	& = \epsilon_{lca} g_{ab}^{(Small)} (\angdiffarg{A} x^c \angdiffarg{B} x^b + \angdiffarg{B} x^c \angdiffarg{A} x^b)
		\notag \\
	& \ \ + \upmu^{-1} \RotRadcomponent{l} \angGdoublearg{A}{B} \Rad \Psi
			+ \angGdoublearg{A}{B} \Rot_{(l)} \Psi.
		\notag 
\end{align}
We now add \eqref{E:DASPHEREPROJECTROTEUCB} and \eqref{E:SECONDFORMDAROTEUCB}
and note the cancellation of the dangerous $\upmu^{-1} \RotRadcomponent{l} \angGdoublearg{A}{B} \Rad \Psi$ terms.
Also using \eqref{E:FIRSTRELATIONROTSPHEREAB} 
and the fact that $g_{ab}^{(Small)}$ is a smooth function of $\Psi$
that vanishes at $\Psi = 0,$ we conclude that
\begin{align} \label{E:ROTDEFORMSPHERE}
	\angdeformarg{\Rot_{(l)}}{A}{B}
	& = 2 \RotRadcomponent{l} \upchi_{AB}
			+ G_{(Frame)} 
			\myarray
				[\Rot_{(l)}]
				{\RotRadcomponent{l}}
			\myarray
				[\Lunit \Psi]
				{\angdiff \Psi}
			+ \smoothfunction(\Psi) 
				\Psi 
				\angdiff x \otimes \angdiff x.
\end{align}
The desired identities
\eqref{E:ROTDEFORMSPHERETRACEFREE} and \eqref{E:ROTDEFORMSPHERETRACE}
now follow from \eqref{E:ROTDEFORMSPHERE}.

As a third example, we now prove \eqref{E:ROTDEFORMLSPHERE}.
To proceed, we use Defs.~\ref{D:ROTATION} and \ref{D:DEFORMTENSDEFINED} to deduce
\begin{align} \label{E:FIRSTRELATIONROTLSPHERE}
	\angdeformarg{\Rot_{(l)}}{\Lunit}{A}
	& = g(\D_{\Lunit} [\sphereproject \Roteucarg{l}], X_A) 
		+ g(\D_A [\sphereproject \Roteucarg{l}], \Lunit) 
		\\
	& \ \ = g([\D_{\Lunit} \sphereproject] \Roteucarg{l}, X_A) 
		+ g([\D_A \sphereproject] \Roteucarg{l}, \Lunit)
		+ g(\D_{\Lunit} \Roteucarg{l}, X_A).
			\notag
\end{align}
From 
\eqref{E:CHIJUNKDEF},
\eqref{E:ROTATIONDECOMPOSITIONINTOEUCLIDEANPLUSRADCOMPONENT},
\eqref{E:ROTATATIONACOMPONENTALTERNATEEXPRESSION},
\eqref{E:DLSPHEREPROJECTAPPLIEDTORADUNIT},
\eqref{E:DLSPHEREPROJECTAPPLIEDTOXA},
\eqref{E:DASPHEREPROJECTAPPLIEDTORADUNIT},
\eqref{E:DASPHEREPROJECTAPPLIEDTOB},
\eqref{E:ZETADECOMPOSED},
\eqref{E:ZETATRANSVERSAL}, 
and \eqref{E:ZETAGOOD},
we deduce that
\begin{align} \label{E:DLSPHEREPROJECTROTEUCAPLUSDASPHEREPROJECTROTEUCL}
g([\D_{\Lunit} \sphereproject] \Roteucarg{l}, X_A)
& + g([\D_A \sphereproject] \Roteucarg{l}, \Lunit)
	\\
& = - \upchi_{AB} \Rot_{(l)}^B
		+ \upmu^{-1} \RotRadcomponent{l} \upzeta_A^{(Trans-\Psi)} 
		+ \RotRadcomponent{l} \upzeta_A^{(Tan-\Psi)}
	\notag	\\
& = - \frac{1}{\rgeo} \epsilon_{lca} g_{ab} x^c \angdiffarg{A} x^b 
		- \upchi_{AB}^{(Small)} \Rot_{(l)}^B
		- \frac{1}{2} \upmu^{-1} \RotRadcomponent{l} \angGdoublearg{\Lunit}{A} \Rad \Psi
		+ \RotRadcomponent{l} 
			G_{(Frame)}
			\myarray
				[\Lunit \Psi]
			{\angdiff \Psi}.
		\notag
\end{align}
Next, we use 
\eqref{E:RECTANGULARCOVARIANTDERIVATIVE},
\eqref{E:RECTANGULARCHRISTOFFEL}, 
Def.~\ref{D:ROTATION}, 
\eqref{E:ROTATIONDECOMPOSITIONINTOEUCLIDEANPLUSRADCOMPONENT},
the identity $\Radunit = \upmu^{-1} \Rad,$
and the identity $\angdiffarg{A} x^i = X_A^i$ to compute 
(relative to rectangular coordinates) that
\begin{align} \label{E:DLROTEUCAPLUSDAROTEUCL}
	g(\D_{\Lunit} \Roteucarg{l}, X_A)
	& = \epsilon_{lca} g_{ab} \Lunit^c \angdiffarg{A} x^b
		+ \frac{1}{2} \angGdoublearg{\Lunit}{A} \Roteucarg{l} \Psi
			\\
	& \ \ 
		+ \frac{1}{2} \angGdoublearg{\Roteucarg{l}}{A} \Lunit \Psi
		- \frac{1}{2} G_{\Lunit \Roteucarg{l}} \angdiffarg{A} \Psi.
		\notag
\end{align}
We now insert the decompositions $\Roteucarg{l} = \Rot_{(l)} + \RotRadcomponent{l} \Radunit$
and $\Lunit^c = \frac{x^c}{\rgeo} + \Lunit_{(Small)}^c$ 
(see \eqref{E:LUNITJUNK})
into \eqref{E:DLROTEUCAPLUSDAROTEUCL} and deduce that
\begin{align} \label{E:SECONDFORMDLROTEUCAPLUSDAROTEUCL}
	g(\D_{\Lunit} \Roteucarg{l}, X_A)
	& = \frac{1}{\rgeo} \epsilon_{lca} g_{ab} x^c \angdiffarg{A} x^b
		+ \epsilon_{lca} g_{ab} \Lunit_{(Small)}^c \angdiffarg{A} x^b
		+ \frac{1}{2} \upmu^{-1} \RotRadcomponent{l} \angGdoublearg{\Lunit}{A} \Rad \Psi
			\\
	& \ \
		+ G_{(Frame)} 
				\myarray
					[\Rot_{(l)}]
					{\RotRadcomponent{l}}
			\myarray
				[\Lunit \Psi]
				{\angdiff \Psi}.
				\notag
\end{align}
We now add \eqref{E:DLSPHEREPROJECTROTEUCAPLUSDASPHEREPROJECTROTEUCL}
and \eqref{E:SECONDFORMDLROTEUCAPLUSDAROTEUCL} and observe
that the first and third terms on the right-hand sides cancel.
The desired identity \eqref{E:ROTDEFORMLSPHERE} thus follows.

As a final example, we prove \eqref{E:ROTDEFORMRADSPHERE}.
This example is somewhat simpler than the previous examples
in the sense that we do not have to observe the cancellation of any 
dangerous $\upmu^{-1} \Rad \Psi-$containing terms
(although we still have to observe some other cancellations).
To proceed, we use Defs.~\ref{D:DEFORMTENSDEFINED} and \ref{D:ROTATION} to deduce
\begin{align} \label{E:FIRSTRELATIONROTRADSPHERE}
	\angdeformarg{\Rot_{(l)}}{\Rad}{A}
	& = g(\D_{\Rad} [\sphereproject \Roteucarg{l}], X_A) 
		+ g(\D_A [\sphereproject \Roteucarg{l}], \Rad) 
		\\
	& \ \ = g([\D_{\Rad} \sphereproject] \Roteucarg{l}, X_A) 
		+ g([\D_A \sphereproject] \Roteucarg{l}, \Rad)
		+ g(\D_{\Rad} \Roteucarg{l}, X_A).
			\notag
\end{align}
Next, using
\eqref{E:CHIJUNKDEF},
\eqref{E:ROTATIONDECOMPOSITIONINTOEUCLIDEANPLUSRADCOMPONENT},
\eqref{E:ROTATATIONACOMPONENTALTERNATEEXPRESSION},
\eqref{E:DRADSPHEREPROJECTAPPLIEDTORADUNIT},
\eqref{E:DRADSPHEREPROJECTAPPLIEDTOA},
\eqref{E:DASPHEREPROJECTAPPLIEDTORADUNIT},
\eqref{E:DASPHEREPROJECTAPPLIEDTOB},
\eqref{E:ANGKDECOMPOSED},
\eqref{E:KABTRANSVERSAL}, 
and \eqref{E:KABGOOD},
we deduce that
\begin{align} \label{E:DRADSPHEREPROJECTROTEUCAPLUSDASPHEREPROJECTROTEUCRAD}
g([\D_{\Rad} \sphereproject] \Roteucarg{l}, X_A)
& + g([\D_A \sphereproject] \Roteucarg{l}, \Rad)
	\\
& = \RotRadcomponent{l} \angdiffarg{A} \upmu
	- \angktriplearg{A}{B}{(Trans-\Psi)} \Rot_{(l)}^B
	- \upmu \angktriplearg{A}{B}{(Tan-\Psi)} \Rot_{(l)}^B
	+ \upmu \upchi_{AB} \Rot_{(l)}^B 
	\notag	\\
& = \frac{1}{\rgeo} \upmu \epsilon_{lca} g_{ab} x^c \angdiffarg{A} x^b 
		+ \upmu \upchi_{AB}^{(Small)} \Rot_{(l)}^B
		+ \RotRadcomponent{l} \angdiffarg{A} \upmu
		- \frac{1}{2} \angGdoublearg{A}{B} \Rot_{(l)}^B \Rad \Psi
		\notag	\\
& \ \ + 
			G_{(Frame)} 
			\myarray
				[\upmu \Lunit \Psi]
				{\upmu \angdiff \Psi}
				\Rot_{(l)}.
		\notag
\end{align}
Next, from 
\eqref{E:RECTANGULARCOVARIANTDERIVATIVE},
\eqref{E:RECTANGULARCHRISTOFFEL}, 
Def.~\ref{D:ROTATION}, 
\eqref{E:ROTATIONDECOMPOSITIONINTOEUCLIDEANPLUSRADCOMPONENT},
the identity $\Radunit = \upmu^{-1} \Rad,$
and the identity $\angdiffarg{A} x^b = X_A^b,$ we compute 
(relative to rectangular coordinates) that
\begin{align} \label{E:DRADROTEUCAPLUSDAROTEUCRAD}
	g(\D_{\Rad} \Roteucarg{l}, X_A) 
	& = \upmu \epsilon_{lca} g_{ab} \Radunit^c \angdiffarg{A} x^b
		+ \frac{1}{2} \angGdoublearg{\Roteucarg{l}}{A} \Rad \Psi
		+ \frac{1}{2} \upmu \angGdoublearg{\Radunit}{A} \Roteucarg{l} \Psi
		- \frac{1}{2} \upmu G_{\Radunit \Roteucarg{l}} \angdiffarg{A} \Psi.
\end{align}
We now insert the decompositions $\Roteucarg{l} = \Rot_{(l)} + \RotRadcomponent{l} \Radunit$
and $\Radunit^c = -\frac{x^c}{\rgeo} + \Radunit_{(Small)}^c$
(see \eqref{E:RADUNITJUNK})
into \eqref{E:DRADROTEUCAPLUSDAROTEUCRAD} and deduce that
\begin{align} \label{E:SECONDFORMDRADROTEUCA}
	g(\D_{\Rad} \Roteucarg{l}, X_A) 
	& = - \frac{1}{\rgeo} \upmu \epsilon_{lca} g_{ab} x^c \angdiffarg{A} x^b
		+ \upmu \epsilon_{lca} g_{bc} \Radunit_{(Small)}^c \angdiffarg{A} x^b
		+ \frac{1}{2} \angGdoublearg{A}{B} \Rot_{(l)}^B \Rad \Psi
			\\
& \ \ + G_{(Frame)}
				\left(
					\upmu \angdiff \Psi
				\right)
				\Rot_{(l)}
		+ G_{(Frame)} \RotRadcomponent{l} \Rad \Psi
		+ \upmu G_{(Frame)} \RotRadcomponent{l} \angdiff \Psi.
				\notag
\end{align}
Adding \eqref{E:DRADSPHEREPROJECTROTEUCAPLUSDASPHEREPROJECTROTEUCRAD} and \eqref{E:SECONDFORMDRADROTEUCA},
noting that the first and third terms on the right-hand side of
\eqref{E:SECONDFORMDRADROTEUCA} are exactly canceled by two terms on the right-hand side of 
\eqref{E:DRADSPHEREPROJECTROTEUCAPLUSDASPHEREPROJECTROTEUCRAD},
and using the fact that $\Radunit^i = - \Lunit_{(Small)}^i$ plus a smooth function of $\Psi$
that vanishes at $\Psi = 0$ (that is, \eqref{E:RADUNITJUNKLIKELMINUSUNITJUNK}),
we arrive at \eqref{E:ROTDEFORMRADSPHERE}.

The remaining identities in the proposition can be proved in a similar fashion
and the computations are much simpler; we leave the details to the reader.

\end{proof}

\begin{corollary}[\textbf{Expression for the spherical covariant of the rotation vectorfields}]
	\label{C:ANGDROTATION}
	The type $\binom{1}{1}$ $S_{t,u}$ tensorfield $\angD \Rot_{(l)}$
	can be expressed as follows, 
	where we are using the notation of 
	Sects.~\ref{S:TRACEANDTRACEFREEPARTS}
	and \ref{S:SCHEMATIC}
	and Remark~\ref{R:CLARIFICATIONOFROT},
	and $\smoothfunction$ is a smooth function of $\Psi:$
	\begin{align} \label{E:ANGDROTATION}
	\angDarg{A} \Rot_{(l)}^B
	& = \epsilon_{lab} (\angdiffarg{A} x^a) \angdiffuparg{B} x^b
			+ \RotRadcomponent{l} \upchi_A^{\ B}
			+ G_{(Frame)}
			\ginversesphere
			\myarray
				[\Rot_{(l)}]
				{\RotRadcomponent{l}}
			\myarray
				[\Lunit \Psi]
				{\angdiff \Psi}
			+ \smoothfunction(\Psi) 
				\Psi 
				\angdiff x \otimes \angdiffuparg{\#} x.
\end{align}
\end{corollary}

\begin{proof}
	Since $\angDarg{A} \Rot_{(l)}^B = (\ginversesphere)^{BC} g(\D_A \Rot_{(l)},X_B),$
	the proof is essentially the same as the proof of \eqref{E:ROTDEFORMSPHERE}.
	The main difference is that in the present case, we do not have the cancellation of the 
	``$\delta_{ab}$ part'' noted below
	equation \eqref{E:DAROTEUCB} because the quantity under consideration 
	is not symmetrized over the indices $A,B.$
	This results in the presence of the ``non-small'' term 
	$\epsilon_{lab} (\angdiffarg{A} x^a) \angdiffuparg{B} x^b$
	on the right-hand side of \eqref{E:ANGDROTATION}.
\end{proof}


\chapter{Geometric Operator Commutator Formulas and Notation for Repeated Differentiation}
\label{C:GEOMETRICFORMULAS}
\thispagestyle{fancy}
In Chapter~\ref{C:GEOMETRICFORMULAS}, we provide geometric expressions for the commutators of some operators that we encounter in our analysis. We also define some shorthand notation connected to repeated differentiation with respect to commutation
vectorfields $Z \in \mathscr{Z}.$

\section{Some differential operators}
In this section, we define various differential operators that play a role in our analysis.

\begin{definition}[\textbf{Angular divergence}]
\label{D:ANGDIV}
\begin{subequations}
If $Y$ is an $S_{t,u}-$tangent vectorfield, then we define
$\angdiv Y$ to be the following function:
\begin{align} \label{E:ANGDIVVECTORFIELD}
	\angdiv Y 
	:= \angDarg{A} Y^A.
\end{align}

If $\xi$ is an $S_{t,u}$ one-form, then we define
$\angdiv \xi$ to be the following function:
\begin{align} \label{E:ANGDIVONEFORM}
	\angdiv \xi  
	:= (\ginversesphere)^{AB}\angDarg{A} \xi_B.
\end{align}

If $\xi$ is a type $\binom{1}{1}$ $S_{t,u}$ tensorfield, then we define
$\angdiv \xi$ to be the following $S_{t,u}$ one-form:
\begin{align} \label{E:ANGDIVYPE11}
	(\angdiv \xi)_A
	:= \angDarg{B} \xi_{\ A}^B.
\end{align}

If $\xi$ is a symmetric type $\binom{0}{2}$ $S_{t,u}$ tensorfield, then we define
$\angdiv \xi$ to be the following $S_{t,u}$ one-form:
\begin{align} \label{D:ANGDIVSYMMETRICTYPE02}
	(\angdiv \xi)_A
	:= (\ginversesphere)^{BC}\angDarg{B} \xi_{AC}.
\end{align}

If $\xi$ is a symmetric type $\binom{2}{0}$ $S_{t,u}$ tensorfield, then we define
$\angdiv \xi$ to be the following $S_{t,u}-$tangent vectorfield:
\begin{align} \label{E:ANGDIVSYMMETRICTYPE20}
	(\angdiv \xi)^A
	:= \angDarg{B} \xi^{AB}.
\end{align}
\end{subequations}
\end{definition}

\begin{definition}[\textbf{Trace-free $S_{t,u}-$projected Lie derivatives}]
\label{D:ANGFREELIEDEF}
If $\xi$ is a type 
$\binom{0}{2}$ $S_{t,u}$ tensorfield
and $V$ is a vectorfield,
then we define the trace-free 
type $\binom{0}{2}$ $S_{t,u}$ tensorfield
$\angfreeLiearg{V} \xi$ as follows:
\begin{align} \label{E:ANGFREELIEDEF}
	\angfreeLiearg{V} \xi_{AB}
	& := \angLie_V \xi_{AB}
		- \frac{1}{2} 
			(\mytr  \angLie_V \xi)
			\gsphere_{AB}.
\end{align}	


\end{definition}

\begin{definition}[\textbf{Trace-free part of} $\angD^2$]
\label{D:ANGFREEDDEF}
If $f$ is any function, then we define the trace-free type 
$\binom{0}{2}$ $S_{t,u}$ tensorfield $\angfreeDsquared f$ as follows:
\begin{align} \label{E:ANGFREEDDEF}
	\angfreeDsquaredarg{A}{B} f
	& := \angDsquaredarg{A}{B} f
		- \frac{1}{2} (\angLap f) \gsphere_{AB}.
\end{align}

\end{definition}

\begin{definition}[\textbf{The operator} $\angcheckD$]
\label{D:ANGDCHECKDEF}
If $\xi$ is a type $\binom{0}{2}$ $S_{t,u}$ tensorfield, then we define
the type $\binom{0}{3}$ $S_{t,u}$ tensorfield $\angcheckD \xi$ as follows:
\begin{align} \label{E:ANGDCHECKDEF}
\angcheckDarg{A} \xi_{BC}
	& := \frac{1}{2} 
		\left\lbrace
			\angDarg{A} \xi_{BC}
			+ \angDarg{B} \xi_{AC}
			- \angDarg{C} \xi_{AB}
		\right\rbrace .
\end{align}

\end{definition}

\section{Operator commutator identities}

In this section, we provide a collection of operator commutator identities
that we use in our analysis.

\begin{definition}[\textbf{Commutator of two operators}]
	\label{D:COMMMUTATOR}
	If $P$ and $Q$ are two operators,
	then $[P,Q]$ denotes their commutator.
\end{definition}

\begin{lemma}[$\Lunit,$ $\rgeo \Lunit,$ $\Rad,$ $\Rot$ \textbf{commute with} $\angdiff$]
\label{L:LANDRADCOMMUTEWITHANGDIFF}
If $f$ is a function and $V \in \lbrace \Lunit, \rgeo \Lunit, \Rad, \Rot_{(1)}, \Rot_{(2)}, \Rot_{(3)}, \rbrace,$
then
\begin{align}
	\angLie_V \angdiff f
	& = \angdiff V f.
		\label{E:ANGLIECOMMUTESWITHANGDIFF}
\end{align}

\end{lemma}

\begin{proof}
	We prove \eqref{E:ANGLIECOMMUTESWITHANGDIFF} only when $V = \Rad$ since the 
	other vectorfields $V$ can be treated by using a similar argument
	and the fact that $\angdiff \rgeo = 0.$ 
	To proceed, we set $V=\Rad,$
	contract the right-hand side of
	\eqref{E:ANGLIECOMMUTESWITHANGDIFF} against the vector $X_A,$ 
	and use the Leibniz rule together with the
	fact that $[X_A, \Rad]$ is $S_{t,u}-$tangent (which follows from \eqref{E:RADCOMMUTETANGENTISTANGENT})
	in order to deduce that
	\begin{align} \label{E:RHSOFANGLIECOMMUTESWITHANGDIFF}
		\angdiffarg{A} \Rad f = X_A(\Rad f) = \Rad(X_A f) + \angdiff f \cdot [X_A, \Rad].
	\end{align}
	On the other hand, we contract the left-hand side of \eqref{E:ANGLIECOMMUTESWITHANGDIFF} against $X_A$
	and use the Leibniz rule together with
	\eqref{E:RADCOMMUTETANGENTISTANGENT} to deduce that
	\begin{align} \label{E:LHSOFANGLIECOMMUTESWITHANGDIFF}
		(\angLie_{\Rad} \angdiff f)\cdot X_A 
		& = (\Lie_{\Rad} \angdiff f)\cdot X_A
			= \Rad(\angdiff f \cdot X_A) + \angdiff f \cdot [X_A, \Rad].
	\end{align}
	We now note that the right-hand sides of 
	\eqref{E:RHSOFANGLIECOMMUTESWITHANGDIFF} and
	\eqref{E:LHSOFANGLIECOMMUTESWITHANGDIFF}
	are equal, which yields the desired identity \eqref{E:ANGLIECOMMUTESWITHANGDIFF}.
	
\end{proof}

\begin{lemma}[\textbf{Leibniz rules with angular differential commutators}]
	\label{L:LEIBNIZRULEWITHANGULARDIFFERENTIALCOMMUTATIONS}
	If $X$ is an $S_{t,u}-$tangent vectorfield, 
	$V \in \lbrace \Lunit, \rgeo \Lunit, \Rad, \Rot_{(1)}, \Rot_{(2)}, \Rot_{(3)}, \rbrace,$
	and $f$ is a function, then the following commutator identity holds:
	\begin{align}
		V (X \cdot \angdiff f)
			& = X \cdot \angdiff \Lunit f
				+ \angLie_V X \cdot \angdiff f.
	\end{align}
\end{lemma}

\begin{proof}
	Lemma~\ref{L:LEIBNIZRULEWITHANGULARDIFFERENTIALCOMMUTATIONS} follows from the Leibniz rule 
	and Lemma~\ref{L:LANDRADCOMMUTEWITHANGDIFF}.
\end{proof}

\begin{lemma}[\textbf{Commutator of the Lie derivative of the trace-free part and the trace-free Lie derivative}]
If $\xi_{AB}$ is a symmetric type $\binom{0}{2}$ $S_{t,u}$ tensorfield
and $Z \in \mathscr{Z},$ then
\begin{subequations}
\begin{align} \label{E:COMMUTINGLIEWITHTRACEFREE}
		\Lie_Z \hat{\xi}_{AB}
		- \angfreeLiearg{Z} \xi_{AB}
		& = \frac{1}{2} \angdeformfreeuparg{Z}{C}{D} \hat{\xi}_{CD} \gsphere_{AB}
			- \frac{1}{2} \mytr \xi \angdeformfreearg{Z}{A}{B}.
\end{align}
	
If in addition $\xi$ is trace-free, then
\begin{align}  \label{E:ANGFREELIEIDFORTRACEFREETENSORS}
		\angLie_Z \xi_{AB}
		- \angfreeLiearg{Z} \xi_{AB}
		& = 
			\frac{1}{2} 
			\angdeformfreeuparg{Z}{C}{D} \xi_{CD}
			\gsphere_{AB}.
\end{align}
\end{subequations}

\begin{proof}
	From the decomposition $\xi = \frac{1}{2} \mytr \xi \gsphere + \hat{\xi},$
	the Leibniz rule, and the identities 
	 $\angLie \gsphere_{AB} = \angdeformarg{Z}{A}{B}$
	 and
	 $(\angLie_Z \ginversesphere)^{AB} = - \angdeformuparg{Z}{A}{B},$
	we compute that
	\begin{align}
		Z \mytr \xi 
		& = - \angdeformuparg{Z}{A}{B} \xi_{AB}
			+ \mytr \angLie_Z \xi,
			\\
		\angLie_Z \hat{\xi}_{AB}
		& = \angLie_Z \xi_{AB} 
			- \frac{1}{2} \mytr \angLie_Z \xi \gsphere_{AB}
			+ \frac{1}{2} \angdeformuparg{Z}{C}{D} \xi_{CD} \gsphere_{AB}
			- \frac{1}{2} \mytr \xi \angdeformarg{Z}{A}{B}
				\label{E:ANGFREELIETRFREECOMPUTATION} \\
		& = \angfreeLiearg{Z} \xi_{AB}
			+ \frac{1}{2} \angdeformuparg{Z}{C}{D} \xi_{CD} \gsphere_{AB}
			- \frac{1}{2} \mytr \xi \angdeformarg{Z}{A}{B}
			\notag \\
		& = \angfreeLiearg{Z} \xi_{AB}
			+ \frac{1}{2} \angdeformfreeuparg{Z}{C}{D} \hat{\xi}_{CD} \gsphere_{AB}
			- \frac{1}{2} \mytr \xi \angdeformfreearg{Z}{A}{B},
			\notag
	\end{align}
	where to deduce the second equality in \eqref{E:ANGFREELIETRFREECOMPUTATION} we used 
	definition \eqref{E:ANGFREELIEDEF}, and to deduce the last equality in \eqref{E:ANGFREELIETRFREECOMPUTATION},
	we used the decompositions 
	$\angdeformarg{Z}{A}{B} = \frac{1}{2} \mytr \angdeform{Z} \gsphere_{AB} 
	+ \angdeformfreearg{Z}{A}{B}$
	and
	$\angdeformuparg{Z}{C}{D} = \frac{1}{2} \mytr \angdeform{Z} (\ginversesphere)^{CD} 
	+ \angdeformfreeuparg{Z}{C}{D}.$
	We have thus proved \eqref{E:COMMUTINGLIEWITHTRACEFREE}.	
	The identity \eqref{E:ANGFREELIEIDFORTRACEFREETENSORS} follows trivially from \eqref{E:COMMUTINGLIEWITHTRACEFREE}.
\end{proof}

\end{lemma}

%

%

\begin{lemma} [\textbf{Commuting $\rgeo \Lunit,$ $\Rad,$ $\Rot$ with $\angD$}]
\label{L:COMMUTINGVEDCTORFIELDSWITHANGD}
For any commutation vectorfield $Z \in \mathscr{Z}$ 
and for $Z = \Lunit,$
and for any $S_{t,u}$
one-form $\xi,$ we have the following commutator identity:
\begin{align} \label{E:ANGDANGLIEZONEFORMCOMMUTATOR}
	([\angDarg{A}, \angLie_Z] \xi)_B
	& = (\angcheckDarg{A} \angdeformmixedarg{Z}{B}{C}) \xi_C.
\end{align}

For any commutation vectorfield $Z \in \mathscr{Z}$ 
and for $Z = \Lunit,$
and for any type $\binom{0}{2}$ $S_{t,u}$  
tensorfield $\xi,$ we have the following commutator identities:
\begin{align} \label{E:ANGDANGLIEZTYPE02COMMUTATOR}
	([\angDarg{A}, \angLie_Z] \xi)_{BC}
	& = (\angcheckDarg{A} \angdeformmixedarg{Z}{B}{D}) \xi_{CD}
		+ (\angcheckDarg{A} \angdeformmixedarg{Z}{C}{D}) \xi_{BD},
			\\
	([\angDarg{A}, \angLie_Z] \angD \xi)_{BCD}
	& = (\angcheckDarg{A} \angdeformmixedarg{Z}{B}{E}) \angDarg{E} \xi_{CD}
		+ (\angcheckDarg{A} \angdeformmixedarg{Z}{C}{E}) \angDarg{B} \xi_{ED}
		+ (\angcheckDarg{A} \angdeformmixedarg{Z}{D}{E}) \angDarg{B} \xi_{CE}.
		\label{E:ANGDANGLIEZANGDTYPE02COMMUTATOR}
\end{align}

\end{lemma}

\begin{proof}
We first prove \eqref{E:ANGDANGLIEZTYPE02COMMUTATOR}
in the case $Z = \Lunit.$ 
We compute relative to the geometric 
coordinates, and we
use the fact that the vectorfield $\Lunit = \frac{\partial}{\partial t}$
commutes with the coordinate vectorfields 
$(\frac{\partial}{\partial \vartheta^1}, \frac{\partial}{\partial \vartheta^2}) = (X_1, X_2).$
Also using Lemma~\ref{L:CONNECTIONBETWEENANGLIEOFGSPHEREANDDEFORMATIONTENSORS},
we deduce the following identities:
$ \frac{\partial}{\partial t} \gsphere_{AB}
= \angLie_{\Lunit} \gsphere_{AB} 
= \angdeformarg{\Lunit}{A}{B},$
$ \frac{\partial}{\partial t} (\ginversesphere)^{AB}
= (\angLie_{\Lunit} \ginversesphere)^{AB} = - \angdeformuparg{\Lunit}{A}{B}.$
To proceed, we note that the Christoffel symbols $\sphereGamma_{A \ B}^{\ C}$ of $\gsphere$ 
(relative to the local coordinates $(\vartheta^1,\vartheta^2)$)
can be expressed as
\begin{align} \label{E:GSHPERECHRISTOFFEL}
	\sphereGamma_{A \ B}^{\ C} 
	& = \frac{1}{2} 
		(\ginversesphere)^{CD}
		\left\lbrace
			\frac{\partial}{\partial \vartheta^A} \gsphere_{DB}
			+ \frac{\partial}{\partial \vartheta^B} \gsphere_{AD}
			- \frac{\partial}{\partial \vartheta^D} \gsphere_{AB}
		\right\rbrace.
\end{align}
Next, we state the following standard covariant derivative identities: 
\begin{align} \label{E:ANGDXICOMPUTED}
	\angDarg{A} \xi_{BC} 
	& = \frac{\partial}{\partial \vartheta^A} \xi_{BC} 
		- \sphereGamma_{A \ B}^{\ D} \xi_{DC} 
		- \sphereGamma_{A \ C}^{\ D} \xi_{BD},
			\\
	\angDarg{A} \angLie_{\Lunit} \xi_{BC}
		& = \angDarg{A} \frac{\partial}{\partial t} \xi_{BC}
			= \frac{\partial}{\partial t} \frac{\partial}{\partial \vartheta^A} \xi_{BC}
			- \sphereGamma_{A \ B}^{\ D} \frac{\partial}{\partial t} \xi_{DC}
			- \sphereGamma_{A \ C}^{\ D} \frac{\partial}{\partial t} \xi_{BD}.
			\label{E:ANGDOFLIEXICOMPUTED}
\end{align}
We now set
$\sphereGamma_{A D B} 
:= \gsphere_{CD} \sphereGamma_{A \ B}^{\ C}
= \left\lbrace
	\partial_A \gsphere_{DB}
	+ \partial_B \gsphere_{AD}
	- \partial_D \gsphere_{AB}
\right\rbrace,$
refer to definition \eqref{E:ANGDCHECKDEF},
and use the above calculations
and the symmetry property 
$\sphereGamma_{A D B} = \sphereGamma_{B D A}$
to compute that
\begin{align} \label{E:COMPUTINGDERIVATIVESOFSPHERECHRISTOFFEL}
\frac{\partial}{\partial t} 
\sphereGamma_{A \ B}^{\ C} 
& = - 
	\angdeformuparg{\Lunit}{C}{D}
	\sphereGamma_{ADB}
	+ 
	\frac{1}{2} 
	(\ginversesphere)^{CD}
	\left\lbrace
		\frac{\partial}{\partial \vartheta^A} \angdeformarg{\Lunit}{D}{B}
		+ \frac{\partial}{\partial \vartheta^B} \angdeformarg{\Lunit}{A}{D}
		- \frac{\partial}{\partial \vartheta^D} \angdeformarg{\Lunit}{A}{B}
	\right\rbrace
	\\
& = (\ginversesphere)^{CD} \angcheckDarg{A}  \angdeformarg{\Lunit}{B}{D}.
\notag
\end{align}
Using \eqref{E:ANGDXICOMPUTED}
and 
\eqref{E:COMPUTINGDERIVATIVESOFSPHERECHRISTOFFEL},
we compute that
\begin{align} \label{E:ANGLIEOFANGDXICOMPUTED}
	\angLie_{\Lunit}
	\angDarg{A} \xi_{BC} 
	= \frac{\partial}{\partial t} 
		\angDarg{A} \xi_{BC} 
	& = \frac{\partial}{\partial \vartheta^A} \frac{\partial}{\partial t} \xi_{BC} 
		- \sphereGamma_{A \ B}^{\ D} \frac{\partial}{\partial t}  \xi_{DC} 
		- \sphereGamma_{A \ C}^{\ D} \frac{\partial}{\partial t} \xi_{BD}
			\\
	& \ \
		 - (\ginversesphere)^{CE} 
		 	 (\angcheckDarg{A} \angdeformarg{\Lunit}{B}{E})
			\xi_{DC}
		- (\ginversesphere)^{DE} 
			(\angcheckDarg{A} \angdeformarg{\Lunit}{C}{E})
			\xi_{BD}.
		\notag
\end{align}
Finally, subtracting
\eqref{E:ANGLIEOFANGDXICOMPUTED}
from \eqref{E:ANGDOFLIEXICOMPUTED} and
noting that the right-hand side of \eqref{E:ANGDOFLIEXICOMPUTED}
cancels the first three terms on the right-hand side of
\eqref{E:ANGLIEOFANGDXICOMPUTED},
we arrive at the desired identity \eqref{E:ANGDANGLIEZTYPE02COMMUTATOR}.

Since $\angdeform{\rgeo \Lunit} = \rgeo \angdeform{\Lunit},$
the identity \eqref{E:ANGDANGLIEZTYPE02COMMUTATOR}
for $\rgeo \Lunit$ follows directly from the identity
for $\Lunit.$
 
A similar argument involving one minor
adjustment yields
\eqref{E:ANGDANGLIEZTYPE02COMMUTATOR}
when $Z = \Rad.$
The minor adjustment is that
we first endow $\Sigma_t^{U^0}$ with local coordinates
$(u, \widetilde{\vartheta}^1, \widetilde{\vartheta}^2)$
(where $u$ is the eikonal function)
in such a way that $\Rad = \frac{\partial}{\partial u}|_{\widetilde{\vartheta}^1, \widetilde{\vartheta}^2}.$
We already carried out such a construction in our proof of Cor.~\ref{C:SPACETIMEVOLUMEFORMWITHUPMU}.
We can then proceed as in the proof of \eqref{E:ANGDANGLIEZTYPE02COMMUTATOR}
in the case $Z = \Lunit,$ but using the coordinates 
$\widetilde{\vartheta}^A$ in place of $\vartheta^A.$

To prove \eqref{E:ANGDANGLIEZTYPE02COMMUTATOR} in the case
$Z = \Rot,$ 
for each real number $\uplambda$ belonging to a small neighborhood of $0,$
we let $\varphi_{(\uplambda)}: \mathbb{R}^4 \rightarrow \mathbb{R}^4$ 
be the flow map of the $S_{t,u}-$tangent vectorfield $\Rot.$
That is, for each spacetime point $p,$
the four rectangular components of $\varphi_{(\uplambda)}(p)$ verify the ODE system
$\frac{d}{d \uplambda} \varphi_{(\uplambda)}^{\nu}(p) = \Rot^{\nu} \circ \varphi_{(\uplambda)}(p)$ 
and $\varphi_{(0)} = I,$ where $I$ is the identity map on $\mathbb{R}^4.$
Let 
$\varphi_{(\uplambda)}^* (\angD \xi),$
$\varphi_{(\uplambda)}^*  \gsphere,$
and
$\varphi_{(\uplambda)}^* \xi$
respectively denote the pullbacks of the 
tensorfields 
$\angD \xi,$
$\gsphere$
and 
$\xi$ by $\varphi_{(\uplambda)},$
and let 
$\varphi_{(\uplambda)}^* \angD$ denote the Levi-Civita connection of $\varphi_{(\uplambda)}^*  \gsphere.$
Then by covariance, we have the identity
$\varphi_{(\uplambda)}^* (\angD \xi)
	= (\varphi_{(\uplambda)}^* \angD) \varphi_{(\uplambda)}^* \xi.$

To simplify the notation, we denote
$\angDleftexp{\uplambda} := \varphi_{(\uplambda)}^* \angD,$
 $\gsphereleftexp{\uplambda} := \varphi_{(\uplambda)}^*  \gsphere,$ 
and $\xileftexp{\uplambda} := \varphi_{(\uplambda)}^* \xi.$
We also let $\sphereGammaleftexparg{\uplambda}{A}{B}{C}$ 
denote the 
Christoffel symbols of $\gsphereleftexp{\uplambda}$ relative to the local coordinates 
$(\vartheta^1,\vartheta^2).$
In particular, the above covariance identity can be written as
\begin{align} \label{E:DIFFEOMORPH}
	\varphi_{(\uplambda)}^* (\angD \xi)
	= \angDleftexp{\uplambda} \xileftexp{\uplambda}.
\end{align}

Using the fact that the Lie derivative of a tensorfield
with respect to a vectorfield
is the derivative (with respect to the flow parameter)
of the pullback of the tensorfield
by the flow map of the vectorfield,
the fact that $\Rot$ and $\gsphere = \sphereproject g$ are $S_{t,u}$ tensorfields,
and Cor.~\ref{C:ANGLIEVFORGOODVOFTENSOREQUALSANGLIEOFSTUPROJECTEDTENSOR},
we deduce the tensorial identities
$\frac{d}{d \uplambda}|_{\uplambda = 0} \gsphereleftexparg{\uplambda}{A}{B}
= \angLie_{\Rot} \gsphere_{AB}
= \angdeformarg{\Rot}{A}{B}$
and 
$\frac{d}{d \uplambda}|_{\uplambda = 0} \ginversesphereleftexparg{\uplambda}{A}{B} 
= - \angdeformuparg{\Rot}{A}{B},$
where $\ginversesphereleftexp{\uplambda}$ denotes the inverse of 
$\gsphereleftexp{\uplambda}.$
From these identities, we deduce that the identity
\eqref{E:COMPUTINGDERIVATIVESOFSPHERECHRISTOFFEL} holds
with the left-hand side replaced by
$\frac{d}{d \uplambda}|_{\uplambda = 0}
\sphereGammaleftexparg{\uplambda}{A}{B}{C}$
and all terms
$\angdeform{\Lunit}$
on the right-hand side replaced by $\angdeform{\Rot}.$ 
That is, we have
\begin{align} \label{E:LIEDERIVATIVEOFSPHERECHRISTOFFEL}
	\frac{d}{d \uplambda}|_{\uplambda = 0}
	\sphereGammaleftexparg{\uplambda}{A}{B}{C}
	& = (\ginversesphere)^{CD} \angcheckDarg{A} \angdeformarg{\Rot}{D}{B}.
\end{align}

We now consider the $ABC$ component of the right-hand side of 
\eqref{E:DIFFEOMORPH},
that is, $\angDleftexparg{\uplambda}{A} \xileftexparg{\uplambda}{B}{C}.$ 
This component is equal to
the right-hand side of \eqref{E:ANGDXICOMPUTED}
except with 
$\angDleftexp{\uplambda},$
$\xileftexp{\uplambda},$
and
$\sphereGammaleftexp{\uplambda}$
respectively 
in place of
$\angD,$
$\xi,$
and $\Gamma.$
Hence, applying the operator
$\frac{d}{d \uplambda}|_{\uplambda = 0}$
to the right-hand side of \eqref{E:DIFFEOMORPH}
and using the identity 
\eqref{E:LIEDERIVATIVEOFSPHERECHRISTOFFEL},
we deduce that 
\begin{align} \label{E:LIEDERIVATIVEOFANGDXIPULLBACKRHS}
	\left(
	\frac{d}{d \uplambda}|_{\uplambda = 0}
	\left\lbrace
		\angDleftexp{\uplambda}
		\xileftexp{\uplambda}
	\right\rbrace
	\right)_{ABC}
	& = (\angD \angLie_{\Rot} \xi)_{ABC}
		+ (\angcheckDarg{A} \angdeformmixedarg{\Rot}{B}{D}) \xi_{DC}
		+ (\angcheckDarg{A} \angdeformmixedarg{\Rot}{C}{D}) \xi_{BD}.
\end{align}
On the other hand, 
applying the operator
$\frac{d}{d \uplambda}|_{\uplambda = 0}$
to the left-hand side of \eqref{E:DIFFEOMORPH}
and considering the $ABC$ component,
we deduce 
(from the representation of Lie differentiation in terms of the derivative of the pullback by the flow map)
that
\begin{align} \label{E:LIEDERIVATIVEOFANGDXIPULLBACKLHS}
	\left(
		\frac{d}{d \uplambda}|_{\uplambda = 0} 
		\left\lbrace
			\varphi_{(\uplambda)}^* (\angD \xi)
		\right\rbrace
	\right)_{ABC}
	= (\angLie_{\Rot} \angD \xi)_{ABC}.
\end{align}
Subtracting 
\eqref{E:LIEDERIVATIVEOFANGDXIPULLBACKLHS}
from
\eqref{E:LIEDERIVATIVEOFANGDXIPULLBACKRHS},
we conclude \eqref{E:ANGDANGLIEZTYPE02COMMUTATOR} in the case $Z = \Rot.$
We have thus proved \eqref{E:ANGDANGLIEZTYPE02COMMUTATOR} in all cases.

The identities 
\eqref{E:ANGDANGLIEZONEFORMCOMMUTATOR}
and
\eqref{E:ANGDANGLIEZANGDTYPE02COMMUTATOR}
can be proved in an analogous fashion,
and we omit the details.

\end{proof}

\begin{lemma}[{$[\left\lbrace \angLie_{\Lunit} + \mytr \upchi \right\rbrace, \angdiv] = 0$} \textbf{for} 
$S_{t,u}-$\textbf{tangent vectorfields}] 
\label{L:LPLUSTRCHIANGDIVCOMMUTE}
If $Y$ is an $S_{t,u}-$tangent vectorfield, then
	\begin{align} \label{E:LPLUSTRCHIANGDIVCOMMUTE}
		\left\lbrace
			\Lunit
			+ \mytr \upchi
		\right\rbrace
		\angdiv Y
		& = \angdiv 
				\left\lbrace
					\angLie_{\Lunit} Y
					+ \mytr \upchi Y
				\right\rbrace.
	\end{align}
	
\end{lemma}
\begin{proof}
	Let $\xi$ be the $S_{t,u}$ one-form that is $\gsphere-$dual to $Y.$
	Using the Leibniz rule 
	and the identity $(\angLie_{\Lunit} \ginversesphere)^{AB} = - \angdeformuparg{\Lunit}{A}{B}$
	and commuting the operators $\angLie_{\Lunit}$ and $\angD,$
	we deduce that
	$\Lunit \angdiv Y = \Lunit\left\lbrace (\ginversesphere)^{AB} \angDarg{A} \xi_B) \right\rbrace
	= - \angdeformuparg{\Lunit}{A}{B} \angDarg{A} Y_B 
	+ (\ginversesphere)^{AB} \angDarg{A} \angLie_{\Lunit} \xi_B
	+ (\ginversesphere)^{AB} [\angLie_{\Lunit}, \angD]_A Y_B.$
	Using the identity $\angdeform{\Lunit}= 2 \upchi$
	and \eqref{E:ANGDANGLIEZONEFORMCOMMUTATOR} with $\Lunit$ in the role of 
	$Z,$ we deduce that 
	$(\ginversesphere)^{AB} [\angLie_{\Lunit}, \angD]_A \xi_B = Y \mytr \upchi - 2 (\angdiv \upchi) \cdot Y.$
	Furthermore, since $\xi_B = \gsphere_{BC} Y^C,$ we deduce from the Leibniz rule that
	$\angLie_{\Lunit} \xi_B = 2 \upchi_{BC} Y^C + \gsphere_{BC} \angLie_{\Lunit} Y^C.$
	Hence, we deduce that 
	$(\ginversesphere)^{AB} \angDarg{A} \angLie_{\Lunit} \xi_B = 
	2 (\angdiv \upchi) \cdot Y + 2 \upchi^{AB} \angDarg{A} Y_B + \angdiv \angLie_{\Lunit} Y.$
	Combining these identities, we deduce that
	\begin{align} \label{E:ALMOLSTDONELPLUSTRCHIANGDIVCOMMUTE}
		\Lunit \angdiv Y
		& = \angdiv \angLie_{\Lunit} Y
			+ Y \mytr \upchi.
	\end{align}
	The desired identity \eqref{E:LPLUSTRCHIANGDIVCOMMUTE}
	now follows easily from \eqref{E:ALMOLSTDONELPLUSTRCHIANGDIVCOMMUTE}.
\end{proof}

\begin{lemma}[\textbf{Commuting $\rgeo \Lunit,$ $\Rad,$ $\Rot$ with $\angD^2$}]
\label{L:COMMUTINGANGDSQUAREDWITHCOMMUTATIONVECTORFIELDS}
For any commutation vectorfield $Z \in \mathscr{Z}$ 
and for $Z = \Lunit,$
and for any scalar-valued function $f,$ we have
the following commutator identities:
\begin{subequations}
\begin{align}
	([\angD^2, \angLie_Z] f)_{AB}
	& = 
		(\angcheckDarg{A} \angdeformmixedarg{Z}{B}{C})
		\angdiffarg{C} f, 
		\label{E:COMMUTINGANGDSQUAREDANDLIEZ}
\end{align}
\begin{align}
	[\angLap, \angLie_Z] f
	& = \angdeformuparg{Z}{A}{B} \angDsquaredarg{A}{B} f
		+ (\angcheckDarg{A}  \angdeformuparg{Z}{A}{B})
			\angdiffarg{B} f.
			\label{E:COMMUTINGANGANGLAPANDLIEZ}
\end{align}
\end{subequations}
%

\end{lemma}

\begin{proof}
The identity \eqref{E:COMMUTINGANGDSQUAREDANDLIEZ}
follows from \eqref{E:ANGDANGLIEZONEFORMCOMMUTATOR}
with $\angdiff f$ in the role of $\xi$
and Lemma~\ref{L:LANDRADCOMMUTEWITHANGDIFF}.

To deduce \eqref{E:COMMUTINGANGANGLAPANDLIEZ},
we take the trace over the $AB$ indices in \eqref{E:COMMUTINGANGDSQUAREDANDLIEZ}
and use the operator commutator identity
$(\ginversesphere)^{AB} \angLie_Z \angD^2_{AB} f
= Z \angLap f + \angdeformuparg{Z}{A}{B} \angDsquaredarg{A}{B} f.$

\end{proof}

\section{Notation for repeated differentiation}
\label{S:REPEATEDDIFF}
In this section, we introduce some compact notation for repeated differentiation.

\begin{definition}[\textbf{Notation for repeated differentiation}]
\label{D:REPEATEDDIFFERENTIATIONSHORTHAND}
We fix a labeling $\mathscr{Z} = \lbrace Z_{(i)} \rbrace_{i=1}^5$ of the five commutation vectorfields
(see Def.~\ref{D:DEFSETOFCOMMUTATORVECTORFIELDS}).
We use the following shorthand notation.
\begin{itemize}
	\item If $\vec{I} = (\iota_{(1)}, \iota_{(2)}, \cdots, \iota_{(N)})$ is a multi-index
	of length $N$
		with $\iota_{(1)}, \iota_{(2)}, \cdots, \iota_{(N)} \in \lbrace 1,2,3,4,5 \rbrace,$
		then $\mathscr{Z}^{\vec{I}} := Z_{\iota_{(1)}} Z_{\iota_{(2)}} \cdots Z_{\iota_{(N)}}$ 
		denotes the corresponding $N^{th}$ order differential operator.
	\item When we are not concerned with the precise structure of 
		the multi-index $\vec{I},$ we abbreviate $\mathscr{Z}^N := Z_{\iota_{(1)}} Z_{\iota_{(2)}} \cdots Z_{\iota_{(N)}}.$ 
	\item Similarly, $\angLie_{\mathscr{Z}}^N := \angLie_{Z_{\iota_{(1)}}} \angLie_{Z_{\iota_{(2)}}} \cdots \angLie_{Z_{\iota_{(N)}}}$	
		denotes an $N^{th}$ order $S_{t,u}-$projected Lie derivative operator 
		(see Def.~\ref{D:PROJECTEDLIE}).
	\item Similarly, $\angfreeLietwoarg{\mathscr{Z}}{N} := \angfreeLiearg{Z_{\iota_{(1)}}} \angfreeLiearg{Z_{\iota_{(2)}}} \cdots \angfreeLiearg{Z_{\iota_{(N)}}}$	
		denotes an $N^{th}$ order trace-free $S_{t,u}-$projected Lie derivative operator
		(see Def.~\ref{D:ANGFREELIEDEF}).
	\item We use similar notation with $\mathscr{S}$ or $\mathscr{O}$ in place of $\mathscr{Z}$ 
			when all derivatives are spatial derivatives or rotation derivatives 
			(see Def.~\ref{D:DEFSETOFSPATIALCOMMUTATORVECTORFIELDS}).
\end{itemize}
\end{definition}

\begin{definition}[\textbf{Shorthand notation for a pointwise norm of a quantity and its derivatives}]
	We use the following shorthand notation:
	\begin{itemize}
		\item If $f$ is a function, then
			$\left|\mathscr{Z}^{\leq N} f \right|$ denotes any term that is
			$\leq \sum_{|\vec{I}| \leq N} c_{\vec{I}} \left| \mathscr{Z}^{\vec{I}} f \right|,$
			where the $c_{\vec{I}}$ are non-negative constants that verify
			$\sum_{|\vec{I}| \leq N} c_{\vec{I}} \leq 1.$
			Similarly, $\left|\mathscr{Z}^N f \right|$ denotes any term that is
			$\leq \sum_{|\vec{I}| = N} c_{\vec{I}} \left| \mathscr{Z}^{\vec{I}} f \right|,$
			where the $c_{\vec{I}}$ are non-negative constants that verify
			$\sum_{|\vec{I}| = N} c_{\vec{I}} \leq 1.$
			Similarly, if $V$ is a vectorfield, then 
			$\left|V \mathscr{Z}^{\leq N} f \right|$ 
			denotes any term that is
			$\leq \sum_{|\vec{I}| \leq N} c_{\vec{I}} \left| V \mathscr{Z}^{\vec{I}} f \right|,$
			where the $c_{\vec{I}}$ are non-negative constants that verify
			$\sum_{|\vec{I}| \leq N} c_{\vec{I}} \leq 1.$
			Similarly,  
			$\left|V \mathscr{Z}^N f \right|$ 
			denotes any term that is
			$\leq \sum_{|\vec{I}| = N} c_{\vec{I}} \left| V \mathscr{Z}^{\vec{I}} f \right|,$
			where the $c_{\vec{I}}$ are non-negative constants that verify
			$\sum_{\vec{I}} c_{\vec{I}} \leq 1.$
			Similarly,
			$\left|\angdiff \mathscr{Z}^{\leq N} f \right|$ 
			denotes any term that is
			$\leq \sum_{|\vec{I}| \leq N} c_{\vec{I}} \left| \angdiff \mathscr{Z}^{\vec{I}} f \right|,$
			where the $c_{\vec{I}}$ are non-negative constants that verify
			$\sum_{\vec{I}} c_{\vec{I}} \leq 1.$
			We use similar notation for other expressions.
		\item If $\xi$ is an $S_{t,u}$ tensor, then 
			$\left| \angLie_{\mathscr{Z}}^{\leq N} \xi \right|$ 
			denotes any term that is
			$\leq \sum_{|\vec{I}| \leq N} c_{\vec{I}} \left| \angLie_{\mathscr{Z}}^{\vec{I}} \xi \right|,$
			where the $c_{\vec{I}}$ are non-negative constants that verify
			$\sum_{|\vec{I}| \leq N} c_{\vec{I}} \leq 1.$
			Similarly, $\left| \angLie_{\mathscr{Z}}^N \xi \right|$ 
			denotes any term that is
			$\leq \sum_{|\vec{I}| = N} c_{\vec{I}} \left| \angLie_{\mathscr{Z}}^{\vec{I}} \xi \right|,$
			where the $c_{\vec{I}}$ are non-negative constants that verify
			$\sum_{|\vec{I}| = N} c_{\vec{I}} \leq 1.$
			Similarly, if $V$ is a vectorfield, then
			$\left| \angLie_V \angLie_{\mathscr{Z}}^{\leq N} \xi \right|$ 
			denotes any term that is
			$\leq \sum_{|\vec{I}| \leq N} c_{\vec{I}} \left| \angLie_V \angLie_{\mathscr{Z}}^{\vec{I}} \xi \right|,$
			where the $c_{\vec{I}}$ are non-negative constants that verify
			$\sum_{|\vec{I}| \leq N} c_{\vec{I}} \leq 1.$
			Similarly,
			$\left| \angLie_V \angLie_{\mathscr{Z}}^N \xi \right|$ 
			denotes any term that is
			$\leq \sum_{|\vec{I}| = N} c_{\vec{I}} \left| \angLie_V \angLie_{\mathscr{Z}}^{\vec{I}} \xi \right|,$
			where the $c_{\vec{I}}$ are non-negative constants that verify
			$\sum_{|\vec{I}| = N} c_{\vec{I}} \leq 1.$
			We use similar notation for other expressions,
			including the case in which trace-free $S_{t,u}-$projected 
			Lie derivatives $\angfreeLie$ are present instead
			of ordinary $S_{t,u}-$projected 
			Lie derivatives $\angLie.$
		\item We use similar notation with $\mathscr{S}$ or $\mathscr{O}$ in place of $\mathscr{Z}$ 
			when all derivatives are spatial derivatives or rotation derivatives 
			(see Def.~\ref{D:DEFSETOFSPATIALCOMMUTATORVECTORFIELDS}).
	\end{itemize}
\end{definition}


\chapter{The Structure of the Wave Equation Inhomogeneous Terms After One Commutation}
\label{C:ONECOMMUTATION}
\thispagestyle{fancy}
In Chapter~\ref{C:ONECOMMUTATION}, we commute the wave equation $\upmu \square_{g(\Psi)} \Psi = 0$
with a vectorfield $Z \in \mathscr{Z}$ (see Def.~\ref{D:DEFSETOFCOMMUTATORVECTORFIELDS})
and decompose the corresponding inhomogeneous terms 
relative to the frame $\lbrace \Lunit, \Rad, X_1, X_2 \rbrace.$ 
In order to prove our sharp classical lifespan theorem, 
we need to know the precise structure,
including numerical constants, of some of the terms.
Furthermore,
the special properties
of $Z$ lead to some exact cancellations 
that are essential for the proof of the theorem.
The main result is contained in Proposition 
\ref{P:COMMUTATIONCURRENTDIVERGENCEFRAMEDECOMP}.

\section{Preliminary calculations}
We start with the following lemma, which provides a preliminary decomposition
of the commutator term $[Z, \upmu \square_{g(\Psi)}] \Psi.$ 

\begin{lemma}[\textbf{Preliminary wave equation commutator expression}] 
 \label{L:BOXZCOM}
	Let $Z$ be any spacetime vectorfield, and let $\deformarg{Z}{\mu}{\nu} := \D_{\mu} Z_{\nu} + \D_{\mu} Z_{\nu}$
	be its deformation tensor. 
	Assume that $\deformarg{Z}{\Lunit}{\Lunit} := \deformarg{Z}{\alpha}{\beta} \Lunit^{\alpha} \Lunit^{\beta} = 0.$ 
	Then the following commutation identity holds:
	\begin{align} \label{E:BOXZCOM}
		\upmu \square_{g(\Psi)} (Z \Psi)
		& =   \upmu 
					\D_{\alpha} 
					\left(\deformuparg{Z}{\alpha}{\beta} \D_{\beta} \Psi
						- \frac{1}{2} \myspacetimetr  \deform{Z} \D^{\alpha} \Psi 
					\right)
						\\
	& \ \ + Z(\upmu \square_{g(\Psi)} \Psi)
				\underbrace{- \upmu^{-1} 
											\left\lbrace 
												\deformarg{Z}{\Lunit}{\Rad}
												+ Z \upmu 
											\right\rbrace
									  	(\upmu \square_{g(\Psi)} \Psi)
									  }_{\mbox{potentially dangerous factor}} 
				+ \frac{1}{2} \mytr  \angdeform{Z} (\upmu \square_{g(\Psi)} \Psi).
				\notag
	\end{align}
	
\end{lemma}

\begin{proof}
	Using the Leibniz rule, we compute that 
	\begin{align} \label{E:FIRSTBOXCOMMUTATORCOMPUTATION}
		\upmu \square_{g(\Psi)} (Z \Psi) 
			& =	Z(\upmu \square_{g(\Psi)} \Psi)
				- (Z \upmu) \square_{g(\Psi)} \Psi
			 	+ \upmu Z^{\alpha} \D^{\beta} \D_{\beta} \D_{\alpha} \Psi 
				- \upmu Z^{\alpha} \D_{\alpha} \D^{\beta} \D_{\beta} \Psi
				+ \upmu \D^{\alpha} Z^{\beta} \D_{\alpha} \D_{\beta} \Psi
					\\
			& \ \ + \frac{1}{2} \upmu \D_{\alpha} \pi^{\alpha \beta} \D_{\beta} \Psi
				+ \frac{1}{2} \upmu \D^{\alpha} (\D_{\alpha} Z^{\beta} - \D^{\beta} Z_{\alpha}) \D_{\beta} \Psi.
				\notag
		\end{align}
	Furthermore, with $\Cur_{\alpha \beta} := (g^{-1})^{\kappa \lambda}
	\Cur_{\alpha \kappa \beta \lambda}$ denoting the Ricci curvature tensor of $g$
	(our sign conventions for the Riemann curvature tensor $\Cur_{\alpha \kappa \beta \lambda}$ are given in Def.~\ref{D:SPACETIMERIEMANN}),
	we compute that the following commutator identities hold:
	\begin{align}
		Z^{\alpha} \D^{\beta} \D_{\beta} \D_{\alpha} \Psi 
			- Z^{\alpha} \D_{\alpha} \D^{\beta} \D_{\beta} \Psi
		& = \Cur_{\alpha}^{\ \beta} Z^{\alpha} \D_{\beta} \Psi, 
			\label{E:FIRSTBOXRICCIRELATION} \\
		\frac{1}{2} \D^{\alpha} (\D_{\alpha} Z^{\beta} - \D^{\beta} Z_{\alpha}) \D_{\beta} \Psi
		& = \frac{1}{2} \D_{\alpha} \deformuparg{Z}{\alpha}{\beta} \D_{\beta} \Psi
			- \frac{1}{2} \D^{\alpha} \deformmixedarg{Z}{\beta}{\beta} \D_{\alpha} \Psi
			- \Cur_{\alpha}^{\ \beta} Z^{\alpha} \D_{\beta} \Psi.
			\label{E:SECONDBOXRICCIRELATION}
	\end{align}
	The desired identity \eqref{E:BOXZCOM} now follows from 
	\eqref{E:FIRSTBOXCOMMUTATORCOMPUTATION},
	\eqref{E:FIRSTBOXRICCIRELATION},
	\eqref{E:SECONDBOXRICCIRELATION},
	Lemma~\ref{L:SPACETIMEMETRICFRAMEVECTORFIELDS},
	the decomposition
	$\frac{1}{2} \myspacetimetr  \deform{Z} = - \upmu^{-1} \deformarg{Z}{\Lunit}{\Rad} + \frac{1}{2} \mytr  \angdeform{Z},$
	and from straightforward computations.
	
\end{proof}

\begin{remark}[\textbf{Some good properties of the commutation vectorfields}]
	\label{R:POTENTIALLYDANGEROUSTERMISNOTDANGEROUS}
	For commutation vectorfields $Z \in \mathscr{Z}$ (see Def.~\ref{D:DEFSETOFCOMMUTATORVECTORFIELDS}),
	we have $\deformarg{Z}{\Lunit}{\Lunit} = 0.$
	Furthermore, the ``potentially dangerous factor'' on the right-hand side of \eqref{E:BOXZCOM} contains no powers
	of $\upmu^{-1}$ and therefore, despite superficial appearances, is not actually dangerous. These properties of $\mathscr{Z}$
	are essential for closing our $L^2$ estimates.
	Specifically, Prop.~\ref{P:DEFORMATIONTENSORFRAMECOMPONENTS} implies that for each $Z \in \mathscr{Z},$
	$\upmu^{-1} \left\lbrace \deform{Z}_{\Lunit \Rad} + Z \upmu \right\rbrace$
	is equal to either $0$ or $-1.$ 
	 
\end{remark}

We now define a family of commutation current vectorfields that will facilitate our analysis of 
the structure of the inhomogeneous terms on the right-hand side of \eqref{E:BOXZCOM}.

\begin{definition}[\textbf{Commutation currents}]
\label{D:COMMUTATIONCURRENT}
Let $Z$ be any vectorfield, and let $\deformarg{Z}{\mu}{\nu} := \D_{\mu} Z_{\nu} + \D_{\mu} Z_{\nu}$
be its deformation tensor. We define the corresponding 
commutation current (vectorfield)
${\Jcurrent{Z}^{\alpha}[\Psi]}$ as follows:
	\begin{align} \label{E:COMMUTATIONCURRENT}
		{\Jcurrent{Z}^{\alpha}[\Psi]}
			:= \deformuparg{Z}{\alpha}{\beta} \D_{\beta} \Psi 
			- \frac{1}{2} \myspacetimetr  \deform{Z} \D^{\alpha} \Psi.
	\end{align}

\end{definition}

We now use commutation currents to provide an alternate description of the 
inhomogeneous terms on the right-hand side of \eqref{E:BOXZCOM}. 
The advantage of the alternate description is that it is easy
to decompose the corresponding expression relative to the frame
$\lbrace \Lunit, \Rad, X_1, X_2 \rbrace.$

\begin{lemma}[\textbf{Basic structure of the inhomogeneous terms in the once-commuted wave equation}]
	\label{L:WAVEONCECOMMUTEDBASICSTRUCTURE}
	Given any $Z \in \mathscr{Z}$ 
	(note that $\deformarg{Z}{\Lunit}{\Lunit} = 0$),
	let $\Jcurrent{Z}[\Psi]$ be the commutation current
	\eqref{E:COMMUTATIONCURRENT}. Then $Z \Psi$ verifies the inhomogeneous wave equation
	\begin{align} \label{E:WAVEONCECOMMUTEDBASICSTRUCTURE}
		\upmu \square_{g(\Psi)} (Z \Psi)
		& = \upmu \D_{\alpha} {\Jcurrent{Z}^{\alpha}[\Psi]}
			+ Z(\upmu \square_{g(\Psi)} \Psi)
				\underbrace{- \upmu^{-1} 
											\left\lbrace 
												\pi_{\Lunit \Rad}
												+ Z \upmu 
											\right\rbrace
									  	(\upmu \square_{g(\Psi)} \Psi)
									  }_{\mbox{potentially dangerous factor}} 
				+ \frac{1}{2} \mytr  \angdeform{Z} (\upmu \square_{g(\Psi)} \Psi).
\end{align}
\end{lemma}

\begin{proof}
	The lemma follows easily from applying $\D_{\alpha}$ to the right-hand side of \eqref{E:COMMUTATIONCURRENT}
	and comparing the resulting expression to the right-hand side of \eqref{E:BOXZCOM}.
\end{proof}

Lemma~\ref{L:WAVEONCECOMMUTEDBASICSTRUCTURE} shows that in order to estimate the inhomogeneous terms in the wave
equation verified by $Z \Psi,$ we have to analyze the structure of the derivatives of the commutation
currents $\upmu \D_{\alpha} {\Jcurrent{Z}^{\alpha}[\Psi]}.$ As a first step, 
in Prop.~\ref{P:COMMUTATIONCURRENTDIVERGENCEFRAMEDECOMP},
we decompose $\upmu {\D_{\alpha} \Jcurrent{Z}^{\alpha}[\Psi]}$ relative to the frame
$\lbrace \Lunit, \Rad, X_1, X_2 \rbrace.$ We use the following preliminary lemma in our proof of
the proposition.

\begin{lemma} [\textbf{The divergence of $\mathscr{J}$ in terms of rescaled frame derivatives}]
\label{L:DIVERGENCEFRAME}
Let $\mathscr{J}$ be any spacetime vectorfield, and consider its decomposition
relative to the frame $\lbrace \Lunit, \Rad, X_1, X_2 \rbrace$
(see Remark~\ref{R:UPDOWNFRAMETRANSLATION}):
\begin{align} \label{E:JFRAMEDECOMP}
		\mathscr{J}
		& = 
			- \mathscr{J}_{\Lunit} \Lunit
			- \upmu^{-1} \mathscr{J}_{\Rad} \Lunit
			- \upmu^{-1} \mathscr{J}_{\Lunit} \Rad
			+ \angJ,
\end{align}
where $\angJ$ is the projection of $\mathscr{J}$ onto the $S_{t,u}.$
Then the following $\upmu-$weighted spacetime covariant divergence identity holds:
\begin{align} \label{E:DIVERGENCEFRAME}
	\upmu \D_{\alpha} \mathscr{J}^{\alpha} 
	& = - \Lunit (\upmu \mathscr{J}_{\Lunit})
		- \Lunit (\mathscr{J}_{\Rad})
		- \Rad (\mathscr{J}_{\Lunit})
		+ \angdiv (\upmu \angJ)
			\\
	& \ \ - \left\lbrace
						\mytr  \angkuparg{(Trans-\Psi)}
						+ \upmu \mytr \angkuparg{(Tan-\Psi)}
					\right\rbrace
					\mathscr{J}_{\Lunit}
		- \mytr \upchi \mathscr{J}_{\Rad}.
		\notag
\end{align}	
The $S_{t,u}$ tensorfields 
$\upchi,$
$\angkuparg{(Trans-\Psi)},$
and
$\angkuparg{(Tan-\Psi)}$
appearing in \eqref{E:DIVERGENCEFRAME}
are defined in 
\eqref{E:CHIDEF},
\eqref{E:KABTRANSVERSAL},
and \eqref{E:KABGOOD}.
\end{lemma}

\begin{proof}
We fix a hypersurface region $\Sigma_t^{U_0}.$ We will show that
\eqref{E:DIVERGENCEFRAME} holds along $\Sigma_t^{U_0}.$
Since the identity \eqref{E:DIVERGENCEFRAME} is invariant
with respect to changes of angular coordinates of the form
$\widetilde{\vartheta}^1 = f^1(u,\vartheta^1, \vartheta^2),$
$\widetilde{\vartheta}^2 = f^2(u,\vartheta^1, \vartheta^2)$
(where the $f^A$ are the coordinate transformation functions),
it suffices to show that
\eqref{E:SPACETIMEVOLUMEFORMWITHUPMU}
holds for a well-chosen version 
$(\widetilde{\vartheta}^1, \widetilde{\vartheta}^2)$
of such angular coordinates.
To this end, we use the angular coordinates constructed in the proof of
Cor.~\ref{C:SPACETIMEVOLUMEFORMWITHUPMU}, which are such that
along $\Sigma_t^{U_0},$ we have
$\Rad = \frac{\partial}{\partial u}|_{t, \widetilde{\vartheta}^1, \widetilde{\vartheta}^2}.$
By transporting $\widetilde{\vartheta}^1, \widetilde{\vartheta}^2$ via the ODEs
$\Lunit \widetilde{\vartheta}^1 = \Lunit \widetilde{\vartheta}^2 = 0,$
we may assume that $\Lunit = \frac{\partial}{\partial t}|_{u,\widetilde{\vartheta}^1, \widetilde{\vartheta}^2}$
along $\Sigma_t^{U_0}.$ 
We now expand $\mathscr{J}$ relative to the coordinate frame 
$\lbrace \frac{\partial}{\partial t}, \frac{\partial}{\partial u}, 
\frac{\partial}{\partial \widetilde{\vartheta}^1}, \frac{\partial}{\partial \widetilde{\vartheta}^2} \rbrace$
	as follows:
		$\mathscr{J} = \mathscr{J}^t \frac{\partial}{\partial t} 
		+ \mathscr{J}^u \frac{\partial}{\partial u}
		+ \angJ^A \frac{\partial}{\partial \widetilde{\vartheta}^A},$
		where $\angJ$ is $S_{t,u}-$tangent.
		Furthermore, 
		from the identities 
		$g(\Lunit, \Lunit) 
		= g(\Lunit, \frac{\partial}{\partial \widetilde{\vartheta}^A}) 
		= g(\Rad, \frac{\partial}{\partial \widetilde{\vartheta}^A})
		= 0,$
		$g(\Lunit, \Rad) = - \upmu,$
		and
		$g(\Rad, \frac{\partial}{\partial u}) = \upmu^2,$
		it is straightforward to compute that
		\begin{align} \label{E:GENERALDIVERGENCERELATIONCOORDINATEVECTROFIELDCOMPONENTSINTERMSOFFRAMECOMPONENTS}
			\mathscr{J}^t & = - \mathscr{J}_{\Lunit} - \upmu^{-1} \mathscr{J}_{\Rad}, \qquad \mathscr{J}^u = - \upmu^{-1} \mathscr{J}_{\Lunit}.
		\end{align}
		We now recall the following standard formula for 
		the divergence $\D_{\alpha} \mathscr{J}^{\alpha}$ of $\mathscr{J}$
		relative to the coordinates $(t,u,\widetilde{\vartheta}^1,\widetilde{\vartheta}^2):$
		\begin{align} \label{E:COORDINATEDIVERGENCEEXPRESSION1}
			\D_{\alpha} \mathscr{J}^{\alpha}
			=
			\frac{1}{\sqrt{|\mbox{\upshape{det}} g|}}
			\left\lbrace
				\frac{\partial}{\partial t} (\sqrt{|\mbox{\upshape{det}} g|} \mathscr{J}^t)
				+ \frac{\partial}{\partial u} (\sqrt{|\mbox{\upshape{det}} g|} \mathscr{J}^u)
				+ \frac{\partial}{\partial \widetilde{\vartheta}^A} (\sqrt{|\mbox{\upshape{det}} g|} \angJ^A)
			\right\rbrace.
		\end{align}
		In \eqref{E:COORDINATEDIVERGENCEEXPRESSION1},
		$\mbox{\upshape{det}} g$ is taken relative to the coordinates
		$(t,u,\widetilde{\vartheta}^1,\widetilde{\vartheta}^2).$
		We next note that the identity \eqref{E:SPACETIMEVOLUMEFORMWITHUPMU} 
		holds relative to the coordinates $(t,u,\widetilde{\vartheta}^1,\widetilde{\vartheta}^2).$
		Inserting the square root of the identity 
		into \eqref{E:COORDINATEDIVERGENCEEXPRESSION1},
		using the fact that
		$\frac{1}{\sqrt{\mbox{$|$\upshape{det}$\gsphere|$}}}
			\frac{\partial}{\partial \widetilde{\vartheta}^A} 
			(\upmu \sqrt{|\mbox{\upshape{det}} \gsphere|} \angJ^A) = \angdiv (\upmu \angJ),$
		and multiplying the resulting equation by $\upmu,$
		we deduce that
		\begin{align} \label{E:COORDINATEDIVERGENCESECONDEXPRESSION2}
			\upmu \D_{\alpha} \mathscr{J}^{\alpha}
			& = 
				\frac{\partial}{\partial t} (\upmu \mathscr{J}^t)
				+ \frac{\partial}{\partial u} (\upmu \mathscr{J}^u)
				+ \angdiv (\upmu \angJ)
					\\
			& \ \ 
				+  	\left(
							\frac{\partial}{\partial t} \ln \sqrt{\mbox{\upshape{det}} \gsphere} 
						\right)
						\upmu \mathscr{J}^t
				+ \left(
						\frac{\partial}{\partial u} \ln \sqrt{\mbox{\upshape{det}} \gsphere}
					\right)
					\upmu \mathscr{J}^u.
					\notag			
		\end{align}
		Next, from the standard matrix identity 
		$\frac{\partial}{\partial t} \ln \sqrt{\mbox{\upshape{det}} \gsphere} = 
		\frac{1}{2} (\ginversesphere)^{AB} \frac{\partial}{\partial t} \gsphere_{AB}$
		and the fact that $\frac{\partial}{\partial t} \gsphere_{AB} = \angLie_{\Lunit} \gsphere_{AB} 
		= \angdeformarg{\Lunit}{A}{B} = 2 \upchi_{AB},$ we deduce that 
		$\frac{\partial}{\partial t} \ln \sqrt{\mbox{\upshape{det}} \gsphere} = \mytr \upchi.$
		Similarly, with the help of \eqref{E:RADDEFORMANG}, 
		we deduce that 
		$\frac{\partial}{\partial u} \ln \sqrt{\mbox{\upshape{det}} \gsphere} 
			= \frac{1}{2} \mytr  \angdeform{\Rad}
			=	- \upmu \mytr \upchi
				+ \mytr  \angkuparg{(Trans-\Psi)}
				+ \upmu \mytr  \angkuparg{(Tan-\Psi)}.$
		Inserting these identities, 
		the identities $\frac{\partial}{\partial t}|_{u,\widetilde{\vartheta}^1,\widetilde{\vartheta}^2} = \Lunit$
		and $\frac{\partial}{\partial u}|_{t,\widetilde{\vartheta}^1,\widetilde{\vartheta}^2} = \Rad$
		(along $\Sigma_t^{U_0}$),
		and \eqref{E:GENERALDIVERGENCERELATIONCOORDINATEVECTROFIELDCOMPONENTSINTERMSOFFRAMECOMPONENTS}
		into \eqref{E:COORDINATEDIVERGENCESECONDEXPRESSION2}, we arrive at \eqref{E:DIVERGENCEFRAME}.

\end{proof}

\section{Frame decomposition of \texorpdfstring{$\upmu \D_{\alpha} {\Jcurrent{Z}^{\alpha}[\Psi]}$}{the commutation current}}
We now use Lemma~\ref{L:DIVERGENCEFRAME} to derive a detailed
frame decomposition of the term 
$\upmu \D_{\alpha} {\Jcurrent{Z}^{\alpha}[\Psi]}$
on the right-hand side of \eqref{E:WAVEONCECOMMUTEDBASICSTRUCTURE}.
When we derive a priori $L^2$ estimates, 
the top-order derivatives of this term are
by far the most difficult quantities that we have to bound.
In fact, in order to close our top-order $L^2$ estimates, 
we must understand the precise structure 
of some of the terms in $\upmu \D_{\alpha} {\Jcurrent{Z}^{\alpha}[\Psi]}$
including the exact numerical constants that appear. In order to
close our below-top-order $L^2$ estimates, 
we need to know much less about their structure.

\begin{proposition} [\textbf{Frame decomposition of the divergence of the commutation current}]
\label{P:COMMUTATIONCURRENTDIVERGENCEFRAMEDECOMP}
	Let $Z$ be any vectorfield verifying $\deformarg{Z}{\Lunit}{\Lunit}=0,$
	and let ${\Jcurrent{Z}[\Psi]}$
	be the commutation current \eqref{E:COMMUTATIONCURRENT} associated to $Z.$
	Then
	\begin{align} \label{E:DIVCOMMUTATIONCURRENTDECOMPOSITION}
		\upmu \D_{\alpha} {\Jcurrent{Z}^{\alpha}[\Psi]}
		& = \mathscr{K}_{(\pi-Danger)}^{(Z)}[\Psi]
			+ \mathscr{K}_{(\pi-Cancel-1)}^{(Z)}[\Psi]
			+ \mathscr{K}_{(\pi-Cancel-2)}^{(Z)}[\Psi]
			+ \mathscr{K}_{(\pi-Elliptic)}^{(Z)}[\Psi]
				\\
		& \ \ 
			+ \mathscr{K}_{(\pi-Good)}^{(Z)}[\Psi]
			+ \mathscr{K}_{(\Psi)}^{(Z)}[\Psi]
			+ \mathscr{K}_{(Low)}^{(Z)}[\Psi],
			\notag
	\end{align}
	where
	\begin{subequations}	
		\begin{align}
			\mathscr{K}_{(\pi-Danger)}^{(Z)}[\Psi]
			& := - (\angdiv \angdeformoneformupsharparg{Z}{\Lunit}) \Rad \Psi,
				\label{E:DIVCURRENTTRANSVERSAL}
				\\
			\mathscr{K}_{(\pi-Cancel-1)}^{(Z)}[\Psi]
			& := \left\lbrace 
						\frac{1}{2} \Rad \mytr  \angdeform{Z}
						- \angdiv \angdeformoneformupsharparg{Z}{\Rad}
						- \upmu \angdiv \angdeformoneformupsharparg{Z}{\Lunit}
					\right\rbrace 
						\Lunit \Psi,
					\label{E:DIVCURRENTCANEL1} \\
			\mathscr{K}_{(\pi-Cancel-2)}^{(Z)}[\Psi]
			& :=	
				\left\lbrace
					- \angLie_{\Rad} \angdeformoneformupsharparg{Z}{\Lunit}
					+ \angdiffuparg{\#} \deformarg{Z}{\Lunit}{\Rad}
				\right\rbrace 
				\cdot
				\angdiff \Psi,
				\label{E:DIVCURRENTCANEL2} \\
			\mathscr{K}_{(\pi-Elliptic)}^{(Z)}[\Psi]
			& := \upmu (\angdiv \angdeformfreeupdoublesharparg{Z}) \cdot \angdiff \Psi, 
				\label{E:DIVCURRENTELLIPTIC} \\
			\mathscr{K}_{(\pi-Good)}^{(Z)}[\Psi] 
			& := \frac{1}{2} \upmu (\Lunit \mytr  \angdeform{Z}) \Lunit \Psi
				+ (\Lunit \deformarg{Z}{\Lunit}{\Rad}) \Lunit \Psi
				+ (\Lunit \deformarg{Z}{\Rad}{\Radunit}) \Lunit \Psi
				\label{E:DIVCURRENTGOOD} \\
			& \ \ + \frac{1}{2} (\Lunit \mytr  \angdeform{Z}) \Rad \Psi
				- \upmu (\angLie_{\Lunit} \angdeformoneformupsharparg{Z}{\Lunit}) \cdot \angdiff \Psi
				- (\angLie_{\Lunit} \angdeformoneformupsharparg{Z}{\Rad}) \cdot \angdiff \Psi,
				\notag 
	\end{align}
	\end{subequations}
	\begin{align}
		\mathscr{K}_{(\Psi)}^{(Z)}[\Psi]
			& := \rgeo^{-1} 
						\left\lbrace
							\frac{1}{2}
							\upmu \mytr  \angdeform{Z}
							+ \deformarg{Z}{\Lunit}{\Rad}
							+ \deformarg{Z}{\Rad}{\Radunit}
						\right\rbrace
						\left\lbrace \Lunit (\rgeo \Lunit \Psi)
							+ \frac{1}{2} \mytr \upchi (\rgeo \Lunit \Psi)
						\right\rbrace  
					\label{E:DIVCURRENTPSI} \\
		& \ \
				+ \mytr  \angdeform{Z}
					\left\lbrace 
						\Lunit \Rad \Psi
						+ \frac{1}{2} \mytr \upchi \Rad \Psi
					\right\rbrace 
					\notag \\
			& \ \ - 2 \rgeo^{-1} \upmu \angdeformoneformupsharparg{Z}{\Lunit} \cdot \angdiff (\rgeo \Lunit \Psi)
					- 2 \angdeformoneformupsharparg{Z}{\Lunit} \cdot \angdiff \Rad \Psi
					- 2 \rgeo^{-1} \angdeformoneformupsharparg{Z}{\Rad} \cdot \angdiff (\rgeo \Lunit \Psi)
				\notag \\
			& \ \ + \deformarg{Z}{\Lunit}{\Rad} \angLap \Psi
				+ \upmu \angdeformfreeupdoublesharparg{Z} \cdot \angfreeDsquared \Psi,
				\notag 
	\end{align}
	and
	\begin{align}
		\mathscr{K}_{(Low)}^{(Z)}[\Psi] 
		& :=  \frac{1}{2} \mytr \upchi^{(Small)} 
						\left\lbrace
							\deformarg{Z}{\Lunit}{\Rad} 
							+ \deformarg{Z}{\Rad}{\Radunit}
							- \frac{1}{2} \upmu \mytr \angdeform{Z}
						\right\rbrace 
						\Lunit \Psi
				\label{E:DIVCURRENTLOW} \\
		& \ \ + \frac{1}{2} 
						\left\lbrace
							\Lunit \upmu
							+ \mytr  \angkuparg{(Trans-\Psi)} 
							+ \upmu \mytr  \angkuparg{(Tan-\Psi)} 
							- 2 \rgeo^{-1} \upmu
						\right\rbrace 
						\mytr  \angdeform{Z} \Lunit \Psi 
						- (\angdeformoneformupsharparg{Z}{\Lunit} \cdot \angdiff \upmu) 
							\Lunit \Psi
				\notag \\
		& \ \	-  
					\left\lbrace
						\Lunit \upmu
						+ \mytr  \angkuparg{(Trans-\Psi)}
						+ \upmu \mytr  \angkuparg{(Tan-\Psi)}
					\right\rbrace
					\angdeformoneformupsharparg{Z}{\Lunit} 
					\cdot
					\angdiff \Psi
					\notag \\
		& \ \ 
				- \left\lbrace2 \rgeo^{-1} + \mytr \upchi^{(Small)}\right\rbrace 
						\angdeformoneformupsharparg{Z}{\Rad} \cdot \angdiff \Psi
					+ \frac{1}{2} 
						\mytr  \angdeform{Z}
						(\angdiffuparg{\#} \upmu) 
						\cdot
						\angdiff \Psi
			\notag \\
		& \ \ + \mytr \angdeform{Z} 
						\left\lbrace
							\upzeta^{(Trans-\Psi)\#}
							+ \upmu \upzeta^{(Tan-\Psi)\#}
						\right\rbrace 
						\cdot
						\angdiff \Psi
					+ (\angdiff \upmu) 
						\cdot
						\angdeformfreeupdoublesharparg{Z} 
						\cdot
						\angdiff \Psi.
			\notag
	\end{align}
	In the above expressions, the $S_{t,u}$ tensorfields
	$\upchi,$ 
	$\upchi^{(Small)},$
	$\upzeta^{(Trans-\Psi)},$
	$\angkuparg{(Trans-\Psi)},$ 
	$\upzeta^{(Tan-\Psi)},$
	and
	$\angkuparg{(Tan-\Psi)}$
	are defined by
	\eqref{E:CHIDEF},
	\eqref{E:CHIJUNKDEF},
	\eqref{E:ZETATRANSVERSAL},
	\eqref{E:KABTRANSVERSAL},
	\eqref{E:ZETAGOOD},
	and \eqref{E:KABGOOD}.
	The type $\binom{0}{2}$ $S_{t,u}$ tensor 
	$\angfreeDsquared \Psi$ is defined by \eqref{E:ANGFREEDDEF}.
	\end{proposition}

\begin{remark}[\textbf{The importance of} $\deformarg{Z}{\Lunit}{\Lunit} = 0$] \label{R:NOPILL}
	The assumption $\deformarg{Z}{\Lunit}{\Lunit} = 0$
	implies the absence of the dangerous quadratic term 
	mentioned in Remark~\ref{R:DEFTENSORCALCULATIONSIMPORTANTASPECTS}.
\end{remark}

\begin{remark}[\textbf{Favorable combinations}]
	Note that we have carefully combined various terms 
	to create sums of the form
	$\Lunit Z \Psi + \frac{1}{2} \mytr \upchi Z \Psi$
	(where $Z \in \mathscr{Z}$)
	and placed them on the right-hand side of in \eqref{E:DIVCURRENTPSI}.
	These combinations have better
	$t-$decay properties than the individual summands.
\end{remark}

\begin{proof}
	With the help of 
	Lemma~\ref{L:SPACETIMEMETRICFRAMEVECTORFIELDS}
	and the decomposition $\frac{1}{2} \myspacetimetr  \deform{Z} = - \upmu^{-1} \deformarg{Z}{\Lunit}{\Rad} + \frac{1}{2} \mytr  \angdeform{Z},$
	we compute that the components of $\Jcurrent{Z}[\Psi]$ relative
	to the frame $\lbrace \Lunit, \Rad, X_1, X_2 \rbrace$ are:
	\begin{align} \label{E:COMMUTATIONJLFRAMEDECOMPOSED}
		\Jcurrent{Z}_{\Lunit}[\Psi] 
		& := g(\Jcurrent{Z}, \Lunit)
			= - \frac{1}{2} \mytr  \angdeform{Z} \Lunit \Psi
			+ \angdeformoneformupsharparg{Z}{\Lunit} \cdot \angdiff \Psi, 
				\\
		\Jcurrent{Z}_{\Rad}[\Psi] 
		& := g(\Jcurrent{Z}, \Rad)		
			= - \deformarg{Z}{\Lunit}{\Rad} \Lunit \Psi 
			- \deformarg{Z}{\Rad}{\Radunit} \Lunit \Psi 
			+ \angdeformoneformupsharparg{Z}{\Rad} \cdot \angdiff \Psi
			- \frac{1}{2} \mytr  \angdeform{Z} \Rad \Psi,
			\label{E:COMMUTATIONJRADFRAMEDECOMPOSED}\\
		\upmu \angJ_A[\Psi]
		& := \upmu g(\Jcurrent{Z}, X_A)
			= - \upmu \angdeformarg{Z}{\Lunit}{A} \Lunit \Psi
			- \angdeformarg{Z}{\Rad}{A} \Lunit \Psi
			- \angdeformarg{Z}{\Lunit}{A} \Rad \Psi
			+ \deformarg{Z}{\Lunit}{\Rad} \angdiffarg{A} \Psi
			+  \upmu \angdeformfreemixedarg{Z}{A}{B} \angdiffarg{B} \Psi.
			\label{E:COMMUTATIONJANGFRAMEDECOMPOSED}
	\end{align}
	The identity \eqref{E:DIVCOMMUTATIONCURRENTDECOMPOSITION} now follows from applying Lemma~\ref{L:DIVERGENCEFRAME}
	to the frame components \eqref{E:COMMUTATIONJLFRAMEDECOMPOSED}-\eqref{E:COMMUTATIONJANGFRAMEDECOMPOSED}
	and carrying out straightforward computations with the help of
	the commutator identity \eqref{E:LRADCOMM},
	Lemma~\ref{L:LANDRADCOMMUTEWITHANGDIFF},
	the identity $\mytr \upchi = 2 \rgeo^{-1} + \mytr \upchi^{(Small)},$
	and the identity $\Lunit \rgeo = 1.$
\end{proof}

\begin{remark}[\textbf{Explanation of the various types of error terms}]
	We now make some remarks concerning the terms on the right-hand side of \eqref{E:DIVCOMMUTATIONCURRENTDECOMPOSITION}
	for commutation vectorfields $Z \in \mathscr{Z}.$
	When we derive top-order $L^2$ estimates for the error terms in the commuted wave equation,
	we have to be very careful in how we estimate the top-order derivatives of 
	$\mathscr{K}_{(\pi-Danger)}^{(Z)}[\Psi],$
	$\mathscr{K}_{(\pi-Cancel-1)}^{(Z)}[\Psi],$
	$\mathscr{K}_{(\pi-Cancel-2)}^{(Z)}[\Psi],$
	$\mathscr{K}_{(\pi-Elliptic)}^{(Z)}[\Psi],$
	and $\mathscr{K}_{(\pi-Good)}^{(Z)}[\Psi],$
	for it is precisely these terms that lead to the presence of the top-order derivatives of the eikonal function quantities;
	as we outlined in Sect.~\ref{SS:INTROTOPORDEERENERGYESTIMATES},
	such terms are difficult to estimate in $L^2$ and lead to $\upmu-$degenerate top-order $L^2$ estimates.
	\begin{itemize}
		\item The most difficult terms are generated by the top-order derivatives of $\mathscr{K}_{(\pi-Danger)}^{(Z)}[\Psi].$
			The reason is that the factor $\Rad \Psi$ decays at the non-integrable rate 
			$\varepsilon (1 + t)^{-1}$ (see Sect.~\ref{SS:PSIBOOTSTRAP}).
			Actually, in closing our the top-order $L^2$ estimates for this term,
			we must find some special structure that allows us to connect the
			product $G_{\Lunit \Lunit} \Rad \Psi$ back to $\Lunit \upmu$
			via the transport equation \eqref{E:UPMUFIRSTTRANSPORT}.
		\item We must exploit some critically important algebraic cancellation that occurs in
			$\mathscr{K}_{(\pi-Cancel-1)}^{(Z)}[\Psi]$
			and $\mathscr{K}_{(\pi-Cancel-2)}^{(Z)}[\Psi];$
			see just below equation \eqref{E:PICANCEL2HARMLESS}.
			The cancellation is based in part on the special structure of the deformation tensors of
			the vectorfields belonging to the commutation set $\mathscr{Z}.$
			In particular, if not for the cancellation present in $\mathscr{K}_{(\pi-Cancel-2)}^{(Z)}[\Psi],$
			our error terms would involve certain top-order derivatives of the eikonal function
			quantities that we would have no means of controlling (that is, our estimates would lose derivatives);
			see Remark~\ref{R:AVOIDDERIVATIVELOSS}.
	\item The top-order derivatives of $\mathscr{K}_{(\pi-Elliptic)}^{(Z)}[\Psi]$ 
			are relatively easy to bound, even though some of the terms must be treated with elliptic estimates.
	\item The top-order derivatives of $\mathscr{K}_{(\pi-Good)}^{(Z)}[\Psi]$
		are relatively easy to bound because all of the top-order derivatives
		of the eikonal function quantities contain at least one $\Lunit$ differentiation;
		for such terms, there is no danger of losing derivatives.
	\item The top-order derivatives of $\mathscr{K}_{(\Psi)}^{(Z)}[\Psi]$
			generate terms involving the top-order derivatives of $\Psi,$ but
			\emph{not of the eikonal function quantities}; such terms are relatively easy to bound.
	\item The term $\mathscr{K}_{(Low)}^{(Z)}[\Psi]$ consists of products that are lower-order
		in the sense of number of derivatives. Hence, $\mathscr{K}_{(Low)}^{(Z)}[\Psi]$ is 
		relatively easy to bound at all derivative levels.
	\end{itemize}
\end{remark}


\chapter{Energy and Cone-Flux Definitions and the Fundamental Divergence Identities} 
\label{C:DIVERGENCETHEOREM}
\thispagestyle{fancy}
In Chapter~\ref{C:DIVERGENCETHEOREM}, we define energies and fluxes, which are $L^2-$type
quantities that we use to control $\Psi.$ 
We then use the vectorfield multiplier method 
to derive divergence theorem-type identities verified by the energies and fluxes. 
To derive the identities, we use
two distinct vectorfield multipliers together with some modifications needed to handle
the lower-order terms. 
We use the identities in Chapter~\ref{C:ERRORTERMSOBOLEV} in order to 
derive, via a lengthy Gronwall argument, a priori
$L^2$ estimates for $\Psi$ and its up-to-top-order derivatives.

\section{Preliminary calculations}

We start with the following standard lemma, which we often use in our analysis.

\begin{remark} [\textbf{The size of the $S_{t,u}$ area form}]
\label{R:SPHEREAREAFORMIMPLICITFACTOR}
	Throughout this monograph,
	the geometric spherical area form $d \spherevol$ on $S_{t,u}$ implicitly contains a factor of
	magnitude $\approx \rgeo^2$ (see inequality \eqref{E:SPHEREVOLUMEFORMCOMPARISON}) compared to the standard
	Euclidean area form on the Euclidean-unit sphere $(\mathbb{S}^2,\Eucsphereunit),$
	where $\Eucsphereunit$ is the Riemannian metric
	induced on the Euclidean unit sphere
	by the Euclidean metric $\Euct$ on $\Sigma_t.$
	Note that relative to the rectangular spatial coordinates on $\Sigma_t,$ 
	we have $\Euct_{ij} = \delta_{ij}.$
\end{remark}

\begin{lemma}[\textbf{Derivatives of $S_{t,u}$ integrals}]
	\label{L:STUDERIVATIVES}
	We have the following identities:	
			\begin{subequations}
			\begin{align}
				\frac{\partial}{\partial t} 
				\left(
					\int_{S_{t,u}} f \, d \spherevol
				\right)
				& = \int_{S_{t,u}} 
							\Lunit f 
							+ \mytr \upchi f
						\, d \spherevol,
				\label{E:OUTGOINGDERIVSTU} \\
			\frac{\partial}{\partial u} 
				\left(
					\int_{S_{t,u}} f \, d \spherevol
				\right)
				& = \int_{S_{t,u}} 
						\Rad f 
						+ \left\lbrace
								- \frac{2}{\rgeo} 
								- \frac{2}{\rgeo} (\upmu - 1)
								- \upmu \mytr \upchi^{(Small)}
								+ \mytr  \angkuparg{(Trans-\Psi)}
								+ \upmu \mytr  \angkuparg{(Tan-\Psi)}
							\right\rbrace
							f
					\, d \spherevol.
					\label{E:UDERIVSTU}
			\end{align}
			\end{subequations}
			In the above expressions, the $S_{t,u}$ tensorfields
	$\upchi^{(Small)},$
	$\angkuparg{(Trans-\Psi)},$ 
	and
	$\angkuparg{(Tan-\Psi)}$
	are defined by
	\eqref{E:CHIJUNKDEF},
	\eqref{E:KABTRANSVERSAL},
	and \eqref{E:KABGOOD}.
			\end{lemma}

		\begin{proof}
			For each real number $\uplambda$ belonging to a small neighborhood of $0,$
			let $\varphi_{(\uplambda)}: \mathbb{R}^4 \rightarrow \mathbb{R}^4$ 
			be the flow map of $\Lunit.$ That is, for each spacetime point $p,$ the four rectangular 
			components of $\varphi_{(\uplambda)}(p)$ verify the ODE system
			$\frac{d}{d \uplambda} \varphi_{(\uplambda)}^{\nu}(p) = \Lunit^{\nu} \circ \varphi_{(\uplambda)}(p)$ 
			and $\varphi_{(0)} = I,$ where $I$ is the identity map on $\mathbb{R}^4.$
			Let $\varphi_{(\uplambda)}^{*}$ denote pullback by $\varphi_{(\uplambda)}$
			so that in particular, $\varphi_{(\uplambda)}^{*} f = f \circ \varphi_{(\uplambda)}$
			for functions $f.$
			Note that since $\Lunit t = 1,$ it follows that
			$\varphi_{(\uplambda)}$ bijectively maps the sphere $S_{t,u}$ onto
			the sphere $S_{t + \uplambda, u}.$ Hence, by covariance, we have
			\begin{align} \label{E:STUPULLBACKINTEGRALIDENTITY}
				\int_{S_{t + \uplambda,u}} 
					f 
				\, d \spherevol
				& = 
				\int_{S_{t,u}}
					\varphi_{(\uplambda)}^{*} f  
				\, d \upsilon_{{\varphi_{(\uplambda)}^{*} \gsphere}}.
			\end{align} 
			We now claim that 
			\begin{align} \label{E:DERIRVATIVEOFSTUVOLUMEFORMPULLBACK}
				\frac{d}{d \uplambda}\left|_{\uplambda = 0}\right. d \upsilon_{{\varphi_{(\uplambda)}^{*} \gsphere}}
				& = \angLie_{\Lunit} d \spherevol
					= \frac{1}{2} \mytr  \angdeform{\Lunit}d \spherevol.
			\end{align}
			The identity \eqref{E:DERIRVATIVEOFSTUVOLUMEFORMPULLBACK} is straightforward to verify
			in the local $S_{t,u}$ 
			coordinates $(\vartheta^1, \vartheta^2),$
			in which
			$d \upsilon_{{\varphi_{(\uplambda)}^{*} \gsphere}}
			= \sqrt{\mbox{\upshape{det}}{\varphi_{(\uplambda)}^{*} \gsphere} } \, d \vartheta^1 \wedge d \vartheta^2.$
			In particular, one uses the identity
			\begin{align}
				\frac{d}{d \uplambda} \left|_{\uplambda = 0}\right. \ln \mbox{\upshape{det}} \varphi_{(\uplambda)}^{*} \gsphere
				= ([\varphi_{(\uplambda)}^{*} \gsphere]^{-1})^{AB}\left|_{\uplambda = 0}\right. 
				\frac{d}{d \uplambda} \left|_{\uplambda = 0}\right. \varphi_{(\uplambda)}^{*} \gsphere_{AB}
				= (\gsphere^{-1})^{AB} \angLie_{\Lunit} \gsphere_{AB}.
			\end{align}
			Thus, differentiating each side of \eqref{E:STUPULLBACKINTEGRALIDENTITY}
			with $\frac{d}{d \uplambda},$
			setting $\uplambda = 0,$ and using \eqref{E:DERIRVATIVEOFSTUVOLUMEFORMPULLBACK}, we deduce
			\begin{align} \label{E:TIMEDERIVATIVEOFSPHEREINTEGRALS}
				\frac{\partial}{\partial t} 
				\left(
					\int_{S_{t,u}} f \, d \spherevol
				\right)
				& 
				= \frac{d}{d \uplambda}\left|_{\uplambda = 0}\right.
				\int_{S_{t + \uplambda,u}} 
					f 
				\, d \spherevol
				= 	\int_{S_{t,u}} 
							(\Lunit f)
						\, d \spherevol
					+ \int_{S_{t,u}} 
							f	
						\, \angLie_{\Lunit} d \spherevol
					\\
			& = 
				\int_{S_{t,u}} 
							(\Lunit f)
						\, d \spherevol
					+ \frac{1}{2}
						\int_{S_{t,u}} 
							f	\mytr  \angdeform{\Lunit}
						\, d \spherevol.
						\notag
			\end{align}
			Similarly, since $\Rad u = 1,$ the identity 
			\eqref{E:TIMEDERIVATIVEOFSPHEREINTEGRALS} holds
			with $\frac{\partial}{\partial t}$
			replaced by $\frac{\partial}{\partial u}$
			and $\Lunit$ replaced by $\Rad.$
			Lemma~\ref{L:STUDERIVATIVES} thus follows from
			the identity \eqref{E:TIMEDERIVATIVEOFSPHEREINTEGRALS}
			and the identities
			\begin{align}
				\frac{1}{2} \mytr  \angdeform{\Lunit} 
				& = \mytr \upchi,	
					\label{E:LDERIVATIVESPHERVOLUMEFORM}  \\
				\frac{1}{2} \mytr  \angdeform{\Rad}
				& = 
					\left\lbrace
						- \frac{2}{\rgeo} \upmu 
						- \upmu \mytr \upchi^{(Small)}
						+ \mytr  \angkuparg{(Trans-\Psi)}
						+ \upmu \mytr  \angkuparg{(Tan-\Psi)}
					\right\rbrace,
				\label{E:RADDERIVATIVESPHERVOLUMEFORM}
			\end{align}	
			which follow from Prop.~\ref{P:DEFORMATIONTENSORFRAMECOMPONENTS}.	
	\end{proof}

The following tensorfield is a core ingredient in our derivation
of $L^2-$type estimates for solutions to the inhomogeneous
wave equation $\upmu \square_{g(\Psi)} \Psi = \waveinhom.$

\begin{definition}[\textbf{Energy-momentum tensorfield}]
\label{D:ENERGYMOMENTUMTENSOR}
We define the energy-momentum tensorfield
$\enmomtensor_{\mu \nu}[\Psi]$
associated to $\Psi$ as follows: 
\begin{align} \label{E:ENERGYMOMENTUMTENSOR}
	\enmomtensor_{\mu \nu}
	=
	\enmomtensor_{\mu \nu}[\Psi]
	& := \D_{\mu} \Psi \D_{\nu} \Psi
	- \frac{1}{2} g_{\mu \nu} (g^{-1})^{\alpha \beta} \D_{\alpha} \Psi \D_{\beta} \Psi.
\end{align}

\end{definition}

In the next lemma, we provide the fundamental divergence property  
of $\enmomtensor_{\mu \nu}[\Psi]$ for solutions to the inhomogeneous wave equation.

\begin{lemma}[\textbf{Fundamental divergence property of $\enmomtensor_{\mu \nu}[\Psi]$}]
	\label{L:ENERGYMOMENTUMTENSORDIVERGNCEFORINHOMOGENEOUSWAVEEQUATION}
	If $\Psi$ verifies the $\upmu-$weighted equation
	\begin{align} \label{E:MUWEIGHTEDINHOMOGENEOUSWAVEEQUATION}
		\upmu \square_g \Psi & = \waveinhom,
	\end{align}
	then
	\begin{align} \label{E:ENERGYMOMENTUMTENSORDIVERGNCEFORINHOMOGENEOUSWAVEEQUATION}
		\upmu
		\D_{\alpha} (\enmomtensor_{\ \nu}^{\alpha}[\Psi])
		& = \waveinhom \D_{\nu} \Psi.
	\end{align}
\end{lemma}
\begin{proof}
To derive \eqref{E:ENERGYMOMENTUMTENSORDIVERGNCEFORINHOMOGENEOUSWAVEEQUATION},
we take the divergence of the right-hand side of \eqref{E:ENERGYMOMENTUMTENSOR},
multiply by $\upmu,$
and replace the factor 
$\upmu \square_{g(\Psi)} \Psi = \upmu (g^{-1})^{\alpha \beta} \D^2_{\alpha \beta} \Psi$
with $\waveinhom$ when it arises. 
\end{proof}

\begin{remark}[\textbf{The role of the $\upmu-$weighted wave operator}]
	We consider the $\upmu-$weighted product $\upmu \square_{g(\Psi)} \Psi$
	in \eqref{E:MUWEIGHTEDINHOMOGENEOUSWAVEEQUATION}
	because this is the product that will appear 
	when we apply the divergence theorem.
	The reason that this product appears rather than
	$\square_{g(\Psi)} \Psi$ is that our spacetime integrals are defined
	through a rescaled spacetime volume $d \vol$ in which the term
	$\upmu$ has been factored out (see Sect.~\ref{S:VOLFORMSANDSOBOLEVNORMS}).
\end{remark}

We use the following two multiplier vectorfields 
to construct the $L^2-$type quantities that control $\Psi.$

\begin{definition}[\textbf{Multiplier vectorfields}]
\label{D:MULTIPLIERVECTORFIELDS}
We define the multiplier vectorfields	$\Mult$ and $\Mor$ by
	\begin{subequations}
	\begin{align}
		\Mult 
		& := (1 + \upmu) \Lunit + \uLgood
			= (1 + 2 \upmu) \Lunit + 2 \Rad, 
			\label{E:DEFINITIONMULT} \\
		\Mor
		& := \rgeo^2 \Lunit,
			\label{E:DEFINITIONMOR}
	\end{align}
	\end{subequations}
	where $\rgeo$ is defined in \eqref{E:GEOMETRICRADIAL}.
\end{definition}

We note that $\Mult$ is $g-$future-directed and $g-$timelike,
while $\Mor$ is $g-$future-directed and $g-$null.
$\Mult$ is an analog of the Minkowskian vectorfield
$2 \partial_t + \Lunit_{(Flat)} = 3 \partial_t + \partial_r,$ 
but it has been carefully
constructed so as to be useful both in regions where 
$\upmu$ is large and in regions where $\upmu$ is small.
$\Mor$ is an analog of the \emph{Morawetz multiplier},
which has played an important role in the global analysis
of various wave-like equations. We note that neither
$\Mult$ nor $\Mor$ are Killing fields, even for the background
solution $\Psi \equiv 0.$ That is, both
$\deform{\Mult}$ 
and
$\deform{\Mor}$
are non-zero.

We now derive expressions for the components of
$\enmomtensor[\Psi]$ relative to the rescaled null frame.

\begin{lemma}[\textbf{Null components of $\enmomtensor_{\mu \nu}$}] \label{L:NULLCOMPONENTSOFENERGYMOMENTUMTENSOR}
	Let $\enmomtensor[\Psi]$ be the energy-momentum tensorfield defined in \eqref{E:ENERGYMOMENTUMTENSOR}.
	The components of $\enmomtensor[\Psi]$ relative to the rescaled null frame $\lbrace \Lunit, \uLgood, X_1, X_2 \rbrace$ are
	\begin{subequations}
	\begin{align}
		\enmomtensor_{\Lunit \Lunit}[\Psi] 
		& = (\Lunit \Psi)^2, 
			\\
		\enmomtensor_{\uLgood \uLgood}[\Psi]
		& = (\uLgood \Psi)^2, 
			\\
		\enmomtensor_{\Lunit \uLgood}[\Psi] 
		& = \upmu |\angdiff \Psi|^2,
			\label{E:TLULGOOD} \\
		\enmomtensor_{\Lunit A}[\Psi] 
		& = (\Lunit \Psi) \angdiffarg{A} \Psi, 
			\\
		\enmomtensor_{\uLgood A}[\Psi] & = (\uLgood \Psi) \angdiffarg{A} \Psi, 
			\\
		\enmomtensor_{AB}[\Psi]
		& = (\angdiffarg{A} \Psi)(\angdiffarg{B} \Psi) 
			- \frac{1}{2} \gsphere_{AB} 
				\left\lbrace
					- \upmu^{-1} (\Lunit \Psi)(\uLgood \Psi)  
					+ |\angdiff \Psi|_{\gsphere}^2 
				\right\rbrace.
	\end{align}
	\end{subequations}
\end{lemma}

\begin{proof}
The proof is a series of computations that can be carried out with the help of 
the decomposition \eqref{E:GINVERSENULLFRAME} of $g^{-1}$ relative to the null frame.
For example, using the identity $g(\Lunit, \uLgood) = - 2 \upmu,$ we compute that
$\enmomtensor_{\Lunit \uLgood}[\Psi] = (\Lunit \Psi)(\uLgood \Psi) 
- \frac{1}{2} g(\Lunit, \uLgood)\left\lbrace - \upmu^{-1} (\Lunit \Psi)(\uLgood \Psi) + |\angdiff \Psi|^2 \right\rbrace 
= \upmu |\angdiff \Psi|^2.$ We have thus proved \eqref{E:TLULGOOD}. We leave the proofs of the remaining identities
to the reader.

\end{proof}

\begin{lemma}[\textbf{Null components of $\deform{\Mult}$ and $\deform{\Mor} - \rgeo^2 \mytr \upchi g$}]
\label{L:NULLCOMPONENTSOFDEFORMATIONTENSORSOFMULTIPLIERS}
Let $\Mult$ and $\Mor$ be the multiplier vectorfields from Def.~\ref{D:MULTIPLIERVECTORFIELDS}.
Let $\deform{\Mult}$ and $\deform{\Mor}$ denote the corresponding deformation tensors as defined in
Def.~\ref{D:DEFORMTENSDEFINED}.
Relative to the rescaled null frame $\lbrace \Lunit, \uLgood, X_1, X_2 \rbrace,$
the components of $\deform{\Mult}$ are:
\begin{subequations}
\begin{align}
	\deformarg{\Mult}{\Lunit}{\Lunit} 
	& = 0, 
		\\
	\deformarg{\Mult}{\uLgood}{\uLgood} 
	& = - 4 \upmu 
		\left\lbrace 
				\uLgood \upmu
				-(1 + \upmu) \Lunit \upmu
		\right\rbrace, 
		\\
	\deformarg{\Mult}{\Lunit}{\uLgood}
	& = - 2 \upmu
		\left\lbrace
			\upmu^{-1} (\uLgood \upmu + \Lunit \upmu)
			+ 3 \Lunit \upmu
		\right\rbrace,
		 \\
	\deformarg{\Mult}{\Lunit}{A} 
	& = - 2 \angdiffarg{A} \upmu 
		- 4 \upzeta_A^{(Trans-\Psi)}
		- 4 \upmu \upzeta_A^{(Tan-\Psi)},
		\\
	\deformarg{\Mult}{\uLgood}{A} 
	& = 2 (1 - \upmu)
			\angdiffarg{A} \upmu 
		+ 4 (1 + \upmu) 
			\left\lbrace
				\upzeta_A^{(Trans-\Psi)}
				+ \upmu \upzeta_A^{(Tan-\Psi)}
			\right\rbrace,
			\\
	\angdeformfreearg{\Mult}{A}{B}
	& = 2 \left\lbrace
					\hat{\upchi}_{AB}
					+ 2 \angkfreetriplearg{A}{B}{(Trans-\Psi)}
					+ 2 \upmu \angkfreetriplearg{A}{B}{(Tan-\Psi)}
				\right\rbrace,
		\label{E:MULTDEFORMSPHERETRACEFREE} \\
	\mytr  \angdeform{\Mult}
	& = 2 \left\lbrace
					\mytr \upchi
					+ 2 \mytr  \angkuparg{(Trans-\Psi)} 
					+ 2 \upmu \mytr  \angkuparg{(Tan-\Psi)}
				\right\rbrace.
				\label{E:MULTDEFORMSPHERETRACE}
\end{align}
\end{subequations}

Furthermore, the components of $\deform{\Mor} - \rgeo^2 \mytr \upchi g$ are:
\begin{subequations}
\begin{align}
	\deformarg{\Mor}{\Lunit}{\Lunit} 
		- \rgeo^2 \mytr \upchi g(\Lunit, \Lunit)
	& = 0, 
	\\
	\deformarg{\Mor}{\uLgood}{\uLgood}
		- \frac{1}{2} \rgeo^2 \mytr \upchi g(\uLgood, \uLgood)
	& = - 4 \upmu 
		\left\lbrace
			2 \rgeo \uLgood \rgeo
			- \rgeo^2 \Lunit \upmu
		\right\rbrace,
		\\
	\deformarg{\Mor}{\Lunit}{\uLgood}
		- \rgeo^2 \mytr \upchi g(\Lunit, \uLgood)
	& = - 2 
			\rgeo^2 
			\upmu
			\left\lbrace
				\frac{\Lunit \upmu}{\upmu}
				-\mytr \upchi^{(Small)}
			\right\rbrace,
			\label{E:RENORMALIZEDMORDEFORMLULGOOD} \\
	\deformarg{\Mor}{\Lunit}{A}
		- \rgeo^2 \mytr \upchi g(\Lunit, A)
	& = 0, \\
	\deformarg{\Mor}{\uLgood}{A}
		- \rgeo^2 \mytr \upchi g(\uLgood, A)
	& = 2 \rgeo^2
				\left\lbrace 
					\angdiffarg{A} \upmu
					+ 2 \upzeta_A^{(Trans-\Psi)}
					+ 2 \upmu \upzeta_A^{(Tan-\Psi)}
				\right\rbrace, 
		\\
	\angdeformfreearg{\Mor}{A}{B}
		- \rgeo^2 \mytr \upchi \gsphere_{AB}
	& = 2 \rgeo^2 \hat{\upchi}_{AB}^{(Small)},
		\\
	\mytr \angdeform{\Mor}
		- 2 \rgeo^2 \mytr \upchi
	& = 0.
	\end{align}
\end{subequations}
In the above expressions, the $S_{t,u}$ tensorfields
	$\upchi,$ 
	$\upchi^{(Small)},$
	$\upzeta^{(Trans-\Psi)},$
	$\angkuparg{(Trans-\Psi)},$ 
	$\upzeta^{(Tan-\Psi)},$
	and
	$\angkuparg{(Tan-\Psi)}$
	are defined by
	\eqref{E:CHIDEF},
	\eqref{E:CHIJUNKDEF},
	\eqref{E:ZETATRANSVERSAL},
	\eqref{E:KABTRANSVERSAL},
	\eqref{E:ZETAGOOD},
	and \eqref{E:KABGOOD}.
\end{lemma}

\begin{proof}
The proof is a series of computations based on the rescaled null frame connection coefficient
identities of Lemma~\ref{L:CONNECTIONCOEFFICENTSOFNULLFRAME}
and the ideas used in the proof of
Lemma~\ref{L:NULLCOMPONENTSOFENERGYMOMENTUMTENSOR}.
We give two examples. We use the notation $\langle V,W \rangle := g(V,W)$ to simplify the presentation.
As a first example, we prove \eqref{E:MULTDEFORMSPHERETRACEFREE} and \eqref{E:MULTDEFORMSPHERETRACE}.
Using 
\eqref{E:DEFORMTENSDEFINED}
and
\eqref{E:DEFINITIONMULT},
we compute that
\begin{align} \label{E:ANGDEFORMMULTABCOMPUTATION}
\angdeformarg{\Mult}{A}{B} 
	& = \left\langle 
				\D_A \left\lbrace (1 + \upmu) \Lunit + \uLgood \right\rbrace, X_B 
			\right\rangle 
		+ \left\langle 
				\D_B \left\lbrace (1 + \upmu) \Lunit + \uLgood \right\rbrace, X_A 
			\right\rangle  
		\\
	& = (1 + \upmu) \langle \D_A \Lunit, X_B \rangle
		+ (1 + \upmu) \langle \D_B \Lunit, X_A \rangle
		+ \langle \D_A \uLgood, X_B \rangle
		+ \langle \D_B \uLgood, X_A \rangle.
		\notag
\end{align}
Also using Lemma
\ref{L:CONNECTIONCOEFFICENTSOFNULLFRAME} and the decomposition \eqref{E:ANGKDECOMPOSED}, 
we deduce that the right-hand side of \eqref{E:ANGDEFORMMULTABCOMPUTATION} is equal to
$2 \upchi_{AB} 
	+ 4 \angktriplearg{A}{B}{(Trans-\Psi)} 
	+ 4 \upmu \angktriplearg{A}{B}{(Tan-\Psi)}.
$
The desired identities \eqref{E:MULTDEFORMSPHERETRACEFREE} and \eqref{E:MULTDEFORMSPHERETRACE}
now follow from splitting this expression
into its trace and trace-free parts.

As a second example, we prove \eqref{E:RENORMALIZEDMORDEFORMLULGOOD}.
Using 
\eqref{E:DEFORMTENSDEFINED}, 
\eqref{E:DEFINITIONMOR}, 
and
the identities 
$\Lunit \rgeo = 1$
and
$\uLgood \rgeo = \upmu - 2,$
we compute that 
\begin{align} \label{E:DEFORMMORLULGOOD}
	\deformarg{\Mor}{\Lunit}{\uLgood}
	& = \langle \D_{\Lunit} (\rgeo^2 \Lunit), \uLgood \rangle 
	+ \langle \D_{\uLgood} (\rgeo^2 \Lunit), \Lunit \rangle
		\\
	& = - 4 \rgeo \upmu
		+ \rgeo^2 \langle \D_{\Lunit} \Lunit, \uLgood \rangle
	+ \rgeo^2 \langle \D_{\uLgood} \Lunit, \Lunit \rangle.
	\notag
\end{align}
Also using Lemma~\ref{L:CONNECTIONCOEFFICENTSOFNULLFRAME}, 
we deduce that the right-hand side of \eqref{E:DEFORMMORLULGOOD} is equal to
$- 4 \rgeo \upmu 
	- 2 \rgeo^2 \Lunit \upmu,$
which is the first product on the left-hand side of \eqref{E:RENORMALIZEDMORDEFORMLULGOOD}.
Furthermore, using the decomposition $\mytr \upchi = 2 \rgeo^{-1} + \mytr \upchi^{(Small)},$
we compute that the second term on the left-hand side of \eqref{E:RENORMALIZEDMORDEFORMLULGOOD}
is 
\begin{align} \label{E:SECONDTERMLHSRENORMALIZEDMORDEFORMLULGOOD}
	- \rgeo^2 \mytr \upchi g(\Lunit, \uLgood) = 4 \rgeo \upmu + 2 \upmu \rgeo^2 \mytr \upchi^{(Small)}.
\end{align}
Adding the previous term 
$- 4 \rgeo \upmu - \rgeo^2 \Lunit \upmu$
to \eqref{E:SECONDTERMLHSRENORMALIZEDMORDEFORMLULGOOD}
yields 
$- \rgeo^2 \Lunit \upmu + 2 \upmu \rgeo^2 \mytr \upchi^{(Small)},$
from which we easily conclude the desired identity \eqref{E:RENORMALIZEDMORDEFORMLULGOOD}.
\end{proof}

We now define vectorfields that we use for bookkeeping
when we integrate by parts to derive energy-flux identities.

\begin{definition}[\textbf{Compatible currents}] \label{D:COMPATIBLECURRENTS}
	Let $\Psi$ be a function, and let $\enmomtensor[\Psi]$ be the energy-momentum tensorfield defined in \eqref{E:ENERGYMOMENTUMTENSOR}.
	Let $\Mult$ and $\Mor$ be the multiplier vectorfields from Def.~\ref{D:MULTIPLIERVECTORFIELDS}.
	We associate to $\Psi$ the following four ``compatible current'' vectorfields:
	\begin{subequations}
	\begin{align}
			\Jenergycurrent{\Mult}^{\nu}[\Psi]
			& := \enmomtensor_{\ \alpha}^{\nu}[\Psi] \Mult^{\alpha},
				\label{E:MULTCURRENTDEF} \\
			\Jenergycurrent{\Mor}^{\nu}[\Psi]
			& := \enmomtensor_{\ \alpha}^{\nu}[\Psi]\Mor^{\alpha},
					\label{E:MORCURRENTDEF} \\
			\Jenergycurrent{Correction}^{\nu}[\Psi]
			& := \frac{1}{2} 
						\left\lbrace
							\rgeo^2 \mytr \upchi \Psi \D^{\nu} \Psi
							- \frac{1}{2} \Psi^2 \D^{\nu}[\rgeo^2 \mytr \upchi]
					  \right\rbrace,
				\label{E:MODCURRENTDEF} \\
			\Jenergycurrent{\Mor+Correction}^{\nu}[\Psi]
			& := {\Jenergycurrent{\Mor}^{\nu}[\Psi]}	
				+ {\Jenergycurrent{Correction}^{\nu}[\Psi]}.
				\label{E:MORPLUSMODCURRENT}
	\end{align}
	\end{subequations}

\end{definition}

In the next lemma, 
as a precursor to applying the divergence theorem,
we compute the divergence of the compatible currents.

\begin{lemma}[\textbf{Preliminary expressions for the $\upmu-$weighted divergence of the currents}]
	\label{L:PRELIMDIV}
	If $\Psi$ verifies the $\upmu-$weighted inhomogeneous wave equation
	\begin{align} \label{E:MUWEIGHTEDINHOMWAVE}
		\upmu \square_{g(\Psi)} \Psi & = \waveinhom,
	\end{align}
	then the following divergence identities hold:
\begin{subequations}
\begin{align}
		\upmu
		\D_{\alpha} {\Jenergycurrent{\Mult}^{\alpha}[\Psi]}
		& =  \left\lbrace
						(1 + 2 \upmu) \Lunit \Psi + 2 \Rad \Psi
					\right\rbrace
					\waveinhom 
				+ \frac{1}{2}
					\upmu
					\enmomtensor^{\alpha \beta}[\Psi] \deformarg{\Mult}{\alpha}{\beta},
				\label{E:MULTCURRENTDIV}	\\
			\upmu
		\D_{\alpha} {\Jenergycurrent{\Mor}^{\alpha}[\Psi]}
		& = \rgeo^2 
				(\Lunit \Psi) 
				\waveinhom 
				+ \frac{1}{2}
					\upmu
					\enmomtensor^{\alpha \beta}[\Psi] \deformarg{\Mor}{\alpha}{\beta},
				\label{E:MORCURRENTDIV}	\\
		\upmu
		\D_{\alpha} {\Jenergycurrent{Correction}^{\alpha}[\Psi]}
		& =  \frac{1}{2} \rgeo^2 \mytr \upchi
					\Psi 
					\waveinhom
			+ \frac{1}{2} 
				\rgeo^2 \mytr \upchi
				\upmu 
				(g^{-1})^{\alpha \beta} 
				\D_{\alpha} 
				\Psi \D_{\beta} \Psi
			- \frac{1}{4} 
				\upmu 
				\left\lbrace 
					\square_{g(\Psi)} [\rgeo^2 \mytr \upchi] 
				\right\rbrace 
				\Psi^2,
				\label{E:MODCURRENTDIV}	\\
		\upmu
		\D_{\alpha} {\Jenergycurrent{\Mor+Correction}^{\alpha}[\Psi]}
		& =  \rgeo^2 
					\waveinhom 
					\left\lbrace
							\Lunit \Psi
							+ \frac{1}{2} \mytr \upchi \Psi
					\right\rbrace
				\label{E:MORMODDIV} \\
		& \ \ + \frac{1}{2}
						\upmu \enmomtensor^{\alpha \beta}[\Psi] 
						\left\lbrace
							\deformarg{\Mor}{\alpha}{\beta}
							- \rgeo^2 \mytr \upchi g_{\alpha \beta}
						\right\rbrace
				- 
				\frac{1}{4} 
				\upmu 
				\left\lbrace 
					\square_{g(\Psi)} [\rgeo^2 \mytr \upchi] 
				\right\rbrace 
				\Psi^2.
				\notag
\end{align}
\end{subequations}

\end{lemma}

\begin{proof}
	Equation \eqref{E:MULTCURRENTDIV} follows	
	applying $\D_{\nu}$ to \eqref{E:MULTCURRENTDEF},
	using the Leibniz rule
	and Lemma~\ref{L:ENERGYMOMENTUMTENSORDIVERGNCEFORINHOMOGENEOUSWAVEEQUATION},
	definition \eqref{E:DEFORMTENSDEFINED},
	definition \eqref{E:DEFINITIONMOR},
	and the symmetry of $\enmomtensor_{\mu \nu}[\Psi]$ and $\deformarg{\Mult}{\mu}{\nu}.$
	The proof of \eqref{E:MORCURRENTDIV}
	is based on definition \eqref{E:MORCURRENTDEF} and is similar; we omit the details.
	
	Equation \eqref{E:MODCURRENTDIV} follows 
	from applying $\D_{\nu}$ to \eqref{E:MODCURRENTDEF},
	using the Leibniz rule, and from
	using equation \eqref{E:MUWEIGHTEDINHOMWAVE}
	to substitute for the term $\upmu \square_{g(\Psi)} \Psi$
	when it arises.
	
	Equation \eqref{E:MORMODDIV}
	follows from 
	definition \eqref{E:MORPLUSMODCURRENT},
	equations
	\eqref{E:MORCURRENTDIV}
	and
	\eqref{E:MODCURRENTDIV},
	and the identity
	$ g_{\alpha \beta} \enmomtensor^{\alpha \beta}[\Psi] = - (g^{-1})^{\alpha \beta} \D_{\alpha} \Psi \D_{\beta} \Psi,$
	which easily follows from definition \eqref{E:ENERGYMOMENTUMTENSOR}.
	
\end{proof}

\section{The energy-flux integral identities}

We are now ready to define the energies and fluxes that we
use in our $L^2$ analysis of solutions.

\begin{definition}[\textbf{Energies and fluxes}]
\label{D:ENERGIESANDFLUXES}
	Let $\Timenormal = \frac{1}{2} \Lunit + \frac{1}{2} \upmu^{-1} \uLgood = \Lunit + \Radunit$
	denote the future-directed unit normal to $\Sigma_t$
	(see Lemma~\ref{L:BASICPROPERTIESOFLUNITRADANDTIMENORMAL}).
	We define the energies $\enzero[\Psi](t,u),$ $\enone[\Psi](t,u)$
	and the cone fluxes $\flzero[\Psi](t,u),$ $\flone[\Psi](t,u)$ as follows,
	where the rescaled volume forms 
	$d \tvol$ 
	and
	$d \conevol$
	are defined in Def.~\ref{D:AREAANDRESCALEDVOLUMEFORM}:
	\begin{subequations}
	\begin{align}
		\enzero[\Psi](t,u)
		& : = \int_{\Sigma_t^u} \upmu \Jenergycurrent{\Mult}_{\Timenormal}[\Psi] \, d \tvol,
			\label{E:E0DEF} \\
		\enone[\Psi](t,u)
		& := \int_{\Sigma_t^u} \upmu \Jenergycurrent{\Mor+Correction}_{\Timenormal}[\Psi] \, d \tvol
			\label{E:E1DEF} \\
		& \ \ - \frac{1}{2} \int_{\Sigma_t^u} \rgeo^2 \mytr \upchi \Psi \Rad \Psi \, d \tvol
			+ \frac{1}{4} \int_{\Sigma_t^u}
				\left\lbrace
					\upmu \Lunit [\rgeo^2 \mytr \upchi]
					+ \Rad [\rgeo^2 \mytr \upchi]
					+ \frac{1}{2} \upmu \rgeo^2 (\mytr \upchi)^2 
				\right\rbrace 
				 \Psi^2 
				\, d \tvol,
				\notag
	\end{align}
	\end{subequations}
	\begin{subequations}
	\begin{align}
		\flzero[\Psi](t,u)
		& := \int_{\mathcal{C}_u^t}  {\Jenergycurrent{\Mult}_{\Lunit}[\Psi]} \, d \conevol,
			\label{E:F0DEF} \\
		\flone[\Psi](t,u)
		& := \int_{\mathcal{C}_u^t} {\Jenergycurrent{\Mor+Correction}_{\Lunit}[\Psi]} \, d \conevol
			\label{E:F1DEF} 
			+ \frac{1}{4} 
				\int_{\mathcal{C}_u^t} 
					\left\lbrace
						\Lunit[\rgeo^2 \mytr \upchi \Psi^2] 
						+ \rgeo^2 (\mytr \upchi)^2 \Psi^2
					\right\rbrace
				\, d \conevol.
	\end{align}
	\end{subequations}
	
	\end{definition}
	The coerciveness properties of the energies and fluxes are established in Lemma~\ref{L:COERCIVEENERGIESANDFLUXES}.
	
	In the next proposition, we derive our main energy-flux identities, 
 	which are the starting points of our $L^2$ analysis of $\Psi$ and its derivatives.
	As a preliminary step, we first prove the following lemma, which 
	provides a version of the divergence theorem on the spacetime region $\mathcal{M}_{t,u}.$
	\emph{Note that we explicitly indicate the factors of} $\upmu$ in the lemma.
	
	\begin{lemma}[\textbf{The divergence theorem on the region} $\mathcal{M}_{t,u}$]
		\label{L:WTUSPACETIMEDIVERGENCETHM} 
		Let $J$ be any spacetime vectorfield.
		Let $\Timenormal = \Lunit + \Radunit$ be the future-directed unit normal to 
		$\Sigma_t^{U_0}.$
		Then the following integral identity holds,
		where the rescaled volume forms are defined in
		Def.~\ref{D:AREAANDRESCALEDVOLUMEFORM}:
		\begin{align} \label{E:WTUSPACETIMEDIVERGENCETHM} 
			\int_{\Sigma_t^u} 
				\upmu J_{\Timenormal} 
			\, d \tvol
			+ 
			\int_{\mathcal{C}_u^t} 
				J_{\Lunit}
			\, d \conevol
			-
			\int_{\Sigma_0^u} 
				\upmu J_{\Timenormal}
			\, d \tvol
			- \int_{\mathcal{C}_0^t} 
					J_{\Lunit}
				\, d \conevol
		& = -
				\int_{\mathcal{M}_{t,u}}
					\upmu \D_{\alpha} J^{\alpha}
				\, d \vol.
		\end{align}
	\end{lemma}
	
	\begin{proof}
		We can expand $J$ relative to the geometric coordinate frame as follows:
		$J = J^t \frac{\partial}{\partial t} 
		+ J^u \frac{\partial}{\partial u}
		+ \angnormalJ^A \frac{\partial}{\partial \vartheta^A},$
		where $\angnormalJ$ is $S_{t,u}-$tangent.
		Furthermore, 
		from the identities 
		$\Lunit = \frac{\partial}{\partial t},$
		$\Rad = \frac{\partial}{\partial u} - \Xi^A \frac{\partial}{\partial \vartheta^A},$
		$g(\Lunit, \Lunit) 
		= g(\Lunit, \frac{\partial}{\partial \vartheta^A}) 
		= g(\Rad, \frac{\partial}{\partial \vartheta^A})
		= 0,$
		$g(\Lunit, \Rad) = - \upmu,$
		and
		$g(\Rad, \Rad) = \upmu^2,$
		it is straightforward to compute that
		\begin{align} \label{E:COORDINATEVECTROFIELDCOMPONENTSINTERMSOFFRAMECOMPONENTS}
			J^t & = - J_{\Timenormal}, && J^u = - \upmu^{-1} J_{\Lunit}.
		\end{align}
		We will use the following standard formula for 
		the divergence $\D_{\alpha} J^{\alpha}$ of $J:$
		\begin{align} \label{E:COORDINATEDIVERGENCEFIRSTEXPRESSION}
			\D_{\alpha} J^{\alpha}
			=
			\frac{1}{\sqrt{|\mbox{\upshape{det}} g|}}
			\left\lbrace
				\frac{\partial}{\partial t} (\sqrt{|\mbox{\upshape{det}} g|} J^t)
				+ \frac{\partial}{\partial u} (\sqrt{|\mbox{\upshape{det}} g|} J^u)
				+ \frac{\partial}{\partial \vartheta^A} (\sqrt{|\mbox{\upshape{det}} g|} \angnormalJ^A)
			\right\rbrace,
		\end{align}
		where the determinant $\mbox{\upshape{det}} g$ is taken relative to the geometric coordinates.
		Inserting the square root of \eqref{E:SPACETIMEVOLUMEFORMWITHUPMU} into \eqref{E:COORDINATEDIVERGENCEFIRSTEXPRESSION},
		we deduce that
		\begin{align} \label{E:COORDINATEDIVERGENCESECONDEXPRESSION}
			\frac{\partial}{\partial t} (\upmu \sqrt{\mbox{\upshape{det}} \gsphere} J^t)
				+ \frac{\partial}{\partial u} (\upmu \sqrt{\mbox{\upshape{det}} \gsphere} J^u)
				+ \frac{\partial}{\partial \vartheta^A} (\upmu \sqrt{|\mbox{\upshape{det}} \gsphere|} \angnormalJ^A)
			= (\upmu \D_{\alpha} J^{\alpha}) \sqrt{\mbox{\upshape{det}} \gsphere}.
		\end{align}
		Integrating \eqref{E:COORDINATEDIVERGENCESECONDEXPRESSION} 
		over $\mathbb{S}^2$
		with respect to $d \vartheta$
		and noting that the integral 
		of the last term on the left-hand side of \eqref{E:COORDINATEDIVERGENCESECONDEXPRESSION} vanishes,
		we deduce that
		\begin{align} \label{E:COORDINATEDIVERGNCESPHEREINTEGRATED}
			\frac{\partial}{\partial t}
				\left(
					\int_{S_{t,u}} 
						\upmu J^t 	
						\, d \spherevol
				\right)
		+ \frac{\partial}{\partial u}
				\left(
					\int_{S_{t,u}} 
						\upmu J^u 
					\, d \spherevol
				\right)
		& = \int_{S_{t,u}}
					\upmu \D_{\alpha} J^{\alpha}
					\, d \spherevol.
		\end{align}
		Integrating \eqref{E:COORDINATEDIVERGNCESPHEREINTEGRATED} 
		with respect to $du' dt'$
		over the coordinate rectangle $[0,t] \times [0,u],$
		we deduce that
		\begin{align} \label{E:SPACETIMEDIVERGENCENEARFINALFORM} 
			\int_{\Sigma_t^u} 
				\upmu J^t 
			\, d \tvol
			-
			\int_{\Sigma_0^u} 
				\upmu J^t 
			\, d \tvol
			+ 
			\int_{\mathcal{C}_u^t} 
				\upmu J^u
			\, d \conevol
			- \int_{\mathcal{C}_0^t} 
				\upmu J^u
			\, d \conevol
		& = \int_{\mathcal{M}_{t,u}}
					\upmu \D_{\alpha} J^{\alpha}
				\, d \vol.
		\end{align}
		Inserting the identities \eqref{E:COORDINATEVECTROFIELDCOMPONENTSINTERMSOFFRAMECOMPONENTS} into \eqref{E:SPACETIMEDIVERGENCENEARFINALFORM},
		we arrive at \eqref{E:WTUSPACETIMEDIVERGENCETHM}.
		
	\end{proof}
	
	We now derive our main energy-flux identities.
	
	\begin{proposition}[\textbf{Divergence theorem with important cancellations}]	
	\label{P:DIVTHMWITHCANCELLATIONS}
		Consider the energies and fluxes of Def.~\ref{D:ENERGIESANDFLUXES}.
		For solutions $\Psi$ to 
		\begin{align*}
			\upmu \square_{g(\Psi)} \Psi & = \waveinhom,
		\end{align*}
		that vanish along the outer null cone $\mathcal{C}_0,$
		we have the following integral identities,
		where the rescaled volume form $d \vol$
		is defined in
		Def.~\ref{D:AREAANDRESCALEDVOLUMEFORM}:
		\begin{subequations}
		\begin{align} \label{E:E0DIVID}
				\enzero[\Psi](t,u)
				+ \flzero[\Psi](t,u)
				& = \enzero[\Psi](0,u)
						\\
				& \ \ - 
							\int_{\mathcal{M}_{t,u}}
								\left\lbrace
									(1 + 2 \upmu) (\Lunit \Psi)
									+ 2 \Rad \Psi 
								\right\rbrace
								\waveinhom 
							\, d \vol
							\notag \\
				& \ \ - \frac{1}{2} 
						\int_{\mathcal{M}_{t,u}}
							\upmu \enmomtensor^{\alpha \beta}[\Psi] \deformarg{\Mult}{\alpha}{\beta}
							\, d \vol,
							\notag
		\end{align}
		\begin{align} \label{E:E1DIVID}
				\enone[\Psi](t,u)
				& + \flone[\Psi](t,u)
					\\
				& = \enone[\Psi](0,u)
					\notag \\
				& \ \ - \frac{1}{4} 
							\int_{\Sigma_0^u} 
								\left\lbrace
									\uLgood [\rgeo^2 \mytr \upchi]
									- \frac{1}{2} \rgeo^2 \upmu (\mytr \upchi)^2
									+ \rgeo^2 \mytr \upchi \mytr \angkuparg{(Trans-\Psi)}
									+ \rgeo^2 \upmu \mytr \upchi \mytr  \angkuparg{(Tan-\Psi)}
								\right\rbrace
								\Psi^2 
							\, d \vol
						\notag \\
				& \ \ + \frac{1}{4} 
							\int_{\Sigma_t^u} 
								\left\lbrace
									\uLgood [\rgeo^2 \mytr \upchi]
									- \frac{1}{2} \rgeo^2  \upmu (\mytr \upchi)^2
									+ \rgeo^2 \mytr \upchi \mytr  \angkuparg{(Trans-\Psi)}
									+ \rgeo^2 \upmu \mytr \upchi \mytr  \angkuparg{(Tan-\Psi)}
								\right\rbrace
								\Psi^2 
							\, d \vol
						\notag \\
			& \ \ - \int_{\mathcal{M}_{t,u}}
								\rgeo^2 
								\left\lbrace
									\Lunit \Psi 
									+ \frac{1}{2} \mytr \upchi \Psi
								\right\rbrace
								\waveinhom
								\,  d \vol
							\notag \\
			& \ \ - \frac{1}{2} \int_{\mathcal{M}_{t,u}}
								\upmu \enmomtensor^{\alpha \beta}[\Psi] 
									\left\lbrace
										\deformarg{\Mor}{\alpha}{\beta} 
										- \rgeo^2 \mytr \upchi g_{\alpha \beta}
									\right\rbrace
								\, d \vol
						\notag \\
			& \ \ + \frac{1}{4} \int_{\mathcal{M}_{t,u}}
								\upmu 
								\left\lbrace 
									\square_{g(\Psi)} [\rgeo^2 \mytr \upchi] 
								\right\rbrace 
								\Psi^2
							\, d \vol.
						\notag 
		\end{align}
	\end{subequations}
	In the above expressions, the $S_{t,u}$ tensorfields
	$\upchi,$ 
	$\angkuparg{(Trans-\Psi)},$ 
	and
	$\angkuparg{(Tan-\Psi)}$
	are defined by
	\eqref{E:CHIDEF},
	\eqref{E:KABTRANSVERSAL},
	and \eqref{E:KABGOOD}.	
		\end{proposition}

	\begin{proof}
			The overall strategy is to apply Lemma~\ref{L:WTUSPACETIMEDIVERGENCETHM} 
			with various currents $J$ from Def.~\ref{D:COMPATIBLECURRENTS}.
			Specifically, \eqref{E:E0DIVID} follows from applying Lemma~\ref{L:WTUSPACETIMEDIVERGENCETHM} 
			with the current
			$\Jenergycurrent{\Mult}_{\Timenormal}[\Psi]$
			defined in \eqref{E:MULTCURRENTDEF}. 
			The proof is simpler than the proof of \eqref{E:E1DIVID}, which we provide below.
			Hence, we omit the details of the proof of \eqref{E:E0DIVID}.
			
			To prove \eqref{E:E1DIVID}, we first apply Lemma~\ref{L:WTUSPACETIMEDIVERGENCETHM} with the current
			$\Jenergycurrent{\Mor+Correction}_{\Timenormal}[\Psi]$ defined in \eqref{E:MORPLUSMODCURRENT}.
			Using the fact that
			$\Psi \equiv 0$ along $\mathcal{C}_0^t,$  
			we deduce the following identity:
			\begin{align} \label{E:DIVTHMMORMOD}
				\int_{\Sigma_t^u} \upmu {\Jenergycurrent{\Mor+Correction}_{\Timenormal}[\Psi]} \, d \tvol
					+ \int_{\mathcal{C}_u^t} {\Jenergycurrent{\Mor}_{\Lunit}[\Psi]} \, d \conevol
				& =  \int_{\Sigma_0^u} \upmu {\Jenergycurrent{\Mor+Correction}_{\Timenormal}[\Psi]} \, d \tvol
					\\
				& \ \ - \int_{\mathcal{M}_{t,u}} \upmu \D_{\alpha} {\Jenergycurrent{\Mor+Correction}^{\alpha}[\Psi]} \, d \vol.
					\notag
			\end{align}	
			Examining \eqref{E:MORMODDIV}, \eqref{E:E1DEF}, \eqref{E:F1DEF},
			\eqref{E:E1DIVID}, and \eqref{E:DIVTHMMORMOD}, 
			we see that the desired
			identity \eqref{E:E1DIVID} will follow once we show that the last two $\Sigma_t^u$ 
			integrals on the right-hand side of \eqref{E:E1DEF} and 
			the last $\mathcal{C}_u^t$ integral on the right-hand side of \eqref{E:F1DEF} sum to
			\begin{align} \label{E:COMBINEDINTEGRALS}
				& \frac{1}{4} 
							\int_{\Sigma_t^u} 
								\left\lbrace
									\uLgood [\rgeo^2 \mytr \upchi]
									- \frac{1}{2} \upmu \rgeo^2 (\mytr \upchi)^2
									+ \rgeo^2 \mytr \upchi \mytr \angkuparg{(Trans-\Psi)}
									+ \upmu \rgeo^2 \mytr \upchi \mytr \angkuparg{(Tan-\Psi)}
								\right\rbrace
								\Psi^2 
							\, d \vol
							\\
			& - \frac{1}{4} 
							\int_{S_{0,u}} 
								\rgeo^2 \mytr \upchi \Psi^2
							\, d \spherevol.
						\notag
			\end{align}
			Note that the integral $- \frac{1}{4} 
							\int_{S_{0,u}} 
								\rgeo^2 \mytr \upchi \Psi^2
							\, d \spherevol$ 
			on the right-hand side of \eqref{E:COMBINEDINTEGRALS} does not depend on $t$ and hence cancels out of
			the expression \eqref{E:E1DIVID}.
			To deduce \eqref{E:COMBINEDINTEGRALS}, we
			first analyze the last two $\Sigma_t^u$ integrals on the right-hand side of \eqref{E:E1DEF}.
			To this end, we set $f := \rgeo^2 \mytr \upchi \Psi^2$ and integrate \eqref{E:UDERIVSTU} 
			in $u$ to compute that
			\begin{align} \label{E:STUDERIVINTEGRATED}
			& \int_{S_{t,u}} 
					\rgeo^2 \mytr \upchi \Psi^2 
				\, d \spherevol
				= \int_0^u \frac{\partial}{\partial u'} 
					\left(
						\int_{S_{t,u'}} 
						\rgeo^2 \mytr \upchi \Psi^2
						\, d \spherevol 
					\right) 
					\, du'
				\\
			& = 	\int_0^u \int_{S_{t,u'}}  
								\rgeo^2 \mytr \upchi \Rad (\Psi^2)
						\, d \spherevol \, du'
						\notag \\
			& \ \ 
				+ \int_0^u \int_{S_{t,u'}}  
							\left\lbrace
									\Rad [\rgeo^2 \mytr \upchi]
									- \upmu \rgeo^2 (\mytr \upchi)^2
									+ \rgeo^2 \mytr \upchi \mytr \angkuparg{(Trans-\Psi)} 
									+ \upmu \rgeo^2 \mytr \upchi \mytr \angkuparg{(Tan-\Psi)} 
							\right\rbrace \Psi^2
						\, d \spherevol \, du'.
						\notag
			\end{align}
			
			We now use \eqref{E:STUDERIVINTEGRATED}, 
			the identity $\upmu \Lunit = \uLgood - 2 \Rad,$
			and straightforward computations 
			to express the last two $\Sigma_t^u$ integrals on the right-hand side of \eqref{E:E1DEF}
			as follows:
			\begin{align} \label{E:E1INTEGRALSID}
				& - \frac{1}{4} 
						\int_{\Sigma_t^u} 
							\rgeo^2 \mytr \upchi \Rad (\Psi^2) 
						\, d \tvol
							\\
				& \ \
					+ \frac{1}{4} 
						\int_{\Sigma_t^u}
						\left\lbrace
							\upmu \Lunit [\rgeo^2 \mytr \upchi]
							+ \Rad [\rgeo^2 \mytr \upchi]
							+ \frac{1}{2} \upmu \rgeo^2 (\mytr \upchi)^2
						\right\rbrace 
				 		\Psi^2 
						\, d \tvol
					\notag	\\
				& = - \frac{1}{4} 
							\int_{S_{t,u}} 
								\rgeo^2 \mytr \upchi \Psi^2 
							\, d \spherevol
							\notag \\
				& \ \ 
						+ \frac{1}{4} 
							\int_{\Sigma_t^u} 
							\left\lbrace
								\uLgood [\rgeo^2 \mytr \upchi]
								- \frac{1}{2} \upmu \rgeo^2 (\mytr \upchi)^2
								+ \rgeo^2 \mytr \upchi \mytr \angkuparg{(Trans-\Psi)}
								+ \upmu \rgeo^2 \mytr \upchi \mytr \angkuparg{(Tan-\Psi)}
							\right\rbrace
							\Psi^2 \, d \vol.
							\notag
			\end{align}
			
			
			Similarly, we use	\eqref{E:OUTGOINGDERIVSTU} 
			with $f := \rgeo^2 \mytr \upchi \Psi^2$ to express
			the last $\mathcal{C}_u^t$ integral on the right-hand side of \eqref{E:F1DEF} as follows:
			\begin{align} \label{E:FLUXINTEGRALSID}
				\frac{1}{4} 
						\int_{\mathcal{C}_u^t} 
							\left\lbrace
								\Lunit[\rgeo^2 \mytr \upchi \Psi^2] 
								+ \rgeo^2 (\mytr \upchi)^2 \Psi^2
							\right\rbrace
						\, d \conevol
				 & = \frac{1}{4} 
				 			\int_{S_{t,u}} 
				 				\rgeo^2 \mytr \upchi \Psi^2 
				 			\, d \spherevol 
							- \frac{1}{4} 
							\int_{S_{0,u}} 
								\rgeo^2 \mytr \upchi \Psi^2
							\, d \spherevol.
			\end{align}	
			
			We now add \eqref{E:E1INTEGRALSID} and \eqref{E:FLUXINTEGRALSID} and note the cancellation of 
			the integrals $\frac{1}{4} \int_{S_{t,u}} \rgeo^2 \mytr \upchi \Psi^2 \, d \spherevol,$ 
			which yields the desired identity \eqref{E:COMBINEDINTEGRALS}.
			
	\end{proof}

\section{Integration by parts identities for the top-order \texorpdfstring{$L^2$}{square integral} estimates}	
In order to estimate some of the top-order $\waveinhom-$containing 
integrals on the right-hand side of \eqref{E:E1DIVID}, we need to 
derive alternate expressions via integration by parts. 
We use the following two simple integration by parts lemmas.
	
	\begin{lemma}[\textbf{Integration by parts on} $\mathcal{M}_{t,u}$] \label{L:KEYIBP}
		If $\Psi$ and $f$ are functions, then we have the following integration by parts identity:
		\begin{align} \label{E:KEYIBP}
			\int_{\mathcal{M}_{t,u}}
				f
				\left\lbrace
					\Lunit + \frac{1}{2} \mytr \upchi
				\right\rbrace
				\Psi
			 \, d \vol 
			 & = 
				- \int_{\mathcal{M}_{t,u}}
						\Psi
						\left\lbrace
							\Lunit + \frac{1}{2} \mytr \upchi 
						\right\rbrace 
						f
			 	 		\, d \vol
			 	- \int_{\Sigma_0^u} 
			 			\Psi f 
			 		\, d \tvol
			 	+ \int_{\Sigma_t^u} 
			 			\Psi f 
			 		\, d \tvol.
			\end{align}
		
	\end{lemma}

\begin{proof}
	From \eqref{E:OUTGOINGDERIVSTU}, we deduce that
	\begin{align} \label{E:FIRSTSTEPKEYIBP}
		\frac{\partial}{\partial t}
		\int_{\Sigma_t^u} 
			\Psi f 
		\, d \tvol
		& = 
		\int_{\Sigma_t^u} 
			\Psi
			\left\lbrace
				\Lunit f + \frac{1}{2} \mytr \upchi f
			\right\rbrace
		\, d \tvol
		+
		\int_{\Sigma_t^u} 
			f
			\left\lbrace
				\Lunit \Psi + \frac{1}{2} \mytr \upchi \Psi
			\right\rbrace
		\, d \tvol.
	\end{align}
	The desired identity \eqref{E:KEYIBP} now follows from integrating
	\eqref{E:FIRSTSTEPKEYIBP} in time from $0$ to $t.$

\end{proof}

\begin{lemma}[\textbf{Integration by parts on $S_{t,u}$}] \label{L:SPHEREIBP}
	If $Y$ is an $S_{t,u}-$tangent vectorfield and
	$\Psi$ and $f$ are functions, then we have the following integration by parts identity:
	\begin{align} \label{E:SPHEREIBP}
		\int_{S_{t,u}} 
			(Y \Psi) f
		\, d \spherevol
		& = - \int_{S_{t,u}} 
					(Y f) \Psi
				 \, d \spherevol
			- \frac{1}{2}
			\int_{S_{t,u}} 
				\mytr  \angdeform{Y} \Psi f
			\, d \spherevol.
	\end{align}

\end{lemma}

\begin{proof}
	By the Leibniz rule, we have
	\begin{align} \label{E:LEIBNIZRULEDIVERGENCEOFSCALARTIMESSTUTANGENTVECTORFIELD}
		\angdiv(f \Psi Y)
		& = (Y \Psi) f
			+ (Y f) \Psi
			+ \frac{1}{2} f \Psi (\ginversesphere)^{AB} 
				\left\lbrace
					\angD_A Y_B
					+ \angD_B Y_A
				\right\rbrace
			= (Y \Psi) f
			+ (Y f) \Psi
			+ \frac{1}{2} f \Psi \mytr  \angdeform{Y}.
	\end{align}
	The desired identity \eqref{E:SPHEREIBP} now follows from
	integrating
	\eqref{E:LEIBNIZRULEDIVERGENCEOFSCALARTIMESSTUTANGENTVECTORFIELD}
	over $S_{t,u}$ and noting that the $S_{t,u}$ integral of the
	angular divergence of an $S_{t,u}-$tangent vectorfield must vanish.
\end{proof}

As we mentioned above, some of the integrands on the right-hand side of \eqref{E:E1DIVID}
are not quite in a form that allows us to close our top-order $L^2$ estimates.
Hence, in order to bound certain top-order error integrals, we need 
to carry out some additional integrations by parts; 
see Sect.~\ref{S:DIFFICULTTOPORDERINTEGRALESTIMATES}.
We carry out these integrations by parts in the next two lemmas. 
The proofs involve lengthy computations 
but are conceptually very simple.

\begin{remark}[\textbf{Important integrals vs. negligible integrals}]
	In the next two lemmas, the
	explicitly written integrals on the right-hand sides 
	are the difficult ones that contribute to the 
	degeneracy of our top-order $L^2$ estimates,
	while the integrals featuring the integrands
	$\IBPerror{i}$
	and $\IBPseconderror{i}$
	have a negligible effect on the dynamics.
\end{remark}
	
	\begin{lemma}[\textbf{Integration by parts identity useful for bounding some top-order error integrals}]
\label{L:TOPORDERMORAWETZREORMALIZEDTRCHIIBP}
If $\newsmoothfunction$ is a function, then we have the following integration by parts identity:
\begin{align}  \label{E:TOPORDERMORAWETZREORMALIZEDTRCHIIBP}
			& - \int_{\mathcal{M}_{t,u}}
						(\Rad \Psi) 
						\left\lbrace
							\Lunit \mathscr{Z}^N \Psi 
							+ \frac{1}{2} \mytr \upchi \mathscr{Z}^N \Psi
						\right\rbrace
						\Rot 
						\newsmoothfunction
						\,  d \vol
						\\
			& =	 - \int_{\mathcal{M}_{t,u}}
								(\Rad \Psi)
								(\Lunit \newsmoothfunction)
								\Rot \mathscr{Z}^N \Psi
							\, d \vol
						\notag \\
			& \ \ + \int_{\Sigma_t^u}
								(\Rad \Psi)
								\newsmoothfunction
								\Rot \mathscr{Z}^N \Psi
							\, d \tvol
							\notag \\
			& \ \ + \int_{\mathcal{M}_{t,u}}
								\IBPerror{1}[\mathscr{Z}^N \Psi]
							\, d \vol
				+ \int_{\Sigma_t^u}
						\IBPerror{2}[\mathscr{Z}^N \Psi]
					\, d \tvol
					\notag \\
			& \ \ + \int_{\Sigma_0^u}
								\IBPerror{3}[\mathscr{Z}^N \Psi]
							\, d \tvol,
				\notag
	\end{align}
	where
	\begin{subequations}
	\begin{align} \label{E:NONDAMAGINGSPACETIMERENORMALZEDTRCHI}
		\IBPerror{1}[\mathscr{Z}^N \Psi]
		& := - 
								(\Rot \Rad \Psi)
								(\Lunit \newsmoothfunction)
								\mathscr{Z}^N \Psi
				- \frac{1}{2}
								\mytr \angdeform{\Rot}
								(\Rad \Psi)
								(\Lunit \newsmoothfunction)
								\mathscr{Z}^N \Psi
							\\
		& \ \     - 
								(\Rad \Psi)
								(\Lunit \newsmoothfunction)
								(\angdeformoneformupsharparg{\Rot}{\Lunit} \cdot \angdiff \mathscr{Z}^N \Psi)
							- 	
								(\angdeformoneformupsharparg{\Rot}{\Lunit} \cdot \angdiff \Rad \Psi)
								(\Lunit \newsmoothfunction)
								\mathscr{Z}^N \Psi
							\notag \\
		& \ \ - 		(\angdiv \angdeformoneformupsharparg{\Rot}{\Lunit})
								(\Rad \Psi)
								(\Lunit \newsmoothfunction)
								\mathscr{Z}^N \Psi
							\notag \\
		& \ \ -		\left\lbrace
								\Lunit (\Rad \Psi) 
								+ \frac{1}{2} \mytr \upchi (\Rad \Psi)
							\right\rbrace
							\newsmoothfunction
							\Rot \mathscr{Z}^N \Psi
							\notag \\
		& \ \ -		\left(
								\Rot 
								\left\lbrace
									\Lunit (\Rad \Psi) 
									+ \frac{1}{2} \mytr \upchi (\Rad \Psi)
								\right\rbrace
							\right)
							\newsmoothfunction
							\mathscr{Z}^N \Psi
							\notag \\
		& \ \ -		\frac{1}{2} 
							\mytr \angdeform{\Rot}
							\left\lbrace
								\Lunit (\Rad \Psi) 
								+ \frac{1}{2} \mytr \upchi (\Rad \Psi)
							\right\rbrace
							\newsmoothfunction
							\mathscr{Z}^N \Psi,
					\notag						
	\end{align}
	\begin{align} \label{E:NONDAMAGINGHYPERSURFACERENORMALZEDTRCHI}
		\IBPerror{2}[\mathscr{Z}^N \Psi]
		& := (\Rot \Rad \Psi) 
				\newsmoothfunction
				\mathscr{Z}^N \Psi
			+ \frac{1}{2} 
				\mytr \angdeform{\Rot}
				(\Rad \Psi)
				\newsmoothfunction
				\mathscr{Z}^N \Psi,
	\end{align}
	\begin{align}   \label{E:NONDAMAGINGINITIALDATARENORMALZEDTRCHI}
		\IBPerror{3}[\mathscr{Z}^N \Psi]
		& := - (\Rad \Psi) 
					\newsmoothfunction
					\Rot \mathscr{Z}^N
					- 
								(\Rad \Psi)
								(\Lunit \newsmoothfunction)
								(\angdeformoneformupsharparg{\Rot}{\Lunit} \cdot \angdiff \mathscr{Z}^N \Psi)
								\\
		& \ \ - \frac{1}{2}
					\mytr \angdeform{\Rot}
					(\Rot \Rad \Psi) 
					\newsmoothfunction
					\mathscr{Z}^N \Psi.
					\notag
	\end{align}
	\end{subequations}
	
\end{lemma}

\begin{proof}
We first use the identity \eqref{E:KEYIBP} with 
$\mathscr{Z}^N \Psi$ in the role of $\Psi$
and 
$(\Rad \Psi) \Rot \newsmoothfunction$
in the role of $f.$
We then split
$\left\lbrace
	\Lunit + \frac{1}{2} \mytr \upchi
\right\rbrace
\left\lbrace
	(\Rad \Psi) \Rot \newsmoothfunction
\right\rbrace$
$
= \left\lbrace
	\Lunit (\Rad \Psi) + \frac{1}{2} \mytr \upchi (\Rad \Psi)
\right\rbrace
\Rot \newsmoothfunction
+ 
(\Rad \Psi)
\Lunit \Rot \newsmoothfunction.
$
Furthermore, we commute the operators
$\Lunit$ and $\Rot$ on the previous term,
and then integrate by parts on the $S_{t,u}$ via Lemma~\ref{L:SPHEREIBP} to pull
the $\Rot$ off of $\newsmoothfunction$ and also the commutator
$[\Lunit, \Rot] = \angdeformoneformupsharparg{\Rot}{\Lunit}$
(see Lemma~\ref{L:VECTORFIELDCOMMUTATORS})
off of $\newsmoothfunction.$ In total, these steps lead to the lemma.

\end{proof}

\begin{lemma}[\textbf{A second integration by parts identity useful for bounding some top-order error integrals}]
\label{L:TOPORDERMORAWETZREORMALIZEDANGLAPUPMUIBP}
Let $Y$ be an $S_{t,u}-$tangent vectorfield and let $w=w(t,u)$ be a ``weight'' function. 
Then we have the following integration by parts identity:
\begin{align}  \label{E:TOPORDERMORAWETZREORMALIZEDANGLAPUPMUIBP}
			& - \int_{\mathcal{M}_{t,u}}
							w
							(\Rad \Psi) 
							\left\lbrace
								\Lunit \mathscr{Z}^N \Psi 
								+ \frac{1}{2} \mytr \upchi \mathscr{Z}^N \Psi
							\right\rbrace
							\angdiv Y
						\,  d \vol
						\\
			& =	 - \int_{\mathcal{M}_{t,u}}
								w
								(\Rad \Psi)
								\left(
									\left\lbrace
										\angLie_{\Lunit}
										+ \mytr \upchi
									\right\rbrace
									Y
								\cdot 
								\angdiff \mathscr{Z}^N \Psi
								\right)
					\, d \vol
						\notag \\
			& \ \ + \int_{\Sigma_t^u}
									w
									(\Rad \Psi) 
									(Y \cdot \angdiff \mathscr{Z}^N \Psi)
							\, d \tvol
							\notag \\
			& \ \ + \int_{\mathcal{M}_{t,u}}
						\IBPseconderror{1}[\mathscr{Z}^N;w]
					\, d \vol
				+ \int_{\Sigma_t^u}
						\IBPseconderror{2}[\mathscr{Z}^N;w]
					\, d \tvol
					\notag \\
			& \ \ + \int_{\Sigma_0^u}
								\IBPseconderror{3}[\mathscr{Z}^N;w]
							\, d \tvol,
				\notag
	\end{align}
	where
	\begin{subequations}
	\begin{align} \label{E:NONDAMAGINGSPACETIMERENORMALZEDANGLAPUPMU}
		\IBPseconderror{1}[\mathscr{Z}^N;w]
		& := - 
				 w
				(\angdiff \Rad \Psi)
				\left\lbrace
					\angLie_{\Lunit} Y
					+ \mytr \upchi Y
				\right\rbrace
				\mathscr{Z}^N \Psi
				\\
	& \ \ 
				-
				\left\lbrace
					\Lunit w 
					- \mytr \upchi w
				\right\rbrace
				(\Rad \Psi)
				(Y \cdot \angdiff \mathscr{Z}^N \Psi)
			\notag \\
		& \ \ 
				-
				\left\lbrace
					\Lunit w 
					- \mytr \upchi w
				\right\rbrace
				(Y \cdot \angdiff \Rad \Psi)
				\mathscr{Z}^N \Psi
			\notag \\
		& \ \ + w
				(\Rad \Psi)
				(Y \cdot \angdiff \mytr \upchi^{(Small)})
				\mathscr{Z}^N \Psi
				\notag \\
		& \ \
			- w
				\left\lbrace
					\Lunit \Rad \Psi 
					+ \frac{1}{2} \mytr \upchi \Rad \Psi
				\right\rbrace
				(Y \cdot \angdiff \mathscr{Z}^N \Psi)
				\notag \\
		& \ \ 
			- w
				\left\lbrace
					Y
					\cdot
					\left(
						\angdiff 
						\left\lbrace
							\Lunit \Rad \Psi 
							+ \frac{1}{2} \mytr \upchi \Rad \Psi
						\right\rbrace
					\right)
				\right\rbrace
				\mathscr{Z}^N \Psi,
				\notag \\
		\IBPseconderror{2}[\mathscr{Z}^N;w]
		& := w
				(Y \cdot \angdiff \Rad \Psi)
				\mathscr{Z}^N \Psi,	
				\label{E:NONDAMAGINGHYPERSURFACERENORMALZEDANGLAPUPMU}
					\\
		\IBPseconderror{3}[\mathscr{Z}^N;w]
		& := - w
					(\Rad \Psi) 
					(Y \cdot \angdiff \mathscr{Z}^N \Psi)
				- w
					(Y \cdot \angdiff \Rad \Psi)
					\mathscr{Z}^N \Psi.
					\label{E:NONDAMAGINGINITIALDATARENORMALZEDANGLAPUPMU}
	\end{align}
	\end{subequations}
	
\end{lemma}

\begin{proof}
We first use \eqref{E:KEYIBP} with
$\mathscr{Z}^N \Psi$ in the role of $\Psi$
and $f = f_1 f_2 f_3 := w (\Rad \Psi) \angdiv Y$
to remove the operator
$\left\lbrace
	\Lunit
		+ \frac{1}{2} \mytr \upchi
\right\rbrace$
off of $\mathscr{Z}^N \Psi$
on the left-hand side of \eqref{E:TOPORDERMORAWETZREORMALIZEDANGLAPUPMUIBP}.
After removal, we have the spacetime integrand factor
\begin{align} \label{E:TRIPLESPLITTING}
	& 
	(\mathscr{Z}^N \Psi)
	\left\lbrace
		\Lunit
		+ \frac{1}{2} \mytr \upchi
	\right\rbrace
(f_1 f_2 f_3)
	\\
& = 
	w
	(\Rad \Psi)
	(\mathscr{Z}^N \Psi)
	\left\lbrace
		\Lunit
			+ \mytr \upchi
	\right\rbrace
	\angdiv Y
	+ 
	(\Rad \Psi)
	(\mathscr{Z}^N \Psi)
	\left\lbrace
		\Lunit w
		- \mytr \upchi w
	\right\rbrace
	\angdiv Y
		\notag \\
	& \ \
	+
	w
	\left\lbrace
		\Lunit \Rad \Psi
		+ \frac{1}{2} \mytr \upchi \Rad \Psi
	\right\rbrace
	(\mathscr{Z}^N \Psi)
	\angdiv Y.
	\notag
\end{align}
We rewrite the first product on the right-hand side of \eqref{E:TRIPLESPLITTING} as
$ w
	(\Rad \Psi)
	(\mathscr{Z}^N \Psi)
	\angdiv 
				\left\lbrace
					\angLie_{\Lunit} Y
					+ \mytr \upchi Y
				\right\rbrace
$
with the help of Lemma~\ref{L:LPLUSTRCHIANGDIVCOMMUTE}.
Finally, we integrate by parts on the $S_{t,u}$ 
to remove all covariant angular derivatives off of
$\left\lbrace
					\angLie_{\Lunit} Y
					+ \mytr \upchi Y
\right\rbrace$
and $Y$ and use the fact that
$\angdiff \mytr \upchi = \angdiff \mytr \upchi^{(Small)}.$
In total, these steps lead to the identity \eqref{E:TOPORDERMORAWETZREORMALIZEDANGLAPUPMUIBP}.
\end{proof}

	\section{Error integrals arising from the deformation tensors of the multiplier vectorfields}	
	In the next definition, we give names to
	the energy error integrands corresponding to the 
	deformation tensors of the multiplier vectorfields
	(that is, the last integrand in \eqref{E:E0DIVID} and the next-to-last integrand in \eqref{E:E1DIVID}).
	
	\begin{definition}[\textbf{The error integrands corresponding to $\Mult$ and $\Mor$}] \label{D:LK0K1ERRORINTEGRANDS}
		Given a function $\Psi,$ we associate to it the following quantities
		$\basicenergyerror{\Mult}[\Psi]$
		and
		$\basicenergyerror{\Mor}[\Psi]:$
		\begin{subequations}
		\begin{align}
				 \basicenergyerror{\Mult}[\Psi] 
				 & := - \frac{1}{2} \upmu \enmomtensor^{\alpha \beta}[\Psi] \deformarg{\Mult}{\alpha}{\beta},
					\label{E:MULTERRORINT} \\
				 \basicenergyerror{\Mor}[\Psi]
				 & := - \frac{1}{2} \upmu \enmomtensor^{\alpha \beta}[\Psi]
								\left\lbrace
									\deformarg{\Mor}{\alpha}{\beta}
								 	- \rgeo^2 \mytr \upchi g
								\right\rbrace.
								\label{E:MORERRORINT}
		\end{align}
		\end{subequations}
	\end{definition}	
	
To derive our main results, we separately analyze
the positive and negative parts of $\Lunit \upmu.$ 
This motivates the following definition.

\begin{definition}[\textbf{$f_+$ and $f_-$}]
\label{D:POSITIVEANDNEGATIVEPARTS}
	Given any real-valued function $f,$ we decompose it into its positive part
	$f_+$ and its negative part $f_-$ as follows:
	\begin{align}
		f = f_+ - f_-,
	\end{align}
	where
	\begin{align}
		f_+ & := \max \lbrace f,0 \rbrace, 
			\qquad
		f_-:= \max \lbrace -f,0 \rbrace.
	\end{align}
\end{definition}

		In the next lemma, decompose the quantities
		$\basicenergyerror{\Mult}[\Psi]$
		and
		$\basicenergyerror{\Mor}[\Psi]$
		from Def.~\ref{D:LK0K1ERRORINTEGRANDS}
		relative to the rescaled null frame
		$\lbrace \Lunit, \uLgood, X_1, X_2 \rbrace.$
	
		\begin{lemma}[\textbf{Null decomposition of the error integrands corresponding to $\Mult$ and $\Mor$}] \label{L:LK0K1ERRORINTEGRANDS}
		Let $[\Lunit \upmu]_-$ and $[\Lunit \upmu]_+$ denote the negative and positive parts of $\Lunit \upmu$
		(see Def.~\ref{D:POSITIVEANDNEGATIVEPARTS}).
		Then the quantities 
		$\basicenergyerror{\Mult}[\Psi]$
		and
		$\basicenergyerror{\Mor}[\Psi]$
		from Def.~\ref{D:LK0K1ERRORINTEGRANDS}
		can be decomposed 
		relative to the rescaled null frame
		$\lbrace \Lunit, \uLgood, X_1, X_2 \rbrace$
		as follows:
		\begin{subequations}
		\begin{align}
				\basicenergyerror{\Mult}[\Psi] 
				& = \sum_{i=1}^6 \basicenergyerrorarg{\Mult}{i}[\Psi],
					\label{E:MULTAWETZENERGYERRORINTEGRANDS} \\
				\basicenergyerror{\Mor}[\Psi] 
				& = - \frac{1}{2}
							\rgeo^2
							|\angdiff \Psi|_{\gsphere}^2 
							[\Lunit \upmu]_-
						+ \sum_{i = 1}^4 \basicenergyerrorarg{\Mor}{i},
						\label{E:MORAWETZENERGYERRORINTEGRANDS}
		\end{align}	
		\end{subequations}
		where
		\begin{subequations}
		\begin{align}
			\basicenergyerrorarg{\Mult}{1}[\Psi] 
				& := \frac{1}{2} (\Lunit \Psi)^2
						\left\lbrace 
							\uLgood \upmu
							-(1 + \upmu) \Lunit \upmu
						\right\rbrace,
						\label{E:MULTERRORINTEGRAND1} \\
			\basicenergyerrorarg{\Mult}{2}[\Psi] 
			& := \frac{1}{2} |\angdiff \Psi|^2
							\left\lbrace
									(\uLgood \upmu + \Lunit \upmu)
								+ 3 \upmu \Lunit \upmu
							\right\rbrace,
						\label{E:MULTERRORINTEGRAND2} \\
			\basicenergyerrorarg{\Mult}{3}[\Psi] 
			& := - (\uLgood \Psi)
						 (\angdiffuparg{\#} \Psi) 
					\cdot	 
					\left\lbrace
						\angdiff \upmu 
						+ 2  \upzeta^{(Trans-\Psi)}
						+ 2\upmu \upzeta^{(Tan-\Psi)}
					\right\rbrace,
					\label{E:MULTERRORINTEGRAND3} \\
			\basicenergyerrorarg{\Mult}{4}[\Psi] 
			& :=  (\Lunit \Psi)
						(\angdiffuparg{\#} \Psi) 
					\cdot	
					\left\lbrace
						(1 - \upmu) \angdiff \upmu 
						+ 2 (1 + \upmu) \upzeta^{(Trans-\Psi)}
						+ 2 \upmu (1 + \upmu) \upzeta^{(Tan-\Psi)}
					\right\rbrace,
					\label{E:MULTERRORINTEGRAND4} \\
			\basicenergyerrorarg{\Mult}{5}[\Psi] 
			& := - (\upmu \angdiffuparg{\#} \Psi \hat{\otimes} \angdiffuparg{\#} \Psi)
				\cdot
				\left\lbrace
					\hat{\upchi}^{(Small)}
					+ 2 \hat{\angk}^{(Trans-\Psi)}
					+ 2 \upmu \hat{\angk}^{(Tan-\Psi)}
				\right\rbrace,
				\label{E:MULTERRORINTEGRAND5} \\
			\basicenergyerrorarg{\Mult}{6}[\Psi] 
			& := - \frac{1}{2} 
				(\Lunit \Psi) (\uLgood \Psi)
				\left\lbrace
					\mytr \upchi
					+ 2 \mytr  \angkuparg{(Trans-\Psi)}
					+ 2 \upmu \mytr  \angkuparg{(Tan-\Psi)}
				\right\rbrace,
				\label{E:MULTERRORINTEGRAND6}
			\end{align}	
		\end{subequations}
		and
		\begin{subequations}
		\begin{align}
			\basicenergyerrorarg{\Mor}{1}[\Psi] 
				& :=(\Lunit \Psi)^2
						\left\lbrace
							2 \rgeo \uLgood \rgeo
							- \rgeo^2 \Lunit \upmu
						\right\rbrace,
						\label{E:MORERRORINTEGRAND1} \\
				\basicenergyerrorarg{\Mor}{2}[\Psi] 
			& := 	\frac{1}{2} 
						\rgeo^2
						(\upmu |\angdiff \Psi|_{\gsphere}^2) 
						\left\lbrace
							\frac{[\Lunit \upmu]_+}{\upmu}
							- \mytr \upchi^{(Small)}
						\right\rbrace,
						\label{E:MORERRORINTEGRAND2} \\
		\basicenergyerrorarg{\Mor}{3}[\Psi] 
			& := \rgeo^2
				(\Lunit \Psi)
				(\angdiffuparg{\#} \Psi)
				\cdot
				\left\lbrace 
					 \angdiff \upmu
					+ 2 \upzeta^{(Trans-\Psi)}
					+ 2 \upmu \upzeta^{(Tan-\Psi)}
				\right\rbrace,
				\label{E:MORERRORINTEGRAND3} \\
			\basicenergyerrorarg{\Mor}{4}[\Psi] 
			& := - \rgeo^2  
						(\upmu \angdiffuparg{\#} \Psi \hat{\otimes} \angdiffuparg{\#} \Psi)
						\cdot
						\hat{\upchi}^{(Small)}.
						\label{E:MORERRORINTEGRAND4}
		\end{align}	
		\end{subequations}
	In the above expressions, the $S_{t,u}$ tensorfields
	$\upchi,$ 
	$\upchi^{(Small)},$
	$\upzeta^{(Trans-\Psi)},$
	$\angkuparg{(Trans-\Psi)},$ 
	$\upzeta^{(Tan-\Psi)},$
	and
	$\angkuparg{(Tan-\Psi)}$
	are defined by
	\eqref{E:CHIDEF},
	\eqref{E:CHIJUNKDEF},
	\eqref{E:ZETATRANSVERSAL},
	\eqref{E:KABTRANSVERSAL},
	\eqref{E:ZETAGOOD},
	and \eqref{E:KABGOOD}.		
	\end{lemma}
	
	\begin{remark}[\textbf{Isolating the Morawetz term}] 
	\label{R:ORIGINOFMORAWETZTERM}
		Note that in \eqref{E:MORAWETZENERGYERRORINTEGRANDS}, we have singled out 
		the non-positive term 
		$-  		\frac{1}{2} \rgeo^2
						(\upmu |\angdiff \Psi|_{\gsphere}^2) 
						\frac{[\Lunit \upmu]_-}{\upmu}$
		by ``removing it'' from $\basicenergyerrorarg{\Mor}{2}[\Psi].$ 
		The reason we have singled out this term is that it generates an important Morawetz-type
		spacetime integral that plays a fundamental role in our analysis;
		see Lemma~\ref{L:MORAWETZSPACETIMECOERCIVITY}.
	\end{remark}
		
	\begin{proof}
		The basic idea is to use the null decomposition \eqref{E:GINVERSENULLFRAME} of $g^{-1}$
		and the null decompositions of Lemmas \ref{L:NULLCOMPONENTSOFDEFORMATIONTENSORSOFMULTIPLIERS}
		and \ref{L:NULLCOMPONENTSOFDEFORMATIONTENSORSOFMULTIPLIERS} to expand
		the products \eqref{E:MULTERRORINT} and \eqref{E:MORERRORINT}
		relative to the frame $\lbrace \Lunit, \uLgood, X_1, X_2 \rbrace.$ 
		For example, to compute \eqref{E:MULTERRORINTEGRAND1}-\eqref{E:MULTERRORINTEGRAND6},
		we first use \eqref{E:GINVERSENULLFRAME} 
		to deduce that the
		left-hand side of \eqref{E:MULTERRORINT} can be expressed as
		\begin{align*}
		&
		- \frac{1}{2} \upmu 
		\left\lbrace
			- \frac{1}{2} 
				\upmu^{-1} 
				\left(
					\Lunit^{\alpha} \uLgood^{\kappa} 
					+ \Lunit^{\kappa} \uLgood^{\alpha}
				\right)
			+ (\ginversesphere)^{AB} X_A^{\alpha} X_B^{\kappa}
		\right\rbrace
			\\
		& \ \ \times
			\left\lbrace
			- \frac{1}{2} \upmu^{-1}
				\left(
					\Lunit^{\beta} \uLgood^{\lambda} 
					+ \Lunit^{\lambda} \uLgood^{\beta}
				\right)
			+ (\ginversesphere)^{AB} X_A^{\beta} X_B^{\lambda}
		\right\rbrace
		\enmomtensor_{\alpha \beta}[\Psi] \deformarg{\Mult}{\kappa}{\lambda}.
		\end{align*}
		We then fully expand the above product in terms of the null components of
		$\enmomtensor[\Psi]$ and $\deform{\Mult}$
		and use Lemmas \ref{L:NULLCOMPONENTSOFDEFORMATIONTENSORSOFMULTIPLIERS}
		and \ref{L:NULLCOMPONENTSOFDEFORMATIONTENSORSOFMULTIPLIERS} to substitute
		for these null components.
		Carrying out straightforward computations, we arrive at
		\eqref{E:MULTERRORINTEGRAND1}-\eqref{E:MULTERRORINTEGRAND6}.
		
		The proofs of \eqref{E:MORERRORINTEGRAND1}-\eqref{E:MORERRORINTEGRAND4}
		are similar. As we noted in Remark~\ref{R:ORIGINOFMORAWETZTERM},
		we removed the important term
		$-  \frac{1}{2} \rgeo^2
						(\upmu |\angdiff \Psi|_{\gsphere}^2) 
						\frac{[\Lunit \upmu]_-}{\upmu}$
		which arises as one of the terms in the product
		\[
			- \frac{1}{4} \upmu^{-1} \enmomtensor_{\Lunit \uLgood}[\Psi] 
			\left\lbrace 
				\deformarg{\Mor}{\Lunit}{\uLgood}
				- \rgeo^2 \mytr \upchi g(\Lunit, \uLgood)
			\right\rbrace,
		\]
		and placed it explicitly on the right-hand side of \eqref{E:MORAWETZENERGYERRORINTEGRANDS}.
	\end{proof}


\chapter{Avoiding Derivative Loss and Other Difficulties Via Modified Quantities}
\label{C:RENORMALIZEDEIKONALFUNCTIONQUANTITIES}
\thispagestyle{fancy}
We use the results of Ch.~\ref{C:RENORMALIZEDEIKONALFUNCTIONQUANTITIES}
only to close our top-order $L^2$ estimates.
Specifically, in this section, we derive transport equations for several
``modified'' versions of  the eikonal function quantities
$\mytr \upchi^{(Small)}$ and $\angdiff \upmu.$
When combined with elliptic estimates on the $S_{t,u},$
the modified quantities allow us to estimate some important top-order derivatives
of $\mytr \upchi^{(Small)}$ \emph{without losing derivatives.}
These steps are necessary even for proving an up-to-top-order 
local well-posedness result relative to the geometric coordinates.
The idea that one should work with modified quantities in order to close the
$L^2$ estimates for the eikonal function $u$ first appeared in Christodoulou-Klainerman's
proof of the global stability of the Minkowski spacetime as a solution to Einstein's equations
\cite{dCsK1993}. A similar idea was later employed by Klainerman-Rodnianski \cite{sKiR2003} in
their proof of low regularity local well-posedness for quasilinear wave equations
of the form $-\partial_t^2 \Psi + (g^{-1})^{ab}(\Psi) \partial_a \partial_b \Psi = \mathcal{N}(\Psi,\partial \Psi);$
their work was also heavily based on using an eikonal function to prove sharp estimates.
A similar idea was later used by Christodoulou in his shock formation monograph \cite{dC2007}.
In the latter two works, the key step is analogous to Cor.~\ref{C:ALPHARENORMALIZED} in the present work.
The main idea is that for solutions to $\square_g \Psi = 0,$
a certain component of the Ricci curvature tensor of $g,$
specifically $(\ginversesphere)^{AB} \Cur_{\Lunit A \Lunit B}$ 
(where $\Cur$ is the Riemann curvature tensor of $g$)
can be written as a perfect $\Lunit$
derivative of the first derivatives of $\Psi$ plus some error terms. The main point is that the 
$\Psi-$dependent error terms involve only $\Psi$ and its first derivatives, and 
\emph{not second-order derivatives of $\Psi$}, which are of principal order from the point of view of differentiability.
The perfect $\Lunit$ derivative structure is what allows us to construct suitably modified versions of
$\mytr \upchi^{(Small)}$ that do not lose derivatives. 

We now provide more details concerning this structure.
Roughly, we have $\Lunit \mytr \upchi^{(Small)} = \D^2 \Psi + Error(\Psi,\D \Psi).$
Integrating this identity along the integral curves of $\Lunit,$ 
we find only that $\mytr \upchi^{(Small)}$ has the same 
degree of differentiability as $\D^2 \Psi.$
We can commute this identity with vectorfields $Z \in \mathcal{Z}$
to obtain a similar result for the higher-order derivatives of $\mytr \upchi^{(Small)}.$
The problem is that the top-order result is $\mathscr{Z}^N \mytr \upchi^{(Small)} \sim \D^2 \mathscr{Z}^N \Psi,$
where the right-hand side involves one more derivative of $\Psi$ than our top-order energy-flux quantities allow us to control.
However, to close the top-order $L^2$ estimates corresponding to the timelike multiplier $\Mult,$
we need to be able to control $\mathscr{Z}^N \mytr \upchi^{(Small)}.$
To overcome the apparent derivative loss, we construct a ``fully modified'' 
version of $\mytr \upchi^{(Small)}$ 
(to be distinguished from the ``partially modified'' quantities discussed below)
by rescaling it and the adding to it certain terms connected to the perfect
$\Lunit$ derivative structure mentioned above. 
The new quantity verifies
$\Lunit (Modified) = Error(\Psi, \D \Psi),$ and we can integrate this identity
to conclude that $Modified$ has the same regularity
as $\D \Psi,$ a gain of one derivative.
Actually, these schematic expressions are somewhat misleading in the following sense:
one of the error terms on the right-hand side of the equation $\Lunit (Modified) = \cdots$
is roughly of the form
$|\hat{\upchi}|^2$ and cannot be eliminated through modification. 
At the top-order, the corresponding term would lead to derivative loss.
To control the top-order derivatives of $|\hat{\upchi}|^2$
we will, as we alluded to above, derive
a family of elliptic estimates on the $S_{t,u}$
(see Ch.~\ref{C:ELLIPTIC}). 
The elliptic estimates allow us to
control the top-order derivatives of 
$\hat{\upchi}$ in terms of quantities that do not lose derivatives.
This strategy was first employed
in \cite{dCsK1993} and later in \cite{sKiR2003} and \cite{dC2007}.
We also stress that in \cite{dC2007} and the present work, the structure of 
the top-order and the lower-order
terms are both important. In particular, 
we have to take care in keeping track of
factors of $\upmu$ and $\upmu^{-1},$ which can greatly
affect the dynamics near the shock.
We derive the transport equation verified by the lowest-order modified quantity 
in Cor.~\ref{C:FIRSTRENORMALIZEDTRANSPORTEQUATION}. 
We provide the higher-order versions of this equation
in Lemma~\ref{L:ONECOMMUTATIONTRANSPORTRENORMALIZEDTRCHIJUNK} and 
Prop.~\ref{P:TOPORDERTRCHIJUNKRENORMALIZEDTRANSPORT}.

In \cite{dC2007}, Christodoulou also encountered a derivative loss difficulty for
the term $\angLap \upmu,$ similar to one for $\mytr \upchi^{(Small)}$ described above.
To circumvent it, he derived an independent transport equation 
and corresponding estimates for a fully modified version of $\angLap \upmu,$
in analogy with his treatment of $\mytr \upchi^{(Small)}.$ 
Although his strategy solved the derivative loss problem, 
implementing it required a large number of complicated calculations 
and estimates. In the present monograph, we adopt an alternate 
strategy to handle this term. Specifically, 
in Chapter~\ref{C:TOPORDEREIKONALFUNCTIONCOMMUTATIONPOINTWISE}, 
we use commutation estimates to control
the relevant top-order derivatives of $\angLap \upmu$
in terms of related top-order derivatives of $\mytr \upchi^{(Small)}.$
In doing so, we avoid having to derive an 
analog of Christodoulou's transport equation
for a modified version of $\angLap \upmu,$
which saves a large amount of work.
Our strategy, when supplemented with elliptic estimates
for $\angLap \upmu,$
allows us to bound all top-order derivatives of $\upmu$ without loss, 
\emph{as long as the derivative 
operator involves an $\Lunit$ derivative 
or two angular derivatives};
see Remark~\ref{R:AVOIDDERIVATIVELOSS}.
Amazingly, the structure of the right-hand
side of \eqref{E:DIVCOMMUTATIONCURRENTDECOMPOSITION} together with the identities
of Prop.~\ref{P:DEFORMATIONTENSORFRAMECOMPONENTS} imply that,
thanks to the special properties of the vectorfields $Z \in \mathscr{Z},$
all top-order derivatives of $\upmu$ that appear in our equations
involve an $\Lunit$ derivative 
or two angular derivatives.

In order to close the top-order $L^2$ estimates
corresponding to the Morawetz multiplier $\Mor,$ we also need
to derive transport equations for ``partially modified'' versions of
$\mytr \upchi^{(Small)}$ and $\angdiff \upmu.$
The partial modifications are not connected to avoiding derivative loss,
but rather to avoiding certain error integrals that
have unfavorable time-growth properties. 
We derive these equations in Sects.~\ref{S:TRCHIJUNKPARTIALRENORMALIZATION} and \ref{S:ANGDIFFUPMUPARTIALRENORMALIZATION}.

\section{Preliminary structural identities}
In order to prove our sharp classical lifespan theorem, 
we must understand the sharp structure
of the transport equations verified by the modified quantities. 
As a first step, we analyze some components of the Riemann curvature tensor of $g.$ 

\begin{definition}[\textbf{The Riemann curvature tensor} $\Cur(W,X,Y,Z)$ \textbf{of} $g$]
\label{D:SPACETIMERIEMANN}
The Riemann curvature tensor $\Cur$ of the spacetime metric $g$ is
the type $\binom{0}{4}$ spacetime tensorfield defined by
\begin{align} \label{E:SPACETIMERIEMANN}
	g(\D_{WX}^2 Y - \D_{XW}^2 Y, Z)
	& = - \Cur(W,X,Y,Z),
\end{align}
where $W,$ $X,$ $Y,$ and $Z$ are arbitrary spacetime vectorfields.
In \eqref{E:SPACETIMERIEMANN}, 
$\D_{WX}^2 Y := W^{\alpha} X^{\beta} \D_{\alpha} \D_{\beta} Y.$
\end{definition}

In the next lemma, we compute the components of $\Cur$ relative to the
rectangular spacetime coordinates. This is a preliminary calculation
that will help us compute some of the components of $\Cur$
relative to the rescaled frame $\lbrace \Lunit, \Rad, X_1, X_2 \rbrace.$

\begin{lemma} [\textbf{Rectangular components of the Riemann curvature tensor of} $g$]
	\label{L:RECTANGULARCURVATURECOMPONENTS}
	Relative to the \textbf{rectangular} coordinate system, the components of $\Cur$
	can be expressed as
	\begin{align}
		\Cur_{\mu \nu \alpha \beta}
			& = \frac{1}{2} \Big\lbrace
					G_{\beta \mu} \D_{\alpha \nu}^2 \Psi
				+ G_{\alpha \nu} \D_{\beta \mu}^2 \Psi
				- G_{\beta \nu} \D_{\alpha \mu}^2 \Psi
				- G_{\alpha \mu} \D_{\beta \nu}^2 \Psi
				\Big \rbrace 
				\label{E:RECTANGULRCURVATURECOMPONENTS} \\
		& \ \ 
				+ 
				\frac{1}{4} 
				G_{\mu \alpha}
				G_{\nu \beta} 
				(g^{-1})^{\kappa \lambda} 
				(\partial_{\kappa} \Psi)
				(\partial_{\lambda} \Psi)
				- 
				\frac{1}{4} 
				G_{\mu \beta} 
				G_{\nu \alpha} 
				(g^{-1})^{\kappa \lambda}
				(\partial_{\kappa} \Psi)
				(\partial_{\lambda} \Psi)
			\notag \\
		& \ \ 
				+ 
				\frac{1}{4} 
				(g^{-1})^{\kappa \lambda} 
				[G_{\lambda (\nu} \partial_{\alpha)} \Psi]
				[G_{\kappa (\mu} \partial_{\beta)} \Psi]
				- 
				\frac{1}{4} 
				(g^{-1})^{\kappa \lambda} 
				[G_{\lambda (\mu} \partial_{\alpha)} \Psi]
				[G_{\kappa (\nu} \partial_{\beta)} \Psi]
			\notag
				\\
			& \ \ + \frac{1}{2} \Big\lbrace
					G_{\beta \mu}' (\partial_{\alpha} \Psi) (\partial_{\nu} \Psi)
				+ G_{\alpha \nu}' (\partial_{\beta} \Psi) (\partial_{\mu} \Psi)
				- G_{\beta \nu}' (\partial_{\alpha} \Psi) (\partial_{\mu} \Psi)
				- G_{\alpha \mu}' (\partial_{\beta} \Psi) (\partial_{\nu} \Psi)
				\Big \rbrace.
				\notag 	
\end{align}	

	In the above formula, $\D$ denotes the Levi-Civita connection of $g,$
	$\D^2 \Psi$ (which is a symmetric type $\binom{0}{2}$ tensorfield)
	denotes the second covariant derivative of $\Psi,$
	$G_{\alpha \beta} = G_{\alpha \beta}(\Psi) = \frac{d g_{\alpha \beta}(\Psi)}{d \Psi},$
	$G_{\alpha \beta}'(\Psi) = \frac{d G_{\alpha \beta}(\Psi)}{d \Psi},$
	$G_{(\mu \nu} \partial_{\lambda)} \Psi 
		: = 
			G_{\mu \nu} \partial_{\lambda} \Psi
		+ G_{\lambda \mu} \partial_{\nu} \Psi
		+ G_{\nu \lambda} \partial_{\mu} \Psi,$
	and
	$G_{\mu (\nu} \partial_{\lambda)} \Psi
		: = G_{\mu \nu} \partial_{\lambda} \Psi
			+ G_{\mu \lambda} \partial_{\nu} \Psi.$
\end{lemma}

\begin{proof}
It is a standard fact (see, for example, \cite{dC2007}*{Ch.4})
that relative to an arbitrary coordinate system, the components of the Riemann curvature tensor of $g$
can be expressed as
	\begin{align}
		\Cur_{\mu \nu \alpha \beta}
			& = \frac{1}{2} \Big\lbrace
					\D_{\alpha \nu}^2 g_{\beta \mu}
				+ \D_{\beta \mu}^2 g_{\alpha \nu}
				- \D_{\alpha \mu}^2 g_{\beta \nu}
				- \D_{\beta \nu}^2 g_{\alpha \mu}
				\Big \rbrace \\
		& \ \ + \frac{1}{8} (g^{-1})^{\kappa \lambda} \Big \lbrace
						\partial_{(\kappa} g_{\beta \mu)} [\partial_{(\alpha} g_{\nu) \lambda} - \partial_{\lambda} g_{\alpha \nu}]
					+ \partial_{(\kappa} g_{\alpha \nu)} [\partial_{(\beta} g_{\mu) \lambda} - \partial_{\lambda} g_{\beta \mu}]
					\Big \rbrace
			\notag \\
		& \ \ - \frac{1}{8} (g^{-1})^{\kappa \lambda} \Big \lbrace 
					\partial_{(\kappa} g_{\beta \nu)} [\partial_{(\alpha} g_{\mu) \lambda} - \partial_{\lambda} g_{\alpha \mu}]
				+ \partial_{(\kappa} g_{\alpha \mu)} [\partial_{(\beta} g_{\nu) \lambda} - \partial_{\lambda} g_{\beta \nu}]
				\Big \rbrace,
			\notag
	\end{align}
	where $\D$ denotes the Levi-Civita connection of $g,$
	the components $g_{\mu \nu}$ are treated as \textbf{scalar-valued functions} for the purpose
	of covariant differentiation,
	$\partial_{(\lambda} g_{\mu \nu)} 
		:=  
		\partial_{\lambda} g_{\mu \nu} 
		+ \partial_{\nu} g_{\lambda \mu} 
		+ \partial_{\mu} g_{\nu \lambda},$
	and 
	$\partial_{(\lambda} g_{\mu) \nu} 
	: = \partial_{\lambda} g_{\mu \nu} + \partial_{\mu} g_{\lambda \nu}.$	
	The desired expression \eqref{E:RECTANGULRCURVATURECOMPONENTS} now follows readily
	from straightforward computations.
	
\end{proof}

We now compute the curvature component $\Cur_{\Lunit A \Lunit B}.$

\begin{corollary}[\textbf{A $\upmu^{-1}-$regular expression for} $\Cur_{\Lunit A \Lunit B}$]
\label{C:ALPHAHASNOBADTERMS}
The curvature component $\Cur_{\Lunit A \Lunit B}$ can be expressed as follows:
\begin{align}
	\Cur_{\Lunit A \Lunit B}
	& = \frac{1}{2}
		\left\lbrace
			- \angGdoublearg{A}{B} \Lunit (\Lunit \Psi)
			+ \angGdoublearg{\Lunit}{A} \angdiffarg{B} \Lunit \Psi
			+ \angGdoublearg{\Lunit}{B} \angdiffarg{A} \Lunit \Psi
			- G_{\Lunit \Lunit} \angDsquaredarg{A}{B} \Psi
		\right\rbrace
			\label{E:ALPHAHASNOBADTERMS} \\
& \ \
		+ \frac{1}{2} \upmu^{-1} G_{\Lunit \Lunit} \upchi_{AB} \Rad \Psi
		- \frac{1}{2} \angGdoublearg{\Lunit}{A} \upchi_B^{\ C} \angdiffarg{C} \Psi
		- \frac{1}{2} \angGdoublearg{\Lunit}{B} \upchi_A^{\ C} \angdiffarg{C} \Psi
		\notag \\
& \ \   + 
				\myarray
				[G_{(Frame)}^2]
				{G_{(Frame)}'}
				\myarray
					[\Lunit \Psi]
					{\angdiff \Psi}^2
				+ 
				G_{(Frame)}^2
					\ginversesphere
					\myarray
					[\Lunit \Psi]
					{\angdiff \Psi}^2.
			\notag
\end{align}

\begin{proof}
We contract each side of \eqref{E:RECTANGULRCURVATURECOMPONENTS}
against $\Lunit^{\mu} X_A^{\nu} \Lunit^{\alpha} X_B^{\beta}$
and use the decompositions (see Lemma~\ref{L:SPACETIMEMETRICFRAMEVECTORFIELDS})
\begin{align} \label{E:GINVERSERECTDECOMP}
	(g^{-1})^{\kappa \lambda}
	& = (\ginversesphere)^{\kappa \lambda}
		- \Lunit^{\kappa} \Lunit^{\lambda}
		- \Lunit^{\kappa} \Radunit^{\lambda}
		- \Radunit^{\kappa} \Lunit^{\lambda}
		= (\ginversesphere)^{\kappa \lambda}
		-  \Lunit^{\kappa} \Lunit^{\lambda}
		- \upmu^{-1} \Lunit^{\kappa} \Rad^{\lambda}
		- \upmu^{-1} \Rad^{\kappa} \Lunit^{\lambda}.
\end{align}
to deduce that
\begin{align} \label{E:ALPHAFIRSTDECOMP}
\Cur_{\mu \nu \alpha \beta} \Lunit^{\mu} X_A^{\nu} \Lunit^{\alpha} X_B^{\beta}
& = \frac{1}{2}
		\left\lbrace
			- \angGdoublearg{A}{B} \D_{\Lunit \Lunit}^2 \Psi
			+ \angGdoublearg{\Lunit}{A} \D_{\Lunit B}^2 \Psi
			+ \angGdoublearg{\Lunit}{B} \D_{\Lunit A}^2 \Psi
			- G_{\Lunit \Lunit} \D_{AB}^2 \Psi
		\right\rbrace
		\\
& \ \ 
		+ \frac{1}{2}
		\upmu^{-1}	
		\left\lbrace 
				- 
				G_{\Lunit \Lunit} 
				\angGdoublearg{A}{B} 
				(\Lunit \Psi)
				\Rad \Psi
				+
				\angGdoublearg{\Lunit}{A} 
				\angGdoublearg{\Lunit}{B} 
				(\Lunit \Psi)
				\Rad \Psi
		\right\rbrace
			\notag \\
	& \ \ +
				\myarray
				[G_{(Frame)}^2]
				{G_{(Frame)}'}
				\myarray
					[\Lunit \Psi]
					{\angdiff \Psi}^2
				+ 
				G_{(Frame)}^2
					\ginversesphere
					\myarray
					[\Lunit \Psi]
					{\angdiff \Psi}^2.
					\notag
\end{align}
In deriving \eqref{E:ALPHAFIRSTDECOMP}, we used the fact that
the terms arising from the last two lines on the right-hand side of \eqref{E:RECTANGULRCURVATURECOMPONENTS} 
generate harmless lower-order terms that are part of the last term on the right-hand side of \eqref{E:ALPHAFIRSTDECOMP}.
The reason is that for the terms in these two lines, the $\lambda, \kappa$ indices are contractions and not differentiations
of $\Psi$
and hence we can use the decomposition 
$(g^{-1})^{\kappa \lambda} = (\ginversesphere)^{\kappa \lambda}
		- \Lunit^{\kappa} \Lunit^{\lambda}
		- \Lunit^{\kappa} \Radunit^{\lambda}
		- \Radunit^{\kappa} \Lunit^{\lambda}.$
Our goal is to show that the $\upmu^{-1}-$containing terms on the second
line of \eqref{E:ALPHAFIRSTDECOMP} are canceled by corresponding terms 
in first line, so that the sum of the first two lines contains no such terms.
To this end, we use Lemmas \ref{L:UPMUFIRSTTRANSPORT} and \ref{L:CONNECTIONLRADFRAME} to deduce
\begin{align}
	\D_{\Lunit \Lunit}^2 \Psi
	& = 	\Lunit \Lunit \Psi
			- \frac{1}{2} \upmu^{-1} G_{\Lunit \Lunit} (\Lunit \Psi)\Rad \Psi
			+ 	G_{(Frame)}
					(\Lunit \Psi)^2,
		\\
	\D_{\Lunit A}^2 \Psi
	& = \angdiffarg{A} \Lunit \Psi
			- \frac{1}{2} \upmu^{-1} \angGdoublearg{\Lunit}{A} (\Lunit \Psi) \Rad \Psi
			- \upchi_A^{\ B} \angdiffarg{B} \Psi
			+ 	G_{(Frame)}
					\myarray
					[\Lunit \Psi]
					{\angdiff \Psi}
					\Lunit \Psi,
		\label{E:DSQUAREDLAPSIINTERMSOFANGDIFFALPSI} \\
	\D_{AB}^2 \Psi
		& = \angDsquaredarg{A}{B} \Psi
			- \frac{1}{2} \upmu^{-1} \angGdoublearg{A}{B} (\Lunit \Psi) \Rad \Psi
			- \upmu^{-1} \upchi_{AB} \Rad \Psi
			+ 	G_{(Frame)}
					\myarray
					[\Lunit \Psi]
					{\angdiff \Psi}
					\Lunit \Psi.
					\label{E:DSQUAREDABPSIINTERMSOFANGDQUAREDAB}
\end{align}
Inserting these identities into \eqref{E:ALPHAFIRSTDECOMP}, we observe the desired
cancellations and thus arrive at \eqref{E:ALPHAHASNOBADTERMS}.

\end{proof}

\end{corollary}

We now show that up to lower-order terms,
$(\ginversesphere)^{AB} \Cur_{\Lunit A \Lunit B}$
can be written as a perfect $\Lunit$ derivative 
of the first derivatives of $\Psi.$ This is the key step
that will allow us to avoid losing derivatives in 
our estimates of the top-order derivatives of the eikonal function $u.$ 
This structure played a fundamental role in Christodoulou's work \cite{dC2007}.
Before \cite{dC2007}, Klainerman and Rodnianski had previously used this 
structure in their proof
of low regularity well-posedness results for quasilinear wave equations
\cite{sKiR2003}.

\begin{corollary}[\textbf{The key identity verified by the curvature component} $(\ginversesphere)^{AB} \Cur_{\Lunit A \Lunit B}$]
\label{C:ALPHARENORMALIZED}
Assume that $\square_{g(\Psi)} \Psi = 0.$
Then the curvature component $\upmu (\ginversesphere)^{AB} \Cur_{\Lunit A \Lunit B}$ can be expressed as 
a perfect $\Lunit$ derivative plus lower order terms with a favorable structure as follows:
\begin{align} \label{E:ALPHARENORMALIZED}
	\upmu
	(\ginversesphere)^{AB} \Cur_{\Lunit A \Lunit B}
	& = \Lunit
			\left\lbrace
				- G_{\Lunit \Lunit} \Rad \Psi
				- \frac{1}{2} \upmu \angGmixedarg{A}{A} \Lunit \Psi
				- \frac{1}{2} \upmu G_{\Lunit \Lunit} \Lunit \Psi 
				+ \upmu \angGmixedarg{\Lunit}{A} \angdiffarg{A} \Psi	
			\right\rbrace
		+ \mathfrak{A},
\end{align}
where $\mathfrak{A}$ has the following schematic structure:
\begin{align} \label{E:RENORMALIZEDALPHAINHOMOGENEOUSTERM}
	\mathfrak{A}
	& =  \sum_{p=0}^1
				\myarray
				[G_{(Frame)}^2 \ginversesphere]
				{G_{(Frame)}'}
				(\ginversesphere)^p
				\threemyarray[\upmu \Lunit \Psi]
					{\Rad \Psi}
					{\upmu \angdiff \Psi}
				\myarray
					[\Lunit \Psi]
					{\angdiff \Psi}.
\end{align}

Furthermore, without assuming $\square_{g(\Psi)} \Psi = 0,$
we have
\begin{align} \label{E:ALPHAPARTIALRENORMALIZED}
	(\ginversesphere)^{AB} \Cur_{\Lunit A \Lunit B}
	& = \frac{(\Lunit \upmu)}{\upmu} \mytr \upchi
		+ \Lunit
			\left\lbrace
				- \frac{1}{2} \angGmixedarg{A}{A} \Lunit \Psi
				- \frac{1}{2} G_{\Lunit \Lunit} \Lunit \Psi 
				+ \angGmixedarg{\Lunit}{A} \angdiffarg{A} \Psi
			\right\rbrace
		- \frac{1}{2} G_{\Lunit \Lunit} \angLap \Psi
		+ \mathfrak{B},
\end{align}
where $\mathfrak{B}$ has the following schematic structure:
\begin{align} \label{E:PARTIALRENORMALIZEDALPHAINHOMOGENEOUSTERM}
	\mathfrak{B}
	& = 	G_{(Frame)}
				\mytr 
				\upchi
				\Lunit \Psi
				+
				\sum_{p=0}^1
				\myarray
				[G_{(Frame)}^2 \ginversesphere]
				{G_{(Frame)}'}
				(\ginversesphere)^p
				\myarray[\Lunit \Psi]
					{\angdiff \Psi}^2.
\end{align}

\end{corollary}

\begin{proof}
We first prove \eqref{E:ALPHARENORMALIZED}.
To proceed, we contract \eqref{E:ALPHAHASNOBADTERMS} against $\upmu (\ginversesphere)^{AB}$
to deduce that
\begin{align} \label{E:ALPHABEFORERENORMALIZING}
	\upmu (\ginversesphere)^{AB} \Cur_{\Lunit A \Lunit B}
	& = - \frac{1}{2} \upmu \angGmixedarg{A}{A} \Lunit (\Lunit \Psi)
	    + \upmu \angGmixedarg{\Lunit}{A} \angdiffarg{A} \Lunit \Psi
			- \frac{1}{2} \upmu G_{\Lunit \Lunit} \angLap \Psi
		\\
& \ \
		+ \frac{1}{2} G_{\Lunit \Lunit} \mytr \upchi \Rad \Psi
		- \upmu \angGmixedarg{\Lunit}{A} \upchi_A^{\ B} \angdiffarg{B} \Psi
		\notag \\
& \ \ + \sum_{p=0}^1
				\myarray
				[G_{(Frame)}^2 \ginversesphere]
				{G_{(Frame)}'}
				(\ginversesphere)^p
				\myarray[\upmu \Lunit \Psi]
					{\upmu \angdiff \Psi}
				\myarray[\Lunit \Psi]
					{\angdiff \Psi}.
		\notag
\end{align}
Using Cor.~\ref{C:GFAMELDERIVATIVE}, 
the identity
$\Lunit \upmu = G_{(Frame)}
				\myarray
				[\upmu \Lunit \Psi]
				{\Rad \Psi}$
(that is, Lemma~\ref{L:UPMUFIRSTTRANSPORT}),				
and the identities $\angLie_{\Lunit} \gsphere_{AB} = 2 \upchi_{AB}$
and
$(\angLie_{\Lunit} \ginversesphere)^{AB} = - 2 \upchi^{AB},$
we rewrite the first two terms on the right-hand side of \eqref{E:ALPHABEFORERENORMALIZING} as
\begin{align} \label{E:ALPHAFIRSTTERMRENORMALIZED}
	- \frac{1}{2} \angGmixedarg{A}{A} \Lunit (\Lunit \Psi)
	& = - \frac{1}{2} 
				\Lunit
				\left\lbrace
					\upmu 
					\angGmixedarg{A}{A} 
					\Lunit \Psi
				\right\rbrace
		+ \sum_{p=0}^1
				\myarray
				[G_{(Frame)}^2 \ginversesphere]
				{G_{(Frame)}'}
				(\ginversesphere)^p
				\threemyarray[\upmu \Lunit \Psi]
					{\Rad \Psi}
					{\upmu \angdiff \Psi}
				\myarray
					[\Lunit \Psi]
					{\angdiff \Psi},
					\\
 \upmu \angGmixedarg{\Lunit}{A} \angdiffarg{A} \Lunit \Psi
	& = \Lunit
			\left\lbrace
				\upmu 
				\angGmixedarg{\Lunit}{A} 
				\angdiffarg{A} \Psi		
			\right\rbrace
	+ \upmu \angGmixedarg{\Lunit}{A} \upchi_A^{\ B} \angdiffarg{B} \Lunit \Psi
		\label{E:ALPHASECONDTERMRENORMALIZED}  
		\\
	& \ \
				+			
				\sum_{p=0}^1
				\myarray
				[G_{(Frame)}^2 \ginversesphere]
				{G_{(Frame)}'}
				(\ginversesphere)^p
				\threemyarray[\upmu \Lunit \Psi]
					{\Rad \Psi}
					{\upmu \angdiff \Psi}
				\myarray
					[\Lunit \Psi]
					{\angdiff \Psi}.
					\notag
\end{align}

Furthermore, from 
Cor.~\ref{C:GFAMELDERIVATIVE} and
\eqref{E:LONOUTSIDEGEOMETRICWAVEOPERATORFRAMEDECOMPOSED}, we deduce
that the third term on the right-hand side of \eqref{E:ALPHABEFORERENORMALIZING}
can be rewritten as
	\begin{align} \label{E:ALPHAWAVEEQUATIONTERMNORMALIZED}
		\frac{1}{2} 
		\upmu G_{\Lunit \Lunit} \angLap \Psi 
		& = \frac{1}{2} 
				\Lunit
				\left\lbrace
					G_{\Lunit \Lunit}
					(\upmu \Lunit \Psi + 2 \Rad \Psi)
				\right\rbrace
			+ \frac{1}{2} G_{\Lunit \Lunit} \mytr \upchi \Rad \Psi
			+ \myarray
				[G_{(Frame)}^2 \ginversesphere]
				{G_{(Frame)}'}
				\threemyarray[\upmu \Lunit \Psi]
					{\Rad \Psi}
					{\upmu \angdiff \Psi}
				\myarray
					[\Lunit \Psi]
					{\angdiff \Psi}.
	\end{align}

Inserting \eqref{E:ALPHAFIRSTTERMRENORMALIZED},
\eqref{E:ALPHASECONDTERMRENORMALIZED},
and \eqref{E:ALPHAWAVEEQUATIONTERMNORMALIZED}
into \eqref{E:ALPHABEFORERENORMALIZING},
we deduce the desired identities
\eqref{E:ALPHARENORMALIZED} and \eqref{E:RENORMALIZEDALPHAINHOMOGENEOUSTERM}.

To prove \eqref{E:ALPHAPARTIALRENORMALIZED}, we repeat the above argument,
but we do not use equation \eqref{E:ALPHAWAVEEQUATIONTERMNORMALIZED} to substitute
for the $\angLap \Psi$ term on the right-hand side of \eqref{E:ALPHABEFORERENORMALIZING}. 
This results in the presence of 
the term $\frac{1}{2} \upmu^{-1} G_{\Lunit \Lunit} \mytr \upchi \Rad \Psi,$
which arises from the first term on the second line of \eqref{E:ALPHABEFORERENORMALIZING},
on the right-hand side of \eqref{E:ALPHAPARTIALRENORMALIZED}.
Using Lemma~\ref{L:UPMUFIRSTTRANSPORT}, we can rewrite this term as
$\upmu^{-1} (\Lunit \upmu) \mytr \upchi$
plus a term that is schematically of the form of the right-hand side of \eqref{E:PARTIALRENORMALIZEDALPHAINHOMOGENEOUSTERM}. 
This explains the origin of the first term on the right-hand side of \eqref{E:ALPHAPARTIALRENORMALIZED}.

\end{proof}

\section{Full modification for \texorpdfstring{$\mytr \upchi^{(Small)}$}{the trace of the re-centered null second fundamental form}}
\label{S:FULLRENORMALIZATIONFORTRCHIJUNK}

We now define our fully modified version of $\mytr \upchi^{(Small)}.$

\begin{definition}[\textbf{Lowest-order fully modified version of} $\mytr \upchi^{(Small)}$]
\label{D:LOWESTORDERTRANSPORTRENORMALIZEDTRCHIJUNK}
We define
\begin{subequations}
\begin{align}
	\chifullmod
	& := \upmu \mytr \upchi^{(Small)} + \mathfrak{X},
		\label{E:LOWESTORDERTRANSPORTRENORMALIZEDTRCHIJUNK} \\
	\mathfrak{X}
	& := - G_{\Lunit \Lunit} \Rad \Psi
			- \frac{1}{2} \upmu \angGmixedarg{A}{A} \Lunit \Psi
			- \frac{1}{2} \upmu G_{\Lunit \Lunit} \Lunit \Psi
			+ \upmu \angGmixedarg{\Lunit}{A} \angdiffarg{A} \Psi.
			\label{E:LOWESTORDERTRANSPORTRENORMALIZEDTRCHIJUNKDISCREPANCY}
\end{align}	
\end{subequations}
\end{definition}

In the next corollary, we derive the transport equation verified by $\chifullmod.$

\begin{corollary}[\textbf{The key ``Raychaudhuri-type'' transport equation verified by the fully modified version of} $\mytr \upchi^{(Small)}$]
\label{C:FIRSTRENORMALIZEDTRANSPORTEQUATION}
Assume that $\square_{g(\Psi)} \Psi = 0.$
Let $\chifullmod$ be the fully modified quantity defined in \eqref{E:LOWESTORDERTRANSPORTRENORMALIZEDTRCHIJUNK}, 
and let $\mathfrak{A}$ be the inhomogeneous term defined in \eqref{E:RENORMALIZEDALPHAINHOMOGENEOUSTERM}.
Then $\chifullmod$ verifies the following transport equation:
\begin{align} \label{E:FIRSTRENORMALIZEDTRANSPORTEQUATION}
	\Lunit \chifullmod	
	& = 2 (\Lunit \upmu) \mytr \upchi
			- \frac{1}{2} \upmu (\mytr \upchi)^2
			- 2 \frac{1}{\rgeo} \Lunit \upmu
			+ 2 \frac{1}{\rgeo^2} \upmu 
			- \upmu |\hat{\upchi}^{(Small)}|^2
			- \mathfrak{A}.
\end{align}

\end{corollary}

\begin{proof}
	From the identity $\upchi_{AB} = g(\D_A \Lunit, X_B),$ 
	the fact that $[\Lunit, X_A]=0,$ the Def.~\ref{D:SPACETIMERIEMANN} of the Riemann curvature tensor, 
	the torsion-free property $\D_{\Lunit} X_B = \D_B \Lunit$ of $\D,$
	and Lemma~\ref{L:CONNECTIONLRADFRAME}, we deduce
	\begin{align} \label{E:CHITRANSPORT}
		\angLie_{\Lunit} \upchi_{AB}
		& =	\Lunit (\upchi_{AB})
			= \Lunit [g(\D_A \Lunit, X_B)]
			= g(\D_{\Lunit} (\D_A \Lunit), X_B)
			+ g(\D_A \Lunit, \D_{\Lunit} X_B)
			\\
		& = g(\D_A (\D_{\Lunit} \Lunit), X_B)
			- \Cur_{\Lunit A \Lunit B}
			+ g(\D_A \Lunit, \D_B \Lunit)
				\notag \\
		& = \frac{(\Lunit \upmu)}{\upmu} \upchi_{AB}
			+ \upchi_A^{\ C} \upchi_{BC}
			- \Cur_{\Lunit A \Lunit B}.
			\notag
	\end{align}
	Contracting the left-hand and right-hand sides of \eqref{E:CHITRANSPORT} with $(\ginversesphere)^{AB},$
	multiplying by $\upmu,$ and using the identity $(\angLie_{\Lunit} \ginversesphere)^{AB} = - 2 \upchi^{AB},$
	we deduce
	\begin{align} \label{E:STANDARDALPHARELATION}
		\upmu \Lunit \mytr \upchi
		& = (\Lunit \upmu) \mytr \upchi
			- \upmu |\upchi|^2
			- \upmu (\ginversesphere)^{AB} \Cur_{\Lunit A \Lunit B}.
	\end{align}
	From Cor.~\ref{C:ALPHARENORMALIZED}, we deduce that the last term on the right-hand side
	of \eqref{E:STANDARDALPHARELATION} can be expressed as
	$\upmu (\ginversesphere)^{AB} \Cur_{\Lunit A \Lunit B} = \Lunit \mathfrak{X} + \mathfrak{A}.$
	Equation \eqref{E:FIRSTRENORMALIZEDTRANSPORTEQUATION} now follows from the decompositions
	$\upchi = \rgeo^{-1} \gsphere + \upchi^{(Small)},$
	$|\upchi|^2 = \frac{1}{2}(\mytr \upchi)^2 + |\hat{\upchi}^{(Small)}|^2,$
	$\mytr \upchi = 2 \rgeo^{-1} + \mytr \upchi^{(Small)},$
	the identity $\Lunit \rgeo = 1,$
	and straightforward calculations.
	
\end{proof}

In the next lemma, we reveal the structure of the equation that arises after 
commuting the modified equation \eqref{E:FIRSTRENORMALIZEDTRANSPORTEQUATION}
with one spatial commutation vectorfield $S \in \mathscr{S}.$ 
The proof is somewhat delicate because
we have to carefully identify some important cancellations.
 
\begin{lemma}[\textbf{The transport equation for the fully modified version of} $S \mytr \upchi^{(Small)}$]
	\label{L:ONECOMMUTATIONTRANSPORTRENORMALIZEDTRCHIJUNK}
	Assume that $\square_{g(\Psi)} \Psi = 0.$
	Let $S \in \mathscr{S}$ be a spatial commutation vectorfield (see Def.~\ref{D:DEFSETOFSPATIALCOMMUTATORVECTORFIELDS}),
	and let $\mathfrak{X}$ be the quantity defined in \eqref{E:LOWESTORDERTRANSPORTRENORMALIZEDTRCHIJUNKDISCREPANCY}.
	We define
	\begin{align} \label{E:ONECOMMUTATIONTRANSPORTRENORMALIZEDTRCHIJUNK}
		\chifullmodarg{S}
		:= \upmu S \mytr \upchi^{(Small)} 
			+ S \mathfrak{X}.
	\end{align}
	Then the fully modified quantity $\chifullmodarg{S}$ verifies the following transport equation:
	\begin{align} \label{E:RENORMALIZEDTRANSPORTEQUATIONAFTERONECOMMUTATION}
	\Lunit \chifullmodarg{S}	
	& = - \mytr \upchi \chifullmodarg{S}
	  + \upmu [\Lunit, S] \mytr \upchi^{(Small)}
		+ 2 (\Lunit \upmu) S \mytr \upchi^{(Small)}
		+ \frac{1}{2} \mytr \upchi S \mathfrak{X}
		- S(\upmu |\hat{\upchi}^{(Small)}|^2)
		+ \chifullmodsourcearg{S},
\end{align}
where the inhomogeneous term $\chifullmodsourcearg{S}$ is equal to 
\begin{align} \label{E:ONECOMMUTATIONRENORMALIZEDTRCHIJUNKINHOMOGENEOUSTERM}
	\chifullmodsourcearg{S}
	& = - S \mathfrak{A}
		+ \frac{1}{2} \mytr \upchi
		S
		\overbrace
		{
		\left\lbrace
			- \frac{3}{2} \upmu G_{\Lunit \Lunit} \Lunit \Psi
			- 2 \upmu G_{\Lunit \Rad} \Lunit \Psi
			- \frac{1}{2} \upmu \angGmixedarg{A}{A} \Lunit \Psi
			+ \upmu \angGmixedarg{\Lunit}{A} \angdiffarg{A} \Psi	
		\right\rbrace
		}^{= 2 \Lunit \upmu + \mathfrak{X}}
		\\
	& \ \ 
		- (S \upmu)
		\left\lbrace 
			\Lunit \mytr \upchi^{(Small)} 
			+ 2 \frac{1}{\rgeo} \mytr \upchi^{(Small)}
			+ \frac{1}{2} (\mytr \upchi^{(Small)})^2
		\right\rbrace
		\notag \\
	& \ \ + [\Lunit, S] \mathfrak{X}
			- 2 \frac{1}{\rgeo^2} (\Lunit \upmu) S \rgeo
			+ 2 \upmu \frac{1}{\rgeo^2} \mytr \upchi^{(Small)} S \rgeo,
			\notag
\end{align}
and $\mathfrak{A}$ is the inhomogeneous term defined in \eqref{E:RENORMALIZEDALPHAINHOMOGENEOUSTERM}.

\end{lemma}


\begin{proof}
We apply $S$ to both sides of \eqref{E:FIRSTRENORMALIZEDTRANSPORTEQUATION}.
Clearly, the last term on the right-hand side of \eqref{E:FIRSTRENORMALIZEDTRANSPORTEQUATION}
gives rise to the first term on the right-hand side of \eqref{E:ONECOMMUTATIONRENORMALIZEDTRCHIJUNKINHOMOGENEOUSTERM}.
In our analysis of the remaining terms, we 
will draw boxes around the important terms that appear 
explicitly on either the left-hand or right-hand sides of equation \eqref{E:RENORMALIZEDTRANSPORTEQUATIONAFTERONECOMMUTATION}; 
the remaining terms are 
error terms that are by definition part of $\chifullmodsourcearg{S}.$

Referring to definitions \eqref{E:LOWESTORDERTRANSPORTRENORMALIZEDTRCHIJUNK} 
and \eqref{E:ONECOMMUTATIONTRANSPORTRENORMALIZEDTRCHIJUNK},
we compute that the term $S \Lunit \chifullmod$ arising from the left-hand
side of \eqref{E:FIRSTRENORMALIZEDTRANSPORTEQUATION} can be written as
\begin{align} \label{E:TERMLEFTDECOMPOSED}
	\boxed{\Lunit \chifullmodarg{S}}	
	+ \boxed{\upmu [S, \Lunit] \mytr \upchi^{(Small)}}
	+ [S, \Lunit] \mathfrak{X}
	+ (S \upmu) \Lunit \mytr \upchi^{(Small)}
	+ (S \Lunit \upmu) \mytr \upchi^{(Small)}.
\end{align}

Using the decomposition $\mytr \upchi = 2 \rgeo^{-1} + \mytr \upchi^{(Small)},$
we compute that the term 
$2 S [(\Lunit \upmu) \mytr \upchi]$ arising from the right-hand
side of \eqref{E:FIRSTRENORMALIZEDTRANSPORTEQUATION} can be written as
\begin{align} \label{E:TERM1DECOMPOSED}
	\boxed{2 (\Lunit \upmu) S \mytr \upchi^{(Small)}}
	- \frac{4}{\rgeo^2} (\Lunit \upmu) S \rgeo
	+ 2 (S \Lunit \upmu) \mytr \upchi.
\end{align}

Referring to definition \eqref{E:ONECOMMUTATIONTRANSPORTRENORMALIZEDTRCHIJUNK}
and using the decomposition $\mytr \upchi = 2 \rgeo^{-1} + \mytr \upchi^{(Small)},$
we compute that the term 
$- \frac{1}{2} S [\upmu (\mytr \upchi)^2]$ arising from the right-hand
side of \eqref{E:FIRSTRENORMALIZEDTRANSPORTEQUATION} can be written as
\begin{align} \label{E:TERM2DECOMPOSED}
& \boxed{- \mytr \upchi \chifullmodarg{S}}
	+ \boxed{\mytr \upchi S \mathfrak{X}}
	+ \frac{4}{\rgeo^3} \upmu S \rgeo
	+ \frac{2}{\rgeo^2} \upmu \mytr \upchi^{(Small)} S \rgeo
		\\
& \ \ 	
	- 2 \frac{1}{\rgeo^2}  S \upmu
	- 2 \frac{1}{\rgeo} \mytr \upchi^{(Small)} S \upmu 
	- \frac{1}{2} (S \upmu) (\mytr \upchi^{(Small)})^2.
	\notag
\end{align}  
We now explicitly place exactly half of the second term $\mytr \upchi S \mathfrak{X}$ in \eqref{E:TERM2DECOMPOSED}
on the right-hand side of \eqref{E:RENORMALIZEDTRANSPORTEQUATIONAFTERONECOMMUTATION}.
The other half is added to half of the last term on the right-hand side of \eqref{E:TERM1DECOMPOSED}
(that is, $(S \Lunit \upmu) \mytr \upchi$)
and the sum is placed on the right-hand side of \eqref{E:ONECOMMUTATIONRENORMALIZEDTRCHIJUNKINHOMOGENEOUSTERM}
as the second term.
We then use equations 
\eqref{E:UPMUFIRSTTRANSPORT} and
\eqref{E:LOWESTORDERTRANSPORTRENORMALIZEDTRCHIJUNKDISCREPANCY}
to put this sum into the form that appears on the right-hand side of 
\eqref{E:ONECOMMUTATIONRENORMALIZEDTRCHIJUNKINHOMOGENEOUSTERM}
(the important point is the cancellation of the factors $G_{\Lunit \Lunit} \Rad \Psi$).
The remaining part $(S \Lunit \upmu) \mytr \upchi$
of the last term on the right-hand side of \eqref{E:TERM1DECOMPOSED}
is completely canceled by the last term on the right-hand
side of \eqref{E:TERMLEFTDECOMPOSED} 
(which gets multiplied by $-1$ when it is moved over to the right-hand side of \eqref{E:RENORMALIZEDTRANSPORTEQUATIONAFTERONECOMMUTATION})
and the last term on the right-hand side of \eqref{E:TERM3DECOMPOSED} below
(the cancellation occurs because of the identity $\mytr \upchi = 2 \rgeo^{-1} + \mytr \upchi^{(Small)}$).

We compute that the term 
$- 2 S[ \frac{1}{\rgeo} \Lunit \upmu]$ arising from the right-hand
side of \eqref{E:FIRSTRENORMALIZEDTRANSPORTEQUATION} can be written as
\begin{align} \label{E:TERM3DECOMPOSED}
	2 \frac{1}{\rgeo^2} (\Lunit \upmu) S \rgeo
	- 2 \frac{1}{\rgeo} S \Lunit \upmu.
\end{align}

We compute that the term 
$2 S[\frac{1}{\rgeo^2} \upmu]$ arising from the right-hand
side of \eqref{E:FIRSTRENORMALIZEDTRANSPORTEQUATION} can be written as
\begin{align} \label{E:TERM4DECOMPOSED}
	2 \frac{1}{\rgeo^2} S \upmu
	- 4 \frac{1}{\rgeo^3} \upmu S \rgeo.
\end{align}

The term 
$- S[\upmu |\hat{\upchi}^{(Small)}|^2]$ arising from the right-hand
side of \eqref{E:FIRSTRENORMALIZEDTRANSPORTEQUATION} 
appears manifestly on the right-hand side of \eqref{E:RENORMALIZEDTRANSPORTEQUATIONAFTERONECOMMUTATION}.

Combining 
\eqref{E:TERMLEFTDECOMPOSED},
\eqref{E:TERM1DECOMPOSED},
\eqref{E:TERM2DECOMPOSED},
and \eqref{E:TERM3DECOMPOSED},
taking into account the remarks made just below
\eqref{E:TERM2DECOMPOSED}
and the previous sentence,
and carrying out straightforward computations,
we arrive at the desired equations 
\eqref{E:RENORMALIZEDTRANSPORTEQUATIONAFTERONECOMMUTATION}
and
\eqref{E:ONECOMMUTATIONRENORMALIZEDTRCHIJUNKINHOMOGENEOUSTERM}.



\end{proof}

We now define higher-order versions of the fully modified quantity $\chifullmod.$
We use the top-order version to close our top-order $L^2$ 
estimates corresponding to the timelike multiplier $\Mult.$

\begin{definition}[\textbf{Fully modified version of the pure spatial derivatives of} $\mytr \upchi^{(Small)}$]
\label{D:TRANSPORTRENORMALIZEDTRCHIJUNK}
Let $\mathscr{S}^N$ be an $N^{th}$ order pure spatial commutation vectorfield operator
(see Def.~\ref{D:DEFSETOFSPATIALCOMMUTATORVECTORFIELDS}),
and let $\mathfrak{X}$ be the quantity defined in
\eqref{E:LOWESTORDERTRANSPORTRENORMALIZEDTRCHIJUNKDISCREPANCY}.
We define the fully modified function $\chifullmodarg{\mathscr{S}^N}$ as follows:
\begin{align}
	\chifullmodarg{\mathscr{S}^N}
	& := \upmu \mathscr{S}^N \mytr \upchi^{(Small)} + \mathscr{S}^N \mathfrak{X}.
		\label{E:TRANSPORTRENORMALIZEDTRCHIJUNK}
\end{align}	

\end{definition}

In the next proposition, we reveal the structure of the equation that arises after 
commuting equation \eqref{E:FIRSTRENORMALIZEDTRANSPORTEQUATION} 
with $\mathscr{S}^N.$

\begin{proposition}[\textbf{The transport equation for the fully modified version of} $\mathscr{S}^N \mytr \upchi^{(Small)}$]
\label{P:TOPORDERTRCHIJUNKRENORMALIZEDTRANSPORT}
Assume that $\square_{g(\Psi)} \Psi = 0.$
Consider an $N^{th}$ order pure spatial commutation vectorfield operator
of the form $\mathscr{S}^N = \mathscr{S}^{N-1} S$
(see Def.~\ref{D:DEFSETOFSPATIALCOMMUTATORVECTORFIELDS}),
and let $\chifullmodarg{\mathscr{S}^N}$ and $\mathfrak{X}$ be the corresponding quantities 
defined in \eqref{E:TRANSPORTRENORMALIZEDTRCHIJUNK} and \eqref{E:LOWESTORDERTRANSPORTRENORMALIZEDTRCHIJUNKDISCREPANCY}.
Then the fully modified quantity $\chifullmodarg{\mathscr{S}^N}$ verifies the following transport equation:
\begin{align} \label{E:TOPORDERTRCHIJUNKRENORMALIZEDTRANSPORT}
\Lunit \chifullmodarg{\mathscr{S}^N}
	- \left\lbrace
			2 \frac{\Lunit \upmu}{\upmu}
			- \mytr \upchi 		
		\right\rbrace	\chifullmodarg{\mathscr{S}^N}
	&= 		\left\lbrace
					\frac{1}{2} \mytr \upchi
					- 2 \frac{\Lunit \upmu}{\upmu} 
				\right\rbrace
				\mathscr{S}^N \mathfrak{X}
		\\
	& \ \ + \upmu [\Lunit, \mathscr{S}^N] \mytr \upchi^{(Small)}
			- 2 \upmu \hat{\upchi}^{(Small)\# \#} \angLie_{\mathscr{S}}^N \hat{\upchi}^{(Small)}
			+ \chifullmodsourcearg{\mathscr{S}^N},
			\notag
\end{align}
where $\chifullmodsourcearg{\mathscr{S}^N}$ is the inhomogeneous term
\begin{align} \label{E:NTIMESCOMMUTATIONRENORMALIZEDTRCHIJUNKINHOMOGENEOUSTERM}
	\chifullmodsourcearg{\mathscr{S}^N}
	& := \mathscr{S}^{N-1} \chifullmodsourcearg{S}
		+ [\Lunit, \mathscr{S}^{N-1}] S \mathfrak{X}
			\\
	& \ \ 
			- [\mathscr{S}^{N-1}, \upmu \mytr \upchi] S \mytr \upchi^{(Small)}
			+ \frac{1}{2} [\mathscr{S}^{N-1}, \mytr \upchi] S \mathfrak{X}
			\notag \\
	& \ \ 
			+ [\mathscr{S}^{N-1}, \Lunit \upmu]
				\left(S \mytr \upchi^{(Small)} \right)
			- [\mathscr{S}^{N-1}, \upmu] 
				\left(S \Lunit \mytr \upchi^{(Small)} \right)
		\notag \\
	& \ \ 
		+ \left\lbrace
				2 \upmu \hat{\upchi}^{(Small) \# \#} \angLie_{\mathscr{S}}^N \hat{\upchi}^{(Small)}
				- \mathscr{S}^N (\upmu \hat{\upchi}^{(Small) \# \#} \hat{\upchi}^{(Small)})
			\right\rbrace,
		\notag
\end{align}
and $\chifullmodsourcearg{S}$ is the term given by \eqref{E:ONECOMMUTATIONRENORMALIZEDTRCHIJUNKINHOMOGENEOUSTERM}.

\end{proposition}

\begin{proof}
We apply $\mathscr{S}^{N-1}$ to the transport equation \eqref{E:RENORMALIZEDTRANSPORTEQUATIONAFTERONECOMMUTATION}
for $\chifullmodarg{S}.$ Clearly, the last term on the right-hand side
of \eqref{E:RENORMALIZEDTRANSPORTEQUATIONAFTERONECOMMUTATION}
gives rise to the first term on the right-hand side of \eqref{E:NTIMESCOMMUTATIONRENORMALIZEDTRCHIJUNKINHOMOGENEOUSTERM}.
In our analysis of the remaining terms, we 
draw boxes around the important terms that appear 
explicitly on either the left-hand or right-hand sides of equation \eqref{E:TOPORDERTRCHIJUNKRENORMALIZEDTRANSPORT}; 
the remaining terms are either canceled by other terms or are
error terms that are by definition part of $\chifullmodsourcearg{\mathscr{S}^N}.$
To begin, we consider definitions 
\eqref{E:ONECOMMUTATIONTRANSPORTRENORMALIZEDTRCHIJUNK}
and \eqref{E:TRANSPORTRENORMALIZEDTRCHIJUNK}
and compute that the corresponding left-hand side of $\mathscr{S}^{N-1}$
applied to \eqref{E:RENORMALIZEDTRANSPORTEQUATIONAFTERONECOMMUTATION}
(that is, $\mathscr{S}^{N-1} \Lunit \chifullmodarg{S}$) 
can be written as
\begin{align} \label{E:LEFTHANDSIDETERMDECOMPOSED}
	& \boxed{\Lunit \chifullmodarg{\mathscr{S}^N}}
		+ \boxed{\upmu  [\mathscr{S}^N, \Lunit] \mytr \upchi^{(Small)}}
		+ [\mathscr{S}^{N-1}, \Lunit \upmu]
			\left(S \mytr \upchi^{(Small)} \right)
		+ [\mathscr{S}^{N-1}, \upmu]
			\left(\Lunit S \mytr \upchi^{(Small)} \right)	
		\\
	& \ \
		+ \upmu \mathscr{S}^{N-1} 
			\left([\Lunit, S] \mytr \upchi^{(Small)} \right)
		+ [\mathscr{S}^{N-1}, \Lunit] S \mathfrak{X}.
		\notag
\end{align}

The term $- \mathscr{S}^{N-1} (\mytr \upchi \chifullmodarg{S})$
arising from the right-hand side of \eqref{E:RENORMALIZEDTRANSPORTEQUATIONAFTERONECOMMUTATION}
can be written as
\begin{align} \label{E:FIRSTTERMDECOMPOSED}
	- \boxed{\mytr \upchi \chifullmodarg{\mathscr{S}^N}}
	- [\mathscr{S}^{N-1}, \upmu \mytr \upchi] S \mytr \upchi^{(Small)}.
\end{align}	
	
The term $\mathscr{S}^{N-1} \left\lbrace \upmu [\Lunit, S] \mytr \upchi^{(Small)} \right\rbrace$
arising from the right-hand side of \eqref{E:RENORMALIZEDTRANSPORTEQUATIONAFTERONECOMMUTATION}
can be written as
\begin{align} \label{E:SECONDTERMDECOMPOSED}
	\upmu \mathscr{S}^{N-1} \left([\Lunit, S] \mytr \upchi^{(Small)}\right)
	+ [\mathscr{S}^{N-1}, \upmu] \left([\Lunit, S] \mytr \upchi^{(Small)}\right).
\end{align}	

The term $2 \mathscr{S}^{N-1} \lbrace (\Lunit \upmu) S \mytr \upchi^{(Small)} \rbrace$
arising from the right-hand side of \eqref{E:RENORMALIZEDTRANSPORTEQUATIONAFTERONECOMMUTATION}
can be written as
\begin{align} \label{E:THIRDTERMDECOMPOSED}
	\boxed{2 \frac{\Lunit \upmu}{\upmu} \chifullmodarg{\mathscr{S}^N}}
	- \boxed{2 \frac{\Lunit \upmu}{\upmu} \mathscr{S}^N \mathfrak{X}}
	+ 2 [\mathscr{S}^{N-1}, \Lunit \upmu] S \mytr \upchi^{(Small)}.
\end{align}

The term $\frac{1}{2} \mathscr{S}^{N-1} (\mytr \upchi S \mathfrak{X})$
arising from the right-hand side of \eqref{E:RENORMALIZEDTRANSPORTEQUATIONAFTERONECOMMUTATION}
can be written as
\begin{align} \label{E:FOURTHTERMDECOMPOSED}
	\boxed{\frac{1}{2} \mytr \upchi \mathscr{S}^N \mathfrak{X}}
	+ \frac{1}{2} [\mathscr{S}^{N-1}, \mytr \upchi] S \mathfrak{X}.
\end{align}

The term $- \mathscr{S}^{N-1} S(\upmu |\hat{\upchi}^{(Small)}|^2)$ 
arising from the right-hand side of \eqref{E:RENORMALIZEDTRANSPORTEQUATIONAFTERONECOMMUTATION}
can be written as
\begin{align} \label{E:FIFTHTERMDECOMPOSED}
	- \boxed{2 \upmu \hat{\upchi}^{(Small) \# \#} \angLie_{\mathscr{S}}^N \hat{\upchi}^{(Small)}}
	+ \left\lbrace
			2 \upmu \hat{\upchi}^{(Small) \# \#} \angLie_{\mathscr{S}}^N \hat{\upchi}^{(Small)}
			-
			\mathscr{S}^N (\upmu \hat{\upchi}^{(Small) \# \#} \hat{\upchi}^{(Small)})
		\right\rbrace.
\end{align}	

Combining
\eqref{E:LEFTHANDSIDETERMDECOMPOSED},
\eqref{E:FIRSTTERMDECOMPOSED},
\eqref{E:SECONDTERMDECOMPOSED},
\eqref{E:THIRDTERMDECOMPOSED},
\eqref{E:FOURTHTERMDECOMPOSED},
and
\eqref{E:FIFTHTERMDECOMPOSED},
we deduce the desired equations 
\eqref{E:TOPORDERTRCHIJUNKRENORMALIZEDTRANSPORT} 
and \eqref{E:NTIMESCOMMUTATIONRENORMALIZEDTRCHIJUNKINHOMOGENEOUSTERM}.

\end{proof}

\section{Partial modification of \texorpdfstring{$\mytr \upchi^{(Small)}$}{the trace of the re-centered null second fundamental form}}
\label{S:TRCHIJUNKPARTIALRENORMALIZATION}

To close our top-order $L^2$ estimates involving the Morawetz multiplier $\Mor,$
we need to use the partially modified version of $\mytr \upchi^{(Small)}$
given in the next definition. The modified quantity allows us to avoid certain spacetime
error integrals with damaging time growth.

\begin{definition}[\textbf{Lowest-order partially modified version of} $\mytr \upchi^{(Small)}$]
\label{D:LOWESTORDERTRANSPORTREPARTIALNORMALIZEDTRCHIJUNK}

We define the scalar-valued functions $\chipartialmod$ and $\chipartialmodinhom$ as follows:
\begin{subequations}
\begin{align}
	\chipartialmod
	& := \mytr \upchi^{(Small)} + \chipartialmodinhom,
		\label{E:LOWESTORDERTRANSPORTPARTIALRENORMALIZEDTRCHIJUNK} \\
	\chipartialmodinhom
	& := - \frac{1}{2} \angGmixedarg{A}{A} \Lunit \Psi
			 - \frac{1}{2} G_{\Lunit \Lunit} \Lunit \Psi
			 + \angGmixedarg{\Lunit}{A} \angdiffarg{A} \Psi.
			\label{E:LOWESTORDERTRANSPORTPARTIALRENORMALIZEDTRCHIJUNKDISCREPANCY}
\end{align}	
\end{subequations}

\end{definition}

In the next lemma, we derive the transport equation verified by $\chipartialmod.$

\begin{lemma}[\textbf{The transport equation verified by the partially modified version of} $\mytr \upchi^{(Small)}$]
\label{L:TRCHIJUNKFIRSTPARTIALRENORMALIZEDTRANSPORTEQUATION}
Let $\chipartialmod$ be the partially modified quantity defined in \eqref{E:LOWESTORDERTRANSPORTPARTIALRENORMALIZEDTRCHIJUNK},
and let $\mathfrak{B}$ be the term defined in \eqref{E:PARTIALRENORMALIZEDALPHAINHOMOGENEOUSTERM}.
Then $\rgeo^2 \chipartialmod$ verifies the following transport equation:
\begin{align} \label{E:TRCHIJUNKFIRSTPARTIALRENORMALIZEDTRANSPORTEQUATION}
	\Lunit (\rgeo^2 \chipartialmod)
	& = - \frac{1}{2} \rgeo^2 G_{\Lunit \Lunit} \angLap \Psi
			- \rgeo^2 |\upchi^{(Small)}|^2
			- \rgeo^2 \mathfrak{B}.
\end{align}

\end{lemma}

\begin{proof}
From equations \eqref{E:ALPHAPARTIALRENORMALIZED} and \eqref{E:STANDARDALPHARELATION}, we deduce that
\begin{align} \label{E:AGAINSTANDARDALPHARELATION}
		\Lunit 
			\left\lbrace
				\mytr \upchi
				- \frac{1}{2} \angGmixedarg{A}{A} \Lunit \Psi
				- \frac{1}{2} G_{\Lunit \Lunit} \Lunit \Psi 
				+ \angGmixedarg{\Lunit}{A} \angdiffarg{A} \Psi
			\right\rbrace
		& = - \frac{1}{2} G_{\Lunit \Lunit} \angLap \Psi
				+ \mathfrak{B}
				- |\upchi|^2.
\end{align}
The desired equation \eqref{E:TRCHIJUNKFIRSTPARTIALRENORMALIZEDTRANSPORTEQUATION}
now follows from \eqref{E:AGAINSTANDARDALPHARELATION},
the decompositions
$\upchi = \rgeo^{-1} \gsphere + \upchi^{(Small)},$
$|\upchi|^2 = \frac{1}{2}(\mytr \upchi)^2 + |\hat{\upchi}^{(Small)}|^2,$
$\mytr \upchi = 2 \rgeo^{-1} + \mytr \upchi^{(Small)},$
the fact that $\Lunit \rgeo = 1,$
and from straightforward computations.

\end{proof}

We now define higher-order versions of $\chipartialmod$ and $\chipartialmodinhom.$

\begin{definition}[\textbf{Partially modified version of} $\mathscr{Z}^N \mytr \upchi^{(Small)}$]
\label{D:TRANSPORTREPARTIALNORMALIZEDTRCHIJUNK}
We define the partially modified function $\chipartialmodarg{\mathscr{Z}^N}$ as follows:
\begin{subequations}
\begin{align}
	\chipartialmodarg{\mathscr{Z}^N}
	& := \mathscr{Z}^N \mytr \upchi^{(Small)} 
		+ \chipartialmodinhomarg{\mathscr{Z}^N},
		\label{E:TRANSPORTPARTIALRENORMALIZEDTRCHIJUNK} \\
	\chipartialmodinhomarg{\mathscr{Z}^N}
	& := - \frac{1}{2} \angGmixedarg{A}{A} \Lunit \mathscr{Z}^N \Psi
			- \frac{1}{2} G_{\Lunit \Lunit} \Lunit \mathscr{Z}^N \Psi
			+ \angGmixedarg{\Lunit}{A} \angdiffarg{A} \mathscr{Z}^N \Psi.
			\label{E:TRANSPORTPARTIALRENORMALIZEDTRCHIJUNKDISCREPANCY}
\end{align}	
\end{subequations}

\end{definition}

In the next lemma, we derive the transport equation verified by $\chipartialmodarg{\mathscr{S}^{N-1}}.$

\begin{lemma}[\textbf{The transport equation verified by the partially modified version of} 
$\mathscr{S}^{N-1} \mytr \upchi^{(Small)}$]
\label{L:COMMUTEDTRCHIJUNKFIRSTPARTIALRENORMALIZEDTRANSPORTEQUATION}
Let $\mathscr{S}^{N-1}$ be an $(N-1)^{st}$ order pure spatial commutation vectorfield operator
and let $\chipartialmodarg{\mathscr{S}^{N-1}}$ 
be the partially modified quantity
defined in \eqref{E:TRANSPORTPARTIALRENORMALIZEDTRCHIJUNK}.
Then $\chipartialmodarg{\mathscr{S}^{N-1}}$ verifies the following transport equation:
\begin{align} \label{E:COMMUTEDTRCHIJUNKFIRSTPARTIALRENORMALIZEDTRANSPORTEQUATION}
	\Lunit (\rgeo^2 \chipartialmodarg{\mathscr{S}^{N-1}})
	& = - \frac{1}{2} \rgeo^2 G_{\Lunit \Lunit} \angLap \mathscr{S}^{N-1} \Psi
			+ \chipartialmodsourcearg{\mathscr{S}^{N-1}},
\end{align}
where the inhomogeneous term $\chipartialmodsourcearg{\mathscr{S}^{N-1}}$ is given by
\begin{align}  \label{E:TRCHIJUNKCOMMUTEDTRANSPORTEQNPARTIALRENORMALIZATIONINHOMOGENEOUSTERM}
\chipartialmodsourcearg{\mathscr{S}^{N-1}}
	& = - \mathscr{S}^{N-1} (\rgeo^2 \mathfrak{B})
			- \mathscr{S}^{N-1} (\rgeo^2 |\upchi^{(Small)}|^2)
				\\
	& \ \ 
			- \frac{1}{2} [\rgeo^2 G_{\Lunit \Lunit}, \mathscr{S}^{N-1}] \angLap \Psi
			- \frac{1}{2} \rgeo^2 G_{\Lunit \Lunit} [\mathscr{S}^{N-1}, \angLap] \Psi
			+ [\Lunit, \mathscr{S}^{N-1}](\rgeo^2 \mytr \upchi^{(Small)})
			  \notag \\
	& \ \ 
				+ [\Lunit, \mathscr{S}^{N-1}](\rgeo^2 \chipartialmodinhom)
				+ 
				\Lunit 
				\left\lbrace
					[\rgeo^2 , \mathscr{S}^{N-1}] \mytr \upchi^{(Small)}
				\right\rbrace
				+ 
				\Lunit	
				\left\lbrace
					\rgeo^2 \chipartialmodinhomarg{\mathscr{S}^{N-1}}
					- \mathscr{S}^{N-1} (\rgeo^2 \chipartialmodinhom)
				\right\rbrace,
				\notag
\end{align}	
$\mathfrak{B}$ is defined in \eqref{E:PARTIALRENORMALIZEDALPHAINHOMOGENEOUSTERM},
$\chipartialmodinhom$ is defined in \eqref{E:LOWESTORDERTRANSPORTPARTIALRENORMALIZEDTRCHIJUNKDISCREPANCY},
and $\chipartialmodinhomarg{\mathscr{S}^{N-1}}$ 
is defined in \eqref{E:TRANSPORTPARTIALRENORMALIZEDTRCHIJUNKDISCREPANCY}.

\end{lemma}

\begin{proof}
	We apply $\mathscr{S}^{N-1}$ to each side of 
	equation \eqref{E:TRCHIJUNKFIRSTPARTIALRENORMALIZEDTRANSPORTEQUATION} 
	and perform several straightforward operator commutations.
\end{proof}

\section{Partial modification of \texorpdfstring{$\angdiff \upmu$}{the angular gradient of the inverse foliation density}}
\label{S:ANGDIFFUPMUPARTIALRENORMALIZATION}
This section complements Sect.~\ref{S:TRCHIJUNKPARTIALRENORMALIZATION}.
Specifically, we define partially modified versions of $\angdiff \upmu$
in order to close our top-order $L^2$ estimates corresponding to the multiplier $\Mor.$
We start by providing the definition of the lowest order partially modified
version of $\angdiff \upmu.$

\begin{definition}[\textbf{Lowest order partially modified version of} $\angdiff \upmu$]
\label{D:LOWESTORDERTRANSPORTREPARTIALNORMALIZEDUPMU}

We define the $S_{t,u}$ one-forms $\mupartialmod$ and $\mupartialmodinhom$ as follows:
\begin{subequations}
\begin{align}
	\mupartialmod
	& := \angdiff \upmu 
				+ \mupartialmodinhom,
		\label{E:LOWESTORDERTRANSPORTPARTIALRENORMALIZEDUPMU} \\
	\mupartialmodinhom
	& := 		\frac{1}{2} \upmu  G_{\Lunit \Lunit} \angdiff \Psi
					+ \upmu G_{\Lunit \Radunit} \angdiff \Psi.
			\label{E:LOWESTORDERTRANSPORTPARTIALRENORMALIZEDUPMUDISCREPANCY}
\end{align}	
\end{subequations}

\end{definition}

In the next lemma, we derive the transport equation verified by $\mupartialmod.$

\begin{lemma}[\textbf{The transport equation verified by the partially modified version of} $\angdiff \upmu$]
\label{L:UPMUFIRSTPARTIALRENORMALIZEDTRANSPORTEQUATION}
Let $\mupartialmod$ be the partially modified $S_{t,u}$ one-form defined in \eqref{E:LOWESTORDERTRANSPORTPARTIALRENORMALIZEDUPMU}. Then $\mupartialmod$
verifies the following transport equation:
\begin{align} \label{E:UPMUFIRSTPARTIALRENORMALIZEDTRANSPORTEQUATION}
	\angLie_{\Lunit} 
		\mupartialmod 
	& = \frac{1}{2} 
			G_{\Lunit \Lunit}
			\angdiff \Rad \Psi
			+ \mathfrak{J}
\end{align}
where the $S_{t,u}$ one-form $\mathfrak{J}$ has the following schematic form:
\begin{align} \label{E:UPMUFIRSTPARTIALRENORMALIZEDTRANSPORTINHOMOGENEOUSTERM}
	\mathfrak{J}
	& = \myarray
				[\upmu \Lunit \Psi]
				{\Rad \Psi}
				\angdiff G_{(Frame)}
			+ \upmu
				(\Lunit G_{(Frame)})
				\angdiff \Psi
			+ G_{(Frame)}^2
				\myarray
					[\upmu \Lunit \Psi]
					{\Rad \Psi}
				\angdiff \Psi
			+ G_{(Frame)}
				(\Lunit \Psi)
				\angdiff \upmu,
\end{align}
and no $S_{t,u}$ tensors such as $\angGarg{\Lunit}$
are present in the $G_{(Frame)}$ term on the right-hand side of \eqref{E:UPMUFIRSTPARTIALRENORMALIZEDTRANSPORTINHOMOGENEOUSTERM}
(that is, only scalar-valued functions such as $G_{\Lunit \Radunit}$ are present).

\end{lemma}

\begin{proof}
	We apply $\angdiff$ to equation \eqref{E:UPMUFIRSTTRANSPORT} and use Lemma~\ref{L:LANDRADCOMMUTEWITHANGDIFF}
	to rewrite the left-hand side as $\angdiff \Lunit \upmu = \angLie_{\Lunit} \angdiff \upmu.$
	We now address the terms that arise when $\angdiff$ falls on the term
	$\frac{1}{2} G_{\Lunit \Lunit}
	\Rad \Psi$
	from the right-hand side of \eqref{E:UPMUFIRSTTRANSPORT}.
	We explicitly place the term 
	$\frac{1}{2} G_{\Lunit \Lunit}
	\angdiff \Rad \Psi$ on the right-hand side of \eqref{E:UPMUFIRSTPARTIALRENORMALIZEDTRANSPORTEQUATION},
	while the term 
	$\frac{1}{2} (\angdiff G_{\Lunit \Lunit}) \Rad \Psi$
	is part of the inhomogeneous term \eqref{E:UPMUFIRSTPARTIALRENORMALIZEDTRANSPORTINHOMOGENEOUSTERM}.
	
	We now address the terms that arise when $\angdiff$ falls on the terms
	$- \frac{1}{2} \upmu G_{\Lunit \Lunit} \Lunit \Psi$
	and
	$- \upmu G_{\Lunit \Radunit} \Lunit \Psi$
	from the right-hand side of \eqref{E:UPMUFIRSTTRANSPORT}.
	When $\angdiff$ falls on $\upmu$ or $G_{(Frame)},$
	we consider these terms to be part of the inhomogeneous term \eqref{E:UPMUFIRSTPARTIALRENORMALIZEDTRANSPORTINHOMOGENEOUSTERM}.
	When $\angdiff$ falls on $\Lunit \Psi,$ we 
	use Lemma~\ref{L:LANDRADCOMMUTEWITHANGDIFF} to commute the $\Lunit$ derivative all the way out,
	and then move the terms
	$- \frac{1}{2} \angLie_{\Lunit}(\upmu G_{\Lunit \Lunit} \angdiff \Psi)$
	and $- \angLie_{\Lunit} (\upmu G_{\Lunit \Radunit} \angdiff \Psi)$
	to the left-hand side; the products underneath the $\angLie_{\Lunit}$ differentiation are 
	part of the quantity $\mupartialmod.$
	The commutator terms are included in 
	the inhomogeneous term \eqref{E:UPMUFIRSTPARTIALRENORMALIZEDTRANSPORTINHOMOGENEOUSTERM},
	and we use the equation \eqref{E:UPMUFIRSTTRANSPORT} to replace all instances of 
	$\Lunit \upmu$ with the right-hand side of \eqref{E:UPMUFIRSTTRANSPORT}.
	Lemma~\ref{L:UPMUFIRSTPARTIALRENORMALIZEDTRANSPORTEQUATION} thus follows.
\end{proof}

We now define higher-order versions of $\mupartialmod$ and $\mupartialmodinhom.$

\begin{definition}[\textbf{Partially modified version of} $\angdiff \mathscr{Z}^N \upmu$]
\label{D:TRANSPORTREPARTIALNORMALIZEDUPMU}
We define the partially modified $S_{t,u}$ one-form $\mupartialmodarg{\mathscr{Z}^N}$ as follows:
\begin{subequations}
\begin{align}
	\mupartialmodarg{\mathscr{Z}^N}
	& := \angdiff \mathscr{Z}^N \upmu 
				+ \mupartialmodinhomarg{\mathscr{Z}^N},
		\label{E:TRANSPORTPARTIALRENORMALIZEDUPMU} \\
	\mupartialmodinhomarg{\mathscr{Z}^N}
	& := 		\frac{1}{2} \upmu G_{\Lunit \Lunit} \angdiff \mathscr{Z}^N \Psi
					+ \upmu G_{\Lunit \Radunit} \angdiff \mathscr{Z}^N \Psi.
			\label{E:TRANSPORTPARTIALRENORMALIZEDUPMUDISCREPANCY}
\end{align}	
\end{subequations}
\end{definition}

In the next lemma, we derive the transport equation verified by 
$\mupartialmodarg{\mathscr{S}^{N-1}}.$

\begin{lemma}[\textbf{The transport equation verified by the partially modified version of} $\angdiff \mathscr{S}^{N-1} \upmu$]
\label{L:COMMUTEDUPMUFIRSTPARTIALRENORMALIZEDTRANSPORTEQUATION}
Let $\mathscr{S}^{N-1}$ be an $(N-1)^{st}$ order pure spatial commutation vectorfield operator
and let $\mupartialmodarg{\mathscr{S}^{N-1}}$ be the partially modified 
$S_{t,u}$ one-form defined in \eqref{E:TRANSPORTPARTIALRENORMALIZEDUPMU}.
Then $\mupartialmodarg{\mathscr{S}^{N-1}}$ verifies the following transport equation:
\begin{align} \label{E:COMMUTEDUPMUFIRSTPARTIALRENORMALIZEDTRANSPORTEQUATION}
	\angLie_{\Lunit} 
		\mupartialmodarg{\mathscr{S}^{N-1}}
	& = \frac{1}{2} 
			G_{\Lunit \Lunit}
			\angdiff \mathscr{S}^{N-1} \Rad \Psi
			+ \inhomleftexparg{\mathscr{S}^{N-1}}{\mathfrak{J}},
\end{align}
where the $S_{t,u}$ one-form $\inhomleftexparg{\mathscr{S}^{N-1}}{\mathfrak{J}}$ is given by
\begin{align} \label{E:COMMUTEDUPMUFIRSTPARTIALRENORMALIZEDTRANSPORTINHOMOGENEOUSTERM}
\mupartialmodsourcearg{\mathscr{S}^{N-1}}
	& = \angLie_{\mathscr{S}}^{N-1} \mathfrak{J}
		+ \frac{1}{2} 
			[\angLie_{\mathscr{S}}^{N-1}, G_{\Lunit \Lunit}]
			\angdiff \Rad \Psi
		+ [\angLie_{\Lunit}, \angLie_{\mathscr{S}}^{N-1}] \angdiff \upmu
		+ [\angLie_{\Lunit}, \angLie_{\mathscr{S}}^{N-1}] \mupartialmodinhom
			\\
	& \ \ 
			+ \angLie_{\Lunit}
			\left\lbrace
				\mupartialmodinhomarg{\mathscr{S}^{N-1}}
				- \angLie_{\mathscr{S}}^{N-1}  \mupartialmodinhom
			\right\rbrace,
			\notag
\end{align}
the $S_{t,u}$ one-form $\mupartialmodinhom$ is defined in \eqref{E:LOWESTORDERTRANSPORTPARTIALRENORMALIZEDUPMUDISCREPANCY},
the $S_{t,u}$ one-form $\mathfrak{J}$ is defined in \eqref{E:UPMUFIRSTPARTIALRENORMALIZEDTRANSPORTINHOMOGENEOUSTERM},
and the $S_{t,u}$ one-form $\mupartialmodinhomarg{\mathscr{S}^{N-1}}$ is defined in \eqref{E:TRANSPORTPARTIALRENORMALIZEDUPMUDISCREPANCY}.
\end{lemma}

\begin{proof}
	We apply $\angLie_{\mathscr{S}}^{N-1}$ to each side of 
	equation \eqref{E:UPMUFIRSTPARTIALRENORMALIZEDTRANSPORTEQUATION} 
	and perform several straightforward operator commutations.
\end{proof}


\chapter{Small Data, \texorpdfstring{$C^0$}{Sup-Norm} Bootstrap Assumptions, and First Pointwise Estimates} 
\label{C:C0BOUNDBOOTSTRAP}
\setcounter{equation}{0}
\thispagestyle{fancy}
Starting in Chapter~\ref{C:C0BOUNDBOOTSTRAP}, we assume that 
$\Psi$ verifies
$\square_{g(\Psi)} \Psi = 0$
on a spacetime region $\mathcal{M}_{\Tboot,U_0}.$
We first define the size of the data.
We then state our fundamental positivity bootstrap assumption for $\upmu,$
our $C^0$ bootstrap assumptions for $\Psi$
and its lower-order derivatives, and our auxiliary $C^0$ bootstrap assumptions for 
$\upmu,$ $\Lunit_{(Small)}^i,$ and $\upchi^{(Small)}$
and their lower-order derivatives
on $\mathcal{M}_{\Tboot,U_0}.$
We then use these bootstrap assumptions to derive 
$C^0$ and pointwise estimates for various quantities and their derivatives
with respect to the commutation vectorfields $Z \in \mathscr{Z}.$
The estimates that we derive in this section are
tedious but not too difficult;
we derive related, but more difficult, estimates 
in Chapters \ref{C:SHARPESTIMATESFORUPMU} and 
\ref{C:TOPORDEREIKONALFUNCTIONCOMMUTATIONPOINTWISE}-\ref{C:POINTWISEESTIMATESDIFFICULTERRORINTEGRANDS}.

\begin{remark}[\textbf{Suppression of the independent variables}]
Throughout the remainder of the monograph, we state many of our pointwise estimates in the form
\begin{align}
	|f_1| \lesssim h(t,u)|f_2|
\end{align}
for some function $h.$
To avoid cluttering the notation,
we use the convention that unless we indicate otherwise,
in such inequalities,
both $f_1$ and $f_2$
are evaluated at the point with geometric coordinates
$(t,u,\vartheta).$
\end{remark}

\section{\texorpdfstring{$\Psi$}{Psi} solves the wave equation}
\label{S:PSISOLVES}
We recall that 
until Sect.~\ref{S:SHOCKFORMINGDATA}, 
$U_0$ denotes a fixed parameter that verifies
$0 < U_0 < 1.$ 
We are primarily interested in the nontrivial portion of the future development
of the portion of the data lying in the exterior of $S_{0,U_0}.$
That is, we are interested in a spacetime region of the form $\mathcal{M}_{\Tboot,U_0},$ 
where $\Tboot > 0$ (see Def.~\ref{D:HYPERSURFACESANDCONICALREGIONS}).
Until Appendix \ref{A:EQUIVALENTPROBLEM}, we assume that on  $\mathcal{M}_{\Tboot,U_0},$ 
$\Psi$ is a solution to the covariant wave equation \eqref{E:WAVEGEO}, that is, that
\[
	\square_{g(\Psi)} \Psi = 0.
\]
We also assume that the metric verifies 
$g_{\mu \nu}(\Psi) = m_{\mu \nu} + g_{\mu \nu}^{(Small)}(\Psi)$
and $g_{\mu \nu}^{(Small)}(0) = 0,$
as stated in
\eqref{E:LITTLEGDECOMPOSED}
and
\eqref{E:METRICNONLINEARITYVANISHESWHENPSIIS0}.

\begin{remark}[\textbf{The constants can depend on} $U_0$]
Throughout our analysis,
some of the explicit constants ``$C$''
and the implicit constants tied to the notation ``$\lesssim$''
depend on $U_0.$
\end{remark}

\section{Small data}
\label{S:SMALLDATA}
In this section, we define the size of the data for the covariant wave equation $\square_{g(\Psi)} \Psi = 0.$

\begin{definition}[\textbf{Definition of the size of the data}] 
\label{D:SMALLDATA}
Let $(\mathring{\Psi} := \Psi|_{\Sigma_0}, \mathring{\Psi}_0 := \partial_t \Psi|_{\Sigma_0})$ 
be initial data for the covariant wave equation \eqref{E:WAVEGEO} that are compactly
supported in the Euclidean unit ball $\Sigma_0^1$ (see \eqref{E:EIKONALFUNCTIONINITIALCONDITIONS} and Def.~\ref{D:HYPERSURFACESANDCONICALREGIONS}).
We define the size of the data as follows:
\begin{align} \label{E:SMALLDATA}
	\mathring{\upepsilon} 
	:= \| \mathring{\Psi} \|_{H_{\Euct}^{25}(\Sigma_0^1)}
		+ \| \mathring{\Psi}_0 \|_{H_{\Euct}^{24}(\Sigma_0^1)}.
\end{align}
In \eqref{E:SMALLDATA}, $H_{\Euct}^N$ is the standard Euclidean Sobolev space involving 
rectangular spatial derivatives along $\Sigma_0^1.$

\end{definition}

\begin{center}
	{\textbf{Preliminary small data assumption}}
\end{center}
We make the following assumption throughout our analysis:
\begin{align} \label{E:SMALLDATASSUMPTION}
	\mathring{\upepsilon} \leq \varepsilon,
\end{align}
where $\varepsilon$ is the number appearing in
bootstrap assumptions of
Sects. 
\ref{SS:PSIBOOTSTRAP}
and
\ref{SS:AUXILIARYBOOTSTRAP}.

\section{Fundamental positivity bootstrap assumption for \texorpdfstring{$\upmu$}{the inverse foliation density}}
\label{S:FUNDAMENTALPOSITIVITY}
We make the following bootstrap assumption
for the solution on $\mathcal{M}_{\Tboot,U_0}:$
\begin{align} \label{E:BAUMPUISPOSITIVE}  \tag{$\mathbf{BA} \upmu > 0$}
	\upmu > 0 \ \mbox{on \ } \mathcal{M}_{\Tboot,U_0}. 
\end{align}
The bootstrap assumption \eqref{E:BAUMPUISPOSITIVE} 
implies that no shock is present in the region $\mathcal{M}_{\Tboot,U_0}.$


\section{\texorpdfstring{$C^0$}{Sup-norm} bootstrap assumptions} \label{S:C0BOUNDBOOTSTRAP}
Our quantitative bootstrap assumptions involve a 
small parameter $\varepsilon > 0.$
Throughout our analysis,
we adjust the smallness of $\varepsilon$
as necessary. 

\subsection{Fundamental \texorpdfstring{$C^0$}{sup-norm} bootstrap assumptions}
\label{SS:PSIBOOTSTRAP}
Our fundamental $C^0$ bootstraps assumption for $\Psi$ are:
\begin{align} \label{E:PSIFUNDAMENTALC0BOUNDBOOTSTRAP} \tag{$\mathbf{BA}\Psi$}
	\| \mathscr{Z}^{\leq 13} \Psi \|_{C^0(\Sigma_t^u)},
		\,
	\| \rgeo  \Lunit \mathscr{Z}^{\leq 12} \Psi \|_{C^0(\Sigma_t^u)},
		\,
	\| \rgeo \angdiff \mathscr{Z}^{\leq 12} \Psi \|_{C^0(\Sigma_t^u)}
		& \leq \frac{\varepsilon}{1 + t},
		&& (t,u) \in [0,\Tboot) \times [0,U_0]. 
\end{align}

\begin{remark}[\textbf{The role of Sobolev embedding}]
	Much later in the monograph, we use Sobolev embedding-type
	estimates to derive an improvement of 
	\eqref{E:PSIFUNDAMENTALC0BOUNDBOOTSTRAP}.
	The improvement is based on
	Cor.~\ref{C:C0BOUNDSOBOLEVINTERMSOFENERGIES}.
\end{remark}

\subsection{Auxiliary bootstrap assumptions}
\label{SS:AUXILIARYBOOTSTRAP}
Let $\upmu,$ $\Lunit_{(Small)}^i,$ $\Radunit_{(Small)}^i,$ and $\upchi^{(Small)}$
be the quantities from Definitions 
\ref{D:UPMU}
and \ref{D:RENORMALIZEDVARIABLES}.
Our auxiliary $C^0$ bootstraps for these quantities are:
\begin{subequations}
\begin{align}
	\| \mathscr{Z}^{\leq 12} (\upmu - 1) \|_{C^0(\Sigma_t^u)}
		& \leq \varepsilon^{1/2} \ln(\myexp + t),  
			&& (t,u) \in [0,\Tboot) \times [0,U_0],
		  \label{E:UPMUBOOT}  \tag{$\mathbf{AUX}\upmu$} \\
	\| \mathscr{Z}^{\leq 12} \Lunit_{(Small)}^i \|_{C^0(\Sigma_t^u)},
		\, \| \mathscr{Z}^{\leq 12} \Radunit_{(Small)}^i \|_{C^0(\Sigma_t^u)} 	
		& \leq \varepsilon^{1/2} \frac{\ln(\myexp + t)}{1 + t},  
			&& (t,u) \in [0,\Tboot) \times [0,U_0], 
	 	 \label{E:FRAMECOMPONENTSBOOT} \tag{$\mathbf{AUX}\Lunit_{(Small)}$} \\
	\| \angLie_{\mathscr{Z}}^{\leq 11} \upchi^{(Small)} \|_{C^0(\Sigma_t^u)}
		& \leq \varepsilon^{1/2} \frac{\ln(\myexp + t)}{(1 + t)^2},  
			&& t \in (t,u) \in [0,\Tboot) \times [0,U_0].
			\label{E:CHIJUNKBOOT} \tag{$\mathbf{AUX}\upchi$}
\end{align}
\end{subequations}

\begin{remark}[\textbf{Improving the auxiliary bootstrap assumptions}]
	Already within Chapter~\ref{C:C0BOUNDBOOTSTRAP}, we 
	derive estimates showing that the bootstrap assumptions
	\eqref{E:UPMUBOOT},
	\eqref{E:FRAMECOMPONENTSBOOT},
	and
	\eqref{E:CHIJUNKBOOT}
	in fact hold with $\sqrt{\varepsilon}$ replaced by $C \varepsilon;$
	see Cor.~\ref{C:AUXBOOTSTRAPIMPROVED}.
	It is in this sense that we consider these bootstrap assumptions to be 
	``auxiliary.''
\end{remark}	

\begin{remark}[\textbf{Conventions for repeated differentiation}]
In this section, we often use the conventions 
for repeated differentiation described in Sect.~\ref{S:REPEATEDDIFF}.	
\end{remark}

\section{Basic estimates for the geometric radial variable}
In this section, we derive some simple pointwise and commutator estimates involving
the geometric radial variable $\rgeo = 1 - u + t.$

\begin{lemma}[\textbf{Basic estimates for the geometric radial variable} $\rgeo = 1 - u + t$]
\label{L:BASICESTIMATESFORRGEO}
There exists a ($U_0-$dependent) constant $C > 1$ such that for 
$0 \leq u \leq U_0,$ we have the following comparison estimates:
\begin{align} \label{E:RGEOMINKOWSKITIMECOMPARISON}
	C^{-1} 
	(1 + t)
	& \leq \rgeo(t,u)
	\leq 1 + t.
\end{align}

Furthermore, the following identities hold, 
where $\mathscr{Z}$ is the set of commutation vectorfields 
from definition \eqref{E:DEFSETOFCOMMUTATORVECTORFIELDS}:
\begin{subequations}
		\begin{align} 
			\Lunit \rgeo 
			&= - \Rad \rgeo = 1, 
			&& \Rot \rgeo = 0,
				\label{E:EXACTRELATIONSZAPPLIEDTORGEO} \\
			\mathscr{Z}^N \rgeo 
			& \in \lbrace 0, - 1, \rgeo \rbrace &&
				\label{E:ZNAPPLIEDTORGEOISNOTTOOLARGE}
		\end{align}
\end{subequations}

Furthermore, if $\xi$ is any $S_{t,u}$ tensorfield, $M$ is an integer (not necessarily positive), and 
$1 \leq N \leq 24$ is an integer, then
we have the following commutator estimates:
\begin{subequations}		
		\begin{align}
			\left|	
				[\angLie_{\mathscr{Z}}^N, \rgeo]
				\xi
			\right|
			& \lesssim 
				\rgeo
				\left|
					\angLie_{\mathscr{Z}}^{\leq N-1} \xi
				\right|,
				\label{E:COMMUTATOROFZNANDRGEOAPPLIEDTOXIPOINTWISE}	\\
			\left|	
				\Lunit
				\left\lbrace
					[\angLie_{\mathscr{Z}}^N, \rgeo^M]
					\xi
				\right\rbrace
			\right|
			& \lesssim 
				\left|
					\angLie_{\Lunit} 
					\left\lbrace
						\rgeo^M 
						\angLie_{\mathscr{Z}}^{\leq N-1} \xi
					\right\rbrace
				\right|
				+ \frac{1}{\rgeo^2}
					\left|
						\rgeo^M \angLie_{\mathscr{Z}}^{\leq N-1} \xi
					\right|.
				\label{E:LAPPLIEDTOCOMMUTATOROFZNANDRGEOAPPLIEDTOXIPOINTWISE}
		\end{align}
	\end{subequations}
	
\end{lemma}

\begin{proof}
	The first two identities in \eqref{E:EXACTRELATIONSZAPPLIEDTORGEO} 
	follow from Lemmas \ref{L:NULLPAIRBASICPROPERTIES} and \ref{L:BASICPROPERTIESOFLUNITRADANDTIMENORMAL}.
	The identity $\Rot \rgeo = 0$ is trivial because $\Rot$ is $S_{t,u}-$tangent.
	The identity \eqref{E:ZNAPPLIEDTORGEOISNOTTOOLARGE} then follows from the definition \eqref{E:DEFSETOFCOMMUTATORVECTORFIELDS}
	of $\mathscr{Z}.$
	Inequalities \eqref{E:COMMUTATOROFZNANDRGEOAPPLIEDTOXIPOINTWISE} and \eqref{E:LAPPLIEDTOCOMMUTATOROFZNANDRGEOAPPLIEDTOXIPOINTWISE}
	follow easily from \eqref{E:EXACTRELATIONSZAPPLIEDTORGEO} and \eqref{E:ZNAPPLIEDTORGEOISNOTTOOLARGE}.
\end{proof}

\begin{remark}[\textbf{Silent use of} $\rgeo(t,u) \approx 1 + t$]
	Throughout the remainder of the monograph, we often use the estimate
	\eqref{E:RGEOMINKOWSKITIMECOMPARISON} without explicitly mentioning it.
\end{remark}

\section{Basic estimates for rectangular spatial coordinate functions}

In this section, we derive some simple pointwise estimates involving the 
rectangular spatial coordinate functions $x^i$ and the radial coordinate $r.$

\begin{lemma}[\textbf{Basic estimates for the rectangular spatial coordinate functions} $x^i$ \textbf{and} $r$]
Let $x^i,$ 
$(i=1,2,3),$
be the rectangular spatial coordinate functions and let $r = \sqrt{\sum_{a=1}^3(x^a)^2}$
be the standard Euclidean radial coordinate.
Under the small-data and bootstrap assumptions 
of Sects.~\ref{S:PSISOLVES}-\ref{S:C0BOUNDBOOTSTRAP},
if $\varepsilon$ is sufficiently small, 
then the following pointwise estimates hold on $\mathcal{M}_{\Tboot,U_0}:$
\begin{subequations}
\begin{align} \label{E:EASYCOORDINATEBOUND}
	|x^i| & \leq 
		\left(
			1+ C \varepsilon^{1/2}\frac{\ln(\myexp + t)}{1 + t}
		\right) \rgeo,
		 \\
	\left|
		\frac{r}{\rgeo} - 1
	\right|	
	& \leq
		C
		\varepsilon^{1/2} \frac{\ln(\myexp + t)}{1 + t}.
		\label{E:EASYEUCLIDEANRADIALVARIABLEBOUND}
\end{align}
\end{subequations}

Furthermore, the following estimates hold on $\mathcal{M}_{\Tboot,U_0}:$
\begin{align} \label{E:FIRSTESTIMATEFORANGULARDIFFERNTIALOFXI}
	|\angdiff x^i|
	& \leq 1 + C \varepsilon \frac{1}{1+t}.
\end{align}

\end{lemma}

\begin{proof}
To prove \eqref{E:EASYEUCLIDEANRADIALVARIABLEBOUND},
we first use \eqref{E:LITTLEGDECOMPOSED},
\eqref{E:HATRADLENGTH},
and
\eqref{E:RADUNITJUNK}
to deduce the identity
	\begin{align} \label{E:USEFULIDEFORRADIALCOORDINATES}
		\frac{r}{\rgeo}
		& =
			\sqrt{ 
			  1 
			  - 2 \delta_{ab} \Radunit^a \Radunit_{(Small)}^b
			  - g_{ab}^{(Small)} \Radunit^a \Radunit^b
				+ \delta_{ab} \Radunit_{(Small)}^a \Radunit_{(Small)}^b
			}.
	\end{align}
Using $|x^i| \leq r$
and the bootstrap assumptions 
\eqref{E:PSIFUNDAMENTALC0BOUNDBOOTSTRAP}
and
\eqref{E:FRAMECOMPONENTSBOOT},
we deduce that 
\begin{align} \label{E:ANNOYINGROVERRGEOESTIMATE}
\frac{r}{\rgeo}
& = \sqrt{
	1 
	+ 
	\mathcal{O}
		\left(
			\varepsilon^{1/2}
			\frac{\ln(\myexp + t)}{1 + t} 
			\frac{r}{\rgeo}
		\right)
	+ 
	\mathcal{O}
		\left(
			\varepsilon
			\frac{1}{1 + t} 
			\frac{r^2}{\rgeo^2}
		\right)
	+ 
	\mathcal{O}
	\left(
		\varepsilon 
		\frac{\ln^2(\myexp + t)}{(1 + t)^2}
	\right)
	}.
\end{align}
The desired estimate \eqref{E:EASYEUCLIDEANRADIALVARIABLEBOUND}
now follows easily from \eqref{E:ANNOYINGROVERRGEOESTIMATE}.


The estimate \eqref{E:EASYCOORDINATEBOUND} then follows from
\eqref{E:EASYEUCLIDEANRADIALVARIABLEBOUND}
and the fact that $|x^i| \leq r.$

To prove \eqref{E:FIRSTESTIMATEFORANGULARDIFFERNTIALOFXI},
we first use \eqref{E:PSIFUNDAMENTALC0BOUNDBOOTSTRAP} 
to deduce that $|\gt_{ij} - \delta_{ij}| = |g_{ij}^{(Small)}| \lesssim \varepsilon (1 + t)^{-1}.$
Hence, Taylor expanding the components of $\gt^{-1}$ in terms of the components of $\gt,$ we deduce that 
$|(\gt^{-1})^{ij} - \delta^{ij}| \lesssim \varepsilon (1 + t)^{-1}.$
Using this estimate and the fact that the magnitude of the 
angular differential of a function as measured by $\gsphere$
is no larger than the magnitude
of its spatial differential 
as measured by $\gt,$
we conclude the desired estimate as follows:
\begin{align}
	|\angdiff x^i|^2 
	& \leq |(\gt^{-1})^{ab} \partial_a x^i \partial_b x^i|
		= |(\gt^{-1})^{ii}|
			\leq 1 + C \varepsilon \frac{1}{1+t}.
\end{align}

\end{proof}

\section{Estimates for the rectangular components of the metrics and the \texorpdfstring{$S_{t,u}$}{spherical} projection}
In this section, we derive pointwise estimates for the rectangular components of the metrics
and the $S_{t,u}$ projection tensorfield.

\begin{lemma}[\textbf{Pointwise estimates for the rectangular components of the metrics and the} $S_{t,u}$ \textbf{projection tensorfield}]
Let $0 \leq N \leq 24$ be an integer.
Let $\gsphere,$ $\ginversesphere,$ and $\sphereproject$
be the $S_{t,u}$ tensors from Defs.~\ref{D:FIRSTFUND} and \ref{D:PROJECTIONS}.
Under the small-data and bootstrap assumptions 
of Sects.~\ref{S:PSISOLVES}-\ref{S:C0BOUNDBOOTSTRAP},
if $\varepsilon$ is sufficiently small, 
then the following pointwise estimates hold on the spacetime domain $\mathcal{M}_{\Tboot,U_0}$
for their non-zero rectangular components $(\mu,\nu = 0,1,2,3$ and $i,j=1,2,3):$
\begin{subequations}
\begin{align}
	\left| 
		\mathscr{Z}^N \overbrace{\left\lbrace g_{\mu \nu} - m_{\mu \nu} \right\rbrace}^{g_{\mu \nu}^{(Small)}} 
	\right|,
	\left| 
		\mathscr{Z}^N \left\lbrace (g^{-1})^{\mu \nu} - (m^{-1})^{\mu \nu} \right\rbrace 
	\right|
	& \lesssim |\mathscr{Z}^{\leq N} \Psi|,
		\label{E:GRECTCOMPONENTINTERMSOFOTHERVARIABLES}	\\
	\left| \mathscr{Z}^N 
			\left\lbrace
				(\gtinverse)^{ij} - \delta^{ij} 
			\right\rbrace 
	\right|
	& \lesssim |\mathscr{Z}^{\leq N} \Psi|,
		\label{E:GTINVERSERECTCOMPONENTSINTERMSOFOTHERVARIABLES} \\
 	\left| \mathscr{Z}^N 
					\left\lbrace
						\gsphere_{ij} 
						- \left(\delta_{ij} - \frac{x^i x^j}{\rgeo^2} \right) 
					\right\rbrace
	\right|
	& \lesssim |\mathscr{Z}^{\leq N} \Psi| 
		+ \sum_{a=1}^3 \left| \mathscr{Z}^{\leq N} \Lunit_{(Small)}^a \right|,
		 \\
	\left| 
			\mathscr{Z}^N
			\left\lbrace
				\gsphere_{ij} 
				- \left(
						\delta_{ij} - \Radunit_i \Radunit_j 
					\right)
			\right\rbrace
	\right|
	& \lesssim |\mathscr{Z}^{\leq N} \Psi| 
		+ \sum_{a=1}^3 \left| \mathscr{Z}^{\leq N} \Lunit_{(Small)}^a \right|,
		\label{E:GSPHERERECTCOMPONENTINTERMSOFRADUNIT}
			\\
	\left| \mathscr{Z}^N 
					\left\lbrace
						(\ginversesphere)^{ij} 
						- \left(\delta^{ij} - \frac{x^i x^j}{\rgeo^2} \right) 
					\right\rbrace
	\right|
	& \lesssim |\mathscr{Z}^{\leq N} \Psi| 
		+ \sum_{a=1}^3 \left| \mathscr{Z}^{\leq N} \Lunit_{(Small)}^a \right|,
		\label{E:GINVERSESPHERERECTCOMPONENTINTERMSOFOTHERVARIABLES} \\
		\left| \mathscr{Z}^N 
					\left\lbrace
						(\ginversesphere)^{ij} 
						- \left(\delta^{ij} - \Radunit^i \Radunit^j \right) 
					\right\rbrace
	\right|
	& \lesssim |\mathscr{Z}^{\leq N} \Psi| 
		+ \sum_{a=1}^3 \left| \mathscr{Z}^{\leq N} \Lunit_{(Small)}^a \right|,
		\label{E:GINVERSESPHERERECTCOMPONENTINTERMSOFRADUNIT} \\
	\left| 
			\mathscr{Z}^N
			\left\lbrace
				\sphereproject_j^{\ i} 
				- \left(
						\delta_j^{\ i} - \frac{x^i x^j}{\rgeo^2} 
					\right)
			\right\rbrace
	\right|
	& \lesssim |\mathscr{Z}^{\leq N} \Psi| 
		+ \sum_{a=1}^3 \left| \mathscr{Z}^{\leq N} \Lunit_{(Small)}^a \right|.
\end{align}
\end{subequations}

Furthermore, the following estimates hold:
\begin{subequations}
\begin{align}
	\left\| 
		\mathscr{Z}^{\leq 13} 	
		\overbrace{\left\lbrace g_{\mu \nu} - m_{\mu \nu} \right\rbrace}^{g_{\mu \nu}^{(Small)}} 
	\right\|_{C^0(\Sigma_t^u)},
	\left\| 
		\mathscr{Z}^{\leq 13} \left\lbrace (g^{-1})^{\mu \nu} - (m^{-1})^{\mu \nu} \right\rbrace 
	\right\|_{C^0(\Sigma_t^u)}
	& \lesssim \varepsilon \frac{1}{1 + t},
			\label{E:GRECTCOMPONENTC0BOUND} \\
	\left\| \mathscr{Z}^{\leq 13} 
			\left\lbrace
				(\gtinverse)^{ij} - \delta^{ij} 
			\right\rbrace 
	\right\|_{C^0(\Sigma_t^u)}
	& \lesssim \varepsilon \frac{1}{1 + t},
		\label{E:GTINVERSERECTCOMPONENTSC0BOUND} \\
	\left\| \mathscr{Z}^{\leq 12} 
					\left\lbrace
						\gsphere_{ij} 
						- \left(\delta_{ij} - \frac{x^i x^j}{\rgeo^2} \right) 
					\right\rbrace
	\right\|_{C^0(\Sigma_t^u)}
	& \lesssim \varepsilon^{1/2} \frac{\ln(\myexp + t)}{1 + t},
		\label{E:GSPHERERECTCOMPONENTC0BOUND} \\
	\left\| 
			\mathscr{Z}^{\leq 12}
			\left\lbrace
				\gsphere_{ij} 
				- \left(
					\delta_{ij} - \Radunit_i \Radunit_j 
				\right)
			\right\rbrace
	\right\|_{C^0(\Sigma_t^u)}
	& \lesssim \varepsilon^{1/2} \frac{\ln(\myexp + t)}{1 + t},
	\label{E:GSPHEREINTERMSOFRADUNITRECTCOMPONENTSC0BOUNDESTIMATE}
		\\
		\left\| \mathscr{Z}^{\leq 12} 
					\left\lbrace
						(\ginversesphere)^{ij} 
						- \left(\delta^{ij} - \frac{x^i x^j}{\rgeo^2} \right) 
					\right\rbrace
	\right\|_{C^0(\Sigma_t^u)}
	& \lesssim \varepsilon^{1/2} \frac{\ln(\myexp + t)}{1 + t},
		\label{E:GINVERSESPHERERECTCOMPONENTSC0BOUNDESTIMATE} \\
	\left\| \mathscr{Z}^{\leq 12} 
					\left\lbrace
						(\ginversesphere)^{ij} 
						- \left(\delta^{ij} - \Radunit^i \Radunit^j \right) 
					\right\rbrace
	\right\|_{C^0(\Sigma_t^u)}
	& \lesssim \varepsilon^{1/2} \frac{\ln(\myexp + t)}{1 + t},
		\label{E:GINVERSESPHEREINTERMSOFRADUNITRECTCOMPONENTSC0BOUNDESTIMATE} \\
	\left\| 
			\mathscr{Z}^{\leq 12}
			\left\lbrace
				\sphereproject_j^{\ i} 
				- \left(
					\delta_j^{\ i} - \frac{x^i x^j}{\rgeo^2} 
				\right)
			\right\rbrace
	\right\|_{C^0(\Sigma_t^u)}
	& \lesssim \varepsilon^{1/2} \frac{\ln(\myexp + t)}{1 + t}.
\end{align}
\end{subequations}

Finally, the following pointwise estimates hold on $\mathcal{M}_{\Tboot,U_0}:$
\begin{align} \label{E:GSPHERERECTCOMPONENTSPOINTWISEBOUND}
	\left|\gsphere_{ij} \right|,
		\, \left|(\ginversesphere)^{ij} \right|,
		\, 	\left|\sphereproject_j^{\ i} \right| 
	& \leq 1 + C \varepsilon^{1/2}.
\end{align}

In the above estimates,
$\delta_{ij},$
$\delta_j^{\ i},$
and $\delta^{ij}$
are all standard Kronecker deltas.
\end{lemma}

\begin{proof}
	Recall that $g_{\mu \nu} = \delta_{\mu \nu} + g_{\mu \nu}^{(Small)}(\Psi),$
	where $g_{\mu \nu}^{(Small)}(\cdot)$ is smooth and verifies $g_{\mu \nu}^{(Small)}(0) = 0.$
	The inequalities in \eqref{E:GRECTCOMPONENTINTERMSOFOTHERVARIABLES} thus follow easily.
	\eqref{E:GRECTCOMPONENTC0BOUND} then follows from 
	\eqref{E:GRECTCOMPONENTINTERMSOFOTHERVARIABLES} and
	\eqref{E:PSIFUNDAMENTALC0BOUNDBOOTSTRAP}.
	The proofs of \eqref{E:GTINVERSERECTCOMPONENTSINTERMSOFOTHERVARIABLES} 
	and \eqref{E:GTINVERSERECTCOMPONENTSC0BOUND} are similar.
	
	To prove \eqref{E:GINVERSESPHERERECTCOMPONENTINTERMSOFOTHERVARIABLES}
	we first use \eqref{E:GINVERSEFRAMEWITHRECTCOORDINATESFORGSPHEREINVERSE}
	to expand 
	$(\ginversesphere)^{ij} 
		= (g^{-1})^{ij} 
		+ \Lunit^i \Lunit^j
		+ \Lunit^i \Radunit^j + \Radunit^i \Lunit^j.$
	We next use Def.~\ref{D:RENORMALIZEDVARIABLES} 
	and \eqref{E:RADUNITJUNKLIKELMINUSUNITJUNK}
	to expand $\Lunit^i = \rgeo^{-1} x^i + \Lunit_{(Small)}^i$
	and $\Radunit^i = -\rgeo^{-1} x^i - \Lunit_{(Small)}^i - (g^{-1})^{0i}.$
	Inequality \eqref{E:GINVERSESPHERERECTCOMPONENTINTERMSOFOTHERVARIABLES} now easily 
	follows from these expansions 
	and the estimates
	\eqref{E:GRECTCOMPONENTINTERMSOFOTHERVARIABLES} and \eqref{E:GTINVERSERECTCOMPONENTSINTERMSOFOTHERVARIABLES}.
	\eqref{E:GINVERSESPHERERECTCOMPONENTSC0BOUNDESTIMATE} then follows from
	\eqref{E:GINVERSESPHERERECTCOMPONENTINTERMSOFOTHERVARIABLES},
	\eqref{E:PSIFUNDAMENTALC0BOUNDBOOTSTRAP},
	and \eqref{E:FRAMECOMPONENTSBOOT}.
	
	The proofs of the remaining inequalities are very similar, and we omit the details.
\end{proof}

\begin{corollary}[\textbf{Comparison between the norms of} $S_{t,u}$ \textbf{tensors and the size of their rectangular components}]
\label{C:NORMOFTHESTUTENSORISCOMPARABLETOTHESUMOFITSCOMPONENTS}
Let be $\xi$ any $S_{t,u}$ tensor.
Then under the small-data and bootstrap assumptions 
of Sects.~\ref{S:PSISOLVES}-\ref{S:C0BOUNDBOOTSTRAP},
the following comparison estimates hold on $\mathcal{M}_{\Tboot,U_0}:$
\begin{align} \label{E:NORMOFTHESTUTENSORISCOMPARABLETOTHESUMOFITSCOMPONENTS}
	|\xi|
	\approx
	\sum |\xi_{Rectangular}|,
\end{align}
where the left-hand side is as in Def.~\ref{D:POINTWISENORMS} 
and the sum on the right-hand side is taken over the components of $\xi$ relative 
to the spatial rectangular coordinate frame.

\end{corollary}
\begin{proof}
 	We give the proof in the case that $\xi$ is an $S_{t,u}$ one-form. The proof
 	will generalize in an obvious fashion to the case of general $S_{t,u}$ tenors.
 	We first recall that $|\xi|^2 = (\ginversesphere)^{ab} \xi_a \xi_b.$
	Hence, by the estimate \eqref{E:GSPHERERECTCOMPONENTSPOINTWISEBOUND}
	for the rectangular components of $\ginversesphere,$ we have
	$|\xi|^2 \lesssim \sum_{a=1}^3 |\xi_a|^2,$
	which easily implies that the left-hand side of \eqref{E:NORMOFTHESTUTENSORISCOMPARABLETOTHESUMOFITSCOMPONENTS} 
	is $\lesssim$ the right-hand side as desired.
	
	To prove the reverse inequality, we first note that
	$\sum_{a=1}^3 |\xi_a|^2 = \delta^{ab} \xi_a \xi_b.$
	Using \eqref{E:GINVERSESPHEREINTERMSOFRADUNITRECTCOMPONENTSC0BOUNDESTIMATE},
	we see that we can replace $\delta^{ab}$ with $(\ginversesphere)^{ab}$
	plus an error term that is in magnitude $\lesssim \varepsilon^{1/2} \ln(\myexp + t)(1 + t)^{-1}$
	plus tensorial products of the vector $\Radunit$ that are completely annihilated because they are paired with
	copies of $\xi.$ We therefore deduce that
	$\left|\delta^{ab} \xi_a \xi_b - |\xi|^2 \right| \lesssim \delta^{ab} \xi_a \xi_b \varepsilon^{1/2} \ln(\myexp + t)(1 + t)^{-1},$
	from which we easily conclude that
	the right-hand side of \eqref{E:NORMOFTHESTUTENSORISCOMPARABLETOTHESUMOFITSCOMPONENTS} 
	is $\lesssim$ the left-hand side as desired.
\end{proof}

\section{The behavior of quantities along  \texorpdfstring{$\Sigma_0$}{the initial data hypersurface}}
\label{S:INITIALBEHAVIOROFQUANTITIES}
In this section, we provide quantitative estimates
showing that when $\mathring{\upepsilon}$ is small
(see \eqref{E:SMALLDATA}), many other quantities such as
$\upmu - 1,$
$\Lunit_{(Small)}^i,$
etc.
are also small along $\Sigma_0^1.$

\begin{lemma}[\textbf{All rectangular derivatives of} $\Psi$ \textbf{are small along} $\Sigma_0^1$]
	\label{L:BEHAVIOROFRECTANGLARDERIVATIVESALONGSIGMA0}
	Assume that $\square_{g(\Psi)} \Psi = 0,$
	and let $\mathring{\upepsilon}$ be the size of the data as defined in 
	Def.~\ref{D:SMALLDATA}.
	Then if $\mathring{\upepsilon}$ is sufficiently small, we have
	\begin{align} \label{E:ALLRECTANGULARDERIVATIVESSMALLALONGSIGMA0}
		\sum_{M=0}^{25}
		\left\|
			\partial_t^M \Psi
		\right\|_{H_{\Euct}^{25 - M}(\Sigma_0^1)}
		& \lesssim \mathring{\upepsilon}.
	\end{align}
	In \eqref{E:ALLRECTANGULARDERIVATIVESSMALLALONGSIGMA0}, 
	$H_{\Euct}^N$ is the standard Euclidean Sobolev space involving 
	order $\leq N$ rectangular spatial derivatives along $\Sigma_0^1.$
\end{lemma}
\begin{proof}
	The cases $M=0,1$ in \eqref{E:ALLRECTANGULARDERIVATIVESSMALLALONGSIGMA0} follow directly from the definition
	of the size of the data. To deduce \eqref{E:ALLRECTANGULARDERIVATIVESSMALLALONGSIGMA0} in the case $M=2,$
	we use the wave equation to solve for $\partial_t^2 \Psi$ as follows:
	\begin{align} \label{E:PARTIALTSQUAREDPSIISOLTATED}
		\partial_t^2 \Psi 
		& = (g^{-1})^{ab} \partial_a \partial_b \Psi
			+ 2 (g^{-1})^{0a} \partial_a \partial_t \Psi
			+ \frac{1}{2}
				(g^{-1})^{\alpha \beta} (g^{-1})^{\kappa \lambda} 
				\left\lbrace
					2 G_{\alpha \kappa} \partial_{\beta} \Psi
				-  G_{\alpha \beta} \partial_{\kappa} \Psi
		\right\rbrace.
	\end{align}
	The desired estimate \eqref{E:ALLRECTANGULARDERIVATIVESSMALLALONGSIGMA0}  
	then follows from setting $t = 0$ in \eqref{E:PARTIALTSQUAREDPSIISOLTATED} and applying the standard Sobolev calculus.
	To deduce \eqref{E:ALLRECTANGULARDERIVATIVESSMALLALONGSIGMA0} in the cases $M \geq 3,$
	we repeatedly differentiate \eqref{E:PARTIALTSQUAREDPSIISOLTATED} with respect to $t$
	and use the equations to express all higher-order time derivatives in terms of spatial derivatives
	of $\Psi$ and $\partial_t \Psi.$ The desired estimate \eqref{E:ALLRECTANGULARDERIVATIVESSMALLALONGSIGMA0}  
	then follows from the standard Sobolev calculus at $t=0$ as in the case $M=2.$
\end{proof}

\begin{lemma}[$\upmu - 1,$ $\Lunit_{(Small)}^i,$ \textbf{and} $\Xi^i$  \textbf{are small along} $\Sigma_0^1$]
\label{L:INITIALEXPRESSIONSFORUPMUANDLJUNKI}
Let 
\begin{align} \label{E:LITTLEZDEF}
	z := \sqrt{(\gtinverse)^{ab} \frac{x^a x^b}{r^2}}
	= \sqrt{1 + \left\lbrace(\gtinverse)^{ab} - \delta^{ab} \right\rbrace \frac{x^a x^b}{r^2}},
\end{align}
where $\gt$ is as defined in Def.~\ref{D:FIRSTFUND}.
Then along the initial data hypersurface region $\Sigma_0^1,$ 
the following identities hold
$(i=1,2,3):$
\begin{subequations}
\begin{align}
	\upmu^{-2} 
	& = 1 
			+
			\left\lbrace 
				(\gtinverse)^{ab} - \delta^{ab}
			\right\rbrace
			\frac{x^a x^b}{r^2},
		\label{E:INITIALMURELATION} \\
	\Lunit_{(Small)}^i
	& = \frac{1}{z}
			\left\lbrace 
					(\gtinverse)^{ia} - \delta^{ia} 
			\right\rbrace
			\frac{x^a}{r}
			+ \left\lbrace
					\frac{1}{z}
					- 1
				\right\rbrace
				\frac{x^i}{r}
			- (g^{-1})^{0i},
			\label{E:INITIALLIJUNKRELATION}
		\\
	\Xi^i
		& = \left\lbrace 
					\frac{1}{z^2} - 1
				\right\rbrace
				(\gtinverse)^{ia} \frac{x^a}{r}
			+ \left\lbrace 
					(\gtinverse)^{ia} - \delta^{ia} 
			\right\rbrace
			\frac{x^a}{r}. 
		\label{E:INITIALXIRELATION}
\end{align}
\end{subequations}
Above, $\Xi$ is the $S_{t,u}-$tangent vectorfield defined by \eqref{E:RADINTERMSOFGEOMETRICCOORDINATEPARTIALDERIVATIVES}.

In addition, under the hypotheses of Lemma~\ref{L:BEHAVIOROFRECTANGLARDERIVATIVESALONGSIGMA0},
the following estimates hold:
\begin{subequations}
\begin{align}
	\sum_{M=0}^{25}
		\| \partial_t^M (\upmu - 1)  \|_{H_{\Euct}^{25 - M}(\Sigma_0^1)} 
	& \lesssim \mathring{\upepsilon},
		\label{E:UPMUMINUSONESMALLDATA} \\
	\sum_{M=0}^{25}
		\| \partial_t^M \Lunit_{(Small)}^i \|_{H_{\Euct}^{25 - M}(\Sigma_0^1)} 
	& \lesssim \mathring{\upepsilon},
		\label{E:LIJUNKSMALLDATA}
		\\
	\sum_{M=0}^{25}
		\| \partial_t^M \Xi^i \|_{H_{\Euct}^{25 - M}(\Sigma_0^1)}
	& \lesssim \mathring{\upepsilon}.
		\label{E:XISMALLDATA}
\end{align}
\end{subequations}
Above, $H_{\Euct}^N$ is the standard Euclidean Sobolev space involving 
order $\leq N$ rectangular spatial derivatives along $\Sigma_0^1.$
\end{lemma}

\begin{proof}
	In this proof, the identities that we state hold along $\Sigma_0^1.$
	In particular,  
	we have $u = 1 - r$ 
	(see \eqref{E:EIKONALFUNCTIONINITIALCONDITIONS}),
	where $r$ is the standard Euclidean radial variables on $\mathbb{R}^3,$
	and hence $\partial_i u = - \frac{x^i}{r}.$
	From this identity and \eqref{E:ALTERNATEEIKONAL},
	we deduce \eqref{E:INITIALMURELATION}.
	Next, to prove \eqref{E:INITIALLIJUNKRELATION}, we first use definition \eqref{E:RADUNITJUNK},
	the identity \eqref{E:RADUNITJUNKLIKELMINUSUNITJUNK},
	the above calculations,
	and the fact that $\rgeo = r$ along $\Sigma_0^1$ to deduce that
	$\Lunit_{(Small)}^i = - \Radunit_{(Small)}^i - (g^{-1})^{0i} 
		= z^{-1} (\gtinverse)^{ia} \frac{x^a}{r} - \frac{x^i}{r} - (g^{-1})^{0i}.$
	The desired identity \eqref{E:INITIALLIJUNKRELATION} follows easily from the previous equation.
	Next, to prove \eqref{E:INITIALXIRELATION}, 
	we note that the geometric coordinates are constructed in such a way that along $\Sigma_0^1,$
	we have $\frac{\partial}{\partial u} = - \partial_r,$ where $\partial_r$ is the standard Euclidean radial
	derivative. Also using 
	\eqref{E:RADINTERMSOFGEOMETRICCOORDINATEPARTIALDERIVATIVES},
	we deduce that along $\Sigma_0^1,$ we have
	\[
	\Xi^i 
	= \Xi x^i 
	= -\partial_r x^i - \Rad^i 
	= - \frac{x^i}{r} - \upmu \Radunit^i 
	= - \frac{x^i}{r} + z^{-2} (\gtinverse)^{ia} \frac{x^a}{r}.
	\]
	The desired identity \eqref{E:INITIALXIRELATION} thus follows.
	
The estimates \eqref{E:UPMUMINUSONESMALLDATA}-\eqref{E:XISMALLDATA}
in the case $M=0$
then follow from the identities \eqref{E:INITIALMURELATION}-\eqref{E:INITIALXIRELATION},
the fact that 
$(\gtinverse)^{ij} - \delta^{ij}$ 
and $(g^{-1})^{0i}$
are smooth scalar-valued functions of $\Psi$
that vanish at $\Psi = 0,$
the estimate \eqref{E:ALLRECTANGULARDERIVATIVESSMALLALONGSIGMA0}, and the standard Sobolev calculus.
To obtain the desired estimates in the cases $M \geq 1,$ we 
inductively use the evolution equations	
\eqref{E:XIEVOLUTION}
(this equation needs to be rewritten relative to rectangular coordinates),
\eqref{E:UPMUFIRSTTRANSPORT},
and
\eqref{E:LUNITJUNKLPROP}
and equations \eqref{E:INITIALMURELATION}-\eqref{E:INITIALXIRELATION}
to solve for all time derivatives (including higher-order ones)
of $\upmu,$ $\Lunit_{(Small)}^i,$ and $\Xi^i$
in terms of the rectangular coordinate derivatives of $\Psi$ and 
the rectangular coordinate derivatives of the eikonal function $u.$ Furthermore, along $\Sigma_0^1,$ 
$u$ is a smooth function of 
the spatial coordinates, while by using the eikonal equation
\eqref{E:OUTGOINGEIKONAL},
we can solve for the time derivatives of $u$ in terms of
the rectangular spatial derivatives of
$\Psi$ and the rectangular spatial derivatives of $u.$
The desired estimates thus follow from the standard Sobolev calculus.
\end{proof}

\begin{lemma}[\textbf{Various Sobolev norms are initially small}]
\label{L:SMALLINITIALSOBOLEVNORMS}
	Under the hypotheses of Lemma~\ref{L:BEHAVIOROFRECTANGLARDERIVATIVESALONGSIGMA0},
	we have the following estimates $(i=1,2,3):$
	\begin{subequations}
	\begin{align} \label{E:SMALLINITIALSOBOLEVNORMS}
		\left\|
			\mathscr{Z}^{\leq 25} \Psi
		\right\|_{L^2(\Sigma_0^1)},
		\,
		\left\|
			\mathscr{Z}^{\leq 25} (\upmu - 1)
		\right\|_{L^2(\Sigma_0^1)},
		\,
		\left\|
			\rgeo \mathscr{Z}^{\leq 25} \Lunit_{(Small)}^i
		\right\|_{L^2(\Sigma_0^1)},
			\,
		\left\|
			\mathscr{Z}^{\leq 25} \Xi^i
		\right\|_{L^2(\Sigma_0^1)}
		& \lesssim \mathring{\upepsilon},	
					\\
		\left\|
			\mathscr{Z}^{\leq 23} \Psi
		\right\|_{C^0(\Sigma_0^1)},
		\,
		\left\|
			\mathscr{Z}^{\leq 23} (\upmu - 1)
		\right\|_{C^0(\Sigma_0^1)},
		\,
		\left\|
			\rgeo \mathscr{Z}^{\leq 23} \Lunit_{(Small)}^i
		\right\|_{C^0(\Sigma_0^1)},
			\,
		\left\|
			\mathscr{Z}^{\leq 23} \Xi^i
		\right\|_{C^0(\Sigma_0^1)}
		& \lesssim \mathring{\upepsilon}.	
		\label{E:SMALLINITIALC0BOUNDNORMS}
	\end{align}
	\end{subequations}
	Above, $\left\| \cdot \right\|_{L^2(\Sigma_0^1)}$ is the Sobolev
	norm from \eqref{E:SIGMATL2NORMDEF}.
	
	Furthermore, we have the following estimates, where 
	$\Xi$ is the $S_{t,u}-$tangent vectorfield defined by \eqref{E:RADINTERMSOFGEOMETRICCOORDINATEPARTIALDERIVATIVES}
	and
	$\Xi_{\flat}$ denotes the
	$\gsphere-$dual of $\Xi:$
	\begin{subequations}
	\begin{align}
		\left\|
			\angLie_{\mathscr{Z}}^{\leq 25} \Xi
		\right\|_{L^2(\Sigma_0^1)},
			\,
		\left\|
			\angLie_{\mathscr{Z}}^{\leq 25} \Xi_{\flat}
		\right\|_{L^2(\Sigma_0^1)}
		& \lesssim \mathring{\upepsilon},
			\label{E:XITENSORIALSMALLSOBOLEVNORM} \\
		\left\|
			\angLie_{\mathscr{Z}}^{\leq 23} \Xi
		\right\|_{C^0(\Sigma_0^1)},
			\, 
		\left\|
			\angLie_{\mathscr{Z}}^{\leq 23} \Xi_{\flat}
		\right\|_{C^0(\Sigma_0^1)}
		& \lesssim \mathring{\upepsilon}.
		\label{E:XIC0SMALLSOBOLEVNORM}
	\end{align}
	\end{subequations}
	
\end{lemma}

\begin{proof}
	Since $d \tvol$ is equal to $\upmu^{-1}$ times the volume form induced by $g$ on $\Sigma_0,$
	it follows that $d \tvol = \upmu^{-1} \mbox{\upshape{det}} \gt dx^1 dx^2 dx^3,$
	where the determinant is taken relative to rectangular coordinates.
	From \eqref{E:UPMUMINUSONESMALLDATA} and the fact that $\gt_{ij} = \delta_{ij} + g_{ij}^{(Small)}(\Psi),$
	we deduce that along $\Sigma_0,$ $\upmu^{-1} \mbox{\upshape{det}} \gt dx^1 dx^2 dx^3$
	is near-Euclidean in the sense that
	\begin{align} \label{E:NEAREUCLIDEANVOLUMEFORM}
		1 - C \mathring{\upepsilon}
		\leq \upmu^{-1} \mbox{\upshape{det}} \gt 
		\leq 1 + C \mathring{\upepsilon}.
	\end{align}
	It follows that the norm
	$\left\| \cdot \right\|_{L^2(\Sigma_0^1)}$
	is comparable to the standard Euclidean Sobolev norm
	$\left\| \cdot \right\|_{L_{\Euct}^2(\Sigma_0^1)}.$
	Hence, to deduce \eqref{E:SMALLINITIALSOBOLEVNORMS}, it suffices to prove the same inequalities
	using the standard Euclidean Sobolev norm 
	$\left\| \cdot \right\|_{L_{\Euct}^2(\Sigma_0^1)}$
	in place of $\left\| \cdot \right\|_{L^2(\Sigma_0^1)}.$
	To this end, we completely expand
	$\mathscr{Z}^{\leq 25} (\upmu - 1)$
	in terms of rectangular coordinate derivatives
	and use the standard Sobolev calculus to deduce that
	\begin{align} \label{E:ZUPMUMINUSONEINITIALEXPANDEDRECTANGLUAR}
		\left\| 
			\mathscr{Z}^{\leq 25} (\upmu - 1)	
		\right\|_{L_{\Euct}^2(\Sigma_0^1)}
		& \lesssim
		\left(
			\sum_{M=0}^{25}
			\left\| 
				\partial_t^M
				(\upmu - 1)
			\right\|_{H_{\Euct}^{25-M}(\Sigma_0^1)}
		\right)
		\left(
		\sum_{M=0}^{24}
		\max_{Z \in \mathscr{Z}}
		\sum_{\alpha = 0}^3
		\left\| 
			\partial_t^M
			Z^{\alpha}
		\right\|_{H_{\Euct}^{24-M}(\Sigma_0^1)}
		\right)^{25}.
	\end{align}	
	Inequality \eqref{E:UPMUMINUSONESMALLDATA} implies that 
	$\left\| 
		\partial_t^M
	(\upmu - 1)
	\right\|_{H_{\Euct}^{25-M}(\Sigma_0^1)} \lesssim \mathring{\upepsilon}.$
	Furthermore, from  the identities 
	$\rgeo \Lunit^0 = \rgeo, 
	\rgeo \Lunit^i = x^i + \rgeo \Lunit_{(Small)}^i,$
	\eqref{E:RADUNITJUNKLIKELMINUSUNITJUNK},
	and
	\eqref{E:ROTATIONRECATNGULARCOMPONENT},
	it follows that all of the rectangular components 
	$Z^{\alpha}$ $(\alpha = 0,1,2,3)$ 
	are smooth functions of the rectangular coordinate functions,
	$\Psi,$ $u,$ $\upmu,$ and the rectangular components $\Lunit_{(Small)}^i.$
	It thus follows from the remarks made in the proof of Lemma
	\ref{L:INITIALEXPRESSIONSFORUPMUANDLJUNKI} about how to estimate
	the rectangular derivatives of $u,$
	the standard Sobolev calculus, 
	\eqref{E:ALLRECTANGULARDERIVATIVESSMALLALONGSIGMA0},
	and
	\eqref{E:UPMUMINUSONESMALLDATA}-\eqref{E:LIJUNKSMALLDATA}
	that $\left\| \partial_t^M Z^{\alpha} \right\|_{H_{\Euct}^{25-M}(\Sigma_0^1)} \lesssim 1.$
	We have thus proved \eqref{E:SMALLINITIALSOBOLEVNORMS} for $\upmu - 1.$ 
	The remaining
	three inequalities in \eqref{E:SMALLINITIALSOBOLEVNORMS} can be proved using a similar argument,
	and we omit the details.
	The estimates \eqref{E:SMALLINITIALC0BOUNDNORMS} follow in a similar fashion
	thanks to the standard Sobolev embedding estimate 
	$\left\| f \right\|_{C^0(\Sigma_0^1)} \lesssim \left\| f \right\|_{H_{\Euct}^2(\Sigma_0^1)}.$
	The estimates \eqref{E:XITENSORIALSMALLSOBOLEVNORM} and \eqref{E:XIC0SMALLSOBOLEVNORM}
	for $\Xi$ follow similarly. In particular, by Lemma~\ref{L:VECTORFIELDCOMMUTATORS}, the projected Lie derivatives
	of $\Xi$ agree with standard Lie derivatives. Hence, 
	$|\angLie_{\mathscr{Z}}^{N} \Xi|^2 = \gsphere_{ab}(\Lie_{\mathscr{Z}}^N \Xi^a)\Lie_{\mathscr{Z}}^N \Xi^b$
	and the computations can be carried out relative
	to the rectangular coordinates. To derive the desired estimates
	\eqref{E:XITENSORIALSMALLSOBOLEVNORM} and \eqref{E:XIC0SMALLSOBOLEVNORM}
	for $\Xi,$ we note the identity
	$|\angLie_{\mathscr{Z}}^N \Xi_{\flat}|^2 
	= (\ginversesphere)^{ab}(\angLie_{\mathscr{Z}}^N (\gsphere_{ac} \Xi^c)) \angLie_{\mathscr{Z}}^N (\gsphere_{bd} \Xi^d),$
	which allows us to carry out the computations relative to the rectangular coordinates.
\end{proof}

\section{Estimates for the derivatives of rectangular components of various vectorfields}
In this section, we derive pointwise estimates for the derivatives of 
the rectangular components of various vectorfields that we use in our analysis.

\begin{lemma}[\textbf{First pointwise estimates for the spatial rectangular components of various vectorfields
and the functions $\RotRadcomponent{l}$}]
	\label{L:DERIVATIVESOFRECTANGULARCOMMUTATORCOMPONENTS}
	Let $0 \leq N \leq 24$ be an integer,
	and let $\Lunit^i$ and $\Radunit^i$ $(i=1,2,3)$
	be the (scalar-valued) rectangular spatial components
	of $\Lunit$ and $\Radunit.$
	Under the small-data and bootstrap assumptions 
	of Sects.~\ref{S:PSISOLVES}-\ref{S:C0BOUNDBOOTSTRAP},
	if $\varepsilon$ is sufficiently small,
	then the following pointwise estimates hold on $\mathcal{M}_{\Tboot,U_0}:$
	\begin{subequations}
	\begin{align} 	\label{E:DERIVATIVESOFRECTANGULARFRAMECOMPONENTS}
		\left|
			\myarray
				[\mathscr{Z}^N \Lunit^i]
				{\mathscr{Z}^N \Radunit^i}
		\right|
		& \lesssim |\mathscr{Z}^{\leq N} \Psi|
			+ \frac{1}{1 + t}
				\left|
						\myarray[\mathscr{Z}^{\leq N-1} (\upmu - 1)]
						{\sum_{a=1}^3 \rgeo |\mathscr{Z}^{\leq N} \Lunit_{(Small)}^a|} 
				\right|
			+ 1,
		   \\
	\left\| 
		\myarray
			[\mathscr{Z}^{\leq 12} \Lunit^i]
			{\mathscr{Z}^{\leq 12} \Radunit^i}
	\right\|_{C^0(\Sigma_t^u)}
		& \lesssim 1.
		\label{E:POINTWISEBOUNDLOWERORDERDERIVATIVESOFRECTANGULARFRAMECOMPONENTS}
	\end{align}
	\end{subequations}
	
	Furthermore, for each $Z \in \mathscr{Z} = \lbrace \rgeo \Lunit, \Rad, \Rot_{(1)}, \Rot_{(2)}, \Rot_{(3)} \rbrace,$ 
	we have the following estimates for the rectangular spatial components $Z^i$ $(i=1,2,3)$
	on $\mathcal{M}_{\Tboot,U_0}:$
	\begin{subequations}
	\begin{align} \label{E:DERIVATIVESOFRECTANGULARCOMMUTATORVECTORFIELDRECTCOMPONENTS}
		\left| \mathscr{Z}^N Z^i \right|
		& \lesssim (1+t) |\mathscr{Z}^{\leq N} \Psi|
			+ \left|
						\myarray[\mathscr{Z}^{\leq N} (\upmu - 1)]
						{\sum_{a=1}^3 \rgeo |\mathscr{Z}^{\leq N} \Lunit_{(Small)}^a|} 
				\right|
			+ 1 + t,
				\\
		\left\| 
			\mathscr{Z}^{\leq 12} Z^i 
		\right\|_{C^0(\Sigma_t^u)}
		& \lesssim 1 + t.
		\label{E:LOWERORDERC0BOUNDDERIVATIVESOFRECTANGULARCOMMUTATORVECTORFIELDRECTCOMPONENTS}
	\end{align}
	\end{subequations}
	
	In addition, we have the following estimates for the scalar-valued functions $\RotRadcomponent{l},$
	$(l = 1,2,3),$ defined in \eqref{E:EUCLIDEANROTATIONRADCOMPONENT} on $\mathcal{M}_{\Tboot,U_0}:$
\begin{subequations}	
\begin{align} \label{E:FIRSTPOINTWISEBOUNDEUCLIDEANROTATIONRADCOMPONENT}
	|\mathscr{Z}^N \RotRadcomponent{l}|
	& \lesssim 
		(1+t) |\mathscr{Z}^{\leq N} \Psi|
		+ \sum_{a=1}^3 \rgeo \left|\mathscr{Z}^{\leq N} \Lunit_{(Small)}^a \right|,
			\\
	\left\|
		\mathscr{Z}^{\leq 12} \RotRadcomponent{l}
	\right\|_{C^0(\Sigma_t^u)}
	& \lesssim \varepsilon^{1/2} \ln(\myexp + t).
		\label{E:LOWERORDERC0BOUNDEUCLIDEANROTATIONRADCOMPONENT}
\end{align}
\end{subequations}

\end{lemma}

\begin{proof}
	We first prove \eqref{E:DERIVATIVESOFRECTANGULARCOMMUTATORVECTORFIELDRECTCOMPONENTS} by induction.
	Throughout this proof, we use the 	
	bootstrap assumptions
	\eqref{E:PSIFUNDAMENTALC0BOUNDBOOTSTRAP},
	\eqref{E:UPMUBOOT},
	and
	\eqref{E:FRAMECOMPONENTSBOOT}.
	
	In the base case $N=0,$ we have to bound the magnitude of
	$\rgeo \Lunit^i = x^i + \rgeo \Lunit_{(Small)}^{j},$
	$\Rad^i = - \upmu \frac{x^i}{\rgeo} + \upmu \Radunit_{(Small)}^i,$
	and
	$\Rot_{(l)}^i = \epsilon_{laj} x^a + \RotRadcomponent{l} \frac{x^i}{\rgeo} - \RotRadcomponent{l} \Radunit_{(Small)}^i$
	by the right-hand side of \eqref{E:DERIVATIVESOFRECTANGULARCOMMUTATORVECTORFIELDRECTCOMPONENTS}.
	Recall that an expression for $\RotRadcomponent{l}$ is given by \eqref{E:EUCLIDEANROTATIONRADCOMPONENT}.
	We bound all factors $x^i$ via \eqref{E:EASYCOORDINATEBOUND}.
	Furthermore, 
	$g_{ab}^{(Small)}$ is a smooth function of
	$\Psi$ that vanishes at $\Psi = 0,$
	and \eqref{E:RADUNITJUNKLIKELMINUSUNITJUNK} implies that
	$\Radunit^i = - \Lunit_{(Small)}^i$ plus a smooth function of $\Psi$
	that vanishes at $\Psi = 0,$ that is, $|\Radunit_{(Small)}^i| \lesssim |\Lunit_{(Small)}^i| + |\Psi|.$
	From these facts and the bootstrap assumptions,
	we deduce the desired estimate \eqref{E:DERIVATIVESOFRECTANGULARCOMMUTATORVECTORFIELDRECTCOMPONENTS} when $N = 0.$	
	
	To carry out the induction, we apply 
	$\mathscr{Z}^N$ to the right-hand side of the above identities for 
	$\rgeo \Lunit^i,$
	$\Rad^i,$
	$\Rot_{(l)}^i.$
	Using the Leibniz rule,
	we can easily bound most terms that arise via 
	the bootstrap assumptions,
	the estimate \eqref{E:EASYCOORDINATEBOUND},
	and the estimate \eqref{E:ZNAPPLIEDTORGEOISNOTTOOLARGE}
	(which we use to bound factors involving derivatives of $\rgeo$]
	as in the previous paragraph (without the need for induction).
	The factors that require the induction hypotheses are those
	of the form $\mathscr{Z}^M x^i,$ $1 \leq M \leq N.$ 
	To bound these terms, we schematically write $\mathscr{Z}^M = \mathscr{Z}^{M-1} Z$
	and hence $\mathscr{Z}^M x^i = \mathscr{Z}^{M-1} Z^i,$
	which can be bounded by the induction hypotheses. We have thus proved 	
	\eqref{E:DERIVATIVESOFRECTANGULARCOMMUTATORVECTORFIELDRECTCOMPONENTS}.
	\eqref{E:LOWERORDERC0BOUNDDERIVATIVESOFRECTANGULARCOMMUTATORVECTORFIELDRECTCOMPONENTS}
	then follows from \eqref{E:DERIVATIVESOFRECTANGULARCOMMUTATORVECTORFIELDRECTCOMPONENTS}
	and the bootstrap assumptions.
	
	The estimate \eqref{E:FIRSTPOINTWISEBOUNDEUCLIDEANROTATIONRADCOMPONENT}
	then follows from \eqref{E:EUCLIDEANROTATIONRADCOMPONENT}, 
	the fact that
	$\Radunit^i = - \Lunit_{(Small)}^i$ plus a smooth function of $\Psi$ that vanishes at $\Psi = 0,$
	the bootstrap assumptions,
	and the estimate \eqref{E:DERIVATIVESOFRECTANGULARCOMMUTATORVECTORFIELDRECTCOMPONENTS}.
	\eqref{E:LOWERORDERC0BOUNDEUCLIDEANROTATIONRADCOMPONENT} then follows from
	\eqref{E:FIRSTPOINTWISEBOUNDEUCLIDEANROTATIONRADCOMPONENT} and the bootstrap assumptions.
	
	The estimate \eqref{E:DERIVATIVESOFRECTANGULARFRAMECOMPONENTS} for $\Lunit^i$ 
	follows similarly from the decomposition $\Lunit^i = \frac{x^i}{\rgeo} + \Lunit_{(Small)}^i$
	and the estimates \eqref{E:ZNAPPLIEDTORGEOISNOTTOOLARGE} 
	and \eqref{E:DERIVATIVESOFRECTANGULARCOMMUTATORVECTORFIELDRECTCOMPONENTS}. 
	The estimate \eqref{E:DERIVATIVESOFRECTANGULARFRAMECOMPONENTS}
	for $\Radunit^i$ then follows from the estimate for $\mathscr{Z}^{\leq N} \Lunit^i,$
	the fact that
	$\Radunit^i = - \Lunit_{(Small)}^i$ plus a smooth function of $\Psi$ that vanishes at $\Psi = 0,$
	and the bootstrap assumptions.
	Inequality \eqref{E:POINTWISEBOUNDLOWERORDERDERIVATIVESOFRECTANGULARFRAMECOMPONENTS}
	then follows from \eqref{E:DERIVATIVESOFRECTANGULARFRAMECOMPONENTS}
	and the bootstrap assumptions.
	
\end{proof}

\section{Estimates for the rectangular components of the \texorpdfstring{$g-$dual}{metric dual} of 
\texorpdfstring{$\Radunit$}{the unit-length radial vectorfield}}
We now derive some simple pointwise estimates for 
$\Radunit_i + \frac{x^i}{\rgeo}$ and $\Radunit_i^{(Small)}.$

We start with the following lemma, 
in which we derive higher-order analogs of inequality \eqref{E:EASYEUCLIDEANRADIALVARIABLEBOUND}

\begin{lemma}[\textbf{Estimates for} $r/\rgeo$]
	Under the small-data and bootstrap assumptions 
	of Sects.~\ref{S:PSISOLVES}-\ref{S:C0BOUNDBOOTSTRAP},
	if $\varepsilon$ is sufficiently small, 
	then the following estimates hold for
	$(t,u) \in [0,\Tboot) \times [0,U_0]:$
	\begin{align}
	\left\|
		\mathscr{Z}^{\leq 12}
			\left(
				\frac{r}{\rgeo} - 1
			\right)
	\right\|_{C^0(\Sigma_t^u)}	
	& \lesssim
		\varepsilon^{1/2} \frac{\ln(\myexp + t)}{1 + t}.
		\label{E:HIGHERORDEEUCLIDEANRADIALVARIABLEBOUND}
	\end{align}

\end{lemma}

\begin{proof}
	Inequality \eqref{E:HIGHERORDEEUCLIDEANRADIALVARIABLEBOUND} follows from
	applying $\mathscr{Z}^{\leq 12}$ to \eqref{E:USEFULIDEFORRADIALCOORDINATES}
	and using the estimate \eqref{E:POINTWISEBOUNDLOWERORDERDERIVATIVESOFRECTANGULARFRAMECOMPONENTS}
	and the bootstrap assumptions
	\eqref{E:PSIFUNDAMENTALC0BOUNDBOOTSTRAP}
	and \eqref{E:FRAMECOMPONENTSBOOT}.
\end{proof}

\begin{lemma}[\textbf{Estimates for} $\Radunit_i + \frac{x^i}{\rgeo}$ and $\Radunit_i^{(Small)}.$]
	Let $0 \leq N \leq 24$ be an integer.
	Under the small-data and bootstrap assumptions 
	of Sects.~\ref{S:PSISOLVES}-\ref{S:C0BOUNDBOOTSTRAP},
	if $\varepsilon$ is sufficiently small, 
	then for $i=1,2,3,$ 
	the following pointwise estimates hold on the spacetime domain $\mathcal{M}_{\Tboot,U_0}:$
	\begin{subequations}
	\begin{align}
	\left|
			\mathscr{Z}^N
			\left(
				\Radunit_i
				+ \frac{x^i}{\rgeo}
			\right)
		\right|
		& \lesssim 
			|\mathscr{Z}^{\leq N} \Psi| 
			+ \sum_{a=1}^3 \left| \mathscr{Z}^{\leq N} \Lunit_{(Small)}^a \right|,
			\label{E:RADUNITLOWEREDBACKGROUNDSUBTRACTEDPOINTWISE}
				\\
		\left|
			\mathscr{Z}^N \Radunit_i^{(Small)}
		\right|
		& \lesssim 
			|\mathscr{Z}^{\leq N} \Psi| 
			+ \sum_{a=1}^3 \left| \mathscr{Z}^{\leq N} \Lunit_{(Small)}^a \right|.
			\label{E:RADUNITSMALLLOWEREDBACKGROUNDSUBTRACTEDPOINTWISE}
	\end{align}
	\end{subequations}
	
	Furthermore, the following estimates hold:
	\begin{subequations}
	\begin{align}
	\left\|
			\mathscr{Z}^{\leq 12}
			\left(
				\Radunit_i
				+ \frac{x^i}{\rgeo}
			\right)
		\right\|_{C^0(\Sigma_t^u)}
		& \lesssim \varepsilon^{1/2} \frac{\ln(\myexp + t)}{1 + t},
		\label{E:RADUNITLOWEREDBACKGROUNDSUBTRACTEDC0}
			\\
		\left\|
			\mathscr{Z}^{\leq 12}
			\Radunit_i^{(Small)}
		\right\|_{C^0(\Sigma_t^u)}
		& \lesssim \varepsilon^{1/2} \frac{\ln(\myexp + t)}{1 + t}.
		\label{E:RADUNITSMALLLOWEREDBACKGROUNDSUBTRACTEDC0}
	\end{align}
	\end{subequations}
	
	In addition,
	\begin{align} \label{E:ALTERNATEVERSIONRADUNITLOWEREDBACKGROUNDSUBTRACTEDC0}
		\mbox{inequality \eqref{E:RADUNITLOWEREDBACKGROUNDSUBTRACTEDC0} 
		holds with the term $\frac{x^i}{\rgeo}$ on the left replaced by $\frac{x^i}{r}.$}
	\end{align}
	
\end{lemma}

\begin{proof}
	To prove \eqref{E:RADUNITLOWEREDBACKGROUNDSUBTRACTEDPOINTWISE},
	we first use 
	\eqref{E:LITTLEGDECOMPOSED},
	\eqref{E:LUNITJUNK},
	and
	\eqref{E:UPSTAIRSSRADUNITPLUSLUNITISAFUNCTIONOFPSI},
	to deduce the identity 
	\begin{align} \label{E:RADUNITILOWEREDEXPANDED}
		\Radunit_i 
		+ \frac{x^i}{\rgeo}
		& = 
			- g_{ia}^{(Small)} \frac{x^a}{\rgeo}
			- g_{ia} \Lunit_{(Small)}^a 
			- g_{ia} (g^{-1})^{0a}.
	\end{align}
	Applying $\mathscr{Z}^N$ to
	\eqref{E:RADUNITILOWEREDEXPANDED},
	using the estimates
	\eqref{E:ZNAPPLIEDTORGEOISNOTTOOLARGE},
	\eqref{E:DERIVATIVESOFRECTANGULARCOMMUTATORVECTORFIELDRECTCOMPONENTS},
	and \eqref{E:LOWERORDERC0BOUNDDERIVATIVESOFRECTANGULARCOMMUTATORVECTORFIELDRECTCOMPONENTS},
	and using the bootstrap assumptions
	\eqref{E:PSIFUNDAMENTALC0BOUNDBOOTSTRAP},
	and \eqref{E:FRAMECOMPONENTSBOOT},
	we conclude the desired estimate
	\eqref{E:RADUNITLOWEREDBACKGROUNDSUBTRACTEDPOINTWISE}.
	The estimate \eqref{E:RADUNITLOWEREDBACKGROUNDSUBTRACTEDC0}
	then follows from
	\eqref{E:RADUNITLOWEREDBACKGROUNDSUBTRACTEDPOINTWISE}
	and the bootstrap assumptions
	\eqref{E:PSIFUNDAMENTALC0BOUNDBOOTSTRAP},
	and \eqref{E:FRAMECOMPONENTSBOOT}.
	
	The proofs of \eqref{E:RADUNITSMALLLOWEREDBACKGROUNDSUBTRACTEDPOINTWISE} and
	\eqref{E:RADUNITSMALLLOWEREDBACKGROUNDSUBTRACTEDC0} are similar,
	and we omit the details.
	
	Inequality \eqref{E:ALTERNATEVERSIONRADUNITLOWEREDBACKGROUNDSUBTRACTEDC0}
	follows from expanding
	$x^i/\rgeo = (x^i/r) + (x^i/r)(r/\rgeo - 1)$
	and using the identity $Z x^i = Z^i$ 
	and the estimates
	\eqref{E:ZNAPPLIEDTORGEOISNOTTOOLARGE},
	\eqref{E:EASYCOORDINATEBOUND},
	\eqref{E:LOWERORDERC0BOUNDDERIVATIVESOFRECTANGULARCOMMUTATORVECTORFIELDRECTCOMPONENTS},
	\eqref{E:HIGHERORDEEUCLIDEANRADIALVARIABLEBOUND},
	and
	\eqref{E:RADUNITLOWEREDBACKGROUNDSUBTRACTEDC0}.
\end{proof}

\section{Precise pointwise estimates for the rotation vectorfields}
\label{S:ROTATIONPRECISEPOINTWISE}
In this section, we derive some sharp pointwise estimates for the rotation vectorfields
and their spherical covariant derivatives. We carefully estimate
the non-small constants such as the ``$1$'' on the right-hand side of \eqref{E:ROTATIONPOINTWISENORMESTIMATE}.
The reason that we are interested in precise constants
is that the same constants appear
in the comparison estimates of Sect.~\ref{S:POINTWISEDIFFERENTIALOPERATORCOMPARISONESTIMATES},
and these comparison estimates ultimately affect the number of derivatives we need to close 
our top-order $L^2$ estimates.

\begin{lemma}[\textbf{Precise pointwise estimates for the rotation vectorfields}]
\label{L:ROTATIONVECTORFIELDSBASICESTIMATES}
Under the small-data and bootstrap assumptions 
of Sects.~\ref{S:PSISOLVES}-\ref{S:C0BOUNDBOOTSTRAP},
if $\varepsilon$ is sufficiently small,
then the following pointwise estimates hold on $\mathcal{M}_{\Tboot,U_0}:$
\begin{subequations}
\begin{align}
	|\Rot_{(l)}| 
	& \leq (1 + C \varepsilon^{1/2}) \rgeo,
		\label{E:ROTATIONPOINTWISENORMESTIMATE} \\
	\left|[\Rot_{(l)}, \Rot_{(m)}] \right| 
	& \leq (1 + C \varepsilon^{1/2}) \rgeo,
		\label{E:POINTWISEESTIMATEFORROTATIONCOMMUTATORS}	
		\\
	\left|
		\angD \Rot_{(l)}
	\right|^2
	& \leq 2(1 + C \varepsilon^{1/2}).
	\label{E:POINTWISEBOUNDSPHERECOVARIANTDERIVATIVEOFROTATIONS}
\end{align}
\end{subequations}

Furthermore, the following pointwise estimates for rectangular components hold on $\mathcal{M}_{\Tboot,U_0},$
($a,b,i,j,m,n=1,2,3$):
\begin{subequations}
\begin{align} \label{E:GINVERSESPHEREISWELLAPPROXIMATEDBYROTATIONTENSORPRODUCTS}
		\left|
			(\ginversesphere)^{mn}
			-
				\rgeo^{-2}
				\sum_{l=1}^3 
				\Rot_{(l)}^m
				\Rot_{(l)}^n
		\right|
		& \lesssim \varepsilon^{1/2} \frac{\ln(\myexp + t)}{1 + t},
			\\
	\left|
			\gsphere_{ij} (\ginversesphere)^{mn} 
			- \sphereproject_i^{\ n} \sphereproject_j^{\ m}
			-
			\sum_{l=1}^3 
				(\angDarg{i} \Rot_{(l)}^m)
				(\angDarg{j} \Rot_{(l)}^n)
		\right|
		& \lesssim \varepsilon^{1/2}
		\frac{\ln(\myexp + t)}{1 + t},
			\label{E:GINVERSEPROJECTIONCOMBINATIONISWELLAPPROXIMATEDBYROTATIONCOVARIANTDERIVATIVETENSORPRODUCTS} 
			\\
		\left|
			(\ginversesphere)^{mn}
			-
			\sum_{l=1}^3 
				(\ginversesphere)^{ab}
				(\angDarg{a} \Rot_{(l)}^m)
				(\angDarg{b} \Rot_{(l)}^n)
		\right|
		& \lesssim \varepsilon^{1/2}
		\frac{\ln(\myexp + t)}{1 + t},
			\label{E:GINVERSESPHEREISWELLAPPROXIMATEDBYROTATIONCOVARIANTDERIVATIVETENSORPRODUCTS} \\
		\left|
			\sum_{l=1}^3 
				\Rot_{(l)}^m
				\angDarg{b} \Rot_{(l)}^n
		\right|
		& \lesssim \varepsilon^{1/2} \ln(\myexp + t).
		\label{E:ROTANGDROTTENSORPRODUCTSUMISSMALL}
	\end{align}	
	\end{subequations}
\end{lemma}

\begin{proof}
We first prove \eqref{E:ROTATIONPOINTWISENORMESTIMATE}. 
Since $\Rot_{(l)}$ is $S_{t,u}-$tangent, the following identities holds
relative to the rectangular spatial coordinates:
$|\Rot_{(l)}|^2 = \gsphere_{ab} \Rot_{(l)}^a \Rot_{(l)}^b = g_{ab} \Rot_{(l)}^a \Rot_{(l)}^b.$
From equations \eqref{E:LITTLEGDECOMPOSED} and 
\eqref{E:ROTATIONDECOMPOSITIONINTOEUCLIDEANPLUSRADCOMPONENT}
and the fact that $g_{ab} \Radunit^a \Radunit^b = 1,$
we deduce that 
$|\Rot_{(l)}|^2 \leq  \delta_{ab} \Roteucarg{l}^a \Roteucarg{l}^b 
+ C |g_{ab}^{(Small)} \Roteucarg{l}^a \Roteucarg{l}^b|
+ C |\RotRadcomponent{l}|^2.$ 
From \eqref{E:EULCIDEANROTATION}
and 
the estimates 
\eqref{E:RGEOMINKOWSKITIMECOMPARISON}
and
\eqref{E:EASYEUCLIDEANRADIALVARIABLEBOUND},
we deduce that
$\delta_{ab} \Roteucarg{l}^a \Roteucarg{l}^b \leq \sum_{a=1}^3 (x^a)^2 = r^2 \leq (1 + C \varepsilon^{1/2}) \rgeo^2.$
Inequality \eqref{E:ROTATIONPOINTWISENORMESTIMATE}
therefore easily follows once we show that
$|g_{ab}^{(Small)} \Roteucarg{l}^a \Roteucarg{l}^b| \leq C \varepsilon^{1/2} \rgeo$
and $|\RotRadcomponent{l}| \leq C \varepsilon^{1/2} \rgeo.$
The first of these two estimates follows from
definition \eqref{E:EULCIDEANROTATION}
and the estimates 
\eqref{E:RGEOMINKOWSKITIMECOMPARISON},
\eqref{E:EASYCOORDINATEBOUND},
and \eqref{E:GRECTCOMPONENTC0BOUND}.
The second estimate follows 
from \eqref{E:LOWERORDERC0BOUNDEUCLIDEANROTATIONRADCOMPONENT}.

We now prove \eqref{E:POINTWISEESTIMATEFORROTATIONCOMMUTATORS}.	
We take the norm of each term on the right-hand side of
\eqref{E:ROTATATIONCOMMUTATORS}. 
By \eqref{E:ROTATIONPOINTWISENORMESTIMATE}, the 
first term is in magnitude $\leq (1 + C \varepsilon^{1/2}) \rgeo.$
It remains for us to show that the terms on the right-hand
side of \eqref{E:ROTATIONCOMMUTATORSERRORONEFORM} are in magnitude
$\lesssim \varepsilon^{1/2} \rgeo.$ To see that the terms on
the first line of the right-hand side of \eqref{E:ROTATIONCOMMUTATORSERRORONEFORM}
are in fact $\lesssim \varepsilon \ln^2(\myexp + t)(1 + t)^{-1},$ we use the bootstrap 
assumption \eqref{E:CHIJUNKBOOT}
and the estimates \eqref{E:LOWERORDERC0BOUNDEUCLIDEANROTATIONRADCOMPONENT} and
\eqref{E:ROTATIONPOINTWISENORMESTIMATE}.
To see that the terms on
the third line of the right-hand side of \eqref{E:ROTATIONCOMMUTATORSERRORONEFORM}
are in fact $\lesssim \varepsilon \ln^2(\myexp + t)(1 + t)^{-1},$
we treat all lowercase Latin-indexed quantities as scalar valued functions and the
$\angdiff x^c$ as $S_{t,u}$ one-forms. The desired bound follows from 
the bootstrap assumption \eqref{E:FRAMECOMPONENTSBOOT}
and the estimates
\eqref{E:FIRSTESTIMATEFORANGULARDIFFERNTIALOFXI},
\eqref{E:GRECTCOMPONENTC0BOUND},
and \eqref{E:LOWERORDERC0BOUNDEUCLIDEANROTATIONRADCOMPONENT}.
To see that the terms on
the second line of the right-hand side of \eqref{E:ROTATIONCOMMUTATORSERRORONEFORM}
are in fact $\lesssim \varepsilon \ln(\myexp + t)(1 + t)^{-1},$
we use 
the estimates
\eqref{E:LOWERORDERC0BOUNDEUCLIDEANROTATIONRADCOMPONENT},
\eqref{E:ROTATIONPOINTWISENORMESTIMATE},
and the estimate 
$|\Theta^{(Tan-\Psi)}| \lesssim \varepsilon (1+t)^{-2}.$
To prove this latter bound, we first use \eqref{E:BIGTHETAGOOD}
to write $\Theta^{(Tan-\Psi)}$ in the schematic form
$\Theta^{(Tan-\Psi)} = G_{(Frame)} \myarray[\Lunit \Psi]{\angdiff \Psi}.$
By the bootstrap assumptions \eqref{E:PSIFUNDAMENTALC0BOUNDBOOTSTRAP},
the terms in the array of $\Psi$ derivatives are $\lesssim \varepsilon (1 + t)^{-2}$ in magnitude.
It remains for us to show that $|G_{(Frame)}| \lesssim 1.$
To this end, we note that the entries of the array $G_{(Frame)}$ can be viewed as smooth functions
of $\Psi,$ of the rectangular components of the vectorfields
$\Lunit$ and $\Radunit,$ and of the rectangular components of the $S_{t,u}$ projection tensorfield $\sphereproject.$
Hence, from the bootstrap assumptions \eqref{E:PSIFUNDAMENTALC0BOUNDBOOTSTRAP}
and \eqref{E:FRAMECOMPONENTSBOOT} and the estimate 
\eqref{E:GSPHERERECTCOMPONENTSPOINTWISEBOUND},
we deduce that $|G_{(Frame)}| \lesssim 1$ as desired.

We now prove \eqref{E:GINVERSESPHEREISWELLAPPROXIMATEDBYROTATIONTENSORPRODUCTS}.
To begin, we first use \eqref{E:GINVERSEFRAMEWITHRECTCOORDINATESFORGSPHEREINVERSE}
to expand 
\begin{align}  \label{E:GINVERSESPHEREEASYEXPANSION}
	(\ginversesphere)^{mn}  
	& = \left\lbrace
				\delta^{mn} - \frac{x^m x^n}{r^2} 
			\right\rbrace
		+ \left\lbrace 
				\frac{1}{r^2} - \frac{1}{\rgeo^2}
			\right\rbrace
			x^m x^n
		+ 
		\left\lbrace
				(\ginversesphere)^{mn} 
				- \left(\delta^{mn} - \frac{x^m x^n}{\rgeo^2} \right) 
		\right\rbrace.
\end{align}
Using \eqref{E:EULCIDEANROTATION}, 
it is straightforward to verify that the first difference on the right-hand side of \eqref{E:GINVERSESPHEREEASYEXPANSION}
can be expressed as
\begin{align} \label{E:EUCLIDEANINVERSESPHEREMETRICINTERMSOFEUCLIDEANROTATIONS}
	\left\lbrace 
		\delta^{mn} - \frac{x^m x^n}{r^2} 
	\right\rbrace = \sum_{l=1}^3 r^{-2} \Roteucarg{l}^m \Roteucarg{l}^n.
\end{align}
From equation \eqref{E:ROTATIONDECOMPOSITIONINTOEUCLIDEANPLUSRADCOMPONENT}
and the estimates
\eqref{E:POINTWISEBOUNDLOWERORDERDERIVATIVESOFRECTANGULARFRAMECOMPONENTS}
and
\eqref{E:LOWERORDERC0BOUNDEUCLIDEANROTATIONRADCOMPONENT},
we deduce that
$|\Roteucarg{l}^m - \Rot_{(l)}^m| = |\RotRadcomponent{l}||\Radunit^m|
\lesssim \varepsilon^{1/2} \ln(\myexp + t).$ Also using the estimates 
\eqref{E:EASYCOORDINATEBOUND} and
\eqref{E:EASYEUCLIDEANRADIALVARIABLEBOUND},
we deduce from \eqref{E:EUCLIDEANINVERSESPHEREMETRICINTERMSOFEUCLIDEANROTATIONS} that
$\left\lbrace \delta^{mn} - \frac{x^m x^n}{\rgeo^2} \right\rbrace = \sum_{l=1}^3 \rgeo^{-2} \Rot_{(l)}^m \Rot_{(l)}^n$
plus an error term that is in magnitude $\lesssim \varepsilon^{1/2} \ln(\myexp + t)(1 + t)^{-1}.$

To bound the second difference $\left\lbrace 
				\frac{1}{r^2} - \frac{1}{\rgeo^2}
			\right\rbrace
			x^m x^n$
on the right-hand side of \eqref{E:GINVERSESPHEREEASYEXPANSION},
we use the estimates \eqref{E:EASYCOORDINATEBOUND} and
\eqref{E:EASYEUCLIDEANRADIALVARIABLEBOUND}
to deduce that it is in magnitude $\lesssim \varepsilon^{1/2} \ln(\myexp + t)(1 + t)^{-1}$
as desired.

To bound the final difference on the right-hand side of \eqref{E:GINVERSESPHEREEASYEXPANSION},
we quote the estimate \eqref{E:GTINVERSERECTCOMPONENTSC0BOUND} to deduce that
its magnitude is $\lesssim \varepsilon^{1/2} \ln(\myexp + t)(1 + t)^{-1}$ as desired.

We now prove \eqref{E:GINVERSEPROJECTIONCOMBINATIONISWELLAPPROXIMATEDBYROTATIONCOVARIANTDERIVATIVETENSORPRODUCTS}.
We use the identity \eqref{E:ANGDROTATION} and as usual, 
we treat all uppercase Latin indices as tensorial $S_{t,u}$ indices,
and all lowercase Latin-indexed quantities as scalar-valued functions.
The estimates proved above imply that 
the norm of the last three products on the right-hand side of \eqref{E:ANGDROTATION}
are (viewed as $S_{t,u}$ tensors) in magnitude $\lesssim \varepsilon^{1/2} \ln(\myexp + t)(1 + t)^{-1}.$
Furthermore, the first term
$\epsilon_{lab} \angdiff x^a \otimes \angdiffuparg{\#} x^b$
is (viewed as an $S_{t,u}$ tensor) in magnitude by $\lesssim 1.$
By Cor.~\ref{C:NORMOFTHESTUTENSORISCOMPARABLETOTHESUMOFITSCOMPONENTS},
the same estimates hold for the rectangular components of these tensors.
Furthermore, using the identities $\angdiffarg{i} x^a = \sphereproject_i^{\ a}$ and $\angdiffuparg{\# i} x^j = (\ginversesphere)^{ij},$
we note that the $i,m$ rectangular components of the type $\binom{1}{1}$ 
$S_{t,u}$ tensor $\epsilon_{lab} \angdiff x^a \otimes \angdiffuparg{\#} x^b$
are $\epsilon_{lab} \sphereproject_i^{\ a} (\ginversesphere)^{mb}.$
Hence, the main contribution to the sum
$\sum_{l=1}^3 (\angDarg{i} \Rot_{(l)}^m)(\angDarg{j} \Rot_{(l)}^n)$ comes from the product term 
$\sum_{l=1}^3 \epsilon_{lab} \epsilon_{lcd} \sphereproject_i^{\ a} (\ginversesphere)^{mb} \sphereproject_j^{\ c}  (\ginversesphere)^{nd}$
$= \delta^{kl}\epsilon_{kab} \epsilon_{lcd} \sphereproject_i^{\ a} (\ginversesphere)^{mb} \sphereproject_j^{\ c}  (\ginversesphere)^{nd},$
and the desired estimate \eqref{E:GINVERSEPROJECTIONCOMBINATIONISWELLAPPROXIMATEDBYROTATIONCOVARIANTDERIVATIVETENSORPRODUCTS}
will follow once we prove the following estimate:
\begin{align}   \label{E:MAINESTIMATEGINVERSEPROJECTIONCOMBINATIONISWELLAPPROXIMATEDBYROTATIONCOVARIANTDERIVATIVETENSORPRODUCTS}
	\left| 
		\gsphere_{ij} (\ginversesphere)^{mn}
		- \sphereproject_i^{\ n} \sphereproject_j^{\ m}
		- \delta^{kl}\epsilon_{kab} \epsilon_{lcd} \sphereproject_i^{\ a} (\ginversesphere)^{mb} \sphereproject_j^{\ c}  (\ginversesphere)^{nd}
	\right|
	\lesssim \varepsilon^{1/2} 
		\frac{\ln(\myexp + t)}{1 + t}.
\end{align}
To proceed, we first use the identity
$\delta^{kl} \epsilon_{kab} \epsilon_{lcd} = \delta_{ac} \delta_{bd} - \delta_{ad} \delta_{bc}$
to deduce that the final product term on the left-hand side of \eqref{E:MAINESTIMATEGINVERSEPROJECTIONCOMBINATIONISWELLAPPROXIMATEDBYROTATIONCOVARIANTDERIVATIVETENSORPRODUCTS}
is equal to
\begin{align} \label{E:ARATHERCOMPLICATEDPRODUCTREDUCEDFORM}
	\left\lbrace
		\delta_{ac} \delta_{bd} - \delta_{ad} \delta_{bc}
	\right\rbrace 
	\sphereproject_i^{\ a} (\ginversesphere)^{mb} \sphereproject_j^{\ c}  (\ginversesphere)^{nd}.
\end{align}
Thus, to conclude \eqref{E:GINVERSEPROJECTIONCOMBINATIONISWELLAPPROXIMATEDBYROTATIONCOVARIANTDERIVATIVETENSORPRODUCTS}, 
it only remains to show that
\eqref{E:ARATHERCOMPLICATEDPRODUCTREDUCEDFORM} is equal to 
$\gsphere_{ij} (\ginversesphere)^{mn} - \sphereproject_i^{\ n} \sphereproject_j^{\ m}$ 
plus error terms that are in magnitude $\lesssim \varepsilon^{1/2} \ln(\myexp + t)(1 + t)^{-1}.$
To this end, we use \eqref{E:GSPHEREINTERMSOFRADUNITRECTCOMPONENTSC0BOUNDESTIMATE}
to replace $\delta_{ac} \delta_{bd} - \delta_{ad} \delta_{bc}$
with $\gsphere_{ac} \gsphere_{bd} - \gsphere_{ad} \gsphere_{bc}$
plus an error term that is in magnitude $\lesssim \varepsilon^{1/2} \ln(\myexp + t)(1 + t)^{-1}$
plus tensorial products of one-forms that are $g-$dual to $\Radunit$ and
hence are annihilated by the
term $\sphereproject_i^{\ a} (\ginversesphere)^{mb} \sphereproject_j^{\ c}  (\ginversesphere)^{nd}$
in \eqref{E:ARATHERCOMPLICATEDPRODUCTREDUCEDFORM}.
Finally, after this-replacement-up-to-errors, we compute that 
\begin{align}
\left\lbrace \gsphere_{ac} \gsphere_{bd} - \gsphere_{ad} \gsphere_{bc} \right\rbrace 
\sphereproject_i^{\ a} (\ginversesphere)^{mb} \sphereproject_j^{\ c}  (\ginversesphere)^{nd}
= \gsphere_{ij} (\ginversesphere)^{mn} - \sphereproject_i^{\ n} \sphereproject_j^{\ m}
\end{align}
as desired.

We now prove \eqref{E:GINVERSESPHEREISWELLAPPROXIMATEDBYROTATIONCOVARIANTDERIVATIVETENSORPRODUCTS}.
To this end, we contract \eqref{E:GINVERSEPROJECTIONCOMBINATIONISWELLAPPROXIMATEDBYROTATIONCOVARIANTDERIVATIVETENSORPRODUCTS}
against $(\ginversesphere)^{ij}$ and use \eqref{E:GSPHERERECTCOMPONENTSPOINTWISEBOUND} to deduce
that the contracted quantity is also $\lesssim$ the right-hand side of
\eqref{E:GINVERSESPHEREISWELLAPPROXIMATEDBYROTATIONCOVARIANTDERIVATIVETENSORPRODUCTS}.
This immediately yields the desired estimate \eqref{E:GINVERSESPHEREISWELLAPPROXIMATEDBYROTATIONCOVARIANTDERIVATIVETENSORPRODUCTS}.

We now prove \eqref{E:POINTWISEBOUNDSPHERECOVARIANTDERIVATIVEOFROTATIONS} with the help of
the identity \eqref{E:ANGDROTATION}.
In our proof of \eqref{E:GINVERSESPHEREISWELLAPPROXIMATEDBYROTATIONCOVARIANTDERIVATIVETENSORPRODUCTS},
we showed that the norms of the last three terms on the right-hand side of \eqref{E:ANGDROTATION}
are (viewed as $S_{t,u}$ tensors) in magnitude $\lesssim \varepsilon \frac{1}{1 + t}.$
Using the identity $\angdiffuparg{A} x^a \angdiffarg{A} x^c = (\ginversesphere)^{ac},$
we deduce that the square of the norm of the first term on the right-hand side of \eqref{E:ANGDROTATION}
is (viewed as an $S_{t,u}$ tensor) equal to
\begin{align} \label{E:MAINTERMSTUCOVARIANTDERIVATIVEOFROTATIONSMAINTERMSQUARENORM}
	|\epsilon_{lab} \epsilon_{lcd} (\ginversesphere)^{ac} (\ginversesphere)^{bd}|,
\end{align}
where there is no summation over $l$ in \eqref{E:MAINTERMSTUCOVARIANTDERIVATIVEOFROTATIONSMAINTERMSQUARENORM}.
The expression \eqref{E:MAINTERMSTUCOVARIANTDERIVATIVEOFROTATIONSMAINTERMSQUARENORM} can be viewed
as the square of the norm of the projection of the type $\binom{0}{2}$ $\Sigma_t$ tensor $\epsilon_{l \cdots}$
onto $S_{t,u}.$ Using the fact that the norm of the projection of such a tensor is no larger than the norm of the tensor itself,
we deduce that \eqref{E:MAINTERMSTUCOVARIANTDERIVATIVEOFROTATIONSMAINTERMSQUARENORM} is 
$\leq \epsilon_{lab} \epsilon_{lcd} (\gtinverse)^{ac} (\gtinverse)^{bd}.$
Using \eqref{E:GTINVERSERECTCOMPONENTSC0BOUND}, 
we see that we can replace $(\gtinverse)^{ac}$ $(\gtinverse)^{bd}$ with
$\delta^{ac}$ and $\delta^{bd}$ up to an error term that is in magnitude $\lesssim \varepsilon (1 + t)^{-1}.$
The remaining term verifies
\begin{align} \label{E:POINTWISEBOUNDMAINTERMSTUCOVARIANTDERIVATIVEOFROTATIONSMAINTERMSQUARENORM}
	|\epsilon_{lab} \epsilon_{lcd} \delta^{ac} \delta^{bd}| = 2.
\end{align}
The desired inequality \eqref{E:POINTWISEBOUNDSPHERECOVARIANTDERIVATIVEOFROTATIONS} thus follows.

We finally prove \eqref{E:ROTANGDROTTENSORPRODUCTSUMISSMALL}.
Using Cor.~\ref{C:NORMOFTHESTUTENSORISCOMPARABLETOTHESUMOFITSCOMPONENTS}, 
we see that it suffices
to prove the following type $\binom{2}{1}$ $S_{t,u}$ tensor bound:
\begin{align} \label{E:GEOMETRICVERSIONROTANGDROTTENSORPRODUCTSUMISSMALL}
	\left| 
		\sum_{l=1}^3 \Rot_{(l)} \otimes \angD \Rot_{(l)}
	\right|
	&
	\lesssim \varepsilon^{1/2} \ln(\myexp + t).
\end{align}
Using the identities 
\eqref{E:ROTATATIONACOMPONENTALTERNATEEXPRESSION}
and \eqref{E:ANGDROTATION} and the line of reasoning
that we used at the beginning of the proof of
\eqref{E:GINVERSEPROJECTIONCOMBINATIONISWELLAPPROXIMATEDBYROTATIONCOVARIANTDERIVATIVETENSORPRODUCTS}, 
we see that the main contribution to the tensor on the left-hand side of 
\eqref{E:GEOMETRICVERSIONROTANGDROTTENSORPRODUCTSUMISSMALL} is
$\delta^{kl}\epsilon_{kab} \epsilon_{lcd} x^a \angdiffuparg{\#} x^b \otimes \angdiff x^c \otimes \angdiffuparg{\#} x^d$
and that all remaining error terms are in magnitude $\lesssim \varepsilon.$
That is, it suffices to prove the following $S_{t,u}$ tensor bound:
\begin{align}  \label{E:MAINTERMROTANGDROTTENSORPRODUCTSUMISSMALL}
	\left|
		\delta^{kl}\epsilon_{kab} \epsilon_{lcd} x^a \angdiffuparg{\#} x^b \otimes \angdiff x^c \otimes \angdiffuparg{\#} x^d
	\right| 
	& \lesssim \varepsilon^{1/2} \ln(\myexp + t).
\end{align}
Using Cor.~\ref{C:NORMOFTHESTUTENSORISCOMPARABLETOTHESUMOFITSCOMPONENTS} 
and the identity $\delta^{kl} \epsilon_{kab} \epsilon_{lcd} = \delta_{ac} \delta_{bd} - \delta_{ad} \delta_{bc},$
we see that it suffices to prove the same bound for the rectangular 
$m,n,p$ components of the $S_{t,u}$ tensor on the left-hand side of \eqref{E:MAINTERMROTANGDROTTENSORPRODUCTSUMISSMALL},
and that the rectangular components of interest can be expressed as
\begin{align} \label{E:SECONDCOMPLICATEDPRODUCTREDUCEDFORM}
	\left\lbrace
		\delta_{ac} \delta_{bd} - \delta_{ad} \delta_{bc}
	\right\rbrace 
	x^a 
	(\ginversesphere)^{mb} 
	\sphereproject_p^c 
	(\ginversesphere)^{nd}
	= x^c \delta_{bd} (\ginversesphere)^{mb} \sphereproject_p^c (\ginversesphere)^{nd}
		- x^d \delta_{bc} (\ginversesphere)^{mb} \sphereproject_p^c (\ginversesphere)^{nd}.
\end{align}
To prove the desired estimate for the rectangular components, 
we first use \eqref{E:RADUNITLOWEREDBACKGROUNDSUBTRACTEDC0}
to replace $x^c$ with $- \rgeo \Radunit_c$
and $x^d$ with $- \rgeo \Radunit_d$
plus an error term with magnitude 
$\lesssim \varepsilon^{1/2} \ln(\myexp + t).$
Furthermore, \eqref{E:GSPHERERECTCOMPONENTSPOINTWISEBOUND} implies
that the remaining ``non-$x$'' factors on the right-hand side of \eqref{E:SECONDCOMPLICATEDPRODUCTREDUCEDFORM}
are in magnitude $\lesssim 1,$ so that the 
total product generated by the error term is
in magnitude $\lesssim \varepsilon^{1/2} \ln(\myexp + t)$ as desired. 
We now note that the $\rgeo \Radunit_c$ and $\rgeo \Radunit_d$ terms 
can be viewed as the components of one-forms that are $S_{t,u}-$orthogonal.
Since the $c$ and $d$ indices in \eqref{E:SECONDCOMPLICATEDPRODUCTREDUCEDFORM}
are paired with the rectangular components of $S_{t,u}$ tensors,
these terms lead to products that completely vanish.
We have thus proved \eqref{E:ROTANGDROTTENSORPRODUCTSUMISSMALL}.

\end{proof}

\section{Precise pointwise differential operator comparison estimates}
\label{S:POINTWISEDIFFERENTIALOPERATORCOMPARISONESTIMATES}
In this section, we derive some pointwise differential operator comparison estimates.
Most of the estimates are sharp in the sense that the explicit order $1$ constants are 
the same constants that our proofs would yield in the case
of Minkowski spacetime.

\begin{lemma}[\textbf{$\angD$ in terms of $\angLie_{\mathscr{O}}$}]
	\label{L:POINTWISEANGDINTERMSOFANGLIEO}
	Let $f$ be a function. 
	Under the small-data and bootstrap assumptions 
	of Sects.~\ref{S:PSISOLVES}-\ref{S:C0BOUNDBOOTSTRAP},
	if $\varepsilon^{1/2}$ is sufficiently small,
	then the following pointwise estimates hold on $\mathcal{M}_{\Tboot,U_0}:$
	\begin{subequations}
	\begin{align}
		\rgeo^2 |\angdiff f|^2
		& \leq (1 + C \varepsilon^{1/2}) 
			\sum_{l=1}^3 |\Rot_{(l)} f|^2,
			\label{E:FUNCTIONPOINTWISEANGDINTERMSOFANGLIEO} 
				\\
		\rgeo^2 |\angD^2 f|^2
			+ |\angdiff f|^2
		& \leq (1 + C \varepsilon^{1/2})
			\sum_{l=1}^3 |\angdiff \Rot_{(l)} f|^2,
			\label{E:ANGDSQUAREDFUNCTIONPOINTWISEINTERMSOFANGDIFFROTATIONS} \\
		\rgeo^2 |\angLap f|^2
		& \leq 2 (1 + C \varepsilon^{1/2})
			\sum_{l=1}^3 |\angdiff \Rot_{(l)} f|^2.
			\label{E:ANGLAPFUNCTIONPOINTWISEINTERMSOFROTATIONS}
	\end{align}
	\end{subequations}
	
	Similarly, if $\xi$ is an $S_{t,u}$ one-form, then
	\begin{align} \label{E:ONEFORMANGDINTERMSOFROTATIONALLIE}
		\rgeo^2 |\angD \xi|^2
		+ |\xi|^2
		& \leq (1 + C \varepsilon^{1/2})
			\sum_{l=1}^3 |\angLie_{\Rot_{(l)}} \xi|^2.
	\end{align}
	
	In addition, if $\xi$ is a symmetric type $\binom{0}{2}$ $S_{t,u}$ tensorfield, then
	\begin{align} \label{E:TYPE02TENSORANGDINTERMSOFROTATIONALLIE}
		\rgeo^2 |\angD \xi|^2 
		+ 4 |\hat{\xi}|^2
		& \leq (1 + C \varepsilon^{1/2})
			\sum_{l=1}^3 |\angLie_{\Rot_{(l)}} \xi|^2
			+ C \varepsilon^{1/2} (\mytr \xi)^2,
	\end{align}
	
	Finally, if $\xi$ is a symmetric type $\binom{0}{2}$ $S_{t,u}$ tensorfield, then
	\begin{align} \label{E:ANGLIEROTOFATENSORINTERMSOFANGD}
	|\angLie_{\Rot} \xi|
	& \leq
		C \rgeo |\angD \xi|
		+ C |\xi|.
	\end{align}
		
\end{lemma}

\begin{proof}
First, we note that \eqref{E:ANGLAPFUNCTIONPOINTWISEINTERMSOFROTATIONS} 
follows from \eqref{E:ANGDSQUAREDFUNCTIONPOINTWISEINTERMSOFANGDIFFROTATIONS}
and the fact that $|\mytr \xi|^2 \leq 2 |\xi|^2$ for symmetric type $\binom{0}{2}$ $S_{t,u}$ tensors.
Next, we note that \eqref{E:ANGDSQUAREDFUNCTIONPOINTWISEINTERMSOFANGDIFFROTATIONS}
follows from \eqref{E:ONEFORMANGDINTERMSOFROTATIONALLIE} 
and the identity $ \angLie_\Rot \angdiff f = \angdiff \angLie_\Rot f = \angdiff \Rot f,$
(that is, Lemma~\ref{L:LANDRADCOMMUTEWITHANGDIFF}).
Hence, to complete the proof of the lemma, we have to prove
\eqref{E:FUNCTIONPOINTWISEANGDINTERMSOFANGLIEO},
\eqref{E:ONEFORMANGDINTERMSOFROTATIONALLIE},
\eqref{E:TYPE02TENSORANGDINTERMSOFROTATIONALLIE},
and \eqref{E:ANGLIEROTOFATENSORINTERMSOFANGD}.

To prove these inequalities, 
we first use Lemma~\ref{L:ALTERNATELIEDERIVATIVEFORSTU} to deduce that
the following identity
holds for any type $\binom{0}{k}$ $S_{t,u}$
tensorfield $\xi$ with rectangular components $\xi_{j_1 \cdots j_k}:$
\begin{align} \label{E:LIEDERIVATIVESINTERMSOFSTUCOVARIANTDERIVATIVES}
\angLie_{\Rot_{(l)}} \xi_{j_1\dots j_k} 
& = \angDarg{\Rot_{(l)}} \xi_{j_1\dots j_k} 
+ \sum_{i = 1}^k \xi_{\cdots j_{i-1} m j_{i+1} \cdots j_k} \angDarg{j_i} \Rot_{(l)}^m.
\end{align}
 
We now note that the desired estimate \eqref{E:ANGLIEROTOFATENSORINTERMSOFANGD}
follows easily from
\eqref{E:LIEDERIVATIVESINTERMSOFSTUCOVARIANTDERIVATIVES},
\eqref{E:ROTATIONPOINTWISENORMESTIMATE},
and \eqref{E:POINTWISEBOUNDSPHERECOVARIANTDERIVATIVEOFROTATIONS}.

To prove the remaining estimates, we first square the norms of both sides of \eqref{E:LIEDERIVATIVESINTERMSOFSTUCOVARIANTDERIVATIVES}
(viewed as $S_{t,u}$ tensors)
and then sum over the index $l$ to deduce that
(where we suppress the $j_{\cdot}$ indices from \eqref{E:LIEDERIVATIVESINTERMSOFSTUCOVARIANTDERIVATIVES})
\begin{align} \label{E:SQUAREOFLIEDERIVATIVESINTERMSOFSTUCOVARIANTDERIVATIVES}
	\sum_{l = 1}^3 |\angLie_{\Rot_{(l)}} \xi|^2 
	& = \sum_{l = 1}^3 \left|\Rot_{(l)}^m \angDarg{m} \xi \right|^2
		+ 2 \sum_{l = 1}^3 
			\sum_{i=1}^k 
			(\Rot_{(l)}^m \angDarg{m} \xi)^{\#}
			\cdot (\angD \Rot_{(l)})^n \xi_{\cdots n\cdots} 
		+ \sum_{l = 1}^3 
			\left|
			\sum_{i=1}^k 
				(\angD \Rot_{(l)})^m \xi_{\cdots m\cdots}
			\right|^2.
\end{align}
From \eqref{E:GINVERSESPHEREISWELLAPPROXIMATEDBYROTATIONTENSORPRODUCTS}, 
we deduce that the first term on the right-hand side of 
\eqref{E:SQUAREOFLIEDERIVATIVESINTERMSOFSTUCOVARIANTDERIVATIVES}
is equal to $\rgeo^2 |\angD \xi|^2$ plus an error term that is in magnitude
$\lesssim  \varepsilon^{1/2} \rgeo^2 |\angD \xi|^2.$
Furthermore, using \eqref{E:ROTANGDROTTENSORPRODUCTSUMISSMALL}
and Cor.~\ref{C:NORMOFTHESTUTENSORISCOMPARABLETOTHESUMOFITSCOMPONENTS},
we deduce that the second term is in magnitude
$\lesssim \varepsilon^{1/2} \rgeo^2 |\angD \xi|^2 + \varepsilon^{1/2} |\xi|^2.$

We split the analysis of the last sum on the right-hand side of \eqref{E:SQUAREOFLIEDERIVATIVESINTERMSOFSTUCOVARIANTDERIVATIVES}
into the cases where $\xi$ is a function, type $\binom{0}{1},$ and symmetric type $\binom{0}{2}$
respectively. 
If $\xi$ is a function $f,$ 
then the left-hand side of \eqref{E:LIEDERIVATIVESINTERMSOFSTUCOVARIANTDERIVATIVES}
is equal to $\sum_{l=1}^3 |\Rot_{(l)} f|^2$ while 
the estimates in the previous paragraph imply that the
right-hand side is equal to $\rgeo^2 |\angdiff f|^2$
plus an error term that is in magnitude $\lesssim  \varepsilon^{1/2} \rgeo^2 |\angdiff f|^2$
(the second and third sums on the right-hand side of \eqref{E:SQUAREOFLIEDERIVATIVESINTERMSOFSTUCOVARIANTDERIVATIVES} are absent).
The desired estimate \eqref{E:FUNCTIONPOINTWISEANGDINTERMSOFANGLIEO} easily follows from these inequalities.

If $\xi$ is type $\binom{0}{1},$
then the last sum 
on the right-hand side of \eqref{E:SQUAREOFLIEDERIVATIVESINTERMSOFSTUCOVARIANTDERIVATIVES}
is equal to 
\begin{align}
	\sum_{l=1}^3	
		(\ginversesphere)^{jj'}(\angDarg{j} \Rot^m_{(l)})(\angDarg{j'} \Rot^n_{(l)}) \xi_m \xi_n.
\end{align}
Hence, using \eqref{E:GINVERSESPHEREISWELLAPPROXIMATEDBYROTATIONCOVARIANTDERIVATIVETENSORPRODUCTS}
and Cor.~\ref{C:NORMOFTHESTUTENSORISCOMPARABLETOTHESUMOFITSCOMPONENTS},
we deduce that the last sum is equal to $|\xi|^2$
plus an error term that is in magnitude $\lesssim \varepsilon^{1/2} |\xi|^2.$
Setting $\xi = \angdiff f,$ we see that the 
desired estimate \eqref{E:ANGDSQUAREDFUNCTIONPOINTWISEINTERMSOFANGDIFFROTATIONS}
now follows easily from this estimate and our prior analysis of the first two sums
on the right-hand side of \eqref{E:SQUAREOFLIEDERIVATIVESINTERMSOFSTUCOVARIANTDERIVATIVES}.

If $\xi$ is symmetric type $\binom{0}{2},$ 
then the last sum on the right-hand side of \eqref{E:SQUAREOFLIEDERIVATIVESINTERMSOFSTUCOVARIANTDERIVATIVES}
is equal to
\begin{align} \label{E:RANKTWOCASETHIRDTERM}
& 
2 \sum_{l=1}^3
			(\ginversesphere)^{jj'}(\ginversesphere)^{kk'} 
			\xi_{mj}
			\xi_{nj'}
			(\angDarg{k}\Rot_{(l)}^m)
			\angDarg{k'}\Rot_{(l)}^n
	+ 2 \sum_{l=1}^3
			(\ginversesphere)^{jj'}(\ginversesphere)^{kk'}
			\xi_{mj}
			\xi_{nk'}
			(\angDarg{k}\Rot_{(l)}^m)
			\angDarg{j'}\Rot_{(l)}^n.
\end{align}
Inequality 
\eqref{E:GINVERSESPHEREISWELLAPPROXIMATEDBYROTATIONCOVARIANTDERIVATIVETENSORPRODUCTS}
and Cor.~\ref{C:NORMOFTHESTUTENSORISCOMPARABLETOTHESUMOFITSCOMPONENTS}
imply that the first sum
\eqref{E:RANKTWOCASETHIRDTERM} is 
equal to $2 |\xi|^2$ plus an error term that is in magnitude $\lesssim \varepsilon^{1/2} |\xi|^2.$
Furthermore, using \eqref{E:GINVERSEPROJECTIONCOMBINATIONISWELLAPPROXIMATEDBYROTATIONCOVARIANTDERIVATIVETENSORPRODUCTS}
and Cor.~\ref{C:NORMOFTHESTUTENSORISCOMPARABLETOTHESUMOFITSCOMPONENTS},
we deduce that
the second sum in \eqref{E:RANKTWOCASETHIRDTERM} is
equal to $2 |\xi|^2 - 2 (\mytr \xi)^2$
plus an error term that is in magnitude $\lesssim \varepsilon^{1/2} |\xi|^2.$
The desired estimate \eqref{E:TYPE02TENSORANGDINTERMSOFROTATIONALLIE}
now follows easily from these two estimates for the two sums in \eqref{E:RANKTWOCASETHIRDTERM}, 
our prior analysis of the first two sums
on the right-hand side of \eqref{E:SQUAREOFLIEDERIVATIVESINTERMSOFSTUCOVARIANTDERIVATIVES},
and the identity $|\xi|^2 = \frac{1}{2} (\mytr \xi)^2 + |\hat{\xi}|^2.$
This completes the proof of Lemma~\ref{L:POINTWISEANGDINTERMSOFANGLIEO}.
\end{proof}

\section{Useful estimates for avoiding detailed commutators}
In this section, we provide some non-optimal commutator estimates
for vectorfields acting on functions.
The estimates are useful when precision is not required.

\begin{lemma}[\textbf{Estimates for avoiding detailed commutators}]
\label{L:AVOIDINGCOMMUTING}
Let $f$ be a function, and let $1 \leq N \leq 24$ be an integer. 
Under the small-data and bootstrap assumptions 
of Sects.~\ref{S:PSISOLVES}-\ref{S:C0BOUNDBOOTSTRAP},
if $\varepsilon$ is sufficiently small,
then the following pointwise estimates hold on $\mathcal{M}_{\Tboot,U_0}:$
\begin{align} \label{E:FUNCTIONDERIVATIVESAVOIDINGCOMMUTING}
	\left| 
		\threemyarray
			[\rgeo \mathscr{Z}^N \Lunit f]
			{\mathscr{Z}^N \Rad f}
			{\rgeo \angLie_{\mathscr{Z}}^N \angdiff f}
	\right|
	& \lesssim
		\left| 
			\threemyarray
				[\rgeo \Lunit \mathscr{Z}^{\leq N} f]
				{\Rad \mathscr{Z}^{\leq N} f}
				{\rgeo \angdiff \mathscr{Z}^{\leq N} f}
		\right|,
			\\
	\left| 
		\mathscr{Z}^N f
	\right|
	& \lesssim
		\left| 
			\fourmyarray
				[\rgeo \Lunit \mathscr{Z}^{\leq N-1} f]
				{\Rad \mathscr{Z}^{\leq N-1} f}
				{\rgeo \angdiff \mathscr{Z}^{\leq N-1} f}
				{\mathscr{Z}^{N-1} f}
			\right|,
		\label{E:FUNCTIONAVOIDINGCOMMUTING}
		\\
	\left| 
		\fourmyarray
			[\rgeo \Lunit \mathscr{Z}^{\leq N-1} f]
			{\Rad \mathscr{Z}^{\leq N-1} f}
			{\rgeo \angdiff \mathscr{Z}^{\leq N-1} f}
			{\mathscr{Z}^{N-1} f}
	\right|
	& \lesssim 
		\left|
			\mathscr{Z}^{\leq N} f
		\right|.
		\label{E:BOUNDINGFRAMEDERIVATIVESINTERMSOFZNF}
\end{align}
\end{lemma}
	
\begin{proof}
	We first prove \eqref{E:FUNCTIONDERIVATIVESAVOIDINGCOMMUTING}.
	The estimate for $\rgeo \angLie_{\mathscr{Z}}^N \angdiff f$
	follows trivially from Lemma~\ref{L:LANDRADCOMMUTEWITHANGDIFF}. 
	
	To deduce the estimate for $\rgeo \mathscr{Z}^N \Lunit f,$
	we express $\Lunit f = \rgeo^{-1} Z f,$ where $Z := \rgeo \Lunit \in \mathscr{Z}.$
	Then using \eqref{E:ZNAPPLIEDTORGEOISNOTTOOLARGE}, we deduce that
	$\left|\rgeo \mathscr{Z}^N \Lunit f \right|
	\lesssim \left|\mathscr{Z}^{\leq N+1} f \right|.$
	If $\mathscr{Z}^{N+1} = \rgeo \Lunit \mathscr{Z}^N,$
	or $\mathscr{Z}^{N+1} = \Rad \mathscr{Z}^N,$
	then the desired bound is obvious.
	If $\mathscr{Z}^{N+1} = \Rot \mathscr{Z}^N,$
	then we use \eqref{E:ROTATIONPOINTWISENORMESTIMATE} to deduce that
	$|\Rot \mathscr{Z}^N f| 
	\lesssim |\Rot| \left|\angdiff \mathscr{Z}^N f \right|
	\lesssim \rgeo \left|\angdiff \mathscr{Z}^N f \right|.$
	
	To deduce the estimate for $\rgeo \mathscr{Z}^N \Rad f,$
	we first note that if $\mathscr{Z}^{N+1} = \rgeo \Lunit \mathscr{Z}^N$
	or $\mathscr{Z}^{N+1} = \Rad \mathscr{Z}^N,$
	then the desired bound is obvious. If $\mathscr{Z}^{N+1} = \Rot \mathscr{Z}^N,$
	then we argue as we did at the end of the previous paragraph.
	We have thus proved \eqref{E:FUNCTIONDERIVATIVESAVOIDINGCOMMUTING}.
	
	Inequalities \eqref{E:FUNCTIONAVOIDINGCOMMUTING} and \eqref{E:BOUNDINGFRAMEDERIVATIVESINTERMSOFZNF}
	can be proved using similar arguments, and we omit the details.
\end{proof}

\section{Estimates for \texorpdfstring{$\angdiff x^i$}{the derivatives of the angular differential of the rectangular coordinate functions}}
In this section, we derive pointwise estimates
for the $S_{t,u}-$projected Lie derivatives of $\angdiff x^i,$ $(i=1,2,3).$

\begin{lemma}
	\label{L:POINTWISEBOUNDPROJECTEDLIEDERIVATIVESANGDIFFCOORDINATEX}
	Let $0 \leq N \leq 24$ be an integer.
	Under the small-data and bootstrap assumptions 
	of Sects.~\ref{S:PSISOLVES}-\ref{S:C0BOUNDBOOTSTRAP},
	if $\varepsilon$ is sufficiently small,
	then the following pointwise estimates hold on $\mathcal{M}_{\Tboot,U_0}:$
\begin{subequations}
\begin{align} \label{E:POINTWISEBOUNDPROJECTEDLIEDERIVATIVESANGDIFFCOORDINATEX}
	\left| 
		\angLie_{\mathscr{Z}}^N \angdiff x^i
	\right|,
	\,
	\left| 
		\angdiff \mathscr{Z}^N x^i
	\right|
	& \lesssim
		|\mathscr{Z}^{\leq N} \Psi|
		+ \frac{1}{1 + t}
			\left|
				\myarray[\mathscr{Z}^{\leq N} (\upmu - 1)]
					{\sum_{a=1}^3 \rgeo |\mathscr{Z}^{\leq N} \Lunit_{(Small)}^a|} 
			\right|
			+ 1,
			\\
	\left\| 
		\myarray
			[\angLie_{\mathscr{Z}}^{\leq 12} \angdiff x^i]
			{\angdiff \mathscr{Z}^{\leq 12}  x^i}
	\right\|_{C^0(\Sigma_t^u)}		
	& \lesssim 1.
		\label{E:LOWERORDERPOINTWISEBOUNDPROJECTEDLIEDERIVATIVESANGDIFFCOORDINATEX}
\end{align}
\end{subequations}
\end{lemma}

\begin{proof}
	To prove \eqref{E:POINTWISEBOUNDPROJECTEDLIEDERIVATIVESANGDIFFCOORDINATEX} 
	and \eqref{E:LOWERORDERPOINTWISEBOUNDPROJECTEDLIEDERIVATIVESANGDIFFCOORDINATEX},
	we use 
	Lemma~\ref{L:LANDRADCOMMUTEWITHANGDIFF} and
	\eqref{E:FUNCTIONPOINTWISEANGDINTERMSOFANGLIEO}
	to deduce that 
	$|\angLie_{\mathscr{Z}}^N \angdiff x^i| \lesssim (1 + t)^{-1} |\mathscr{Z}^{N+1} x^i| 
		= (1 + t)^{-1} |\mathscr{Z}^N Z^i|.$ The desired estimates thus follow from
	\eqref{E:DERIVATIVESOFRECTANGULARFRAMECOMPONENTS}
	and \eqref{E:LOWERORDERC0BOUNDDERIVATIVESOFRECTANGULARCOMMUTATORVECTORFIELDRECTCOMPONENTS}.
\end{proof}

\section{Pointwise estimates for the Lie derivatives of \texorpdfstring{$G_{(Frame)}$ and $G_{(Frame)}'$}{the frame components of the derivative of the rectangular components of the metric with respect to the solution}}
In this section, we derive pointwise estimates for the $S_{t,u}-$projected Lie derivatives of 
the arrays of $S_{t,u}$ tensorfields $G_{(Frame)}$ and $G_{(Frame)}'.$

\begin{lemma}[\textbf{Pointwise estimates for the Lie derivatives of} $G_{(Frame)}$ \textbf{and} $G_{(Frame)}'$] 
\label{L:POINTWISEESTIMATESGFRAMEINTERMSOFOTHERQUANTITIES}
 	Let $G_{(Frame)}$ and $G_{(Frame)}'$ be the arrays of $S_{t,u}$ tensorfields
 	given by Defs.~\ref{D:BIGGANDBIGGPRIME} and \ref{D:GFRAMECOMPONENTS}.
	Let $0 \leq N \leq 24$ be an integer. 
	Under the small-data and bootstrap assumptions 
	of Sects.~\ref{S:PSISOLVES}-\ref{S:C0BOUNDBOOTSTRAP},
	the following pointwise estimates hold on $\mathcal{M}_{\Tboot,U_0}:$
	\begin{subequations}
	\begin{align}
		\left| 
			\angLie_{\mathscr{Z}}^{\leq N} G_{(Frame)}
		\right|
		& \lesssim |\mathscr{Z}^{\leq N} \Psi|
			+ 
			\frac{1}{1 + t}
			\left| 
				\myarray
					[\mathscr{Z}^{\leq N} (\upmu - 1)]
						{\rgeo \sum_{a=1}^3 |\mathscr{Z}^{\leq N} \Lunit_{(Small)}^a|}
			\right|
			+ 1,
				\label{E:LIEDERIVATIVESOFGRAMEINTERMSOFOTHERVARIABLES} \\
		\left\| 
			\angLie_{\mathscr{Z}}^{\leq 12} G_{(Frame)}
		\right\|_{C^0(\Sigma_t^u)}
		& \lesssim 1.
		\label{E:LOWERORDERC0BOUNDLIEDERIVATIVESOFGRAME}
	\end{align}
	\end{subequations}	
	Furthermore, the same estimates hold for 
	$G_{(Frame)}'.$
	
\end{lemma}

\begin{proof}
We prove only the desired estimate for the $S_{t,u}$ one-form $\angGarg{\Lunit};$
the proofs of the estimates for the remaining elements of $G_{(Frame)}$
and for $G_{(Frame)}'$ are essentially the same.
We use the notation of Lemma~\ref{L:PROJECTEDLIEDERIVATIVESINTERMSOFOTHERVARIABLES}.
In the following argument, for convenience, we denote by $G'$
any derivative (including higher-order ones)
of the rectangular component functions $G_{\mu \nu}(\Psi)$ with respect to $\Psi.$
Since the $G'$ are smooth functions of $\Psi$, the bootstrap assumptions \eqref{E:PSIFUNDAMENTALC0BOUNDBOOTSTRAP}
imply that $|G'| \lesssim 1.$ Using in addition the estimates \eqref{E:DERIVATIVESOFRECTANGULARFRAMECOMPONENTS}
and the estimates \eqref{E:GSPHERERECTCOMPONENTSPOINTWISEBOUND} for the rectangular components of 
$\ginversesphere$ and the $S_{t,u}$ projections $\sphereproject,$ 
we conclude that the scalar functions $G_{Va}$ and the $S_{t,u}$ one-forms $\angGarg{V}$
are $\lesssim 1$ in magnitude, and similarly with $G'$ in place of $G.$
Furthermore, Lemma~\ref{L:LANDRADCOMMUTEWITHANGDIFF},
inequality \eqref{E:FUNCTIONPOINTWISEANGDINTERMSOFANGLIEO},
and the estimates \eqref{E:POINTWISEBOUNDLOWERORDERDERIVATIVESOFRECTANGULARFRAMECOMPONENTS}
and \eqref{E:LOWERORDERC0BOUNDDERIVATIVESOFRECTANGULARCOMMUTATORVECTORFIELDRECTCOMPONENTS}
imply that 
$\left\| \mathscr{Z}^{\leq 12} \Lunit^j \right\|_{C^0(\Sigma_t^u)},$
$\left\| \mathscr{Z}^{\leq 12} \Radunit^j \right\|_{C^0(\Sigma_t^u)} \lesssim 1$
and that for any $Z \in \mathscr{Z},$ 
$\left\| \angLie_{\mathscr{Z}}^{\leq 11} \angdiff Z^j \right\|_{C^0(\Sigma_t^u)},$
$\left\| \angLie_{\mathscr{Z}}^{\leq 11} \angdiff Z^j \right\|_{C^0(\Sigma_t^u)}
\lesssim 1.$ Using all of these estimates and the bootstrap assumptions \eqref{E:PSIFUNDAMENTALC0BOUNDBOOTSTRAP},
we now repeatedly $\angLie-$differentiate the identity \eqref{E:GFRAMEONEFORMDERIVATIVE}
with respect to vectorfields $Z \in \mathscr{Z}$ and use  
Lemma~\ref{L:LANDRADCOMMUTEWITHANGDIFF} and inequality \eqref{E:FUNCTIONPOINTWISEANGDINTERMSOFANGLIEO}  
to deduce that 
\begin{align} \label{E:FIRSTESTIMATELIEDERIVATIVESOFGRAMEINTERMSOFOTHERVARIABLES}
	\left| \angLie_{\mathscr{Z}}^N \angGarg{V} \right|
	& \lesssim \left|\mathscr{Z}^{\leq N} \Psi \right|
		+ \sum_{a=1}^3 |\mathscr{Z}^{\leq N} \Lunit^a|
		+ \sum_{Z \in \mathscr{Z}} \sum_{a=1}^3 |\angdiff \mathscr{Z}^{\leq N-1} Z^a|
		\\
	& \lesssim 
		\left|\mathscr{Z}^{\leq N} \Psi \right|
		+ \sum_{a=1}^3 |\mathscr{Z}^{\leq N} \Lunit^a|
		+ \frac{1}{\rgeo} \sum_{Z \in \mathscr{Z}} \sum_{a=1}^3 |\mathscr{Z}^{\leq N} Z^a|.
		\notag
\end{align}
The desired estimate
\eqref{E:LIEDERIVATIVESOFGRAMEINTERMSOFOTHERVARIABLES}
now follows from \eqref{E:FIRSTESTIMATELIEDERIVATIVESOFGRAMEINTERMSOFOTHERVARIABLES},
\eqref{E:DERIVATIVESOFRECTANGULARFRAMECOMPONENTS},
and \eqref{E:DERIVATIVESOFRECTANGULARCOMMUTATORVECTORFIELDRECTCOMPONENTS}.
Inequality \eqref{E:LOWERORDERC0BOUNDLIEDERIVATIVESOFGRAME} then follows from
\eqref{E:LIEDERIVATIVESOFGRAMEINTERMSOFOTHERVARIABLES} and the bootstrap assumptions.

\end{proof}

\section{Crude pointwise estimates for 
\texorpdfstring{$\angLie_{\mathscr{Z}}^N \angdeform{V}$}
{the Lie derivatives of the angular components of the deformation tensors}}
In this section, 
for various vectorfields $V,$
we provide some crude pointwise estimates for the $S_{t,u}-$projected Lie derivatives of the
$S_{t,u}$ tensorfields $\angdeform{V}$ and 
the Lie derivatives of some related tensorfields. 

\begin{lemma}[\textbf{Crude pointwise estimates for the Lie derivatives of the angular components of the deformation tensors}]
\label{L:POINTWISEBOUNDSDERIVATIVESOFANGULARDEFORMATIONTENSORS}
	Let $0 \leq N \leq 23$ be an integer, and let $Z \in \mathscr{Z}$ be a commutation vectorfield.
	Under the small-data and bootstrap assumptions 
	of Sects.~\ref{S:PSISOLVES}-\ref{S:C0BOUNDBOOTSTRAP},
	if $\varepsilon$ is sufficiently small,
	then the following pointwise estimates hold on $\mathcal{M}_{\Tboot,U_0}:$
	\begin{subequations}
	\begin{align} \label{E:CRUDEPOINTWISEBOUNDSDERIVATIVESOFANGULARDEFORMATIONTENSORS}
		\left|
		\ninemyarray
			[\rgeo \angLie_{\mathscr{Z}}^N \upchi]
			{\rgeo \angLie_{\mathscr{Z}}^N \upchi^{\#}}
			{\rgeo \mathscr{Z}^N \mytr \upchi}
			{\angLie_{\mathscr{Z}}^N \angdeform{Z}}
			{\angLie_{\mathscr{Z}}^N \angdeformupsharparg{Z}}
			{\angLie_{\mathscr{Z}}^N \angdeformupdoublesharparg{Z}}
			{\mathscr{Z}^N \mytr \angdeform{Z}}
			{\angLie_{\mathscr{Z}}^{N+1} \gsphere}
			{\angLie_{\mathscr{Z}}^{N+1} \ginversesphere}
		\right|
		& \lesssim 
			\left|
				\fourmyarray[\rgeo \Lunit \mathscr{Z}^{\leq N} \Psi]
					{\Rad \mathscr{Z}^{\leq N} \Psi}
					{\rgeo \angdiff \mathscr{Z}^{\leq N} \Psi} 
					{\mathscr{Z}^{\leq N} \Psi}
			\right|
		+ \frac{1}{1 + t}
			\left| 
				\myarray
					[\mathscr{Z}^{\leq N+1} (\upmu - 1)]
						{\rgeo \sum_{a=1}^3 |\mathscr{Z}^{\leq N+1} \Lunit_{(Small)}^a|}
				\right|
			+ 1,	
			\\
		\left\|
		\ninemyarray
			[\rgeo \angLie_{\mathscr{Z}}^{\leq 11} \upchi]
			{\rgeo \angLie_{\mathscr{Z}}^{\leq 11} \upchi^{\#}}
			{\rgeo \mathscr{Z}^{\leq 11} \mytr \upchi}
			{\angLie_{\mathscr{Z}}^{\leq 11} \angdeform{Z}}
			{\angLie_{\mathscr{Z}}^{\leq 11} \angdeformupsharparg{Z}}
			{\angLie_{\mathscr{Z}}^{\leq 11} \angdeformupdoublesharparg{Z}}
			{\mathscr{Z}^{\leq 11} \mytr  \angdeform{Z}}
			{\angLie_{\mathscr{Z}}^{\leq 12} \gsphere}
			{\angLie_{\mathscr{Z}}^{\leq 12} \ginversesphere}
		\right\|_{C^0(\Sigma_t^u)}
		& \lesssim 1.
		\label{E:CRUDELOWERORDERC0BOUNDDERIVATIVESOFANGULARDEFORMATIONTENSORS}
\end{align}
\end{subequations}

\end{lemma}

\begin{proof}
	Throughout this proof, we use the 	
	bootstrap assumptions
	\eqref{E:PSIFUNDAMENTALC0BOUNDBOOTSTRAP},
	\eqref{E:UPMUBOOT},
	and
	\eqref{E:FRAMECOMPONENTSBOOT}.
We first prove \eqref{E:CRUDEPOINTWISEBOUNDSDERIVATIVESOFANGULARDEFORMATIONTENSORS}	
for $\gsphere.$ We differentiate both sides of the identity
$\gsphere_{AB} = g_{ab} \angdiffarg{A} x^a \angdiffarg{B} x^b$
with the operator $\angLie_{\mathscr{Z}}^N$ and apply the Leibniz
rule to the right-hand side.
We treat all uppercase Latin indices as tensorial $S_{t,u}$ indices,
and all lowercase Latin-indexed quantities as scalar-valued functions.
We use \eqref{E:ANGLIECOMMUTESWITHANGDIFF} to commute the operator
$\angLie_{\mathscr{Z}}^N$ under $\angdiff.$ It follows that
\begin{align} \label{E:GSPHERELIEDIFFERNTIATEDLEIBNIZ}
	\left|
		\angLie_{\mathscr{Z}}^{N+1} \gsphere
	\right| 
	& \lesssim 
	\sum_{N_1 + N_2 + N_3 = N+1}
		\left|
			\mathscr{Z}^{N_1} g_{ab} 
		\right|
		\left|
			\angdiff \mathscr{Z}^{N_2} x^a
		\right|
		\left|
			\angdiff \mathscr{Z}^{N_3} x^b
		\right|.
\end{align}	
The desired estimate \eqref{E:CRUDEPOINTWISEBOUNDSDERIVATIVESOFANGULARDEFORMATIONTENSORS} for
$\angLie_{\mathscr{Z}}^{N+1} \gsphere$ now follows from 
\eqref{E:GSPHERELIEDIFFERNTIATEDLEIBNIZ}
and the estimates
\eqref{E:GRECTCOMPONENTINTERMSOFOTHERVARIABLES},
\eqref{E:GRECTCOMPONENTC0BOUND},
\eqref{E:POINTWISEBOUNDPROJECTEDLIEDERIVATIVESANGDIFFCOORDINATEX},
\eqref{E:LOWERORDERPOINTWISEBOUNDPROJECTEDLIEDERIVATIVESANGDIFFCOORDINATEX},
and 
\eqref{E:FUNCTIONAVOIDINGCOMMUTING}.
The desired estimate \eqref{E:CRUDELOWERORDERC0BOUNDDERIVATIVESOFANGULARDEFORMATIONTENSORS}
for $\angLie_{\mathscr{Z}}^{N+1} \gsphere$ then follows from
the estimate \eqref{E:CRUDEPOINTWISEBOUNDSDERIVATIVESOFANGULARDEFORMATIONTENSORS} for
$\angLie_{\mathscr{Z}}^{N+1} \gsphere$ and the bootstrap assumptions.

To prove the estimate \eqref{E:CRUDEPOINTWISEBOUNDSDERIVATIVESOFANGULARDEFORMATIONTENSORS}	
for $\angdeform{Z},$ we simply note that by Lemma~\ref{L:CONNECTIONBETWEENANGLIEOFGSPHEREANDDEFORMATIONTENSORS},
we have $\angdeform{Z} = \angLie_Z \gsphere.$
Hence the desired estimate for $\angLie_{\mathscr{Z}}^N \angdeform{Z}$
follows from the previously proven estimate for $\angLie_{\mathscr{Z}}^{N+1} \gsphere.$
The desired estimate \eqref{E:CRUDELOWERORDERC0BOUNDDERIVATIVESOFANGULARDEFORMATIONTENSORS}
for $\angdeform{Z}$ then follows from
the estimate \eqref{E:CRUDEPOINTWISEBOUNDSDERIVATIVESOFANGULARDEFORMATIONTENSORS} for
$\angdeform{Z}$ and the bootstrap assumptions.

To prove the estimate \eqref{E:CRUDEPOINTWISEBOUNDSDERIVATIVESOFANGULARDEFORMATIONTENSORS}
for $\angLie_{\mathscr{Z}}^N \upchi,$ we set $Z = \rgeo \Lunit$ and note that 
by Lemma~\ref{L:CONNECTIONBETWEENANGLIEOFGSPHEREANDDEFORMATIONTENSORS},
\eqref{E:RGEOLDEFORMTRFREESPHERE}, and \eqref{E:RGEOLDEFORMTRSPHERE}, we have
$2 \rgeo \upchi = \angLie_Z \gsphere.$ Using Lemma~\ref{L:BASICESTIMATESFORRGEO}, we see that
the desired estimate follows from the previously proven estimate for $\angLie_{\mathscr{Z}}^{N+1} \gsphere.$
The desired estimate \eqref{E:CRUDELOWERORDERC0BOUNDDERIVATIVESOFANGULARDEFORMATIONTENSORS}
for $\angLie_{\mathscr{Z}}^N \upchi$ then follows from
the estimate \eqref{E:CRUDEPOINTWISEBOUNDSDERIVATIVESOFANGULARDEFORMATIONTENSORS} for
$\angLie_{\mathscr{Z}}^N \upchi$ and the bootstrap assumptions.

To prove the estimate \eqref{E:CRUDEPOINTWISEBOUNDSDERIVATIVESOFANGULARDEFORMATIONTENSORS}	
for $\ginversesphere,$ we start with the identity
$\angLie_{Z} \ginversesphere = - \angdeformupdoublesharparg{Z}$
(see Lemma~\ref{L:CONNECTIONBETWEENANGLIEOFGSPHEREANDDEFORMATIONTENSORS}).
Applying $\angLie_{\mathscr{Z}}^N$ to this identity and arguing inductively, we see
that the resulting expression involves only $S_{t,u}-$tensorial products of $\ginversesphere$ and 
the Lie derivatives of the $\angpi.$ From the previously proven estimate
\eqref{E:CRUDELOWERORDERC0BOUNDDERIVATIVESOFANGULARDEFORMATIONTENSORS} for the lower-order derivatives of the
$\angpi,$ we infer that all terms that are quadratic in the derivatives of the $\angpi$ can be pointwise bounded by
a constant times a term that is linear in the derivatives of the $\angpi.$ It thus follows that
\begin{align}
	\left| \angLie_{\mathscr{Z}}^{N+1} \ginversesphere \right|
	& \lesssim \max_{Z \in \mathscr{Z}} \left| \angLie_{\mathscr{Z}}^{\leq N} \angdeform{Z} \right|,
\end{align}
and hence the desired estimate \eqref{E:CRUDEPOINTWISEBOUNDSDERIVATIVESOFANGULARDEFORMATIONTENSORS}
for $\angLie_{\mathscr{Z}}^{N+1} \ginversesphere$ follows from the 
previously proven estimate \eqref{E:CRUDEPOINTWISEBOUNDSDERIVATIVESOFANGULARDEFORMATIONTENSORS}	
for $\angdeform{Z}.$ 
The desired estimate \eqref{E:CRUDELOWERORDERC0BOUNDDERIVATIVESOFANGULARDEFORMATIONTENSORS}
for $\angLie_{\mathscr{Z}}^{N+1} \ginversesphere$ then follows from
the estimate \eqref{E:CRUDEPOINTWISEBOUNDSDERIVATIVESOFANGULARDEFORMATIONTENSORS} for
$\angLie_{\mathscr{Z}}^{N+1} \ginversesphere$ and the bootstrap assumptions.
Finally, since $\upchi^{\#} = \ginversesphere \upchi,$
$\angdeformupsharparg{Z} = \ginversesphere \angdeform{Z},$
and $\angdeformupdoublesharparg{Z} = (\ginversesphere)^2 \angdeform{Z},$
the estimates 
\eqref{E:CRUDEPOINTWISEBOUNDSDERIVATIVESOFANGULARDEFORMATIONTENSORS}
and
\eqref{E:CRUDELOWERORDERC0BOUNDDERIVATIVESOFANGULARDEFORMATIONTENSORS}
for $\angLie_{\mathscr{Z}}^N \upchi^{\#}$ 
and 
$\angLie_{\mathscr{Z}}^N \angdeformupsharparg{Z}$
follow easily from the Leibniz rule, the
previously proven estimates
for $\ginversesphere,$ 
$\upchi,$ and 
$\angdeform{Z}$
and the bootstrap assumptions.
\end{proof}

\section{Two additional crude differential operator comparison estimates}
In this section, we provide two additional differential operator comparison estimates.
We do not bother to derive sharp constants in the estimates.

\begin{lemma}[\textbf{A crude differential operator comparison estimate involving $\angLie_X$ in terms of $\angLie_{\mathscr{O}}$}]
	\label{L:POINTWISEANGLIEXINTERMSOFANGLIEO}
	Let $X$ be an $S_{t,u}$ vectorfield and let
	$\xi$ be an $S_{t,u}$ one-form or
	a symmetric type $\binom{0}{2}$ $S_{t,u}$ tensorfield.
	Under the small-data and bootstrap assumptions 
	of Sects.~\ref{S:PSISOLVES}-\ref{S:C0BOUNDBOOTSTRAP},
	if $\varepsilon$ is sufficiently small,	
	then the following pointwise estimate holds on $\mathcal{M}_{\Tboot,U_0}:$
	\begin{align} \label{E:POINTWISEANGLIEXINTERMSOFANGLIEO}
		|\angLie_X \xi| 
		& \leq C \frac{1}{1 + t}
			\left(
				|X| |\angLie_{\mathscr{O}}^{\leq 1} \xi|
				+ |\xi| |\angLie_{\mathscr{O}}^{\leq 1} X| 
			\right).
	\end{align}
	
\end{lemma}

\begin{proof}
	Let $X_{\flat}$ be the $S_{t,u}$ one-form that is $\gsphere-$dual to $X.$
	Using the identity \eqref{E:LIEDERIVATIVESINTERMSOFSTUCOVARIANTDERIVATIVES}
	with $X$ in place of $\Rot_{(l)}$ 
	to express $S_{t,u}-$projected Lie derivatives in terms of $S_{t,u}-$covariant derivatives,
	we deduce the schematic identity 
	$\angLie_X \xi = X \angD \xi + \xi^{\#} \angD X_{\flat}.$
	From this identity, \eqref{E:ONEFORMANGDINTERMSOFROTATIONALLIE}, and
	\eqref{E:TYPE02TENSORANGDINTERMSOFROTATIONALLIE}, we deduce that
	$|\angLie_X \xi| \lesssim 
		(1+t)^{-1}
			\left(
				|X| |\angLie_{\mathscr{O}}^{\leq 1} \xi|
				+ |\xi| |\angLie_{\mathscr{O}}^1 X_{\flat}| 
			\right).$
	To bound the Lie derivatives of $X_{\flat},$ we first note that
	$|\angLie_{\Rot} X_{\flat}| = |\angLie_{\Rot}(\gsphere X)| \lesssim |\angLie_{\Rot} X| 
	+ |\angdeform{\Rot}||X|.$
	Using this inequality and
	the estimate \eqref{E:CRUDELOWERORDERC0BOUNDDERIVATIVESOFANGULARDEFORMATIONTENSORS}
	(to bound the magnitude of $\angdeform{\Rot}$),
	we deduce that $|\angLie_{\Rot} X_{\flat}| \lesssim |\angLie_{\mathscr{O}}^{\leq 1} X|,$
	and the desired estimate \eqref{E:POINTWISEANGLIEXINTERMSOFANGLIEO} readily follows.
	
\end{proof}

\begin{lemma}[\textbf{A crude differential operator comparison estimate involving second spherical covariant derivatives}]
	\label{L:TYPE02TENSORANGDTWOTIMESINTERMSOFROTATIONALLIETWOTIME}
	Let $\xi$ be a symmetric type $\binom{0}{2}$ $S_{t,u}$ tensorfield.
 	Under the small-data and bootstrap assumptions 
	of Sects.~\ref{S:PSISOLVES}-\ref{S:C0BOUNDBOOTSTRAP},
	if $\varepsilon$ is sufficiently small,
 	then the following pointwise estimate holds on $\mathcal{M}_{\Tboot,U_0}:$
	\begin{align} \label{E:TYPE02TENSORANGDTWOTIMESINTERMSOFROTATIONALLIETWOTIME}
		 \rgeo^4 |\angD^2 \xi|^2 
		& \leq C |\angLie_{\mathscr{O}}^{\leq 2} \xi|^2.
	\end{align}
	
\end{lemma}

\begin{proof}
We first use the identity \eqref{E:SQUAREOFLIEDERIVATIVESINTERMSOFSTUCOVARIANTDERIVATIVES}
with $\angD \xi$ in the role of $\xi,$
the line of reasoning just below it, 
and inequality \eqref{E:POINTWISEBOUNDSPHERECOVARIANTDERIVATIVEOFROTATIONS}
(to bound the final product on the right-hand side of \eqref{E:SQUAREOFLIEDERIVATIVESINTERMSOFSTUCOVARIANTDERIVATIVES})
to deduce that
\begin{align} \label{E:CRUDESECONDCOVARIANTDERIVATIVEBOUNDINTERMSOFLIEDERIVATIVES}
	\rgeo^4 |\angD^2 \xi|^2 
	& \lesssim 
		\rgeo^2 \sum_{l = 1}^3 |\angLie_{\Rot_{(l)}} \angD \xi|^2
		+  \rgeo^2 |\angD \xi|^2.
\end{align}
Using the commutator identity \eqref{E:ANGDANGLIEZTYPE02COMMUTATOR},
the estimate \eqref{E:TYPE02TENSORANGDINTERMSOFROTATIONALLIE}, 
and the estimate \eqref{E:CRUDELOWERORDERC0BOUNDDERIVATIVESOFANGULARDEFORMATIONTENSORS}
(to bound the magnitude of the Lie derivatives of the $\angdeform{\Rot_{(l)}}$),
we deduce that the right-hand side of \eqref{E:CRUDESECONDCOVARIANTDERIVATIVEBOUNDINTERMSOFLIEDERIVATIVES}
is
\begin{align}
	& \lesssim
		\rgeo^2 \sum_{l = 1}^3 |\angD (\angLie_{\Rot_{(l)}} \xi)|^2
		+ \rgeo^2 \sum_{l=1}^3 |\angD \angdeform{\Rot_{(l)}}|^2 |\xi|^2
		+  \rgeo^2 |\angD \xi|^2
			\\
	& \lesssim
		|\angLie_{\mathscr{O}}^{\leq 2} \xi|^2
		+ \sum_{l=1}^3 |\angLie_{\mathscr{O}}^{\leq 1} \angdeform{\Rot_{(l)}}|^2 |\xi|^2
		\lesssim |\angLie_{\mathscr{O}}^{\leq 2} \xi|^2.
		\notag
\end{align}

\end{proof}

\section{Pointwise estimates for \texorpdfstring{$\angLie_{\mathscr{Z}}^N \upchi^{(Small)}$}
{the derivatives of the re-centered null second fundamental form} 
in terms of other quantities}
In this section, we derive pointwise estimates for the $S_{t,u}-$projected Lie derivatives of $\upchi^{(Small)}.$

\begin{lemma}[\textbf{Pointwise estimates for} $\upchi^{(Small)}$ \textbf{in terms of other quantities}] 
\label{L:POINTWISEESTIMATESFORCHIJUNKINTERMSOFOTHERVARIABLES}
	Let $0 \leq N \leq 23$ be an integer.
	Then under the small-data and bootstrap assumptions 
	of Sects.~\ref{S:PSISOLVES}-\ref{S:C0BOUNDBOOTSTRAP},
	the following pointwise estimates hold on $\mathcal{M}_{\Tboot,U_0}:$
	\begin{subequations}
	\begin{align} \label{E:POINTWISEESTIMATESFORCHIJUNKINTERMSOFOTHERVARIABLES}
		\left|
			\sevenmyarray
				[\angLie_{\mathscr{Z}}^N \upchi^{(Small)}]
				{\mathscr{Z}^N \mytr \upchi^{(Small)}}
				{\angLie_{\mathscr{Z}}^N \hat{\upchi}^{(Small)}}
				{\angLie_{\mathscr{Z}}^N \upchi^{(Small) \#}}
				{\angLie_{\mathscr{Z}}^N \hat{\upchi}^{(Small) \#}}
				{\angLie_{\mathscr{Z}}^N \upchi^{(Small) \# \#}}
				{\angLie_{\mathscr{Z}}^N \hat{\upchi}^{(Small) \# \#}}
		 \right|
		& \lesssim
			\frac{1}{1 + t}
			\left|
				\fourmyarray[\rgeo \Lunit \mathscr{Z}^{\leq N} \Psi]
					{\Rad \mathscr{Z}^{\leq N} \Psi}
					{\rgeo \angdiff \mathscr{Z}^{\leq N} \Psi} 
					{\mathscr{Z}^{\leq N} \Psi}
			\right|
		+ \frac{1}{(1 + t)^2}
			\left| 
				\myarray
					[\mathscr{Z}^{\leq N} (\upmu - 1)]
						{\sum_{a=1}^3 \rgeo |\mathscr{Z}^{\leq N+1} \Lunit_{(Small)}^a|}
			\right|,
			\\
		\left|
			\sevenmyarray
				[\angLie_{\mathscr{Z}}^N \upchi^{(Small)}]
				{\mathscr{Z}^N \mytr \upchi^{(Small)}}
				{\angLie_{\mathscr{Z}}^N \hat{\upchi}^{(Small)}}
				{\angLie_{\mathscr{Z}}^N \upchi^{(Small) \#}}
				{\angLie_{\mathscr{Z}}^N \hat{\upchi}^{(Small) \#}}
				{\angLie_{\mathscr{Z}}^N \upchi^{(Small) \# \#}}
				{\angLie_{\mathscr{Z}}^N \hat{\upchi}^{(Small) \# \#}}
		 \right\|_{C^0(\Sigma_t^u)}
		& \lesssim
			\varepsilon^{1/2} \frac{\ln(\myexp + t)}{(1 + t)^2}.
			\label{E:LOWERORDERC0BOUNDCHIJUNK} 
	\end{align}	
	\end{subequations}
\end{lemma}

\begin{proof}
By \eqref{E:CHIJUNKINTERMSOFOTHERVARIABLES}, we have the following identity, where
the last product on the right-hand side is written schematically:
\begin{align} \label{E:CHIJUNKSCHEMATIC}
\upchi_{AB}^{(Small)}
& = g_{ab} (\angdiffarg{A} x^a) \angdiffarg{B} \Lunit_{(Small)}^b
	+ G_{(Frame)} 
		\myarray
			[\Lunit \Psi]
			{\angdiff \Psi}.
\end{align}
We now apply $\angLie_{\mathscr{Z}}^N$ to both sides of \eqref{E:CHIJUNKSCHEMATIC}.
As in our proof of \eqref{E:CRUDEPOINTWISEBOUNDSDERIVATIVESOFANGULARDEFORMATIONTENSORS},
we treat all uppercase Latin indices in the term $g_{ab} (\angdiffarg{A} x^a) \angdiffarg{B} \Lunit_{(Small)}^b$
as tensorial $S_{t,u}$ indices, and all lowercase Latin-indexed quantities as scalar-valued functions.
The Leibniz rule, the identity \eqref{E:ANGLIECOMMUTESWITHANGDIFF}, 
and inequality \eqref{E:FUNCTIONPOINTWISEANGDINTERMSOFANGLIEO} thus yield
\begin{align} \label{E:CHIJUNKDIFFERNTIATEDLEIBNIZ}
	\left|
		\angLie_{\mathscr{Z}}^N \upchi^{(Small)}
	\right| 
	& \lesssim 
	\frac{1}{1 + t}
	\mathop{\sum_{N_1 + N_2 + N_3 \leq N+1}}_{N_1, N_2 \leq N}
		\left|
			\mathscr{Z}^{N_1} g_{ab} 
		\right|
		\left|
			\angdiff \mathscr{Z}^{N_2} x^a
		\right|
		\left|
			\mathscr{Z}^{N_3} \Lunit_{(Small)}^b
		\right|
		\\
	& \ \ 
	+
	\sum_{N_1 + N_2 = N}
		\left|
			\angLie_{\mathscr{Z}}^{N_1} G_{(Frame)} 
		\right| 
		\left|
			\myarray
				[\mathscr{Z}^{N_2} \Lunit \Psi]
				{\angdiff \mathscr{Z}^{N_2} \Psi}
		\right|.
			\notag
\end{align}	
The desired estimate \eqref{E:POINTWISEESTIMATESFORCHIJUNKINTERMSOFOTHERVARIABLES}
for $\angLie_{\mathscr{Z}}^N \upchi^{(Small)}$ now follows from
inequality 
\eqref{E:CHIJUNKDIFFERNTIATEDLEIBNIZ}, 
the estimates 
\eqref{E:GRECTCOMPONENTINTERMSOFOTHERVARIABLES},
\eqref{E:GRECTCOMPONENTC0BOUND},
\eqref{E:LIEDERIVATIVESOFGRAMEINTERMSOFOTHERVARIABLES},
\eqref{E:LOWERORDERC0BOUNDLIEDERIVATIVESOFGRAME},
\eqref{E:POINTWISEBOUNDPROJECTEDLIEDERIVATIVESANGDIFFCOORDINATEX},
\eqref{E:LOWERORDERPOINTWISEBOUNDPROJECTEDLIEDERIVATIVESANGDIFFCOORDINATEX},
Lemma~\ref{L:AVOIDINGCOMMUTING},
and the bootstrap assumptions.
The desired estimate \eqref{E:LOWERORDERC0BOUNDCHIJUNK}
for $\angLie_{\mathscr{Z}}^N \upchi^{(Small)}$ then follows from the
estimate \eqref{E:POINTWISEESTIMATESFORCHIJUNKINTERMSOFOTHERVARIABLES}
for $\angLie_{\mathscr{Z}}^N \upchi^{(Small)}$ and the bootstrap assumptions.

The desired estimates
\eqref{E:LOWERORDERC0BOUNDCHIJUNK}
and
\eqref{E:POINTWISEESTIMATESFORCHIJUNKINTERMSOFOTHERVARIABLES}
for $\upchi^{(Small) \#}$ and $\upchi^{(Small) \# \#}$
then follow from the estimates for $\upchi^{(Small)}$
in the same way that the 
estimates for $\#$ quantities in Lemma~\ref{L:POINTWISEBOUNDSDERIVATIVESOFANGULARDEFORMATIONTENSORS}
follow from the estimates for the non-$\#$ quantities.

The desired estimates
\eqref{E:POINTWISEESTIMATESFORCHIJUNKINTERMSOFOTHERVARIABLES}
and
\eqref{E:LOWERORDERC0BOUNDCHIJUNK}
for $\mytr \upchi^{(Small)}$
then follow from the corresponding estimates for 
$\upchi^{(Small) \#},$
the fact that 
$\mathscr{Z}^{\leq N} \mytr \upchi^{(Small)}$ is the pure trace of
the type $\binom{1}{1}$ tensorfield $\angLie_{\mathscr{Z}}^N \upchi^{(Small) \#},$
and the fact that magnitude of the pure trace of a type $\binom{1}{1}$ tensorfield
is $\lesssim$ the magnitude of the type $\binom{1}{1}$ tensorfield itself.

The desired estimates
\eqref{E:POINTWISEESTIMATESFORCHIJUNKINTERMSOFOTHERVARIABLES}
and
\eqref{E:LOWERORDERC0BOUNDCHIJUNK}
for $\hat{\upchi}^{(Small)}$
then follow from the identity
$\hat{\upchi}^{(Small)} = \upchi^{(Small)} - \frac{1}{2} \mytr \upchi^{(Small)} \gsphere,$
the estimates
\eqref{E:POINTWISEESTIMATESFORCHIJUNKINTERMSOFOTHERVARIABLES}
and
\eqref{E:LOWERORDERC0BOUNDCHIJUNK}
for 
$\upchi^{(Small)}$
and
$\mytr \upchi^{(Small)},$
and the estimates of Lemma~\ref{L:POINTWISEBOUNDSDERIVATIVESOFANGULARDEFORMATIONTENSORS}
for the quantities $\angLie_{\mathscr{Z}}^M \gsphere.$

The desired estimates
\eqref{E:LOWERORDERC0BOUNDCHIJUNK}
and
\eqref{E:POINTWISEESTIMATESFORCHIJUNKINTERMSOFOTHERVARIABLES}
for $\hat{\upchi}^{(Small) \#}$
and $\hat{\upchi}^{(Small) \# \#}$
then follow from the estimates for $\hat{\upchi}^{(Small)}$
in the same way that the 
estimates for $\#$ quantities in Lemma~\ref{L:POINTWISEBOUNDSDERIVATIVESOFANGULARDEFORMATIONTENSORS}
follow from the estimates for the non-$\#$ quantities.

\end{proof}

\begin{corollary}[\textbf{Pointwise estimates for the Lie derivatives of} $G_{(Frame)}^{\#}$]
\label{C:POINTWISEESTIMATESGFRAMESHARPINTERMSOFOTHERQUANTITIES}
The estimates 
\eqref{E:LIEDERIVATIVESOFGRAMEINTERMSOFOTHERVARIABLES} 
and \eqref{E:LOWERORDERC0BOUNDLIEDERIVATIVESOFGRAME}
also hold with $G_{(Frame)}^{\#}$ in place of $G_{(Frame)}.$
\end{corollary}

\begin{proof}
	We have the schematic identity $G_{(Frame)}^{\#} = \ginversesphere G_{(Frame)}.$
	We now differentiate this identity
	with $\angLie_{\mathscr{Z}}^N,$ apply the Leibniz rule,
	and using the estimates
	\eqref{E:LIEDERIVATIVESOFGRAMEINTERMSOFOTHERVARIABLES}
	and
	\eqref{E:LOWERORDERC0BOUNDLIEDERIVATIVESOFGRAME}
	for $G_{(Frame)}$
	and 
	\eqref{E:CRUDEPOINTWISEBOUNDSDERIVATIVESOFANGULARDEFORMATIONTENSORS}
	and 
	\eqref{E:CRUDELOWERORDERC0BOUNDDERIVATIVESOFANGULARDEFORMATIONTENSORS} 
	for $\ginversesphere.$ The corollary thus follows.
\end{proof}

\section{Pointwise estimates for the Lie derivatives of the rotation vectorfields}
In this section, we derive pointwise estimates for the
$S_{t,u}-$projected Lie derivatives of the rotation vectorfields 
$\lbrace \Rot_{(1)}, \Rot_{(2)}, \Rot_{(3)} \rbrace.$

\begin{lemma}[\textbf{Pointwise estimates for the Lie derivatives of the rotation vectorfields}]
\label{L:ROTATIONSPOINTWISEESTIMATES}
	Let $0 \leq N \leq 24$ be an integer.
	Under the small-data and bootstrap assumptions 
	of Sects.~\ref{S:PSISOLVES}-\ref{S:C0BOUNDBOOTSTRAP},
	if $\varepsilon$ is sufficiently small,
	then for each $\Rot \in \lbrace \Rot_{(1)}, \Rot_{(2)}, \Rot_{(3)} \rbrace,$
	the following pointwise estimates hold on $\mathcal{M}_{\Tboot,U_0}:$
\begin{subequations}
\begin{align}
	\left|
		\angLie_{\mathscr{Z}}^{\leq N} \Rot
	\right|
	& \lesssim
		(1 + t)
		\left|
			\mathscr{Z}^{\leq N} \Psi
		\right|
			+ 
			\left|
				\myarray[\mathscr{Z}^{\leq N} (\upmu - 1)]
					{\sum_{a=1}^3 \rgeo |\mathscr{Z}^{\leq N} \Lunit_{(Small)}^a|} 
			\right|
		+ 1 + t,
		\label{E:LIEDERIVATIVESOFROTATIONSPOINTWISEESTIMATE}
			\\
	\left\|
		\angLie_{\mathscr{Z}}^{\leq 12} \Rot
	\right\|_{C^0(\Sigma_t^u)}
	& \lesssim 1 + t.
		\label{E:LOWERORDERLIEDERIVATIVESOFROTATIONSC0BOUND}
\end{align}
\end{subequations}
\end{lemma}

\begin{proof}
To prove \eqref{E:LIEDERIVATIVESOFROTATIONSPOINTWISEESTIMATE}, we decompose
$\Rot^A = g_{ab} (\ginversesphere)^{AB} (\angdiffarg{B} x^a) \Rot^b.$
We now apply $\angLie_{\mathscr{Z}}^N$ to both sides of this identity.
Arguing as in our proof of \eqref{E:CRUDEPOINTWISEBOUNDSDERIVATIVESOFANGULARDEFORMATIONTENSORS},
we apply the Leibniz rule to the right-hand side and 
treat all uppercase Latin indices as tensorial $S_{t,u}$ indices 
and all lowercase Latin-indexed quantities as scalar-valued functions.
Also using Lemma~\ref{L:LANDRADCOMMUTEWITHANGDIFF}, we see that
\begin{align} \label{E:ROTATIONVECTORFIELDSDIFFERNTIATEDLEIBNIZ}
	\left|
		\angLie_{\mathscr{Z}}^N \Rot
	\right| 
	& \lesssim 
		\sum_{N_1 + N_2 + N_3 + N_4 \leq N}
		\left|
			\mathscr{Z}^{N_1} g_{ab} 
		\right|
		\left|
			\angLie_{\mathscr{Z}}^{N_2} \ginversesphere
		\right|
		\left|
			\angdiff \mathscr{Z}^{N_3} x^a
		\right|
		\left|
			\mathscr{Z}^{N_4} \Rot^b
		\right|.
\end{align}
The desired bounds 
\eqref{E:LIEDERIVATIVESOFROTATIONSPOINTWISEESTIMATE}
and \eqref{E:LOWERORDERLIEDERIVATIVESOFROTATIONSC0BOUND}
now follow from the estimates
\eqref{E:ROTATIONVECTORFIELDSDIFFERNTIATEDLEIBNIZ},
\eqref{E:GRECTCOMPONENTINTERMSOFOTHERVARIABLES},
\eqref{E:GRECTCOMPONENTC0BOUND},
\eqref{E:POINTWISEBOUNDPROJECTEDLIEDERIVATIVESANGDIFFCOORDINATEX},
\eqref{E:LOWERORDERPOINTWISEBOUNDPROJECTEDLIEDERIVATIVESANGDIFFCOORDINATEX},
\eqref{E:DERIVATIVESOFRECTANGULARCOMMUTATORVECTORFIELDRECTCOMPONENTS},
\eqref{E:LOWERORDERC0BOUNDDERIVATIVESOFRECTANGULARCOMMUTATORVECTORFIELDRECTCOMPONENTS},
\eqref{E:CRUDEPOINTWISEBOUNDSDERIVATIVESOFANGULARDEFORMATIONTENSORS},
\eqref{E:CRUDELOWERORDERC0BOUNDDERIVATIVESOFANGULARDEFORMATIONTENSORS},
and \eqref{E:BOUNDINGFRAMEDERIVATIVESINTERMSOFZNF},
and the bootstrap assumptions
\eqref{E:PSIFUNDAMENTALC0BOUNDBOOTSTRAP},
\eqref{E:UPMUBOOT},
and
\eqref{E:FRAMECOMPONENTSBOOT}.

\end{proof}

\section{Pointwise estimates for
\texorpdfstring{$\angdeformoneformarg{V}{\Lunit} \ \mbox{and} \ \angdeformoneformarg{V}{\Lunit}$}{the angular one-forms and vectorfields corresponding to the commutation vectorfield
deformation tensors}
}
In this section, for various vectorfields $V,$ we derive pointwise estimates for the 
$S_{t,u}-$projected Lie derivatives of the
$S_{t,u}$ tensorfields
$\angdeformoneformarg{V}{\Lunit},$ 
$\angdeformoneformarg{V}{\Rad},$
and their $\gsphere-$duals.

\begin{lemma}[\textbf{Pointwise estimates for the angular one-forms and vectorfields corresponding to the commutation vectorfield
	deformation tensors}]
	\label{L:STUONEFORMANDVECTORFIELDSCORRESPONDINGTOCOMMUTATORDEFORMATIONTENSORS}
	Let $0 \leq N \leq 23$ be an integer
	and let $\Rot \in \lbrace \Rot_{(1)}, \Rot_{(2)}, \Rot_{(3)} \rbrace.$
	Under the small-data and bootstrap assumptions 
	of Sects.~\ref{S:PSISOLVES}-\ref{S:C0BOUNDBOOTSTRAP},
	if $\varepsilon$ is sufficiently small,
	then the following pointwise estimates hold on $\mathcal{M}_{\Tboot,U_0}:$
	\begin{subequations}
	\begin{align}	
	\left|
		\myarray
			[\angLie_{\mathscr{Z}}^N \angdeformoneformarg{\Rad}{\Lunit}]
			{\angLie_{\mathscr{Z}}^N \angdeformoneformupsharparg{\Rad}{\Lunit}}
	\right|
	& \lesssim
		\left|
				\fourmyarray[\rgeo \Lunit \mathscr{Z}^{\leq N} \Psi]
					{\Rad \mathscr{Z}^{\leq N} \Psi}
					{\rgeo \angdiff \mathscr{Z}^{\leq N} \Psi} 
					{\mathscr{Z}^{\leq N} \Psi}
		\right|
		+ 
			\frac{1}{1 + t}
			\left|
				\myarray[\mathscr{Z}^{\leq N+1} (\upmu - 1)]
					{\sum_{a=1}^3 \rgeo |\mathscr{Z}^{\leq N} \Lunit_{(Small)}^a|} 
			\right|,
		 \label{E:RADDEFORMSPHERELPOINTWISE} \\
		\left\|
			\myarray
				[\angLie_{\mathscr{Z}}^{\leq 11} \angdeformoneformarg{\Rad}{\Lunit}]
				{\angLie_{\mathscr{Z}}^{\leq 11} \angdeformoneformupsharparg{\Rad}{\Lunit}}
		\right\|_{C^0(\Sigma_t^u)}
	& \lesssim \varepsilon^{1/2} \frac{\ln(\myexp + t)}{1 + t},
		\label{E:LOWERORDERC0BOUNDRADDEFORMSPHEREL}
	\end{align}
	\end{subequations}

\begin{subequations}
\begin{align}
	\left|
		\myarray
			[\angLie_{\mathscr{Z}}^N \angdeformoneformarg{\rgeo \Lunit}{\Rad}]
			{\angLie_{\mathscr{Z}}^N \angdeformoneformupsharparg{\rgeo \Lunit}{\Rad}}
	\right|
	& \lesssim
		(1 + t)
		\left|
				\fourmyarray[\rgeo \Lunit \mathscr{Z}^{\leq N} \Psi]
					{\Rad \mathscr{Z}^{\leq N} \Psi}
					{\rgeo \angdiff \mathscr{Z}^{\leq N} \Psi} 
					{\mathscr{Z}^{\leq N} \Psi}
		\right|
		+ 
			\left|
				\myarray[\mathscr{Z}^{\leq N+1} (\upmu - 1)]
					{\sum_{a=1}^3 \rgeo |\mathscr{Z}^{\leq N} \Lunit_{(Small)}^a|} 
			\right|,
				\label{E:RGEOLDEFORMSPHERERADPOINTWISE} \\
	\left\|
		\myarray
			[\angLie_{\mathscr{Z}}^{\leq 11} \angdeformoneformarg{\rgeo \Lunit}{\Rad}]
			{\angLie_{\mathscr{Z}}^{\leq 11} \angdeformoneformupsharparg{\rgeo \Lunit}{\Rad}}
	\right\|_{C^0(\Sigma_t^u)}
	& \lesssim \varepsilon^{1/2} \ln(\myexp + t),
			\label{E:LOWERORDERC0BOUNDRGEOLDEFORMSPHERERAD}
\end{align}
\end{subequations}

\begin{subequations}
\begin{align}	
	\left|
		\myarray
			[\angLie_{\mathscr{Z}}^N \angdeformoneformarg{\Rot}{\Lunit}]
			{\angLie_{\mathscr{Z}}^N \angdeformoneformupsharparg{\Rot}{\Lunit}}
	\right|
	& \lesssim
		\left|
				\fourmyarray[\rgeo \Lunit \mathscr{Z}^{\leq N} \Psi]
					{\Rad \mathscr{Z}^{\leq N} \Psi}
					{\rgeo \angdiff \mathscr{Z}^{\leq N} \Psi} 
					{\mathscr{Z}^{\leq N} \Psi}
			\right|
		+ 
			\frac{1}{1 + t}
			\left|
				\myarray[\mathscr{Z}^{\leq N} (\upmu - 1)]
					{\sum_{a=1}^3 \rgeo |\mathscr{Z}^{\leq N+1} \Lunit_{(Small)}^a|} 
			\right|,
		 \label{E:ROTDEFORMSPHERELPOINTWISE} \\
	\left\|
		\myarray
			[\angLie_{\mathscr{Z}}^{\leq 11} \angdeformoneformarg{\Rot}{\Lunit}]
			{\angLie_{\mathscr{Z}}^{\leq 11} \angdeformoneformupsharparg{\Rot}{\Lunit}}
	\right\|_{C^0(\Sigma_t^u)}
	& \lesssim \varepsilon^{1/2} \frac{\ln(\myexp + t)}{1 + t},
		\label{E:LOWERORDERC0BOUNDROTDEFORMSPHEREL}
	\end{align}
	\end{subequations}

\begin{subequations}
\begin{align}
	\left|
		\myarray
			[\angLie_{\mathscr{Z}}^N \angdeformoneformarg{\Rot}{\Rad}]
			{\angLie_{\mathscr{Z}}^N \angdeformoneformupsharparg{\Rot}{\Rad}}
	\right|
	& \lesssim
		\ln(\myexp + t)
		\left|
				\fourmyarray[\rgeo \Lunit \mathscr{Z}^{\leq N} \Psi]
					{\Rad \mathscr{Z}^{\leq N} \Psi}
					{\rgeo \angdiff \mathscr{Z}^{\leq N} \Psi} 
					{\mathscr{Z}^{\leq N} \Psi}
		\right|
		+ 
			\frac{\ln(\myexp + t)}{1 + t}
			\left|
				\myarray[\mathscr{Z}^{\leq N+1} (\upmu - 1)]
					{\sum_{a=1}^3 \rgeo |\mathscr{Z}^{\leq N+1} \Lunit_{(Small)}^a|} 
			\right|,
				\label{E:ROTDEFORMSPHERERADPOINTWISE} \\
	\left\|
		\myarray
			[\angLie_{\mathscr{Z}}^{\leq 11} \angdeformoneformarg{\Rot}{\Rad}]
			{\angLie_{\mathscr{Z}}^{\leq 11} \angdeformoneformupsharparg{\Rot}{\Rad}}
	\right\|_{C^0(\Sigma_t^u)}
	& \lesssim \varepsilon^{1/2} \frac{\ln^2(\myexp + t)}{1 + t}.
			\label{E:LOWERORDERC0BOUNDROTDEFORMSPHERERAD}
\end{align}
\end{subequations}

\end{lemma}

\begin{proof}
To prove \eqref{E:ROTDEFORMSPHERERADPOINTWISE} for $\angLie_{\mathscr{Z}}^{\leq N} \angdeformoneformarg{\Rot}{\Rad},$ 
we first use \eqref{E:ROTDEFORMRADSPHERE} and \eqref{E:ERRORROTDEFORMRADSPHERE}
to decompose (see Remark~\ref{R:CLARIFICATIONOFROT})
\begin{align} \label{E:ROTDEFORMRADSPHERESCHEMATIC}
	\angdeformarg{\Rot_{(l)}}{\Rad}{A}
	& = 	\upmu \upchi_{AB}^{(Small)} \Rot_{(l)}^B 
				+ \RotRadcomponent{l} \angdiffarg{A} \upmu
					\\
	& \ \ +
				G_{(Frame)} 
				\myarray
					[\upmu 	\Rot_{(l)}]
					{\upmu \RotRadcomponent{l}}
			\myarray
				[\Lunit \Psi]
				{\angdiff \Psi}
		+ G_{(Frame)} 
			\RotRadcomponent{l} 
			\Rad \Psi
		+ \upmu
			\smoothfunction(\Psi) 
			\myarray[\Psi]
				{\Lunit_{(Small)}}	 
				\angdiff x.
				\notag
\end{align}
We now apply $\angLie_{\mathscr{Z}}^N$ to the terms on the right-hand
side of \eqref{E:ROTDEFORMRADSPHERESCHEMATIC} and apply the Leibniz rule.
To bound $\mathscr{Z}^M (\upmu - 1)$ for $M \leq 12,$ we use
the bootstrap assumptions \eqref{E:UPMUBOOT}, 
while the higher-order derivatives of $\upmu - 1$ appear explicitly on the right-hand side of
\eqref{E:ROTDEFORMSPHERERADPOINTWISE}.
To bound $\angLie_{\mathscr{Z}}^M \angdiff \upmu$ for $M \leq 11,$ 
we use Lemma~\ref{L:LANDRADCOMMUTEWITHANGDIFF},
inequality \eqref{E:FUNCTIONPOINTWISEANGDINTERMSOFANGLIEO},
and the bootstrap assumptions \eqref{E:UPMUBOOT}, 
while we bound the higher-order derivatives
by $\left|\angLie_{\mathscr{Z}}^M \angdiff \upmu \right| 
\lesssim (1+t)^{-1} \mathscr{Z}^{M+1} (\upmu - 1),$ 
and the latter term appears
explicitly on the right-hand side of
\eqref{E:ROTDEFORMSPHERERADPOINTWISE}.
We similarly bound
$\angLie_{\mathscr{Z}}^M \Lunit_{(Small)}$
with the help of the bootstrap assumptions \eqref{E:FRAMECOMPONENTSBOOT}.
To bound $\angLie_{\mathscr{Z}}^M \upchi^{(Small)},$
we use Lemma~\ref{L:POINTWISEESTIMATESFORCHIJUNKINTERMSOFOTHERVARIABLES}.
To bound $\angLie_{\mathscr{Z}}^M \Rot_{(l)},$
we use Lemma~\ref{L:ROTATIONSPOINTWISEESTIMATES}.
To bound $\angLie_{\mathscr{Z}}^M \RotRadcomponent{l},$
we use inequalities 
\eqref{E:FIRSTPOINTWISEBOUNDEUCLIDEANROTATIONRADCOMPONENT}
and \eqref{E:LOWERORDERC0BOUNDEUCLIDEANROTATIONRADCOMPONENT}.
We bound the terms 
$\angLie_{\mathscr{Z}}^M G_{(Frame)}$ with Lemma~\ref{L:POINTWISEESTIMATESGFRAMEINTERMSOFOTHERQUANTITIES}.
To bound the terms
$\myarray[\mathscr{Z}^M \Lunit \Psi] {\angLie_{\mathscr{Z}}^M \angdiff \Psi},$
$\angLie_{\mathscr{Z}}^M \smoothfunction(\Psi),$ 
and $\mathscr{Z}^M \Psi,$
we use Lemma~\ref{L:AVOIDINGCOMMUTING}
and the bootstrap assumptions \eqref{E:PSIFUNDAMENTALC0BOUNDBOOTSTRAP}.
To bound the terms $\angLie_{\mathscr{Z}}^M \angdiff x,$
we use the estimates 
\eqref{E:POINTWISEBOUNDPROJECTEDLIEDERIVATIVESANGDIFFCOORDINATEX}
and
\eqref{E:LOWERORDERPOINTWISEBOUNDPROJECTEDLIEDERIVATIVESANGDIFFCOORDINATEX}.
In total, these estimates yield inequality \eqref{E:ROTDEFORMSPHERERADPOINTWISE}
for $\angLie_{\mathscr{Z}}^{\leq N} \angdeformoneformarg{\Rot}{\Rad}.$
We can prove the estimate \eqref{E:ROTDEFORMSPHERERADPOINTWISE}
for 
$\angLie_{\mathscr{Z}}^{\leq N} \angdeformoneformupsharparg{\Rot}{\Rad} = \angLie_{\mathscr{Z}}^{\leq N} (\ginversesphere \angdeformoneformarg{\Rot}{\Rad})$
in a similar fashion, but we also need to use the estimates 
\eqref{E:CRUDEPOINTWISEBOUNDSDERIVATIVESOFANGULARDEFORMATIONTENSORS}
and \eqref{E:CRUDELOWERORDERC0BOUNDDERIVATIVESOFANGULARDEFORMATIONTENSORS}
to bound the factors $\angLie_{\mathscr{Z}}^M \ginversesphere$ that arise in the estimates.

Inequality \eqref{E:LOWERORDERC0BOUNDROTDEFORMSPHERERAD} then follows from
\eqref{E:ROTDEFORMSPHERERADPOINTWISE} and the bootstrap assumptions
\eqref{E:PSIFUNDAMENTALC0BOUNDBOOTSTRAP},
\eqref{E:UPMUBOOT},
and \eqref{E:FRAMECOMPONENTSBOOT}.

Using the identities
\eqref{E:ROTDEFORMLSPHERE}
and 
\eqref{E:ERRORROTDEFORMLSPHERE},
we similarly deduce
inequalities 
\eqref{E:ROTDEFORMSPHERELPOINTWISE}
and \eqref{E:LOWERORDERC0BOUNDROTDEFORMSPHEREL}.

A similar but simpler proof
based on the identity \eqref{E:RADDEFORMLA}
yields inequalities
\eqref{E:RADDEFORMSPHERELPOINTWISE}
and
\eqref{E:LOWERORDERC0BOUNDRADDEFORMSPHEREL}.

A similar but simpler proof
based on the identity \eqref{E:RGEOLDEFORMRADA}
yields inequalities
\eqref{E:RGEOLDEFORMSPHERERADPOINTWISE} 
and \eqref{E:LOWERORDERC0BOUNDRGEOLDEFORMSPHERERAD}.

\end{proof}

\section{Preliminary Lie commutator estimates}

In this section, we derive some preliminary commutator estimates
involving $S_{t,u}-$projected Lie derivatives.

\begin{lemma}[\textbf{Preliminary quantitative estimate involving Lie commutators}] \label{L:PRELIMINARYQUANTITATIVELIECOMMUTATIONS}
	Let $1 \leq N \leq 23$ be an integer, and let $\xi$ be an $S_{t,u}$ one-form 
	or a symmetric type $\binom{0}{2}$ $S_{t,u}$ tensorfield. 
	Let $\Rot \in \lbrace \Rot_{(1)}, \Rot_{(2)}, \Rot_{(3)} \rbrace.$
	Under the small-data and bootstrap assumptions 
	of Sects.~\ref{S:PSISOLVES}-\ref{S:C0BOUNDBOOTSTRAP},
	if $\varepsilon$ is sufficiently small, then the following
	pointwise estimates hold on $\mathcal{M}_{\Tboot,U_0}:$
	\begin{subequations}
	\begin{align} 
		\left| 
			[\angLie_{\rgeo \Lunit}, \angLie_{\mathscr{Z}}^N] \xi 
		\right|
		& \lesssim
				\left| \angLie_{\Lunit} \angLie_{\mathscr{Z}}^{\leq N-1} \xi \right|
				+
				\frac{1}{1 + t}
				\sum_{N_1 + N_2 \leq N}
				\max_{l=1,2,3}	
					\left|
						\angLie_{\mathscr{Z}}^{N_1}
						\myarray
							[\rgeo \angdeformoneformupsharparg{\Rad}{\Lunit}]
								{\rgeo \angdeformoneformupsharparg{\Rot_{(l)}}{\Lunit}}
					 \right|	
					\left| \angLie_{\mathscr{Z}}^{N_2} \xi \right|,
			\label{E:PRELIMINARYQUANTITATIVERGEOLLIECOMMUTATIONSTENSORS}
		\\
		\left| 
			[\angLie_{\Rad}, \angLie_{\mathscr{Z}}^N] \xi 
		\right|
		& \lesssim
				\left| \angLie_{\Lunit} \angLie_{\mathscr{Z}}^{\leq N-1} \xi \right|
				+
				\frac{1}{1 + t}
				\sum_{N_1 + N_2 \leq N}
				\max_{l=1,2,3}	
					\left|
						\angLie_{\mathscr{Z}}^{N_1}
						\threemyarray
							[\rgeo \angdeformoneformupsharparg{\Rad}{\Lunit}]
								{\angdeformoneformupsharparg{\Rot_{(l)}}{\Lunit}}
								{\angdeformoneformupsharparg{\Rot_{(l)}}{\Rad}}
					\right|	
					\left| \angLie_{\mathscr{Z}}^{N_2} \xi \right|,
			\label{E:PRELIMINARYQUANTITATIVERADLIECOMMUTATIONSTENSORS}
		\\
		\left| 
			[\angLie_{\Rot}, \angLie_{\mathscr{Z}}^N] \xi 
		\right|
		& \lesssim
				\frac{1}{1 + t}
				\sum_{N_1 + N_2 \leq N}
				\max_{l,m=1,2,3}	
					\left|
						\angLie_{\mathscr{Z}}^{N_1}
						\threemyarray
							[\rgeo \angdeformoneformupsharparg{\Rot_{(l)}}{\Lunit}]
								{\angdeformoneformupsharparg{\Rot_{(l)}}{\Rad}}
								{\angLie_{\Rot_{(l)}} \Rot_{(m)}}
						\right|	
					\left| \angLie_{\mathscr{Z}}^{N_2} \xi \right|,
			\label{E:PRELIMINARYQUANTITATIVEROTLIECOMMUTATIONSTENSORS}
		\\
		\left| 
			[\angLie_{\Lunit}, \angLie_{\mathscr{Z}}^N] \xi 
		\right|
		& \lesssim
				\left| \angLie_{\Lunit} \angLie_{\mathscr{Z}}^{\leq N-1} \xi \right|
				+
				\frac{1}{1 + t}
					\sum_{N_1 + N_2 \leq N}
				\max_{l=1,2,3}	
					\left|
						\angLie_{\mathscr{Z}}^{N_1}
						\myarray
							[\angdeformoneformupsharparg{\Rad}{\Lunit}]
								{\angdeformoneformupsharparg{\Rot_{(l)}}{\Lunit}}
					 \right|	
					\left| \angLie_{\mathscr{Z}}^{N_2} \xi \right|,
				\label{E:PRELIMINARYQUANTITATIVELIECOMMUTATIONS}
					\\
		\left| 
			[\angLie_{\Lunit}, \angLie_{\mathscr{S}}^N] \xi 
		\right|
		& \lesssim
				\frac{1}{1 + t}
					\sum_{N_1 + N_2 \leq N}
				\max_{l=1,2,3}	
					\left|
						\angLie_{\mathscr{Z}}^{N_1}
						\myarray
							[\angdeformoneformupsharparg{\Rad}{\Lunit}]
								{\angdeformoneformupsharparg{\Rot_{(l)}}{\Lunit}}
					 \right|	
					\left| \angLie_{\mathscr{Z}}^{N_2} \xi \right|.
				\label{E:SPATIALANDLPRELIMINARYQUANTITATIVELIECOMMUTATIONS}
	\end{align}
	\end{subequations}
	
	Furthermore, if $1 \leq N \leq 24$ and $f$ is a function, then
	\begin{align} \label{E:PRELIMINARYQUANTITATIVELIEDERIVATIVEBOUNDFORFUNCTIONS}
		& \mbox{the estimates \ }
		\eqref{E:PRELIMINARYQUANTITATIVERGEOLLIECOMMUTATIONSTENSORS}-\eqref{E:SPATIALANDLPRELIMINARYQUANTITATIVELIECOMMUTATIONS}
		\mbox{\ also hold with \ } 
		\xi
		\mbox{\ replaced by \ }f, 
			\\
		& \mbox{but with the sums \ }
		\sum_{N_1 + N_2 \leq N}
		\mbox{\ on the right-hand sides replaced by \ }
		\mathop{\sum_{N_1 + N_2 \leq N}}_{N_1 \leq N-1}.
		\notag
	\end{align}	
	
\end{lemma}

\begin{proof}
In this proof, we sometimes use 
Lemma~\ref{L:VECTORFIELDCOMMUTATORS}
and
Cor.~\ref{C:STUPROJECTEDDERIVATIVESCOMMUTATOR} silently.
We begin by proving \eqref{E:PRELIMINARYQUANTITATIVELIECOMMUTATIONS}.
As a first step, let $Z \in \mathscr{Z}.$ Then either $Z = \rgeo \Lunit,$ in which case
$[\Lunit,Z] = \Lunit,$ or $Z = S \in \mathscr{S}$ is a spatial commutation vectorfield, 
in which case $[\Lunit,S]$ is $S_{t,u}-$tangent 
by Lemma~\ref{L:VECTORFIELDCOMMUTATORS}.
In the former case, we have 
$[\angLie_{\Lunit}, \angLie_Z] \xi = \angLie_{\Lunit} \xi,$
while in the latter case, we have 
$[\angLie_{\Lunit}, \angLie_S] \xi = \angLie_{[\Lunit, S]} \xi = \angLie_{\angdeformoneformupsharparg{S}{\Lunit}} \xi.$
Iterating these identities, we deduce the following schematic identity,
where some terms on the right-hand side may be absent:
\begin{align} \label{E:FIRSTSCHEMATICPRELIMINARYQUANTITATIVELIECOMMUTATIONS}
	[\angLie_{\Lunit}, \angLie_{\mathscr{Z}}^N] \xi
	& = \sum_{M \leq N-1} \angLie_{\Lunit} \angLie_{\mathscr{Z}}^M \xi
		+ \sum_{S \in \mathscr{S}} \sum_{N_1 + N_2 \leq N-1} 
				\angLie_{\angLie_{\mathscr{Z}}^{N_1} \angdeformoneformupsharparg{S}{\Lunit}} \angLie_{\mathscr{Z}}^{N_2} \xi.
\end{align}
The desired estimate \eqref{E:PRELIMINARYQUANTITATIVELIECOMMUTATIONS}
now follows 
from \eqref{E:FIRSTSCHEMATICPRELIMINARYQUANTITATIVELIECOMMUTATIONS}
and
inequality \eqref{E:POINTWISEANGLIEXINTERMSOFANGLIEO}.

To prove the estimate \eqref{E:PRELIMINARYQUANTITATIVELIEDERIVATIVEBOUNDFORFUNCTIONS}
in the case of $[\Lunit, \mathscr{Z}^N] f,$ we first use an 
argument similar to the one used to derive \eqref{E:FIRSTSCHEMATICPRELIMINARYQUANTITATIVELIECOMMUTATIONS}
together with Lemma~\ref{L:LANDRADCOMMUTEWITHANGDIFF} 
in order to deduce that
\begin{align} \label{E:FIRSTSCHEMATICPRELIMINARYQUANTITATIVELIECOMMUTATIONSFORFUNCTIONS}
	[\Lunit, \mathscr{Z}^N] f
	& = \sum_{M \leq N-1} \Lunit \mathscr{Z}^M f
		+ \sum_{S \in \mathscr{S}} \sum_{N_1 + N_2 \leq N-1} 
				(\angLie_{\mathscr{Z}}^{N_1} \angdeformoneformupsharparg{S}{\Lunit}) \cdot \angdiff \mathscr{Z}^{N_2} f.
\end{align}
The desired estimate \eqref{E:PRELIMINARYQUANTITATIVELIEDERIVATIVEBOUNDFORFUNCTIONS} in this case
now follows from \eqref{E:FIRSTSCHEMATICPRELIMINARYQUANTITATIVELIECOMMUTATIONSFORFUNCTIONS}
and inequality \eqref{E:FUNCTIONPOINTWISEANGDINTERMSOFANGLIEO}.
Note that inequality \eqref{E:FIRSTSCHEMATICPRELIMINARYQUANTITATIVELIECOMMUTATIONSFORFUNCTIONS}
involves one fewer derivatives of $\angdeformoneformupsharparg{S}{\Lunit}$  
compared to \eqref{E:PRELIMINARYQUANTITATIVELIECOMMUTATIONS} because we 
do not need to use inequality \eqref{E:POINTWISEANGLIEXINTERMSOFANGLIEO} 
in the proof of \eqref{E:FIRSTSCHEMATICPRELIMINARYQUANTITATIVELIECOMMUTATIONSFORFUNCTIONS}.

The proofs of 
\eqref{E:PRELIMINARYQUANTITATIVERGEOLLIECOMMUTATIONSTENSORS}-\eqref{E:PRELIMINARYQUANTITATIVEROTLIECOMMUTATIONSTENSORS}
and \eqref{E:SPATIALANDLPRELIMINARYQUANTITATIVELIECOMMUTATIONS}
are similar to the proof of \eqref{E:PRELIMINARYQUANTITATIVELIECOMMUTATIONS}.
To proceed, we first note that by 
Lemma~\ref{L:VECTORFIELDCOMMUTATORS}
and \eqref{E:ZNAPPLIEDTORGEOISNOTTOOLARGE},
if $Z \in \mathscr{Z},$ then
$[\rgeo \Lunit, Z]$ is either 
$\rgeo \angdeformoneformupsharparg{\Rad}{\Lunit} + \Lunit$
or
$\rgeo \angdeformoneformupsharparg{\Rot_{(l)}}{\Lunit};$
$[\Rad, Z]$ is either 
$-(\rgeo \angdeformoneformupsharparg{\Rad}{\Lunit} + \Lunit)$
or
$ \angdeformoneformupsharparg{\Rot_{(l)}}{\Rad};$
$[\Rot, Z]$ is either 
$\pm \rgeo \angdeformoneformupsharparg{\Rot}{\Lunit},$
$- \angdeformoneformupsharparg{\Rot}{\Rad},$
or $\angLie_{\Rot} \Rot_{(m)};$
and if $S \in \mathscr{S},$ then
$[\Lunit,S] = \angdeformoneformupsharparg{S}{\Lunit}.$
Based on these commutator identities and \eqref{E:ZNAPPLIEDTORGEOISNOTTOOLARGE},
the remainder of the proofs of
\eqref{E:PRELIMINARYQUANTITATIVERGEOLLIECOMMUTATIONSTENSORS}-\eqref{E:PRELIMINARYQUANTITATIVEROTLIECOMMUTATIONSTENSORS}
and \eqref{E:SPATIALANDLPRELIMINARYQUANTITATIVELIECOMMUTATIONS}
essentially mirrors the proof of \eqref{E:PRELIMINARYQUANTITATIVELIECOMMUTATIONS}.

The proofs of the remaining estimates in \eqref{E:PRELIMINARYQUANTITATIVELIEDERIVATIVEBOUNDFORFUNCTIONS}
are similarly connected to the proof in the case of $[\Lunit, \mathscr{Z}^N] f.$

\end{proof}

\section{Commutator estimates for vectorfields acting on functions and \texorpdfstring{$S_{t,u}$}{spherical} covariant tensorfields}
In this section, we derive pointwise commutator estimates for vectorfields acting on functions and $S_{t,u}$
covariant tensorfields.

\begin{remark}[\textbf{Floor and ceiling functions}]
	\label{R:FLOORFUNCTION}
	In what follows, $\lfloor \cdot \rfloor$
	and $\lceil \cdot \rceil$
	respectively denote the floor and ceiling functions. 
	That is, if $M$ is a non-negative integer,
	then $\lfloor M/2 \rfloor = M/2$ for $M$ even and 
	$\lfloor M/2 \rfloor = (M-1)/2$ for $M$ odd,
	while
	$\lceil M/2 \rceil = M/2$ for $M$ even and 
	$\lceil M/2 \rceil = (M+1)/2$ for $M$ odd.
\end{remark}

\begin{lemma}[\textbf{Commutator estimates for vectorfields acting on functions}]
\label{L:COMMUTATORESTIMATESVECTORFIELDSACTINGONFUNCTIONS}
	Let $1 \leq N \leq 24$ be an integer, and let $f$ be a function. 
	Let $\Rot \in \lbrace \Rot_{(1)}, \Rot_{(2)}, \Rot_{(3)} \rbrace.$
	Under the small-data and bootstrap assumptions 
	of Sects.~\ref{S:PSISOLVES}-\ref{S:C0BOUNDBOOTSTRAP},
	if $\varepsilon$ is sufficiently small,
	then the following pointwise estimates hold on $\mathcal{M}_{\Tboot,U_0}:$
\begin{subequations}
\begin{align}
	\left|
		\myarray[{[\mathscr{Z}^N, \rgeo \Lunit] f}]
			{{[\mathscr{Z}^N, \Rad] f}}
	\right|
	& \lesssim 
			\left|
				\Lunit \mathscr{Z}^{\leq N-1} f
			\right|
			+
			\varepsilon^{1/2}
			\frac{\ln(\myexp + t)}{1 + t}
			\left| 
				\fourmyarray[\rgeo \Lunit \mathscr{Z}^{\leq N-1} f]
					{\Rad \mathscr{Z}^{\leq N-1} f}
					{\rgeo \angdiff \mathscr{Z}^{\leq N-1} f}
					{\mathscr{Z}^{\leq N-1} f}
			\right|
			\label{E:RGEOLORRADZNCOMMUTATORACTINGONFUNCTIONSPOINTWISE} \\
	& \ \ 
			+
			\left\|
				\mathscr{Z}^{\leq \lfloor N/2 \rfloor} f
			\right\|_{C^0(\Sigma_t^u)}
			\left| 
				\fourmyarray[ \rgeo \Lunit \mathscr{Z}^{\leq N-1} \Psi]
					{\Rad \mathscr{Z}^{\leq N-1} \Psi}
					{\rgeo \angdiff \mathscr{Z}^{\leq N-1} \Psi}
					{\mathscr{Z}^{\leq N-1} \Psi}
			\right|
			\notag \\
	&  \ \ 
			+ 
			\frac{1}{1 + t}
			\left\|
				\mathscr{Z}^{\leq \lfloor N/2 \rfloor} f
			\right\|_{C^0(\Sigma_t^u)}
			\left|
				\myarray[\mathscr{Z}^{\leq N} (\upmu - 1)]
					{\sum_{a=1}^3 \rgeo |\mathscr{Z}^{\leq N} \Lunit_{(Small)}^a|} 
			\right|,
		\notag \\
	\left|
		[\mathscr{Z}^N, \Rot] f
	\right|
	& \lesssim 
		 \left| 
				\fourmyarray[\rgeo \Lunit \mathscr{Z}^{\leq N-1} f]
					{\Rad \mathscr{Z}^{\leq N-1} f}
					{\rgeo \angdiff \mathscr{Z}^{\leq N-1} f}
					{\mathscr{Z}^{\leq N-1} f}
			\right|
			+
			(1 + t)
			\left\|
				\mathscr{Z}^{\leq \lfloor N/2 \rfloor} f
			\right\|_{C^0(\Sigma_t^u)}
			\left| 
				\fourmyarray[ \rgeo \Lunit \mathscr{Z}^{\leq N-1} \Psi]
					{\Rad \mathscr{Z}^{\leq N-1} \Psi}
					{\rgeo \angdiff \mathscr{Z}^{\leq N-1} \Psi}
					{\mathscr{Z}^{\leq N-1} \Psi}
			\right|
			\label{E:ROTZNCOMMUTATORACTINGONFUNCTIONSPOINTWISE}  \\
	&  \ \ 
			+ 
			\left\|
				\mathscr{Z}^{\leq \lfloor N/2 \rfloor} f
			\right\|_{C^0(\Sigma_t^u)}
			\left|
				\myarray[\mathscr{Z}^{\leq N} (\upmu - 1)]
					{\sum_{a=1}^3 \rgeo |\mathscr{Z}^{\leq N} \Lunit_{(Small)}^a|} 
			\right|,
		\notag \\	
	\left|
		[\mathscr{Z}^N, \Lunit] f
	\right|
	& \lesssim 
			\left|
				\Lunit \mathscr{Z}^{\leq N-1} f
			\right|
			+
			\varepsilon^{1/2}
			\frac{\ln(\myexp + t)}{(1 + t)^2}
			\left| 
				\fourmyarray[ \rgeo \Lunit \mathscr{Z}^{\leq N-1} f]
					{\Rad \mathscr{Z}^{\leq N-1} f}
					{\rgeo \angdiff \mathscr{Z}^{\leq N-1} f}
					{\mathscr{Z}^{\leq N-1} f}
			\right|
			\label{E:LZNCOMMUTATORACTINGONFUNCTIONSPOINTWISE} \\
	& \ \ 
			+ 
			\frac{1}{1 + t}
			\left\|
				\mathscr{Z}^{\leq \lfloor N/2 \rfloor} f
			\right\|_{C^0(\Sigma_t^u)}
			\left| 
				\fourmyarray[ \rgeo \Lunit \mathscr{Z}^{\leq N-1} \Psi]
					{\Rad \mathscr{Z}^{\leq N-1} \Psi}
					{\rgeo \angdiff \mathscr{Z}^{\leq N-1} \Psi}
					{\mathscr{Z}^{\leq N-1} \Psi}
			\right|
			\notag \\
	&  \ \ 
			+ 
			\frac{1}{(1 + t)^2}
			\left\|
				\mathscr{Z}^{\leq \lfloor N/2 \rfloor} f
			\right\|_{C^0(\Sigma_t^u)}
			\left|
				\myarray[\mathscr{Z}^{\leq N} (\upmu - 1)]
					{\sum_{a=1}^3 \rgeo |\mathscr{Z}^{\leq N} \Lunit_{(Small)}^a|} 
			\right|.
				\notag 
\end{align}
\end{subequations}

In addition,
\begin{align} \label{E:ORDERSWITCHEDLZNCOMMUTATORACTINGONFUNCTIONSPOINTWISE} 
	\eqref{E:LZNCOMMUTATORACTINGONFUNCTIONSPOINTWISE} 
	\mbox{\ also holds with the first term on the right replaced by}
	 	\left|
			 \mathscr{Z}^{\leq N-1} \Lunit f
		\right|.
\end{align}

In addition,
\begin{align} \label{E:LSNCOMMUTATORACTINGONFUNCTIONSPOINTWISE} 
	& \eqref{E:LZNCOMMUTATORACTINGONFUNCTIONSPOINTWISE} 
	\mbox{\ also holds without the first term on the right if the left-hand side is equal to \ }
		\\
	& \left|
			[\mathscr{S}^N, \Lunit] f
		\right|,
	\mbox{\ where }	
	\mathscr{S}^N 
	\mbox{\ is an $N^{th}$ order pure spatial commutation vectorfield operator}.
		\notag
\end{align}

\end{lemma}

\begin{proof}
	The two estimates in \eqref{E:RGEOLORRADZNCOMMUTATORACTINGONFUNCTIONSPOINTWISE} 
	and the estimate \eqref{E:LZNCOMMUTATORACTINGONFUNCTIONSPOINTWISE}
	follow from
	\eqref{E:PRELIMINARYQUANTITATIVELIEDERIVATIVEBOUNDFORFUNCTIONS},
	\eqref{E:ZNAPPLIEDTORGEOISNOTTOOLARGE},	
	\eqref{E:FUNCTIONAVOIDINGCOMMUTING},
	\eqref{E:FUNCTIONPOINTWISEANGDINTERMSOFANGLIEO},
	\eqref{E:RADDEFORMSPHERELPOINTWISE},
	\eqref{E:LOWERORDERC0BOUNDRADDEFORMSPHEREL}
	\eqref{E:RGEOLDEFORMSPHERERADPOINTWISE},
	\eqref{E:LOWERORDERC0BOUNDRGEOLDEFORMSPHERERAD},
	\eqref{E:ROTDEFORMSPHERELPOINTWISE},
	\eqref{E:LOWERORDERC0BOUNDROTDEFORMSPHEREL},
	\eqref{E:ROTDEFORMSPHERERADPOINTWISE},
	\eqref{E:LOWERORDERC0BOUNDROTDEFORMSPHERERAD},
	and the fact that
	$\rgeo \angdeformoneformarg{\Rad}{\Lunit} = - \angdeformoneformarg{\rgeo \Lunit}{\Rad}.$
	The estimate \eqref{E:ORDERSWITCHEDLZNCOMMUTATORACTINGONFUNCTIONSPOINTWISE}  
	then follows inductively from \eqref{E:LZNCOMMUTATORACTINGONFUNCTIONSPOINTWISE}.
	
	The proof of 
	\eqref{E:LSNCOMMUTATORACTINGONFUNCTIONSPOINTWISE} 
	is similar, but relies on inequality
	\eqref{E:SPATIALANDLPRELIMINARYQUANTITATIVELIECOMMUTATIONS}.
	
	The proof of \eqref{E:ROTZNCOMMUTATORACTINGONFUNCTIONSPOINTWISE} is similar, 
	but we also need use the estimates
	\eqref{E:LIEDERIVATIVESOFROTATIONSPOINTWISEESTIMATE}
	and
	\eqref{E:LOWERORDERLIEDERIVATIVESOFROTATIONSC0BOUND},
	which results in a less favorable right-hand side.
	
	
\end{proof}

\begin{lemma}[\textbf{Pointwise estimates involving commutations with the operator} $\Lunit + \frac{1}{2} \mytr \upchi$]
\label{L:COMMUTINGFUNCTIONSWITHLPLUSHALFTRACECHI}
	Let $1 \leq N \leq 23$ be an integer and let $f$ be a function.
	Under the small-data and bootstrap assumptions 
	of Sects.~\ref{S:PSISOLVES}-\ref{S:C0BOUNDBOOTSTRAP},
	if $\varepsilon$ is sufficiently small,
	then the following pointwise estimates hold on $\mathcal{M}_{\Tboot,U_0}:$
	\begin{align} \label{E:COMMUTINGFUNCTIONSWITHLPLUSHALFTRACECHI}
		\left|
			\left[
				\mathscr{Z}^N,
				\left\lbrace 
					\Lunit 
					+ \frac{1}{2} \mytr \upchi 
				\right\rbrace
			\right]
		\right|
		& \lesssim
			\left|
				\left\lbrace
					\Lunit + \frac{1}{2} \mytr \upchi 
				\right\rbrace
				\mathscr{Z}^{\leq N-1} f
			\right|
		+ \frac{\ln(\myexp + t)}{(1 + t)^2}
			\left| 
				\threemyarray
					[\Rad \mathscr{Z}^{\leq N-1} f]
					{\rgeo \angdiff \mathscr{Z}^{\leq N-1} f}
					{\mathscr{Z}^{\leq N-1} f}
			\right|	
				\\
	& \ \	
			+ 
			\frac{1}{1 + t}
			\left\|
				\mathscr{Z}^{\leq \lfloor N/2 \rfloor} f
			\right\|_{C^0(\Sigma_t^u)}
			\left| 
				\fourmyarray[ \rgeo \Lunit \mathscr{Z}^{\leq N-1} \Psi]
					{\Rad \mathscr{Z}^{\leq N-1} \Psi}
					{\rgeo \angdiff \mathscr{Z}^{\leq N-1} \Psi}
					{\mathscr{Z}^{\leq N-1} \Psi}
			\right|
			\notag	\\
		& \ \
				+ 
				\frac{1}{(1 + t)^2}
				\left\|
					\mathscr{Z}^{\leq \lfloor N/2 \rfloor} f
				\right\|_{C^0(\Sigma_t^u)}
				\left|	
					\myarray
						[\mathscr{Z}^{\leq N} (\upmu - 1)]
						{\sum_{a=1}^3 \rgeo |\mathscr{Z}^{\leq N} \Lunit_{(Small)}^a|} 
				\right|.
			\notag
	\end{align}	
\end{lemma}

\begin{proof}
We first recall the decomposition $\frac{1}{2} \mytr \upchi = \rgeo^{-1} + \frac{1}{2} \mytr \upchi^{(Small)}.$
We separate the proof into case \textbf{i)}
$\mathscr{Z}^N = \mathscr{Z}^{N-1} S,$ where $S \in \lbrace \Rad, \Rot_{(1)}, \Rot_{(2)}, \Rot_{(3)} \rbrace$
is a spatial commutation vectorfield, 
and case \textbf{ii)} $\mathscr{Z}^N = \mathscr{Z}^{N-1} Z,$
where $Z = \rgeo \Lunit.$ 
Using the aforementioned decomposition, the identity $\Lunit \rgeo = 1,$
and the fact that $[\Lunit, S] =  \angdeformoneformupsharparg{S}{\Lunit}$
(see Lemma~\ref{L:CONNECTIONBETWEENANGLIEOFGSPHEREANDDEFORMATIONTENSORS}),
we compute that in cases \textbf{i)} and \textbf{ii)}, we respectively have
\begin{align}  \label{E:FIRSTSTEPCOMMUTINGFUNCTIONSWITHLPLUSHALFTRACECHI}
	\mathscr{Z}^N
		\left\lbrace 
			\Lunit f
			+ \frac{1}{2} \mytr \upchi f
		\right\rbrace
	& = 	
		\mathscr{Z}^{N-1}
		\left\lbrace 
			\Lunit S f
			+ \frac{1}{2} \mytr \upchi S f
		\right\rbrace
			\\
	& \ \ 
		- \mathscr{Z}^{N-1} (\angdeformoneformupsharparg{S}{\Lunit} \cdot \angdiff f)
		+ \frac{1}{2} 
			\mathscr{Z}^{N-1} 
			\left\lbrace
				(S \rgeo^{-1} + S \mytr \upchi^{(Small)}) f
			\right\rbrace,
				\notag \\
	\mathscr{Z}^N
		\left\lbrace 
			\Lunit f
			+ \frac{1}{2} \mytr \upchi f
		\right\rbrace
	& = \mathscr{Z}^{N-1}
		\left\lbrace 
			\Lunit Z f
			+ \frac{1}{2} \mytr \upchi Z f
		\right\rbrace
		 \label{E:RGEOLCASEFIRSTSTEPCOMMUTINGFUNCTIONSWITHLPLUSHALFTRACECHI}	
		 \\
	& \ \ 
		- \mathscr{Z}^{N-1}
			\left\lbrace
				\Lunit f
				+ \frac{1}{2} \mytr \upchi f	
			\right\rbrace
		+ \frac{1}{2}
			\mathscr{Z}^{N-1}  
			\left\lbrace
				(Z \mytr \upchi^{(Small)})
				f
			\right\rbrace.
			\notag
\end{align}
We now use the identities
\eqref{E:FIRSTSTEPCOMMUTINGFUNCTIONSWITHLPLUSHALFTRACECHI}
and
\eqref{E:RGEOLCASEFIRSTSTEPCOMMUTINGFUNCTIONSWITHLPLUSHALFTRACECHI}
to inductively commute the operators $\mathscr{Z}^M$
through the operator
$\left\lbrace 
			\Lunit 
			+ \frac{1}{2} \mytr \upchi 
		\right\rbrace.
$ 
Using also the Leibniz rule, Lemma~\ref{L:LANDRADCOMMUTEWITHANGDIFF}, 
the fact that $S \rgeo^{-1}$ is equal to either $\rgeo^{-2}$ or $0,$
and \eqref{E:ZNAPPLIEDTORGEOISNOTTOOLARGE},
we deduce that
\begin{align}  \label{E:FIRSTQUANTITATIVEESTIMATECOMMUTINGFUNCTIONSWITHLPLUSHALFTRACECHI}
	&
	\left|
			\mathscr{Z}^N
			\left\lbrace 
				\Lunit f
				+ \frac{1}{2} \mytr \upchi f
			\right\rbrace
		- \left\lbrace
				\Lunit + \frac{1}{2} \mytr \upchi 
			\right\rbrace
			\mathscr{Z}^N f
		\right|
			\\
	& \lesssim
		\left|
			\left\lbrace
				\Lunit + \frac{1}{2} \mytr \upchi 
			\right\rbrace
			\mathscr{Z}^{\leq N-1} f
		\right|
			+ \sum_{N_1 + N_2 \leq N-1}
				\sum_{S \in \mathscr{S}}
				\left|
					\angLie_{\mathscr{Z}}^{N_1} \angdeformoneformupsharparg{S}{\Lunit} 
				\right|
				\left| 
					\angdiff \mathscr{Z}^{N_2} f
				\right|
				\notag	\\
		& \ \ 
			+ \mathop{\sum_{N_1 + N_2 \leq N}}_{N_2 \leq N-1}
				\left|
					\mathscr{Z}^{N_1} \mytr \upchi^{(Small)}
				\right|
				\left| 
					\mathscr{Z}^{N_2} f
				\right|
			+ \frac{1}{(1 + t)^2}
				\left|
					\mathscr{Z}^{\leq N-1} f
				\right|.
				\notag
\end{align}
The desired bound \eqref{E:COMMUTINGFUNCTIONSWITHLPLUSHALFTRACECHI}
now follows from \eqref{E:FIRSTQUANTITATIVEESTIMATECOMMUTINGFUNCTIONSWITHLPLUSHALFTRACECHI},
inequality \eqref{E:FUNCTIONPOINTWISEANGDINTERMSOFANGLIEO},
the estimates 
\eqref{E:POINTWISEESTIMATESFORCHIJUNKINTERMSOFOTHERVARIABLES},
\eqref{E:RADDEFORMSPHERELPOINTWISE},
\eqref{E:LOWERORDERC0BOUNDRADDEFORMSPHEREL},
\eqref{E:ROTDEFORMSPHERELPOINTWISE},
\eqref{E:LOWERORDERC0BOUNDROTDEFORMSPHEREL},
and \eqref{E:C0BOUNDCRUCIALEIKONALFUNCTIONQUANTITIES},
and the bootstrap assumptions \eqref{E:PSIFUNDAMENTALC0BOUNDBOOTSTRAP}.

\end{proof}

\begin{corollary}[\textbf{Pointwise estimates involving commutations with the operator} $\Lunit + \frac{1}{2} \mytr \upchi$]
\label{C:COMMUTINGFUNCTIONSWITHLPLUSHALFTRACECHI}
Let $0 \leq N \leq 23$ be an integer and let $Z \in \mathscr{Z}.$
	Under the small-data and bootstrap assumptions 
	of Sects.~\ref{S:PSISOLVES}-\ref{S:C0BOUNDBOOTSTRAP},
	if $\varepsilon$ is sufficiently small,
	then the following pointwise estimates hold on $\mathcal{M}_{\Tboot,U_0}:$
	\begin{subequations}
	\begin{align} \label{E:PSIVERSIONCOMMUTINGFUNCTIONSWITHLPLUSHALFTRACECHI}
		\left|
			\mathscr{Z}^N
			\left\lbrace 
				\Lunit Z \Psi
				+ \frac{1}{2} \mytr \upchi \Psi
			\right\rbrace
		\right|
		& \lesssim
			\left|
				\left\lbrace
					\Lunit + \frac{1}{2} \mytr \upchi 
				\right\rbrace
				\mathscr{Z}^{\leq N+1} \Psi
			\right|
		+ \frac{\ln(\myexp + t)}{(1 + t)^2}
			\left| 
				\threemyarray
					[\Rad \mathscr{Z}^{\leq N} \Psi]
					{\rgeo \angdiff \mathscr{Z}^{\leq N} \Psi}
					{\mathscr{Z}^{\leq N} \Psi}
			\right|	
				\\
	& \ \ +   \varepsilon
				\frac{1}{(1 + t)^3}
				\left|	
					\myarray
						[\mathscr{Z}^{\leq N} (\upmu - 1)]
						{\sum_{a=1}^3 \rgeo |\mathscr{Z}^{\leq N} \Lunit_{(Small)}^a|} 
				\right|,
				\notag \\
		\left\|
			\mathscr{Z}^{\leq 11}
			\left\lbrace
				\Lunit Z \Psi
				+ \frac{1}{2} \mytr \upchi Z \Psi
			\right\rbrace
		\right \|_{C^0(\Sigma_t^u)}
		& \lesssim
			\varepsilon
			\frac{1}{(1 + t)^2}.
			\label{E:PSIVERSIONLOWERORDERC0BOUNDCOMMUTINGFUNCTIONSWITHLPLUSHALFTRACECHI}
\end{align}
\end{subequations}	
\end{corollary}

\begin{proof}
	We simply use Lemma~\ref{L:COMMUTINGFUNCTIONSWITHLPLUSHALFTRACECHI} with 
	$Z \Psi$ in the role of $f$ and the bootstrap assumptions \eqref{E:PSIFUNDAMENTALC0BOUNDBOOTSTRAP}.
\end{proof}

\begin{lemma}[\textbf{Commutator estimates for vectorfields acting on tensorfields}]
\label{L:COMMUTATORESTIMATESVECTORFIELDSACTINGONTENSORS}
Let $1 \leq N \leq 23$ be an integer, and let $\xi$ be an $S_{t,u}$ one-form 
or a symmetric type $\binom{0}{2}$ $S_{t,u}$ tensorfield.
Let $Z \in \mathscr{Z}$ be a commutation vectorfield.
Under the small-data and bootstrap assumptions 
of Sects.~\ref{S:PSISOLVES}-\ref{S:C0BOUNDBOOTSTRAP},
if $\varepsilon$ is sufficiently small,
then the following commutator estimates hold on $\mathcal{M}_{\Tboot,U_0}:$
\begin{subequations}
\begin{align}
	\left|
		\myarray[{[\angLie_{\mathscr{Z}}^N, \angLie_{\rgeo \Lunit}] \xi}]
			{{[\angLie_{\mathscr{Z}}^N, \angLie_{\Rad}] \xi}}
	\right|
	& \lesssim 
			\left|
				\angLie_{\Lunit} \angLie_{\mathscr{Z}}^{N-1} \xi
			\right|
			+
			\varepsilon^{1/2}
			\frac{\ln(\myexp + t)}{1 + t}
			\left| 
				\angLie_{\mathscr{Z}}^{\leq N} \xi
			\right|
			\label{E:RGEOLORRADZNCOMMUTATORACTINGONTENSORFIELDSPOINTWISE} \\
	& \ \ 
			+
			\left\|
				\angLie_{\mathscr{Z}}^{\leq \lfloor N/2 \rfloor} \xi
			\right\|_{C^0(\Sigma_t^u)}
			\left| 
				\fourmyarray[ \rgeo \Lunit \mathscr{Z}^{\leq N} \Psi]
					{\Rad \mathscr{Z}^{\leq N} \Psi}
					{\rgeo \angdiff \mathscr{Z}^{\leq N} \Psi}
					{\mathscr{Z}^{\leq N} \Psi}
			\right|
			\notag \\
	&  \ \ 
			+ 
			\frac{1}{1 + t}
			\left\|
				\angLie_{\mathscr{Z}}^{\leq \lfloor N/2 \rfloor} \xi
			\right\|_{C^0(\Sigma_t^u)}
			\left|
				\myarray[\mathscr{Z}^{\leq N+1} (\upmu - 1)]
					{\sum_{a=1}^3 \rgeo |\mathscr{Z}^{\leq N+1} \Lunit_{(Small)}^a|} 
			\right|,
		\notag \\
	\left|
		[\angLie_{\mathscr{Z}}^N, \angLie_{\Rot}] \xi
	\right|
	& \lesssim 
		 \left| 
				\angLie_{\mathscr{Z}}^{\leq N} \xi
			\right|
			+
			(1 + t)
			\left\|
				\angLie_{\mathscr{Z}}^{\leq \lfloor N/2 \rfloor} \xi
			\right\|_{C^0(\Sigma_t^u)}
			\left| 
				\fourmyarray[ \rgeo \Lunit \mathscr{Z}^{\leq N} \Psi]
					{\Rad \mathscr{Z}^{\leq N} \Psi}
					{\rgeo \angdiff \mathscr{Z}^{\leq N} \Psi}
					{\mathscr{Z}^{\leq N} \Psi}
			\right|
			\label{E:ROTZNCOMMUTATORACTINGONTENSORFIELDSPOINTWISE}  \\
	&  \ \ 
			+ 
			\left\|
				\angLie_{\mathscr{Z}}^{\leq \lfloor N/2 \rfloor} \xi
			\right\|_{C^0(\Sigma_t^u)}
			\left|
				\myarray[\mathscr{Z}^{\leq N+1} (\upmu - 1)]
					{\sum_{a=1}^3 \rgeo |\mathscr{Z}^{\leq N+1} \Lunit_{(Small)}^a|} 
			\right|,
		\notag \\	
	\left|
		[\angLie_{\mathscr{Z}}^N, \angLie_{\Lunit}] \xi
	\right|
	& \lesssim 
			\left|
				\angLie_{\Lunit} \angLie_{\mathscr{Z}}^{N-1} \xi
			\right|
			+
			\varepsilon^{1/2}
			\frac{\ln(\myexp + t)}{(1 + t)^2}
			\left| 
				\angLie_{\mathscr{Z}}^{\leq N} \xi
			\right|
			\label{E:LZNCOMMUTATORACTINGONTENSORFIELDSPOINTWISE} \\
	& \ \ 
			+ 
			\frac{1}{1 + t}
			\left\|
				\mathscr{Z}^{\leq \lfloor N/2 \rfloor} \xi
			\right\|_{C^0(\Sigma_t^u)}
			\left| 
				\fourmyarray[ \rgeo \Lunit \mathscr{Z}^{\leq N} \Psi]
					{\Rad \mathscr{Z}^{\leq N} \Psi}
					{\rgeo \angdiff \mathscr{Z}^{\leq N} \Psi}
					{\mathscr{Z}^{\leq N} \Psi}
			\right|
			\notag \\
	&  \ \ 
			+ 
			\frac{1}{(1 + t)^2}
			\left\|
				\mathscr{Z}^{\leq \lfloor N/2 \rfloor} \xi
			\right\|_{C^0(\Sigma_t^u)}
			\left|
				\myarray[\mathscr{Z}^{\leq N+1} (\upmu - 1)]
					{\sum_{a=1}^3 \rgeo |\mathscr{Z}^{\leq N+1} \Lunit_{(Small)}^a|} 
			\right|.
		\notag 
\end{align}
\end{subequations}

In addition,
\begin{align} \label{E:ORDERSWITCHEDLZNCOMMUTATORACTINGONTENSORFIELDSPOINTWISE} 
	\eqref{E:LZNCOMMUTATORACTINGONTENSORFIELDSPOINTWISE} 
	\mbox{\ also holds with the first term on the right replaced by}
	 	\left|
			 \angLie_{\mathscr{Z}}^{\leq N-1} \angLie_{\Lunit} \xi
		\right|.
\end{align}

In addition,
\begin{align} \label{E:LSNCOMMUTATORACTINGONTENSORFIELDSPOINTWISE} 
	& \eqref{E:LZNCOMMUTATORACTINGONTENSORFIELDSPOINTWISE} 
	\mbox{\ also holds without the first term on the right if the left-hand side is equal to \ }
		\\
	& \left|
		[\angLie_{\mathscr{S}}^N, \angLie_{\Lunit}] \xi
		\right|,
	\mbox{\ where }	
	\mathscr{S}^N 
	\mbox{\ is an $N^{th}$ order pure spatial commutation vectorfield operator}.
		\notag
\end{align}

\end{lemma}

\begin{proof}
	The proof of Lemma~\ref{L:COMMUTATORESTIMATESVECTORFIELDSACTINGONTENSORS} 
	is essentially the same as the proof of Lemma~\ref{L:COMMUTATORESTIMATESVECTORFIELDSACTINGONFUNCTIONS}.
	but we use the estimates 
	\eqref{E:PRELIMINARYQUANTITATIVERGEOLLIECOMMUTATIONSTENSORS}-\eqref{E:SPATIALANDLPRELIMINARYQUANTITATIVELIECOMMUTATIONS}
	in place of 
	\eqref{E:PRELIMINARYQUANTITATIVELIEDERIVATIVEBOUNDFORFUNCTIONS}
	(the former estimates involve one additional derivative of deformation tensors compared to the latter).
\end{proof}

\section{Commutator estimates for vectorfields acting on \texorpdfstring{$\angD \xi$}{the covariant angular derivative of a spherical tensorfield}}
In this section, we derive pointwise commutator estimates for vectorfields acting on 
$\angD \xi,$ where $\xi$ is a covariant $S_{t,u}$ tensorfield.

\begin{lemma}[\textbf{Commutator estimates for vectorfields acting on $\angD \xi,$ where $\xi$ is an $S_{t,u}$ tensorfield}]
\label{L:COMMUTATORESTIMATESVECTORFIELDSACTINGONANGDTENSORS}
Let $1 \leq N \leq 23$ be an integer, and let $\xi$ 
be an $S_{t,u}$ one-form
or a symmetric type $\binom{0}{2}$ $S_{t,u}$ tensorfield.
Under the small-data and bootstrap assumptions 
	of Sects.~\ref{S:PSISOLVES}-\ref{S:C0BOUNDBOOTSTRAP},
if $\varepsilon$ is sufficiently small,
then the following commutator estimates hold on $\mathcal{M}_{\Tboot,U_0}:$
\begin{align} \label{E:COMMUTATORESTIMATESVECTORFIELDSACTINGONANGDTENSORS}
	\left|
		[\angD, \angLie_{\mathscr{Z}}^N] \xi
	\right|
	& \lesssim 
			\frac{1}{1 + t}
			\left|
				\angLie_{\mathscr{Z}}^{\leq N-1} \xi
			\right|
			+ 
			\frac{1}{1 + t}
			\left\|
				\angLie_{\mathscr{Z}}^{\leq \lceil N/2 \rceil} \xi
			\right\|_{C^0(\Sigma_t^u)}
			\left| 
				\fourmyarray[\rgeo \Lunit \mathscr{Z}^{\leq N} \Psi]
					{\Rad \mathscr{Z}^{\leq N} \Psi}
					{\rgeo \angdiff \mathscr{Z}^{\leq N} \Psi}
					{\mathscr{Z}^{\leq N} \Psi}
			\right|
		\\
	&  \ \ 
			+ 
			\frac{1}{(1 + t)^2}
			\left\|
				\angLie_{\mathscr{Z}}^{\leq \lceil N/2 \rceil} \xi
			\right\|_{C^0(\Sigma_t^u)}
			\left|
				\myarray[\mathscr{Z}^{\leq N+1} (\upmu - 1)]
					{\sum_{a=1}^3 \rgeo |\mathscr{Z}^{\leq N+1} \Lunit_{(Small)}^a|} 
			\right|.
			\notag
\end{align}

\end{lemma}

\begin{proof}
	We first prove \eqref{E:COMMUTATORESTIMATESVECTORFIELDSACTINGONANGDTENSORS} in the case
	that $\xi$ is a symmetric type $\binom{0}{2}$ $S_{t,u}$ tensor.
	Let $Z \in \mathscr{Z}.$
	From \eqref{E:ANGDANGLIEZTYPE02COMMUTATOR}, we deduce	the schematic identity
	\begin{align}  \label{E:SCHEMATICRELATIONCOMMUTINGANGDANDLIEDERIVATIVESONTENSORS}
		[\angD, \angLie_Z] \xi
		& = \ginversesphere (\angD \angdeform{Z}) \xi.
	\end{align}
	Iterating \eqref{E:SCHEMATICRELATIONCOMMUTINGANGDANDLIEDERIVATIVESONTENSORS},
	using the schematic identity $\angLie_Z \ginversesphere = - (\ginversesphere)^2 \angdeform{Z}$
	(see Lemma~\ref{L:CONNECTIONBETWEENANGLIEOFGSPHEREANDDEFORMATIONTENSORS}),
	and using the fact that $\angdeform{Z}$ is also a symmetric
	type $\binom{0}{2}$ $S_{t,u}$ tensor in order to commute Lie derivatives with
	the operator $\angD$ in $\angD \angdeform{Z},$
	we deduce the schematic identity
	\begin{align} \label{E:ITERATEDSCHEMATICRELATIONCOMMUTINGANGDANDLIEDERIVATIVESONTENSORS}
		& [\angD, \angLie_{\mathscr{Z}}^N] \xi
			\\
	& = 
		\sum_{M=1}^N
		\sum_{N_1 + N_2 + \cdots + N_{M+2} = N-M} 
			(\angLie_{\mathscr{Z}}^{N_1} \ginversesphere)
			\cdots 
			(\angLie_{\mathscr{Z}}^{N_k} \ginversesphere)
			(\angD \angLie_{\mathscr{Z}}^{N_{k+1}} \angpi) 
			(\angLie_{\mathscr{Z}}^{N_{k+2}} \angpi)
			\cdots
			(\angLie_{\mathscr{Z}}^{N_{M+1}} \angpi)
			\angLie_{\mathscr{Z}}^{N_{M+2}} \xi,
			\notag
	\end{align}
	where the $\angpi$ are the
	$S_{t,u}$ projections of deformation tensors of vectorfields in $\mathscr{Z}.$
	Inequality \eqref{E:COMMUTATORESTIMATESVECTORFIELDSACTINGONANGDTENSORS}
	now follows from the identity \eqref{E:ITERATEDSCHEMATICRELATIONCOMMUTINGANGDANDLIEDERIVATIVESONTENSORS}
	and inequalities 
	\eqref{E:CRUDEPOINTWISEBOUNDSDERIVATIVESOFANGULARDEFORMATIONTENSORS},
	\eqref{E:CRUDELOWERORDERC0BOUNDDERIVATIVESOFANGULARDEFORMATIONTENSORS},
	and \eqref{E:TYPE02TENSORANGDINTERMSOFROTATIONALLIE}.
	
	Thanks to the identity \eqref{E:ANGDANGLIEZONEFORMCOMMUTATOR},
	the proof of \eqref{E:COMMUTATORESTIMATESVECTORFIELDSACTINGONANGDTENSORS} 
	in the case that $\xi$ is an $S_{t,u}$ one-form
	is essentially the same; we omit the details.

\end{proof}

\section{Commutator estimates for vectorfields acting on \texorpdfstring{$\angD^2 f$}{the angular Hessian of a function}}
In this section, we derive pointwise commutator estimates for vectorfields acting on 
$\angD^2 f,$ where $f$ is a function.

\begin{lemma}[\textbf{Commutator estimates for vectorfields acting on $\angD^2 f$}]
\label{L:ANGDSQUAREDLIEZNCOMMUTATORACTINGONFUNCTIONSSPOINTWISE} 
Let $1 \leq N \leq 23$ be an integer, and let $f$ be a function.
	Under the small-data and bootstrap assumptions 
	of Sects.~\ref{S:PSISOLVES}-\ref{S:C0BOUNDBOOTSTRAP},
	if $\varepsilon$ is sufficiently small,
	then the following commutator estimates hold on $\mathcal{M}_{\Tboot,U_0}:$
\begin{align} \label{E:ANGDSQUAREDLIEZNCOMMUTATORACTINGONFUNCTIONSSPOINTWISE} 
	\left|
		[\angD^2, \angLie_{\mathscr{Z}}^N] f
	\right|
	& \lesssim 
		\frac{1}{(1 + t)^2}
		\left|
			\mathscr{Z}^{\leq N} f
		\right|
		+ \frac{1}{(1 + t)^2}
			\left\|
				\mathscr{Z}^{\leq \lceil N/2 \rceil} f
			\right\|_{C^0(\Sigma_t^u)}
			\left| 
				\fourmyarray[ \rgeo \Lunit \mathscr{Z}^{\leq N} \Psi]
					{\Rad \mathscr{Z}^{\leq N} \Psi}
					{\rgeo \angdiff \mathscr{Z}^{\leq N} \Psi}
					{\mathscr{Z}^{\leq N} \Psi}
			\right|	
					\\
	& \ \ 
		+ \frac{1}{(1 + t)^3}
			\left\|
				\mathscr{Z}^{\leq \lceil N/2 \rceil} f
			\right\|_{C^0(\Sigma_t^u)}
			\left|
				\myarray[\mathscr{Z}^{\leq N+1} (\upmu - 1)]
					{\sum_{a=1}^3 \rgeo |\mathscr{Z}^{\leq N+1} \Lunit_{(Small)}^a|} 
			\right|,
			\notag
\end{align}

\begin{align} \label{E:ANGFREEDSQUAREDLIEZNCOMMUTATORACTINGONFUNCTIONSSPOINTWISE} 
	\left|
		([\angfreeDsquared, \angLie_{\mathscr{Z}}^N] f)
	\right|
	& \lesssim 
		\frac{1}{(1 + t)^2}
		\left|
			\mathscr{Z}^{\leq N+1} f
		\right|
		+ \frac{1}{(1 + t)^2}
			\left\|
				\mathscr{Z}^{\leq \lceil N/2 \rceil} f
			\right\|_{C^0(\Sigma_t^u)}
			\left| 
				\fourmyarray[ \rgeo \Lunit \mathscr{Z}^{\leq N} \Psi]
					{\Rad \mathscr{Z}^{\leq N} \Psi}
					{\rgeo \angdiff \mathscr{Z}^{\leq N} \Psi}
					{\mathscr{Z}^{\leq N} \Psi}
			\right|	
					\\
	& \ \ 
		+ \frac{1}{(1 + t)^3}
			\left\|
				\mathscr{Z}^{\leq \lceil N/2 \rceil} f
			\right\|_{C^0(\Sigma_t^u)}
			\left|
				\myarray[\mathscr{Z}^{\leq N+1} (\upmu - 1)]
					{\sum_{a=1}^3 \rgeo |\mathscr{Z}^{\leq N+1} \Lunit_{(Small)}^a|} 
			\right|,
			\notag
\end{align}

\begin{align} \label{E:ANGLAPZNCOMMUTATORACTINGONFUNCTIONSSPOINTWISE} 
	\left|
		([\angLap, \mathscr{Z}^N] f)
	\right|
	& \lesssim 
		\frac{1}{(1 + t)^2}
		\left|
			\mathscr{Z}^{\leq N+1} f
		\right|
		+ \frac{1}{(1 + t)^2}
			\left\|
				\mathscr{Z}^{\leq \lceil N/2 \rceil} f
			\right\|_{C^0(\Sigma_t^u)}
			\left| 
				\fourmyarray[ \rgeo \Lunit \mathscr{Z}^{\leq N} \Psi]
					{\Rad \mathscr{Z}^{\leq N} \Psi}
					{\rgeo \angdiff \mathscr{Z}^{\leq N} \Psi}
					{\mathscr{Z}^{\leq N} \Psi}
			\right|	
					\\
	& \ \ 
		+ \frac{1}{(1 + t)^3}
			\left\|
				\mathscr{Z}^{\leq \lceil N/2 \rceil} f
			\right\|_{C^0(\Sigma_t^u)}
			\left|
				\myarray[\mathscr{Z}^{\leq N+1} (\upmu - 1)]
					{\sum_{a=1}^3 \rgeo |\mathscr{Z}^{\leq N+1} \Lunit_{(Small)}^a|} 
			\right|.
			\notag
\end{align}

\end{lemma}

\begin{proof}
	We first prove \eqref{E:ANGDSQUAREDLIEZNCOMMUTATORACTINGONFUNCTIONSSPOINTWISE}.
	To this end, let $Z \in \mathscr{Z}.$
	From \eqref{E:COMMUTINGANGDSQUAREDANDLIEZ}, we deduce	the schematic identity
	\begin{align}  \label{E:ANGDSQUAREDSCHEMATICRELATIONCOMMUTINGANGDANDLIEDERIVATIVESONTENSORS}
		[\angD^2, \angLie_Z] f
		& = \ginversesphere (\angD \angdeform{Z}) \angdiff f.
	\end{align}
	Iterating \eqref{E:ANGDSQUAREDSCHEMATICRELATIONCOMMUTINGANGDANDLIEDERIVATIVESONTENSORS},
	using the schematic identity $\angLie_Z \ginversesphere = - (\ginversesphere)^2 \angdeform{Z}$
	(see Lemma~\ref{L:CONNECTIONBETWEENANGLIEOFGSPHEREANDDEFORMATIONTENSORS}),
	commuting Lie derivatives with the operator $\angD$ in $\angD \angdeform{Z}$
	as in our proof of \eqref{E:ITERATEDSCHEMATICRELATIONCOMMUTINGANGDANDLIEDERIVATIVESONTENSORS},
	and using Lemma~\ref{L:LANDRADCOMMUTEWITHANGDIFF},
	we deduce the schematic identity
	\begin{align} \label{E:ITERATEDANGDSQUAREDSCHEMATICRELATIONCOMMUTINGANGDANDLIEDERIVATIVESONTENSORS}
		& [\angD^2, \angLie_{\mathscr{Z}}^N] f
			\\
		& = 
		\sum_{M=1}^N
		\sum_{N_1 + N_2 + \cdots + N_{M+2} = N-M} 
			(\angLie_{\mathscr{Z}}^{N_1} \ginversesphere)
			\cdots 
			(\angLie_{\mathscr{Z}}^{N_k} \ginversesphere)
			(\angD \angLie_{\mathscr{Z}}^{N_{k+1}} \angpi) 
			(\angLie_{\mathscr{Z}}^{N_{k+2}} \angpi)
			\cdots
			(\angLie_{\mathscr{Z}}^{N_{M+1}} \angpi)
			\angdiff \angLie_{\mathscr{Z}}^{N_{M+2}} f,
			\notag
	\end{align}
	where the $\angpi$ are 
	the $S_{t,u}$ projections of	deformation tensors of vectorfields in $\mathscr{Z}.$
	Inequality \eqref{E:ANGDSQUAREDLIEZNCOMMUTATORACTINGONFUNCTIONSSPOINTWISE}
	now follows from the identity \eqref{E:ITERATEDANGDSQUAREDSCHEMATICRELATIONCOMMUTINGANGDANDLIEDERIVATIVESONTENSORS}
	and inequalities 
	\eqref{E:CRUDEPOINTWISEBOUNDSDERIVATIVESOFANGULARDEFORMATIONTENSORS},
	\eqref{E:CRUDELOWERORDERC0BOUNDDERIVATIVESOFANGULARDEFORMATIONTENSORS},
	\eqref{E:FUNCTIONPOINTWISEANGDINTERMSOFANGLIEO} ,
	and
	\eqref{E:TYPE02TENSORANGDINTERMSOFROTATIONALLIE}.
	
	Inequality \eqref{E:ANGLAPZNCOMMUTATORACTINGONFUNCTIONSSPOINTWISE}
	can be proved in a similar fashion 
	with the help of
	\eqref{E:COMMUTINGANGANGLAPANDLIEZ} and 
	\eqref{E:ANGDSQUAREDFUNCTIONPOINTWISEINTERMSOFANGDIFFROTATIONS}.
	Note that the right-hand side of
	\eqref{E:COMMUTINGANGANGLAPANDLIEZ}
	depends on one extra derivative of 
	$f$ compared to \eqref{E:COMMUTINGANGDSQUAREDANDLIEZ}
	and hence the same is true for 
	\eqref{E:ANGLAPZNCOMMUTATORACTINGONFUNCTIONSSPOINTWISE}
	compared to \eqref{E:ANGDSQUAREDLIEZNCOMMUTATORACTINGONFUNCTIONSSPOINTWISE}.
	
	To prove \eqref{E:ANGFREEDSQUAREDLIEZNCOMMUTATORACTINGONFUNCTIONSSPOINTWISE}, we
	first use the decomposition 
	$\angfreeDsquared f 
	= \angD^2 f 
		- \frac{1}{2}(\angLap f) \gsphere$
	and the Leibniz rule to deduce the following identity:
	\begin{align} \label{E:FIRSTPOINTWISEANGFREEDSQUAREDLIEZNCOMMUTATORACTINGONFUNCTIONSSPOINTWISE}
		[\angfreeDsquared, \angLie_{\mathscr{Z}}^N] f
		& = [\angD^2, \angLie_{\mathscr{Z}}^N] f
			- \frac{1}{2}([\angLap,  \angLie_{\mathscr{Z}}^N] f) \gsphere
			+ \mathop{\sum_{N_1 + N_2 \leq N}}_{N_1 \leq N-1}
					(\angLap \angLie_{\mathscr{Z}}^{N_1} f)
					\angLie_{\mathscr{Z}}^{N_2} \gsphere
			+ \mathop{\sum_{N_1 + N_2 \leq N}}_{N_1 \leq N-1}
					([\angLie_{\mathscr{Z}}^{N_1}, \angLap] f)
					\angLie_{\mathscr{Z}}^{N_2} \gsphere.
	\end{align}	
	To bound the magnitude first term on the right-hand side of 	
	\eqref{E:FIRSTPOINTWISEANGFREEDSQUAREDLIEZNCOMMUTATORACTINGONFUNCTIONSSPOINTWISE}
	by the right-hand side of \eqref{E:ANGFREEDSQUAREDLIEZNCOMMUTATORACTINGONFUNCTIONSSPOINTWISE},
	we simply use the estimate \eqref{E:ANGDSQUAREDLIEZNCOMMUTATORACTINGONFUNCTIONSSPOINTWISE}.
	Similarly, to bound the second term, we use \eqref{E:ANGLAPZNCOMMUTATORACTINGONFUNCTIONSSPOINTWISE}.
	To bound the magnitude of the third term on the right-hand side of 
	\eqref{E:FIRSTPOINTWISEANGFREEDSQUAREDLIEZNCOMMUTATORACTINGONFUNCTIONSSPOINTWISE}
	by the right-hand side of \eqref{E:ANGFREEDSQUAREDLIEZNCOMMUTATORACTINGONFUNCTIONSSPOINTWISE},
	we use 
	\eqref{E:CRUDEPOINTWISEBOUNDSDERIVATIVESOFANGULARDEFORMATIONTENSORS},
	\eqref{E:CRUDELOWERORDERC0BOUNDDERIVATIVESOFANGULARDEFORMATIONTENSORS},
	\eqref{E:FUNCTIONPOINTWISEANGDINTERMSOFANGLIEO},
	and
	\eqref{E:ANGLAPFUNCTIONPOINTWISEINTERMSOFROTATIONS}.
	To bound the magnitude of the last term on the right-hand side of 
	\eqref{E:FIRSTPOINTWISEANGFREEDSQUAREDLIEZNCOMMUTATORACTINGONFUNCTIONSSPOINTWISE}
	by the right-hand side of \eqref{E:ANGFREEDSQUAREDLIEZNCOMMUTATORACTINGONFUNCTIONSSPOINTWISE},
	we use 
	\eqref{E:CRUDEPOINTWISEBOUNDSDERIVATIVESOFANGULARDEFORMATIONTENSORS},
	\eqref{E:CRUDELOWERORDERC0BOUNDDERIVATIVESOFANGULARDEFORMATIONTENSORS},
	\eqref{E:ANGLAPZNCOMMUTATORACTINGONFUNCTIONSSPOINTWISE},
	and the bootstrap assumptions \eqref{E:PSIFUNDAMENTALC0BOUNDBOOTSTRAP}.
	
\end{proof}

\begin{corollary} [\textbf{Commutator estimates for vectorfields acting on $\angD^2 \Psi$}]
\label{C:POINTWISECOMMUTINGANGDSQUAREDPSIWITHLIEZN}
Let $0 \leq N \leq 23$ be an integer.
Under the small-data and bootstrap assumptions 
of Sects.~\ref{S:PSISOLVES}-\ref{S:C0BOUNDBOOTSTRAP},
if $\varepsilon$ is sufficiently small,
then the following estimates hold on $\mathcal{M}_{\Tboot,U_0}:$	
\begin{subequations}
\begin{align}
		\label{E:POINTWISECOMMUTINGANGDSQUAREDPSIWITHLIEZN}
		\left|
			\threemyarray
				[\angLie_{\mathscr{Z}}^N \angD^2 \Psi]
				{\mathscr{Z}^N \angLap \Psi}
				{\angLie_{\mathscr{Z}}^N \angfreeDsquared \Psi}
			\right|
		& 
		\lesssim 
		\frac{1}{(1 + t)^2}
		\left| 
			\fourmyarray[ \rgeo \Lunit \mathscr{Z}^{\leq N+1} \Psi]
				{\Rad \mathscr{Z}^{\leq N+1} \Psi}
				{\rgeo \angdiff \mathscr{Z}^{\leq N+1} \Psi}
				{\mathscr{Z}^{\leq N+1} \Psi}
			\right|	
		+ \varepsilon
			\frac{1}{(1 + t)^4}
			\left|
				\myarray[\mathscr{Z}^{\leq N+1} (\upmu - 1)]
					{\sum_{a=1}^3 \rgeo |\mathscr{Z}^{\leq N+1} \Lunit_{(Small)}^a|} 
			\right|,
			\\
		\left|
			\threemyarray
				[\angLie_{\mathscr{Z}}^{\leq 11} \angD^2 \Psi]
				{\mathscr{Z}^{\leq 11} \angLap \Psi}
				{\angLie_{\mathscr{Z}}^{\leq 1} \angfreeDsquared \Psi}
			\right\|_{C^0(\Sigma_t^u)}
		& \lesssim 
			\varepsilon
			\frac{1}{(1 + t)^3}.
			\label{E:LINFFTYLOWERORDERCOMMUTINGANGDSQUAREDPSIWITHLIEZN}
	\end{align}
	\end{subequations}
\end{corollary}	

\begin{proof}
The estimates \eqref{E:POINTWISECOMMUTINGANGDSQUAREDPSIWITHLIEZN}
and \eqref{E:LINFFTYLOWERORDERCOMMUTINGANGDSQUAREDPSIWITHLIEZN}
follow easily from Lemma~\ref{L:ANGDSQUAREDLIEZNCOMMUTATORACTINGONFUNCTIONSSPOINTWISE} 
with $\Psi$ in the role of $f$
and the bootstrap assumptions \eqref{E:PSIFUNDAMENTALC0BOUNDBOOTSTRAP}.
\end{proof}

\section{Commutator estimates involving the trace and trace-free parts}
In this section, we derive pointwise commutator estimates involving the trace and trace-free parts of
symmetric type $\binom{0}{2}$ $S_{t,u}$ tensorfields.

\begin{lemma}[\textbf{Commutator estimates involving the trace and trace-free parts of symmetric type $\binom{0}{2}$ $S_{t,u}$ tensorfields}]
Let $0 \leq N \leq 24$ be an integer, and let $\xi$ 
be a symmetric type $\binom{0}{2}$ $S_{t,u}$ tensorfield.	
Under the small-data and bootstrap assumptions 
of Sects.~\ref{S:PSISOLVES}-\ref{S:C0BOUNDBOOTSTRAP},
if $\varepsilon$ is sufficiently small,
then the following pointwise estimates hold on $\mathcal{M}_{\Tboot,U_0}:$
\begin{align} \label{E:LIEDERIVATIVESOFTRACEFREEPARTINTERMSOFLIEDERIVATIVESOFFULLTENSOR}
	\left|
		\angLie_{\mathscr{Z}}^N \hat{\xi}
	\right|
	& \lesssim 
		\left|
			\angLie_{\mathscr{Z}}^{\leq N} \xi
		\right|
		+
		\left\|
			\angLie_{\mathscr{Z}}^{\leq \lfloor N/2 \rfloor} \xi
		\right\|_{C^0(\Sigma_t^u)}
		\left|
				\fourmyarray[\rgeo \Lunit \mathscr{Z}^{\leq N-1} \Psi]
					{\Rad \mathscr{Z}^{\leq N-1} \Psi}
					{\rgeo \angdiff \mathscr{Z}^{\leq N-1} \Psi} 
					{\mathscr{Z}^{\leq N-1} \Psi}
			\right|
				\\
	& \ \ 
			+ 
			\frac{1}{1 + t}
			\left\|
				\angLie_{\mathscr{Z}}^{\leq \lfloor N/2 \rfloor} \xi
			\right\|_{C^0(\Sigma_t^u)}
			\left| 
				\myarray
					[\mathscr{Z}^{\leq N} (\upmu - 1)]
						{\rgeo \sum_{a=1}^3 |\mathscr{Z}^{\leq N} \Lunit_{(Small)}^a|}
		\right|.
		\notag
\end{align}

Furthermore, if $1 \leq N \leq 24$ is an integer, 
then the following pointwise estimates hold on $\mathcal{M}_{\Tboot,U_0}:$
\begin{align} \label{E:DIFFERENCEBETWEENTRACEFREELIEDERIVATIVESANDLIEDERIVATIVESACTINGONATRACEFREETENSOR}
	\left|
		\angLie_{\mathscr{Z}}^N \hat{\xi}
		- \angfreeLietwoarg{\mathscr{Z}}{N} \hat{\xi}
	\right|
	& \lesssim 
		\left|
			\angLie_{\mathscr{Z}}^{\leq N-1} \xi
		\right|
		+
		\left\|
				\angLie_{\mathscr{Z}}^{\leq \lfloor N/2 \rfloor} \xi
		\right\|_{C^0(\Sigma_t^u)}
		\left|
				\fourmyarray[\rgeo \Lunit \mathscr{Z}^{\leq N-1} \Psi]
					{\Rad \mathscr{Z}^{\leq N-1} \Psi}
					{\rgeo \angdiff \mathscr{Z}^{\leq N-1} \Psi} 
					{\mathscr{Z}^{\leq N-1} \Psi}
			\right|
				\\
	& \ \ + \frac{1}{1 + t}
			\left\|
				\angLie_{\mathscr{Z}}^{\leq \lfloor N/2 \rfloor} \xi
			\right\|_{C^0(\Sigma_t^u)}
			\left| 
				\myarray
					[\mathscr{Z}^{\leq N} (\upmu - 1)]
						{\rgeo \sum_{a=1}^3 |\mathscr{Z}^{\leq N} \Lunit_{(Small)}^a|}
			\right|,
			\notag
\end{align}

\begin{align} \label{E:COMMUTINGTRACEWITHLIEDERIVATIVES}
		\left|
			\mytr  \angLie_{\mathscr{Z}}^N \xi
			- \mathscr{Z}^N \mytr \xi
		\right|	
		& \lesssim 
		\left|
			\angLie_{\mathscr{Z}}^{\leq N-1} \xi
		\right|
		+
		\left\|
				\angLie_{\mathscr{Z}}^{\leq \lfloor N/2 \rfloor} \xi
		\right\|_{C^0(\Sigma_t^u)}
		\left|
				\fourmyarray[\rgeo \Lunit \mathscr{Z}^{\leq N-1} \Psi]
					{\Rad \mathscr{Z}^{\leq N-1} \Psi}
					{\rgeo \angdiff \mathscr{Z}^{\leq N-1} \Psi} 
					{\mathscr{Z}^{\leq N-1} \Psi}
			\right|
				\\
	& \ \ + \frac{1}{1 + t}
			\left\|
				\angLie_{\mathscr{Z}}^{\leq \lfloor N/2 \rfloor} \xi
			\right\|_{C^0(\Sigma_t^u)}
			\left| 
				\myarray
					[\mathscr{Z}^{\leq N} (\upmu - 1)]
						{\rgeo \sum_{a=1}^3 |\mathscr{Z}^{\leq N} \Lunit_{(Small)}^a|}
				\right|,
			\notag
	\end{align}
	
	\begin{align} \label{E:COMMUTINGTRACEFREEWITHLIEDERIVATIVES}
		\left|
			\Lie_{\mathscr{Z}}^N \hat{\xi}
			- \widehat{(\angLie_{\mathscr{Z}}^N \xi)}
		\right|	
		& \lesssim 
			\left|
				\angLie_{\mathscr{Z}}^{\leq N-1} \xi
			\right|
			+
		\left\|
			\angLie_{\mathscr{Z}}^{\leq \lfloor N/2 \rfloor} \xi
		\right\|_{C^0(\Sigma_t^u)}
		\left|
				\fourmyarray[\rgeo \Lunit \mathscr{Z}^{\leq N-1} \Psi]
					{\Rad \mathscr{Z}^{\leq N-1} \Psi}
					{\rgeo \angdiff \mathscr{Z}^{\leq N-1} \Psi} 
					{\mathscr{Z}^{\leq N-1} \Psi}
			\right|
				\\
	& \ \ + \frac{1}{1 + t}
			\left\|
				\angLie_{\mathscr{Z}}^{\leq \lfloor N/2 \rfloor} \xi
			\right\|_{C^0(\Sigma_t^u)}
			\left| 
				\myarray
					[\mathscr{Z}^{\leq N} (\upmu - 1)]
						{\rgeo \sum_{a=1}^3 |\mathscr{Z}^{\leq N} \Lunit_{(Small)}^a|}
				\right|,
			\notag
	\end{align}
	
	\begin{align}
		\left|
			\angLie_{\mathscr{Z}}^N \xi
				-
			\left\lbrace
				\angfreeLietwoarg{\mathscr{Z}}{N} \hat{\xi}
				+ \frac{1}{2} (\mathscr{Z}^N \mytr \xi) \gsphere
			\right\rbrace
		\right|
		& \lesssim 
			\left|
				\angLie_{\mathscr{Z}}^{\leq N-1} \xi
			\right|
			+
		\left\|
			\angLie_{\mathscr{Z}}^{\leq \lfloor N/2 \rfloor} \xi
		\right\|_{C^0(\Sigma_t^u)}
		\left|
				\fourmyarray[\rgeo \Lunit \mathscr{Z}^{\leq N-1} \Psi]
					{\Rad \mathscr{Z}^{\leq N-1} \Psi}
					{\rgeo \angdiff \mathscr{Z}^{\leq N-1} \Psi} 
					{\mathscr{Z}^{\leq N-1} \Psi}
			\right|
				\label{E:COMMUTINGANGLIEZANDTRACETRACEFREESPLITTING} \\
	& \ \ + \frac{1}{1 + t}
			\left\|
				\angLie_{\mathscr{Z}}^{\leq \lfloor N/2 \rfloor} \xi
			\right\|_{C^0(\Sigma_t^u)}
			\left| 
				\myarray
					[\mathscr{Z}^{\leq N} (\upmu - 1)]
						{\rgeo \sum_{a=1}^3 |\mathscr{Z}^{\leq N} \Lunit_{(Small)}^a|}
			\right|.
			\notag
	\end{align}
	
	Furthermore, if $0 \leq N \leq 23$ is an integer,
	then the following estimates hold on $\mathcal{M}_{\Tboot,U_0}:$
	\begin{align} \label{E:COMMUTINGANGDANDTRACETRACEFREESPLITTING}
		\left|
			\angD \angLie_{\mathscr{Z}}^N \xi
			-
			\left\lbrace
				\angD \angfreeLietwoarg{\mathscr{Z}}{N} \hat{\xi}
				+ \frac{1}{2} (\angdiff \mathscr{Z}^N \mytr \xi) \gsphere
			\right\rbrace
		\right|
		& \lesssim 
		\frac{1}{1 + t}	
		\left|
			\angLie_{\mathscr{Z}}^{\leq N} \xi
		\right|
		+
		\frac{1}{1 + t}	
		\left\|
			\angLie_{\mathscr{Z}}^{\leq \lfloor N/2 \rfloor} \xi
		\right\|_{C^0(\Sigma_t^u)}
			\left|
				\fourmyarray[\rgeo \Lunit \mathscr{Z}^{\leq N} \Psi]
					{\Rad \mathscr{Z}^{\leq N} \Psi}
					{\rgeo \angdiff \mathscr{Z}^{\leq N} \Psi} 
					{\mathscr{Z}^{\leq N} \Psi}
			\right|
				\\
	& \ \ + \frac{1}{(1 + t)^2}
			\left\|
				\angLie_{\mathscr{Z}}^{\leq \lfloor N/2 \rfloor} \xi
			\right\|_{C^0(\Sigma_t^u)}
			\left| 
				\myarray
					[\mathscr{Z}^{\leq N+1} (\upmu - 1)]
						{\rgeo \sum_{a=1}^3 |\mathscr{Z}^{\leq N+1} \Lunit_{(Small)}^a|}
			\right|,
			\notag
	\end{align}
	
	\begin{align} \label{E:COMMUTINGTRACEFREELIEDERIVATIVESTHROUGHANGD}
		\left|
			\angD \angLie_{\mathscr{Z}}^N \hat{\xi}
			- \angD \angfreeLietwoarg{\mathscr{Z}}{N} \hat{\xi}
		\right|
		& \lesssim 
		\frac{1}{1 + t}	
		\left|
			\angLie_{\mathscr{Z}}^{\leq N} \xi
		\right|
		+
		\frac{1}{1 + t}	
		\left\|
			\angLie_{\mathscr{Z}}^{\leq \lfloor N/2 \rfloor} \xi
		\right\|_{C^0(\Sigma_t^u)}
			\left|
				\fourmyarray[\rgeo \Lunit \mathscr{Z}^{\leq N} \Psi]
					{\Rad \mathscr{Z}^{\leq N} \Psi}
					{\rgeo \angdiff \mathscr{Z}^{\leq N} \Psi} 
					{\mathscr{Z}^{\leq N} \Psi}
			\right|
				\\
	& \ \ + \frac{1}{(1 + t)^2}
			\left\|
				\angLie_{\mathscr{Z}}^{\leq \lfloor N/2 \rfloor} \xi
			\right\|_{C^0(\Sigma_t^u)}
			\left| 
				\myarray
					[\mathscr{Z}^{\leq N+1} (\upmu - 1)]
						{\rgeo \sum_{a=1}^3 |\mathscr{Z}^{\leq N+1} \Lunit_{(Small)}^a|}
			\right|,
			\notag
	\end{align}

\begin{align} \label{E:COMMUTINGLIEDERIVATIVESTHROUGHANGDIV}
		\left|
			\angLie_{\mathscr{Z}}^N \angdiv \xi
			-
			\angdiv \angLie_{\mathscr{Z}}^N \xi
		\right|	
		& \lesssim 
		\frac{1}{1 + t}	
		\left|
			\angLie_{\mathscr{Z}}^{\leq N} \xi
		\right|
		+
		\frac{1}{1 + t}	
		\left\|
			\angLie_{\mathscr{Z}}^{\leq \lfloor N/2 \rfloor} \xi
		\right\|_{C^0(\Sigma_t^u)}
		\left|
				\fourmyarray[\rgeo \Lunit \mathscr{Z}^{\leq N} \Psi]
					{\Rad \mathscr{Z}^{\leq N} \Psi}
					{\rgeo \angdiff \mathscr{Z}^{\leq N} \Psi} 
					{\mathscr{Z}^{\leq N} \Psi}
			\right|
				\\
	& \ \ + \frac{1}{(1 + t)^2}
			\left\|
				\angLie_{\mathscr{Z}}^{\leq \lfloor N/2 \rfloor} \xi
			\right\|_{C^0(\Sigma_t^u)}
			\left| 
				\myarray
					[\mathscr{Z}^{\leq N+1} (\upmu - 1)]
						{\rgeo \sum_{a=1}^3 |\mathscr{Z}^{\leq N+1} \Lunit_{(Small)}^a|}
				\right|,
			\notag
	\end{align}	
	
\begin{align} \label{E:COMMUTINGTRACEFREELIEDERIVATIVESTHROUGHANGDIV}
		\left|
			\angLie_{\mathscr{Z}}^N \angdiv \hat{\xi}
			-
			\angdiv \angfreeLietwoarg{\mathscr{Z}}{N} \hat{\xi}
		\right|	
		& \lesssim 
		\frac{1}{1 + t}	
		\left|
			\angLie_{\mathscr{Z}}^{\leq N} \xi
		\right|
		+
		\frac{1}{1 + t}	
		\left\|
			\angLie_{\mathscr{Z}}^{\leq \lfloor N/2 \rfloor} \xi
		\right\|_{C^0(\Sigma_t^u)}
		\left|
				\fourmyarray[\rgeo \Lunit \mathscr{Z}^{\leq N} \Psi]
					{\Rad \mathscr{Z}^{\leq N} \Psi}
					{\rgeo \angdiff \mathscr{Z}^{\leq N} \Psi} 
					{\mathscr{Z}^{\leq N} \Psi}
			\right|
				\\
	& \ \ + \frac{1}{(1 + t)^2}
			\left\|
				\angLie_{\mathscr{Z}}^{\leq \lfloor N/2 \rfloor} \xi
			\right\|_{C^0(\Sigma_t^u)}
			\left| 
				\myarray
					[\mathscr{Z}^{\leq N+1} (\upmu - 1)]
						{\rgeo \sum_{a=1}^3 |\mathscr{Z}^{\leq N+1} \Lunit_{(Small)}^a|}
				\right|.
			\notag
	\end{align}
	
\end{lemma}

\begin{proof}
To prove \eqref{E:LIEDERIVATIVESOFTRACEFREEPARTINTERMSOFLIEDERIVATIVESOFFULLTENSOR}, we 
first decompose $\hat{\xi}_{AB} = \xi_{AB} - \frac{1}{2} (\ginversesphere)^{CD} \xi_{CD} \gsphere_{AB}.$
Applying $\angLie_{\mathscr{Z}}^N$ to both sides of this identity and applying the Leibniz rule to the
product on the right-hand side, we deduce that
\begin{align} \label{E:FIRSTPOINTWISEESTIMATELIEDERIVATIVESOFTRACEFREEPARTINTERMSOFLIEDERIVATIVESOFFULLTENSOR}
	\left|
		\angLie_{\mathscr{Z}}^N \hat{\xi}
	\right|
	& \lesssim 
		\left|
			\angLie_{\mathscr{Z}}^N \xi
		\right|
		+
		\sum_{N_1 + N_2 + N_3 = N}
			\left|
				\angLie_{\mathscr{Z}}^{N_1} \xi
			\right|
			\left|
				\angLie_{\mathscr{Z}}^{N_2} \gsphere
			\right|
			\left|
				\angLie_{\mathscr{Z}}^{N_3} \ginversesphere
			\right|.
\end{align}	
The desired estimate \eqref{E:LIEDERIVATIVESOFTRACEFREEPARTINTERMSOFLIEDERIVATIVESOFFULLTENSOR} now follows from
\eqref{E:FIRSTPOINTWISEESTIMATELIEDERIVATIVESOFTRACEFREEPARTINTERMSOFLIEDERIVATIVESOFFULLTENSOR}
and the estimates
\eqref{E:CRUDEPOINTWISEBOUNDSDERIVATIVESOFANGULARDEFORMATIONTENSORS}
and
\eqref{E:CRUDELOWERORDERC0BOUNDDERIVATIVESOFANGULARDEFORMATIONTENSORS}.

To prove \eqref{E:COMMUTINGTRACEWITHLIEDERIVATIVES}, we first deduce from the Leibniz rule that\begin{align} \label{E:FIRSTPOINTWISECOMMUTINGTRACEWITHLIEDERIVATIVES}
	\left|
		\mytr  \angLie_{\mathscr{Z}}^N \xi
		- \mathscr{Z}^N \mytr \xi
	\right|
	& \lesssim
	\mathop{\sum_{N_1 + N_2 = N}}_{N_2 \leq N-1}
			\left|
				\angLie_{\mathscr{Z}}^{N_1} \gsphere
			\right|
			\left|
				\angLie_{\mathscr{Z}}^{N_2} \xi
			\right|.
\end{align}
The desired estimate \eqref{E:COMMUTINGTRACEWITHLIEDERIVATIVES}
now follows from
\eqref{E:FIRSTPOINTWISECOMMUTINGTRACEWITHLIEDERIVATIVES},
\eqref{E:CRUDEPOINTWISEBOUNDSDERIVATIVESOFANGULARDEFORMATIONTENSORS}
and
\eqref{E:CRUDELOWERORDERC0BOUNDDERIVATIVESOFANGULARDEFORMATIONTENSORS}.

To prove \eqref{E:COMMUTINGTRACEFREEWITHLIEDERIVATIVES},
we first decompose $\hat{\xi}_{AB} = \xi_{AB} - \frac{1}{2} (\ginversesphere)^{CD} \xi_{CD} \gsphere_{AB}.$
Applying $\angLie_{\mathscr{Z}}^N$ to both sides of this identity
and applying the Leibniz rule to the
product on the right-hand side,
we deduce that
\begin{align} \label{E:FIRSTPOINTWISETRACELIEZNOFTRACEFREEISLOWERORDER}
	\left|
		\angLie_{\mathscr{Z}}^N \hat{\xi}
		- \left\lbrace
				\angLie_{\mathscr{Z}}^N \xi
				- \frac{1}{2} (\mytr  \angLie_{\mathscr{Z}}^N \xi) \gsphere
			\right\rbrace
	\right|
	& \lesssim 
	\mathop{\sum_{N_1 + N_2 + N_3 = N}}_{N_1 \leq N-1}
			\left|
				\angLie_{\mathscr{Z}}^{N_1} \xi
			\right|
			\left|
				\angLie_{\mathscr{Z}}^{N_2} \gsphere
			\right|
			\left|
				\angLie_{\mathscr{Z}}^{N_3} \ginversesphere
			\right|.
\end{align}
We now conclude the desired estimate \eqref{E:COMMUTINGTRACEFREEWITHLIEDERIVATIVES}
by noting that the term in braces on the left-hand side of
\eqref{E:FIRSTPOINTWISETRACELIEZNOFTRACEFREEISLOWERORDER}
is precisely $\widehat{(\angLie_{\mathscr{Z}}^N \xi)}$
and by
using the estimates
\eqref{E:CRUDEPOINTWISEBOUNDSDERIVATIVESOFANGULARDEFORMATIONTENSORS}
and
\eqref{E:CRUDELOWERORDERC0BOUNDDERIVATIVESOFANGULARDEFORMATIONTENSORS}
to bound the magnitude of the right-hand sides of
\eqref{E:FIRSTPOINTWISETRACELIEZNOFTRACEFREEISLOWERORDER}
by the right-hand side of \eqref{E:COMMUTINGTRACEFREEWITHLIEDERIVATIVES}.

To prove \eqref{E:DIFFERENCEBETWEENTRACEFREELIEDERIVATIVESANDLIEDERIVATIVESACTINGONATRACEFREETENSOR},
we let $Z \in \mathscr{Z}.$ From \eqref{E:ANGFREELIEIDFORTRACEFREETENSORS}
and the definition of the trace-free part of a tensor,
we deduce the schematic identity
\begin{align} \label{E:RESTATEDANGFREELIEIDFORTRACEFREETENSORS}
	\Lie_Z \hat{\xi}
	- \angfreeLiearg{Z} \hat{\xi}
	& = \left\lbrace
				(\ginversesphere)^2 \angdeform{Z} \xi 
				+ (\ginversesphere) \angdeform{Z} \ginversesphere \xi
			\right\rbrace
		\gsphere.
\end{align}
Iterating \eqref{E:RESTATEDANGFREELIEIDFORTRACEFREETENSORS},
using the schematic identities 
$\angLie_Z \gsphere = \angdeform{Z}$
and
$\angLie_Z \ginversesphere = - (\ginversesphere)^2 \angdeform{Z}$
(see Lemma~\ref{L:CONNECTIONBETWEENANGLIEOFGSPHEREANDDEFORMATIONTENSORS}),
and using \eqref{E:ANGFREELIEDEF} to express the
right-hand side of \eqref{E:ITERATEDSCHEMATICRESTATEDANGFREELIEIDFORTRACEFREETENSORS} in
terms of the Lie derivatives of $\xi,$ $\gsphere$ and $\ginversesphere,$
we deduce the following schematic identity, in which not all factors on the right-hand side
appear in the ``correct order'' (the order is irrelevant from the point of view of our estimates): 
	\begin{align} \label{E:ITERATEDSCHEMATICRESTATEDANGFREELIEIDFORTRACEFREETENSORS}
		\Lie_{\mathscr{Z}}^N \hat{\xi}
		- \angfreeLietwoarg{\mathscr{Z}}{N} \hat{\xi}
		& = 
		\sum_{M=1}^N
		\sum_{N_1 + N_2 + \cdots + N_k = N-M} 
			(\angLie_{\mathscr{Z}}^{N_1} \ginversesphere)
			\cdots 
			(\angLie_{\mathscr{Z}}^{N_a} \ginversesphere)
			(\angLie_{\mathscr{Z}}^{N_{a+1}} \gsphere)
			\cdots
			(\angLie_{\mathscr{Z}}^{N_{k-2}} \gsphere)
			(\angLie_{\mathscr{Z}}^{N_{k-1}} \angpi)
			(\angLie_{\mathscr{Z}}^{N_k} \xi).
	\end{align}
	Above, the $\angpi$ are the $S_{t,u}$ projections of deformation tensors of vectorfields in $\mathscr{Z}.$
	We now conclude the desired estimate 
	\eqref{E:DIFFERENCEBETWEENTRACEFREELIEDERIVATIVESANDLIEDERIVATIVESACTINGONATRACEFREETENSOR}	
	by using the estimates
	\eqref{E:CRUDEPOINTWISEBOUNDSDERIVATIVESOFANGULARDEFORMATIONTENSORS}
	and
	\eqref{E:CRUDELOWERORDERC0BOUNDDERIVATIVESOFANGULARDEFORMATIONTENSORS}
	to bound the magnitude of the right-hand side of
	\eqref{E:ITERATEDSCHEMATICRESTATEDANGFREELIEIDFORTRACEFREETENSORS}
	by the right-hand side of \eqref{E:DIFFERENCEBETWEENTRACEFREELIEDERIVATIVESANDLIEDERIVATIVESACTINGONATRACEFREETENSOR}.
	
	To prove \eqref{E:COMMUTINGTRACEFREELIEDERIVATIVESTHROUGHANGD}, we first apply
	$\angD$ to both sides of the schematic identity \eqref{E:ITERATEDSCHEMATICRESTATEDANGFREELIEIDFORTRACEFREETENSORS}.
	We then use inequality \eqref{E:TYPE02TENSORANGDINTERMSOFROTATIONALLIE}
	and the same reasoning use just after equation \eqref{E:ITERATEDSCHEMATICRESTATEDANGFREELIEIDFORTRACEFREETENSORS}
	to bound the magnitude of the right-hand side 
	by the right-hand side of \eqref{E:COMMUTINGTRACEFREELIEDERIVATIVESTHROUGHANGD} as desired.
	
	To prove \eqref{E:COMMUTINGANGLIEZANDTRACETRACEFREESPLITTING},
	we first decompose $\xi_{AB} = \hat{\xi}_{AB} + \frac{1}{2} (\ginversesphere)^{CD} \xi_{CD} \gsphere_{AB}.$
	Applying $\angLie_{\mathscr{Z}}^N$ to this identity and using the Leibniz rule,
	we deduce that
	\begin{align} \label{E:FIRSTPOINTWISEBOUNDCOMMUTINGANGLIEZANDTRACETRACEFREESPLITTING}
		\left|
			\angLie_{\mathscr{Z}}^N \xi
			-
			\left\lbrace
				\angLie_{\mathscr{Z}}^N \hat{\xi}
				+ \frac{1}{2} (\mathscr{Z}^N \mytr \xi) \gsphere
			\right\rbrace
		\right|
		& \lesssim	
			\mathop{\sum_{N_1 + N_2 + N_3 \leq N}}_{N_2 \leq N-1}
				\left|
					\angLie_{\mathscr{Z}}^{N_1} \ginversesphere
				\right|
				\left|
					\angLie_{\mathscr{Z}}^{N_2} \xi
				\right|
				\left|
					\angLie_{\mathscr{Z}}^{N_3} \gsphere
				\right|.
	\end{align}
	We then use the estimates
	\eqref{E:CRUDEPOINTWISEBOUNDSDERIVATIVESOFANGULARDEFORMATIONTENSORS}
	and
	\eqref{E:CRUDELOWERORDERC0BOUNDDERIVATIVESOFANGULARDEFORMATIONTENSORS}
	to bound the magnitude of the right-hand side of
	\eqref{E:FIRSTPOINTWISEBOUNDCOMMUTINGANGLIEZANDTRACETRACEFREESPLITTING}
	by the right-hand side of \eqref{E:COMMUTINGANGLIEZANDTRACETRACEFREESPLITTING}
	as desired.
	
	To prove \eqref{E:COMMUTINGANGDANDTRACETRACEFREESPLITTING}, 
	we first decompose $\xi_{AB} = \hat{\xi}_{AB} + \frac{1}{2} (\ginversesphere)^{CD} \xi_{CD} \gsphere_{AB}.$
	Applying $\angD \angLie_{\mathscr{Z}}^N$ to this identity and using the Leibniz rule
	and inequality \eqref{E:TYPE02TENSORANGDINTERMSOFROTATIONALLIE},
	we deduce that
	\begin{align} \label{E:FIRSTPOINTWISEBOUNDCOMMUTINGANGDANDTRACETRACEFREESPLITTING}
		\left|
			\angD \angLie_{\mathscr{Z}}^N \xi
			-
			\left\lbrace
				\angD \angLie_{\mathscr{Z}}^N \hat{\xi}
				+ \frac{1}{2} (\angdiff \mathscr{Z}^N \mytr \xi) \gsphere
			\right\rbrace
		\right|
		& \lesssim	
			\frac{1}{1 + t}
			\mathop{\sum_{N_1 + N_2 + N_3 \leq N+1}}_{N_2 \leq N}
				\left|
					\angLie_{\mathscr{Z}}^{N_1} \ginversesphere
				\right|
				\left|
					\angLie_{\mathscr{Z}}^{N_2} \xi
				\right|
				\left|
					\angLie_{\mathscr{Z}}^{N_3} \gsphere
				\right|.
	\end{align}
	Next, we use the estimates
	\eqref{E:CRUDEPOINTWISEBOUNDSDERIVATIVESOFANGULARDEFORMATIONTENSORS}
	and
	\eqref{E:CRUDELOWERORDERC0BOUNDDERIVATIVESOFANGULARDEFORMATIONTENSORS}
	to bound the magnitude of the right-hand side of
	\eqref{E:FIRSTPOINTWISEBOUNDCOMMUTINGANGDANDTRACETRACEFREESPLITTING}
	by the right-hand side of \eqref{E:COMMUTINGANGDANDTRACETRACEFREESPLITTING}.
	Finally, we use the already proven estimate \eqref{E:COMMUTINGTRACEFREELIEDERIVATIVESTHROUGHANGD}
	to replace the term $\angD \angLie_{\mathscr{Z}}^N \hat{\xi}$
	on the left-hand side of \eqref{E:FIRSTPOINTWISEBOUNDCOMMUTINGANGDANDTRACETRACEFREESPLITTING}
	with the term $\angD \angfreeLietwoarg{\mathscr{Z}}{N} \hat{\xi}$
	up to error terms that are bounded by in magnitude by
	the right-hand side of \eqref{E:COMMUTINGANGDANDTRACETRACEFREESPLITTING}.
	
	To prove \eqref{E:COMMUTINGTRACEFREELIEDERIVATIVESTHROUGHANGDIV}, 
	we consider the $ABC$ components of the difference of the two
	type $\binom{0}{3}$ $S_{t,u}$ tensors on the left-hand side
	of \eqref{E:COMMUTINGTRACEFREELIEDERIVATIVESTHROUGHANGD},
	where the $A$ index corresponds to $\angD.$ The
	$A,B$ traces of these tensors are also in magnitude $\lesssim$ the
	right-hand side of \eqref{E:COMMUTINGTRACEFREELIEDERIVATIVESTHROUGHANGD}.
	Furthermore, the $AB$ trace of the second tensor on the left-hand side of 
	\eqref{E:COMMUTINGTRACEFREELIEDERIVATIVESTHROUGHANGD} is precisely
	$\angdiv \angfreeLietwoarg{\mathscr{Z}}{N} \hat{\xi}.$ Hence, the
	desired estimate \eqref{E:COMMUTINGTRACEFREELIEDERIVATIVESTHROUGHANGDIV}
	will follow once we show that the
	difference between the $AB$ trace of the first tensor on the left-hand
	side of \eqref{E:COMMUTINGTRACEFREELIEDERIVATIVESTHROUGHANGD}
	and $\angLie_{\mathscr{Z}}^N \angdiv \hat{\xi}$ is in magnitude
	$\lesssim$ the right-hand side of \eqref{E:COMMUTINGTRACEFREELIEDERIVATIVESTHROUGHANGDIV}.
	To this end, we use the Leibniz rule to commute 	
	the operators $\angD$ and $\ginversesphere \angLie_{\mathscr{Z}}^{N_2},$
	thereby deducing that the difference under consideration can be bounded as follows:
	\begin{align} \label{E:FIRSTPOINTWISECOMMUTINGTRACEFREELIEDERIVATIVESTHROUGHANGDIV}
		\left|
			\angLie_{\mathscr{Z}}^N 
			\left\lbrace
				(\ginversesphere)^{AB} \angDarg{A} \hat{\xi}_{B \cdot})
			\right\rbrace
			- (\ginversesphere)^{AB} \angLie_{\mathscr{Z}}^N \angDarg{A} \hat{\xi}_{B \cdot}
		\right|
		& \lesssim
			\mathop{\sum_{N_1 + N_2 \leq N}}_{N_2 \leq N-1}
				\left|
					\angLie_{\mathscr{Z}}^{N_1} \ginversesphere
				\right|
				\left|
					\angD \angLie_{\mathscr{Z}}^{N_2} \hat{\xi}
				\right|
			\\
	& \ \ 
				+
				\mathop{\sum_{N_1 + N_2 \leq N}}_{N_2 \leq N-1}
				\left|
					\angLie_{\mathscr{Z}}^{N_1} \ginversesphere
				\right|
				\left|
					[\angD, \angLie_{\mathscr{Z}}^{N_2}] \hat{\xi}
				\right|,
				\notag
	\end{align}
	where the $\cdot'$s in \eqref{E:FIRSTPOINTWISECOMMUTINGTRACEFREELIEDERIVATIVESTHROUGHANGDIV} 
	signify that the contracted quantities on the left-hand side are $S_{t,u}$ one-forms.
	To bound the first product on the right-hand side of
	\eqref{E:FIRSTPOINTWISECOMMUTINGTRACEFREELIEDERIVATIVESTHROUGHANGDIV}
	by the right-hand side of \eqref{E:COMMUTINGTRACEFREELIEDERIVATIVESTHROUGHANGDIV},
	we use the estimates
	\eqref{E:CRUDEPOINTWISEBOUNDSDERIVATIVESOFANGULARDEFORMATIONTENSORS},
	\eqref{E:CRUDELOWERORDERC0BOUNDDERIVATIVESOFANGULARDEFORMATIONTENSORS},
	\eqref{E:TYPE02TENSORANGDINTERMSOFROTATIONALLIE},
	and
	\eqref{E:LIEDERIVATIVESOFTRACEFREEPARTINTERMSOFLIEDERIVATIVESOFFULLTENSOR},
	and the bootstrap assumptions 
	\eqref{E:PSIFUNDAMENTALC0BOUNDBOOTSTRAP},
	\eqref{E:UPMUBOOT},
	and \eqref{E:FRAMECOMPONENTSBOOT}.
	To bound the second product on the right-hand side of
	\eqref{E:FIRSTPOINTWISECOMMUTINGTRACEFREELIEDERIVATIVESTHROUGHANGDIV}
	by the right-hand side of \eqref{E:COMMUTINGTRACEFREELIEDERIVATIVESTHROUGHANGDIV},
	we use 
	the estimates
	\eqref{E:CRUDEPOINTWISEBOUNDSDERIVATIVESOFANGULARDEFORMATIONTENSORS},
	\eqref{E:CRUDELOWERORDERC0BOUNDDERIVATIVESOFANGULARDEFORMATIONTENSORS},
	\eqref{E:COMMUTATORESTIMATESVECTORFIELDSACTINGONANGDTENSORS},
	and
	\eqref{E:LIEDERIVATIVESOFTRACEFREEPARTINTERMSOFLIEDERIVATIVESOFFULLTENSOR},
	and the bootstrap assumptions 
	\eqref{E:PSIFUNDAMENTALC0BOUNDBOOTSTRAP},
	\eqref{E:UPMUBOOT},
	and \eqref{E:FRAMECOMPONENTSBOOT}.
	We have thus proved inequality \eqref{E:COMMUTINGTRACEFREELIEDERIVATIVESTHROUGHANGDIV}.
	Inequality \eqref{E:COMMUTINGLIEDERIVATIVESTHROUGHANGDIV}
	can be proved by using essentially the same argument, and
	we omit the details.
\end{proof}

\section{Pointwise estimates for \texorpdfstring{$\angLie_{\Lunit} \angLie_{\mathscr{Z}}^N \upchi^{(Small)}$}{
the Lie derivatives of the re-centered null second fundamental form 
involving an outgoing null differentiation} in terms of other quantities}
In this section, we derive pointwise estimates for $\angLie_{\Lunit} \angLie_{\mathscr{Z}}^N  \upchi_{(Small)}$
and for related quantities. We start with the following lemma, which provides an expression for
$\angLie_{\Lunit} \upchi^{(Small)}$ in terms of $\Psi$ and $\Lunit^i,$ $(i=1,2,3).$

\begin{lemma}[\textbf{Expressions for} 
$\angLie_{\Lunit}\upchi_{(Small)}$ 
\textbf{in terms of other quantities}]
\label{L:LIELANDLIERADCHIJUNKINTERMSOFOTHERVARIABLES}
The symmetric type $\binom{0}{2}$ $S_{t,u}$ tensorfield $\upchi^{(Small)}$ 
defined in \eqref{E:CHIJUNKDEF}
verifies the following transport equation:
\begin{align} \label{E:LIELCHIJUNKINTERMSOFOTHERVARIABLES}
	\angLie_{\Lunit} \upchi_{AB}^{(Small)}
	& =	\rgeo^{-1} 
			g_{ab} (\angdiffarg{A} x^a)
	 		\angdiffarg{B} \Lunit (\rgeo \Lunit_{(Small)}^b) 
			+ 
			\rgeo^{-2} 	
			g_{ab} (\angdiffarg{A} (\rgeo \Lunit_{(Small)}^a))
	 		\angdiffarg{B} (\rgeo \Lunit_{(Small)}^b)
	 				\\
	 & \ \ 
	 	+ \rgeo^{-1} 
	 		G_{ab} 
	 		(\angdiffarg{A} x^a) 
	 		(\angdiffarg{B} (\rgeo \Lunit_{(Small)}^b))
	 		\Lunit \Psi
	 	- \angLie_{\Lunit} 
	 		\Lambda_{AB}^{(Tan-\Psi)},
	 	\notag 
\end{align}
where $\Lambda_{AB}^{(Tan-\Psi)}$ is the type $\binom{0}{2}$ $S_{t,u}$ tensorfield given by \eqref{E:BIGLAMBDAGOOD}.
\end{lemma}

\begin{proof}
We first use the decomposition 
$\Lunit^a = \rgeo^{-1} x^a + \Lunit_{(Small)}^a$
(that is, \eqref{E:LUNITJUNK}),
the identity $\Lunit x^a = \Lunit^a,$
and Lemma~\ref{L:LANDRADCOMMUTEWITHANGDIFF} to deduce that
\begin{align} \label{E:ANGLIELOFANGDIFFRECTCOORDINATE}
	\angLie_{\Lunit} \angdiffarg{A} x^a = \angdiffarg{A} \Lunit^a = \rgeo^{-1} \angdiffarg{A} x^a + \angdiffarg{A} \Lunit_{(Small)}^a.
\end{align}
We now apply $\angLie_{\Lunit}$ to both sides of equation \eqref{E:CHIJUNKINTERMSOFOTHERVARIABLES}
and treat all uppercase Latin indices as tensorial $S_{t,u}$ indices,
and all lowercase Latin-indexed quantities as scalar-valued functions.
Using \eqref{E:ANGLIELOFANGDIFFRECTCOORDINATE}, 
the identity $\Lunit \rgeo = 1,$
and the chain rule identity $\Lunit g_{ab} = G_{ab} \Lunit \Psi,$
we deduce that
\begin{align} \label{E:ANGLIELCHIJUNKFIRSTRELATION}
\angLie_{\Lunit} \upchi^{(Small)}_{AB} 
 & =
	\rgeo^{-1} 
	g_{ab} (\angdiffarg{A} x^a)
	\angdiffarg{B} \Lunit (\rgeo \Lunit_{(Small)}^b)
	+
	g_{ab} (\angdiffarg{A} \Lunit_{(Small)}^a)
	\angdiffarg{B} \Lunit_{(Small)}^b
		\\
	& \ \ 
		+ G_{ab} 
	 	(\angdiffarg{A} x^a) 
	 	(\angdiffarg{B} (\rgeo \Lunit_{(Small)}^b))
	 	\Lunit \Psi
		- \angLie_{\Lunit} \Lambda_{AB}^{(Tan-\Psi)}.
		\notag
\end{align}
The desired identity \eqref{E:LIELCHIJUNKINTERMSOFOTHERVARIABLES} 
now easily follows from \eqref{E:ANGLIELCHIJUNKFIRSTRELATION}
and the fact that $\angdiff \rgeo = 0.$

\end{proof}

\begin{lemma} [\textbf{Pointwise estimates for} 
$\angLie_{\Lunit} \angLie_{\mathscr{Z}}^N  \upchi_{(Small)}$ 
\textbf{in terms of other variables}]
\label{L:POINTWISEESTIMATESFORLDERIVATIVESOFCHIJUNKINTERMSOFOTHERVARIABLES}
Let $0 \leq N \leq 23$ be an integer. 
Under the small-data and bootstrap assumptions 
of Sects.~\ref{S:PSISOLVES}-\ref{S:C0BOUNDBOOTSTRAP},
if $\varepsilon$ is sufficiently small,
then the following estimates hold on $\mathcal{M}_{\Tboot,U_0}:$
\begin{align} \label{E:POINTWISEESTIMATESFORLDERIVATIVESOFCHIJUNKINTERMSOFOTHERVARIABLES}
	\left|
		\fivemyarray
			[\rgeo^2 \angLie_{\Lunit} \angLie_{\mathscr{Z}}^N \upchi^{(Small)}]
			{\angLie_{\Lunit} (\rgeo^2 \angLie_{\mathscr{Z}}^N \upchi^{(Small)\#})}
			{\Lunit (\rgeo^2 \mathscr{Z}^N \mytr \upchi^{(Small)})}
			{\rgeo^2 \angLie_{\Lunit} \angLie_{\mathscr{Z}}^N \hat{\upchi}^{(Small)}}
			{\angLie_{\Lunit} (\rgeo^2 \angLie_{\mathscr{Z}}^N \hat{\upchi}^{(Small)\#})}
	\right|
	& \lesssim
		\left|
				\threemyarray[\rgeo \Lunit \mathscr{Z}^{\leq N+1} \Psi]
					{\rgeo \angdiff \mathscr{Z}^{\leq N+1} \Psi} 
					{\mathscr{Z}^{\leq N+1} \Psi}
		\right|
		+	\sum_{a=1}^3 
				\left|
						\Lunit (\rgeo \mathscr{Z}^{\leq N+1} \Lunit_{(Small)}^a)
				\right|
				\\
	& \ \ 
		+ \frac{\ln(\myexp + t)}{(1+t)^2}
			\left|
				\myarray[\mathscr{Z}^{\leq N+1} (\upmu - 1)]
					{\sum_{a=1}^3 \rgeo |\mathscr{Z}^{\leq N+1} \Lunit_{(Small)}^a|} 
			\right|.
			\notag
\end{align}

\end{lemma}

\begin{proof}
We first prove the estimate
\eqref{E:POINTWISEESTIMATESFORLDERIVATIVESOFCHIJUNKINTERMSOFOTHERVARIABLES}
for $\rgeo^2 \angLie_{\Lunit} \angLie_{\mathscr{Z}}^N \upchi^{(Small)}.$
As a first step, we bound $\rgeo^2 \angLie_{\mathscr{Z}}^N \angLie_{\Lunit} \upchi^{(Small)}$
by the right-hand side of \eqref{E:POINTWISEESTIMATESFORLDERIVATIVESOFCHIJUNKINTERMSOFOTHERVARIABLES}.
To this end,
we write equation \eqref{E:LIELCHIJUNKINTERMSOFOTHERVARIABLES} in the following schematic form
(recall that $\rgeo \Lunit \in \mathscr{Z}$):
\begin{align} \label{E:SCHEMATICAFORMANGLIELCHIJUNK}
	\angLie_{\Lunit} \upchi^{(Small)}
	& = \frac{1}{\rgeo^2} \smoothfunction(\Psi) (\angdiff x) \angdiff (\rgeo \Lunit (\rgeo \Lunit_{(Small)}))
		+ \frac{1}{\rgeo^2} \smoothfunction(\Psi) (\angdiff (\rgeo \Lunit_{(Small)}))^2
		+ \frac{1}{\rgeo} \smoothfunction(\Psi) (\angdiff x) (\angdiff (\rgeo \Lunit_{(Small)})) \Lunit \Psi
			\\
	& \ \
		+ \frac{1}{\rgeo}
			(\angLie_{\rgeo \Lunit} G_{(Frame)})
			\myarray
				[\Lunit \Psi]
				{\angdiff \Psi}
		+ \frac{1}{\rgeo}
			G_{(Frame)}
			\angLie_{\rgeo \Lunit}
			\myarray
				[\Lunit \Psi]
				{\angdiff \Psi},
				\notag
\end{align}
where the $\smoothfunction$ are smooth scalar-valued functions of $\Psi.$
We then apply the operator $\angLie_{\mathscr{Z}}^N$ to
the right-hand side of \eqref{E:SCHEMATICAFORMANGLIELCHIJUNK} 
and apply the Leibniz rule to the products on the right-hand side.
We bound the terms $\mathscr{Z}^M \rgeo$ with \eqref{E:ZNAPPLIEDTORGEOISNOTTOOLARGE}.
We bound the terms $\mathscr{Z}^M \smoothfunction(\Psi)$ with
the bootstrap assumptions \eqref{E:PSIFUNDAMENTALC0BOUNDBOOTSTRAP}.
The factors $\angLie_{\mathscr{Z}}^M \angdiff (\rgeo \Lunit (\rgeo \Lunit_{(Small)}))$ 
are the only ones that contribute to the second term on the right-hand side of \eqref{E:POINTWISEESTIMATESFORLDERIVATIVESOFCHIJUNKINTERMSOFOTHERVARIABLES};
we bound them with 
\eqref{E:ZNAPPLIEDTORGEOISNOTTOOLARGE},
Lemma~\ref{L:LANDRADCOMMUTEWITHANGDIFF},
inequality \eqref{E:FUNCTIONPOINTWISEANGDINTERMSOFANGLIEO},
and the bootstrap assumptions \eqref{E:FRAMECOMPONENTSBOOT}.
We bound the terms
$\angLie_{\mathscr{Z}}^M G_{(Frame)}$
and
$\angLie_{\mathscr{Z}}^M (\angLie_{\rgeo \Lunit} G_{(Frame)})$ 
with Lemma~\ref{L:POINTWISEESTIMATESGFRAMEINTERMSOFOTHERQUANTITIES}.
We bound the terms $\angLie_{\mathscr{Z}}^M \angdiff x$
with the estimates 
\eqref{E:POINTWISEBOUNDPROJECTEDLIEDERIVATIVESANGDIFFCOORDINATEX}
and
\eqref{E:LOWERORDERPOINTWISEBOUNDPROJECTEDLIEDERIVATIVESANGDIFFCOORDINATEX}.
We bound the terms
$\myarray[\mathscr{Z}^M \Lunit \Psi] {\angLie_{\mathscr{Z}}^M \angdiff \Psi}$
with Lemma~\ref{L:AVOIDINGCOMMUTING}
and the bootstrap assumptions \eqref{E:PSIFUNDAMENTALC0BOUNDBOOTSTRAP}.
Also using inequality \eqref{E:FUNCTIONAVOIDINGCOMMUTING},
we see that in total,
these estimates imply that 
$\rgeo^2$ times
the magnitude
of the right-hand side of \eqref{E:SCHEMATICAFORMANGLIELCHIJUNK}
is $\lesssim$ the right-hand side of
\eqref{E:POINTWISEESTIMATESFORLDERIVATIVESOFCHIJUNKINTERMSOFOTHERVARIABLES}
as desired. To complete the proof, we must bound the commutator term
$\rgeo^2 [\angLie_{\Lunit}, \angLie_{\mathscr{Z}}^N] \upchi^{(Small)}$
by the right-hand side \eqref{E:POINTWISEESTIMATESFORLDERIVATIVESOFCHIJUNKINTERMSOFOTHERVARIABLES}.
To derive the desired estimate, we use inequality
\eqref{E:LZNCOMMUTATORACTINGONTENSORFIELDSPOINTWISE}
with $\upchi^{(Small)}$ in the role of $\xi$
(and the first term 
$\left|
	\angLie_{\Lunit} \angLie_{\mathscr{Z}}^{N-1} \upchi^{(Small)}
\right|$
on the right-hand side of \eqref{E:LZNCOMMUTATORACTINGONTENSORFIELDSPOINTWISE}
is bounded by induction)
and the estimates of Lemma~\ref{L:POINTWISEESTIMATESFORCHIJUNKINTERMSOFOTHERVARIABLES}.

To bound the magnitude of 
$\angLie_{\Lunit} (\rgeo^2 \angLie_{\mathscr{Z}}^N \upchi^{(Small)\#})$
by the right-hand side of
\eqref{E:POINTWISEESTIMATESFORLDERIVATIVESOFCHIJUNKINTERMSOFOTHERVARIABLES},
we first note the following identity, which holds
for any type $\binom{0}{2}$ $S_{t,u}$ tensorfield $\xi:$
\begin{align} \label{E:ANGLIELRGEOSQUAREDXIINTERMSOFRGEOSQUAREDANGLIELXI}
	\angLie_{\Lunit}
	(\rgeo^2 \xi^{\#})
	& = - 2 \rgeo^2 \upchi_{AB}^{(Small)} \xi^{AB}
	+ \rgeo^2 (\angLie_{\Lunit} \xi)^{\#}.
\end{align}
The identity \eqref{E:ANGLIELRGEOSQUAREDXIINTERMSOFRGEOSQUAREDANGLIELXI}
follows easily from the identity $(\angLie_{\Lunit} \ginversesphere)^{AB} = - 2 \upchi^{AB},$
\eqref{E:CHIJUNKDEF}, and the identity $\Lunit \rgeo = 1.$
We now set $\xi = \angLie_{\mathscr{Z}}^N \upchi^{(Small)}$ 
in \eqref{E:ANGLIELRGEOSQUAREDXIINTERMSOFRGEOSQUAREDANGLIELXI}.
The second term on the right-hand side of 
\eqref{E:ANGLIELRGEOSQUAREDXIINTERMSOFRGEOSQUAREDANGLIELXI} has already been suitably 
bounded by previously proven estimate
\eqref{E:POINTWISEESTIMATESFORLDERIVATIVESOFCHIJUNKINTERMSOFOTHERVARIABLES}
for  $\rgeo^2 \angLie_{\Lunit} \angLie_{\mathscr{Z}}^N \upchi^{(Small)}.$
To bound the first term on the right-hand side of
\eqref{E:ANGLIELRGEOSQUAREDXIINTERMSOFRGEOSQUAREDANGLIELXI}
in magnitude by the right-hand side of
\eqref{E:POINTWISEESTIMATESFORLDERIVATIVESOFCHIJUNKINTERMSOFOTHERVARIABLES},
we use the estimates of 
Lemma~\ref{L:POINTWISEESTIMATESFORCHIJUNKINTERMSOFOTHERVARIABLES}.

To bound the magnitude of 
$\Lunit (\rgeo^2 \mathscr{Z}^N \mytr \upchi^{(Small)})$
by the right-hand side of
\eqref{E:POINTWISEESTIMATESFORLDERIVATIVESOFCHIJUNKINTERMSOFOTHERVARIABLES},
we note that 
$\Lunit (\rgeo^2 \mathscr{Z}^N \mytr \upchi^{(Small)})$
is the pure trace of the type 
$\binom{1}{1}$ $S_{t,u}$ tensorfield
$\angLie_{\Lunit} (\rgeo^2 \angLie_{\mathscr{Z}}^N \upchi^{(Small)\#}).$
Hence, the desired estimate follows from
the previously proven estimate
\eqref{E:POINTWISEESTIMATESFORLDERIVATIVESOFCHIJUNKINTERMSOFOTHERVARIABLES}
for 
$\angLie_{\Lunit} (\rgeo^2 \angLie_{\mathscr{Z}}^N \upchi^{(Small)\#}).$

To bound the magnitude of
$\rgeo^2 \angLie_{\Lunit} \angLie_{\mathscr{Z}}^N \hat{\upchi}^{(Small)}$
by the right-hand side of
\eqref{E:POINTWISEESTIMATESFORLDERIVATIVESOFCHIJUNKINTERMSOFOTHERVARIABLES},
we first use the Leibniz rule to
derive the following analog of \eqref{E:FIRSTPOINTWISEESTIMATELIEDERIVATIVESOFTRACEFREEPARTINTERMSOFLIEDERIVATIVESOFFULLTENSOR}:
\begin{align} \label{E:FIRSTPOINTWISEESTIMATELIEDERIVATIVESWITHLOFTRACEFREEPARTINTERMSOFLIEDERIVATIVESOFFULLTENSOR}
	\left|
		\angLie_{\Lunit} \angLie_{\mathscr{Z}}^N \hat{\upchi}^{(Small)}
	\right|
	& \lesssim 
		\left|
			\angLie_{\Lunit} \angLie_{\mathscr{Z}}^N \upchi^{(Small)}
		\right|
		+
		\sum_{i_1 + i_2 + i_3 = 1}
		\sum_{N_1 + N_2 + N_3 = N}
			\left|
				\angLie_{\Lunit}^{i_1}
				\angLie_{\mathscr{Z}}^{N_1} \upchi^{(Small)}
			\right|
			\left|
				\angLie_{\Lunit}^{i_2}
				\angLie_{\mathscr{Z}}^{N_2} \gsphere
			\right|
			\left|
				\angLie_{\Lunit}^{i_3}
				\angLie_{\mathscr{Z}}^{N_3} \ginversesphere
			\right|.
\end{align}	
The desired bound now follows from
\eqref{E:FIRSTPOINTWISEESTIMATELIEDERIVATIVESWITHLOFTRACEFREEPARTINTERMSOFLIEDERIVATIVESOFFULLTENSOR},
the previously proven estimate \eqref{E:POINTWISEESTIMATESFORLDERIVATIVESOFCHIJUNKINTERMSOFOTHERVARIABLES}
for  
$\rgeo^2 \angLie_{\Lunit} \angLie_{\mathscr{Z}}^N \upchi^{(Small)},$
the estimates of Lemmas 
\ref{L:POINTWISEBOUNDSDERIVATIVESOFANGULARDEFORMATIONTENSORS} and
\ref{L:POINTWISEESTIMATESFORCHIJUNKINTERMSOFOTHERVARIABLES},
and the fact that $\rgeo \Lunit \in \mathscr{Z}.$

Finally, the 
desired estimate \eqref{E:POINTWISEESTIMATESFORLDERIVATIVESOFCHIJUNKINTERMSOFOTHERVARIABLES}
for
$\angLie_{\Lunit} (\rgeo^2 \angLie_{\mathscr{Z}}^N \hat{\upchi}^{(Small)\#})$
follows from the estimate 
\eqref{E:POINTWISEESTIMATESFORLDERIVATIVESOFCHIJUNKINTERMSOFOTHERVARIABLES}
for $\rgeo^2 \angLie_{\Lunit} \angLie_{\mathscr{Z}}^N \hat{\upchi}^{(Small)}$
in the same way that the estimate
for $\angLie_{\Lunit} (\rgeo^2 \angLie_{\mathscr{Z}}^N \upchi^{(Small)\#})$
followed from the estimate for
$\rgeo^2 \angLie_{\Lunit} \angLie_{\mathscr{Z}}^N \upchi^{(Small)}.$

\end{proof}

\section{Improvement of the auxiliary bootstrap assumptions}
In this section, we use previously derived 
transport equations to derive improvements of 
the $C^0$ bootstrap assumptions for
$\upmu,$ $\Lunit_{(Small)}^i,$ and their lower-order derivatives.
These estimates \emph{involve a loss of one order of differentiability} because the right-hand 
sides of the transport equations depend on higher derivatives of $\Psi.$

\begin{proposition}[\textbf{Estimates for} $\upmu,$ $\Lunit_{(Small)}^i,$ \textbf{and} $\upchi^{(Small)}$
\textbf{derived from transport equations in the direction} $\Lunit$]
\label{P:CRUCICALTRANSPORTINTEQUALITIES}
Let $0 \leq N \leq 24$ be an integer.
Under the 
small data and
bootstrap assumptions 
of Sects.~\ref{S:PSISOLVES}-\ref{S:C0BOUNDBOOTSTRAP},
if $\varepsilon$ is sufficiently small, 
then the following pointwise estimates hold on $\mathcal{M}_{\Tboot,U_0}:$
\begin{subequations}
\begin{align} \label{E:POINTWISEEIKONALFUNCTIONTRANSPORTBASEDINTEQUALITIES}
	\left| 
		\fivemyarray
			[\rgeo^2 \angLie_{\mathscr{Z}}^N \upchi^{(Small)}]
			{\rgeo^2 \angLie_{\mathscr{Z}}^N \upchi^{(Small) \#}}
			{\rgeo^2 \mathscr{Z}^N \mytr \upchi^{(Small)}}
			{\rgeo^2 \angLie_{\mathscr{Z}}^N \hat{\upchi}^{(Small)} }
			{\rgeo^2 \angLie_{\mathscr{Z}}^N \hat{\upchi}^{(Small) \#}}
	\right|
	& \lesssim
		(1 + t)
		\left| 
				\fourmyarray[ \rgeo \Lunit \mathscr{Z}^{\leq N} \Psi]
					{\Rad \mathscr{Z}^{\leq N} \Psi}
					{\rgeo \angdiff \mathscr{Z}^{\leq N} \Psi}
					{\mathscr{Z}^{\leq N} \Psi}
		\right|
		+ 			\left| 
							\myarray
								[\mathscr{Z}^{\leq N} (\upmu - 1)]
								{\sum_{a=1}^3 \rgeo \left|\mathscr{Z}^{\leq N+1} \Lunit_{(Small)}^a \right|}
						\right|,
			\\
		&
	\left| 
		\sevenmyarray
			[\Lunit \mathscr{Z}^{\leq N} \upmu]
			{\sum_{a=1}^3 \left|\Lunit (\rgeo \mathscr{Z}^{\leq N} \Lunit_{(Small)}^a) \right|}
			{\rgeo^2 \angLie_{\Lunit} \angLie_{\mathscr{Z}}^{\leq N-1} \upchi^{(Small)} }
			{\angLie_{\Lunit} (\rgeo^2 \angLie_{\mathscr{Z}}^{\leq N-1} \upchi^{(Small) \#})}
			{\Lunit(\rgeo^2 \mathscr{Z}^{\leq N-1} \mytr \upchi^{(Small)})}
			{\rgeo^2 \angLie_{\Lunit} \angLie_{\mathscr{Z}}^{\leq N-1} \hat{\upchi}^{(Small)} }
			{\angLie_{\Lunit} (\rgeo^2 \angLie_{\mathscr{Z}}^{\leq N-1} \hat{\upchi}^{(Small) \#})}
		\right|,
		\, 
		\left|
		\sevenmyarray
			[\mathscr{Z}^{\leq N} \Lunit \upmu]
			{\sum_{a=1}^3 \left|\mathscr{Z}^{\leq N} (\Lunit (\rgeo \Lunit_{(Small)}^a)) \right|}
			{\angLie_{\mathscr{Z}}^{\leq N-1} (\rgeo^2 \angLie_{\Lunit} \upchi^{(Small)}) }
			{ \angLie_{\mathscr{Z}}^{\leq N-1} (\angLie_{\Lunit} (\rgeo^2 \upchi^{(Small) \#}))}
			{\mathscr{Z}^{\leq N-1} (\Lunit(\rgeo^2 \mytr \upchi^{(Small)}))}
			{\angLie_{\mathscr{Z}}^{\leq N-1} (\rgeo^2 \angLie_{\Lunit} \hat{\upchi}^{(Small)}) }
			{\angLie_{\mathscr{Z}}^{\leq N-1} (\angLie_{\Lunit} (\rgeo^2 \hat{\upchi}^{(Small) \#}))}
		\right|
			\label{E:LDERIVATIVECRUCICALTRANSPORTINTEQUALITIES} \\
	& \lesssim
		\left| 
				\fourmyarray[ \rgeo \Lunit \mathscr{Z}^{\leq N} \Psi]
					{\Rad \mathscr{Z}^{\leq N} \Psi}
					{\rgeo \angdiff \mathscr{Z}^{\leq N} \Psi}
					{\mathscr{Z}^{\leq N} \Psi}
		\right|
		+ 			\frac{\ln(\myexp + t)}{(1 + t)^2}
						\left| 
							\myarray
								[\mathscr{Z}^{\leq N} (\upmu - 1)]
								{\sum_{a=1}^3 \rgeo \left|\mathscr{Z}^{\leq N} \Lunit_{(Small)}^a \right|}
						\right|.
						\notag
\end{align}
\end{subequations}

Furthermore, the following estimates hold for $t \in [0,\Tboot):$
\begin{align}  \label{E:C0BOUNDCRUCIALEIKONALFUNCTIONQUANTITIES}
	\left\| 
		\sevenmyarray
			[\mathscr{Z}^{\leq 12} (\upmu - 1)]
			{\sum_{a=1}^3 \left|\rgeo \mathscr{Z}^{\leq 12} \Lunit_{(Small)}^a \right|}
			{\rgeo^2 \angLie_{\mathscr{Z}}^{\leq 11} \upchi^{(Small)} }
			{\rgeo^2 \angLie_{\mathscr{Z}}^{\leq 11} \upchi^{(Small) \#}}
			{\rgeo^2 \mathscr{Z}^{\leq 11} \mytr \upchi^{(Small)}}
			{\rgeo^2 \angLie_{\mathscr{Z}}^{\leq 11} \hat{\upchi}^{(Small)} }
			{\rgeo^2 \angLie_{\mathscr{Z}}^{\leq 11} \hat{\upchi}^{(Small) \#}}
		\right\|_{C^0(\Sigma_t^u)}
	& \lesssim \varepsilon \ln(\myexp + t).
\end{align}

Finally, the following estimates hold for $t \in [0,\Tboot):$
\begin{align}   \label{E:C0BOUNDLDERIVATIVECRUCICALEIKONALFUNCTIONQUANTITIES}
	\left\|
		\sevenmyarray
			[\Lunit \mathscr{Z}^{\leq 12} \upmu]
			{\sum_{a=1}^3 \left|\Lunit (\rgeo \mathscr{Z}^{\leq 12} \Lunit_{(Small)}^a) \right|}
			{\rgeo^2 \angLie_{\Lunit} \angLie_{\mathscr{Z}}^{\leq 11} \upchi^{(Small)} }
			{\angLie_{\Lunit} (\rgeo^2 \angLie_{\mathscr{Z}}^{\leq 11} \upchi^{(Small) \#})}
			{\Lunit(\rgeo^2 \mathscr{Z}^{\leq 11} \mytr \upchi^{(Small)})}
			{\rgeo^2 \angLie_{\Lunit} \angLie_{\mathscr{Z}}^{\leq 11} \hat{\upchi}^{(Small)} }
			{\angLie_{\Lunit} (\rgeo^2 \angLie_{\mathscr{Z}}^{\leq 11} \hat{\upchi}^{(Small) \#})}
		\right\|_{C^0(\Sigma_t^u)},
		\, 
		\left\|
		\sevenmyarray
			[\mathscr{Z}^{\leq 12} \Lunit \upmu]
			{\sum_{a=1}^3 \left|\mathscr{Z}^{\leq 12} (\Lunit (\rgeo \Lunit_{(Small)}^a)) \right|}
			{\angLie_{\mathscr{Z}}^{\leq 11} (\rgeo^2 \angLie_{\Lunit} \upchi^{(Small)}) }
			{ \angLie_{\mathscr{Z}}^{\leq 11} (\angLie_{\Lunit} (\rgeo^2 \upchi^{(Small) \#}))}
			{\mathscr{Z}^{\leq 11} (\Lunit(\rgeo^2 \mytr \upchi^{(Small)}))}
			{\angLie_{\mathscr{Z}}^{\leq 11} (\rgeo^2 \angLie_{\Lunit} \hat{\upchi}^{(Small)}) }
			{\angLie_{\mathscr{Z}}^{\leq 11} (\angLie_{\Lunit} (\rgeo^2 \hat{\upchi}^{(Small) \#}))}
		\right\|_{C^0(\Sigma_t^u)}
	& \lesssim \varepsilon \frac{1}{1 + t}.
\end{align}

\end{proposition}

\begin{proof}
We first prove \eqref{E:LDERIVATIVECRUCICALTRANSPORTINTEQUALITIES}. We prove only the estimates for the first
array on the left-hand side because the proofs of the estimates for the second array do not involve 
commutations and hence are simpler. To proceed with our estimates for the first array, we note that
by Lemma~\ref{L:POINTWISEESTIMATESFORLDERIVATIVESOFCHIJUNKINTERMSOFOTHERVARIABLES}, the estimates
for the $\upchi^{(Small)}$ terms in the array will follow once we prove inequality 
\eqref{E:LDERIVATIVECRUCICALTRANSPORTINTEQUALITIES} for only the first
two terms in the array.

We now use induction in $N$ to prove the desired estimate for the term
$\left|\Lunit (\rgeo \mathscr{Z}^{\leq N} \Lunit_{(Small)}^i) \right|$
on the left-hand side of \eqref{E:LDERIVATIVECRUCICALTRANSPORTINTEQUALITIES}.
To proceed, we write equation \eqref{E:LUNITJUNKLPROP} in the schematic form
\begin{align} \label{E:SCHEMATICLUNITJUNKLPROP}
	\Lunit (\rgeo \Lunit_{(Small)}^i)
	& = G_{(Frame)}
			\myarray[\Lunit \Psi]
				{\angdiff \Psi}
			\fourmyarray
				[\rgeo \Lunit_{(Small)}^i]
				{x^i}
				{\rgeo f^i(\Psi)}
				{\rgeo \angdiffuparg{\#} x^i}
		:= \mathfrak{I}^i,
\end{align}	
where $f^i(\Psi)$ is smooth function of $\Psi$ that vanishes at $\Psi = 0$
(specifically, $f^i(\Psi) =  (g^{-1})^{0i}$).
Commuting \eqref{E:SCHEMATICLUNITJUNKLPROP} with $\mathscr{Z}^N,$ we deduce that
\begin{align} \label{E:COMMUTEDSCHEMATICLUNITJUNKLPROP}
	\Lunit (\rgeo \mathscr{Z}^N \Lunit_{(Small)}^i)
	& = \mathscr{Z}^N \mathfrak{I}^i
		+ [\Lunit, \mathscr{Z}^N] (\rgeo \Lunit_{(Small)}^i)
		+ \Lunit 
			\left\lbrace
				[\rgeo, \mathscr{Z}^N]
				\Lunit_{(Small)}^i
			\right\rbrace.
\end{align}
We claim that the following bound holds for the first term 
$\mathscr{Z}^N \mathfrak{I}^i$ on the right-hand side of
\eqref{E:COMMUTEDSCHEMATICLUNITJUNKLPROP}:
\begin{align} \label{E:SCHEMATICLUNITJUNKLPROPINHOMOGENEOUSTERMPOINTWISE}
	\left|
		\mathscr{Z}^N \mathfrak{I}^i
	\right|
	& \lesssim 
		\left| 
				\fourmyarray[ \rgeo \Lunit \mathscr{Z}^{\leq N} \Psi]
					{\rgeo \angdiff \mathscr{Z}^{\leq N} \Psi}
					{\Rad \mathscr{Z}^{\leq N} \Psi}
					{\mathscr{Z}^{\leq N} \Psi}
		\right|
	+ 			\frac{\ln(\myexp + t)}{(1 + t)^2}
						\left| 
							\myarray
								[\mathscr{Z}^{\leq N} (\upmu - 1)]
								{\sum_{a=1}^3 \rgeo \left|\mathscr{Z}^{\leq N} \Lunit_{(Small)}^a \right|}
						\right|.
\end{align}
Let us accept \eqref{E:SCHEMATICLUNITJUNKLPROPINHOMOGENEOUSTERMPOINTWISE} for the moment;
we will independently prove \eqref{E:SCHEMATICLUNITJUNKLPROPINHOMOGENEOUSTERMPOINTWISE} at the end of the argument
without using induction.
We note that the right-hand side of \eqref{E:SCHEMATICLUNITJUNKLPROPINHOMOGENEOUSTERMPOINTWISE}
is $\lesssim$ the right-hand side of \eqref{E:LDERIVATIVECRUCICALTRANSPORTINTEQUALITIES} as desired.
This immediately yields the desired estimate in the base case $N=0.$ 
To carry out the induction, we assume that the estimate \eqref{E:LDERIVATIVECRUCICALTRANSPORTINTEQUALITIES}
for $\left|\Lunit (\rgeo \mathscr{Z}^{\leq N-1} \Lunit_{(Small)}^i) \right|$
has been proved.

To bound the term 
$\Lunit 
			\left\lbrace
				[\rgeo, \mathscr{Z}^N]
				\Lunit_{(Small)}^i
			\right\rbrace$
on the right-hand side of
\eqref{E:COMMUTEDSCHEMATICLUNITJUNKLPROP},			
we use \eqref{E:LAPPLIEDTOCOMMUTATOROFZNANDRGEOAPPLIEDTOXIPOINTWISE} to deduce that
\begin{align}  \label{E:SCHEMATICLUNITJUNKLPROPFIRSTCOMMUTATORPOINTWISE}
	\left|
		\Lunit 
			\left\lbrace
				[\rgeo, \mathscr{Z}^N]
				\Lunit_{(Small)}^i
			\right\rbrace
	\right|
	& \lesssim 
		\left|
			\Lunit (\rgeo \mathscr{Z}^{\leq N-1} \Lunit_{(Small)}^i)
		\right|
	+ \frac{1}{(1 + t)^2}
		\left|
			\rgeo \mathscr{Z}^{\leq N-1} \Lunit_{(Small)}^i
		\right|.
\end{align}
The second term on the right-hand side of \eqref{E:SCHEMATICLUNITJUNKLPROPFIRSTCOMMUTATORPOINTWISE}
is manifestly $\lesssim$ the right-hand side of \eqref{E:LDERIVATIVECRUCICALTRANSPORTINTEQUALITIES} as desired,
while to bound the first, we use the induction hypothesis.

To bound the term $[\Lunit, \mathscr{Z}^N] (\rgeo \Lunit_{(Small)}^i)$
on the right-hand side of \eqref{E:COMMUTEDSCHEMATICLUNITJUNKLPROP},
we use the commutator estimate \eqref{E:LZNCOMMUTATORACTINGONFUNCTIONSPOINTWISE}
with $f = \rgeo \Lunit_{(Small)}^i,$
inequalities \eqref{E:ZNAPPLIEDTORGEOISNOTTOOLARGE}
and \eqref{E:BOUNDINGFRAMEDERIVATIVESINTERMSOFZNF},
and the bootstrap assumption \eqref{E:FRAMECOMPONENTSBOOT} to deduce that
\begin{align} \label{E:SCHEMATICLUNITJUNKLPROPSECONDCOMMUTATORPOINTWISE}
	\left|
		[\Lunit, \mathscr{Z}^N] (\rgeo \Lunit_{(Small)}^i)
	\right|
	& \lesssim 
		\left|
			\Lunit \mathscr{Z}^{\leq N-1} (\rgeo \Lunit_{(Small)}^i)
		\right|
		+ 
			\frac{\ln(\myexp + t)}{1 + t}
			\left|
				\fourmyarray[\rgeo \Lunit \mathscr{Z}^{\leq N} \Psi]
					{\Rad \mathscr{Z}^{\leq N} \Psi}
					{\rgeo \angdiff \mathscr{Z}^{\leq N} \Psi} 
					{\mathscr{Z}^{\leq N} \Psi}
			\right|
			\\
	& \ \ 
			+  \frac{\ln(\myexp + t)}{(1 + t)^2}
			\left|
				\myarray[\mathscr{Z}^{\leq N} (\upmu - 1)]
					{\sum_{a=1}^3 \rgeo |\mathscr{Z}^{\leq N} \Lunit_{(Small)}^a|} 
			\right|.
			\notag
\end{align}
The last two terms on the right-hand side of \eqref{E:SCHEMATICLUNITJUNKLPROPSECONDCOMMUTATORPOINTWISE}
are manifestly $\lesssim$ the right-hand side of \eqref{E:LDERIVATIVECRUCICALTRANSPORTINTEQUALITIES} as desired,
while to bound the first, we use the induction hypothesis.

It remains for us to prove \eqref{E:SCHEMATICLUNITJUNKLPROPINHOMOGENEOUSTERMPOINTWISE}. We
apply the operator $\angLie_{\mathscr{Z}}^N$ to the terms in \eqref{E:SCHEMATICLUNITJUNKLPROP}
and use the Leibniz rule. 
We bound the terms 
$\angLie_{\mathscr{Z}}^M G_{(Frame)}$ with the estimates of Lemma~\ref{L:POINTWISEESTIMATESGFRAMEINTERMSOFOTHERQUANTITIES}.
To bound
$\myarray[\mathscr{Z}^M \Lunit \Psi] {\angLie_{\mathscr{Z}}^M\angdiff \Psi},$
we use Lemma~\ref{L:AVOIDINGCOMMUTING}
and the bootstrap assumptions \eqref{E:PSIFUNDAMENTALC0BOUNDBOOTSTRAP}.
To bound $\mathscr{Z}^M (\rgeo \Lunit_{(Small)}^i)$ when $M \leq 12,$ we use 
inequality \eqref{E:ZNAPPLIEDTORGEOISNOTTOOLARGE}
and the bootstrap 
assumptions \eqref{E:FRAMECOMPONENTSBOOT} 
(when $M \geq 13,$ the term $\rgeo \mathscr{Z}^M \Lunit_{(Small)}^i$
appears on the right-hand side of \eqref{E:SCHEMATICLUNITJUNKLPROPINHOMOGENEOUSTERMPOINTWISE}).
To bound $\mathscr{Z}^M x^i,$ we use inequalities 
\eqref{E:EASYCOORDINATEBOUND},
\eqref{E:DERIVATIVESOFRECTANGULARCOMMUTATORVECTORFIELDRECTCOMPONENTS},
and \eqref{E:LOWERORDERC0BOUNDDERIVATIVESOFRECTANGULARCOMMUTATORVECTORFIELDRECTCOMPONENTS}
as well as the observation that $\mathscr{Z}^M x^i = \mathscr{Z}^{M-1} Z x^i = \mathscr{Z}^{M-1} Z^i.$ 
To bound $\mathscr{Z}^M  (\rgeo f^i(\Psi)),$ we use
inequality \eqref{E:ZNAPPLIEDTORGEOISNOTTOOLARGE}
and the bootstrap assumptions \eqref{E:PSIFUNDAMENTALC0BOUNDBOOTSTRAP},
which in particular imply that $|\mathscr{Z}^M f^i(\Psi)| \lesssim |\mathscr{Z}^{\leq M} \Psi|.$
To bound $\angLie_{\mathscr{Z}}^M (\rgeo \angdiffuparg{\#} x^i),$
we use inequalities 
\eqref{E:ZNAPPLIEDTORGEOISNOTTOOLARGE},
\eqref{E:CRUDEPOINTWISEBOUNDSDERIVATIVESOFANGULARDEFORMATIONTENSORS},
\eqref{E:CRUDELOWERORDERC0BOUNDDERIVATIVESOFANGULARDEFORMATIONTENSORS},
\eqref{E:POINTWISEBOUNDPROJECTEDLIEDERIVATIVESANGDIFFCOORDINATEX},
and \eqref{E:LOWERORDERPOINTWISEBOUNDPROJECTEDLIEDERIVATIVESANGDIFFCOORDINATEX}.
In total, these estimates yield inequality \eqref{E:SCHEMATICLUNITJUNKLPROPINHOMOGENEOUSTERMPOINTWISE},
and our proof of the bound for $\left|\Lunit (\rgeo \mathscr{Z}^{\leq N} \Lunit_{(Small)}^i) \right|$ is complete.

\ \\

\noindent \textbf{Proof of the bound for} $\left|\Lunit \mathscr{Z}^{\leq N} \upmu \right|:$
We now use induction in $N$ to prove the desired estimate for the term
$\left|\Lunit \mathscr{Z}^{\leq N} \upmu \right|$
on the left-hand side of \eqref{E:LDERIVATIVECRUCICALTRANSPORTINTEQUALITIES}.
The proof is simpler than our proof of the bound for 
$\left|\Lunit (\rgeo \mathscr{Z}^{\leq N} \Lunit_{(Small)}^i) \right|.$
To proceed, we write equation \eqref{E:UPMUFIRSTTRANSPORT} in the schematic form
\begin{align} \label{E:SCHEMATICLUPMPROP}
	\Lunit \upmu 
	& = G_{(Frame)}
			\myarray
				[\upmu \Lunit \Psi]
				{\Rad \Psi}
		:= \mathfrak{I}.
\end{align}
Commuting \eqref{E:SCHEMATICLUPMPROP} with $\mathscr{Z}^N,$ we have
\begin{align} \label{E:COMMUTEDSCHEMATICLUPMPROP}
	\left|
		\Lunit \mathscr{Z}^N \upmu
	\right|
	& \lesssim
		\left|
			[\Lunit, \mathscr{Z}^N] \upmu 
		\right|
		+ \left|
				\mathscr{Z}^N \mathfrak{I}
			\right|.
\end{align}
We claim that
\begin{align} \label{E:SCHEMATICLUNITUPMUPROPINHOMOGENEOUSTERMPOINTWISE}
	 \left|
			\mathscr{Z}^N \mathfrak{I}
	 \right|
	 & \lesssim 
	 	\left| 
				\fourmyarray[\rgeo \Lunit \mathscr{Z}^{\leq N} \Psi]
					{\Rad \mathscr{Z}^{\leq N} \Psi}
					{\rgeo \angdiff \mathscr{Z}^{\leq N} \Psi}
					{\mathscr{Z}^{\leq N} \Psi}
		\right|
		+ 			\frac{\ln(\myexp + t)}{(1 + t)^2}
						\left| 
							\myarray
								[\mathscr{Z}^{\leq N} (\upmu - 1)]
								{\sum_{a=1}^3 \rgeo \left|\mathscr{Z}^{\leq N} \Lunit_{(Small)}^a \right|}
						\right|.
\end{align}
To derive the estimate \eqref{E:SCHEMATICLUNITUPMUPROPINHOMOGENEOUSTERMPOINTWISE}, we apply
the Leibniz rule. We bound the terms 
$\mathscr{Z}^M G_{(Frame)}$ with the estimates of Lemma~\ref{L:POINTWISEESTIMATESGFRAMEINTERMSOFOTHERQUANTITIES}.
To bound the terms $\angLie_{\mathscr{Z}}^M 
			\myarray
				[\Lunit \Psi]
				{\Rad \Psi},$
we use Lemma~\ref{L:AVOIDINGCOMMUTING}
and the bootstrap assumptions \eqref{E:PSIFUNDAMENTALC0BOUNDBOOTSTRAP}.
To bound $\mathscr{Z}^M \upmu$ when $M \leq 12,$ we use the bootstrap 
assumptions \eqref{E:UPMUBOOT}
(when $M \geq 13,$ the term $\mathscr{Z}^M (\upmu - 1)$
appears on the right-hand side of \eqref{E:SCHEMATICLUNITUPMUPROPINHOMOGENEOUSTERMPOINTWISE}).

We now note that the right-hand side of \eqref{E:SCHEMATICLUNITUPMUPROPINHOMOGENEOUSTERMPOINTWISE}
is manifestly $\lesssim$ the right-hand side of \eqref{E:LDERIVATIVECRUCICALTRANSPORTINTEQUALITIES} as desired.
The estimate \eqref{E:SCHEMATICLUNITUPMUPROPINHOMOGENEOUSTERMPOINTWISE} therefore
immediately yields the desired estimate for $\left|\Lunit \mathscr{Z}^{\leq N} \upmu \right|$
in the base case $N=0.$ To carry out the induction, we assume that the estimate \eqref{E:LDERIVATIVECRUCICALTRANSPORTINTEQUALITIES}
for $\left|\Lunit \mathscr{Z}^{\leq N-1} \upmu \right|$
has been proved.

To bound the commutator term on the right-hand side of \eqref{E:COMMUTEDSCHEMATICLUPMPROP},
we use the commutator estimate \eqref{E:LZNCOMMUTATORACTINGONFUNCTIONSPOINTWISE} 
with $f = \upmu,$
inequality \eqref{E:BOUNDINGFRAMEDERIVATIVESINTERMSOFZNF},
and the bootstrap assumptions \eqref{E:UPMUBOOT}
to deduce that
\begin{align}  \label{E:SCHEMATICLUNITUPMUPROPCOMMUTATORTERMPOINTWISE}
\left|
		[\Lunit, \mathscr{Z}^N] \upmu 
\right|	
& \lesssim
		\left|
			\Lunit \mathscr{Z}^{\leq N-1} \upmu
		\right|
	+ 
		\frac{\ln(\myexp + t)}{1 + t}
			\left|
				\mathscr{Z}^{\leq N} \Psi
			\right|
	+ \frac{\ln(\myexp + t)}{(1 + t)^2}
			\left|
				\myarray[\mathscr{Z}^{\leq N} (\upmu - 1)]
					{\sum_{a=1}^3 \rgeo |\mathscr{Z}^{\leq N} \Lunit_{(Small)}^a|} 
			\right|.
\end{align}
The last two terms on the right-hand side of 
\eqref{E:SCHEMATICLUNITUPMUPROPCOMMUTATORTERMPOINTWISE}
are manifestly bounded by the right-hand side of
\eqref{E:LDERIVATIVECRUCICALTRANSPORTINTEQUALITIES} as desired,
while to bound the first, we use the induction hypothesis.

\noindent \textbf{Proof of \eqref{E:C0BOUNDCRUCIALEIKONALFUNCTIONQUANTITIES}-\eqref{E:C0BOUNDLDERIVATIVECRUCICALEIKONALFUNCTIONQUANTITIES}}:
Let $0 \leq N \leq 12$ be an integer.
To prove \eqref{E:C0BOUNDCRUCIALEIKONALFUNCTIONQUANTITIES}, 
we first simplify the notation by defining
\begin{align}
			q_N(t,u,\vartheta) 
			&:= 
				\left|
					\myarray
						[\mathscr{Z}^{\leq N} (\upmu - 1)]
						{\sum_{a=1}^3 \left|\rgeo \mathscr{Z}^{\leq N} \Lunit_{(Small)}^a \right|}
				\right|,
			\qquad
			\mathring{q}_N(u,\vartheta)
			:= q_N(0,u,\vartheta).
	\end{align}
The bootstrap assumptions \eqref{E:PSIFUNDAMENTALC0BOUNDBOOTSTRAP} 
imply that the array of $\Psi$ quantities on the right-hand side of \eqref{E:LDERIVATIVECRUCICALTRANSPORTINTEQUALITIES}
is bounded in the norm $\| \cdot \|_{C^0(\Sigma_t^u)}$ by $\lesssim \varepsilon (1+t)^{-1}.$
Hence, integrating inequality \eqref{E:LDERIVATIVECRUCICALTRANSPORTINTEQUALITIES} in time along the integral 
curves of $\Lunit,$ we deduce that
\begin{align}  \label{E:GRONWALLREADYCRUCICALLDERIVATIVEOFEIKONALFUNCTIONQUANTITIESEQUALITIES}
	q_N(t,u,\vartheta) 		
	& \leq \mathring{q}_N(u,\vartheta)
		+ C \varepsilon \ln(\myexp + t)
		+ C
			\int_{t'=0}^t
				\frac{\ln(\myexp + t')}{(1 + t')^2}|q_N|(t',u,\vartheta)			
			\, dt'.
\end{align}
Inequality \eqref{E:SMALLINITIALC0BOUNDNORMS} 
and the small data assumption $\mathring{\upepsilon} \leq \varepsilon$
together imply that $\mathring{q}_N(u,\vartheta) \lesssim \varepsilon.$
Hence, applying Gronwall's inequality to \eqref{E:GRONWALLREADYCRUCICALLDERIVATIVEOFEIKONALFUNCTIONQUANTITIESEQUALITIES},
we deduce that $q_N(t,u,\vartheta) \lesssim \varepsilon \ln(\myexp + t).$ This yields the desired
inequalities for the first two terms in the array on the left-hand side of \eqref{E:C0BOUNDCRUCIALEIKONALFUNCTIONQUANTITIES}.
The estimates for the remaining quantities on the left-hand side of \eqref{E:C0BOUNDCRUCIALEIKONALFUNCTIONQUANTITIES}
then follow from these estimates, inequality \eqref{E:POINTWISEESTIMATESFORCHIJUNKINTERMSOFOTHERVARIABLES},
and the bootstrap assumptions \eqref{E:PSIFUNDAMENTALC0BOUNDBOOTSTRAP}.
	
Inequality \eqref{E:C0BOUNDLDERIVATIVECRUCICALEIKONALFUNCTIONQUANTITIES} then	
follows from inequalities 
\eqref{E:LDERIVATIVECRUCICALTRANSPORTINTEQUALITIES} 
and 
\eqref{E:C0BOUNDCRUCIALEIKONALFUNCTIONQUANTITIES}
and the bootstrap assumptions \eqref{E:PSIFUNDAMENTALC0BOUNDBOOTSTRAP}.
\end{proof}

\begin{corollary}[\textbf{Improvement of the auxiliary bootstrap assumptions}]
\label{C:AUXBOOTSTRAPIMPROVED}
Under the 
small data and
bootstrap assumptions 
of Sects.~\ref{S:PSISOLVES}-\ref{S:C0BOUNDBOOTSTRAP},
if $\varepsilon$ is sufficiently small, 
then the auxiliary bootstrap assumptions \eqref{E:UPMUBOOT}-\eqref{E:CHIJUNKBOOT}
hold with $\varepsilon^{1/2}$ replaced by $C \varepsilon.$
\end{corollary}

\begin{proof}
	The corollary follows directly from
	\eqref{E:C0BOUNDCRUCIALEIKONALFUNCTIONQUANTITIES}.
\end{proof}

Using the estimates of Cor.~\ref{C:AUXBOOTSTRAPIMPROVED} in place of the 
auxiliary bootstrap assumptions, we can repeat the proofs of all of the previous
estimates in this section to arrive at the following corollary.

\begin{corollary}[$\varepsilon^{1/2}$ \textbf{can be replaced with} $C \varepsilon$]
\label{C:SQRTEPSILONREPLCEDWITHCEPSILON}
	Under the 
	small data and
	bootstrap assumptions 
	of Sects.~\ref{S:PSISOLVES}-\ref{S:C0BOUNDBOOTSTRAP},
	all of the estimates in Chapter~\ref{C:C0BOUNDBOOTSTRAP} 
	that we proved before
	Cor.~\ref{C:AUXBOOTSTRAPIMPROVED}
	and that involve an explicit factor of precisely $\varepsilon^{1/2}$ 
	in fact hold with $\varepsilon^{1/2}$ replaced by $C \varepsilon.$
\end{corollary}

\section{Sharp pointwise estimates for \texorpdfstring{$G_{\Lunit \Lunit}$}{a frame component of the derivative of the metric with respect to the solution}}
In this section, we provide sharp pointwise estimates for 
$G_{\Lunit \Lunit}$ and $\Lunit G_{\Lunit \Lunit}.$ In particular, we show that
$G_{\Lunit \Lunit}$ is well-approximated by the future null condition failure factor $\FutFailFac.$

\begin{lemma}[\textbf{Bounds for $G_{\Lunit \Lunit}$}] \label{L:NULLCONDFACT}
Let $\FutFailFac$ be the future null condition failure factor from Def.~\ref{D:FAILUREFACTOR}.
Under the 
small data and
bootstrap assumptions 
of Sects.~\ref{S:PSISOLVES}-\ref{S:C0BOUNDBOOTSTRAP},
if $\varepsilon$ is sufficiently small,
then the following pointwise estimates hold on $\mathcal{M}_{\Tboot,U_0}:$
\begin{align} \label{E:IMPORTANTLDERIVATIVEOFGFRAMEC0BOUNDESTIMATES}
	\left|
		\Lunit G_{\Lunit \Lunit}
	\right|
	& \lesssim \varepsilon \frac{1}{(1 + t)^2},
\end{align}
\begin{align}  \label{E:NULLCONDFACTANNOYINGSTIMATE}
		\left| 
			G_{\Lunit \Lunit}
			- \FutFailFac
		\right|
		& \lesssim \varepsilon \frac{\ln(\myexp + t)}{1 + t}.
\end{align}

Furthermore, if $0 \leq s \leq t,$ then
\begin{align} \label{E:IMPORTANTGLLDIFFERENCEESTIMATE}
	\left|
		G_{\Lunit \Lunit}(t,u,\vartheta)
		- G_{\Lunit \Lunit}(s,u,\vartheta)
	\right|
	& \lesssim \varepsilon \frac{t-s}{(1 + s)(1+t)}
		\lesssim \varepsilon \frac{1}{1 + s}.
\end{align}
\end{lemma}

\begin{proof}
	From \eqref{E:LDERIVATIVEGLL}, the estimate \eqref{E:LOWERORDERC0BOUNDLIEDERIVATIVESOFGRAME},
	and the bootstrap assumptions \eqref{E:PSIFUNDAMENTALC0BOUNDBOOTSTRAP}, we 
	deduce that
	\begin{align} \label{E:LOFGLLBOUND}
		\left|
			\Lunit G_{\Lunit \Lunit}
		\right| 
		& \lesssim \varepsilon \frac{1}{(1 + t)^2},
	\end{align}
	which is the desired estimate \eqref{E:IMPORTANTLDERIVATIVEOFGFRAMEC0BOUNDESTIMATES}.
	
	Inequality \eqref{E:IMPORTANTGLLDIFFERENCEESTIMATE} then follows from
	integrating \eqref{E:IMPORTANTLDERIVATIVEOFGFRAMEC0BOUNDESTIMATES} 
	along the integral curves of $\Lunit$ from time $s$ to time $t.$
	
	
	To prove \eqref{E:NULLCONDFACTANNOYINGSTIMATE}, we first note that 
	by the mean value theorem, the estimate
	\eqref{E:POINTWISEBOUNDLOWERORDERDERIVATIVESOFRECTANGULARFRAMECOMPONENTS},
	and the bootstrap assumptions \eqref{E:PSIFUNDAMENTALC0BOUNDBOOTSTRAP}, we
	have
	\begin{align} \label{E:GLLNEARGATZEROPSILL}
		\left|
			G_{\Lunit \Lunit}
			- G_{\alpha \beta}(\Psi = 0)\Lunit^{\alpha} \Lunit^{\beta}
		\right|
		& \lesssim \frac{\varepsilon}{1 + t}.
	\end{align}
	Next, we recall that $ \Lunit_{(Flat)} = \partial_t + (x^a/r) \partial_a.$
	Hence, using definition \eqref{E:LUNITJUNK},
	the estimates 
	\eqref{E:EASYEUCLIDEANRADIALVARIABLEBOUND}
	and \eqref{E:C0BOUNDCRUCIALEIKONALFUNCTIONQUANTITIES},
	Cor.~\ref{C:SQRTEPSILONREPLCEDWITHCEPSILON},
	and the fact that $\Lunit^0 = \Lunit_{(Flat)}^0 = 1,$ we deduce that
	the following estimate holds for $\nu = 0,1,2,3:$
	\begin{align} \label{E:LCONVERGESTOLFLAT}
		\left|
			\Lunit^{\nu} - \Lunit_{(Flat)}^{\nu}
		\right|
		& \lesssim \frac{\varepsilon \ln(\myexp + t)}{1 + t}.
	\end{align}
	Using the estimates 
	\eqref{E:POINTWISEBOUNDLOWERORDERDERIVATIVESOFRECTANGULARFRAMECOMPONENTS}, 
	\eqref{E:GLLNEARGATZEROPSILL}, 
	and \eqref{E:LCONVERGESTOLFLAT}, 
	we conclude the desired bound \eqref{E:NULLCONDFACTANNOYINGSTIMATE}.
	
\end{proof}

\section{Pointwise estimates for \texorpdfstring{$\angLap \mathscr{Z}^{\leq N-1} \Lunit_{(small)}^i$}{the angular Laplacian of the derivatives of the re-centered version of the outgoing null vectorfield}}

In this section, we derive pointwise estimates for the scalar-valued functions
$\angLap \mathscr{Z}^{\leq N-1} \Lunit_{(Small)}^i,$ $(i=1,2,3).$

\begin{lemma}[\textbf{Pointwise estimate for} $\angLap \mathscr{Z}^{\leq N-1} \Lunit_{(Small)}^i$ \textbf{in terms of other variables}]
	\label{L:POINTWISEBOUNDFORANGLAPLJUNKICOMMUTED}
	Let $1 \leq N \leq 24$ be an integer, and let $\Lunit_{(Small)}^i,$ $(i=1,2,3),$
	be the scalar-valued function defined in \eqref{E:LUNITJUNK}.
	Under the 
	small data and
	bootstrap assumptions 
	of Sects.~\ref{S:PSISOLVES}-\ref{S:C0BOUNDBOOTSTRAP},
	if $\varepsilon$ is sufficiently small, then
	the following pointwise estimate holds on $\mathcal{M}_{\Tboot,U_0}:$
	\begin{align} \label{E:POINTWISEBOUNDFORANGLAPLJUNKICOMMUTED}
		\left|
			\angLap \mathscr{Z}^{\leq N-1} \Lunit_{(Small)}^i
		\right|
		& \lesssim 
			\frac{1}{(1 + t)^2}
		\left| 
			\fourmyarray[ \rgeo \Lunit \mathscr{Z}^{\leq N} \Psi]
				{\Rad \mathscr{Z}^{\leq N} \Psi}
				{\rgeo \angdiff \mathscr{Z}^{\leq N} \Psi}
				{\mathscr{Z}^{\leq N} \Psi}
			\right|	
			\\
	& \ \ + \frac{1}{(1 + t)^3}
			\left|
				\myarray[\mathscr{Z}^{\leq N} (\upmu - 1)]
					{\sum_{a=1}^3 \rgeo |\mathscr{Z}^{\leq N} \Lunit_{(Small)}^a|} 
			\right|
			+
			\left|
				\angdiv \angLie_{\mathscr{Z}}^{\leq N-1} \upchi^{(Small)}
			\right|.
			\notag
	\end{align}
\end{lemma} 

\begin{proof}
		We first apply $\angLie_{\mathscr{Z}}^{\leq N-1}$ and then $\angD^A$ 
		to each size of equation \eqref{E:ANGDIFFLJUNKI}. By Lemma~\ref{L:LANDRADCOMMUTEWITHANGDIFF},
		the left-hand side of the resulting expression is equal to $\angLap \mathscr{Z}^{\leq N-1} \Lunit_{(Small)}^i.$
		Hence, Lemma~\ref{L:POINTWISEBOUNDFORANGLAPLJUNKICOMMUTED} will follow once 
		we show that the right-hand side of the resulting expression is bounded in magnitude by 
		$\lesssim$ the right-hand side of \eqref{E:POINTWISEBOUNDFORANGLAPLJUNKICOMMUTED}.
		We begin by addressing the first product on the right-hand side of \eqref{E:ANGDIFFLJUNKI}:
		$\angD^A 
		\angLie_{\mathscr{Z}}^{\leq N-1} 
		\left\lbrace 
			(\ginversesphere)^{BC} \upchi^{(Small)}_{AB} \angdiffarg{C} x^i 
		\right\rbrace.$
		When all derivatives fall on $\upchi^{(Small)}_{AB},$ we use the estimate
		\eqref{E:LOWERORDERPOINTWISEBOUNDPROJECTEDLIEDERIVATIVESANGDIFFCOORDINATEX} 
		to bound this term by
		$\lesssim 
			\left|
				\angdiv \angLie_{\mathscr{Z}}^{\leq N-1} \upchi^{(Small)}
			\right|$
		as desired.
		Using 
		the schematic identity
		$\angLie_Z \ginversesphere = - (\ginversesphere)^2 \angdeform{Z},$
		Lemma~\ref{L:LANDRADCOMMUTEWITHANGDIFF} 
		and
		inequalities \eqref{E:ONEFORMANGDINTERMSOFROTATIONALLIE}
		and
		\eqref{E:TYPE02TENSORANGDINTERMSOFROTATIONALLIE},
		we see that the remaining terms that arise in the Leibniz expansion are bounded in magnitude by
		\begin{align} \label{E:ANGLAPLJUNKIFIRSTPRODUCTLOWERORDERTERMESTIMATE}
			\lesssim 
			\frac{1}{1 + t}
			\sum_{k=0}^N
			\mathop{\sum_{N_1 + \cdots + N_{k+2} \leq N-k}}_{N_{k+1} \leq N-1}
			\left|
				\angLie_{\mathscr{Z}}^{N_1} \angpi
			\right|
			\left|
				\angLie_{\mathscr{Z}}^{N_2} \angpi
			\right|
			\cdots
			\left|
				\angLie_{\mathscr{Z}}^{N_k} \angpi
			\right|	
			\left|
				\angLie_{\mathscr{Z}}^{N_{k+1}} \upchi^{(Small)}
			\right|
			\left|
				\angdiff \mathscr{Z}^{N_{k+2}}x^i
			\right|,
		\end{align}
		where the $\angpi$ are the $S_{t,u}-$projections of deformation tensors of vectorfields belonging to $\mathscr{Z}.$
		From the estimates
		\eqref{E:POINTWISEBOUNDPROJECTEDLIEDERIVATIVESANGDIFFCOORDINATEX},
		\eqref{E:LOWERORDERPOINTWISEBOUNDPROJECTEDLIEDERIVATIVESANGDIFFCOORDINATEX},	
		\eqref{E:CRUDEPOINTWISEBOUNDSDERIVATIVESOFANGULARDEFORMATIONTENSORS},
		\eqref{E:CRUDELOWERORDERC0BOUNDDERIVATIVESOFANGULARDEFORMATIONTENSORS},
		\eqref{E:POINTWISEESTIMATESFORCHIJUNKINTERMSOFOTHERVARIABLES},
		and
		\eqref{E:LOWERORDERC0BOUNDCHIJUNK},
		we deduce that the right-hand side of \eqref{E:ANGLAPLJUNKIFIRSTPRODUCTLOWERORDERTERMESTIMATE} is 
		$\lesssim$
		the right-hand side of \eqref{E:POINTWISEBOUNDFORANGLAPLJUNKICOMMUTED} as desired.
		
		We now address the second product on the right-hand side of \eqref{E:ANGDIFFLJUNKI}:
		$\angD^A 
		\angLie_{\mathscr{Z}}^{\leq N-1} 
		\left\lbrace 
			\uplambda_A^{(Tan-\Psi)} \Radunit_{(Small)}^i 
		\right\rbrace.$ From the identity 
		\eqref{E:RADUNITJUNKLIKELMINUSUNITJUNK}, definition \eqref{E:LITTLELAMBDAGOOD},
		Lemma~\ref{L:LANDRADCOMMUTEWITHANGDIFF},
		inequalities
		\eqref{E:FUNCTIONPOINTWISEANGDINTERMSOFANGLIEO} 
		and
		\eqref{E:ONEFORMANGDINTERMSOFROTATIONALLIE}, and the Leibniz rule, it follows that
		this term is bounded in magnitude by
		\begin{align} \label{E:ANGLAPLJUNKISECONDPRODUCTESTIMATE}
			& \lesssim 
				\frac{1}{1 + t}
				\sum_{a=1}^3
				\sum_{N_1 + N_2 + N_3 \leq N}
				\left|
					\angLie_{\mathscr{Z}}^{N_1} G_{(Frame)}
				\right|
				\left|
					\angdiff \mathscr{Z}^{N_2} \Psi
				\right|
				\left|
					\mathscr{Z}^{N_3} \Lunit_{(Small)}^a
				\right|
					\\
				&
				\ \ + \frac{1}{1 + t}
				\sum_{N_1 + N_2 + N_3 \leq N}
				\left|
					\angLie_{\mathscr{Z}}^{N_1} G_{(Frame)}
				\right|
				\left|
					\angdiff \mathscr{Z}^{N_2} \Psi
				\right|
				\left|
					\mathscr{Z}^{N_3} \Psi
				\right|.
				\notag
		\end{align}
		From the bootstrap assumptions \eqref{E:PSIFUNDAMENTALC0BOUNDBOOTSTRAP}
		and the estimates
		\eqref{E:LIEDERIVATIVESOFGRAMEINTERMSOFOTHERVARIABLES},
		\eqref{E:LOWERORDERC0BOUNDLIEDERIVATIVESOFGRAME},
		and
		\eqref{E:C0BOUNDCRUCIALEIKONALFUNCTIONQUANTITIES},
		we deduce that the right-hand side of \eqref{E:ANGLAPLJUNKISECONDPRODUCTESTIMATE} is $\lesssim$
		the right-hand side of \eqref{E:POINTWISEBOUNDFORANGLAPLJUNKICOMMUTED} as desired.	
	
		We now address the third product on the right-hand side of \eqref{E:ANGDIFFLJUNKI}:
		$\angD^A 
		\angLie_{\mathscr{Z}}^{\leq N-1} 
		\left\lbrace 
			\frac{x^i}{\rgeo} \uplambda_A^{(Tan-\Psi)} 
		\right\rbrace.$ From definition \eqref{E:LITTLELAMBDAGOOD},
		inequality \eqref{E:ONEFORMANGDINTERMSOFROTATIONALLIE}, and the Leibniz rule, it follows that
		this term is bounded in magnitude by
		\begin{align} \label{E:ANGLAPLJUNKITHIRDPRODUCTESTIMATE}
			& \lesssim 
				\frac{1}{1 + t}
				\sum_{N_1 + N_2 + N_3 \leq N}
				\left|
					\mathscr{Z}^{N_1} \left(\frac{x^i}{\rgeo}\right)
				\right|
				\left|
					\angLie_{\mathscr{Z}}^{N_2} G_{(Frame)}
				\right|
				\left|
					\angdiff \mathscr{Z}^{N_3} \Psi
				\right|.
		\end{align}
		From the bootstrap assumptions \eqref{E:PSIFUNDAMENTALC0BOUNDBOOTSTRAP},
		the identity \eqref{E:ZNAPPLIEDTORGEOISNOTTOOLARGE},
		the estimates
		\eqref{E:EASYCOORDINATEBOUND},
		\eqref{E:DERIVATIVESOFRECTANGULARCOMMUTATORVECTORFIELDRECTCOMPONENTS},
		\eqref{E:LOWERORDERC0BOUNDDERIVATIVESOFRECTANGULARCOMMUTATORVECTORFIELDRECTCOMPONENTS},		
		\eqref{E:LIEDERIVATIVESOFGRAMEINTERMSOFOTHERVARIABLES},
		\eqref{E:LOWERORDERC0BOUNDLIEDERIVATIVESOFGRAME},
		and
		\eqref{E:C0BOUNDCRUCIALEIKONALFUNCTIONQUANTITIES},
		and the fact that $Z x^i = Z^i,$
		we deduce that the right-hand side of \eqref{E:ANGLAPLJUNKITHIRDPRODUCTESTIMATE} is $\lesssim$
		the right-hand side of \eqref{E:POINTWISEBOUNDFORANGLAPLJUNKICOMMUTED} as desired.
		
		We now address the final product on the right-hand side of \eqref{E:ANGDIFFLJUNKI}:
		$\angD^A 
		\angLie_{\mathscr{Z}}^{\leq N-1} 
		\left\lbrace 
			(\ginversesphere)^{BC} \Lambda_{AB}^{(Tan-\Psi)} \angdiffarg{C} x^i 
		\right\rbrace.$ From definition \eqref{E:BIGLAMBDAGOOD},
		Lemma~\ref{L:LANDRADCOMMUTEWITHANGDIFF},
		inequality \eqref{E:ONEFORMANGDINTERMSOFROTATIONALLIE},
		and the Leibniz rule, we argue as in the proof of
		\eqref{E:ANGLAPLJUNKIFIRSTPRODUCTLOWERORDERTERMESTIMATE} to bound this term in magnitude by
		\begin{align} \label{E:ANGLAPLJUNKIFOURTHPRODUCTESTIMATE}
			& \lesssim 
				\frac{1}{1 + t}
				\sum_{k=0}^N
				\sum_{N_1 + \cdots + N_{k+3} \leq N-k}
				\left|
				\angLie_{\mathscr{Z}}^{N_1} \angpi
			\right|
			\left|
				\angLie_{\mathscr{Z}}^{N_2} \angpi
			\right|
			\cdots
			\left|
				\angLie_{\mathscr{Z}}^{N_k} \angpi
			\right|	
			\left|
				\angLie_{\mathscr{Z}}^{N_{k+1}} G_{(Frame)}
			\right|
				\left|
					\angLie_{\mathscr{Z}}^{N_{k+2}}
					\myarray[\Lunit \Psi]
						{\angdiff \Psi}
				\right|
				\left|
					\angdiff \mathscr{Z}^{N_{k+3}} x^i
				\right|.
		\end{align}
		From the bootstrap assumptions \eqref{E:PSIFUNDAMENTALC0BOUNDBOOTSTRAP}
		and the estimates
		\eqref{E:FUNCTIONDERIVATIVESAVOIDINGCOMMUTING},
		\eqref{E:POINTWISEBOUNDPROJECTEDLIEDERIVATIVESANGDIFFCOORDINATEX},
		\eqref{E:LOWERORDERPOINTWISEBOUNDPROJECTEDLIEDERIVATIVESANGDIFFCOORDINATEX},	
		\eqref{E:CRUDEPOINTWISEBOUNDSDERIVATIVESOFANGULARDEFORMATIONTENSORS},
		\eqref{E:CRUDELOWERORDERC0BOUNDDERIVATIVESOFANGULARDEFORMATIONTENSORS},
		\eqref{E:LOWERORDERC0BOUNDLIEDERIVATIVESOFGRAME},
		and
		\eqref{E:C0BOUNDCRUCIALEIKONALFUNCTIONQUANTITIES},
		we deduce that the right-hand side of \eqref{E:ANGLAPLJUNKIFOURTHPRODUCTESTIMATE} is $\lesssim$
		the right-hand side of \eqref{E:POINTWISEBOUNDFORANGLAPLJUNKICOMMUTED} as desired.
\end{proof}	

\section{Estimates related to integrals over \texorpdfstring{$S_{t,u}$}{the spheres}}
In this section, we provide a collection of estimates that are connected to 
integrals over $S_{t,u}.$

\begin{lemma}[\textbf{Estimates for the error terms in} $\angLie_{\Rad} d \spherevol$]
	\label{L:LIERADSTUAREAFORMERRORTERMSPOINTWISEESTIMATE}
	Under the 
	small data and
	bootstrap assumptions 
	of Sects.~\ref{S:PSISOLVES}-\ref{S:C0BOUNDBOOTSTRAP},
	if $\varepsilon$ is sufficiently small, 
	then the last four terms in braces in equation \eqref{E:UDERIVSTU}	
	verify the following pointwise estimate on $\mathcal{M}_{\Tboot,U_0}:$	
	\begin{align} \label{E:LIERADSTUAREAFORMERRORTERMSPOINTWISEESTIMATE}
		\left|
			- \frac{2}{\rgeo} (\upmu - 1)
			- \upmu \mytr \upchi^{(Small)}
			+ \mytr  \angkuparg{(Trans-\Psi)}
			+ \upmu \mytr  \angkuparg{(Tan-\Psi)}
		\right|
		& \lesssim \varepsilon \frac{\ln(\myexp + t)}{(1 + t)}.
	\end{align}
	
\end{lemma}

\begin{proof}
The last two terms on the left-hand side of \eqref{E:LIERADSTUAREAFORMERRORTERMSPOINTWISEESTIMATE} have the schematic form
$G_{(Frame)}^{\#} \threemyarray[\upmu \Lunit \Psi]{\Rad \Psi}{ \upmu \angdiff \Psi}.$
Hence, the desired estimate \eqref{E:LIERADSTUAREAFORMERRORTERMSPOINTWISEESTIMATE}
follows from 
Cor.~\ref{C:POINTWISEESTIMATESGFRAMESHARPINTERMSOFOTHERQUANTITIES},
the estimates 
\eqref{E:LOWERORDERC0BOUNDLIEDERIVATIVESOFGRAME}
and
\eqref{E:C0BOUNDCRUCIALEIKONALFUNCTIONQUANTITIES},
and the bootstrap assumptions \eqref{E:PSIFUNDAMENTALC0BOUNDBOOTSTRAP}.
\end{proof}

\begin{lemma}[\textbf{Comparison of the area forms} $d \spherevol$ \textbf{and} $d \Eucspherevol$]
\label{L:STUINTEGRALCOMPARISON}
Let $p=p(\vartheta)$ be a non-negative function of the geometric angular coordinates. 
Under the 
small data and
bootstrap assumptions 
of Sects.~\ref{S:PSISOLVES}-\ref{S:C0BOUNDBOOTSTRAP},
if $\varepsilon$ is sufficiently small, 
then the following estimate holds for $(t,u) \in [0,\Tboot) \times [0,U_0]:$
\begin{align} \label{E:SPHEREVOLUMEFORMCOMPARISON}
	(1 - C \varepsilon) \rgeo^2(t,u) \int_{\vartheta \in \mathbb{S}^2} p(\vartheta) d \argEucspherevol{(\vartheta)}
	\leq
	\int_{S_{t,u}} p(\vartheta) d \argspherevol{(t,u,\vartheta)}
	& \leq (1 + C \varepsilon) \rgeo^2(t,u) \int_{\vartheta \in \mathbb{S}^2} p(\vartheta) d \argEucspherevol{(\vartheta)},
\end{align}
where $d \argEucspherevol{(\vartheta)}$ denotes the standard Euclidean area form
on the \textbf{Euclidean-unit} sphere $\mathbb{S}^2.$
\end{lemma}

\begin{proof}
We set $P(t,u) := \int_{S_{t,u}} p(\vartheta) d \argspherevol{(t,u,\vartheta)}.$
From Lemma~\ref{L:STUDERIVATIVES} and the decomposition $\mytr \upchi = 2 \rgeo^{-1} + \mytr \upchi^{(Small)},$
 we have
\begin{align}
				\frac{\partial}{\partial t} 
				P(t,u)
				& = \int_{S_{t,u}} \left\lbrace \frac{2}{\rgeo} + \mytr \upchi^{(Small)} \right\rbrace p
					\, d \spherevol,
				\label{E:PTDERIVATIVEIDENTITY} \\
			\frac{\partial}{\partial u} 
				P(t,u)
				& = \int_{S_{t,u}} 
						 \left\lbrace
								- \frac{2}{\rgeo} 
								- \frac{2}{\rgeo} (\upmu - 1)
								- \upmu \mytr \upchi^{(Small)}
								+ \mytr  \angkuparg{(Trans-\Psi)}
								+ \upmu \mytr  \angkuparg{(Tan-\Psi)}
							\right\rbrace
						p
					\, d \spherevol.
					\label{E:PUDERIVATIVEIDENTITY}
		\end{align}
Inserting the estimates $|\mytr \upchi^{(Small)}| \lesssim \varepsilon \ln(\myexp + t) (1+t)^{-2}$
(that is, \eqref{E:C0BOUNDCRUCIALEIKONALFUNCTIONQUANTITIES}) 
and
\eqref{E:LIERADSTUAREAFORMERRORTERMSPOINTWISEESTIMATE} 
at $t = 0$
into
\eqref{E:PTDERIVATIVEIDENTITY} and \eqref{E:PUDERIVATIVEIDENTITY}
respectively,
we deduce that 
\begin{align}
	\left(
		1 - C \varepsilon \frac{\ln(\myexp + t)}{1 + t}
	\right)
	\frac{2}{\rgeo(t,u)} P(t,u)
	& \leq \frac{\partial}{\partial t} P(t,u)
	\leq 
		\left(
			1 + C \varepsilon \frac{\ln(\myexp + t)}{1 + t}
		\right) 
		\frac{2}{\rgeo(t,u)} P(t,u),
		\label{E:PTDERIVATIVEESTIMATE} \\
	- (1 + C \varepsilon)\frac{2}{\rgeo(0,u)} P(0,u)
	& \leq \frac{\partial}{\partial u} P(0,u)
	\leq - (1 - C \varepsilon) \frac{2}{\rgeo(0,u)} P(0,u).
		\label{E:PUDERIVATIVEESTIMATE}
\end{align}

We now apply Gronwall's inequality to 
\eqref{E:PTDERIVATIVEESTIMATE} and \eqref{E:PUDERIVATIVEESTIMATE} 
and recall that $\rgeo(t,u) = 1 - u + t,$
which leads to the following inequalities:
\begin{align}
	(1 - C \varepsilon) \frac{\rgeo^2(t,u)}{\rgeo^2(0,u)}
	& \leq \frac{P(t,u)}{P(0,u)}
	\leq  (1 + C \varepsilon) \frac{\rgeo^2(t,u)}{\rgeo^2(0,u)},
		\label{E:PATTUOVERPAT0UINEQUALITIES} \\
	(1 - C \varepsilon) \frac{\rgeo^2(0,u)}{\rgeo^2(0,0)}
	& \leq \frac{P(0,u)}{P(0,0)}
	\leq  (1 + C \varepsilon) \frac{\rgeo^2(0,u)}{\rgeo^2(0,0)}.
		\label{E:PAT0UOVERPAT00INEQUALITIES}
\end{align}

The desired estimate \eqref{E:SPHEREVOLUMEFORMCOMPARISON} now follows
easily from inequalities \eqref{E:PATTUOVERPAT0UINEQUALITIES} and \eqref{E:PAT0UOVERPAT00INEQUALITIES}
and the identities $\rgeo(0,0)=1$
and  
$P(0,0) = \int_{\vartheta \in \mathbb{S}^2} p(\vartheta) d \argEucspherevol{(\vartheta)}.$
\end{proof}

\begin{corollary}[\textbf{Estimates for integrals of constants}]
\label{C:ESTIMTATESFORINTEGRALSOFCONSTANTS}
Under the 
small data and
bootstrap assumptions 
of Sects.~\ref{S:PSISOLVES}-\ref{S:C0BOUNDBOOTSTRAP},
if $\varepsilon$ is sufficiently small, 
then the following estimates hold for $(t,u) \in [0,\Tboot) \times [0,U_0]:$
\begin{subequations}
\begin{align}
	4 \pi (1 - C \varepsilon) \rgeo^2(t,u)
	& \leq \int_{S_{t,u}} 1 \, d \spherevol
		\leq 4 \pi(1 + C \varepsilon) \rgeo^2(t,u),
		\\
	C^{-1} \rgeo^2(t,u)
	& \leq \int_{\Sigma_t^u} 1 \, d \tvol	
		\leq C \rgeo^2(t,u),
			\\
	C^{-1} \rgeo(t,u)
	& \leq \left\| 1 \right\|_{L^2(\Sigma_t^u)}	
		\leq C \rgeo(t,u).
		\label{E:ESTIMATEFORSPATIALL2NORMOF1}
\end{align}
\end{subequations}

\end{corollary}

\begin{proof}
Cor.~\ref{C:ESTIMTATESFORINTEGRALSOFCONSTANTS} now follows
easily from the definitions of the quantities involved, 
Lemma~\ref{L:STUINTEGRALCOMPARISON}, and the identity
$\int_{\vartheta \in \mathbb{S}^2} 1 d \argEucspherevol{(\vartheta)} = 4 \pi.$
\end{proof}

\begin{lemma}[\textbf{Pointwise estimates for} $\sqrt{\mbox{\upshape{det}} \gsphere}$]
\label{L:SQRTDETANGSPHEREMETRICPOINTWISEESTIMATE}
Let $\sqrt{\mbox{\upshape{det}} \gsphere}(t,u,\vartheta)$ denote the area form factor
of $\gsphere$ relative to the geometric coordinates,
and let $\sqrt{\mbox{\upshape{det}} \Eucsphereunit}(\vartheta)$ denote the area form 
factor of $\Eucsphereunit$ relative to the geometric coordinates.
Here, $\Eucsphereunit$ is the Riemannian metric
induced on the Euclidean unit sphere
by the Euclidean metric\footnote{Relative to the rectangular spatial coordinates on $\Sigma_t,$ 
we have $\Euct_{ij} = \delta_{ij}.$} 
$\Euct$ on $\Sigma_t.$
Under the 
small data and
bootstrap assumptions 
of Sects.~\ref{S:PSISOLVES}-\ref{S:C0BOUNDBOOTSTRAP},
if $\varepsilon$ is sufficiently small, 
then the following estimate holds on the spacetime domain 
$\mathcal{M}_{\Tboot,U_0}:$
\begin{align}  \label{E:SQRTDETANGSPHEREMETRICPOINTWISEESTIMATE}
	\left|
		\frac{\rgeo^{-2} \sqrt{\mbox{\upshape{det}} \gsphere}}
			{\sqrt{\mbox{\upshape{det}} \Eucsphereunit}}
		- 1
	\right|
	& \lesssim \varepsilon \frac{\ln(\myexp + t)}{1 + t}.
\end{align}
\end{lemma}

\begin{proof}
	Let $p \in \Sigma_t^{U_0}$ 
	and let $\gamma: [0,U_0] \rightarrow \Sigma_t^{U_0}$ be the integral curve of $\Rad$ 
	that passes through $p$ and that is parametrized by the eikonal function $u.$
	From the proof of Lemma~\ref{L:STUDERIVATIVES} and in particular
	\eqref{E:RADDERIVATIVESPHERVOLUMEFORM},
	and the identities
	$\frac{d}{du} \rgeo \circ \gamma = \Rad \rgeo = - 1$
	and
	$\Rad u = 1,$
	we deduce that
	\begin{align} \label{E:UDERIVATIVEOFRENORMALIZEDSPHEREAREAFORM}
		\frac{d}{du}
		\left\lbrace
			\ln [(\rgeo^{-2} \sqrt{\mbox{\upshape{det}} \gsphere})] \circ \gamma
		\right\rbrace
		& = 
			- \frac{2}{\rgeo} (\upmu - 1)
			- \upmu \mytr \upchi^{(Small)}
			+ \mytr  \angkuparg{(Trans-\Psi)}
			+ \upmu \mytr  \angkuparg{(Tan-\Psi)}.
	\end{align}
	Inserting the estimate \eqref{E:LIERADSTUAREAFORMERRORTERMSPOINTWISEESTIMATE} into the right-hand 
	side of \eqref{E:UDERIVATIVEOFRENORMALIZEDSPHEREAREAFORM}, we deduce that
	\begin{align} \label{E:UDERIVATIVEOFRENORMALIZEDSPHEREAREAFORMESTIMATED}
		\left|
			\frac{d}{du}	
			\left\lbrace
				\ln [(\rgeo^{-2} \sqrt{\mbox{\upshape{det}} \gsphere})] \circ \gamma
			\right\rbrace	
		\right|
		& \lesssim \varepsilon \frac{\ln(\myexp + t)}{1 + t}.
	\end{align}
	Integrating \eqref{E:UDERIVATIVEOFRENORMALIZEDSPHEREAREAFORMESTIMATED} $du$ from the point $p$
	until the initial sphere $S_{t,0}$ (where $u = 0$) and using the fact that
	$\rgeo^{-2} \sqrt{\mbox{\upshape{det}} \gsphere}|_{S_{t,0}} = \sqrt{\mbox{\upshape{det}} \Eucsphereunit},$ 
	we conclude the desired estimate \eqref{E:SQRTDETANGSPHEREMETRICPOINTWISEESTIMATE}.
\end{proof}

\begin{lemma}[\textbf{Pointwise estimates for the time integral of} $\mytr \upchi$]
\label{L:POINTWISEESTIMATEFORTRCHIANTIDERIVATIVE}
Under the 
small data and
bootstrap assumptions 
of Sects.~\ref{S:PSISOLVES}-\ref{S:C0BOUNDBOOTSTRAP},
if $\varepsilon$ is sufficiently small, 
then the following estimate holds for $(t,u) \in [0,\Tboot) \times [0,U_0]:$
\begin{align}  \label{E:POINTWISEESTIMATEFORTRCHIANTIDERIVATIVE}
\ln \left( \frac{\rgeo^2(t,u)}{\rgeo^2(0,u)} \right)
- C \varepsilon
\leq
\int_{s = 0}^t
	\mytr \upchi (s,u,\vartheta)
\, ds
\leq 
\ln \left( \frac{\rgeo^2(t,u)}{\rgeo^2(0,u)} \right)
+ C \varepsilon.
\end{align}

\end{lemma}

\begin{proof}
	The inequalities in \eqref{E:POINTWISEESTIMATEFORTRCHIANTIDERIVATIVE} follow easily
	from the decomposition $\mytr \upchi = 2 \rgeo^{-1} + \mytr \upchi^{(Small)}$
	and the pointwise estimate $|\mytr \upchi^{(Small)}| \leq C \varepsilon \ln(\myexp + t) (1+t)^{-2}$
	(that is, \eqref{E:C0BOUNDCRUCIALEIKONALFUNCTIONQUANTITIES}).
\end{proof}

\begin{lemma}[\textbf{Estimate for the norm} $\| \cdot \|_{L^2(\Sigma_t^u)}$ \textbf{of time-integrated functions}] 
\label{L:L2NORMSOFTIMEINTEGRATEDFUNCTIONS}
Let $f$ be a function on spacetime, and let
\begin{align}
	F(t,u,\vartheta) := \int_{t'=0}^t f(t',u,\vartheta) \, dt'.
\end{align}

Under the 
small data and
bootstrap assumptions 
of Sects.~\ref{S:PSISOLVES}-\ref{S:C0BOUNDBOOTSTRAP},
if $\varepsilon$ is sufficiently small, 
then the following estimate holds for $(t,u) \in [0,\Tboot) \times [0,U_0]:$
\begin{align} \label{E:L2NORMSOFTIMEINTEGRATEDFUNCTIONS}
	\| F \|_{L^2(\Sigma_t^u)} 
	& \leq (1 + C \varepsilon) \rgeo(t,u) 
		\int_{t'=0}^t 
			\frac{\| f \|_{L^2(\Sigma_{t'}^u)}}{\rgeo(t',u)} 
		\, dt'.
\end{align}

\end{lemma}

\begin{proof}
	By the definitions of $F$ and $\| F \|_{L^2(\Sigma_t^u)}^2,$ 
	the estimate \eqref{E:SPHEREVOLUMEFORMCOMPARISON}, 
	Minkowski's inequality for integrals, 
	and the fact that $\rgeo(t',u') \geq \rgeo(t',u)$ when $0 \leq u' \leq u,$
	we have
	\begin{align} \label{E:PROOFOFL2NORMSOFTIMEINTEGRATEDFUNCTIONS} 
		\| F \|_{L^2(\Sigma_t^u)}^2
		& = \int_{u'=0}^u
					\int_{S_{t,u}}
						F^2(t,u',\vartheta) 
					\, d \spherevol
				\, du'
				\\
	& \leq
				(1 + C \varepsilon) \rgeo^2(t,u)
				\int_{u'=0}^u
					\int_{\vartheta \in \mathbb{S}^2}
						F^2(t,u',\vartheta) 
					\, d \Eucspherevol
				\, du'
			\notag \\
		& =
			(1 + C \varepsilon) 
			\rgeo^2(t,u)
				\int_{u'=0}^u
				\int_{\vartheta \in \mathbb{S}^2}
				\left(
					\int_{t'=0}^t		
						f(t',u',\vartheta) 
					\, dt'
				\right)^2
				\, d \Eucspherevol
				\, du'
			\notag 
			\\
		& \leq
			(1 + C \varepsilon) 
			\rgeo^2(t,u)
			\left\lbrace
			\int_{t'=0}^t
				\left(
				\int_{u'=0}^u
					\int_{\vartheta \in \mathbb{S}^2}
						f^2(t',u',\vartheta) 
					\, d \Eucspherevol
				\, du'
				\right)^{1/2}	
			\, dt'
			\right\rbrace^2
			\notag 
			\\
			& \leq
				(1 + C \varepsilon) \rgeo^2(t,u)
				\left\lbrace
				\int_{t'=0}^t
					\left(
					\int_{u'=0}^u
						\int_{S_{t',u'}}
							f^2(t',u',\vartheta)
						\, \frac{d \spherevol}{\rgeo^2(t',u')}
						\, du'
						\right)^{1/2}	
				\, dt'
				\right\rbrace^2
				\notag \\
			& \leq
				(1 + C \varepsilon) \rgeo^2(t,u)
				\left\lbrace
				\int_{t'=0}^t
					\left(
					\frac{1}{\rgeo^2(t',u)}
					\int_{u'=0}^u
						\int_{S_{t',u'}}
							f^2(t',u',\vartheta)
						\, d \spherevol
						\, du'
						\right)^{1/2}	
				\, dt'
				\right\rbrace^2
				\notag \\
			& = 
				(1 + C \varepsilon) 
				\rgeo^2(t,u)
				\left\lbrace
				\int_{t'=0}^t
						\frac{1}{\rgeo(t',u)}
						\left\| 
							f
						\right\|_{L^2(\Sigma_{t'}^u)}
				\, dt'
				\right\rbrace^2.
				\notag
	\end{align}
	Inequality \eqref{E:L2NORMSOFTIMEINTEGRATEDFUNCTIONS} now follows from taking 
	the square root of both sides of \eqref{E:PROOFOFL2NORMSOFTIMEINTEGRATEDFUNCTIONS}.
\end{proof}

\section{Faster than expected decay for certain \texorpdfstring{$\Psi-$}{wave solution-}related quantities}

We use the estimates in the following lemma
in order to bound some of the error terms that appear 
in the integration by parts identities
of Lemmas \ref{L:TOPORDERMORAWETZREORMALIZEDTRCHIIBP} and
\ref{L:TOPORDERMORAWETZREORMALIZEDANGLAPUPMUIBP}.

\begin{lemma}[\textbf{Faster than expected decay for certain} $\Psi-$\textbf{related quantities}]
\label{L:FASTERTHANEXPECTEDPSIDECAY}
Recall that $\Psi$ verifies the covariant wave equation 
$\upmu \square_{g(\Psi)} \Psi = 0$ and that
$\Psi$ vanishes in the exterior of the outgoing null cone $\mathcal{C}_0.$
Under the 
small data and
bootstrap assumptions 
of Sects.~\ref{S:PSISOLVES}-\ref{S:C0BOUNDBOOTSTRAP},
if $\varepsilon$ is sufficiently small, then 
the following estimate holds for $t \in [0,\Tboot):$
\begin{align} \label{E:LOFRGEOWEIGHTEDZPSISTRONGC0BOUND}
	\left\| 
		\Lunit (\rgeo \Rad \Psi)
	\right\|_{C^0(\Sigma_t^u)},
	\qquad
	\left\| 
		\Lunit (\rgeo \Lunit \Psi)
	\right\|_{C^0(\Sigma_t^u)}
	& \lesssim \varepsilon \frac{\ln(\myexp + t)}{(1 + t)^2},
\end{align}

\begin{align} \label{E:LPSIPLUSHALFTRCHIPSISTRONGC0BOUND}
	\left\| 
		\mathscr{O}^{\leq 1} 
		\left(
			\Lunit \Psi 
			+ \frac{1}{2} \mytr \upchi \Psi 
		\right)
	\right\|_{C^0(\Sigma_t^u)},
	\qquad
	\left\| 
			\mathscr{O}^{\leq 1}
			\left(
				\Lunit \Rad \Psi + \frac{1}{2} \mytr \upchi \Rad \Psi 
			\right)
	\right\|_{C^0(\Sigma_t^u)}
		& \lesssim \varepsilon \frac{\ln(\myexp + t)}{(1 + t)^3}.
\end{align}

\end{lemma}

\begin{proof}
We first prove \eqref{E:LPSIPLUSHALFTRCHIPSISTRONGC0BOUND} for
$\Lunit \Psi + \frac{1}{2} \mytr \upchi \Psi.$ 
Writing the right-hand side of
the wave operator decomposition \eqref{E:WAVEOPDECOMPRADRGEOLPSIPLUSHALFTRCHIPSI}
(with $f = \Psi$)
in schematic form, we have
\begin{align} \label{E:SCHEMATICWAVEOPDECOMPRADRGEOLPSIPLUSHALFTRCHIPSI}
\Rad 
\left\lbrace
	\rgeo 
	\left[
		\Lunit \Psi
		+ \frac{1}{2} \mytr \upchi \Psi
	\right]
\right\rbrace				
& = \myarray
			[\upmu \Lunit (\rgeo \Lunit \Psi)]
			{\upmu \rgeo \angLap \Psi}
		+ \rgeo
			G_{(Frame)}
			\ginversesphere
			\threemyarray
				[\upmu \Lunit \Psi]
				{\Rad \Psi}	
				{\upmu \angdiff \Psi}	
			\myarray
				[\Lunit \Psi]
				{\angdiff \Psi}
				\\
& \ \ +
				\myarray
					[\upmu]
					{1}
					\Lunit \Psi
				+ \rgeo \ginversesphere (\angdiff \upmu) \angdiff \Psi
				+ \myarray
						[\rgeo \Rad \mytr \upchi^{(Small)}]
						{\mytr \upchi^{(Small)}}	
					\Psi.
					\notag
\end{align}
From \eqref{E:SCHEMATICWAVEOPDECOMPRADRGEOLPSIPLUSHALFTRCHIPSI}, 
the inequalities
\eqref{E:FUNCTIONPOINTWISEANGDINTERMSOFANGLIEO}
and
\eqref{E:ANGLAPFUNCTIONPOINTWISEINTERMSOFROTATIONS},
the estimates \eqref{E:LOWERORDERC0BOUNDLIEDERIVATIVESOFGRAME},
\eqref{E:CRUDELOWERORDERC0BOUNDDERIVATIVESOFANGULARDEFORMATIONTENSORS},
and
\eqref{E:C0BOUNDCRUCIALEIKONALFUNCTIONQUANTITIES},
and the bootstrap assumptions \eqref{E:PSIFUNDAMENTALC0BOUNDBOOTSTRAP}, 
we deduce that
\begin{align} \label{E:RADRGEOLUNITPSIPLUSHAFLTRCHIPSIPOINTWISEBOUND}
\left|
	\Rad 
		\left\lbrace
			\rgeo 
			\left[
				\Lunit \Psi
				+ \frac{1}{2} \mytr \upchi \Psi
			\right]
		\right\rbrace	
\right|
& \lesssim \varepsilon \frac{\ln(\myexp + t)}{(1 + t)^2}.
\end{align}
Fixing $t,$ integrating \eqref{E:RADRGEOLUNITPSIPLUSHAFLTRCHIPSIPOINTWISEBOUND} along the integral curves of
$- \Rad$ until reaching the sphere $S_{t,0},$ recalling that $\Rad u = 1,$
and recalling that $\Psi \equiv 0$ in the exterior 
of the outgoing null cone portion $\mathcal{C}_0^{\Tboot},$ we deduce that
\begin{align}
	\left| \label{E:RGEOLUNITPSIPLUSHAFLTRCHIPSIPOINTWISEBOUND}
			\rgeo 
			\left[
				\Lunit \Psi
				+ \frac{1}{2} \mytr \upchi \Psi
			\right]
		\right|
& \lesssim \varepsilon \frac{\ln(\myexp + t)}{(1 + t)^2}.
\end{align}
The desired estimate \eqref{E:LPSIPLUSHALFTRCHIPSISTRONGC0BOUND} for
$\Lunit \Psi + \frac{1}{2} \mytr \upchi \Psi$ now follows from
dividing inequality \eqref{E:RGEOLUNITPSIPLUSHAFLTRCHIPSIPOINTWISEBOUND} by $\rgeo.$

To prove the desired estimate 
\eqref{E:LPSIPLUSHALFTRCHIPSISTRONGC0BOUND} for
$\Rot(\Lunit \Psi + \frac{1}{2} \mytr \upchi \Psi),$
we note that $\Rot(\Lunit \Psi + \frac{1}{2} \mytr \upchi \Psi)$
verifies an equation similar to the equation
\eqref{E:SCHEMATICWAVEOPDECOMPRADRGEOLPSIPLUSHALFTRCHIPSI}
verified by $\Lunit \Psi + \frac{1}{2} \mytr \upchi \Psi,$
but with the right-hand side of \eqref{E:SCHEMATICWAVEOPDECOMPRADRGEOLPSIPLUSHALFTRCHIPSI}
replaced by $\Rot$ applied to the right-hand side of \eqref{E:SCHEMATICWAVEOPDECOMPRADRGEOLPSIPLUSHALFTRCHIPSI},
and with the additional commutator term (see \eqref{E:RADCOMMUTETANGENTISTANGENT})
$\rgeo \angdeformoneformupsharparg{\Rot}{\Rad} \cdot \angdiff (\Lunit \Psi + \frac{1}{2} \mytr \upchi \Psi)$ on the right-hand side.
We now claim that the magnitudes of the terms arising from $\Rot$ applied to the right-hand side of \eqref{E:SCHEMATICWAVEOPDECOMPRADRGEOLPSIPLUSHALFTRCHIPSI}
are $\lesssim \varepsilon \ln(\myexp + t) (1 + t)^{-2}.$
To prove this claim, we use 
Lemma~\ref{L:LANDRADCOMMUTEWITHANGDIFF},
the same estimates that we used to deduce \eqref{E:RADRGEOLUNITPSIPLUSHAFLTRCHIPSIPOINTWISEBOUND}, 
and, to bound the term $\Rot \angLap \Psi,$ the estimate \eqref{E:LINFFTYLOWERORDERCOMMUTINGANGDSQUAREDPSIWITHLIEZN}. 
Furthermore, using the decomposition $\mytr \upchi = 2 \rgeo^{-1} + \mytr \upchi^{(Small)},$
inequality \eqref{E:FUNCTIONPOINTWISEANGDINTERMSOFANGLIEO},
the estimates \eqref{E:LOWERORDERC0BOUNDROTDEFORMSPHERERAD}
and \eqref{E:C0BOUNDCRUCIALEIKONALFUNCTIONQUANTITIES},
and the bootstrap assumptions \eqref{E:PSIFUNDAMENTALC0BOUNDBOOTSTRAP}, we bound the commutator term as follows:
\begin{align}
	\rgeo 
	\left|
		 \angdeformoneformupsharparg{\Rot}{\Rad}
		 \cdot
		 \angdiff
		(\Lunit \Psi + \frac{1}{2} \mytr \upchi \Psi)
	\right|
	& \lesssim 
		\left| \angdeformoneformarg{\Rot}{\Rad} \right|
		\left\lbrace
			\sum_{l=1}^3 
			\left(
				|\Rot_{(l)} \Lunit \Psi|
				+ |\Rot_{(l)} \mytr \upchi^{(Small)}| |\Psi|
				+ \frac{1}{\rgeo} |\Rot_{(l)} \Psi| 
			\right)
		\right\rbrace
			\\
	& \ \ \lesssim \varepsilon \frac{\ln^2(\myexp + t)}{(1 + t)^3}.
		\notag
\end{align}
Combining these estimates, we deduce that inequality \eqref{E:RADRGEOLUNITPSIPLUSHAFLTRCHIPSIPOINTWISEBOUND} holds with
$			\rgeo 
			\left[
				\Lunit \Psi
				+ \frac{1}{2} \mytr \upchi \Psi
			\right]
$
replaced by 
$\Rot \left\lbrace
			\rgeo 
			\left[
				\Lunit \Psi
				+ \frac{1}{2} \mytr \upchi \Psi
			\right]
		\right\rbrace,	
$
and the desired estimate \eqref{E:LPSIPLUSHALFTRCHIPSISTRONGC0BOUND}
for 
$\Rot \left\lbrace
			\left[
				\Lunit \Psi
				+ \frac{1}{2} \mytr \upchi \Psi
			\right]
		\right\rbrace$
follows by integrating along the integral curves of
$- \Rad$ as above.

To prove estimate \eqref{E:LPSIPLUSHALFTRCHIPSISTRONGC0BOUND} for
$\Lunit \Rad \Psi + \frac{1}{2} \mytr \upchi \Rad \Psi,$
we first use the wave operator decomposition 
\eqref{E:LONOUTSIDEGEOMETRICWAVEOPERATORFRAMEDECOMPOSED}
with $f = \Psi,$
Lemma~\ref{L:UPMUFIRSTTRANSPORT},
and the simple identity 
$\Lunit \Lunit \Psi = \frac{1}{\rgeo} \Lunit (\rgeo \Lunit \Psi) - \frac{1}{\rgeo} \Lunit \Psi$
to deduce the schematic equation
\begin{align} \label{E:SCHEMATICWAVEOPDECOMPLRADPSIPLUSHALFTRCHIRADPSI}
\Lunit \Rad \Psi + \frac{1}{2} \mytr \upchi \Rad \Psi			
& = \myarray
			[\frac{1}{\rgeo} \upmu \Lunit (\rgeo \Lunit \Psi)]
			{\upmu \angLap \Psi}
		+  \frac{1}{\rgeo} \upmu \Lunit \Psi
		+ G_{(Frame)}
			\ginversesphere
			\threemyarray
				[\upmu \Lunit \Psi]
				{\Rad \Psi}	
				{\upmu \angdiff \Psi}	
			\myarray
				[\Lunit \Psi]
				{\angdiff \Psi}.
\end{align}
Using the same estimates we used to prove \eqref{E:RADRGEOLUNITPSIPLUSHAFLTRCHIPSIPOINTWISEBOUND}
and also \eqref{E:ZNAPPLIEDTORGEOISNOTTOOLARGE},
we bound the right-hand side 
of \eqref{E:SCHEMATICWAVEOPDECOMPLRADPSIPLUSHALFTRCHIRADPSI}
in magnitude by $\lesssim \varepsilon \ln(\myexp + t) (1 + t)^{-3}.$
We have thus proved the desired estimate \eqref{E:LPSIPLUSHALFTRCHIPSISTRONGC0BOUND} for
$\Lunit \Rad \Psi + \frac{1}{2} \mytr \upchi \Rad \Psi.$

To prove the estimate \eqref{E:LPSIPLUSHALFTRCHIPSISTRONGC0BOUND} for
$\Rot \left\lbrace \Lunit \Rad \Psi + \frac{1}{2} \mytr \upchi \Rad \Psi \right\rbrace,$
we apply $\Rot$ to both sides of \eqref{E:SCHEMATICWAVEOPDECOMPLRADPSIPLUSHALFTRCHIRADPSI}.
Using 
Lemma~\ref{L:LANDRADCOMMUTEWITHANGDIFF},
the same estimates we used to bound the right-hand side 
of \eqref{E:SCHEMATICWAVEOPDECOMPLRADPSIPLUSHALFTRCHIRADPSI},
and the estimate \eqref{E:LINFFTYLOWERORDERCOMMUTINGANGDSQUAREDPSIWITHLIEZN} to bound the term $\Rot \angLap \Psi,$ 
we conclude that
$\left| 
	\Rot 
	\left\lbrace 
		\Lunit \Rad \Psi 
		+ \frac{1}{2} \mytr \upchi \Rad \Psi 
	\right\rbrace
\right| 
\lesssim \varepsilon \ln(\myexp + t) (1 + t)^{-3}$
as desired.

The desired estimate \eqref{E:LOFRGEOWEIGHTEDZPSISTRONGC0BOUND} for 
$\Lunit (\rgeo \Rad \Psi)$ follows from the identity
$\Lunit (\rgeo \Rad \Psi) = \rgeo \left\lbrace \Lunit \Rad \Psi + \frac{1}{2} \mytr \upchi \Rad \Psi \right\rbrace 
- \frac{1}{2} \rgeo \mytr \upchi^{(Small)} \Rad \Psi,$
the estimate \eqref{E:LPSIPLUSHALFTRCHIPSISTRONGC0BOUND} for
$\Lunit \Rad \Psi + \frac{1}{2} \mytr \upchi \Rad \Psi,$
the estimate \eqref{E:C0BOUNDCRUCIALEIKONALFUNCTIONQUANTITIES} for $\mytr \upchi^{(Small)},$ 
and the bootstrap assumption $\| \Rad \Psi \|_{C^0(\Sigma_t^u)} \leq \varepsilon (1 + t)^{-1}.$

The desired estimate \eqref{E:LOFRGEOWEIGHTEDZPSISTRONGC0BOUND} for 
$\Lunit (\rgeo \Lunit \Psi)$ follows easily from the bootstrap assumption
$|\rgeo \Lunit (\rgeo \Lunit \Psi)| \leq \varepsilon (1 + t)^{-1}$
(that is, \eqref{E:PSIFUNDAMENTALC0BOUNDBOOTSTRAP}).
\end{proof}

\section{Pointwise estimates for \texorpdfstring{$\Xi$}{the vectorfield Xi}}
In this section, we derive pointwise estimates for the 
vectorfield $\Xi$ from the decomposition $\Rad = \frac{\partial}{\partial u} - \Xi.$

\begin{lemma}[\textbf{Estimates for} $\Xi$]
\label{L:XIESTIMATES}
Let $\Xi$ be the $S_{t,u}-$tangent vectorfield from \eqref{E:RADINTERMSOFGEOMETRICCOORDINATEPARTIALDERIVATIVES},
and let $\Xi_{\flat}$ be the corresponding $\gsphere-$dual one-form.
Let $0 \leq N \leq 23$ be an integer.
Under the 
small data and
bootstrap assumptions 
of Sects.~\ref{S:PSISOLVES}-\ref{S:C0BOUNDBOOTSTRAP},
if $\varepsilon$ is sufficiently small, then 
the following pointwise estimates hold on the spacetime domain 
$\mathcal{M}_{\Tboot,U_0}:$
\begin{subequations}
\begin{align}  \label{E:POINTWISEGEOMETRICXI}
	\left|
		\angLie_{\Lunit} 
		(\rgeo^{-2} \angLie_{\mathscr{Z}}^N \Xi_{\flat})
	\right|
		& \lesssim
			\frac{1}{(1 + t)^2}
			\left|
				\fourmyarray[\rgeo \Lunit \mathscr{Z}^{\leq N} \Psi]
					{\Rad \mathscr{Z}^{\leq N} \Psi}
					{\rgeo \angdiff \mathscr{Z}^{\leq N} \Psi} 
					{\mathscr{Z}^{\leq N} \Psi}
			\right|
		+ 
			\frac{1}{(1 + t)^3}
			\left|
				\myarray[\mathscr{Z}^{\leq N+1} (\upmu - 1)]
					{\sum_{a=1}^3 \rgeo |\mathscr{Z}^{\leq N+1} \Lunit_{(Small)}^a|} 
			\right|
				\\
		& \ \ 
			+  
			\frac{\ln(\myexp + t)}{(1 + t)^2}
			\left|
				\rgeo^{-2}
				\angLie_{\mathscr{Z}}^{\leq N} 
				\Xi_{\flat}
			\right|.
			\notag
\end{align}

Furthermore, the following estimate holds for $(t,u) \in [0,\Tboot) \times [0,U_0]:$
	\begin{align} \label{E:C0BOUNDSGEOMETRICXI}
		\left\|
			\angLie_{\mathscr{Z}}^{\leq 11} \Xi
		\right\|_{C^0(\Sigma_t^u)}
		& \lesssim
			\varepsilon (1 + t).
	\end{align}	

Furthermore, the following estimate holds for the rectangular spatial components $\Xi^i,$ $(i=1,2,3):$
	\begin{align} \label{E:C0BOUNDSXIRECTANGULARCOMPONENTS}
		\left\|
			\mathscr{Z}^{\leq 11} \Xi^i
		\right\|_{C^0(\Sigma_t^u)}
			& \lesssim 
				\varepsilon (1 + t).
	\end{align}
\end{subequations}
\end{lemma}

\begin{proof}
	We first lower the index in equation \eqref{E:XIEVOLUTION} with $\gsphere,$
	use the identities $\angLie_{\Lunit} \gsphere = 2 \rgeo^{-1} \gsphere + 2 \upchi^{(Small)}$
	(see \eqref{E:CHIALTDEF} and \eqref{E:CHIJUNKDEF})
	and
	$\Lunit \rgeo = 1,$
	use the fact that $[\Lunit, \Rad]_{\flat} = \angdeformoneformarg{\Rad}{\Lunit}$
	(see \eqref{E:LCOMMUTERADISSTUTANGENT}),
	and commute with $\rgeo^{-2}$ 
	followed by $\angLie_{\mathscr{Z}}^N$ 
	to deduce that
	\begin{align} \label{E:XICOMMUTEDEQN}
		\angLie_{\Lunit}
		(\rgeo^{-2} \angLie_{\mathscr{Z}}^N \Xi_{\flat})
		& = 
			2
			\angLie_{\mathscr{Z}}^N
			\left\lbrace
				\rgeo^{-2} \upchi^{(Small) \#} \cdot \Xi_{\flat}
			\right\rbrace
			- \angLie_{\mathscr{Z}}^N (\rgeo^{-2} \angdeformoneformarg{\Rad}{\Lunit})
			+ [\angLie_{\Lunit}, \angLie_{\mathscr{Z}}^N] (\rgeo^{-2} \Xi_{\flat}) 
			+
			\angLie_{\Lunit}
			\left\lbrace
				[\rgeo^{-2}, \angLie_{\mathscr{Z}}^N]
				\Xi_{\flat}
			\right\rbrace	.
	\end{align}
	
	We now derive bounds for the right-hand side of \eqref{E:XICOMMUTEDEQN}.
	To simplify the proof, we make the following temporary bootstrap assumptions:
	\begin{align} \label{E:TEMPBOOTXI} \tag{$\mathbf{TEMP-BA}\Xi_{\flat}$}
			\left\|
				\angLie_{\mathscr{Z}}^{\leq 11} \Xi_{\flat}
			\right\|_{C^0(\Sigma_t^u)}
			& \leq 1 + t,
			&& (t,u) \in [0,\Tboot) \times [0,U_0].
	\end{align}
	
	From the Leibniz rule,
	the estimates 
	\eqref{E:ZNAPPLIEDTORGEOISNOTTOOLARGE},
	\eqref{E:POINTWISEESTIMATESFORCHIJUNKINTERMSOFOTHERVARIABLES},
	and
	\eqref{E:C0BOUNDCRUCIALEIKONALFUNCTIONQUANTITIES},
	and the bootstrap assumptions \eqref{E:PSIFUNDAMENTALC0BOUNDBOOTSTRAP}
	and \eqref{E:TEMPBOOTXI},
	we deduce that 
	\begin{align} \label{E:XIESTIMATELOWERORDERJUNKTERM}
		\left|
			\angLie_{\mathscr{Z}}^N
			\left\lbrace
				\rgeo^{-2} \upchi^{(Small) \#} \cdot \Xi_{\flat}
			\right\rbrace
		\right|
		& \lesssim
			\frac{1}{(1 + t)^2}
			\left|
				\fourmyarray[\rgeo \Lunit \mathscr{Z}^{\leq N} \Psi]
					{\Rad \mathscr{Z}^{\leq N} \Psi}
					{\rgeo \angdiff \mathscr{Z}^{\leq N} \Psi} 
					{\mathscr{Z}^{\leq N} \Psi}
			\right|
			+ 
			\frac{1}{(1 + t)^3}
			\left|
				\myarray[\mathscr{Z}^{\leq N} (\upmu - 1)]
					{\sum_{a=1}^3 \rgeo |\mathscr{Z}^{\leq N+1} \Lunit_{(Small)}^a|} 
			\right|
				\\
		& \ \ 
			+
			\varepsilon
			\frac{\ln(\myexp + t)}{(1 + t)^2}
			\left|
				\rgeo^{-2} \angLie_{\mathscr{Z}}^{\leq N} \Xi_{\flat}
			\right|.
			\notag
	\end{align}
	From the Leibniz rule,
	the identity \eqref{E:RADDEFORMLA},
	Lemma~\ref{L:LANDRADCOMMUTEWITHANGDIFF},
	inequalities 
	\eqref{E:FUNCTIONPOINTWISEANGDINTERMSOFANGLIEO}
	and
	\eqref{E:FUNCTIONDERIVATIVESAVOIDINGCOMMUTING},
	the estimates
	\eqref{E:ZNAPPLIEDTORGEOISNOTTOOLARGE},
	\eqref{E:LIEDERIVATIVESOFGRAMEINTERMSOFOTHERVARIABLES},
	\eqref{E:LOWERORDERC0BOUNDLIEDERIVATIVESOFGRAME},
	and \eqref{E:C0BOUNDCRUCIALEIKONALFUNCTIONQUANTITIES},
	and the bootstrap assumptions
	\eqref{E:PSIFUNDAMENTALC0BOUNDBOOTSTRAP},
	it follows that
	\begin{align} \label{E:POINTWISEANGLIEZNAPPLIEDTOANGLIEXIINHOMOGENEOUSTERM}
		\left|
			\angLie_{\mathscr{Z}}^N 
			(\rgeo^{-2} \angdeformoneformarg{\Rad}{\Lunit})
		\right|
		& \lesssim
	 	\frac{1}{(1 + t)^2}
	 	\left|
				\fourmyarray[\rgeo \Lunit \mathscr{Z}^{\leq N} \Psi]
					{\Rad \mathscr{Z}^{\leq N} \Psi}
					{\rgeo \angdiff \mathscr{Z}^{\leq N} \Psi} 
					{\mathscr{Z}^{\leq N} \Psi}
		\right|
		+ 
			\frac{1}{(1 + t)^3}
			\left|
				\myarray[\mathscr{Z}^{\leq N+1} (\upmu - 1)]
					{\sum_{a=1}^3 \rgeo |\mathscr{Z}^{\leq N} \Lunit_{(Small)}^a|} 
			\right|.
	\end{align}
	In addition, using 
	inequality \eqref{E:LAPPLIEDTOCOMMUTATOROFZNANDRGEOAPPLIEDTOXIPOINTWISE},
	inequality \eqref{E:LZNCOMMUTATORACTINGONTENSORFIELDSPOINTWISE}
	with $\rgeo^{-2} \Xi_{\flat}$ in the role of $\xi$
	(and the first term 
$\left|
	\angLie_{\Lunit} \angLie_{\mathscr{Z}}^{N-1} (\rgeo^{-2} \Xi_{\flat})
\right|$
on the right-hand side of \eqref{E:LZNCOMMUTATORACTINGONTENSORFIELDSPOINTWISE}
is bounded by induction),
	Cor.~\ref{C:SQRTEPSILONREPLCEDWITHCEPSILON},
	the estimates \eqref{E:ZNAPPLIEDTORGEOISNOTTOOLARGE}
	and
	\eqref{E:C0BOUNDCRUCIALEIKONALFUNCTIONQUANTITIES},
	and the bootstrap assumptions
	\eqref{E:PSIFUNDAMENTALC0BOUNDBOOTSTRAP}
	and \eqref{E:TEMPBOOTXI},
	we deduce that
  \begin{align}  \label{E:POINTWISEANGLIEXICOMMUTATORTERM}
		\left|
			[\angLie_{\Lunit}, \angLie_{\mathscr{Z}}^N] 
			(\rgeo^{-2} \Xi_{\flat})
		\right|
		& \lesssim
		\frac{1}{(1 + t)^2}
		\left|
				\fourmyarray[\rgeo \Lunit \mathscr{Z}^{\leq N} \Psi]
					{\Rad \mathscr{Z}^{\leq N} \Psi}
					{\rgeo \angdiff \mathscr{Z}^{\leq N} \Psi} 
					{\mathscr{Z}^{\leq N} \Psi}
		\right|
		+ \frac{1}{(1 + t)^3}
			\left|
				\myarray[\mathscr{Z}^{\leq N+1} (\upmu - 1)]
					{\sum_{a=1}^3 \rgeo |\mathscr{Z}^{\leq N} \Lunit_{(Small)}^a|} 
			\right|
				\\
		& \ \ +
			\left|
				\angLie_{\Lunit} 
				(\rgeo^{-2} \angLie_{\mathscr{Z}}^{\leq N-1} \Xi_{\flat})
			\right|
			+  
			\frac{\ln(\myexp + t)}{(1 + t)^2}
			\left|
				\rgeo^{-2}
				\angLie_{\mathscr{Z}}^{\leq N} \Xi_{\flat}
			\right|.
			\notag
	\end{align}
	Furthermore, from \eqref{E:LAPPLIEDTOCOMMUTATOROFZNANDRGEOAPPLIEDTOXIPOINTWISE}, we deduce that
	\begin{align} \label{E:POINTWISEXITORSIONANNOYINGCOMMUTATORTERM}
		\left|
			\angLie_{\Lunit}
			\left\lbrace
				[\rgeo^{-2}, \angLie_{\mathscr{Z}}^N]
				\Xi_{\flat}
			\right\rbrace
		\right|
		& \lesssim 
			\left|
				\angLie_{\Lunit}
				(\rgeo^{-2} \angLie_{\mathscr{Z}}^{\leq N-1} \Xi_{\flat})
			\right|
			+ \frac{1}{(1 + t)^2}
			\left|
				\rgeo^{-2} \angLie_{\mathscr{Z}}^{\leq N} \Xi_{\flat}
			\right|.
	\end{align}
	Combining \eqref{E:XICOMMUTEDEQN}, 
	\eqref{E:XIESTIMATELOWERORDERJUNKTERM},
	\eqref{E:POINTWISEANGLIEZNAPPLIEDTOANGLIEXIINHOMOGENEOUSTERM},
	\eqref{E:POINTWISEANGLIEXICOMMUTATORTERM}, 
	and \eqref{E:POINTWISEXITORSIONANNOYINGCOMMUTATORTERM},
	and applying induction in $N$
	to handle the next-to-last term on the right-hand side of \eqref{E:POINTWISEANGLIEXICOMMUTATORTERM}
	and the first term on the right-hand side of \eqref{E:POINTWISEXITORSIONANNOYINGCOMMUTATORTERM},
	we deduce that
	\begin{align} \label{E:ALMOSTDONEPOINTWISEGEOMETRICXI}
		\left|
			\angLie_{\Lunit} 
			(\rgeo^{-2} \angLie_{\mathscr{Z}}^N \Xi_{\flat})
		\right|
		& \lesssim
			\frac{1}{(1 + t)^2}
			\left|
				\fourmyarray[\rgeo \Lunit \mathscr{Z}^{\leq N} \Psi]
					{\Rad \mathscr{Z}^{\leq N} \Psi}
					{\rgeo \angdiff \mathscr{Z}^{\leq N} \Psi} 
					{\mathscr{Z}^{\leq N} \Psi}
			\right|
		+ 
			\frac{1}{(1 + t)^3}
			\left|
				\myarray[\mathscr{Z}^{\leq N+1} (\upmu - 1)]
					{\sum_{a=1}^3 \rgeo |\mathscr{Z}^{\leq N+1} \Lunit_{(Small)}^a|} 
			\right|
				\\
		& \ \ 
			+  
			\frac{\ln(\myexp + t)}{(1 + t)^2}
			\left|
				\rgeo^{-2} \angLie_{\mathscr{Z}}^{\leq N} \Xi_{\flat}
			\right|.
			\notag
	\end{align}
	In particular, we have proved the desired bound \eqref{E:POINTWISEGEOMETRICXI}.
	
	We now use the identities 
	$\Lunit \rgeo = 1$
	and $\angLie_{\Lunit} \ginversesphere = - 2 \rgeo^{-1} \ginversesphere - 2 \upchi^{(Small)\# \#}$
	to deduce that for any $S_{t,u}$ one-form $\xi,$ we have
	$\Lunit |\rgeo \xi|^2 
		= 2 \rgeo^2 \upchi_{AB}^{(Small)} \xi^A \xi^B 
		+ 2 \rgeo^2 \xi^A \angLie_{\Lunit} \xi_A.$
	Using this identity with $\xi = \rgeo^{-2} \angLie_{\mathscr{Z}}^N \Xi_{\flat},$ 
	the Cauchy-Schwarz inequality, 
	the estimates
	\eqref{E:C0BOUNDCRUCIALEIKONALFUNCTIONQUANTITIES}
	and \eqref{E:ALMOSTDONEPOINTWISEGEOMETRICXI},
	and the bootstrap assumptions
	\eqref{E:PSIFUNDAMENTALC0BOUNDBOOTSTRAP},
	we deduce that
	\begin{align} \label{E:POINTWISEGEOMETRICXIGRONWALLREADY}
		\left|
		\Lunit 
		\left|
			\rgeo^{-1} \angLie_{\mathscr{Z}}^{\leq N} \Xi_{\flat}
		\right|
		\right|
		& \lesssim 
			\varepsilon \frac{\ln(\myexp + t)}{(1 + t)^2}
		+ \frac{\ln(\myexp + t)}{(1 + t)^2}
			\left\|
				\rgeo^{-1} \angLie_{\mathscr{Z}}^{\leq N} \Xi_{\flat}
			\right\|_{C^0(\Sigma_t^u)}.
	\end{align}
	We now set $N=11,$ 
	use the small data estimate \eqref{E:XIC0SMALLSOBOLEVNORM},
	and apply Gronwall's inequality to the quantity 
	$\left|
			\rgeo^{-1} \angLie_{\mathscr{Z}}^{\leq 11} \Xi_{\flat}
	\right|
	$
	in \eqref{E:POINTWISEGEOMETRICXIGRONWALLREADY}
	along the integral curves of $\Lunit.$
	We thus conclude that
	$\left\|
		\rgeo^{-1} \angLie_{\mathscr{Z}}^{\leq 11} \Xi_{\flat}
	\right\|_{C^0(\Sigma_t^u)}  \lesssim \varepsilon.$
	In particular, we have improved the bootstrap assumption \eqref{E:TEMPBOOTXI}
	and proved \eqref{E:C0BOUNDSGEOMETRICXI} with 
	$\Xi_{\flat}$ in place of $\Xi.$
	To conclude the desired bound \eqref{E:C0BOUNDSGEOMETRICXI}
	for $\Xi,$ we apply the Leibniz rule to the right-hand side of the identity
	$\angLie_{\mathscr{Z}}^N \Xi = \angLie_{\mathscr{Z}}^N(\ginversesphere \Xi_{\flat})$
	and use the estimate for 
	$\left\|
		\rgeo^{-1} \angLie_{\mathscr{Z}}^{\leq 11} \Xi_{\flat}
	\right\|_{C^0(\Sigma_t^u)}$
	together with
	\eqref{E:CRUDELOWERORDERC0BOUNDDERIVATIVESOFANGULARDEFORMATIONTENSORS}.	
	
	Inequality \eqref{E:C0BOUNDSXIRECTANGULARCOMPONENTS} then follows from
	the identity $\mathscr{Z}^N \Xi^i = \angLie_{\mathscr{Z}}^N (\Xi \cdot \angdiff x^i),$
	the Leibniz rule,
	and inequalities 
	\eqref{E:LOWERORDERPOINTWISEBOUNDPROJECTEDLIEDERIVATIVESANGDIFFCOORDINATEX}
	and \eqref{E:C0BOUNDSGEOMETRICXI}.
	
\end{proof}

\section{Estimates for the components of \texorpdfstring{$Z$}{the commutation vectorfields} relative to the geometric coordinates}
In this section, 
we prove a lemma that provides estimates for the components of the commutation vectofields 
$Z \in \mathscr{Z}$ relative to the geometric coordinates.
As a simple corollary, we deduce that if $f$ is a function, 
$0 \leq N \leq 12,$
and $\mathscr{Z}^{\leq N} f$ is a 
continuous function of the geometric coordinates $(t,u,\vartheta^1,\vartheta^2),$ 
then $f$ is a $C^N$ function of $(t,u,\vartheta^1,\vartheta^2).$

\begin{lemma}[\textbf{Estimates for the components of the commutation vectorfields relative to the geometric coordinates}]
	\label{L:GEOMETRICCOMPONENTS}
		Let $0 \leq N \leq 11$ be an integer
		and let $Z \in \mathscr{Z}.$
		Let $\lbrace (\mathbb{D}_i,\vartheta_i^1, \vartheta_i^2) \rbrace_{i=1,2}$
		be the atlas from Def.~\ref{D:ATLASONS2}
		and let $\mathbb{D}_i' \subset \mathbb{D}_i$
		be compact subsets such that $\mathbb{S}^2 = \mathbb{D}_1' \cup \mathbb{D}_2'.$
		Under the small-data and bootstrap assumptions 
		of Sects.~\ref{S:PSISOLVES}-\ref{S:C0BOUNDBOOTSTRAP},
		if $\varepsilon$ is sufficiently small,
		then for $\mathbb{D}' \in \lbrace \mathbb{D}_1', \mathbb{D}_2' \rbrace$
		and $A=1,2,$ the following pointwise estimates hold
		for $(t,u) \in [0,\Tboot) \times [0,U_0]:$
	\begin{subequations}
	\begin{align} 
		\left\|
			\Lunit \mathscr{Z}^{\leq N} Z \vartheta^A
		\right\|_{C^0(\lbrace t \rbrace \times [0,u] \times \mathbb{D}')}
		& \overset{\mathbb{D}'}{\lesssim} 
			\frac{\ln(\myexp + t)}{(1 + t)^2},
			\label{E:TIMEDERIVATIVEGEOMETRICCOMPONENTS}	\\
		\left\|
			\mathscr{Z}^{\leq N} Z \vartheta^A
		\right\|_{C^0(\lbrace t \rbrace \times [0,u] \times \mathbb{D}')}
		& \overset{\mathbb{D}'}{\lesssim} 
			1.
		\label{E:GEOMETRICCOMPONENTS}
	\end{align}
	\end{subequations}
	Furthermore, the following estimates hold:
	\begin{subequations}
	\begin{align} 
		\left\|
			\Lunit \mathscr{Z}^{\leq N} Z t
		\right\|_{C^0(\Sigma_t^u)}
		& \lesssim
			1,
			\label{E:TIMEDERIVATIVETIMECOMPONENTS}	\\
		\left\|
			\mathscr{Z}^{\leq N} Z t
		\right\|_{C^0(\Sigma_t^u)}
		& \lesssim
			1 + t,
		\label{E:TIMECOMPONENTS}
	\end{align}
	\end{subequations}
	\begin{subequations}
	\begin{align} 
		\left\|
			\Lunit \mathscr{Z}^{\leq N} Z u
		\right\|_{C^0(\Sigma_t^u)}
		& = 0,
			\label{E:TIMEDERIVATIVEEIKONALCOMPONENTS}	\\
		\left\|
			\mathscr{Z}^{\leq N} Z u
		\right\|_{C^0(\Sigma_t^u)}
		& \lesssim
			1.
		\label{E:EIKONALCOMPONENTS}
	\end{align}
	\end{subequations}
\end{lemma}	

\begin{proof}
	We first note that since $\Lunit \vartheta^A = 0,$
	it follows that 
	\begin{align}
		\Lunit \mathscr{Z}^N Z \vartheta^A
		& = 
			[\Lunit, \mathscr{Z}^N Z] \vartheta^A
			= 
			[\Lunit, \mathscr{Z}^{N+1}] \vartheta^A.
	\end{align}
	Hence, by making straightforward changes to the proof of \eqref{E:LZNCOMMUTATORACTINGONFUNCTIONSPOINTWISE},
	using the bootstrap assumptions 
	\eqref{E:PSIFUNDAMENTALC0BOUNDBOOTSTRAP}
	and \eqref{E:UPMUBOOT}-\eqref{E:CHIJUNKBOOT},
	and arguing inductively to bound the term analogous to the first term on the right-hand side of 	\eqref{E:LZNCOMMUTATORACTINGONFUNCTIONSPOINTWISE},
	we deduce the following estimate 
	on $[0,\Tboot) \times [0,U_0] \times \mathbb{D}':$
	\begin{align} \label{E:GEOMETRICCOMPONENTSGRONWALLREADY}
		\left|
			\Lunit \mathscr{Z}^{\leq N} Z \vartheta^A
		\right|
		& \lesssim
			\left|
			\Lunit \mathscr{Z}^{\leq N + 1} \vartheta^A
		\right|
		\lesssim 
			\frac{\ln(\myexp + t)}{(1 + t)^2}
			\left|
				\mathscr{Z}^{\leq N + 1} \vartheta^A
			\right|.
	\end{align}
	Applying Gronwall's inequality to the quantity $\left|\mathscr{Z}^{\leq N + 1} \vartheta^A \right|$ 
	in \eqref{E:GEOMETRICCOMPONENTSGRONWALLREADY}
	and using the small-data estimate
	$\left| \mathscr{Z}^{\leq N + 1} \vartheta^A \right|(0,u,\vartheta) \overset{\mathbb{D}'}{\lesssim} 1,$
	which is straightforward to verify using the estimates of Sect.~\ref{S:INITIALBEHAVIOROFQUANTITIES},
	we conclude the desired estimate \eqref{E:GEOMETRICCOMPONENTS}.
	The estimate \eqref{E:TIMEDERIVATIVEGEOMETRICCOMPONENTS} 
	then follows from \eqref{E:GEOMETRICCOMPONENTSGRONWALLREADY} and \eqref{E:GEOMETRICCOMPONENTS}.
	Next, since $Z t \in \lbrace 0, \rgeo \rbrace,$
	the estimates
	\eqref{E:TIMEDERIVATIVETIMECOMPONENTS}-\eqref{E:TIMECOMPONENTS}
	follow easily from Lemma~\ref{L:BASICESTIMATESFORRGEO}.
	Finally, since $Z u \in \lbrace 0, 1 \rbrace,$
	the desired estimates
	\eqref{E:TIMEDERIVATIVEEIKONALCOMPONENTS}-\eqref{E:EIKONALCOMPONENTS}
	follow easily.
\end{proof}

\begin{corollary}[\textbf{Geometric coordinate regularity follows from vectorfield differential operator regularity}]
	\label{C:VECTORFIEDLDIFFERENTIABILITYVSGEOCOORDINATEDIFFERENTIABILITY}
	Let $f$ be a function 
	and let $0 \leq N \leq 12$ be an integer.
	Assume that for $0 \leq M \leq N,$
	$\mathscr{Z}^M f$ is a continuous function
	of $(t,u,\vartheta).$ 
	Under the small-data and bootstrap assumptions 
	of Sects.~\ref{S:PSISOLVES}-\ref{S:C0BOUNDBOOTSTRAP},
	if $\varepsilon$ is sufficiently small,
	the $f$ is $N-$times continuously differentiable with respect to the geometric coordinates 
	$(t,u,\vartheta^1,\vartheta^2)$
	on $\mathcal{M}_{\Tboot,U_0}.$ 
\end{corollary}

\begin{proof}
	We first show that if 
	$Z \in \mathscr{Z},$ then the components
	of $Z$ relative to the geometric coordinates
	are $C^{11}$ functions of $(t,u,\vartheta).$
	To this end, we first note that the commutation set $\mathscr{Z}$ has span equal to the span of the geometric coordinate vectorfields
	$\lbrace \frac{\partial}{\partial t}, \frac{\partial}{\partial u}, X_1, X_2 \rbrace.$
	From Lemma~\ref{L:GEOMETRICCOMPONENTS} with $N=0,$
	it follows that each coordinate vectorfield can (locally) be written as a linear combination
	of the vectorfields in $\mathscr{Z}$ with coefficients that are $C^0$ functions of 
	$(t,u,\vartheta^1,\vartheta^2).$
	Hence, we can use Lemma~\ref{L:GEOMETRICCOMPONENTS} with $N=1$ to conclude that
	these coefficients are in fact $C^1$ functions of $(t,u,\vartheta^1,\vartheta^2)$
	and that the same regularity statement holds for the components 
	$Zt, Zu, Z \vartheta^1, Z \vartheta^2$
	of the vectorfields $Z \in \mathscr{Z}$ 
	relative to the geometric coordinates.
	Continuing by induction with the help of Lemma~\ref{L:GEOMETRICCOMPONENTS}, 
	we conclude that the components are $C^{11}$ as desired.
	Cor.~\ref{C:VECTORFIEDLDIFFERENTIABILITYVSGEOCOORDINATEDIFFERENTIABILITY}
	now follows from expressing the geometric coordinate derivatives of $f$
	in terms of the aforementioned linear combinations of 
	vectorfields $Z \in \mathscr{Z}$ and noting that 
	in the resulting expressions,
	no more than $11$ derivatives
	fall on the components of the $Z.$ 
\end{proof}

\section{Estimates for the rectangular spatial derivatives of \texorpdfstring{$u$}{the eikonal function}}

\begin{lemma}[\textbf{Estimates for the rectangular spatial derivatives of $u$}]
	\label{L:EIKONALFUNCTIONRECTANGULARDERIVATIVES}
	Under the small-data and bootstrap assumptions 
	of Sects.~\ref{S:PSISOLVES}-\ref{S:C0BOUNDBOOTSTRAP},
	if $\varepsilon$ is sufficiently small,
	then for $i=1,2,3,$ 
	the following estimates hold for the rectangular spatial
	derivatives of the eikonal function
	for $(t,u) \in [0,\Tboot) \times [0,U_0]:$
	\begin{align} \label{E:EIKONALFUNCTIONRECTANGULARDERIVATIVES}
		\left\|
			\mathscr{Z}^{\leq 12}
			\left(
				\upmu \partial_i u
				+ \frac{x^i}{r}
			\right)
		\right\|_{C^0(\Sigma_t^u)}
		& \leq C \varepsilon \frac{\ln(\myexp + t)}{1 + t}.
	\end{align}
\end{lemma}

\begin{proof}
	Inequality \eqref{E:EIKONALFUNCTIONRECTANGULARDERIVATIVES}
	follows from the identity \eqref{E:PARTIALIUINTERMSOFRADUNIT},
	the estimate \eqref{E:ALTERNATEVERSIONRADUNITLOWEREDBACKGROUNDSUBTRACTEDC0},
	and Cor.~\ref{C:SQRTEPSILONREPLCEDWITHCEPSILON}.
\end{proof}


\chapter{Sharp Estimates for \texorpdfstring{$\upmu$}{the Inverse Foliation Density}}
\label{C:SHARPESTIMATESFORUPMU}
\thispagestyle{fancy}
In Chapter~\ref{C:SHARPESTIMATESFORUPMU}, we derive a variety of sharp pointwise estimates for the inverse foliation
density $\upmu$ and some of its derivatives. We then derive corresponding estimates
for time integrals of some related quantities. 
The time integral estimates play a critical role in our
proof of the main Gronwall-type lemma 
(see Lemma~\ref{L:FUNDAMENTALGRONWALL}),
which we use to deduce a priori $L^2$ estimates for $\Psi.$
Our time integral estimates are sensitive and affect 
the degree of degeneracy-in-$\upmu^{-1}$ of our $L^2$
estimates, which in turn affects the number of derivatives we
need to close the proof of the sharp classical lifespan theorem.
The most difficult estimates in this section
are based on the fact
that $\Lunit(\rgeo \Lunit \upmu)$ is integrable 
in time along the integral curves of $\Lunit.$
This fact leads to the approximate monotonicity of $\upmu$ along the integral curves,
where along a fixed integral curve, the sign of the monotonicity up to some ``late'' time $t$ 
is determined by the sign of $\Lunit \upmu$ at time $t;$
see Sect.~\ref{SS:INTROLUPMUOVERUPMU} for an overview.

\begin{remark}[\textbf{The role of a posteriori estimates}]
	The most difficult estimates involving $\upmu$ 
	in this section are not a priori nature, but rather a posteriori. 
	Hence, many of our estimates are derived on time intervals of
	the form $0 \leq s \leq t,$ and the arguments are
	based, on a case-by-case basis, on the possible (but not known in advance) behavior 
	of various quantities at time $t.$
\end{remark}

\section{Basic ingredients in the analysis}
We begin by recalling that $\upmu$ verifies the transport equation (see \eqref{E:UPMUFIRSTTRANSPORT})
\begin{align}
	\Lunit \upmu & = \upomega.
\end{align}
In this section, we often view $\upmu = \upmu(s,u,\vartheta)$
and similarly for other quantities.
In doing so, we are referring to the moving time 
coordinate as ``$s$'' in order to distinguish it from a
``fixed later time'' $t$ verifying $t \geq s.$ 
We recall that relative to these coordinates, we
have $\Lunit = \frac{\partial}{\partial s} := \frac{\partial}{\partial s}|_{u,\vartheta}.$
In summary, we have
\begin{align}
	\frac{\partial}{\partial s} \upmu(s,u,\vartheta) & = \upomega(s,u,\vartheta).
\end{align}

If a function $f$ depends on the geometric coordinates $(s,u,\vartheta)$ as well as the 
late time $t \geq s,$ then we indicate this by writing $f(s,u,\vartheta;t).$

We now define some auxiliary quantities that we use to analyze $\upmu.$

\begin{definition}[\textbf{Auxiliary quantities used to analyze $\upmu$}]
	\label{D:AUXQUANTITIES}
	We define the following quantities, 
	where we assume that $0 \leq s \leq t$
	for those quantities that depend on both $s$ and $t:$
	\begin{subequations}
	\begin{align}
		\Omega(s,u,\vartheta)
		& := \rgeo(s,u) \upomega(s,u,\vartheta)
			= \rgeo(s,u) \Lunit \upmu(s,u,\vartheta),
		\label{E:BIGOMEGADEF} \\
	m(s,u,\vartheta;t) 
	& := \frac{\Omega(t,u,\vartheta) - \Omega(s,u,\vartheta)}{\rgeo(s,u)},
		\label{E:LITTLEMDEF} \\
	M(s,u,\vartheta;t) 
	& := \int_{s'=s}^{s'=t} m(s',u,\vartheta;t) \, ds',
		\label{E:BIGMDEF} \\
	\mathring{\upmu}(u,\vartheta)
	& := \upmu(s=0,u,\vartheta),
		\\
	\widetilde{\Omega}(s,u,\vartheta;t)
	& := \frac{\Omega(s,u,\vartheta)}{\mathring{\upmu}(u,\vartheta) - M(0,u,\vartheta;t)},
		\label{E:WIDETILDEBIGOMEGADEF} \\
	\widetilde{m}(s,u,\vartheta;t)
	& := \frac{m(s,u,\vartheta;t)}{\mathring{\upmu}(u,\vartheta) - M(0,u,\vartheta;t)},
		\label{E:WIDETILDELITTLEMDEF} \\
	\widetilde{M}(s,u,\vartheta;t)
	& := \frac{M(s,u,\vartheta;t)}{\mathring{\upmu}(u,\vartheta) - M(0,u,\vartheta;t)},
		\label{E:WIDETILDEBIGMDEF} \\
	\upmu_{(Approx)}(s,u,\vartheta;t)
		& := 1
			+  \widetilde{\Omega}(t,u,\vartheta;t) \ln \left(\frac{\rgeo(s,u)}{\rgeo(0,u)} \right)
			+ \widetilde{M}(s,u,\vartheta;t).
		\label{E:MUAPPROXDEF}
	\end{align}
	\end{subequations}
\end{definition}

\begin{remark}[\textbf{The role of} $\upmu_{(Approx)}(s,u,\vartheta;t)$]
	The main point of Def.~\ref{D:AUXQUANTITIES} is that
	$\upmu_{(Approx)}(s,u,\vartheta;t)$ is a good approximation to 
	$\upmu(s,u,\vartheta)$ that is monotonic 
	along each integral curve of
	$\Lunit$ for $0 \leq s \leq t$
	up to the error term $\widetilde{M}(s,u,\vartheta;t).$
	Furthermore, the sign of the approximate monotonicity is determined by 
	$\widetilde{\Omega}(t,u,\vartheta;t).$
\end{remark}

The identities stated in the next lemma follow easily from Def.~\ref{D:AUXQUANTITIES};
we omit the simple proofs.

\begin{lemma}[\textbf{Some identities verified by the auxiliary quantities}]
The following identities hold,
where we assume that $0 \leq s \leq t$
for those quantities that depend on both $s$ and $t:$
\begin{subequations}
\begin{align}
	\frac{\partial}{\partial s} \upmu(s,u,\vartheta)
	& = \frac{\Omega(s,u,\vartheta)}{\rgeo(s,u)}
		= \frac{\Omega(t,u,\vartheta)}{\rgeo(s,u)}
			- m(s,u,\vartheta;t),
		\label{E:LMURELATION} \\
	\frac{\partial}{\partial s} M(s,u,\vartheta;t)
	& = - m(s,u,\vartheta;t),
		\\	
	\frac{\partial}{\partial s} \upmu_{(Approx)}(s,u,\vartheta;t)
	& = \frac{\widetilde{\Omega}(t,u,\vartheta;t)}{\rgeo(s,u)}
		- \widetilde{m}(s,u,\vartheta;t),
		\label{E:LMUAPPROXRELATION}
\end{align}

\begin{align} \label{E:EQUALRATIOS}
	\frac{\frac{\partial}{\partial s} \upmu_{(Approx)}(s,u,\vartheta;t)}{\upmu_{(Approx)}(s,u,\vartheta;t)}
	& = \frac{\frac{\partial}{\partial s} \upmu(s,u,\vartheta)}{\upmu(s,u,\vartheta)},
\end{align}

\begin{align} \label{E:LMUAPPROXINTERMSOFMUAPPROX}
	\frac{\partial}{\partial s} \upmu_{(Approx)}(s,u,\vartheta;t)
	& = \frac{\upmu_{(Approx)}(s,u,\vartheta;t) - 1}{\rgeo(s,u) \left\lbrace 1 + \ln \left(\frac{\rgeo(s,u)}{\rgeo(0,u)} \right) \right\rbrace}
		+ \frac{\widetilde{\Omega}(t,u,\vartheta;t) - \widetilde{M}(s,u,\vartheta;t)}	
					{\rgeo(s,u) \left\lbrace 1 + \ln \left(\frac{\rgeo(s,u)}{\rgeo(0,u)} \right) \right\rbrace}
		- \widetilde{m}(s,u,\vartheta;t),
\end{align}

\begin{align}
	\frac{\partial}{\partial s} \upmu(s,u,\vartheta)
	& = \frac{\Omega(t,u,\vartheta)}{\rgeo(s,u)}
		- m(s,u,\vartheta;t),
\end{align}

\begin{align}
	\frac{\partial}{\partial s} \upmu(s,u,\vartheta)
	& = \left\lbrace
				\mathring{\upmu}(u,\vartheta) - M(0,u,\vartheta;t)
			\right\rbrace 
			\left\lbrace
				\frac{\widetilde{\Omega}(t,u,\vartheta;t)}{\rgeo(s,u)}
				- \widetilde{m}(s,u,\vartheta;t)
			\right\rbrace,
\end{align}

\begin{align} \label{E:MUSPLIT}	
	\upmu(s,u,\vartheta) 
	& = \left\lbrace
				\mathring{\upmu}(u,\vartheta) - M(0,u,\vartheta;t)
			\right\rbrace
			\upmu_{(Approx)}(s,u,\vartheta;t) .
\end{align}
\end{subequations}
\end{lemma}

$\hfill \qed$

We now define some quantities that play a role in our analysis of $\upmu.$
\begin{definition}[$\upmu_{(Min)},$ $\upmu_{\star},$ \textbf{and} $\widetilde{\Omega}_{(Min)}$]
	We define the following quantities,
	where we assume that $0 \leq s \leq t$
	in \eqref{E:TILDEOMEGAMIN}:	
	\begin{align}
		\upmu_{(Min)}(s,u)
		& := \min_{\Sigma_s^u} \upmu = \min_{u' \in [0,u], \vartheta \in \mathbb{S}^2} \upmu(s,u',\vartheta),
			\label{E:MUMINDEF} \\
		\upmu_{\star}(s,u)
		& := \min \left\lbrace 1, \min_{\Sigma_s^u} \upmu \right\rbrace
			= \min \left\lbrace 1, \upmu_{(Min)}(s,u) \right\rbrace,
			\label{E:MUSTARDEF} \\
		\widetilde{\Omega}_{(Min)}(s,u;t)
		& := \min_{u' \in [0,u], \vartheta \in \mathbb{S}^2} \widetilde{\Omega}(s,u',\vartheta;t).
			\label{E:TILDEOMEGAMIN}
\end{align}
	
\end{definition}

In the next lemma, 
we provide some basic estimates for the auxiliary quantities.
Below, we use these basic estimates to prove more intricate estimates
(see Prop.~\ref{P:SHARPMU}).

\begin{lemma}[\textbf{First estimates for the auxiliary quantities used to analyze $\upmu$}]
\label{E:FIRSTESTIMATESFORAUXILIARYUPMUQUANTITIES}
Under the 
small data and
bootstrap assumptions 
of Sects.~\ref{S:PSISOLVES}-\ref{S:C0BOUNDBOOTSTRAP},
if $\varepsilon$ is sufficiently small, 
then the following pointwise estimates hold 
for $(t,u,\vartheta) \in [0,\Tboot) \times [0,U_0] \times \mathbb{S}^2$
and $0 \leq s \leq t:$ 
\begin{align}
	|\mathring{\upmu}(u,\vartheta) - 1|
	& \lesssim \varepsilon,
		\label{E:MUINITIALDATAESTIMATE}
		\\
	|\mathring{\upmu}(u,\vartheta) - M(0,u,\vartheta;t) - 1|
	& \lesssim \varepsilon,
		\label{E:MUAMPLITUDENEARONE}
\end{align}
\begin{align}
	|\Omega(s,u,\vartheta)|, |\widetilde{\Omega}(s,u,\vartheta;t)| 
	& \lesssim \varepsilon,
	\label{E:BIGOMEGATRIVIALBOUND}
	\\
	 \left|
			\Omega(t,u,\vartheta)
			- \frac{1}{2} [\rgeo G_{\Lunit \Lunit} \Rad \Psi](t,u,\vartheta)
		\right|
		& \lesssim 
			\varepsilon \frac{1}{1+t}.
	\label{E:BIGOMEGAISWELLAPPROXIMATEDBYNULLCONDITIONFAILUREFACTORTERM}
\end{align}

 In addition, the following pointwise estimates hold:
\begin{align}
	|\Omega(s,u,\vartheta) - \Omega(t,u,\vartheta)|,
	|\widetilde{\Omega}(s,u,\vartheta;t) - \widetilde{\Omega}(t,u,\vartheta;t)|
	& \lesssim \varepsilon
		\frac{\ln(\myexp + s)}{1 + s}
		\left(
			\frac{t - s}{1 + t}
		\right),
	\label{E:BIGOMEAGDIFFERENCEEST}
\end{align}
\begin{align}
	|m(s,u,\vartheta;t)|, |\widetilde{m}(s,u,\vartheta;t)|
	& \lesssim \varepsilon 
		\frac{\ln(\myexp + s)}{(1 + s)^2}
		\left(
			\frac{t - s}{1 + t}
		\right),
		\label{E:LITTLEMEST} \\
	|M(s,u,\vartheta;t)|, |\widetilde{M}(s,u,\vartheta;t)|
	& \lesssim \varepsilon 
		\frac{\ln(\myexp + s)}{1 + s}
		\left(
			\frac{t - s}{1 + t}
		\right),
		\label{E:BIGMEST} 
\end{align}
\begin{align} \label{E:MUAPPROXMINUSONEEASYBOUND}
	|\upmu_{(Approx)}(s,u,\vartheta;t) - 1| & \leq C \varepsilon \ln(\myexp + s),
\end{align}
\begin{align} \label{E:MUAPPROXMISLIKEMU}
	(1 - C \varepsilon) \upmu_{(Approx)}(s,u,\vartheta;t)
	& \leq \upmu(s,u,\vartheta) 
	\leq (1 + C \varepsilon)\upmu_{(Approx)}(s,u,\vartheta;t),
		\\
	(1 - C \varepsilon) \frac{\partial}{\partial s} \upmu_{(Approx)}(s,u,\vartheta;t)
	& \leq \Lunit \upmu(s,u,\vartheta) 
	\leq (1 + C \varepsilon) \frac{\partial}{\partial s} \upmu_{(Approx)}(s,u,\vartheta;t),
	\label{E:LMUAPPROXMISLIKELMU}
\end{align}
\begin{align}
	\left|
		\upmu(s,u,\vartheta)
		-
		\left\lbrace	 
			\mathring{\upmu}(u,\vartheta)
			+ \frac{1}{2}
				\ln \left(\frac{\rgeo(s,u)}{\rgeo(0,u)} \right)
				[\rgeo G_{\Lunit \Lunit} \Rad \Psi](t,u,\vartheta)
		\right\rbrace
		\right|	
		& \lesssim \varepsilon.
		\label{E:UPMUALMOSTSHARPBOUND}
\end{align}

\end{lemma}

\begin{proof}
The estimate \eqref{E:MUINITIALDATAESTIMATE} was already proved in \eqref{E:SMALLINITIALC0BOUNDNORMS}.

The estimate \eqref{E:BIGOMEGATRIVIALBOUND} for $\Omega$ follows from definition \eqref{E:BIGOMEGADEF} and 
the estimate \eqref{E:C0BOUNDLDERIVATIVECRUCICALEIKONALFUNCTIONQUANTITIES}.

To prove \eqref{E:BIGOMEGAISWELLAPPROXIMATEDBYNULLCONDITIONFAILUREFACTORTERM},
we first use 
\eqref{E:UPMUFIRSTTRANSPORT}
and \eqref{E:BIGOMEGADEF}
to deduce that
\begin{align} \label{E:BIGOMEGADECOMPOSED}
	\Omega
	-  \frac{1}{2} \rgeo G_{\Lunit \Lunit} \Rad \Psi
	& = \rgeo G_{(Frame)} \Lunit \Psi.
\end{align}
The desired estimate 
\eqref{E:BIGOMEGAISWELLAPPROXIMATEDBYNULLCONDITIONFAILUREFACTORTERM}
now follows from \eqref{E:BIGOMEGADECOMPOSED},
the estimate
\eqref{E:LOWERORDERC0BOUNDLIEDERIVATIVESOFGRAME},
and the bootstrap assumptions
\eqref{E:PSIFUNDAMENTALC0BOUNDBOOTSTRAP}. 

We now prove the estimate \eqref{E:BIGOMEAGDIFFERENCEEST} for $\Omega.$
We first use 
\eqref{E:UPMUFIRSTTRANSPORT}
and \eqref{E:BIGOMEGADEF}
to deduce that
\begin{align} \label{E:LOFBIGOMEGA}
	\Lunit \Omega 
	& = 
		\Lunit
		\left\lbrace
			\frac{1}{2} G_{\Lunit \Lunit} \rgeo \Rad \Psi
			- \frac{1}{2} \upmu G_{\Lunit \Lunit} \rgeo \Lunit \Psi
			- \upmu G_{\Lunit \Radunit} \rgeo \Lunit \Psi
		\right\rbrace,
\end{align}
where we view the left-hand and right-hand sides as functions of 
$(s',u,\vartheta)$ and $\Lunit = \frac{\partial}{\partial s'}.$
From the estimates
\eqref{E:LOWERORDERC0BOUNDLIEDERIVATIVESOFGRAME},
\eqref{E:C0BOUNDCRUCIALEIKONALFUNCTIONQUANTITIES},
\eqref{E:C0BOUNDLDERIVATIVECRUCICALEIKONALFUNCTIONQUANTITIES},
\eqref{E:IMPORTANTLDERIVATIVEOFGFRAMEC0BOUNDESTIMATES},
and
\eqref{E:LOFRGEOWEIGHTEDZPSISTRONGC0BOUND}
and the bootstrap assumptions
\eqref{E:PSIFUNDAMENTALC0BOUNDBOOTSTRAP},
it follows that the magnitude of the right-hand side of \eqref{E:LOFBIGOMEGA} is
$\lesssim \varepsilon \ln(\myexp + s')(\myexp + s')^{-2}.$ Hence, we have
\begin{align}
	\left| 
		\Omega(s,u,\vartheta) 
		- \Omega(t,u,\vartheta)
	\right|
	& \lesssim 
	\varepsilon
	\int_{s'=s}^t
		\frac{\ln(\myexp + s')}{(\myexp + s')^2}
	\, ds'
		\\
	& =
		\varepsilon
		\left\lbrace
			\frac{\ln(\myexp + s) + 1}{\myexp + s}
		\right\rbrace
		-
		\varepsilon
		\left\lbrace
			\frac{\ln(\myexp + t) + 1}{\myexp + t}
		\right\rbrace
		\notag \\
	& \leq 
		\varepsilon
		\left\lbrace
			\frac{\ln(\myexp + s) + 1}{\myexp + s}
		\right\rbrace
		-
		\varepsilon
		\left\lbrace
			\frac{\ln(\myexp+s) + 1}{\myexp + t}
		\right\rbrace
		\notag \\
	& \lesssim 
		\varepsilon	
		\frac{\ln(\myexp+s)}{\myexp + s} 
		\left(
			\frac{t - s}{\myexp + t}
		\right)
	\leq
		\varepsilon	
		\frac{\ln(\myexp+s)}{1 + s} 
		\left(
			\frac{t - s}{1 + t}
		\right).
		\notag
\end{align}
We have thus proved the desired inequality \eqref{E:BIGOMEAGDIFFERENCEEST} that involves $\Omega.$ 

The estimate \eqref{E:LITTLEMEST} involving $m$ then follows easily from
definition \eqref{E:LITTLEMDEF} and the inequality \eqref{E:BIGOMEAGDIFFERENCEEST} involving $\Omega.$ 

The estimate \eqref{E:BIGMEST} involving $M$ then follows easily from 
definition \eqref{E:BIGMDEF} and the inequality \eqref{E:LITTLEMEST} involving $m.$ 

The estimate \eqref{E:MUAMPLITUDENEARONE} then follows from \eqref{E:MUINITIALDATAESTIMATE}
and the estimate \eqref{E:BIGMEST} for $|M(0,u,\vartheta;t)|.$

The estimates 
\eqref{E:BIGOMEGATRIVIALBOUND},
\eqref{E:BIGOMEAGDIFFERENCEEST}, 
\eqref{E:LITTLEMEST}, 
and \eqref{E:BIGMEST}
involving $\widetilde{\Omega},$ $\widetilde{m},$ and $\widetilde{M}$ then follow from
definitions \eqref{E:WIDETILDEBIGOMEGADEF}, \eqref{E:WIDETILDELITTLEMDEF}, \eqref{E:WIDETILDEBIGMDEF}, 
the corresponding estimates for the non-tilded quantities, and inequality \eqref{E:MUAMPLITUDENEARONE}.

The estimate \eqref{E:MUAPPROXMINUSONEEASYBOUND} then follows from
definition \eqref{E:MUAPPROXDEF} and the estimates \eqref{E:BIGOMEGATRIVIALBOUND} and \eqref{E:BIGMEST}
involving $\widetilde{\Omega}$ and $\widetilde{M}.$ 
	
The estimate \eqref{E:MUAPPROXMISLIKEMU} follows from the identity \eqref{E:MUSPLIT}
and the estimate \eqref{E:MUAMPLITUDENEARONE}.

The estimate \eqref{E:LMUAPPROXMISLIKELMU} follows from the identity \eqref{E:EQUALRATIOS}
and the estimate \eqref{E:MUAPPROXMISLIKEMU}.

		To prove \eqref{E:UPMUALMOSTSHARPBOUND},
		we first use definition \eqref{E:BIGOMEGADEF} and the estimates 
		\eqref{E:BIGOMEGAISWELLAPPROXIMATEDBYNULLCONDITIONFAILUREFACTORTERM}
		and
		\eqref{E:BIGOMEAGDIFFERENCEEST}
		to deduce that for $0 \leq s' \leq t,$ we have
		\begin{align}  \label{E:MUTRANSPORTINEQUALITYWITHNULLCONDITIONFAILUREFACTOR}
		\left|
			\Lunit \upmu(s',u,\vartheta)
			- \frac{1}{2} \frac{1}{\rgeo(s',u)} [\rgeo G_{\Lunit \Lunit} \Rad \Psi](t,u,\vartheta) 
		\right|
		& \lesssim 
			\varepsilon \frac{\ln(\myexp + s')}{(1 + s')^2}.
	\end{align}
	Integrating \eqref{E:MUTRANSPORTINEQUALITYWITHNULLCONDITIONFAILUREFACTOR} 
	along the integral curves of $\Lunit = \frac{\partial}{\partial s'}$ 
	from $s'=0$ to $s'=s,$ we conclude \eqref{E:UPMUALMOSTSHARPBOUND}.
 
\end{proof}

\section{Sharp pointwise estimates for \texorpdfstring{$\upmu,$ $\Lunit \upmu,$ and $\Rad \upmu$}{the inverse foliation density}}
When we derive our a priori $L^2$ estimates, it will be essential that we
treat the region where $\upmu$ is not appreciably shrinking in a different fashion than
we treat the region where $\upmu$ is appreciably shrinking. 
The estimates of Lemma~\ref{E:FIRSTESTIMATESFORAUXILIARYUPMUQUANTITIES}
imply that along the portion of the integral curve of $\Lunit$
corresponding to $0 \leq s \leq t$ and fixed $(u,\vartheta),$
the behavior of $\upmu$ is essentially determined by 
$\widetilde{\Omega}(t,u,\vartheta;t).$
This discussion motivates the following definitions.

\begin{definition}[\textbf{Regions of distinct $\upmu$ behavior}]
\label{D:REGIONSOFDISTINCTUPMUBEHAVIOR}
Let $\widetilde{\Omega}$ be the function defined in \eqref{E:WIDETILDEBIGOMEGADEF}.
For each $s \in [0,t]$ and $u \in [0,U_0],$ we partition 
\begin{subequations}
\begin{align}
	[0,u] \times \mathbb{S}^2 
	& = \Vplus{t}{u} \cup \Vminus{t}{u},
		\label{E:OUINTERVALCROSSS2SPLIT} \\
	\Sigma_s^u
	& = \Sigmaplus{s}{t}{u} \cup \Sigmaminus{s}{t}{u},
	\label{E:SIGMASSPLIT}
\end{align}
\end{subequations}
where
\begin{subequations}
\begin{align}
	\Vplus{t}{u}
	& := 
	\left\lbrace
		(u',\vartheta) \in [0,u] \times \mathbb{S}^2 \ | \ \widetilde{\Omega}(t,u',\vartheta;t) \geq 0
	\right\rbrace,
		\label{E:ANGLESANDUWITHNONDECAYINUPMUGBEHAVIOR} \\
	\Vminus{t}{u}
	& := 
	\left\lbrace
		(u',\vartheta) \in [0,u] \times \mathbb{S}^2 \ | \ \widetilde{\Omega}(t,u',\vartheta;t) < 0
	\right\rbrace,
		\label{E:ANGLESANDUWITHDECAYINUPMUGBEHAVIOR} \\
	\Sigmaplus{s}{t}{u}
	& := 
	\left\lbrace
		(s,u',\vartheta) \in \Sigma_s^u \ | \ \widetilde{\Omega}(t,u',\vartheta;t) \geq 0
	\right\rbrace,
		\label{E:SIGMAPLUS} \\
	\Sigmaminus{s}{t}{u}
	& := 
	\left\lbrace
		(s,u',\vartheta) \in \Sigma_s^u \ | \ \widetilde{\Omega}(t,u',\vartheta;t) < 0
	\right\rbrace.
	\label{E:SIGMAMINUS}
\end{align}
\end{subequations}
	
\end{definition}

In the next proposition, we provide the
sharp pointwise estimates for the inverse foliation
density $\upmu$ and some of its derivatives.

\begin{proposition}[\textbf{Sharp pointwise estimates for $\upmu,$ $\Lunit \upmu,$ and $\Rad \upmu$}]
\label{P:SHARPMU} 
Assume that the 
small data and
bootstrap assumptions 
of Sects.~\ref{S:PSISOLVES}-\ref{S:C0BOUNDBOOTSTRAP}
hold for $(t,u,\vartheta) \in \mathcal{M}_{\Tboot,U_0}.$
If $\varepsilon$ is sufficiently small,
$0 \leq s \leq t < \Tboot,$ 
and $0 \leq u \leq U_0,$
then the following estimates hold.

\medskip

\noindent \underline{\textbf{Weak upper bound for $\Lunit \upmu$ under all conditions.}}
\begin{align} \label{E:MUWEAKUPPERBOUND}
	\| \Lunit \upmu \|_{C^0(\Sigma_s^u)}
	& \leq C \varepsilon \frac{1}{1 + s}.
\end{align}

\medskip

\noindent \underline{\textbf{Upper bounds for $\frac{[\Lunit \upmu]_+}{\upmu}$ under all conditions.}}
\begin{align} \label{E:POSITIVEPARTOFLMUOVERMUISSMALL}
	\left\|
		\frac{[\Lunit \upmu]_+}{\upmu}
	\right\|_{C^0(\Sigma_s^u)}
	& \leq (1 + C \varepsilon) \frac{1}{\rgeo(s,u) \left\lbrace 1 + \ln \left(\frac{\rgeo(s,u)}{\rgeo(0,u)} \right) \right\rbrace}.
\end{align}

Furthermore,
\begin{align} \label{E:ALTERNATEPOSITIVEPARTOFLMUOVERMUISSMALL}
	\left\|
		\frac{[\Lunit \upmu]_+}{\upmu}
	\right\|_{C^0(\Sigma_s^u)}
	& \leq C \varepsilon \frac{1}{1 + s}.
\end{align}

\medskip

\noindent \underline{\textbf{Small $\upmu$ implies $\Lunit \upmu$ is negative.}}
\begin{align} \label{E:SMALLMUIMPLIESLMUISNEGATIVE}
	\upmu(s,u,\vartheta) \leq \frac{1}{4}
	\implies
	\Lunit \upmu(s,u,\vartheta) \leq - \frac{2}{\rgeo(s,u) \left\lbrace 1 + \ln \left(\frac{\rgeo(s,u)}{\rgeo(0,u)} \right) \right\rbrace}.
\end{align}

\medskip

\noindent \underline{\textbf{Upper bound for $\frac{[\Rad \upmu]_+}{\upmu}$ under all conditions.}}
\begin{align} \label{E:UNIFORMBOUNDFORMRADMUOVERMU}
	\left\|
		\frac{[\Rad \upmu]_+}{\upmu}
	\right\|_{C^0(\Sigma_s^u)}
	& \leq C \varepsilon^{1/2} 	
					\frac{\ln(\myexp + s)}{\sqrt{\ln(\myexp + t) - \ln(\myexp + s)}}
			+ C \varepsilon \ln(\myexp + s).
\end{align}

\medskip

\noindent \underline{\textbf{Upper bound for $\frac{[\Lunit \upmu + \uLgood \upmu]_+}{\upmu}$ under all conditions.}}
\begin{align} \label{E:POSITIVEPARTOFLMUPLUSRADMUOVERMUISSMALL}
	\left\|
		\frac{[\Lunit \upmu + \uLgood \upmu]_+}{\upmu}
	\right\|_{C^0(\Sigma_s^u)}
	& \leq C \varepsilon^{1/2} 	
					\frac{\ln(\myexp + s)}{\sqrt{\ln(\myexp + t) - \ln(\myexp + s)}}
				+ C \varepsilon \ln(\myexp + s).
\end{align}

\medskip

\noindent \underline{\textbf{Sharp spatially uniform estimates when $\upmu$ is not decaying.}}
Consider a time interval $s \in [0,t],$ and assume the right endpoint inequality
\begin{align}
	\widetilde{\Omega}_{(Min)}(t,u;t) \geq 0,
\end{align}
where $\widetilde{\Omega}_{(Min)}$ is defined in \eqref{E:TILDEOMEGAMIN}. Let 
\begin{align*}
	\Contwo > 0
\end{align*}
be a real number. There exists a constant $C > 0$ 
such that if $\Contwo \sqrt{\varepsilon} < 1,$ 
then the following estimates hold for $s \in [0,t]:$
\begin{subequations}
\begin{align}
		1 - C \varepsilon
	& \leq \upmu_{\star}(s,u)
		\leq 1,
		\label{E:EASYCASEMUSTARTABOUNDS}
		\\
	1 - C \sqrt{\varepsilon}
	& \leq \upmu_{\star}^{\Contwo}(s,u)
		\leq 1,
		\label{E:APOWEREASYCASEMUSTARTABOUNDS} 
		\\
		\| \rgeo [\Lunit \upmu]_- \|_{C^0(\Sigma_s^u)}
	& \leq C \varepsilon \frac{\ln(\myexp + s)}{1 + s}.
		\label{E:EASYCASEOMEGAMINUSBOUNDS} 
\end{align}
\end{subequations}

\medskip

\noindent \underline{\textbf{Sharp spatially uniform estimates when $\upmu$ is decaying.}}
Consider a time interval $s \in [0,t],$ and assume the right endpoint inequality
\begin{align} \label{E:CRUCIALDELTADEF}
	\updelta & := - \widetilde{\Omega}_{(Min)}(t,u;t) > 0,
\end{align}
where $\widetilde{\Omega}_{(Min)}$ is defined in \eqref{E:TILDEOMEGAMIN}.
Then the following short-time and large-time estimates hold.
\ \\

\noindent \underline{\textbf{Short-time estimates.}}
Let 
\begin{align*}
	\Contwo > 0
\end{align*}
be a real number. There exists
a constant $C > 0$
such that if $\Contwo \sqrt{\varepsilon} |\ln \varepsilon| < 1$
and ${s \in [0, (1 - u) \frac{\varepsilon}{\updelta^2} - (1 - u))},$ 
then
\begin{subequations}
\begin{align} 
	1 - C \varepsilon |\ln \varepsilon|
	& \leq \upmu_{\star}(s,u) 
		\leq 1,
		\label{E:SHORTTIMEHARDCASEMUSTARBOUNDS} 
		\\
	1 - C \sqrt{\varepsilon}
	& \leq \upmu_{\star}^{\Contwo}(s,u) 
		\leq 1,
		\label{E:POWERASHORTTIMEHARDCASEMUSTARBOUNDS} 
		\\
		\| \rgeo [\Lunit \upmu]_-\|_{C^0(\Sigma_s^u)}
		= \| \Omega_-\|_{C^0(\Sigma_s^u)}
	& \leq 
		\left(
			1 + C \varepsilon
		\right)
		\left\lbrace
			\updelta 
			+ C \varepsilon \frac{\ln(\myexp + s)}{1 + s}
		\right\rbrace.
		\label{E:SHORTTIMEHARDCASEOMEGAMINUSBOUND}
\end{align}
\end{subequations}

\medskip

\noindent \underline{\textbf{Large-time estimates.}}
Under the short-time estimate hypotheses,
there exists a constant $C > 0$ 
such that if $\Contwo \sqrt{\varepsilon} |\ln \varepsilon| < 1$ 
and $s \in [(1 - u) \frac{\varepsilon}{\updelta^2} - (1 - u),t],$ 
then
\begin{subequations}
\begin{align}
	\left(
				1 
				- C \varepsilon |\ln \varepsilon|
			\right)
			\left\lbrace 1 - \updelta \ln \left( \frac{\rgeo(s,u)}{\rgeo(0,u)} \right) \right\rbrace
	& \leq
	\upmu_{\star}(s,u)
	\leq 
			\left(
				1 
				+ C \varepsilon |\ln \varepsilon|
			\right)
			\left\lbrace 1 - \updelta \ln \left( \frac{\rgeo(s,u)}{\rgeo(0,u)} \right) \right\rbrace,
	\label{E:LARGETIMEHARDCASEMUSTARBOUNDS}  
		\\
	\left(
		1 - C \sqrt{\varepsilon}
	\right)
	\left\lbrace 1 - \updelta \ln \left( \frac{\rgeo(s,u)}{\rgeo(0,u)} \right) \right\rbrace^{\Contwo}
	& \leq
	\upmu_{\star}^{\Contwo}(s,u)
	\leq 
	\left(
		1 + C \sqrt{\varepsilon}
	\right)
	\left\lbrace 1 - \updelta \ln \left( \frac{\rgeo(s,u)}{\rgeo(0,u)} \right) \right\rbrace^{\Contwo},
	\label{E:POWERALARGETIMEHARDCASEMUSTARBOUNDS}
		\\
	\| \rgeo [\Lunit \upmu]_- \|_{C^0(\Sigma_s^u)}
	= \| \Omega_-\|_{C^0(\Sigma_s^u)}
	& \leq 
		\left(1	
					+
					C \sqrt{\varepsilon}
		\right)
		\updelta.
	\label{E:LARGETIMEHARDCASEOMEGAMINUSBOUND}
\end{align}
\end{subequations}

\medskip

\noindent \underline{\textbf{Sharp estimates when $(u',\vartheta) \in \Vplus{t}{u}$}}
For each $t \geq 0,$ we define the ($t-$dependent) constant $\upgamma \geq 0$ by
\begin{align} \label{E:CRUCIALMAXIMALMUNOTDECAYINGPARAMETER}
	\upgamma & := \sup_{(u',\vartheta) \in \Vplus{t}{u}} \widetilde{\Omega}(t,u',\vartheta;t),
\end{align}
where $\widetilde{\Omega}$ is defined in \eqref{E:WIDETILDEBIGOMEGADEF}
and $ \Vplus{t}{u}$ is defined in \eqref{E:ANGLESANDUWITHNONDECAYINUPMUGBEHAVIOR}.

The following estimate holds:
\begin{align} \label{E:GAMMACONSTANTVERYSMALL}
	\upgamma \leq C \varepsilon.
\end{align}
If $0 \leq s_1 \leq s_2 \leq t,$ then the following estimates hold:
\begin{subequations}
\begin{align} \label{E:SHARPLOCALIZEDMUCANTGROWTOOFAST}
	\sup_{(u',\vartheta) \in \Vplus{t}{u}}
	\frac{\upmu(s_2,u',\vartheta)}{\upmu(s_1,u',\vartheta)}
	& \leq (1 + C \varepsilon) \frac{1 + \upgamma \ln(\myexp + s_2)}{1 + \upgamma \ln(\myexp + s_1)},
		\\
	\sup_{(u',\vartheta) \in \Vplus{t}{u}}
	\frac{\upmu(s_2,u',\vartheta)}{\upmu(s_1,u',\vartheta)}
	& \leq (1 + C \varepsilon) \ln(\myexp + s_2).
	\label{E:LOCALIZEDMUCANTGROWTOOFAST}
\end{align}
\end{subequations}
Furthermore, 
there exists a constant $C > 0$ 
such that if $0 < \Littleconone \leq 1$ is any constant, then
\begin{align} \label{E:LARGETIMELOCALIZEDMUNOTDECAYINGRATIOESTIMATE}
	\mathop{\sup_{1 \leq s_2^{\Littleconone} \leq s_1 \leq s_2 \leq t}}_{(u',\vartheta) \in \Vplus{t}{u}}
	\frac{\upmu(s_2,u',\vartheta)}{\upmu(s_1,u',\vartheta)}
	& \leq (1 + C \varepsilon) \frac{1}{\Littleconone}.
\end{align}
In addition, if $s \in [0,t]$ and $\Sigmaplus{s}{t}{u}$ is as defined in \eqref{E:SIGMAPLUS}, then we have
\begin{align} 
	\left\| \frac{\Lunit \upmu}{\upmu} \right\|_{C^0(\Sigmaplus{s}{t}{u})}
	& \leq (1 + C \varepsilon) \frac{\upgamma}{(1 + s)[1 + \upgamma \ln(\myexp + s)]}
	+ C \varepsilon \frac{\ln(\myexp + s)}{(1 + s)^2}.
		\label{E:KEYMUNOTDECAYINGLMUOVERMUBOUND}
\end{align}
Furthermore, if $s \in [0,t],$ then we have
\begin{align} 
	\left\| \frac{[\Lunit \upmu]_-}{\upmu} \right\|_{C^0(\Sigmaplus{s}{t}{u})}
	& \leq C \varepsilon \frac{\ln(\myexp + s)}{(1 + s)^2}.
		\label{E:KEYMUNOTDECAYINGMINUSPARTLMUOVERMUBOUND}
\end{align}




\medskip

\noindent \underline{\textbf{Sharp estimates when $(u',\vartheta) \in \Vminus{t}{u}$}}
Assume that the set $\Vminus{t}{u}$ defined in \eqref{E:ANGLESANDUWITHDECAYINUPMUGBEHAVIOR} 
is non-empty, and
consider a time interval $s \in [0,t].$ Let $\updelta > 0$
be as in \eqref{E:CRUCIALDELTADEF} and note that since $\Vminus{t}{u}$ is non-empty, we have
$\updelta = - \min_{(u',\vartheta) \in \Vminus{t}{u}} \widetilde{\Omega}(t,u',\vartheta;t).$
Then the following estimate holds:
\begin{align} \label{E:LOCALIZEDMUMUSTSHRINK}
	\mathop{\sup_{0 \leq s_1 \leq s_2 \leq t}}_{(u',\vartheta) \in \Vminus{t}{u}}
	\frac{\upmu(s_2,u',\vartheta)}{\upmu(s_1,u',\vartheta)}
	& \leq 1 + C \sqrt{\varepsilon}.
\end{align}
Furthermore, if $s \in [0,t]$ and $\Sigmaminus{s}{t}{u}$
is as defined in \eqref{E:SIGMAMINUS}, then the following estimate holds:
\begin{align} \label{E:LMUPLUSNEGLIGIBLEINSIGMAMINUS}
	\left\| [\Lunit \upmu]_+ \right\|_{C^0(\Sigmaminus{s}{t}{u})}
	& \leq C \varepsilon \frac{\ln(\myexp + s)}{(1 + s)^2}.
\end{align}
In addition, there exists
a constant $C > 0$ such that if
$0 \leq s \leq (1 - u) \frac{\varepsilon}{\updelta^2} - (1 - u) \leq t,$
then
\begin{subequations}
\begin{align} 
	\| \rgeo [\Lunit \upmu]_-\|_{C^0(\Sigmaminus{s}{t}{u})}
	& \leq 
		\left(
			1 + C \varepsilon
		\right)
		\left\lbrace
			\updelta 
			+ C \varepsilon \frac{\ln(\myexp + s)}{1 + s}
		\right\rbrace.
		\label{E:HYPERSURFACESHORTTIMEHARDCASEOMEGAMINUSBOUND}
\end{align}
Finally, there exists a constant $C > 0$ 
such that if $(1 - u) \frac{\varepsilon}{\updelta^2} - (1 - u) \leq s \leq t,$ then
\begin{align}		\label{E:HYPERSURFACELARGETIMEHARDCASEOMEGAMINUSBOUND}
	\| \rgeo [\Lunit \upmu]_- \|_{C^0(\Sigmaminus{s}{t}{u})}
	& \leq 
		\left(1	
					+
					C \sqrt{\varepsilon}
		\right)
		\updelta.
\end{align}
\end{subequations}




\noindent \underline{\textbf{Approximate time-monotonicity of $\upmu_{\star}^{-\Contwo}(s,u)$ under all conditions.}}
Let 
\begin{align*}
	\Contwo > 0
\end{align*}
be a real number.
There exists a constant $C > 0$ 
such that if $\Contwo \sqrt{\varepsilon} |\ln \varepsilon| \leq 1$ and
$0 \leq s_1 \leq s_2 \leq t,$ then
\begin{align} \label{E:MUSTARINVERSETOPOWERAMUSTGROWUPTOACONSTANT}
	\upmu_{\star}^{-\Contwo}(s_1,u) & \leq (1 + C \sqrt{\varepsilon}) \upmu_{\star}^{-\Contwo}(s_2,u).
\end{align}

\medskip
\noindent \underline{\textbf{A localized-in$-u-$and-$\vartheta$ signed upper bound for $\frac{\Lunit \upmu}{\upmu}$ under all conditions.}}
The following estimate holds:
\begin{align} \label{E:LUMUOVERMUSIGNEDUPPERBOUND}
	\frac{\Lunit \upmu(s,u,\vartheta)}
			 {\upmu(s,u,\vartheta)} & \leq C \varepsilon \frac{1}{1 + s}.
\end{align}

\end{proposition}

\begin{remark}[\textbf{A choice made out of convenience}]
	There is nothing special about the value $1/4$ in inequality \eqref{E:SMALLMUIMPLIESLMUISNEGATIVE}; 
	we could prove a similar estimate with ``$1/4$'' replaced by any positive constant less than $1$
	(but the smallness of $\varepsilon$ would have to be adapted to that constant).
\end{remark}

\begin{proof}

\noindent \textbf{Proof of \eqref{E:MUWEAKUPPERBOUND}:}
The estimate \eqref{E:MUWEAKUPPERBOUND} has already been proved as \eqref{E:C0BOUNDLDERIVATIVECRUCICALEIKONALFUNCTIONQUANTITIES}.

\ \\

\noindent \textbf{Proof of \eqref{E:POSITIVEPARTOFLMUOVERMUISSMALL}-\eqref{E:ALTERNATEPOSITIVEPARTOFLMUOVERMUISSMALL}:}
To prove \eqref{E:POSITIVEPARTOFLMUOVERMUISSMALL}, we note that by \eqref{E:EQUALRATIOS},
it suffices to bound the ratio
\begin{align} \label{E:SUFFICIENTRATIO}
	\frac{[\frac{\partial}{\partial s} \upmu_{(Approx)}(s,u,\vartheta;t)]_+}{\upmu_{(Approx)}(s,u,\vartheta;t)}
\end{align}
by the right-hand side of \eqref{E:POSITIVEPARTOFLMUOVERMUISSMALL}.
To bound the ratio \eqref{E:SUFFICIENTRATIO}, 
we first divide both sides of the identity \eqref{E:LMUAPPROXINTERMSOFMUAPPROX} 
by $\upmu_{(Approx)}(s,u,\vartheta;t)$ to deduce
\begin{align} \label{E:KEYLMUAPPROXMUAPPROXRATIO}
	\frac{\frac{\partial}{\partial s} \upmu_{(Approx)}(s,u,\vartheta;t)}{\upmu_{(Approx)}(s,u,\vartheta;t)}
	& = \frac{1}{\rgeo(s,u) \left\lbrace 1 + \ln \left(\frac{\rgeo(s,u)}{\rgeo(0,u)} \right) \right\rbrace}
		- \frac{1}{\upmu_{(Approx)}(s,u,\vartheta;t) \rgeo(s,u) \left\lbrace 1 + \ln \left(\frac{\rgeo(s,u)}{\rgeo(0,u)} \right) \right\rbrace}
		\\
	& \ \ + \frac{\widetilde{\Omega}(t,u,\vartheta;t) - \widetilde{M}(s,u,\vartheta;t)}	
					{\upmu_{(Approx)}(s,u,\vartheta;t) \rgeo(s,u) \left\lbrace 1 + \ln \left(\frac{\rgeo(s,u)}{\rgeo(0,u)} \right) \right\rbrace}
		- \frac{\widetilde{m}(s,u,\vartheta;t)}
				{\upmu_{(Approx)}(s,u,\vartheta;t)}.
		\notag
\end{align}
We now observe that the right-hand side of \eqref{E:KEYLMUAPPROXMUAPPROXRATIO} is negative 
(and hence the quantity \eqref{E:SUFFICIENTRATIO} is $0$) unless
\begin{align} \label{E:MUAPPROXNEARONEFORPOSITIVITY}
	\upmu_{(Approx)}(s,u,\vartheta;t) - 1
	\geq \frac{- \widetilde{\Omega}(t,u,\vartheta;t) + \widetilde{M}(s,u,\vartheta;t)}	
					{\rgeo(s,u) \left\lbrace 1 + \ln \left(\frac{\rgeo(s,u)}{\rgeo(0,u)} \right) \right\rbrace}
		+ \widetilde{m}(s,u,\vartheta;t) \rgeo(s,u) \left\lbrace 1 + \ln \left(\frac{\rgeo(s,u)}{\rgeo(0,u)} \right) \right\rbrace.
\end{align}
We may thus assume that inequality \eqref{E:MUAPPROXNEARONEFORPOSITIVITY} holds, 
in which case the estimates 
\eqref{E:BIGOMEGATRIVIALBOUND},
\eqref{E:LITTLEMEST}, 
and \eqref{E:BIGMEST} 
imply that
\begin{align} \label{E:ONEOVERMUAPPROXLMUPOSITIVEPARTARGUMENT}
	\frac{1}{\upmu_{(Approx)}(s,u,\vartheta;t)} 
	& \leq 1 
		+ C \varepsilon 
			\frac{1}
				{\rgeo(s,u) \left\lbrace 1 + \ln \left(\frac{\rgeo(s,u)}{\rgeo(0,u)} \right) \right\rbrace}.
\end{align}
Now since the positive part of the left-hand side of \eqref{E:KEYLMUAPPROXMUAPPROXRATIO} is bounded by the 
sum of the first term on the right-hand side and the sum of the magnitude of the last two terms on the right-hand side,
the estimates 
\eqref{E:BIGOMEGATRIVIALBOUND},
\eqref{E:LITTLEMEST}, 
\eqref{E:BIGMEST},
and \eqref{E:ONEOVERMUAPPROXLMUPOSITIVEPARTARGUMENT}
thus yield that 
\eqref{E:SUFFICIENTRATIO} is bounded by the right-hand side of 
\eqref{E:POSITIVEPARTOFLMUOVERMUISSMALL} as desired.

To prove \eqref{E:ALTERNATEPOSITIVEPARTOFLMUOVERMUISSMALL}, we use a slightly different argument that
takes into account the bound \eqref{E:MUAPPROXMINUSONEEASYBOUND}. 
Inserting this bound and 
the estimates 
\eqref{E:BIGOMEGATRIVIALBOUND},
\eqref{E:LITTLEMEST}, 
\eqref{E:BIGMEST},
and \eqref{E:ONEOVERMUAPPROXLMUPOSITIVEPARTARGUMENT}
into \eqref{E:KEYLMUAPPROXMUAPPROXRATIO}, 
we deduce that
\begin{align} \label{E:MUAPPROXRATIOBOUND}
	\frac{[\frac{\partial}{\partial s} \upmu_{(Approx)}(s,u,\vartheta;t)]_+}{\upmu_{(Approx)}(s,u,\vartheta;t)}
	& \leq \frac{1}{\rgeo(s,u) \left\lbrace 1 + \ln \left(\frac{\rgeo(s,u)}{\rgeo(0,u)} \right) \right\rbrace} 
				\frac{|\upmu_{(Approx)}(s,u,\vartheta;t) - 1|}{\upmu_{(Approx)}(s,u,\vartheta;t)} 
				\\
	& \ \ + \frac{1}{\rgeo(s,u) \left\lbrace 1 + \ln \left(\frac{\rgeo(s,u)}{\rgeo(0,u)} \right) \right\rbrace}
				\frac{1}{\upmu_{(Approx)}(s,u,\vartheta;t)}
				\left|
					\widetilde{\Omega}(t,u,\vartheta;t) - \widetilde{M}(s,u,\vartheta;t)
				\right|
			\notag \\
	& \ \ 
		+ \frac{1}{\upmu_{(Approx)}(s,u,\vartheta;t)}
			\left|
				\widetilde{m}(s,u,\vartheta;t)
			\right|
		\notag \\
	& \leq C \varepsilon \frac{1}{1 + s}.
		\notag
\end{align}
The desired estimate 
\eqref{E:ALTERNATEPOSITIVEPARTOFLMUOVERMUISSMALL}
now follows from
\eqref{E:EQUALRATIOS} 
and 
\eqref{E:MUAPPROXRATIOBOUND}.

\ \\

\noindent \textbf{Proof of \eqref{E:SMALLMUIMPLIESLMUISNEGATIVE}:}
To prove \eqref{E:SMALLMUIMPLIESLMUISNEGATIVE}, we first note 
that the identity \eqref{E:KEYLMUAPPROXMUAPPROXRATIO}
and the estimates 
\eqref{E:BIGOMEGATRIVIALBOUND},
\eqref{E:LITTLEMEST}, 
and
\eqref{E:BIGMEST} 
imply that
whenever
\begin{align}
	\upmu_{(Approx)}(s,u,\vartheta;t) \leq \frac{1}{3.5},
\end{align}
the following inequality necessarily holds:
\begin{align} \label{E:MUAPPROXISNEGATIVE}
	\frac{\partial}{\partial s} \upmu_{(Approx)}(s,u,\vartheta;t)
	\leq \frac{1 - 3.5 + C \varepsilon}{\rgeo(s,u) \left\lbrace 1 + \ln \left(\frac{\rgeo(s,u)}{\rgeo(0,u)} \right) \right\rbrace}.
\end{align}
The desired conclusion \eqref{E:SMALLMUIMPLIESLMUISNEGATIVE}
follows from the estimates 
\eqref{E:MUAMPLITUDENEARONE},
\eqref{E:MUAPPROXMISLIKEMU}, 
and \eqref{E:MUAPPROXISNEGATIVE},
and the identity
\begin{align}
	\Lunit \upmu(s,u,\vartheta) 
	& = \left\lbrace
				\mathring{\upmu}(u,\vartheta) - M(0,u,\vartheta;t)
			\right\rbrace
		\frac{\partial}{\partial s} \upmu_{(Approx)}(s,u,\vartheta;t),
\end{align}
which follows from \eqref{E:MUSPLIT}.

\ \\

\noindent \textbf{Proof of \eqref{E:EASYCASEMUSTARTABOUNDS}  - \eqref{E:EASYCASEOMEGAMINUSBOUNDS}:}	
To prove \eqref{E:EASYCASEMUSTARTABOUNDS} and
\eqref{E:APOWEREASYCASEMUSTARTABOUNDS}, we first note that since 
$\widetilde{\Omega}(t,u,\vartheta;t) \geq \widetilde{\Omega}_{(Min)}(t,u;t) \geq 0,$
\eqref{E:MUAPPROXDEF} and \eqref{E:BIGMEST} imply that
the following estimate holds for $s \in [0,t]:$ 
\begin{align} \label{E:MUAPPROXEASYBOUND}
	\upmu_{(Approx)}(s,u,\vartheta;t) \geq 1 - C \varepsilon.
\end{align}
Furthermore, 
\eqref{E:MUSPLIT}, 
\eqref{E:MUSTARDEF},
\eqref{E:MUAMPLITUDENEARONE},
and the estimate \eqref{E:MUAPPROXEASYBOUND} 
imply that the following estimate holds for $s \in [0,t]:$
\begin{align} \label{E:APOWERMUAPPROXEASYBOUND}
	\upmu_{\star}^{\Contwo}(s,u) \geq (1 - C \varepsilon)^{\Contwo}.
\end{align}
In particular, we have proved \eqref{E:EASYCASEMUSTARTABOUNDS}
(the upper bound is trivial).
Furthermore, since the function $h(y) := (1 - y)^{1/y}$ 
monotonically increases to
$e^{-1}$ as $y$ decreases to $0,$ we can set 
$y = C \varepsilon$ and use \eqref{E:APOWERMUAPPROXEASYBOUND}
and the trivial bound $\exp(- C \Contwo \varepsilon) \geq 1 - C \Contwo \varepsilon$
to conclude the desired estimate 
\eqref{E:APOWEREASYCASEMUSTARTABOUNDS} (the upper bound is trivial).

To prove \eqref{E:EASYCASEOMEGAMINUSBOUNDS}, we first decompose
\begin{align} \label{E:BIGOMEGASPLIT}
	\widetilde{\Omega}(s,u,\vartheta;t) 
	& = \widetilde{\Omega}(t,u,\vartheta;t) 
		+ \left\lbrace 
				\widetilde{\Omega}(s,u,\vartheta;t) - \widetilde{\Omega}(t,u,\vartheta;t) 
			\right\rbrace.
\end{align}
Then since $\widetilde{\Omega}(t,u,\vartheta;t) \geq 0,$ it follows from 
\eqref{E:BIGOMEAGDIFFERENCEEST} and \eqref{E:BIGOMEGASPLIT} that
\begin{align} \label{E:EASYCASEOMEGAMINUSFIRSTBOUND}
	\widetilde{\Omega}_-(s,u,\vartheta;t)
	& \leq \left| \widetilde{\Omega}(s,u,\vartheta;t) - \widetilde{\Omega}(t,u,\vartheta;t) \right|
		\leq C \varepsilon
				\frac{\ln(\myexp + s)}{1 + s}
				\left(
					\frac{t - s}{1 + t}
				\right)
		\leq C \varepsilon \frac{\ln(\myexp + s)}{1 + s}.
\end{align}
The desired estimate \eqref{E:EASYCASEOMEGAMINUSBOUNDS} 
now easily follows from 
definitions 
\eqref{E:BIGOMEGADEF}
and
\eqref{E:WIDETILDEBIGOMEGADEF} 
and
inequalities 
\eqref{E:MUAMPLITUDENEARONE}
and
\eqref{E:EASYCASEOMEGAMINUSFIRSTBOUND}.

\ \\

\noindent \textbf{Proof of \eqref{E:SHORTTIMEHARDCASEMUSTARBOUNDS}-\eqref{E:SHORTTIMEHARDCASEOMEGAMINUSBOUND}:}
We first prove \eqref{E:SHORTTIMEHARDCASEMUSTARBOUNDS}.
In order to simplify the presentation, we define 
\begin{align} \label{E:SMALLTIMEDEF}
	t_1 
	& := (1 - u) \frac{\varepsilon}{\updelta^2} - (1 - u).
\end{align}
Note that by \eqref{E:BIGOMEGATRIVIALBOUND}, the following simple inequality holds:
\begin{align} \label{E:DELTATRIVIALBLOUND}
	\updelta & \leq C \varepsilon.
\end{align}
Hence, it follows from 
\eqref{E:SMALLTIMEDEF}
and
\eqref{E:DELTATRIVIALBLOUND}
that if $\varepsilon$ is sufficiently small, then 
\begin{align} \label{E:T1ISPOSITIVE}
	t_1 \geq \frac{1}{C \varepsilon}.
\end{align}
We may thus assume that \eqref{E:T1ISPOSITIVE} holds.

For the estimates under consideration, we are assuming that
\begin{align}
	0 \leq s \leq t_1.
\end{align}
From \eqref{E:SMALLTIMEDEF}, \eqref{E:DELTATRIVIALBLOUND}, 
and the identities 
$\ln \left(\frac{\rgeo(t_1,u)}{\rgeo(0,u)} \right) = 
\ln \left( \frac{\varepsilon}{\updelta^2} \right)
\leq |\ln \updelta| + \left|\ln \frac{\updelta}{\varepsilon}\right| ,$
we deduce that for $0 \leq s \leq t_1,$ we have
\begin{align} \label{E:ONEMINUSDELTALOGSBOUND}
	1 - \updelta \ln \left(\frac{\rgeo(s,u)}{\rgeo(0,u)} \right) 
	& \geq 1 - \updelta \ln \left(\frac{\rgeo(t_1,u)}{\rgeo(0,u)} \right)
		\geq 1 - \updelta 
							\left\lbrace 
								|\ln \updelta| 
									+ \left|\ln \frac{\updelta}{\varepsilon}\right| 
							\right\rbrace
							\\
	& \geq 1 - C \varepsilon |\ln \varepsilon|.
		\notag
\end{align}	
It then follows from 
\eqref{E:MUAPPROXDEF},
\eqref{E:MUSPLIT},
\eqref{E:MUAMPLITUDENEARONE},
\eqref{E:BIGMEST},
\eqref{E:CRUCIALDELTADEF},
and \eqref{E:ONEMINUSDELTALOGSBOUND} 
that for $0 \leq s \leq t_1,$ we have
\begin{align}  \label{E:SHORTTIMEHARDCASEMULOWER}
	\upmu(s,u,\vartheta)
	& \geq  1 - C \varepsilon |\ln \varepsilon|.
\end{align}
Since inequality \eqref{E:SHORTTIMEHARDCASEMULOWER} is independent of $u$ and $\vartheta,$
the same estimate holds for $\upmu_{\star}(s,u).$ Hence, we have
\begin{align}  \label{E:POWERASHORTTIMEHARDCASEMULOWER}
	\upmu_{\star}^{\Contwo}(s,u)
	& \geq (1 - C \varepsilon |\ln \varepsilon|)^{\Contwo}.
\end{align}
In particular, we have proved \eqref{E:SHORTTIMEHARDCASEMUSTARBOUNDS}
(the upper bound is trivial). Furthermore, by arguing
as in our proof of \eqref{E:APOWEREASYCASEMUSTARTABOUNDS} and in particular
using \eqref{E:POWERASHORTTIMEHARDCASEMULOWER}, we 
deduce that
\begin{align} \label{E:MUSTARPOWERHARDCASEALMOSTCRUCIALLOWERBOUND}
	\upmu_{\star}^{\Contwo}(s,u)
	& \geq \exp(-C \Contwo \varepsilon |\ln \varepsilon|) 
	\geq 1 - C \Contwo \varepsilon |\ln \varepsilon| \geq 1 - C \sqrt{\varepsilon}.
\end{align}
We have thus shown that
the desired bounds \eqref{E:POWERASHORTTIMEHARDCASEMUSTARBOUNDS} hold for $0 \leq s \leq t_1$
(the upper bound is trivial).

To prove \eqref{E:SHORTTIMEHARDCASEOMEGAMINUSBOUND}, we
first use the splitting \eqref{E:BIGOMEGASPLIT} and the estimate \eqref{E:BIGOMEAGDIFFERENCEEST}
to deduce that
\begin{align} \label{E:SHORTTIMEHARDCASEWIDETILDEOMEGAMINUSUPPERBOUND}
	\widetilde{\Omega}_-(s,u,\vartheta;t) 
	& \leq
		\updelta 
		+ C \varepsilon \frac{\ln(\myexp + s)}{1 + s}.
\end{align}
Hence, the desired bound \eqref{E:SHORTTIMEHARDCASEOMEGAMINUSBOUND}
follows from 
definitions 
\eqref{E:BIGOMEGADEF}
and
\eqref{E:WIDETILDEBIGOMEGADEF}
and inequalities
\eqref{E:MUAMPLITUDENEARONE}
and
\eqref{E:SHORTTIMEHARDCASEWIDETILDEOMEGAMINUSUPPERBOUND}.

\ \\

\noindent \textbf{Proof of \eqref{E:LARGETIMEHARDCASEMUSTARBOUNDS}-\eqref{E:LARGETIMEHARDCASEOMEGAMINUSBOUND}:}
Let $t_1$ be as in \eqref{E:SMALLTIMEDEF}. We introduce the logarithmic variables
\begin{align}
	\uptau 
	&:= \ln \left(\frac{\rgeo(t,u)}{\rgeo(0,u)} \right),
		\label{E:TAUCHOV} \\
	\uptau_1 
	&:= \ln \left(\frac{\rgeo(t_1,u)}{\rgeo(0,u)} \right)
		= \ln \left(\frac{\varepsilon}{\updelta^2}\right),
	\label{E:TAU1CHOV} \\
	\upsigma 
	&:= \ln \left(\frac{\rgeo(s,u)}{\rgeo(0,u)} \right).
	\label{E:SIGMACHOV}
\end{align}

To prove \eqref{E:LARGETIMEHARDCASEMUSTARBOUNDS},
we first note that \eqref{E:MUAPPROXDEF} implies 
that for $(t,u) \in [0,\Tboot) \times [0,U_0],$ we have
\begin{align} \label{E:UAPPROXATTT}
	\upmu_{(Approx)}(t,u,\vartheta;t) = 1 - \updelta \ln \left( \frac{\rgeo(t,u)}{\rgeo(0,u)} \right).
\end{align}
Since 
\eqref{E:BAUMPUISPOSITIVE},
\eqref{E:MUSPLIT},
and
\eqref{E:MUAMPLITUDENEARONE}
imply that $0 < \upmu_{(Approx)},$ 
it follows from \eqref{E:UAPPROXATTT} that
\begin{align} \label{E:TNOTTOOBIG}
	 \updelta \ln \left(\frac{\rgeo(t,u)}{\rgeo(0,u)} \right) < 1.
\end{align}



We now use
\eqref{E:MUAPPROXDEF},
\eqref{E:MUSPLIT},
and
\eqref{E:BIGMEST}
to deduce that for $s \in [t_1, t],$ we have
\begin{align} \label{E:MULARGETIMEFIRSTLOWERBOUND}
	\upmu(s,u,\vartheta)
	& \geq 
		\left\lbrace
				\mathring{\upmu}(u,\vartheta) - M(0,u,\vartheta;t)
		\right\rbrace
		\left\lbrace	
			1 
			- \updelta \ln \left(\frac{\rgeo(s,u)}{\rgeo(0,u)} \right) 
			- C \varepsilon
			\frac{\ln(\myexp + s)}{1 + s}
			\left(
				\frac{t - s}{1 + t}
			\right) 
		\right\rbrace.
\end{align}
Using \eqref{E:MULARGETIMEFIRSTLOWERBOUND} and the inequalities
\begin{align} 
	\frac{\ln(\myexp + s)}{1 + s}
	& \leq C \frac{1 + \upsigma}{\exp(\upsigma)},
		\label{E:ASIMPLELOGINEQUALITY} \\
	\frac{t - s}{1 + t}
	& \leq \frac{\rgeo(t,u) - \rgeo(s,u)}{\rgeo(t,u)},
		\\
	0 \leq \frac{\rgeo(t,u) - \rgeo(s,u)}{\rgeo(t,u)}
	& = \frac{\exp(\uptau) - \exp(\upsigma)}{\exp(\uptau)} 
		= 1 - \exp(\upsigma - \uptau) 
	\leq \uptau - \upsigma,
	\label{E:SIMPLEEXPINEQUALITY}	
\end{align}
we deduce that for $s \in [t_1, t],$ we have
\begin{align} \label{E:MULARGETIMESECONDLOWERBOUND}
	\upmu(s,u,\vartheta)
	& \geq 
		\left\lbrace
				\mathring{\upmu}(u,\vartheta) - M(0,u,\vartheta;t)
		\right\rbrace
		\left\lbrace	
			1 
			- \updelta \upsigma
			- C \varepsilon
			\frac{1 + \upsigma}{\exp(\upsigma)}
			(\uptau - \upsigma)  
		\right\rbrace
		\\
	& =
		\left\lbrace
				\mathring{\upmu}(u,\vartheta) - M(0,u,\vartheta;t)
		\right\rbrace
		\left\lbrace	
			1 
			- \updelta \uptau
			+ \updelta
			\left(
				1
				- C \frac{\varepsilon}{\updelta}
				\frac{1 + \upsigma}{\exp(\upsigma)}
			\right)
			(\uptau - \upsigma)  
		\right\rbrace.
		\notag
\end{align}
We now estimate the crucially important term 
$\frac{\varepsilon}{\updelta} \frac{1 + \upsigma}{\exp(\upsigma)}$
appearing in the second line of \eqref{E:MULARGETIMESECONDLOWERBOUND}.
To this end, we use the fact that $s \in [t_1, t]$ 
and inequalities 
\eqref{E:DELTATRIVIALBLOUND}
and
\eqref{E:T1ISPOSITIVE} 
to deduce that
\begin{align} \label{E:KEYPRODUCTBOUND}
	\frac{\varepsilon}{\updelta}
		\frac{1 + \upsigma}{\exp(\upsigma)}
	& \leq \frac{\varepsilon}{\updelta}
			\frac{1 + \uptau_1}{\exp(\uptau_1)}
		= 
		\updelta
		\left\lbrace
			1
			+
			\ln \left(\frac{\varepsilon}{\updelta^2}\right)
		\right\rbrace
		\leq C \varepsilon |\ln \varepsilon|.
\end{align}
Inserting the bound \eqref{E:KEYPRODUCTBOUND}
into inequality \eqref{E:MULARGETIMESECONDLOWERBOUND}
and using the bound \eqref{E:MUAMPLITUDENEARONE},
we deduce the following lower bound for $s \in [t_1, t]:$
\begin{align} \label{E:MULARGETIMETHIRDLOWERBOUND}
	\upmu(s,u,\vartheta)
	& \geq 	
		\left\lbrace
				\mathring{\upmu}(u,\vartheta) - M(0,u,\vartheta;t)
		\right\rbrace
		\left\lbrace	
			1 
			- \updelta \uptau
			+ \updelta
			\left(
				1
				- C \frac{\varepsilon}{\updelta}
				\frac{1 + \upsigma}{\exp(\upsigma)}
			\right)
			(\uptau - \upsigma)  
		\right\rbrace
		\\
		& \geq 	
		\left\lbrace
				\mathring{\upmu}(u,\vartheta) - M(0,u,\vartheta;t)
		\right\rbrace
		\left\lbrace	
			\left(
				1 
				- C \varepsilon |\ln \varepsilon|
			\right)
			\left(
				1 
				- \updelta \uptau
			\right)
			+ 
			\updelta
			\left(
				1 
				- C \varepsilon |\ln \varepsilon|
			\right)
			(\uptau - \upsigma)  
		\right\rbrace
		\notag \\
		& =
			\left\lbrace
				\mathring{\upmu}(u,\vartheta) - M(0,u,\vartheta;t)
			\right\rbrace
			\left(
				1 
				- C \varepsilon |\ln \varepsilon|
			\right)
			(1 - \updelta \upsigma)
			\notag
			\\
		& \geq
			\left(
				1 
				- C \varepsilon |\ln \varepsilon|
			\right)
			(1 - \updelta \upsigma).
			\notag
\end{align}
Note that by \eqref{E:TNOTTOOBIG}, the right-hand side of \eqref{E:MULARGETIMETHIRDLOWERBOUND}
is positive.
Since the final inequality in \eqref{E:MULARGETIMETHIRDLOWERBOUND}  
is independent of $\vartheta$ 
and since $- \updelta \upsigma = - \updelta \upsigma(s,u)$ is decreasing in $u,$
the same estimate holds for $\upmu_{\star}(s,u).$
We have thus deduced the desired lower bound in \eqref{E:LARGETIMEHARDCASEMUSTARBOUNDS}.
A similar argument yields the upper bound. Then by using an
argument similar to the one we used just after inequality \eqref{E:APOWERMUAPPROXEASYBOUND}, 
we also deduce from the bounds \eqref{E:LARGETIMEHARDCASEMUSTARBOUNDS} that
the upper and lower bounds in \eqref{E:POWERALARGETIMEHARDCASEMUSTARBOUNDS} hold.

To prove \eqref{E:LARGETIMEHARDCASEOMEGAMINUSBOUND}, we first decompose
$\widetilde{\Omega}(s,u,\vartheta;t)$ as in \eqref{E:BIGOMEGASPLIT} and use 
inequalities
\eqref{E:BIGOMEAGDIFFERENCEEST}
and
\eqref{E:ASIMPLELOGINEQUALITY}
to deduce
\begin{align} \label{E:BIGOMEGAFIRSTHARDESTIMATE}
	\widetilde{\Omega}_-(s,u,\vartheta;t) 
	& \leq \updelta
			+ C \varepsilon
			\frac{1 + \upsigma}{\exp(\upsigma)}
			= \updelta
			\left\lbrace
				1
				+ C \frac{\varepsilon}{\updelta}
					\frac{1 + \upsigma}{\exp(\upsigma)}
			\right\rbrace.
\end{align}
Arguing as in our proof of 
\eqref{E:KEYPRODUCTBOUND},
we deduce from inequality \eqref{E:BIGOMEGAFIRSTHARDESTIMATE} that 
\begin{align} \label{E:BIGOMEGASECONDHARDESTIMATE}
	\widetilde{\Omega}_-(s,u,\vartheta;t) 
	& \leq 
		\left(
			1 
			+ C \varepsilon |\ln \varepsilon|
		\right)
		\updelta
		\leq 
		\left(
			1 
			+ C \sqrt{\varepsilon}
		\right)
		\updelta.
\end{align}
From
definition
\eqref{E:WIDETILDEBIGOMEGADEF}
and the estimates \eqref{E:MUAMPLITUDENEARONE}
and \eqref{E:BIGOMEGASECONDHARDESTIMATE},
we deduce the following bound for $s \in [t_1, t]:$
\begin{align}\label{E:LARGETIMEHARDCASEOMEGAMINUSKEYBOUNDPROOF}
	\Omega_-(s,u,\vartheta) 
	& \leq 
		\left(
			1 
			+ C \sqrt{\varepsilon}
		\right)
		\updelta.
\end{align}	
Since inequality \eqref{E:LARGETIMEHARDCASEOMEGAMINUSKEYBOUNDPROOF}
is independent of $u$ and $\vartheta,$
the desired bound \eqref{E:LARGETIMEHARDCASEOMEGAMINUSBOUND} now 
follows from definition \eqref{E:BIGOMEGADEF}.

\ \\

\noindent \textbf{Proof of \eqref{E:GAMMACONSTANTVERYSMALL}:}
Inequality \eqref{E:GAMMACONSTANTVERYSMALL} follows from \eqref{E:BIGOMEGATRIVIALBOUND}.

\ \\
\noindent \textbf{Proof of \eqref{E:SHARPLOCALIZEDMUCANTGROWTOOFAST}-\eqref{E:LARGETIMELOCALIZEDMUNOTDECAYINGRATIOESTIMATE}:}
To prove \eqref{E:SHARPLOCALIZEDMUCANTGROWTOOFAST}, 
we first note that by \eqref{E:MUSPLIT},
it suffices to show that 
for each $(u',\vartheta) \in \Vplus{t}{u}$
and each $0 \leq s_1 \leq s_2 \leq t,$
the ratio $\frac{\upmu_{(Approx)}(s_2,u',\vartheta)}{\upmu_{(Approx)}(s_1,u',\vartheta)}$
is $\leq$ the right-hand side of \eqref{E:SHARPLOCALIZEDMUCANTGROWTOOFAST}.
To this end, we use \eqref{E:MUAPPROXDEF} and the estimates 
\eqref{E:BIGOMEGATRIVIALBOUND} and \eqref{E:BIGMEST} to deduce that
\begin{align} \label{E:MUOFS2OVERMUOFS1FIRSTBOUND}
	\frac{\upmu_{(Approx)}(s_2,u',\vartheta)}{\upmu_{(Approx)}(s_1,u',\vartheta)} 
	& = \frac{1 + \mathcal{O}(\varepsilon) + \widetilde{\Omega}(t,u',\vartheta;t) \ln(\myexp + s_2)}
					 {1 + \mathcal{O}(\varepsilon) + \widetilde{\Omega}(t,u',\vartheta;t) \ln(\myexp + s_1)}
					 \\
	& \leq (1 + C \varepsilon)
					\frac{1 + \widetilde{\Omega}(t,u',\vartheta;t) \ln(\myexp + s_2)}
					 {1 + \widetilde{\Omega}(t,u',\vartheta;t) \ln(\myexp + s_1)}.
					 \notag
\end{align}
Since the right-hand side of \eqref{E:MUOFS2OVERMUOFS1FIRSTBOUND}
is increasing when viewed as a function of $\widetilde{\Omega}(t,u',\vartheta;t),$
it follows that the right-hand side of \eqref{E:MUOFS2OVERMUOFS1FIRSTBOUND} is
\begin{align}
	\leq (1 + C \varepsilon)
					\frac{1 + \upgamma \ln(\myexp + s_2)}
					 {1 + \upgamma \ln(\myexp + s_1)}
\end{align}
as desired.

The estimate \eqref{E:LOCALIZEDMUCANTGROWTOOFAST} follows from 
\eqref{E:SHARPLOCALIZEDMUCANTGROWTOOFAST}
and \eqref{E:GAMMACONSTANTVERYSMALL}.

To deduce \eqref{E:LARGETIMELOCALIZEDMUNOTDECAYINGRATIOESTIMATE}, 
we use \eqref{E:SHARPLOCALIZEDMUCANTGROWTOOFAST}
and the following simple estimate, which holds when 
$0 < \Littleconone \leq 1$
and
$1 \leq s_2^{\Littleconone} \leq s_1 \leq s_2:$
$\ln(\myexp + s_1) \geq \Littleconone \ln(\myexp + s_2).$

\ \\
\noindent \textbf{Proof of \eqref{E:KEYMUNOTDECAYINGLMUOVERMUBOUND}:}
To prove \eqref{E:KEYMUNOTDECAYINGLMUOVERMUBOUND}, we first note that by \eqref{E:EQUALRATIOS},
it suffices to bound the ratio on the left-hand side of \eqref{E:KEYLMUAPPROXMUAPPROXRATIO} 
on the domain $(s,u',\vartheta) \in \Sigmaplus{s}{t}{u}$
by the right-hand side of \eqref{E:KEYMUNOTDECAYINGLMUOVERMUBOUND}.
To this end, for each $(u',\vartheta) \in \Vplus{t}{u},$
we use \eqref{E:MUAPPROXDEF} and \eqref{E:LMUAPPROXRELATION} 
and the estimates \eqref{E:BIGOMEGATRIVIALBOUND}, \eqref{E:LITTLEMEST}, and \eqref{E:BIGMEST} to deduce that
\begin{align} \label{E:KEYINITIALINEQUALITYMUNOTDECAYINGLMUOVERMU}
	\frac{\frac{\partial}{\partial s} \upmu_{(Approx)}(s,u',\vartheta;t)}{\upmu_{(Approx)}(s,u',\vartheta;t)}
	& = \frac{1}{\rgeo(s,u')}
			\frac{\widetilde{\Omega}(t,u',\vartheta;t)}
					 {\left[1 + \widetilde{\Omega}(t,u',\vartheta;t) \ln \left(\frac{\rgeo(s,u')}{\rgeo(0,u')} \right) \right]}
					+ \mathcal{O}\left(\varepsilon \frac{\ln(\myexp + s)}{(1 + s)^2} \frac{(t-s)}{(1+t)}\right).
\end{align}
Since the ratio 
$\frac{\widetilde{\Omega}(t,u',\vartheta;t)}
{\left[1 + \widetilde{\Omega}(t,u',\vartheta;t) \ln \left(\frac{\rgeo(s,u')}{\rgeo(0,u')} \right) \right]}$ 
is increasing when
viewed as a function of $\widetilde{\Omega}(t,u',\vartheta;t),$
it follows that the right-hand side of \eqref{E:KEYINITIALINEQUALITYMUNOTDECAYINGLMUOVERMU} is bounded 
by
\begin{align} \label{E:ALMOSTPROVEDKEYMUNOTDECAYINGLMUOVERMUBOUND}
	\leq
	\frac{1}{\rgeo(s,u')}
			\frac{\upgamma}
					 {\left[1 + \upgamma \ln \left(\frac{\rgeo(s,u')}{\rgeo(0,u')} \right) \right]}
					 + \mathcal{O}\left(\varepsilon \frac{\ln(\myexp + s)}{(1 + s)^2} \right).
\end{align}
The desired estimate \eqref{E:KEYMUNOTDECAYINGLMUOVERMUBOUND} now follows 
easily from 
\eqref{E:GAMMACONSTANTVERYSMALL}
and
\eqref{E:ALMOSTPROVEDKEYMUNOTDECAYINGLMUOVERMUBOUND}.

\ \\
\noindent \textbf{Proof of \eqref{E:KEYMUNOTDECAYINGMINUSPARTLMUOVERMUBOUND}:}
To prove \eqref{E:KEYMUNOTDECAYINGMINUSPARTLMUOVERMUBOUND}, we
first use the splitting \eqref{E:BIGOMEGASPLIT} and the estimate \eqref{E:BIGOMEAGDIFFERENCEEST}
to deduce that for $(s,u',\vartheta) \in \Sigmaplus{s}{t}{u},$ we have
\begin{align} \label{E:AGAINSHORTTIMEHARDCASEWIDETILDEOMEGAMINUSUPPERBOUND}
	\widetilde{\Omega}_-(s,u',\vartheta;t) 
	& \leq
		 C \varepsilon \frac{\ln(\myexp + s)}{1 + s}.
\end{align}
It thus follows from
\eqref{E:BIGOMEGADEF},
\eqref{E:WIDETILDEBIGOMEGADEF},
\eqref{E:MUAMPLITUDENEARONE},
and
\eqref{E:AGAINSHORTTIMEHARDCASEWIDETILDEOMEGAMINUSUPPERBOUND} 
that
\begin{align}
	[\Lunit \upmu]_-(s,u',\vartheta) 
	& \leq C \varepsilon \frac{\ln(\myexp + s)}{(1 + s)^2},
		\label{E:LUNITMINUSNEGLIGIBLEINDECAYINGREGION} \\
	\Lunit \upmu(s,u',\vartheta)  
	& \geq - C \varepsilon \frac{\ln(\myexp + s)}{(1 + s)^2}.
	\label{E:LUNITUPMUNOTVERYNEGATIVEINGOWINGREGION}
\end{align}
Integrating \eqref{E:LUNITUPMUNOTVERYNEGATIVEINGOWINGREGION} along the integral
curves of $\Lunit = \frac{\partial}{\partial s}$ emanating from $\Sigma_0$ and using the small-data estimate \eqref{E:MUINITIALDATAESTIMATE},
we deduce that 
\begin{align} \label{E:UPMUNOTSMALLINGROWINGREGION}
	\upmu(s,u,\vartheta)
	& \leq 1 - C \varepsilon.
\end{align}
The desired estimate \eqref{E:KEYMUNOTDECAYINGMINUSPARTLMUOVERMUBOUND}
now follows from
\eqref{E:LUNITMINUSNEGLIGIBLEINDECAYINGREGION}
and 
\eqref{E:UPMUNOTSMALLINGROWINGREGION}.






\ \\
\noindent \textbf{Proof of \eqref{E:LOCALIZEDMUMUSTSHRINK}:}
Note that we are investigating points $(s,u',\vartheta)$ such that
$(u',\vartheta) \in \Vminus{t}{u}$
(see definition \eqref{E:ANGLESANDUWITHDECAYINUPMUGBEHAVIOR}).
We start by repeating the proof of inequality \eqref{E:SHORTTIMEHARDCASEMULOWER},
but with $\updelta$ in the proof
everywhere replaced by $\widetilde{\Omega}(t,u,\vartheta;t),$
which changes the value of $t_1$ as defined in \eqref{E:SMALLTIMEDEF}.
The net effect is that for this new value of $t_1,$ 
if $s \in [0,t_1],$ then we have
\begin{align}  \label{E:AGAINSHORTTIMEHARDCASEMULOWER}
	\upmu(s,u,\vartheta)
	& \geq  1 - C \sqrt{\varepsilon}.
\end{align}
In addition, by making straightforward modifications to the proof
of \eqref{E:AGAINSHORTTIMEHARDCASEMULOWER}, we 
find that if $s \in [0,t_1],$ then
\begin{align}  \label{E:UPPERAGAINSHORTTIMEHARDCASEMULOWER}
	\upmu(s,u,\vartheta)
	& \leq  1 + C \sqrt{\varepsilon}.
\end{align}

Similarly, we repeat the proof of inequality \eqref{E:MULARGETIMETHIRDLOWERBOUND},
but with $\updelta$ in the proof
everywhere replaced by $\widetilde{\Omega}(t,u,\vartheta;t).$
As above, this replacement 
changes the value of 
$t_1$ as defined in \eqref{E:SMALLTIMEDEF}
and $\uptau_1$ as defined in \eqref{E:TAU1CHOV}.
In total, we find that for this new value of $t_1,$
if $s \in [t_1,t],$
then we have
\begin{align} \label{E:SHRINKINGLOWERBOUND}
\upmu(s,u',\vartheta)
& \geq
			\left(
				1 
				- C \sqrt{\varepsilon}
			\right)
			\left\lbrace
				1 
				+ 
				\widetilde{\Omega}(t,u',\vartheta;t)
				\ln 
					\left(
						\frac{\rgeo(s,u')}
							{\rgeo(0,u')}
					\right)
			\right\rbrace.
\end{align}
In addition, by making straightforward modifications
to the proof of \eqref{E:SHRINKINGLOWERBOUND}, 
we find that for $s \in [t_1,t],$ we have
\begin{align} \label{E:SHRINKINGUPPERBOUND}
\upmu(s,u',\vartheta)
& \leq
			\left(
				1 
				+ C \sqrt{\varepsilon}
			\right)
			\left\lbrace
				1 
				+ 
				\widetilde{\Omega}(t,u',\vartheta;t)
				\ln 
					\left(
						\frac{\rgeo(s,u')}
							{\rgeo(0,u')}
					\right)
			\right\rbrace.
\end{align}

We now separately investigate the short-time regime $s_2 \leq t_1$ and the large-time regime
$s_2 \geq t_1,$ where we are referring to the new value of $t_1$ defined above.
In the short-time regime,
we easily bound the ratio
$\frac{\upmu(s_2,u',\vartheta)}{\upmu(s_1,u',\vartheta)}$
by $\leq$ the right-hand side of
\eqref{E:LOCALIZEDMUMUSTSHRINK}
with the help of 
\eqref{E:AGAINSHORTTIMEHARDCASEMULOWER}
and \eqref{E:UPPERAGAINSHORTTIMEHARDCASEMULOWER}.

We now consider the large-time regime $s_2 \geq t_1.$ We split this regime into two sub-cases:
$s_1 \geq t_1$ and $s_1 \leq t_1.$ In the first sub-case,
we bound the ratio
$\frac{\upmu(s_2,u',\vartheta)}{\upmu(s_1,u',\vartheta)}$
by $\leq$ the right-hand side of
\eqref{E:LOCALIZEDMUMUSTSHRINK}
by using the estimates
\eqref{E:SHRINKINGLOWERBOUND} and \eqref{E:SHRINKINGUPPERBOUND},
the assumption $\widetilde{\Omega}(t,u',\vartheta;t) < 0,$
and the fact that
$\ln \left( \frac{\rgeo(s,u')}{\rgeo(0,u')} \right)$ is increasing in $s.$
In the second sub-case $s_1 \leq t_1,$ 
we bound the ratio
$\frac{\upmu(s_2,u',\vartheta)}{\upmu(s_1,u',\vartheta)}$
by $\leq$ the right-hand side of
\eqref{E:LOCALIZEDMUMUSTSHRINK}
by using the estimates
\eqref{E:AGAINSHORTTIMEHARDCASEMULOWER}
and
\eqref{E:SHRINKINGUPPERBOUND}.

\ \\
 
\noindent \textbf{Proof of \eqref{E:LMUPLUSNEGLIGIBLEINSIGMAMINUS}:}
To prove \eqref{E:LMUPLUSNEGLIGIBLEINSIGMAMINUS}, we
first use the splitting \eqref{E:BIGOMEGASPLIT} and the estimate \eqref{E:BIGOMEAGDIFFERENCEEST}
to deduce that for $(s,u',\vartheta) \in \Sigmaminus{s}{t}{u},$ we have
\begin{align} \label{E:EASYCASEWIDETILDEOMEGAPLUSUPPERBOUND}
	\widetilde{\Omega}_+(s,u',\vartheta;t) 
	& \leq
			C \varepsilon \frac{\ln(\myexp + s)}{1 + s}.
\end{align}
From \eqref{E:BIGOMEGADEF}
and
\eqref{E:WIDETILDEBIGOMEGADEF}
and the estimates
\eqref{E:MUAMPLITUDENEARONE}
and
\eqref{E:EASYCASEWIDETILDEOMEGAPLUSUPPERBOUND},
we conclude the desired estimate \eqref{E:EASYCASEWIDETILDEOMEGAPLUSUPPERBOUND}.

\ \\ 

\noindent \textbf{Proof of \eqref{E:HYPERSURFACESHORTTIMEHARDCASEOMEGAMINUSBOUND}-\eqref{E:HYPERSURFACELARGETIMEHARDCASEOMEGAMINUSBOUND}:}
Inequality \eqref{E:HYPERSURFACESHORTTIMEHARDCASEOMEGAMINUSBOUND}
can be proved by using essentially the same argument that we used to prove \eqref{E:SHORTTIMEHARDCASEOMEGAMINUSBOUND}.

Inequality \eqref{E:HYPERSURFACELARGETIMEHARDCASEOMEGAMINUSBOUND} can be proved by using
essentially the same argument that we used to prove \eqref{E:LARGETIMEHARDCASEOMEGAMINUSBOUND}.

\ \\

\noindent \textbf{Proof of \eqref{E:MUSTARINVERSETOPOWERAMUSTGROWUPTOACONSTANT}:}
We separately investigate the cases 
$\widetilde{\Omega}_{(Min)}(t,u;t) \geq 0$
and $\widetilde{\Omega}_{(Min)}(t,u;t) < 0.$
In the former case, the desired estimate
\eqref{E:MUSTARINVERSETOPOWERAMUSTGROWUPTOACONSTANT}
follows from
\eqref{E:APOWEREASYCASEMUSTARTABOUNDS}.

We now investigate the remaining case in which 
(as in \eqref{E:CRUCIALDELTADEF}) $\updelta := - \widetilde{\Omega}_{(Min)}(t,u;t) > 0.$
We separately investigate
the short-time regime 
$s_2 \leq t_1$ and the large-time regime
$s_2 \geq t_1,$ where $t_1$ is defined in \eqref{E:SMALLTIMEDEF}.
In the short-time regime $s_2 \leq t_1,$ the desired 
estimate \eqref{E:MUSTARINVERSETOPOWERAMUSTGROWUPTOACONSTANT}
follows from \eqref{E:POWERASHORTTIMEHARDCASEMUSTARBOUNDS}.

We now consider the large-time regime $s_2 \geq t_1.$ This regime splits into two sub-cases:
$s_1 \geq t_1$ and $s_1 \leq t_1.$ In the first sub-case, the desired estimate \eqref{E:MUSTARINVERSETOPOWERAMUSTGROWUPTOACONSTANT} 
follows from \eqref{E:POWERALARGETIMEHARDCASEMUSTARBOUNDS} and the fact that 
$1 - \updelta \ln \left( \frac{\rgeo(s,u)}{\rgeo(0,u)} \right)$ is decreasing in $s.$
In the second sub-case, the desired estimate \eqref{E:MUSTARINVERSETOPOWERAMUSTGROWUPTOACONSTANT} follows from 
\eqref{E:POWERASHORTTIMEHARDCASEMUSTARBOUNDS}
and \eqref{E:POWERALARGETIMEHARDCASEMUSTARBOUNDS}.

\ \\

\noindent \textbf{Proof of \eqref{E:LUMUOVERMUSIGNEDUPPERBOUND}:}
If $\Lunit \upmu(s,u,\vartheta) < 0,$ then the bound \eqref{E:LUMUOVERMUSIGNEDUPPERBOUND} is trivial.
If $\Lunit \upmu(s,u,\vartheta) \geq 0,$ then the bound \eqref{E:LUMUOVERMUSIGNEDUPPERBOUND} follows from
\eqref{E:ALTERNATEPOSITIVEPARTOFLMUOVERMUISSMALL}.

\ \\

\noindent \textbf{Proof of \eqref{E:UNIFORMBOUNDFORMRADMUOVERMU}:}

\noindent \emph{Case 1: $\widetilde{\Omega}_{(Min)}(t,U_0;t) \geq 0$:}
In this (easy) case, \eqref{E:MUAPPROXDEF} and \eqref{E:MUSPLIT}
and the estimates \eqref{E:MUAMPLITUDENEARONE} and \eqref{E:BIGMEST} imply that
\begin{align} \label{E:UPMUEASYLOWERBOUND}
	\upmu(s,u,\vartheta)
	& \geq 1 - C \varepsilon.
\end{align}
Furthermore, from the estimate \eqref{E:C0BOUNDCRUCIALEIKONALFUNCTIONQUANTITIES}, we have
\begin{align} \label{E:RADUPMUEASYBOUND}
	\| \Rad \upmu \|_{C^0(\Sigma_s^u)}
	& \leq C \varepsilon \ln(\myexp + s).
\end{align}
Hence, in this case, 
we can use
\eqref{E:UPMUEASYLOWERBOUND} and \eqref{E:RADUPMUEASYBOUND}
to deduce that
$\frac{[\Rad \upmu(s,u,\vartheta)]_+}{\upmu(s,u,\vartheta)}$
is $\leq$ the second term on the right-hand side
of \eqref{E:UNIFORMBOUNDFORMRADMUOVERMU} as desired.

\ \\

\noindent \emph{Case 2: Mild decay rate:}
In this case, we assume that $\widetilde{\Omega}_{(Min)}(t,U_0;t) < 0.$
Using the notation $\updelta := - \widetilde{\Omega}_{(Min)}(t,U_0;t) > 0$ 
much as in \eqref{E:CRUCIALDELTADEF}, we define the \emph{open} sets
\begin{align}
	\mathcal{V}(U_0;t) 
	& := \left\lbrace
				(u,\vartheta) \in [0,U_0] \times \mathbb{S}^2 \ | \ \widetilde{\Omega}(t,u,\vartheta;t) < - \frac{\updelta}{2}
				\right\rbrace,
				\\
	\mathcal{U}(s,U_0;t) 
	& := \left\lbrace
				(s,u,\vartheta) \in \Sigma_s^{U_0} \ | \ (u,\vartheta) \in \mathcal{V}(U_0;t)
				\right\rbrace.
\end{align}
Roughly, $\mathcal{U}(s,U_0;t)$ is the ``dangerous'' subset of $\Sigma_s^{U_0}$
where $\upmu$ can decay at a relatively rapid rate.

Now if $(s,u,\vartheta) \in \Sigma_s^{U_0} \backslash \mathcal{U}(s,U_0;t),$ then 
we can slightly modify our proof of \eqref{E:UPMUEASYLOWERBOUND}
(which was based in part on \eqref{E:MUAPPROXDEF})
and also use the bound \eqref{E:TNOTTOOBIG} 
to deduce the following estimate:
\begin{align}  \label{E:UPMUCASE2LOWERBOUND}
	\upmu(s,u,\vartheta)
	& \geq (1 - C \varepsilon)
		\left\lbrace
			1
			- \frac{1}{2} \updelta \ln \left( \frac{\rgeo(s,u)}{\rgeo(0,u)} \right) 
		\right\rbrace
		\geq (1 - C \varepsilon)
		\left\lbrace
			1
			- \frac{1}{2} \updelta \ln \left( \frac{\rgeo(t,u)}{\rgeo(0,u)} \right) 
		\right\rbrace
			\\
	& \geq \frac{1}{2} (1 - C \varepsilon).
	\notag
\end{align}
From \eqref{E:RADUPMUEASYBOUND} and \eqref{E:UPMUCASE2LOWERBOUND}, 
we see that 
$\frac{[\Rad \upmu(s,u,\vartheta)]_+}{\upmu(s,u,\vartheta)}$
is $\leq$ the second term on the right-hand side
of \eqref{E:UNIFORMBOUNDFORMRADMUOVERMU} as desired,
as in Case 1.

\ \\

\noindent \emph{Case 3: Dangerous decay rate, but short times:}
In this case, we assume that $\updelta := - \widetilde{\Omega}_{(Min)}(t,U_0;t) > 0,$
that $(s,u,\vartheta) \in \mathcal{U}(s,U_0;t),$ 
and that $s \leq t_1,$ where we set
\begin{align}
	t_1 & := (1-u)\exp(\uptau_1) - (1 - u), 
		\label{E:SUBCASE3T1} \\
	\uptau_1 & := \frac{1}{2 \updelta}.
\end{align}
Using the estimate \eqref{E:DELTATRIVIALBLOUND}, we conclude that
if $\varepsilon$ is sufficiently small, then
\begin{align} \label{E:T1ISPOSITIVECASE4}
	\uptau_1 \geq \frac{1}{C \varepsilon}.
\end{align}

As in our proof of \eqref{E:UPMUCASE2LOWERBOUND}, 
we can slightly modify our proof of \eqref{E:UPMUEASYLOWERBOUND}
(which was based in part on \eqref{E:MUAPPROXDEF})
and also use the assumption $s \leq t_1$ to deduce the following estimate:
\begin{align}  \label{E:UPMUCASE3LOWERBOUND}
	\upmu(s,u,\vartheta)
	\geq (1 - C \varepsilon)
		\left\lbrace
			1
			- \updelta \ln \left( \frac{\rgeo(s,u)}{\rgeo(0,u)} \right)
		\right\rbrace
	& \geq (1 - C \varepsilon)
		\left\lbrace
			1
			- \updelta \uptau_1
		\right\rbrace
			\\
	& = \frac{1}{2} (1 - C \varepsilon).
	\notag
\end{align}
Also using the bound \eqref{E:RADUPMUEASYBOUND}, we see that
$\frac{[\Rad \upmu(s,u,\vartheta)]_+}{\upmu(s,u,\vartheta)}$
is $\leq$ the second term on the right-hand side
of \eqref{E:UNIFORMBOUNDFORMRADMUOVERMU} as desired,
as in Cases 1 and 2.
\ \\

\noindent \emph{Case 4: The most difficult case involving dangerous decay rates and large times:}
In this case, we assume that $\updelta := - \widetilde{\Omega}_{(Min)}(t,U_0;t) > 0,$
that $(s,u,\vartheta) \in \mathcal{U}(s,U_0;t),$ 
and that $t_1 \leq s \leq t,$ where $t_1$ is defined in \eqref{E:SUBCASE3T1}.
In particular, in terms of the logarithmic change of variable
$\upsigma := \ln \left( \frac{\rgeo(s,u)}{\rgeo(0,u)} \right),$ 
the following estimate holds for $t_1 \leq s \leq t:$
\begin{align} \label{E:FINALHARDSUBCASEDELTAINEQUALITY}
	2 \updelta 
	& \geq \frac{1}{\upsigma}
		\geq \frac{1}{\upsigma + 1}.
\end{align}
Furthermore, with $\uptau := \ln \left( \frac{\rgeo(t,u)}{\rgeo(0,u)} \right) ,$
then as in \eqref{E:TNOTTOOBIG}, the following inequality holds:
\begin{align} \label{E:TNOTTOOBIGAGAIN}
	 \updelta \uptau < 1.
\end{align}

We may assume that $\Rad \upmu(s,u,\vartheta) > 0,$
since otherwise the estimate \eqref{E:UNIFORMBOUNDFORMRADMUOVERMU} is trivial.
To prove \eqref{E:UNIFORMBOUNDFORMRADMUOVERMU} under the assumption $\Rad \upmu(s,u,\vartheta) > 0,$ 
we study the integral curves of $\Rad.$ Thus, let 
\begin{align}
	\upgamma: [0,U_0] 
	& \rightarrow \Sigma_s^{U_0}
\end{align}	
denote the integral curve of $\Rad$ passing through
the point $p :=(s,u,\vartheta).$ We parametrize $\upgamma$
via the eikonal function $u',$ that is, $\upgamma = \upgamma(u'),$
where $u' \in [0,U_0],$ and $u$ is the value of $u'$ that corresponds
to the point $p$ of interest. We use notation such as
\begin{align}
	\dot{\upgamma}(u') & := \frac{d}{d u'} \upgamma(u'),
		\qquad
	\ddot{\upgamma}(u') := \frac{d^2}{d (u')^2} \upgamma(u').
\end{align}
We also set
\begin{align}
	F(u') & := \upmu \circ \upgamma(u'),
\end{align}
which implies that 
\begin{align}
	\dot{F}(u') 
	& = \Rad \upmu|_{\upgamma(u')},
		\\
	\ddot{F}(u') 
	& = \Rad \Rad \upmu|_{\upgamma(u')}.
		\label{E:SECONDERIVATIVENOTATION}
\end{align}	
Note that by assumption, we have
\begin{align} \label{E:DERIVATIVEATUISPOSITIVE}
	\dot{F}(u'=u) > 0.
\end{align}

We now set
\begin{align}
	u_* 
	& := \inf 
				\left\lbrace
					u' \in [0,u] \ : \ \upgamma(u') \in \mathcal{U}(s,U_0;t)
				\right\rbrace. 
\end{align}
Note that $\upgamma(u) \in \mathcal{U}(s,U_0;t)$ by assumption,
while $\upgamma(0) \notin \mathcal{U}(s,U_0;t)$ because
$\widetilde{\Omega}(t,0,\vartheta;t) \equiv 0$
(recall that $u'=0$ corresponds to a point on the null cone $\mathcal{C}_0,$ where the solution is 
trivial). Also using the fact that $\mathcal{U}(s,U_0;t)$ is open, we deduce that
\begin{align}
	0 & < u_* < u \leq U_0,
		\\
	\upgamma(u_*) & \in \Sigma_s^{U_0} \backslash \mathcal{U}(s,U_0;t).
	\label{E:POINTINGOODSET}
\end{align}

We consider two sub-cases. In the first sub-case, we assume that
\begin{align} \label{E:FINCREASING}
		\dot{F}(u') & > 0 \ \mbox{for} \ u' \in [u_*, u].
\end{align}
In this first sub-case, it follows from \eqref{E:FINCREASING} that
\begin{align} \label{E:FORUISBIGGER}
	F(u_*) < F(u).
\end{align}
Furthermore, the identity \eqref{E:POINTINGOODSET} and the argument
used in Case 2 imply that the bound
\eqref{E:UPMUCASE2LOWERBOUND} holds for
$F(u_*) = \upmu \circ \upgamma(u_*).$ Also using the identity \eqref{E:FORUISBIGGER},
we deduce that
\begin{align}
	\upmu(s,u,\vartheta)
	= F(u) 
	> F(u_*)
	\geq \frac{1}{2} (1 - C \varepsilon).
\end{align}
Also using the bound \eqref{E:RADUPMUEASYBOUND}, we see that
$\frac{[\Rad \upmu(s,u,\vartheta)]_+}{\upmu(s,u,\vartheta)}$
is $\leq$ the second term on the right-hand side
of \eqref{E:UNIFORMBOUNDFORMRADMUOVERMU} as desired,
as in the other cases we have previously considered.

In the final (most difficult) sub-case, we assume that the statement \eqref{E:FINCREASING} is false. Hence, 
there exists a point $u_{**} \in (u_*, u]$ such that
\begin{align}
	\dot{F}(u_{**}) & = 0,
	\label{E:DERIVATIVEATUSTARSTARIS0}
		\\
	\dot{F}(u') & > 0 \ \mbox{for} \ u' \in (u_{**}, u].
		\label{E:DERIVATIVEABOVEUSTARSTARISPOSITIVE}
\end{align}
In this case, our analysis involves the following function $H(\cdot)$ of a single real variable:
\begin{align} \label{E:KEYHESSIANDEF}
	H(s) & := \sup_{\Sigma_s^{U_0}} \Rad \Rad \upmu.
\end{align}

By \eqref{E:C0BOUNDCRUCIALEIKONALFUNCTIONQUANTITIES}, we have the following bound:
\begin{align} \label{E:RADRADMUUPPERBOUND}
	H(s) & \leq C \varepsilon \ln(\myexp + s).
\end{align}
By carefully studying the integral curves of $\Rad,$ we will
prove that the following two inequalities also hold in this case:
\begin{align} \label{E:HESSIANPOSITIVE}
	H(s) & > 0,
\end{align}
\begin{align} \label{E:RADMUOVERMUKEYRATIO3}
	\frac{[\Rad \upmu(s,u,\vartheta)]_+}{\upmu(s,u,\vartheta)}
	& \leq \sqrt{\frac{H(s)}{\upmu_{(Min)}(s,u)}},
\end{align}
where $\upmu_{(Min)}(s,u)$ is defined in \eqref{E:MUMINDEF}. Let us accept \eqref{E:HESSIANPOSITIVE} and \eqref{E:RADMUOVERMUKEYRATIO3}
for the moment; we show that these inequalities hold at the end of the proof.

In order to derive an upper bound for the right-hand side of \eqref{E:RADMUOVERMUKEYRATIO3}, 
we now derive a lower bound for $\upmu_{(Min)}(s,u).$ We use
the logarithmic change of variables 
$\upsigma := \ln \left(\frac{\rgeo(s,u)}{\rgeo(0,u)} \right),$ 
$\uptau := \ln \left(\frac{\rgeo(t,u)}{\rgeo(0,u)} \right),$ 
and
$\uptau_1 := \ln \left(\frac{\rgeo(t_1,u)}{\rgeo(0,u)} \right) = \frac{1}{2 \updelta}.$ 
Using essentially
the same reasoning that led to the first line of \eqref{E:MULARGETIMESECONDLOWERBOUND}
together with the estimates 
\eqref{E:MUAMPLITUDENEARONE}
and \eqref{E:TNOTTOOBIGAGAIN},
we deduce that
\begin{align} \label{E:MUMINOFSFIRSTLOWERBOUNDANNOYINGCASE}
	\upmu_{(Min)}(s,u)
	& \geq (1 - C \varepsilon)
		\left\lbrace
			1 - \updelta \upsigma - C \varepsilon \frac{1 + \upsigma}{\exp(\upsigma)}(\uptau - \upsigma)
		\right\rbrace
		\\
	& \geq (1 - C \varepsilon)
		\left\lbrace
			\updelta (\uptau - \upsigma) - C \varepsilon \frac{1 + \upsigma}{\exp(\upsigma)}(\uptau - \upsigma)
		\right\rbrace
		\notag \\
	& \geq \frac{1}{2} \updelta
		\left\lbrace
			1 - C \frac{\varepsilon}{\updelta} \frac{1 + \upsigma}{\exp(\upsigma)}
		\right\rbrace
		(\uptau - \upsigma).
		\notag
\end{align}

We now estimate the crucially important term $C \frac{\varepsilon}{\updelta} \frac{1 + \upsigma}{\exp(\upsigma)}$
on the right-hand side of \eqref{E:MUMINOFSFIRSTLOWERBOUNDANNOYINGCASE}. Since $\upsigma \geq \uptau_1,$ we have
\begin{align} \label{E:TMUOVERMUANNOYINGRATIO}
	& C \frac{\varepsilon}{\updelta}
			\frac{1 + \upsigma}{\exp(\upsigma)}
		\leq  
		C \frac{\varepsilon}{\updelta}
		\frac{1+ \uptau_1}{\exp(\uptau_1)}
	= 2 C \varepsilon 
		\frac{1}{2 \updelta}
		\left\lbrace
	    1 +	
			\frac{1}{2 \updelta}
		\right\rbrace
		\exp \left(- \frac{1}{2 \updelta} \right).
\end{align}
Since 
\begin{align}
	\lim_{y \to \infty} y (1+y) \exp(-y) = 0,
\end{align}
it follows from 
\eqref{E:MUMINOFSFIRSTLOWERBOUNDANNOYINGCASE}
(where we view $y = \frac{1}{2 \updelta}$),
\eqref{E:TMUOVERMUANNOYINGRATIO}, 
the estimate \eqref{E:FINALHARDSUBCASEDELTAINEQUALITY},
and the simple bound \eqref{E:DELTATRIVIALBLOUND}
that if $\varepsilon$ is sufficiently small, then
the quantity in braces on the right-hand side of 
\eqref{E:MUMINOFSFIRSTLOWERBOUNDANNOYINGCASE} is $\geq 1/2,$ and hence
\begin{align} \label{E:UPMUMINANNOYINGCASEFINALLOWERBOUND}
	\upmu_{(Min)}(s,u) 
	& \geq \frac{1}{8} \frac{(\uptau - \upsigma)}{1 + \upsigma}.
\end{align}

Inserting the bounds 
\eqref{E:RADRADMUUPPERBOUND}
and \eqref{E:UPMUMINANNOYINGCASEFINALLOWERBOUND}
into the right-hand side of inequality \eqref{E:RADMUOVERMUKEYRATIO3}, 
we deduce that
\begin{align} \label{E:TUMUOVERMUPOSITIVEPARTANNOYINGCASEFINALLOWERBOUND}
	\frac{[\Rad \upmu(s,u,\vartheta)]_+}{\upmu(s,u,\vartheta)}
	& \leq C \varepsilon^{1/2} \frac{(1 + \upsigma)}{\sqrt{\uptau - \upsigma}}.
\end{align}
It is straightforward to check using the basic properties of the function $f(x) = \ln(x)$
that the right-hand side of \eqref{E:TUMUOVERMUPOSITIVEPARTANNOYINGCASEFINALLOWERBOUND}
is $\leq$ the first term on the right-hand side of \eqref{E:UNIFORMBOUNDFORMRADMUOVERMU}
as desired. Hence, in order to conclude the desired estimate \eqref{E:UNIFORMBOUNDFORMRADMUOVERMU} in this final sub-case,
it remains only for us to prove inequalities \eqref{E:HESSIANPOSITIVE} and \eqref{E:RADMUOVERMUKEYRATIO3}.

To prove \eqref{E:HESSIANPOSITIVE},
we first use \eqref{E:DERIVATIVEATUISPOSITIVE} and \eqref{E:DERIVATIVEATUSTARSTARIS0}
to deduce that
\begin{align}	\label{E:SECONDRADDERIVATIVEINTEGRATEDPOSITIVITY}
	\int_{u'=u_{**}}^{u'= u} \ddot{F}(u') \, du'
	& = \dot{F}(u) > 0. 
\end{align}
Hence, it follows from
the definition \eqref{E:KEYHESSIANDEF} of $H(s),$
\eqref{E:SECONDERIVATIVENOTATION},
and
\eqref{E:SECONDRADDERIVATIVEINTEGRATEDPOSITIVITY}
that $H(s) > 0$ as desired.

We now derive inequality \eqref{E:RADMUOVERMUKEYRATIO3}.
To this end, we
first use the fundamental theorem of calculus to deduce that
if $\hat{u} \in [u_{**},u],$ then
\begin{align} \label{E:FUNDTHMCALCINEQUALITY}
	\dot{F}(u) - \dot{F}(\hat{u})
	& = \int_{u' = \hat{u}}^{u' = u}
		\ddot{F}(u') \, du'
	\leq H(s) (u - \hat{u}).
\end{align}
From inequality \eqref{E:FUNDTHMCALCINEQUALITY}, it follows that whenever
$\hat{u} \in [u_{**},u]$ and
\begin{align} \label{E:KEYDOMAIN}
	|\hat{u} - u|
	& \leq \frac{1}{2} \frac{\dot{F}(u)}{H(s)},
\end{align}
we have
\begin{align} \label{E:FIRSTDERIVATVELARGE}
	\dot{F}(\hat{u}) 
	& \geq \frac{1}{2} \dot{F}(u) > 0.
\end{align}
In particular, by choosing $\hat{u} = u_1 := u - \frac{1}{2} \frac{\dot{F}(u)}{H(s)}$
and using the lower bound \eqref{E:FIRSTDERIVATVELARGE},
we deduce that
\begin{align} \label{E:KEYHARDMUCASEBOUND}
	F(u) - F(u_1) 
	= \int_{u'=u_1}^{u' = u} \dot{F}(u') \, du'
	& \geq \frac{1}{2} \dot{F}(u) (u - u_1)
		\\
	& = \frac{1}{4} \frac{\dot{F}^2(u)}{H(s)}.
	\notag
\end{align}
Moreover, since $\dot{F}(u_{**}) = 0,$ the bound \eqref{E:FIRSTDERIVATVELARGE} 
for the domain of $\hat{u}$ values defined by \eqref{E:KEYDOMAIN} 
implies that $u_{**}$ cannot belong to the domain, 
that is, that
\begin{align} \label{E:USTARSTARSMALLERTHANU1}
	u_{**} < u_1.
\end{align}
Since \eqref{E:DERIVATIVEABOVEUSTARSTARISPOSITIVE} implies that $\dot{F}(u') > 0$ for $u' \in (u_{**},u_1],$ 
we also have
\begin{align}
	F(u_{**}) < F(u_1).
\end{align}
Hence, using the trivial bound $F(u_1) \geq \upmu_{(Min)}(s,u)$
(which follows from the definitions of $F(u_1)$ and $\upmu_{(Min)}(s,u)$)
and the estimate \eqref{E:KEYHARDMUCASEBOUND}, we deduce that
\begin{align} \label{E:RADMUOVERMUKEYRATIO1}
	\upmu(s,u,\vartheta) - \upmu_{(Min)}(s,u)
		= F(u) - \upmu_{(Min)}(s,u)
	& \geq F(u) - F(u_1)
		\\
	& \geq \frac{1}{4} \frac{\dot{F}^2(u)}{H(s)}
		\geq \frac{1}{4} \frac{[\Rad \upmu(s,u,\vartheta)]_+^2}{H(s)}.
		\notag
\end{align}
Solving for $\upmu(s,u,\vartheta)$ in inequality \eqref{E:RADMUOVERMUKEYRATIO1}
and then using the resulting bound to deduce a bound for
$\frac{[\Rad \upmu(s,u,\vartheta)]_+}{\upmu(s,u,\vartheta)},$
we find that
\begin{align} \label{E:RADMUOVERMUKEYRATIO2}
	\frac{[\Rad \upmu(s,u,\vartheta)]_+}{\upmu(s,u,\vartheta)}
	& \leq \frac{[\Rad \upmu(s,u,\vartheta)]_+}{\frac{1}{4} \frac{[\Rad \upmu(s,u,\vartheta)]_+^2}{H(s)} + \upmu_{(Min)}(s,u)}.
\end{align}

We now deduce an ``algebraic'' bound from inequality \eqref{E:RADMUOVERMUKEYRATIO2}.
To this end, we note that a standard calculus exercise
yields that for any constant $a > 0,$ 
the function
\begin{align}
	h(x,y) := \frac{x}{a^{-1} x^2 + y}
\end{align}
on the domain $[0,\infty) \times [0,\infty)$ 
is bounded by
\begin{align} \label{E:LITTLEFBOUND}
	|h(x,y)| & \leq \frac{\sqrt{a}}{2 \sqrt{y}}.
\end{align}

Applying inequality \eqref{E:LITTLEFBOUND} 
to the right-hand side of \eqref{E:RADMUOVERMUKEYRATIO2}
with $a= 4 H(s),$
$x= [\Rad \upmu(s,u,\vartheta)]_+,$ and $y = \upmu_{(Min)}(s,u),$ we arrive at
the desired inequality \eqref{E:RADMUOVERMUKEYRATIO3}.
This completes our proof of \eqref{E:UNIFORMBOUNDFORMRADMUOVERMU} in this final sub-case.

\ \\

\noindent \textbf{Proof of \eqref{E:POSITIVEPARTOFLMUPLUSRADMUOVERMUISSMALL}:}
From the identity $\uLgood \upmu = \upmu \Lunit + 2 \Rad$
and the estimates
\eqref{E:MUWEAKUPPERBOUND},
\eqref{E:ALTERNATEPOSITIVEPARTOFLMUOVERMUISSMALL},
and \eqref{E:UNIFORMBOUNDFORMRADMUOVERMU},
we deduce the desired estimate as follows:
\begin{align}
	\frac{[\Lunit \upmu(s,u,\vartheta) + \uLgood \upmu(s,u,\vartheta)]_+}{\upmu(s,u,\vartheta)}
	& \leq |\Lunit \upmu(s,u,\vartheta)|
		+ \frac{[\Lunit \upmu(s,u,\vartheta)]_+}{\upmu(s,u,\vartheta)}
		+ 2 \frac{[\Rad \upmu(s,u,\vartheta)]_+}{\upmu(s,u,\vartheta)}
		\\
	& \leq 	
			C \varepsilon^{1/2} \frac{\ln(\myexp + s)}{\sqrt{\ln(\myexp + t) - \ln(\myexp + s)}}
			+ C \varepsilon \frac{1}{1 + s}.
			\notag
\end{align}

\end{proof}

The following simple corollary plays a role in our analysis of the transport inequalities of 
Sect.~\ref{S:INVERTINGFULLYRENORMALIZEDTRANSPORT}.

\begin{corollary}[\textbf{Analysis of} $\frac{1}{2} \mytr \upchi - 2 \frac{\Lunit \upmu}{\upmu}$]
	\label{C:RENORMALIZATIONCOEFFICIENTISPOSITIVE}
	Under the small-data and bootstrap assumptions 
	of Sects.~\ref{S:PSISOLVES}-\ref{S:C0BOUNDBOOTSTRAP},
	if $\varepsilon$ is sufficiently small, 
	then the following pointwise estimate holds on $\mathcal{M}_{\Tboot,U_0}:$
	\begin{align}  \label{E:RENORMALIZATIONCOEFFICIENTISPOSITIVE}
		(1 - C \varepsilon) \rgeo^{-1}
		& \leq \frac{1}{2} \mytr \upchi - 2 \frac{\Lunit \upmu}{\upmu}
		\leq \frac{1}{2} \mytr \upchi + 2 \frac{[\Lunit \upmu]_-}{\upmu}.
	\end{align}
\end{corollary}

\begin{proof}
	From the decomposition $\mytr \upchi = 2 \rgeo^{-1} + \mytr \upchi^{(Small)}$
	and inequality \eqref{E:C0BOUNDCRUCIALEIKONALFUNCTIONQUANTITIES},
	we have that $\frac{1}{2} \mytr \upchi \geq (1 - C \varepsilon) \rgeo^{-1}.$
	Furthermore, by \eqref{E:ALTERNATEPOSITIVEPARTOFLMUOVERMUISSMALL}, 
	we have $- 2 \frac{\Lunit \upmu}{\upmu} \geq - C \varepsilon \rgeo^{-1}.$
	Adding these two estimates, we conclude
	\eqref{E:RENORMALIZATIONCOEFFICIENTISPOSITIVE}.
\end{proof}

\section{Fundamental estimates for time integrals involving \texorpdfstring{$\upmu^{-1}$}{the foliation density}}
We now use the estimates of Prop.~\ref{P:SHARPMU} to derive estimates for 
time integrals that involve powers of $\upmu^{-1}.$ These estimates
play a crucial role in our derivation of a priori $L^2$ estimates for $\Psi.$
Specifically, the time integrals appear in Chapter~\ref{C:ERRORTERMSOBOLEV}, 
when we bound the error integrals
of Prop.~\ref{P:DIVTHMWITHCANCELLATIONS}.

\begin{proposition}[\textbf{Fundamental estimates for time integrals involving $\upmu^{-1}$}] \label{P:MUINVERSEINTEGRALESTIMATES}
	Let $\upmu_{\star}(t,u)$ be as defined in \eqref{E:MUSTARDEF}.
	Let	
	\begin{align*}
		\Contwo > 1
	\end{align*}
	be a real number. 
	Under the small-data and bootstrap assumptions 
	of Sects.~\ref{S:PSISOLVES}-\ref{S:C0BOUNDBOOTSTRAP},
	if $\varepsilon$ is sufficiently small, then the following
	estimates hold for $0 \leq t < \Tboot$ and $0 \leq u \leq U_0.$
	
	\medskip
	
	\noindent \underline{\textbf{Estimates relevant for Gronwall estimates for borderline top-order spacetime integrals.}}
	There exists a constant $C > 0$ 
	such that if $\Contwo \sqrt{\varepsilon} |\ln \varepsilon| \leq 1,$ then
	\begin{align} \label{E:KEYMUTOAPOWERINTEGRALBOUND}
		\int_{s=0}^t 
			\frac{\left\| [\Lunit \upmu]_- \right\|_{C^0(\Sigma_s^u)}} 
					 {\upmu_{\star}^{\Contwo}(s,u)}
		\, ds 
		& \leq \frac{1 + C \sqrt{\varepsilon}}{\Contwo-1} 
			\upmu_{\star}^{1-\Contwo}(t,u).
	\end{align}	
	In addition, there exists a constant $C > 0$ such that
	if $\sqrt{\varepsilon} < \Littlecontwo \leq 1/2$ is a constant, then
	\begin{align} \label{E:ANNLOYINGSQRTMUOVERMUINTEGRATEDBOUND}
		\frac{1}{\rgeo(t,u)}
		\int_{s=0}^t
			\left\|
				\left(
					\frac{\upmu(t,\cdot)}{\upmu}
				\right)^2
			\right\|_{C^0(\Sigma_s^u)}
		\, ds
		& \leq 1 + C \Littlecontwo
			+ C \frac{\ln^2(\myexp + t)}{(1 + t)^{\Littlecontwo}}.
	\end{align}
	
	
	
	\medskip
		
	\noindent \underline{\textbf{Estimates relevant for Gronwall estimates for borderline top-order hypersurface integrals.}}
	Let $\Sigmaplus{s}{t}{u} \subset \Sigma_s^u$ 
	and
	$\Sigmaminus{s}{t}{u} \subset \Sigma_s^u$ be as defined in
	\eqref{E:SIGMAPLUS}
	and
	\eqref{E:SIGMAMINUS}.
	There exists a constant $C > 0$ 
	such that if $\Contwo \sqrt{\varepsilon} |\ln \varepsilon| \leq 1$ and 
	if $\Sigmaminus{t}{t}{u}$ is non-empty, then we have
	\begin{subequations}
	\begin{align} \label{E:KEYMUDECAYHYPERSURFACEMUTOAPOWERINTEGRALBOUND}
		\| \rgeo \Lunit \upmu \|_{C^0(\Sigmaminus{t}{t}{u})} 
		\int_{s=0}^t 
			\frac{1} 
				{\rgeo(s,u) \upmu_{\star}^{\Contwo}(s,u)}
			\, ds 
		& \leq \frac{1 + C \sqrt{\varepsilon}}{\Contwo-1} 
			\upmu_{\star}^{1-\Contwo}(t,u).
	\end{align}
	In addition, if $\Conone \geq 0$ is a real number and
	$\Sigmaplus{t}{t}{u}$ is non-empty, we have
	\begin{align} 	\label{E:KEYMUNODECAYHYPERSURFACEMUTOAPOWERINTEGRALBOUND}
		\left\| 
			\frac{\Lunit \upmu}{\upmu}
		\right\|_{C^0(\Sigmaplus{t}{t}{u})}
		\int_{s=0}^t 
			\left\|
					\sqrt{
						\frac{\upmu(t,\cdot)} 
					       {\upmu}
					     }
			\right\|_{C^0(\Sigmaplus{s}{t}{u})}
			\ln^{\Conone}(\myexp + s)
		\, ds 
		& \leq \frac{1 + C \varepsilon}{\Conone+1/2}
			\ln^{\Conone}(\myexp + t).
	\end{align}	
	\end{subequations}
	
	\medskip
	
	\noindent \underline{\textbf{Estimates relevant for Gronwall estimates for less dangerous top-order spacetime integrals.}}
	There exists a constant $C > 0$ depending on 
	\textbf{upper} bounds for  $\Littleconone$ and $\Conone$ but \textbf{independent of} $\Contwo$
	such that if 
	$\Littleconone \geq \varepsilon^{1/4},$
	$\Conone \varepsilon \leq 1,$ 
	and $\Contwo \sqrt{\varepsilon} |\ln \varepsilon| \leq 1,$ then
	\begin{align} \label{E:LESSSINGULARTERMSMUINTEGRALBOUND}
		\int_{s=0}^t \frac{\ln^{\Conone}(\myexp + s)} 
									{(1 + s)^{1 + \Littleconone} \upmu_{\star}^{\Contwo}(s,u)}
									\, ds 
		& \leq C \left\lbrace 1 + \frac{1}{\Contwo-1} \right\rbrace \frac{1}{\Littleconone^{\Conone+1}} \upmu_{\star}^{1-\Contwo}(t,u).
	\end{align}
	In addition, there exists 
	a constant $C > 0$
	such that if $\Contwo \sqrt{\varepsilon} |\ln \varepsilon| \leq 1,$  
	then
	\begin{align} \label{E:LOGLOSSKEYMUINTEGRALBOUND}
		\int_{s=0}^t \frac{1} 
									{(1 + s) \upmu_{\star}^{\Contwo}(s,u)}
									\, ds 
		& \leq C \left\lbrace 1 + \frac{1}{\Contwo-1} \right\rbrace \ln(\myexp + t) \upmu_{\star}^{1-\Contwo}(t,u).
	\end{align} 
	
	\medskip
	
	\noindent \underline{\textbf{Estimates for integrals that lead to only $\ln \upmu_{\star}^{-1}$ degeneracy.}} 
	There exists a constant $C > 0$ such that
	\begin{align} \label{E:KEYMUINVERSEINTEGRALBOUND}
		\int_{s=0}^t 
			\frac{\left\| [\Lunit \upmu]_- \right\|_{C^0(\Sigma_s^u)}} 
					 {\upmu_{\star}(s,u)}
		\, ds 
		& \leq (1 + C \sqrt{\varepsilon}) \ln \upmu_{\star}^{-1}(t,u) + C \sqrt{\varepsilon}.
	\end{align}
	Furthermore, 
	there exists a constant $C > 0$ depending on 
	\textbf{upper} bounds for $\Littleconone$ and $\Conone$
	such that if 
	$\Littleconone \geq \varepsilon^{1/4}$
	and
	$\Conone \varepsilon \leq 1,$
	then
	\begin{align} \label{E:LESSSINGULARLOGUPMULOSSTERMSMUINTEGRALBOUND}
		\int_{s=0}^t \frac{\ln^{\Conone}(\myexp + s)} 
									{(1 + s)^{1 + \Littleconone} \upmu_{\star}(s,u)}
								\, ds 
		& \leq C \frac{1}{\Littleconone^{\Conone+1}} 
			\left\lbrace
				\ln \upmu_{\star}^{-1}(t,u) 
				+ 1
			\right\rbrace.
	\end{align}
	In addition, there exists a constant $C > 0$ such that
	\begin{align} \label{E:LOGLOSSMUINVERSEINTEGRALBOUND}
		\int_{s=0}^t 
			\frac{1}{(1 + s)}
			\frac{1}{\upmu_{\star}(s,u)}
		\, ds 
		& \leq  C \ln(\myexp + t) \left\lbrace \ln \upmu_{\star}^{-1}(t,u) + 1 \right\rbrace.
	\end{align}
	
	\medskip
	
	\noindent \underline{\textbf{Estimates for integrals that break the $\upmu_{\star}^{-1}$ degeneracy.}} 
	There exists a constant $C > 0$ depending on 
	\textbf{upper} bounds for $\Littleconone$ and $\Conone$
	such that if 
	$\Littleconone \geq \varepsilon^{1/4}$
	and
	$\Conone \varepsilon \leq 1,$
	then
	\begin{align} \label{E:LESSSINGULARTERMSMUTHREEFOURTHSINTEGRALBOUND}
		\int_{s=0}^t \frac{\ln^{\Conone}(\myexp + s)} 
									{(1 + s)^{1 + \Littleconone} \upmu_{\star}^{3/4}(s,u)}
									\, ds 
		& \leq \frac{C}{\Littleconone^{\Conone}}.
	\end{align}
	In addition, there exists a constant $C > 0$ such that
	\begin{align} \label{E:LOGLOSSLESSSINGULARTERMSMTHREEFOURTHSINTEGRALBOUND}
		\int_{s=0}^t \frac{1} 
									{(1 + s) \upmu_{\star}^{3/4}(s,u)}
									\, ds 
		& \leq C \ln(\myexp + t).
	\end{align}

\end{proposition}

\begin{proof}
	
\noindent \textbf{Proof of \eqref{E:KEYMUTOAPOWERINTEGRALBOUND} and \eqref{E:KEYMUINVERSEINTEGRALBOUND}:}
We prove only \eqref{E:KEYMUTOAPOWERINTEGRALBOUND}. 
Inequality \eqref{E:KEYMUINVERSEINTEGRALBOUND} can be proved 
by making simple modifications to our proof of \eqref{E:KEYMUTOAPOWERINTEGRALBOUND}.
To prove \eqref{E:KEYMUTOAPOWERINTEGRALBOUND}, we first recall that
\begin{align}
	[\Lunit \upmu(s,u,\vartheta)]_- = \frac{[\Omega(s,u,\vartheta)]_-}{\rgeo(s,u)}.
\end{align}
We consider separately the cases $\widetilde{\Omega}_{(Min)}(t,u;t) \geq 0$
and $\widetilde{\Omega}_{(Min)}(t,u;t) < 0.$ In the case  $\widetilde{\Omega}_{(Min)}(t,u;t) \geq 0,$
the estimates
\eqref{E:APOWEREASYCASEMUSTARTABOUNDS}
and
\eqref{E:EASYCASEOMEGAMINUSBOUNDS}
imply that
\begin{align} \label{E:EASYCASEINTEGRALINTERMEDIATEDESCENT}
	\int_{s = 0}^t \frac{\left \| [\Lunit \upmu]_- \right \|_{C^0(\Sigma_s^u)}}
											{\upmu_{\star}^{\Contwo}(s,u)} \, ds
	& \leq \int_{s = 0}^t \frac{\left \| \Omega_- \right \|_{C^0(\Sigma_s^u)}}
											{\rgeo(s,u) \upmu_{\star}^{\Contwo}(s,u)} \, ds
			\\
	& \leq C \varepsilon 
		\int_{s = 0}^t \frac{\ln(\myexp + s)}{(1 + s)^2} \, ds
	\leq C \varepsilon \upmu_{\star}^{1-\Contwo}(s,u).
	\notag
\end{align}
We have thus proved \eqref{E:KEYMUTOAPOWERINTEGRALBOUND} in this case.

In the remaining case, which is $\updelta := - \widetilde{\Omega}_{(Min)}(t,u;t) > 0,$
we split the integral under consideration into the two portions
$\int_{s=0}^{t_1} \cdots \, ds$ and $\int_{s=t_1}^t \cdots \, ds,$
where as in \eqref{E:SMALLTIMEDEF},
\begin{align*}
	t_1 
	& = (1 - u) \frac{\varepsilon}{\updelta^2} - (1 - u), 
			\\
	\uptau_1 
	& = \ln \left(\frac{\rgeo(t_1,u)}{\rgeo(0,u)} \right)
		= \ln \left( \frac{\varepsilon}{\updelta^2} \right).
\end{align*}
We recall that if $\varepsilon$ is sufficiently small,
then we have $t_1 > 0$ (see \eqref{E:T1ISPOSITIVE}).
To estimate the integral portions, we use the change of variables
\eqref{E:TAUCHOV}-\eqref{E:SIGMACHOV}. In particular,
$d \upsigma = \rgeo^{-1}(s,u) \, ds.$ 
We first estimate the 
integral portion $\int_{s=0}^{t_1} \cdots \, ds.$
Using the change of variables,
the short-time estimates \eqref{E:POWERASHORTTIMEHARDCASEMUSTARBOUNDS} and \eqref{E:SHORTTIMEHARDCASEOMEGAMINUSBOUND},
the estimate \eqref{E:DELTATRIVIALBLOUND},
and the estimate $\int_{\upsigma = 0}^{\infty} (1 + \upsigma) \exp(- \upsigma) \leq C,$
we bound the integral of interest as follows:
\begin{align} \label{E:SHORTTIMEINTEGRALINTERMEDIATEDESCENT}
		\int_{s=0}^{t_1}  \frac{\left \| [\Lunit \upmu]_- \right \|_{C^0(\Sigma_s^u)}}
											     {\upmu_{\star}^{\Contwo}(s,u)}
		\, d s
		& \leq \int_{s=0}^{t_1} \frac{\left \| \Omega_- \right \|_{C^0(\Sigma_s^u)}}
											           {\upmu_{\star}^{\Contwo}(s,u)}
			     \, \frac{d s}{\rgeo(s,u)}
					\\
		& \leq (1 + C \sqrt{\varepsilon})
						\int_{\upsigma = 0}^{\uptau_1} 
							\left\lbrace
								\updelta + C \varepsilon (1 + \upsigma) \exp(- \upsigma)
							\right\rbrace
							\, d \upsigma
					\notag		\\
		& \leq 
			(1 + C \sqrt{\varepsilon})
			\left\lbrace
				\uptau_1 \updelta + C \varepsilon 
			\right\rbrace
			\notag \\
		& \leq 
			(1 + C \sqrt{\varepsilon})
			\left\lbrace
				2 \updelta |\ln \updelta| 
				+ \updelta |\ln \varepsilon|
			+ C \varepsilon
			\right\rbrace
			\notag \\
		& \leq C \varepsilon |\ln \varepsilon| 
			\notag \\
		& \leq 
			 C \frac{1}{\Contwo-1} 
			(\Contwo-1) \varepsilon |\ln \varepsilon|
			\upmu_{\star}^{1-\Contwo}(t,u).
			\notag 
\end{align}
Using our upper bound assumptions on the size of $\Contwo$ in terms of $\varepsilon,$
it is straightforward to see that the 
right-hand side of \eqref{E:SHORTTIMEINTEGRALINTERMEDIATEDESCENT} is
$\leq$ the right-hand side of
\eqref{E:KEYMUTOAPOWERINTEGRALBOUND} as desired.

In the final step, we bound the integral portion
$\int_{s=t_1}^t \cdots \, ds$ by $\leq \frac{1 + C \sqrt{\varepsilon}}{\Contwo-1} \upmu_{\star}^{1-\Contwo}(t,u).$ 
To this end, we use the large-time estimates \eqref{E:POWERALARGETIMEHARDCASEMUSTARBOUNDS} and \eqref{E:LARGETIMEHARDCASEOMEGAMINUSBOUND}
and the change of variables \eqref{E:TAUCHOV} and \eqref{E:SIGMACHOV}
to bound the integral of interest as follows:
\begin{align} \label{E:LARGETIMEINTEGRALINTERMEDIATEDESCENT}
		\int_{s=t_1}^{t}  \frac{\left \| [\Lunit \upmu]_- \right \|_{C^0(\Sigma_s^u)}}
											     {\upmu_{\star}^{\Contwo}(s,u)}
		\, d s
		& \leq 
		\int_{s=t_1}^t \frac{\left \| \Omega_- \right \|_{C^0(\Sigma_s^u)}}
											  {\upmu_{\star}^{\Contwo}(s,u)}
		\, \frac{d s}{\rgeo(s,u)}
			\\
		& \leq (1 + C \sqrt{\varepsilon}) 
					 \updelta	
					\int_{\upsigma = \uptau_1}^{\uptau} 
						\frac{1}{(1 - \updelta \upsigma)^{\Contwo}}
					\, d \upsigma
				\notag \\
		& \leq \frac{1}{\Contwo-1} 
				(1 + C \sqrt{\varepsilon})
				\frac{1}{(1 - \updelta \uptau)^{\Contwo-1}}
				\notag \\
		& \leq \frac{1}{\Contwo-1} 
			(1 + C \sqrt{\varepsilon})
			\upmu_{\star}^{1-\Contwo}(t,u).
		\notag
\end{align}
Clearly, the 
right-hand side of \eqref{E:LARGETIMEINTEGRALINTERMEDIATEDESCENT} is
$\leq$ the right-hand side of
\eqref{E:KEYMUTOAPOWERINTEGRALBOUND} as desired.
Combining \eqref{E:EASYCASEINTEGRALINTERMEDIATEDESCENT}
and the sum of
\eqref{E:SHORTTIMEINTEGRALINTERMEDIATEDESCENT}
and \eqref{E:LARGETIMEINTEGRALINTERMEDIATEDESCENT},
we conclude the desired estimate \eqref{E:KEYMUTOAPOWERINTEGRALBOUND}
in all cases.

\ \\

\noindent \textbf{Proof of \eqref{E:LESSSINGULARTERMSMUINTEGRALBOUND}, \eqref{E:LESSSINGULARLOGUPMULOSSTERMSMUINTEGRALBOUND},
and \eqref{E:LESSSINGULARTERMSMUTHREEFOURTHSINTEGRALBOUND}:}
We prove only \eqref{E:LESSSINGULARTERMSMUINTEGRALBOUND}; the estimates 
\eqref{E:LESSSINGULARLOGUPMULOSSTERMSMUINTEGRALBOUND}
and
\eqref{E:LESSSINGULARTERMSMUTHREEFOURTHSINTEGRALBOUND}
can be proved by making straightforward modifications to our proof of \eqref{E:LESSSINGULARTERMSMUINTEGRALBOUND}.
To prove \eqref{E:LESSSINGULARTERMSMUINTEGRALBOUND}, we modify our proof of \eqref{E:KEYMUTOAPOWERINTEGRALBOUND}.
We use the fact that the function 
\begin{align}
	f(\upsigma) := (1 + \upsigma)^{\Conone} \exp(-\Littleconone (1 + \upsigma))
\end{align}
defined on the domain $[-1,\infty)$
achieves its maximum at $\upsigma_{(Max)} := \frac{\Conone}{\Littleconone} - 1$
and that
\begin{align} \label{E:POINT}
	f'(\upsigma) & \geq 0, & & \upsigma \in [-1, \upsigma_{(Max)}),
		\\
	f(\upsigma_{(Max)}) & = \left(\frac{\Conone}{e}\right)^{\Conone} \frac{1}{\Littleconone^{\Conone}}, & &
		\\
	f'(\upsigma) & < 0, & & \upsigma \in (\upsigma_{(Max)}, \infty).
		\label{E:POINTDECREASING}
\end{align}
Furthermore, $f$ verifies the integral bound
\begin{align} \label{E:INTBOUND}
	\int_{\upsigma = 0}^{\infty} f(\upsigma) \, d \upsigma 
		\leq \frac{\Conone!}{\Littleconone^{\Conone+1}}.
\end{align}
We consider separately the cases $\widetilde{\Omega}_{(Min)}(t,u;t) \geq 0$
and $\widetilde{\Omega}_{(Min)}(t,u;t) < 0.$ In the case $\widetilde{\Omega}_{(Min)}(t,u;t) \geq 0,$
we use the estimates
\eqref{E:APOWEREASYCASEMUSTARTABOUNDS}
and
\eqref{E:INTBOUND}
and the change of variables 
$\upsigma := \ln \left( \frac{\rgeo(s,u)}{\rgeo(0,u)} \right)$
to deduce 
\begin{align}  \label{E:EASYCASEUGLYINTEGRAL}
	\int_{s=0}^t 
		\frac{\ln^{\Conone}(\myexp + s)} 
				 {(1 + s)^{1 + \Littleconone} \upmu_{\star}^{\Contwo}(s,u)}
	\, ds
	& \leq
		C
		\int_{s=0}^t 
			\frac{\left\lbrace 1 + \ln \left(\frac{\rgeo(s,u)}{\rgeo(0,u)} \right) \right\rbrace^{\Conone}} 
				   {\left\lbrace 1 + \frac{\rgeo(s,u)}{\rgeo(0,u)} \right\rbrace^{\Littleconone} \upmu_{\star}^{\Contwo}(s,u)}
		\, \frac{ds}{\rgeo(s,u)}
		\\
	& \leq 
		C
		\int_{\upsigma = 0}^{\infty} (1 + \upsigma)^{\Conone} \exp\left\lbrace - \Littleconone (1 + \upsigma) \right\rbrace \, d \upsigma
		\leq C \frac{1}{\Littleconone^{\Conone+1}}
			\notag \\
	& \leq C \frac{1}{\Littleconone^{\Conone+1}} \upmu_{\star}^{1-\Contwo}(t,u).
	\notag
\end{align}
We have thus proved the desired bound \eqref{E:LESSSINGULARTERMSMUINTEGRALBOUND} in this case.

In the remaining case, which is $\updelta := - \widetilde{\Omega}_{(Min)}(t,u;t) > 0,$
we split the time integration into the two portions 
$\int_{s=0}^{t_1} \cdots \, ds$ and 
$\int_{s=t_1}^t \cdots \, ds,$
where we now set $t_1 := (1 - u) \exp \left(\frac{1}{2 \Contwo \updelta} \right) - (1 - u),$
$\uptau_1 := \ln \left( \frac{\rgeo(t_1,u)}{\rgeo(0,u)} \right) = \frac{1}{2 \Contwo \updelta},$
$\uptau := \ln \left( \frac{\rgeo(t,u)}{\rgeo(0,u)} \right).$
We first bound the integral portion $\int_{s=0}^{t_1} \cdots \, ds.$
To this end, we first use the estimates
\eqref{E:POWERASHORTTIMEHARDCASEMUSTARBOUNDS}
and
\eqref{E:POWERALARGETIMEHARDCASEMUSTARBOUNDS},
the fact that $1- \updelta \ln \left( \frac{\rgeo(s,u)}{\rgeo(0,u)} \right) \geq 1 - \frac{1}{2\Contwo}$
for $0 \leq s \leq t_1,$
and the fact that there exists a uniform constant $C > 1$ such that
$C^{-1} \leq \lbrace 1 - \frac{1}{2 \Contwo} \rbrace^{\Contwo} \leq C$ 
for $\Contwo > 1$ 
to deduce that there exists a uniform constant $C > 1$ 
\emph{independent of} $\Contwo > 1$
such that the following estimates hold:
\begin{align}
	C^{-1} & \leq \upmu_{\star}^{\Contwo}(s,u) \leq C, && s \in [0, t_1],
		 \label{E:SHORTTIMEMUSTARBOUNDPROOF} \\ 
	C^{-1} & \leq \upmu_{\star}^{1-\Contwo}(s,u) \leq C, && s \in [0, t_1].
		\label{E:SHORTTIMESECONDMUSTARBOUNDPROOF}
\end{align}
Also using \eqref{E:INTBOUND}, we can estimate the integral portion of interest 
by using essentially the same argument that we used in deriving inequality 
\eqref{E:EASYCASEUGLYINTEGRAL}:
\begin{align} \label{E:SHORTTIMELESSDANGEROUS}
	\int_{s=0}^{t_1}
		\frac{\ln^{\Conone}(\myexp + s)} 
				{(1 + s)^{1 + \Littleconone} \upmu_{\star}^{\Contwo}(s,u)}
	\, ds
	& \leq C \frac{1}{\Littleconone^{\Conone + 1}} \upmu_{\star}^{1-\Contwo}(t,u).
\end{align}
Clearly, the right-hand side of \eqref{E:SHORTTIMELESSDANGEROUS} is $\leq$
the right-hand side of \eqref{E:LESSSINGULARTERMSMUINTEGRALBOUND} 
as desired.

We now bound the integral portion $\int_{s= t_1}^t \cdots \, ds.$
We first note that the estimate \eqref{E:DELTATRIVIALBLOUND}
and our assumptions on 
$\Conone,$
$\Contwo,$
and $\Littleconone$ imply that $\uptau_1 > 1 + \upsigma_{(Max)}$
and thus by \eqref{E:POINTDECREASING}, 
$f(\upsigma)$ is decreasing for $\upsigma \geq \uptau_1 - 1.$
Hence, using the estimate
\eqref{E:POWERALARGETIMEHARDCASEMUSTARBOUNDS},
the fact that $f(\upsigma) < f(\upsigma_{(Max)})$ over the integration range,
and the definition of $\uptau_1,$
we deduce that
\begin{align}
	\label{E:LARGETIMELESSDANGEROUS}
	\int_{s=t_1}^t 
		\frac{\ln^{\Conone}(\myexp + s)} 
				 {(1 + s)^{1 + \Littleconone} \upmu_{\star}^{\Contwo}(s,u)}
	\, ds
	& \leq C
		\int_{\upsigma= \uptau_1}^{\uptau} \frac{(1 + \upsigma)^{\Conone} \exp(-\Littleconone (1 + \upsigma))}{\upmu_{\star}^{\Contwo}} \, d \upsigma
		\\
	& \leq 
		C
		\uptau_1^{\Conone} \exp(-\Littleconone \uptau_1)
		\int_{\upsigma = \uptau_1}^{\uptau}
			\frac{1}{\left\lbrace 1 - \updelta \upsigma \right\rbrace^{\Contwo}}
		\, ds
		\notag \\
	& \leq 
		C
		\uptau_1^{\Conone} \exp(-\Littleconone \uptau_1) \frac{1}{\updelta} \frac{1}{(\Contwo-1)} 
		(1 - \updelta \uptau)^{1-\Contwo}
		\notag \\
	& \leq C \frac{\Contwo}{\Contwo-1} 
				\left(\frac{\Littleconone}{2 \Contwo \updelta} \right)^{\Conone+1} 
				\exp\left(- \frac{\Littleconone}{2 \Contwo \updelta} \right) 
				\frac{1}{\Littleconone^{\Conone+1}} \upmu_{\star}^{1-\Contwo}(t,u).
		\notag
\end{align}
Using the fact that the function $h(x):= x^{\Conone+1} \exp(-x)$ 
is $\leq \left\lbrace (A+1)/e \right\rbrace^{A+1}$ 
on the domain $x \in [0,\infty),$
we deduce (with $x := \frac{\Littleconone}{2 \Contwo \updelta}$)
that the product
$
\left(\frac{\Littleconone}{2 \Contwo \updelta} \right)^{\Conone+1} 
\exp\left(- \frac{\Littleconone}{2 \Contwo \updelta} \right)
$ 
on the right-hand side of
\eqref{E:LARGETIMELESSDANGEROUS} 
is 
$\leq \left\lbrace (A+1)/e \right\rbrace^{A+1}.$ 
The desired estimate \eqref{E:LESSSINGULARTERMSMUINTEGRALBOUND} thus follows in this case.

\ \\

\noindent \textbf{Proof of \eqref{E:LOGLOSSKEYMUINTEGRALBOUND}, \eqref{E:LOGLOSSMUINVERSEINTEGRALBOUND}, and 
\eqref{E:LOGLOSSLESSSINGULARTERMSMTHREEFOURTHSINTEGRALBOUND}:}
The proof of \eqref{E:LOGLOSSKEYMUINTEGRALBOUND} is similar to the proof of \eqref{E:KEYMUTOAPOWERINTEGRALBOUND} but requires a few simple changes. To prove \eqref{E:LOGLOSSKEYMUINTEGRALBOUND}, we consider separately the cases $\widetilde{\Omega}_{(Min)}(t,u;t) \geq 0$ and $\widetilde{\Omega}_{(Min)}(t,u;t) < 0.$
In the case  $\widetilde{\Omega}_{(Min)}(t,u;t) \geq 0,$
we use the estimate
\eqref{E:APOWEREASYCASEMUSTARTABOUNDS}
to deduce that
\[
\int_{s=0}^t \frac{1} 
									{(1 + s) \upmu_{\star}^{\Contwo}(s,u)}
\, ds 	
\leq 
C
\int_{s=0}^t \frac{1} 
									{1 + s} 
\, ds 	
\leq C \ln(\myexp + t) \upmu_{\star}^{1-\Contwo}(t,u)
\]
as desired.	

In the remaining case, which is $\updelta := - \widetilde{\Omega}_{(Min)}(t,u;t) > 0,$
we set $t_1 := (1 - u) \exp\left(\frac{1}{2 \Contwo \updelta} \right) - (1 - u),$ 
$\uptau_1 := \ln \left( \frac{\rgeo(t_1,u)}{\rgeo(0,u)} \right) = \frac{1}{2 \Contwo \updelta}.$
We use the change of variables \eqref{E:TAUCHOV} and \eqref{E:SIGMACHOV}.
We first consider the sub-case $t \leq t_1.$ In this sub-case, we use the 
estimates 
\eqref{E:SHORTTIMEMUSTARBOUNDPROOF}
and \eqref{E:SHORTTIMESECONDMUSTARBOUNDPROOF}
to bound the integral of interest as follows:
\begin{align} \label{E:SHORTTIMELOGLOSS}
	\int_{s=0}^t \frac{1} 
									  {(1 + s) \upmu_{\star}^{\Contwo}(s,u)}
	\, ds 	
	& \leq 
		C
		\int_{s=0}^t \frac{1} 
									{\upmu_{\star}^{\Contwo}(s,u)}
		\, \frac{ds}{\rgeo(s,u)} 	
		\\
& \leq 
	C
	\int_{\upsigma = 0}^{\uptau} 
		1
	\, d\upsigma 	
	\leq 	C 
		\uptau
		\upmu_{\star}^{1-\Contwo}(t,u)
		\notag \\
& \leq C \ln(\myexp + t) \upmu_{\star}^{1-\Contwo}(t,u).
\notag
\end{align}

We now consider the sub-case $\updelta := - \widetilde{\Omega}_{(Min)}(t,u;t) > 0,$ $t_1 < t.$
We split the integral under consideration into the two portions
$\int_{s=0}^{t_1} \cdots \, ds$ and $\int_{s=t_1}^t \cdots \, ds.$
We first estimate the integral portion $\int_{s=0}^{t_1} \cdots \, ds.$
We use the same argument used to prove \eqref{E:SHORTTIMELOGLOSS}
together with the approximate monotonicity inequality \eqref{E:MUSTARINVERSETOPOWERAMUSTGROWUPTOACONSTANT}
to bound the integral of interest as follows:
\begin{align}
	\int_{s=0}^{t_1} \frac{1} 
									{(1 + s) \upmu_{\star}^{\Contwo}(s,u)}
	\, ds 	
& \leq 	C \ln(\myexp + t_1) 
		\upmu_{\star}^{1-\Contwo}(t_1,u)
		\\
& \leq C \ln(\myexp + t) \upmu_{\star}^{1-\Contwo}(t,u).
\notag
\end{align}

It remains for us to bound the integral portion
$\int_{s=t_1}^t \cdots \, ds$ by $\leq C \left\lbrace 1 + \frac{1}{\Contwo-1} \right\rbrace \ln(\myexp + t) \upmu_{\star}^{1-\Contwo}(t,u).$ 
To this end, we use 
the short-time estimate
\eqref{E:POWERASHORTTIMEHARDCASEMUSTARBOUNDS},
large-time estimate \eqref{E:POWERALARGETIMEHARDCASEMUSTARBOUNDS},
the change of variables \eqref{E:TAUCHOV} and \eqref{E:SIGMACHOV},
and the inequalities
\begin{align}
	\frac{1}{(\Contwo-1) \updelta}
	= \frac{2\Contwo}{\Contwo-1} \uptau_1 
	\leq C \left\lbrace 1 + \frac{1}{\Contwo-1} \right\rbrace \uptau
	\leq C \left\lbrace 1 + \frac{1}{\Contwo-1} \right\rbrace \ln (\myexp + t)
\end{align}
to bound the integral of interest as follows:
\begin{align}
	\int_{s=t_1}^t \frac{1} 
									{(1 + s) \upmu_{\star}^{\Contwo}(s,u)}
	\, ds 	
	& \leq 
		C
		\int_{s=t_1}^t \frac{1} 
									      {\upmu_{\star}^{\Contwo}(s,u)}
		\, \frac{ds}{\rgeo(s,u)}
		\\
	& \leq C
			\int_{\upsigma = \uptau_1}^{\uptau} 
				\frac{1} 
			       {\left\lbrace 1 - \updelta \upsigma \right\rbrace^{\Contwo}} 
			\, d \upsigma 	
			\notag \\
   & \leq \frac{C}{(\Contwo-1)\updelta} \left\lbrace 1 - \updelta \uptau \right\rbrace^{1-\Contwo}
   	\notag \\
	 & \leq C \left\lbrace 1 + \frac{1}{\Contwo-1} \right\rbrace \ln(\myexp + t) \upmu_{\star}^{1-\Contwo}(t,u).
	 \notag
\end{align}
We have thus proved the desired estimate \eqref{E:LOGLOSSKEYMUINTEGRALBOUND} in all cases.

The estimate \eqref{E:LOGLOSSMUINVERSEINTEGRALBOUND}
can be proved in a similar fashion by taking
$t_1 := (1 - u) \exp\left(\frac{1}{2 \updelta} \right) - (1 - u).$

The estimate \eqref{E:LOGLOSSLESSSINGULARTERMSMTHREEFOURTHSINTEGRALBOUND} can
be proved in a similar fashion by taking 
$t_1 := (1 - u) \exp\left(\frac{1}{(3/2) \updelta} \right) - (1 - u).$ 

\ \\

\noindent \textbf{Proof of \eqref{E:ANNLOYINGSQRTMUOVERMUINTEGRATEDBOUND}:}
To prove \eqref{E:ANNLOYINGSQRTMUOVERMUINTEGRATEDBOUND}, we first consider the 
case $0 \leq t \leq 1.$ In this case, we use
the bound
$
\left\|
	\frac{\upmu(t,\cdot)}{\upmu}
\right\|_{C^0(\Sigma_s^u)}
\leq 1 + C \sqrt{\varepsilon}
$
implied by
Def.~\ref{D:REGIONSOFDISTINCTUPMUBEHAVIOR}
and the estimates
\eqref{E:LOCALIZEDMUCANTGROWTOOFAST} 
and
\eqref{E:LOCALIZEDMUMUSTSHRINK}
(with $s_1:=s$ and $s_2:=t$)
in order to deduce that
\begin{align} \label{E:EASYCASEEARLYTIMEANNLOYINGSQRTMUOVERMUINTEGRATEDBOUND}
\frac{1}{1 - u + t}
		\int_{s=0}^t
			\left\|
				\left(
					\frac{\upmu(t,\cdot)}{\upmu}
				\right)^2
			\right\|_{C^0(\Sigma_s^u)}
		\, ds
	& \leq (1 + C \sqrt{\varepsilon})\frac{t}{1 - u + t}
		\leq C \frac{\ln^2(\myexp + t)}{(1 + t)^{\Littlecontwo}}
\end{align}
as desired.

We now consider the case $t \geq 1.$ To proceed, 
we split the integral into the 
two portions $\int_{s=0}^{t^{1 - \Littlecontwo}} \cdots \, ds$ 
and
$\int_{s= t^{1 - \Littlecontwo}}^t \cdots \, ds.$
To bound the first portion, we use
the bound
$
\left\|
	\frac{\upmu(t,\cdot)}{\upmu}
\right\|_{C^0(\Sigma_s^u)}
\leq (1 + C \sqrt{\varepsilon}) \ln(\myexp + t)
$
implied by
Def.~\ref{D:REGIONSOFDISTINCTUPMUBEHAVIOR} 
and the estimates
\eqref{E:LOCALIZEDMUCANTGROWTOOFAST} 
and
\eqref{E:LOCALIZEDMUMUSTSHRINK}
(with $s_1:=s$ and $s_2:=t$)
in order to deduce that
\begin{align} \label{E:EARLYTIMEANNLOYINGSQRTMUOVERMUINTEGRATEDBOUND}
\frac{1}{1 - u + t}
		\int_{s=0}^{t^{1-\Littlecontwo}}
			\left\|
				\left(
					\frac{\upmu(t,\cdot)}{\upmu}
				\right)^2
			\right\|_{C^0(\Sigma_s^u)}
		\, ds
	& \leq C \frac{\ln^2(\myexp + t)}{1 + t} 
		(t^{1 - \Littlecontwo})
		\leq C \frac{\ln^2(\myexp + t)}{(1 + t)^{\Littlecontwo}}
\end{align}
as desired.

To bound the second portion, we use
the bound
$
\sup_{s \in [t^{1 - \Littlecontwo}, t]}
\left\|
	\frac{\upmu(t,\cdot)}{\upmu}
\right\|_{C^0(\Sigma_s^u)}
\leq 1 + C \sqrt{\varepsilon} + C \Littlecontwo
$
implied by
Def.~\ref{D:REGIONSOFDISTINCTUPMUBEHAVIOR} 
and the estimates
\eqref{E:LARGETIMELOCALIZEDMUNOTDECAYINGRATIOESTIMATE}
and
\eqref{E:LOCALIZEDMUMUSTSHRINK}
(with $\Littleconone:= 1 - \Littlecontwo,$ $s_1:=s,$ and $s_2:=t$)
in order to deduce that
\begin{align} \label{E:LATETIMEANNLOYINGSQRTMUOVERMUINTEGRATEDBOUND}
\frac{1}{1 - u + t}
		\int_{s= t^{1 - \Littlecontwo}}^{t}
			\left\|
				\left(
					\frac{\upmu(t,\cdot)}{\upmu}
				\right)^2
			\right\|_{C^0(\Sigma_s^u)}
		\, ds
	& \leq (1 + C \sqrt{\varepsilon} + C \Littlecontwo) \frac{t(1 - t^{- \Littlecontwo})}{1 - u + t}
		\\
	& \leq 1+ C \sqrt{\varepsilon} + C \Littlecontwo
		+ C \frac{1}{1 + t}.
		\notag 
\end{align}
We now observe that since $\sqrt{\varepsilon} \leq \Littlecontwo,$
the sum of \eqref{E:EARLYTIMEANNLOYINGSQRTMUOVERMUINTEGRATEDBOUND} and 
\eqref{E:LATETIMEANNLOYINGSQRTMUOVERMUINTEGRATEDBOUND}
is $\leq$ the right-hand side of \eqref{E:ANNLOYINGSQRTMUOVERMUINTEGRATEDBOUND} as desired.

\ \\

\noindent \textbf{Proof of \eqref{E:KEYMUNODECAYHYPERSURFACEMUTOAPOWERINTEGRALBOUND}:}
We start by considering inequality \eqref{E:KEYMUNOTDECAYINGLMUOVERMUBOUND} with $s = t.$
We now derive a suitable bound for the term on the left-hand side of
\eqref{E:KEYMUNODECAYHYPERSURFACEMUTOAPOWERINTEGRALBOUND}
generated by the first term on the right-hand side of
\eqref{E:KEYMUNOTDECAYINGLMUOVERMUBOUND}.
To this end, we use
\eqref{E:GAMMACONSTANTVERYSMALL},
\eqref{E:SHARPLOCALIZEDMUCANTGROWTOOFAST}
with $s_1 = s$ and $s_2 = t,$
the change of variables $\upsigma = \ln(\myexp + s),$ 
and the trivial bound $(\myexp + t)^{-1} \leq (\myexp+s)^{-1}$ for $0 \leq s \leq t$
to deduce that
\begin{align} \label{E:FIRSTKEYMUNOTDECAYHYPERSURFACETIMEINTEGRALLOGESTIMATE}
	& \leq
			(1 + C \varepsilon)
			\frac{\upgamma}{\sqrt{1 + \upgamma \ln(\myexp + t)}}
			\int_{s=0}^t
				\frac{\ln^{\Conone}(\myexp+s)}{(\myexp + s) \sqrt{1 + \upgamma \ln(\myexp + s)}}
			\, ds
			\\
		& =
			(1 + C \varepsilon)
			\frac{\upgamma}{\sqrt{1 + \upgamma \ln(\myexp + t)}}
			\int_{\upsigma = 1}^{\ln(\myexp + t)}
				\frac{\upsigma^{\Conone}}{\sqrt{1 + \upgamma \upsigma}}
			\, d \upsigma.
			\notag
\end{align}
Using the trivial bound $\upsigma^{\Conone} / \sqrt{1 + \upgamma \upsigma} \leq \upsigma^{\Conone-1/2} / \sqrt{\upgamma}$
to pointwise bound the integrand,
we deduce that the right-hand side of \eqref{E:FIRSTKEYMUNOTDECAYHYPERSURFACETIMEINTEGRALLOGESTIMATE} is
\begin{align} \label{E:FINALKEYMUNOTDECAYHYPERSURFACETIMEINTEGRALLOGESTIMATE}
	& \leq 
	(1 + C \varepsilon)
	\frac{\upgamma}{(\myexp + t)\sqrt{1 + \upgamma \ln(\myexp + t)}}
	\frac{1}
		{\sqrt{\upgamma}}
	\int_{\upsigma = 1}^{\ln(\myexp + t)}
		\upsigma^{\Conone-1/2}
	\, d \upsigma
	\\
	& \leq 
		(1 + C \varepsilon)
		\frac{\sqrt{\upgamma}}{\sqrt{1 + \upgamma \ln(\myexp + t)}}
		\frac{\ln^{\Conone+1/2}(\myexp + t)}{(\Conone+1/2)}.
			\notag
\end{align}

It now easily follows that the right-hand side of
\eqref{E:FINALKEYMUNOTDECAYHYPERSURFACETIMEINTEGRALLOGESTIMATE}
is $\leq$ the right-hand side of \eqref{E:KEYMUNODECAYHYPERSURFACEMUTOAPOWERINTEGRALBOUND} 
as desired.

It remains for us to address 
the term on the left-hand side of
\eqref{E:KEYMUNODECAYHYPERSURFACEMUTOAPOWERINTEGRALBOUND}
generated by the second term on the right-hand side of
\eqref{E:KEYMUNOTDECAYINGLMUOVERMUBOUND}
(where we set $s=t$ in \eqref{E:KEYMUNOTDECAYINGLMUOVERMUBOUND}).
Using \eqref{E:LOCALIZEDMUCANTGROWTOOFAST} 
with $s_1 = s$ and $s_2 = t,$
we deduce that the product under consideration is 
\begin{align} \label{E:SECONDTERMFIRSTKEYMUNOTDECAYHYPERSURFACETIMEINTEGRALLOGESTIMATE}
	& \leq
			C \varepsilon
			\frac{\ln(\myexp + t)}{(1 + t)^2}
			\int_{s=0}^t
				\ln^{\Conone + 1/2}(\myexp+s)
			\, ds
		\leq
			C \varepsilon
			\frac{\ln^{\Conone + 3/2}(\myexp+t)}{1 + t}.
\end{align}

Clearly, the right-hand side of 
\eqref{E:SECONDTERMFIRSTKEYMUNOTDECAYHYPERSURFACETIMEINTEGRALLOGESTIMATE} is
$\leq$ the right-hand side of \eqref{E:KEYMUNODECAYHYPERSURFACEMUTOAPOWERINTEGRALBOUND} 
as desired.

\ \\

\noindent \textbf{Proof of \eqref{E:KEYMUDECAYHYPERSURFACEMUTOAPOWERINTEGRALBOUND}:}
The proof of \eqref{E:KEYMUDECAYHYPERSURFACEMUTOAPOWERINTEGRALBOUND} is very similar to
that of \eqref{E:KEYMUTOAPOWERINTEGRALBOUND},
so we are somewhat terse concerning the details. We begin with the trivial estimate 
\begin{align} \label{E:TRIVIALBOUNDFORRGEOLUNITUMPU}
	\| \rgeo \Lunit \upmu \|_{C^0(\Sigmaminus{t}{t}{u})}
	\leq \| \rgeo [\Lunit \upmu]_- \|_{C^0(\Sigmaminus{t}{t}{u})}
		+ \| \rgeo [\Lunit \upmu]_+ \|_{C^0(\Sigmaminus{t}{t}{u})}.
\end{align}
We separately bound the terms corresponding to 
the two terms in \eqref{E:TRIVIALBOUNDFORRGEOLUNITUMPU}.
We first prove the easier bound involving $\| \rgeo [\Lunit \upmu]_+ \|_{C^0(\Sigmaminus{t}{t}{u})}.$
To this end, we use the estimates 
\eqref{E:LMUPLUSNEGLIGIBLEINSIGMAMINUS}
(with $s=t$)
and
\eqref{E:LOGLOSSKEYMUINTEGRALBOUND}
to deduce the bound
\begin{align} \label{E:POSITIVEPARTKEYMUDECAYHYPERSURFACEMUTOAPOWERINTEGRALBOUND}
		\| \rgeo [\Lunit \upmu]_+ \|_{C^0(\Sigmaminus{t}{t}{u})} 
		\int_{s=0}^t 
			\frac{1} 
				{\rgeo(s,u) \upmu_{\star}^{\Contwo}(s,u)}
			\, ds 
		& \leq 
			C \varepsilon \frac{\ln(\myexp + t)}{1 + t}
		\int_{s=0}^t 
			\frac{1} 
				{(\myexp + s) \upmu_{\star}^{\Contwo}(s,u)}
			\, ds 
			\\
	& \leq
		C \varepsilon \frac{1}{\Contwo-1} \frac{\ln^2(\myexp + t)}{1 + t} \upmu_{\star}^{1-\Contwo}(t,u).
		\notag
	\end{align}
Clearly, the right-hand side of \eqref{E:POSITIVEPARTKEYMUDECAYHYPERSURFACEMUTOAPOWERINTEGRALBOUND} 
is $\leq$ the right-hand side of \eqref{E:KEYMUDECAYHYPERSURFACEMUTOAPOWERINTEGRALBOUND} 
as desired.

To prove the bound corresponding to the second term $\| \rgeo [\Lunit \upmu]_- \|_{C^0(\Sigmaminus{t}{t}{u})}$
on the right-hand side of \eqref{E:POSITIVEPARTKEYMUDECAYHYPERSURFACEMUTOAPOWERINTEGRALBOUND},
we first note that since $\Sigmaminus{t}{t}{u}$ is non-empty,
the positivity assumption on $\updelta$ stated in \eqref{E:CRUCIALDELTADEF} holds.
We separately consider the cases $0 \leq t \leq t_1$ 
and $t \geq t_1,$ where $t_1 := (1 - u) \frac{\varepsilon}{\updelta^2} - (1 - u).$ 
In the first case $0 \leq t \leq t_1,$ 
we use the estimates
\eqref{E:POWERASHORTTIMEHARDCASEMUSTARBOUNDS},
\eqref{E:HYPERSURFACESHORTTIMEHARDCASEOMEGAMINUSBOUND},
and 
\eqref{E:DELTATRIVIALBLOUND}
and the bounds
$\int_{s=0}^t
				\frac{1}{\myexp + s}
			\, ds \leq \ln \left(\myexp + \frac{\varepsilon}{\updelta^2} \right),$
$\int_{s=0}^t
				\frac{1}{\myexp + s}
			\, ds \leq \ln(\myexp + t)$			
to deduce the desired estimate as follows:
\begin{align} \label{E:NEGATIVEPARTSHORTTIMEKEYMUDECAYHYPERSURFACEMUTOAPOWERINTEGRALBOUND}
		\| \rgeo [\Lunit \upmu]_- \|_{C^0(\Sigmaminus{t}{t}{u})} 
		\int_{s=0}^t 
			\frac{1} 
				{\rgeo(s,u) \upmu_{\star}^{\Contwo}(s,u)}
			\, ds 
		& \leq 
			 C 
			\left\lbrace
				\updelta 
				+ C \varepsilon \frac{\ln(\myexp + t)}{1 + t}
			\right\rbrace
			\int_{s=0}^t
				\frac{1}{\myexp + s}
			\, ds
		\\
	& \leq
			C
			\left\lbrace
				\updelta \ln \left(\myexp + \frac{\varepsilon}{\updelta^2} \right)
				+ C \varepsilon \frac{\ln^2(\myexp + t)}{1 + t}
			\right\rbrace
		\notag
		\\
	& \leq C \varepsilon |\ln \varepsilon|
		\leq C \frac{1}{\Contwo-1} (\Contwo-1) \varepsilon |\ln \varepsilon| 
		\upmu_{\star}^{1-\Contwo}(t,u)
		\notag \\
	& \leq C \sqrt{\varepsilon} \frac{1}{\Contwo-1} \upmu_{\star}^{1-\Contwo}(t,u).
		\notag
	\end{align}

In the remaining case $t \geq t_1,$ we use the estimates
\eqref{E:POWERASHORTTIMEHARDCASEMUSTARBOUNDS},
\eqref{E:POWERALARGETIMEHARDCASEMUSTARBOUNDS},
\eqref{E:HYPERSURFACESHORTTIMEHARDCASEOMEGAMINUSBOUND},
\eqref{E:HYPERSURFACELARGETIMEHARDCASEOMEGAMINUSBOUND},
and
\eqref{E:DELTATRIVIALBLOUND}
to deduce the desired bound as follows:
\begin{align} \label{E:NEGATIVEPARTLARGETIMEKEYMUDECAYHYPERSURFACEMUTOAPOWERINTEGRALBOUND}
		& \| \rgeo [\Lunit \upmu]_- \|_{C^0(\Sigmaminus{t}{t}{u})} 
		\int_{s=0}^t 
			\frac{1} 
				{\rgeo(s,u) \upmu_{\star}^{\Contwo}(s,u)}
			\, ds 
				\\
		& \leq 
			\left(1	
					+
					C \sqrt{\varepsilon} 
			\right)
			\updelta
			\int_{s=0}^{t_1}
				\frac{1}
				{\myexp + s}
			\, ds
			+
			\left(1	
					+
					C \sqrt{\varepsilon}
			\right)
			\updelta
			\int_{s=t_1}^t
				\frac{1}
				{\rgeo(s,u) \left\lbrace 1 - \updelta \ln \left( \frac{\rgeo(s,u)}{\rgeo(0,u)} \right) \right\rbrace^{\Contwo}}
			\, ds
			\notag 
			\\
		& \leq 
			C
			\varepsilon |\ln \varepsilon|
			+ 
			\frac{1	+ C \sqrt{\varepsilon}}
					 {\Contwo-1}
			\left\lbrace 1 - \updelta \ln \left( \frac{\rgeo(t,u)}{\rgeo(0,u)} \right) \right\rbrace^{1-\Contwo}
			\notag \\
		& \leq 
			\frac{1	+ C \sqrt{\varepsilon}}
					 {\Contwo-1}
			\upmu_{\star}^{1-\Contwo}(t,u).
			\notag
\end{align}
We have thus proved the desired estimate \eqref{E:KEYMUDECAYHYPERSURFACEMUTOAPOWERINTEGRALBOUND}.

\end{proof}


\chapter{\texorpdfstring{$L^2$}{Square Integral} Coerciveness and the Fundamental \texorpdfstring{$L^2-$}{Square Integral-}Controlling Quantities}
\label{C:L2COERCIVENESS}
\thispagestyle{fancy}
In Chapter~\ref{C:L2COERCIVENESS}, we reveal the coercive nature of the energies and fluxes
defined in Def.~\ref{D:ENERGIESANDFLUXES}.
We then define a related family of coercive $L^2-$based 
quantities that we later use to derive a priori $L^2$ estimates for 
$\Psi$ and its derivatives.

\section{Coerciveness of the energies and fluxes}
In this section, we show that the energy-flux quantities from Def.~\ref{D:ENERGIESANDFLUXES} are coercive.

\begin{lemma}[\textbf{Coerciveness of the energies and fluxes}]
		\label{L:COERCIVEENERGIESANDFLUXES}
		Under the 
		small data and
		bootstrap assumptions 
		of Sects.~\ref{S:PSISOLVES}-\ref{S:C0BOUNDBOOTSTRAP},
		if $\varepsilon$ is sufficiently small, 
		then the energies $\enzero[\Psi](t,u),$ $\enone[\Psi](t,u)$ and the cone fluxes 
		$\flzero[\Psi](t,u),$ $\flone[\Psi](t,u)$ from Def.~\ref{D:ENERGIESANDFLUXES}
		have the following coerciveness properties (see Remark~\ref{R:SPHEREAREAFORMIMPLICITFACTOR}):
		\begin{subequations}
		\begin{align} \label{E:MULTENERGYCOERCIVITY}
				\enzero[\Psi](t,u)
				& = \frac{1}{2} 
						\int_{\Sigma_t^u}
							\left\lbrace
							\upmu(1 + \upmu) (\Lunit \Psi)^2
							+ (\uLgood \Psi)^2
							+ \upmu (1 + 2\upmu) |\angdiff \Psi|^2
							\right\rbrace
						\, d \tvol
					\\
					& \geq 
					\max 
						\left\lbrace
							C^{-1} \| \Psi \|_{L^2(S_{t,u})}^2,
							C^{-1} \| \Psi \|_{L^2(\Sigma_t^u)}^2,
							\frac{1}{2} \| \sqrt{\upmu} \Lunit \Psi \|_{L^2(\Sigma_t^u)}^2,
							\frac{1}{2} \| \upmu \Lunit \Psi \|_{L^2(\Sigma_t^u)}^2,
							\| \Rad \Psi \|_{L^2(\Sigma_t^u)}^2
						\right\rbrace
						\notag	\\
					& \ \ + \frac{1}{2} \| \sqrt{\upmu} \angdiff \Psi \|_{L^2(\Sigma_t^u)}^2
							+  \| \upmu \angdiff \Psi \|_{L^2(\Sigma_t^u)}^2,
						\notag \\
				\flzero[\Psi](t,u)
					& = \int_{\mathcal{C}_u^t}
								\left\lbrace
								(1 + \upmu)(\Lunit \Psi)^2 
							+ \upmu |\angdiff \Psi|^2
							\right\rbrace
						\, d \conevol
						\label{E:MULTCONEFLUXCOERCIVITY} \\
					& =
						\| \Lunit \Psi \|_{L^2(\mathcal{C}_u^t)}^2
						+ \| \sqrt{\upmu} \Lunit \Psi \|_{L^2(\mathcal{C}_u^t)}^2
						+ \| \sqrt{\upmu} \angdiff \Psi \|_{L^2(\mathcal{C}_u^t)}^2,
						\notag
		\end{align}	
		\end{subequations}
		
		\begin{subequations}
		\begin{align}
				\enone[\Psi](t,u)
				& = \frac{1}{2}
					\int_{\Sigma_t^u}
						\left\lbrace
						\rgeo^2 
						\upmu
						\left(
							\Lunit \Psi 
							+ \frac{1}{2} \mytr \upchi \Psi 
						\right)^2 
						+ \rgeo^2
							\upmu |\angdiff \Psi|^2
						\right\rbrace
					\, d \tvol
					\label{E:ENONECOERCIVENESS} \\
				& \approx (1 + t)^2 
									\left\| \sqrt{\upmu} 
										\left(
											\Lunit \Psi 
											+ \frac{1}{2} \mytr \upchi \Psi 
										\right)
									\right\|_{L^2(\Sigma_t^u)}^2
							+ (1 + t)^2 
								\| \sqrt{\upmu} 
									\angdiff \Psi	
								\|_{L^2(\Sigma_t^u)}^2,
						\notag \\
				\flone[\Psi](t,u)
				& =
					\int_{\mathcal{C}_u^t}
						\rgeo^2
						\left(
							\Lunit \Psi 
							+ \frac{1}{2} \mytr \upchi \Psi 
						\right)^2 
					\, d \conevol
					\label{E:FLUXONECOERCIVENESS} \\
				& \approx 
					\left\| 
						(1 + t') 
						\left(
							\Lunit \Psi 
							+ \frac{1}{2} \mytr \upchi \Psi 
						\right)
					\right\|_{L^2(\mathcal{C}_u^t)}^2.
					\notag
		\end{align}	
		\end{subequations}
		
	\end{lemma}

\begin{proof}
To prove \eqref{E:ENONECOERCIVENESS}, we first use
Lemma~\ref{L:NULLCOMPONENTSOFENERGYMOMENTUMTENSOR}
and definition \eqref{E:MORPLUSMODCURRENT} to compute that
(recall that $\Timenormal = \frac{1}{2} \Lunit + \frac{1}{2} \uLunit = \Lunit + \upmu^{-1} \Rad$)
\begin{align} \label{E:ENONEJMORPLUSMODEXPANSION}
	\upmu {\Jenergycurrent{\Mor+Correction}_{\Timenormal}[\Psi]}
	& = \frac{1}{2}
			\upmu
			\rgeo^2
			\left\lbrace
				\Lunit \Psi
				+ \frac{1}{2} \mytr \upchi \Psi
			\right\rbrace^2
		+ \frac{1}{2} 
			\upmu
			\rgeo^2
			|\angdiff \Psi|^2
			\\
	& \ \ 
		+ \frac{1}{2} \rgeo^2 \mytr \upchi \Psi \Rad \Psi
		- \frac{1}{8} \upmu \rgeo^2 (\mytr \upchi)^2 \Psi^2
		- \frac{1}{4} \upmu \Psi^2 \Lunit [\rgeo^2 \mytr \upchi]
		- \frac{1}{4} \Psi^2 \Rad [\rgeo^2 \mytr \upchi].
		\notag
\end{align}
From \eqref{E:ENONEJMORPLUSMODEXPANSION}, it follows that the 
right-hand side of \eqref{E:E1DEF}
is equal to the integral lying to the right of the equal sign in 
\eqref{E:ENONECOERCIVENESS} as desired.
The final inequality the right-hand side of \eqref{E:ENONECOERCIVENESS}
follows trivially from the estimate $\rgeo(t',u) \approx 1 + t'.$

To prove \eqref{E:FLUXONECOERCIVENESS}, we first use Lemma~\ref{L:NULLCOMPONENTSOFENERGYMOMENTUMTENSOR}
and definition \eqref{E:MORPLUSMODCURRENT} to compute that 
\begin{align} \label{E:JMORPLUSMODEXPANSION}
	{\Jenergycurrent{\Mor+Correction}_{\Lunit}[\Psi]}
	& =
		\rgeo^2 \left\lbrace \Lunit \Psi + \frac{1}{2} \mytr \upchi \Psi \right\rbrace^2
		- \frac{1}{4} \rgeo^2 (\mytr \upchi)^2 \Psi^2
		- \frac{1}{2} \rgeo^2 \mytr \upchi \Psi \Lunit \Psi
		- \frac{1}{4} \Psi^2 \Lunit [\rgeo^2 \mytr \upchi]
			\\
	& =
		\rgeo^2 \left\lbrace \Lunit \Psi + \frac{1}{2} \mytr \upchi \Psi \right\rbrace^2
		- \frac{1}{4} \rgeo^2 (\mytr \upchi)^2 \Psi^2
		- \frac{1}{4} \Psi^2 \Lunit [\rgeo^2 \mytr \upchi \Psi^2].
			\notag
\end{align}
From \eqref{E:JMORPLUSMODEXPANSION}, it follows that the
right-hand side of \eqref{E:F1DEF} is equal to the integral 
lying to the right of the equal sign in \eqref{E:FLUXONECOERCIVENESS} as desired.
The final inequality the right-hand side of \eqref{E:FLUXONECOERCIVENESS}
follows trivially from the estimate $\rgeo(t',u) \approx 1 + t'.$

With two exceptions, the proofs of estimates \eqref{E:MULTENERGYCOERCIVITY} 
and \eqref{E:MULTCONEFLUXCOERCIVITY} are similar, and we omit the details.
The two exceptions are the estimates for the terms 
$C^{-1} \| \Psi \|_{L^2(S_{t,u})}^2$ and $C^{-1} \| \Psi \|_{L^2(\Sigma_t^u)}^2$
in inequality \eqref{E:MULTENERGYCOERCIVITY}. The latter follows from the former by integrating in $u.$
Hence, we only prove the former estimate.
To this end, we first use the identity \eqref{E:UDERIVSTU} 
and the estimate \eqref{E:LIERADSTUAREAFORMERRORTERMSPOINTWISEESTIMATE} to deduce that
\begin{align} \label{E:PARTIALUINTEGRALOFPSISQUREDOVERSPHERESFIRSTESTIMATE}
	\left|
		\frac{\partial}{\partial u}
		\int_{S_{t,u}} 
			\Psi^2 
		\, d \argspherevol{(t,u,\vartheta)}
	\right|
	& \leq \int_{S_{t,u}}
				2 |\Psi \Rad \Psi|	
				+ \left\lbrace
						\frac{2}{\rgeo}
						+ C \varepsilon \frac{\ln(\myexp + t)}{1 + t}	
					\right\rbrace
							\Psi^2
			\, d \argspherevol{(t,u,\vartheta)}.
\end{align}
From \eqref{E:PARTIALUINTEGRALOFPSISQUREDOVERSPHERESFIRSTESTIMATE} and Cauchy-Schwarz, 
we deduce that
\begin{align} \label{E:GRONWALLREADYPARTIALUINTEGRALOFPSISQUREDOVERSPHERESID}
	\left|
		\frac{\partial}{\partial u}
		\| \Psi \|_{L^2(S_{t,u})}
	\right|
	& \leq \| \Rad \Psi \|_{L^2(S_{t,u})}
		+ C \| \Psi \|_{L^2(S_{t,u})}.
\end{align}
Applying Gronwall's inequality to \eqref{E:GRONWALLREADYPARTIALUINTEGRALOFPSISQUREDOVERSPHERESID},
using the fact that $\| \Psi \|_{L^2(S_{t,0})} = 0,$
using the fact that $0 \leq u \leq U_0,$ 
and using Cauchy-Schwarz,
we deduce that
\begin{align} \label{E:PSIL2ONTHESPHERESGRONWALLED}
	\| \Psi \|_{L^2(S_{t,u})} 
	& \leq C 
				\int_{u'=0}^u
			 		\| \Rad \Psi \|_{L^2(S_{t,u'})}
			 \, du'
			 \\
	& \leq 	
			C 
			\left(
				\int_{u'=0}^u
			 		\| \Rad \Psi \|_{L^2(S_{t,u'})}^2
				\, du'	
			\right)^{1/2}
			= C \| \Rad \Psi \|_{L^2(\Sigma_t^u)}.
			\notag
\end{align}
Since the already proven inequality \eqref{E:MULTENERGYCOERCIVITY} for $\| \Rad \Psi \|_{L^2(\Sigma_t^u)}$
implies that the right-hand side of \eqref{E:PSIL2ONTHESPHERESGRONWALLED} is
$\lesssim \enzero^{1/2}[\Psi](t,u),$ the desired estimate \eqref{E:MULTENERGYCOERCIVITY}
for $\| \Psi \|_{L^2(S_{t,u})}$ thus follows.

\end{proof}

\section{The fundamental \texorpdfstring{$L^2-$}{square integral-}controlling quantities}
In this section, we define the main quantities that we use
to control $\Psi$ and its derivatives in $L^2.$
In particular, we define a family of Morawetz spacetime integrals
and exhibit their key coerciveness properties.

\begin{definition}[\textbf{The main coercive quantities used for deriving} $L^2$ \textbf{estimates}]
\label{D:MAINCOERCIVEQUANT}
In terms of the energies and fluxes of Def.~\ref{D:ENERGIESANDFLUXES}, we define
\begin{subequations}
\begin{align}
	\totzeromax{N}(t,u)
	& := \max_{|\vec{I}| = N} \sup_{(t',u') \in [0,t] \times [0,u]} 
		\left\lbrace
			\enzero[\mathscr{Z}^{\vec{I}} \Psi](t',u')	
			+ \flzero[\mathscr{Z}^{\vec{I}} \Psi](t',u')
		\right\rbrace,
			\label{E:Q0LEQNDEF} \\
	\totonemax{N}(t,u)
	& := \max_{|\vec{I}| = N} \sup_{(t',u') \in [0,t] \times [0,u]} 
		\left\lbrace
			\enone[\mathscr{Z}^{\vec{I}} \Psi](t',u')	
			+ \flone[\mathscr{Z}^{\vec{I}} \Psi](t',u')
		\right\rbrace,
		\label{E:Q1LEQNDEF}
\end{align}

\begin{align}
	\totzeromax{\leq N}(t,u)
	& := \max_{0 \leq M \leq N} \totzeromax{M}(t,u),
		\label{E:MAXEDQ0LEQNDEF} \\
	\totonemax{\leq N}(t,u)
	& := \max_{0 \leq M \leq N} \totonemax{M}(t,u).
		\label{E:MAXEDQ1LEQNDEF}
\end{align}

\end{subequations}

\end{definition}

The following coercive spacetime integrals play a fundamental role
in controlling some of the error integrals.

\begin{definition}[\textbf{The coercive Morawetz spacetime integral}]
\label{D:COERCIVEMORAWETZINTEGRAL}
	To a function $\Psi,$ we associate the following integrals:
	\begin{subequations}
	\begin{align} \label{E:COERCIVEMORDEF} 
		\Morint[\Psi](t,u)
		& :=
	 	\frac{1}{2}
	 	\int_{\mathcal{M}_{t,u}}
			\rgeo^2
			[\Lunit \upmu]_{-}
			|\angdiff \Psi|^2
		\, d \vol, 
				\\
		\totMormax{N}(t,u) 
		& := \max_{|\vec{I}|=N} \Morint[\mathscr{Z}^{\vec{I}} \Psi](t,u),
			\\
		\totMormax{\leq N}(t,u) 
		& := \max_{0 \leq M \leq N} \totMormax{M}(t,u).
	\end{align}
	\end{subequations}
\end{definition}

\begin{remark}[\textbf{The origin of the Morawetz integral}]
	Recall that the integral \eqref{E:COERCIVEMORDEF} appears 
	in the energy-flux identities for $\Psi$ 
	generated by the Morawetz multiplier $\Mor;$
	see \eqref{E:MORAWETZENERGYERRORINTEGRANDS}.
\end{remark}

We now quantify the $L^2$ coerciveness of the integrals from Def.~\ref{D:COERCIVEMORAWETZINTEGRAL}.

\begin{lemma}[\textbf{Quantification of the $L^2$ coerciveness of $\Morint[\Psi](t,u)$}] 
	\label{L:MORAWETZSPACETIMECOERCIVITY}
	Let $\mathbf{1}_{\lbrace \upmu \leq 1/4 \rbrace}(t,u,\vartheta)$ 
	denote the characteristic function of the spacetime set $\lbrace \upmu \leq 1/4 \rbrace.$ 
	Under the small-data and bootstrap assumptions 
	of Sects.~\ref{S:PSISOLVES}-\ref{S:C0BOUNDBOOTSTRAP},
	if $\varepsilon$ is sufficiently small, 
	then the following inequality holds for 
	the Morawetz spacetime integral \eqref{E:COERCIVEMORDEF} 
	whenever $(t,u) \in [0,\Tboot) \times [0,U_0]:$
	\begin{align} \label{E:MORAWETZSPACETIMECOERCIVITY}
		\Morint[\Psi](t,u) 
		& \geq 
		\int_{\mathcal{M}_{t,u}}
			\mathbf{1}_{\lbrace \upmu \leq 1/4 \rbrace}
			\frac{\rgeo(t',u')}{1 + \ln \left(\frac{\rgeo(t',u')}{\rgeo(0,u')} \right)}
			|\angdiff \Psi|^2
		\, d \vol.
	\end{align}
\end{lemma}

\begin{proof}
	Inequality \eqref{E:MORAWETZSPACETIMECOERCIVITY} follows from
	the estimate \eqref{E:SMALLMUIMPLIESLMUISNEGATIVE}.
\end{proof}

In the next proposition, we explicitly quantify the 
$L^2$ coerciveness of $\totzeromax{N}(t,u)$ and $\totonemax{N}(t,u)$
in a manner that will help us estimate the error integrals on the right-hand
side of the identities of Prop.~\ref{P:DIVTHMWITHCANCELLATIONS}.

\begin{remark}[\textbf{The importance of some sharp constants}]
Some of the sharp constants, such as the $1$ on the right-hand side of \eqref{E:RADPSIL2INTERMSOFZERO}
and the $\sqrt{2}$ on the right-hand side of \eqref{E:SQRTUPMUANGDIFFPSIL2INTERMSOFONE},
are important because they affect the number of derivatives
we need to close our estimates; see, for example, inequalities
\eqref{E:KEYPROOFSTEPSIGMAMINUSANGDIFFUPMURENORMALIZEDSHARPL2INTERMSOFQANDWIDETILDEQ},
\eqref{E:KEYFIRSTTERMWHERETHECONSTANTMATTERS},
and \eqref{E:SIMPLEENERGYNORMCOMPARISON}.
\end{remark}

\begin{proposition}[\textbf{Quantitative $L^2$ coerciveness of $\totzeromax{N}(t,u)$ and $\totonemax{N}(t,u)$}]
		\label{P:L2NORMSOFPSIINTERMSOFTHECOERCIVEQUANTITIES}
		Let $0 \leq N \leq 24$ be an integer.
		Under the 
		small data and
		bootstrap assumptions 
		of Sects.~\ref{S:PSISOLVES}-\ref{S:C0BOUNDBOOTSTRAP},
		if $\varepsilon$ is sufficiently small, 
		then the quantities $\totzeromax{N}(t,u)$ and $\totonemax{N}(t,u)$
		from Def.~\ref{D:MAINCOERCIVEQUANT}
		have the following coerciveness properties.
		
		\medskip
		
		\noindent \underline{\textbf{Coerciveness along $\Sigma_t^u$.}}
		\begin{subequations}
		\begin{align}
			\left\| 
				\mathscr{Z}^N \Psi 
			\right\|_{L^2(S_{t,u})}
			& \leq C \totzeromax{N}^{1/2}(t,u),
				\label{E:PSIL2SPHERESINTERMSOFZERO} \\
			\left\| 
				\mathscr{Z}^N \Psi 
			\right\|_{L^2(\Sigma_t^u)}
			& \leq C \totzeromax{N}^{1/2}(t,u),
				\label{E:PSIL2INTERMSOFZERO} \\
			\left\| 
				\sqrt{\upmu} \Lunit \mathscr{Z}^N \Psi 
			\right\|_{L^2(\Sigma_t^u)}
			& \leq \sqrt{2} \totzeromax{N}^{1/2}(t,u),
				\label{E:SQRTUPMULPSIL2INTERMSOFZERO} \\
		\left\| 
			\upmu \Lunit \mathscr{Z}^N \Psi 
		\right\|_{L^2(\Sigma_t^u)}
			& \leq \sqrt{2} \totzeromax{N}^{1/2}(t,u),
				\label{E:UPMULPSIL2INTERMSOFZERO} \\
		\left\| 
			\Rad \mathscr{Z}^N \Psi 
		\right\|_{L^2(\Sigma_t^u)}
		& \leq \totzeromax{N}^{1/2}(t,u),
				\label{E:RADPSIL2INTERMSOFZERO} \\
		\left\| 
			\sqrt{\upmu} \angdiff \mathscr{Z}^N \Psi 
		\right\|_{L^2(\Sigma_t^u)}
		& \leq \sqrt{2} \totzeromax{N}^{1/2}(t,u),
				\label{E:SQRTUPMUANGDIFFPSIL2INTERMSOFZERO} \\
		\left\| 
			\upmu \angdiff \mathscr{Z}^N \Psi 
		\right\|_{L^2(\Sigma_t^u)}
		& \leq  \totzeromax{N}^{1/2}(t,u),
		\label{E:UPMUANGDIFFPSIL2INTERMSOFZERO}
		\end{align}
		\end{subequations}

		\begin{subequations}
		\begin{align}
			(1 + t)
			\left\| 
				\Lunit \mathscr{Z}^N \Psi 
			\right\|_{L^2(\Sigma_t^u)}
			& \leq C \totzeromax{N}^{1/2}(t,u)
				+ C \frac{1}{\upmu_{\star}^{1/2}(t,u)} \totonemax{N}^{1/2}(t,u),
				\label{E:LPSIL2INTERMSOFZEROANDONE} \\
			(1 + t)
			\left\|
				\sqrt{\upmu} \Lunit \mathscr{Z}^N \Psi 
			\right\|_{L^2(\Sigma_t^u)}
			& \leq C \ln^{1/2}(\myexp + t)\totzeromax{N}^{1/2}(t,u)
				+ C \totonemax{N}^{1/2}(t,u),
				\label{E:SQRTUPMULPSIL2INTERMSOFZEROANDONE} \\
		  \left\| 
				\rgeo
				\sqrt{\upmu}
				\left\lbrace 
					\Lunit + \frac{1}{2} \mytr \upchi 
				\right\rbrace	
				\mathscr{Z}^N \Psi 
			\right\|_{L^2(\Sigma_t^u)}
			& \leq \sqrt{2} \totonemax{N}^{1/2}(t,u),
				\label{E:SQRTUPMULPLUSTRCHIPSIL2INTERMSOFONE} \\
			\left\| 
				\rgeo \sqrt{\upmu} \angdiff \mathscr{Z}^N \Psi 
			\right\|_{L^2(\Sigma_t^u)}
			& \leq \sqrt{2} \totonemax{N}^{1/2}(t,u),
				\label{E:SQRTUPMUANGDIFFPSIL2INTERMSOFONE}
		\end{align}
		\end{subequations}
		
		\begin{align}
			\left\| 
				\mathscr{Z}^N \Psi 
			\right\|_{L^2(\Sigma_t^u)}
			& \leq C \totzeromax{N-1}^{1/2}(t,u)
				+ C \frac{1}{\upmu_{\star}^{1/2}(t,u)} \totonemax{N-1}^{1/2}(t,u),
			\label{E:PSIL2INTERMSOFLOWERORDERZEROANDONE}
		\end{align}
		
		\begin{subequations}
		\begin{align}
			(1+t)
			\left\| 
				\Lunit \mathscr{Z}^N \Psi 
			\right\|_{L^2(\Sigma_t^u)}
			& \leq C \totzeromax{N+1}^{1/2}(t,u),
				\label{E:LPSIL2INTERMSOFHIGHERORDERZERO} \\
			(1 + t)
			\left\| 
				\angdiff \mathscr{Z}^N \Psi 
			\right\|_{L^2(\Sigma_t^u)}
			& \leq C \totzeromax{N+1}^{1/2}(t,u).
			\label{E:ANGDIFFPSIL2INTERMSOFHIGHERORDERZERO}
		\end{align}
		\end{subequations}
		
		\medskip
		
		\noindent \underline{\textbf{Coerciveness along $\mathcal{C}_u^t$.}}
		\begin{subequations}
		\begin{align}
			\left\| 
				\Lunit \mathscr{Z}^N \Psi 
			\right\|_{L^2(\mathcal{C}_u^t)}
			& \leq \totzeromax{N}^{1/2}(t,u),
				\label{E:LPSICONEL2INTERMSOFZERO} \\
			\left\| 
				\sqrt{\upmu} \Lunit \mathscr{Z}^N \Psi 
			\right\|_{L^2(\mathcal{C}_u^t)}
			& \leq \totzeromax{N}^{1/2}(t,u),
				\label{E:SQRTUPMULPSICONEL2INTERMSOFZERO} \\
			\left\| 
				\sqrt{\upmu} \angdiff \mathscr{Z}^N \Psi 
			\right\|_{L^2(\mathcal{C}_u^t)}
			& \leq \totzeromax{N}^{1/2}(t,u),
			\label{E:SQRTUPMUANGDIFFCONEL2INTERMSOFZERO}
		\end{align}
		\end{subequations}
		
		\begin{align}
			\left\| 
				(1 + t')
				\left\lbrace 
					\Lunit + \frac{1}{2} \mytr \upchi 
				\right\rbrace	
				\mathscr{Z}^N \Psi 
			\right\|_{L^2(\mathcal{C}_u^t)}
			& \leq C \totonemax{N}^{1/2}(t,u).
			\label{E:ONEPLUSTTIMESLPLUSTRCHIPSICONEL2INTERMSOFONE}
		\end{align}
		
\end{proposition}

\begin{proof}
	Most of the inequalities in
	Prop.~\ref{P:L2NORMSOFPSIINTERMSOFTHECOERCIVEQUANTITIES}
	follow directly from Lemma~\ref{L:COERCIVEENERGIESANDFLUXES}
	and Def.~\ref{D:MAINCOERCIVEQUANT}.
	We sketch the proofs of the estimates that require some additional ingredients.
	To deduce \eqref{E:LPSIL2INTERMSOFZEROANDONE},
	we split
	$\Lunit \mathscr{Z}^N \Psi 
	= \left\lbrace \Lunit + \frac{1}{2} \mytr \upchi \right\rbrace \mathscr{Z}^N \Psi
	- \frac{1}{2} \mytr \upchi \mathscr{Z}^N \Psi$
	and use
	the estimate $\rgeo |\mytr \upchi| \lesssim 1$ (that is, \eqref{E:CRUDELOWERORDERC0BOUNDDERIVATIVESOFANGULARDEFORMATIONTENSORS}).
	To deduce \eqref{E:SQRTUPMULPSIL2INTERMSOFZEROANDONE}, we also use
	the estimate $\upmu \lesssim \ln(\myexp + t)$ (that is, \eqref{E:C0BOUNDCRUCIALEIKONALFUNCTIONQUANTITIES}).
	The estimate \eqref{E:PSIL2INTERMSOFLOWERORDERZEROANDONE}
	follows from separately considering the three cases
	$\mathscr{Z}^N = \rgeo \Lunit \mathscr{Z}^{N-1},$
	$\mathscr{Z}^N = \Rad \mathscr{Z}^{N-1},$
	and $\mathscr{Z}^N = \Rot \mathscr{Z}^{N-1},$
	where in last case, we first use 
	the estimate \eqref{E:ROTATIONPOINTWISENORMESTIMATE} to bound
	$|\Rot \mathscr{Z}^{N-1} \Psi| \lesssim C (1 + t) |\angdiff \mathscr{Z}^{N-1} \Psi|.$
	To deduce \eqref{E:LPSIL2INTERMSOFHIGHERORDERZERO}, we
	only need to note that since $\rgeo \Lunit \in \mathscr{Z},$ we have $\Lunit \mathscr{Z}^N = \rgeo^{-1} \mathscr{Z}^{N+1}.$
	To deduce \eqref{E:ANGDIFFPSIL2INTERMSOFHIGHERORDERZERO}, we only need to use
	inequality \eqref{E:FUNCTIONPOINTWISEANGDINTERMSOFANGLIEO}
	to deduce that $\left| \angdiff \mathscr{Z}^N \Psi \right| \lesssim \rgeo^{-1} \sum_{l=1}^3 \left| \Rot_{(l)} \mathscr{Z}^N \Psi  \right|.$
\end{proof}

In the next lemma, we quantify the smallness of $\totzeromax{N}$ and $\totonemax{N}$ along $\Sigma_0^1$
under our small data and bootstrap assumptions.

\begin{lemma}[\textbf{Smallness of $\totzeromax{N}$ and $\totonemax{N}$ along $\Sigma_0^1$}]
		\label{L:CONTROLLINGQUANTITIESAREINTIALLYSMALL}
		Let $0 \leq N \leq 24$ be an integer, 
		and let $\mathring{\upepsilon}$ be the size of the data as defined by \eqref{E:SMALLDATA}.
		Under the 
		small data and
		bootstrap assumptions 
		of Sects.~\ref{S:PSISOLVES}-\ref{S:C0BOUNDBOOTSTRAP},
		if $\varepsilon$ is sufficiently small
		and $u \in [0,U_0],$ then
		\begin{align} \label{E:TOTONEMAXINITIALLYBOUNDEDBYTOZEROMAX}
			\totonemax{N}^{1/2}(0,u) \leq C \totzeromax{N}^{1/2}(0,u)
		\end{align}
		and
		\begin{align} \label{E:TOTZEROMAXDATABOUND}
			\totzeromax{N}^{1/2}(0,u)
			& \leq C \mathring{\upepsilon}.
		\end{align}
\end{lemma}

\begin{proof}
	Lemma~\ref{L:CONTROLLINGQUANTITIESAREINTIALLYSMALL} follows
	easily from Lemma~\ref{L:SMALLINITIALSOBOLEVNORMS},
	Lemma~\ref{L:COERCIVEENERGIESANDFLUXES},
	Def.~\ref{D:MAINCOERCIVEQUANT},
	and the fact that 
	$\rgeo |\mytr \upchi| \lesssim 1$
	(that is, \eqref{E:CRUDELOWERORDERC0BOUNDDERIVATIVESOFANGULARDEFORMATIONTENSORS}).
\end{proof}


\chapter{Top-Order Pointwise Commutator Estimates Involving the Eikonal Function}
\label{C:TOPORDEREIKONALFUNCTIONCOMMUTATIONPOINTWISE}
\thispagestyle{fancy}
In Chapter~\ref{C:TOPORDEREIKONALFUNCTIONCOMMUTATIONPOINTWISE}, we derive a collection of pointwise commutator estimates
for various quantities that are constructed out of the eikonal function.
In Chapter~\ref{C:ERRORTERMSOBOLEV}, 
we use these commutator estimates to help us derive top-order
$L^2$ estimates without losing derivatives.

\section{Top-order pointwise commutator estimates 
connecting \texorpdfstring{$\angD^2 \upmu$}{the angular Hessian of the inverse foliation density} to 
\texorpdfstring{$\angLie_{\Rad} \upchi^{(Small)}$}{the radial Lie derivative of the re-centered null second fundamental form}}
Our main goal in this section is to show that 
all top-order derivatives of 
$\angD^2 \upmu$ 
can be controlled in terms of
the top-order derivatives of 
$\angLie_{\Rad} \upchi^{(Small)}$
and vice versa. This will save us a great deal of redundant effort
by allowing us to avoid deriving
a fully modified equation
and separate fully modified
estimates for $\angLap \upmu.$
Specifically, by the last inequality in 
\eqref{E:TOPORDERDERIVATIVESOFANGDSQUAREDUPMUINTERMSOFCONTROLLABLE},
the top-order derivatives 
$\angLap \mathscr{Z}^{N-1}\upmu$ are determined 
by the top-order derivatives
$\mathscr{Z}^{N-1} \Rad \mytr \upchi^{(Small)}$
up to error terms; 
the error terms can be controlled without the need to 
use modified quantities,
and furthermore, we already derived the necessary fully modified
equation for $\mathscr{Z}^{N-1} \Rad \mytr \upchi^{(Small)}$
in Prop.~\ref{P:TOPORDERTRCHIJUNKRENORMALIZEDTRANSPORT}
(in the case $\mathscr{Z}^{N-1} = \mathscr{S}^{N-1},$ 
where full modification is needed).

\begin{proposition}[\textbf{Top-order pointwise commutator estimates connecting} $\angD^2 \upmu$ \textbf{to} $\angLie_{\Rad} \upchi^{(Small)}$]
\label{P:TOPORDERDERIVATIVESOFANGLAPUPMUINTERMSOFCONTROLLABLE}
	Let $1 \leq N \leq 24$ be an integer.
	Under the small-data and bootstrap assumptions 
	of Sects.~\ref{S:PSISOLVES}-\ref{S:C0BOUNDBOOTSTRAP},
	if $\varepsilon$ is sufficiently small, 
	then the following pointwise estimates hold on $\mathcal{M}_{\Tboot,U_0}:$
	\begin{align}
		&
		\left|
			\angD^2 \mathscr{Z}^{N-1} \upmu
			- \angLie_{\mathscr{Z}}^{N-1} \angLie_{\Rad} \upchi^{(Small)}
		\right|,
		  \label{E:TOPORDERDERIVATIVESOFANGDSQUAREDUPMUINTERMSOFCONTROLLABLE} \\
		& 
		\left|
			\angfreeDsquared \mathscr{Z}^{N-1} \upmu
			- \angLie_{\Rad} \angLie_{\mathscr{Z}}^{N-1} \hat{\upchi}^{(Small)}
		\right|,
			\notag \\
		& 
		\left|
			\angLap \mathscr{Z}^{N-1} \upmu
			- \mathscr{Z}^{N-1} \Rad \mytr \upchi^{(Small)}
		\right|,
			\notag \\
		& 
		\left|
			\angLap \mathscr{Z}^{N-1} \upmu
			- \Rad \mathscr{Z}^{N-1} \mytr \upchi^{(Small)}
		\right|
			\notag \\
		& \lesssim 
			\frac{1}{1 + t}
			\left| 
					\fourmyarray
					 [\rgeo \Lunit \mathscr{Z}^{\leq N} \Psi]
					 {\Rad \mathscr{Z}^{\leq N} \Psi}
					 {\rgeo \angdiff \mathscr{Z}^{\leq N} \Psi}
					{\mathscr{Z}^{\leq N} \Psi}
			\right|
			+ 
			\frac{1}{(1 + t)^2}
	 		\left| 
				\myarray
					[\mathscr{Z}^{\leq N} (\upmu - 1)]
					{\rgeo \sum_{a=1}^3 |\mathscr{Z}^{\leq N} \Lunit_{(Small)}^a|}
			\right|.
			\notag
	\end{align}
	
\end{proposition}

Our proof of Prop.~\ref{P:TOPORDERDERIVATIVESOFANGLAPUPMUINTERMSOFCONTROLLABLE} is based on the following
commutator-type lemma.

\begin{lemma}[\textbf{An expression for} $\angLie_{\Rad} \upchi^{(Small)}$ \textbf{in terms of other variables}]
\label{L:LIERADCHIJUNKTRANSPORT}
The symmetric type $\binom{0}{2}$ $S_{t,u}$ tensorfield $\upchi^{(Small)}$ 
defined in \eqref{E:CHIJUNKDEF}
verifies the following transport equation,
where the capital Latin-indexed terms are exact and
the remaining ones are schematic:
\begin{align} \label{E:LIERADCHIJUNKTRANSPORT}
	\angLie_{\Rad} \upchi_{AB}^{(Small)}
	& = \angDsquaredarg{A}{B} \upmu
			+ \frac{\upmu - 1}{\rgeo^2} \gsphere_{AB}
				 \\		
	& \ \ + G_{(Frame)}
				\threemyarray
					[\upmu \Lunit \Rad \Psi]
					{\angdiff \Rad \Psi}
					{\upmu \angD^2 \Psi}
			\notag	\\
	& \ \ + \upmu \upchi^{(Small) \# \#} \upchi_{(Small)}
				+ 	G_{(Frame)}
					\upchi^{(Small)}
					\myarray
						[\upmu \Lunit \Psi]
						{\Rad \Psi}
			+ 	G_{(Frame)}
					\ginversesphere
					\upchi^{(Small)}
					\threemyarray
						[\upmu \Lunit \Psi]
						{\Rad \Psi}
						{\upmu \angdiff \Psi}
				\notag \\
	& \ \ 
			+ 	G_{(Frame)} 
					(\angdiff \upmu)	
					\myarray
						[\Lunit \Psi]
						{\angdiff \Psi}
			\notag
			\\
	& \ \
		+ 		\frac{1}{\rgeo}
					G_{(Frame)}
					\gsphere
					\myarray
						[\upmu \Lunit \Psi]
						{\Rad \Psi}
			+ 	\frac{1}{\rgeo}
					G_{(Frame)}
					\threemyarray
						[\upmu \Lunit \Psi]
						{\Rad \Psi}
						{\upmu \angdiff \Psi}
				\notag \\
	& \ \
		+ (\angD G_{(Frame)})
				\threemyarray
					[\upmu \Lunit \Psi]
					{\Rad \Psi}
					{\upmu \angdiff \Psi}
					\notag \\
	& \ \ 
			+ \myarray
				[G_{(Frame)}^2]
				{G_{(Frame)}'}
				\threemyarray
					[\upmu \Lunit \Psi]
					{\Rad \Psi}
					{\upmu \angdiff \Psi}
				\myarray
					[\Lunit \Psi]
					{\angdiff \Psi}
				+ 
				G_{(Frame)}^2
					\ginversesphere
					\threemyarray
						[\upmu \Lunit \Psi]
						{\Rad \Psi}
						{\upmu \angdiff \Psi}
					\myarray
						[\Lunit \Psi]
						{\angdiff \Psi}.
						\notag
	\end{align}

\end{lemma}

Before proving Lemma~\ref{L:LIERADCHIJUNKTRANSPORT}, we first use it to prove the proposition.

\begin{proof}[Proof of Prop.~\ref{P:TOPORDERDERIVATIVESOFANGLAPUPMUINTERMSOFCONTROLLABLE}]
	We first prove the bound for the first term
	$\left|
			\angD^2 \mathscr{Z}^{N-1} \upmu
			- \angLie_{\Rad} \angLie_{\mathscr{Z}}^{N-1} \upchi^{(Small)}
		\right|$
	on the left-hand side of \eqref{E:TOPORDERDERIVATIVESOFANGDSQUAREDUPMUINTERMSOFCONTROLLABLE}.
	To this end, we write equation \eqref{E:LIERADCHIJUNKTRANSPORT} in the form
	\begin{align} \label{E:SCHEMATICLIERADCHIJUNK}
		\angLie_{\Rad} \upchi_{AB}^{(Small)}
		& = \angDsquaredarg{A}{B} \upmu
			+ \mathfrak{I}_{AB},
	\end{align}
	where $\mathfrak{I}$ denotes all of the terms
	on the right-hand side of \eqref{E:LIERADCHIJUNKTRANSPORT} except for the first one.
	Applying $\angLie_{\mathscr{Z}}^{N-1}$ to equation \eqref{E:SCHEMATICLIERADCHIJUNK},
	we deduce that
	\begin{align} \label{E:COMMUTEDSCHEMATICLIERADCHIJUNK}
		\angLie_{\mathscr{Z}}^{N-1} \angLie_{\Rad}  \upchi_{AB}^{(Small)}
		- \angDsquaredarg{A}{B} \mathscr{Z}^{N-1} \upmu
		& = [\angLie_{\Rad}, \angLie_{\mathscr{Z}}^{N-1}] \upchi_{AB}^{(Small)}
			+ \angLie_{\mathscr{Z}}^{N-1} \mathfrak{I}_{AB}.
	\end{align}
	The desired estimate will follow once we show that the 
	terms on the right-hand side of \eqref{E:COMMUTEDSCHEMATICLIERADCHIJUNK} are in
	magnitude $\lesssim$ the right-hand side of \eqref{E:TOPORDERDERIVATIVESOFANGDSQUAREDUPMUINTERMSOFCONTROLLABLE}.
	To bound the term $([\angD^2, \angLie_{\mathscr{Z}}^{N-1}] \upmu)_{AB},$
	we use inequality \eqref{E:ANGDSQUAREDLIEZNCOMMUTATORACTINGONFUNCTIONSSPOINTWISE}
	with $N-1$ in the role of $N$ and
	$\upmu - 1$ in the role of $f$
	and inequality \eqref{E:C0BOUNDCRUCIALEIKONALFUNCTIONQUANTITIES} 
	to deduce that
	\begin{align} \label{E:SECONDCOMMUTATORTERMESTIMATECOMMUTEDSCHEMATICLIERADCHIJUNK}
		\left|
			[\angD^2, \angLie_{\mathscr{Z}}^{N-1}] (\upmu - 1)
		\right|
		& \lesssim
			\frac{\ln(\myexp + t)}{(1 + t)^2}
			\left| 
				\fourmyarray[\rgeo \Lunit \mathscr{Z}^{\leq N-1} \Psi]
					{\Rad \mathscr{Z}^{\leq N-1} \Psi}
					{\rgeo \angdiff \mathscr{Z}^{\leq N-1} \Psi}
					{\mathscr{Z}^{\leq N-1} \Psi}
			\right|
		+ 
			\frac{1}{(1 + t)^2}
			\left|
				\myarray[\mathscr{Z}^{\leq N} (\upmu - 1)]
					{\sum_{a=1}^3 \rgeo |\mathscr{Z}^{\leq N} \Lunit_{(Small)}^a|} 
			\right|.
	\end{align}


	It remains for us to bound the magnitude of 
	$\angLie_{\mathscr{Z}}^{N-1} \mathfrak{I}$ by the right-hand
	side of \eqref{E:TOPORDERDERIVATIVESOFANGDSQUAREDUPMUINTERMSOFCONTROLLABLE}.
	To this end, we apply $\angLie_{\mathscr{Z}}^{N-1}$ to 
	\eqref{E:LIERADCHIJUNKTRANSPORT}
	and apply the Leibniz rule
	to the terms on the right-hand side of \eqref{E:LIERADCHIJUNKTRANSPORT}
	(except for $\angD^2 \upmu,$ which we have already handled).
	We bound the terms $\mathscr{Z}^M \rgeo$ with \eqref{E:ZNAPPLIEDTORGEOISNOTTOOLARGE}.
	We bound the terms 
	$\mathscr{Z}^M \upmu,$ 
	$\mathscr{Z}^M (\upmu-1),$
	and $\angLie_{\mathscr{Z}}^M \angdiff \upmu$
	with	
	Lemma~\ref{L:LANDRADCOMMUTEWITHANGDIFF}, 
	\eqref{E:FUNCTIONPOINTWISEANGDINTERMSOFANGLIEO},
	and
	\eqref{E:C0BOUNDCRUCIALEIKONALFUNCTIONQUANTITIES}.
	We bound the terms $\angLie_{\mathscr{Z}}^M \gsphere$ 
	and $\angLie_{\mathscr{Z}}^M \ginversesphere$
	with Lemma~\ref{L:POINTWISEBOUNDSDERIVATIVESOFANGULARDEFORMATIONTENSORS}.
	We bound the terms
	$\angLie_{\mathscr{Z}}^M G_{(Frame)}$
	and
	$\angLie_{\mathscr{Z}}^M G_{(Frame)}'$
	with Lemma~\ref{L:POINTWISEESTIMATESGFRAMEINTERMSOFOTHERQUANTITIES}.
	We bound the terms
	$\angLie_{\mathscr{Z}}^M \angD G_{(Frame)}$
	with 
	Lemma~\ref{L:LANDRADCOMMUTEWITHANGDIFF},
	\eqref{E:FUNCTIONPOINTWISEANGDINTERMSOFANGLIEO},
	\eqref{E:ONEFORMANGDINTERMSOFROTATIONALLIE},
	\eqref{E:TYPE02TENSORANGDINTERMSOFROTATIONALLIE},
	\eqref{E:COMMUTATORESTIMATESVECTORFIELDSACTINGONANGDTENSORS},
	Lemma~\ref{L:POINTWISEESTIMATESGFRAMEINTERMSOFOTHERQUANTITIES},
	\eqref{E:C0BOUNDCRUCIALEIKONALFUNCTIONQUANTITIES},
	and the bootstrap assumptions \eqref{E:PSIFUNDAMENTALC0BOUNDBOOTSTRAP}.
	We bound the terms 
	$\angLie_{\mathscr{Z}}^M \upchi^{(Small)}$
	and
	$\angLie_{\mathscr{Z}}^M \upchi^{(Small) \# \#}$
	with \eqref{E:POINTWISEESTIMATESFORCHIJUNKINTERMSOFOTHERVARIABLES} and
	\eqref{E:C0BOUNDCRUCIALEIKONALFUNCTIONQUANTITIES}.
	We bound the terms
	$\threemyarray[\mathscr{Z}^M \Lunit \Psi]{\mathscr{Z}^M \Rad \Psi} {\angLie_{\mathscr{Z}}^M \angdiff \Psi}$
	with Lemma~\ref{L:AVOIDINGCOMMUTING}
	and the bootstrap assumptions \eqref{E:PSIFUNDAMENTALC0BOUNDBOOTSTRAP}.
	In total, these estimates imply that
	$\left|\angLie_{\mathscr{Z}}^{N-1} \mathfrak{I}\right|$ 
	is $\lesssim$
	the right-hand
	side of \eqref{E:TOPORDERDERIVATIVESOFANGDSQUAREDUPMUINTERMSOFCONTROLLABLE} as desired,
	which completes the proof of the bound for
	$\left|
			\angD^2 \mathscr{Z}^{N-1} \upmu
			- \angLie_{\Rad} \angLie_{\mathscr{Z}}^{N-1} \upchi^{(Small)}
		\right|.$
	
To prove \eqref{E:TOPORDERDERIVATIVESOFANGDSQUAREDUPMUINTERMSOFCONTROLLABLE} for 
$\left|
			\angLap \mathscr{Z}^{N-1} \upmu
			-  \Rad \mathscr{Z}^{N-1} \mytr \upchi^{(Small)}
		\right|,$
we first note that it follows trivially from the previously proven estimate for
$\left|
			\angD^2 \mathscr{Z}^{N-1} \upmu
			- \angLie_{\mathscr{Z}}^{N-1} \angLie_{\Rad} \upchi^{(Small)}
		\right|$
that
$\left| \angLap \mathscr{Z}^{N-1} \upmu - \mytr  \angLie_{\mathscr{Z}}^{N-1} \angLie_{\Rad} \upchi^{(Small)} \right|$
is $\lesssim$ the right-hand side of \eqref{E:TOPORDERDERIVATIVESOFANGDSQUAREDUPMUINTERMSOFCONTROLLABLE}.
Hence, to prove the desired estimate, it suffices to show that
\begin{align} \label{E:FIRSTCOMMUTATORTERMESTIMATECOMMUTEDSCHEMATICLIERADCHIJUNK}
		&
		\left|
			[\angLie_{\Rad}, \angLie_{\mathscr{Z}}^{N-1}] \upchi^{(Small)}
		\right|, \\
		&
		\left|
			\mytr  \angLie_{\Rad} \angLie_{\mathscr{Z}}^{N-1} \upchi^{(Small)}
			-  \Rad \mathscr{Z}^{N-1} \mytr \upchi^{(Small)}
		\right|
			\notag \\
		& \lesssim
			\frac{1}{1 + t}
			\left| 
				\fourmyarray[\rgeo \Lunit \mathscr{Z}^{\leq N-1} \Psi]
					{\Rad \mathscr{Z}^{\leq N-1} \Psi}
					{\rgeo \angdiff \mathscr{Z}^{\leq N-1} \Psi}
					{\mathscr{Z}^{\leq N-1} \Psi}
			\right|
		+ 
			\frac{1}{(1 + t)^2}
			\left|
				\myarray[\mathscr{Z}^{\leq N} (\upmu - 1)]
					{\sum_{a=1}^3 \rgeo |\mathscr{Z}^{\leq N} \Lunit_{(Small)}^a|} 
			\right|.
			\notag
	\end{align}
	To deduce the desired estimate \eqref{E:FIRSTCOMMUTATORTERMESTIMATECOMMUTEDSCHEMATICLIERADCHIJUNK} 
	for the first term on the left-hand side of 
	\eqref{E:FIRSTCOMMUTATORTERMESTIMATECOMMUTEDSCHEMATICLIERADCHIJUNK}, 
	we use inequality \eqref{E:RGEOLORRADZNCOMMUTATORACTINGONTENSORFIELDSPOINTWISE} with
	$N-1$ in the role of $N$ and
	$\upchi^{(Small)}$ in the role of $\xi,$ 
	inequality \eqref{E:POINTWISEESTIMATESFORCHIJUNKINTERMSOFOTHERVARIABLES}
	and inequality \eqref{E:C0BOUNDCRUCIALEIKONALFUNCTIONQUANTITIES}.
	To deduce the desired estimate \eqref{E:FIRSTCOMMUTATORTERMESTIMATECOMMUTEDSCHEMATICLIERADCHIJUNK}
	for the second term on the left-hand side of 
	\eqref{E:FIRSTCOMMUTATORTERMESTIMATECOMMUTEDSCHEMATICLIERADCHIJUNK},  
	we use inequality \eqref{E:COMMUTINGTRACEWITHLIEDERIVATIVES} with 
	$\upchi^{(Small)}$ in the role of $\xi,$
	\eqref{E:C0BOUNDCRUCIALEIKONALFUNCTIONQUANTITIES},
	and \eqref{E:POINTWISEESTIMATESFORCHIJUNKINTERMSOFOTHERVARIABLES}
	with $N-1$ in the role of $N$ (to bound the first term on the right-hand side of \eqref{E:COMMUTINGTRACEWITHLIEDERIVATIVES}).
	We have thus proved the desired bound for
	$\left|
			\angLap \mathscr{Z}^{N-1} \upmu
			- \Rad \mathscr{Z}^{N-1} \mytr \upchi^{(Small)}
		\right|.$	
 The proof of the estimate \eqref{E:FIRSTCOMMUTATORTERMESTIMATECOMMUTEDSCHEMATICLIERADCHIJUNK} for
 $\left|
		\angLap \mathscr{Z}^{N-1} \upmu
		- \mathscr{Z}^{N-1} \Rad \mytr \upchi^{(Small)}
\right|$	
is similar; we omit the details.
	
Finally, since the trace-free part of a tensor is in magnitude $\lesssim$ the magnitude of the tensor itself,
we infer from the previously proven estimate for
$\left|
			\angD^2 \mathscr{Z}^{N-1} \upmu
			- \angLie_{\mathscr{Z}}^{N-1} \angLie_{\Rad} \upchi^{(Small)}
		\right|$
that in order to prove that
$\left| \angfreeDsquared \mathscr{Z}^{N-1} \upmu - \angLie_{\Rad} \angLie_{\mathscr{Z}}^{N-1} \hat{\upchi}^{(Small)} \right|$
is $\lesssim$ the right-hand side of \eqref{E:TOPORDERDERIVATIVESOFANGDSQUAREDUPMUINTERMSOFCONTROLLABLE}, 
it suffices to show that
\begin{align} \label{E:COMMUTINGTRACEFREEANDLIEDIFFERENTIATIONRADCHIJUNK}
&
\left|
	\angLie_{\Rad} \angLie_{\mathscr{Z}}^{N-1} \hat{\upchi}^{(Small)}			
	- 
	\left\lbrace
		(\angLie_{\Rad} \angLie_{\mathscr{Z}}^{N-1} \upchi^{(Small)})
		- \frac{1}{2} \mytr  (\angLie_{\Rad} \angLie_{\mathscr{Z}}^{N-1} \upchi^{(Small)}) \gsphere
	\right\rbrace
\right|
	\\
& \lesssim 
			\frac{1}{1 + t}
			\left| 
					\fourmyarray
					 [\rgeo \Lunit \mathscr{Z}^{\leq N} \Psi]
					 {\Rad \mathscr{Z}^{\leq N} \Psi}
					 {\rgeo \angdiff \mathscr{Z}^{\leq N} \Psi}
					{\mathscr{Z}^{\leq N} \Psi}
			\right|
			+ 
			\frac{1}{(1 + t)^2}
	 		\left| 
				\myarray
				[\mathscr{Z}^{\leq N} (\upmu - 1)]
				{\rgeo \sum_{a=1}^3 |\mathscr{Z}^{\leq N} \Lunit_{(Small)}^a|}
			\right|.
			\notag
\end{align}
The estimate \eqref{E:COMMUTINGTRACEFREEANDLIEDIFFERENTIATIONRADCHIJUNK} follows from
inequality \eqref{E:COMMUTINGTRACEFREEWITHLIEDERIVATIVES} with 
$\upchi^{(Small)}$ in the role of $\xi,$
\eqref{E:C0BOUNDCRUCIALEIKONALFUNCTIONQUANTITIES},
and \eqref{E:POINTWISEESTIMATESFORCHIJUNKINTERMSOFOTHERVARIABLES}
with $N-1$ in the role of $N$ (to bound the first term on the right-hand side of \eqref{E:COMMUTINGTRACEFREEWITHLIEDERIVATIVES}).

\end{proof}

It remains for us to prove Lemma~\ref{L:LIERADCHIJUNKTRANSPORT}. Our proof
is based in part on the 
following lemma, which is an analog of Cor.~\ref{C:ALPHAHASNOBADTERMS}.
Among the most important aspects of the lemma is our identification of the
$\upmu^{-1}-$singular terms, which are on the second line of the right-hand side
of \eqref{E:GAMMACURVATURECOMPONENTEXPRESSION}.

\begin{lemma}[\textbf{An expression for the curvature component} $\Cur_{\Rad A \Lunit B}$]
\label{L:GAMMACURVATURECOMPONENTEXPRESSION}
Let $\Cur$ be the Riemann curvature tensor of the spacetime metric $g$
from Def.~\ref{D:SPACETIMERIEMANN}.
Then the curvature component $\Cur_{\Rad A \Lunit B}$
can be expressed as follows:
\begin{align}  \label{E:GAMMACURVATURECOMPONENTEXPRESSION}
	\Cur_{\Rad A \Lunit B}
	& = \frac{1}{2} 
			\left\lbrace
				\upmu \angGdoublearg{\Radunit}{A} \angdiffarg{B} \Lunit \Psi
					+ \angGdoublearg{\Lunit}{A} \angdiffarg{B} \Rad \Psi
					- \angGdoublearg{A}{B} \Lunit \Rad \Psi
					- \upmu G_{\Lunit \Radunit} \angDsquaredarg{A}{B} \Psi
			\right\rbrace
				\\
& \ \ + \frac{1}{4} 
				\upmu^{-1}
				\angGdoublearg{\Lunit}{A}
				\angGdoublearg{\Lunit}{B}
				(\Rad \Psi)^2
		- \frac{1}{2} 
				\upmu^{-1}	
				\angGdoublearg{\Lunit}{A}
				(\angdiffarg{B} \upmu) 
				\Rad \Psi
		\notag \\& \ \ + \frac{1}{2} 
				\left\lbrace
					G_{\Lunit \Radunit} \upchi_{AB} \Rad \Psi 
					+ \upmu \angGdoublearg{\Lunit}{A} \upchi_B^{\ C} \angdiffarg{C} \Psi		
					- \upmu \angGdoublearg{\Radunit}{B} \upchi_A^{\ C} \angdiffarg{C} \Psi
				\right\rbrace
			\notag \\
& \ \   +  
				\myarray
				[G_{(Frame)}^2]
				{G_{(Frame)}'}
				\threemyarray
					[\upmu \Lunit \Psi]
					{\Rad \Psi}
					{\upmu \angdiff \Psi}
				\myarray
					[\Lunit \Psi]
					{\angdiff \Psi}
				+ 
				G_{(Frame)}^2
					\ginversesphere
					\threemyarray
						[\upmu \Lunit \Psi]
						{\Rad \Psi}
						{\upmu \angdiff \Psi}
					\myarray
						[\Lunit \Psi]
						{\angdiff \Psi}.
			\notag
\end{align}

\end{lemma}

\begin{proof}
We claim that by contracting \eqref{E:RECTANGULRCURVATURECOMPONENTS} against
$\Rad^{\mu} X_A^{\nu} \Lunit^{\alpha} X_B^{\beta},$
we can deduce the following identity:
\begin{align} \label{E:FIRSTRELATIONGAMMACURVATURECOMPONENTEXPRESSION}
	\Cur_{\Rad A \Lunit B}
	& = \frac{1}{2} 
		\left\lbrace
			\angGdoublearg{\Lunit}{A} \D_{\Rad B}^2 \Psi
			+ \upmu \angGdoublearg{\Radunit}{B} \D_{\Lunit A}^2 \Psi
			- \angGdoublearg{A}{B} \D_{\Lunit \Rad}^2 \Psi
			- \upmu G_{\Lunit \Radunit} \D_{AB}^2 \Psi
		\right\rbrace
			\\
	& \ \ +  
				\myarray
				[G_{(Frame)}^2]
				{G_{(Frame)}'}
				\threemyarray
					[\upmu \Lunit \Psi]
					{\Rad \Psi}
					{\upmu \angdiff \Psi}
				\myarray
					[\Lunit \Psi]
					{\angdiff \Psi}
				+ 
				G_{(Frame)}^2
					\ginversesphere
					\myarray
						[\upmu \Lunit \Psi]
						{\upmu \angdiff \Psi}
					\myarray
						[\Lunit \Psi]
						{\angdiff \Psi}.
				\notag
\end{align}
Clearly, the terms on the first line of the right-hand side of \eqref{E:FIRSTRELATIONGAMMACURVATURECOMPONENTEXPRESSION}
arise from the terms in the first line on the right-hand side of \eqref{E:RECTANGULRCURVATURECOMPONENTS}.
To analyze the terms that arise from the terms 
in the second line of the right-hand side of \eqref{E:RECTANGULRCURVATURECOMPONENTS}, 
we first factor the contraction vector $\Rad$ as follows: $\Rad^{\mu} = \upmu \Radunit^{\mu}.$
We then pair this factor $\upmu$ with the factor $(g^{-1})^{\kappa \lambda}$
and use the decomposition \eqref{E:GINVERSERECTDECOMP} to deduce that the terms of interest 
are a product of $G_{(Frame)} G_{(Frame)}^{\#}$ and
$\upmu (g^{-1})^{\kappa \lambda} (\partial_{\kappa}\Psi) (\partial_{\lambda} \Psi) 
= \upmu |\angdiff \Psi|^2
- \upmu (\Lunit \Psi)^2
- 2 (\Lunit \Psi) (\Rad \Psi).$
Hence, these terms are of the form of the terms on the second line of the right-hand
side of \eqref{E:FIRSTRELATIONGAMMACURVATURECOMPONENTEXPRESSION}.
In the terms that arise from the terms on the third line of the right-hand side of \eqref{E:RECTANGULRCURVATURECOMPONENTS},
the $\lambda, \kappa$ indices are contractions and not differentiations
of $\Psi.$ Hence we can use the decomposition $(g^{-1})^{\kappa \lambda} = (\ginversesphere)^{\kappa \lambda}
		- \Lunit^{\kappa} \Lunit^{\lambda}
		- \Lunit^{\kappa} \Radunit^{\lambda}
		- \Radunit^{\kappa} \Lunit^{\lambda}$ 
and thus conclude that the terms that arise
are a product of $\upmu G_{(Frame)}^2$ and two distinct elements of  
the set $\lbrace \Radunit \Psi, \Lunit \Psi, \angdiffarg{A} \Psi, \angdiffarg{B} \Psi \rbrace.$
Hence, these terms are of the form of the terms on the second line of the right-hand
side of \eqref{E:FIRSTRELATIONGAMMACURVATURECOMPONENTEXPRESSION}. Similarly,
the terms that arise from the terms on the last line of the right-hand side of \eqref{E:RECTANGULRCURVATURECOMPONENTS}
are a product of $\upmu G_{(Frame)}'$ and two distinct elements of  
the set $\lbrace \Radunit \Psi, \Lunit \Psi, \angdiffarg{A} \Psi, \angdiffarg{B} \Psi \rbrace$
and thus are of the form of the terms on the second line of the right-hand
side of \eqref{E:FIRSTRELATIONGAMMACURVATURECOMPONENTEXPRESSION}.

We now further analyze the terms on the first line of the right-hand side of \eqref{E:FIRSTRELATIONGAMMACURVATURECOMPONENTEXPRESSION}.
To this end, we use Lemma~\ref{L:CONNECTIONLRADFRAME},
the decompositions \eqref{E:ZETADECOMPOSED}
and \eqref{E:ANGKDECOMPOSED}, 
and the identity $\Lunit \upmu = G_{(Frame)} \myarray [\upmu \Lunit \Psi] {\Rad \Psi}$ 
(which follows from Lemma~\ref{L:UPMUFIRSTTRANSPORT}) to deduce that
\begin{align}
	\D_{\Rad B}^2 \Psi
	& = \angdiffarg{B} \Rad \Psi
		+ \frac{1}{2} \upmu^{-1} \angGdoublearg{\Lunit}{B} (\Rad \Psi)^2
		- \upmu^{-1} (\angdiffarg{B} \upmu) \Rad \Psi
		+ \upmu \upchi_B^{\ C} \angdiffarg{C} \Psi 
			\\
	& \ \ 
		+   G_{(Frame)}
				\myarray
					[\upmu \Lunit \Psi]
					{\Rad \Psi}
				\myarray
					[\Lunit \Psi]
					{\angdiff \Psi}
		+   G_{(Frame)}
				\ginversesphere
				\threemyarray
					[\upmu \Lunit \Psi]
					{\Rad \Psi}
					{\upmu \angdiff \Psi}
				\angdiff \Psi,
					\notag 	\\
	\D_{\Lunit \Rad}^2 \Psi
	&	= \Lunit \Rad \Psi
		+ G_{(Frame)}
				\myarray
					[\upmu \Lunit \Psi]
					{\Rad \Psi}
				\Lunit \Psi
		+ G_{(Frame)}
			\ginversesphere
				\threemyarray
					[\upmu \Lunit \Psi]
					{\Rad \Psi}
					{\upmu \angdiff \Psi}
				\angdiff \Psi.
\end{align}
Substituting these two identities and also the identities 
\eqref{E:DSQUAREDLAPSIINTERMSOFANGDIFFALPSI}
and
\eqref{E:DSQUAREDABPSIINTERMSOFANGDQUAREDAB}
into the terms on the first line of the right-hand side of
\eqref{E:FIRSTRELATIONGAMMACURVATURECOMPONENTEXPRESSION},
we arrive at the desired decomposition \eqref{E:GAMMACURVATURECOMPONENTEXPRESSION}.

\end{proof}

The next lemma provides an analog of the identity \eqref{E:CHITRANSPORT}
but with the derivative $\angLie_{\Rad}$ in place of $\angLie_{\Lunit}.$
Among the most important aspects of the lemma is our identification of the
$\upmu^{-1}-$singular terms, which are on the second line of the right-hand side
of \eqref{E:PRELIMINARYRADDIRECTIONCHITRANSPORT}.

\begin{lemma}[\textbf{A preliminary expression for} $\angLie_{\Rad} \upchi$ \textbf{in terms of other variables}]
\label{L:PRELIMINARYRADDIRECTIONCHITRANSPORT}
The type $\binom{0}{2}$ $S_{t,u}$ tensor $\upchi$ verifies the following transport equation,
where the capital Latin-indexed terms are exact and the remaining ones are schematic:
\begin{align}  \label{E:PRELIMINARYRADDIRECTIONCHITRANSPORT}
	\angLie_{\Rad} \upchi_{AB}
	& = \angDsquaredarg{A}{B} \upmu
			+ 
			\frac{1}{2}
			\left\lbrace
				\angDarg{A} (\upmu \upzeta)_B
				+ \angDarg{B} (\upmu \upzeta)_A	
			\right\rbrace
		-	\frac{1}{2}
			\left\lbrace
				\Cur_{\Rad A \Lunit B}
				+ \Cur_{\Rad B \Lunit A}
			\right\rbrace
			\\
	& \ \
			+ \frac{1}{4} 
				\upmu^{-1} 
				\angGdoublearg{\Lunit}{A} 
				\angGdoublearg{\Lunit}{B} 
				(\Rad \Psi)^2
			- \frac{1}{4} 
				\upmu^{-1} 
				\angGdoublearg{\Lunit}{A} 
				\angdiffarg{B} \upmu
			- \frac{1}{4} 
				\upmu^{-1} 
				\angGdoublearg{\Lunit}{B} 
				\angdiffarg{A} \upmu
				\notag \\
	& \ \ 
			- \upmu \upchi_A^{\ C} \upchi_{BC}
		\notag \\ 
	& \ \ +
			\left\lbrace
				- \frac{1}{2} G_{\Lunit \Lunit} \Rad \Psi
				+ \frac{1}{2} \upmu G_{\Lunit \Lunit} \Lunit \Psi
				+ \upmu G_{\Lunit \Radunit} \Lunit \Psi
			\right\rbrace
			\upchi_{AB}
		+ \frac{1}{2}
			\left\lbrace
				 \upmu \upchi_A^{\ C} \angkdoublearg{B}{C}
				+ \upmu \upchi_B^{\ C} \angkdoublearg{A}{C}	
			\right\rbrace
			\notag \\
	& \ \ 
			+ \myarray
				[G_{(Frame)}^2]
				{G_{(Frame)}'}
				\threemyarray
					[\upmu \Lunit \Psi]
					{\Rad \Psi}
					{\upmu \angdiff \Psi}
				\myarray
					[\Lunit \Psi]
					{\angdiff \Psi}
				+ 
				G_{(Frame)}^2
					\ginversesphere
					\threemyarray
						[\upmu \Lunit \Psi]
						{\Rad \Psi}
						{\upmu \angdiff \Psi}
					\myarray
						[\Lunit \Psi]
						{\angdiff \Psi}
						\notag \\
	& \ \ 
			+ 	G_{(Frame)} 
					(\angdiff \upmu)	
					\myarray
						[\Lunit \Psi]
						{\angdiff \Psi}.
					\notag
\end{align}

\end{lemma}

\begin{proof}
Using the torsion-free property 
$\Lie_{\Rad} X_A 
= [\Rad, X_A] 
= \D_{\Rad} X_A - \D_A \Rad,$
we compute that
\begin{align} \label{E:GAMMACURVATUREFIRSTRELATION}
		\angLie_{\Rad} \upchi_{AB}
		= \Lie_{\Rad} \upchi_{AB}
		& = \Rad(\upchi_{AB})
				- \upchi(\Lie_{\Rad} X_A, X_B)
				- \upchi(X_A, \Lie_{\Rad} X_B)
				\\
		& = \Rad(\upchi_{AB})
			+ \upchi(\D_A \Rad, X_B)
			+ \upchi(\D_B \Rad, X_A)
			- \upchi(\D_{\Rad} X_A, X_B)
			- \upchi(\D_{\Rad} X_B, X_A).
			\notag
	\end{align}
	We next use the Def.~\ref{D:SPACETIMERIEMANN} of the Riemann curvature tensor
	to express the first term on the right-hand side of \eqref{E:GAMMACURVATUREFIRSTRELATION} as follows:
		\begin{align} \label{E:GAMMACURVATURESECONDRELATION}
		\Rad(\upchi_{AB})
		& = \Rad(g(\D_A \Lunit, X_B))
			= g(\D_{\Rad} (\D_A \Lunit), X_B)
			+ g(\D_A \Lunit, \D_{\Rad} X_B)
			\\
		& = g(\D_A (\D_{\Rad} \Lunit), X_B)
			+ g(\D_{[\Rad, X_A]} \Lunit, X_B)
			+ g(\D_A \Lunit, \D_{\Rad} X_B)
			- \Cur_{\Rad A \Lunit B}.
			\notag
	\end{align}	
	Using the above torsion-free property,
	Lemma~\ref{L:VECTORFIELDCOMMUTATORS},
	and Lemma~\ref{L:CONNECTIONLRADFRAME},
	we compute that
	\begin{align}
		g(\D_A (\D_{\Rad} \Lunit), X_B)
		& = \angDsquaredarg{A}{B} \upmu
			+ \angDarg{A} (\upmu \upzeta)_B
			- (\Lunit \upmu) \upchi_{AB},
				\label{E:GAMMACURVATURETHIRDRELATION} \\
		g(\D_{[\Rad, X_A]} \Lunit, X_B)
		& = \upchi([\Rad, X_A], X_B)
			= \upchi(\D_{\Rad} X_A, X_B)
				- \upchi(\D_A \Rad, X_B),
				\label{E:GAMMACURVATUREFOURTHRELATION} \\
		g(\D_A \Lunit, \D_{\Rad} X_B)
		& = - \upzeta_A g(\Lunit, \D_{\Rad} X_B)
			+ \upchi(X_A, \D_{\Rad} X_B)
			= \upzeta_A g(\D_{\Rad} \Lunit, X_B)
			+ \upchi(X_A, \D_{\Rad} X_B)
			\label{E:GAMMACURVATUREFIFTHRELATION} \\
		& = \upmu \upzeta_A \upzeta_B
			+ \upzeta_A \angdiffarg{B} \upmu 
			+ \upchi(X_A, \D_{\Rad} X_B),
			\notag \\
		\upchi(\D_A \Rad, X_B)
		& = - \upmu \upchi_A^{\ C} \upchi_{BC}
				+ \upmu \upchi_A^{\ C} \angkdoublearg{B}{C}.
					\label{E:GAMMACURVATURESIXTHRELATION}
	\end{align}
	
	The first two terms on the second line of \eqref{E:GAMMACURVATUREFIFTHRELATION}
	are the only $\upmu^{-1}-$singular terms. More precisely, using the decomposition \ref{E:ZETADECOMPOSED},
	we deduce that
	\begin{align}
		\upmu \upzeta_A \upzeta_B
		& = \frac{1}{4} 
				\upmu^{-1} 
				\angGdoublearg{\Lunit}{A} 
				\angGdoublearg{\Lunit}{B} 
				(\Rad \Psi)^2
			+  G_{(Frame)}^2
				\threemyarray
					[\upmu \Lunit \Psi]
					{\Rad \Psi}
					{\upmu \angdiff \Psi}
				\myarray
					[\Lunit \Psi]
					{\angdiff \Psi},
			 \label{E:LIERADCHIFIRSTSINGULARTERM} \\
		\upzeta_A \angdiffarg{B} \upmu
		& = - \frac{1}{2} 
					\upmu^{-1} 
					\angGdoublearg{\Lunit}{A} 
					\angdiffarg{B} \upmu
			 +   G_{(Frame)} 
					(\angdiff \upmu)	
					\myarray
						[\Lunit \Psi]
						{\angdiff \Psi}.
					\label{E:LIERADCHISECONDSINGULARTERM}
	\end{align}
	Substituting 
	\eqref{E:GAMMACURVATURESECONDRELATION}-\eqref{E:GAMMACURVATURESIXTHRELATION}
	into the right-hand side of \eqref{E:GAMMACURVATUREFIRSTRELATION},
	using \eqref{E:LIERADCHIFIRSTSINGULARTERM}-\eqref{E:LIERADCHISECONDSINGULARTERM},
	using Lemma~\ref{L:UPMUFIRSTTRANSPORT} to express the factor $\Lunit \upmu$ 
	in \eqref{E:GAMMACURVATURETHIRDRELATION}
	in terms of derivatives of $\Psi,$ 
	and symmetrizing over the indices $AB$ (since $\angLie_{\Rad} \upchi$ is symmetric),
	we arrive at the desired identity \eqref{E:PRELIMINARYRADDIRECTIONCHITRANSPORT}.
	
\end{proof}

\begin{proof}[Proof of Lemma~\ref{L:LIERADCHIJUNKTRANSPORT}]
	We first claim that $\angLie_{\Rad} \upchi_{AB}$ can be expressed as follows:
	\begin{align}  \label{E:LIERADCHINOTYETCHIJUNK}
	\angLie_{\Rad} \upchi_{AB}
	& = \angDsquaredarg{A}{B} \upmu
			+ 
			\frac{1}{2}
			\left\lbrace
				\angGdoublearg{A}{B} \Lunit \Rad \Psi
				- \angGdoublearg{\Lunit}{A} \angdiffarg{B} \Rad \Psi
				- \angGdoublearg{\Lunit}{B} \angdiffarg{A} \Rad \Psi
				- \upmu G_{\Radunit \Radunit} \angDsquaredarg{A}{B} \Psi
			\right\rbrace
			\\
	& \ \
			- \upmu \upchi_A^{\ C} \upchi_{BC}
			+ 	G_{(Frame)}
					\upchi
					\myarray
						[\upmu \Lunit \Psi]
						{\Rad \Psi}
			+ 	G_{(Frame)}
					\ginversesphere
					\upchi
					\threemyarray
						[\upmu \Lunit \Psi]
						{\Rad \Psi}
						{\upmu \angdiff \Psi}
				\notag \\
	& \ \ 
			+ 	G_{(Frame)} 
					(\angdiff \upmu)	
					\myarray
						[\Lunit \Psi]
						{\angdiff \Psi}
			+ (\angD G_{(Frame)})
				\threemyarray
					[\upmu \Lunit \Psi]
					{\Rad \Psi}
					{\upmu \angdiff \Psi}
			\notag
			\\
	& \ \ 
			+ \myarray
				[G_{(Frame)}^2]
				{G_{(Frame)}'}
				\threemyarray
					[\upmu \Lunit \Psi]
					{\Rad \Psi}
					{\upmu \angdiff \Psi}
				\myarray
					[\Lunit \Psi]
					{\angdiff \Psi}
				+ 
				G_{(Frame)}^2
					\ginversesphere
					\threemyarray
						[\upmu \Lunit \Psi]
						{\Rad \Psi}
						{\upmu \angdiff \Psi}
					\myarray
						[\Lunit \Psi]
						{\angdiff \Psi}.
						\notag 
\end{align}
	To derive \eqref{E:LIERADCHINOTYETCHIJUNK}, we first substitute
	the right-hand side of \eqref{E:GAMMACURVATURECOMPONENTEXPRESSION}
	for the curvature terms in the first line on the right-hand side of 
	\eqref{E:PRELIMINARYRADDIRECTIONCHITRANSPORT}. We observe
	that the $\upmu^{-1}-$singular terms on the second line of \eqref{E:GAMMACURVATURECOMPONENTEXPRESSION}
	exactly cancel the $\upmu^{-1}-$singular terms on the second line of
	\eqref{E:PRELIMINARYRADDIRECTIONCHITRANSPORT}.
	
	Next, we use the decomposition \eqref{E:ZETADECOMPOSED} to compute that
	\begin{align} \label{E:ANGDOFMUZETADECOMPOSED}
		\angDarg{A}(\upmu \upzeta)_B
		& = \frac{1}{2}
				\left\lbrace
					- \angGdoublearg{\Lunit}{B} \angdiffarg{A} \Rad \Psi
					+ \upmu \angGdoublearg{\Radunit}{B} \angdiffarg{A} \Lunit \Psi
					- \upmu G_{\Lunit \Radunit} \angDsquaredarg{A}{B} \Psi
					- \upmu G_{\Radunit \Radunit} \angDsquaredarg{A}{B} \Psi
				\right\rbrace
					\\
	& \ \ + \threemyarray
					[\upmu \Lunit \Psi]
					{\Rad \Psi}
					{\upmu \angdiff \Psi}
					\angD G_{(Frame)}
		+ G_{(Frame)} 
				(\angdiff \upmu)	
				\myarray
					[\Lunit \Psi]
					{\angdiff \Psi}.
					\notag
	\end{align}
	Substituting the right-hand side of \eqref{E:ANGDOFMUZETADECOMPOSED}
	for the first terms in braces on the right-hand
	side of \eqref{E:PRELIMINARYRADDIRECTIONCHITRANSPORT},
	we arrive at the desired identity \eqref{E:LIERADCHINOTYETCHIJUNK}.
	
	It remains for us to use the equation \eqref{E:LIERADCHINOTYETCHIJUNK}
	verified by $\angLie_{\Rad} \upchi$ to derive
	an analogous equation for $\angLie_{\Rad} \upchi^{(Small)}.$
	To this end, we use the identities 
	$\upchi_{AB}^{(Small)} = \upchi_{AB} - \rgeo^{-1} \gsphere_{AB},$
	$\Rad \rgeo = -1,$
	$\angLie_{\Rad} \gsphere_{AB} 
	= \angdeformarg{\Rad}{A}{B}
	= - 2 \upmu \upchi_{AB} 
		+ G_{(Frame)}
			\threemyarray
				[\upmu \Lunit \Psi]
				{\Rad \Psi}
				{\upmu \angdiff \Psi}$
	(see 
	\eqref{E:CHIJUNKDEF},
	Lemma~\ref{L:CONNECTIONBETWEENANGLIEOFGSPHEREANDDEFORMATIONTENSORS},
	\eqref{E:RADDEFORMTRFREEANG}-\eqref{E:RADDEFORMANG},
	and \eqref{E:EXACTRELATIONSZAPPLIEDTORGEO})
	to compute that
	\begin{align} \label{E:LIERADOFCHIJUNKINTERMSOFLIERADCHI}
		\angLie_{\Rad} \upchi_{AB}^{(Small)}
		& = \angLie_{\Rad} \upchi_{AB}
			- \frac{1}{\rgeo^2} \gsphere_{AB}
			+ 2 \frac{1}{\rgeo^2} \upmu \gsphere_{AB} 
			+ 2 \frac{1}{\rgeo} \upmu \upchi_{AB}^{(Small)}
			+ 
			\frac{1}{\rgeo}
			G_{(Frame)}
			\threemyarray
				[\upmu \Lunit \Psi]
				{\Rad \Psi}
				{\upmu \angdiff \Psi},
					\\
		\upmu \upchi_A^{\ C} \upchi_{BC}
			& = \upmu \upchi_A^{(Small)C} \upchi_{BC}^{(Small)}
				+ 2 \frac{1}{\rgeo} \upmu \upchi_{AB}^{(Small)}
				+ \frac{1}{\rgeo^2} \upmu \gsphere_{AB}.
				\label{E:CHITIMESCHIDECOMPOSED}
	\end{align}
	Using equation \eqref{E:LIERADOFCHIJUNKINTERMSOFLIERADCHI}
	to substitute for
	$\angLie_{\Rad} \upchi_{AB}$
	in equation \eqref{E:LIERADCHINOTYETCHIJUNK},
	and substituting the right-hand side of \eqref{E:CHITIMESCHIDECOMPOSED}
	for the first product on the second line of \eqref{E:LIERADCHINOTYETCHIJUNK},
	we finally arrive at the desired decomposition \eqref{E:LIERADCHIJUNKTRANSPORT}.
	
\end{proof}

\section{Top-order pointwise commutator estimates corresponding to the \texorpdfstring{$S_{t,u}$}{spherical} Codazzi equations}
From 
Propositions \ref{P:DEFORMATIONTENSORFRAMECOMPONENTS}
and \ref{P:COMMUTATIONCURRENTDIVERGENCEFRAMEDECOMP}
and Lemma~\ref{L:WAVEONCECOMMUTEDBASICSTRUCTURE},
it follows that the terms
$\angdiv \upchi^{(Small)}$
and $\angdiv \hat{\upchi}^{(Small)}$
are present in the commuted wave equation.
If we estimated the top-order derivatives of these
quantities in a naive fashion via transport estimates,
then our estimates would lose a derivative.
However, the next lemma shows that these quantities can be
expressed in terms of 
$\angdiff \mytr \upchi^{(Small)}$ plus
terms that do not lose derivatives.
The main point is that later in the monograph,
we will show how to use
the fully modified
quantities defined in Chapter~\ref{C:RENORMALIZEDEIKONALFUNCTIONQUANTITIES}
in conjunction with some elliptic estimates on the $S_{t,u}$
to estimate the top-order derivatives of 
$\angdiff \mytr \upchi^{(Small)}$
without losing derivatives.
The lemma is essentially a version of the 
Codazzi equations for the Riemannian manifolds
$(S_{t,u},\gsphere)$ viewed as embedded submanifolds
of the Lorentzian manifold-with-boundary $(\mathcal{M}_{\Tboot,U_0},g).$ 
In order to close our estimates, we do not need to know much
about the precise structure of the
terms on the right-hand side of
\eqref{E:CODAZZITYPERELATIONFORCHIJUNK},
except for the first one.

\begin{lemma}[\textbf{Codazzi-type identities involving} 
$\angdiv \upchi^{(Small)},$
$\angdiv \hat{\upchi}^{(Small)},$
 \textbf{and} $\angdiff \mytr \upchi_{(Small)}$]
\label{L:CODAZZITYPERELATION}
The $S_{t,u}$ one-forms 
$\angdiv \upchi^{(Small)}$ and 
$\angdiv \hat{\upchi}^{(Small)}$ verify the following equations,
where the terms on the left-hand side and the first term on the right-hand side are
exact, and the remaining terms are schematic:
\begin{align} 
	\angdiv \upchi^{(Small)},
	\,
	2 \angdiv \hat{\upchi}^{(Small)}
	& = \overbrace{\angdiff \mytr \upchi^{(Small)}}^{\mbox{precise term}}
		+ \sum_{i_1 + i_2 =1}
			\ginversesphere
			(\angD^{i_1} G_{(Frame)})
			\angD^{i_2}
			\myarray
				[\Lunit \Psi]
				{\angdiff \Psi}
		\label{E:CODAZZITYPERELATIONFORCHIJUNK} \\
		& \ \ 
		+ \sum_{i_1 + i_2 =1}
				(\angD^{i_1} \smoothfunction(\Psi))
				(\angD^{i_2} \angdiffuparg{\#} x) 
				\angdiff \Lunit_{(Small)},
		\notag 
\end{align}
and the quantities $\smoothfunction$ are smooth scalar-valued functions of $\Psi.$
\end{lemma}

\begin{remark}[\textbf{Absence of} $\upmu^{-1}$]
An important structural feature of the right-hand side of
\eqref{E:CODAZZITYPERELATIONFORCHIJUNK} 
is that there are no factors $\upmu^{-1}.$
\end{remark}

\begin{proof}
	We apply $\angD^A$ to both sides of \eqref{E:CHIJUNKINTERMSOFOTHERVARIABLES}.
	When carrying out spherical covariant differentiation with $\angD,$
	we view all lowercase Latin-indexed quantities as scalar-valued functions on $S_{t,u}.$
	Using the fact that $\angDarg{B} \angdiffuparg{A} = \angD^A \angdiffarg{B}$ when applied to functions
	and the chain rule identity $\angdiffarg{B} g_{ab} = G_{ab} \angdiffarg{B} \Psi,$
	we see that the first term on the right-hand side of the resulting identity can be written as
	\begin{align} \label{E:ANGDIVCHIJUNKMAINTERMEXPRESSION}
		& 
		\angdiffarg{B}
		\Big(
			\overbrace{
			g_{ab} (\angdiffarg{A} x^a)
	 		\angdiffuparg{A} \Lunit_{(Small)}^b
	 		}^{\angdiffarg{A} \Lunit_{(Small)}^a}
	 	\Big)
	 	+ G_{ab} (\angdiffuparg{A} \Psi) (\angdiffarg{A} x^a) \angdiffarg{B} \Lunit_{(Small)}^b
	 	- G_{ab} (\angdiffarg{B} \Psi) (\angdiffarg{A} x^a) \angdiffuparg{A} \Lunit_{(Small)}^b
	 		\\
	 	& + g_{ab} (\angLap x^a) \angdiffarg{B} \Lunit_{(Small)}^b
	 	- g_{ab} (\angD^2_{AB} x^a) \angdiffuparg{A} \Lunit_{(Small)}^b.
	 		\notag
	\end{align}
	Furthermore, the second term on the right-hand side of the resulting identity is
	\begin{align}
		\angD^A \Lambda_{AB}^{(Tan-\Psi)}
		& = \ginversesphere
				G_{(Frame)}
				\angD 
				\myarray
					[\Lunit \Psi]
					{\angdiff \Psi}
			+ \ginversesphere
				(\angD G_{(Frame)})
				\myarray
					[\Lunit \Psi]
					{\angdiff \Psi}.
	\end{align}
	From \eqref{E:TRCHIJUNKINTERMSOFOTHERVARIABLES} and \eqref{E:BIGLAMBDAGOOD}, 
	we see that the first product in parentheses in
	\eqref{E:ANGDIVCHIJUNKMAINTERMEXPRESSION} is equal to
	$\mytr \upchi^{(Small)} + \mytr  \Lambda^{(Tan-\Psi)} = \mytr \upchi^{(Small)} - \frac{1}{2} 
	\angGmixedarg{A}{A} \Lunit \Psi.$
	Combining these identities, we arrive at the desired identity \eqref{E:CODAZZITYPERELATIONFORCHIJUNK}
	for $\angdiv \upchi^{(Small)}.$ The corresponding identity for
	$\angdiv \hat{\upchi}^{(Small)}$ follows easily from the identity
	$\angdiv \hat{\upchi}^{(Small)} = \angdiv \upchi^{(Small)} - \frac{1}{2} \angdiff \mytr \upchi^{(Small)}.$
\end{proof}

We now provide a commuted version of Lemma~\ref{L:CODAZZITYPERELATION}.

\begin{lemma}[\textbf{Commutator estimates for the Codazzi equation}]
Let $1 \leq N \leq 23$ be an integer. 
Under the small-data and bootstrap assumptions 
of Sects.~\ref{S:PSISOLVES}-\ref{S:C0BOUNDBOOTSTRAP},
if $\varepsilon$ is sufficiently small, 
then the following pointwise commutator estimates hold
on $\mathcal{M}_{\Tboot,U_0}:$
\begin{align} \label{E:DIVCHIJUNKANGDIFFTRCHIJUNKHIGHERORDERCOMMMUTOR}
	& \left|
		\angLie_{\mathscr{Z}}^{N-1} \angdiv \upchi^{(Small)}
		- \angdiff \mathscr{Z}^{N-1} \mytr \upchi^{(Small)}
	\right|,
		\\	
	& \left|
			\angLie_{\mathscr{Z}}^{N-1} \angdiv \hat{\upchi}^{(Small)}
			- \frac{1}{2} \angdiff \mathscr{Z}^{N-1} \mytr \upchi^{(Small)}
		\right|,
		\notag \\
	& \left|
			\angdiv \angLie_{\mathscr{Z}}^{N-1} \upchi^{(Small)}
			- \angdiff \mathscr{Z}^{N-1} \mytr \upchi^{(Small)}
		\right|,
		\notag \\	
	&	\left|
			\angdiv \angfreeLietwoarg{\mathscr{Z}}{N-1} \hat{\upchi}^{(Small)}
			- \frac{1}{2} \angdiff \mathscr{Z}^{N-1} \mytr \upchi^{(Small)}
		\right|
		\notag \\
	& \lesssim 
		\frac{1}{(1 + t)^2} 
		\left| 
			\fourmyarray[\rgeo \Lunit \mathscr{Z}^{\leq N} \Psi ]
				{\Rad \mathscr{Z}^{\leq N} \Psi}
				{\rgeo \angdiff \mathscr{Z}^{\leq N} \Psi}
				{\mathscr{Z}^{\leq N} \Psi}
		\right|	
		+ \frac{1}{(1 + t)^3}
			\left|
				\myarray[\mathscr{Z}^{\leq N} (\upmu - 1)]
					{\sum_{a=1}^3 \rgeo |\mathscr{Z}^{\leq N} \Lunit_{(Small)}^a|} 
			\right|.
			\notag
\end{align}

\end{lemma}

\begin{proof}
	We prove the estimate \eqref{E:DIVCHIJUNKANGDIFFTRCHIJUNKHIGHERORDERCOMMMUTOR}
	for $\angLie_{\mathscr{Z}}^{N-1} \angdiv \upchi^{(Small)}$ in detail.
	The estimate for $\angLie_{\mathscr{Z}}^{N-1} \angdiv \hat{\upchi}^{(Small)}$
	can be proved by a nearly identical argument.
	The estimate \eqref{E:DIVCHIJUNKANGDIFFTRCHIJUNKHIGHERORDERCOMMMUTOR}
	for $\angdiv \angLie_{\mathscr{Z}}^{N-1} \upchi^{(Small)}$
	then follows from the estimate for $\angLie_{\mathscr{Z}}^{N-1} \angdiv \upchi^{(Small)},$
	the commutator estimate \eqref{E:COMMUTINGLIEDERIVATIVESTHROUGHANGDIV}
	with $\upchi^{(Small)}$ in the role of $\xi$
	and $N-1$ in the role of $N,$
	and inequalities \eqref{E:POINTWISEESTIMATESFORCHIJUNKINTERMSOFOTHERVARIABLES} 
	and \eqref{E:C0BOUNDCRUCIALEIKONALFUNCTIONQUANTITIES}.
	The final estimate \eqref{E:DIVCHIJUNKANGDIFFTRCHIJUNKHIGHERORDERCOMMMUTOR}
	for $\angdiv \angfreeLietwoarg{\mathscr{Z}}{N-1} \hat{\upchi}^{(Small)}$
	can be proved similarly with the help of 
	the commutator estimate \eqref{E:COMMUTINGTRACEFREELIEDERIVATIVESTHROUGHANGDIV}.
	
	To derive the estimate for $\angLie_{\mathscr{Z}}^{N-1} \angdiv \upchi^{(Small)},$
	we first use equation \eqref{E:CODAZZITYPERELATIONFORCHIJUNK} and Lemma~\ref{L:LANDRADCOMMUTEWITHANGDIFF}
	to deduce that
	\begin{align} 
	\left| 
		\angLie_{\mathscr{Z}}^{N-1} \angdiv \upchi^{(Small)}
		- \angdiff \mathscr{Z}^{N-1} \mytr \upchi^{(Small)}
	\right|
	& \lesssim 
		\left|
				\angLie_{\mathscr{Z}}^{N-1}
				\sum_{i_1 + i_2 =1}
				\ginversesphere
				(\angD^{i_1} G_{(Frame)})
				\angD^{i_2}
				\myarray
					[\Lunit \Psi]
					{\angdiff \Psi}
			\right|
		\label{E:CODAZZICOMMUTATORTERMFIRSTESTIMATE} \\
		& \ \ 
		+ \left|
				\angLie_{\mathscr{Z}}^{N-1}
				\sum_{i_1 + i_2 =1}
				\rgeo^{-1}
				(\angD^{i_1} \smoothfunction(\Psi))
				(\angD^{i_2} \angdiffuparg{\#} x) 
				\angdiff (\rgeo \Lunit_{(Small)})
			\right|.
		\notag
		\end{align}
We now apply the Leibniz rule to the terms on the right-hand side of \eqref{E:CODAZZICOMMUTATORTERMFIRSTESTIMATE}.
We bound the terms $\angLie_{\mathscr{Z}}^M \ginversesphere$
by using the estimates of
Lemma~\ref{L:POINTWISEBOUNDSDERIVATIVESOFANGULARDEFORMATIONTENSORS}.
We bound the terms $\angLie_{\mathscr{Z}}^M G_{(Frame)}$ and
$\angLie_{\mathscr{Z}}^M \angD G_{(Frame)}$ by using
Lemma~\ref{L:LANDRADCOMMUTEWITHANGDIFF},
\eqref{E:COMMUTATORESTIMATESVECTORFIELDSACTINGONANGDTENSORS},		
\eqref{E:FUNCTIONPOINTWISEANGDINTERMSOFANGLIEO},
\eqref{E:ONEFORMANGDINTERMSOFROTATIONALLIE},
\eqref{E:TYPE02TENSORANGDINTERMSOFROTATIONALLIE},		
and the estimates of Lemma~\ref{L:POINTWISEESTIMATESGFRAMEINTERMSOFOTHERQUANTITIES}.
To bound $\myarray[\mathscr{Z}^M \Lunit \Psi] {\angLie_{\mathscr{Z}}^M \angdiff \Psi}$
and $\angLie_{\mathscr{Z}}^M \angD \myarray[\Lunit \Psi] {\angdiff \Psi},$
we use Lemma~\ref{L:LANDRADCOMMUTEWITHANGDIFF},
\eqref{E:ANGDSQUAREDFUNCTIONPOINTWISEINTERMSOFANGDIFFROTATIONS},
\eqref{E:ANGDSQUAREDLIEZNCOMMUTATORACTINGONFUNCTIONSSPOINTWISE},
and
Lemma~\ref{L:AVOIDINGCOMMUTING}. 
We bound the terms $\angLie_{\mathscr{Z}}^M \rgeo^{-1}$
by using \eqref{E:ZNAPPLIEDTORGEOISNOTTOOLARGE}.
We bound the terms
$\mathscr{Z}^M \smoothfunction(\Psi)$ by using 
the bootstrap assumptions \eqref{E:PSIFUNDAMENTALC0BOUNDBOOTSTRAP}
to deduce that $\left| \mathscr{Z}^M \smoothfunction(\Psi) \right| \lesssim \left| \mathscr{Z}^{\leq M} \Psi \right| + 1.$
We bound the terms
$\angLie_{\mathscr{Z}}^M \angdiff \smoothfunction(\Psi)$ by 
using Lemma~\ref{L:LANDRADCOMMUTEWITHANGDIFF}
to deduce that $\left| \angLie_{\mathscr{Z}}^M \angdiff \smoothfunction(\Psi) \right| \lesssim \left| \angdiff \mathscr{Z}^{\leq M} \Psi \right|.$
We bound the terms $\angLie_{\mathscr{Z}}^M  \angdiffuparg{\#} x$ by using
Lemma~\ref{L:LANDRADCOMMUTEWITHANGDIFF} and the estimates
\eqref{E:POINTWISEBOUNDPROJECTEDLIEDERIVATIVESANGDIFFCOORDINATEX}, 
\eqref{E:LOWERORDERPOINTWISEBOUNDPROJECTEDLIEDERIVATIVESANGDIFFCOORDINATEX},
\eqref{E:CRUDEPOINTWISEBOUNDSDERIVATIVESOFANGULARDEFORMATIONTENSORS},
and 
\eqref{E:CRUDELOWERORDERC0BOUNDDERIVATIVESOFANGULARDEFORMATIONTENSORS}.
We bound the terms 
$\angLie_{\mathscr{Z}}^M \angD \angdiffuparg{\#} x
= \angLie_{\mathscr{Z}}^M (\ginversesphere \angD \angdiff x)$ 
by using
Lemma~\ref{L:LANDRADCOMMUTEWITHANGDIFF},
\eqref{E:COMMUTATORESTIMATESVECTORFIELDSACTINGONANGDTENSORS}
and		
\eqref{E:ONEFORMANGDINTERMSOFROTATIONALLIE}
with $\angdiff x^i$ in the role of $\xi,$
\eqref{E:POINTWISEBOUNDPROJECTEDLIEDERIVATIVESANGDIFFCOORDINATEX}, 
\eqref{E:LOWERORDERPOINTWISEBOUNDPROJECTEDLIEDERIVATIVESANGDIFFCOORDINATEX},
\eqref{E:CRUDEPOINTWISEBOUNDSDERIVATIVESOFANGULARDEFORMATIONTENSORS},
and 
\eqref{E:CRUDELOWERORDERC0BOUNDDERIVATIVESOFANGULARDEFORMATIONTENSORS}.
We bound the terms
$\angLie_{\mathscr{Z}}^M \angdiff (\rgeo \Lunit_{(Small)})$
by using Lemma~\ref{L:LANDRADCOMMUTEWITHANGDIFF},
\eqref{E:ZNAPPLIEDTORGEOISNOTTOOLARGE},
\eqref{E:FUNCTIONPOINTWISEANGDINTERMSOFANGLIEO},
and
\eqref{E:C0BOUNDCRUCIALEIKONALFUNCTIONQUANTITIES}.

\end{proof}


\chapter[Pointwise Estimates for the Easy Error Integrands]{Pointwise Estimates for the Easy Error Integrands and 
Identification of the Difficult Error Integrands Corresponding to the Commuted Wave Equation}
\label{C:POINTWISEBOUNDSFOREASYERRORINTEGRANDS}
\thispagestyle{fancy}
Recall that we derived our energy-flux identities for solutions 
to $\upmu \square_{g(\Psi)} \Psi = 0$
in Prop.~\ref{P:DIVTHMWITHCANCELLATIONS}.
Furthermore, similar identities hold for its higher-order analogs $\mathscr{Z}^N \Psi,$
which verify $\upmu \square_{g(\Psi)} (\mathscr{Z}^N \Psi) = \inhomarg{\mathscr{Z}^N},$ 
where $\mathscr{Z}=\lbrace \rgeo \Lunit, \Rad, \Rot_{(1)}, \Rot_{(2)}, \Rot_{(3)}\rbrace$ is the set of commutation vectorfields.
In Chapter~\ref{C:POINTWISEBOUNDSFOREASYERRORINTEGRANDS}, our main goal is to derive pointwise estimates for the ``easy factors'' 
in the integrands on the right-hand side of the identities of Prop.~\ref{P:DIVTHMWITHCANCELLATIONS}.
These pointwise estimates play an important role in Chapter~\ref{C:ERRORTERMSOBOLEV}, 
where we use them to bound the corresponding error integrals 
in terms of the fundamental $L^2-$controlling quantities defined in Chapter~\ref{C:L2COERCIVENESS}.
Most of the integrand factors are easy to pointwise bound,
but some terms found in $\waveinhom := \inhomarg{\mathscr{Z}^N}$ 
(see Prop.~\ref{P:DIVTHMWITHCANCELLATIONS})
are not; we postpone the analysis
of the difficult terms until Ch.~\ref{C:POINTWISEESTIMATESDIFFICULTERRORINTEGRANDS}.
More precisely, in Sects.~\ref{S:PRELIMINARYANALYSIS}-\ref{S:PROOFOFCOROLLARYPOINTWISEESTIMATESFOREASYCOMMUTATORTERMS},
we decompose $\inhomarg{\mathscr{Z}^N}$
and identify those terms that lead to error integrals that are easy to bound.
This analysis is among the most important of the entire monograph, 
\emph{for our pointwise estimates for the few important terms in} 
$\inhomarg{\mathscr{Z}^N},$ 
\emph{which are provided in Chapter~\ref{C:POINTWISEESTIMATESDIFFICULTERRORINTEGRANDS}, affect the degree to which our a priori $L^2$ estimates can degenerate with respect to} 
$\upmu_{\star}^{-1}.$ 
In Sect.~\ref{S:POINTWISEESTIMATESFORBASICENERGYERRORINTEGRANDS},
we derive pointwise estimates for the remaining
non-$\waveinhom$ integrand factors appearing in 
Prop.~\ref{P:DIVTHMWITHCANCELLATIONS}, which are relatively easy to bound.
Finally, in Sect.~\ref{S:POINTWISEFORCLOSINGTHEELLIPTICESTIMATES}, 
we derive some simple pointwise estimates that we need to close 
our top-order elliptic estimates.

\section{Preliminary analysis and the definition of \texorpdfstring{$Harmless^{\leq N}$}{harmless} terms}
\label{S:PRELIMINARYANALYSIS}
We start by proving a preliminary lemma. Roughly speaking, 
when combined with Cor.~\ref{C:POINTWISEESTIMATESFOREASYCOMMUTATORTERMS},
the lemma shows that difficult inhomogeneous terms in the $N-$times-commuted wave equation
can only arise from repeatedly differentiating some of the inhomogeneous terms that arise after 
a single commutation. In particular, the difficult terms involve certain top-order derivatives
of the eikonal function quantities
$\upmu$
and
$\upchi^{(Small)}.$

\begin{lemma}[\textbf{Basic structure of the inhomogeneous terms in the $N-$times commuted wave equation}]
	\label{L:WAVENTIMESCOMMUTEDBASICSTRUCTURE}
	Assume that $\Psi$ verifies the wave equation	
	\begin{align}
		\square_{g(\Psi)} \Psi = 0,
	\end{align}
	and let $0 \leq N \leq 24$ be an integer.
	Let $\mathscr{Z}^N$ be an $N^{th}$ order commutation vectorfield operator
	consisting of iterated vectorfields belonging to the commutation set
	$\mathscr{Z} = \lbrace \rgeo \Lunit, \Rad, \Rot_{(1)}, \Rot_{(2)}, \Rot_{(3)} \rbrace.$
	Assume that $\mathscr{Z}^N$ is of the form $\mathscr{Z}^N = \mathscr{Z}^{N-1} Z,$
	and let $\Jcurrent{Z}[\Psi]$ be the commutation current \eqref{E:COMMUTATIONCURRENT}. 
	Under the small-data and bootstrap assumptions 
	of Sects.~\ref{S:PSISOLVES}-\ref{S:C0BOUNDBOOTSTRAP},
	if $\varepsilon$ is sufficiently small, then
	$\mathscr{Z}^N \Psi$ verifies the inhomogeneous wave equation
	\begin{align} \label{E:WAVENTIMESCOMMUTEDBASICSTRUCTURE}
		\upmu \square_{g(\Psi)} (\mathscr{Z}^N \Psi)
		& = \inhomarg{\mathscr{Z}^N},
			\\
		\inhomarg{\mathscr{Z}^N} 
		& = \mathscr{Z}^{N-1}(\upmu \D_{\alpha} {\Jcurrent{Z}^{\alpha}[\Psi]})
			+ \inhomarg{\mathscr{Z}^N}_{(Eikonal-Low)},
			\label{E:NTIMESCOMMUTEDWAVEINHOMOGENEOUSTERMFIRSTSPLITTING}
	\end{align}
	where $\inhomarg{\mathscr{Z}^N}_{(Eikonal-Low)} = 0$ if $N=1$ and otherwise 
	$\inhomarg{\mathscr{Z}^N}_{(Eikonal-Low)}$ verifies the following 
	pointwise inequality on $\mathcal{M}_{\Tboot,U_0}:$
	\begin{align}
		\left|\inhomarg{\mathscr{Z}^N}_{(Eikonal-Low)} \right|
		& \lesssim 	\mathop{\sum_{N_1 + N_2 + N_3 \leq N-1}}_{N_1, N_2 \leq N - 2} 
								\sum_{Z_1, Z_2 \in \mathscr{Z}}
									\left(
										\left| 
											\mathscr{Z}^{N_1} \mytr  \angdeform{Z_1} 
										\right|
									+ 1	
									\right)
									\left| 
										\mathscr{Z}^{N_2} (\upmu \D_{\alpha} \Jcurrent{Z_2}^{\alpha}[\mathscr{Z}^{N_3} \Psi]) 
									\right|.
			\label{E:LOWERORDERINHOMOGENEOUSTERMSFIRSTPOINTWISE}
\end{align}
\end{lemma}

\begin{proof}
We iterate the identity \eqref{E:WAVEONCECOMMUTEDBASICSTRUCTURE}. 
Clearly the first term on the right-hand side of \eqref{E:NTIMESCOMMUTEDWAVEINHOMOGENEOUSTERMFIRSTSPLITTING}
arises when $\mathscr{Z}^{N-1}$ falls on the first term on the right-hand side of \eqref{E:WAVEONCECOMMUTEDBASICSTRUCTURE}.
Using the observation made in Remark~\ref{R:POTENTIALLYDANGEROUSTERMISNOTDANGEROUS}
and the estimate $\left\| \mathscr{Z}^{\leq 11} \mytr  \angdeform{Z} \right\|_{C^0(\Sigma_t^u)} \lesssim 1$
(that is, \eqref{E:CRUDELOWERORDERC0BOUNDDERIVATIVESOFANGULARDEFORMATIONTENSORS}),
we deduce that all of the remaining terms that arise are in magnitude
$\lesssim$ the right-hand side of \eqref{E:LOWERORDERINHOMOGENEOUSTERMSFIRSTPOINTWISE} as desired.

\end{proof}

The term $\mathscr{Z}^{N-1}(\upmu \D_{\alpha} \Jcurrent{Z}^{\alpha}[\Psi])$
on the right-hand side of \eqref{E:NTIMESCOMMUTEDWAVEINHOMOGENEOUSTERMFIRSTSPLITTING}
is difficult to bound
The two sums on the right-hand side of \eqref{E:LOWERORDERINHOMOGENEOUSTERMSFIRSTPOINTWISE} are 
lower-order in terms of number of derivatives of the eikonal function quantities and are relatively easy to bound. 

In our analysis of 
$\mathscr{Z}^{N-1}(\upmu \D_{\alpha} {\Jcurrent{Z}^{\alpha}[\Psi]})$
and
$\inhomarg{\mathscr{Z}^N}_{(Eikonal-Low)},$
the vast majority of the terms that we encounter have a negligible effect on the dynamics. 
We call such terms ``harmless.'' We now give a precise definition of what we mean by ``harmless.''

\begin{definition}[\textbf{Harmless terms}] 
\label{D:HARMLESSTERMS}
Let $0 \leq N \leq 24$ be an integer.
A ``$Harmless^{\leq N}$'' term is any term 
such that under the small-data and bootstrap assumptions 
of Sects.~\ref{S:PSISOLVES}-\ref{S:C0BOUNDBOOTSTRAP},
if $\varepsilon$ is sufficiently small,
then on the spacetime domain $\mathcal{M}_{\Tboot,U_0},$
it verifies the following bound:
	\begin{align} \label{E:HARMLESSTERMS}
		\left|
			Harmless^{\leq N} 
		\right|
		& \lesssim
			\left|\left\lbrace \Lunit + \frac{1}{2} \mytr \upchi \right\rbrace \mathscr{Z}^{\leq N} \Psi \right|
			+ \frac{\ln(\myexp + t)}{(1 + t)^2}
				\left|
					\fourmyarray
						[\rgeo (1 + \upmu) \Lunit \mathscr{Z}^{\leq N} \Psi]
						{\Rad \mathscr{Z}^{\leq N} \Psi}
						{\rgeo (1 + \upmu) \angdiff \mathscr{Z}^{\leq N} \Psi}
					 	{\mathscr{Z}^{\leq N} \Psi}
					\right|	
				\\
		& \ \ + \varepsilon 
				\frac{1}{(1 + t)^2}
				\left| 
					\sevenmyarray
						[\Lunit \mathscr{Z}^{\leq N} \upmu]
							{\rgeo \sum_{a=1}^3 |\Lunit \mathscr{Z}^{\leq N} \Lunit_{(Small)}^a|}
							{\rgeo^2 \angLie_{\Lunit} \angLie_{\mathscr{Z}}^{\leq N-1} \upchi^{(Small)}}
							{\angLie_{\Lunit} (\rgeo^2 \angLie_{\mathscr{Z}}^{\leq N-1} \upchi^{(Small)\#})}
							{\Lunit (\rgeo^2 \mathscr{Z}^{\leq N-1} \mytr \upchi^{(Small)})}
							{\rgeo^2 \angLie_{\Lunit} \angLie_{\mathscr{Z}}^{\leq N-1} \hat{\upchi}^{(Small)}}
							{\angLie_{\Lunit} (\rgeo^2 \angLie_{\mathscr{Z}}^{\leq N-1} \hat{\upchi}^{(Small)\#})}
				\right|
				\notag
				\\
		& \ \ + 
			\varepsilon 
				\frac{\ln(\myexp + t)}{(1 + t)^3}
				\left| 
					\sevenmyarray
						[\mathscr{Z}^{\leq N} (\upmu - 1)]
						{\rgeo \sum_{a=1}^3 |\mathscr{Z}^{\leq N-1} \Lunit_{(Small)}^a|}
						{\rgeo^2 \angLie_{\mathscr{Z}}^{\leq N-1} \upchi^{(Small)}}
						{\rgeo^2 \angLie_{\mathscr{Z}}^{\leq N-1} \upchi^{(Small)\#}}
						{\rgeo^2 \mathscr{Z}^{\leq N-1} \mytr \upchi^{(Small)}}
						{\rgeo^2 \angLie_{\mathscr{Z}}^{\leq N-1} \hat{\upchi}^{(Small)}}
						{\rgeo^2 \angLie_{\mathscr{Z}}^{\leq N-1} \hat{\upchi}^{(Small)\#}}
				\right|
				\notag 
	\end{align}
	
\end{definition}

\begin{remark}[\textbf{The main features of the terms} $Harmless^{\leq N}$]
	The important features of $Harmless^{\leq N}$ terms are the following.
	\begin{itemize}
		\item Various $L^2$ norms of the terms 
		on the right-hand side of \eqref{E:HARMLESSTERMS}
		can be controlled
		\emph{without the need for introducing modified versions and/or elliptic estimates}.
		\item The $t$ weight factors on the right-hand side of \eqref{E:HARMLESSTERMS}
		are strong enough that these terms have 
		a negligible effect on the dynamics in $\mathcal{M}_{\Tboot,U_0}.$
	\end{itemize}
\end{remark}

In the next proposition, we identify the difficult
products in the main term $\mathscr{Z}^{N-1}(\upmu \D_{\alpha} \Jcurrent{Z}^{\alpha}[\Psi])$ 
on the right-hand side of \eqref{E:NTIMESCOMMUTEDWAVEINHOMOGENEOUSTERMFIRSTSPLITTING}.
The difficult products depend on the innermost commutation vectorfield $Z.$

\begin{proposition}[\textbf{Identification of the key difficult error term factors}]
\label{P:IDOFKEYDIFFICULTENREGYERRORTERMS}
Let $\Jcurrent{Z}[\Psi]$ 
be the commutation current \eqref{E:COMMUTATIONCURRENT}.
Under the small-data and bootstrap assumptions 
of Sects.~\ref{S:PSISOLVES}-\ref{S:C0BOUNDBOOTSTRAP},
if $\varepsilon$ is sufficiently small,
then the following pointwise estimates hold on $\mathcal{M}_{\Tboot,U_0}:$
\begin{subequations}
\begin{align}
	\mathscr{Z}^{N-1}(\upmu \D_{\alpha} {\Jcurrent{\Rad}^{\alpha}[\Psi]})
	& = (\Rad \Psi) \angLap \mathscr{Z}^{N-1} \upmu
			- (\upmu \angdiffuparg{\#} \Psi) \cdot (\upmu \angdiff \mathscr{Z}^{N-1} \mytr \upchi^{(Small)})
			+ Harmless^{\leq N},
			\label{E:RADISTHEFIRSTCOMMUTATORIMPORTANTTERMS} \\
	\mathscr{Z}^{N-1}(\upmu \D_{\alpha} {\Jcurrent{\rgeo \Lunit}^{\alpha}[\Psi]})
	& = \rgeo (\angdiffuparg{\#} \Psi) \cdot (\upmu \angdiff \mathscr{Z}^{N-1} \mytr \upchi^{(Small)})
			+ Harmless^{\leq N},
			\label{E:RGEOLISTHEFIRSTCOMMUTATORIMPORTANTTERMS} \\
	\mathscr{Z}^{N-1}(\upmu \D_{\alpha} {\Jcurrent{\Rot_{(l)}}^{\alpha}[\Psi]})
	& = (\Rad \Psi) \Rot_{(l)} \mathscr{Z}^{N-1} \mytr \upchi^{(Small)}
			+ \RotRadcomponent{l}(\angdiffuparg{\#} \Psi) \cdot (\upmu \angdiff \mathscr{Z}^{N-1} \mytr \upchi^{(Small)})
				\label{E:ROTISTHEFIRSTCOMMUTATORIMPORTANTTERMS} \\
	& \ \ + Harmless^{\leq N}.
		\notag
\end{align}
\end{subequations}
\end{proposition}

We provide the proof of Prop.~\ref{P:IDOFKEYDIFFICULTENREGYERRORTERMS} 
in Sect.~\ref{S:PROOFOFPROPOSITIONPIDOFKEYDIFFICULTENREGYERRORTERMS}.

In the next corollary,
we show that the second term on the right-hand side of
\eqref{E:NTIMESCOMMUTEDWAVEINHOMOGENEOUSTERMFIRSTSPLITTING}
is $Harmless^{\leq N}.$

\begin{corollary}[\textbf{Pointwise estimates for} $\inhomarg{\mathscr{Z}^N}_{(Eikonal-Low)}$]
	\label{C:POINTWISEESTIMATESFOREASYCOMMUTATORTERMS}
	Let $0 \leq N \leq 24$ be an integer.
	Under the small-data and bootstrap assumptions 
	of Sects.~\ref{S:PSISOLVES}-\ref{S:C0BOUNDBOOTSTRAP},
	if $\varepsilon$ is sufficiently small, then
	the following pointwise estimates hold 
	for the term $\inhomarg{\mathscr{Z}^N}_{(Eikonal-Low)}$
	from \eqref{E:LOWERORDERINHOMOGENEOUSTERMSFIRSTPOINTWISE}
	on $\mathcal{M}_{\Tboot,U_0}:$
 \begin{align}
		\inhomarg{\mathscr{Z}^N}_{(Eikonal-Low)}
		& = Harmless^{\leq N},
			\label{E:LOWERORDERDERIVATIVESOFCOMMUTATIONCURRENTSAREHARMLESS}
	\end{align}
	where $Harmless^{\leq N}$ terms are defined in
	Def.~\ref{D:HARMLESSTERMS}.
\end{corollary}

We provide the proof of Cor.~\ref{C:POINTWISEESTIMATESFOREASYCOMMUTATORTERMS}
in Sect.~\ref{S:PROOFOFCOROLLARYPOINTWISEESTIMATESFOREASYCOMMUTATORTERMS}.

The main point of the next corollary is that when a few factors of $\rgeo \Lunit$ are present
in the differential operator $\mathscr{Z}^N$ from \eqref{E:WAVENTIMESCOMMUTEDBASICSTRUCTURE},
we can completely avoid the use of modified quantities and elliptic estimates in deriving our 
a priori $L^2$ estimates.
Thus, in order to close our $L^2$ estimates, 
we only have to perform a detailed analysis for a 
handful of cases in which $\mathscr{Z}^N$ contains at most one factor of $\rgeo \Lunit.$
Furthermore, corresponding to each of these handful of cases, there are only one or two terms 
that are difficult to estimate.

\begin{corollary}[\textbf{Reduction of the $L^2$ analysis to essentially estimates for 
pure spatial commutation vectorfield operators}]
\label{C:REDUCTIONOFPROOFTOPURESPATIALCOMMUTATORS}
Assume that the differential operator $\mathscr{Z}^N$ contains precisely one factor of $\rgeo \Lunit.$
Under the small-data and bootstrap assumptions 
of Sects.~\ref{S:PSISOLVES}-\ref{S:C0BOUNDBOOTSTRAP},
if $\varepsilon$ is sufficiently small, then
the term $\inhomarg{\mathscr{Z}^N}$ from the right-hand side of \eqref{E:WAVENTIMESCOMMUTEDBASICSTRUCTURE}
verifies one of the following two estimates on $\mathcal{M}_{\Tboot,U_0}:$
\begin{align} \label{E:ONEGOODCOMMUTATORESTIMATEEASIEST}
	\inhomarg{\mathscr{Z}^N} 
	& = Harmless^{\leq N},
		\\
	\inhomarg{\mathscr{Z}^N}
	& = \rgeo 
			(\angdiff \Psi^{\#}) 
			\cdot(\upmu \angdiff \mathscr{S}^{N-1} \mytr \upchi^{(Small)})
		+ Harmless^{\leq N},
		\label{E:ONEGOODCOMMUTATORESTIMATE}
\end{align}
where $\mathscr{S}^{N-1}$ is an $(N-1)^{st}$ order pure spatial commutation vectorfield operator
(see Definition \eqref{E:DEFSETOFSPATIALCOMMUTATORVECTORFIELDS}).

Furthermore, if the differential operator $\mathscr{Z}^N$ contains at least two factors of $\rgeo \Lunit,$ then
\begin{align} \label{E:TWOGOODCOMMUTATORSESTIMATE}
	\inhomarg{\mathscr{Z}^N}
	& = Harmless^{\leq N}.
\end{align}

\end{corollary}

\begin{proof}
	We first assume that $\mathscr{Z}^N$ contains precisely one factor of $\rgeo \Lunit$
	and prove \eqref{E:ONEGOODCOMMUTATORESTIMATEEASIEST} and \eqref{E:ONEGOODCOMMUTATORESTIMATE}.
	We split the argument into three cases:
	\textbf{i)} $\mathscr{Z}^N = \mathscr{Z}^{N-1} \rgeo \Lunit,$
	\textbf{ii)} $\mathscr{Z}^N = \mathscr{Z}^{N-1} \Rad,$
	\textbf{iii)} $\mathscr{Z}^N = \mathscr{Z}^{N-1} \Rot_{(l)}.$
	In case \textbf{i)}, the estimate \eqref{E:ONEGOODCOMMUTATORESTIMATE} follows from
	\eqref{E:NTIMESCOMMUTEDWAVEINHOMOGENEOUSTERMFIRSTSPLITTING},
	\eqref{E:RGEOLISTHEFIRSTCOMMUTATORIMPORTANTTERMS},
	and Cor.~\ref{C:POINTWISEESTIMATESFOREASYCOMMUTATORTERMS}.
	In case \textbf{ii)}, $\mathscr{Z}^{N-1}$
	must contain a factor of $\rgeo \Lunit.$ We will show that 
	$\inhomarg{\mathscr{Z}^N} = Harmless^{\leq N}$ in this case.
	By \eqref{E:NTIMESCOMMUTEDWAVEINHOMOGENEOUSTERMFIRSTSPLITTING},
	\eqref{E:RADISTHEFIRSTCOMMUTATORIMPORTANTTERMS},
	and Cor.~\ref{C:POINTWISEESTIMATESFOREASYCOMMUTATORTERMS},
	$\inhomarg{\mathscr{Z}^N}$ is equal to the first two terms on the right-hand
	side of \eqref{E:RADISTHEFIRSTCOMMUTATORIMPORTANTTERMS} $+ Harmless^{\leq N}.$
	Using inequalities \eqref{E:FUNCTIONPOINTWISEANGDINTERMSOFANGLIEO} 
	and \eqref{E:ANGLAPFUNCTIONPOINTWISEINTERMSOFROTATIONS}
	and the bootstrap assumption $\| \Rad \Psi \|_{C^0(\Sigma_t^u)} \leq \varepsilon (1 + t)^{-1}$
	(that is, \eqref{E:PSIFUNDAMENTALC0BOUNDBOOTSTRAP}),
	we deduce that the first term on the right-hand
	side of \eqref{E:RADISTHEFIRSTCOMMUTATORIMPORTANTTERMS}
	is bounded in magnitude by
	\begin{align} \label{E:RADFIRSTWITHANRGEOLFACTORWANTSTOBEHARMLESS}
		\lesssim
		\varepsilon \frac{1}{(1 + t)^3}
		\sum_{l,m=1}^3
		\left|
			\Rot_{(l)} \Rot_{(m)} \mathscr{Z}^{N-1} \upmu
		\right|.
	\end{align}
	Next, we repeatedly use inequalities
	\eqref{E:RGEOLORRADZNCOMMUTATORACTINGONFUNCTIONSPOINTWISE}
	and \eqref{E:ROTZNCOMMUTATORACTINGONFUNCTIONSPOINTWISE}
	with $f = \upmu - 1,$
	\eqref{E:C0BOUNDCRUCIALEIKONALFUNCTIONQUANTITIES}, 
	and \eqref{E:BOUNDINGFRAMEDERIVATIVESINTERMSOFZNF}
	to commute the factor $Z := \rgeo \Lunit$ in $\mathscr{Z}^{N-1}$ all the way in front, 
	which allows us to bound
	the right-hand side of \eqref{E:RADFIRSTWITHANRGEOLFACTORWANTSTOBEHARMLESS} by
	\begin{align} \label{E:COMMUTINGDONERADFIRSTWITHANRGEOLFACTORWANTSTOBEHARMLESS}
		& 
		\lesssim
		\varepsilon \frac{1}{(1 + t)^2}
		\left|
			\Lunit \mathscr{Z}^{\leq N} \upmu
		\right|
		+ \varepsilon
			\frac{\ln(\myexp + t)}{(1 + t)^2}
			\left| 
				\fourmyarray[\rgeo \Lunit \mathscr{Z}^{\leq N-1} \Psi]
					{\Rad \mathscr{Z}^{\leq N-1} \Psi}
					{\rgeo \angdiff \mathscr{Z}^{\leq N-1} \Psi}
					{\mathscr{Z}^{\leq N-1} \Psi}
			\right|
				\\
& \ \ + 
			\varepsilon
			\frac{\ln(\myexp + t)}{(1 + t)^3}
			\left|
				\myarray[\mathscr{Z}^{\leq N} (\upmu - 1)]
					{\sum_{a=1}^3 \rgeo |\mathscr{Z}^{\leq N} \Lunit_{(Small)}^a|} 
			\right|.
			\notag
	\end{align}
	Referring to Def.~\ref{D:HARMLESSTERMS}, we see that
	the right-hand side of \eqref{E:COMMUTINGDONERADFIRSTWITHANRGEOLFACTORWANTSTOBEHARMLESS}
	is $=Harmless^{\leq N}$ as desired. 
	Using a similar argument, we also deduce that the second term on the right-hand
	side of \eqref{E:RADISTHEFIRSTCOMMUTATORIMPORTANTTERMS} is $=Harmless^{\leq N}$ as desired. 
	This completes the proof in \textbf{ii)}.
	In case \textbf{iii)}, 
	we use
	arguments similar to the ones we used in the case \textbf{ii)},
	except that we use \eqref{E:ROTISTHEFIRSTCOMMUTATORIMPORTANTTERMS}
	in place of \eqref{E:RADISTHEFIRSTCOMMUTATORIMPORTANTTERMS}.
	We conclude that $\inhomarg{\mathscr{Z}^N} = Harmless^{\leq N}$ in case \textbf{iii)}.
	We have thus proved \eqref{E:ONEGOODCOMMUTATORESTIMATEEASIEST} and \eqref{E:ONEGOODCOMMUTATORESTIMATE}.
	
	We now prove \eqref{E:TWOGOODCOMMUTATORSESTIMATE}. Under the present assumptions, in 
	the cases \textbf{ii)} and \textbf{iii)} from above, 
	$\mathscr{Z}^{N-1}$ must contain two factors of $\rgeo \Lunit.$
	The desired estimate \eqref{E:TWOGOODCOMMUTATORSESTIMATE} 
	follows from the previously proven results for cases \textbf{ii)} and \textbf{iiii)}
	since it was enough for $\mathscr{Z}^{N-1}$ to contain a single factor of $\rgeo \Lunit.$
	In case \textbf{i)}, the estimate
	\eqref{E:ONEGOODCOMMUTATORESTIMATE} implies that
	$\inhomarg{\mathscr{Z}^N} 
	= \rgeo 
		(\angdiff \Psi^{\#}) 
		\cdot
		(\upmu \angdiff \mathscr{Z}^{N-1} \mytr \upchi^{(Small)})
	+ Harmless^{\leq N}.$ The main point is that $\mathscr{Z}^{N-1}$ must contain at least
	one factor of $\rgeo \Lunit.$ Hence, essentially the same arguments given 
	for case \textbf{ii)} in the previous paragraph
	imply that
	$\rgeo 
		(\angdiff \Psi^{\#}) 
		\cdot
		(\upmu \angdiff \mathscr{Z}^{N-1} \mytr \upchi^{(Small)}) = Harmless^{\leq N},$
	which yields the desired estimate \eqref{E:TWOGOODCOMMUTATORSESTIMATE} in this case.

\end{proof}

\section{The important terms in the top-order derivatives of 
\texorpdfstring{$\deform{Z}$}{the deformation tensors of the commutation vectorfields}}
	\label{S:IDENTIFICATIONOFDIFFICULTTOPORDERTERMS}
	As we noted above, the most difficult terms in the commuted wave equation 
	\eqref{E:WAVENTIMESCOMMUTEDBASICSTRUCTURE}
	are contained
	in the first term 
	$\mathscr{Z}^{N-1}(\upmu \D_{\alpha} \Jcurrent{Z}^{\alpha}[\Psi])$
	on the right-hand side of equation \eqref{E:NTIMESCOMMUTEDWAVEINHOMOGENEOUSTERMFIRSTSPLITTING}.
	Specifically, the difficult terms arise from the top-order derivatives of 
	the deformation tensors of the commutation vectorfields,
	whose first derivatives appear after a single commutation;
	see Prop.~\ref{P:COMMUTATIONCURRENTDIVERGENCEFRAMEDECOMP}.
	The terms are difficult precisely because at top order, we need to use modified quantities
	and elliptic estimates to bound them in $L^2.$
	In the next three lemmas, 
	in the three cases 
	$Z = \rgeo \Lunit,$
	$Z = \Rad,$
	$Z= \Rot,$
	we identify these ``important'' difficult terms
	in $\mathscr{Z}^{N-1}(\upmu \D_{\alpha} \Jcurrent{Z}^{\alpha}[\Psi]).$
	The important terms appear on the left-hand sides of
	the inequalities of the lemmas.

\begin{lemma} [\textbf{Identification of the important top-order terms in the derivatives of} $\deform{\Rad}$]
\label{L:RADIMPORTANTDEFTENSORTERMS}
	Let $1 \leq N \leq 24$ be an integer. 
	Under the small-data and bootstrap assumptions 
	of Sects.~\ref{S:PSISOLVES}-\ref{S:C0BOUNDBOOTSTRAP},
	if $\varepsilon$ is sufficiently small, then
	the following estimates hold on $\mathcal{M}_{\Tboot,U_0}:$
\begin{subequations}
\begin{align}  \label{E:RADDEFIMPORTANTLIERADSPHERELANDRADTRACESPHERETERMS}
		& \left|
				\angLie_{\mathscr{Z}}^{N-1} \angLie_{\Rad} \angdeformoneformupsharparg{\Rad}{\Lunit}
				+ \angdiffuparg{\#} \mathscr{Z}^{N-1} \Rad \upmu
			\right|,
			\\
		& \left|
			\mathscr{Z}^{N-1} \Rad \mytr  \angdeform{\Rad}
			+ 2 \upmu \mathscr{Z}^{N-1} \angLap \upmu 
		\right|
			\notag \\
		& \lesssim
			\left|
				\fourmyarray
					[\rgeo \Lunit \mathscr{Z}^{\leq N} \Psi]
					{\Rad \mathscr{Z}^{\leq N} \Psi}
					{\rgeo \angdiff \mathscr{Z}^{\leq N} \Psi}
					{\mathscr{Z}^{\leq N} \Psi}
			\right|
			+ \frac{1}{1 + t}
				\left| 
					\myarray
						[\mathscr{Z}^{\leq N} (\upmu - 1)]
						{\sum_{a=1}^3 \rgeo |\mathscr{Z}^{\leq N} \Lunit_{(Small)}^a|}
				\right|
			+ \frac{1}{1 + t},
				\notag
	\end{align}
	
\begin{align}  \label{E:RADDEFIMPORTANTANGDIVSPHERELANDANGDIFFPILRADTERMS}
		& 
		(1+t)
		\left|
			\mathscr{Z}^{N-1} \angdiv \angdeformoneformupsharparg{\Rad}{\Lunit}
				+ \angLap \mathscr{Z}^{N-1} \upmu
		\right|,
 			\\
 		& 
			\left| 
				\angLie_{\mathscr{Z}}^{N-1} \angdiffuparg{\#} \deformarg{\Rad}{\Lunit}{\Rad}
				+ \angdiffuparg{\#} \mathscr{Z}^{N-1} \Rad \upmu
			\right|,
			\notag \\
		& (1+t)
			\left|
				\angLie_{\mathscr{Z}}^{N-1} \angdiv \angdeformfreeupdoublesharparg{\Rad}
				+ \upmu \angdiffuparg{\#} \mathscr{Z}^{N-1} \mytr \upchi^{(Small)}
			\right|
			\notag \\
		& \lesssim
			\left|
				\fourmyarray
					[\rgeo \Lunit \mathscr{Z}^{\leq N} \Psi]
					{\Rad \mathscr{Z}^{\leq N} \Psi}
					{\rgeo \angdiff \mathscr{Z}^{\leq N} \Psi}
					{\mathscr{Z}^{\leq N} \Psi}
			\right|
			+ \frac{1}{1 + t}
				\left| 
					\myarray
						[\mathscr{Z}^{\leq N} (\upmu - 1)]
						{\sum_{a=1}^3 \rgeo |\mathscr{Z}^{\leq N} \Lunit_{(Small)}^a|}
				\right|.
				\notag
	\end{align}
	\end{subequations}
\end{lemma}

\begin{proof}
	We first prove the bound \eqref{E:RADDEFIMPORTANTLIERADSPHERELANDRADTRACESPHERETERMS}
	for 
	$\angLie_{\mathscr{Z}}^{N-1} \angLie_{\Rad} \angdeformoneformupsharparg{\Rad}{\Lunit} + \angdiffuparg{\#} \mathscr{Z}^{N-1} \Rad \upmu.$	
	We apply $\angLie_{\mathscr{Z}}^{N-1} \angLie_{\Rad}$ to the $\gsphere-$dual of the right-hand side of  
	\eqref{E:RADDEFORMLA}. By Lemma~\ref{L:LANDRADCOMMUTEWITHANGDIFF} and the Leibniz rule, 
	the first term on right-hand side is equal to the principal term
	$- \angdiffuparg{\#} \mathscr{Z}^{N-1} \Rad \upmu$ plus an error term 
	that is bounded in magnitude by
	$\lesssim 
	\mathop{\sum_{N_1 + N_2 \leq N}}_{N_2 \leq N-1}
	\left| 
		\angLie_{\mathscr{Z}}^{N_1} \ginversesphere 
	\right| 
	\left| \angdiff 
		\mathscr{Z}^{N_2} \upmu 
	\right|.$
	We bring the principal term over to the left, while by
	inequality \eqref{E:FUNCTIONPOINTWISEANGDINTERMSOFANGLIEO}
	and the estimates
	\eqref{E:CRUDEPOINTWISEBOUNDSDERIVATIVESOFANGULARDEFORMATIONTENSORS},
	\eqref{E:CRUDELOWERORDERC0BOUNDDERIVATIVESOFANGULARDEFORMATIONTENSORS},
	and
	\eqref{E:C0BOUNDCRUCIALEIKONALFUNCTIONQUANTITIES},
	the error term is $\lesssim$
	the right-hand side of \eqref{E:RADDEFIMPORTANTLIERADSPHERELANDRADTRACESPHERETERMS} as desired. 
	Hence, to conclude the desired bound, it only remains for us to bound
	the magnitude of $\angLie_{\mathscr{Z}}^{N-1} \angLie_{\Rad} = \angLie_{\mathscr{Z}}^N$ applied to the 
	$\gsphere-$dual of the
	second product on the right-hand side of \eqref{E:RADDEFORMLA} by the right-hand side
	of \eqref{E:RADDEFIMPORTANTLIERADSPHERELANDRADTRACESPHERETERMS}.
	This product is of the form 
	$ \xi :=	
		G_{(Frame)}
		\ginversesphere
			\threemyarray
				[\upmu \Lunit \Psi]
				{\Rad \Psi}
				{\upmu \angdiff \Psi}.
	$
	To bound $\angLie_{\mathscr{Z}}^N \xi$
	in magnitude by the right-hand side of \eqref{E:RADDEFIMPORTANTLIERADSPHERELANDRADTRACESPHERETERMS},
	we use the Leibniz rule
	and the estimates
	\eqref{E:FUNCTIONAVOIDINGCOMMUTING}
	(with $\Psi$ in the role of $f$),
	\eqref{E:LIEDERIVATIVESOFGRAMEINTERMSOFOTHERVARIABLES},
	\eqref{E:LOWERORDERC0BOUNDLIEDERIVATIVESOFGRAME},	
	\eqref{E:CRUDEPOINTWISEBOUNDSDERIVATIVESOFANGULARDEFORMATIONTENSORS},
	\eqref{E:CRUDELOWERORDERC0BOUNDDERIVATIVESOFANGULARDEFORMATIONTENSORS},
  \eqref{E:C0BOUNDCRUCIALEIKONALFUNCTIONQUANTITIES},
	and the bootstrap assumptions \eqref{E:PSIFUNDAMENTALC0BOUNDBOOTSTRAP}
	to deduce that
	\begin{subequations}
	\begin{align} \label{E:GFRAMEGINVERSESPHEREUPMUPSIDERIVATIVESLEIBNIZEXPANDED}
	\left|
		\angLie_{\mathscr{Z}}^N \xi
	\right|
	& \lesssim
			\sum_{N_1 + N_2 + N_3 + N_4 \leq N}
			\left|
				\angLie_{\mathscr{Z}}^{N_1} G_{(Frame)}
			\right|
			\left|
				\angLie_{\mathscr{Z}}^{N_2} \ginversesphere
			\right|
			\left|
				\mathscr{Z}^{N_3} \upmu
			\right|
			\left|
			\angLie_{\mathscr{Z}}^{N_4}
				\myarray 
					[\Lunit \Psi]
					{\angdiff \Psi}
			\right|
				\\
		& \ \
			+
			\sum_{N_1 + N_2 + N_3 \leq N}
			\left|
				\angLie_{\mathscr{Z}}^{N_1} G_{(Frame)}
			\right|
			\left|
				\angLie_{\mathscr{Z}}^{N_2} \ginversesphere
			\right|
			\left|
				\mathscr{Z}^{N_3} \Rad \Psi
			\right|
			\notag
				\\
		& \lesssim
		\frac{\ln(\myexp + t)}{1 + t}
			\left|
				\threemyarray[\rgeo \Lunit \mathscr{Z}^{\leq N} \Psi]
					{\rgeo \angdiff \mathscr{Z}^{\leq N} \Psi} 
					{\mathscr{Z}^{\leq N} \Psi}
			\right|
			+ \left|
					\Rad \mathscr{Z}^{\leq N} \Psi
				\right|
			+
			\varepsilon
			\frac{1}{(1 + t)^2}
			\left| 
				\myarray
					[\mathscr{Z}^{\leq N} (\upmu - 1)]
						{\rgeo \sum_{a=1}^3 |\mathscr{Z}^{\leq N} \Lunit_{(Small)}^a|}
			\right|,
				\notag \\
		\left\|
			\angLie_{\mathscr{Z}}^{\leq 12} \xi
		\right\|_{C^0(\Sigma_t^u)}
		& \lesssim \varepsilon \frac{1}{1+t}.
			\label{E:C0BOUNDGFRAMEGINVERSESPHEREUPMUPSIDERIVATIVESLEIBNIZEXPANDED}
	\end{align}	
	\end{subequations}
	as desired (we use the bound \eqref{E:C0BOUNDGFRAMEGINVERSESPHEREUPMUPSIDERIVATIVESLEIBNIZEXPANDED} below).
	
	To prove the bound \eqref{E:RADDEFIMPORTANTLIERADSPHERELANDRADTRACESPHERETERMS} for 
	$\mathscr{Z}^{N-1} \Rad \mytr  \angdeform{\Rad} + 2 \upmu \mathscr{Z}^{N-1} \angLap \upmu,$ we apply
	$\mathscr{Z}^{N-1} \Rad$ to the right-hand side of 
	\eqref{E:RADDEFORMANG}. The last term is of the form 
	$\xi := 
		G_{(Frame)}
		\ginversesphere
			\threemyarray
				[\upmu \Lunit \Psi]
				{\Rad \Psi}
				{\upmu \angdiff \Psi}$
	and hence the proof of \eqref{E:GFRAMEGINVERSESPHEREUPMUPSIDERIVATIVESLEIBNIZEXPANDED}
	yields the desired bound \eqref{E:GFRAMEGINVERSESPHEREUPMUPSIDERIVATIVESLEIBNIZEXPANDED}
	for $\xi$.
	To bound the first term on the right-hand side of \eqref{E:RADDEFORMANG}, 
	we use \eqref{E:ZNAPPLIEDTORGEOISNOTTOOLARGE}
	to deduce that $\left| \mathscr{Z}^{N-1} \Rad (\rgeo^{-1} \upmu )\right| \lesssim 
	(1+t)^{-1}
	\left|
		\mathscr{Z}^{\leq N} (\upmu - 1)
	\right|
	+ (1+t)^{-1}$ as desired.
	It remains for us to address the important term
	$- 2 \mathscr{Z}^{N-1} \Rad (\upmu \mytr \upchi^{(Small)})$
	and to show that it is equal to $-2 \upmu \mathscr{Z}^{N-1} \angLap \upmu$
	plus an error term that is in magnitude $\lesssim$
	the right-hand side of \eqref{E:RADDEFIMPORTANTLIERADSPHERELANDRADTRACESPHERETERMS}.
	To see that this is the case, we apply the Leibniz rule. 
	We bound all terms except the one in which all derivatives fall on $\mytr \upchi^{(Small)}$
	by using 
	\eqref{E:C0BOUNDCRUCIALEIKONALFUNCTIONQUANTITIES}
	and \eqref{E:POINTWISEESTIMATESFORCHIJUNKINTERMSOFOTHERVARIABLES}.
	To handle the one term in which all derivatives fall on $\mytr \upchi^{(Small)},$
	we use in addition the third inequality in \eqref{E:TOPORDERDERIVATIVESOFANGDSQUAREDUPMUINTERMSOFCONTROLLABLE}.
	We have thus proved the desired bound.
	
	We now prove the bound \eqref{E:RADDEFIMPORTANTANGDIVSPHERELANDANGDIFFPILRADTERMS}
	for $\mathscr{Z}^{N-1} \angdiv \angdeformoneformupsharparg{\Rad}{\Lunit} + \angLap \mathscr{Z}^{N-1} \upmu.$
	We apply $\mathscr{Z}^{N-1} \angdiv$
	to the $\gsphere-$dual of the right-hand side of \eqref{E:RADDEFORMLA}.
	We begin by addressing the second term on the right-hand side, 
	which is equal to the $\gsphere-$dual of an $S_{t,u}$ one-form $\xi$ that is of the form 
	$\xi := G_{(Frame)}
			\threemyarray
				[\upmu \Lunit \Psi]
				{\Rad \Psi}
				{\upmu \angdiff \Psi}.$
	We now note that a similar proof to that of \eqref{E:GFRAMEGINVERSESPHEREUPMUPSIDERIVATIVESLEIBNIZEXPANDED}
	and \eqref{E:C0BOUNDGFRAMEGINVERSESPHEREUPMUPSIDERIVATIVESLEIBNIZEXPANDED} yields that
	$\xi$ verifies the bounds
	\eqref{E:GFRAMEGINVERSESPHEREUPMUPSIDERIVATIVESLEIBNIZEXPANDED}
	and \eqref{E:C0BOUNDGFRAMEGINVERSESPHEREUPMUPSIDERIVATIVESLEIBNIZEXPANDED}.
	Next, to derive the desired bound for $\mathscr{Z}^{N-1} \angdiv \xi^{\#},$
	we use the Leibniz rule to deduce that
	\begin{align} \label{E:LIEZNANGDIVXISHARPLEIBNIZEXPANDED}
		\left|
			\mathscr{Z}^{N-1} \angdiv \xi^{\#}
		\right|
		& \lesssim
			\sum_{N_1 + N_2 \leq N-1}
			\left|
				\angLie_{\mathscr{Z}}^{N_1} \ginversesphere
			\right|
			\left|
				\angD \angLie_{\mathscr{Z}}^{N_2} \xi
			\right|
			+
			\sum_{N_1 + N_2 \leq N-1}
			\left|
				\angLie_{\mathscr{Z}}^{N_1} \ginversesphere
			\right|
			\left|
				[\angLie_{\mathscr{Z}}^{N_2}, \angD] \xi
			\right|.
	\end{align}	
	To deduce that the first sum on the right-hand side 
	of \eqref{E:LIEZNANGDIVXISHARPLEIBNIZEXPANDED} is
	$\lesssim$ 
	the product of $(1 + t)^{-1}$
	and the right-hand side of \eqref{E:RADDEFIMPORTANTANGDIVSPHERELANDANGDIFFPILRADTERMS},
	we use the estimates 
	\eqref{E:TYPE02TENSORANGDINTERMSOFROTATIONALLIE},
	\eqref{E:CRUDEPOINTWISEBOUNDSDERIVATIVESOFANGULARDEFORMATIONTENSORS},
	\eqref{E:CRUDELOWERORDERC0BOUNDDERIVATIVESOFANGULARDEFORMATIONTENSORS},
	and the aforementioned bounds
	\eqref{E:GFRAMEGINVERSESPHEREUPMUPSIDERIVATIVESLEIBNIZEXPANDED}
	and \eqref{E:C0BOUNDGFRAMEGINVERSESPHEREUPMUPSIDERIVATIVESLEIBNIZEXPANDED}.
	To deduce 
	that the second sum on the right-hand side of 
	\eqref{E:SECONDPOINTWISEBOUNDJUNKTERMSDIVERGENCEOFROTDEFORMRADSPHERE}
	is
	$\lesssim$ 
	the product of $(1 + t)^{-1}$
	and the right-hand side of \eqref{E:RADDEFIMPORTANTANGDIVSPHERELANDANGDIFFPILRADTERMS},
	we combine similar reasoning with 
	the commutator estimate
	\eqref{E:COMMUTATORESTIMATESVECTORFIELDSACTINGONANGDTENSORS},
	where $N_2$ is in the role of $N.$
	We now address the term	$-\mathscr{Z}^{N-1} \angLap \upmu$ 
	arising from the first term on the right-hand side of \eqref{E:RADDEFORMLA}.
	We commute $\mathscr{Z}^{N-1}$ and $\angLap$ and move the principal
	term $\angLap \mathscr{Z}^{N-1} \upmu$ to the left-hand side.
	To complete the proof
	of the bound for $\mathscr{Z}^{N-1} \angdiv \angdeformoneformupsharparg{\Rad}{\Lunit} + \angLap \mathscr{Z}^{N-1} \upmu,$
	it only remains for us to bound the magnitude of the commutator term
	$[\mathscr{Z}^{N-1}, \angLap] \upmu$ 
	by the product of $(1 + t)^{-1}$ and the right-hand side of \eqref{E:RADDEFIMPORTANTANGDIVSPHERELANDANGDIFFPILRADTERMS}.
	To this end, we use the estimate
	\eqref{E:ANGLAPZNCOMMUTATORACTINGONFUNCTIONSSPOINTWISE} with $N-1$ in the role of $N$
	and $\upmu$ in the role of $f$ and inequality \eqref{E:C0BOUNDCRUCIALEIKONALFUNCTIONQUANTITIES}.
	
	We now prove the bound \eqref{E:RADDEFIMPORTANTANGDIVSPHERELANDANGDIFFPILRADTERMS} for
	$\angLie_{\mathscr{Z}}^{N-1} \angdiffuparg{\#} \deformarg{\Rad}{\Lunit}{\Rad}
				+ \angdiffuparg{\#} \mathscr{Z}^{N-1} \Rad \upmu.$
	We apply $\angLie_{\mathscr{Z}}^{N-1} \angdiffuparg{\#}$ 
	to the right-hand side of \eqref{E:RADDEFORMLRADUNIT}.
	The quantity we must estimate
	is $- \angLie_{\mathscr{Z}}^{N-1} \angdiffuparg{\#} \Rad \upmu.$
	From Lemma~\ref{L:LANDRADCOMMUTEWITHANGDIFF}
	and the Leibniz rule,
	it follows that this quantity is equal to the principal term
	$- \angdiffuparg{\#} \mathscr{Z}^{N-1} \Rad \upmu$
	plus an error term that is bounded in magnitude by
	\begin{align} \label{E:FIRSTPOINTWISELIEZNANGDIFFRADDEFORMLRAD}
		\lesssim
		\sum_{N_1 + N_2 \leq N}
		\left|
			\angLie_{\mathscr{Z}}^{N_1}
			\ginversesphere
		\right|
		\left|
			\angdiff \angLie_{\mathscr{Z}}^{N_2} \upmu
		\right|.
	\end{align}
	We bring the principal term over to the left-hand side. 
	To complete the proof, we need to show that the 
	right-hand side of \eqref{E:FIRSTPOINTWISELIEZNANGDIFFRADDEFORMLRAD} is 
	$\lesssim$ the right-hand side of \eqref{E:RADDEFIMPORTANTANGDIVSPHERELANDANGDIFFPILRADTERMS}.
	To deduce the desired bound, 
	we use the estimates
	\eqref{E:CRUDEPOINTWISEBOUNDSDERIVATIVESOFANGULARDEFORMATIONTENSORS},
	\eqref{E:CRUDELOWERORDERC0BOUNDDERIVATIVESOFANGULARDEFORMATIONTENSORS},
	\eqref{E:FUNCTIONPOINTWISEANGDINTERMSOFANGLIEO},
	and
	\eqref{E:C0BOUNDCRUCIALEIKONALFUNCTIONQUANTITIES}.
		
	Finally, we prove \eqref{E:RADDEFIMPORTANTANGDIVSPHERELANDANGDIFFPILRADTERMS} for
	the term
	$\angLie_{\mathscr{Z}}^{N-1} \angdiv \angdeformfreeupdoublesharparg{\Rad}
				+ \upmu \angdiffuparg{\#} \mathscr{Z}^{N-1} \mytr \upchi^{(Small)}
	.$	
	We apply $\angLie_{\mathscr{Z}}^{N-1} \angdiv$ 
	to the double $\gsphere-$dual of the right-hand side of \eqref{E:RADDEFORMTRFREEANG}.
	We begin by addressing the first term, which is the difficult one:
	$-2 \angLie_{\mathscr{Z}}^{N-1} 
	\angDarg{B} 
	\left\lbrace
		\upmu (\ginversesphere)^{AC} (\ginversesphere)^{BD} \hat{\upchi}_{CD}^{(Small)}
	\right\rbrace.$
	When all derivatives fall on $(\ginversesphere)^{BD} \hat{\upchi}_{CD}^{(Small)},$ we use
	\eqref{E:C0BOUNDCRUCIALEIKONALFUNCTIONQUANTITIES}
	and the second inequality in \eqref{E:DIVCHIJUNKANGDIFFTRCHIJUNKHIGHERORDERCOMMMUTOR}
	to rewrite it as the principal term
	$- \upmu \angdiffuparg{\#} \mathscr{Z}^{N-1} \mytr \upchi^{(Small)}$
	plus an error term with magnitude $\lesssim$ 
	the product of $(1+t)^{-1}$
	and the right-hand side
	of \eqref{E:RADDEFIMPORTANTANGDIVSPHERELANDANGDIFFPILRADTERMS}.
	We bring the principal term over to the left-hand side.
	Using Lemma~\ref{L:LANDRADCOMMUTEWITHANGDIFF}, we see that the remaining terms arising from the 
	Leibniz expansion of 
 	$-2 \angLie_{\mathscr{Z}}^{N-1} 
	\angDarg{B} 
	\left\lbrace
		\upmu (\ginversesphere)^{AC} (\ginversesphere)^{BD} \hat{\upchi}_{CD}^{(Small)}
	\right\rbrace$
	are bounded in magnitude by
	\begin{align} \label{E:RADDEFORMSPHERETRFREENONPRINCIPALCHIJUNKPRODUCT}
		&
		\lesssim
		\mathop{\sum_{N_1 + N_2 + N_3 + N_4 \leq N-1}}_{N_4 \leq N - 2}
			\left|
				\angLie_{\mathscr{Z}}^{N_1} \upmu
			\right|
			\left|
				\angLie_{\mathscr{Z}}^{N_2} \ginversesphere
			\right|
			\left|
				\angLie_{\mathscr{Z}}^{N_3} \ginversesphere
			\right|
			\left|
				\angLie_{\mathscr{Z}}^{N_4} \angD \hat{\upchi}^{(Small)}
			\right|
				\\
	& \ \ 
		+ \sum_{N_1 + N_2 + N_3 + N_4 \leq N-1}
			\left|
				\angdiff \angLie_{\mathscr{Z}}^{N_1} \upmu
			\right|
			\left|
				\angLie_{\mathscr{Z}}^{N_2} \ginversesphere
			\right|
			\left|
				\angLie_{\mathscr{Z}}^{N_3} \ginversesphere
			\right|
			\left|
				\angLie_{\mathscr{Z}}^{N_4} \hat{\upchi}^{(Small)}
			\right|.
				\notag 
		\end{align}
		Using inequalities
		\eqref{E:FUNCTIONPOINTWISEANGDINTERMSOFANGLIEO} 
		and \eqref{E:TYPE02TENSORANGDINTERMSOFROTATIONALLIE},
		the estimates
		\eqref{E:CRUDEPOINTWISEBOUNDSDERIVATIVESOFANGULARDEFORMATIONTENSORS},
		\eqref{E:CRUDELOWERORDERC0BOUNDDERIVATIVESOFANGULARDEFORMATIONTENSORS},
		\eqref{E:POINTWISEESTIMATESFORCHIJUNKINTERMSOFOTHERVARIABLES},
		and
		\eqref{E:C0BOUNDCRUCIALEIKONALFUNCTIONQUANTITIES},
		the commutator estimate \eqref{E:COMMUTATORESTIMATESVECTORFIELDSACTINGONANGDTENSORS}
		with $\xi = \hat{\upchi}^{(Small)},$
		and \eqref{E:LIEDERIVATIVESOFTRACEFREEPARTINTERMSOFLIEDERIVATIVESOFFULLTENSOR}
		with $\xi = \upchi^{(Small)},$
		we deduce that the right-hand side of 
		\eqref{E:RADDEFORMSPHERETRFREENONPRINCIPALCHIJUNKPRODUCT}
		is $\lesssim$ 
		the product of $(1+t)^{-1}$
		and the right-hand side of \eqref{E:RADDEFIMPORTANTANGDIVSPHERELANDANGDIFFPILRADTERMS}	
		as desired.
		To complete the proof, it remains for us to show that
		the $\angLie_{\mathscr{Z}}^{N-1} \angdiv$ derivatives of 
		the double $\gsphere-$dual of the
		remaining terms on the right-hand side of \eqref{E:RADDEFORMTRFREEANG}
		have magnitudes that are
		$\lesssim$ 
		the product of $(1+t)^{-1}$
		and the right-hand side of \eqref{E:RADDEFIMPORTANTANGDIVSPHERELANDANGDIFFPILRADTERMS}.	
		To this end, we let $\hat{\xi}^{\# \#}$ denote the remaining terms,
		where we note that $\xi$ is of the form
		$\xi := G_{(Frame)} 
			\otimes
			\threemyarray
				[\upmu \Lunit \Psi]
				{\Rad \Psi}
				{\upmu \angdiff \Psi}.
				$	
	We now note that a similar proof to that of \eqref{E:GFRAMEGINVERSESPHEREUPMUPSIDERIVATIVESLEIBNIZEXPANDED}
	and \eqref{E:C0BOUNDGFRAMEGINVERSESPHEREUPMUPSIDERIVATIVESLEIBNIZEXPANDED} yields that
	$\xi$ verifies the bounds
	\eqref{E:GFRAMEGINVERSESPHEREUPMUPSIDERIVATIVESLEIBNIZEXPANDED}
	and \eqref{E:C0BOUNDGFRAMEGINVERSESPHEREUPMUPSIDERIVATIVESLEIBNIZEXPANDED}.
	We now recall that our goal is to bound
	$\angLie_{\mathscr{Z}}^{N-1} \angdiv \hat{\xi}^{\# \#}$
	in magnitude by
	$\lesssim$ 
	the product of $(1+t)^{-1}$
	and the right-hand side of \eqref{E:ROTDEFIMPORTANTANGDIVSPHERELANDANGDIFFPILRADTERMS}. 
	To proceed, we use the Leibniz rule to deduce that
	inequality \eqref{E:LIEZNANGDIVXISHARPLEIBNIZEXPANDED} holds
	with $\angLie_{\mathscr{Z}}^{N-1} \angdiv \hat{\xi}^{\# \#}$ 
	in place of $\mathscr{Z}^{N-1} \angdiv \hat{\xi}^{\# \#}$ on the left-hand side
	and $\hat{\xi}$ in place of $\xi$ on the right-hand side.
	Finally, to bound the right-hand side of 
	\eqref{E:LIEZNANGDIVXISHARPLEIBNIZEXPANDED} 
	(with $\hat{\xi}$ in place of $\xi$)
	by
	$\lesssim$ 
	the product of $(1+t)^{-1}$
	and the right-hand side of \eqref{E:ROTDEFIMPORTANTANGDIVSPHERELANDANGDIFFPILRADTERMS},
	we use the estimates 
	\eqref{E:TYPE02TENSORANGDINTERMSOFROTATIONALLIE},
	\eqref{E:CRUDEPOINTWISEBOUNDSDERIVATIVESOFANGULARDEFORMATIONTENSORS},
	\eqref{E:CRUDELOWERORDERC0BOUNDDERIVATIVESOFANGULARDEFORMATIONTENSORS},
	\eqref{E:COMMUTATORESTIMATESVECTORFIELDSACTINGONANGDTENSORS},
	\eqref{E:LIEDERIVATIVESOFTRACEFREEPARTINTERMSOFLIEDERIVATIVESOFFULLTENSOR},
	and the aforementioned bounds
	\eqref{E:GFRAMEGINVERSESPHEREUPMUPSIDERIVATIVESLEIBNIZEXPANDED}
	and \eqref{E:C0BOUNDGFRAMEGINVERSESPHEREUPMUPSIDERIVATIVESLEIBNIZEXPANDED}.
		
	\end{proof}

\begin{lemma} [\textbf{Identification of the important top-order terms in the derivatives of} $\deform{\Rot_{(l)}}$]
	\label{L:ROTIMPORTANTDEFTENSORTERMS}
	Let $1 \leq N \leq 24$ be an integer. 
	Under the small-data and bootstrap assumptions 
	of Sects.~\ref{S:PSISOLVES}-\ref{S:C0BOUNDBOOTSTRAP},
	if $\varepsilon$ is sufficiently small, then
	the following estimates hold on $\mathcal{M}_{\Tboot,U_0}:$
	\begin{subequations}
	\begin{align}  \label{E:ROTDEFIMPORTANTLIERADSPHERELANDRADTRACESPHERETERMS}
		& \left|
			 \angLie_{\mathscr{Z}}^{N-1}	\angLie_{\Rad} \angdeformoneformupsharparg{\Rot_{(l)}}{\Lunit}
				+ (\angD^{2 \#} \mathscr{Z}^{N-1} \upmu) \cdot \Rot_{(l)} 
			\right|,
			\\
		& \left|
			\mathscr{Z}^{N-1} \Rad \mytr  \angdeform{\Rot_{(l)}}
			- 2 \RotRadcomponent{l} \angLap \mathscr{Z}^{N-1} \upmu
		\right|
			\notag \\
		& \lesssim
			\left|
				\fourmyarray
					[\rgeo \Lunit \mathscr{Z}^{\leq N} \Psi]
					{\Rad \mathscr{Z}^{\leq N} \Psi}
					{\rgeo \angdiff \mathscr{Z}^{\leq N} \Psi}
					{\mathscr{Z}^{\leq N} \Psi}
			\right|
			+ \frac{1}{1 + t}
				\left| 
					\myarray
						[\mathscr{Z}^{\leq N} (\upmu - 1)]
						{\sum_{a=1}^3 \rgeo |\mathscr{Z}^{\leq N} \Lunit_{(Small)}^a|}
				\right|,
				\notag
	\end{align}

	\begin{align} \label{E:ROTDEFIMPORTANTANGDIVSPHERELANDANGDIFFPILRADTERMS}
		& 
		(1+t)
		\left|
			\mathscr{Z}^{N-1} \angdiv \angdeformoneformupsharparg{\Rot_{(l)}}{\Lunit}
				+ \Rot_{(l)} \mathscr{Z}^{N-1} \mytr \upchi^{(Small)}
		\right|,
 			\\
		& \frac{(1+t)}{\ln(\myexp + t)}
			\left|
			\mathscr{Z}^{N-1} \angdiv \angdeformoneformupsharparg{\Rot_{(l)}}{\Rad}
			- \left \lbrace
					\upmu \Rot_{(l)} \mathscr{Z}^{N-1} \mytr \upchi^{(Small)}
					+ \RotRadcomponent{l} \angLap \mathscr{Z}^{N-1} \upmu
				\right\rbrace
		\right|,
			\notag \\
		& \left| 
				\angLie_{\mathscr{Z}}^{N-1} \angdiffuparg{\#} \deformarg{\Rot_{(l)}}{\Lunit}{\Rad}
				+ (\angD^{2 \#} \mathscr{Z}^{N-1} \upmu) \cdot \Rot_{(l)}
			\right|,
			\notag \\
		& 
		(1+t)
		\left|
			\angLie_{\mathscr{Z}}^{N-1} \angdiv \angdeformfreeupdoublesharparg{\Rot_{(l)}}
			- \RotRadcomponent{l} \angdiffuparg{\#} \mathscr{Z}^{N-1} \mytr \upchi^{(Small)}
		\right|
			\notag \\
		& \lesssim
			\left|
				\fourmyarray
					[\rgeo \Lunit \mathscr{Z}^{\leq N} \Psi]
					{\Rad \mathscr{Z}^{\leq N} \Psi}
					{\rgeo \angdiff \mathscr{Z}^{\leq N} \Psi}
					{\mathscr{Z}^{\leq N} \Psi}
			\right|
			+ \frac{1}{1+t}
				\left| 
					\myarray
						[\mathscr{Z}^{\leq N} (\upmu - 1)]
						{\sum_{a=1}^3 \rgeo |\mathscr{Z}^{\leq N} \Lunit_{(Small)}^a|}
				\right|.
				\notag
	\end{align}
	\end{subequations}	
\end{lemma}	

\begin{proof}
	We first prove \eqref{E:ROTDEFIMPORTANTLIERADSPHERELANDRADTRACESPHERETERMS} for
	$\angLie_{\mathscr{Z}}^{N-1}	\angLie_{\Rad} \angdeformoneformupsharparg{\Rot_{(l)}}{\Lunit}
	+ (\angD^{2 \#} \mathscr{Z}^{N-1} \upmu) \cdot \Rot_{(l)}.$	We apply
	$\angLie_{\mathscr{Z}}^{N-1}	\angLie_{\Rad}$ to the $\gsphere-$dual of the right-hand side of
	\eqref{E:ROTDEFORMLSPHERE}. We begin by addressing the first term, which is the difficult one:
	$- \angLie_{\mathscr{Z}}^{N-1} \angLie_{\Rad} \left\lbrace (\ginversesphere)^{AB} \upchi_{BC}^{(Small)} \Rot_{(l)}^C \right\rbrace.$ 
	When all derivatives fall on $\upchi_{BC}^{(Small)},$ 
	we use \eqref{E:ROTATIONPOINTWISENORMESTIMATE}
	and the first inequality in \eqref{E:TOPORDERDERIVATIVESOFANGDSQUAREDUPMUINTERMSOFCONTROLLABLE} 
	to rewrite
	this term as the principal eikonal function term 
	$-(\ginversesphere)^{AB} (\angDsquaredarg{B}{C} \mathscr{Z}^{N-1} \upmu)\Rot_{(l)}^C$
	plus an error term with magnitude $\lesssim$ the right-hand side of \eqref{E:ROTDEFIMPORTANTLIERADSPHERELANDRADTRACESPHERETERMS}.
	We then bring the principal term over to the left-hand side. To conclude the desired inequality,
	it remains for us to show that all remaining terms arising from the right-hand side of
	\eqref{E:ROTDEFORMLSPHERE} have magnitudes that are
	$\lesssim$ the right-hand side of \eqref{E:ROTDEFIMPORTANTLIERADSPHERELANDRADTRACESPHERETERMS}.
	The remaining terms arising in the Leibniz expansion of
	$- \angLie_{\mathscr{Z}}^{N-1} \angLie_{\Rad} \left\lbrace (\ginversesphere)^{AB} \upchi_{BC}^{(Small)} \Rot_{(l)}^C \right\rbrace$ 
	are bounded in magnitude by
	\begin{align} \label{E:ROTDEFORMSPHERELNONPRINCIPALCHIJUNKPRODUCT}
		\lesssim
		\mathop{\sum_{N_1 + N_2 + N_3 \leq N}}_{N_2 \leq N-1}
			\left|
				\angLie_{\mathscr{Z}}^{N_1} \ginversesphere
			\right|
			\left|
				\angLie_{\mathscr{Z}}^{N_2} \upchi^{(Small)}
			\right|
			\left|
				\angLie_{\mathscr{Z}}^{N_3} \Rot_{(l)}
			\right|.
	\end{align}
	From the estimates
	\eqref{E:CRUDEPOINTWISEBOUNDSDERIVATIVESOFANGULARDEFORMATIONTENSORS},
	\eqref{E:CRUDELOWERORDERC0BOUNDDERIVATIVESOFANGULARDEFORMATIONTENSORS},
	\eqref{E:POINTWISEESTIMATESFORCHIJUNKINTERMSOFOTHERVARIABLES},
	\eqref{E:LIEDERIVATIVESOFROTATIONSPOINTWISEESTIMATE},
	\eqref{E:LOWERORDERLIEDERIVATIVESOFROTATIONSC0BOUND},
	and \eqref{E:C0BOUNDCRUCIALEIKONALFUNCTIONQUANTITIES},
	we conclude that the right-hand side of \eqref{E:ROTDEFORMSPHERELNONPRINCIPALCHIJUNKPRODUCT}
	is $\lesssim$ the right-hand side of \eqref{E:ROTDEFIMPORTANTLIERADSPHERELANDRADTRACESPHERETERMS} 
	as desired. To bound the first term
	$\angLie_{\mathscr{Z}}^{N-1}	\angLie_{\Rad}
		\left\lbrace
		G_{(Frame)} 
		\ginversesphere
			\myarray
				[\Rot_{(l)}]
				{\RotRadcomponent{l}}
			\myarray
				[\Lunit \Psi]
				{\angdiff \Psi}
	\right\rbrace
	$
	arising from the right-hand side of \eqref{E:ERRORROTDEFORMLSPHERE}
	by the right-hand side of \eqref{E:ROTDEFIMPORTANTLIERADSPHERELANDRADTRACESPHERETERMS},
	we use the estimates
	\eqref{E:FUNCTIONDERIVATIVESAVOIDINGCOMMUTING},
	\eqref{E:FIRSTPOINTWISEBOUNDEUCLIDEANROTATIONRADCOMPONENT},
	\eqref{E:LOWERORDERC0BOUNDEUCLIDEANROTATIONRADCOMPONENT},
	\eqref{E:LIEDERIVATIVESOFGRAMEINTERMSOFOTHERVARIABLES},
	\eqref{E:LOWERORDERC0BOUNDLIEDERIVATIVESOFGRAME},
	\eqref{E:CRUDEPOINTWISEBOUNDSDERIVATIVESOFANGULARDEFORMATIONTENSORS},
	\eqref{E:CRUDELOWERORDERC0BOUNDDERIVATIVESOFANGULARDEFORMATIONTENSORS},
	\eqref{E:LIEDERIVATIVESOFROTATIONSPOINTWISEESTIMATE},
	and \eqref{E:LOWERORDERLIEDERIVATIVESOFROTATIONSC0BOUND}
	and the bootstrap assumptions \eqref{E:PSIFUNDAMENTALC0BOUNDBOOTSTRAP}.
	To bound the second term 
	$\angLie_{\mathscr{Z}}^{N-1}	\angLie_{\Rad}
		\left\lbrace
			\smoothfunction(\Psi) \Lunit_{(Small)} \angdiffuparg{\#} x
		\right\rbrace
	$
	arising from the 
	right-hand side of \eqref{E:ERRORROTDEFORMLSPHERE}
	by the right-hand side of \eqref{E:ROTDEFIMPORTANTLIERADSPHERELANDRADTRACESPHERETERMS},
	we use the estimates
	\eqref{E:FUNCTIONDERIVATIVESAVOIDINGCOMMUTING},
	\eqref{E:POINTWISEBOUNDPROJECTEDLIEDERIVATIVESANGDIFFCOORDINATEX}
	\eqref{E:LOWERORDERPOINTWISEBOUNDPROJECTEDLIEDERIVATIVESANGDIFFCOORDINATEX},
	\eqref{E:CRUDEPOINTWISEBOUNDSDERIVATIVESOFANGULARDEFORMATIONTENSORS},
	\eqref{E:CRUDELOWERORDERC0BOUNDDERIVATIVESOFANGULARDEFORMATIONTENSORS},
	and
	\eqref{E:C0BOUNDCRUCIALEIKONALFUNCTIONQUANTITIES},
	and the bootstrap assumptions \eqref{E:PSIFUNDAMENTALC0BOUNDBOOTSTRAP}.

	We now prove \eqref{E:ROTDEFIMPORTANTLIERADSPHERELANDRADTRACESPHERETERMS} for
	$\mathscr{Z}^{N-1} \Rad \mytr  \angdeform{\Rot_{(l)}}
			- 2 \RotRadcomponent{l} \angLap \mathscr{Z}^{N-1} \upmu.$
	We apply $\mathscr{Z}^{N-1} \Rad$ to the right-hand side of
	\eqref{E:ROTDEFORMSPHERETRACE} and apply the Leibniz rule.
	We begin by addressing the first term, which is the difficult one:
	$2 \mathscr{Z}^{N-1} \Rad(\RotRadcomponent{l} \mytr \upchi^{(Small)}).$
	When all derivatives fall on $\mytr \upchi^{(Small)},$
	we use \eqref{E:LOWERORDERC0BOUNDEUCLIDEANROTATIONRADCOMPONENT}
	and the third inequality in \eqref{E:TOPORDERDERIVATIVESOFANGDSQUAREDUPMUINTERMSOFCONTROLLABLE}
	to rewrite this term as the principal term
	$2 \RotRadcomponent{l} \angLap \mathscr{Z}^{N-1} \upmu$
	plus an error term with magnitude $\lesssim$
	the right-hand side of \eqref{E:ROTDEFIMPORTANTLIERADSPHERELANDRADTRACESPHERETERMS}.
	We then bring the principal term over to the left-hand side.
	To complete the proof, it remains for us to show that
	the $\mathscr{Z}^{N-1} \Rad$ derivatives of the
	remaining terms on the right-hand side of \eqref{E:ROTDEFORMSPHERETRACE}
	have magnitudes that are
	$\lesssim$ the right-hand side of \eqref{E:ROTDEFIMPORTANTLIERADSPHERELANDRADTRACESPHERETERMS}.
	To bound the remaining terms arising from the Leibniz expansion of
	$2 \mathscr{Z}^{N-1} \Rad(\RotRadcomponent{l} \mytr \upchi^{(Small)}),$
	we use 
	inequalities
	\eqref{E:FIRSTPOINTWISEBOUNDEUCLIDEANROTATIONRADCOMPONENT},
	\eqref{E:LOWERORDERC0BOUNDEUCLIDEANROTATIONRADCOMPONENT},
	\eqref{E:POINTWISEESTIMATESFORCHIJUNKINTERMSOFOTHERVARIABLES},
	and \eqref{E:C0BOUNDCRUCIALEIKONALFUNCTIONQUANTITIES}.
	To bound the magnitude of the term
	$\mathscr{Z}^{N-1} \Rad \left(\frac{\RotRadcomponent{l}}{\rgeo}\right)$
	arising from the right-hand side of \eqref{E:ERRORROTDEFORMSPHERETRACE}
	by the right-hand side of \eqref{E:ROTDEFIMPORTANTLIERADSPHERELANDRADTRACESPHERETERMS},
	we use
	\eqref{E:ZNAPPLIEDTORGEOISNOTTOOLARGE}
	and \eqref{E:FIRSTPOINTWISEBOUNDEUCLIDEANROTATIONRADCOMPONENT}.
	To bound the magnitude of the term
	$ \angLie_{\mathscr{Z}}^{N-1} \angLie_{\Rad}
		\left\lbrace	
			\smoothfunction(\Psi) 
			 	\Psi 
			 	(\angdiffuparg{\#} x) \angdiff x
		\right\rbrace$
	arising from the right-hand side of \eqref{E:ERRORROTDEFORMSPHERETRACE}
	by the right-hand side of \eqref{E:ROTDEFIMPORTANTLIERADSPHERELANDRADTRACESPHERETERMS},
	we use the estimates
	\eqref{E:POINTWISEBOUNDPROJECTEDLIEDERIVATIVESANGDIFFCOORDINATEX},
	\eqref{E:LOWERORDERPOINTWISEBOUNDPROJECTEDLIEDERIVATIVESANGDIFFCOORDINATEX},
	\eqref{E:CRUDEPOINTWISEBOUNDSDERIVATIVESOFANGULARDEFORMATIONTENSORS},
	\eqref{E:CRUDELOWERORDERC0BOUNDDERIVATIVESOFANGULARDEFORMATIONTENSORS},
	and the bootstrap assumptions \eqref{E:PSIFUNDAMENTALC0BOUNDBOOTSTRAP}.
	To bound the magnitude of the remaining terms 
	$\angLie_{\mathscr{Z}}^{N-1} \angLie_{\Rad}
		\left\lbrace
		G_{(Frame)}
		\ginversesphere 
			\myarray
				[	\Rot_{(l)}]
				{\RotRadcomponent{l}}
			\myarray
				[\Lunit \Psi]
				{\angdiff \Psi}
		\right\rbrace$
	arising from the right-hand side of \eqref{E:ERRORROTDEFORMSPHERETRACE}
	by the right-hand side of \eqref{E:ROTDEFIMPORTANTLIERADSPHERELANDRADTRACESPHERETERMS},
	we note that these terms have essentially the same structure as the 
	$\gsphere-$dual of some of 
	the terms in \eqref{E:ERRORROTDEFORMLSPHERE}.
	Hence, the analysis in the previous paragraph yields the desired bound.
	
	We now prove \eqref{E:ROTDEFIMPORTANTANGDIVSPHERELANDANGDIFFPILRADTERMS}	
	for the term 
	$\mathscr{Z}^{N-1} \angdiv \angdeformoneformupsharparg{\Rot_{(l)}}{\Rad}
			- \left \lbrace
					\upmu \Rot_{(l)} \mathscr{Z}^{N-1} \mytr \upchi^{(Small)}
					+ \RotRadcomponent{l} \angLap \mathscr{Z}^{N-1} \upmu
				\right\rbrace.$
	We apply 
	$\mathscr{Z}^{N-1} \angdiv$
	to the $\gsphere-$dual of 
	\eqref{E:ROTDEFORMRADSPHERE}
	and apply the Leibniz rule to the terms on the right-hand side.
	We begin by addressing the first term, which is the difficult one:
	$\mathscr{Z}^{N-1} \left\lbrace(\angDarg{A} (\upmu (\ginversesphere)^{AB} \upchi_{BC}^{(Small)} \Rot_{(l)}^C) \right\rbrace.$
	When all derivatives fall on $(\ginversesphere)^{AB} \upchi_{BC}^{(Small)},$
	we use \eqref{E:ROTATIONPOINTWISENORMESTIMATE},
	\eqref{E:C0BOUNDCRUCIALEIKONALFUNCTIONQUANTITIES},
	and the first inequality in \eqref{E:DIVCHIJUNKANGDIFFTRCHIJUNKHIGHERORDERCOMMMUTOR}
	to deduce that this term is equal to
	the principal term	
	$\upmu \Rot_{(l)} \mathscr{Z}^{N-1} \mytr \upchi^{(Small)}$
	plus an error term with magnitude $\lesssim$
	the product of $\ln(\myexp + t) (1 + t)^{-1}$ and
	the right-hand side of \eqref{E:ROTDEFIMPORTANTANGDIVSPHERELANDANGDIFFPILRADTERMS}.
	We move the principal term to the left-hand side.
	From the Leibniz rule, we deduce that 
	the remaining terms in the expansion of
	$\mathscr{Z}^{N-1} \left\lbrace(\angDarg{A} (\upmu (\ginversesphere)^{AB} \upchi_{BC}^{(Small)} \Rot_{(l)}^C) \right\rbrace$
	are bounded in magnitude 
	(we are also using the identity $\angdiv \Rot_{(l)} = \frac{1}{2} \mytr  \angdeform{\Rot_{(l)}}$)
	as follows:
	\begin{align} \label{E:ROTDEFORMSPHERERADNONPRINCIPALCHIJUNKPRODUCT}
		&
		\lesssim
		\mathop{\sum_{N_1 + N_2 + N_3 + N_4 \leq N-1}}_{N_3 \leq N - 2}
			\left|
				\angLie_{\mathscr{Z}}^{N_1} \upmu
			\right|
			\left|
				\angLie_{\mathscr{Z}}^{N_2} \ginversesphere
			\right|
			\left|
				\angLie_{\mathscr{Z}}^{N_3} \angD \upchi^{(Small)}
			\right|
			\left|
				\angLie_{\mathscr{Z}}^{N_4} \Rot_{(l)}
			\right|
			\\
	& \ \ 
		+ \sum_{N_1 + N_2 + N_3 + N_4 \leq N-1}
			\left|
				\angdiff \angLie_{\mathscr{Z}}^{N_1} \upmu
			\right|
			\left|
				\angLie_{\mathscr{Z}}^{N_2} \ginversesphere
			\right|
			\left|
				\angLie_{\mathscr{Z}}^{N_3} \upchi^{(Small)}
			\right|
			\left|
				\angLie_{\mathscr{Z}}^{N_4} \Rot_{(l)}
			\right|
				\notag \\
		& \ \
		+ \sum_{N_1 + N_2 + N_3 + N_4 \leq N-1}
			\left|
				\angLie_{\mathscr{Z}}^{N_1} \upmu
			\right|
			\left|
				\angLie_{\mathscr{Z}}^{N_2} \ginversesphere
			\right|
			\left|
				\angLie_{\mathscr{Z}}^{N_3} \upchi^{(Small)}
			\right|
			\left|
				\mathscr{Z}^{N_4} \mytr  \angdeform{\Rot_{(l)}}
			\right|.
			\notag
	\end{align}
	Using inequalities
	\eqref{E:CRUDEPOINTWISEBOUNDSDERIVATIVESOFANGULARDEFORMATIONTENSORS},
		\eqref{E:CRUDELOWERORDERC0BOUNDDERIVATIVESOFANGULARDEFORMATIONTENSORS},
	\eqref{E:LIEDERIVATIVESOFROTATIONSPOINTWISEESTIMATE},
	\eqref{E:LOWERORDERLIEDERIVATIVESOFROTATIONSC0BOUND},
	\eqref{E:POINTWISEESTIMATESFORCHIJUNKINTERMSOFOTHERVARIABLES},
	\eqref{E:FUNCTIONPOINTWISEANGDINTERMSOFANGLIEO},
	\eqref{E:TYPE02TENSORANGDINTERMSOFROTATIONALLIE},
	\eqref{E:C0BOUNDCRUCIALEIKONALFUNCTIONQUANTITIES},
	and
	\eqref{E:COMMUTATORESTIMATESVECTORFIELDSACTINGONANGDTENSORS},
	we deduce that the right-hand side of 
	\eqref{E:ROTDEFORMSPHERERADNONPRINCIPALCHIJUNKPRODUCT} 
	is $\lesssim$ 
	the product of $\ln(\myexp + t) (1 + t)^{-1}$ and
	the right-hand side of 
	\eqref{E:ROTDEFIMPORTANTLIERADSPHERELANDRADTRACESPHERETERMS} 
	as desired.
	We now address the estimates corresponding to the second term on the right-hand side of \eqref{E:ROTDEFORMRADSPHERE}, 
	that is, corresponding to 
	$\mathscr{Z}^{N-1} \left\lbrace \RotRadcomponent{l} \angLap \upmu + (\angdiffuparg{\#} \RotRadcomponent{l}) \cdot \angdiff \upmu \right\rbrace.$
	Using Lemma~\ref{L:LANDRADCOMMUTEWITHANGDIFF}, 
	we see that these terms are equal to the principal term $\RotRadcomponent{l} \angLap \mathscr{Z}^{N-1} \upmu,$
	which we move to the left, plus an error term with magnitude
	\begin{align}  \label{E:ROTDEFORMSPHERERADNONPRINCIPALUPMUPRODUCT}
		&
		\lesssim
		|\RotRadcomponent{l}|
			\left|
				[\mathscr{Z}^{N-1}, \angLap] \upmu
			\right|
		+	
		\mathop{\sum_{N_1 + N_2 \leq N-1}}_{N_2 \leq N - 2}
		\left|
			\mathscr{Z}^{N_1} \RotRadcomponent{l}
		\right|
		\left|
			\mathscr{Z}^{N_2} \angLap \upmu
		\right|
			\\
	& \ \
		+ 
		\sum_{N_1 + N_2 + N_3 \leq N-1}
		\left|
			\angLie_{\mathscr{Z}}^{N_1} \ginversesphere
		\right|
		\left|
			\angdiff \mathscr{Z}^{N_2} \RotRadcomponent{l}
		\right|
		\left|
			\angdiff \mathscr{Z}^{N_3} \upmu
		\right|.
		\notag
	\end{align}
	Using inequalities
	\eqref{E:FUNCTIONPOINTWISEANGDINTERMSOFANGLIEO} 
	and \eqref{E:ANGLAPFUNCTIONPOINTWISEINTERMSOFROTATIONS},
	the estimates
	\eqref{E:FIRSTPOINTWISEBOUNDEUCLIDEANROTATIONRADCOMPONENT},
	\eqref{E:LOWERORDERC0BOUNDEUCLIDEANROTATIONRADCOMPONENT},
	\eqref{E:CRUDEPOINTWISEBOUNDSDERIVATIVESOFANGULARDEFORMATIONTENSORS},
	\eqref{E:CRUDELOWERORDERC0BOUNDDERIVATIVESOFANGULARDEFORMATIONTENSORS},
	\eqref{E:C0BOUNDCRUCIALEIKONALFUNCTIONQUANTITIES},
	and the commutator estimate \eqref{E:ANGLAPZNCOMMUTATORACTINGONFUNCTIONSSPOINTWISE}
	with $\upmu$ in the role of $f$
	and $N_2$ in the role of $N,$
	we deduce that the right-hand 
	side of \eqref{E:ROTDEFORMSPHERERADNONPRINCIPALUPMUPRODUCT}
	is $\lesssim$ 
	the product of $\ln(\myexp + t) (1 + t)^{-1}$ and
	the right-hand side of 
	\eqref{E:ROTDEFIMPORTANTLIERADSPHERELANDRADTRACESPHERETERMS}.
	To complete the proof, it remains for us to show that
	the $\mathscr{Z}^{N-1} \angdiv$ derivatives of the
	remaining terms on the right-hand side of \eqref{E:ROTDEFORMRADSPHERE}
	have magnitudes that are
	$\lesssim$ 
	the product of $\ln(\myexp + t) (1 + t)^{-1}$ and
	the right-hand side of \eqref{E:ROTDEFIMPORTANTANGDIVSPHERELANDANGDIFFPILRADTERMS}.
	To this end, we let $\xi^{\#}$ denote the remaining terms,
	which are the $\gsphere-$dual of the terms on the second line of
	\eqref{E:ROTDEFORMRADSPHERESCHEMATIC}. Our goal is to bound
	$\mathscr{Z}^{N-1} \angdiv \xi^{\#}.$ 
	To this end, we first use the Leibniz rule to deduce that
	inequality \eqref{E:LIEZNANGDIVXISHARPLEIBNIZEXPANDED} holds.
	Next, we note that the argument given in the discussion following equation \eqref{E:ROTDEFORMRADSPHERESCHEMATIC}
	implies that 
	\begin{align}  \label{E:SECONDPOINTWISEBOUNDJUNKTERMSDIVERGENCEOFROTDEFORMRADSPHERE}
	\left| 
		\angLie_{\mathscr{Z}}^{\leq N} \xi
	\right| 
		& \lesssim 
	\ln(\myexp + t)
		\left|
				\fourmyarray[\rgeo \Lunit \mathscr{Z}^{\leq N} \Psi]
					{\Rad \mathscr{Z}^{\leq N} \Psi}
					{\rgeo \angdiff \mathscr{Z}^{\leq N} \Psi} 
					{\mathscr{Z}^{\leq N} \Psi}
		\right|
		+ 
			\frac{\ln(\myexp + t)}{1 + t}
			\left|
				\myarray[\mathscr{Z}^{\leq N} (\upmu - 1)]
					{\sum_{a=1}^3 \rgeo |\mathscr{Z}^{\leq N} \Lunit_{(Small)}^a|} 
			\right|,
			\\
	\left| 
		\angLie_{\mathscr{Z}}^{\leq 12} \xi
	\right| 
	& \lesssim \frac{\ln^2(\myexp + t)}{1 + t}.
		\label{E:C0BOUNDSECONDPOINTWISEBOUNDJUNKTERMSDIVERGENCEOFROTDEFORMRADSPHERE}
	\end{align}
		From 
		\eqref{E:ONEFORMANGDINTERMSOFROTATIONALLIE},
		\eqref{E:CRUDEPOINTWISEBOUNDSDERIVATIVESOFANGULARDEFORMATIONTENSORS},
		\eqref{E:CRUDELOWERORDERC0BOUNDDERIVATIVESOFANGULARDEFORMATIONTENSORS},		
		\eqref{E:SECONDPOINTWISEBOUNDJUNKTERMSDIVERGENCEOFROTDEFORMRADSPHERE},
		and
		\eqref{E:C0BOUNDSECONDPOINTWISEBOUNDJUNKTERMSDIVERGENCEOFROTDEFORMRADSPHERE},
		we deduce that
		the first sum on the right-hand side of \eqref{E:LIEZNANGDIVXISHARPLEIBNIZEXPANDED}
		is
		$\lesssim$ 
		the product of $\ln(\myexp + t) (1 + t)^{-1}$ and
		the right-hand side of \eqref{E:ROTDEFIMPORTANTANGDIVSPHERELANDANGDIFFPILRADTERMS} as desired.
		To show that the second sum on the right-hand side of \eqref{E:LIEZNANGDIVXISHARPLEIBNIZEXPANDED}
		is
		$\lesssim$
		the product of $\ln(\myexp + t) (1 + t)^{-1}$ and
		the right-hand side of \eqref{E:ROTDEFIMPORTANTANGDIVSPHERELANDANGDIFFPILRADTERMS},
		we combine similar reasoning with 
		the commutator estimate
		\eqref{E:COMMUTATORESTIMATESVECTORFIELDSACTINGONANGDTENSORS},
		where $N_2$ is in the role of $N.$
	
	The proof of the estimate \eqref{E:ROTDEFIMPORTANTANGDIVSPHERELANDANGDIFFPILRADTERMS}	
	for the term $\mathscr{Z}^{N-1} \angdiv \angdeformoneformupsharparg{\Rot_{(l)}}{\Lunit}
	+ \Rot_{(l)} \mathscr{Z}^{N-1} \mytr \upchi^{(Small)}$
	is similar to the one
	for $\mathscr{Z}^{N-1} \angdiv \angdeformoneformupsharparg{\Rot_{(l)}}{\Rad}
			- \left \lbrace
					\upmu \Rot_{(l)} \mathscr{Z}^{N-1} \mytr \upchi^{(Small)}
					+ \RotRadcomponent{l} \angLap \mathscr{Z}^{N-1} \upmu
				\right\rbrace$
	but easier because the terms have a similar but slightly simpler structure.	
	We apply 
	$\mathscr{Z}^{N-1} \angdiv$
	to the $\gsphere-$dual of  
	right-hand side of \eqref{E:ROTDEFORMLSPHERE}
	and argue as in the previous paragraph; we omit the details.
	
	We now prove the estimate \eqref{E:ROTDEFIMPORTANTANGDIVSPHERELANDANGDIFFPILRADTERMS}	
	for 
	$\angLie_{\mathscr{Z}}^{N-1} \angdiv \angdeformfreeupdoublesharparg{\Rot_{(l)}}
			- \RotRadcomponent{l} \angdiffuparg{\#} \mathscr{Z}^{N-1} \mytr \upchi^{(Small)}.$
	We apply 
	$\angLie_{\mathscr{Z}}^{N-1} \angdiv$
	to the double $\gsphere-$dual of 
	\eqref{E:ROTDEFORMSPHERETRACEFREE}
	and apply the Leibniz rule to the terms on the right-hand side.
 	We begin by addressing the first term, which is the difficult one:
	$2 \angLie_{\mathscr{Z}}^{N-1} 
	\angDarg{B} 
	\left\lbrace 
		\RotRadcomponent{l} 
		(\ginversesphere)^{AC} (\ginversesphere)^{BD} 
		\hat{\upchi}_{CD}^{(Small)} 
	\right\rbrace.$
	When all derivatives fall on 
	$(\ginversesphere)^{BD} \hat{\upchi}_{CD}^{(Small)},$
	we use \eqref{E:LOWERORDERC0BOUNDEUCLIDEANROTATIONRADCOMPONENT}
	and the second inequality in \eqref{E:DIVCHIJUNKANGDIFFTRCHIJUNKHIGHERORDERCOMMMUTOR}
	to rewrite this term as the principal eikonal function term
	$\RotRadcomponent{l} \angdiffuparg{\#} \mathscr{Z}^{N-1} \mytr \upchi^{(Small)}$
	plus an error term with magnitude $\lesssim$
	the product of $(1 + t)^{-1}$ and
	the right-hand side of \eqref{E:ROTDEFIMPORTANTANGDIVSPHERELANDANGDIFFPILRADTERMS}.
	We move the principal term the left-hand side.
	Using Lemma~\ref{L:LANDRADCOMMUTEWITHANGDIFF}, 
	we see that the remaining terms arising from the Leibniz expansion of
	$2 \angLie_{\mathscr{Z}}^{N-1} 
	\angDarg{B} 
	\left\lbrace (\ginversesphere)^{AC} (\ginversesphere)^{BD} 
		\RotRadcomponent{l} \hat{\upchi}_{CD}^{(Small)} 
	\right\rbrace$
	are bounded in magnitude as follows:
	\begin{align} \label{E:ROTDEFORMSPHERETRFREENONPRINCIPALCHIJUNKPRODUCT}
		&
		\lesssim
		\mathop{\sum_{N_1 + N_2 + N_3 + N_4 \leq N-1}}_{N_4 \leq N - 2}
			\left|
				\angLie_{\mathscr{Z}}^{N_1} \RotRadcomponent{l} 
			\right|
			\left|
				\angLie_{\mathscr{Z}}^{N_2} \ginversesphere
			\right|
			\left|
				\angLie_{\mathscr{Z}}^{N_3} \ginversesphere
			\right|
			\left|
				\angLie_{\mathscr{Z}}^{N_4} \angD \hat{\upchi}^{(Small)}
			\right|
				\\
	& \ \ 
		+ \sum_{N_1 + N_2 + N_3 + N_4 \leq N-1}
			\left|
				\angdiff \angLie_{\mathscr{Z}}^{N_1} \RotRadcomponent{l} 
			\right|
			\left|
				\angLie_{\mathscr{Z}}^{N_2} \ginversesphere
			\right|
			\left|
				\angLie_{\mathscr{Z}}^{N_3} \ginversesphere
			\right|
			\left|
				\angLie_{\mathscr{Z}}^{N_4} \hat{\upchi}^{(Small)}
			\right|.
				\notag 
		\end{align}
		Using inequalities
		\eqref{E:FUNCTIONPOINTWISEANGDINTERMSOFANGLIEO} 
		and \eqref{E:TYPE02TENSORANGDINTERMSOFROTATIONALLIE},
		the estimates
		\eqref{E:FIRSTPOINTWISEBOUNDEUCLIDEANROTATIONRADCOMPONENT},
		\eqref{E:LOWERORDERC0BOUNDEUCLIDEANROTATIONRADCOMPONENT},
		\eqref{E:CRUDEPOINTWISEBOUNDSDERIVATIVESOFANGULARDEFORMATIONTENSORS},
		\eqref{E:CRUDELOWERORDERC0BOUNDDERIVATIVESOFANGULARDEFORMATIONTENSORS},
		\eqref{E:POINTWISEESTIMATESFORCHIJUNKINTERMSOFOTHERVARIABLES},
		and
		\eqref{E:C0BOUNDCRUCIALEIKONALFUNCTIONQUANTITIES},
		the commutator estimate \eqref{E:COMMUTATORESTIMATESVECTORFIELDSACTINGONANGDTENSORS}
		with $\xi = \hat{\upchi}^{(Small)}$ and $N_4$ in the role of $N,$
		and \eqref{E:LIEDERIVATIVESOFTRACEFREEPARTINTERMSOFLIEDERIVATIVESOFFULLTENSOR}
		with $\xi = \upchi^{(Small)},$
		we deduce that the right-hand side of 
		\eqref{E:ROTDEFORMSPHERETRFREENONPRINCIPALCHIJUNKPRODUCT}
		is $\lesssim$ 
		the product of $(1 + t)^{-1}$ and
		the right-hand side of \eqref{E:ROTDEFIMPORTANTANGDIVSPHERELANDANGDIFFPILRADTERMS}	
		as desired.
		To complete the proof, it remains for us to show that
		the $\angLie_{\mathscr{Z}}^{N-1} \angdiv$ derivatives of 
		the $\gsphere-$dual of the
		remaining terms on the right-hand side of \eqref{E:ROTDEFORMSPHERETRACEFREE}
		have magnitudes that are
		$\lesssim$ 
		the product of $(1 + t)^{-1}$ and
		the right-hand side of \eqref{E:ROTDEFIMPORTANTANGDIVSPHERELANDANGDIFFPILRADTERMS}.
	To this end, we let $\hat{\xi}^{\# \#}$ denote the remaining terms,
	which are the $\gsphere-$dual of the terms 
	in equation \eqref{E:ERRORROTDEFORMSPHERETRACEFREE}. 
	We first claim that the following bounds
	hold for 
	\begin{align} \label{E:XIFORESTIMATINGERRORROTDEFORMSPHERETRACEFREE}
		\xi := G_{(Frame)} 
			\otimes
			\myarray
				[	\Rot_{(l)}]
				{\RotRadcomponent{l}}
			\myarray
				[\Lunit \Psi]
				{\angdiff \Psi}
			+ \smoothfunction(\Psi) 
				\Psi 
				\angdiff x \otimes \angdiff x:
	\end{align}
	\begin{align}  \label{E:JUNKTERMSROTDEFANGDIVSPHERELANDANGDIFFPILRADTERMS}
		\left|
			\angLie_{\mathscr{Z}}^{\leq N} \xi
		\right|
		& \lesssim
			\left|
				\fourmyarray
					[\rgeo \Lunit \mathscr{Z}^{\leq N} \Psi]
					{\Rad \mathscr{Z}^{\leq N} \Psi}
					{\rgeo \angdiff \mathscr{Z}^{\leq N} \Psi}
					{\mathscr{Z}^{\leq N} \Psi}
			\right|
			+ \frac{1}{(1 + t)^2}
				\left| 
					\myarray
						[\mathscr{Z}^{\leq N} (\upmu - 1)]
						{\sum_{a=1}^3 \rgeo |\mathscr{Z}^{\leq N} \Lunit_{(Small)}^a|}
				\right|, 
					\\
		\left\|
			\angLie_{\mathscr{Z}}^{\leq 12} \xi
		\right\|_{C^0(\Sigma_t^u)}
		& \lesssim
			\varepsilon
			\frac{1}{1 + t}.
			\label{E:C0BOUNDLOWERORDERJUNKTERMSROTDEFANGDIVSPHERELANDANGDIFFPILRADTERMS}
	\end{align}
	To prove \eqref{E:JUNKTERMSROTDEFANGDIVSPHERELANDANGDIFFPILRADTERMS},
	we apply $\angLie_{\mathscr{Z}}^{\leq N}$ to $\xi$ and apply the 
	Leibniz rule.
	We bound the terms 
	$\angLie_{\mathscr{Z}}^M G_{(Frame)}$ with the estimates of Lemma~\ref{L:POINTWISEESTIMATESGFRAMEINTERMSOFOTHERQUANTITIES}.
	We bound the terms
	$\threemyarray[\mathscr{Z}^M \Lunit \Psi] {\angLie_{\mathscr{Z}}^M \angdiff \Psi}{\mathscr{Z}^M \Psi}$
	with Lemma~\ref{L:AVOIDINGCOMMUTING}
	and the bootstrap assumptions \eqref{E:PSIFUNDAMENTALC0BOUNDBOOTSTRAP}.
	We bound the terms
	$\mathscr{Z}^M \smoothfunction(\Psi)$ by using 
	the bootstrap assumptions \eqref{E:PSIFUNDAMENTALC0BOUNDBOOTSTRAP}
	to deduce that $\left| \mathscr{Z}^M \smoothfunction(\Psi) \right| \lesssim \left| \mathscr{Z}^{\leq M} \Psi \right| + 1.$
	We bound the terms $\angLie_{\mathscr{Z}}^M  \angdiff x$ with
	\eqref{E:POINTWISEBOUNDPROJECTEDLIEDERIVATIVESANGDIFFCOORDINATEX}, 
	and
	\eqref{E:LOWERORDERPOINTWISEBOUNDPROJECTEDLIEDERIVATIVESANGDIFFCOORDINATEX}.
	We bound the terms 
	$\mathscr{Z}^M \RotRadcomponent{l}$		
	with \eqref{E:FIRSTPOINTWISEBOUNDEUCLIDEANROTATIONRADCOMPONENT} and 
	\eqref{E:LOWERORDERC0BOUNDEUCLIDEANROTATIONRADCOMPONENT}.
	We bound the terms 
	$\angLie_{\mathscr{Z}}^M \Rot_{(l)}$
	with \eqref{E:LIEDERIVATIVESOFROTATIONSPOINTWISEESTIMATE} and 
	\eqref{E:LOWERORDERLIEDERIVATIVESOFROTATIONSC0BOUND}.
	In total, these estimates yield the desired bound \eqref{E:JUNKTERMSROTDEFANGDIVSPHERELANDANGDIFFPILRADTERMS}.
	The bound \eqref{E:C0BOUNDLOWERORDERJUNKTERMSROTDEFANGDIVSPHERELANDANGDIFFPILRADTERMS} then follows 
	from \eqref{E:JUNKTERMSROTDEFANGDIVSPHERELANDANGDIFFPILRADTERMS},
	\eqref{E:C0BOUNDCRUCIALEIKONALFUNCTIONQUANTITIES},
	and the bootstrap assumptions \eqref{E:PSIFUNDAMENTALC0BOUNDBOOTSTRAP}.
	We now recall that our goal is to bound
	$\angLie_{\mathscr{Z}}^{N-1} \angdiv \hat{\xi}^{\# \#}$
	in magnitude by
	$\lesssim$ 
	the product of $(1 + t)^{-1}$ and
	the right-hand side of \eqref{E:ROTDEFIMPORTANTANGDIVSPHERELANDANGDIFFPILRADTERMS}. 
	To proceed, we use the Leibniz rule to deduce that
	\begin{align} \label{E:LIEZNANGDIVXIDOUBLESHARPLEIBNIZEXPANDED}
		\left|
			\angLie_{\mathscr{Z}}^{N-1} \angdiv \hat{\xi}^{\# \#}
		\right|
		& \lesssim
			\sum_{N_1 + N_2 + N_3 \leq N-1}
			\left|
				\angLie_{\mathscr{Z}}^{N_1} \ginversesphere
			\right|
			\left|
				\angLie_{\mathscr{Z}}^{N_2} \ginversesphere
			\right|
			\left|
				\angD \angLie_{\mathscr{Z}}^{N_3} \hat{\xi}
			\right|	
				\\
		& \ \
			+
			\sum_{N_1 + N_2 + N_3 \leq N-1}
			\left|
				\angLie_{\mathscr{Z}}^{N_1} \ginversesphere
			\right|
			\left|
				\angLie_{\mathscr{Z}}^{N_2} \ginversesphere
			\right|
			\left|
				[\angLie_{\mathscr{Z}}^{N_3}, \angD] \hat{\xi}
			\right|.
			\notag
	\end{align}	
	Finally, to bound the right-hand side of 
	\eqref{E:LIEZNANGDIVXIDOUBLESHARPLEIBNIZEXPANDED}
	by $\lesssim$ 
	the product of $(1 + t)^{-1}$ and
	the right-hand side of \eqref{E:ROTDEFIMPORTANTANGDIVSPHERELANDANGDIFFPILRADTERMS},
	we use the estimates 
	\eqref{E:TYPE02TENSORANGDINTERMSOFROTATIONALLIE},
	\eqref{E:CRUDEPOINTWISEBOUNDSDERIVATIVESOFANGULARDEFORMATIONTENSORS},
	\eqref{E:CRUDELOWERORDERC0BOUNDDERIVATIVESOFANGULARDEFORMATIONTENSORS},
	the commutator estimate
	\eqref{E:COMMUTATORESTIMATESVECTORFIELDSACTINGONANGDTENSORS}
	with $\hat{\xi}$ in the role of $\xi$
	and $N_3$ in the role of $N,$
	\eqref{E:LIEDERIVATIVESOFTRACEFREEPARTINTERMSOFLIEDERIVATIVESOFFULLTENSOR},
	\eqref{E:JUNKTERMSROTDEFANGDIVSPHERELANDANGDIFFPILRADTERMS},
	and	
	\eqref{E:C0BOUNDLOWERORDERJUNKTERMSROTDEFANGDIVSPHERELANDANGDIFFPILRADTERMS}.		
\end{proof}

\begin{lemma} [\textbf{Identification of the important top-order terms in the derivatives of} $\deform{\rgeo \Lunit}$]
\label{L:RGEOLIMPORTANTDEFTENSORTERMS}
	Let $1 \leq N \leq 24$ be an integer. 
	Under the small-data and bootstrap assumptions 
	of Sects.~\ref{S:PSISOLVES}-\ref{S:C0BOUNDBOOTSTRAP},
	if $\varepsilon$ is sufficiently small, then
	the following estimates hold on $\mathcal{M}_{\Tboot,U_0}:$
	\begin{subequations}
	\begin{align} \label{E:RGEOLDEFIMPORTANTRADTRACEANGTERMS}
	& \left|
			\mathscr{Z}^{N-1} \Rad \mytr  \angdeform{\rgeo \Lunit}
			- 2 \rgeo \angLap \mathscr{Z}^{N-1}  \upmu 
		\right|
			\\
		& \lesssim 
		\left| 
			\fourmyarray[\rgeo \Lunit \mathscr{Z}^{\leq N} \Psi ]
				{\Rad \mathscr{Z}^{\leq N} \Psi}
				{\rgeo \angdiff \mathscr{Z}^{\leq N} \Psi}
				{\mathscr{Z}^{\leq N} \Psi}
		\right|	
		+ \frac{1}{1 + t}
			\left|
				\myarray[\mathscr{Z}^{\leq N} (\upmu - 1)]
					{\sum_{a=1}^3 \rgeo |\mathscr{Z}^{\leq N} \Lunit_{(Small)}^a|} 
			\right|,
			\notag
	\end{align}
	
	\begin{align}  \label{E:RGEOLDEFIMPORTANTLIERADSPHERERADANDAGNDIVSPHERETERMS}
		&
		\left| 
			\angLie_{\mathscr{Z}}^{N-1} \angdiffuparg{\#} \deformarg{\rgeo \Lunit}{\Lunit}{\Rad}
		\right|,
			\\
		& \left|
				\mathscr{Z}^{N-1} \angdiv \angdeformoneformupsharparg{\rgeo \Lunit}{\Rad}
				- \rgeo \angLap \mathscr{Z}^{N-1} \upmu
			\right|,
			\notag \\
		& (1 + t)
			\left|
				\angLie_{\mathscr{Z}}^{N-1} \angdiv \angdeformfreeupdoublesharparg{\rgeo \Lunit}
				-  \rgeo \angdiffuparg{\#} \mathscr{Z}^{N-1} \mytr \upchi^{(Small)}
			\right|
			\notag \\
		& \lesssim 
		\left| 
			\fourmyarray[\rgeo \Lunit \mathscr{Z}^{\leq N} \Psi ]
				{\Rad \mathscr{Z}^{\leq N} \Psi}
				{\rgeo \angdiff \mathscr{Z}^{\leq N} \Psi}
				{\mathscr{Z}^{\leq N} \Psi}
		\right|	
		+ \frac{1}{1 + t}
			\left|
				\myarray[\mathscr{Z}^{\leq N} (\upmu - 1)]
					{\sum_{a=1}^3 \rgeo |\mathscr{Z}^{\leq N} \Lunit_{(Small)}^a|} 
			\right|.
			\notag
	\end{align}
	\end{subequations}
\end{lemma}

\begin{proof}
	We first prove \eqref{E:RGEOLDEFIMPORTANTRADTRACEANGTERMS}.
	From \eqref{E:RGEOLDEFORMTRSPHERE}, we see that we have to 
	estimate $\mathscr{Z}^{N-1} \Rad(\rgeo \mytr \upchi^{(Small)}).$
	When all derivatives fall on $\mytr \upchi^{(Small)},$
	we use the third inequality in \eqref{E:TOPORDERDERIVATIVESOFANGDSQUAREDUPMUINTERMSOFCONTROLLABLE}	
	to rewrite this term as the principal eikonal function term $2 \rgeo \angLap \mathscr{Z}^{N-1}  \upmu$
	plus an error term with magnitude $\lesssim$
	the right-hand side of \eqref{E:RGEOLDEFIMPORTANTRADTRACEANGTERMS}.	
	We move the principal term to the left-hand side.
	We claim that the remaining terms 
	arising in the Leibniz expansion of
	$\mathscr{Z}^{N-1} \Rad(\rgeo \mytr \upchi^{(Small)})$
	are error terms
	with magnitudes that are  $\lesssim$
	the right-hand side of \eqref{E:RGEOLDEFIMPORTANTRADTRACEANGTERMS} as desired.
	To prove the claim,	we use
	\eqref{E:ZNAPPLIEDTORGEOISNOTTOOLARGE},
	\eqref{E:POINTWISEESTIMATESFORCHIJUNKINTERMSOFOTHERVARIABLES},
	and \eqref{E:C0BOUNDCRUCIALEIKONALFUNCTIONQUANTITIES}.
	
	We now prove \eqref{E:RGEOLDEFIMPORTANTLIERADSPHERERADANDAGNDIVSPHERETERMS} for
	$\angLie_{\mathscr{Z}}^{N-1} \angdiffuparg{\#} \deformarg{\rgeo \Lunit}{\Lunit}{\Rad}.$
	We apply $\angLie_{\mathscr{Z}}^{N-1} \angdiffuparg{\#}$ to the terms on
	the right-hand side of \eqref{E:RGEOLDEFORMLRAD}.
	From the Leibniz rule, 
	Lemma~\ref{L:LANDRADCOMMUTEWITHANGDIFF},
	\eqref{E:ZNAPPLIEDTORGEOISNOTTOOLARGE}, 
	\eqref{E:FUNCTIONPOINTWISEANGDINTERMSOFANGLIEO},
	we deduce that
	the terms of interest are bounded in magnitude by
	\begin{align}	\label{E:LEIBNIZEXPANDEDLIEZANGDIFFSHARPNRGEOLDEFORMLRAD}
		\lesssim
		\mathop{\sum_{N_1 + N_2 \leq N}}_{N_1 \leq N-1}
			\left|
				\angLie_{\mathscr{Z}}^{N_1} \ginversesphere
			\right|
			\left|
				\angLie_{\mathscr{Z}}^{N_2} \Lunit \upmu
			\right|
		+ \frac{1}{1 + t}
			\mathop{\sum_{N_1 + N_2 \leq N}}_{N_1 \leq N-1}
			\left|
				\angLie_{\mathscr{Z}}^{N_1} \ginversesphere
			\right|
			\left|
				\angLie_{\mathscr{Z}}^{N_2} (\upmu - 1)
			\right|.
	\end{align}
	To bound the right-hand side of \eqref{E:LEIBNIZEXPANDEDLIEZANGDIFFSHARPNRGEOLDEFORMLRAD} 
	by $\lesssim$ the right-hand side of \eqref{E:RGEOLDEFIMPORTANTLIERADSPHERERADANDAGNDIVSPHERETERMS},
	we use
	\eqref{E:CRUDEPOINTWISEBOUNDSDERIVATIVESOFANGULARDEFORMATIONTENSORS},
	\eqref{E:CRUDELOWERORDERC0BOUNDDERIVATIVESOFANGULARDEFORMATIONTENSORS},
	\eqref{E:LDERIVATIVECRUCICALTRANSPORTINTEQUALITIES},
	\eqref{E:C0BOUNDCRUCIALEIKONALFUNCTIONQUANTITIES},
	and \eqref{E:C0BOUNDLDERIVATIVECRUCICALEIKONALFUNCTIONQUANTITIES}.
	
	We now prove \eqref{E:RGEOLDEFIMPORTANTLIERADSPHERERADANDAGNDIVSPHERETERMS}
	for 
	$\mathscr{Z}^{N-1} \angdiv \angdeformoneformupsharparg{\rgeo \Lunit}{\Rad}
	- \rgeo \angLap \mathscr{Z}^{N-1} \upmu.$	
	We apply $\mathscr{Z}^{N-1} \angdiv$
	to the $\gsphere-$dual of the right-hand side of \eqref{E:RGEOLDEFORMRADA}.
	We begin by addressing the first term, which is the difficult one:
	$\mathscr{Z}^{N-1} (\rgeo \angLap \upmu).$
	By the Leibniz rule and \eqref{E:ZNAPPLIEDTORGEOISNOTTOOLARGE}, 
	we see that this term is equal to the principal eikonal function term
	$\rgeo \angLap \mathscr{Z}^{N-1} \upmu$ plus an error term that 
	is bounded in magnitude by 
	\begin{align}  \label{E:JUNKTERMRGEOLDEFORMLIERADSPHERERADANDAGNDIVSPHERETERMS}
		\lesssim
		(1 + t)
		\left|
			[\mathscr{Z}^{\leq N-1}, \angLap] \upmu
		\right|
		+
		(1 + t) 
		\left|
			\angLap \mathscr{Z}^{\leq N - 2} \upmu
		\right|.
	\end{align}
	From inequalities 
	\eqref{E:FUNCTIONPOINTWISEANGDINTERMSOFANGLIEO}
	and \eqref{E:ANGLAPFUNCTIONPOINTWISEINTERMSOFROTATIONS},
	the commutator estimate \eqref{E:ANGLAPZNCOMMUTATORACTINGONFUNCTIONSSPOINTWISE}
	with $N-1$ in the role of $N$ and $\upmu$ in the role of $f,$
	and the estimate \eqref{E:C0BOUNDCRUCIALEIKONALFUNCTIONQUANTITIES},
	we deduce 
	that the right-hand side of \eqref{E:JUNKTERMRGEOLDEFORMLIERADSPHERERADANDAGNDIVSPHERETERMS}
	is $\lesssim$ the right-hand side of \eqref{E:RGEOLDEFIMPORTANTLIERADSPHERERADANDAGNDIVSPHERETERMS}.
	We then bring the principal term over to the left-hand side.
	To complete the proof, it remains for us to show 
	that the $\mathscr{Z}^{N-1} \angdiv$ derivatives of 
	the $\gsphere-$dual of the
	remaining terms on the right-hand side of \eqref{E:RGEOLDEFORMRADA}
	have magnitudes that 
	are $\lesssim$ the right-hand side of \eqref{E:RGEOLDEFIMPORTANTLIERADSPHERERADANDAGNDIVSPHERETERMS}.
	To this end, we let $\xi^{\#}$ denote the remaining terms,
	which are the $\gsphere-$dual of $S_{t,u}$ one-forms of the form
	$	\xi : =	\rgeo
	 		G_{(Frame)}
			\threemyarray
				[\upmu \Lunit \Psi]
				{\Rad \Psi}
				{\upmu \angdiff \Psi}.$	
	We now claim that $\xi$ verifies the following bounds:
	\begin{align}  \label{E:JUNKTERMSBOUNDRGEOLDEFRADSPHERE}
		\left|
			\angLie_{\mathscr{Z}}^{\leq N} \xi
		\right|
		& \lesssim
			\left|
				\fourmyarray
					[\rgeo \Lunit \mathscr{Z}^{\leq N} \Psi]
					{\Rad \mathscr{Z}^{\leq N} \Psi}
					{\rgeo \angdiff \mathscr{Z}^{\leq N} \Psi}
					{\mathscr{Z}^{\leq N} \Psi}
			\right|
			+ \frac{1}{1 + t}
				\left| 
					\myarray
						[\mathscr{Z}^{\leq N} (\upmu - 1)]
						{\sum_{a=1}^3 \rgeo |\mathscr{Z}^{\leq N} \Lunit_{(Small)}^a|}
				\right|, 
					\\
		\left\|
			\angLie_{\mathscr{Z}}^{\leq 12} \xi
		\right\|_{C^0(\Sigma_t^u)}
		& \lesssim
			\varepsilon
			\frac{1}{1 + t}.
			\label{E:C0BOUNDLOWERORDERJUNKTERMSBOUNDRGEOLDEFRADSPHERE}
	\end{align}
	To prove \eqref{E:JUNKTERMSBOUNDRGEOLDEFRADSPHERE},
	we apply $\angLie_{\mathscr{Z}}^{\leq N}$ to $\xi$ and apply the 
	Leibniz rule.
	We bound the terms $\mathscr{Z}^M \rgeo$
	with \eqref{E:ZNAPPLIEDTORGEOISNOTTOOLARGE}.
	We bound the terms 
	$\angLie_{\mathscr{Z}}^M G_{(Frame)}$ with the estimates of Lemma~\ref{L:POINTWISEESTIMATESGFRAMEINTERMSOFOTHERQUANTITIES}.
	We bound the terms
	$\threemyarray
		[\mathscr{Z}^M \Lunit \Psi] 
		{\mathscr{Z}^M \Rad \Psi}
		{\angLie_{\mathscr{Z}}^M \angdiff \Psi}
	$
	with Lemma~\ref{L:AVOIDINGCOMMUTING}
	and the bootstrap assumptions \eqref{E:PSIFUNDAMENTALC0BOUNDBOOTSTRAP}.			
	We bound the terms $\mathscr{Z}^M \upmu$
	with \eqref{E:C0BOUNDCRUCIALEIKONALFUNCTIONQUANTITIES}.
	In total, these estimates yield the desired bound
	\eqref{E:JUNKTERMSBOUNDRGEOLDEFRADSPHERE}.
	The bound \eqref{E:C0BOUNDLOWERORDERJUNKTERMSBOUNDRGEOLDEFRADSPHERE} then follows 
	from \eqref{E:JUNKTERMSBOUNDRGEOLDEFRADSPHERE},
	\eqref{E:C0BOUNDCRUCIALEIKONALFUNCTIONQUANTITIES},
	and the bootstrap assumptions \eqref{E:PSIFUNDAMENTALC0BOUNDBOOTSTRAP}.			
	We now recall that our goal is to bound
	$\mathscr{Z}^{N-1} \angdiv \xi^{\#}$
	in magnitude by
	$\lesssim$ the right-hand side of \eqref{E:ROTDEFIMPORTANTANGDIVSPHERELANDANGDIFFPILRADTERMS}. 
	To proceed, we use the Leibniz rule to deduce that
	inequality \eqref{E:LIEZNANGDIVXISHARPLEIBNIZEXPANDED} holds.	
	Finally, to bound 
	the right-hand side of 
	\eqref{E:LIEZNANGDIVXISHARPLEIBNIZEXPANDED} 
	by
	$\lesssim$ the right-hand side of \eqref{E:RGEOLDEFIMPORTANTLIERADSPHERERADANDAGNDIVSPHERETERMS},
	we use the estimates 
	\eqref{E:FUNCTIONPOINTWISEANGDINTERMSOFANGLIEO},
	\eqref{E:CRUDEPOINTWISEBOUNDSDERIVATIVESOFANGULARDEFORMATIONTENSORS},
	\eqref{E:CRUDELOWERORDERC0BOUNDDERIVATIVESOFANGULARDEFORMATIONTENSORS},
	\eqref{E:COMMUTATORESTIMATESVECTORFIELDSACTINGONANGDTENSORS}
	with $N_2$ in the role of $N,$
	\eqref{E:JUNKTERMSBOUNDRGEOLDEFRADSPHERE},
	and	
	\eqref{E:C0BOUNDLOWERORDERJUNKTERMSBOUNDRGEOLDEFRADSPHERE}.
	
	Finally, we prove \eqref{E:RGEOLDEFIMPORTANTLIERADSPHERERADANDAGNDIVSPHERETERMS}
	for the term
	$\angLie_{\mathscr{Z}}^{N-1} \angdiv \angdeformfreeupdoublesharparg{\rgeo \Lunit}
				-  \rgeo \angdiffuparg{\#} \mathscr{Z}^{N-1} \mytr \upchi^{(Small)}.
	$
	We apply $\angLie_{\mathscr{Z}}^{N-1} \angdiv$ to 
	the double $\gsphere-$dual of the product on
	the right-hand side of \eqref{E:RGEOLDEFORMTRFREESPHERE}, 
	that is, 
	we apply $\angLie_{\mathscr{Z}}^{N-1} \angD_B$
	to $2 \rgeo (\ginversesphere)^{AC}(\ginversesphere)^{BD} \hat{\upchi}^{(Small)}_{CD}.$
	We then apply the Leibniz rule.
	When all derivatives fall on $(\ginversesphere)^{BD} \hat{\upchi}^{(Small)}_{CD},$
	we use the second inequality in \eqref{E:DIVCHIJUNKANGDIFFTRCHIJUNKHIGHERORDERCOMMMUTOR}
	to deduce that the product is equal to the principal eikonal function term
	$\rgeo \angdiff^A \mathscr{Z}^{N-1} \mytr \upchi^{(Small)}$ plus an error 
	term that is bounded in magnitude by $\lesssim$
	$(1 + t)^{-1}$ times the right-hand side of \eqref{E:RGEOLDEFIMPORTANTLIERADSPHERERADANDAGNDIVSPHERETERMS}
	as desired. Using \eqref{E:ZNAPPLIEDTORGEOISNOTTOOLARGE}, we deduce that
	the remaining terms arising from the Leibniz expansion of
	$\angLie_{\mathscr{Z}}^{N-1} 
	\angDarg{B} 
	\left\lbrace 
		\rgeo (\ginversesphere)^{AC} (\ginversesphere)^{BD} \hat{\upchi}_{CD}^{(Small)} 
	\right\rbrace$ 
	are bounded in magnitude as follows:
	\begin{align} \label{E:RGEOLDEFORMSPHERETRFREENONPRINCIPALCHIJUNKPRODUCT}
		&
		\lesssim
		(1 + t)
		\mathop{\sum_{N_1 + N_2 + N_3 \leq N-1}}_{N_3 \leq N - 2}
			\left|
				\angLie_{\mathscr{Z}}^{N_1} \ginversesphere
			\right|
			\left|
				\angLie_{\mathscr{Z}}^{N_2} \ginversesphere
			\right|
			\left|
				\angLie_{\mathscr{Z}}^{N_3} \angD \hat{\upchi}^{(Small)}
			\right|.
	\end{align}
		Using
		\eqref{E:TYPE02TENSORANGDINTERMSOFROTATIONALLIE},
		\eqref{E:CRUDEPOINTWISEBOUNDSDERIVATIVESOFANGULARDEFORMATIONTENSORS},
		\eqref{E:CRUDELOWERORDERC0BOUNDDERIVATIVESOFANGULARDEFORMATIONTENSORS},
		\eqref{E:POINTWISEESTIMATESFORCHIJUNKINTERMSOFOTHERVARIABLES},
		\eqref{E:C0BOUNDCRUCIALEIKONALFUNCTIONQUANTITIES},
		and the commutator estimate \eqref{E:COMMUTATORESTIMATESVECTORFIELDSACTINGONANGDTENSORS}
		with $\xi = \hat{\upchi}^{(Small)}$ and $N_3$ in the role of $N,$
		we deduce that the right-hand side of 
		\eqref{E:RGEOLDEFORMSPHERETRFREENONPRINCIPALCHIJUNKPRODUCT}
		is $\lesssim$ $(1 + t)^{-1}$ times
		the right-hand side of \eqref{E:RGEOLDEFIMPORTANTLIERADSPHERERADANDAGNDIVSPHERETERMS}	
		as desired.
	
\end{proof}

\section{Crude pointwise estimates for the below-top-order derivatives of 
\texorpdfstring{$\deform{Z}$}{the deformation tensors of the commutation vectorfields}}
\label{S:CRUDEPOINTWISEESTIMATESFORSOMEBELOWTOPORDERDEFORMATIONTENSORDERIVATIVES}

In Lemmas \ref{L:CRUDEPOINTWISEBELOWTOPORDERDEFORMATIONTENSORESTIMATES}
and \ref{L:ONEANGDERIVLOWERDEFORMATIONTENSORDERIVATIVESPOINTWISE}, 
we establish some rather crude pointwise estimates
for the frame components of the deformation tensors that appear in the 
wave equation error terms \eqref{E:DIVCURRENTPSI} and \eqref{E:DIVCURRENTLOW}
and their below-top-order derivatives. The lemmas collectively show that there are no difficult terms present in the 
below-top-order derivatives of the frame components of 
$\deform{Z}$ for $Z \in \mathscr{Z}.$ The quantities that we bound in 
Lemma~\ref{L:ONEANGDERIVLOWERDEFORMATIONTENSORDERIVATIVESPOINTWISE} 
involve covariant angular differentiation. Hence, 
in proving the corresponding estimates, we use a few additional ingredients
that are not needed in the proof of Lemma~\ref{L:CRUDEPOINTWISEBELOWTOPORDERDEFORMATIONTENSORESTIMATES}.

We now state and prove Lemma~\ref{L:CRUDEPOINTWISEBELOWTOPORDERDEFORMATIONTENSORESTIMATES}.
From the point of view of proving our sharp classical lifespan theorem,
the most important aspect of these estimates is that
the quantities in \eqref{E:LRADANDRADRADUNITC0BOUNDLOWERORDERCRUDEPOINTWISEBELOWTOPORDERDEFORMATIONTENSORESTIMATES} can experience some logarithmic growth in time while those in 
\eqref{E:C0BOUNDLOWERORDERCRUDEPOINTWISEBELOWTOPORDERDEFORMATIONTENSORESTIMATES} do not.

\begin{lemma}[\textbf{Crude pointwise estimates for the below-top-order derivatives of some frame components of}
$\deform{Z}$]
\label{L:CRUDEPOINTWISEBELOWTOPORDERDEFORMATIONTENSORESTIMATES}
Let $0 \leq N \leq 23$ be an integer
and let $Z \in \mathscr{Z}.$
Under the small-data and bootstrap assumptions 
of Sects.~\ref{S:PSISOLVES}-\ref{S:C0BOUNDBOOTSTRAP},
if $\varepsilon$ is sufficiently small, then
the following estimates hold on $\mathcal{M}_{\Tboot,U_0}:$
\begin{subequations}
\begin{align} \label{E:LRADANDRADRADUINTCRUDEPOINTWISEBELOWTOPORDERDEFORMATIONTENSORESTIMATES}
		\left|
			\threemyarray
				[\mathscr{Z}^N \deformarg{Z}{\Lunit}{\Rad}]
				{\mathscr{Z}^N \deformarg{Z}{\Rad}{\Radunit}}
				{\angLie_{\mathscr{Z}}^N \angdeformoneformupsharparg{Z}{\Rad}}
			\right|
		& \lesssim
			(1 + t)
			\left|
				\fourmyarray[\rgeo \Lunit \mathscr{Z}^{\leq N} \Psi]
					{\Rad \mathscr{Z}^{\leq N} \Psi}
					{\rgeo \angdiff \mathscr{Z}^{\leq N} \Psi} 
					{\mathscr{Z}^{\leq N} \Psi}
			\right|
		+ \left|
				\myarray[\mathscr{Z}^{\leq N+1} (\upmu-1)]
					{\sum_{a=1}^3 \rgeo |\mathscr{Z}^{\leq N+1} \Lunit_{(Small)}^a|} 
			\right|
		+ 1,
			\\
		\left\|
			\threemyarray
				[\mathscr{Z}^{\leq 11} \deformarg{Z}{\Lunit}{\Rad}]
				{\mathscr{Z}^{\leq 11} \deformarg{Z}{\Rad}{\Radunit}}
				{\angLie_{\mathscr{Z}}^{\leq 11} \angdeformoneformupsharparg{Z}{\Rad}}
		\right\|_{C^0(\Sigma_t^u)}
		& \lesssim \ln(\myexp + t),
			\label{E:LRADANDRADRADUNITC0BOUNDLOWERORDERCRUDEPOINTWISEBELOWTOPORDERDEFORMATIONTENSORESTIMATES}
\end{align}
\end{subequations}

\begin{subequations}
\begin{align} \label{E:CRUDEPOINTWISEBELOWTOPORDERDEFORMATIONTENSORESTIMATES}
		\left|
			\fivemyarray
				[\angLie_{\mathscr{Z}}^N \angdeformoneformupsharparg{Z}{\Lunit}]
				{\angLie_{\mathscr{Z}}^N \angdeformupsharparg{Z}}
				{\mathscr{Z}^N \mytr  \angdeform{Z}}
				{\angLie_{\mathscr{Z}}^N \angdeform{Z}}
				{\angLie_{\mathscr{Z}}^N \angdeformfreeupdoublesharparg{Z}}
			\right|
		& \lesssim
			\left|
				\fourmyarray[\rgeo \Lunit \mathscr{Z}^{\leq N} \Psi]
					{\Rad \mathscr{Z}^{\leq N} \Psi}
					{\rgeo \angdiff \mathscr{Z}^{\leq N} \Psi} 
					{\mathscr{Z}^{\leq N} \Psi}
			\right|
		+ \frac{1}{1 + t}
			\left|
				\myarray[\mathscr{Z}^{\leq N+1} (\upmu - 1)]
					{\sum_{a=1}^3 \rgeo |\mathscr{Z}^{\leq N+1} \Lunit_{(Small)}^a|} 
			\right|
		+ 1,
			\\
	\left\|
			\fivemyarray
				[\angLie_{\mathscr{Z}}^{\leq 11} \angdeformoneformupsharparg{Z}{\Lunit}]
				{\angLie_{\mathscr{Z}}^{\leq 11} \angdeformupsharparg{Z}}
				{\mathscr{Z}^{\leq 11} \mytr  \angdeform{Z}}
				{\angLie_{\mathscr{Z}}^N \angdeform{Z}}
				{\angLie_{\mathscr{Z}}^{\leq 11} \angdeformfreeupdoublesharparg{Z}}
			\right\|_{C^0(\Sigma_t^u)}
		& \lesssim 1.
			\label{E:C0BOUNDLOWERORDERCRUDEPOINTWISEBELOWTOPORDERDEFORMATIONTENSORESTIMATES}
\end{align}	
\end{subequations}

\end{lemma}

\begin{proof}
	The desired estimates for 
	$\angLie_{\mathscr{Z}}^N \angdeformoneformupsharparg{Z}{\Lunit},$
	$\angLie_{\mathscr{Z}}^N \angdeformoneformupsharparg{Z}{\Rad},$
	$\angLie_{\mathscr{Z}}^N \angdeform{Z},$
	and
	$\mathscr{Z}^N \mytr \angdeform{Z}$
	have already been established in
	Lemma~\ref{L:POINTWISEBOUNDSDERIVATIVESOFANGULARDEFORMATIONTENSORS}
	and Lemma~\ref{L:STUONEFORMANDVECTORFIELDSCORRESPONDINGTOCOMMUTATORDEFORMATIONTENSORS}.
	The desired estimates for $\angLie_{\mathscr{Z}}^N \angdeformfreeupdoublesharparg{Z} = (\ginversesphere)^2 \angdeformfree{Z}$
	then follow from these bounds together with
	Lemma~\ref{L:POINTWISEBOUNDSDERIVATIVESOFANGULARDEFORMATIONTENSORS},
	\eqref{E:LIEDERIVATIVESOFTRACEFREEPARTINTERMSOFLIEDERIVATIVESOFFULLTENSOR},
	\eqref{E:C0BOUNDCRUCIALEIKONALFUNCTIONQUANTITIES},
	and the bootstrap assumptions 
	\eqref{E:PSIFUNDAMENTALC0BOUNDBOOTSTRAP}.
	
	The desired estimates \eqref{E:LRADANDRADRADUINTCRUDEPOINTWISEBELOWTOPORDERDEFORMATIONTENSORESTIMATES}
	for
	$\deformarg{Z}{\Lunit}{\Rad}$
	and
	$\deformarg{Z}{\Rad}{\Radunit}$
	follow easily 
	from the identities
	\eqref{E:RADDEFORMRADRADUNIT},
	\eqref{E:RADDEFORMLRADUNIT},	
	\eqref{E:RGEOLDEFORMRADRADUNIT},
	\eqref{E:RGEOLDEFORMLRAD},	
	\eqref{E:ROTDEFORMRADRADUNIT},
	and
	\eqref{E:ROTDEFORMLRAD}
	and the estimates
	\eqref{E:ZNAPPLIEDTORGEOISNOTTOOLARGE},
	and 
	\eqref{E:LDERIVATIVECRUCICALTRANSPORTINTEQUALITIES}.
	
	The desired estimates \eqref{E:LRADANDRADRADUNITC0BOUNDLOWERORDERCRUDEPOINTWISEBELOWTOPORDERDEFORMATIONTENSORESTIMATES}
	for
	$\deformarg{Z}{\Lunit}{\Rad}$
	and
	$\deformarg{Z}{\Rad}{\Radunit}$
	then follow from the bounds \eqref{E:LRADANDRADRADUINTCRUDEPOINTWISEBELOWTOPORDERDEFORMATIONTENSORESTIMATES}
	and
	\eqref{E:C0BOUNDCRUCIALEIKONALFUNCTIONQUANTITIES},
	and the bootstrap assumptions \eqref{E:PSIFUNDAMENTALC0BOUNDBOOTSTRAP}.
	
\end{proof}

\begin{lemma} [\textbf{Pointwise estimates for the below-top-order derivatives of $\deform{Z}$ involving at least one angular derivative}]
\label{L:ONEANGDERIVLOWERDEFORMATIONTENSORDERIVATIVESPOINTWISE}
Let $2 \leq N \leq 24$ be an integer. 
Under the small-data and bootstrap assumptions 
	of Sects.~\ref{S:PSISOLVES}-\ref{S:C0BOUNDBOOTSTRAP},
	if $\varepsilon$ is sufficiently small, then 
	the following estimates hold on $\mathcal{M}_{\Tboot,U_0}:$
	\begin{subequations}
	\begin{align}  \label{E:LOWERDEFORMATIONTENSORDERIVATIVESANGDIVSPHERELPOINTWISE}
		&
		\left|
			\mathscr{Z}^{N-2} \angdiv \angdeformoneformupsharparg{\Rad}{\Lunit} 
		\right|,
			\,
		\left|
			\mathscr{Z}^{N-2} \angdiv \angdeformoneformupsharparg{\Rot_{(l)}}{\Lunit}
		\right|
			\\
		& \lesssim
			\frac{1}{1 + t}
			\left|
				\fourmyarray
					[\rgeo \Lunit \mathscr{Z}^{\leq N-1} \Psi]
					{\Rad \mathscr{Z}^{\leq N-1} \Psi}
					{\rgeo \angdiff \mathscr{Z}^{\leq N-1} \Psi}
					{\mathscr{Z}^{\leq N-1} \Psi}
			\right|
			+ \frac{1}{(1 + t)^2}
				\left| 
					\myarray
						[\mathscr{Z}^{\leq N} (\upmu - 1)]
						{\sum_{a=1}^3 \rgeo |\mathscr{Z}^{\leq N} \Lunit_{(Small)}^a|}
				\right|,
				\notag
	\end{align}
	
	\begin{align}  \label{E:C0BOUNDLOWERDEFORMATIONTENSORDERIVATIVESANGDIVSPHEREL}
		&
		\left\|
			\mathscr{Z}^{\leq 10} \angdiv \angdeformoneformupsharparg{\Rad}{\Lunit}
		\right\|_{C^0(\Sigma_t^u)},
			\,
		\left\|
			\mathscr{Z}^{\leq 10} \angdiv \angdeformoneformupsharparg{\Rot_{(l)}}{\Lunit}
		\right\|_{C^0(\Sigma_t^u)}
			\\
		& \lesssim
			\varepsilon
			\frac{\ln(\myexp + t)}{(1 + t)^2}.
			\notag
	\end{align}
	\end{subequations}
	
	Furthermore, the following estimates hold on $\mathcal{M}_{\Tboot,U_0}:$
	\begin{subequations}
	\begin{align} \label{E:LOWERDEFORMATIONTENSORDERIVATIVESANGDIVTRACEFREESPHEREPOINTWISE}
		&
		\left|
			\angLie_{\mathscr{Z}}^{N-2} \angdiv \angdeformfreeupdoublesharparg{\rgeo \Lunit}
		\right|,
			\,
		\left|
		 	\angLie_{\mathscr{Z}}^{N-2} \angdiv \angdeformfreeupdoublesharparg{\Rad}
		\right|,
			\, 
		\left|
			\angLie_{\mathscr{Z}}^{N-2} \angdiv \angdeformfreeupdoublesharparg{\Rot_{(l)}}
		\right|
			\\
		& \lesssim
			\frac{1}{1 + t}
			\left|
				\fourmyarray
					[\rgeo \Lunit \mathscr{Z}^{\leq N-1} \Psi]
					{\Rad \mathscr{Z}^{\leq N-1} \Psi}
					{\rgeo \angdiff \mathscr{Z}^{\leq N-1} \Psi}
					{\mathscr{Z}^{\leq N-1} \Psi}
			\right|
			+ \frac{1}{(1 + t)^2}
				\left| 
					\myarray
						[\mathscr{Z}^{\leq N} (\upmu - 1)]
						{\sum_{a=1}^3 \rgeo |\mathscr{Z}^{\leq N} \Lunit_{(Small)}^a|}
				\right|,
				\notag
		\end{align}
		
		\begin{align} \label{E:C0BOUNDLOWERDEFORMATIONTENSORDERIVATIVESANGDIVTRACEFREESPHEREORDERIVATIVES}
		\left\|
			\angLie_{\mathscr{Z}}^{\leq 10} \angdiv \angdeformfreeupdoublesharparg{\Rad}
		\right\|_{C^0(\Sigma_t^u)},
			\,
		\left\|
			\angLie_{\mathscr{Z}}^{\leq 10} \angdiv \angdeformfreeupdoublesharparg{\rgeo \Lunit}
		\right\|_{C^0(\Sigma_t^u)},
			\,
		\left\|
			\angLie_{\mathscr{Z}}^{\leq 10} \angdiv \angdeformfreeupdoublesharparg{\Rot_{(l)}}
		\right\|_{C^0(\Sigma_t^u)}
		& \lesssim
			\varepsilon 
			\frac{\ln(\myexp + t)}{(1 + t)^2}.
	\end{align}
	\end{subequations}
\end{lemma}

\begin{proof}
	We first prove \eqref{E:LOWERDEFORMATIONTENSORDERIVATIVESANGDIVSPHERELPOINTWISE} 
	for $\mathscr{Z}^{N-2} \angdiv \angdeformoneformupsharparg{\Rad}{\Lunit}.$
	From inequality \eqref{E:RADDEFIMPORTANTANGDIVSPHERELANDANGDIFFPILRADTERMS} with
	$N-1$ in the role of $N,$ we deduce that
	$\left|\mathscr{Z}^{N-2} \angdiv \angdeformoneformupsharparg{\Rad}{\Lunit} \right|
	\leq \left|\angLap \mathscr{Z}^{N-2} \upmu \right|$ 
	plus error terms that are $\lesssim$ the right-hand side of
	\eqref{E:LOWERDEFORMATIONTENSORDERIVATIVESANGDIVSPHERELPOINTWISE} as desired.
	To bound $\left|\angLap \mathscr{Z}^{N-2} \upmu \right|$
	by $\lesssim$ the right-hand side of
	\eqref{E:LOWERDEFORMATIONTENSORDERIVATIVESANGDIVSPHERELPOINTWISE},
	we use inequality \eqref{E:ANGLAPFUNCTIONPOINTWISEINTERMSOFROTATIONS}.
	The inequality \eqref{E:LOWERDEFORMATIONTENSORDERIVATIVESANGDIVSPHERELPOINTWISE}
	for $\mathscr{Z}^{N-2} \angdiv \angdeformoneformupsharparg{\Rot_{(l)}}{\Lunit}$
	follows similarly from the first inequality in  
	\eqref{E:ROTDEFIMPORTANTANGDIVSPHERELANDANGDIFFPILRADTERMS}
	and the fact that $\left| \Rot_{(l)} \mathscr{Z}^{N-2} \mytr \upchi^{(Small)} \right|$
	is $\lesssim$ the right-hand side of
	\eqref{E:LOWERDEFORMATIONTENSORDERIVATIVESANGDIVSPHERELPOINTWISE}, 
	which follows from \eqref{E:POINTWISEESTIMATESFORCHIJUNKINTERMSOFOTHERVARIABLES}.
	
	The inequality \eqref{E:LOWERDEFORMATIONTENSORDERIVATIVESANGDIVTRACEFREESPHEREPOINTWISE}
	for $\angLie_{\mathscr{Z}}^{N-2} \angdiv \angdeformfreeupdoublesharparg{\rgeo \Lunit}$
	follows similarly from the last inequality in \eqref{E:RGEOLDEFIMPORTANTLIERADSPHERERADANDAGNDIVSPHERETERMS}
	and the fact that $\left|\rgeo \angdiffuparg{\#} \mathscr{Z}^{N-2} \mytr \upchi^{(Small)} \right|$
	is $\lesssim$ the right-hand side of
	\eqref{E:LOWERDEFORMATIONTENSORDERIVATIVESANGDIVSPHERELPOINTWISE},
	which follows from \eqref{E:FUNCTIONPOINTWISEANGDINTERMSOFANGLIEO} 
	and
	\eqref{E:POINTWISEESTIMATESFORCHIJUNKINTERMSOFOTHERVARIABLES}.
	The inequality \eqref{E:LOWERDEFORMATIONTENSORDERIVATIVESANGDIVTRACEFREESPHEREPOINTWISE}
	for $\angLie_{\mathscr{Z}}^{N-2} \angdiv \angdeformfreeupdoublesharparg{\Rad}$
	follows similarly from the last inequality in \eqref{E:RADDEFIMPORTANTANGDIVSPHERELANDANGDIFFPILRADTERMS},
	the bound for $\left|\angdiffuparg{\#} \mathscr{Z}^{N-2} \mytr \upchi^{(Small)} \right|$ mentioned just above,
	and the inequality $\upmu \lesssim \ln(\myexp + t)$ (that is, \eqref{E:C0BOUNDCRUCIALEIKONALFUNCTIONQUANTITIES}). 
	The inequality \eqref{E:LOWERDEFORMATIONTENSORDERIVATIVESANGDIVTRACEFREESPHEREPOINTWISE}
	for $\angLie_{\mathscr{Z}}^{N-2} \angdiv \angdeformfreeupdoublesharparg{\Rot_{(l)}}$
	follows similarly from the last inequality in \eqref{E:ROTDEFIMPORTANTANGDIVSPHERELANDANGDIFFPILRADTERMS},
	the bound for $\left|\angdiffuparg{\#} \mathscr{Z}^{N-2} \mytr \upchi^{(Small)} \right|$ mentioned just above,
	and the inequality $|\RotRadcomponent{l}| \lesssim \ln(\myexp + t)$ (that is, \eqref{E:LOWERORDERC0BOUNDEUCLIDEANROTATIONRADCOMPONENT}). 
	
	The estimate \eqref{E:C0BOUNDLOWERDEFORMATIONTENSORDERIVATIVESANGDIVSPHEREL}
	then follows from 
	\eqref{E:LOWERDEFORMATIONTENSORDERIVATIVESANGDIVSPHERELPOINTWISE},
	\eqref{E:C0BOUNDCRUCIALEIKONALFUNCTIONQUANTITIES},
	and the bootstrap assumptions \eqref{E:PSIFUNDAMENTALC0BOUNDBOOTSTRAP}.
	Similarly, the estimate 
	\eqref{E:C0BOUNDLOWERDEFORMATIONTENSORDERIVATIVESANGDIVTRACEFREESPHEREORDERIVATIVES}
	then follows from
	\eqref{E:LOWERDEFORMATIONTENSORDERIVATIVESANGDIVTRACEFREESPHEREPOINTWISE},
	\eqref{E:C0BOUNDCRUCIALEIKONALFUNCTIONQUANTITIES},
	and the bootstrap assumptions \eqref{E:PSIFUNDAMENTALC0BOUNDBOOTSTRAP}.

\end{proof}

\section{Pointwise estimates for the top-order derivatives of
\texorpdfstring{$\angLie_{\Lunit} \deform{Z}$}
{the outgoing null derivative of the commutation vectorfield deformation tensors}}

In the next lemma, we establish pointwise estimates for 	
the top-order derivatives of some deformation tensor components
where one of the derivatives is in the direction $\Lunit.$
Unlike the other top-order derivatives of the deformation tensors,
top-order derivatives involving an $\Lunit$ derivative
are relatively easy to bound. In particular, to bound them in $L^2,$ 
we do not need to use the modified quantities of \ref{C:RENORMALIZEDEIKONALFUNCTIONQUANTITIES}.
	
\begin{lemma}[\textbf{Pointwise estimates for the top-order derivatives, 
involving an $\Lunit$ derivative, of $\deform{Z}$}]
\label{L:POINTWISEESTIMATESDERIVATIVESOFLDERIVATIVESOFDEFORMATIONTENSORCOMPONENTS}
Let $1 \leq N \leq 24$ be an integer. 
Under the small-data and bootstrap assumptions 
of Sects.~\ref{S:PSISOLVES}-\ref{S:C0BOUNDBOOTSTRAP},
if $\varepsilon$ is sufficiently small, 
and $Z \in \mathscr{Z}$ is a commutation vectorfield, 
then the following estimates hold on $\mathcal{M}_{\Tboot,U_0}:$
\begin{subequations}
\begin{align} \label{E:LDERIVATIVESOFDEFORMATIONTENSORCOMPONENTSPOINTWISE}
	\left|
		\threemyarray
			[\mathscr{Z}^{N-1} \Lunit \deformarg{Z}{\Lunit}{\Rad}]
			{\mathscr{Z}^{N-1} \Lunit \deformarg{Z}{\Rad}{\Radunit}}
			{\angLie_{\mathscr{Z}}^{N-1} \angLie_{\Lunit} \angdeformoneformupsharparg{Z}{\Rad}}
		\right|
	& \lesssim 
		\left| 
			\fourmyarray[\rgeo \Lunit \mathscr{Z}^{\leq N} \Psi]
				{\Rad \mathscr{Z}^{\leq N} \Psi}
				{\rgeo \angdiff \mathscr{Z}^{\leq N} \Psi}
				{\mathscr{Z}^{\leq N} \Psi}
		\right|	
		+ \frac{1}{1 + t}
			\left|
				\myarray[\mathscr{Z}^{\leq N} (\upmu - 1)]
					{\sum_{a=1}^3 \rgeo |\mathscr{Z}^{\leq N} \Lunit_{(Small)}^a|} 
			\right|,
			\\
	\threemyarray
		[\left\|
				\mathscr{Z}^{\leq 11} \Lunit \deformarg{Z}{\Lunit}{\Rad}
		 \right\|_{C^0(\Sigma_t^u)}]
		{\left\|
				\mathscr{Z}^{\leq 11} \Lunit \deformarg{Z}{\Rad}{\Radunit}
			\right\|_{C^0(\Sigma_t^u)}}
		{\left\|
				\angLie_{\mathscr{Z}}^{\leq 11} \angLie_{\Lunit} \angdeformoneformupsharparg{Z}{\Rad}
			\right\|_{C^0(\Sigma_t^u)}}
& \lesssim 
	\varepsilon
	\frac{\ln(\myexp + t)}{1 + t}.
	\label{E:LOWERORDERC0BOUNDLDERIVATIVESOFDEFORMATIONTENSORCOMPONENTS}
\end{align}
\end{subequations}
Furthermore, the following estimates hold on $\mathcal{M}_{\Tboot,U_0}:$
\begin{subequations}
\begin{align}
	\left|
		\angLie_{\mathscr{Z}}^{N-1} \angLie_{\Lunit} \angdeformoneformupsharparg{Z}{\Lunit}
	\right|
& \lesssim 
		\frac{1}{1 + t}
		\left| 
			\fourmyarray[\rgeo \Lunit \mathscr{Z}^{\leq N} \Psi ]
				{\Rad \mathscr{Z}^{\leq N} \Psi}
				{\rgeo \angdiff \mathscr{Z}^{\leq N} \Psi}
				{\mathscr{Z}^{\leq N} \Psi}
		\right|	
		+ \frac{1}{(1 + t)^2}
			\left|
				\myarray[\mathscr{Z}^{\leq N} (\upmu - 1)]
					{\sum_{a=1}^3 \rgeo |\mathscr{Z}^{\leq N} \Lunit_{(Small)}^a|} 
			\right|,
		\label{E:LDERIVATIVESOFDEFORMATIONTENSORLSPHEREPOINTWISE} \\
	\left\|
		\angLie_{\mathscr{Z}}^{\leq 11} \angLie_{\Lunit} \angdeformoneformupsharparg{Z}{\Lunit}
	\right\|_{C^0(\Sigma_t^u)}
	& \lesssim
		\varepsilon
		\frac{\ln(\myexp + t)}{(1 + t)^2}.
		\label{E:LOWERORDERC0BOUNDLDERIVATIVESOFDEFORMATIONTENSORLSPHEREPOINTWISE}
\end{align}
\end{subequations}
Finally, the following estimates hold on $\mathcal{M}_{\Tboot,U_0}:$
\begin{subequations}
\begin{align}
	\left|
		\mathscr{Z}^{N-1} \Lunit \mytr  \angdeform{Z}
	\right|
& \lesssim 
		\frac{1}{1 + t}
		\left| 
			\fourmyarray[\rgeo \Lunit \mathscr{Z}^{\leq N} \Psi ]
				{\Rad \mathscr{Z}^{\leq N} \Psi}
				{\rgeo \angdiff \mathscr{Z}^{\leq N} \Psi}
				{\mathscr{Z}^{\leq N} \Psi}
		\right|	
		+ \frac{1}{(1 + t)^2}
			\left|
				\myarray[\mathscr{Z}^{\leq N} (\upmu - 1)]
					{\sum_{a=1}^3 \rgeo |\mathscr{Z}^{\leq N} \Lunit_{(Small)}^a|} 
			\right|
	+ \frac{\ln(\myexp + t)}{(1 + t)^2},
		\label{E:LDERIVATIVESOFDEFORMATIONTENSORTRACESPHEREPOINTWISE} \\
	\left\|
		\mathscr{Z}^{\leq 11} \Lunit \mytr  \angdeform{Z}
	\right\|_{C^0(\Sigma_t^u)}
	& \lesssim
		\frac{\ln(\myexp + t)}{(1 + t)^2}.
		\label{E:LOWERORDERC0BOUNDLDERIVATIVESOFDEFORMATIONTENSORTRACESPHERE}
\end{align}
\end{subequations}

\end{lemma}	

\begin{proof}
	Once we have proved
	\eqref{E:LDERIVATIVESOFDEFORMATIONTENSORCOMPONENTSPOINTWISE},
	\eqref{E:LDERIVATIVESOFDEFORMATIONTENSORLSPHEREPOINTWISE},
	and \eqref{E:LDERIVATIVESOFDEFORMATIONTENSORTRACESPHEREPOINTWISE},
	the remaining estimates 
	\eqref{E:LOWERORDERC0BOUNDLDERIVATIVESOFDEFORMATIONTENSORCOMPONENTS},
	\eqref{E:LOWERORDERC0BOUNDLDERIVATIVESOFDEFORMATIONTENSORLSPHEREPOINTWISE},
	and \eqref{E:LOWERORDERC0BOUNDLDERIVATIVESOFDEFORMATIONTENSORTRACESPHERE}
	then follow from 
	\eqref{E:C0BOUNDCRUCIALEIKONALFUNCTIONQUANTITIES}
	and the bootstrap assumptions \eqref{E:PSIFUNDAMENTALC0BOUNDBOOTSTRAP}.
	
	It remains for us to prove
	\eqref{E:LDERIVATIVESOFDEFORMATIONTENSORCOMPONENTSPOINTWISE},
	\eqref{E:LDERIVATIVESOFDEFORMATIONTENSORLSPHEREPOINTWISE},
	and \eqref{E:LDERIVATIVESOFDEFORMATIONTENSORTRACESPHEREPOINTWISE}.
	The proofs are very similar
	to the proofs of the estimates for
	$\mathscr{Z}^N \deformarg{Z}{\Lunit}{\Rad},$
	$\mathscr{Z}^N \deformarg{Z}{\Radunit}{\Rad},$
	$\angLie_{\mathscr{Z}}^N \angdeformoneformupsharparg{Z}{\Lunit},$
	$\angLie_{\mathscr{Z}}^N \angdeformoneformupsharparg{Z}{\Rad},$
	and $\mathscr{Z}^N \mytr \angdeform{Z}$
	that we derived in Lemma~\ref{L:CRUDEPOINTWISEBELOWTOPORDERDEFORMATIONTENSORESTIMATES},
	so we only highlight the three important differences. 
	First, the fact $\rgeo \Lunit \in \mathscr{Z}$
	and inequality \eqref{E:ZNAPPLIEDTORGEOISNOTTOOLARGE}
	together imply that the estimates
	\eqref{E:LDERIVATIVESOFDEFORMATIONTENSORCOMPONENTSPOINTWISE},
	\eqref{E:LDERIVATIVESOFDEFORMATIONTENSORLSPHEREPOINTWISE},
	and \eqref{E:LDERIVATIVESOFDEFORMATIONTENSORTRACESPHEREPOINTWISE},
	which involve the operator $\mathscr{Z}^{N-1} \Lunit,$
	are better by a factor of $(1 + t)^{-1}$ than their counterpart
	estimates from Lemma~\ref{L:CRUDEPOINTWISEBELOWTOPORDERDEFORMATIONTENSORESTIMATES}.
	Second, and by far most importantly, 
	the right-hand sides of the estimates 
	\eqref{E:LDERIVATIVESOFDEFORMATIONTENSORCOMPONENTSPOINTWISE},
	\eqref{E:LDERIVATIVESOFDEFORMATIONTENSORLSPHEREPOINTWISE},
	and \eqref{E:LDERIVATIVESOFDEFORMATIONTENSORTRACESPHEREPOINTWISE}	
	involve \emph{one fewer derivatives} of $\upmu$ and $\Lunit_{(Small)}^i$
	compared to the corresponding estimates of 
	Lemma~\ref{L:CRUDEPOINTWISEBELOWTOPORDERDEFORMATIONTENSORESTIMATES}.
	This crucially important fact essentially follows from
	inequality \eqref{E:LDERIVATIVECRUCICALTRANSPORTINTEQUALITIES},
	whose right-hand side involves one fewer derivatives of
	the eikonal function quantities
	$\upmu$ and $\Lunit_{(Small)}^i$ 
	compared to the number of derivatives 
	that one would obtain by using
	\eqref{E:POINTWISEEIKONALFUNCTIONTRANSPORTBASEDINTEQUALITIES}.
	Finally, unlike the estimate \eqref{E:LRADANDRADRADUINTCRUDEPOINTWISEBELOWTOPORDERDEFORMATIONTENSORESTIMATES},
	the right-hand side of the estimates 
	\eqref{E:LDERIVATIVESOFDEFORMATIONTENSORCOMPONENTSPOINTWISE}
	and \eqref{E:LDERIVATIVESOFDEFORMATIONTENSORLSPHEREPOINTWISE}
	do not involve a constant ``$+1$'' type term.
	The reason is that the $+1$ term on the right-hand side of
	\eqref{E:LRADANDRADRADUINTCRUDEPOINTWISEBELOWTOPORDERDEFORMATIONTENSORESTIMATES}
	arises from the $2 \upmu =2 + 2 (\upmu - 1)$
	and $\upmu = 1 + (\upmu - 1)$
	terms on the right-hand sides of 
	\eqref{E:RGEOLDEFORMRADRADUNIT} 
	and \eqref{E:RGEOLDEFORMLRAD};
	the constants are therefore annihilated by the $\Lunit$ derivative on the left-hand side of
	\eqref{E:LDERIVATIVESOFDEFORMATIONTENSORCOMPONENTSPOINTWISE}.
	
\end{proof}

\section{Proof of Prop.~\ref{P:IDOFKEYDIFFICULTENREGYERRORTERMS}}
\label{S:PROOFOFPROPOSITIONPIDOFKEYDIFFICULTENREGYERRORTERMS}
We now prove Prop.~\ref{P:IDOFKEYDIFFICULTENREGYERRORTERMS}.
Our goal is to analyze $\mathscr{Z}^{N-1} (\upmu \D_{\alpha} \Jcurrent{Z}^{\alpha}[\Psi])$
for each of the $5$ commutation vectorfields $Z \in \mathscr{Z} = \lbrace \rgeo \Lunit, \Rad, \Rot_{(1)}, \Rot_{(2)}, \Rot_{(3)} \rbrace.$ 
To this end, we first use the identity \eqref{E:DIVCOMMUTATIONCURRENTDECOMPOSITION} to decompose
\begin{align*}
	\upmu \D_{\alpha} {\Jcurrent{Z}^{\alpha}[\Psi]}
	& = \mathscr{K}_{(\pi-Danger)}^{(Z)}[\Psi]
			+ \mathscr{K}_{(\pi-Cancel-1)}^{(Z)}[\Psi]
			+ \mathscr{K}_{(\pi-Cancel-2)}^{(Z)}[\Psi]
			+ \mathscr{K}_{(\pi-Elliptic)}^{(Z)}[\Psi]
				\\
	& \ \ 
			+ \mathscr{K}_{(\pi-Good)}^{(Z)}[\Psi]
			+ \mathscr{K}_{(\Psi)}^{(Z)}[\Psi]
			+ \mathscr{K}_{(Low)}^{(Z)}[\Psi].
\end{align*}	
We now separately analyze $\mathscr{Z}^{N-1}$ applied to each 
of the above $7$ terms. Throughout our analysis, we implicitly use 
the definition of $Harmless^{\leq N}$ terms; see Def.~\ref{D:HARMLESSTERMS}.

\noindent \textbf{Analysis of $\mathscr{Z}^{N-1} \mathscr{K}_{(\pi-Danger)}^{(Z)}[\Psi]$}:
We first consider the case $Z = \Rot_{(l)}.$
We will show that 
\[
\mathscr{Z}^{N-1} \mathscr{K}_{(\pi-Danger)}^{(\Rot_{(l)})}[\Psi] 
=  (\Rad \Psi) \Rot_{(l)} \mathscr{Z}^{N-1} \mytr \upchi^{(Small)} + Harmless^{\leq N},
\]
where $\mathscr{K}_{(\pi-Danger)}^{(\Rot_{(l)})}[\Psi]$ is defined by \eqref{E:DIVCURRENTTRANSVERSAL}.
By the Leibniz rule, we have to bound the following terms:
\begin{align} \label{E:DANGEROOUSTERMLEIBNIZ}
	&
	- \left(\mathscr{Z}^{N-1} \angdiv \angdeformoneformupsharparg{\Rot_{(l)}}{\Lunit} \right)
		\Rad \Psi
	-
	\mathop{\sum_{N_1 + N_2 \leq N-1}}_{N_1 \leq N - 2}
	\left(
		 \mathscr{Z}^{N_1} \angdiv \angdeformoneformupsharparg{\Rot_{(l)}}{\Lunit}
	\right) 
	\mathscr{Z}^{N_2}
	\Rad \Psi.
\end{align}
Using the first inequality in \eqref{E:ROTDEFIMPORTANTANGDIVSPHERELANDANGDIFFPILRADTERMS}
and the bootstrap assumption $\|\Rad \Psi\|_{C^0(\Sigma_t^u)} \leq \varepsilon (1 + t)^{-1}$
(that is, \eqref{E:PSIFUNDAMENTALC0BOUNDBOOTSTRAP}),
we deduce that the first term in
\eqref{E:DANGEROOUSTERMLEIBNIZ} is equal to the top-order-eikonal function-containing product
$(\Rad \Psi) \Rot_{(l)} \mathscr{Z}^{N-1} \mytr \upchi^{(Small)}$
plus a quadratic error term that is bounded in magnitude by
\begin{align}  \label{E:JUNKESTIMATEDANGEROUSROTWAVEEQUATIONINHOMOGENEOUSTERM}
	& \lesssim
			\varepsilon
			\frac{1}{(1 + t)^2}
			\left|
				\fourmyarray
					[\rgeo \Lunit \mathscr{Z}^{\leq N} \Psi]
					{\Rad \mathscr{Z}^{\leq N} \Psi}
					{\rgeo \angdiff \mathscr{Z}^{\leq N} \Psi}
					{\mathscr{Z}^{\leq N} \Psi}
			\right|
			+ \varepsilon
				\frac{1}{(1 + t)^3}
				\left| 
					\myarray
						[\mathscr{Z}^{\leq N} (\upmu - 1)]
						{\sum_{a=1}^3 \rgeo |\mathscr{Z}^{\leq N} \Lunit_{(Small)}^a|}
				\right|
				= Harmless^{\leq N}
\end{align}
as desired.
To bound the magnitude of the remaining sum in \eqref{E:DANGEROOUSTERMLEIBNIZ}
by \eqref{E:JUNKESTIMATEDANGEROUSROTWAVEEQUATIONINHOMOGENEOUSTERM},
we use the second inequalities in 
\eqref{E:LOWERDEFORMATIONTENSORDERIVATIVESANGDIVSPHERELPOINTWISE} 
and \eqref{E:C0BOUNDLOWERDEFORMATIONTENSORDERIVATIVESANGDIVSPHEREL}
and the bootstrap assumptions \eqref{E:PSIFUNDAMENTALC0BOUNDBOOTSTRAP}.

Similarly, in the case $Z = \Rad,$
we use the first inequality in 
\eqref{E:RADDEFIMPORTANTANGDIVSPHERELANDANGDIFFPILRADTERMS}
and the first inequalities in 
\eqref{E:LOWERDEFORMATIONTENSORDERIVATIVESANGDIVSPHERELPOINTWISE} 
and \eqref{E:C0BOUNDLOWERDEFORMATIONTENSORDERIVATIVESANGDIVSPHEREL}
to deduce that the analog of the expression
\eqref{E:DANGEROOUSTERMLEIBNIZ}
is equal to the top-order-eikonal function-containing product 
$(\Rad \Psi)\angLap \mathscr{Z}^{N-1} \upmu$
plus a quadratic error term that is bounded in magnitude by
the right-hand side of \eqref{E:JUNKESTIMATEDANGEROUSROTWAVEEQUATIONINHOMOGENEOUSTERM}.

In the case $Z = \rgeo \Lunit,$ we have
$\angdeformoneformupsharparg{\rgeo \Lunit}{\Lunit} = 0$
(see \eqref{E:RGEOLDEFORMLA})
and the desired estimate is trivial.

\ \\
\noindent \textbf{Analysis of $\mathscr{Z}^{N-1} \mathscr{K}_{(\pi-Cancel-1)}^{(Z)}[\Psi]$}:
We first consider the case $Z = \Rot_{(l)}.$
We will show that 
\[
	\mathscr{Z}^{N-1} \mathscr{K}_{(\pi-Cancel-1)}^{(\Rot_{(l)})}[\Psi] =  Harmless^{\leq N},
\]
where $\mathscr{K}_{(\pi-Cancel-1)}^{(\Rot_{(l)})}[\Psi]$ is defined by \eqref{E:DIVCURRENTCANEL1}.
By the Leibniz rule, we have to analyze the following sum:
\begin{align} \label{E:ROTCANCELLATIONTERMLEIBNIZ}
	\sum_{N_1 + N_2 \leq N-1}
	\left(
	\mathscr{Z}^{N_1}
	\left\lbrace 
		\frac{1}{2} \Rad \mytr  \angdeform{\Rot_{(l)}}
		- \angdiv \angdeformoneformupsharparg{\Rot_{(l)}}{\Rad}
		- \upmu \angdiv \angdeformoneformupsharparg{\Rot_{(l)}}{\Lunit}
	\right\rbrace 
	\right)
	\mathscr{Z}^{N_2}
	\Lunit \Psi.
\end{align}

From the second inequality in 
\eqref{E:ROTDEFIMPORTANTLIERADSPHERELANDRADTRACESPHERETERMS}
and the first two inequalities in
\eqref{E:ROTDEFIMPORTANTANGDIVSPHERELANDANGDIFFPILRADTERMS}
we deduce that the dangerous terms 
$\RotRadcomponent{l} \angLap \mathscr{Z}^{N_1} \upmu$
and 
$\upmu \Rot_{(l)} \mathscr{Z}^{N_1} \mytr \upchi^{(Small)}$
completely cancel out of the sum in braces in \eqref{E:ROTCANCELLATIONTERMLEIBNIZ}.
Using in addition the estimate $\upmu \lesssim \ln(\myexp + t)$
(that is, \eqref{E:C0BOUNDCRUCIALEIKONALFUNCTIONQUANTITIES}),
we deduce that for $N_1 \leq N-1,$ we have
\begin{align}   \label{E:ROTCANCEL1WAVEEQNERRORTERMPOINTWISE}
	&
	\left|
	\mathscr{Z}^{N_1}
	\left\lbrace 
		\frac{1}{2} \Rad \mytr  \angdeform{\Rot_{(l)}}
		- \angdiv \angdeformoneformupsharparg{\Rot_{(l)}}{\Rad}
		- \upmu \angdiv \angdeformoneformupsharparg{\Rot_{(l)}}{\Lunit}
	\right\rbrace 
	\right|
		\\
	& \lesssim
			\left|
				\fourmyarray
					[\rgeo \Lunit \mathscr{Z}^{\leq N_1+1} \Psi]
					{\Rad \mathscr{Z}^{\leq N_1+1} \Psi}
					{\rgeo \angdiff \mathscr{Z}^{\leq N_1+1} \Psi}
					{\mathscr{Z}^{\leq N_1+1} \Psi}
			\right|
			+ \frac{1}{1 + t}
				\left| 
					\myarray
						[\mathscr{Z}^{\leq N_1+1} (\upmu - 1)]
						{\sum_{a=1}^3 \rgeo |\mathscr{Z}^{\leq N_1+1} \Lunit_{(Small)}^a|}
				\right|.
				\notag 
\end{align}

Furthermore, 
from \eqref{E:ROTCANCEL1WAVEEQNERRORTERMPOINTWISE},
\eqref{E:C0BOUNDCRUCIALEIKONALFUNCTIONQUANTITIES},
and the bootstrap assumptions \eqref{E:PSIFUNDAMENTALC0BOUNDBOOTSTRAP},
we deduce that
\begin{align}  \label{E:LINFITYLOWERORDERROTCANCEL1WAVEEQNERRORTERMPOINTWISE}
	\left\|
	\mathscr{Z}^{\leq 11}
	\left\lbrace 
		\frac{1}{2} \Rad \mytr  \angdeform{\Rot_{(l)}}
		- \angdiv \angdeformoneformupsharparg{\Rot_{(l)}}{\Rad}
		- \upmu \angdiv \angdeformoneformupsharparg{\Rot_{(l)}}{\Lunit}
	\right\rbrace 
	\right\|_{C^0(\Sigma_t^u)}
	& \lesssim
		\varepsilon \frac{\ln(\myexp + t)}{1 + t}.
\end{align}
From 
\eqref{E:ROTCANCEL1WAVEEQNERRORTERMPOINTWISE},
\eqref{E:LINFITYLOWERORDERROTCANCEL1WAVEEQNERRORTERMPOINTWISE},
Lemma~\ref{L:AVOIDINGCOMMUTING},
and the bootstrap assumptions \eqref{E:PSIFUNDAMENTALC0BOUNDBOOTSTRAP},
we deduce that the sum \eqref{E:ROTCANCELLATIONTERMLEIBNIZ}
$= Harmless^{\leq N}$ as desired.
	
Similarly, in the case $Z = \rgeo \Lunit,$ 
we use 
\eqref{E:RGEOLDEFORMLA},
\eqref{E:RGEOLDEFIMPORTANTRADTRACEANGTERMS},
and
the first inequality in
\eqref{E:RGEOLDEFIMPORTANTLIERADSPHERERADANDAGNDIVSPHERETERMS}
to deduce that the terms
$\rgeo \angLap \mathscr{Z}^{N-1} \upmu$
cancel from the analog of the term in braces in
\eqref{E:ROTCANCELLATIONTERMLEIBNIZ}
and that
\begin{align}
	\mathscr{Z}^{N-1} \mathscr{K}_{(\pi-Cancel-1)}^{(\rgeo \Lunit)}[\Psi]
	& = Harmless^{\leq N}
\end{align}
as desired.

Similarly, in the case $Z = \Rad,$
we use \eqref{E:RADDEFORMRADA},
the second inequality in \eqref{E:RADDEFIMPORTANTLIERADSPHERELANDRADTRACESPHERETERMS},
and the first inequality in \eqref{E:RADDEFIMPORTANTANGDIVSPHERELANDANGDIFFPILRADTERMS}
to deduce that the terms
$\upmu \angLap \mathscr{Z}^{N-1} \upmu$
cancel from the analog of the term in braces in
\eqref{E:ROTCANCELLATIONTERMLEIBNIZ} and that
\begin{align}
	\mathscr{Z}^{N-1} \mathscr{K}_{(\pi-Cancel-1)}^{(\Rad)}[\Psi]
	& = Harmless^{\leq N}
\end{align}
as desired.

\ \\
\noindent \textbf{Analysis of $\mathscr{Z}^{N-1} \mathscr{K}_{(\pi-Cancel-2)}^{(Z)}[\Psi]$}
The analysis is similar to our analysis of $\mathscr{Z}^{N-1} \mathscr{K}_{(\pi-Cancel-1)}^{(Z)}[\Psi].$
We first consider the case $Z = \Rot_{(l)}.$
We will show that 
\[
	\mathscr{Z}^{N-1} \mathscr{K}_{(\pi-Cancel-2)}^{(\Rot_{(l)})}[\Psi] =  Harmless^{\leq N},
\]
where $\mathscr{K}_{(\pi-Cancel-2)}^{(\Rot_{(2)})}[\Psi]$ is defined by \eqref{E:DIVCURRENTCANEL2}.
By the Leibniz rule and Lemma~\ref{L:LEIBNIZRULEWITHANGULARDIFFERENTIALCOMMUTATIONS}, 
we have to analyze the following sum:
\begin{align} \label{E:ROTCANCELLATIONSECONDTERMLEIBNIZ}
	\sum_{N_1 + N_2 \leq N-1}
	\left(
	\angLie_{\mathscr{Z}}^{N_1}
	\left\lbrace
		- \angLie_{\Rad} \angdeformoneformupsharparg{\Rot_{(l)}}{\Lunit}
		+ \angdiffuparg{\#} \deformarg{\Rot_{(l)}}{\Lunit}{\Rad}
		\right\rbrace  
	\right)
	\cdot
	\angdiff \mathscr{Z}^{N_2} \Psi.
\end{align}
Using the first inequality in \eqref{E:ROTDEFIMPORTANTLIERADSPHERELANDRADTRACESPHERETERMS}
and the third inequality in \eqref{E:ROTDEFIMPORTANTANGDIVSPHERELANDANGDIFFPILRADTERMS},
we deduce that the dangerous terms 
$(\angD^{2 \#} \mathscr{Z}^{N-1} \upmu) \cdot \Rot_{(l)}$
completely cancel out of the sum in braces in \eqref{E:ROTCANCELLATIONSECONDTERMLEIBNIZ}.
Thanks to this cancellation, we can argue as in our analysis of the sum \eqref{E:ROTCANCELLATIONTERMLEIBNIZ}
to deduce that $\mathscr{Z}^{N-1} \mathscr{K}_{(\pi-Cancel-2)}^{(\Rot_{(l)})}[\Psi] =  Harmless^{\leq N}$
as desired.

Similarly, in the case $Z = \Rad,$ 
we use the first inequality in \eqref{E:RADDEFIMPORTANTLIERADSPHERELANDRADTRACESPHERETERMS}
and the second inequality in \eqref{E:RADDEFIMPORTANTANGDIVSPHERELANDANGDIFFPILRADTERMS}
to deduce that the terms
$\angdiffuparg{\#} \mathscr{Z}^{N-1} \Rad \upmu$
cancel from the analog of the term in braces in
\eqref{E:ROTCANCELLATIONSECONDTERMLEIBNIZ},
and a similar argument yields
\begin{align} \label{E:PICANCEL2HARMLESS}
	\mathscr{Z}^{N-1} \mathscr{K}_{(\pi-Cancel-2)}^{(\Rad)}[\Psi]
	& = Harmless^{\leq N}
\end{align}
as desired. \emph{This cancellation is critically important because
we do not have any way to estimate the top-order terms} $\angdiffuparg{\#} \mathscr{Z}^{23} \Rad \upmu.$

In the case of the vectorfield $Z = \rgeo \Lunit,$
\eqref{E:RGEOLDEFORMLA}
and the first inequality in \eqref{E:RGEOLDEFIMPORTANTLIERADSPHERERADANDAGNDIVSPHERETERMS}
together imply that there are no dangerous top-order 
eikonal function quantities
present in the analog of the term in braces in 
\eqref{E:ROTCANCELLATIONSECONDTERMLEIBNIZ},
and a similar argument yields the desired estimate
\begin{align}
	\mathscr{Z}^{N-1} \mathscr{K}_{(\pi-Cancel-2)}^{(\rgeo \Lunit)}[\Psi]
	& = Harmless^{\leq N}.
\end{align}

\ \\
\noindent \textbf{Analysis of $\mathscr{Z}^{N-1} \mathscr{K}_{(\pi-Elliptic)}^{(Z)}[\Psi]$}
We first consider the case $Z = \Rot_{(l)}.$
We will show that 
\[
	\mathscr{Z}^{N-1} \mathscr{K}_{(\pi-Elliptic)}^{(\Rot_{(l)})}[\Psi] 
	= \RotRadcomponent{l} (\angdiffuparg{\#} \Psi) \cdot (\upmu \angdiff \mathscr{Z}^{N-1} \mytr \upchi^{(Small)}) + Harmless^{\leq N},
\]
where $\mathscr{K}_{(\pi-Elliptic)}^{(\Rot_{(l)})}$ is defined by \eqref{E:DIVCURRENTELLIPTIC}.
By the Leibniz rule and Lemma~\ref{L:LEIBNIZRULEWITHANGULARDIFFERENTIALCOMMUTATIONS}, 
we have to analyze the following sum:
\begin{align} \label{E:ROTELLIPTICLEIBNIZ}
	&
	\upmu
	(\angLie_{\mathscr{Z}}^{N-1} \angdiv \angdeformfreeupdoublesharparg{\Rot_{(l)}}) 
	\cdot
	\angdiff \Psi
	+ 
	\mathop{\sum_{N_1 + N_2 + N_3 \leq N-1}}_{N_2 \leq N - 2}
	(\mathscr{Z}^{N_1} \upmu)
	(\angLie_{\mathscr{Z}}^{N_2} \angdiv \angdeformfreeupdoublesharparg{\Rot_{(l)}}) 
	\cdot
	\angdiff \mathscr{Z}^{N_3} \Psi.
\end{align}
Using the last inequality in \eqref{E:ROTDEFIMPORTANTANGDIVSPHERELANDANGDIFFPILRADTERMS},
the estimate $|\angdiff \Psi| \lesssim \varepsilon (1 + t)^{-2}$ (see \eqref{E:PSIFUNDAMENTALC0BOUNDBOOTSTRAP}),
and the bound $\upmu \lesssim \ln(\myexp + t)$ (which follows from \eqref{E:C0BOUNDCRUCIALEIKONALFUNCTIONQUANTITIES}),
we deduce that the first term in \eqref{E:ROTELLIPTICLEIBNIZ}
is equal to the top-order-eikonal function-containing product
$\RotRadcomponent{l} (\angdiff \Psi)\cdot(\upmu \angdiffuparg{\#} \mathscr{Z}^{N-1} \mytr \upchi^{(Small)})$
plus a quadratic error term that $= Harmless^{\leq N}$ as desired.

To bound the sum in \eqref{E:ROTELLIPTICLEIBNIZ},
we use the bootstrap assumptions
and the estimates \eqref{E:C0BOUNDCRUCIALEIKONALFUNCTIONQUANTITIES},
\eqref{E:LOWERDEFORMATIONTENSORDERIVATIVESANGDIVTRACEFREESPHEREPOINTWISE},
and \eqref{E:C0BOUNDLOWERDEFORMATIONTENSORDERIVATIVESANGDIVTRACEFREESPHEREORDERIVATIVES}
to deduce that it $=Harmless^{\leq N}$ as desired.

Similarly, in the case $Z = \Rad,$ 
we use the third inequality in \eqref{E:RADDEFIMPORTANTANGDIVSPHERELANDANGDIFFPILRADTERMS},
\eqref{E:C0BOUNDCRUCIALEIKONALFUNCTIONQUANTITIES},
\eqref{E:LOWERDEFORMATIONTENSORDERIVATIVESANGDIVTRACEFREESPHEREPOINTWISE},
and \eqref{E:C0BOUNDLOWERDEFORMATIONTENSORDERIVATIVESANGDIVTRACEFREESPHEREORDERIVATIVES}
to deduce that the analog of the sum
\eqref{E:ROTELLIPTICLEIBNIZ}
is equal to the top-order-eikonal function-containing product
$(\angdiffuparg{\#} \Psi) \cdot (\upmu \angdiff \mathscr{Z}^{N-1} \mytr \upchi^{(Small)})$
plus a quadratic error term that 
$= Harmless^{\leq N}$ as desired.

Similarly, in the case $Z = \rgeo \Lunit,$ 
we use the second inequality in
\eqref{E:RGEOLDEFIMPORTANTLIERADSPHERERADANDAGNDIVSPHERETERMS},
\eqref{E:C0BOUNDCRUCIALEIKONALFUNCTIONQUANTITIES},
\eqref{E:LOWERDEFORMATIONTENSORDERIVATIVESANGDIVTRACEFREESPHEREPOINTWISE},
and \eqref{E:C0BOUNDLOWERDEFORMATIONTENSORDERIVATIVESANGDIVTRACEFREESPHEREORDERIVATIVES}
to deduce that 
the analog of the sum
\eqref{E:ROTELLIPTICLEIBNIZ}
is equal to the top-order-eikonal function-containing product
$(\rgeo \angdiffuparg{\#} \Psi) \cdot (\upmu \angdiff \mathscr{Z}^{\leq N-1} \mytr \upchi^{(Small)})$
plus a quadratic error term that 
$= Harmless^{\leq N}$ as desired.

\ \\

\noindent \textbf{Analysis of $\mathscr{Z}^{N-1} \mathscr{K}_{(\pi-Good)}^{(Z)}[\Psi]$}:
We will show that 
\[	
	\mathscr{Z}^{N-1} \mathscr{K}_{(\pi-Good)}^{(Z)}[\Psi] = Harmless^{\leq N},
\]
where $\mathscr{K}_{(\pi-Good)}^{(Z)}[\Psi]$ is given by \eqref{E:DIVCURRENTGOOD}.
To begin, we use \eqref{E:DIVCURRENTGOOD} to deduce the following schematic identity:
\begin{align} \label{E:SCHEMATICWAVEEQUATIONINHOMGENEOUSTERMLDERIVATIVESOFDEFORMATIONTENSORS}
		\mathscr{K}_{(\pi-Good)}^{(Z)}[\Psi]
		& = 
		\threemyarray
			[\Lunit \mytr  \angdeform{Z}]
			{\Lunit \deformarg{Z}{\Lunit}{\Rad}}
			{\Lunit \deformarg{Z}{\Rad}{\Radunit}}
			\upmu \Lunit \Psi
		+ (\Lunit \mytr  \angdeform{Z})
			\Rad \Psi
		+ \myarray
				[\upmu \angLie_{\Lunit} \angdeformoneformupsharparg{Z}{\Lunit}] 
				{\angLie_{\Lunit} \angdeformoneformupsharparg{Z}{\Rad}}
			\angdiff \Psi.
\end{align}
We now apply $\angLie_{\mathscr{Z}}^{N-1}$ to 
\eqref{E:SCHEMATICWAVEEQUATIONINHOMGENEOUSTERMLDERIVATIVESOFDEFORMATIONTENSORS}
and then apply the Leibniz rule to the products on the right-hand side.
Using 
Lemma~\ref{L:LANDRADCOMMUTEWITHANGDIFF},
inequality \eqref{E:FUNCTIONDERIVATIVESAVOIDINGCOMMUTING} with $f = \Psi,$
the estimates \eqref{E:C0BOUNDCRUCIALEIKONALFUNCTIONQUANTITIES},
\eqref{E:LDERIVATIVESOFDEFORMATIONTENSORCOMPONENTSPOINTWISE},
\eqref{E:LOWERORDERC0BOUNDLDERIVATIVESOFDEFORMATIONTENSORCOMPONENTS},
\eqref{E:LDERIVATIVESOFDEFORMATIONTENSORLSPHEREPOINTWISE},
\eqref{E:LOWERORDERC0BOUNDLDERIVATIVESOFDEFORMATIONTENSORLSPHEREPOINTWISE},
\eqref{E:LDERIVATIVESOFDEFORMATIONTENSORTRACESPHEREPOINTWISE},
and
\eqref{E:LOWERORDERC0BOUNDLDERIVATIVESOFDEFORMATIONTENSORTRACESPHERE},
and the bootstrap assumptions \eqref{E:PSIFUNDAMENTALC0BOUNDBOOTSTRAP},
we deduce that
\begin{align} \label{E:HARMLESSBOUNDSCHEMATICWAVEEQUATIONINHOMGENEOUSTERMLDERIVATIVESOFDEFORMATIONTENSORS}
	\left|
		\mathscr{Z}^{N-1} \mathscr{K}_{(\pi-Good)}^{(Z)}[\Psi]
	\right|
	& \lesssim
		\frac{\ln(\myexp + t)}{(1 + t)^2}
		\left| 
			\fourmyarray[\rgeo \Lunit \mathscr{Z}^{\leq N} \Psi ]
				{\Rad \mathscr{Z}^{\leq N} \Psi}
				{\rgeo \angdiff \mathscr{Z}^{\leq N} \Psi}
				{\mathscr{Z}^{\leq N} \Psi}
		\right|
		+ \varepsilon
			\frac{\ln(\myexp + t)}{(1 + t)^3}
			\left|
				\myarray[\mathscr{Z}^{\leq N} (\upmu - 1)]
					{\sum_{a=1}^3 \rgeo |\mathscr{Z}^{\leq N} \Lunit_{(Small)}^a|} 
			\right|
			\notag
			\\
& = Harmless^{\leq N} \notag
\end{align}
as desired. 

\ \\

\noindent \textbf{Analysis of $\mathscr{Z}^{N-1} \mathscr{K}_{(\Psi)}^{(Z)}[\Psi]$}:
We will show that 
\[ 
	\mathscr{Z}^{N-1} \mathscr{K}_{(\Psi)}^{(Z)}[\Psi] = Harmless^{\leq N},
\]
where $\mathscr{K}_{(\Psi)}^{(Z)}[\Psi]$ is given by \eqref{E:DIVCURRENTPSI}.
To this end, we apply $\mathscr{Z}^{N-1}$ to both sides of \eqref{E:DIVCURRENTPSI}.
We first address the terms that arise from the first line of the right-hand 
side of \eqref{E:DIVCURRENTPSI}. Using the Leibniz rule,
\eqref{E:ZNAPPLIEDTORGEOISNOTTOOLARGE},
the estimates
\eqref{E:C0BOUNDCRUCIALEIKONALFUNCTIONQUANTITIES},
\eqref{E:LRADANDRADRADUINTCRUDEPOINTWISEBELOWTOPORDERDEFORMATIONTENSORESTIMATES},
\eqref{E:LRADANDRADRADUNITC0BOUNDLOWERORDERCRUDEPOINTWISEBELOWTOPORDERDEFORMATIONTENSORESTIMATES},
\eqref{E:CRUDEPOINTWISEBELOWTOPORDERDEFORMATIONTENSORESTIMATES},
\eqref{E:C0BOUNDLOWERORDERCRUDEPOINTWISEBELOWTOPORDERDEFORMATIONTENSORESTIMATES},
\eqref{E:PSIVERSIONCOMMUTINGFUNCTIONSWITHLPLUSHALFTRACECHI},
and \eqref{E:PSIVERSIONLOWERORDERC0BOUNDCOMMUTINGFUNCTIONSWITHLPLUSHALFTRACECHI},
and the bootstrap assumptions \eqref{E:PSIFUNDAMENTALC0BOUNDBOOTSTRAP},
we bound the terms of interest by
\begin{align} \label{E:PSIHIGHTERMSREQUIRINGMORAWETZBOUND}
	& \lesssim 
	\left|
		\left\lbrace
			\Lunit + \frac{1}{2} \mytr \upchi 
		\right\rbrace
		\mathscr{Z}^{\leq N} \Psi
	\right|
	+ 	\frac{\ln(\myexp + t)}{(1 + t)^2}
			\left|
				\fourmyarray[\rgeo \Lunit \mathscr{Z}^{\leq N} \Psi]
					{\Rad \mathscr{Z}^{\leq N} \Psi}
					{\rgeo \angdiff \mathscr{Z}^{\leq N} \Psi} 
					{\mathscr{Z}^{\leq N} \Psi}
			\right|
		+ \varepsilon
			\frac{1}{(1 + t)^3}
			\left|
				\myarray[\mathscr{Z}^{\leq N} (\upmu - 1)]
					{\sum_{a=1}^3 \rgeo |\mathscr{Z}^{\leq N} \Lunit_{(Small)}^a|} 
			\right| 
			\\
		& = Harmless^{\leq N}
			\notag
\end{align}
as desired.

We now address the terms that arise from the second and third lines of the right-hand 
side of \eqref{E:DIVCURRENTPSI}. Using the Leibniz rule,
\eqref{E:ZNAPPLIEDTORGEOISNOTTOOLARGE},
the estimates
\eqref{E:C0BOUNDCRUCIALEIKONALFUNCTIONQUANTITIES},
\eqref{E:LRADANDRADRADUINTCRUDEPOINTWISEBELOWTOPORDERDEFORMATIONTENSORESTIMATES},
\eqref{E:LRADANDRADRADUNITC0BOUNDLOWERORDERCRUDEPOINTWISEBELOWTOPORDERDEFORMATIONTENSORESTIMATES},
\eqref{E:CRUDEPOINTWISEBELOWTOPORDERDEFORMATIONTENSORESTIMATES},
\eqref{E:C0BOUNDLOWERORDERCRUDEPOINTWISEBELOWTOPORDERDEFORMATIONTENSORESTIMATES},
\eqref{E:POINTWISECOMMUTINGANGDSQUAREDPSIWITHLIEZN},
and
\eqref{E:LINFFTYLOWERORDERCOMMUTINGANGDSQUAREDPSIWITHLIEZN},
and the bootstrap assumptions \eqref{E:PSIFUNDAMENTALC0BOUNDBOOTSTRAP},
we bound the terms of interest by
\begin{align} \label{E:PSIHIGHTERMSSLIGHTLYEASIER}
	&
	\lesssim 
		\frac{\ln(\myexp + t)}{(1 + t)^2}
			\left|
				\fourmyarray[\rgeo \Lunit \mathscr{Z}^{\leq N} \Psi]
					{\Rad \mathscr{Z}^{\leq N} \Psi}
					{\rgeo \angdiff \mathscr{Z}^{\leq N} \Psi} 
					{\mathscr{Z}^{\leq N} \Psi}
			\right|
		+ \varepsilon
			\frac{1}{(1 + t)^3}
			\left|
				\myarray[\mathscr{Z}^{\leq N} (\upmu - 1)]
					{\sum_{a=1}^3 \rgeo |\mathscr{Z}^{\leq N} \Lunit_{(Small)}^a|} 
			\right| 
			\\
		& = Harmless^{\leq N}
			\notag
\end{align}
as desired.

\ \\

\noindent \textbf{Analysis of $\mathscr{Z}^{N-1} \mathscr{K}_{(Low)}^{(Z)}[\Psi]$}:
We will show that 
\[
	\mathscr{Z}^{N-1} \mathscr{K}_{(Low)}^{(Z)}[\Psi] = Harmless^{\leq N},
\]
where $\mathscr{K}_{(Low)}^{(Z)}[\Psi]$ is given by \eqref{E:DIVCURRENTLOW}.
To this end, we apply $\angLie_{\mathscr{Z}}^{N-1}$ to both sides of \eqref{E:DIVCURRENTLOW}
and then apply the Leibniz rule to terms on the right-hand side.
We claim that all resulting products are in magnitude $\lesssim$ the right-hand
side of \eqref{E:PSIHIGHTERMSSLIGHTLYEASIER} and hence are $Harmless^{\leq N}$
as desired. We bound all derivatives of all deformation tensors with
the estimates of Lemma~\ref{L:CRUDEPOINTWISEBELOWTOPORDERDEFORMATIONTENSORESTIMATES}.
We bound the quantities $\mathscr{Z}^M \rgeo^{-1}$
with \eqref{E:ZNAPPLIEDTORGEOISNOTTOOLARGE}.
To bound
$\myarray[\mathscr{Z}^M \Lunit \Psi] {\angLie_{\mathscr{Z}}^M\angdiff \Psi}$
and $\mathscr{Z}^M \Psi,$
we use Lemma~\ref{L:AVOIDINGCOMMUTING}
and the bootstrap assumptions \eqref{E:PSIFUNDAMENTALC0BOUNDBOOTSTRAP}.
To bound $\mathscr{Z}^M \upmu$ and
$\angLie_{\mathscr{Z}}^M \angdiff \upmu,$
we use Lemma~\ref{L:LANDRADCOMMUTEWITHANGDIFF},
inequality \eqref{E:FUNCTIONPOINTWISEANGDINTERMSOFANGLIEO},
and \eqref{E:C0BOUNDCRUCIALEIKONALFUNCTIONQUANTITIES}.
To bound $\angLie_{\mathscr{Z}}^M \angdiffuparg{\#} \upmu,$
we use Lemma~\ref{L:LANDRADCOMMUTEWITHANGDIFF},
inequality \eqref{E:FUNCTIONPOINTWISEANGDINTERMSOFANGLIEO},
\eqref{E:C0BOUNDCRUCIALEIKONALFUNCTIONQUANTITIES},
and the estimates of Lemma~\ref{L:POINTWISEBOUNDSDERIVATIVESOFANGULARDEFORMATIONTENSORS}
(to bound the derivatives of $\ginversesphere$).
To bound $\mathscr{Z}^M \mytr \upchi^{(Small)},$
we use \eqref{E:POINTWISEESTIMATESFORCHIJUNKINTERMSOFOTHERVARIABLES}
and \eqref{E:C0BOUNDCRUCIALEIKONALFUNCTIONQUANTITIES}.
The quantities
$\Lunit \upmu,$
$\mytr  \angkuparg{(Trans-\Psi)},$
$\upmu \mytr  \angkuparg{(Tan-\Psi)},$
$\upzeta^{(Trans-\Psi)\#}$
and $\upmu \upzeta^{(Tan-\Psi)\#}$
are all schematically of the form 
$G_{(Frame)}
\ginversesphere
\threemyarray[\upmu \Lunit \Psi]
{\Rad \Psi}
{\upmu \angdiff \Psi}$
(see 
\eqref{E:UPMUFIRSTTRANSPORT},
\eqref{E:KABTRANSVERSAL},
\eqref{E:KABGOOD},
\eqref{E:ZETATRANSVERSAL},
and
\eqref{E:ZETAGOOD}).
Hence, we can bound their derivatives by using
the estimates of Lemma~\ref{L:POINTWISEESTIMATESGFRAMEINTERMSOFOTHERQUANTITIES},
the previously mentioned estimates for $\mathscr{Z}^M \upmu,$
$\angLie_{\mathscr{Z}}^M \ginversesphere,$
and
$\myarray[\mathscr{Z}^M \Lunit \Psi] {\angLie_{\mathscr{Z}}^M \angdiff \Psi},$
and the bootstrap assumptions \eqref{E:PSIFUNDAMENTALC0BOUNDBOOTSTRAP}.
In total, these estimates yield that the terms of interest are
$\lesssim$ the right-hand
side of \eqref{E:PSIHIGHTERMSSLIGHTLYEASIER}
as desired.

$\hfill \qed$

\section{Proof of Cor.~\ref{C:POINTWISEESTIMATESFOREASYCOMMUTATORTERMS}}
\label{S:PROOFOFCOROLLARYPOINTWISEESTIMATESFOREASYCOMMUTATORTERMS}
We now prove Cor.~\ref{C:POINTWISEESTIMATESFOREASYCOMMUTATORTERMS}.
We first consider inequality \eqref{E:LOWERORDERINHOMOGENEOUSTERMSFIRSTPOINTWISE}
in the case $N_1 \leq 11.$ Then from the estimate 
$\left\| \mathscr{Z}^{\leq 11} \mytr  \angdeform{Z} \right\|_{C^0(\Sigma_t^u)} \lesssim 1$
(that is, \eqref{E:CRUDELOWERORDERC0BOUNDDERIVATIVESOFANGULARDEFORMATIONTENSORS}),
we deduce that
\begin{align} \label{E:EASIERTERMLOWERORDERINHOMOGENEOUSTERMSSECONDPOINTWISE}
	\left| 
		\inhomarg{\mathscr{Z}^N}_{(Eikonal-Low)}
	\right|
		& \lesssim 
			\mathop{\sum_{N_2 + N_3 \leq N-1}}_{N_2 \leq N - 2} 
								\sum_{Z \in \mathscr{Z}}
								\left| 
									\mathscr{Z}^{N_2} (\upmu \D_{\alpha} \Jcurrent{Z}^{\alpha}[\mathscr{Z}^{N_3} \Psi]) 
								\right|.
\end{align}
Our goal is to show that 
\begin{align} \label{E:LOWERTERMSAREHARMLESSLEQN}
	\mbox{the right-hand side of \eqref{E:EASIERTERMLOWERORDERINHOMOGENEOUSTERMSSECONDPOINTWISE} is} = Harmless^{\leq N}.
\end{align}
To prove \eqref{E:LOWERTERMSAREHARMLESSLEQN}, we begin by repeating the proof of 
Prop.~\ref{P:IDOFKEYDIFFICULTENREGYERRORTERMS} with 
$N-1$ replaced by $N_2$ 
(where $N_2 \leq N - 2$)
and $\Psi$ replaced by $\mathscr{Z}^{N_3} \Psi.$
The same proof yields that all terms in \eqref{E:EASIERTERMLOWERORDERINHOMOGENEOUSTERMSSECONDPOINTWISE} 
are $Harmless^{\leq N_2}$ (and hence they are also $Harmless^{\leq N}$)
except for the terms 
corresponding to the explicitly written terms in
\eqref{E:RADISTHEFIRSTCOMMUTATORIMPORTANTTERMS}-\eqref{E:ROTISTHEFIRSTCOMMUTATORIMPORTANTTERMS},
that is, 
for terms of the form
$(\Rad \mathscr{Z}^{N_3} \Psi) \angLap \mathscr{Z}^{N_2} \upmu,$
$(\Rad \mathscr{Z}^{N_3} \Psi) \Rot \mathscr{Z}^{N_2} \mytr \upchi^{(Small)},$
$(\upmu \angdiffuparg{\#} \mathscr{Z}^{N_3} \Psi) \cdot (\upmu \angdiff \mathscr{Z}^{N_2} \mytr \upchi^{(Small)}),$
and
$\RotRadcomponent{l}(\angdiffuparg{\#} \mathscr{Z}^{N_3} \Psi) \cdot (\upmu \angdiff \mathscr{Z}^{N_2} \mytr \upchi^{(Small)}).$
To handle the first two of these terms, we use
the bootstrap assumptions \eqref{E:PSIFUNDAMENTALC0BOUNDBOOTSTRAP},
inequalities \eqref{E:FUNCTIONPOINTWISEANGDINTERMSOFANGLIEO} and \eqref{E:ANGLAPFUNCTIONPOINTWISEINTERMSOFROTATIONS},
the pointwise estimate \eqref{E:POINTWISEESTIMATESFORCHIJUNKINTERMSOFOTHERVARIABLES}, 
inequality \eqref{E:C0BOUNDCRUCIALEIKONALFUNCTIONQUANTITIES},
and the fact that $N_2 \leq N - 2$
to deduce that
\begin{align} \label{E:LOWERORDERHARMLESS1}
	\mathop{\sum_{N_2 + N_3 \leq N-1}}_{N_2 \leq N - 2} 
	\left|
		(\Rad \mathscr{Z}^{\leq N_3} \Psi) \angLap \mathscr{Z}^{\leq N_2} \upmu
	\right|
	& \lesssim 
		\varepsilon
		\frac{\ln(\myexp + t)}{(1 + t)^2} 
		\left|
			\Rad \mathscr{Z}^{\leq N} \Psi
		\right|
		+
		\varepsilon
		\frac{1}{(1 + t)^3} 
		\left|
			\mathscr{Z}^{\leq N} (\upmu - 1)
		\right|,
	\\
	\mathop{\sum_{N_2 + N_3 \leq N-1}}_{N_2 \leq N - 2}
	\left|
		(\Rad \mathscr{Z}^{N_3} \Psi) \Rot \mathscr{Z}^{N_2} \mytr \upchi^{(Small)}
	\right|
	& \lesssim
		\varepsilon
		\frac{\ln(\myexp + t)}{(1+t)^2}
		\left|
				\fourmyarray[\rgeo \Lunit \mathscr{Z}^{\leq N} \Psi]
					{\Rad \mathscr{Z}^{\leq N} \Psi}
					{\rgeo \angdiff \mathscr{Z}^{\leq N} \Psi} 
					{\mathscr{Z}^{\leq N} \Psi}
		\right|
				\label{E:LOWERORDERHARMLESS2} \\
	& \ \ 
		+ \varepsilon
			\frac{1}{(1+t)^3}
			\left|
				\myarray[\mathscr{Z}^{\leq N-1} (\upmu - 1)]
					{\sum_{a=1}^3 \rgeo |\mathscr{Z}^{\leq N} \Lunit_{(Small)}^a|} 
			\right|.
			\notag
\end{align}
Hence, by Def.~\ref{D:HARMLESSTERMS},
the terms in \eqref{E:LOWERORDERHARMLESS1}-\eqref{E:LOWERORDERHARMLESS2}
are $= Harmless^{\leq N}.$ 
Also using the pointwise estimates 
$\upmu \lesssim \ln(\myexp + t)$ 
and 
$|\RotRadcomponent{l}| \lesssim \varepsilon \ln(\myexp + t),$
we use a similar argument to deduce that
\begin{align}
\mathop{\sum_{N_2 + N_3 \leq N-1}}_{N_2 \leq N - 2}
\left| (\upmu \angdiffuparg{\#} \mathscr{Z}^{N_3} \Psi) (\upmu \angdiff \mathscr{Z}^{N_2} \mytr \upchi^{(Small)}) \right|
	& = Harmless^{\leq N}, \\
\mathop{\sum_{N_2 + N_3 \leq N-1}}_{N_2 \leq N - 2}
\left| \RotRadcomponent{l}(\angdiffuparg{\#} \mathscr{Z}^{N_3} \Psi) (\upmu \angdiff \mathscr{Z}^{N_2} \mytr \upchi^{(Small)}) \right|
	& = Harmless^{\leq N}.
\end{align}
We have thus proved \eqref{E:LOWERTERMSAREHARMLESSLEQN}.

It remains for us to consider inequality \eqref{E:LOWERORDERINHOMOGENEOUSTERMSFIRSTPOINTWISE}
in the case that $N_1 \geq 12$ and thus $N_2 + N_3 \leq 11.$
The proof in the previous paragraph yields that
$\mathscr{Z}^{N_2} (\upmu \D_{\alpha} \Jcurrent{Z}^{\alpha}[\mathscr{Z}^{N_3} \Psi])
= Harmless^{\leq 12}.$ Hence, by Def.~\ref{D:HARMLESSTERMS},
the bootstrap assumptions \eqref{E:PSIFUNDAMENTALC0BOUNDBOOTSTRAP},
the estimate $|\mytr \upchi| \lesssim (1 + t)^{-1}$ (that is, \eqref{E:CRUDELOWERORDERC0BOUNDDERIVATIVESOFANGULARDEFORMATIONTENSORS}),
and inequality \eqref{E:C0BOUNDCRUCIALEIKONALFUNCTIONQUANTITIES},
we have that
\begin{align} \label{E:LOWERORDERDERIVATIVESUPMUWEIGHTEDDIVCURRENTC0BOUND}
	\sum_{N_2 + N_3 \leq 11}
	\left\|
		\mathscr{Z}^{N_2} (\upmu \D_{\alpha} \Jcurrent{Z}^{\alpha}[\mathscr{Z}^{N_3} \Psi])
	\right\|_{C^0(\Sigma_t^u)}
	& \lesssim \varepsilon \frac{1}{(1 + t)^2}.
\end{align}
Actually, by treating the first term on the right-hand side of \eqref{E:HARMLESSTERMS} with 
more refined arguments like those used in the proof of Lemma~\ref{L:FASTERTHANEXPECTEDPSIDECAY},
we could show that the right-hand side of \eqref{E:LOWERORDERDERIVATIVESUPMUWEIGHTEDDIVCURRENTC0BOUND} 
can be improved to $\varepsilon \ln(\myexp + t)(1 + t)^{-3},$ but we have no need for the improvement here.
Then from \eqref{E:LOWERORDERDERIVATIVESUPMUWEIGHTEDDIVCURRENTC0BOUND}, 
the estimates 
\eqref{E:CRUDEPOINTWISEBOUNDSDERIVATIVESOFANGULARDEFORMATIONTENSORS}
and \eqref{E:CRUDELOWERORDERC0BOUNDDERIVATIVESOFANGULARDEFORMATIONTENSORS}
for the derivatives of $\mytr \angdeform{Z},$ 
\eqref{E:LOWERTERMSAREHARMLESSLEQN},
and Def.~\ref{D:HARMLESSTERMS},
we deduce that
\begin{align} \label{E:MOREANNOYINGTERMLOWERORDERINHOMOGENEOUSTERMSSECONDPOINTWISE}
	& \mathop{\mathop{\sum_{N_1 + N_2 + N_3 \leq N-1}}_{N_1, N_2 \leq N - 2}}_{N_2 + N_3 \leq 11}
	\sum_{Z_1, Z_2 \in \mathscr{Z}}
		\left| 
		 \mathscr{Z}^{N_1} \mytr  \angdeform{Z} 
		\right|
		\left| 
			\mathscr{Z}^{N_2} (\upmu \D_{\alpha} \Jcurrent{Z}^{\alpha}[\mathscr{Z}^{N_3} \Psi]) 
		\right|
			\\
	& \lesssim 
			\varepsilon
			\frac{1}{(1 + t)^2}
			\left|
				\fourmyarray[\rgeo \Lunit \mathscr{Z}^{\leq N-1} \Psi]
					{\Rad \mathscr{Z}^{\leq N-1} \Psi}
					{\rgeo \angdiff \mathscr{Z}^{\leq N-1} \Psi} 
					{\mathscr{Z}^{\leq N-1} \Psi}
			\right|
		+ \varepsilon
			\frac{1}{(1 + t)^3}
			\left| 
				\myarray
					[\mathscr{Z}^{\leq N} (\upmu - 1)]
						{\rgeo \sum_{a=1}^3 |\mathscr{Z}^{\leq N} \Lunit_{(Small)}^a|}
				\right|
	 + Harmless^{\leq N}
				\notag \\
	& = Harmless^{\leq N}
		\notag
\end{align}
as desired.
$\hfill \qed$

\section{Pointwise estimates for the error integrands involving 
\texorpdfstring{$\deform{\Mult} \ \mbox{and} \ \deform{\Mor}$}
{the deformation tensors of the multiplier vectorfields}}
\label{S:POINTWISEESTIMATESFORBASICENERGYERRORINTEGRANDS}
We have now derived all of the necessary pointwise estimates corresponding
to the inhomogeneous terms in the commuted wave equation.
However, in order to derive our desired $L^2$ estimates, 
we also have to derive pointwise estimates for the error integrands
from Def.~\ref{D:LK0K1ERRORINTEGRANDS} and Lemma~\ref{L:LK0K1ERRORINTEGRANDS}.
We derive the desired pointwise estimates in the next proposition.

\begin{proposition}[\textbf{Pointwise estimates for ${\basicenergyerror{\Mult}[\Psi]}$ and ${\basicenergyerror{\Mor}[\Psi]}$}]
	\label{P:POINTWISEESTIMATESBASICERRORINTEGRANDS}
	Let $\basicenergyerror{\Mult}[\Psi]$ and $\basicenergyerror{\Mor}[\Psi]$ be the error integrands from
	Def.~\ref{D:LK0K1ERRORINTEGRANDS} and Lemma~\ref{L:LK0K1ERRORINTEGRANDS}.	
	There exists a constant $C > 0$ 
	such that under the small-data and bootstrap assumptions 
	of Sects.~\ref{S:PSISOLVES}-\ref{S:C0BOUNDBOOTSTRAP},
	if $\varepsilon$ is sufficiently small, 
	then we have the following upper bound for $(t',u',\vartheta) \in \mathcal{M}_{\Tboot,U_0}$
	(without taking the absolute value on the left-hand side),
	where we view both the left and right-hand sides as functions of
	$(t',u',\vartheta)$ and $t$ is any fixed time verifying $t' \leq t < \Tboot:$
	\begin{subequations}
	\begin{align} \label{E:MULTERRORINTEGRANDPOINTWISE}
		\basicenergyerror{\Mult}[\Psi]
		& \leq	
			C \frac{\ln(\myexp + t')}{(1 + t')^2} \Psi^2
			+ C \ln(\myexp + t')\left\lbrace\Lunit \Psi + \frac{1}{2} \mytr \upchi \Psi \right\rbrace^2
			+ C \frac{1}{(1 + t')^2} (\Rad \Psi)^2
			\\
		& \ \ + C \ln^2(\myexp + t') \upmu |\angdiff \Psi|^2
				+ C
					\varepsilon^{1/2}
					\frac{\ln(\myexp + t')}{\sqrt{\ln(\myexp + t) - \ln(\myexp + t')}}
					\upmu |\angdiff \Psi|^2
					\notag \\
		& \ \ + C \varepsilon 
						\ln^2(\myexp + t')
						\mathbf{1}_{\lbrace \upmu \leq 1/4 \rbrace} 
						|\angdiff \Psi|^2.
				\notag
	\end{align}	
	
	In addition, the following pointwise estimate holds
	for $(t',u',\vartheta) \in \mathcal{M}_{\Tboot,U_0}:$
	\begin{align} \label{E:MORAWETZERRORINTEGRANDPOINTWISE}
			\left|
				\basicenergyerror{\Mor}[\Psi]
				+ \frac{1}{4} 
					\rgeo^2
					|\angdiff \Psi|^2 
					[\Lunit \upmu]_-
				\right| 
		& \leq C \frac{\ln^3(\myexp + t')}{1 + t'} \Psi^2
					+ C \ln^3(\myexp + t') (1 + t') 
					\left\lbrace
						\Lunit \Psi + \frac{1}{2} \mytr \upchi \Psi 
					\right\rbrace^2
					\\
		& \ \ + \frac{1}{2} 
						(1 + C \varepsilon)
						\left\lbrace
							\frac{1}{\rgeo(t',u') \left\lbrace 1 + \ln \left(\frac{\rgeo(t',u')}{\rgeo(0,u')} \right) \right\rbrace}
						\right\rbrace
						\rgeo^2(t',u') 
						\upmu |\angdiff \Psi|^2
						\notag \\
		& \ \ + C \varepsilon 
						\frac{(1 + t')}{\ln(\myexp + t')} 
						\mathbf{1}_{\lbrace \upmu \leq 1/4 \rbrace} 
						|\angdiff \Psi|^2.
			\notag
	\end{align}
	\end{subequations}
		
\end{proposition}

\begin{proof}
Throughout this proof, we freely
use the trivial bound $|\uLgood f| \lesssim \upmu |\Lunit f| + |\Rad f|$
without mentioning it each time. 
We also silently use the simple inequality
$(\Lunit f)^2 \leq \left\lbrace\Lunit f + \frac{1}{2} \mytr \upchi f \right\rbrace^2
+ C (1 + t')^{-2} f^2,$
which follows easily from the estimate 
$|\rgeo \mytr \upchi| \lesssim 1$
(that is, \eqref{E:CRUDELOWERORDERC0BOUNDDERIVATIVESOFANGULARDEFORMATIONTENSORS}).
Furthermore, we silently use simple inequalities of the form
$|f_1 f_2| \lesssim h f_1^2 + h^{-1} f_2^2,$
where $h > 0$ is allowed to depend on time.

We first prove \eqref{E:MULTERRORINTEGRANDPOINTWISE}.
We separately bound each term $\basicenergyerrorarg{\Mult}{i}[\Psi]$ in the sum \eqref{E:MULTAWETZENERGYERRORINTEGRANDS}.
From \eqref{E:MULTERRORINTEGRAND1}, 
\eqref{E:C0BOUNDCRUCIALEIKONALFUNCTIONQUANTITIES},
and \eqref{E:C0BOUNDLDERIVATIVECRUCICALEIKONALFUNCTIONQUANTITIES},
we deduce that 
$\left|\basicenergyerrorarg{\Mult}{1}[\Psi]\right| \lesssim 
\varepsilon \ln(\myexp + t')(\Lunit \Psi)^2,$
which is easily seen to be $\leq$ the right-hand side of
\eqref{E:MULTERRORINTEGRANDPOINTWISE}. 
To bound the term $\basicenergyerrorarg{\Mult}{2}[\Psi]$ from
\eqref{E:MULTERRORINTEGRAND2}, we first use
\eqref{E:C0BOUNDLDERIVATIVECRUCICALEIKONALFUNCTIONQUANTITIES}, 
to deduce that
$|\angdiff \Psi|^2 |3 \upmu \Lunit \upmu| \lesssim \varepsilon (1+t')^{-1} \upmu |\angdiff \Psi|^2$ as desired.
We then use \eqref{E:POSITIVEPARTOFLMUPLUSRADMUOVERMUISSMALL} to deduce that the remaining product in \eqref{E:MULTERRORINTEGRAND2}
verifies 
$\frac{1}{2} (\uLgood \upmu + \Lunit \upmu) |\angdiff \Psi|^2
\leq C \varepsilon^{1/2}	
			\frac{\ln(\myexp + t')}{\sqrt{\ln(\myexp + t) - \ln(\myexp + t')}}
			 \upmu |\angdiff \Psi|^2 
		+ C \varepsilon \ln(\myexp + t') \upmu |\angdiff \Psi|^2$
as desired
(we do \emph{not} take the absolute value on the left-hand side because of the nature of the estimate \eqref{E:POSITIVEPARTOFLMUPLUSRADMUOVERMUISSMALL}).
To bound the term $\basicenergyerrorarg{\Mult}{4}[\Psi]$ from
\eqref{E:MULTERRORINTEGRAND4}, we first note that the terms in braces multiplying
$(\Lunit \Psi) (\angdiffuparg{\#} \Psi)$ are schematically of the form 
$\angdiff \upmu,$ 
$\upmu \angdiff \upmu,$
or $(1 + \upmu) G_{(Frame)} \threemyarray[\upmu \Lunit \Psi]{\Rad \Psi}{\upmu \angdiff \Psi}.$
Hence, from 
inequality \eqref{E:FUNCTIONPOINTWISEANGDINTERMSOFANGLIEO},
Lemma~\ref{L:POINTWISEESTIMATESGFRAMEINTERMSOFOTHERQUANTITIES},
the estimate \eqref{E:C0BOUNDCRUCIALEIKONALFUNCTIONQUANTITIES},
and the bootstrap assumptions \eqref{E:PSIFUNDAMENTALC0BOUNDBOOTSTRAP},
we deduce that the terms in braces multiplying $(\Lunit \Psi) (\angdiffuparg{\#} \Psi)$ 
are in magnitude 
$\lesssim \varepsilon \ln^2(\myexp + t') (1 + t')^{-1}.$
When $\upmu \leq 1/4,$ it thus
follows that 
$\left|\basicenergyerrorarg{\Mult}{4}[\Psi]\right| 
\lesssim \varepsilon \mathbf{1}_{\lbrace \upmu \leq 1/4 \rbrace} \ln^2(\myexp + t') (1 + t')^{-1} (\Lunit \Psi)^2
+ \varepsilon \mathbf{1}_{\lbrace \upmu \leq 1/4 \rbrace} \ln^2(\myexp + t') (1 + t')^{-1} |\angdiff \Psi|^2$
as desired. When $\upmu > 1/4,$ it thus follows that
$\left|\basicenergyerrorarg{\Mult}{4}[\Psi]\right| 
\lesssim \varepsilon \ln^2(\myexp + t') (1 + t')^{-1}(\Lunit \Psi)^2 + \varepsilon \ln^2(\myexp + t') (1 + t')^{-1} \upmu |\angdiff \Psi|$
as desired.
To bound the term $\basicenergyerrorarg{\Mult}{3}[\Psi]$ from
\eqref{E:MULTERRORINTEGRAND3}, we first note that 
the same reasoning we used in analyzing
$\basicenergyerrorarg{\Mult}{4}[\Psi]$ 
yields that the terms in braces multiplying 
$-(\uLgood \Psi) (\angdiffuparg{\#} \Psi)$ are in magnitude 
$\lesssim \varepsilon \ln(\myexp + t') (1 + t')^{-1}.$ 
Hence, using the bound $\upmu \lesssim \ln(\myexp + t')$ (that is, \eqref{E:C0BOUNDCRUCIALEIKONALFUNCTIONQUANTITIES}),
we deduce that 
$\left|\basicenergyerrorarg{\Mult}{3}[\Psi]\right| 
\lesssim \varepsilon \ln^2(\myexp + t') (1 + t')^{-1}(\Lunit \Psi)^2 
+ \varepsilon \ln(\myexp + t') (1 + t')^{-1} \upmu |\angdiff \Psi|^2
+ \varepsilon \ln(\myexp + t') (1 + t')^{-1} |\Rad \Psi| |\angdiff \Psi|.$
The first two terms are $\lesssim$ the right-hand side of \eqref{E:MULTERRORINTEGRANDPOINTWISE} as desired.
Furthermore, by separately considering the cases $\upmu \leq 1/4$ and $\upmu > 1/4$
(as we did at end of the proof of the bound for $\basicenergyerrorarg{\Mult}{4}[\Psi]$),
we deduce that
$\varepsilon \ln(\myexp + t') (1 + t')^{-1} |\Rad \Psi| |\angdiff \Psi|
	\lesssim 
		\varepsilon (1 + t')^{-2} (\Rad \Psi)^2
	+ \varepsilon \mathbf{1}_{\lbrace \upmu \leq 1/4 \rbrace} \ln^2(\myexp + t') |\angdiff \Psi|^2
	+ \ln^2(\myexp + t') \upmu |\angdiff \Psi|^2
$
as desired.
To bound the term $\basicenergyerrorarg{\Mult}{5}[\Psi]$ from
\eqref{E:MULTERRORINTEGRAND5}, we first 
note that the terms in braces are schematically of the form
$\hat{\upchi}^{(Small)}$
or the trace-free part of
$G_{(Frame)} \threemyarray[\upmu \Lunit \Psi]{\Rad \Psi}{\upmu \angdiff \Psi}.$
Hence, using the estimate $|\hat{\xi}| \lesssim |\xi|$ for symmetric type $\binom{0}{2}$ tensors,
the argument used to bound the terms in braces in $\basicenergyerrorarg{\Mult}{4}[\Psi],$
and \eqref{E:C0BOUNDCRUCIALEIKONALFUNCTIONQUANTITIES}, we deduce that the terms in
braces are in magnitude $\lesssim \varepsilon (1+t')^{-1}.$
It follows that $\left|\basicenergyerrorarg{\Mult}{5}[\Psi]\right| \lesssim \varepsilon (1+t')^{-1} \upmu |\angdiff \Psi|^2$
as desired.
To bound the term $\basicenergyerrorarg{\Mult}{6}[\Psi]$ from
\eqref{E:MULTERRORINTEGRAND6}, we first 
note that the terms in braces 
multiplying
$
- \frac{1}{2}
	(\Lunit \Psi)
	(\uLgood \Psi)
$
are schematically of the form
$\mytr \upchi$
plus the trace (with respect to $\gsphere$) of
$G_{(Frame)} \threemyarray[\upmu \Lunit \Psi]{\Rad \Psi}{\upmu \angdiff \Psi}.$
Hence, using the fact that $|\mytr \xi| \lesssim |\xi|$ for symmetric type $\binom{0}{2}$ tensors,
the estimate $|\rgeo \mytr \upchi| \lesssim 1$ mentioned at the beginning of the proof,
and the above bound for the terms in braces in $\basicenergyerrorarg{\Mult}{5}[\Psi],$
we deduce that the terms in braces are in magnitude $\lesssim (1+t')^{-1}.$
Also using the estimate $\upmu \lesssim \ln(\myexp + t')$ (that is, \eqref{E:C0BOUNDCRUCIALEIKONALFUNCTIONQUANTITIES}),
we deduce that
$\left|\basicenergyerrorarg{\Mult}{6}[\Psi] \right| 
\lesssim (\Lunit \Psi)^2 + (1+t')^{-2} (\Rad \Psi)^2,$
which can easily be seen to be $\leq$ the right-hand side of \eqref{E:MULTERRORINTEGRANDPOINTWISE} 
as desired. 

We now prove \eqref{E:MORAWETZERRORINTEGRANDPOINTWISE}.
We separately bound each term $\basicenergyerrorarg{\Mor}{i}[\Psi]$ in the sum \eqref{E:MORERRORINT}.
To bound the term $\basicenergyerrorarg{\Mor}{1}[\Psi]$
from \eqref{E:MORERRORINTEGRAND1}, we first use
the identity $\rgeo \uLgood \rgeo = \rgeo (\upmu - 2)$ and
the estimates \eqref{E:C0BOUNDCRUCIALEIKONALFUNCTIONQUANTITIES} 
and \eqref{E:C0BOUNDLDERIVATIVECRUCICALEIKONALFUNCTIONQUANTITIES}
to deduce
that the terms in braces multiplying $(\Lunit \Psi)^2$ are in magnitude $\lesssim \ln(\myexp + t')(1 + t').$
The desired bound now readily follows.
To bound the term $\basicenergyerrorarg{\Mor}{2}[\Psi]$
from \eqref{E:MORERRORINTEGRAND2},
we first use
\eqref{E:C0BOUNDCRUCIALEIKONALFUNCTIONQUANTITIES}
and \eqref{E:POSITIVEPARTOFLMUOVERMUISSMALL} to deduce that the terms in 
braces multiplying
$
\frac{1}{2}
\rgeo^2
\upmu 
(
\left|
	\angdiff \Psi 
\right|^2
)
$
are in magnitude
$\leq (1 + C \varepsilon) \frac{1}{\rgeo(t',u') \left\lbrace 1 + \ln \left(\frac{\rgeo(t',u')}{\rgeo(0,u')} \right) \right\rbrace}.$
It follows that 
$\left| \basicenergyerrorarg{\Mor}{2}[\Psi] \right|
\leq \frac{1}{2}  (1 + C \varepsilon) \frac{1}{\rgeo(t',u') \left\lbrace 1 + \ln \left(\frac{\rgeo(t',u')}{\rgeo(0,u')} \right) \right\rbrace}
\rgeo^2(t',u') \upmu |\angdiff \Psi|^2$ as desired.
To bound the term $\basicenergyerrorarg{\Mor}{3}[\Psi]$
from \eqref{E:MORERRORINTEGRAND3},
we first use essentially the same analysis 
that we used in analyzing $\basicenergyerrorarg{\Mult}{4}[\Psi]$
to deduce that the terms in braces multiplying
$\rgeo^2 (\Lunit \Psi)(\angdiffuparg{\#} \Psi$) are in magnitude $\lesssim \varepsilon \ln(\myexp + t') (1 + t')^{-1}.$
It follows that
$\left| 
	\basicenergyerrorarg{\Mor}{2}[\Psi] 
\right|
\lesssim
\varepsilon \ln^3(\myexp + t') (1 + t') (\Lunit \Psi)^2
+ \varepsilon \frac{1 + t'}{\ln(\myexp + t')} |\angdiff \Psi|^2.
$
The first term is easily seen to be 
$\leq$ the right-hand side of \eqref{E:MORAWETZERRORINTEGRANDPOINTWISE}.
To bound the second term 
$\varepsilon \frac{1 + t'}{\ln(\myexp + t')} |\angdiff \Psi|^2$ 
by the right-hand side of \eqref{E:MORAWETZERRORINTEGRANDPOINTWISE},
we use the same reasoning that we used at end of the proof of the bound for $\basicenergyerrorarg{\Mult}{4}[\Psi],$
where we separately considered the cases $\upmu \leq 1/4$ and $\upmu > 1/4$
(in the case $\upmu > 1/4,$ 
$\varepsilon \frac{1 + t'}{\ln(\myexp + t')} |\angdiff \Psi|^2$
is bounded by the term on the second line of the
right-hand side of \eqref{E:MORAWETZERRORINTEGRANDPOINTWISE}).
To bound the last term $\basicenergyerrorarg{\Mor}{4}[\Psi]$
from \eqref{E:MORERRORINTEGRAND4}, 
we first use 
\eqref{E:C0BOUNDCRUCIALEIKONALFUNCTIONQUANTITIES}
to deduce that 
$|\hat{\upchi}^{(Small)}| \lesssim \varepsilon \ln(\myexp + t') (1+t')^{-2}.$
Using this bound and the fact that $|\hat{\xi}| \lesssim |\xi|$ for symmetric type $\binom{0}{2}$ tensors,
we conclude that 
$\left| \basicenergyerrorarg{\Mor}{4}[\Psi] \right| 
\lesssim \varepsilon \ln(\myexp + t') \upmu |\angdiff \Psi|^2,$
which is $\leq$ the term on the second line of the
right-hand side of \eqref{E:MORAWETZERRORINTEGRANDPOINTWISE}
as desired. 

\end{proof}

In order to derive our desired a priori estimates for the energy quantities
$\totonemax{\leq N},$ we still have to derive pointwise estimates for
two more terms: the remaining two error integrands
present on the right-hand side of \eqref{E:E1DIVID}. 
We derive these simple estimates in the next lemma.

\begin{lemma}[\textbf{Pointwise estimates for the remaining easy energy error integrands}]
	Under the small-data and bootstrap assumptions 
	of Sects.~\ref{S:PSISOLVES}-\ref{S:C0BOUNDBOOTSTRAP},
	if $\varepsilon$ is sufficiently small, then
	the following estimates hold on $\mathcal{M}_{\Tboot,U_0}:$
	\begin{subequations}
	\begin{align} \label{E:EASYERRORINTEGRANDONEPOINTWISE}
		\left|
			\uLgood [\rgeo^2 \mytr \upchi]
				- \frac{1}{2} \rgeo^2 \upmu (\mytr \upchi)^2
				+ \rgeo^2 \mytr \upchi \mytr  \angkuparg{(Trans-\Psi)}
				+ \rgeo^2 \upmu \mytr \upchi \angkuparg{(Tan-\Psi)}
		\right|
		\Psi^2
		& \lesssim \ln(\myexp + t) \Psi^2,	
			\\
		\left| \upmu \square_{g(\Psi)} [\rgeo^2 \mytr \upchi] \right| \Psi^2 
		& \lesssim
			\frac{1}{1 + t}
			\Psi^2.
			\label{E:EASYERRORINTEGRANDTWOPOINTWISE}
	\end{align}
	\end{subequations}
	
\end{lemma}

\begin{proof}
	We first prove \eqref{E:EASYERRORINTEGRANDONEPOINTWISE}.
	Using \eqref{E:CHIJUNKDEF} and the identity 
	$\uLgood \rgeo = \upmu - 2,$ we compute that
	$ \uLgood [\rgeo^2 \mytr \upchi]
		= 2(\upmu - 2)
		+ \upmu \Lunit [\rgeo^2 \mytr \upchi^{(Small)}]
	- 4 \rgeo \mytr \upchi^{(Small)}
	+ 2 \rgeo^2 \Rad \mytr \upchi^{(Small)}.$
	Hence, from \eqref{E:C0BOUNDCRUCIALEIKONALFUNCTIONQUANTITIES}
	and \eqref{E:C0BOUNDLDERIVATIVECRUCICALEIKONALFUNCTIONQUANTITIES},
	we deduce that the magnitude of this term is $\lesssim \ln(\myexp + t)$ as desired.
	From 
	\eqref{E:CRUDELOWERORDERC0BOUNDDERIVATIVESOFANGULARDEFORMATIONTENSORS}
	and
	\eqref{E:C0BOUNDCRUCIALEIKONALFUNCTIONQUANTITIES},
	we deduce that the second term in the absolute value
	verifies 
	$\left|
		\rgeo^2 \upmu (\mytr \upchi)^2
	 \right| 
	 \lesssim \ln(\myexp + t)$
	as desired. The third and fourth terms are schematically of the form
	$\rgeo^2 \mytr \upchi G_{(Frame)} \ginversesphere \threemyarray[\upmu \Lunit \Psi]{\Rad \Psi}{\upmu \angdiff \Psi}.$
	Hence, 
	by Lemma~\ref{L:POINTWISEESTIMATESGFRAMEINTERMSOFOTHERQUANTITIES},
	the estimates \eqref{E:CRUDELOWERORDERC0BOUNDDERIVATIVESOFANGULARDEFORMATIONTENSORS}
	and \eqref{E:C0BOUNDCRUCIALEIKONALFUNCTIONQUANTITIES},
	and the bootstrap assumptions \eqref{E:PSIFUNDAMENTALC0BOUNDBOOTSTRAP},
	these terms are in magnitude $\lesssim \varepsilon$ as desired.
	
	To deduce \eqref{E:EASYERRORINTEGRANDTWOPOINTWISE}, we use	
	the wave operator decomposition identity \eqref{E:LONOUTSIDEGEOMETRICWAVEOPERATORFRAMEDECOMPOSED} with
	$\rgeo^2 \mytr \upchi = 2 \rgeo + \rgeo^2 \mytr \upchi^{(Small)}$ in the role of $f$
	and hence $\Lunit f = 2 + \Lunit(\rgeo^2 \mytr \upchi^{(Small)})$
	and $\Rad f = - 2 - 2 \rgeo \mytr \upchi^{(Small)} + \rgeo^2 \Rad \mytr \upchi^{(Small)}.$
	Up to an overall minus sign, the first term on the right-hand side of 
	\eqref{E:LONOUTSIDEGEOMETRICWAVEOPERATORFRAMEDECOMPOSED} 
	is equal to $\upmu \Lunit \Lunit f + 2 \Lunit \Rad f + (\Lunit \upmu) \Lunit f,$
	which in the present case can be expressed as
	\begin{align} \label{E:FIRSTTERMWAVEOPERATORDCOMPERGEOSQUAREDTRCHISMALL}
		&
		\upmu \Lunit \Lunit(\rgeo^2 \mytr \upchi^{(Small)})
		- 4 \rgeo^{-1} \Lunit (\rgeo^2 \mytr \upchi^{(Small)})
		+ 4 \mytr \upchi^{(Small)}
			\\
		& \ \ 
			+ 2 \Lunit(\rgeo^2 \mytr  \Rad \upchi^{(Small)})
			+ (\Lunit \upmu) \left\lbrace 2 + \Lunit(\rgeo^2 \mytr \upchi^{(Small)}) \right\rbrace.
			\notag
	\end{align}
	From \eqref{E:C0BOUNDCRUCIALEIKONALFUNCTIONQUANTITIES} 
	and \eqref{E:C0BOUNDLDERIVATIVECRUCICALEIKONALFUNCTIONQUANTITIES},
	we deduce that all terms in \eqref{E:FIRSTTERMWAVEOPERATORDCOMPERGEOSQUAREDTRCHISMALL}
	are in magnitude $\lesssim \varepsilon (1+t)^{-1}$
	as desired.
	The second term on the first line of the right-hand side of 
	\eqref{E:LONOUTSIDEGEOMETRICWAVEOPERATORFRAMEDECOMPOSED} is equal to
	$\rgeo^2 \upmu \angLap \mytr \upchi^{(Small)}.$
	Hence, from \eqref{E:FUNCTIONPOINTWISEANGDINTERMSOFANGLIEO},
	\eqref{E:ANGLAPFUNCTIONPOINTWISEINTERMSOFROTATIONS},
	and \eqref{E:C0BOUNDCRUCIALEIKONALFUNCTIONQUANTITIES},
	we deduce that its magnitude is $\lesssim \varepsilon \ln^2(\myexp + t) (1 + t)^{-2}$
	as desired.
	The term $\mytr \upchi \Rad f$ on the second line of \eqref{E:LONOUTSIDEGEOMETRICWAVEOPERATORFRAMEDECOMPOSED}
	is equal to $(2 \rgeo^{-1} + \mytr \upchi^{(Small)})(-2 - 2 \rgeo \mytr \upchi^{(Small)} + \rgeo^2 \Rad \mytr \upchi^{(Small)}).$
	Hence, from \eqref{E:C0BOUNDCRUCIALEIKONALFUNCTIONQUANTITIES}, we deduce that
	this term is in magnitude $\lesssim (1+t)^{-1}$ as desired.
	The last four terms on the right-hand side of \eqref{E:LONOUTSIDEGEOMETRICWAVEOPERATORFRAMEDECOMPOSED}
	are schematically of the form
	\[
	G_{(Frame)} \ginversesphere \threemyarray[\upmu \Lunit \Psi]{\Rad \Psi}{\angdiff \Psi} 
	\myarray[2 + \Lunit (\rgeo^2 \mytr \upchi^{(Small)})]
		{\rgeo^2 \angdiff \mytr \upchi^{(Small)}}. 
	\]
	Hence, by 
	inequality \eqref{E:FUNCTIONPOINTWISEANGDINTERMSOFANGLIEO},
	Lemma~\ref{L:POINTWISEESTIMATESGFRAMEINTERMSOFOTHERQUANTITIES},
	the estimates \eqref{E:CRUDELOWERORDERC0BOUNDDERIVATIVESOFANGULARDEFORMATIONTENSORS}
	\eqref{E:C0BOUNDCRUCIALEIKONALFUNCTIONQUANTITIES},
	and \eqref{E:C0BOUNDLDERIVATIVECRUCICALEIKONALFUNCTIONQUANTITIES},
	and the bootstrap assumptions \eqref{E:PSIFUNDAMENTALC0BOUNDBOOTSTRAP},
	these terms are in magnitude $\lesssim \varepsilon (1+t)^{-1}$ as desired.	
	
\end{proof}

\section{Pointwise estimates needed to close the elliptic estimates}
\label{S:POINTWISEFORCLOSINGTHEELLIPTICESTIMATES}

In Lemma~\ref{L:TOPORDERELLIPTICRECOVERY}, we derive our main top-order elliptic estimates
for $\upmu$ and $\upchi^{(Small)}.$ In the next lemma, we provide a collection of 
preliminary pointwise estimates that play a role in our derivation of the elliptic estimates.

\begin{lemma}[\textbf{{Pointwise estimates needed to close the elliptic estimates}}]
\label{L:POINTWISECOMMUTEDANGDCHIJUNKINTERMSOFTRACEFREEANDTRACEPART}
	Let $1 \leq N \leq 24$ be an integer. 
	Under the small-data and bootstrap assumptions 
of Sects.~\ref{S:PSISOLVES}-\ref{S:C0BOUNDBOOTSTRAP},
if $\varepsilon$ is sufficiently small, then
the following estimates hold on $\mathcal{M}_{\Tboot,U_0}:$
	\begin{align} \label{E:POINTWISECOMMUTEDANGDCHIJUNKINTERMSOFTRACEFREEANDTRACEPART}
		&
		\left|
			\angD \angLie_{\mathscr{Z}}^{N-1} \upchi^{(Small)}
			-
			\left\lbrace
				\angD \angLie_{\mathscr{Z}}^{N-1} \hat{\upchi}^{(Small)}
				+ \frac{1}{2} (\angdiff \mathscr{Z}^{N-1} \mytr \upchi^{(Small)}) \gsphere 	
			\right\rbrace
		\right|
		\\
		& \lesssim 
			\frac{1}{(1 + t)^2}
			\left|
				\fourmyarray[\rgeo \Lunit \mathscr{Z}^{\leq N} \Psi]
					{\Rad \mathscr{Z}^{\leq N} \Psi}
					{\rgeo \angdiff \mathscr{Z}^{\leq N} \Psi} 
					{\mathscr{Z}^{\leq N} \Psi}
			\right|
			+ 
			\frac{1}{(1 + t)^3}
	 		\left| 
				\myarray
					[\mathscr{Z}^{\leq N} (\upmu - 1)]
					{\rgeo \sum_{a=1}^3 |\mathscr{Z}^{\leq N} \Lunit_{(Small)}^a|}
			\right|,
			\notag
	\end{align}
	
	\begin{align}  \label{E:POINTWISECOMMUTEDANGDTRACEFREECHIJUNKINTERMSOFFREECOMMUTEDANGDTRACEFREECHI}
		& 
		\left|
			\angD \angLie_{\mathscr{Z}}^{N-1} \hat{\upchi}^{(Small)}
			- \angD \angfreeLietwoarg{\mathscr{Z}}{N-1} \hat{\upchi}^{(Small)}	
		\right|	
		\\
		& \lesssim 
			\frac{1}{(1 + t)^2}
			\left|
				\fourmyarray[\rgeo \Lunit \mathscr{Z}^{\leq N} \Psi]
					{\Rad \mathscr{Z}^{\leq N} \Psi}
					{\rgeo \angdiff \mathscr{Z}^{\leq N} \Psi} 
					{\mathscr{Z}^{\leq N} \Psi}
			\right|
		+ \frac{1}{(1 + t)^3}
	 		\left| 
				\myarray
					[\mathscr{Z}^{\leq N} (\upmu - 1)]
					{\rgeo \sum_{a=1}^3 |\mathscr{Z}^{\leq N} \Lunit_{(Small)}^a|}
			\right|,
			\notag
	\end{align}

	\begin{align} \label{E:POINTWISEESTIMATEROTTRCHIJUNKINTERMSOFANGDCHIJUNKPLUSJUNK}
		\left| \Rot \mathscr{Z}^{N-1} \mytr \upchi^{(Small)} \right|
		& \lesssim
			\rgeo \left| \angD \angLie_{\mathscr{Z}}^{N-1} \upchi^{(Small)} \right|
			\\
		& \ \ 
			+ 
			\frac{1}{1 + t}
				\left| 
					\fourmyarray
					 [\rgeo \Lunit \mathscr{Z}^{\leq N} \Psi]
					 {\Rad \mathscr{Z}^{\leq N} \Psi}
					 {\rgeo \angdiff \mathscr{Z}^{\leq N} \Psi}
					{\mathscr{Z}^{\leq N} \Psi}
				\right|
			+ 
			\frac{1}{(1 + t)^2}
	 		\left| 
				\myarray
					[\mathscr{Z}^{\leq N} (\upmu - 1)]
					{\rgeo \sum_{a=1}^3 |\mathscr{Z}^{\leq N} \Lunit_{(Small)}^a|}
				\right|,
			\notag
	\end{align}

	\begin{align} \label{E:POINTWISELOWERORDERANGREELIETRACEFREECHIJUNKISJUNK}
		\left| 
			\angfreeLietwoarg{\mathscr{Z}}{N-1} \hat{\upchi}^{(Small)} 
		\right|
		& \lesssim
			\frac{1}{1 + t}
				\left| 
					\fourmyarray
					 [\rgeo \Lunit \mathscr{Z}^{\leq N} \Psi]
					 {\Rad \mathscr{Z}^{\leq N} \Psi}
					 {\rgeo \angdiff \mathscr{Z}^{\leq N} \Psi}
					{\mathscr{Z}^{\leq N} \Psi}
				\right|
			+ 
			\frac{1}{(1 + t)^2}
	 		\left| 
				\myarray
					[\mathscr{Z}^{\leq N} (\upmu - 1)]
					{\rgeo \sum_{a=1}^3 |\mathscr{Z}^{\leq N} \Lunit_{(Small)}^a|}
			\right|.
		\end{align}
	
	\end{lemma}

\begin{proof}
	The estimate \eqref{E:POINTWISECOMMUTEDANGDCHIJUNKINTERMSOFTRACEFREEANDTRACEPART} 
	follows from
	\eqref{E:COMMUTINGANGDANDTRACETRACEFREESPLITTING}
	with $\upchi_{(Small)}$ in the role of $\xi$
	and $N-1$ in the role of $N,$
	\eqref{E:POINTWISEESTIMATESFORCHIJUNKINTERMSOFOTHERVARIABLES},
	and \eqref{E:C0BOUNDCRUCIALEIKONALFUNCTIONQUANTITIES}.	
	
	The estimate \eqref{E:POINTWISECOMMUTEDANGDTRACEFREECHIJUNKINTERMSOFFREECOMMUTEDANGDTRACEFREECHI}
	follows from
	\eqref{E:COMMUTINGTRACEFREELIEDERIVATIVESTHROUGHANGD}
	with $\upchi_{(Small)}$ in the role of $\xi$
	and $N-1$ in the role of $N,$
	\eqref{E:POINTWISEESTIMATESFORCHIJUNKINTERMSOFOTHERVARIABLES},
	and \eqref{E:C0BOUNDCRUCIALEIKONALFUNCTIONQUANTITIES}.	
	
	To prove \eqref{E:POINTWISEESTIMATEROTTRCHIJUNKINTERMSOFANGDCHIJUNKPLUSJUNK},
	we first use inequality \eqref{E:COMMUTINGTRACEWITHLIEDERIVATIVES}
	with $\upchi^{(Small)}$ in the role of $\xi$ and $\Rot \mathscr{Z}^{N-1}$
	in the role of $\mathscr{Z}^N$ together with
	the fact that the trace of a type $\binom{0}{2}$ $S_{t,u}$ tensor is
	in magnitude $\lesssim$ the magnitude of the tensor itself
	in order to deduce that
	$\left| \Rot \mathscr{Z}^{N-1} \mytr \upchi^{(Small)} \right|$
	is $\lesssim \left| \angLie_{\Rot} \angLie_{\mathscr{Z}}^{N-1} \mytr \upchi^{(Small)} \right|$
	plus the right-hand side of \eqref{E:COMMUTINGTRACEWITHLIEDERIVATIVES}.
	Using \eqref{E:C0BOUNDCRUCIALEIKONALFUNCTIONQUANTITIES}
	and the estimate
	\eqref{E:POINTWISEESTIMATESFORCHIJUNKINTERMSOFOTHERVARIABLES} with $N-1$
	in the role of $N,$ 
	we deduce that the right-hand side of \eqref{E:COMMUTINGTRACEWITHLIEDERIVATIVES} is
	$\lesssim$ the right-hand side of \eqref{E:POINTWISEESTIMATEROTTRCHIJUNKINTERMSOFANGDCHIJUNKPLUSJUNK}
	as desired. Furthermore, using \eqref{E:ANGLIEROTOFATENSORINTERMSOFANGD}, we deduce that
	$\left| \angLie_{\Rot} \angLie_{\mathscr{Z}}^{N-1} \mytr \upchi^{(Small)} \right|
	\lesssim \rgeo \left| \angD \angLie_{\mathscr{Z}}^{N-1} \upchi^{(Small)} \right|
	+ \left|\angLie_{\mathscr{Z}}^{N-1} \upchi^{(Small)} \right|.$
	The first term on the right-hand side of the previous inequality is 
	the first term on the right-hand side of \eqref{E:POINTWISEESTIMATEROTTRCHIJUNKINTERMSOFANGDCHIJUNKPLUSJUNK}, while
	the second is $\lesssim$ the remaining terms on the right-hand side of 	
	\eqref{E:POINTWISEESTIMATEROTTRCHIJUNKINTERMSOFANGDCHIJUNKPLUSJUNK}
	by virtue of the estimate
	\eqref{E:POINTWISEESTIMATESFORCHIJUNKINTERMSOFOTHERVARIABLES} with $N-1$
	in the role of $N.$

To prove \eqref{E:POINTWISELOWERORDERANGREELIETRACEFREECHIJUNKISJUNK},
we first use inequalities
\eqref{E:LIEDERIVATIVESOFTRACEFREEPARTINTERMSOFLIEDERIVATIVESOFFULLTENSOR}
and
\eqref{E:DIFFERENCEBETWEENTRACEFREELIEDERIVATIVESANDLIEDERIVATIVESACTINGONATRACEFREETENSOR}
with $\upchi^{(Small)}$ in the role of $\xi$ 
and $N-1$ in the role of $N$
to deduce that
$\left| \angfreeLietwoarg{\mathscr{Z}}{N-1} \hat{\upchi}^{(Small)} \right|$
is $\lesssim$ $\left| \angLie_{\mathscr{Z}}^{\leq N-1} \upchi^{(Small)} \right|$
plus the terms on the 
right-hand sides of
\eqref{E:LIEDERIVATIVESOFTRACEFREEPARTINTERMSOFLIEDERIVATIVESOFFULLTENSOR}
and
\eqref{E:DIFFERENCEBETWEENTRACEFREELIEDERIVATIVESANDLIEDERIVATIVESACTINGONATRACEFREETENSOR}.
Using inequalities \eqref{E:POINTWISEESTIMATESFORCHIJUNKINTERMSOFOTHERVARIABLES}
and \eqref{E:C0BOUNDCRUCIALEIKONALFUNCTIONQUANTITIES},
we deduce that all of these terms are
$\lesssim$ the right-hand side of \eqref{E:POINTWISELOWERORDERANGREELIETRACEFREECHIJUNKISJUNK} as desired.

\end{proof}


\chapter
[Pointwise Estimates for the Difficult Error Integrands]
{Pointwise Estimates for the Difficult Error Integrands Corresponding to the Commuted Wave Equation} 
\label{C:POINTWISEESTIMATESDIFFICULTERRORINTEGRANDS}
\thispagestyle{fancy}
Recall that for solutions to the commuted wave equation
$\upmu \square_{g(\Psi)} (\mathscr{Z}^N \Psi) = \inhomarg{\mathscr{Z}^N},$
the basic energy-flux identities verified by $\mathscr{Z}^N \Psi$
were established Prop.~\ref{P:DIVTHMWITHCANCELLATIONS},
with $\mathscr{Z}^N \Psi$ in the role of $\Psi$ and 
$\inhomarg{\mathscr{Z}^N}$ in the role of $\waveinhom.$
In Chapter~\ref{C:POINTWISEESTIMATESDIFFICULTERRORINTEGRANDS}, we complete the task,
initiated in Chapter~\ref{C:POINTWISEBOUNDSFOREASYERRORINTEGRANDS},
of deriving pointwise estimates
for the error integrands on the right-hand sides of the identities of the proposition;
these pointwise estimates are a preliminary ingredient 
in our derivation of a priori $L^2$ estimates, 
which we carry out in Chapter~\ref{C:ERRORTERMSOBOLEV}.
In particular, thanks to the decomposition \eqref{E:NTIMESCOMMUTEDWAVEINHOMOGENEOUSTERMFIRSTSPLITTING},
Prop.~\ref{P:IDOFKEYDIFFICULTENREGYERRORTERMS}, and
Cor.~\ref{C:POINTWISEESTIMATESFOREASYCOMMUTATORTERMS}, 
we have already derived pointwise
estimates for all terms in $\inhomarg{\mathscr{Z}^N}$ except for a handful,
and in fact, aside from the handful,
the terms are $Harmless^{\leq N}$ in the sense of Def.~\ref{D:HARMLESSTERMS}.
Our goal in this section 
is to derive pointwise estimates for the remaining terms 
in $\inhomarg{\mathscr{Z}^N}$
that are not $Harmless^{\leq N}.$ These difficult terms would cause 
derivative loss and other problems
if they were not properly handled, and to derive suitable pointwise estimates for them, 
we need to use the modified quantities of Ch.~\ref{C:RENORMALIZEDEIKONALFUNCTIONQUANTITIES}.

The pointwise estimates that we derive in Chapter~\ref{C:POINTWISEESTIMATESDIFFICULTERRORINTEGRANDS}
can be split into two classes.
In Chapter~\ref{C:ERRORTERMSOBOLEV}, we use the first class of estimates to bound the difficult top-order
error integrals corresponding to the timelike multiplier $\Mult.$
We derive the first class estimates with the help of the fully modified quantities of Ch.~\ref{C:RENORMALIZEDEIKONALFUNCTIONQUANTITIES}.
We provide the main estimates in Prop.~\ref{P:MAINCOMMUTEDWAVEEQNINHOMOGENEOUSTERMPOINTWISEESTIMATES}
of Sect.~\ref{S:FULLYRENORMALIZEDPOINTWISEESTIMATESFORDIFFICULTENERGYERRORINTEGRAND},
but in order to prove the proposition, we first prove a series of auxiliary lemmas
in Sects.~\ref{S:PRELIMINARYPOINTWISEFULLYMOD}
and \ref{S:INVERTINGFULLYRENORMALIZEDTRANSPORT}.
Similarly, in Chapter~\ref{C:ERRORTERMSOBOLEV}, we use the second class of estimates to bound the difficult top-order
error integrals corresponding to the Morawetz multiplier $\Mor.$
We derive the second class estimates with the help of the partially modified quantities of Ch.~\ref{C:RENORMALIZEDEIKONALFUNCTIONQUANTITIES}.
These estimates are easier to derive than those of the first class. We prove
some auxiliary lemmas in Sect.~\ref{S:PRELIMINARYPOINTWISEPARTIALLYMODIFIED},
and we then provide the main estimates in
Sect.~\ref{S:POINTWISEPARTIALLYMODIFIED}.

\section{Preliminary pointwise estimates for 
\texorpdfstring{$\mathscr{Z}^N \mathfrak{X} \ \mbox{and} \ \mathscr{Z}^{N-1} \chifullmodsourcearg{S}$}
{the derivatives of the inhomogeneous terms in the transport equations for the fully modified quantities}}
\label{S:PRELIMINARYPOINTWISEFULLYMOD}
Recall that Prop.~\ref{P:TOPORDERTRCHIJUNKRENORMALIZEDTRANSPORT} 
provides the transport equation verified by the fully modified quantity
$\chifullmodarg{\mathscr{S}^N}$ and that
$\mathscr{S}^N \mathfrak{X}$ 
and $\mathscr{S}^{N-1} \chifullmodsourcearg{S}$
are two source terms in the equation
(see the right-hand side of \eqref{E:NTIMESCOMMUTATIONRENORMALIZEDTRCHIJUNKINHOMOGENEOUSTERM}).
In the next lemma, we derive pointwise estimates 
for these source terms.

\begin{lemma}[\textbf{Pointwise estimates for} $\mathscr{Z}^N \mathfrak{X}$ \textbf{and} $\mathscr{Z}^{N-1} \chifullmodsourcearg{S}$]
	\label{L:POINTWISEESTIMATEZNAPPLIEDTOONECOMMUTATIONRENORMALIZEDTRCHIJUNKINHOMOGENEOUSTERM}
	Let $0 \leq N \leq 24$ be an integer. Let $\mathfrak{X}$ be the quantity 
	defined in \eqref{E:LOWESTORDERTRANSPORTRENORMALIZEDTRCHIJUNKDISCREPANCY}.
	Under the small-data and bootstrap assumptions 
	of Sects.~\ref{S:PSISOLVES}-\ref{S:C0BOUNDBOOTSTRAP},
	if $\varepsilon$ is sufficiently small, then
	the following estimates hold on $\mathcal{M}_{\Tboot,U_0}:$
	\begin{subequations}
	\begin{align} \label{E:POINTWISEESTIMATEZNAPPLIEDTOLOWESTTRCHIJUNKDISCREPANCYTERM}
		\left|
			\mathscr{Z}^N \mathfrak{X}
		\right|
		& \lesssim
		\frac{\ln(\myexp + t)}{1 + t}
			\left|
				\threemyarray[\rgeo \Lunit \mathscr{Z}^{\leq N} \Psi]
					{\rgeo \angdiff \mathscr{Z}^{\leq N} \Psi} 
					{\mathscr{Z}^{\leq N} \Psi}
			\right|
			+ \left|
					\Rad \mathscr{Z}^{\leq N} \Psi
				\right|
			+
			\varepsilon
			\frac{1}{(1 + t)^2}
			\left| 
				\myarray
					[\mathscr{Z}^{\leq N} (\upmu - 1)]
						{\rgeo \sum_{a=1}^3 |\mathscr{Z}^{\leq N} \Lunit_{(Small)}^a|}
			\right|,
			\\
		\left\|
			\mathscr{Z}^{\leq 12} \mathfrak{X}
		\right\|_{C^0(\Sigma_t^u)}
		& \lesssim 
			\varepsilon
			\frac{1}{1 + t}.
			\label{E:LOWERORDERC0BOUNDZNAPPLIEDTOLOWESTTRCHIJUNKDISCREPANCYTERM}
	\end{align}
	\end{subequations}
	
	Furthermore, the following estimates hold on $\mathcal{M}_{\Tboot,U_0}:$
	\begin{align} \label{E:POINTWISEESTIMATECOMMUTATORINVOLVINGZNAPPLIEDTOLOWESTTRCHIJUNKDISCREPANCYTERM}
		&
		\left|
			\mathscr{Z}^N \mathfrak{X}
			- 
			\left\lbrace
				- G_{\Lunit \Lunit} \Rad \mathscr{Z}^N \Psi
				- \frac{1}{2} \upmu \angGmixedarg{A}{A} \Lunit \mathscr{Z}^N  \Psi
				- \frac{1}{2} \upmu G_{\Lunit \Lunit} \Lunit \mathscr{Z}^N \Psi
				+ \upmu \angGmixedarg{\Lunit}{A} \angdiffarg{A} \mathscr{Z}^N \Psi	
			\right\rbrace
		\right|	
			\\
		& \lesssim
			\frac{\ln(\myexp + t)}{1 + t}
			\left|
				\threemyarray[\rgeo \Lunit \mathscr{Z}^{\leq N-1} \Psi]
					{\rgeo \angdiff \mathscr{Z}^{\leq N-1} \Psi} 
					{\mathscr{Z}^{\leq N-1} \Psi}
			\right|
			+ \left|
					\Rad \mathscr{Z}^{\leq N-1} \Psi
				\right|	
			+
			\varepsilon
			\frac{1}{(1 + t)^2}
			\left| 
				\myarray
					[\mathscr{Z}^{\leq N} (\upmu - 1)]
						{\rgeo \sum_{a=1}^3 |\mathscr{Z}^{\leq N} \Lunit_{(Small)}^a|}
			\right|.
			\notag 
	\end{align}
	
	Furthermore, let $S \in \mathscr{S}$ be a spatial commutation vectorfield, and
	let $\chifullmodsourcearg{S}$ be the inhomogeneous term \eqref{E:ONECOMMUTATIONRENORMALIZEDTRCHIJUNKINHOMOGENEOUSTERM}.
	Then the following estimates hold on $\mathcal{M}_{\Tboot,U_0}:$
	\begin{subequations}
	\begin{align} \label{E:POINTWISEESTIMATEZNAPPLIEDTOONECOMMUTATIONRENORMALIZEDTRCHIJUNKINHOMOGENEOUSTERM}
		\left|
			\mathscr{Z}^{N-1} \chifullmodsourcearg{S}
		\right|
		& \lesssim
			\frac{\ln(\myexp + t)}{(1 + t)^2}
			\left|
				\fourmyarray[\rgeo \Lunit \mathscr{Z}^{\leq N} \Psi]
					{\Rad \mathscr{Z}^{\leq N} \Psi}
					{\rgeo \angdiff \mathscr{Z}^{\leq N} \Psi} 
					{\mathscr{Z}^{\leq N} \Psi}
			\right|
			+
			\varepsilon
			\frac{\ln(\myexp + t)}{(1 + t)^3}
			\left| 
				\myarray
					[\mathscr{Z}^{\leq N} (\upmu - 1)]
						{\rgeo \sum_{a=1}^3 |\mathscr{Z}^{\leq N} \Lunit_{(Small)}^a|}
			\right|,
			\\
		\left\|
			\mathscr{Z}^{\leq 11} \chifullmodsourcearg{S}
		\right\|_{C^0(\Sigma_t^u)}
		& \lesssim 
			\varepsilon
			\frac{\ln^2(\myexp + t)}{(1 + t)^3}.
			\label{E:C0BOUNDESTIMTELOWERORDERZNAPPLIEDTOONECOMMUTATIONRENORMALIZEDTRCHIJUNKINHOMOGENEOUSTERM}
	\end{align}
	\end{subequations}
\end{lemma}

\begin{proof}
\noindent{\textbf{Proof of} \eqref{E:POINTWISEESTIMATEZNAPPLIEDTOLOWESTTRCHIJUNKDISCREPANCYTERM}
	\textbf{and} \eqref{E:LOWERORDERC0BOUNDZNAPPLIEDTOLOWESTTRCHIJUNKDISCREPANCYTERM}:}	
  We note the schematic identity 
	$\mathfrak{X}
	= G_{(Frame)}^{\#} 
				\threemyarray
					[\upmu \Lunit \Psi]
					{\Rad \Psi}
					{\upmu \angdiff \Psi}.
	$					
	Hence, our proofs of \eqref{E:GFRAMEGINVERSESPHEREUPMUPSIDERIVATIVESLEIBNIZEXPANDED}
	and \eqref{E:C0BOUNDGFRAMEGINVERSESPHEREUPMUPSIDERIVATIVESLEIBNIZEXPANDED}
	also yield the desired bounds
	\eqref{E:POINTWISEESTIMATEZNAPPLIEDTOLOWESTTRCHIJUNKDISCREPANCYTERM}
	and \eqref{E:LOWERORDERC0BOUNDZNAPPLIEDTOLOWESTTRCHIJUNKDISCREPANCYTERM}.

\ \\
\noindent{\textbf{Proof of} \eqref{E:POINTWISEESTIMATECOMMUTATORINVOLVINGZNAPPLIEDTOLOWESTTRCHIJUNKDISCREPANCYTERM}}:
We first use definition \eqref{E:LOWESTORDERTRANSPORTRENORMALIZEDTRCHIJUNKDISCREPANCY}
and Lemma~\ref{L:LANDRADCOMMUTEWITHANGDIFF}
to deduce that the term in the absolute value on the left-hand side of 
\eqref{E:POINTWISEESTIMATECOMMUTATORINVOLVINGZNAPPLIEDTOLOWESTTRCHIJUNKDISCREPANCYTERM}
can be written as
\begin{align} \label{E:COMMUTATORDECOMPOSED}
& - [\mathscr{Z}^N, G_{\Lunit \Lunit}] \Rad \Psi
	- \frac{1}{2} [\angLie_{\mathscr{Z}}^N, \upmu \angGmixedarg{A}{A}] \Lunit \Psi
	- \frac{1}{2} [\mathscr{Z}^N, \upmu G_{\Lunit \Lunit}] \Lunit \Psi
	+ [\angLie_{\mathscr{Z}}^N, \upmu \angGmixedarg{\Lunit}{A}] \angdiffarg{A} \Psi	
		\\
	& \ \ 
		- G_{\Lunit \Lunit} [\mathscr{Z}^N, \Rad] \Psi
		- \frac{1}{2} \upmu \angGmixedarg{A}{A} [\angLie_{\mathscr{Z}}^N, \Lunit] \Psi
		- \frac{1}{2} \upmu G_{\Lunit \Lunit} [\mathscr{Z}^N, \Lunit] \Psi.
			\notag
\end{align}
Hence, using the Leibniz rule and the reasoning used to prove 
\eqref{E:GFRAMEGINVERSESPHEREUPMUPSIDERIVATIVESLEIBNIZEXPANDED}
and
\eqref{E:POINTWISEESTIMATEZNAPPLIEDTOLOWESTTRCHIJUNKDISCREPANCYTERM},
we deduce that the terms on the first line of \eqref{E:COMMUTATORDECOMPOSED} are in magnitude 
$\lesssim$ the right-hand side of \eqref{E:GFRAMEGINVERSESPHEREUPMUPSIDERIVATIVESLEIBNIZEXPANDED},
\emph{with $N-1$ in place of $N$ in the $\Psi-$containing terms on the right-hand side}.
It follows that these terms are in magnitude $\lesssim$ the right-hand side of
\eqref{E:POINTWISEESTIMATECOMMUTATORINVOLVINGZNAPPLIEDTOLOWESTTRCHIJUNKDISCREPANCYTERM} as desired.
To deduce that the terms on the second line of \eqref{E:COMMUTATORDECOMPOSED} are in magnitude 
$\lesssim$ the right-hand side of 
\eqref{E:POINTWISEESTIMATECOMMUTATORINVOLVINGZNAPPLIEDTOLOWESTTRCHIJUNKDISCREPANCYTERM},
we use the estimates 
\eqref{E:LOWERORDERC0BOUNDLIEDERIVATIVESOFGRAME}
and \eqref{E:C0BOUNDCRUCIALEIKONALFUNCTIONQUANTITIES},
inequalities \eqref{E:RGEOLORRADZNCOMMUTATORACTINGONFUNCTIONSPOINTWISE} 
and \eqref{E:LZNCOMMUTATORACTINGONFUNCTIONSPOINTWISE}
with $f = \Psi,$
and the bootstrap assumptions \eqref{E:PSIFUNDAMENTALC0BOUNDBOOTSTRAP}.
We have thus proved the desired bound \eqref{E:POINTWISEESTIMATECOMMUTATORINVOLVINGZNAPPLIEDTOLOWESTTRCHIJUNKDISCREPANCYTERM}.

\ \\
\noindent{\textbf{Proof of} \eqref{E:POINTWISEESTIMATEZNAPPLIEDTOONECOMMUTATIONRENORMALIZEDTRCHIJUNKINHOMOGENEOUSTERM}
	\textbf{and} \eqref{E:C0BOUNDESTIMTELOWERORDERZNAPPLIEDTOONECOMMUTATIONRENORMALIZEDTRCHIJUNKINHOMOGENEOUSTERM}:}	
	It suffices to prove \eqref{E:POINTWISEESTIMATEZNAPPLIEDTOONECOMMUTATIONRENORMALIZEDTRCHIJUNKINHOMOGENEOUSTERM}
	since \eqref{E:C0BOUNDESTIMTELOWERORDERZNAPPLIEDTOONECOMMUTATIONRENORMALIZEDTRCHIJUNKINHOMOGENEOUSTERM}
	follows from 
	\eqref{E:POINTWISEESTIMATEZNAPPLIEDTOONECOMMUTATIONRENORMALIZEDTRCHIJUNKINHOMOGENEOUSTERM},
	\eqref{E:C0BOUNDCRUCIALEIKONALFUNCTIONQUANTITIES},
	and the bootstrap assumptions \eqref{E:PSIFUNDAMENTALC0BOUNDBOOTSTRAP}.
	To proceed, we apply $\angLie_{\mathscr{Z}}^{N-1}$ to the right-hand side of 	
	\eqref{E:ONECOMMUTATIONRENORMALIZEDTRCHIJUNKINHOMOGENEOUSTERM}.
	We begin by bounding the term $\angLie_{\mathscr{Z}}^{N-1} \angLie_S \mathfrak{A}$ arising from the first term on the right-hand side of
	\eqref{E:ONECOMMUTATIONRENORMALIZEDTRCHIJUNKINHOMOGENEOUSTERM}, where $\mathfrak{A}$ is given by
	\eqref{E:RENORMALIZEDALPHAINHOMOGENEOUSTERM}. To this end, we apply the Leibniz rule to the terms on the right-hand side of
	\eqref{E:RENORMALIZEDALPHAINHOMOGENEOUSTERM}. The terms of interest can be bounded by using
	essentially the same argument used to deduce 
	\eqref{E:POINTWISEESTIMATEZNAPPLIEDTOLOWESTTRCHIJUNKDISCREPANCYTERM}
	and \eqref{E:LOWERORDERC0BOUNDZNAPPLIEDTOLOWESTTRCHIJUNKDISCREPANCYTERM},
	except that we use Lemma~\ref{L:POINTWISEESTIMATESGFRAMEINTERMSOFOTHERQUANTITIES}
	in place of Cor.~\ref{C:POINTWISEESTIMATESGFRAMESHARPINTERMSOFOTHERQUANTITIES}
	and we also use the estimates 
	\eqref{E:CRUDEPOINTWISEBOUNDSDERIVATIVESOFANGULARDEFORMATIONTENSORS} and
	\eqref{E:CRUDELOWERORDERC0BOUNDDERIVATIVESOFANGULARDEFORMATIONTENSORS}
	to bound the Lie derivatives of $\ginversesphere.$
	We note that the resulting bounds are better by a factor of $(1 + t)^{-2}$ 
	compared to 
	\eqref{E:POINTWISEESTIMATEZNAPPLIEDTOLOWESTTRCHIJUNKDISCREPANCYTERM}
	and \eqref{E:LOWERORDERC0BOUNDZNAPPLIEDTOLOWESTTRCHIJUNKDISCREPANCYTERM}
	because the right-hand side
	of \eqref{E:RENORMALIZEDALPHAINHOMOGENEOUSTERM} features an additional factor of
	$\myarray
		[\Lunit \Psi]
		{\angdiff \Psi}$
	compared to the 
	term $\mathfrak{X}
	= G_{(Frame)}^{\#} 
				\threemyarray
					[\upmu \Lunit \Psi]
					{\Rad \Psi}
					{\upmu \angdiff \Psi}
	$
	from \eqref{E:POINTWISEESTIMATEZNAPPLIEDTOLOWESTTRCHIJUNKDISCREPANCYTERM}
	and \eqref{E:LOWERORDERC0BOUNDZNAPPLIEDTOLOWESTTRCHIJUNKDISCREPANCYTERM}.
	
	We now bound the terms arising from the
	last products on the first line of the right-hand side of \eqref{E:ONECOMMUTATIONRENORMALIZEDTRCHIJUNKINHOMOGENEOUSTERM}
	(that is, the product involving the terms in braces).
	This term can be schematically written as
	\begin{align}  \label{E:SCHEMATICONECOMMUTATIONRENORMALIZEDTRCHIJUNKINHOMOGENEOUSSECONDTERM}
				\mytr \upchi
				\angLie_S
				\left\lbrace
					G_{(Frame)}^{\#} 
					\myarray
						[\upmu \Lunit \Psi]
						{\upmu \angdiff \Psi}
				\right\rbrace.
		\end{align}
	Applying $\angLie_{\mathscr{Z}}^{N-1}$ 
	to the right-hand side of \eqref{E:SCHEMATICONECOMMUTATIONRENORMALIZEDTRCHIJUNKINHOMOGENEOUSSECONDTERM}
	and using the Leibniz rule, we bound the resulting terms by
	\begin{align}	 \label{E:FIRSTPOINTWISEESTIMATEZNAPPLIEDTOSECONDTERMLOWESTTRCHIJUNKINHOMOGENEOUS}
		& \lesssim  
			\mathop{\sum_{N_1 + N_2 + N_3 + N_4 \leq N}}_{N_1 \leq N - 1}
			\left|
				\mathscr{Z}^{\leq N_1} \mytr \upchi
			\right|
			\left|
				\angLie_{\mathscr{Z}}^{\leq N_2} G_{(Frame)}^{\#} 
			\right|
			\left|
				\mathscr{Z}^{\leq N_3} \upmu
			\right|
			\left|	
				\angLie_{\mathscr{Z}}^{\leq N_4}
				\myarray
					[\Lunit \Psi]
					{\angdiff \Psi}
			\right|.
	\end{align}
	Using \eqref{E:FIRSTPOINTWISEESTIMATEZNAPPLIEDTOSECONDTERMLOWESTTRCHIJUNKINHOMOGENEOUS},	
	Lemma~\ref{L:AVOIDINGCOMMUTING} with $f=\Psi$,
	Lemma~\ref{L:POINTWISEBOUNDSDERIVATIVESOFANGULARDEFORMATIONTENSORS} with $N_1 \leq N-1$ in the role of $N,$
	Cor.~\ref{C:POINTWISEESTIMATESGFRAMESHARPINTERMSOFOTHERQUANTITIES},
	\eqref{E:C0BOUNDCRUCIALEIKONALFUNCTIONQUANTITIES},
	and the bootstrap assumptions \eqref{E:PSIFUNDAMENTALC0BOUNDBOOTSTRAP},
	we deduce that the right-hand side of 
	\eqref{E:FIRSTPOINTWISEESTIMATEZNAPPLIEDTOSECONDTERMLOWESTTRCHIJUNKINHOMOGENEOUS} 
	is $\lesssim$ the right-hand side of \eqref{E:POINTWISEESTIMATEZNAPPLIEDTOONECOMMUTATIONRENORMALIZEDTRCHIJUNKINHOMOGENEOUSTERM} 
	as desired.
	
	We now bound the first product $(S \upmu) \Lunit \mytr \upchi^{(Small)}$
	on the second line of
	the right-hand side of 	
	\eqref{E:ONECOMMUTATIONRENORMALIZEDTRCHIJUNKINHOMOGENEOUSTERM}.
	We apply $\mathscr{Z}^{N-1}$ to this product
	and use the Leibniz rule, which leads to the following bound:
	\begin{align} \label{E:FIRSTPOINTWISEBOUNDLIEZNSPATIALUPMUTIMESLTRCHIJUNK}
		\left|
			\mathscr{Z}^{N-1}
			\left\lbrace
				(S \upmu)
				\Lunit \mytr \upchi^{(Small)} 
			\right\rbrace
		\right|
			& \lesssim 
				\left|
					\mathscr{Z}^{\leq 12} (\upmu - 1)
				\right|
				\left|
					\mathscr{Z}^{\leq N-1} \Lunit \mytr \upchi^{(Small)}
				\right|
				\\
		& \ \
				+
				\left|
					\mathscr{Z}^{\leq 11} \Lunit \mytr \upchi^{(Small)}
				\right|
				\left|
					\mathscr{Z}^{\leq N} (\upmu - 1)
				\right|.
				\notag
	\end{align}	
	We then use the estimates 
	\eqref{E:ZNAPPLIEDTORGEOISNOTTOOLARGE},
	\eqref{E:POINTWISEESTIMATESFORCHIJUNKINTERMSOFOTHERVARIABLES},
	\eqref{E:LDERIVATIVECRUCICALTRANSPORTINTEQUALITIES},
	\eqref{E:C0BOUNDCRUCIALEIKONALFUNCTIONQUANTITIES},
	and
	\eqref{E:C0BOUNDLDERIVATIVECRUCICALEIKONALFUNCTIONQUANTITIES}
	to deduce that the right-hand side of 
	\eqref{E:FIRSTPOINTWISEBOUNDLIEZNSPATIALUPMUTIMESLTRCHIJUNK} is
	$\lesssim$ the right-hand side of \eqref{E:POINTWISEESTIMATEZNAPPLIEDTOONECOMMUTATIONRENORMALIZEDTRCHIJUNKINHOMOGENEOUSTERM}
	as desired. The remaining products 
	on the second line of the right-hand side of \eqref{E:ONECOMMUTATIONRENORMALIZEDTRCHIJUNKINHOMOGENEOUSTERM}
	can be bounded in a similar fashion.
	The last two terms on the last line of the
	right-hand side of \eqref{E:ONECOMMUTATIONRENORMALIZEDTRCHIJUNKINHOMOGENEOUSTERM}
	can also be bounded in a similar fashion. For these estimates, it is important
	that $S \rgeo \in \lbrace 0,-1 \rbrace$ for spatial commutation vectorfields $S \in \mathscr{S}.$
	
	It remains for us to bound the first term on the last line of the right-hand side of 	
	\eqref{E:ONECOMMUTATIONRENORMALIZEDTRCHIJUNKINHOMOGENEOUSTERM}.
	To this end, we note that the identity $[\Lunit,S] = \angdeformoneformupsharparg{S}{\Lunit}$ 
	(see \eqref{E:LCOMMUTETANGENTISTANGENT})
	implies that the 
	term can be written as $\angdeformoneformupsharparg{S}{\Lunit} \cdot \angdiff \mathfrak{X}.$
	Applying $\angLie_{\mathscr{Z}}^{N-1}$ to this product and using 
	the Leibniz rule, Lemma~\ref{L:LANDRADCOMMUTEWITHANGDIFF}, and inequality \eqref{E:FUNCTIONPOINTWISEANGDINTERMSOFANGLIEO},
	we see that it suffices to bound
	\begin{align} \label{E:FIRSTPOINTWISEESTIMATEZNAPPLIEDTOLOWESTTRCHIJUNKINHOMOGENEOUSCOMMUTATORTERM}
		\frac{1}{1 + t}
		\left|
			\angLie_{\mathscr{Z}}^{\leq 11} \angdeformoneformupsharparg{S}{\Lunit}
		\right|
		\left|
			\mathscr{Z}^{\leq N} \mathfrak{X}
		\right|
		+ 
		\frac{1}{1 + t}
		\left|
			\mathscr{Z}^{\leq 12} \mathfrak{X}
		\right|
		\left|
			\angLie_{\mathscr{Z}}^{\leq N-1} \angdeformoneformupsharparg{S}{\Lunit}
		\right|.
	\end{align}
	Using \eqref{E:POINTWISEESTIMATEZNAPPLIEDTOLOWESTTRCHIJUNKDISCREPANCYTERM},
	\eqref{E:LOWERORDERC0BOUNDZNAPPLIEDTOLOWESTTRCHIJUNKDISCREPANCYTERM},
	\eqref{E:RADDEFORMSPHERELPOINTWISE},
	\eqref{E:LOWERORDERC0BOUNDRADDEFORMSPHEREL},
	\eqref{E:ROTDEFORMSPHERELPOINTWISE},
	and \eqref{E:LOWERORDERC0BOUNDROTDEFORMSPHEREL},
	we see that 
	\eqref{E:FIRSTPOINTWISEESTIMATEZNAPPLIEDTOLOWESTTRCHIJUNKINHOMOGENEOUSCOMMUTATORTERM}
	is $\lesssim$ the right-hand side of \eqref{E:POINTWISEESTIMATEZNAPPLIEDTOONECOMMUTATIONRENORMALIZEDTRCHIJUNKINHOMOGENEOUSTERM}
	as desired.
	
\end{proof}

In the next lemma, we derive pointwise estimates for the inhomogeneous term
$\chifullmodsourcearg{\mathscr{S}^N}$ appearing on the right-hand side of
\eqref{E:TOPORDERTRCHIJUNKRENORMALIZEDTRANSPORT}.

\begin{lemma}[\textbf{Pointwise estimate for the inhomogeneous term} $\chifullmodsourcearg{\mathscr{S}^N}$]
\label{L:TOPORDERRENORMALIZEDTRCHIJUNKTRANSPORTINHOMOGENEOUSTERMPOINTWISEESTIMATE}
Let $1 \leq N \leq 24$ be an integer, and let
$\chifullmodsourcearg{\mathscr{S}^N}$ be the inhomogeneous term 
defined in \eqref{E:NTIMESCOMMUTATIONRENORMALIZEDTRCHIJUNKINHOMOGENEOUSTERM}.
Under the small-data and bootstrap assumptions 
of Sects.~\ref{S:PSISOLVES}-\ref{S:C0BOUNDBOOTSTRAP},
if $\varepsilon$ is sufficiently small, then
the following estimates hold on $\mathcal{M}_{\Tboot,U_0}:$
\begin{align} \label{E:TOPORDERRENORMALIZEDTRCHIJUNKTRANSPORTINHOMOGENEOUSTERMPOINTWISEESTIMATE}
	\left|
		\chifullmodsourcearg{\mathscr{S}^N}
	\right|
	& \lesssim 
		\frac{\ln(\myexp + t)}{(1 + t)^2}
			\left|
				\fourmyarray[\rgeo \Lunit \mathscr{Z}^{\leq N} \Psi]
					{\Rad \mathscr{Z}^{\leq N} \Psi}
					{\rgeo \angdiff \mathscr{Z}^{\leq N} \Psi} 
					{\mathscr{Z}^{\leq N} \Psi}
			\right|
		+ \frac{\ln(\myexp + t)}{(1 + t)^3}
			\left| 
				\myarray
					[\mathscr{Z}^{\leq N} (\upmu - 1)]
						{\rgeo \sum_{a=1}^3 |\mathscr{Z}^{\leq N} \Lunit_{(Small)}^a|}
				\right|.
\end{align}

\end{lemma}

\begin{proof}
	We begin by noting that many of the commutator terms that we estimate in this proof are absent when $N=1.$
	To bound the first term $\mathscr{S}^{N-1} \chifullmodsourcearg{S}$
	on the first line of \eqref{E:NTIMESCOMMUTATIONRENORMALIZEDTRCHIJUNKINHOMOGENEOUSTERM}
	by the right-hand side of \eqref{E:TOPORDERRENORMALIZEDTRCHIJUNKTRANSPORTINHOMOGENEOUSTERMPOINTWISEESTIMATE},
	we simply quote the estimate \eqref{E:POINTWISEESTIMATEZNAPPLIEDTOONECOMMUTATIONRENORMALIZEDTRCHIJUNKINHOMOGENEOUSTERM}.
	
	To bound the second term $[\Lunit, \mathscr{S}^{N-1}] S \mathfrak{X}$
	on the first line of \eqref{E:NTIMESCOMMUTATIONRENORMALIZEDTRCHIJUNKINHOMOGENEOUSTERM}
	by the right-hand side of \eqref{E:TOPORDERRENORMALIZEDTRCHIJUNKTRANSPORTINHOMOGENEOUSTERMPOINTWISEESTIMATE},
	we first use the commutator estimate \eqref{E:LSNCOMMUTATORACTINGONFUNCTIONSPOINTWISE} 
	with $N-1$ in the role of $N$ and $S \mathfrak{X}$ in the role of $f$ to deduce that
	\begin{align} \label{E:FIRSTPOINTWISEESTIMATEDERIVATIVESOFCOMMUTATORAPPLIEDTOLOWERSTORDERDISCREPANCYTERM}
		\left|
			[\Lunit, \mathscr{S}^{N-1}] S \mathfrak{X}
		\right|
		& \lesssim
		\varepsilon^{1/2}
			\frac{\ln(\myexp + t)}{(1 + t)^2}
			\left| 
				\fourmyarray[ \rgeo \Lunit \mathscr{Z}^{\leq N} \mathfrak{X}]
					{\Rad \mathscr{Z}^{\leq N} \mathfrak{X}}
					{\rgeo \angdiff \mathscr{Z}^{\leq N} \mathfrak{X}}
					{\mathscr{Z}^{\leq N} \mathfrak{X}}
			\right|
			+ \frac{1}{1 + t}
			\left\|
				\mathscr{Z}^{\leq 12} \mathfrak{X}
			\right\|_{C^0(\Sigma_t^u)}
			\left| 
				\fourmyarray[ \rgeo \Lunit \mathscr{Z}^{\leq N-1} \Psi]
					{\Rad \mathscr{Z}^{\leq N-1} \Psi}
					{\rgeo \angdiff \mathscr{Z}^{\leq N-1} \Psi}
					{\mathscr{Z}^{\leq N-1} \Psi}
			\right|
			\\
	&  \ \ 
			+ 
			\frac{1}{(1 + t)^2}
			\left\|
				\mathscr{Z}^{\leq 12} \mathfrak{X}
			\right\|_{C^0(\Sigma_t^u)}
			\left|
				\myarray[\mathscr{Z}^{\leq N} (\upmu - 1)]
					{\sum_{a=1}^3 \rgeo |\mathscr{Z}^{\leq N} \Lunit_{(Small)}^a|} 
			\right|.
			\notag
	\end{align}
	The desired bound now follows from
	\eqref{E:FIRSTPOINTWISEESTIMATEDERIVATIVESOFCOMMUTATORAPPLIEDTOLOWERSTORDERDISCREPANCYTERM},
	\eqref{E:POINTWISEESTIMATEZNAPPLIEDTOLOWESTTRCHIJUNKDISCREPANCYTERM}, 
	and \eqref{E:LOWERORDERC0BOUNDZNAPPLIEDTOLOWESTTRCHIJUNKDISCREPANCYTERM}.
	
	To bound the first term $- [\mathscr{S}^{N-1}, \upmu \mytr \upchi] S \mytr \upchi^{(Small)}$ on the second line
	of \eqref{E:NTIMESCOMMUTATIONRENORMALIZEDTRCHIJUNKINHOMOGENEOUSTERM} 
	by the right-hand side of \eqref{E:TOPORDERRENORMALIZEDTRCHIJUNKTRANSPORTINHOMOGENEOUSTERMPOINTWISEESTIMATE},
	we first use the Leibniz rule, 
	the decomposition $\mytr \upchi = 2 \rgeo^{-1} + \mytr \upchi^{(Small)},$
	and the bound $\left|\mathscr{S}^{\leq N-1} \rgeo^{-1} \right| \lesssim \frac{1}{1 + t}$
	(which follows from \eqref{E:EXACTRELATIONSZAPPLIEDTORGEO})
	to deduce that
	\begin{align} \label{E:FIRSTPOINTWISEESTIMATETRCHITIMESCOMMUTINGVECTORFIELDSWITHPUMPUTIMESSTRICHIJUNK}
		\left|
				[\mathscr{S}^{N-1}, \upmu \mytr \upchi] S \mytr \upchi^{(Small)}
		\right|	
		& \lesssim 
			\frac{1}{1 + t}
			\mathop{\sum_{N_1 + N_2 \leq N}}_{N_1, N_2 \leq N-1}
			\left|
				\mathscr{Z}^{N_1} \upmu
			\right|
			\left|
				\mathscr{Z}^{N_2} \mytr \upchi^{(Small)}
			\right|
				\\
	& \ \ 
			+
			\mathop{\sum_{N_1 + N_2 + N_3 \leq N}}_{N_1, N_2, N_3 \leq N-1}
			\left|
				\mathscr{Z}^{N_1} \upmu
			\right|
			\left|
				\mathscr{Z}^{N_2} \mytr \upchi^{(Small)}
			\right|
			\left|
				\mathscr{Z}^{N_3} \mytr \upchi^{(Small)}
			\right|.
			\notag
	\end{align}
	The desired bound now follows from \eqref{E:FIRSTPOINTWISEESTIMATETRCHITIMESCOMMUTINGVECTORFIELDSWITHPUMPUTIMESSTRICHIJUNK},
	\eqref{E:POINTWISEESTIMATESFORCHIJUNKINTERMSOFOTHERVARIABLES}
	and \eqref{E:C0BOUNDCRUCIALEIKONALFUNCTIONQUANTITIES}.
	
	To bound the second term 
	$\frac{1}{2} [\mathscr{S}^{N-1}, \mytr \upchi] S \mathfrak{X} 
	= \frac{1}{2} [\mathscr{S}^{N-1}, 2 \rgeo^{-1} + \mytr \upchi^{(Small)}] S \mathfrak{X}$ 
	on the second line
	of \eqref{E:NTIMESCOMMUTATIONRENORMALIZEDTRCHIJUNKINHOMOGENEOUSTERM} 
	by the right-hand side of \eqref{E:TOPORDERRENORMALIZEDTRCHIJUNKTRANSPORTINHOMOGENEOUSTERMPOINTWISEESTIMATE},
	we first use the Leibniz rule 
	and the bound $\left|\mathscr{S}^M \rgeo^{-1} \right| \lesssim \frac{1}{(1 + t)^2}$
	(which holds for $M \geq 1$ and follows from \eqref{E:EXACTRELATIONSZAPPLIEDTORGEO})
	to bound its magnitude by
	\begin{align} \label{E:LOWERORDERTERMSINVOLVINGSPATIALDERIVATIVESOFRGEO}
		& \lesssim
			\left|
				\mathscr{Z}^{\leq 11} \mytr \upchi^{(Small)}
			\right|
			\left|
				\mathscr{Z}^{N-1} \mathfrak{X}
			\right|
			+
			\frac{1}{(1 + t)^2}
			\left|
				\mathscr{Z}^{N-1} \mathfrak{X}
			\right|
			+ 
			\left|
				\mathscr{Z}^{11} \mathfrak{X}
			\right|
			\left|
				\mathscr{Z}^{\leq N-1} \mytr \upchi^{(Small)}
			\right|.
	\end{align}
	The desired bound now follows
	from \eqref{E:LOWERORDERTERMSINVOLVINGSPATIALDERIVATIVESOFRGEO},
	\eqref{E:CRUDEPOINTWISEBOUNDSDERIVATIVESOFANGULARDEFORMATIONTENSORS}
	with $N-1$ in the role of $N,$
	\eqref{E:CRUDELOWERORDERC0BOUNDDERIVATIVESOFANGULARDEFORMATIONTENSORS}
	\eqref{E:POINTWISEESTIMATEZNAPPLIEDTOLOWESTTRCHIJUNKDISCREPANCYTERM}
	with $N-1$ in the role of $N,$
	and
	\eqref{E:LOWERORDERC0BOUNDZNAPPLIEDTOLOWESTTRCHIJUNKDISCREPANCYTERM}.
	
	We omit the proof of the desired bound for the first term on the third line 
	of \eqref{E:NTIMESCOMMUTATIONRENORMALIZEDTRCHIJUNKINHOMOGENEOUSTERM}
	because it can be bounded by using the same argument that we now use
	to bound the last term on the third line.
	To bound this last term 
	$- [\mathscr{S}^{N-1}, \upmu] \left(S \Lunit \mytr \upchi^{(Small)} \right)$	
	by the right-hand side of \eqref{E:TOPORDERRENORMALIZEDTRCHIJUNKTRANSPORTINHOMOGENEOUSTERMPOINTWISEESTIMATE},
	we first use the Leibniz rule to deduce that its magnitude is
	\begin{align} \label{E:HIGHERCOMMUTATORSMUANDSAPPLIEDTOLTRCHIJUNK}
		\lesssim 
		\left|
			\mathscr{Z}^{\leq 12} (\upmu - 1)
		\right|
		\left|
			\mathscr{Z}^{\leq N-1} \Lunit \mytr \upchi^{(Small)}
		\right|
		+ \lesssim 
		\left|
			\mathscr{Z}^{\leq 11} \Lunit \mytr \upchi^{(Small)}
		\right|
		\left|
			\mathscr{Z}^{\leq N-1} (\upmu - 1)
		\right|.
	\end{align}
	The desired bound now follows from \eqref{E:HIGHERCOMMUTATORSMUANDSAPPLIEDTOLTRCHIJUNK},
	\eqref{E:ZNAPPLIEDTORGEOISNOTTOOLARGE},
	and the estimates of Prop.~\ref{P:CRUCICALTRANSPORTINTEQUALITIES}.
	
	To bound the difference on the last line 
	of \eqref{E:NTIMESCOMMUTATIONRENORMALIZEDTRCHIJUNKINHOMOGENEOUSTERM}
	by the right-hand side of \eqref{E:TOPORDERRENORMALIZEDTRCHIJUNKTRANSPORTINHOMOGENEOUSTERMPOINTWISEESTIMATE}, 
	we use the Leibniz rule,
	the estimates 
	\eqref{E:CRUDEPOINTWISEBOUNDSDERIVATIVESOFANGULARDEFORMATIONTENSORS},
	\eqref{E:CRUDELOWERORDERC0BOUNDDERIVATIVESOFANGULARDEFORMATIONTENSORS},
	\eqref{E:POINTWISEEIKONALFUNCTIONTRANSPORTBASEDINTEQUALITIES},
	\eqref{E:C0BOUNDCRUCIALEIKONALFUNCTIONQUANTITIES} 
	and the bootstrap assumptions \eqref{E:PSIFUNDAMENTALC0BOUNDBOOTSTRAP}
	to deduce that
	\begin{align} \label{E:LOWERORDERUPMUTRFRECHISQUAREDTERMTOPORDERRENORMALIZEDTRCHIJUNKTRANSPORTINHOMOGENEOUSTERM}
		&
		\left|
			2 \upmu \hat{\upchi}^{(Small) \# \#} \angLie_{\mathscr{S}}^N \hat{\upchi}^{(Small)}
			- \mathscr{S}^N (\upmu \hat{\upchi}^{(Small) \# \#} \hat{\upchi}^{(Small)})
		\right|
			\\
		& \lesssim
			\mathop{\sum_{N_1 + N_2 + N_3 + N_4 + N_5 \leq N}}_{N_4, N_5 \leq N-1}
			\left|
				\mathscr{Z}^{N_1} \upmu 
			\right|
			\left|
				\angLie_{\mathscr{Z}}^{N_2} \ginversesphere
			\right|
			\left|
				\angLie_{\mathscr{Z}}^{N_3} \ginversesphere
			\right|
			\left|
				\angLie_{\mathscr{Z}}^{N_4} \hat{\upchi}^{(Small)}
			\right|
			\left|
				\angLie_{\mathscr{Z}}^{N_5} \hat{\upchi}^{(Small)}
			\right|
			\notag \\
		& \lesssim 
			\frac{\ln^2(\myexp + t)}{(1 + t)^3} 
			\left|
				\fourmyarray[\rgeo \Lunit \mathscr{Z}^{\leq N} \Psi]
					{\Rad \mathscr{Z}^{\leq N} \Psi}
					{\rgeo \angdiff \mathscr{Z}^{\leq N} \Psi} 
					{\mathscr{Z}^{\leq N} \Psi}
			\right|
		+ \frac{\ln(\myexp + t)}{(1 + t)^4}
			\left| 
				\myarray
					[\mathscr{Z}^{\leq N} (\upmu - 1)]
						{\rgeo \sum_{a=1}^3 |\mathscr{Z}^{\leq N} \Lunit_{(Small)}^a|}
				\right|.
			\notag
	\end{align}
	We have thus proved the desired estimate \eqref{E:TOPORDERRENORMALIZEDTRCHIJUNKTRANSPORTINHOMOGENEOUSTERMPOINTWISEESTIMATE}.
\end{proof}

\section{Preliminary pointwise estimates for 
\texorpdfstring{$\chipartialmodsourcearg{\mathscr{S}^{N-1}}$}
{the inhomogeneous terms in the partially modified equations}}
\label{S:PRELIMINARYPOINTWISEPARTIALLYMODIFIED}
Recall that Lemma~\ref{L:COMMUTEDTRCHIJUNKFIRSTPARTIALRENORMALIZEDTRANSPORTEQUATION}
provides the transport equation verified by the partially modified quantity
$\rgeo^2 \chipartialmodarg{\mathscr{S}^{N-1}}$
and that $\chipartialmodsourcearg{\mathscr{S}^{N-1}}$
is a source term in the equation. In the next lemma, we derive pointwise estimates
for this source term.

\begin{lemma}[\textbf{Pointwise estimates for the inhomogeneous terms corresponding to the partially modified version of}
$\mathscr{S}^{N-1} \mytr \upchi^{(Small)}$]
\label{L:POINTWISEESTIMATESTRCHIJUNKREFINEDPARTIALLYRENORMALIZDINHOMOGENEOUSTERM}
Let $1 \leq N \leq 24$ be an integer 
and let $ \mathscr{S}^{N-1}$ be an $(N-1)^{st}$ order pure spatial commutation vectorfield operator
(see definition \eqref{E:DEFSETOFSPATIALCOMMUTATORVECTORFIELDS}).
Let $\chipartialmodsourcearg{\mathscr{S}^{N-1}}$ be the inhomogeneous term
defined in \eqref{E:TRCHIJUNKCOMMUTEDTRANSPORTEQNPARTIALRENORMALIZATIONINHOMOGENEOUSTERM}.
Under the small-data and bootstrap assumptions 
of Sects.~\ref{S:PSISOLVES}-\ref{S:C0BOUNDBOOTSTRAP},
if $\varepsilon$ is sufficiently small, then
the following estimates hold on $\mathcal{M}_{\Tboot,U_0}:$
\begin{align}  \label{E:POINTWISEESTIMATESTRCHIJUNKREFINEDPARTIALLYRENORMALIZDINHOMOGENEOUSTERM}
	\left|
	\chipartialmodsourcearg{\mathscr{S}^{N-1}}
	\right|
	& \lesssim
		\left| 
				\fourmyarray[ \rgeo \Lunit \mathscr{Z}^{\leq N-1} \Psi]
					{\Rad \mathscr{Z}^{\leq N-1} \Psi}
					{\rgeo \angdiff \mathscr{Z}^{\leq N-1} \Psi}
					{\mathscr{Z}^{\leq N-1} \Psi}
		\right|
		+ 
					\frac{\ln(\myexp + t)}{(1 + t)^2}
						\left| 
							\myarray
								[\mathscr{Z}^{\leq N} (\upmu - 1)]
								{\rgeo \sum_{a=1}^3 |\mathscr{Z}^{\leq N} \Lunit_{(Small)}^a|}
						\right|.
\end{align}
\end{lemma}

\begin{proof}
	We start by bounding the first term $\mathscr{S}^{N-1} (\rgeo^2 \mathfrak{B})$ on the right-hand
	side of \eqref{E:TRCHIJUNKCOMMUTEDTRANSPORTEQNPARTIALRENORMALIZATIONINHOMOGENEOUSTERM}.
	To this end, we apply $\angLie_{\mathscr{S}}^{N-1}$ to $\rgeo^2$ times the
	right-hand side of \eqref{E:PARTIALRENORMALIZEDALPHAINHOMOGENEOUSTERM} and apply the Leibniz rule. We bound 
	$\mathscr{S}^M \rgeo$ with \eqref{E:EXACTRELATIONSZAPPLIEDTORGEO}.
	We bound the terms $\angLie_{\mathscr{S}}^M \myarray[G_{(Frame)}]{G_{(Frame)}'}$ 
	with the estimates of Lemma~\ref{L:POINTWISEESTIMATESGFRAMEINTERMSOFOTHERQUANTITIES}.
	We bound $\angLie_{\mathscr{S}}^M \ginversesphere$ and $\angLie_{\mathscr{S}}^M \mytr \upchi$
	with Lemma~\ref{L:POINTWISEBOUNDSDERIVATIVESOFANGULARDEFORMATIONTENSORS}.
	We bound the terms
	$\angLie_{\mathscr{S}}^M
				\myarray
					[\Lunit \Psi]
					{\angdiff \Psi}$
	with Lemma~\ref{L:AVOIDINGCOMMUTING}
	and the bootstrap assumptions \eqref{E:PSIFUNDAMENTALC0BOUNDBOOTSTRAP}.	In total, these estimates yield that
	$\left|\mathscr{S}^{N-1} (\rgeo^2 \mathfrak{B}) \right|$ is $\lesssim$
	the right-hand side of \eqref{E:POINTWISEESTIMATESTRCHIJUNKREFINEDPARTIALLYRENORMALIZDINHOMOGENEOUSTERM}
	as desired.
	
	We next bound the second term $\mathscr{S}^{N-1} (\rgeo^2 |\upchi^{(Small)}|^2)$
	on the right-hand side of \eqref{E:TRCHIJUNKCOMMUTEDTRANSPORTEQNPARTIALRENORMALIZATIONINHOMOGENEOUSTERM}.
	From \eqref{E:ZNAPPLIEDTORGEOISNOTTOOLARGE} and \eqref{E:C0BOUNDCRUCIALEIKONALFUNCTIONQUANTITIES},
	we deduce that the magnitude of this term is
	\begin{align} \label{E:PRELIMINARYPOINTWISEQUADRATICCHIJUNKTERMNOTTOPORDER}
		\lesssim \varepsilon \ln(\myexp + t) 
		\left|
			\angLie_{\mathscr{S}}^{\leq N-1}
			\upchi^{(Small)}
		\right|.
	\end{align}
	Using inequality \eqref{E:POINTWISEESTIMATESFORCHIJUNKINTERMSOFOTHERVARIABLES}
	with $N-1$ in the role of $N,$
	we deduce that \eqref{E:PRELIMINARYPOINTWISEQUADRATICCHIJUNKTERMNOTTOPORDER}
	is $\lesssim$
	the right-hand side of \eqref{E:POINTWISEESTIMATESTRCHIJUNKREFINEDPARTIALLYRENORMALIZDINHOMOGENEOUSTERM}
	as desired.
	
	We next bound the first term $\frac{1}{2} [\rgeo^2 G_{\Lunit \Lunit}, \mathscr{S}^{N-1}] \angLap \Psi$
	on the second line of the right-hand side of \eqref{E:TRCHIJUNKCOMMUTEDTRANSPORTEQNPARTIALRENORMALIZATIONINHOMOGENEOUSTERM}.
	The second term on this line can be bounded by using a similar argument, and we omit those details.
	From 
	\eqref{E:EXACTRELATIONSZAPPLIEDTORGEO}
	and the Leibniz rule,
	we deduce that the magnitude of the first term is
	\begin{align}  \label{E:PRELIMINARYPOINTWISEPRODUCTESTIMATEINVOLVINGANGLAP}
		\lesssim 
		\rgeo^2
		\left|
			\mathscr{S}^{\leq 12} G_{\Lunit \Lunit}
		\right|
		\left|
			\mathscr{S}^{\leq N-2} \angLap \Psi
		\right|
		+ \rgeo^2
			\left|
				\mathscr{S}^{\leq 11} \angLap \Psi
			\right|
			\left|
				\mathscr{S}^{\leq N-1} G_{\Lunit \Lunit}
			\right|.
	\end{align}
	From inequalities 
	\eqref{E:LIEDERIVATIVESOFGRAMEINTERMSOFOTHERVARIABLES},
	\eqref{E:LOWERORDERC0BOUNDLIEDERIVATIVESOFGRAME},
	\eqref{E:POINTWISECOMMUTINGANGDSQUAREDPSIWITHLIEZN},
	and
	\eqref{E:LINFFTYLOWERORDERCOMMUTINGANGDSQUAREDPSIWITHLIEZN},
	and the bootstrap assumptions \eqref{E:PSIFUNDAMENTALC0BOUNDBOOTSTRAP},
	we deduce that the magnitude of this 
	right-hand side of \eqref{E:PRELIMINARYPOINTWISEPRODUCTESTIMATEINVOLVINGANGLAP}
	is $\lesssim$
	the right-hand side of \eqref{E:POINTWISEESTIMATESTRCHIJUNKREFINEDPARTIALLYRENORMALIZDINHOMOGENEOUSTERM}
	as desired.
	
	We next bound the last term $[\Lunit, \mathscr{S}^{N-1}](\rgeo^2 \mytr \upchi^{(Small)})$
	on the second line of the right-hand side of \eqref{E:TRCHIJUNKCOMMUTEDTRANSPORTEQNPARTIALRENORMALIZATIONINHOMOGENEOUSTERM}.
	Using the commutator estimate \eqref{E:LSNCOMMUTATORACTINGONFUNCTIONSPOINTWISE}
	with $N-1$ in the role of $N$ and $\rgeo^2 \mytr \upchi^{(Small)}$ in the role of $f,$
	and inequalities \eqref{E:ZNAPPLIEDTORGEOISNOTTOOLARGE},
	\eqref{E:POINTWISEESTIMATESFORCHIJUNKINTERMSOFOTHERVARIABLES},
	\eqref{E:LDERIVATIVECRUCICALTRANSPORTINTEQUALITIES},
	and
	\eqref{E:C0BOUNDCRUCIALEIKONALFUNCTIONQUANTITIES},
	we deduce that the magnitude of this term is
	$\lesssim$
	the right-hand side of \eqref{E:POINTWISEESTIMATESTRCHIJUNKREFINEDPARTIALLYRENORMALIZDINHOMOGENEOUSTERM}
	as desired.
		
	We next bound the first term $[\Lunit, \mathscr{S}^{N-1}](\rgeo^2 \chipartialmodinhom)$
	on the last line of the right-hand side of \eqref{E:TRCHIJUNKCOMMUTEDTRANSPORTEQNPARTIALRENORMALIZATIONINHOMOGENEOUSTERM}.	
	Referring to definition \eqref{E:LOWESTORDERTRANSPORTPARTIALRENORMALIZEDTRCHIJUNKDISCREPANCY},
	we see that we have to bound
	\begin{align} \label{E:BSECONDLINEBOUNDLASTTERM}
	[\angLie_{\Lunit}, \angLie_{\mathscr{S}}^{N-1}]
	\left\lbrace
		- \frac{1}{2} \rgeo^2 \angGmixedarg{A}{A} \Lunit \Psi
		- \frac{1}{2} \rgeo^2 G_{\Lunit \Lunit} \Lunit \Psi
		+ \rgeo^2 \angGmixedarg{\Lunit}{A} \angdiffarg{A} \Psi
	\right\rbrace.
	\end{align}
	From Lemma~\ref{L:LANDRADCOMMUTEWITHANGDIFF},
	inequality \eqref{E:ZNAPPLIEDTORGEOISNOTTOOLARGE},
	inequality \eqref{E:FUNCTIONPOINTWISEANGDINTERMSOFANGLIEO},
	the estimates of Cor.~\ref{C:POINTWISEESTIMATESGFRAMESHARPINTERMSOFOTHERQUANTITIES},
	inequality \eqref{E:LSNCOMMUTATORACTINGONFUNCTIONSPOINTWISE} 
	with $N-1$ in the role of $N$ and $f$ equal to the terms in braces in \eqref{E:BSECONDLINEBOUNDLASTTERM},
	and the bootstrap assumptions \eqref{E:PSIFUNDAMENTALC0BOUNDBOOTSTRAP},	
	we deduce that
	\begin{align}
		\left|
			[\angLie_{\Lunit}, \angLie_{\mathscr{S}}^{N-1}]
			\left\lbrace
				- \frac{1}{2} \rgeo^2 \angGmixedarg{A}{A} \Lunit \Psi
				- \frac{1}{2} \rgeo^2 G_{\Lunit \Lunit} \Lunit \Psi
				+ \rgeo^2 \angGmixedarg{\Lunit}{A} \angdiffarg{A} \Psi
			\right\rbrace
		\right|
		& \lesssim
			\frac{\ln(\myexp + t)}{1 + t}
			\left| 
				\fourmyarray[ \rgeo \Lunit \mathscr{Z}^{\leq N-1} \Psi]
					{\Rad \mathscr{Z}^{\leq N-1} \Psi}
					{\rgeo \angdiff \mathscr{Z}^{\leq N-1} \Psi}
					{\mathscr{Z}^{\leq N-1} \Psi}
			\right|
			\\
	&  \ \ 
			+ 
			\frac{1}{(1 + t)^2}
			\left|
			\myarray[\mathscr{Z}^{\leq N} (\upmu - 1)]
				{\sum_{a=1}^3 \rgeo |\mathscr{Z}^{\leq N} \Lunit_{(Small)}^a|} 
			\right|,
		\notag 
	\end{align}
	which is 
	$\lesssim$
	the right-hand side of \eqref{E:POINTWISEESTIMATESTRCHIJUNKREFINEDPARTIALLYRENORMALIZDINHOMOGENEOUSTERM}
	as desired.
	
	We now bound the second term 
	on the last line of the right-hand side of \eqref{E:TRCHIJUNKCOMMUTEDTRANSPORTEQNPARTIALRENORMALIZATIONINHOMOGENEOUSTERM}.
	We first use \eqref{E:EXACTRELATIONSZAPPLIEDTORGEO} 
	and the fact that $\rgeo \Lunit \in \mathscr{Z}$
	to deduce that
	\begin{align} \label{E:NEXTTOLASTPARTIALLMODIFIEDESTIMATE}
		\left|
			\Lunit 
			\left\lbrace
				[\rgeo^2 , \mathscr{S}^{N-1}] \mytr \upchi^{(Small)}
			\right\rbrace
		\right|
		& \lesssim 
			\left|
				\mathscr{Z}^{\leq N-1} \mytr \upchi^{(Small)}
			\right|.
	\end{align}
	To conclude that the right-hand side of \eqref{E:NEXTTOLASTPARTIALLMODIFIEDESTIMATE}
	is 
	$\lesssim$
	the right-hand side of \eqref{E:POINTWISEESTIMATESTRCHIJUNKREFINEDPARTIALLYRENORMALIZDINHOMOGENEOUSTERM}
	as desired,
	we use inequality \eqref{E:POINTWISEESTIMATESFORCHIJUNKINTERMSOFOTHERVARIABLES}
	with $N-1$ in the role of $N.$
	
	We now bound the final term 
	on the last line of the right-hand side of \eqref{E:TRCHIJUNKCOMMUTEDTRANSPORTEQNPARTIALRENORMALIZATIONINHOMOGENEOUSTERM}.
	We first use definitions \eqref{E:LOWESTORDERTRANSPORTPARTIALRENORMALIZEDTRCHIJUNKDISCREPANCY}
	and \eqref{E:TRANSPORTPARTIALRENORMALIZEDTRCHIJUNKDISCREPANCY}
	to deduce that the term can be rewritten as follows:
	\begin{align} \label{E:FINALTERMSPLITINTOTHREETERMS}
	\Lunit	
				\left\lbrace
					\rgeo^2 \chipartialmodinhomarg{\mathscr{S}^{N-1}}
					- \mathscr{S}^{N-1} (\rgeo^2 \chipartialmodinhom)
				\right\rbrace
	& =	
		- \frac{1}{2} \Lunit [\rgeo \angGmixedarg{A}{A}, \angLie_{\mathscr{S}}^{N-1}] (\rgeo \Lunit \Psi)
		- \frac{1}{2} \Lunit [\rgeo G_{\Lunit \Lunit}, \angLie_{\mathscr{S}}^{N-1}] (\rgeo \Lunit \Psi)
			\\
	& \ \ 
		+ \angLie_{\Lunit} [\rgeo^2 \angGmixedarg{\Lunit}{A}, \angLie_{\mathscr{S}}^{N-1}] \angdiffarg{A} \Psi.
		\notag
	\end{align}
	We show how to bound the last term on the right-hand side of \eqref{E:FINALTERMSPLITINTOTHREETERMS};
	the first two terms can be bounded using essentially the same reasoning.
	From 
	the fact that $\rgeo \Lunit \in \mathscr{Z},$
	Lemma~\ref{L:LANDRADCOMMUTEWITHANGDIFF},
	inequality \eqref{E:ZNAPPLIEDTORGEOISNOTTOOLARGE},
	the estimates of Cor.~\ref{C:POINTWISEESTIMATESGFRAMESHARPINTERMSOFOTHERQUANTITIES},
	and the bootstrap assumptions \eqref{E:PSIFUNDAMENTALC0BOUNDBOOTSTRAP},
	we deduce that
	\begin{align} \label{E:FIRSTESTIMATEFINALTERMCOMMUTEDTRCHIJUNKFIRSTPARTIALRENORMALIZEDTRANSPORTINHOMOGENEOUSTERM}
		\rgeo
		\left|
			\angLie_{\Lunit}
			[\rgeo^2 \angGmixedarg{\Lunit}{A}, \angLie_{\mathscr{S}}^{N-1}] 
			\angdiffarg{A} \Psi
		\right|
		& \lesssim 
			\frac{1}{(1 + t)^2} 
			\left| \mathscr{Z}^{\leq N} (\upmu - 1) \right|
			+ \ln(\myexp + t)
				\left|
					\angdiff \mathscr{Z}^{\leq N-1} \Psi 
				\right|
				\\
		& \ \ 
			+ \ln(\myexp + t)
			\left|
				\angdiff \mathscr{Z}^{\leq 12} \Psi
			\right|
			\left|	
				\angLie_{\mathscr{Z}}^{\leq N} \angGmixedarg{\Lunit}{\#}
			\right|.
			\notag
	\end{align}
	Cor.~\ref{C:POINTWISEESTIMATESGFRAMESHARPINTERMSOFOTHERQUANTITIES},
	inequality \eqref{E:LIEDERIVATIVESOFGRAMEINTERMSOFOTHERVARIABLES},
	and the bootstrap assumptions \eqref{E:PSIFUNDAMENTALC0BOUNDBOOTSTRAP}
	then yield that the right-hand side of 
	\eqref{E:FIRSTESTIMATEFINALTERMCOMMUTEDTRCHIJUNKFIRSTPARTIALRENORMALIZEDTRANSPORTINHOMOGENEOUSTERM}
	is $\lesssim$
	the right-hand side of \eqref{E:POINTWISEESTIMATESTRCHIJUNKREFINEDPARTIALLYRENORMALIZDINHOMOGENEOUSTERM}
	as desired.

\end{proof}

Recall that Lemma~\ref{L:COMMUTEDUPMUFIRSTPARTIALRENORMALIZEDTRANSPORTEQUATION}
provides the transport equation verified by the partially modified quantity
$\angdiff \mathscr{S}^{N-1} \upmu$
and that $\mupartialmodsourcearg{\mathscr{S}^{N-1}}$
is a source term in the equation. In the next lemma, we derive pointwise estimates
for this source term.

\begin{lemma}[\textbf{Pointwise estimates for the inhomogeneous terms corresponding to the partially modified version of}
$\angdiff \mathscr{S}^{N-1} \upmu$]
\label{L:POINTWISEESTIMATESUPMUPARTIALLYRENORMALIZDINHOMOGENEOUSTERM}
Let $1 \leq N \leq 24$ be an integer 
and let $ \mathscr{S}^{N-1}$ be an $(N-1)^{st}$ order pure spatial commutation vectorfield operator
(see definition \eqref{E:DEFSETOFSPATIALCOMMUTATORVECTORFIELDS}).
Let $\mupartialmodsourcearg{\mathscr{S}^{N-1}}$ be the $S_{t,u}$ one-form inhomogeneous term
\eqref{E:COMMUTEDUPMUFIRSTPARTIALRENORMALIZEDTRANSPORTINHOMOGENEOUSTERM}.
Under the small-data and bootstrap assumptions 
of Sects.~\ref{S:PSISOLVES}-\ref{S:C0BOUNDBOOTSTRAP},
if $\varepsilon$ is sufficiently small, then
the following estimates hold on $\mathcal{M}_{\Tboot,U_0}:$
\begin{align}  \label{E:POINTWISEESTIMATESUPMUREFINEDPARTIALLYRENORMALIZDINHOMOGENEOUSTERM}
	\rgeo
	\left|
		\mupartialmodsourcearg{\mathscr{S}^{N-1}}
	\right|
	& \lesssim
		\left| 
				\fourmyarray[ \rgeo \Lunit \mathscr{Z}^{\leq N-1} \Psi]
					{\Rad \mathscr{Z}^{\leq N-1} \Psi}
					{\rgeo \angdiff \mathscr{Z}^{\leq N-1} \Psi}
					{\mathscr{Z}^{\leq N-1} \Psi}
		\right|
		+ 
					\frac{\ln(\myexp + t)}{(1 + t)^2}
						\left| 
							\myarray
								[\mathscr{Z}^{\leq N} (\upmu - 1)]
								{\rgeo \sum_{a=1}^3 |\mathscr{Z}^{\leq N} \Lunit_{(Small)}^a|}
						\right|.
\end{align}
\end{lemma}

\begin{proof}
	We start by bounding the
	product $\rgeo \angLie_{\mathscr{S}}^{N-1} \mathfrak{J}$ arising from the first term
	on the right-hand side of \eqref{E:COMMUTEDUPMUFIRSTPARTIALRENORMALIZEDTRANSPORTINHOMOGENEOUSTERM},
	where the $S_{t,u}$ one-form $\mathfrak{J}$ is defined in \eqref{E:UPMUFIRSTPARTIALRENORMALIZEDTRANSPORTINHOMOGENEOUSTERM}.
	We apply $\rgeo$ times $\angLie_{\mathscr{S}}^{N-1}$ to the right-hand side of 	
	\eqref{E:UPMUFIRSTPARTIALRENORMALIZEDTRANSPORTINHOMOGENEOUSTERM}
	and apply the Leibniz rule. 
	We bound the terms $\mathscr{S}^M G_{(Frame)},$ 
	$\mathscr{S}^M \Lunit G_{(Frame)},$ 
	and $\angLie_{\mathscr{S}}^M \angdiff G_{(Frame)}$
	by using the fact that $\rgeo \Lunit \in \mathscr{Z},$
	inequality \eqref{E:ZNAPPLIEDTORGEOISNOTTOOLARGE},
	Lemma~\ref{L:LANDRADCOMMUTEWITHANGDIFF},
	Lemma~\ref{L:AVOIDINGCOMMUTING},
	inequality \eqref{E:FUNCTIONPOINTWISEANGDINTERMSOFANGLIEO},
	and the estimates of
	Lemma~\ref{L:POINTWISEESTIMATESGFRAMEINTERMSOFOTHERQUANTITIES}.
	We bound the terms
	$\mathscr{S}^M \upmu$ 
	and $\angLie_{\mathscr{S}}^M \angdiff \upmu$
	with Lemma~\ref{L:LANDRADCOMMUTEWITHANGDIFF} 
	and inequality \eqref{E:C0BOUNDCRUCIALEIKONALFUNCTIONQUANTITIES}.
	We bound the terms
	$\angLie_{\mathscr{S}}^M
				\threemyarray
					[\Lunit \Psi]
					{\Rad \Psi}
					{\angdiff \Psi}$
	with Lemma~\ref{L:AVOIDINGCOMMUTING}
	and the bootstrap assumptions \eqref{E:PSIFUNDAMENTALC0BOUNDBOOTSTRAP}.
	In total, these estimates yield that
	$\rgeo \left|\angLie_{\mathscr{S}}^{N-1} \mathfrak{J} \right|$ is $\lesssim$
	the right-hand side of \eqref{E:POINTWISEESTIMATESUPMUREFINEDPARTIALLYRENORMALIZDINHOMOGENEOUSTERM}
	as desired.
	
	To bound $\rgeo$ times the second term on the right-hand side of 		
	\eqref{E:COMMUTEDUPMUFIRSTPARTIALRENORMALIZEDTRANSPORTINHOMOGENEOUSTERM},
	we first use the Leibniz rule and Lemma~\ref{L:LANDRADCOMMUTEWITHANGDIFF}
	to deduce that
	\begin{align} \label{E:FIRSTESTIMATESECONDTERMCOMMUTEDUPMUFIRSTPARTIALRENORMALIZEDTRANSPORTINHOMOGENEOUSTERM}
		\rgeo
		\left|
				[\angLie_{\mathscr{S}}^{N-1}, G_{\Lunit \Lunit}]
				\angdiff \Rad \Psi
		\right|
		& \lesssim 
			\left| 
				\angLie_{\mathscr{S}}^{\leq 12} G_{\Lunit \Lunit}
			\right|
			\left| 
				\rgeo \angdiff \mathscr{S}^{\leq N-2} \Rad \Psi
			\right|
		+ \rgeo
			\left|
				\angdiff \mathscr{S}^{\leq 12} \Psi 
			\right|
			\left| 
				\angLie_{\mathscr{S}}^{\leq N-1} G_{\Lunit \Lunit}
			\right|.
	\end{align}
	Inequalities \eqref{E:LIEDERIVATIVESOFGRAMEINTERMSOFOTHERVARIABLES}
	and \eqref{E:LOWERORDERC0BOUNDLIEDERIVATIVESOFGRAME}
	and the bootstrap assumptions \eqref{E:PSIFUNDAMENTALC0BOUNDBOOTSTRAP}
	then yield that the right-hand side of 
	\eqref{E:FIRSTESTIMATESECONDTERMCOMMUTEDUPMUFIRSTPARTIALRENORMALIZEDTRANSPORTINHOMOGENEOUSTERM}
	is $\lesssim$
	the right-hand side of \eqref{E:POINTWISEESTIMATESUPMUREFINEDPARTIALLYRENORMALIZDINHOMOGENEOUSTERM}
	as desired.
	
	To bound $\rgeo$ times the third term on the right-hand side of 		
	\eqref{E:COMMUTEDUPMUFIRSTPARTIALRENORMALIZEDTRANSPORTINHOMOGENEOUSTERM},
	we use inequality \eqref{E:LSNCOMMUTATORACTINGONTENSORFIELDSPOINTWISE} 
	with $N-1$ in the role of $N$ and
	$\angdiff \upmu$ in the role of $\xi,$ Lemma~\ref{L:LANDRADCOMMUTEWITHANGDIFF}, and inequalities
	\eqref{E:FUNCTIONPOINTWISEANGDINTERMSOFANGLIEO}
	and \eqref{E:LDERIVATIVECRUCICALTRANSPORTINTEQUALITIES}
	to deduce the following desired bound:
	\begin{align}
		\rgeo
		\left|
			[\angLie_{\Lunit}, \angLie_{\mathscr{S}}^{N-1}] 
				\angdiff \upmu
		\right|
	 & \lesssim 
	 			\left| 
				\fourmyarray[ \rgeo \Lunit \mathscr{Z}^{\leq N-1} \Psi]
					{\rgeo \angdiff \mathscr{Z}^{\leq N-1} \Psi}
					{\Rad \mathscr{Z}^{\leq N-1} \Psi}
					{\mathscr{Z}^{\leq N-1} \Psi}
				\right|
		+ \frac{\ln(\myexp + t)}{(1 + t)^2}
			\left|
				\myarray[\mathscr{Z}^{\leq N} (\upmu - 1)]
					{\sum_{a=1}^3 \rgeo |\mathscr{Z}^{\leq N} \Lunit_{(Small)}^a|} 
			\right|.
	\end{align}	
	
	To bound $\rgeo$ times the fourth term on the right-hand side of 		
	\eqref{E:COMMUTEDUPMUFIRSTPARTIALRENORMALIZEDTRANSPORTINHOMOGENEOUSTERM},
	we refer to definition \eqref{E:LOWESTORDERTRANSPORTPARTIALRENORMALIZEDUPMUDISCREPANCY}
	to see that we have to bound 
	$\rgeo [\angLie_{\Lunit}, \angLie_{\mathscr{S}}^{N-1}] 
			\left\lbrace
					\frac{1}{2} \upmu  G_{\Lunit \Lunit} \angdiff \Psi
					+ \upmu G_{\Lunit \Radunit} \angdiff \Psi
			\right\rbrace.$
	We now address the term $\upmu G_{\Lunit \Lunit} \angdiff \Psi$ 
	(the term $\upmu G_{\Lunit \Radunit} \angdiff \Psi$ can be handled using the same arguments).
	The estimates \eqref{E:LOWERORDERC0BOUNDLIEDERIVATIVESOFGRAME},
	\eqref{E:C0BOUNDCRUCIALEIKONALFUNCTIONQUANTITIES},
	the bootstrap assumptions \eqref{E:PSIFUNDAMENTALC0BOUNDBOOTSTRAP},
	and Lemma~\ref{L:LANDRADCOMMUTEWITHANGDIFF} together imply that 
	$\left|\angLie_{\mathscr{Z}}^{\leq 12} (\upmu  G_{\Lunit \Lunit} \angdiff \Psi) \right| 
	\lesssim \varepsilon \ln(\myexp + t)(1 + t)^{-2}.$
	Combining this estimate with inequality 
	\eqref{E:LSNCOMMUTATORACTINGONTENSORFIELDSPOINTWISE},
	where $\upmu  G_{\Lunit \Lunit} \angdiff \Psi$ plays the role of $\xi$
	and $N-1$ plays the role of $N,$
	we deduce that
	\begin{align} \label{E:FIRSTESTIMATESTHIRDTERMCOMMUTEDUPMUFIRSTPARTIALRENORMALIZEDTRANSPORTINHOMOGENEOUSTERM}
		\rgeo
		\left|
			[\angLie_{\Lunit}, \angLie_{\mathscr{S}}^{N-1}] 
			(\upmu  G_{\Lunit \Lunit} \angdiff \Psi)
		\right|
		& \lesssim 
			\frac{\ln(\myexp + t)}{1 + t}
			\left|
				\angLie_{\mathscr{Z}}^{\leq N-1}
				(\upmu  G_{\Lunit \Lunit} \angdiff \Psi)
			\right|
			\\
		& \ \
			+
			\frac{\ln(\myexp + t)}{(1 + t)^2}
				\left| 
				\fourmyarray[ \rgeo \Lunit \mathscr{Z}^{\leq N-1} \Psi]
					{\rgeo \angdiff \mathscr{Z}^{\leq N-1} \Psi}
					{\Rad \mathscr{Z}^{\leq N-1} \Psi}
					{\mathscr{Z}^{\leq N-1} \Psi}
				\right|
			+ 
			\frac{\ln(\myexp + t)}{(1 + t)^3}
			\left|
				\myarray[\mathscr{Z}^{\leq N} (\upmu - 1)]
					{\sum_{a=1}^3 \rgeo |\mathscr{Z}^{\leq N} \Lunit_{(Small)}^a|} 
			\right|.
		\notag	
	\end{align}
	The last two terms in \eqref{E:FIRSTESTIMATESTHIRDTERMCOMMUTEDUPMUFIRSTPARTIALRENORMALIZEDTRANSPORTINHOMOGENEOUSTERM}
	are manifestly $\lesssim$
	the right-hand side of \eqref{E:POINTWISEESTIMATESUPMUREFINEDPARTIALLYRENORMALIZDINHOMOGENEOUSTERM}
	as desired. To bound the first term by
	the right-hand side of \eqref{E:UPMUANDTRCHIJUNKREFINEDPARTIALLYRENORMALIZEDTRANSPORTINEQUALITIES} 
	we use the Leibniz rule, 
	Lemma~\ref{L:LANDRADCOMMUTEWITHANGDIFF}, 
	the estimates
	\eqref{E:LIEDERIVATIVESOFGRAMEINTERMSOFOTHERVARIABLES},
	\eqref{E:LOWERORDERC0BOUNDLIEDERIVATIVESOFGRAME},
	\eqref{E:C0BOUNDCRUCIALEIKONALFUNCTIONQUANTITIES},
	and the bootstrap assumptions \eqref{E:PSIFUNDAMENTALC0BOUNDBOOTSTRAP}.
	
	To bound $\rgeo$ times the last term on the right-hand side of 		
	\eqref{E:COMMUTEDUPMUFIRSTPARTIALRENORMALIZEDTRANSPORTINHOMOGENEOUSTERM},
	we use definitions \eqref{E:LOWESTORDERTRANSPORTPARTIALRENORMALIZEDUPMUDISCREPANCY} 
	and \eqref{E:TRANSPORTPARTIALRENORMALIZEDUPMUDISCREPANCY} 
	and Lemma~\ref{L:LANDRADCOMMUTEWITHANGDIFF}
	to deduce that the term can be rewritten as follows:
	\begin{align} \label{E:LASTTERMSPLITINTOTWOTERMS}
	\angLie_{\Lunit}
			\left\lbrace
				\mupartialmodinhomarg{\mathscr{S}^{N-1}}
				- \angLie_{\mathscr{S}}^{N-1}  \widetilde{\mathfrak{M}}
			\right\rbrace
		& = \frac{1}{2} 
			\angLie_{\Lunit}
			[\angLie_{\mathscr{S}}^{N-1}, \upmu G_{\Lunit \Lunit}] 
			\angdiff \Psi
			\\
	& \ \ + 
			\angLie_{\Lunit}	
			[\angLie_{\mathscr{S}}^{N-1}, \upmu G_{\Lunit \Radunit}] 
			\angdiff \Psi.
			\notag
	\end{align}
	We show how to bound the first term on the right-hand
	side of \eqref{E:LASTTERMSPLITINTOTWOTERMS};
	the proof of the bound for the second one is identical.	
	From the fact that $\rgeo \Lunit \in \mathscr{Z},$
	Lemma~\ref{L:LANDRADCOMMUTEWITHANGDIFF},
	the estimates 
	\eqref{E:LOWERORDERC0BOUNDLIEDERIVATIVESOFGRAME}
	and
	\eqref{E:C0BOUNDCRUCIALEIKONALFUNCTIONQUANTITIES},
	and the bootstrap assumptions \eqref{E:PSIFUNDAMENTALC0BOUNDBOOTSTRAP},
	we deduce that
	\begin{align} \label{E:FIRSTESTIMATEFOURTHTERMCOMMUTEDUPMUFIRSTPARTIALRENORMALIZEDTRANSPORTINHOMOGENEOUSTERM}
		\rgeo
		\left|
			\angLie_{\Lunit} 
			[\angLie_{\mathscr{S}}^{N-1}, \upmu G_{\Lunit \Lunit}] 
			\angdiff \Psi
		\right|
		& \lesssim 
			\frac{1}{(1 + t)^2} 
			\left| \mathscr{Z}^{\leq N} (\upmu - 1) \right|
			+ \ln(\myexp + t)
				\left|
					\angdiff \mathscr{Z}^{\leq N-1} \Psi 
				\right|
				\\
		& \ \ 
			+ \ln(\myexp + t)
			\left|
				\angdiff \mathscr{Z}^{\leq 12} \Psi
			\right|
			\left|	
				\mathscr{Z}^{\leq N} G_{\Lunit \Lunit} 
			\right|.
			\notag
	\end{align}
	Inequality \eqref{E:LIEDERIVATIVESOFGRAMEINTERMSOFOTHERVARIABLES}
	and the bootstrap assumptions \eqref{E:PSIFUNDAMENTALC0BOUNDBOOTSTRAP}
	then yield that the right-hand side of 
	\eqref{E:FIRSTESTIMATEFOURTHTERMCOMMUTEDUPMUFIRSTPARTIALRENORMALIZEDTRANSPORTINHOMOGENEOUSTERM}
	is $\lesssim$
	the right-hand side of \eqref{E:POINTWISEESTIMATESUPMUREFINEDPARTIALLYRENORMALIZDINHOMOGENEOUSTERM}
	as desired.
	
\end{proof}

\section{Inverting the transport equation verified by the fully modified version of 
\texorpdfstring{$\mathscr{S}^N \mytr \upchi^{(Small)}$}{the spatial derivatives of the trace of
the re-centered null second fundamental form}}
\label{S:INVERTINGFULLYRENORMALIZEDTRANSPORT}
In Sect.~\ref{S:FULLYRENORMALIZEDPOINTWISEESTIMATESFORDIFFICULTENERGYERRORINTEGRAND},
we derive pointwise estimates for the most difficult top-order
wave equation error integrands
that arise from the multiplier $\Mult.$
Our derivation of these pointwise estimates is based on a detailed analysis 
of the fully modified quantities 
defined in Sect.~\ref{S:FULLRENORMALIZATIONFORTRCHIJUNK}.
In this section, we carry out a preliminary step.
Specifically, we invert the 
transport equation \eqref{E:TOPORDERTRCHIJUNKRENORMALIZEDTRANSPORT} 
verified by the fully modified quantities $\chifullmodarg{\mathscr{S}^N}$
in order to obtain preliminary pointwise estimates for them, 
where the $\mathscr{S}^N$ are $N^{th}$ order pure spatial commutation vectorfield operators. 
One minor difficulty is that the estimates for distinct 
$\chifullmodarg{\mathscr{S}^N}$
are weakly coupled due to the presence of the important top-order term
$\upmu [\Lunit, \mathscr{S}^N] \mytr \upchi^{(Small)}$
on the right-hand side of \eqref{E:TOPORDERTRCHIJUNKRENORMALIZEDTRANSPORT}.
We overcome this difficulty by a deriving Gronwall-type estimate that 
simultaneously involves all of the $N^{th}$ order quantities $\chifullmodarg{\mathscr{S}^N}.$
The details are provided in the next lemma.

\begin{lemma}[\textbf{Inverting the transport equation verified by} $\chifullmodarg{\mathscr{S}^N}$]
	\label{L:RENORMALIZEDTOPORDERTRCHIJUNKTRANSPORTINVERTED}
	Let $0 \leq N \leq 24$ be an integer.
	Let $\chifullmodarg{\mathscr{S}^N},$ $\mathfrak{X},$ and 
	$\chifullmodsourcearg{\mathscr{S}^N}$ be as in
	Prop.~\ref{P:TOPORDERTRCHIJUNKRENORMALIZEDTRANSPORT}.
	Under the small-data and bootstrap assumptions 
	of Sects.~\ref{S:PSISOLVES}-\ref{S:C0BOUNDBOOTSTRAP},
	if $\varepsilon$ is sufficiently small, then
	the following estimates hold on $\mathcal{M}_{\Tboot,U_0}:$
	\begin{align} \label{E:RENORMALIZEDTOPORDERTRCHIJUNKTRANSPORTINVERTED}
		& \rgeo^2(t,u) \left|\chifullmodarg{\mathscr{S}^N} \right|(t,u,\vartheta)
			\\
		& \leq 
			(1 + C \varepsilon)
			\sup_{0 \leq s \leq t}
			\left(\frac{\upmu(s,u,\vartheta)}{\upmu(0,u,\vartheta)}\right)^2 \rgeo^2(0,u) \left|\chifullmod^{[N]} \right|(0,u,\vartheta)
			\notag \\
		& \ \ + 	2 (1 + C \varepsilon)
							\int_{s=0}^t 
								\left(
									\frac{\upmu(t,u,\vartheta)}{\upmu(s,u,\vartheta)}
								\right)^2 
								\rgeo^2(s,u)
								\frac{[\Lunit \upmu(s,u,\vartheta)]_-}{\upmu(s,u,\vartheta)}
								\left| \mathfrak{X}^{[N]} \right|(s,u,\vartheta)
							\, ds
					\notag \\
		& \ \ + \frac{1}{2} (1 + C \varepsilon) 
											 \int_{s=0}^t 
													\left(\frac{\upmu(t,u,\vartheta)}{\upmu(s,u,\vartheta)}\right)^2 \rgeo^2(s,u)
													\mytr \upchi(s,u,\vartheta)
													\left| \mathfrak{X}^{[N]} \right|(s,u,\vartheta)
												\, ds
							\notag \\
		& \ \ + 2 (1 + C \varepsilon)
					\int_{s=0}^t 
							\left(
								\frac{\upmu(t,u,\vartheta)}{\upmu(s,u,\vartheta)}
							\right)^2 
							\rgeo^2(s,u)
							\left|
								\hat{\upchi}^{(Small)}
							\right|
							(s,u,\vartheta)
							\left| 
								\upmu \angLie_{\mathscr{S}}^{\leq N} \hat{\upchi}^{(Small)}
							\right|(s,u,\vartheta)
						\, ds
				\notag \\
		& \ \ + (1 + C \varepsilon)
						\int_{s=0}^t 
							\left(
								\frac{\upmu(t,u,\vartheta)}{\upmu(s,u,\vartheta)}
							\right)^2 
							\rgeo^2(s,u)
							\left\lbrace
								\left| \mathfrak{I}^{[N]} \right|(s,u,\vartheta)
								+ \left|\widetilde{\mathfrak{I}}^{N} \right|(s,u,\vartheta)
							\right\rbrace
						\, ds.
						\notag 
\end{align}
Above, 
$\left|\chifullmod^{[N]} \right|$
denotes a term 
that is $\leq \sum_{|\vec{I}| = N} c_{\vec{I}} \left|\chifullmodarg{\mathscr{S}^{\vec{I}}} \right|,$
where the $c_{\vec{I}}$ are non-negative constants verifying 
$\sum_{|\vec{I}| = N} c_{\vec{I}} \leq 1.$
Similarly, 
$\left| \mathfrak{X}^{[N]} \right|$ denotes a term
that is $\leq \sum_{|\vec{I}| = N} c_{\vec{I}} \left| \mathscr{S}^{\vec{I}} \mathfrak{X} \right|,$
where the $c_{\vec{I}}$ are non-negative constants verifying 
$\sum_{|\vec{I}| = N} c_{\vec{I}} \leq 1,$
and $\left| \mathfrak{I}^{[N]} \right|$ denotes a term that is
$\leq \sum_{|\vec{I}| = N} c_{\vec{I}} \left| \inhomleftexparg{\mathscr{S}^{\vec{I}}}{\mathfrak{I}} \right|,$
where the $c_{\vec{I}}$ are non-negative constants verifying 
$\sum_{|\vec{I}| = N} c_{\vec{I}} \leq 1.$
Furthermore,
$\widetilde{\mathfrak{I}}^{N}$ 
is an inhomogeneous term that verifies the
pointwise estimate
\begin{align} \label{E:NEEDTOTAKEMAXANDTHENGRONWALLADDITIONALINHOMOGENEOUSTERM}
\left|\widetilde{\mathfrak{I}}^{N} \right|(t,u,\vartheta)
& \lesssim 
		\frac{\ln(\myexp + t)}{(1 + t)^2}
			\left|
				\fourmyarray[\rgeo \Lunit \mathscr{Z}^{\leq N} \Psi]
					{\Rad \mathscr{Z}^{\leq N} \Psi}
					{\rgeo \angdiff \mathscr{Z}^{\leq N} \Psi} 
					{\mathscr{Z}^{\leq N} \Psi}
			\right|
		+ \frac{\ln(\myexp + t)}{(1 + t)^3}
			\left| 
				\myarray
					[\mathscr{Z}^{\leq N} (\upmu - 1)]
						{\rgeo \sum_{a=1}^3 |\mathscr{Z}^{\leq N} \Lunit_{(Small)}^a|}
				\right|.
\end{align}		
\end{lemma}

\begin{proof}
We view the terms in the transport equation \eqref{E:TOPORDERTRCHIJUNKRENORMALIZEDTRANSPORT}
as functions of $(s,u,\vartheta).$ We multiply both sides of \eqref{E:TOPORDERTRCHIJUNKRENORMALIZEDTRANSPORT}
by the integrating factor
$\iota(s,u,\vartheta)
	:=
	\exp 
	\left(
		\int_{t' = 0}^s
			\left\lbrace
					\mytr \upchi 
				- 2 \frac{\Lunit \upmu}{\upmu}
			\right\rbrace
			(t',u,\vartheta)
		\, dt'
		\right)$
(corresponding to the coefficient of $\chifullmodarg{\mathscr{S}^N}$ on the left-hand side of \eqref{E:TOPORDERTRCHIJUNKRENORMALIZEDTRANSPORT})
and then integrate the resulting transport equation for $\iota \chifullmodarg{\mathscr{S}^N}$
along the integral curves of $\Lunit = \frac{\partial}{\partial s}$
from $s=0$ to $s=t.$
With the help of the estimate \eqref{E:POINTWISEESTIMATEFORTRCHIANTIDERIVATIVE}, 
we see that 
\begin{align}
	\iota(s,u,\vartheta)
	& :=
	\exp 
	\left(
		\int_{t' = 0}^s
			\left\lbrace
					\mytr \upchi 
				- 2 \frac{\Lunit \upmu}{\upmu}
			\right\rbrace
			(t',u,\vartheta)
		\, dt'
	\right)
	= \frac{\upmu^2(0,u,\vartheta)}{\upmu^2(s,u,\vartheta)}
	\exp 
	\left(
		\int_{t' = 0}^s
			\mytr \upchi (t',u,\vartheta)
		\, dt'
	\right)
		\\
	& = (1 + \mathcal{O}(\varepsilon))
			\frac{\rgeo^2(s,u)}{\rgeo^2(0,u)}
			\frac{\upmu^2(0,u,\vartheta)}{\upmu^2(s,u,\vartheta)}.
			\notag
	\end{align}
We also use the inequality $0 < - 2 \frac{\Lunit \upmu}{\upmu} + \frac{1}{2} \mytr \upchi \leq 
2 \frac{[\Lunit \upmu]_-}{\upmu} + \frac{1}{2} \mytr \upchi,$ 
which follows from \eqref{E:RENORMALIZATIONCOEFFICIENTISPOSITIVE}, in order to bound the
magnitude of the coefficient of the factor $\mathscr{S}^N \mathfrak{X}$ on the right-hand side of \eqref{E:TOPORDERTRCHIJUNKRENORMALIZEDTRANSPORT}.
This leads to the pointwise estimate \eqref{E:RENORMALIZEDTOPORDERTRCHIJUNKTRANSPORTINVERTED}
except that on the right-hand side of \eqref{E:RENORMALIZEDTOPORDERTRCHIJUNKTRANSPORTINVERTED},
in place of the terms 
$\chifullmod^{[N]},$
$\mathfrak{X}^{[N]},$
and $\mathfrak{I}^{[N]},$
there respectively appears
$\chifullmodarg{\mathscr{S}^N},$
$\mathscr{S}^N \mathfrak{X},$
and $\chifullmodsourcearg{\mathscr{S}^N},$
and in place of the
integrand term
$\left| \widetilde{\mathfrak{I}}^N \right|,$
there appears the commutator term $\left| \upmu [\Lunit, \mathscr{S}^N] \mytr \upchi^{(Small)} \right|$
from the right-hand side of \eqref{E:TOPORDERTRCHIJUNKRENORMALIZEDTRANSPORT}. 

We now claim that the following estimate holds:
\begin{align} \label{E:TOPORDERRENORMALIZEDTRCHIJUNKANNOYINGCOMMUTATORTERMESTIMATE}
	\left|
		\upmu [\Lunit, \mathscr{S}^N] \mytr \upchi^{(Small)}
	\right|
	(t,u,\vartheta)
	& \leq 
		C \varepsilon 
		\frac{\ln(\myexp + t)}{(1 + t)^2} \sum_{|\vec{I}|=N} \chifullmodarg{\mathscr{S}^{\vec{I}}}
		+ C
			\frac{\ln(\myexp + t)}{(1 + t)^2}
			\left|
				\fourmyarray[\rgeo \Lunit \mathscr{Z}^{\leq N} \Psi]
					{\Rad \mathscr{Z}^{\leq N} \Psi}
					{\rgeo \angdiff \mathscr{Z}^{\leq N} \Psi} 
					{\mathscr{Z}^{\leq N} \Psi}
			\right|
				\\
	& \ \ 
			+ C
			\frac{\ln(\myexp + t)}{(1 + t)^3}
			\left| 
				\myarray
					[\mathscr{Z}^{\leq N} (\upmu - 1)]
						{\rgeo \sum_{a=1}^3 |\mathscr{Z}^{\leq N} \Lunit_{(Small)}^a|}
				\right|.
				\notag 
\end{align}
Once we have shown \eqref{E:TOPORDERRENORMALIZEDTRCHIJUNKANNOYINGCOMMUTATORTERMESTIMATE},
in order to deduce the desired estimate \eqref{E:RENORMALIZEDTOPORDERTRCHIJUNKTRANSPORTINVERTED}, 
we first sum the integral inequalities verified by each fixed $\left|\chifullmodarg{\mathscr{S}^N}\right|$
in order to deduce an integral inequality for 
$\sum_{|\vec{I}|=N} \left|\chifullmodarg{\mathscr{S}^{\vec{I}}} \right|.$
We then use Gronwall's inequality to derive a pointwise estimate for the quantity
$\sum_{|\vec{I}|=N} \left| \chifullmodarg{\mathscr{S}^{\vec{I}}} \right|$
in which the quantity $\sum_{|\vec{I}|=N} \left|\chifullmodarg{\mathscr{S}^{\vec{I}}} \right|$
no longer appears on the right-hand side. 
In deriving the Gronwall estimate, 
we use the time-integrability of the factor $\varepsilon \ln(\myexp + t) (1+t)^{-2}$
that multiplies $\sum_{|\vec{I}|=N} \left| \chifullmodarg{\mathscr{S}^{\vec{I}}} \right|$ 
on the right-hand side of 
\eqref{E:TOPORDERRENORMALIZEDTRCHIJUNKANNOYINGCOMMUTATORTERMESTIMATE}.
We then revisit the original 
integral inequality verified by 
the fixed $\left|\chifullmodarg{\mathscr{S}^N} \right|$
and insert the Gronwall estimate for
$\sum_{|\vec{I}|=N} \left| \chifullmodarg{\mathscr{S}^{\vec{I}}} \right|$
into the first term on the right-hand side of
\eqref{E:TOPORDERRENORMALIZEDTRCHIJUNKANNOYINGCOMMUTATORTERMESTIMATE}.
The remaining terms on the right-hand side of
\eqref{E:TOPORDERRENORMALIZEDTRCHIJUNKANNOYINGCOMMUTATORTERMESTIMATE}
result in the presence of the term
$\left|\widetilde{\mathfrak{I}}^{N} \right|$
on the right-hand side of \eqref{E:RENORMALIZEDTOPORDERTRCHIJUNKTRANSPORTINVERTED}.


It remains for us to prove \eqref{E:TOPORDERRENORMALIZEDTRCHIJUNKANNOYINGCOMMUTATORTERMESTIMATE}.
As a first step, we use the pointwise estimates 
$\upmu \lesssim \ln(\myexp + t)$
and $\left|\mathscr{Z}^{\leq 11} \mytr \upchi^{(Small)} \right| \lesssim \varepsilon \ln(\myexp + t)(1 + t)^{-2}$
(which follow from \eqref{E:C0BOUNDCRUCIALEIKONALFUNCTIONQUANTITIES}),
the commutator estimate
\eqref{E:LSNCOMMUTATORACTINGONFUNCTIONSPOINTWISE} with 
$\mytr \upchi^{(Small)}$ in the role of $f,$
the bound \eqref{E:LDERIVATIVECRUCICALTRANSPORTINTEQUALITIES}
for 
$\left|\Lunit(\rgeo^2 \mathscr{Z}^{\leq N-1} \mytr \upchi^{(Small)})\right|,$
the identity $\Lunit \rgeo = 1,$
the bound \eqref{E:POINTWISEESTIMATESFORCHIJUNKINTERMSOFOTHERVARIABLES} for
$\left|\mathscr{Z}^{\leq N-1} \mytr \upchi^{(Small)} \right|$ in terms of other variables,
Cor.~\ref{C:SQRTEPSILONREPLCEDWITHCEPSILON}, 
and inequality \eqref{E:FUNCTIONPOINTWISEANGDINTERMSOFANGLIEO}
to deduce that
\begin{align} \label{E:FIRSTTOPORDERRENORMALIZEDTRCHIJUNKANNOYINGCOMMUTATORTERMESTIMATE}
	\left|
		\upmu [\Lunit, \mathscr{S}^N] \mytr \upchi^{(Small)}
	\right|
	& \lesssim 
			\varepsilon
			\frac{\ln(\myexp + t)}{(1 + t)^2}
			\sum_{l=1}^3
			\left| 
				\myarray
					[\upmu \Rad \mathscr{Z}^{\leq N-1} \mytr \upchi^{(Small)}]
					{\upmu \Rot_{(l)} \mathscr{Z}^{\leq N-1} \mytr \upchi^{(Small)}}
			\right|
		+ 
			\varepsilon \frac{\ln^2(\myexp + t)}{(1 + t)^3}
			\left| 
				\fourmyarray[ \rgeo \Lunit \mathscr{Z}^{\leq N} \Psi]
					{\Rad \mathscr{Z}^{\leq N} \Psi}
					{\rgeo \angdiff \mathscr{Z}^{\leq N} \Psi}
					{\mathscr{Z}^{\leq N-} \Psi}
			\right|
				\\
	& \ \ 
			+ 
			\varepsilon
			\frac{\ln^2(\myexp + t)}{(1 + t)^4}
			\left|
				\myarray[\mathscr{Z}^{\leq N} (\upmu - 1)]
					{\sum_{a=1}^3 \rgeo |\mathscr{Z}^{\leq N} \Lunit_{(Small)}^a|} 
			\right|.
			\notag	
\end{align}
	Next, we consider the first array on the 
	right-hand side of \eqref{E:FIRSTTOPORDERRENORMALIZEDTRCHIJUNKANNOYINGCOMMUTATORTERMESTIMATE}. 
	We repeatedly use inequalities
	\eqref{E:RGEOLORRADZNCOMMUTATORACTINGONFUNCTIONSPOINTWISE}
	and \eqref{E:ROTZNCOMMUTATORACTINGONFUNCTIONSPOINTWISE}
	with $f = \mytr \upchi^{(Small)},$
	\eqref{E:C0BOUNDCRUCIALEIKONALFUNCTIONQUANTITIES}, 
	and \eqref{E:BOUNDINGFRAMEDERIVATIVESINTERMSOFZNF}
	to commute any factor $Z := \rgeo \Lunit$ in 
	$\Rad \mathscr{Z}^{\leq N-1}$ or $\Rot \mathscr{Z}^{\leq N-1}$ 
	all the way in front.
	We then use the estimates described just above
	inequality \eqref{E:FIRSTTOPORDERRENORMALIZEDTRCHIJUNKANNOYINGCOMMUTATORTERMESTIMATE}
	to bound all derivatives of $\mytr \upchi^{(Small)}$ except the pure spatial top-order ones
	$\mathscr{S}^N$ in terms of other variables, thus deducing that
\begin{align} \label{E:ONLYPURETOPORDERSPATIALDERIVATIVESREMAIN}
			\varepsilon
			\frac{\ln(\myexp + t)}{(1 + t)^2}
			\sum_{l=1}^3
			\left| 
				\myarray
					[\upmu \Rad \mathscr{Z}^{\leq N-1} \mytr \upchi^{(Small)}]
					{\upmu \Rot_{(l)} \mathscr{Z}^{\leq N-1} \mytr \upchi^{(Small)}}
			\right|
			& \lesssim 
			\varepsilon
			\frac{\ln(\myexp + t)}{(1 + t)^2}
			\left|
				\upmu \mathscr{S}^N \mytr \upchi^{(Small)}
			\right|
				\\
		& \ \ 
			+
			\varepsilon \frac{\ln^2(\myexp + t)}{(1 + t)^3}
			\left| 
				\fourmyarray[ \rgeo \Lunit \mathscr{Z}^{\leq N} \Psi]
					{\Rad \mathscr{Z}^{\leq N} \Psi}
					{\rgeo \angdiff \mathscr{Z}^{\leq N} \Psi}
					{\mathscr{Z}^{\leq N} \Psi}
			\right|
				\notag \\
	& \ \ 
			+ 
			\varepsilon
			\frac{\ln^2(\myexp + t)}{(1 + t)^4}
			\left|
				\myarray[\mathscr{Z}^{\leq N} (\upmu - 1)]
					{\sum_{a=1}^3 \rgeo |\mathscr{Z}^{\leq N} \Lunit_{(Small)}^a|} 
			\right|.
			\notag
\end{align}
Using the identity
$\upmu \mathscr{S}^{\vec{I}} \mytr \upchi^{(Small)} = \chifullmodarg{\mathscr{S}^{\vec{I}}} - \mathscr{S}^{\vec{I}} \mathfrak{X},$
and the estimate $\upmu \lesssim \ln(\myexp + t),$
we bound the first term on the right-hand side of 
\eqref{E:ONLYPURETOPORDERSPATIALDERIVATIVESREMAIN} as follows:
\begin{align} \label{E:TOPORDERPOURESPATIALDERIVATIVEBOUND}
	\left|
		\upmu \mathscr{S}^N \mytr \upchi^{(Small)}
	\right|
	& \lesssim 
		\sum_{|\vec{I}|=N} \chifullmodarg{\mathscr{S}^{\vec{I}}}
		+ \sum_{|\vec{I}|=N} \left| \mathscr{S}^{\vec{I}} \mathfrak{X} \right|.
	\end{align}
Combining
\eqref{E:FIRSTTOPORDERRENORMALIZEDTRCHIJUNKANNOYINGCOMMUTATORTERMESTIMATE},
\eqref{E:ONLYPURETOPORDERSPATIALDERIVATIVESREMAIN},
and
\eqref{E:TOPORDERPOURESPATIALDERIVATIVEBOUND},
and using \eqref{E:POINTWISEESTIMATEZNAPPLIEDTOLOWESTTRCHIJUNKDISCREPANCYTERM}
to bound the terms $\left| \mathscr{S}^{\vec{I}} \mathfrak{X} \right|$  
on the right-hand side of \eqref{E:TOPORDERPOURESPATIALDERIVATIVEBOUND},
we conclude the desired bound \eqref{E:TOPORDERRENORMALIZEDTRCHIJUNKANNOYINGCOMMUTATORTERMESTIMATE}.
We have thus proved Lemma~\ref{L:RENORMALIZEDTOPORDERTRCHIJUNKTRANSPORTINVERTED}.

\end{proof}

\section{Pointwise estimates for the difficult error integrands requiring full modification}
\label{S:FULLYRENORMALIZEDPOINTWISEESTIMATESFORDIFFICULTENERGYERRORINTEGRAND}

In the next proposition, we derive precise pointwise estimates for the 
difficult error integrand products
$(\Rad \Psi) \Rot \mathscr{S}^{N-1} \mytr \upchi^{(Small)}$
and
$(\Rad \Psi) \angLap \mathscr{S}^{N-1} \upmu$
from Prop.~\ref{P:IDOFKEYDIFFICULTENREGYERRORTERMS},
where $\mathscr{S}^{N-1}$ is an $(N-1)^{st}$ order pure spatial commutation vectorfield operator.
More precisely, 
the last inequality in \eqref{E:TOPORDERDERIVATIVESOFANGDSQUAREDUPMUINTERMSOFCONTROLLABLE}
allows us to estimate
$(\Rad \Psi) \Rad \mathscr{S}^{N-1} \mytr \upchi^{(Small)}$
instead of
$(\Rad \Psi) \angLap \mathscr{S}^{N-1} \upmu.$
We must use the fully modified quantities 
defined in Sect.~\ref{S:FULLRENORMALIZATIONFORTRCHIJUNK}
to derive these estimates,
which play an important role in our derivation of 
a priori estimates for the top-order $L^2$ quantity $\totzeromax{\leq 24}.$
\textbf{We stress that the large $\upmu_{\star}^{-1}-$degeneracy of our a priori estimates 
for $\totzeromax{\leq 24}$ 
(see inequality \eqref{E:Q0TOPORDERKEYBOUND})
stems in part from the ``boxed'' constants appearing on the right-hand sides
of the estimates of Prop.~\ref{P:MAINCOMMUTEDWAVEEQNINHOMOGENEOUSTERMPOINTWISEESTIMATES}.}	

\begin{proposition}[\textbf{Pointwise estimates for} 
$(\Rad \Psi) \Rot \mathscr{S}^{N-1} \mytr \upchi^{(Small)})$
\textbf{and}
$(\Rad \Psi) \Rad \mathscr{S}^{N-1} \mytr \upchi^{(Small)})$
]
\label{P:MAINCOMMUTEDWAVEEQNINHOMOGENEOUSTERMPOINTWISEESTIMATES}
	Let $1 \leq N \leq 24$ be an integer 
	and let $\mathscr{S}^{N-1}$ be an $(N-1)^{st}$ order pure spatial commutation vectorfield operator
	(see definition \eqref{E:DEFSETOFSPATIALCOMMUTATORVECTORFIELDS}).
	Under the small-data and bootstrap assumptions 
	of Sects.~\ref{S:PSISOLVES}-\ref{S:C0BOUNDBOOTSTRAP},
	if $\varepsilon$ is sufficiently small, then
	the following pointwise estimates hold on $\mathcal{M}_{\Tboot,U_0},$
	where 
	$\left|\chifullmod^{[N]} \right|$
	denotes a term 
	that is $\leq \sum_{|\vec{I}| = N} c_{\vec{I}} \left|\chifullmodarg{\mathscr{S}^{\vec{I}}} \right|,$
	where $\chifullmodarg{\mathscr{S}^{\vec{I}}}$ is defined by \eqref{E:TRANSPORTRENORMALIZEDTRCHIJUNK} 
	and the $c_{\vec{I}}$ are non-negative constants verifying 
	$\sum_{|\vec{I}| = N} c_{\vec{I}} \leq 1:$
	\begin{align} 
		& 
		\left| (\Rad \Psi) \Rot \mathscr{S}^{N-1} \mytr \upchi^{(Small)} \right|,
			\, \left| (\Rad \Psi) \Rad \mathscr{S}^{N-1} \mytr \upchi^{(Small)} \right|
			\label{E:TOPORDERTRCHIJUNKENERGYERRORTERMKEYPOINTWISEESTIMATE}  
					\\
		& \leq \boxed{2} \left \|
								\frac{\Lunit \upmu}{\upmu}
							\right \|_{C^0(\Sigma_t^u)}
						\left| \Rad \mathscr{S}^N \Psi \right|
				\notag	\\
		&  \ \ + 
						\boxed{4} (1 + C \sqrt{\varepsilon})\frac{1}{\rgeo(t,u)}
							\left \|
								\frac{\Lunit \upmu}{\upmu}
							\right \|_{C^0(\Sigma_t^u)}
					\int_{t'=0}^t 
					\frac{
						\left \|
							[\Lunit \upmu]_-
						\right \|_{C^0(\Sigma_{t'}^u)}}
							{\upmu_{\star}(t',u)}
					\rgeo(t',u)
					\left| \Rad \mathscr{S}^N \Psi \right|(t',u,\vartheta)
				\, dt'
				\notag \\
			& \ \ 
					+ 
						\boxed{2}(1 + C \varepsilon)
						\frac{1}{\rgeo^2(t,u)}
						\left \|
								\frac{\Lunit \upmu}{\upmu}
						\right \|_{C^0(\Sigma_t^u)}
						\int_{t'=0}^t 
							\left\|
								\left(\frac{\upmu(t,\cdot)}{\upmu}\right)^2
							\right\|_{C^0(\Sigma_{t'}^u)}
							\rgeo(t',u)
							\left| \Rad \mathscr{S}^N \Psi \right|(t',u,\vartheta)
						\, dt'
					\notag \\
		& \ \	+ 
			C \varepsilon 
			\upmu_{\star}^{-1}(t,u)
			\frac{\ln^2(\myexp + t)}{(1 + t)^3}
			|\chifullmod^{[N]}|(0,u,\vartheta)
				\notag \\
		& \ \ + C \varepsilon 
							\upmu_{\star}^{-1/2}(t,u)
							\frac{1}{(1 + t)^2}
							\left|
								\threemyarray
									[\rgeo \sqrt{\upmu} \left\lbrace \Lunit + \frac{1}{2} \mytr \upchi \right\rbrace \mathscr{S}^{\leq N} \Psi]
									{\Rad \mathscr{S}^{\leq N} \Psi}
									{\rgeo \sqrt{\upmu} \angdiff \mathscr{S}^{\leq N} \Psi}
							\right|
					 \notag \\
		& \ \ + 
					C \varepsilon 
					\upmu_{\star}^{-1}(t,u)
					\frac{1}{(1 + t)}
					\left|
						\myarray
							[\Rad \mathscr{Z}^{\leq N-1} \Psi]
							{\mathscr{Z}^{\leq N-1} \Psi}
					\right|
					+ 
					C \varepsilon 
							\upmu_{\star}^{-3/2}(t,u)
							\frac{\ln(\myexp + t)}{(1 + t)^2}
							\left|
								\myarray
									[\rgeo \sqrt{\upmu} \left\lbrace \Lunit + \frac{1}{2} \mytr \upchi \right\rbrace \mathscr{Z}^{\leq N-1} \Psi]
										{ \rgeo \sqrt{\upmu} \angdiff \mathscr{Z}^{\leq N-1} \Psi}
								\right|
			\notag \\
	& \ \ 
		+ C \varepsilon
			\upmu_{\star}^{-1}(t,u)
			\frac{1}{(1 + t)^3}
			\left| 
				\myarray
					[\mathscr{Z}^{\leq N} (\upmu - 1)]
					{\rgeo \sum_{a=1}^3 |\mathscr{Z}^{\leq N} \Lunit_{(Small)}^a|}
			\right|
		\notag \\
	& \ \ + C \varepsilon 
						\upmu_{\star}^{-1}(t,u)
						\frac{\ln^3(\myexp + t)}{(1 + t)^3}
						\int_{t'=0}^t 
							\frac{1}{\upmu_{\star}(t',u)}
							\left| 
								\threemyarray
									[ \rgeo \sqrt{\upmu} \left\lbrace \Lunit + \frac{1}{2} \mytr \upchi \right\rbrace \mathscr{Z}^{\leq N} \Psi]
									{\Rad \mathscr{Z}^{\leq N} \Psi}
									{ \rgeo \sqrt{\upmu} \angdiff \mathscr{Z}^{\leq N} \Psi}
							 \right|
							(t',u,\vartheta)
						\, dt'
						\notag \\
	& \ \ + C \varepsilon 
						\upmu_{\star}^{-1}(t,u)
						\frac{\ln^2(\myexp + t)}{(1 + t)^3}
						\int_{t'=0}^t 
							\rgeo(t',u)
							\frac{1}{\upmu_{\star}(t',u)}
							\left| 
								\myarray
									[\Rad \mathscr{Z}^{\leq N-1} \Psi]
									{\mathscr{Z}^{\leq N-1} \Psi}
							\right|
							(t',u,\vartheta)
						\, dt'
						\notag \\
	& \ \ + C \varepsilon 
						\upmu_{\star}^{-1}(t,u)
						\frac{\ln^3(\myexp + t)}{(1 + t)^3}
						\int_{t'=0}^t 
							\frac{1}{\upmu_{\star}^{3/2}(t',u)}
							\left| 
								\myarray
									[\rgeo \sqrt{\upmu} \left\lbrace \Lunit + \frac{1}{2} \mytr \upchi \right\rbrace \mathscr{Z}^{\leq N-1} \Psi]
									{\rgeo \sqrt{\upmu} \angdiff \mathscr{Z}^{\leq N-1} \Psi}
							\right|
							(t',u,\vartheta)
						\, dt'
						\notag \\
		& \ \ + C \varepsilon 
						\upmu_{\star}^{-1}(t,u)
						\frac{\ln^3(\myexp + t)}{(1 + t)^3}
						\int_{t'=0}^t 
							\left| 
								\myarray
									[\upmu \hat{\angD}^2 \mathscr{S}^{\leq N-1} \upmu]
									{\rgeo \upmu \angD \angLie_{\mathscr{S}}^{\leq N-1} \hat{\upchi}^{(Small)}}
							\right|
							(t',u,\vartheta)
							\, dt'
					\notag \\
		& \ \ + C \varepsilon
					\upmu_{\star}^{-1}(t,u)
					\frac{\ln^2(\myexp + t)}{(1 + t)^3}	
					\int_{t'=0}^t 
						\frac{1}{(1 + t')}
						\frac{1}{\upmu_{\star}(t',u)}
						\left| 
							\myarray
								[\mathscr{Z}^{\leq N} (\upmu - 1)]
								{\rgeo \sum_{a=1}^3 |\mathscr{Z}^{\leq N} \Lunit_{(Small)}^a|}
						\right|
						(t',u,\vartheta)
					\, dt'.
					\notag 
	\end{align}
	
	
\end{proposition}

\begin{proof}
For definiteness, we prove inequality \eqref{E:TOPORDERTRCHIJUNKENERGYERRORTERMKEYPOINTWISEESTIMATE} only for the term
$(\Rad \Psi) \Rot \mathscr{S}^{N-1} \mytr \upchi^{(Small)};$
the proof for the term $(\Rad \Psi) \Rad \mathscr{S}^{N-1} \mytr \upchi^{(Small)}$
is identical. We begin by setting $\mathscr{S}^N := \Rot \mathscr{S}^{N-1}$
and using Def.~\ref{D:TRANSPORTRENORMALIZEDTRCHIJUNK}
(see also \eqref{E:LOWESTORDERTRANSPORTRENORMALIZEDTRCHIJUNKDISCREPANCY})
to deduce that
\begin{align} \label{E:PRELIMINARYREWRITINGOFTOPORDERRENORMALIZEDTHRCHIJUNKTERMINTERMSOFSTUFFICANESTIMATE}
	&
	\left|
		\upmu \Rot \mathscr{S}^{N-1} \mytr \upchi^{(Small)}
		- \chifullmodarg{\mathscr{S}^N}
	\right|
		\\
	& \leq 
		\left| 	
			G_{\Lunit \Lunit} \Rad \mathscr{S}^N \Psi
		\right|
		+
		\frac{1}{2} 
		\left|
			\upmu \angGmixedarg{A}{A} \Lunit \mathscr{S}^N \Psi
		\right|
		+ 
		\frac{1}{2} 
		\left|
			\upmu G_{\Lunit \Lunit} \Lunit \mathscr{S}^N \Psi
		\right|
	  +
		\left|
			\upmu \angGmixedarg{\Lunit}{A} \angdiffarg{A} \mathscr{S}^N \Psi
		\right| 
			\notag \\
		& \ \ 
		+
		\left|
			\mathscr{S}^N \mathfrak{X}
			- 
				\left\lbrace
					- G_{\Lunit \Lunit} \Rad \mathscr{S}^N \Psi
					- \frac{1}{2} \upmu \angGmixedarg{A}{A} \Lunit \mathscr{S}^N  \Psi
					- \frac{1}{2} \upmu G_{\Lunit \Lunit} \Lunit \mathscr{S}^N \Psi
					+ \upmu \angGmixedarg{\Lunit}{A} \angdiffarg{A} \mathscr{S}^N \Psi
				\right\rbrace
			\right|.
			\notag
\end{align}
From \eqref{E:PRELIMINARYREWRITINGOFTOPORDERRENORMALIZEDTHRCHIJUNKTERMINTERMSOFSTUFFICANESTIMATE},
the bootstrap assumption $\|\Rad \Psi\|_{C^0(\Sigma_t^u)} \leq \varepsilon (1 + t)^{-1}$
(that is, \eqref{E:PSIFUNDAMENTALC0BOUNDBOOTSTRAP}),
inequalities \eqref{E:LOWERORDERC0BOUNDLIEDERIVATIVESOFGRAME}
and \eqref{E:POINTWISEESTIMATECOMMUTATORINVOLVINGZNAPPLIEDTOLOWESTTRCHIJUNKDISCREPANCYTERM},
and the bound $\left|\upmu - 1 \right| \lesssim \varepsilon \ln(\myexp + t)$
(that is, \eqref{E:C0BOUNDCRUCIALEIKONALFUNCTIONQUANTITIES}),
we deduce that
\begin{align} \label{E:FIRSTSPLITTINGDIFFICULTOPORDERTRCHIJUNKENERGYTERM}
	&
	\left| 
		(\Rad \Psi) \Rot \mathscr{S}^{N-1} \mytr \upchi^{(Small)}
	\right|
		\\
	&  \leq \frac{1}{\rgeo^2} 
		\left|\frac{\Rad \Psi}{\upmu} 
			(\rgeo^2 \chifullmodarg{\mathscr{S}^N})
		\right|
		+ \left| 	
				G_{\Lunit \Lunit} \frac{\Rad \Psi}{\upmu} \Rad \mathscr{S}^N \Psi
			\right|
		+ C \varepsilon \frac{1}{(1 + t)^2}
				\left|
					\myarray
						[\rgeo \Lunit \mathscr{S}^N \Psi]
						{\rgeo \angdiff \mathscr{S}^N \Psi}		
				\right|
			\notag	\\
	& \ \
		+ C \varepsilon \frac{\ln(\myexp + t)}{(1 + t)^2}
		\frac{1}{\upmu}
		 \left|
				\threemyarray[\rgeo \Lunit \mathscr{Z}^{\leq N-1} \Psi]
					{\rgeo \angdiff \mathscr{Z}^{\leq N-1} \Psi} 
					{\mathscr{Z}^{\leq N-1} \Psi}
			\right|
			+ C \varepsilon 
					\frac{1}{(1 + t)}
					\frac{1}{\upmu}
				\left| 
					\Rad \mathscr{Z}^{\leq N-1} \Psi
				\right|
				\notag \\
		& \ \ 
		+ C \varepsilon
			\frac{1}{(1 + t)^3}
			\frac{1}{\upmu}
			\left| 
				\myarray
					[\mathscr{Z}^{\leq N} (\upmu - 1)]
						{\rgeo \sum_{a=1}^3 |\mathscr{Z}^{\leq N} \Lunit_{(Small)}^a|}
			\right|.
			\notag
\end{align}
	We now bound the right-hand side of \eqref{E:FIRSTSPLITTINGDIFFICULTOPORDERTRCHIJUNKENERGYTERM}
	by the right-hand side of \eqref{E:TOPORDERTRCHIJUNKENERGYERRORTERMKEYPOINTWISEESTIMATE}
	by using a separate argument for each term. 
	Throughout the remainder of this proof, we use
	the definition of $\upmu_{\star}$ and in particular the estimate $\upmu_{\star} \leq 1$
	without mentioning it each time.
	We also silently use the simple inequality
	$|\Lunit f| 
		\leq 
		\left|
			\left\lbrace
				\Lunit f + \frac{1}{2} \mytr \upchi f 
			\right\rbrace
		\right|
	+ C (1 + t)^{-1} |f|,$
which follows easily from the estimate 
$|\rgeo \mytr \upchi| \lesssim 1$
(that is, \eqref{E:CRUDELOWERORDERC0BOUNDDERIVATIVESOFANGULARDEFORMATIONTENSORS}).
	The first term on the right-hand side of \eqref{E:FIRSTSPLITTINGDIFFICULTOPORDERTRCHIJUNKENERGYTERM}
	is the most difficult, and we handle it at the end of the proof.
	To bound the second term on the right-hand side of \eqref{E:FIRSTSPLITTINGDIFFICULTOPORDERTRCHIJUNKENERGYTERM},
	we first use the transport equation \eqref{E:UPMUFIRSTTRANSPORT} for $\Lunit \upmu,$ 
	the bootstrap assumption $ \| \rgeo \Lunit \Psi \|_{C^0(\Sigma_t^u)} \leq \varepsilon (1 + t)^{-1}$
	(that is, \eqref{E:PSIFUNDAMENTALC0BOUNDBOOTSTRAP}),
	and inequality \eqref{E:LOWERORDERC0BOUNDLIEDERIVATIVESOFGRAME}
	to deduce that
	\begin{align} \label{E:KEYGLLRADPSIOVERUPMUTIMESRADZNPSIBOUND}
		\left| 
			G_{\Lunit \Lunit} \frac{\Rad \Psi}{\upmu} \Rad \mathscr{S}^N \Psi 
		\right|
		& \leq 
			2
			\left\| \frac{\Lunit \upmu}{\upmu} \right\|_{C^0(\Sigma_t^u)}
			\left|
				\Rad \mathscr{S}^N \Psi 
			\right|
			+ C \varepsilon
				\frac{1}{(1 + t)^2}
				\left| 
					\Rad \mathscr{S}^N \Psi
				\right|.
	\end{align}
	The first term on the right-hand side of \eqref{E:KEYGLLRADPSIOVERUPMUTIMESRADZNPSIBOUND}
	accounts for the first term on the right-hand side of \eqref{E:TOPORDERTRCHIJUNKENERGYERRORTERMKEYPOINTWISEESTIMATE},
	and the second term is also easily seen to be $\leq$ the right-hand side of 
	\eqref{E:TOPORDERTRCHIJUNKENERGYERRORTERMKEYPOINTWISEESTIMATE}.
	The last term on the first line of the right-hand side of \eqref{E:FIRSTSPLITTINGDIFFICULTOPORDERTRCHIJUNKENERGYTERM}
	and the terms on the second and third lines of the right-hand side of \eqref{E:FIRSTSPLITTINGDIFFICULTOPORDERTRCHIJUNKENERGYTERM}
	are also easily seen to be $\leq$ the right-hand side of \eqref{E:TOPORDERTRCHIJUNKENERGYERRORTERMKEYPOINTWISEESTIMATE} as desired.
	
	\noindent{\textbf{Bounding the difficult term piece by piece:}}
	To bound the (difficult) first term on the right-hand side of \eqref{E:FIRSTSPLITTINGDIFFICULTOPORDERTRCHIJUNKENERGYTERM}
	by the right-hand side of \eqref{E:TOPORDERTRCHIJUNKENERGYERRORTERMKEYPOINTWISEESTIMATE},
	we use the preliminary pointwise estimates provided by Lemma~\ref{L:RENORMALIZEDTOPORDERTRCHIJUNKTRANSPORTINVERTED},
	which provide a pointwise bound for $\rgeo^2 \chifullmodarg{\mathscr{S}^N}.$
	Specifically, we see that it suffices to bound $\frac{1}{\rgeo^2} \frac{\Rad \Psi}{\upmu}$
	times the right-hand side of \eqref{E:RENORMALIZEDTOPORDERTRCHIJUNKTRANSPORTINVERTED};
	we derive the desired bounds by arguing term by term.
	To bound $\frac{1}{\rgeo^2} \frac{\Rad \Psi}{\upmu}$ 
	times the first term on the right-hand side of \eqref{E:RENORMALIZEDTOPORDERTRCHIJUNKTRANSPORTINVERTED},
	we use the bootstrap assumption $\|\Rad \Psi\|_{C^0(\Sigma_t^u)} \leq \varepsilon (1 + t)^{-1}$
	(that is, \eqref{E:PSIFUNDAMENTALC0BOUNDBOOTSTRAP})
	and the aforementioned bound $\left|\upmu - 1 \right| \lesssim \varepsilon \ln(\myexp + t)$
	to deduce that
	\begin{align} \label{E:DATATERMDIFFICULTMULTPOINTWISEENERGYERRORPRODUCT}
		&
		\left|
			\frac{1}{\rgeo^2(t,u)} \frac{\Rad \Psi(t,u,\vartheta)}{\upmu(t,u,\vartheta)}
			\max_{0 \leq s \leq t}
			\left(\frac{\upmu(s,u,\vartheta)}{\upmu(0,u,\vartheta)}\right)^2 \rgeo^2(0,u) 
		\right|
		|\chifullmod^{[N]}|(0,u,\vartheta)
			\\
		& \leq
			C
			\varepsilon \frac{1}{\upmu_{\star}(t,u)} 
			\frac{\ln^2(\myexp + t)}{(1 + t)^3}
			|\chifullmod^{[N]}|(0,u,\vartheta),
			\notag
	\end{align}
	where $|\chifullmod^{[N]}|$
	is defined just below \eqref{E:RENORMALIZEDTOPORDERTRCHIJUNKTRANSPORTINVERTED}.
	Clearly the right-hand side of \eqref{E:DATATERMDIFFICULTMULTPOINTWISEENERGYERRORPRODUCT} 
	is $\leq$ the right-hand side of \eqref{E:TOPORDERTRCHIJUNKENERGYERRORTERMKEYPOINTWISEESTIMATE} as desired.

\noindent{\textbf{Bound for the term arising from the first time integral on the right-hand side of} 
\eqref{E:RENORMALIZEDTOPORDERTRCHIJUNKTRANSPORTINVERTED}}:	
	We now address the first time integral on the right-hand side of \eqref{E:RENORMALIZEDTOPORDERTRCHIJUNKTRANSPORTINVERTED}.
	To this end, we first note that by the definition 
	of $\left|\mathfrak{X}^{[N]} \right|,$ 
	which is defined just below \eqref{E:RENORMALIZEDTOPORDERTRCHIJUNKTRANSPORTINVERTED},
	the estimate \eqref{E:FUNCTIONAVOIDINGCOMMUTING} (with $\Psi$ in the role of $f$),
	the estimate \eqref{E:LOWERORDERC0BOUNDLIEDERIVATIVESOFGRAME},
	and inequality \eqref{E:POINTWISEESTIMATECOMMUTATORINVOLVINGZNAPPLIEDTOLOWESTTRCHIJUNKDISCREPANCYTERM}, 
	we have
	\begin{align} \label{E:MAXEDRENORMALIZEDTRCHIJUNKDISCREPANCYPOINTWISE}
		\left| \mathfrak{X}^{[N]} \right|
		& \leq |G_{\Lunit \Lunit} \Rad \mathscr{S}^N \Psi|
			+ C 
			\frac{1}{1+t}
			 \upmu 
				\left|
				\myarray
					[\rgeo \left\lbrace \Lunit + \frac{1}{2} \mytr \upchi \right\rbrace \mathscr{S}^{\leq N} \Psi]
					{\rgeo \angdiff \mathscr{S}^{\leq N} \Psi}
				\right|
				\\
	& \ \ 
		+ C
			\frac{\ln(\myexp + t)}{1 + t}
			\left|
				\threemyarray[\rgeo \left\lbrace \Lunit + \frac{1}{2} \mytr \upchi \right\rbrace \mathscr{Z}^{\leq N-1} \Psi]
					{\rgeo \angdiff \mathscr{Z}^{\leq N-1} \Psi} 
					{\mathscr{Z}^{\leq N-1} \Psi}
			\right|
			+ C
				\left|
					\Rad \mathscr{Z}^{\leq N-1} \Psi
				\right|		
			\notag	\\
		& \ \
			+
			C
			\varepsilon
			\frac{1}{(1 + t)^2}
			\left| 
				\myarray
					[\mathscr{Z}^{\leq N} (\upmu - 1)]
						{\rgeo \sum_{a=1}^3 |\mathscr{Z}^{\leq N} \Lunit_{(Small)}^a|}
			\right|.
			\notag
	\end{align}
	Next, we pointwise bound the product of the first and third integrand terms 
	in the first time integral on the right-hand side of 
	\eqref{E:RENORMALIZEDTOPORDERTRCHIJUNKTRANSPORTINVERTED}
	by using Def.~\ref{D:REGIONSOFDISTINCTUPMUBEHAVIOR}
	and the estimates
	\eqref{E:GAMMACONSTANTVERYSMALL},
	\eqref{E:SHARPLOCALIZEDMUCANTGROWTOOFAST},
	\eqref{E:KEYMUNOTDECAYINGMINUSPARTLMUOVERMUBOUND},
	and
	\eqref{E:LOCALIZEDMUMUSTSHRINK}
	to deduce that for $0 \leq s \leq t,$ we have
	\begin{align} \label{E:ANNOYINGINTEGRATINGFACTORBOUND}
		\left(
			\frac{\upmu(t,u,\vartheta)}{\upmu(s,u,\vartheta)}
		\right)^2 
		\frac{[\Lunit \upmu(s,u,\vartheta)]_-}{\upmu(s,u,\vartheta)}
		& \leq (1 + C \sqrt{\varepsilon}) 
				\frac{[\Lunit \upmu(s,u,\vartheta)]_-}{\upmu(s,u,\vartheta)}
			+ C  \varepsilon \frac{\ln^3(\myexp + s)}{(1 + s)^2}.
	\end{align}
	We now use 
	\eqref{E:MAXEDRENORMALIZEDTRCHIJUNKDISCREPANCYPOINTWISE},
	\eqref{E:ANNOYINGINTEGRATINGFACTORBOUND}, 
	\eqref{E:LOWERORDERC0BOUNDLIEDERIVATIVESOFGRAME},
	the estimate
	$|[\Lunit \upmu]_-|(s,u,\vartheta) \leq |\Lunit \upmu|(s,u,\vartheta) \lesssim \varepsilon (1 + s)^{-1},$
	(which follows from \eqref{E:C0BOUNDLDERIVATIVECRUCICALEIKONALFUNCTIONQUANTITIES}),
	the bootstrap assumption $\|\Rad \Psi\|_{C^0(\Sigma_t^u)} \leq \varepsilon (1 + t)^{-1}$
	(that is, \eqref{E:PSIFUNDAMENTALC0BOUNDBOOTSTRAP}),
	and the definition of $\upmu_{\star}$ to deduce that
	$\frac{1}{\rgeo^2} \frac{\Rad \Psi}{\upmu}$ times the 
	first time integral on the right-hand side of 
	\eqref{E:RENORMALIZEDTOPORDERTRCHIJUNKTRANSPORTINVERTED}
	is bounded in magnitude by
	\begin{align} \label{E:FIRSTDIFFICULTTIMEINTEGRALDIFFICULTMULTPOINTWISEENERGYERRORPRODUCT}
		& \leq 
			2 (1 + C \varepsilon)
			\frac{1}{\upmu}
			\frac{1}{\rgeo^2}
			|\Rad \Psi|
			\int_{s=0}^t 
								\left(
									\frac{\upmu(t,u,\vartheta)}{\upmu(s,u,\vartheta)}
								\right)^2 
								\rgeo^2(s,u)
								\frac{[\Lunit \upmu(s,u,\vartheta)]_-}{\upmu(s,u,\vartheta)}
								\left| \mathfrak{X}^{[N]} \right|(s,u,\vartheta)
			\, ds
				\\
		& \leq 	2 (1 + C \sqrt{\varepsilon})
			\frac{1}{\upmu}
			\frac{1}{\rgeo^2}
			|\Rad \Psi|
			\int_{s=0}^t 
								\rgeo^2(s,u)
								\frac
									{
									\left \|
										[\Lunit \upmu]_-
									\right \|_{C^0(\Sigma_s^u)}
									}
									{\upmu_{\star}(s,u)}
								\left| G_{\Lunit \Lunit} \Rad \mathscr{S}^N \Psi \right|(s,u,\vartheta)
			\, ds
				\notag
				\\
		& \ \ 	
			+ C \varepsilon 
					\frac{1}{\upmu_{\star}}
					\frac{\ln^3(\myexp + t)}{(1 + t)^3}
					\int_{s=0}^t 
						\left| 
							\Rad \mathscr{S}^N \Psi
						\right|(s,u,\vartheta)
					\, ds
			\notag \\
		& \ \
			+ C \varepsilon
			\frac{1}{\upmu_{\star}}
			\frac{1}{(1 + t)^3}
			\int_{s=0}^t 
				\left| 
					\myarray
						[\rgeo \left\lbrace \Lunit + \frac{1}{2} \mytr \upchi \right\rbrace \mathscr{S}^{\leq N} \Psi]
						{\rgeo \angdiff \mathscr{S}^{\leq N} \Psi}
				\right|(s,u,\vartheta)
			\, ds
				\notag
				\\
		& \ \
			+ C \varepsilon
			\frac{1}{\upmu_{\star}}
			\frac{1}{(1 + t)^3}
			\int_{s=0}^t 
				\rgeo(s,u)
				\frac{1}{\upmu_{\star}(s,u)}
				\left| \Rad \mathscr{Z}^{\leq N-1} \Psi \right|(s,u,\vartheta)
			\, ds
				\notag
				\\
		& \ \
			+ C \varepsilon
			\frac{1}{\upmu_{\star}}
			\frac{\ln(\myexp + t)}{(1 + t)^3}
			\int_{s=0}^t 
									\frac
									{1}
									{\upmu_{\star}(s,u)}
								\left| 
									\threemyarray[\rgeo \left\lbrace \Lunit + \frac{1}{2} \mytr \upchi \right\rbrace \mathscr{Z}^{\leq N-1} \Psi]
										{\rgeo \angdiff \mathscr{Z}^{\leq N-1} \Psi}
										{\mathscr{Z}^{\leq N-1} \Psi}
								\right|(s,u,\vartheta)
			\, ds
				\notag
				\\
		& \ \
			+ C \varepsilon
			\frac{1}{\upmu_{\star}}
			\frac{1}{(1 + t)^3}
			\int_{s=0}^t 
								\frac{1}{(1 + s)}
								\frac{1}{\upmu_{\star}(s,u)}
								\left| 
									\myarray
									[\mathscr{Z}^{\leq N} (\upmu - 1)]
										{\rgeo \sum_{a=1}^3 |\mathscr{Z}^{\leq N} \Lunit_{(Small)}^a|}
								\right|(s,u,\vartheta)
			\, ds.
				\notag
	\end{align}
	The main challenge is to bound the (difficult) first product 
	$2 (1 + C \sqrt{\varepsilon}) \cdots$
	on the right-hand side of 
	\eqref{E:FIRSTDIFFICULTTIMEINTEGRALDIFFICULTMULTPOINTWISEENERGYERRORPRODUCT}
	by the right-hand side of \eqref{E:TOPORDERTRCHIJUNKENERGYERRORTERMKEYPOINTWISEESTIMATE};
	it is straightforward to see that the remaining products on the right-hand side of 
	\eqref{E:FIRSTDIFFICULTTIMEINTEGRALDIFFICULTMULTPOINTWISEENERGYERRORPRODUCT}
	are $\leq$ the right-hand side of \eqref{E:TOPORDERTRCHIJUNKENERGYERRORTERMKEYPOINTWISEESTIMATE} 
	as desired.
	To bound the difficult first product on the right-hand side of \eqref{E:FIRSTDIFFICULTTIMEINTEGRALDIFFICULTMULTPOINTWISEENERGYERRORPRODUCT},
	we first use the estimate \eqref{E:IMPORTANTGLLDIFFERENCEESTIMATE},
	the estimate $\left \|
										[\Lunit \upmu]_-
							\right \|_{C^0(\Sigma_s^u)} \lesssim \varepsilon (1 + s)^{-1}$
	(that is, \eqref{E:C0BOUNDLDERIVATIVECRUCICALEIKONALFUNCTIONQUANTITIES})							
	and the bootstrap assumption $\|\Rad \Psi\|_{C^0(\Sigma_t^u)} \leq \varepsilon (1 + t)^{-1}$ 
	(that is, \eqref{E:PSIFUNDAMENTALC0BOUNDBOOTSTRAP})
	to replace the integrand factor $G_{\Lunit \Lunit}(s,u,\vartheta)$ 
	with $G_{\Lunit \Lunit}(t,u,\vartheta)$
	up to error terms, which allows us 
	to bound the first product by
	\begin{align}  \label{E:SECONDESTIMATEFIRSTDIFFICULTTIMEINTEGRALDIFFICULTMULTPOINTWISEENERGYERRORPRODUCT}
		& \leq 
			2(1 + C \sqrt{\varepsilon}) 
			\frac{1}{\rgeo}
			\frac{1}{\upmu}
			|\Rad \Psi|
			|G_{\Lunit \Lunit}|
			\int_{s=0}^t 
								\rgeo(s,u)
								\frac
									{
									\left \|
										[\Lunit \upmu]_-
									\right \|_{C^0(\Sigma_s^u)}
									}
									{\upmu_{\star}(s,u)}
								\left|
									\Rad \mathscr{S}^N \Psi 
								\right|(s,u,\vartheta)
			\, ds
				\\
	& \ \
		+ C \varepsilon
			\frac{1}{(1 + t)^3}
			\frac{1}{\upmu_{\star}}
			\int_{s=0}^t 
				 \frac{1}{\upmu_{\star}(s,u)}
				 \left|
				  	\Rad \mathscr{S}^N \Psi 
				 \right|(s,u,\vartheta)
			\, ds.
			\notag 
	\end{align}
	Clearly, the second product on the right-hand side of 
	\eqref{E:SECONDESTIMATEFIRSTDIFFICULTTIMEINTEGRALDIFFICULTMULTPOINTWISEENERGYERRORPRODUCT}
	is $\leq$ the right-hand side of \eqref{E:TOPORDERTRCHIJUNKENERGYERRORTERMKEYPOINTWISEESTIMATE} as desired.
	To bound the first product on the right-hand side of 
	\eqref{E:SECONDESTIMATEFIRSTDIFFICULTTIMEINTEGRALDIFFICULTMULTPOINTWISEENERGYERRORPRODUCT},
	we use essentially the same reasoning that we used to deduce \eqref{E:KEYGLLRADPSIOVERUPMUTIMESRADZNPSIBOUND}
	in order to bound it by
	\begin{align} \label{E:THIRDESTIMATEFIRSTDIFFICULTTIMEINTEGRALDIFFICULTMULTPOINTWISEENERGYERRORPRODUCT}
		& \leq 
			4(1 + C \sqrt{\varepsilon}) 
			\frac{1}{\rgeo}
			\left \|
				\frac{\Lunit \upmu}{\upmu}
			\right \|_{C^0(\Sigma_t^u)}
			\int_{s=0}^t 
								\rgeo(s,u)
								\frac
									{
									\left \|
										[\Lunit \upmu]_-
									\right \|_{C^0(\Sigma_s^u)}
									}
									{\upmu_{\star}(s,u)}
								\left|
									\Rad \mathscr{S}^N \Psi 
								\right|(s,u,\vartheta)
			\, ds
				\\
		& \ \
			+
			C \varepsilon
			\frac{1}{(1 + t)^3}
			\int_{s=0}^t 
								\frac
									{1}{\upmu_{\star}(s,u)}
								\left|
									\Rad \mathscr{S}^N \Psi 
								\right|(s,u,\vartheta)
			\, ds.
			\notag
	\end{align}
	Finally, we observe that both terms on the right-hand side of 
	\eqref{E:THIRDESTIMATEFIRSTDIFFICULTTIMEINTEGRALDIFFICULTMULTPOINTWISEENERGYERRORPRODUCT}
	are $\leq$ the right-hand side of \eqref{E:TOPORDERTRCHIJUNKENERGYERRORTERMKEYPOINTWISEESTIMATE} as desired.
	
	\noindent{\textbf{Bound for the term arising from the second time integral on the right-hand side of} 
	\eqref{E:RENORMALIZEDTOPORDERTRCHIJUNKTRANSPORTINVERTED}}:	
	We now bound 
	$\frac{1}{\rgeo^2} \frac{\Rad \Psi}{\upmu}$ 
	times the 
	second time integral on the right-hand side of 
	\eqref{E:RENORMALIZEDTOPORDERTRCHIJUNKTRANSPORTINVERTED}.
	We begin by using the estimate \eqref{E:C0BOUNDCRUCIALEIKONALFUNCTIONQUANTITIES}
	to pointwise bound the integrand factor as follows:
	$\mytr \upchi(s,u,\vartheta) = 2 \rgeo^{-1}(s,u) + \mytr \upchi^{(Small)}(s,u,\vartheta)
	\leq 2 (1 + C \varepsilon) \rgeo^{-1}(s,u).$
	Next, we substitute this bound and the bound
	\eqref{E:MAXEDRENORMALIZEDTRCHIJUNKDISCREPANCYPOINTWISE}
	for the integrand factor $\left| \mathfrak{X}^{[N]} \right|$
	into the second time integral. We give a separate argument for each term on the right-hand side of
	\eqref{E:MAXEDRENORMALIZEDTRCHIJUNKDISCREPANCYPOINTWISE}.
	First, using the above estimate for $\mytr \upchi,$ 
	we bound $\frac{1}{\rgeo^2} \frac{\Rad \Psi}{\upmu}$ 
	times the part of the second time integral that arises from
	the first term on the right-hand side of \eqref{E:MAXEDRENORMALIZEDTRCHIJUNKDISCREPANCYPOINTWISE} by
 	\begin{align} \label{E:DIFFICULTPARTOFTHESECONDDIFFICULTTIMEINTEGRALDIFFICULTMULTPOINTWISEENERGYERRORPRODUCT}
		& \leq
			(1 + C \varepsilon)
			\frac{1}{\rgeo^2}
			\frac{1}{\upmu}
			|\Rad \Psi|
			\int_{s=0}^t 
				\left(
					\frac{\upmu(t,u,\vartheta)}{\upmu(s,u,\vartheta)}
				\right)^2
				\rgeo(s,u)
				\left|
					G_{\Lunit \Lunit} \Rad \mathscr{S}^N \Psi
				\right|
				(s,u,\vartheta)
			\, ds.
	\end{align}
 	We now use arguments similar to the ones we used above in order 
 	to replace the integrand factor $G_{\Lunit \Lunit}(s,u,\vartheta)$ 
	with $G_{\Lunit \Lunit}(t,u,\vartheta)$
	in \eqref{E:DIFFICULTPARTOFTHESECONDDIFFICULTTIMEINTEGRALDIFFICULTMULTPOINTWISEENERGYERRORPRODUCT}
 	up to error terms and to then bound the product $\frac{1}{\upmu} |G_{\Lunit \Lunit}| |\Rad \Psi|$
 	by $2 \left|\frac{\Lunit \upmu}{\upmu} \right|$ up to error terms.
 	This line of reasoning allows us to bound the right-hand side of 
 	\eqref{E:DIFFICULTPARTOFTHESECONDDIFFICULTTIMEINTEGRALDIFFICULTMULTPOINTWISEENERGYERRORPRODUCT}
 	by the product $\boxed{2}(1 + C \varepsilon) \frac{1}{\rgeo^2(t,u)} \cdots$
 	on the right-hand side of \eqref{E:TOPORDERTRCHIJUNKENERGYERRORTERMKEYPOINTWISEESTIMATE}
 	plus other non-boxed-constant-multiplied error terms
 	on the right-hand side of \eqref{E:TOPORDERTRCHIJUNKENERGYERRORTERMKEYPOINTWISEESTIMATE}.
	
	Next, we use the estimate $|\rgeo \mytr \upchi| \lesssim 1$ noted above
	and the bootstrap assumption $\|\Rad \Psi\|_{C^0(\Sigma_t^u)} \leq \varepsilon (1 + t)^{-1}$
	to bound $\frac{1}{\rgeo^2} \frac{\Rad \Psi}{\upmu}$ 
	times the part of the second time integral that arises from
	the second product on the right-hand side of \eqref{E:MAXEDRENORMALIZEDTRCHIJUNKDISCREPANCYPOINTWISE} by
	\begin{align} \label{E:EASYTOPORDERPARTPARTOFTHESECONDDIFFICULTTIMEINTEGRALDIFFICULTMULTPOINTWISEENERGYERRORPRODUCT}
		& \leq
			C \varepsilon
			\frac{1}{(1 + t)^3}
			\int_{s=0}^t 
					\frac{\upmu(t,u,\vartheta)}{\upmu(s,u,\vartheta)}
				\left|
					\myarray
						[\rgeo \left\lbrace \Lunit + \frac{1}{2} \mytr \upchi \right\rbrace \mathscr{S}^{\leq N} \Psi]
						{\rgeo \angdiff \mathscr{S}^{\leq N} \Psi}
				\right|
				(s,u,\vartheta)
			\, ds.
	\end{align}
	Furthermore, from Definition \eqref{D:REGIONSOFDISTINCTUPMUBEHAVIOR}
	and the estimates 
	\eqref{E:LOCALIZEDMUCANTGROWTOOFAST}
	and	
	\eqref{E:LOCALIZEDMUMUSTSHRINK},
	it follows that 
	the factor $\frac{\upmu(t,u,\vartheta)}{\upmu(s,u,\vartheta)}$ in
	\eqref{E:EASYTOPORDERPARTPARTOFTHESECONDDIFFICULTTIMEINTEGRALDIFFICULTMULTPOINTWISEENERGYERRORPRODUCT}
	is $\lesssim \ln(\myexp + t).$ In total, we see
	that \eqref{E:EASYTOPORDERPARTPARTOFTHESECONDDIFFICULTTIMEINTEGRALDIFFICULTMULTPOINTWISEENERGYERRORPRODUCT}
	is $\leq$ the right-hand side of \eqref{E:TOPORDERTRCHIJUNKENERGYERRORTERMKEYPOINTWISEESTIMATE} as desired.
	
	Using similar reasoning,
	we bound $\frac{1}{\rgeo^2} \frac{\Rad \Psi}{\upmu}$ 
	times the part of the second time integral that arises from
	the third, fourth, and fifth products on the right-hand side of \eqref{E:MAXEDRENORMALIZEDTRCHIJUNKDISCREPANCYPOINTWISE} 
	respectively by
	\begin{align} \label{E:LOWERORDERPARTPARTOFTHESECONDDIFFICULTTIMEINTEGRALDIFFICULTMULTPOINTWISEENERGYERRORPRODUCT}
		& \leq
			C \varepsilon
			\frac{1}{\upmu_{\star}}
			\frac{\ln^3(\myexp + t)}{(1 + t)^3}
			\int_{s=0}^t 
				\left|
					\threemyarray
						[\rgeo \left\lbrace \Lunit + \frac{1}{2} \mytr \upchi \right\rbrace \mathscr{Z}^{\leq N-1} \Psi]
						{\rgeo \angdiff \mathscr{Z}^{\leq N-1} \Psi}
						{\mathscr{Z}^{\leq N-1} \Psi}
				\right|
				(s,u,\vartheta)
			\, ds,
	\end{align}
	\begin{align} \label{E:RADPARTLOWERORDERPARTPARTOFTHESECONDDIFFICULTTIMEINTEGRALDIFFICULTMULTPOINTWISEENERGYERRORPRODUCT}
		& \leq
			C \varepsilon
			\frac{1}{\upmu_{\star}}
			\frac{\ln^2(\myexp + t)}{(1 + t)^3}
			\int_{s=0}^t 
				\rgeo(s,u,\vartheta)
				\left|
					\Rad \mathscr{Z}^{N-1} \Psi
				\right|
				(s,u,\vartheta)
			\, ds,
	\end{align}
	and
	\begin{align} \label{E:EIKONALTERMSPARTOFTHESECONDDIFFICULTTIMEINTEGRALDIFFICULTMULTPOINTWISEENERGYERRORPRODUCT}
		& \leq
			C \varepsilon
			\frac{1}{\upmu_{\star}}
			\frac{\ln^2(\myexp + t)}{(1 + t)^3}
			\int_{s=0}^t 
			\frac{1}{1 + s}
			\left| 
				\myarray
					[\mathscr{Z}^{\leq N} (\upmu - 1)]
					{\rgeo \sum_{a=1}^3 |\mathscr{Z}^{\leq N} \Lunit_{(Small)}^a|}
			\right|
				(s,u,\vartheta)
			\, ds
	\end{align}
	as desired.

\noindent{\textbf{Bound for the term arising from the third time integral on the right-hand side of} \eqref{E:RENORMALIZEDTOPORDERTRCHIJUNKTRANSPORTINVERTED}}:	 
We now pointwise bound 
$\frac{1}{\rgeo^2} \frac{\Rad \Psi}{\upmu}$ times the 
third time integral on the right-hand side of 
\eqref{E:RENORMALIZEDTOPORDERTRCHIJUNKTRANSPORTINVERTED}.
We first note that \eqref{E:C0BOUNDCRUCIALEIKONALFUNCTIONQUANTITIES}
implies the pointwise bound
$\left| \hat{\upchi}^{(Small)} \right|(s,u,\vartheta) \lesssim \epsilon \ln(\myexp + s)(1 + s)^{-2}.$
Also using the bootstrap assumption $\|\Rad \Psi\|_{C^0(\Sigma_t^u)} \leq \varepsilon (1 + t)^{-1}$
(that is, \eqref{E:PSIFUNDAMENTALC0BOUNDBOOTSTRAP}),
we deduce that $\frac{1}{\rgeo^2} \frac{\Rad \Psi}{\upmu}$ times the 
third time integral on the right-hand side of 
\eqref{E:RENORMALIZEDTOPORDERTRCHIJUNKTRANSPORTINVERTED}
is bounded in magnitude by 
\begin{align} \label{E:THIRDDIFFICULTTIMEINTEGRALDIFFICULTMULTPOINTWISEENERGYERRORPRODUCT}
	\leq 
	C
	\varepsilon \frac{\ln(\myexp + t)}{(1 + t)^3}
	\frac{1}{\upmu(t,u,\vartheta)}
	\int_{s=0}^t
		\frac{\upmu^2(t,u,\vartheta)}{\upmu^2(s,u,\vartheta)}
		\left| \upmu \angLie_{\mathscr{S}}^{\leq N} \hat{\upchi}^{(Small)} \right|(s,u,\vartheta)	
	\, ds.
\end{align}

We now notice that the analysis 
of the right-hand side of \eqref{E:THIRDDIFFICULTTIMEINTEGRALDIFFICULTMULTPOINTWISEENERGYERRORPRODUCT} 
easily reduces to the case of the 
top-order operators $\angLie_{\mathscr{S}}^N;$
in the cases of the 
lower-order operators $\angLie_{\mathscr{S}}^{\leq N-1},$
we use inequality \eqref{E:POINTWISEESTIMATESFORCHIJUNKINTERMSOFOTHERVARIABLES}
with $\leq N-1$ in the role of $N$
and the pointwise bound $\left| \frac{\upmu(t,u,\vartheta)}{\upmu(s,u,\vartheta)} \right| \lesssim \ln(\myexp + t),$
(which follows from the line of reasoning just below  \eqref{E:EASYTOPORDERPARTPARTOFTHESECONDDIFFICULTTIMEINTEGRALDIFFICULTMULTPOINTWISEENERGYERRORPRODUCT})
to bound the below-top-order terms
on the right-hand side of \eqref{E:THIRDDIFFICULTTIMEINTEGRALDIFFICULTMULTPOINTWISEENERGYERRORPRODUCT}  
by $\leq$ the right-hand side of \eqref{E:TOPORDERTRCHIJUNKENERGYERRORTERMKEYPOINTWISEESTIMATE} as desired. 

We now consider two cases depending on the structure of 
the top-order operators $\angLie_{\mathscr{S}}^N$
in \eqref{E:THIRDDIFFICULTTIMEINTEGRALDIFFICULTMULTPOINTWISEENERGYERRORPRODUCT}:
\textbf{i)} $\mathscr{S}^N = \Rad \mathscr{S}^{N-1},$
and 
\textbf{ii)} $\mathscr{S}^N = \Rot \mathscr{S}^{N-1}.$
If $\mathscr{S}^N = \Rad \mathscr{S}^{N-1},$ then we use 
inequality \eqref{E:FUNCTIONAVOIDINGCOMMUTING} with $\Psi$ in the role of $f,$
the second inequality in 
\eqref{E:TOPORDERDERIVATIVESOFANGDSQUAREDUPMUINTERMSOFCONTROLLABLE},
and the pointwise bound $\frac{\upmu(t,u,\vartheta)}{\upmu(s,u,\vartheta)} \lesssim \ln(\myexp + t)$
noted above to deduce that the right-hand side of \eqref{E:THIRDDIFFICULTTIMEINTEGRALDIFFICULTMULTPOINTWISEENERGYERRORPRODUCT} is 
\begin{align} \label{E:RADCASETHIRDDIFFICULTTIMEINTEGRALDIFFICULTMULTPOINTWISEENERGYERRORPRODUCT}
	&
	\leq
	C
	\varepsilon \frac{\ln^3(\myexp + t)}{(1 + t)^3}
	\frac{1}{\upmu_{\star}(t,u)}
	\int_{s=0}^t
		\left| 
			\upmu \angfreeDsquared \mathscr{S}^{\leq N-1} \upmu 
		\right|(s,u,\vartheta)	
	\, ds
		\\
	& \ \ 
	+ C
		\varepsilon \frac{\ln(\myexp + t)}{(1 + t)^3}
		\int_{s=0}^t
			\frac{1}{1 + s}
			\left| 
					\fourmyarray
					 [\rgeo \left\lbrace \Lunit + \frac{1}{2} \mytr \upchi \right\rbrace \mathscr{Z}^{\leq N} \Psi]
					 {\Rad \mathscr{Z}^{\leq N} \Psi}
					 {\rgeo \angdiff \mathscr{Z}^{\leq N} \Psi}
					 {\mathscr{Z}^{\leq N} \Psi}
			\right|
			(s,u,\vartheta)	
		\, ds
		\notag	\\
	& \ \ 
		+
		C
		\varepsilon \frac{\ln(\myexp + t)}{(1 + t)^3}
		\int_{s=0}^t
			\frac{1}{(1 + s)^2}
			\left| 
				\myarray
				[\mathscr{Z}^{\leq N} (\upmu - 1)]
				{\rgeo \sum_{a=1}^3 |\mathscr{Z}^{\leq N} \Lunit_{(Small)}^a|}
			\right|
			(s,u,\vartheta)	
		\, ds.
		\notag
\end{align}
It is easy to see that the right-hand side of \eqref{E:RADCASETHIRDDIFFICULTTIMEINTEGRALDIFFICULTMULTPOINTWISEENERGYERRORPRODUCT}
is $\leq$ the right-hand side of \eqref{E:TOPORDERTRCHIJUNKENERGYERRORTERMKEYPOINTWISEESTIMATE} as desired. 

If $\mathscr{S}^N = \Rot \mathscr{S}^{N-1},$ then 
we first use inequality \eqref{E:ANGLIEROTOFATENSORINTERMSOFANGD}
to deduce that
\begin{align}  \label{E:ROTCASETHIRDDIFFICULTTIMEINTEGRALDIFFICULTMULTPOINTWISEENERGYERRORPRODUCT}
	\left| 
		\upmu \angLie_{\mathscr{S}}^N \hat{\upchi}^{(Small)} 
	\right|(s,u,\vartheta)	
	& \lesssim 
		\rgeo
		\left| 
			\upmu \angD \angLie_{\mathscr{S}}^{\leq N-1} \hat{\upchi}^{(Small)} 
		\right|(s,u,\vartheta)
		+ \left| 
				\upmu \angLie_{\mathscr{S}}^{\leq N-1} \hat{\upchi}^{(Small)} 
			\right|(s,u,\vartheta).
\end{align}
Using inequality \eqref{E:ROTCASETHIRDDIFFICULTTIMEINTEGRALDIFFICULTMULTPOINTWISEENERGYERRORPRODUCT},
reasoning as in the previous case $\mathscr{S}^N = \Rad \mathscr{S}^{N-1},$ 
and using inequality \eqref{E:POINTWISEESTIMATESFORCHIJUNKINTERMSOFOTHERVARIABLES}
with $N-1$ in the role of $N$ to pointwise bound the second (lower-order) term on the right-hand side of 
\eqref{E:ROTCASETHIRDDIFFICULTTIMEINTEGRALDIFFICULTMULTPOINTWISEENERGYERRORPRODUCT},
we conclude that in the present case $\mathscr{S}^N = \Rot \mathscr{S}^{N-1},$
the right-hand side of \eqref{E:THIRDDIFFICULTTIMEINTEGRALDIFFICULTMULTPOINTWISEENERGYERRORPRODUCT} is
$\leq$ the right-hand side of \eqref{E:TOPORDERTRCHIJUNKENERGYERRORTERMKEYPOINTWISEESTIMATE} as desired.
 

\noindent{\textbf{Bound for the term arising from the final time integral on the right-hand side of} 
\eqref{E:RENORMALIZEDTOPORDERTRCHIJUNKTRANSPORTINVERTED}}:	
We now pointwise bound $\frac{1}{\rgeo^2} \frac{\Rad \Psi}{\upmu}$ times the 
final time integral on the right-hand side of 
\eqref{E:RENORMALIZEDTOPORDERTRCHIJUNKTRANSPORTINVERTED}.
The integrand terms
$\left| \mathfrak{I}^{[N]} \right|(s,u,\vartheta)$
and
$\left|\widetilde{\mathfrak{I}}^{N} \right|(s,u,\vartheta)$
have already been suitably bounded by 
virtue of the definition of $\left| \mathfrak{I}^{[N]} \right|(s,u,\vartheta),$
which is defined just below \eqref{E:RENORMALIZEDTOPORDERTRCHIJUNKTRANSPORTINVERTED},
and the estimates
\eqref{E:TOPORDERRENORMALIZEDTRCHIJUNKTRANSPORTINHOMOGENEOUSTERMPOINTWISEESTIMATE}
and
\eqref{E:NEEDTOTAKEMAXANDTHENGRONWALLADDITIONALINHOMOGENEOUSTERM}
(with $s$ in the role of $t$ in these estimates).
The bootstrap assumption $\|\Rad \Psi\|_{C^0(\Sigma_t^u)} \leq \varepsilon (1 + t)^{-1}$
(that is, \eqref{E:PSIFUNDAMENTALC0BOUNDBOOTSTRAP})
yields the pointwise estimate
$\left|\frac{1}{\rgeo^2} \frac{\Rad \Psi}{\upmu}\right|
\lesssim \varepsilon \frac{1}{\upmu} (1+t)^{-3}.$
Furthermore, we have already proved the integrand factor estimate
$\frac{\upmu(t,u,\vartheta)}{\upmu(s,u,\vartheta)} \lesssim \ln(\myexp + t).$
Also accounting for the integrand factor $\rgeo^2(s,u),$
we see that these estimates imply that
$\frac{1}{\rgeo^2} \frac{\Rad \Psi}{\upmu}$ times the 
final time integral on the right-hand side of 
\eqref{E:RENORMALIZEDTOPORDERTRCHIJUNKTRANSPORTINVERTED}
is in magnitude $\leq$ the right-hand side of
\eqref{E:TOPORDERTRCHIJUNKENERGYERRORTERMKEYPOINTWISEESTIMATE} as desired.
\end{proof}

In the next lemma, we provide a cruder version of
Prop.~\ref{P:MAINCOMMUTEDWAVEEQNINHOMOGENEOUSTERMPOINTWISEESTIMATES}
that is useful for bounding some of the top-order terms
that we encounter in our analysis.

\begin{lemma}[\textbf{Cruder pointwise estimates for}  
$\upmu \Rot \mathscr{S}^{N-1} \mytr \upchi^{(Small)},$
$\upmu \Rad \mathscr{S}^{N-1} \mytr \upchi^{(Small)},$
$\upmu \angLap \mathscr{S}^{N-1} \upmu$]
\label{L:SLIGHTLYLESSSHARPPOINTWISEFORANLGAPUPMUANDROTTRCHI}
Let $1 \leq N \leq 24$ be an integer, and let
$\mathscr{S}^{N-1}$ be an $(N-1)^{st}$ order pure spatial commutation vectorfield operator.
Under the small-data and bootstrap assumptions 
of Sects.~\ref{S:PSISOLVES}-\ref{S:C0BOUNDBOOTSTRAP},
if $\varepsilon$ is sufficiently small, 
then the following pointwise estimates hold on $\mathcal{M}_{\Tboot,U_0}:$
\begin{align} \label{E:SLIGHTLYLESSSHARPPOINTWISEANLGAPUPMUANDROTTRCHI}
	&
	\left|\upmu \Rot \mathscr{S}^{N-1} \mytr \upchi^{(Small)} \right|,
		\\
	&
	\left|\upmu \Rad \mathscr{S}^{N-1} \mytr \upchi^{(Small)} \right|,
		\notag \\
	& \left|\upmu \angLap \mathscr{S}^{N-1} \upmu \right|
		 \notag \\
	& \lesssim
		 \left| 
		 	\fourmyarray
		 		[\upmu \Lunit \mathscr{Z}^{\leq N} \Psi]
				{\Rad \mathscr{Z}^{\leq N} \Psi}
				{\upmu \angdiff \mathscr{Z}^{\leq N} \Psi}
				{\mathscr{Z}^{\leq N} \Psi}
			\right|	
		+ 
			\frac{1}{(1 + t)^2}
	 		\left| 
				\myarray
					[\mathscr{Z}^{\leq N} (\upmu - 1)]
					{\rgeo \sum_{a=1}^3 |\mathscr{Z}^{\leq N} \Lunit_{(Small)}^a|}
			\right|	
			\notag	\\
	& \ \ 
		+ \frac{1}{\rgeo}
			\int_{t'=0}^t 
					\frac{
						\left \|
							[\Lunit \upmu]_-
						\right \|_{C^0(\Sigma_{t'}^u)}}
							{\upmu_{\star}(t',u)}
					\rgeo(t',u)
					\left| 
						\myarray
							[\Rad \mathscr{Z}^{\leq N} \Psi]
							{\mathscr{Z}^{\leq N} \Psi}
					\right|(t',u,\vartheta)
			\, dt'
				\notag \\
		& \ \ 
			+ 
			\frac{1}{\rgeo^2}
			\int_{t'=0}^t 
				\left\|
					\left(\frac{\upmu(t,\cdot)}{\upmu}\right)^2
				\right\|_{C^0(\Sigma_{t'}^u)}
				\rgeo(t',u)
				\left| 
					\myarray
						[\Rad \mathscr{Z}^{\leq N} \Psi]
						{\mathscr{Z}^{\leq N} \Psi}
				\right|(t',u,\vartheta)
			\, dt'
			+ 
			\frac{\ln^2(\myexp + t)}{(1 + t)^2}
			\left|\chifullmod^{[N]} \right|(0,u,\vartheta)
			\notag \\
	 & \ \ + \frac{\ln^3(\myexp + t)}{(1 + t)^2}
						\int_{t'=0}^t 
							\left| 
								\fourmyarray
								[\rgeo \Lunit \mathscr{Z}^{\leq N} \Psi]
								{\Rad \mathscr{Z}^{\leq N} \Psi}
								{\rgeo \angdiff \mathscr{Z}^{\leq N} \Psi}
								{\mathscr{Z}^{\leq N} \Psi}
							\right|
							(t',u,\vartheta)
						\, dt'
						\notag \\
		& \ \ + 
					\varepsilon
					\frac{\ln^3(\myexp + t)}{(1 + t)^2}
					\int_{t'=0}^t 
						\left| 
						\myarray
							[\upmu \hat{\angD}^2 \mathscr{S}^{\leq N-1} \upmu]
							{\rgeo \upmu \angD \angLie_{\mathscr{S}}^{\leq N-1} \hat{\upchi}^{(Small)}}
						\right|
						(t',u,\vartheta)
					\, dt'
					\notag \\
		& \ \ +
					\frac{\ln^3(\myexp + t)}{(1 + t)^2}	
					\int_{t'=0}^t 
						\frac{1}{(1 + t')}
						\frac{1}{\upmu_{\star}(t',u)}
						\left| 
							\myarray
								[\mathscr{Z}^{\leq N} (\upmu - 1)]
								{\rgeo \sum_{a=1}^3 |\mathscr{Z}^{\leq N} \Lunit_{(Small)}^a|}
						\right|
						(t',u,\vartheta)
					\, dt'.
					\notag 
\end{align}
The term $\left|\chifullmod^{[N]} \right|$ 
appearing on the right-hand side of \eqref{E:SLIGHTLYLESSSHARPPOINTWISEANLGAPUPMUANDROTTRCHI}
is defined just below \eqref{E:RENORMALIZEDTOPORDERTRCHIJUNKTRANSPORTINVERTED}.
\end{lemma}

\begin{proof}
The estimate \eqref{E:SLIGHTLYLESSSHARPPOINTWISEANLGAPUPMUANDROTTRCHI}	
for $\left|\upmu \Rad \mathscr{S}^{N-1} \mytr \upchi^{(Small)}\right|$ 
and 
$\left|\upmu \Rad \mathscr{S}^{N-1} \mytr \upchi^{(Small)}\right|$
can be proved by using arguments similar to the
ones we used in the proof of \eqref{E:TOPORDERTRCHIJUNKENERGYERRORTERMKEYPOINTWISEESTIMATE},
along with inequality \eqref{E:BOUNDINGFRAMEDERIVATIVESINTERMSOFZNF} with $\Psi$ in the role of $f.$
We omit the details
(note however that the present estimates contain an extra factor of $\upmu$ 
compared to \eqref{E:TOPORDERTRCHIJUNKENERGYERRORTERMKEYPOINTWISEESTIMATE} and that the factor $\Rad \Psi$ from
\eqref{E:TOPORDERTRCHIJUNKENERGYERRORTERMKEYPOINTWISEESTIMATE} is not present).
The argument is in fact significantly simpler because we are no longer seeking to pair
$\Rad \Psi$ with $G_{\Lunit \Lunit}$ in an effort to generate the factor $\Lunit \upmu.$ 
Hence, we do not bother to replace the integrand factor $G_{\Lunit \Lunit}(t',u,\vartheta)$ 
with $G_{\Lunit \Lunit}(t,u,\vartheta)$ up to error terms 
as we did in \eqref{E:SECONDESTIMATEFIRSTDIFFICULTTIMEINTEGRALDIFFICULTMULTPOINTWISEENERGYERRORPRODUCT}.
Moreover, to derive \eqref{E:SLIGHTLYLESSSHARPPOINTWISEANLGAPUPMUANDROTTRCHI}, 
it suffices to use the following cruder version of
inequality \eqref{E:MAXEDRENORMALIZEDTRCHIJUNKDISCREPANCYPOINTWISE}:
\begin{align} \label{E:CRUDERMAXEDRENORMALIZEDTRCHIJUNKDISCREPANCYPOINTWISE}
		\left| \mathfrak{X}^{[N]} \right|
		& \lesssim
			 \left|
			 	\myarray[\Rad \mathscr{Z}^{\leq N} \Psi]
			 		{\mathscr{Z}^{\leq N} \Psi}
			 	\right|
			+ 
			\frac{1}{1+t}
			 \upmu 
				\left|
				\myarray
					[\rgeo \Lunit \mathscr{Z}^{\leq N} \Psi]
					{\rgeo \angdiff \mathscr{Z}^{\leq N} \Psi}
				\right|
		+
			\varepsilon
			\frac{1}{(1 + t)^2}
			\left| 
				\myarray
					[\mathscr{Z}^{\leq N} (\upmu - 1)]
						{\rgeo \sum_{a=1}^3 |\mathscr{Z}^{\leq N} \Lunit_{(Small)}^a|}
			\right|.
\end{align}

Inequality \eqref{E:SLIGHTLYLESSSHARPPOINTWISEANLGAPUPMUANDROTTRCHI} for 
$\left|\upmu \angLap \mathscr{S}^{N-1} \upmu \right|$
then follows from the bound \eqref{E:SLIGHTLYLESSSHARPPOINTWISEANLGAPUPMUANDROTTRCHI} for 
$\left|\upmu \Rad \mathscr{S}^{N-1} \mytr \upchi^{(Small)}\right|,$
the final inequality in
\eqref{E:TOPORDERDERIVATIVESOFANGDSQUAREDUPMUINTERMSOFCONTROLLABLE},
and the bound $\upmu \lesssim \ln(\myexp + t)$ (that is, \eqref{E:C0BOUNDCRUCIALEIKONALFUNCTIONQUANTITIES}).

\end{proof}

\section{Pointwise estimates for the difficult error integrands requiring partial modification}
\label{S:POINTWISEPARTIALLYMODIFIED}
In order to derive suitable a priori estimates for the top-order $L^2$ quantity 
$\totonemax{\leq 24},$
we need to derive sharpened pointwise and $L^2$ estimates for 
the partially modified versions of 
$\mathscr{S}^{23} \mytr \upchi^{(Small)}$
and
$\angdiff \mathscr{S}^{23} \upmu$
defined in 
Sects.~\ref{S:TRCHIJUNKPARTIALRENORMALIZATION}
and \ref{S:ANGDIFFUPMUPARTIALRENORMALIZATION};
the estimates that would follow from
Prop.~\ref{P:CRUCICALTRANSPORTINTEQUALITIES} are
not sufficient to close the top-order estimates.
We provide the sharpened pointwise estimates in the next two lemmas.

\begin{lemma}[\textbf{Sharp transport inequalities for partially modified versions of}
$\angdiff \mathscr{S}^{N-1} \upmu,$
$\mathscr{S}^{N-1} \mytr \upchi^{(Small)}$]
\label{L:REFINEDRENORMALIZEDTRANSPORTINEQUALITIES}
Let $1 \leq N \leq 24$ be an integer 
and let $\mathscr{S}^{N-1}$ be an $(N-1)^{st}$ order pure spatial commutation vectorfield operator
(see definition \eqref{E:DEFSETOFSPATIALCOMMUTATORVECTORFIELDS}).
Let $\chipartialmodarg{\mathscr{S}^{N-1}}$
be the partially modified quantity defined in \eqref{E:TRANSPORTPARTIALRENORMALIZEDTRCHIJUNK}
and let $\mupartialmodarg{\mathscr{S}^{N-1}}$ be the partially modified 
$S_{t,u}$ one-form defined in \eqref{E:TRANSPORTPARTIALRENORMALIZEDUPMU}.
Under the small-data and bootstrap assumptions 
of Sects.~\ref{S:PSISOLVES}-\ref{S:C0BOUNDBOOTSTRAP},
if $\varepsilon$ is sufficiently small, then
the following estimates hold on $\mathcal{M}_{\Tboot,U_0}:$
\begin{align}  \label{E:UPMUANDTRCHIJUNKREFINEDPARTIALLYRENORMALIZEDTRANSPORTINEQUALITIES}
 & \left|
		\left|
		\Lunit
		\left|
			\rgeo \mupartialmodarg{\mathscr{S}^{N-1}}
		\right|
		\right|
			- \frac{1}{2} 
				\left|
					\rgeo G_{\Lunit \Lunit} 
					\angdiff \mathscr{S}^{N-1} \Rad \Psi
				\right|
		\right|,	
		\\
 & \left|
		\rgeo
		\left\lbrace
			\angLie_{\Lunit}
			+ \mytr \upchi 
		\right\rbrace
		\mupartialmodarg{\mathscr{S}^{N-1}}^{\#}
			- \frac{1}{2} \rgeo G_{\Lunit \Lunit} \angdiffuparg{\#} \mathscr{S}^{N-1} \Rad \Psi
		\right|,
		\notag \\
		& \left|
			\Lunit
		  \left\lbrace
		 		\rgeo^2 \chipartialmodarg{\mathscr{S}^{N-1}}
			\right\rbrace
			+ \frac{1}{2} \rgeo^2 G_{\Lunit \Lunit} \angLap \mathscr{S}^{N-1} \Psi
	  \right|,
	  	\notag \\
	& \left|
			\rgeo^2
			\left\lbrace
				\Lunit 
				+ \mytr \upchi
			\right\rbrace	
			\chipartialmodarg{\mathscr{S}^{N-1}}
			+ 
			\frac{1}{2} \rgeo^2 G_{\Lunit \Lunit} \angLap \mathscr{S}^{N-1} \Psi
	  \right|
	  \notag \\
	& \lesssim
		\left| 
				\fourmyarray[ \rgeo \left\lbrace \Lunit + \frac{1}{2} \mytr \upchi \right\rbrace \mathscr{Z}^{\leq N-1} \Psi]
					{\Rad \mathscr{Z}^{\leq N-1} \Psi}
					{\rgeo \angdiff \mathscr{Z}^{\leq N-1} \Psi}
					{\mathscr{Z}^{\leq N-1} \Psi}
		\right|
		+ 
					\frac{\ln(\myexp + t)}{(1 + t)^2}
						\left| 
							\myarray
								[\mathscr{Z}^{\leq N} (\upmu - 1)]
								{\rgeo \sum_{a=1}^3 |\mathscr{Z}^{\leq N} \Lunit_{(Small)}^a|}
						\right|.
		\notag 
\end{align}

\end{lemma}

\begin{proof}
	\noindent \textbf{Proof of the second inequality in \eqref{E:UPMUANDTRCHIJUNKREFINEDPARTIALLYRENORMALIZEDTRANSPORTINEQUALITIES}}:
	We begin by proving the inequality \eqref{E:UPMUANDTRCHIJUNKREFINEDPARTIALLYRENORMALIZEDTRANSPORTINEQUALITIES} for the second
	term on the left-hand side. 
	From equation \eqref{E:COMMUTEDUPMUFIRSTPARTIALRENORMALIZEDTRANSPORTEQUATION}
	and the identity
	$\left\lbrace
			\angLie_{\Lunit}
			+ \mytr \upchi 
	\right\rbrace \ginversesphere = - 2 \hat{\upchi}^{(Small) \# \#},$ 
	which follows from 
	\eqref{E:CHIALTDEF},
	\eqref{E:CHIJUNKDEF},
	and
	\eqref{E:CONNECTIONBETWEENANGLIEOFGINVERSESPHEREANDDEFORMATIONTENSORS},
	we see that the desired inequality \eqref{E:UPMUANDTRCHIJUNKREFINEDPARTIALLYRENORMALIZEDTRANSPORTINEQUALITIES}
	will follow once we show that
	$\rgeo \left|\hat{\upchi}^{(Small)} \right| \left| \mupartialmodarg{\mathscr{S}^{N-1}} \right|$
	and $\rgeo \left|\mupartialmodsourcearg{\mathscr{S}^{N-1}} \right|$ are each bounded in magnitude by $\lesssim$
	the right-hand side of \eqref{E:UPMUANDTRCHIJUNKREFINEDPARTIALLYRENORMALIZEDTRANSPORTINEQUALITIES},
	where $\inhomleftexparg{\mathscr{S}^{N-1}}{\mathfrak{J}}$ is the inhomogeneous term
	\eqref{E:COMMUTEDUPMUFIRSTPARTIALRENORMALIZEDTRANSPORTINHOMOGENEOUSTERM}.
	We have already proved the desired bound for
	$\rgeo \left|\mupartialmodsourcearg{\mathscr{S}^{N-1}}\right|$
	in Lemma~\ref{L:POINTWISEESTIMATESUPMUPARTIALLYRENORMALIZDINHOMOGENEOUSTERM}.
	To bound the remaining product, we first use
	\eqref{E:FUNCTIONPOINTWISEANGDINTERMSOFANGLIEO},
	\eqref{E:LOWERORDERC0BOUNDLIEDERIVATIVESOFGRAME},
	and the bound $\upmu \lesssim \ln(\myexp + t)$
	(that is, \eqref{E:C0BOUNDCRUCIALEIKONALFUNCTIONQUANTITIES}) to deduce
	that $\rgeo \left|\hat{\upchi}^{(Small)}\right| \lesssim \varepsilon \ln(\myexp + t) (1 + t)^{-1}$
	and
	$\left| \mupartialmodarg{\mathscr{S}^{N-1}} \right| \lesssim 
	(1+t)^{-1} \left|\mathscr{Z}^{\leq N} (\upmu - 1) \right|
	+ \ln(\myexp + t) \left|\angdiff \mathscr{S}^{N-1} \Psi \right|.$
	It follows that $\rgeo \left|\hat{\upchi}^{(Small)}\right| \left| \mupartialmodarg{\mathscr{S}^{N-1}} \right|$
	is $\lesssim$
	the right-hand side of \eqref{E:UPMUANDTRCHIJUNKREFINEDPARTIALLYRENORMALIZEDTRANSPORTINEQUALITIES}
	as desired.
	
\ \\
\noindent \textbf{Proof of the first inequality in \eqref{E:UPMUANDTRCHIJUNKREFINEDPARTIALLYRENORMALIZEDTRANSPORTINEQUALITIES}}:
	We first note the following simple inequality, which holds for any $S_{t,u}$ one-form $\xi:$
	\begin{align} \label{E:LUNITRGEOXISIMPLEINEQUALITY}
		\left|
			\Lunit(\rgeo |\xi|)
		\right|
		& \leq |\upchi^{(Small)}| |\rgeo \xi|
			+ \rgeo 
				\left|
					(\angLie_{\Lunit} + \mytr \upchi)
					\xi^{\#}
				\right|.
	\end{align}
	Inequality \eqref{E:LUNITRGEOXISIMPLEINEQUALITY} follows easily from the 
	identity $\angLie_{\Lunit} \ginversesphere = - 2 \upchi^{(Small) \# \#},$
	the decomposition $\upchi = \rgeo^{-1} \gsphere + \upchi^{(Small)},$
	and the fact that $\Lunit \rgeo = 1.$
	We now set $\xi = \mupartialmodarg{\mathscr{S}^{N-1}}.$ 
	A suitable estimate for the term 
	$\rgeo 
				\left|
					(\angLie_{\Lunit} + \mytr \upchi)
					\xi^{\#}
				\right|$
	is provided by the already proven second inequality in 
	\eqref{E:UPMUANDTRCHIJUNKREFINEDPARTIALLYRENORMALIZEDTRANSPORTINEQUALITIES}.
	Furthermore, the same reasoning that we used
	in our proof of the second inequality in 
	\eqref{E:UPMUANDTRCHIJUNKREFINEDPARTIALLYRENORMALIZEDTRANSPORTINEQUALITIES}
	also implies that
	$|\upchi^{(Small)}| |\rgeo \xi| = |\upchi^{(Small)}| |\rgeo \xi^{\#}|$
	is $\lesssim$ the right-hand side of \eqref{E:UPMUANDTRCHIJUNKREFINEDPARTIALLYRENORMALIZEDTRANSPORTINEQUALITIES}.
	We have thus proved the desired bound for  
	the first term on the left-hand side of \eqref{E:UPMUANDTRCHIJUNKREFINEDPARTIALLYRENORMALIZEDTRANSPORTINEQUALITIES}.
	
	\ \\
	\noindent \textbf{Proof of the third inequality in \eqref{E:UPMUANDTRCHIJUNKREFINEDPARTIALLYRENORMALIZEDTRANSPORTINEQUALITIES}}:
	We now prove inequality \eqref{E:UPMUANDTRCHIJUNKREFINEDPARTIALLYRENORMALIZEDTRANSPORTINEQUALITIES} for the third
	term on the left-hand side. To derive the desired bound, we have to show that
	$\left|\chipartialmodsourcearg{\mathscr{S}^{N-1}}\right|$
	is $\lesssim$ the right-hand side of \eqref{E:UPMUANDTRCHIJUNKREFINEDPARTIALLYRENORMALIZEDTRANSPORTINEQUALITIES},
	where $\chipartialmodsourcearg{\mathscr{S}^{N-1}}$ is the inhomogeneous term defined in
	\eqref{E:TRCHIJUNKCOMMUTEDTRANSPORTEQNPARTIALRENORMALIZATIONINHOMOGENEOUSTERM}.
	The desired bound follows from
	Lemma
	\ref{L:POINTWISEESTIMATESTRCHIJUNKREFINEDPARTIALLYRENORMALIZDINHOMOGENEOUSTERM}
	and the estimate $\rgeo \left|\mytr \upchi \right| \lesssim 1$
	(that is, \eqref{E:CRUDELOWERORDERC0BOUNDDERIVATIVESOFANGULARDEFORMATIONTENSORS}).
	
	\ \\
	\noindent \textbf{Proof of the fourth inequality in \eqref{E:UPMUANDTRCHIJUNKREFINEDPARTIALLYRENORMALIZEDTRANSPORTINEQUALITIES}}:
	We begin by noting that if $f$ is any function, then
	\begin{align} \label{E:SIMPLELDERIVATIVETRCHIFUNCTIONRELATION}
		\rgeo^2
			\left\lbrace
				\Lunit 
				+ \mytr \upchi
			\right\rbrace
			f
			- 
			\Lunit (\rgeo^2 f)
		& = \rgeo^2 \mytr \upchi^{(Small)} f.
	\end{align}
	The identity \eqref{E:SIMPLELDERIVATIVETRCHIFUNCTIONRELATION} follows easily from the 
	decomposition $\mytr \upchi = 2 \rgeo^{-1} + \mytr \upchi^{(Small)}$
	and the fact that $\Lunit \rgeo = 1.$
	We now set $f = \chipartialmodarg{\mathscr{S}^{N-1}},$ 
	that is, equal to the relevant factor in the fourth term 
	on the left-hand side of \eqref{E:UPMUANDTRCHIJUNKREFINEDPARTIALLYRENORMALIZEDTRANSPORTINEQUALITIES}.
	A suitable estimate for the term $\Lunit (\rgeo^2 f)$
	is provided by the already proven third inequality in 
	\eqref{E:UPMUANDTRCHIJUNKREFINEDPARTIALLYRENORMALIZEDTRANSPORTINEQUALITIES}.
	By inequality \eqref{E:LDERIVATIVECRUCICALTRANSPORTINTEQUALITIES}, the magnitude of
	the remaining term $\rgeo^2 \mytr \upchi^{(Small)} f$
	is $\lesssim \varepsilon \ln(\myexp + t) |f|.$ From
	\eqref{E:LOWERORDERC0BOUNDLIEDERIVATIVESOFGRAME},
	Cor.~\ref{C:POINTWISEESTIMATESGFRAMESHARPINTERMSOFOTHERQUANTITIES},
	the estimate $\rgeo \left|\mytr \upchi \right| \lesssim 1$
	(that is, \eqref{E:CRUDELOWERORDERC0BOUNDDERIVATIVESOFANGULARDEFORMATIONTENSORS}),
	and inequality \eqref{E:POINTWISEESTIMATESFORCHIJUNKINTERMSOFOTHERVARIABLES}
	with $N-1$ in the role of $N,$
	we deduce that $\varepsilon \ln(\myexp + t) |f|$ is 
	$\lesssim$
	the right-hand side of \eqref{E:UPMUANDTRCHIJUNKREFINEDPARTIALLYRENORMALIZEDTRANSPORTINEQUALITIES}
	as desired.
	
\end{proof}	

In the next lemma, we provide
time-integrated versions of the estimates of Lemma~\ref{L:REFINEDRENORMALIZEDTRANSPORTINEQUALITIES}.

\begin{lemma}[\textbf{Sharp pointwise estimates for partially modified versions of} $\angdiff \mathscr{S}^{N-1} \upmu,$
$\mathscr{S}^{N-1} \mytr \upchi^{(Small)}$] \label{L:RENORMALIZEDSHARPPOINTWISE}
	Let $1 \leq N \leq 24$ be an integer.
	Let $\chipartialmodarg{\mathscr{S}^{N-1}}$
	be the partially modified quantity defined in \eqref{E:TRANSPORTPARTIALRENORMALIZEDTRCHIJUNK}
	and let $\mupartialmodarg{\mathscr{S}^{N-1}}$ be the partially modified 
	$S_{t,u}$ one-form defined in \eqref{E:TRANSPORTPARTIALRENORMALIZEDUPMU}.
	Under the small-data and bootstrap assumptions 
	of Sects.~\ref{S:PSISOLVES}-\ref{S:C0BOUNDBOOTSTRAP},
	if $\varepsilon$ is sufficiently small, then
	the following estimates hold on $\mathcal{M}_{\Tboot,U_0}:$	
	\begin{align} \label{E:ANGDIFFUPMUSHARPPOINTWISE}
	\left|
		\rgeo
		\mupartialmodarg{\mathscr{S}^{N-1}}
		\right|
		& \leq \frac{1}{2} 
				\left|
					G_{\Lunit \Lunit}
				\right|
					\int_{t'=0}^t
						\left|\rgeo \angdiff \mathscr{S}^N \Psi \right|(t',u,\vartheta)
					\, dt'
				\\
		& \ \ 
			+ C \varepsilon
					\int_{t'=0}^t
					 \left|
					 	\angdiff \mathscr{S}^N \Psi 
					 \right|(t',u,\vartheta)	
					\, dt'
			+ C \varepsilon
			\notag	\\
		& \ \ + C 
					\int_{t'=0}^t
					\left| 
							\fourmyarray[\rgeo \left\lbrace \Lunit + \frac{1}{2} \mytr \upchi \right\rbrace \mathscr{Z}^{\leq N-1} \Psi]
								{\Rad \mathscr{Z}^{\leq N-1} \Psi}
								{\rgeo \angdiff \mathscr{Z}^{\leq N-1} \Psi}
								{\mathscr{Z}^{\leq N-1} \Psi}
					\right|(t',u,\vartheta)
					\, dt'
					\notag \\
		& \ \ + C
					\int_{t'=0}^t
						\frac{\ln(\myexp + t')}{(1 + t')^2}
						\left| 
							\myarray
								[\mathscr{Z}^{\leq N} (\upmu - 1)]
								{\rgeo \sum_{a=1}^3 |\mathscr{Z}^{\leq N} \Lunit_{(Small)}^a|}
						\right|(t',u,\vartheta)
					\, dt',
					\notag
	\end{align}
	\begin{align} \label{E:TRCHIRENORMALIZEDSHARPPOINTWISE}
		\left|
			\rgeo^2 \chipartialmodarg{\mathscr{S}^{N-1}}
		\right|
		& \leq \frac{1}{2} 
				\left|
					G_{\Lunit \Lunit}
				\right|
					\int_{t'=0}^t
						\left|
							\rgeo^2 \angLap \mathscr{S}^{N-1} \Psi 
					\right|(t',u,\vartheta)
					\, dt'
				\\
		& \ \ 
			+ C \varepsilon
					\int_{t'=0}^t
				   \left|
					   \angdiff \mathscr{S}^N \Psi 
					 \right|(t',u,\vartheta)	
					\, dt'
			+ C \varepsilon
			\notag	\\
		& \ \ + C 
					\int_{t'=0}^t
					\left| 
							\fourmyarray[\rgeo \left\lbrace \Lunit + \frac{1}{2} \mytr \upchi \right\rbrace \mathscr{Z}^{\leq N-1} \Psi]
								{\Rad \mathscr{Z}^{\leq N-1} \Psi}
								{\rgeo \angdiff \mathscr{Z}^{\leq N-1} \Psi}
								{\mathscr{Z}^{\leq N-1} \Psi}
					\right|(t',u,\vartheta)
					\, dt'
					\notag \\
		& \ \ + C 
					\int_{t'=0}^t
						\frac{\ln(\myexp + t')}{(1 + t')^2}
						\left| 
							\myarray
								[\mathscr{Z}^{\leq N} (\upmu - 1)]
								{\rgeo \sum_{a=1}^3 |\mathscr{Z}^{\leq N} \Lunit_{(Small)}^a|}
							\right|(t',u,\vartheta)
					\, dt'.
					\notag
	\end{align}
\end{lemma}

\begin{proof}
	To deduce \eqref{E:ANGDIFFUPMUSHARPPOINTWISE}, we integrate
	the first inequality of \eqref{E:UPMUANDTRCHIJUNKREFINEDPARTIALLYRENORMALIZEDTRANSPORTINEQUALITIES} along
	the integral curves of $\Lunit.$ We use 
	the estimate \eqref{E:IMPORTANTGLLDIFFERENCEESTIMATE} 
	to replace the integrand factor $G_{\Lunit \Lunit}(t',u,\vartheta)$ 
	with $G_{\Lunit \Lunit}(t,u,\vartheta)$
	up to an error term
	that is bounded by the second time integral on the right-hand side of \eqref{E:ANGDIFFUPMUSHARPPOINTWISE}.
	We also use Lemma~\ref{L:SMALLINITIALSOBOLEVNORMS} to
	deduce that the initial data term is $\lesssim \varepsilon$ in magnitude.
	The estimate \eqref{E:ANGDIFFUPMUSHARPPOINTWISE} thus follows.
	
	The estimate \eqref{E:TRCHIRENORMALIZEDSHARPPOINTWISE} follows from
	the third inequality of \eqref{E:UPMUANDTRCHIJUNKREFINEDPARTIALLYRENORMALIZEDTRANSPORTINEQUALITIES},
	inequality \eqref{E:ANGLAPFUNCTIONPOINTWISEINTERMSOFROTATIONS} with $f=\mathscr{S}^{N-1} \Psi,$
	and the same argument.
\end{proof}

In the next lemma, we derive pointwise estimates for some negligible error terms that arise when
we a priori derive estimates for the top-order quantity $\totonemax{\leq 24}.$
Specifically, these terms arise in Sect.~\ref{S:PROOFOFLEMMADANGEROUSTOPORDERMORERRORINTEGRAL}
when, to avoid error integrals with unfavorable time-growth,
we substitute the partially modified quantities
$\Rot \chipartialmodarg{\mathscr{S}^{N-1}}$
and
$\angdiv \mupartialmodarg{\mathscr{S}^{N-1}}$
in place of $\Rot \mathscr{S}^{N-1} \mytr \upchi^{(Small)}$
and $\angLap \mathscr{S}^{N-1} \upmu.$

\begin{lemma}[\textbf{Pointwise estimates for some harmless error terms that arise from the top-order Morawetz multiplier estimates}]
\label{L:MORAWETZREPLACEMENTTERMSAREHARMLESS}
Let $1 \leq N \leq 23$ be an integer and let $\mathscr{S}^{N-1}$ be an $(N-1)^{st}$ order 
pure spatial commutation vectorfield operator
(see Def.~\ref{D:DEFSETOFSPATIALCOMMUTATORVECTORFIELDS}).
Let $\chipartialmodinhomarg{\mathscr{S}^{N-1}}$
be the function defined in \eqref{E:TRANSPORTPARTIALRENORMALIZEDTRCHIJUNKDISCREPANCY}
and let $\mupartialmodinhomarg{\mathscr{S}^{N-1}}^{\#}$
be the $S_{t,u}-$tangent vectorfield that is $\gsphere-$dual to the one-form
defined in \eqref{E:TRANSPORTPARTIALRENORMALIZEDUPMUDISCREPANCY}.
Under the small-data and bootstrap assumptions 
of Sects.~\ref{S:PSISOLVES}-\ref{S:C0BOUNDBOOTSTRAP},
if $\varepsilon$ is sufficiently small, then
the following estimates hold on $\mathcal{M}_{\Tboot,U_0}:$
\begin{align} \label{E:TRCHIJUNKMORAWETZREPLACEMENTTERMSAREHARMLESS}
	(\Rad \Psi) \Rot \chipartialmodinhomarg{\mathscr{S}^{N-1}}
	& = Harmless^{\leq N},
		\\
	(\Rad \Psi) \angdiv \mupartialmodinhomarg{\mathscr{S}^{N-1}}^{\#}
	& = Harmless^{\leq N}.
		\label{E:UPMUMORAWETZREPLACEMENTTERMSAREHARMLESS}
\end{align}
\end{lemma}

\begin{proof}
	 We first prove \eqref{E:TRCHIJUNKMORAWETZREPLACEMENTTERMSAREHARMLESS}.
	 From definition \eqref{E:TRANSPORTPARTIALRENORMALIZEDTRCHIJUNKDISCREPANCY}, 
	 the Leibniz rule,
	 the bootstrap assumptions \eqref{E:PSIFUNDAMENTALC0BOUNDBOOTSTRAP},
	 inequality \eqref{E:FUNCTIONDERIVATIVESAVOIDINGCOMMUTING} with 
	 $\Psi$ in the role of $f,$
	 the estimate \eqref{E:LOWERORDERC0BOUNDLIEDERIVATIVESOFGRAME},
	 and Cor.~\ref{C:POINTWISEESTIMATESGFRAMESHARPINTERMSOFOTHERQUANTITIES},
	 we deduce that
	 \begin{align}  \label{E:ALMOSTTHERETRCHIJUNKMORAWETZREPLACEMENTTERMSAREHARMLESS}
	 	\left|
	 		(\Rad \Psi)
	 		\Rot \chipartialmodinhomarg{\mathscr{S}^{N-1}}
	 	\right|
	 	& \lesssim 
	 		\varepsilon
	 		\frac{1}{(1 + t)^2}
	 		\left|
				\fourmyarray[\rgeo \Lunit \mathscr{Z}^{\leq N} \Psi]
					{\Rad \mathscr{Z}^{\leq N} \Psi}
					{\rgeo \angdiff \mathscr{Z}^{\leq N} \Psi} 
					{\mathscr{Z}^{\leq N} \Psi}
			\right|.
	 \end{align}	
	 In view of Def.~\ref{D:HARMLESSTERMS}, 
	 	we see that the desired result \eqref{E:TRCHIJUNKMORAWETZREPLACEMENTTERMSAREHARMLESS}
		follows from \eqref{E:ALMOSTTHERETRCHIJUNKMORAWETZREPLACEMENTTERMSAREHARMLESS}.
	 
		We now prove \eqref{E:UPMUMORAWETZREPLACEMENTTERMSAREHARMLESS}. 
		From definition \eqref{E:TRANSPORTPARTIALRENORMALIZEDUPMUDISCREPANCY}, 
	 	the Leibniz rule,
	 	inequalities \eqref{E:FUNCTIONPOINTWISEANGDINTERMSOFANGLIEO} 
	 	and \eqref{E:ANGLAPFUNCTIONPOINTWISEINTERMSOFROTATIONS},
	 	the estimate\eqref{E:LOWERORDERC0BOUNDLIEDERIVATIVESOFGRAME},
	 	the estimate $\upmu \lesssim \ln(\myexp + t)$
	 	(that is, \eqref{E:C0BOUNDCRUCIALEIKONALFUNCTIONQUANTITIES}),
	 	and the bootstrap assumptions \eqref{E:PSIFUNDAMENTALC0BOUNDBOOTSTRAP},
	 	we deduce that
		\begin{align}  \label{E:ALMOSTTHEREUPMUMORAWETZREPLACEMENTTERMSAREHARMLESS}
	 	\left|
	 		(\Rad \Psi)
	 		\angdiv \mupartialmodinhomarg{\mathscr{S}^{N-1}}^{\#}
	 	\right|
	 	& \lesssim 
	 		\varepsilon
	 		\frac{\ln(\myexp + t)}{(1 + t)^2}
	 		\left|
				\fourmyarray[\rgeo \Lunit \mathscr{Z}^{\leq N} \Psi]
					{\Rad \mathscr{Z}^{\leq N} \Psi}
					{\rgeo \angdiff \mathscr{Z}^{\leq N} \Psi} 
					{\mathscr{Z}^{\leq N} \Psi}
			\right|.
	 \end{align}			
		In view of Def.~\ref{D:HARMLESSTERMS}, 
	 	we see that the desired result \eqref{E:UPMUMORAWETZREPLACEMENTTERMSAREHARMLESS}
		follows from \eqref{E:ALMOSTTHEREUPMUMORAWETZREPLACEMENTTERMSAREHARMLESS}.
	
	\end{proof}


\chapter{Elliptic Estimates and Sobolev Embedding on \texorpdfstring{the Riemannian manifolds 
$(S_{t,u},\gsphere)$}{the Spheres}} \label{C:ELLIPTIC}
\thispagestyle{fancy}
In Chapter~\ref{C:ELLIPTIC}, we provide some elliptic estimates
corresponding to various elliptic operators on the Riemannian manifolds 
$(S_{t,u}, \gsphere).$
We use these estimates in Chapter~\ref{C:ERRORTERMSOBOLEV}
in order to close our top-order a priori $L^2$ estimates.
Specifically, we must use the elliptic estimates when 
we bound the top-order derivatives of
$\angD^2 \upmu$ and $\angD \upchi^{(Small)}.$
We also prove a standard Sobolev embedding-type proposition,
which we eventually use to improve the fundamental $C^0$ bootstrap assumptions 
\eqref{E:PSIFUNDAMENTALC0BOUNDBOOTSTRAP}.

\section{Elliptic estimates} \label{S:ELLIPTICESTIMATES}
To derive the desired elliptic estimates, we first have to obtain pointwise estimates for the
Gaussian curvature $\Gauss$ of $\gsphere.$
We start with the following basic lemma, which provides a connection between
$\Gauss$ and the Riemann curvature tensor
$\Rsphere_{ABCD}$ of $\gsphere.$ 

\begin{lemma}[\textbf{Connection between the Gaussian and Riemann curvatures of $\gsphere$}]
	\label{L:RICCIISGAUSSTIMEMETRIC}
	Let $\Gauss$ and $\Rsphere_{ABCD}$ 
	respectively denote the Gaussian curvature and Riemann curvature tensor of the
	Riemannian metric $\gsphere$ on the two-dimensional manifold $S_{t,u}.$
	Then
	\begin{align} \label{E:RIEMANNSGAUSSTUMETRIC}
		\Rsphere_{ABCD} 
		& = \Gauss 
				 \left\lbrace 
				 	\gsphere_{AC} \gsphere_{BD} 
				 	- \gsphere_{AD} \gsphere_{BC} 
				 \right\rbrace.
	\end{align}
\end{lemma}

\begin{proof}
	This is a standard result.
	See, for example, \cite{jJ2008}*{Sect. 3.3} for a proof.
\end{proof}

We now state a second basic lemma that provides the well-known
formula for the first variation of the scalar curvature $\Rsphere$ of $\gsphere.$

\begin{lemma}[\textbf{First variation formula for the first variation of $\Rsphere$}]
\label{L:FIRSTVARIATIONOFSCALARCURVATURE}
Let $\gsphereleftexp{\uplambda}$ be a one-parameter family of
Riemannian metrics on $S_{t,u}$ verifying 
$\gsphereleftexp{0} = \gsphere,$ and let
$\dot{\gsphere} := \frac{d}{d \uplambda}|_{\uplambda = 0} \gsphereleftexp{\uplambda}.$
Similarly, let $\Scalarsphereleftexp{\uplambda}$ denote the scalar curvature of 
$\gsphereleftexp{\uplambda},$ and let
$\dot{\Rsphere} := \frac{d}{d \uplambda}|_{\uplambda = 0} \Scalarsphereleftexp{\uplambda}.$ 
Then we have the following identity for the first variation of the scalar curvature:
\begin{align} \label{E:SCALARCURVATUREFIRSTVAR}
		\dot{\Rsphere}
		& = 
		- \Rsphere^{AB} \dot{\gsphere}_{AB}
		- \angLap \mytr \dot{\gsphere}
		+ (\ginversesphere)^{AB} (\ginversesphere)^{CD} \angDsquaredarg{A}{C} \dot{\gsphere}_{BD},
\end{align}
where $\Rsphere_{AB} := (\ginversesphere)^{CD} \Rsphere_{ACBD}$
is the Ricci curvature tensor of $\gsphere.$

\end{lemma}

\begin{proof}
	This is another standard result.
	See, for example, \cite{bCpLnL2006}*{Ch. 2}.
\end{proof}

Next, we use Lemma \eqref{L:FIRSTVARIATIONOFSCALARCURVATURE} 
to derive an expression for $\Rad \Rsphere.$

\begin{corollary}[\textbf{An expression $\Rad \Rsphere$}]
\label{C:FIRSTVARIATIONOFSCALARCURVATURE}
Let $\Rsphere_{AB} := (\ginversesphere)^{CD} \Rsphere_{ACBD}$
and $\Rsphere := (\ginversesphere)^{AB} \Rsphere_{AB}$
respectively denote the Ricci curvature tensor and scalar curvature of 
$\gsphere.$ Then we have the following evolution equation:
\begin{align} \label{E:RADSCALARCURVATUREFIRSTVAR}
		\Rad \Rsphere 
		& = 
		- \Rsphere^{AB} \angdeformarg{\Rad}{A}{B}
		- \angLap \mytr  \angdeform{\Rad}
		+ \angDsquaredarg{A}{B} \angdeformuparg{\Rad}{A}{B}.
\end{align}

\end{corollary}

\begin{proof}
We use Lemma~\ref{L:FIRSTVARIATIONOFSCALARCURVATURE} 
in the case that $\gsphereleftexp{\uplambda}$ is the pullback of 
$\gsphere$ by the flow map $\varphi_{(\uplambda)}$ of $\Rad$
(cf. the proof of Lemma~\ref{L:COMMUTINGVEDCTORFIELDSWITHANGD} in the case $Z=\Rot$).
Then using Lemma~\ref{L:CONNECTIONBETWEENANGLIEOFGSPHEREANDDEFORMATIONTENSORS}
and the fact that the Lie derivative of a tensorfield
with respect to a vectorfield
is the derivative (with respect to the flow parameter)
of the pullback of the tensorfield
by the flow map of the vectorfield,
we have $\dot{\gsphere} = \angLie_{\Rad} \gsphere = \angdeform{\Rad}$
and $\dot{\Rsphere} = \Rad \Rsphere.$
Substituting for $\dot{\gsphere}$ and $\dot{\Rsphere}$ in 
\eqref{E:SCALARCURVATUREFIRSTVAR}, we arrive at the desired identity
\eqref{E:RADSCALARCURVATUREFIRSTVAR}.
\end{proof}

The following corollary provides an evolution equation for $\Gauss$
along the integral curves of $\Rad.$ 

\begin{corollary}[\textbf{Evolution equation for the Gaussian curvature} $\Gauss$ \textbf{of} $\gsphere$]
\label{C:GAUSSIANCURVATUREEVOLUTION}
Let $\Gauss$ denote the Gaussian curvature of the Riemannian metric $\gsphere$ on $S_{t,u}.$
Then $\rgeo^2 \Gauss$ verifies the following transport equation:
\begin{align} 
	\Rad (\rgeo^2 \Gauss)
	& = (\rgeo^2 \Gauss)
		\left\lbrace
			\frac{\upmu - 1}{\rgeo}
			+ \upmu \mytr \upchi^{(Small)}
			- \mytr  \angkuparg{(Trans-\Psi)}
			+  \upmu \mytr  \angkuparg{(Tan-\Psi)}
		\right\rbrace	
			\label{E:RADGAUSSCURVATURE} \\
		& 
		\ \	
		+ \rgeo \angLap \upmu
		+ \rgeo^2 \angLap (\upmu \mytr \upchi^{(Small)})
		- \rgeo^2 \angLap \mytr  \angkuparg{(Trans-\Psi)}
		- \rgeo^2 \angLap (\upmu \mytr  \angkuparg{(Tan-\Psi)})
			\notag \\
		& \ \ 
		- \rgeo^2 (\angD^2)^{AB} (\upmu \upchi_{AB}^{(Small)})
		+ \rgeo^2 (\angD^2)^{AB} \angktriplearg{A}{B}{(Trans-\Psi)}
		+ \rgeo^2 (\angD^2)^{AB} (\upmu \angktriplearg{A}{B}{(Tan-\Psi)}).
			\notag
\end{align}
	In \eqref{E:RADGAUSSCURVATURE}, the $S_{t,u}$ tensors 
	$\upchi^{(Small)},$
	$\angkuparg{(Trans-\Psi)},$ 
	and
	$\angkuparg{(Tan-\Psi)}$
	are defined by
	\eqref{E:CHIJUNKDEF},
	\eqref{E:KABTRANSVERSAL},
	and \eqref{E:KABGOOD}.

\end{corollary}

\begin{proof}
 	Contracting the identity \eqref{E:RIEMANNSGAUSSTUMETRIC} 
 	against 
 	$(\ginversesphere)^{AC}$
 	and $(\ginversesphere)^{AC}(\ginversesphere)^{BD},$
	we see that $\Rsphere_{AB} = \Gauss \gsphere_{AB}$
	and $\Rsphere = 2 \Gauss,$
	where $\Rsphere_{AB}$ and $\Rsphere$  
	respectively denote the Ricci curvature and
	the scalar curvature of $\gsphere.$
	The desired identity \eqref{E:RADGAUSSCURVATURE}
	now follows from these identities,
	\eqref{E:RADSCALARCURVATUREFIRSTVAR},
	the identity $\Rad \rgeo = - 1,$
	and the identity
	$	\angdeformarg{\Rad}{A}{B}
		= - 2 \rgeo^{-1} \upmu \gsphere_{AB}
			- 2 \upmu \upchi_{AB}^{(Small)}
			+ 2 \angktriplearg{A}{B}{(Trans-\Psi)} 
			+ 2 \upmu \angktriplearg{A}{B}{(Tan-\Psi)}$
	(see \eqref{E:RADDEFORMTRFREEANG}-\eqref{E:RADDEFORMANG}).
\end{proof}

We proceed to the next lemma, which, based on the
bootstrap assumptions and some previously proven $C^0$ estimates, 
provides pointwise estimates for the Gaussian curvature $\Gauss$ of $\gsphere.$ 
This lemma is important because in order to derive our desired elliptic estimates,
we use the following consequence of the lemma: $\Gauss \geq 0.$ 
Actually, we prove a much stronger estimate in the lemma.

\begin{lemma}[\textbf{Pointwise estimate for the Gaussian curvature} $\Gauss$ \textbf{of} $\gsphere$]
\label{L:POINTWISEESTIMATERESCALEDGAUSSIANCURVATURE}
Under the small-data and bootstrap assumptions 
of Sects.~\ref{S:PSISOLVES}-\ref{S:C0BOUNDBOOTSTRAP},
if $\varepsilon$ is sufficiently small, 
then the following pointwise estimate holds for the Gaussian curvature 
$\Gauss$ of $\gsphere$ on $\mathcal{M}_{\Tboot,U_0}:$
\begin{align} \label{E:POINTWISEESTIMATERESCALEDGAUSSIANCURVATURE}
	\left|
		\rgeo^2 \Gauss  - 1 
	\right|(t,u,\vartheta)
	& \lesssim \varepsilon \frac{\ln(\myexp + t)}{1 + t}.
\end{align}

\end{lemma}

\begin{proof}
We fix the time $t.$
We will apply a Gronwall argument to the quantity 
\begin{align}
	q: = \rgeo^2 \Gauss - 1
\end{align}
along the integral curves of $\Rad,$ which are contained
in $\Sigma_t^{U_0}.$
Recall that $\Rad u = 1$ and hence we can view 
$\Rad = \frac{d}{du}$ along the integral curves. 
We first derive suitable pointwise estimates for 
the non-$\Gauss$ terms on the right-hand side of 
\eqref{E:RADGAUSSCURVATURE} at time $t.$
To begin, we use 
\eqref{E:KABTRANSVERSAL}
and \eqref{E:KABGOOD}
to deduce the following schematic identities:
$\angkuparg{(Trans-\Psi)} = G_{(Frame)} \Rad \Psi,$
$\angkuparg{(Tan-\Psi)} = G_{(Frame)} 
\myarray
	[\Lunit \Psi]
	{\angdiff \Psi}.$
Furthermore, the same identities hold for
$\mytr  \angkuparg{(Trans-\Psi)}$ and $\mytr  \angkuparg{(Tan-\Psi)}$
but with $G_{(Frame)} \ginversesphere$ in place of $G_{(Frame)}.$
We then use these schematic identities,
inequalities 
\eqref{E:FUNCTIONPOINTWISEANGDINTERMSOFANGLIEO},
\eqref{E:ANGLAPFUNCTIONPOINTWISEINTERMSOFROTATIONS},
and
\eqref{E:TYPE02TENSORANGDTWOTIMESINTERMSOFROTATIONALLIETWOTIME},
the $C^0$ estimates 
\eqref{E:LOWERORDERC0BOUNDLIEDERIVATIVESOFGRAME},
\eqref{E:CRUDELOWERORDERC0BOUNDDERIVATIVESOFANGULARDEFORMATIONTENSORS},
and
\eqref{E:C0BOUNDCRUCIALEIKONALFUNCTIONQUANTITIES},
and the bootstrap assumptions \eqref{E:PSIFUNDAMENTALC0BOUNDBOOTSTRAP}
to deduce that at time $t,$ all terms on the 
second and third lines of the right-hand side of
\eqref{E:RADGAUSSCURVATURE} have magnitudes bounded by
$\lesssim \varepsilon \ln(\myexp + t)(1+t)^{-1}.$
Similarly, we deduce that the terms in braces multiplying 
$\rgeo^2 \Gauss$
in the first line on the right-hand side of \eqref{E:RADGAUSSCURVATURE}
are in magnitude bounded by $\lesssim \varepsilon \ln(\myexp + t)(1+t)^{-1}.$
Inserting these estimates into \eqref{E:RADGAUSSCURVATURE}, 
we deduce that
\begin{align} \label{E:RESCALEDGAUSSIANCURVATUREGRONWALLREADY}
	\left|
		\frac{d}{d u} q
	\right|
	& \lesssim
		\varepsilon \frac{\ln(\myexp + t)}{1 + t} |q|
		+ \varepsilon \frac{\ln(\myexp + t)}{1 + t}.
\end{align}
Integrating \eqref{E:RESCALEDGAUSSIANCURVATUREGRONWALLREADY}
(along the integral curves of $\Rad$)
from $0$ to $u$ 
and using the fact that $q(u=0) \equiv 0$ 
(that is, that $S_{t,0}$ is a round Euclidean sphere of radius $\rgeo(t,0)= 1 + t$
and Gaussian curvature $(1+t)^{-2}$),
we deduce that
\begin{align} \label{E:RESCALEDGAUSSIANCURVATUREGRONWALLREADYINTEGRALFORM}
	\left|
		q(u)
	\right|
	& \leq 
			C \varepsilon \frac{\ln(\myexp + t)}{1 + t}
			+
			C \varepsilon \frac{\ln(\myexp + t)}{1 + t}
			\int_{u'=0}^u
				\left|
					q(u')
				\right|
			\, du'.
\end{align}
Finally, applying Gronwall's inequality to \eqref{E:RESCALEDGAUSSIANCURVATUREGRONWALLREADYINTEGRALFORM}
and using the fact that $0 \leq u \leq U_0 < 1,$
we conclude that $\left|q(u) \right| \lesssim \varepsilon \ln(\myexp + t)(1+t)^{-1},$
which is precisely the desired inequality \eqref{E:POINTWISEESTIMATERESCALEDGAUSSIANCURVATURE}.

\end{proof}





We now use Lemma~\ref{L:POINTWISEESTIMATERESCALEDGAUSSIANCURVATURE} to 
derive the desired elliptic estimates. We split the estimates into two lemmas.
The first lemma is for scalar-valued functions $f,$ and we later apply it
with $f$ equal to a high-order derivative of $\upmu.$ The second is 
for symmetric, trace-free type $\binom{0}{2}$ $S_{t,u}$ tensorfields
$\xi,$ and we later apply it with $\xi$ 
equal to a high-order derivative of $\hat{\upchi}^{(Small)}.$

\begin{lemma}[\textbf{Elliptic estimates for solutions to Poisson's equation on} $S_{t,u}$]
\label{L:POISSONELLIPTIC}
	Let $f$ be a function.
	Under the small-data and bootstrap assumptions 
	of Sects.~\ref{S:PSISOLVES}-\ref{S:C0BOUNDBOOTSTRAP},
	if $\varepsilon$ is sufficiently small, then
	the following integral inequalities hold for $(t,u) \in [0,\Tboot) \times [0,U_0]:$
	\begin{align} \label{E:POISSONELLIPTIC}
		\int_{S_{t,u}}  
			\upmu^2 |\angD^2 f|^2 
		\, d \upsilon_{\gsphere}
		& \leq 4 
					\int_{S_{t,u}} 
						\upmu^2 (\angLap f)^2
					\, d \upsilon_{\gsphere}
			+ C \varepsilon^2 
					\frac{\ln^2(\myexp + t)}{(1 + t)^2}
					\int_{S_{t,u}} 
						|\angdiff f|^2 
					\, d \upsilon_{\gsphere}.
	\end{align}
	
\end{lemma}

\begin{proof}
Integrating by parts  
on the Riemannian manifold $(S_{t,u}, \gsphere)$
and using the commutator identity
$\angDarg{A} \angLap f = \angD^{B} \angDsquaredarg{A}{B} f - \Gauss \angdiff_A f$
(where $\Gauss$ is the Gaussian curvature of $\gsphere$),
we deduce that
\begin{align} \label{E:POISSONELLIPTICINTEGRATEDBYPARTS}
	\int_{S_{t,u}}  
		\upmu^2 
		\left\lbrace
			|\angD^2 f|^2 
			+ \Gauss |\angdiff f|^2
		\right\rbrace
	\, d \upsilon_{\gsphere}
	& = 
			- 2 
		\int_{S_{t,u}}  
			\upmu
			(\angDsquaredarg{A}{B} f)(\angdiffuparg{A} f) \angdiffuparg{B} \upmu
		\, d \upsilon_{\gsphere}
		\\
	& \ \ 
		+
		\int_{S_{t,u}}  
			\upmu^2 (\angLap f)^2
			+ 2 \upmu (\angLap f) (\angdiffuparg{A} f) \angdiffarg{A} \upmu
		\, d \upsilon_{\gsphere}.
		\notag
\end{align}
Applying Cauchy-Schwarz to the last integral on the right-hand
side of \eqref{E:POISSONELLIPTICINTEGRATEDBYPARTS}, we deduce that
\begin{align} \label{E:BASICSCALARANGLUARLAPLACIANELLIPTICESTIMATE}
			\int_{S_{t,u}}  
			\upmu^2 
			\left\lbrace
				\frac{1}{2} |\angD^2 f|^2 
				+ \Gauss |\angdiff f|^2
			\right\rbrace
			\, d \upsilon_{\gsphere}
		& \leq 2 	\int_{S_{t,u}}  
								\upmu^2
								(\angLap f)^2 
							\, d \upsilon_{\gsphere}
			+ 3 	\int_{S_{t,u}}  
						|\angdiff \upmu|^2 |\angdiff f|^2 
					\, d \upsilon_{\gsphere}.
	\end{align}
	By \eqref{E:POINTWISEESTIMATERESCALEDGAUSSIANCURVATURE}, we have $\Gauss > 0.$
	Also using the inequality $|\angdiff \upmu| \lesssim \varepsilon \ln(\myexp + t)(1 + t)^{-1},$
	which follows from \eqref{E:FUNCTIONPOINTWISEANGDINTERMSOFANGLIEO}
	and \eqref{E:C0BOUNDCRUCIALEIKONALFUNCTIONQUANTITIES},
	to estimate the last integral on the right-hand side of \eqref{E:BASICSCALARANGLUARLAPLACIANELLIPTICESTIMATE}, we
	arrive at the desired inequality \eqref{E:POISSONELLIPTIC}.
\end{proof}

\begin{lemma}[\textbf{Elliptic estimates for symmetric, trace-free type $\binom{0}{2}$ tensorfields on the} $S_{t,u}$]
	\label{L:TYPE02ELLIPTICESTIMATESPHERES}
	Let $\xi$ be a \textbf{symmetric, trace-free} type $\binom{0}{2}$ $S_{t,u}$ tensorfield.
	Under the small-data and bootstrap assumptions 
	of Sects.~\ref{S:PSISOLVES}-\ref{S:C0BOUNDBOOTSTRAP},
	if $\varepsilon$ is sufficiently small, 
	then the following integral inequalities hold for $(t,u) \in [0,\Tboot) \times [0,U_0]:$
	\begin{align} \label{E:TYPE02ELLIPTICESTIMATESPHERES}
		\int_{S_{t,u}} \upmu^2 
			|\angD \xi|^2 
		\, d \spherevol
		& \leq 6 \int_{S_{t,u}} \upmu |\angdiv \xi|^2 \, d \spherevol
			+ C \varepsilon^2 \frac{\ln^2(\myexp + t)}{(1 + t)^2} 
			\int_{S_{t,u}} |\xi|^2 \, d \spherevol.
	\end{align}
	
\end{lemma}

\begin{proof}
For $\xi$ verifying the hypotheses of the lemma, 
it is straightforward to verify 
using the identity
\eqref{E:RIEMANNSGAUSSTUMETRIC}
that we have the commutator identity
$\angDsquaredarg{C}{A} \xi_B^{\ C} - \angDsquaredarg{A}{C} \xi_B^{\ C} = 2 \Gauss \xi_{AB},$
where $\Riemannsphere_{ABCD}$ is the 
Riemann curvature tensor of the metric $\gsphere_{AB}$ on $S_{t,u},$
$\Gauss$ is its Gaussian curvature,
and our curvature sign conventions are as in Def.~\ref{D:SPACETIMERIEMANN}.
Using this identity
and the identity $|\angD \xi|^2 = |\angdiv \xi|^2 + (\angDarg{A} \xi_{BC})\angD^B \xi^{AC}$ 
(which holds for $\xi$ verifying the hypotheses of the lemma),
we carry out a series of tedious but straightforward calculations 
to deduce that
\begin{align} \label{E:XISYMMTRACEFREEDIVIDENTITY}
	|\angD \xi|^2
	+ 2 \Gauss |\xi|^2
	& = 2 |\angdiv \xi|^2
		+ \angdiv Y,
		\\
	Y^A & := \xi_{BC} \angD^B \xi^{AC} - \xi^{AB} (\angdiv \xi)_B.
\end{align}
Multiplying \eqref{E:XISYMMTRACEFREEDIVIDENTITY} by $\upmu^2$ and integrating
by parts over $S_{t,u},$ we deduce that
\begin{align} \label{E:TYPE02ELLIPTICESTIMATESPHERESINTEGRATEDBYPARTS}
	\int_{S_{t,u}} 
		\upmu^2 
		\left\lbrace
			|\angD \xi|^2 + 2 \Gauss |\xi|^2
		\right\rbrace
	\, d \spherevol
	& = 2 
		\int_{S_{t,u}} 
			\upmu^2 
			|\angdiv \xi|^2
		\, d \spherevol
		- 2 
		\int_{S_{t,u}} 
		\upmu Y^A \angdiffarg{A} \upmu
		\, d \spherevol.
\end{align}
Applying Cauchy-Schwarz to the integrals on the right-hand side of 
\eqref{E:TYPE02ELLIPTICESTIMATESPHERESINTEGRATEDBYPARTS}, we deduce that
\begin{align} \label{E:JUSTNEEDTOINSERTESTIMATESTYPE02ELLIPTIC}
		\int_{S_{t,u}} \upmu^2 
			\left\lbrace
				\frac{1}{2} |\angD \xi|^2 
				+ 2 \Gauss |\xi|^2
			\right\rbrace
			\, d \spherevol
		& \leq 3 \int_{S_{t,u}} \upmu^2| \angdiv \xi|^2 \, d \spherevol
			+ 3 \int_{S_{t,u}} |\angdiff \upmu|^2 |\xi|^2 \, d \spherevol.
	\end{align}
	The desired estimate \eqref{E:TYPE02ELLIPTICESTIMATESPHERES} now follows 
	from \eqref{E:JUSTNEEDTOINSERTESTIMATESTYPE02ELLIPTIC},
	as in the end of our proof of Lemma~\ref{L:POISSONELLIPTIC}.
\end{proof}

\section{Sobolev embedding}
\label{S:SOBOLEVEMBEDDING}
After we have derived suitable a priori $L^2$ estimates, we then use the following
Sobolev embedding-type proposition and the subsequent corollary
to improve the fundamental $C^0$ bootstrap assumptions \eqref{E:PSIFUNDAMENTALC0BOUNDBOOTSTRAP}
for $\left\|\mathscr{Z}^{\leq 13} \Psi \right\|_{C^0(\Sigma_t^u)}.$
For an alternate proof, see \cite{dC2007}*{Lemma 13.1}.

\begin{proposition} [\textbf{$C^0$ bounds for $f$ in terms of the 
norm $\| \cdot \|_{L^2(S_{t,u})}$ of angular derivatives of $f$}]
\label{P:SOBOLEVONSPHERES}
	Let $\mathscr{O} := \lbrace \Rot_{(1)}, \Rot_{(2)}, \Rot_{(3)} \rbrace$ denote the set of
	$S_{t,u}-$tangent rotation vectorfields. Let $f$ be any function.
	Under the small-data and bootstrap assumptions 
	of Sects.~\ref{S:PSISOLVES}-\ref{S:C0BOUNDBOOTSTRAP},
	if $\varepsilon$ is sufficiently small, then
	the following Sobolev embedding estimate holds for $(t,u) \in [0,\Tboot) \times [0,U_0]:$
	\begin{align} \label{E:SOBOLEVONSPHERES}
		\left\| f \right\|_{C^0(S_{t,u})}
		& \lesssim \frac{1}{1 + t} 
			\left\| \mathscr{O}^{\leq 2} f \right\|_{L^2(S_{t,u})}.
	\end{align}
	
\end{proposition}

\begin{proof}
	We will prove the following Sobolev embedding result:
	\begin{align} \label{E:GOODENOUGHFORTHESOBOLEVESTIMATEWEWANT}
		\rgeo^2(t,u) \left\| f \right\|_{C^0(S_{t,u})}^2
		& \leq C 
			\int_{S_{t,u}} |f|^2 + \rgeo^2 |\angdiff f|^2 + \rgeo^4 |\angD^2 f|^2 \, d \spherevol.
	\end{align}
	Inequality \eqref{E:SOBOLEVONSPHERES} then follows from 
	\eqref{E:FUNCTIONPOINTWISEANGDINTERMSOFANGLIEO},
	\eqref{E:ANGDSQUAREDFUNCTIONPOINTWISEINTERMSOFANGDIFFROTATIONS},
	and
	\eqref{E:GOODENOUGHFORTHESOBOLEVESTIMATEWEWANT}.
	
	To prove \eqref{E:GOODENOUGHFORTHESOBOLEVESTIMATEWEWANT}, 
	we fix $t,u$ and define the rescaled metric $\rescaledmetricsphere := \rgeo^{-2} \gsphere$
	on $S_{t,u}.$ From scaling considerations, we deduce that
	\eqref{E:GOODENOUGHFORTHESOBOLEVESTIMATEWEWANT} is equivalent to
	\begin{align} \label{E:RESCALEDGOODENOUGHFORTHESOBOLEVESTIMATEWEWANT}
		\left\| f \right\|_{C^0(S_{t,u})}^2
		& \leq C 
			\int_{S_{t,u}} |f|^2 
				+ |\angdiff f|_{\rescaledmetricsphere}^2 
				+ |\angD^2 f|_{\rescaledmetricsphere}^2 
			\, d \rescaledspherevol
			=  
			C
			\int_{\vartheta \in \mathbb{S}^2}
					|f|^2 
				+ |\angdiff f|_{\rescaledmetricsphere}^2 
				+ |\angD^2 f|_{\rescaledmetricsphere}^2 
			\, d \rescaledargspherevol{(t,u,\vartheta)},
	\end{align}
	where $|\cdot|_{\rescaledmetricsphere}$ denotes the pointwise norm of a tensor as measured by $\rescaledmetricsphere,$
	and there is no ambiguity in the meaning of $\angD^2$ since the Levi-Civita connections of 
	$\gsphere$ and $\rescaledmetricsphere$ agree.
	
	To prove \eqref{E:RESCALEDGOODENOUGHFORTHESOBOLEVESTIMATEWEWANT}, we let
	$i_{(inject)}(\rescaledmetricsphere)(t,u)$ denote the 
	radius of injectivity of $(\mathbb{S}^2,\rescaledmetricsphere(t,u,\cdot)).$
	By \cite{tA1982}*{Theorem 2.20}, an estimate of the form \eqref{E:RESCALEDGOODENOUGHFORTHESOBOLEVESTIMATEWEWANT} holds,
	and furthermore, the constant $C$ can be bounded from above by a continuous function
	that depends only on a lower bound for $i_{(inject)}(\rescaledmetricsphere)(t,u)$ 
	and the absolute value of the Gaussian curvature $\Gauss(\rescaledmetricsphere)(t,u,\vartheta)$
	of $\rescaledmetricsphere(t,u,\vartheta)$
	(this follows from the proof of \cite{tA1982}*{Theorem 2.20}).
	To complete the proof of the proposition, it only remains for us to derive these bounds.
	To proceed, we first use scaling considerations and 
	the estimate
	\eqref{E:POINTWISEESTIMATERESCALEDGAUSSIANCURVATURE}
	for $\Gauss(\gsphere)$ to deduce
	\begin{align} \label{E:RESCALEDGAUSSCURVATUREESTIMATE}
		\left| 
			\Gauss(\rescaledmetricsphere)
			- 1
		\right|(t,u,\vartheta)
		& \lesssim \varepsilon \frac{\ln(\myexp + t)}{1 + t}.
	\end{align}
	The desired bounds for $|\Gauss(\rescaledmetricsphere)|$ clearly follow from 	
	\eqref{E:RESCALEDGAUSSCURVATUREESTIMATE}.
	Furthermore, since \eqref{E:RESCALEDGAUSSCURVATUREESTIMATE}
	implies in particular that $\Gauss(\rescaledmetricsphere)$ is positive, 
	an argument of Klingenberg 
	(see \cite{jCdE2008}*{Theorem 5.9} or the proof 
	given on \cite{dCsK1993}*{pg. 34})
	implies that 
	\begin{align} \label{E:RADOFINJECTLOWERBOUND}
		i_{(inject)}(\rescaledmetricsphere)(t,u)
		\geq \frac{\pi}{\sqrt{\sup_{\vartheta \in \mathbb{S}^2} \Gauss(\rescaledmetricsphere)(t,u,\vartheta)}}
		\geq \pi - C \varepsilon \frac{\ln(\myexp + t)}{1 + t}.
	\end{align}
	The desired lower bound for $i_{(inject)}(\rescaledmetricsphere)(t,u)$
	clearly follows from \eqref{E:RADOFINJECTLOWERBOUND}. We have thus proved the proposition.
	
\end{proof}

\begin{corollary}[\textbf{$C^0$ bounds for $\mathscr{Z}^N \Psi$ in terms of the energies}] 
\label{C:C0BOUNDSOBOLEVINTERMSOFENERGIES}
	Let $\enzero[\cdot](t,u)$ be the energy functional defined in \eqref{E:E0DEF}
	and let $f$ be any function.
	Under the small-data and bootstrap assumptions 
	of Sects.~\ref{S:PSISOLVES}-\ref{S:C0BOUNDBOOTSTRAP},
	if $\varepsilon$ is sufficiently small, then
	the following estimates hold for $(t,u) \in [0,\Tboot) \times [0,U_0]:$
	\begin{align} \label{E:C0BOUNDSOBOLEVINTERMSOFENERGIES}
		\left\| f \right\|_{C^0(S_{t,u})}
		& \leq C \frac{1}{1 + t} 
			\sup_{(t',u') \in [0,t] \times [0,u]} 
			\enzero^{1/2}[\mathscr{O}^{\leq 2} f](t',u').
	\end{align}
	\end{corollary}
\begin{proof}
		Inequality \eqref{E:C0BOUNDSOBOLEVINTERMSOFENERGIES}
		follows from
		\eqref{E:MULTENERGYCOERCIVITY} and
		\eqref{E:SOBOLEVONSPHERES}.
\end{proof}


\chapter[\texorpdfstring{$L^2$}{Square Integral} Estimates for \texorpdfstring{$u$}{the Eikonal Function}
That Do Not Rely on Modified Quantities]{\texorpdfstring{$L^2$}{Square Integral} Estimates for the Eikonal Function Quantities
That Do Not Rely on Modified Quantities}
\label{C:EIKONALBELOWTOPORDERSOBOLEVESTIMATES}
\thispagestyle{fancy}
In Chapter~\ref{C:EIKONALBELOWTOPORDERSOBOLEVESTIMATES}, 
we use the pointwise estimates of Prop.~\ref{P:CRUCICALTRANSPORTINTEQUALITIES} to derive
corresponding $L^2$ estimates
for the below-top-order derivatives of $\upmu,$ $\Lunit_{(Small)}^i,$ and $\upchi^{(Small)}$
in terms of the fundamental $L^2-$controlling quantities.
We also estimate the top-order derivatives involving at least one $\Lunit$ differentiation.
They are easy to derive because we allow them to lose one derivative relative to
$\Psi$ and hence we do not need to invoke the modified quantities of Ch.~\ref{C:RENORMALIZEDEIKONALFUNCTIONQUANTITIES}. 

\section[\texorpdfstring{$L^2$}{Square integral} estimates for \texorpdfstring{$u$}{the eikonal function}
that do not rely on modified quantities]{\texorpdfstring{$L^2$}{Square Integral} estimates for the eikonal function quantities
that do not rely on modified quantities}
\label{S:EIKONALBELOWTOPORDERSOBOLEVESTIMATES}
We provide the desired $L^2$ estimates in the next lemma.

\begin{lemma}[\textbf{Preliminary $L^2$ estimates for the below-top-order derivatives and the top-order $\Lunit$ derivatives
of $\upmu,$ $\Lunit_{(Small)}^i,$ and $\upchi^{(Small)}$}]
\label{L:PRELIMINARYBELOWTOPORDERSOBOLEVESTIMATESFORMUANDL}
Let $0 \leq N \leq 23$ be an integer.
Under the small-data and bootstrap assumptions 
of Sects.~\ref{S:PSISOLVES}-\ref{S:C0BOUNDBOOTSTRAP},
if $\varepsilon$ is sufficiently small, then
we have the following estimates for $(t,u) \in [0,\Tboot) \times [0,U_0]:$
\begin{subequations}
\begin{align} \label{E:EXTRADERIVATIVELOSSEIKONALFUNCTIONQUANTITIESL2BOUNDSINTERMSOFQ0ANDQ1}
			\left\| 
					\sevenmyarray
						[\mathscr{Z}^{\leq N} (\upmu - 1)]
						{\rgeo \sum_{a=1}^3 |\mathscr{Z}^{\leq N} \Lunit_{(Small)}^a|}
						{\rgeo^2 \angLie_{\mathscr{Z}}^{\leq N-1} \upchi^{(Small)}}
						{\rgeo^2 \angLie_{\mathscr{Z}}^{\leq N-1} \upchi^{(Small)\#}}
						{\rgeo^2 \mathscr{Z}^{\leq N-1} \mytr \upchi^{(Small)}}
						{\rgeo^2 \angLie_{\mathscr{Z}}^{\leq N-1} \hat{\upchi}^{(Small)}}
						{\rgeo^2 \angLie_{\mathscr{Z}}^{\leq N-1} \hat{\upchi}^{(Small)\#}}
				\right\|_{L^2(\Sigma_t^u)}
		& \lesssim	
			(1 + t) \varepsilon
			+ (1 + t)
				\int_{s=0}^t
					\frac{\totzeromax{\leq N+1}^{1/2}(s,u)}{1 + s}
				\, ds,
				\\
	\left\| 
						\sevenmyarray
							[\Lunit \mathscr{Z}^{\leq N} \upmu]
							{\rgeo \sum_{a=1}^3 |\Lunit \mathscr{Z}^{\leq N} \Lunit_{(Small)}^a|}
							{\rgeo^2 \angLie_{\Lunit} \angLie_{\mathscr{Z}}^{\leq N-1} \upchi^{(Small)}}
							{\angLie_{\Lunit} (\rgeo^2 \angLie_{\mathscr{Z}}^{\leq N-1} \upchi^{(Small)\#})}
							{\Lunit (\rgeo^2 \mathscr{Z}^{\leq N-1} \mytr \upchi^{(Small)})}
							{\rgeo^2 \angLie_{\Lunit} \angLie_{\mathscr{Z}}^{\leq N-1} \hat{\upchi}^{(Small)}}
							{\angLie_{\Lunit} (\rgeo^2 \angLie_{\mathscr{Z}}^{\leq N-1} \hat{\upchi}^{(Small)\#})}
				\right\|_{L^2(\Sigma_t^u)}
	& \lesssim
		\varepsilon \frac{\ln(\myexp + t)}{1 + t} 
		+ \totzeromax{\leq N+1}^{1/2}(t,u).
		\label{E:EXTRADERIVATIVELOSSLDERIVATIVEEIKONALFUNCTIONQUANTITIESL2BOUNDSINTERMSOFQ0ANDQ1}
\end{align}
\end{subequations}

Furthermore, if $0 \leq N \leq 24$ is an integer, we have the following estimates:
\begin{subequations}
\begin{align} \label{E:EIKONALFUNCTIONQUANTITIESL2BOUNDSINTERMSOFQ0ANDQ1}
		\left\| 
					\sevenmyarray
						[\mathscr{Z}^{\leq N} (\upmu - 1)]
						{\rgeo \sum_{a=1}^3 |\mathscr{Z}^{\leq N} \Lunit_{(Small)}^a|}
						{\rgeo^2 \angLie_{\mathscr{Z}}^{\leq N-1} \upchi^{(Small)}}
						{\rgeo^2 \angLie_{\mathscr{Z}}^{\leq N-1} \upchi^{(Small)\#}}
						{\rgeo^2 \mathscr{Z}^{\leq N-1} \mytr \upchi^{(Small)}}
						{\rgeo^2 \angLie_{\mathscr{Z}}^{\leq N-1} \hat{\upchi}^{(Small)}}
						{\rgeo^2 \angLie_{\mathscr{Z}}^{\leq N-1} \hat{\upchi}^{(Small)\#}}
		\right\|_{L^2(\Sigma_t^u)}
		& \lesssim
		(1 + t) \varepsilon
		+ 
		(1 + t) 
		\int_{s=0}^t 
			\left\lbrace
				\frac{ \totzeromax{\leq N}^{1/2}(s,u)}{1 + s} 
				+ \frac{\totonemax{\leq N}^{1/2}(s,u)}{(1 + s)\upmu_{\star}^{1/2}(s,u)} 
			\right\rbrace
		\, ds,
			\\
		\left\| 
						\sevenmyarray
							[\Lunit \mathscr{Z}^{\leq N} \upmu]
							{\rgeo \sum_{a=1}^3 |\Lunit \mathscr{Z}^{\leq N} \Lunit_{(Small)}^a|}
							{\rgeo^2 \angLie_{\Lunit} \angLie_{\mathscr{Z}}^{\leq N-1} \upchi^{(Small)}}
							{\angLie_{\Lunit} (\rgeo^2 \angLie_{\mathscr{Z}}^{\leq N-1} \upchi^{(Small)\#})}
							{\Lunit (\rgeo^2 \mathscr{Z}^{\leq N-1} \mytr \upchi^{(Small)})}
							{\rgeo^2 \angLie_{\Lunit} \angLie_{\mathscr{Z}}^{\leq N-1} \hat{\upchi}^{(Small)}}
							{\angLie_{\Lunit} (\rgeo^2 \angLie_{\mathscr{Z}}^{\leq N-1} \hat{\upchi}^{(Small)\#})}
				\right\|_{L^2(\Sigma_t^u)}	
		& \lesssim
			\varepsilon \frac{\ln(\myexp + t)}{1 + t}  
			+ \totzeromax{\leq N}^{1/2}(t,u)
			+ \upmu_{\star}^{-1/2}(t,u) \totonemax{\leq N}^{1/2}(t,u).
		\label{E:LDERIVATIVEEIKONALFUNCTIONQUANTITIESL2BOUNDSINTERMSOFQ0ANDQ1}
\end{align}
\end{subequations}

\end{lemma}

\begin{proof}	
	To simplify the notation, we set
	\begin{align}
			q_N(t,u,\vartheta) 
			&:= 
			\left| 
					\threemyarray
						[\mathscr{Z}^{\leq N} (\upmu - 1)]
						{\rgeo \sum_{a=1}^3 |\mathscr{Z}^{\leq N} \Lunit_{(Small)}^a|}
						{\mathscr{Z}^{\leq N-1} (\Lunit(\rgeo^2 \mytr \upchi^{(Small)}))}
			\right|(t,u,\vartheta),
			\qquad
			\mathring{q}_N(u,\vartheta)
			:= q_N(0,u,\vartheta).
	\end{align}
	Integrating inequality \eqref{E:LDERIVATIVECRUCICALTRANSPORTINTEQUALITIES} along the integral 
	curves of $\Lunit$ to derive a pointwise bound for $q_N$
	and then using Lemma~\ref{L:L2NORMSOFTIMEINTEGRATEDFUNCTIONS} to obtain an $L^2$ bound
	for $q_N,$
	we deduce that
	\begin{align} \label{E:EIKONALFUNCTIONLOWERORDERL2ALMOSTGRONWALLREADY}
		\frac{1}{1 + t} 
		\| q_N \|_{L^2(\Sigma_t^u)} 
		& \leq
		\frac{1}{1 + t}
		\| \mathring{q}_N \|_{L^2(\Sigma_t^u)}
		+ 
		C 
		\int_{s=0}^t
			\frac{1}{1 + s}
			\left\| 
				\fourmyarray[ \rgeo \Lunit \mathscr{Z}^{\leq N} \Psi]
					{\rgeo \angdiff \mathscr{Z}^{\leq N} \Psi}
					{\Rad \mathscr{Z}^{\leq N} \Psi}
					{\mathscr{Z}^{\leq N} \Psi}
			\right\|_{L^2(\Sigma_{s}^u)}
		\, ds
			\\
	& \ \ 
			+ C 
			\int_{s=0}^t
				\frac{\ln(\myexp + s)}{(1 + s)^3}
				\left\|
					q_N 
				\right\|_{L^2(\Sigma_s^u)}
			\, ds.
			\notag
	\end{align}
	To bound the first term on the right-hand side of \eqref{E:EIKONALFUNCTIONLOWERORDERL2ALMOSTGRONWALLREADY} 
	by $\lesssim \varepsilon,$ we use
	Lemma~\ref{L:STUINTEGRALCOMPARISON} and the definition \eqref{E:SIGMATL2NORMDEF} of 
	the norm $\| \cdot \|_{L^2(\Sigma_t^u)}$ to deduce that
	$(1+t)^{-1} \| \mathring{q}_N \|_{L^2(\Sigma_t^u)}
	\lesssim \| \mathring{q}_N \|_{L^2(\Sigma_0^u)},$
	and then Lemma~\ref{L:SMALLINITIALSOBOLEVNORMS}, 
	Lemma~\ref{L:POINTWISEESTIMATESFORCHIJUNKINTERMSOFOTHERVARIABLES},
	and \eqref{E:SMALLDATASSUMPTION}
	to deduce that
	$\left\| \mathring{q}_N \right \|_{L^2(\Sigma_0^u)} \lesssim \mathring{\upepsilon} \leq \varepsilon.$
	
	To bound the second term on the right-hand side of \eqref{E:EIKONALFUNCTIONLOWERORDERL2ALMOSTGRONWALLREADY},
	we first use Prop.~\ref{P:L2NORMSOFPSIINTERMSOFTHECOERCIVEQUANTITIES}
	to deduce that the array under the time integral is bounded in the norm $\| \cdot \|_{L^2(\Sigma_s^u)}$ 
	by $\lesssim \totzeromax{\leq N}^{1/2}(s,u) + \upmu_{\star}^{-1/2}(s,u) \totonemax{\leq N}^{1/2}(s,u).$
	Hence, the second term is
	\begin{align}
		\leq 
		C 
		\int_{s=0}^t
			\frac{1}{1 + s} \totzeromax{\leq N}^{1/2}(s,u)
		\, ds
		+
		C 
		\int_{s=0}^t
			\frac{1}{(1 + s)\upmu_{\star}^{1/2}(s,u)} \totonemax{\leq N}^{1/2}(s,u)
		\, ds.
	\end{align}	
	
	We now insert these two bounds into the right-hand side of \eqref{E:EIKONALFUNCTIONLOWERORDERL2ALMOSTGRONWALLREADY} 
	and apply Gronwall's inequality to the quantity $(1+t)^{-1} \| q_N \|_{L^2(\Sigma_t^u)},$
	which thus leads to the desired bound
	\eqref{E:EIKONALFUNCTIONQUANTITIESL2BOUNDSINTERMSOFQ0ANDQ1}
	for the first two terms and the fifth term in the array on the left-hand side.
	
	To obtain the desired bounds \eqref{E:EIKONALFUNCTIONQUANTITIESL2BOUNDSINTERMSOFQ0ANDQ1}
	for
	$\rgeo^2 \angLie_{\mathscr{Z}}^{\leq N-1} \upchi^{(Small)}$
	and
	$\rgeo^2 \angLie_{\mathscr{Z}}^{\leq N-1} \hat{\upchi}^{(Small)},$
	we first note the following inequality, which holds
	for any symmetric type $\binom{0}{2}$ $S_{t,u}$ tensorfield:
	\begin{align} \label{E:TYPE02LDERIVATIVEOFNORMIDENTITY}
		\left| 
			\Lunit (\rgeo^2 |\xi|) 
		\right|
		& \leq 
			\left|
				\rgeo^2 \angLie_{\Lunit} \xi		
			\right|
			+ 2 |\upchi^{(Small)}|
					|\rgeo^2 \xi|
			\\
		& \leq 
			\left|
				\rgeo^2 \angLie_{\Lunit} \xi		
			\right|
			+ C \varepsilon 
					\frac{\ln(\myexp + t)}{(1 + t)^2}
					|\rgeo^2 \xi|.
			\notag
	\end{align}
	The inequality \eqref{E:TYPE02LDERIVATIVEOFNORMIDENTITY}
	follows easily from the decomposition $\mytr \upchi = \rgeo^{-1} \gsphere + \upchi^{(Small)},$
	the fact that $\Lunit \rgeo = 1,$ and the estimate \eqref{E:C0BOUNDCRUCIALEIKONALFUNCTIONQUANTITIES}.
	We now apply \eqref{E:TYPE02LDERIVATIVEOFNORMIDENTITY} with 
	$\angLie_{\mathscr{Z}}^{\leq N-1} \upchi^{(Small)}$
	and
	$\angLie_{\mathscr{Z}}^{\leq N-1} \hat{\upchi}^{(Small)}$
	in the role of $\xi$
	and also use inequality \eqref{E:LDERIVATIVECRUCICALTRANSPORTINTEQUALITIES} 
	to bound the term $\left| \rgeo^2 \angLie_{\Lunit} \xi \right|.$
	A Gronwall argument similar to the one used in the previous paragraph
	now yields the desired estimates.
	We remark that the smallness of 
	the data of
	$\rgeo^2 \angLie_{\mathscr{Z}}^{\leq N-1} \upchi^{(Small)}$
	and
	$\rgeo^2 \angLie_{\mathscr{Z}}^{\leq N-1} \hat{\upchi}^{(Small)}$
	follows easily from \eqref{E:POINTWISEEIKONALFUNCTIONTRANSPORTBASEDINTEQUALITIES} at $t=0.$
	
	To obtain the desired bounds \eqref{E:EIKONALFUNCTIONQUANTITIESL2BOUNDSINTERMSOFQ0ANDQ1}
	for
	$\rgeo^2 \angLie_{\mathscr{Z}}^{\leq N-1} \upchi^{(Small) \#}$
	and
	$\rgeo^2 \angLie_{\mathscr{Z}}^{\leq N-1} \hat{\upchi}^{(Small) \#},$
	we first note the following inequality, which holds
	any type $\binom{1}{1}$ $S_{t,u}$ tensorfield:
	\begin{align} \label{E:TYPE11LDERIVATIVEOFNORMIDENTITY}
		\left| 
			\Lunit |\xi| 
		\right|
		& \leq 
			\left|
				\angLie_{\Lunit} \xi		
			\right|
			+ 2 |\upchi^{(Small)}|
					|\xi|
					\\
		& \leq 
			\left|
				\angLie_{\Lunit} \xi		
			\right|
			+ C \varepsilon 
					\frac{\ln(\myexp + t)}{(1 + t)^2}
					|\xi|.
			\notag
	\end{align}
	The inequality \eqref{E:TYPE11LDERIVATIVEOFNORMIDENTITY}
	follows easily from the decomposition $\chi = \rgeo^{-1} \gsphere + \upchi^{(Small)},$
	the fact that $\Lunit \rgeo = 1,$ and the estimate \eqref{E:C0BOUNDCRUCIALEIKONALFUNCTIONQUANTITIES}.
	We now apply \eqref{E:TYPE11LDERIVATIVEOFNORMIDENTITY} with 
	$\rgeo^2 \angLie_{\mathscr{Z}}^{\leq N-1} \upchi^{(Small) \#}$
	and
	$\rgeo^2 \angLie_{\mathscr{Z}}^{\leq N-1} \hat{\upchi}^{(Small) \#}$
	in the role of $\xi$
	and also use inequality \eqref{E:LDERIVATIVECRUCICALTRANSPORTINTEQUALITIES}
	to bound the term $\left| \angLie_{\Lunit} \xi \right|.$
	A Gronwall argument similar to the one used above now
	yields the desired estimates. We remark that the smallness of 
	the data of
	$\rgeo^2 \angLie_{\mathscr{Z}}^{\leq N-1} \upchi^{(Small) \#}$
	and
	$\rgeo^2 \angLie_{\mathscr{Z}}^{\leq N-1} \hat{\upchi}^{(Small) \#}$
	follows easily from \eqref{E:POINTWISEEIKONALFUNCTIONTRANSPORTBASEDINTEQUALITIES} at $t=0.$

	The proof of \eqref{E:EXTRADERIVATIVELOSSEIKONALFUNCTIONQUANTITIESL2BOUNDSINTERMSOFQ0ANDQ1} is 
	very similar. The only difference is that we 
	use Prop.~\ref{P:L2NORMSOFPSIINTERMSOFTHECOERCIVEQUANTITIES}
	to bound the array under the time
	integral in the second term on the right-hand side of \eqref{E:EIKONALFUNCTIONLOWERORDERL2ALMOSTGRONWALLREADY}
	in the norm $\| \cdot \|_{L^2(\Sigma_s^u)}$ by
	$\lesssim \totzeromax{\leq N+1}^{1/2}(s,u).$
	
	Inequality \eqref{E:LDERIVATIVEEIKONALFUNCTIONQUANTITIESL2BOUNDSINTERMSOFQ0ANDQ1} follows from
	the estimates \eqref{E:LDERIVATIVECRUCICALTRANSPORTINTEQUALITIES}
	and \eqref{E:EIKONALFUNCTIONQUANTITIESL2BOUNDSINTERMSOFQ0ANDQ1},
	inequality \eqref{E:LOGLOSSLESSSINGULARTERMSMTHREEFOURTHSINTEGRALBOUND},
	the fact that 
	the first array on the right-hand side of \eqref{E:LDERIVATIVECRUCICALTRANSPORTINTEQUALITIES}
	is bounded in the norm $\| \cdot \|_{L^2(\Sigma_t^u)}$ by 
	$\lesssim \totzeromax{\leq N}^{1/2}(t,u)
		+ \upmu_{\star}^{-1/2}(t,u) \totonemax{\leq N}^{1/2}(t,u)$
	(this fact follows from Prop.~\ref{P:L2NORMSOFPSIINTERMSOFTHECOERCIVEQUANTITIES}),
	and the fact that 
	$\totzeromax{\leq N}$
	and $\totonemax{\leq N}$
	are increasing in their arguments.
	
	Inequality \eqref{E:EXTRADERIVATIVELOSSLDERIVATIVEEIKONALFUNCTIONQUANTITIESL2BOUNDSINTERMSOFQ0ANDQ1}
	similarly follows from the estimates \eqref{E:LDERIVATIVECRUCICALTRANSPORTINTEQUALITIES}
	and \eqref{E:EXTRADERIVATIVELOSSEIKONALFUNCTIONQUANTITIESL2BOUNDSINTERMSOFQ0ANDQ1}
	and the fact that 
	the first array on the right-hand side of \eqref{E:LDERIVATIVECRUCICALTRANSPORTINTEQUALITIES}
	is bounded in the norm $\| \cdot \|_{L^2(\Sigma_t^u)}$ by 
	$\lesssim \totzeromax{\leq N+1}^{1/2}(t,u)$
	(this fact follows from Prop.~\ref{P:L2NORMSOFPSIINTERMSOFTHECOERCIVEQUANTITIES}).	
		
\end{proof}


\chapter{A Priori Estimates for the Fundamental \texorpdfstring{$L^2-$}{Square Integral-}Controlling Quantities} \label{C:ERRORTERMSOBOLEV}
\thispagestyle{fancy}
In Chapter~\ref{C:ERRORTERMSOBOLEV},
we derive the most important estimates of the monograph:
a priori estimates for the fundamental $L^2-$controlling quantities
$\totzeromax{N}(t,u),$ 
$\totonemax{N}(t,u),$ 
and
$\totMormax{N}(t,u).$
To this end, we first carry out the main precursor step: we use the pointwise estimates 
from Chapters 
\ref{C:POINTWISEBOUNDSFOREASYERRORINTEGRANDS}
and
\ref{C:POINTWISEESTIMATESDIFFICULTERRORINTEGRANDS} 
to bound the error integrals that arise in the energy-flux identities for the commuted wave equation.
By ``energy-flux identities,'' we mean the identities of Prop.~\ref{P:DIVTHMWITHCANCELLATIONS}
for $\mathscr{Z}^{\leq 24} \Psi.$
To make the analysis more tractable, 
we devote different sections to 
bounding the different kinds of error integrals.
We carry out this analysis in Sects.
\ref{S:EASYERRORINTEGRALESTIMATES}-\ref{S:PROOFOFLEMMADANGEROUSTOPORDERMORERRORINTEGRAL},
and in particular, we bound the most difficult error integrals 
in Sects.~\ref{S:PROOFOFLEMMADANGEROUSTOPORDERMULTERRORINTEGRAL}-\ref{S:PROOFOFLEMMADANGEROUSTOPORDERMORERRORINTEGRAL}.
We then combine the large number of error integral estimates into 
our two main propositions of this section, which are stated in Sect.~\ref{S:STATEMENTOFTWOPROPSANDGRONWALLLEMMA}:
Propositions \ref{P:MAINTOPORDERENERGYANDFLUXINTEGRALINEQUALITIES}
and \ref{P:MAINBELOWTOPORDERENERGYANDFLUXINTEGRALINEQUALITIES}.
The first proposition 
provides inequalities verified by the up-to-top-order $(N \leq 24)$
$L^2-$controlling quantities,
while the second provides less degenerate inequalities for the lower-order 
$(N \leq 23)$
$L^2-$controlling quantities.
The main point is that the estimates of the propositions
imply that the hypotheses of the Gronwall-type Lemma~\ref{L:FUNDAMENTALGRONWALL} are verified
by
$\totzeromax{N}(t,u),$ 
$\totonemax{N}(t,u),$ 
and
$\totMormax{N}(t,u).$ 
Hence, the lemma yields the desired a priori estimates.
We provide the proofs of
Propositions \ref{P:MAINTOPORDERENERGYANDFLUXINTEGRALINEQUALITIES}
and \ref{P:MAINBELOWTOPORDERENERGYANDFLUXINTEGRALINEQUALITIES}
in Sects.~\ref{S:PROOFOFPROPOSITIONMAINTOPORDERENERGYANDFLUXINTEGRALINEQUALITIES}
and \ref{S:PROOFOFPROPOSITIONMAINBELOWTOPORDERENERGYANDFLUXINTEGRALINEQUALITIES}.
In Sect.~\ref{S:PROOFOFLEMMAFUNDAMENTALGRONWALL}, we provide the proof of the Gronwall-type lemma.

\section{Bootstrap assumptions for the fundamental \texorpdfstring{$L^2-$}{square integral-}controlling quantities}
\label{S:LTWOBOOTSTRAPASSUMPTIONS}
In deriving our a priori estimates,
we find it convenient to make bootstrap assumptions 
on the $L^2-$controlling quantities
$\totzeromax{N}(t,u),$ 
$\totonemax{N}(t,u),$ 
and
$\totMormax{N}(t,u)$
from Definitions 
\ref{D:MAINCOERCIVEQUANT}
and \ref{D:COERCIVEMORAWETZINTEGRAL}.
To state the bootstrap assumptions, we first recall that $\upmu_{\star}(t,u) := \min\lbrace 1, \inf_{\Sigma_t^u} \upmu \rbrace$
and that $\varepsilon$ is the small constant from Ch.~\ref{C:C0BOUNDBOOTSTRAP}.
Our $L^2$ bootstrap assumptions are 
that the following inequalities hold for $(t,u) \in [0,\Tboot) \times [0,U_0]:$
\begin{subequations}
\begin{align}
	\totzeromax{N}^{1/2}(t,u)
	& \leq \varepsilon, 
		&& (0 \leq N \leq 15), 
		\label{E:Q0LOWBOOT} \\  
	\totonemax{N}^{1/2}(t,u)
		+ \totMormax{N}^{1/2}(t,u)
	& \leq \varepsilon \ln^2 (\myexp + t), 
		&& (0 \leq N \leq 15),
		\label{E:Q1LOWBOOT} \\
	\totzeromax{16 + M}^{1/2}(t,u)
	& \leq \varepsilon \upmu_{\star}^{- .75 - M}(t,u), 
		&& (0 \leq M \leq 7), 
		  \label{E:Q0MIDBOOT} \\  
	\totonemax{16 + M}^{1/2}(t,u)
		+ \totMormax{16 + M}^{1/2}(t,u)
	& \leq \varepsilon \ln^2 (\myexp + t) \upmu_{\star}^{-.75 - M}(t,u), 
		&& (0 \leq M \leq 7),
		\label{E:Q1MIDBOOT} \\
		\totzeromax{24}^{1/2}(t,u)
	& \leq \varepsilon \ln^{\Cononestar} (\myexp + t) \upmu_{\star}^{-8.75}(t,u), 
		&&  
				\label{E:Q0TOPBOOT} \\  
	\totonemax{24}^{1/2}(t,u)
		+ \totMormax{24}^{1/2}(t,u)
	& \leq \varepsilon \ln^{\Cononestar + 2} (\myexp + t) \upmu_{\star}^{-8.75}(t,u),
		\label{E:Q1TOPBOOT}
\end{align}
\end{subequations}
where $\Cononestar > 4$ is a large constant that we choose below.
Our goal in Chapter~\ref{C:ERRORTERMSOBOLEV} 
is to use the results from the previous sections
to prove that for sufficiently small data,
the above bootstrap assumptions 
in fact hold with
$C \mathring{\upepsilon}$ 
in place of $\varepsilon,$
where $\mathring{\upepsilon}$
is the size of the data (see Def.~\ref{D:SMALLDATA}). 

\begin{remark}[\textbf{We do not need the top-order bootstrap assumptions}]
	In our proof of the fundamental Gronwall lemma, 
	Lemma~\ref{L:FUNDAMENTALGRONWALL}, 
	we do not rely on the top-order bootstrap assumptions
	\eqref{E:Q0TOPBOOT}-\eqref{E:Q1TOPBOOT}.
	However, we make them anyway for the sake illustrating
	the expected behavior of 
	$\totzeromax{24}(t,u),$
	$\totonemax{24}(t,u),$
	and $\totMormax{24}(t,u).$
\end{remark}

\section{Statement of the two main propositions and the fundamental Gronwall lemma}
\label{S:STATEMENTOFTWOPROPSANDGRONWALLLEMMA}
We begin by defining four classes of quantities that appear in
the statements of the two main propositions. 
These quantities collectively provide 
suitable upper bounds for the various error integrals
that arise in our analysis.

\begin{remark}[\textbf{The importance of the boxed constants}] \label{R:BOXEDCONSTANTS}
The ``boxed'' constants that appear in the below definitions
(and on the right-hand sides of the top-order estimates of Prop.~\ref{P:MAINTOPORDERENERGYANDFLUXINTEGRALINEQUALITIES})
are the ones that force us to prove estimates for the top-order quantities
$\totzeromax{\leq 24},$ 
$\totonemax{\leq 24},$
and $\totMormax{\leq 24}$
that feature a somewhat large degeneracy with respect to powers of $\upmu_{\star}^{-1}.$ 
The size of the boxed constants is intimately connected to the powers of $\upmu_{\star}^{-1}$ that appear in our
a priori estimates for 
$\totzeromax{\leq N},$ 
$\totonemax{\leq N},$
and $\totMormax{\leq N}.$
The most important feature of the top-order estimates is that 
in the solution regime that we study in this monograph, 
the boxed constants can be chosen to be
\textbf{independent of the number of times the equations are commuted with vectorfields}
$Z \in \mathscr{Z}.$
\end{remark}

\begin{definition}[\textbf{Quantities that bound difficult top-order error integrals generated by the multiplier $\Mult$}] 
\label{D:TOPORDERMULTBOUNDINGQUANTITIES}
We define the classes of terms 
$\topboxedmulterrorone{\leq N}(t,u),$
$\topmulterrortwo{\leq N}(t,u),$
$\cdots,$
and
$\topmulterroreight{\leq N}(t,u)$
to be any functions of $(t,u)$ 
such that there exist 
a constant $0 < \Littleconone < 1/2$ and a constant
$C > 0$ such that 
under the small-data and bootstrap assumptions 
of Sects.~\ref{S:PSISOLVES}-\ref{S:C0BOUNDBOOTSTRAP},
if $\varepsilon$ is sufficiently small,
then the following estimates hold
for $(t,u) \in [0,\Tboot) \times [0, U_0]:$
\begin{subequations}
	\begin{align}
	\topboxedmulterrorone{\leq N}
		& \leq  \boxed{9}
				\int_{t'=0}^t
						\frac{\| [\Lunit \upmu]_- \|_{C^0(\Sigma_{t'}^u)}} 
								 {\upmu_{\star}(t',u)} 
						\totzeromax{\leq N}^{1/2}(t',u) 
						\int_{s=0}^{t'}
							\frac{\| [\Lunit \upmu]_- \|_{C^0(\Sigma_s^u)}} 
									{\upmu_{\star}(s,u)} 
							\totzeromax{\leq N}^{1/2}(s,u) 
						\, ds
				\, dt',
				\label{E:FIRSTHARDMULTTERM} \\
	\topmulterrortwo{\leq N}
		& \leq C \varepsilon
				\int_{t'=0}^t
						\frac{1} 
								 {(1 + t')^{1 + \Littleconone} \upmu_{\star}(t',u)} 
						\totzeromax{\leq N}^{1/2}(t',u) 
						\int_{s=0}^{t'}
							\frac{1}{(1 + s)} 
							\frac{1}{\upmu_{\star}(s,u)} 
							\totzeromax{\leq N}^{1/2}(s,u) 
						\, ds
				\, dt',
				\label{E:TOPORDERMULTTERMSMALLA} \\		
		\topmulterrorthree{\leq N}
		& \leq C \varepsilon
				\int_{t'=0}^t
						\frac{1} 
								 {(1 + t')^{1 + \Littleconone} \upmu_{\star}(t',u)} 
						\totzeromax{\leq N}^{1/2}(t',u) 
						\int_{s=0}^{t'}
							\frac{1}{(1 + s)} 
							\frac{1}{\upmu_{\star}(s,u)} 
							\totonemax{\leq N}^{1/2}(s,u) 
						\, ds
				\, dt',
				\label{E:TOPORDERMULTTERMSECONDSMALLA} \\
		\topmulterrorfour{\leq N}
		& \leq  9
				\int_{t'=0}^t
						\frac{1}{\rgeo(t',u) \left\lbrace 1 + \ln \left(\frac{\rgeo(t',u)}{\rgeo(0,u)} \right) \right\rbrace}
						\totzeromax{\leq N}^{1/2}(t',u) 
						\int_{s=0}^{t'}
							\frac{\| [\Lunit \upmu]_- \|_{C^0(\Sigma_s^u)}} 
									{\upmu_{\star}(s,u)} 
							\totzeromax{\leq N}^{1/2}(s,u) 
						\, ds
				\, dt',
				\label{E:MIXEDTYPEHARDMULTTERM} \\
	\topboxedmulterrorfive{\leq N}
	& \leq \boxed{9} 
				\int_{t'=0}^t
					\frac{\| [\Lunit \upmu]_- \|_{C^0(\Sigma_{t'}^u)}} 
							 {\upmu_{\star}(t',u)} 
				  \totzeromax{\leq N}(t',u)
				\, dt',
				\label{E:SECONDHARDMULTTERM} \\
	\topmulterrorsix{\leq N}
	& \leq  9 
					\int_{t'=0}^t 
						\frac{1}{\rgeo(t',u) \left\lbrace 1 + \ln \left(\frac{\rgeo(t',u)}{\rgeo(0,u)} \right) \right\rbrace} 
						\totzeromax{\leq N}(t',u)
					\, dt',
			\label{E:TOPORDERMULTTERCREATESLOGGROWTH} \\		
	\topmulterrorseven{\leq N}
	& \leq C \varepsilon
				\int_{t'=0}^t
					\frac{1} 
							 {(1 + t')^{1+\Littleconone} \upmu_{\star}(t',u)} 
				  \totzeromax{\leq N}(t',u)
				\, dt',
		\label{E:TOPORDERQ0LARGETIMEDECAYCOUNTERSUPMULOSSTERM} \\
	\topmulterroreight{\leq N}
	& \leq C \varepsilon
				\int_{t'=0}^t
					\frac{1} 
							 {(1 + t')^{1+\Littleconone} \upmu_{\star}(t',u)} 
				  \totonemax{\leq N}(t',u)
				\, dt'.
				\label{E:TOPORDERQ0INVOLVINGQ1LARGETIMEDECAYCOUNTERSUPMULOSSTERM} 
\end{align}
\end{subequations}


\end{definition}

\begin{definition}[\textbf{Quantities that bound easy top-order and below-top-order error integrals generated by the multiplier $\Mult$}] 
\label{D:BELOWORDERMULTBOUNDINGQUANTITIES}
We define the classes of terms 
$\multerrorzero{\leq N+1}(t,u),$
$\multerrorone{\leq N}(t,u),$
$\multerrortwo{\leq N}(t,u),$
$\cdots,$
$\multerrorsix{\leq N}(t,u),$
to be any functions of $(t,u)$ such that there exists 
a constant
$C > 0$ such that 
under the small-data and bootstrap assumptions 
of Sects.~\ref{S:PSISOLVES}-\ref{S:C0BOUNDBOOTSTRAP},
if $\varepsilon$ is sufficiently small,
then the following estimates hold
for $(t,u) \in [0,\Tboot) \times [0, U_0]:$
\begin{subequations}
	\begin{align}
	\multerrorzero{\leq N+1}
		& \leq  
			\underbrace{
				C \varepsilon
				\int_{t'=0}^t
					\frac{1}{(1 + t')^{3/2} \upmu_{\star}^{1/2}(t',u)}
					\totzeromax{\leq N}^{1/2}(t',u)	
					\int_{s=0}^{t'}
						\frac{1}{1 + s}
						\totzeromax{\leq N+1}^{1/2}(s,u)
					\, ds
				\, dt'	
			}_{\mbox{depends on an order $N+1$ quantity}}
	\label{E:MULTLOSSOFONEDERIVERRORINTEGRALBOUND} 
		\\
	& \ \ 
	+	\underbrace{
			C \varepsilon
			\int_{t'=0}^t
				\frac{1}{(1 + t')^{3/2} \upmu_{\star}^{1/2}(t',u)}
				\totzeromax{\leq N}^{1/2}(t',u)	
				\int_{s=0}^{t'}
					\frac{1}{(1 + s)\upmu_{\star}^{1/2}(s,u)}
					\totonemax{\leq N+1}^{1/2}(s,u)
				\, ds
			\, dt'
		}_{\mbox{depends on an order $N+1$ quantity}},
		\notag \\
	\multerrorone{\leq N} 
	& \leq	C 
			\int_{t'=0}^t 
				\frac{1}{(1 + t')^{3/2}} 
				\totzeromax{\leq N}(t',u) 
			\, dt',
			\label{E:BELOWTOPORDEREASYQ0LARGETIMEDECAYTERM} \\
	\multerrortwo{\leq N}
	& \leq  C
				\varepsilon^{1/2}
				\int_{t'=0}^t 
					\frac{\ln(\myexp + t')}{(1 + t')^2 \sqrt{\ln(\myexp + t) - \ln(\myexp + t')}}
					\totonemax{\leq N}(t',u) 
				\, dt',
			 \label{E:STRANGEINTEGRATINGFACTORINTEGRALBOUND}	\\
	\multerrorthree{\leq N} 
			& \leq
					 C
					\int_{t'=0}^t 
						\frac{1}{(1 + t')^{3/2}} \totonemax{\leq N}(t',u) 
					\, dt',
					\label{E:EASYBELOWTOPORDERQ0LARGETIMETERMTHATINVOLVESQ1} \\
	\multerrorfour{\leq N}
	& \leq C 
			\int_{u'=0}^u
				\totzeromax{\leq N}(t,u')		
			\, du',
			 \label{E:EASYBELOWTOPORDERQ0UTERM} \\
	\multerrorfive{\leq N}
	& \leq 	C 
				\frac{1}{(1 + t)^{1/2}} 
				\int_{u'=0}^u
					\totonemax{\leq N}(t,u')		
				\, du',
				\label{E:EASYBELOWTOPORDERQ0UTERMTHATINVOLVESQ1} \\
	\multerrorsix{\leq N}
	& \leq C 
		\sup_{t' \in [0,t]}\frac{\Morint_{(\leq N)}(t',u)}{(1 + t')^{1/2}}.
		\label{E:BELOWTOPORDEREASYMORAWETZINTEGRALTERM}
\end{align}
\end{subequations}

\end{definition}

\begin{definition}[\textbf{Quantities that bound difficult top-order error integrals generated by the multiplier $\Mor$}]
\label{D:TOPORDERMORBOUNDINGQUANTITIES}
We define the classes of terms 
$\topboxedmorerrorone{\leq N}(t,u),$
$\topmorerrortwo{\leq N}(t,u),$
$\cdots,$
and
$\topmorerroreight{\leq N}(t,u)$
to be any functions of $(t,u)$ such that there exists 
a constant
$C > 0$ such that 
under the small-data and bootstrap assumptions 
of Sects.~\ref{S:PSISOLVES}-\ref{S:C0BOUNDBOOTSTRAP},
if $\varepsilon$ is sufficiently small,
then the following estimates hold
for $(t,u) \in [0,\Tboot) \times [0, U_0]:$
\begin{subequations}
	\begin{align}
		\topboxedmorerrorone{\leq N}
		& \leq \boxed{5} 
				\int_{t'=0}^t
					\frac{\| [\Lunit \upmu]_- \|_{C^0(\Sigma_{t'}^u)}} 
							 {\upmu_{\star}(t',u)} 
					\totonemax{\leq N}(t',u)
				\, dt',
				\label{E:FIRSTHARDMORTERM} \\
	\topmorerrortwo{\leq N} 
	& \leq 5 
			\int_{t'=0}^t
				\frac{1}{\rgeo(t',u) \left\lbrace 1 + \ln \left(\frac{\rgeo(t',u)}{\rgeo(0,u)} \right) \right\rbrace}
				\totonemax{\leq N}(t',u)
			\, dt',
			\label{E:TOPORDERMORTERCREATESLOGGROWTH} \\
		\topmorerrorthree{\leq N}
		& \leq C \varepsilon
					 \int_{t'=0}^t
	  			 	 \frac{1}{(1 + t')^{3/2} \upmu_{\star}^{1/2}(t',u)} 
	  			 	 \totzeromax{\leq N}(t',u)
	  			 \, dt',
				\label{E:TOPORDERMORDANGEROUSUPMUBUTLOTSOFTIMEDECAYINTERMSOFQ0} \\
				\topmorerrorfour{\leq N}
		& \leq C \varepsilon
					 \int_{t'=0}^t
	  			 	 \frac{1}{(1 + t')^{3/2} \upmu_{\star}(t',u)} 
	  			 	 \totonemax{\leq N}(t',u)
	  			 \, dt',
	  		\label{E:TOPORDERMORDANGEROUSUPMUBUTLOTSOFTIMEDECAYINTERMSOFQ1} \\
		\topboxedmorerrorfive{\leq N}
		& \leq \boxed{5}
				\frac{\| \Lunit \upmu \|_{C^0(\Sigmaminus{t}{t}{u})}} 
			  	 	 {\upmu_{\star}^{1/2}(t,u)}
				\totonemax{\leq N}^{1/2}(t,u)
				\int_{t' = 0}^t
					\frac{1}{\upmu_{\star}^{1/2}(t',u)}
					\totonemax{\leq N}^{1/2}(t',u)
		 		\, dt',
				\label{E:SECONDHARDMORTERM} \\
	\topmorerrorsix{\leq N}
	& \leq  5
			\left\| 
				\frac{\Lunit \upmu}{\upmu} 
			\right \|_{C^0(\Sigmaplus{t}{t}{u})} 
			\totonemax{\leq N}^{1/2}(t,u)
			\int_{t' = 0}^t
				\left\|
					\sqrt{
					\frac{\upmu(t,\cdot)} 
					     {\upmu}
					     }
				\right\|_{C^0(\Sigmaplus{t'}{t}{u}}
			\totonemax{\leq N}^{1/2}(t',u)
		 	\, dt',
		 	\label{E:ANNOYINGTOPORDERMORTERCREATESLOGGROWTH} \\
	\topmorerrorseven{\leq N}
	& \leq C \varepsilon
			\frac{1}{\upmu_{\star}^{1/2}(t,u)}
			\totonemax{\leq N}^{1/2}(t,u)	
			\int_{t'=0}^t
				\frac{1}{(1 + t')^{3/2}} \totzeromax{\leq N}^{1/2}(t',u)
			\, dt',
			\label{E:TOPORDERMORHYPERSURFACEDANGEROUSUPMUBUTLOTSOFTIMEDECAYINTERMSOFQ0} 	\\
	\topmorerroreight{\leq N}
	& \leq C \varepsilon
				\frac{1}{\upmu_{\star}^{1/2}(t,u)}
				\totonemax{\leq N}^{1/2}(t,u)
				\int_{t'=0}^t
					\frac{1}{(1 + t')^{3/2} \upmu_{\star}^{1/2}(t',u)} \totonemax{\leq N}^{1/2}(t',u)
				\, dt'.
				\label{E:TOPORDERMORHYPERSURFACEDANGEROUSUPMUBUTLOTSOFTIMEDECAYINTERMSOFQ1}
\end{align}
\end{subequations}


\end{definition}

\begin{definition}[\textbf{Quantities that bound easy top-order and below-top-order error integrals generated by the multiplier 
$\Mor$}] 
\label{D:BELOWORDERMORBOUNDINGQUANTITIES}
We define the classes of terms 
$\morerrorzero{\leq N+1}(t,u),$
$\morerrorone{\leq N}(t,u),$
$\morerrortwo{\leq N}(t,u),$
$\cdots, 
\morerrorfive{\leq N}(t,u),$
to be any functions of $(t,u)$ such that there exists 
a constant
$C > 0$ such that 
under the small-data and bootstrap assumptions 
of Sects.~\ref{S:PSISOLVES}-\ref{S:C0BOUNDBOOTSTRAP},
if $\varepsilon$ is sufficiently small,
then the following estimates hold
for $(t,u) \in [0,\Tboot) \times [0,U_0]:$
	\begin{subequations}
	\begin{align}
	\morerrorzero{\leq N+1}
		& \leq
			\underbrace{
				C \varepsilon
				\int_{t'=0}^t
					\frac{1}{(1 + t')^2}
					\left(
					\int_{s=0}^{t'}
						\frac{ \totzeromax{\leq N+1}^{1/2}(s,u)}{1 + s} 
					\, ds
					\right)^2
				\, dt'
			}_{\mbox{depends on an order $N+1$ quantity}}
	 		\label{E:MORLOSSOFONEDERIVERRORINTEGRALTERM}  \\
	& \ \ 
			\underbrace{
			+ C \varepsilon
				\int_{t'=0}^t
					\frac{1}{(1 + t')^2}
					\left(
					\int_{s=0}^{t'}
						\frac{\totonemax{\leq N+1}^{1/2}(s,u)}{(1 + s)\upmu_{\star}^{1/2}(s,u)} 
					\, ds
				\right)^2
				\, dt'	
					}_{\mbox{depends on an order $N+1$ quantity}},
					\notag \\
	 \morerrorone{\leq N} 
	& \leq C \ln^4(\myexp + t) \totzeromax{\leq N}(t,u),
			\label{E:EASYBELOWTOPORDERQ1THATINVOLVESQ0TERM} \\
	\morerrortwo{\leq N}
	& \leq 2 
					\int_{t'=0}^t 
						\frac{1}{\rgeo(t',u) \left\lbrace 1 + \ln \left(\frac{\rgeo(t',u)}{\rgeo(0,u)} \right) \right\rbrace} 
						\totonemax{\leq N}(t',u) 
					\, dt',
		\label{E:EASYBELOWTOPORDERQ1LOGLOSSTERM} \\
	\morerrorthree{\leq N}
	& \leq C 
			\int_{t'=0}^t
				\frac{1}{(1 + t')^{3/2}}
				\totonemax{\leq N}(t',u)		
			\, dt',
			\label{E:EASYBELOWTOPORDERLARGETIMEDECAYTERM} \\
	\morerrorfour{\leq N}			
	& \leq C
			\int_{u'=0}^u
				\totonemax{\leq N}(t,u')		
			\, du',
		\label{E:EASYBELOWTOPORDERQ1UTERM}
			\\
	\morerrorfive{\leq N}			
	& \leq C
		\Morint_{(\leq N)}(t,u).
		\label{E:MORAWETZINTEGRALQ1TERM}
	\end{align}
\end{subequations}
\end{definition}

We now state the two propositions. 
We provide their proofs in 
Sects.~\ref{S:PROOFOFPROPOSITIONMAINTOPORDERENERGYANDFLUXINTEGRALINEQUALITIES}
and \ref{S:PROOFOFPROPOSITIONMAINBELOWTOPORDERENERGYANDFLUXINTEGRALINEQUALITIES}.

\begin{proposition} [\textbf{The main top-order energy-flux integral inequalities}]
\label{P:MAINTOPORDERENERGYANDFLUXINTEGRALINEQUALITIES}
	Assume that $\square_{g(\Psi)} \Psi = 0$
	and consider the quantities defined in 
	Def.~\ref{D:TOPORDERMULTBOUNDINGQUANTITIES}-Def.~\ref{D:BELOWORDERMORBOUNDINGQUANTITIES}.
	Let $0 \leq N \leq 24$ be an integer.
	There exists a constant $C > 0$ such that
	under the small-data and bootstrap assumptions 
	of Sects.~\ref{S:PSISOLVES}-\ref{S:C0BOUNDBOOTSTRAP},
	if $\varepsilon$ is sufficiently small
	and $\varsigma > 0,$ $\widetilde{\varsigma} > 0$ are any real numbers,
	then the following inequalities hold for $(t,u) \in [0,\Tboot) \times [0,U_0]:$
\begin{subequations}	
\begin{align}
	\totzeromax{\leq N}(t,u)
	& \leq \totzeromax{\leq N}(0,u)
		+ C \varepsilon^3 \upmu_{\star}^{-1}(t,u)
		+ C \varepsilon \totzeromax{\leq N}(t,u)
			\label{E:Q0TOPORDERGRONWALLREADYINEQUALITY} \\
	& \ \ 
		+ C \varepsilon \upmu_{\star}^{-1}(t,u) \ln^8(\myexp + t) \totzeromax{\leq N-1}(t,u)
		+ C \varepsilon \upmu_{\star}^{-1}(t,u) \totonemax{\leq N-1}(t,u)
			\notag  \\
	& \ \ + \topboxedmulterrorone{\leq N}(t,u)
		+ \topmulterrortwo{\leq N}(t,u)
		+ \cdots
		+ \topmulterroreight{\leq N}(t,u)
			\notag
			\\
	& \ \ 
		+ \left\lbrace 1 + \varsigma^{-1} \right\rbrace \multerrorone{\leq N}(t,u)
		+ \multerrortwo{\leq N}(t,u)
		+ \left\lbrace 1 + \varsigma^{-1} \right\rbrace \multerrorthree{\leq N}(t,u)
			\notag \\
	& \ \ 
		+ \multerrorfour{\leq N}(t,u)
		+ \left\lbrace 1 + \varsigma^{-1} \right\rbrace \multerrorfive{\leq N}(t,u)
		+ \left\lbrace \varepsilon + \varsigma \right\rbrace \multerrorsix{\leq N}(t,u),
			 \notag \\
	\totonemax{\leq N}(t,u)
		+ \totMormax{\leq N}(t,u)
	& \leq 
			C \totzeromax{\leq N}(0,u)
		+ C \varepsilon^3 \upmu_{\star}^{-1}(t,u)
		+ C \varepsilon \totzeromax{\leq N}(t,u)
		+ C \varepsilon \totonemax{\leq N}(t,u)
		\label{E:Q1TOPORDERGRONWALLREADYINEQUALITY}  \\
	& \ \ 
		+ C \varepsilon \upmu_{\star}^{-1}(t,u) \ln^2(\myexp + t) \totzeromax{\leq N-1}(t,u)
		+ C \varepsilon \upmu_{\star}^{-1}(t,u) \ln^2(\myexp + t) \totonemax{\leq N-1}(t,u)
			\notag  \\
	& \ \ + \topboxedmorerrorone{\leq N}(t,u)
		+ \topmorerrortwo{\leq N}(t,u)
		+ \cdots 
		+ \topmorerroreight{\leq N}(t,u)
		\notag \\
	& \ \  
		+ \left\lbrace 1 + \widetilde{\varsigma}^{-1} \right\rbrace \morerrorone{\leq N}(t,u)
		+ \morerrortwo{\leq N}(t,u)
		+ \left\lbrace 1 + \widetilde{\varsigma}^{-1} \right\rbrace \morerrorthree{\leq N}(t,u)
			\notag \\
	& \ \ 
		+ \left\lbrace 1 + \widetilde{\varsigma}^{-1} \right\rbrace \morerrorfour{\leq N}(t,u)
		+ \left\lbrace \varepsilon + \widetilde{\varsigma} \right\rbrace \morerrorfive{\leq N}(t,u).
		\notag
\end{align}
\end{subequations}

\end{proposition}

\begin{proposition} [\textbf{The main below-top-order energy-flux integral inequalities}]
\label{P:MAINBELOWTOPORDERENERGYANDFLUXINTEGRALINEQUALITIES}
	Assume that $\square_{g(\Psi)} \Psi = 0$
	and let $0 \leq N \leq 23$ be an integer.
	Consider the quantities defined in
	Def.~\ref{D:TOPORDERMULTBOUNDINGQUANTITIES}-Def.~\ref{D:BELOWORDERMORBOUNDINGQUANTITIES}.
	There exists a constant $C > 0$ such that
	under the small-data and bootstrap assumptions 
	of Sects.~\ref{S:PSISOLVES}-\ref{S:C0BOUNDBOOTSTRAP},
	if $\varepsilon$ is sufficiently small
	and $\varsigma > 0,$ $\widetilde{\varsigma} > 0$ are any real numbers,
	then the following inequalities hold for $(t,u) \in [0,\Tboot) \times [0,U_0]:$
	\begin{subequations}	
\begin{align}
	\totzeromax{\leq N}(t,u)
	& \leq 
			\totzeromax{\leq N}(0,u)
		+ C \varepsilon^3 
			\label{E:Q0BELOWTOPORDERGRONWALLREADYINEQUALITY} \\
	& \ \ + \mathbf{0})_{\leq N+1}(t,u)
			\notag \\
& \ \ 
		+ \left\lbrace 1 + \varsigma^{-1} \right\rbrace \multerrorone{\leq N}(t,u)
		+ \multerrortwo{\leq N}(t,u)
		+ \left\lbrace 1 + \varsigma^{-1} \right\rbrace \multerrorthree{\leq N}(t,u)
			\notag \\
& \ \ 
		+ \multerrorfour{\leq N}(t,u)
		+ \left\lbrace 1 + \varsigma^{-1} \right\rbrace \multerrorfive{\leq N}(t,u)
		+ \left\lbrace \varepsilon + \varsigma \right\rbrace \multerrorsix{\leq N}(t,u),
			 \notag \\
	\totonemax{\leq N}(t,u)
	+ \totMormax{\leq N}(t,u)
	& \leq  
			\totonemax{\leq N}(0,u)
		+ C \varepsilon^3
		\label{E:Q1BELOWTOPORDERGRONWALLREADYINEQUALITY} \\
	& \ \ + \mathbf{\widetilde{0}})_{\leq N+1}(t,u)
		\notag \\
	& \ \  
		+ \left\lbrace 1 + \widetilde{\varsigma}^{-1}\right\rbrace \morerrorone{\leq N}(t,u)
		+ \morerrortwo{\leq N}(t,u)
		+ \left\lbrace 1 + \widetilde{\varsigma}^{-1}\right\rbrace \morerrorthree{\leq N}(t,u)
			\notag \\
	& \ \ 
		+ \left\lbrace 1 + \widetilde{\varsigma}^{-1}\right\rbrace \morerrorfour{\leq N}(t,u)
		+ \left\lbrace \varepsilon + \widetilde{\varsigma} \right\rbrace \morerrorfive{\leq N}(t,u).
		\notag
\end{align}
\end{subequations}

\end{proposition}

Before proving Prop.~\ref{P:MAINTOPORDERENERGYANDFLUXINTEGRALINEQUALITIES}
and Prop.~\ref{P:MAINBELOWTOPORDERENERGYANDFLUXINTEGRALINEQUALITIES},
we first state our main Gronwall-type lemma. The lemma 
allows us to derive suitable a priori estimates for 
$\totzeromax{\leq N}(t,u),$ 
$\totonemax{\leq N}(t,u),$
and $\totMormax{\leq N}(t,u).$ 
The previous two propositions imply that 
the hypotheses of the Gronwall lemma are satisfied.
We provide the proof of the lemma in Sect.~\ref{S:PROOFOFLEMMAFUNDAMENTALGRONWALL}.
Its proof is lengthy because of the sheer number of terms involved.
However, it is not difficult; 
we all already established the difficult estimates
in Prop.~\ref{P:MUINVERSEINTEGRALESTIMATES}.

\begin{lemma}[\textbf{The fundamental Gronwall lemma}] 
\label{L:FUNDAMENTALGRONWALL}
Let $\mathring{\upepsilon}$ be the size of the data as defined in Def.~\ref{D:SMALLDATA}.
Assume that on the domain $(t,u) \in [0,\Tboot) \times [0,U_0],$
the quantities $\totzeromax{\leq N}(t,u),$
$\totonemax{\leq N}(t,u),$
$\totMormax{\leq N}(t,u)$
verify the estimates of Prop.~\ref{P:MAINTOPORDERENERGYANDFLUXINTEGRALINEQUALITIES}
for $0 \leq N \leq 24$
and the estimates of Prop.~\ref{P:MAINBELOWTOPORDERENERGYANDFLUXINTEGRALINEQUALITIES}
for $0 \leq N \leq 23.$ 
Assume in addition that the bootstrap assumptions \eqref{E:Q0LOWBOOT}-\eqref{E:Q1TOPBOOT} hold on the same domain.
Then if $\Cononestar > 4$ is sufficiently large, 
there exists a large constant $C > 0$
(depending on $\Cononestar$) such that
if $\varepsilon$ is sufficiently small,
then the following estimates hold for $(t,u) \in [0,\Tboot) \times [0,U_0]:$
\begin{subequations}
\begin{align}
	\totzeromax{N}^{1/2}(t,u)
	& \leq C 
				\left\lbrace
					\mathring{\upepsilon}
					+ \varepsilon^{3/2}
				\right\rbrace, 
		&& (0 \leq N \leq 15), 
		\label{E:Q0LOWESTORDERKEYBOUND} \\  
	\totonemax{N}^{1/2}(t,u)
		+ \totMormax{N}^{1/2}(t,u)
	& \leq C 
				\left\lbrace
					\mathring{\upepsilon}
					+ \varepsilon^{3/2}
				\right\rbrace
			\ln^2(\myexp + t), 
		&& (0 \leq N \leq 15),
		\label{E:Q1LOWESTORDERKEYBOUND} \\
	\totzeromax{16 + M}^{1/2}(t,u)
	& \leq C 
				\left\lbrace
					\mathring{\upepsilon}
					+ \varepsilon^{3/2}
				\right\rbrace 
				\upmu_{\star}^{- .75 - M}(t,u), 
		&& (0 \leq M \leq 7), 
		\label{E:Q0MIDORDERKEYBOUND} \\  
	\totonemax{16 + M}^{1/2}(t,u)
		+ \totMormax{16 + M}^{1/2}(t,u)
	& \leq C 
				\left\lbrace
					\mathring{\upepsilon}
					+ \varepsilon^{3/2}
				\right\rbrace 
			\ln^2(\myexp + t)
			\upmu_{\star}^{-.75 - M}(t,u), 
		&& (0 \leq M \leq 7),
		\label{E:Q1MIDORDERKEYBOUND} \\
		\totzeromax{24}^{1/2}(t,u)
	& \leq C 
				\left\lbrace
					\mathring{\upepsilon}
					+ \varepsilon^{3/2}
				\right\rbrace 
				\ln^{\Cononestar}(\myexp + t)
			\upmu_{\star}^{-8.75}(t,u), 
		&&  
		\label{E:Q0TOPORDERKEYBOUND} \\  
	\totonemax{24}^{1/2}(t,u)
		+ \totMormax{24}^{1/2}(t,u)
	& \leq C 
				\left\lbrace
					\mathring{\upepsilon}
					+ \varepsilon^{3/2}
				\right\rbrace
				\ln^{\Cononestar + 2}(\myexp + t) 
				\upmu_{\star}^{-8.75}(t,u).
				\label{E:Q1TOPORDERKEYBOUND}
\end{align}
\end{subequations}
\end{lemma}

\section{Estimates for all but the most difficult error integrals}
\label{S:EASYERRORINTEGRALESTIMATES}
In our proofs of Prop.~\ref{P:MAINTOPORDERENERGYANDFLUXINTEGRALINEQUALITIES}
and Prop.~\ref{P:MAINBELOWTOPORDERENERGYANDFLUXINTEGRALINEQUALITIES},
we encounter many quadratic spacetime error integrals that we must bound
in terms of the $L^2-$type quantities
$\totzeromax{\leq N},$ 
$\totonemax{\leq N},$ 
and $\totMormax{\leq N}.$
In this section, we derive estimates
for all such error integrals except the most difficult top-order ones.

We begin with the next three lemmas,
in which we derive estimates for 
most of the error integrals appearing our analysis, 
including the ones
corresponding to the $Harmless^{\leq N}$ inhomogeneous terms
(see Def.~\ref{D:HARMLESSTERMS})
appearing in the commuted wave equation.

\begin{lemma}[\textbf{Bounds for spacetime integrals in terms of the fundamental $L^2-$controlling quantities}]
\label{L:SPACETIMEINTEGRALBOUNDSINTERMSOFQ0ANDQ1}
Let $\Conone$ be a constant. Under the small-data and bootstrap assumptions 
of Sects.~\ref{S:PSISOLVES}-\ref{S:C0BOUNDBOOTSTRAP},
if $\varepsilon$ is sufficiently small, 
then the following inequalities hold for $(t,u) \in [0,\Tboot) \times [0,U_0]:$
	\begin{subequations}
	\begin{align}  \label{E:STANDARDPSISPACETIMEINTEGRALSWITHNOSMALLMU}
		\int_{\mathcal{M}_{t,u}}
		 	\frac{\ln^{\Conone}(\myexp + t')}{(1 + t')^2}
		 	& \left|
				\fourmyarray[\rgeo (1 + \upmu) \Lunit \mathscr{Z}^{\leq N} \Psi]
					{\Rad \mathscr{Z}^{\leq N} \Psi}
					{\rgeo (\sqrt{\upmu} + \upmu) \angdiff \mathscr{Z}^{\leq N} \Psi}
					{\mathscr{Z}^{\leq N} \Psi}
				\right|
			\left|
				\myarray[(1 + \upmu)\Lunit \mathscr{Z}^{\leq N} \Psi]
					{\Rad \mathscr{Z}^{\leq N} \Psi}
			\right|
		 \, d \vol
		 	\\
		 & \overset{\Conone}{\lesssim}
		 	\int_{t'=0}^t
				\frac{1}{(1 + t')^{3/2}}
				\totzeromax{\leq N}(t',u)
			\, dt'
			+ 
		 	\int_{t'=0}^t
				\frac{1}{(1 + t')^{3/2}}
				\totonemax{\leq N}(t',u)
			\, dt'
				\notag \\
		& \ \
			+
			\frac{1}{(1 + t)^{1/2}}
			\int_{u'=0}^u
				\totzeromax{\leq N}(t,u')
			\, du',
			\notag \\
	\int_{\mathcal{M}_{t,u}}
		 	\frac{\ln^{\Conone}(\myexp + t)}{(1 + t')^2}
				& \left| 
						\sevenmyarray
							[\Lunit \mathscr{Z}^{\leq N} \upmu]
							{\rgeo \sum_{a=1}^3 |\Lunit \mathscr{Z}^{\leq N} \Lunit_{(Small)}^a|}
							{\rgeo^2 \angLie_{\Lunit} \angLie_{\mathscr{Z}}^{\leq N-1} \upchi^{(Small)}}
							{\angLie_{\Lunit} (\rgeo^2 \angLie_{\mathscr{Z}}^{\leq N-1} \upchi^{(Small)\#})}
							{\Lunit (\rgeo^2 \mathscr{Z}^{\leq N-1} \mytr \upchi^{(Small)})}
							{\rgeo^2 \angLie_{\Lunit} \angLie_{\mathscr{Z}}^{\leq N-1} \hat{\upchi}^{(Small)}}
							{\angLie_{\Lunit} (\rgeo^2 \angLie_{\mathscr{Z}}^{\leq N-1} \hat{\upchi}^{(Small)\#})}
				\right|
				\left|
					\myarray[(1 + \upmu)\Lunit \mathscr{Z}^{\leq N} \Psi]
						{\Rad \mathscr{Z}^{\leq N} \Psi}
				\right|
		 \, d \vol
		   \label{E:STANDARDPSISPACETIMEINTEGRALSINVOLVINGLDERIVATESOFEIKONALWITHNOSMALLMU} \\
		& \overset{\Conone}{\lesssim}
		 	\int_{t'=0}^t
				\frac{1}{(1 + t')^{3/2}}
				\totzeromax{\leq N}(t',u)
			\, dt'
			+ 
			\int_{t'=0}^t
				\frac{1}{(1 + t')^{3/2}}
				\totonemax{\leq N}(t',u)
			\, dt'
			\notag	\\
		& \ \
			+
			\int_{u'=0}^u
				\totzeromax{\leq N}(t,u')
			\, du'
			+ \varepsilon^2,
			\notag 		
			\\
	\int_{\mathcal{M}_{t,u}}
		 	\frac{\ln^{\Conone}(\myexp + t)}{(1 + t')^3}
				& 	
				\left| 
					\sevenmyarray
						[\mathscr{Z}^{\leq N} (\upmu - 1)]
						{\rgeo \sum_{a=1}^3 |\mathscr{Z}^{\leq N} \Lunit_{(Small)}^a|}
						{\rgeo^2 \angLie_{\mathscr{Z}}^{\leq N-1} \upchi^{(Small)}}
						{\rgeo^2 \angLie_{\mathscr{Z}}^{\leq N-1} \upchi^{(Small)\#}}
						{\rgeo^2 \mathscr{Z}^{\leq N-1} \mytr \upchi^{(Small)}}
						{\rgeo^2 \angLie_{\mathscr{Z}}^{\leq N-1} \hat{\upchi}^{(Small)}}
						{\rgeo^2 \angLie_{\mathscr{Z}}^{\leq N-1} \hat{\upchi}^{(Small)\#}}
				\right|
				\left|
					\myarray[(1 + \upmu)\Lunit \mathscr{Z}^{\leq N} \Psi]
						{\Rad \mathscr{Z}^{\leq N} \Psi}
				\right|
		 \, d \vol
		    \label{E:STANDARDPSISPACETIMEINTEGRALSTIMESEASYEIKONALFUNCTIONQUANTITIES} \\
		& \overset{\Conone}{\lesssim}
		 	\int_{t'=0}^t
				\frac{1}{(1 + t')^{3/2}}
				\totzeromax{\leq N}(t',u)
			\, dt'
			+ 
			\int_{t'=0}^t
				\frac{1}{(1 + t')^{3/2}}
				\totonemax{\leq N}(t',u)
			\, dt'
			\notag	\\
		& \ \
			+
			\int_{u'=0}^u
				\totzeromax{\leq N}(t,u')
			\, du'
			+ \varepsilon^2.
			\notag
	\end{align}
	\end{subequations}
\end{lemma}

\begin{remark}[\textbf{We need the Morawetz integral to control non $\upmu-$weighted angular derivatives of $\Psi$}]
Note that Lemma~\ref{L:SPACETIMEINTEGRALBOUNDSINTERMSOFQ0ANDQ1}
does not provide control of integrals involving the non $\upmu-$weighted quantity 
$|\angdiff \mathscr{Z}^{\leq N} \Psi|.$
We obtain control of this term in Lemma~\ref{L:MORAWETZSPACETIMEINTEGRALESTIMATES}
with the help of the Morawetz spacetime integral of Def.~\ref{D:COERCIVEMORAWETZINTEGRAL}.
\end{remark}

\begin{proof}
The proof is essentially a tedious exercise in Cauchy-Schwarz.
Since the estimates are all very similar, we bound 
only one representative product from \eqref{E:STANDARDPSISPACETIMEINTEGRALSWITHNOSMALLMU}.
Specifically, we use spacetime Cauchy-Schwarz
and Prop.~\ref{P:L2NORMSOFPSIINTERMSOFTHECOERCIVEQUANTITIES} to deduce that
\begin{align} \label{E:REPTERMSTANDARDPSISPACETIMEINTEGRALSWITHNOSMALLMU}
	&
	\left|
	\int_{\mathcal{M}_{t,u}}
		 	\frac{\ln^{\Conone}(\myexp + t')}{(1 + t')^2}
		 	\Big[\rgeo \sqrt{\upmu} \angdiff \mathscr{Z}^{\leq N} \Psi \Big]
		 	\Lunit \mathscr{Z}^{\leq N} \Psi
	\, d \vol
	\right|
		\\
	& \lesssim 
	\int_{\mathcal{M}_{t,u}}
		\frac{\ln^{\Conone}(\myexp + t')}{(1 + t')^2}
		\rgeo^2 \upmu |\angdiff \mathscr{Z}^{\leq N} \Psi|^2
	\, d \vol
	+ 
	\int_{\mathcal{M}_{t,u}}
		\frac{\ln^{\Conone}(\myexp + t')}{(1 + t')^2}
		 |\Lunit \mathscr{Z}^{\leq N} \Psi|^2
	\, d \vol
		\notag
		\\
	& \lesssim 
	\int_{t'=0}^t 
		\frac{\ln^{\Conone}(\myexp + t')}{(1 + t')^2}
		\totonemax{\leq N}(t',u)
	\, dt'	
	+ \frac{\ln^{\Conone}(\myexp + t)}{(1 + t)^2}
		\int_{u'=0}^u 
			\int_{t'=0}^t \int_{S_{t',u'}}
			|\Lunit \mathscr{Z}^{\leq N} \Psi|^2
		\, d \spherevol \, dt' \, du'
		\notag \\
	& \ \ 
		+ \int_{t'=0}^t
			\int_{u'=0}^u 
			\int_{s=0}^{t'}
			\int_{S_{s,u'}}
			|\Lunit \mathscr{Z}^{\leq N} \Psi|^2
			\, d \argspherevol{(s,u',\vartheta)} 
			\, ds
			\, du' 
			\left|
				\frac{d}{dt'}
				\frac{\ln^{\Conone}(\myexp + t')}{(1 + t')^2}
			\right|
			\, dt'
		\notag 
		\\
	& \lesssim 
		\int_{t'=0}^t 
			\frac{\ln^{\Conone}(\myexp + t')}{(1 + t')^2}
			\totonemax{\leq N}(t',u)
		\, dt'
		+ \frac{\ln^{\Conone}(\myexp + t)}{(1 + t)^2}
			\int_{u'=0}^u 
				\totzeromax{\leq N}(t,u')
		\, du'
			\notag \\
	& \ \
		+ \int_{t'=0}^t 
				\left| \frac{d}{dt'} \frac{\ln^{\Conone}(\myexp + t')}{(1 + t')^2} \right|
				\totzeromax{\leq N}(t',u)
		\, dt',
		\notag
\end{align}	
where to conclude the next-to-last inequality, we integrated by parts in $t',$
and to conclude the last inequality, we used the fact that
$\totzeromax{\leq N}$ and $\totonemax{\leq N}$ are increasing in their arguments
and the fact that $u \leq U_0.$ Clearly, the right-hand side of \eqref{E:REPTERMSTANDARDPSISPACETIMEINTEGRALSWITHNOSMALLMU}
is bounded by the right-hand side of \eqref{E:STANDARDPSISPACETIMEINTEGRALSWITHNOSMALLMU} as desired.

The proof of \eqref{E:STANDARDPSISPACETIMEINTEGRALSTIMESEASYEIKONALFUNCTIONQUANTITIES} is similar, but
after using Cauchy-Schwarz and suitably distributing the $t'$ weights, we have to estimate integrals such as
\begin{align} \label{E:REPEASYEIKONALINTEGRAL}
	\int_{\mathcal{M}_{t,u}}
		\frac{1}{(1 + t')^4}
		 |\mathscr{Z}^{\leq N} (\upmu - 1)|^2
	\, d \vol
	& = 
	\int_{t'=0}^t 
		\frac{1}{(1 + t')^4}
	\int_{\Sigma_{t'}^u}
		|\mathscr{Z}^{\leq N} (\upmu - 1)|^2
	\, d \tvol \, dt'.
\end{align}
We use \eqref{E:EIKONALFUNCTIONQUANTITIESL2BOUNDSINTERMSOFQ0ANDQ1} to deduce that the right-hand side of \eqref{E:REPEASYEIKONALINTEGRAL} is
\begin{align}
	& 
	\leq
	\int_{t'=0}^t 
		\frac{1}{(1 + t')^2}
		\varepsilon^2
	\, dt'
	+ 
	\int_{t'=0}^t 
		\frac{1}{(1 + t')^2}
		\left(
		\int_{s=0}^{t'} 
			\frac{\totzeromax{\leq N}^{1/2}(s,u)}{1 + s} 
		\, ds
		\right)^2
	\, dt'
	\\
	& \ \
	+ \int_{t'=0}^t 
		\frac{1}{(1 + t')^2}
		\left(
		\int_{s=0}^{t'} 
			\frac{\totonemax{\leq N}^{1/2}(s,u)}{(1 + s)\upmu_{\star}^{1/2}(s,u)} 
		\, ds
		\right)^2
		\, dt'
		\notag \\
	& \lesssim
	\varepsilon^2
	+
		\int_{t'=0}^t 
			\frac{1}{(1 + t')^{3/2}}
				\totzeromax{\leq N}(t',u)
		\, dt'
	+ \int_{t'=0}^t 
			\frac{1}{(1 + t')^{3/2}}
				\totonemax{\leq N}(t',u)
		\, dt',
		\notag
\end{align}
where in the last step, we used the estimate \eqref{E:LOGLOSSLESSSINGULARTERMSMTHREEFOURTHSINTEGRALBOUND} 
to destroy the factor $\upmu_{\star}^{1/2}(s,u)$ in the denominator.

The proof of \eqref{E:STANDARDPSISPACETIMEINTEGRALSINVOLVINGLDERIVATESOFEIKONALWITHNOSMALLMU} is similar, thanks
to the estimate \eqref{E:LDERIVATIVEEIKONALFUNCTIONQUANTITIESL2BOUNDSINTERMSOFQ0ANDQ1}.

\end{proof}

\begin{lemma}[\textbf{Bounds for spacetime integrals involving} 
$\lbrace \Lunit + \frac{1}{2} \mytr \upchi \rbrace \mathscr{Z}^{\leq N} \Psi$
\textbf{in terms of the fundamental $L^2-$controlling quantities}]
	\label{L:SPACETIMEINTEGRALSWITHLPSIPLUSHALFTRCHIPSI}
	Let $\Conone$ be a constant.
	Under the small-data and bootstrap assumptions 
	of Sects.~\ref{S:PSISOLVES}-\ref{S:C0BOUNDBOOTSTRAP},
	if $\varepsilon$ is sufficiently small, 
	then the following inequalities hold for $(t,u) \in [0,\Tboot) \times [0,U_0]:$
	\begin{subequations}
	\begin{align} \label{E:RENORMALIZEDLTIMESFRAMEDERIVATIVESINTERMSOFQ0ANDQ1}
		\int_{\mathcal{M}_{t,u}}
		 	\ln^{\Conone}(\myexp + t') 
		 	& \left|
		 		\left\lbrace
		 			\Lunit 
		 			+ \frac{1}{2} \mytr \upchi
		 		\right\rbrace	
		 		\mathscr{Z}^{\leq N} \Psi
		 	\right|
		 \left|
				\fourmyarray[\rgeo (1 + \upmu) \Lunit \mathscr{Z}^{\leq N} \Psi]
					{\Rad \mathscr{Z}^{\leq N} \Psi}
					{\rgeo (\sqrt{\upmu} + \upmu) \angdiff \mathscr{Z}^{\leq N} \Psi}
					{\mathscr{Z}^{\leq N} \Psi}
			\right|
		 \, d \vol
		 \\
		& \overset{\Conone}{\lesssim}
		 	\int_{t'=0}^t
				\frac{1}{(1 + t')^{3/2}}
				\totzeromax{\leq N}(t',u)
			\, dt'
			+ 
			\int_{t'=0}^t
				\frac{1}{(1 + t')^{3/2}}
				\totonemax{\leq N}(t',u)
			\, dt'
			\notag	\\
		& \ \
			+
			\frac{1}{(1 + t)^{1/2}}
			\int_{u'=0}^u
				\totonemax{\leq N}(t,u')
			\, du',
			\notag
			\\
	\int_{\mathcal{M}_{t,u}}
		 	\rgeo^2 
		 	& \left|
		 		\left\lbrace
		 			\Lunit 
		 			+ \frac{1}{2} \mytr \upchi
		 		\right\rbrace	
		 		\mathscr{Z}^{\leq N} \Psi
		 	\right|^2
		 \, d \vol
		 \label{E:RENORMALIZEDLSQUAREDINTERMSOFQ0ANDQ1} \\
		 & \lesssim
		 	\int_{u'=0}^u
				\totonemax{\leq N}(t,u')
				\, du',
				\notag
	\end{align}

\begin{align}  \label{E:RENORMALIZEDLTIMESLDERIVATIVEOFEIKONALFUNCTIONINTERMSOFQ0ANDQ1}
	\int_{\mathcal{M}_{t,u}}
		 	\ln^{\Conone}(\myexp + t')
		 	& \left|
		 		\left\lbrace
		 			\Lunit 
		 			+ \frac{1}{2} \mytr \upchi
		 		\right\rbrace	
		 		\mathscr{Z}^{\leq N} \Psi
		 	\right|
		 	\left| 
						\sevenmyarray
							[\Lunit \mathscr{Z}^{\leq N} \upmu]
							{\rgeo \sum_{a=1}^3 |\Lunit \mathscr{Z}^{\leq N} \Lunit_{(Small)}^a|}
							{\rgeo^2 \angLie_{\Lunit} \angLie_{\mathscr{Z}}^{\leq N-1} \upchi^{(Small)}}
							{\angLie_{\Lunit} (\rgeo^2 \angLie_{\mathscr{Z}}^{\leq N-1} \upchi^{(Small)\#})}
							{\Lunit (\rgeo^2 \mathscr{Z}^{\leq N-1} \mytr \upchi^{(Small)})}
							{\rgeo^2 \angLie_{\Lunit} \angLie_{\mathscr{Z}}^{\leq N-1} \hat{\upchi}^{(Small)}}
							{\angLie_{\Lunit} (\rgeo^2 \angLie_{\mathscr{Z}}^{\leq N-1} \hat{\upchi}^{(Small)\#})}
			\right|
		 \, d \vol
		 \\
		& \overset{\Conone}{\lesssim}
		 	\int_{t'=0}^t
				\frac{1}{(1 + t')^{3/2}}
				\totzeromax{\leq N}(t',u)
			\, dt'
			+ 
			\int_{t'=0}^t
				\frac{1}{(1 + t')^{3/2}}
				\totonemax{\leq N}(t',u)
			\, dt'
			\notag	\\
		& \ \
			+
			\int_{u'=0}^u
				\totonemax{\leq N}(t,u')
			\, du'
			+ \varepsilon^2,
			\notag
			\\
	\int_{\mathcal{M}_{t,u}}
		 	\frac{\ln^{\Conone}(\myexp + t')}{1+t'}
		 	& \left|
		 		\left\lbrace
		 			\Lunit 
		 			+ \frac{1}{2} \mytr \upchi
		 		\right\rbrace	
		 		\mathscr{Z}^{\leq N} \Psi
		 	\right|
		 	\left| 
					\sevenmyarray
						[\mathscr{Z}^{\leq N} (\upmu - 1)]
						{\rgeo \sum_{a=1}^3 |\mathscr{Z}^{\leq N} \Lunit_{(Small)}^a|}
						{\rgeo^2 \angLie_{\mathscr{Z}}^{\leq N-1} \upchi^{(Small)}}
						{\rgeo^2 \angLie_{\mathscr{Z}}^{\leq N-1} \upchi^{(Small)\#}}
						{\rgeo^2 \mathscr{Z}^{\leq N-1} \mytr \upchi^{(Small)}}
						{\rgeo^2 \angLie_{\mathscr{Z}}^{\leq N-1} \hat{\upchi}^{(Small)}}
						{\rgeo^2 \angLie_{\mathscr{Z}}^{\leq N-1} \hat{\upchi}^{(Small)\#}}
				\right|
		 \, d \vol
		 \label{E:RENORMALIZEDLTIMESEIKONALFUNCTIONINTERMSOFQ0ANDQ1} \\
		& \overset{\Conone}{\lesssim}
		 	\int_{t'=0}^t
				\frac{1}{(1 + t')^{3/2}}
				\totzeromax{\leq N}(t',u)
			\, dt' 
		+ 
			\int_{t'=0}^t
				\frac{1}{(1 + t')^{3/2}}
				\totonemax{\leq N}(t',u)
			\, dt'
			\notag	\\
		& \ \
			+
			\int_{u'=0}^u
				\totonemax{\leq N}(t,u')
			\, du'
			+ \varepsilon^2.
			\notag 
	\end{align}
  \end{subequations}
\end{lemma}

\begin{proof}
The proof is very similar to that of Lemma~\ref{L:SPACETIMEINTEGRALBOUNDSINTERMSOFQ0ANDQ1}.
Since the estimates are all very similar, we estimate 
only one representative product from \eqref{E:RENORMALIZEDLTIMESFRAMEDERIVATIVESINTERMSOFQ0ANDQ1}.
For example, decomposing
$\Lunit \Psi =\left \lbrace \Lunit \Psi + \frac{1}{2} \mytr \upchi \Psi \right\rbrace - \frac{1}{2} \mytr \upchi \Psi,$
using the bound $|\mytr \upchi|(t,u,\vartheta) \lesssim (1 + t)^{-1}$
(that is, \eqref{E:CRUDELOWERORDERC0BOUNDDERIVATIVESOFANGULARDEFORMATIONTENSORS}),
using Prop.~\ref{P:L2NORMSOFPSIINTERMSOFTHECOERCIVEQUANTITIES},
and arguing as in our proof of \eqref{E:REPTERMSTANDARDPSISPACETIMEINTEGRALSWITHNOSMALLMU}, 
we deduce that
\begin{align} \label{E:REPTERMRENORMALIZEDLDERIVATIVESPACETIMEINTEGRALSWITHNOSMALLMU}
&
\int_{\mathcal{M}_{t,u}}
		 	\ln^{\Conone}(\myexp + t') 
		 	\left|
		 		\left\lbrace
		 			\Lunit 
		 			+ \frac{1}{2} \mytr \upchi
		 		\right\rbrace
		 		 \mathscr{Z}^{\leq N} \Psi
		 	\right|
		 	\left|
		 		\rgeo \Lunit \mathscr{Z}^{\leq N} \Psi
		 	\right|
		 \, d \vol
		\\
& \lesssim 
\int_{u'=0}^u \int_{t'=0}^t \int_{S_{t',u'}}
		 	\frac{\ln^{\Conone}(\myexp + t')}{1 + t'}
		 	\left|
		 		\rgeo
		 		\left\lbrace
		 			\Lunit 
		 			+ \frac{1}{2} \mytr \upchi 
		 		\right\rbrace
		 		\mathscr{Z}^{\leq N} \Psi
		 	\right|^2
\, d \spherevol \, dt' \, du'
	\notag \\
& \ \ 
+ 
\int_{t'=0}^t 
	\frac{1}{(1 + t')^{3/2}}
	\int_{\Sigma_{t'}^u} 
	|\mathscr{Z}^{\leq N} \Psi|^2
	\, d \tvol
\, dt'
\notag \\
& \lesssim 
		\frac{\ln^{\Conone}(\myexp + t)}{1 + t}
		\int_{u'=0}^u 
			\totonemax{\leq N}(t,u')
		\, du'
		+ \int_{t'=0}^t 
				\left| \frac{d}{dt'} \frac{\ln^{\Conone}(\myexp + t')}{1 + t'} \right|
				\totonemax{\leq N}(t',u)
		\, dt'	
			\notag \\
	& \ \
		+ \int_{t'=0}^t
				\frac{1}{(1 + t')^{3/2}}
				\totzeromax{\leq N}(t',u)
			\, dt',
		\notag
\end{align}
where to deduce the last inequality, 
we integrated by parts in $t'$ and used the facts that
$\totzeromax{\leq N}$ and $\totonemax{\leq N}$ are increasing in their arguments
and that $u \leq U_0$ (as in our proof of \eqref{E:REPTERMSTANDARDPSISPACETIMEINTEGRALSWITHNOSMALLMU}). 
Clearly, the right-hand side of \eqref{E:REPTERMRENORMALIZEDLDERIVATIVESPACETIMEINTEGRALSWITHNOSMALLMU}
is bounded by the right-hand side of \eqref{E:RENORMALIZEDLTIMESFRAMEDERIVATIVESINTERMSOFQ0ANDQ1} as desired.

The bound \eqref{E:RENORMALIZEDLSQUAREDINTERMSOFQ0ANDQ1} follows easily from 
\eqref{E:ONEPLUSTTIMESLPLUSTRCHIPSICONEL2INTERMSOFONE}.
The bounds 
\eqref{E:RENORMALIZEDLTIMESLDERIVATIVEOFEIKONALFUNCTIONINTERMSOFQ0ANDQ1}
and
\eqref{E:RENORMALIZEDLTIMESEIKONALFUNCTIONINTERMSOFQ0ANDQ1} 
can be proved by using arguments similar to those used in proving
\eqref{E:STANDARDPSISPACETIMEINTEGRALSINVOLVINGLDERIVATESOFEIKONALWITHNOSMALLMU}
and
\eqref{E:STANDARDPSISPACETIMEINTEGRALSTIMESEASYEIKONALFUNCTIONQUANTITIES}.

\end{proof}

The following key lemma,
which is based on the availability of the coercive
Morawetz spacetime integral \eqref{E:MORAWETZSPACETIMECOERCIVITY}, 
provides us with control over
error integrals involving non-$\upmu-$weighted angular derivatives
of $\Psi.$

\begin{lemma}[\textbf{Bounds for spacetime integrals in terms of the coercive Morawetz integral $\Morint$}]
	\label{L:MORAWETZSPACETIMEINTEGRALESTIMATES}
	Let $\Conone$ be a constant. 
	Under the small-data and bootstrap assumptions 
	of Sects.~\ref{S:PSISOLVES}-\ref{S:C0BOUNDBOOTSTRAP},
	if $\varepsilon$ is sufficiently small, 	
	then the following inequalities hold for $(t,u) \in [0,\Tboot) \times [0,U_0]$
	(and the constants
	implicit in ``$\overset{\Conone}{\lesssim}$'' are independent of $\varsigma$):
	\begin{subequations}
	\begin{align}  \label{E:MORAWETZANGULARDERIVATIVESFRAMEDERIVATIVESINTERMSOFQ0ANDQ1}
		&
		\int_{\mathcal{M}_{t,u}}
			\mathbf{1}_{\lbrace \upmu \leq 1/4 \rbrace} 
			\frac{\ln^{\Conone}(\myexp + t')}{1 + t'} 
			\left|
				\angdiff \mathscr{Z}^{\leq N} \Psi
			\right|	
			\left|
				\fourmyarray[(1 + \upmu)\Lunit \mathscr{Z}^{\leq N} \Psi]
					{\angdiff \mathscr{Z}^{\leq N} \Psi}
					{\Rad \mathscr{Z}^{\leq N} \Psi}
					{\mathscr{Z}^{\leq N} \Psi}
			\right|
		 \, d \vol
		 	\\
		 &  \overset{\Conone}{\lesssim}
		 		\varsigma^{-1}
		 		\int_{t'=0}^t
					\frac{1}{(1 + t')^{3/2}}
					\totzeromax{\leq N}(t',u)
				 \, dt'
		 		+ \varsigma^{-1}
		 			\int_{u'=0}^u
					\totzeromax{\leq N}(t,u')
				\, du'
			\notag \\
		& \ \ + 
				\varsigma 
				\sup_{t' \in [0,t)} \frac{\totMormax{\leq N}(t',u)}{(1 + t')^{1/2}},
			\notag 	
				\\
	&
		\int_{\mathcal{M}_{t,u}}
			\mathbf{1}_{\lbrace \upmu \leq 1/4 \rbrace} 
			\ln^{\Conone}(\myexp + t')
			\left| 
				\angdiff \mathscr{Z}^{\leq N} \Psi 
			\right|	
			\left| 
				\myarray
					[\rgeo 
					\left\lbrace
		 				\Lunit 
		 				+ \frac{1}{2} \mytr \upchi
		 			 \right\rbrace	
		 			 \mathscr{Z}^{\leq N} \Psi
		 			 ]
		 			{\mathscr{Z}^{\leq N} \Psi}
			\right|
		  \, d \vol
		  		\label{E:MORAWETZANGULARDERIVATIVESRENORMALIZEDLANDPSIINTERMSOFQ0ANDQ1}	\\
		 & \overset{\Conone}{\lesssim}
		 	\varsigma^{-1}
			\int_{t'=0}^t
				\frac{1}{(1 + t')^{3/2}}
				\totzeromax{\leq N}(t',u)
			 \, dt'
			+
			\varsigma^{-1}
			\int_{u'=0}^u
				\totonemax{\leq N}(t,u')
			\, du'
			\notag \\
		& \ \ + \varsigma \totMormax{\leq N}(t,u).
			 \notag
	\end{align}
	\end{subequations}

\end{lemma}

\begin{proof}[Proof of Lemma~\ref{L:MORAWETZSPACETIMEINTEGRALESTIMATES}]
We prove the estimate for the term
$\ln^{\Conone}(\myexp + t')(1 + t')^{-1} |\angdiff \mathscr{Z}^{\leq N} \Psi||\Lunit \mathscr{Z}^{\leq N} \Psi|$
in \eqref{E:MORAWETZANGULARDERIVATIVESFRAMEDERIVATIVESINTERMSOFQ0ANDQ1};
the proofs of the remaining estimates in \eqref{E:MORAWETZANGULARDERIVATIVESFRAMEDERIVATIVESINTERMSOFQ0ANDQ1} are similar.
The proof essentially involves Cauchy-Schwarz and a suitable partitioning of $[0,t].$ To proceed,
we partition $[0,t]$ into sub-intervals
$[t_1,t_2],$ $[t_2,t_3],$ $\cdots,$ $[t_{M-1}, t_M],$ where $t_m$ is the sequence of times 
defined by $t_m = m^4 - 1$ and the last sub-interval may be a partial one.
In particular, $m = (1 + t_m)^{1/4},$ and $\frac{(1 + t_{m+1})^{1/4}}{(1 + t_m)^{1/4}} \leq 2.$
From this last inequality, spacetime Cauchy-Schwarz, 
Prop.~\ref{P:L2NORMSOFPSIINTERMSOFTHECOERCIVEQUANTITIES},
and Lemma~\ref{L:MORAWETZSPACETIMECOERCIVITY}, 
we deduce that
	\begin{align}
		& \int_{\mathcal{M}_{t,u}}
			\frac{\ln^{\Conone}(\myexp + t')}{1 + t'} 
			\mathbf{1}_{\lbrace \upmu \leq 1/4 \rbrace} 
			|\Lunit \mathscr{Z}^{\leq N} \Psi|
			|\angdiff \mathscr{Z}^{\leq N} \Psi|	
		 \, d \vol
		 \\
		& \leq 
		C
		\sum_{m=1}^{M-1}
		\left(
			\int_{t'=t_m}^{t_{m+1}} \int_{\Sigma_{t'}^u}
				\frac{(1+t')}{\ln(\myexp + t')} 
				\mathbf{1}_{\lbrace \upmu \leq 1/4 \rbrace} 
				|\angdiff \mathscr{Z}^{\leq N} \Psi|^2	
			\, d \tvol \, dt'
		\right)^{1/2}
				\notag \\
	& \ \ \ \ \ \ \ \ \ \ \
		\times	
		\left(
			\int_{t'=t_m}^{t_{m+1}} \int_{\Sigma_{t'}^u}
				\frac{\ln^{2\Conone+1} (\myexp + t')}{(1+t')^3} 
				|\Lunit \mathscr{Z}^{\leq N} \Psi|^2	
			\, d \tvol \, dt'
		\right)^{1/2}
		\notag
			\\
		& \leq 
			C
			\sum_{m=1}^{M-1} \frac{\totMormax{\leq N}^{1/2}(t_{m+1},u)}{(1 +t_{m+1})^{1/4}}
				\frac{(1 +t_{m+1})^{1/4}}{(1 +t_m)^{1/4}}
					\left(
			\int_{t'=t_m}^{t_{m+1}} \int_{\Sigma_{t'}^u}
				\frac{\ln^{2\Conone+1}(\myexp + t')}{(1+t')^{5/2}} 
				|\Lunit \mathscr{Z}^{\leq N} \Psi|^2	
			\, d \tvol \, dt'
		\right)^{1/2}
			\notag \\
		& \leq C \sup_{s \in [0,t)} \frac{\totMormax{\leq N}^{1/2}(s,u)}{(1+s)^{1/4}}
			\sum_{m=1}^{M-1}
			\left(
				\int_{t'=t_m}^{t_{m+1}} \int_{\Sigma_{t'}^u}
					\frac{\ln^{2\Conone+1}(\myexp + t')}{(1+t')^{5/2}} 
					|\Lunit \mathscr{Z}^{\leq N} \Psi|^2	
				\, d \tvol \, dt'
			\right)^{1/2}	
				\notag \\
		& \leq C
			\left(
				\sup_{s \in [0,t)} \frac{\totMormax{\leq N}^{1/2}(s,u)}{(1+s)^{1/4}}
			\right)
			\left( 
				\sum_{m=1}^{\infty}
					m^{-2}
			\right)^{1/2}
			\left(
			\sum_{m=1}^{M-1}
				\int_{t'=t_m}^{t_{m+1}} \int_{\Sigma_{t'}^u}
					m^2
					\frac{\ln^{2\Conone+1}(\myexp + t')}{(1+t')^{5/2}} 
					|\Lunit \mathscr{Z}^{\leq N} \Psi|^2	
				\, d \tvol \, dt'
			\right)^{1/2}	
			\notag \\
		& \leq C
			\left(
				\sup_{s \in [0,t)} \frac{\totMormax{\leq N}^{1/2}(s,u)}{(1+s)^{1/4}}
			\right)
			\left(
				\int_{u'=0}^u 
				\int_{\mathcal{C}_{u'}^t}
					|\Lunit \mathscr{Z}^{\leq N} \Psi|^2	
				\, d \conevol \, du'
			\right)^{1/2}	
			\notag \\
		& \leq C
			\left(
				\sup_{s \in [0,t)} \frac{\totMormax{\leq N}^{1/2}(s,u)}{(1+s)^{1/4}}
			\right)
			\left(
			\int_{u'=0}^u
				\totzeromax{\leq N}(t,u')
			\, du'
			\right)^{1/2}
			\notag \\
		& \leq C \varsigma
				\sup_{s \in [0,t)} \frac{\totMormax{\leq N}(s,u)}{(1 + s)^{1/2}}
			+ C \varsigma^{-1}
				\int_{u'=0}^u
					\totzeromax{\leq N}(t,u')
				\, du',
			\notag			
	\end{align}
	which is in turn bounded by the right-hand side of \eqref{E:MORAWETZANGULARDERIVATIVESFRAMEDERIVATIVESINTERMSOFQ0ANDQ1}
	as desired.
	
	The estimates in \eqref{E:MORAWETZANGULARDERIVATIVESRENORMALIZEDLANDPSIINTERMSOFQ0ANDQ1} can similarly be proved
	by taking $t_m = m^B-1$ for a suitable choice of $B.$
	
\end{proof}

We summarize the content of the previous three lemmas in the next two corollaries.

\begin{corollary}[\textbf{Integral estimates for the Harmless spacetime error integrals corresponding to the multiplier $\Mult$}]
\label{C:EASYHARMLESSMULTERRORINTEGRAL}
Let $0 \leq N \leq 24$ be an integer and let $\varsigma > 0$ be a number. 
Under the small-data and bootstrap assumptions 
of Sects.~\ref{S:PSISOLVES}-\ref{S:C0BOUNDBOOTSTRAP},
if $\varepsilon$ is sufficiently small, then the following inequalities hold 
for the $Harmless^{\leq N}$ terms from Def.~\ref{D:HARMLESSTERMS}
for $(t,u) \in [0,\Tboot) \times [0,U_0]$
(and the implicit constants are \textbf{independent of} $\varsigma$):
\begin{align} \label{E:EASYHARMLESSMULTERRORINTEGRAL}
	& 	\int_{\mathcal{M}_{t,u}}
				\left\lbrace
				\left|
					(1 + 2 \upmu) \Lunit \mathscr{Z}^{\leq N} \Psi
				\right|  
				+ 
				\left|	
					2 \Rad \mathscr{Z}^{\leq N} \Psi
				\right|
				\right\rbrace
				\left|
					Harmless^{\leq N}
				\right|
				\, d \vol
	 	\\
	& \lesssim 
		(1 + \varsigma^{-1})
		\int_{t'=0}^t 
			\frac{1}{(1 + t')^{3/2}}
			\totzeromax{\leq N}(t',u) 
		\, dt'	 
		+ (1 + \varsigma^{-1})
			\int_{t'=0}^t 
				\frac{1}{(1 + t')^{3/2}}
				\totonemax{\leq N}(t',u) 
			\, dt'
		\notag \\
	& \ \
		+ (1 + \varsigma^{-1})
			\frac{1}{(1 + t)^{1/2}}
			\int_{u'=0}^u
				\totonemax{\leq N}(t,u') 
			\, du'
		+ \varsigma 
			\sup_{t' \in [0,t]} 
			\frac{\totMormax{\leq N}(t',u)}{(1 + t')^{1/2}}
		+ \varepsilon^3.
			\notag
\end{align}

\end{corollary}

\begin{proof}
	Referring to Def.~\ref{D:HARMLESSTERMS}
	and using the trivial bound 
	$|\angdiff \mathscr{Z}^N \Psi| \lesssim \mathbf{1}_{\lbrace \upmu \leq 1/4 \rbrace} |\angdiff \mathscr{Z}^N \Psi| + \sqrt{\upmu} |\angdiff \mathscr{Z}^N \Psi|,$
	we see that all of the terms on the left-hand side of \eqref{E:EASYHARMLESSMULTERRORINTEGRAL}
	can be bounded by the right-hand side of \eqref{E:EASYHARMLESSMULTERRORINTEGRAL}
	with the help of the estimates of
	Lemma~\ref{L:SPACETIMEINTEGRALBOUNDSINTERMSOFQ0ANDQ1},
	Lemma~\ref{L:SPACETIMEINTEGRALSWITHLPSIPLUSHALFTRCHIPSI},
	and
	Lemma~\ref{L:MORAWETZSPACETIMEINTEGRALESTIMATES}.
\end{proof}

\begin{corollary}[\textbf{Estimates for the Harmless error integrals corresponding to the multiplier $\Mor$}]
\label{C:EASYHARMLESSMORERRORINTEGRAL}
Let $0 \leq N \leq 24$ be an integer and let $\widetilde{\varsigma} > 0$ be a number. 
Under the small-data and bootstrap assumptions 
of Sects.~\ref{S:PSISOLVES}-\ref{S:C0BOUNDBOOTSTRAP},
if $\varepsilon$ is sufficiently small, 
then the following inequalities hold 
for the $Harmless^{\leq N}$ terms from Def.~\ref{D:HARMLESSTERMS}
for $(t,u) \in [0,\Tboot) \times [0,U_0]$
(and the implicit constants are \textbf{independent of} $\widetilde{\varsigma}$):
\begin{align} \label{E:EASYHARMLESSMORERRORINTEGRAL}
	& 	\int_{\mathcal{M}_{t,u}}
				\rgeo^2 
					\left|
						\left\lbrace
							\Lunit 
							+ \frac{1}{2} \mytr \upchi 
						\right\rbrace
						\mathscr{Z}^{\leq N} \Psi
					\right|
				\left|
					Harmless^{\leq N}
				\right|
				\, d \vol
	 		\\
	 & \lesssim 
	 		(1 + \widetilde{\varsigma}^{-1})
	 		\int_{t'=0}^t
				\frac{1}{(1 + t')^{3/2}}
				\totzeromax{\leq N}(t',u)
			\, dt'
			+ 
			(1 + \widetilde{\varsigma}^{-1})
			\int_{t'=0}^t
				\frac{1}{(1 + t')^{3/2}}
				\totonemax{\leq N}(t',u)
			\, dt'
			\notag	\\
	 	& \ \ + 
	 				(1 + \widetilde{\varsigma}^{-1})
	 				\int_{u'=0}^u
						\totonemax{\leq N}(t,u') 
					\, du'
			+ \widetilde{\varsigma} \totMormax{\leq N}(t,u)
			+ \varepsilon^3.
			\notag
\end{align}

\end{corollary}

\begin{proof}
	The proof is essentially the same as that of Cor.~\ref{C:EASYHARMLESSMULTERRORINTEGRAL}.
\end{proof}

As we will see, when deriving some of our $L^2$ estimates 
for the below-top-order derivatives of $\Psi,$
we must allow some of the error integrals to lose one derivative.
This is essential, for if we instead chose to avoid losing derivatives, 
then we would have to use the modified quantities of Ch.~\ref{C:RENORMALIZEDEIKONALFUNCTIONQUANTITIES},
and the resulting $L^2$ estimates would be as degenerate with
respect to powers of $\upmu_{\star}^{-1}$ as they
at the top order; such degeneracy would prevent us from 
deriving improvements of the lower-order bootstrap assumptions of Sect.~\ref{S:LTWOBOOTSTRAPASSUMPTIONS}.
In the next two lemmas, we derive estimates that are sufficient for controlling
the derivative-losing error integrals.

\begin{lemma}[\textbf{Estimates for the error integrals corresponding to the multiplier $\Mult$ involving a loss of one derivative}]
\label{L:MULTERRORINTEGRALESTIMATELOSSOFONEDERIVATIVE}
Let $1 \leq N \leq 23$ be an integer. 
Under the small-data and bootstrap assumptions 
of Sects.~\ref{S:PSISOLVES}-\ref{S:C0BOUNDBOOTSTRAP},
if $\varepsilon$ is sufficiently small, 
then the following inequalities hold for $(t,u) \in [0,\Tboot) \times [0,U_0]:$
\begin{align} \label{E:MULTERRORINTEGRALESTIMATELOSSOFONEDERIVATIVE}
	& 	\int_{\mathcal{M}_{t,u}}
				\frac{1}{1 + t'}
				\left|
					(1 + 2 \upmu) \Lunit \mathscr{Z}^N \Psi  + \Rad \mathscr{Z}^N \Psi
				\right|
				\left|
					\myarray
						[\angLap \mathscr{Z}^{N-1} \upmu]
						{\Rot \mathscr{Z}^{N-1} \mytr \upchi^{(Small)}}
				\right|
				\, d \vol
	 	\\
& \lesssim
\int_{t'=0}^t
	\frac{1}{(1 + t')^{3/2} \upmu_{\star}^{1/2}(t',u)}
		\totzeromax{\leq N}^{1/2}(t',u)	
		\int_{s=0}^{t'}
			\frac{1}{1 + s}
			\totzeromax{\leq N+1}^{1/2}(s,u)
		\, ds
\, dt'	
\notag 
\\
& \ \ 
+ \int_{t'=0}^t
	\frac{1}{(1 + t')^{3/2} \upmu_{\star}^{1/2}(t',u)}
		\totzeromax{\leq N}^{1/2}(t',u)	
		\int_{s=0}^{t'}
			\frac{1}{(1 + s)\upmu_{\star}^{1/2}(s,u)}
			\totonemax{\leq N+1}^{1/2}(s,u)
		\, ds
\, dt'
	\notag \\
& \ \
+ 
\int_{t'=0}^t
	\frac{1}{(1 + t')^{3/2} \upmu_{\star}^{1/2}(t',u)}
		\totzeromax{\leq N}(t',u)	
\, dt'
+ \varepsilon^2.
	\notag 
\end{align}
\end{lemma}

\begin{proof}
Since the proofs of both of the estimates are similar, we provide only the proof of one of them.
Specifically, using 
\eqref{E:FUNCTIONPOINTWISEANGDINTERMSOFANGLIEO}
and
\eqref{E:ANGLAPFUNCTIONPOINTWISEINTERMSOFROTATIONS}, 
we see that 
$|\angLap \mathscr{Z}^{N-1} \upmu| \lesssim (1 + t)^{-2} |\mathscr{Z}^{\leq N+1} \upmu|.$
Hence, using this bound,
Cauchy-Schwarz on the $\Sigma_{t'}^u,$
Prop.~\ref{P:L2NORMSOFPSIINTERMSOFTHECOERCIVEQUANTITIES},
and
\eqref{E:EIKONALFUNCTIONQUANTITIESL2BOUNDSINTERMSOFQ0ANDQ1}, 
we have
\begin{align} \label{E:REPTERMMULTERRORINTEGRALESTIMATELOSSOFONEDERIVATIVE}
& \int_{\mathcal{M}_{t,u}}
		\frac{1}{1 + t'}
		|\Lunit \mathscr{Z}^N \Psi|
		|\angLap \mathscr{Z}^{N-1} \upmu|
\, d \vol		
	\\
& \lesssim
\int_{t'=0}^t
	\frac{1}{(1 + t')^3}
	\int_{\Sigma_{t'}^u}
		|\Lunit \mathscr{Z}^N \Psi|	
		|\mathscr{Z}^{\leq N+1} \upmu|
	\, d \tvol
\, dt'
	\notag \\
& \lesssim
\varepsilon
\int_{t'=0}^t
	\frac{1}{(1 + t')^2 \upmu_{\star}^{1/2}(t',u)}
		\totzeromax{\leq N}^{1/2}(t',u)	
\, dt'
	\notag \\
&	\ \
+
\int_{t'=0}^t
	\frac{1}{(1 + t')^2 \upmu_{\star}^{1/2}(t',u)}
		\totzeromax{\leq N}^{1/2}(t',u)	
		\int_{s=0}^{t'}
			\frac{1}{1 + s}
			\totzeromax{\leq N+1}^{1/2}(s,u)
		\, ds
\, dt'	
\notag 
\\
& \ \ 
+
\int_{t'=0}^t
	\frac{1}{(1 + t')^2 \upmu_{\star}^{1/2}(t',u)}
		\totzeromax{\leq N}^{1/2}(t',u)	
		\int_{s=0}^{t'}
			\frac{1}{(1 + s)\upmu_{\star}^{1/2}(s,u)}
			\totonemax{\leq N+1}^{1/2}(s,u)
		\, ds
\, dt'.
\notag
\end{align}			
To complete the proof of inequality \eqref{E:MULTERRORINTEGRALESTIMATELOSSOFONEDERIVATIVE},
it remains only for us to bound the first integral on the right-hand side of 
\eqref{E:REPTERMMULTERRORINTEGRALESTIMATELOSSOFONEDERIVATIVE}
by the last two terms on the right-hand side of \eqref{E:MULTERRORINTEGRALESTIMATELOSSOFONEDERIVATIVE}.
To derive this bound, we use the simple inequality
$\varepsilon \totzeromax{\leq N}^{1/2}(t',u) \lesssim \varepsilon^2 + \totzeromax{\leq N}(t',u)$
together with inequality \eqref{E:LESSSINGULARTERMSMUTHREEFOURTHSINTEGRALBOUND}.

\end{proof}

\begin{lemma}[\textbf{Estimates for the error integrals corresponding to the multiplier $\Mor$ involving a loss of one derivative}]
\label{L:MORERRORINTEGRALESTIMATELOSSOFONEDERIVATIVE}
Let $1 \leq N \leq 23$ be an integer. 
Under the small-data and bootstrap assumptions 
of Sects.~\ref{S:PSISOLVES}-\ref{S:C0BOUNDBOOTSTRAP},
if $\varepsilon$ is sufficiently small, 
then the following inequalities hold for $(t,u) \in [0,\Tboot) \times [0,U_0]:$
\begin{align} \label{E:MORERRORINTEGRALESTIMATELOSSOFONEDERIVATIVE}
	& 	\int_{\mathcal{M}_{t,u}}
				\rgeo
				\left|
					\Lunit \mathscr{Z}^N \Psi 
						+ \frac{1}{2} \mytr \upchi \mathscr{Z}^N \Psi
				\right|
				\left|
					\myarray
						[\angLap \mathscr{Z}^{N-1} \upmu]
						{\rgeo \angdiff \mathscr{Z}^{N-1} \mytr \upchi^{(Small)}}
				\right|
				\, d \vol
	 	\\
	 & \lesssim 
	 	\int_{t'=0}^t
					\frac{1}{(1 + t')^2}
					\left(
					\int_{s=0}^{t'}
						\frac{ \totzeromax{\leq N+1}^{1/2}(s,u)}{1 + s} 
					\, ds
					\right)^2
		\, dt'
		+ 
		\int_{t'=0}^t
			\frac{1}{(1 + t')^2}
			\left(
				\int_{s=0}^{t'}
					\frac{\totonemax{\leq N+1}^{1/2}(s,u)}{(1 + s)\upmu_{\star}^{1/2}(s,u)} 
				\, ds
			\right)^2
		\, dt'
	 	\notag \\
	 & \ \ 
		+ \int_{u' = 0}^u
				\totonemax{\leq N}(t,u')
			\, du' 
		+ \varepsilon^2.
		\notag
\end{align}
\end{lemma}

\begin{proof}
	Since the proofs of both of the estimates are similar, we provide only the proof of one of them.
	Specifically, using the bound 
	$|\angLap \mathscr{Z}^{N-1} \upmu| \lesssim (1 + t)^{-2} |\mathscr{Z}^{\leq N+1} \upmu|$
	noted in the proof of Lemma~\ref{L:MULTERRORINTEGRALESTIMATELOSSOFONEDERIVATIVE},
	Cauchy-Schwarz,
	and Prop.~\ref{P:L2NORMSOFPSIINTERMSOFTHECOERCIVEQUANTITIES}, 
	we have
	\begin{align} \label{E:REPTERMMORERRORINTEGRALESTIMATELOSSOFONEDERIVATIVE}
		&
		\int_{\mathcal{M}_{t,u}}
			\rgeo
			\left|
				\Lunit \mathscr{Z}^N \Psi 
				+ \frac{1}{2} \mytr \upchi \mathscr{Z}^N \Psi
			\right|
			|\angLap \mathscr{Z}^{N-1} \upmu|
		\, d \vol
			\\
		& \lesssim 
			\int_{\mathcal{M}_{t,u}}
				\rgeo^2
				\left|
					\Lunit \mathscr{Z}^N \Psi 
					+ \frac{1}{2} \mytr \upchi \mathscr{Z}^N \Psi
				\right|^2
			\, d \vol
			+
			\int_{t'=0}^t
				\frac{1}{(1 + t')^4}
				\int_{\Sigma_{t'}^u}
				 |\mathscr{Z}^{\leq N+1} (\upmu - 1)|^2
				 \, d \tvol
			\, dt'
				\notag \\
			& \lesssim 
			\int_{u'=0}^u
				\totonemax{\leq N}(t,u')
			\, d u'
			+
			\int_{t'=0}^t
				\frac{1}{(1 + t')^4}
				\int_{\Sigma_{t'}^u}
				 |\mathscr{Z}^{\leq N+1} (\upmu - 1)|^2
				 \, d \tvol
			\, dt'.
				\notag
	\end{align}
	To bound the second spacetime integral on the right-hand side of 
	\eqref{E:REPTERMMORERRORINTEGRALESTIMATELOSSOFONEDERIVATIVE},
	we use the estimate
	\eqref{E:EIKONALFUNCTIONQUANTITIESL2BOUNDSINTERMSOFQ0ANDQ1}
	to deduce that
	\begin{align} \label{E:ANNOYINGMORTRIPLETIMEINTEGRATIONBELOWTOPORDERTERM}
		& \int_{t'=0}^t
			\frac{1}{(1 + t')^4}
			\int_{\Sigma_{t'}^u}
				 |\mathscr{Z}^{\leq N+1} (\upmu - 1)|^2
			\, d \tvol
		\, dt'
			\\
		& \lesssim 
			\varepsilon^2
			\int_{t'=0}^t
				\frac{1}{(1 + t')^2}	
			\, dt'
			\notag \\
		& \ \ 
			+ \int_{t'=0}^t
					\frac{1}{(1 + t')^2}
					\left(
					\int_{s=0}^{t'}
						\frac{ \totzeromax{\leq N+1}^{1/2}(s,u)}{1 + s} 
					\, ds
					\right)^2
				\, dt'
				\notag \\
		& \ \ + 
				\int_{t'=0}^t
					\frac{1}{(1 + t')^2}
					\left(
					\int_{s=0}^{t'}
						\frac{\totonemax{\leq N+1}^{1/2}(s,u)}{(1 + s)\upmu_{\star}^{1/2}(s,u)} 
					\, ds
					\right)^2
				\, dt'.
				\notag
\end{align}
Combining \eqref{E:REPTERMMORERRORINTEGRALESTIMATELOSSOFONEDERIVATIVE}	
and \eqref{E:ANNOYINGMORTRIPLETIMEINTEGRATIONBELOWTOPORDERTERM},
we deduce the desired estimate
\eqref{E:MORERRORINTEGRALESTIMATELOSSOFONEDERIVATIVE}.	

\end{proof}

The next lemma, which is based on a combination 
of transport equation estimates and elliptic estimates,
will allow us to derive non-sharp (but suitable for proving our sharp classical lifespan theorem)
top-order Sobolev estimates
for $\upmu \angD^2 \mathscr{Z}^{N-1} \upmu$ and $\upmu \angD \mathscr{Z}^{N-1} \upchi^{(Small)}.$
The main point is that the tensorfields 
$\upmu \angfreeDsquared \mathscr{Z}^{N-1} \upmu$ and
$\upmu \angD \mathscr{Z}^{N-1} \hat{\upchi}^{(Small)}$
cannot be estimated via pure energy estimates. The obstacle
is that we have not been able to derive a transport equation for a modified version of them
with a consistent number of derivatives on the right-hand side. 
To overcome this difficulty, we combine the elliptic estimates of Ch.~\ref{C:ELLIPTIC} with a Gronwall 
argument, which will allow us to bound these tensorfields in terms of other quantities that do not lose derivatives.

\begin{lemma}[\textbf{Non-sharp top-order $L^2$ estimates for 
$\angD^2 \mathscr{Z}^{N-1} \upmu,$ 
$\angD^2 \mathscr{Z}^{N-1} \Lunit_{(Small)}^i,$
and 
$\angD \mathscr{Z}^{N-1} \upchi^{(Small)}$ in terms of $\totzeromax{N}$ and $\totonemax{N}$}]
\label{L:TOPORDERELLIPTICRECOVERY}
Let $1 \leq N \leq 24$ be an integer.
Under the small-data and bootstrap assumptions 
of Sects.~\ref{S:PSISOLVES}-\ref{S:C0BOUNDBOOTSTRAP},
if $\varepsilon$ is sufficiently small, 
then for $(t,u) \in [0,\Tboot) \times [0,U_0],$ we have the following estimates:
\begin{subequations}
\begin{align} \label{E:ELLIPTICRECOVERY}
	& \fivemyarray
		[\| \upmu \angD^2 \mathscr{Z}^{\leq N-1} \upmu \|_{L^2(\Sigma_t^u)}]
		{\sum_{a=1}^3 \|\rgeo \upmu \angD^2 \mathscr{Z}^{\leq N-1} \Lunit_{(Small)}^a \|_{L^2(\Sigma_t^u)}}
		{\| \rgeo \upmu \angD \angLie_{\mathscr{Z}}^{\leq N-1} \upchi^{(Small)} \|_{L^2(\Sigma_t^u)}}
		{\| \rgeo \upmu \angdiff \mathscr{Z}^{\leq N-1} \mytr \upchi^{(Small)} \|_{L^2(\Sigma_t^u)}}
		{\| \rgeo \upmu \angD \angLie_{\mathscr{Z}}^{\leq N-1} \hat{\upchi}^{(Small)} \|_{L^2(\Sigma_t^u)}}
		\\
	& \lesssim
			\int_{t'=0}^t 
						\frac{\| [\Lunit \upmu]_- \|_{C^0(\Sigma_{t'}^u)}} 
							   {\upmu_{\star}({t'},u)} 
						\totzeromax{\leq N}^{1/2}({t'},u) 
			\, dt'
			+
			\int_{t'=0}^t 
					\frac{1}{(1 + t')^{3/2} \upmu_{\star}(t',u)}
					\totzeromax{\leq N}^{1/2}(t',u)
				\, dt'
			+ \totzeromax{\leq N}^{1/2}(t,u) 
			\notag	\\
	& \ \ + \int_{t'=0}^t 
						\frac{1}{(1 + t')^{3/2} \upmu_{\star}(t',u)} 
						\totonemax{\leq N}^{1/2}(t',u)
					\, dt'
				+ \varepsilon
					\left\lbrace
						\ln \upmu_{\star}^{-1}(t,u) + 1
					\right\rbrace,
				\notag
\end{align}
\begin{align}				
	&
	\left\|
	\int_{t'=0}^t 
			\left| 
				\fivemyarray
					[\upmu \angD^2 \mathscr{Z}^{\leq N-1} \upmu]
					{\sum_{a=1}^3 \|\rgeo \upmu \angD^2 \mathscr{Z}^{\leq N-1} \Lunit_{(Small)}^a \|_{L^2(\Sigma_t^u)}}
					{\rgeo \upmu \angD \angLie_{\mathscr{Z}}^{\leq N-1} \upchi^{(Small)} }
					{\rgeo \upmu \angdiff \mathscr{Z}^{\leq N-1} \mytr \upchi^{(Small)} }
					{\rgeo \upmu \angD \angLie_{\mathscr{Z}}^{\leq N-1} \hat{\upchi}^{(Small)} }
			\right|
			(t',\cdot)
		\, dt'
	\right\|_{L^2(\Sigma_t^u)}
		\label{E:NONSHARPTOPORDERELLIPTICRECOVERY} \\
	& \lesssim 
			\ln(\myexp + t) (1 + t) \totzeromax{\leq N}^{1/2}(t,u)
			+ \ln(\myexp + t) (1 + t) \totonemax{\leq N}^{1/2}(t,u)
			+ \varepsilon \ln(\myexp + t) (1 + t).
		\notag
\end{align}
\end{subequations}
\end{lemma}
We provide the proof of Lemma~\ref{L:TOPORDERELLIPTICRECOVERY} in
Sect.~\ref{S:PROOFOFLEMMATOPORDERELLIPTICRECOVERY}.

In the next two lemmas, we estimate some top-order error integrals.
The integrals are nontrivial to bound
in the sense that in order to avoid losing derivatives,
we need to use the transport-elliptic estimates
of Lemma~\ref{L:TOPORDERELLIPTICRECOVERY}.
However, because of the favorable time decay factors present in the integrands,
the integrals are not too difficult to estimate.

\begin{lemma}[\textbf{Estimates for the easy eikonal top-order error integrals corresponding to the multiplier $\Mult$}]
\label{L:HARMLESSELLIPTICTOPORDERMULTERRORINTEGRAL}
Let $1 \leq N \leq 24$ be an integer. 
Under the small-data and bootstrap assumptions 
of Sects.~\ref{S:PSISOLVES}-\ref{S:C0BOUNDBOOTSTRAP},
if $\varepsilon$ is sufficiently small, 
then the following inequalities hold for $(t,u) \in [0,\Tboot) \times [0,U_0]:$
\begin{align}  \label{E:HARMLESSELLIPTICTOPORDERMULTERRORINTEGRAL}
	& 	\int_{\mathcal{M}_{t,u}}
				\left|
					(1 + 2 \upmu) \Lunit \mathscr{Z}^{\leq N} \Psi  + 2 \Rad \mathscr{Z}^{\leq N} \Psi
				\right|
				\left|
					\threemyarray
						[(\upmu \angdiffuparg{\#} \Psi) \cdot (\upmu \angdiff \mathscr{Z}^{N-1} \mytr \upchi^{(Small)})]
						{\rgeo (\angdiff \Psi^{\#}) \cdot (\upmu \angdiff \mathscr{Z}^{N-1} \mytr \upchi^{(Small)})}
						{\RotRadcomponent{l}(\angdiffuparg{\#} \Psi) \cdot (\upmu \angdiff \mathscr{Z}^{N-1} \mytr \upchi^{(Small)})}
				\right|
				\, d \vol
	 	\\
	 & \lesssim 
	 	\varepsilon \totzeromax{\leq N}(t,u) 
		+ \varepsilon
		\int_{t'=0}^t
			\frac{1}{(1 + t')^2}
			\totonemax{\leq N}(t',u)
		\, dt'
		+ \varepsilon^3.
			\notag
\end{align}
\end{lemma}

\begin{proof}
Since the proofs of all of the estimates
involving $\Lunit \mathscr{Z}^{\leq N} \Psi,$ 
$\upmu \Lunit \mathscr{Z}^{\leq N} \Psi,$
and $\Rad \mathscr{Z}^{\leq N} \Psi$ are similar,
we provide only the proof for one representative product.
First, using the pointwise estimates
$\upmu \lesssim \ln(\myexp + t)$
(that is, \eqref{E:C0BOUNDCRUCIALEIKONALFUNCTIONQUANTITIES}),
$|\RotRadcomponent{l}| \lesssim \varepsilon \ln(\myexp + t)$
(that is, \eqref{E:LOWERORDERC0BOUNDEUCLIDEANROTATIONRADCOMPONENT} and Cor.~\ref{C:SQRTEPSILONREPLCEDWITHCEPSILON}),
and $|\angdiff \Psi| \lesssim \varepsilon (1 + t)^{-2}$
(that is, the bootstrap assumptions \eqref{E:PSIFUNDAMENTALC0BOUNDBOOTSTRAP}),
we deduce that
\begin{align} \label{E:EASYEIKONALMULTTERMPOINTWISEBOUND}
	\left|
		\threemyarray
		[(\upmu \angdiffuparg{\#} \Psi) \cdot (\upmu \angdiff \mathscr{Z}^{N-1} \mytr \upchi^{(Small)})]
		{\rgeo (\angdiff \Psi^{\#}) \cdot (\upmu \angdiff \mathscr{Z}^{N-1} \mytr \upchi^{(Small)})}
		{\RotRadcomponent{l}(\angdiffuparg{\#} \Psi) \cdot (\upmu \angdiff \mathscr{Z}^{N-1} \mytr \upchi^{(Small)})}
	\right|
	& \lesssim 
		\varepsilon
		\frac{1}{(1 + t)^2}
		\left|
			\rgeo \upmu \angdiff \mathscr{Z}^{N-1} \mytr \upchi^{(Small)}
		\right|.
\end{align}

Hence, using spacetime Cauchy-Schwarz, 
inequality \eqref{E:KEYMUINVERSEINTEGRALBOUND},
inequality \eqref{E:LESSSINGULARTERMSMUTHREEFOURTHSINTEGRALBOUND},
inequality \eqref{E:ELLIPTICRECOVERY}, 
inequality \eqref{E:EASYEIKONALMULTTERMPOINTWISEBOUND},
and Prop.~\ref{P:L2NORMSOFPSIINTERMSOFTHECOERCIVEQUANTITIES},
we have
\begin{align} \label{E:REPTERMHARMLESSELLIPTICTOPORDERMULTERRORINTEGRAL}
	&
	\int_{\mathcal{M}_{t,u}}
				\left|
					\Lunit \mathscr{Z}^{\leq N} \Psi
				\right|
				\left|
					\threemyarray
						[(\upmu \angdiffuparg{\#} \Psi) \cdot (\upmu \angdiff \mathscr{Z}^{N-1} \mytr \upchi^{(Small)})]
						{\rgeo(\angdiff \Psi^{\#}) \cdot (\upmu \angdiff \mathscr{Z}^{N-1} \mytr \upchi^{(Small)})}
						{\RotRadcomponent{l}(\angdiffuparg{\#} \Psi) \cdot (\upmu \angdiff \mathscr{Z}^{N-1} \mytr \upchi^{(Small)})}
				\right|
				\, d \vol
	 \\
	& \lesssim 
		\varepsilon
		\int_{u'=0}^u
			\int_{\mathcal{C}_{u'}^t}
				\frac{1}{(1 + t')^2}
				\left|
					\Lunit \mathscr{Z}^{\leq N} \Psi
				\right|^2
			\, d \conevol
		\, du'
		+ 
		\varepsilon
		\int_{t'=0}^t
			\frac{1}{(1 + t')^2}
			\int_{\Sigma_{t'}^u}	
				\left|
					\rgeo \upmu \angdiff \mathscr{Z}^{N-1} \mytr \upchi
				\right|^2	
			\, d \tvol
		\, dt'
				\notag \\
& \lesssim
		\varepsilon
		\int_{u'=0}^u
			\totzeromax{\leq N}(t,u')
		\, du'
	+ \varepsilon
		\int_{t'=0}^t
			\frac{1}{(1 + t')^2}
			\left(
				\int_{s=0}^{t'}
						\frac{\| [\Lunit \upmu]_- \|_{C^0(\Sigma_s^u)}} 
							   {\upmu_{\star}(s,u)} 
						\totzeromax{\leq N}^{1/2}(s,u) 
				\, ds
			\right)^2
		\, dt'
		\notag \\
& \ \ 
	+ \varepsilon
		\int_{t'=0}^t
			\frac{1}{(1 + t')^2}
			\totzeromax{\leq N}(t',u)
		\, dt'
	+ \varepsilon
		\int_{t'=0}^t
			\frac{1}{(1 + t')^2}
			\totonemax{\leq N}(t',u)
		\, dt'	
	+ \int_{t'=0}^t
			\frac{1}{(1 + t')^2}
			\varepsilon^3
		\, dt'
		\notag
		\\
& \lesssim
	\varepsilon
	\totzeromax{\leq N}(t,u)
	+ \int_{t'=0}^t
		\frac{1}{(1 + t')^{3/2}}
			\totzeromax{\leq N}(t',u)
			(\ln \upmu_{\star}^{-1}
			+ \sqrt{\varepsilon}
			)^2
		\, dt'
	+ \varepsilon
		\int_{t'=0}^t
			\frac{1}{(1 + t')^2}
			\totonemax{\leq N}(t',u)
		\, dt'
	+ \varepsilon^3
		\notag \\
	& \lesssim 
		\varepsilon \totzeromax{\leq N}(t,u) 
		+ \varepsilon
		\int_{t'=0}^t
			\frac{1}{(1 + t')^2}
			\totonemax{\leq N}(t',u)
		\, dt'
		+ \varepsilon^3
		\notag
\end{align}
as desired. 

\end{proof}

\begin{lemma}[\textbf{Estimates for the easy eikonal top-order error integrals corresponding to the multiplier $\Mor$}]
\label{L:HARMLESSELLIPTICTOPORDERMORERRORINTEGRAL}
Let $0 \leq N \leq 24$ be an integer.
Under the small-data and bootstrap assumptions 
of Sects.~\ref{S:PSISOLVES}-\ref{S:C0BOUNDBOOTSTRAP},
if $\varepsilon$ is sufficiently small, 
then the following inequalities hold for $(t,u) \in [0,\Tboot) \times [0,U_0]:$
\begin{align}
	& 	\int_{\mathcal{M}_{t,u}}
				\rgeo^2
				\left|
					\left\lbrace
						\Lunit 
						+ \frac{1}{2} \mytr \upchi 
					\right\rbrace
					\mathscr{Z}^{\leq N} \Psi
				\right|
				\left|
					\threemyarray
						[(\upmu \angdiffuparg{\#} \Psi) \cdot (\upmu \angdiff \mathscr{Z}^{N-1} \mytr \upchi^{(Small)})]
						{\rgeo(\angdiff \Psi^{\#}) \cdot (\upmu \angdiff \mathscr{Z}^{N-1} \mytr \upchi^{(Small)})}
						{\RotRadcomponent{l}(\angdiffuparg{\#} \Psi) \cdot (\upmu \angdiff \mathscr{Z}^{N-1} \mytr \upchi^{(Small)})}
				\right|
				\, d \vol
	 	\\
	 & \lesssim 
			\varepsilon \totzeromax{\leq N}(t,u)
	 		+ \varepsilon \totonemax{\leq N}(t,u)
	 		+ \varepsilon^3.
		\notag
\end{align}
\end{lemma}

\begin{proof}
	Using spacetime Cauchy-Schwarz, 
	inequality \eqref{E:EASYEIKONALMULTTERMPOINTWISEBOUND},
	and Prop.~\ref{P:L2NORMSOFPSIINTERMSOFTHECOERCIVEQUANTITIES},
	we have
	\begin{align}
		&
		\int_{\mathcal{M}_{t,u}}
				\rgeo^2
				\left|
					\left\lbrace
						\Lunit 
						+ \frac{1}{2} \mytr \upchi 
					\right\rbrace
					\mathscr{Z}^{\leq N} \Psi
				\right|
				\left|
					\threemyarray
						[(\upmu \angdiffuparg{\#} \Psi) \cdot (\upmu \angdiff \mathscr{Z}^{N-1} \mytr \upchi^{(Small)})]
						{\rgeo(\angdiff \Psi^{\#}) \cdot (\upmu \angdiff \mathscr{Z}^{N-1} \mytr \upchi^{(Small)})}
						{\RotRadcomponent{l}(\angdiffuparg{\#} \Psi) \cdot (\upmu \angdiff \mathscr{Z}^{N-1} \mytr \upchi^{(Small)})}
				\right|
			\,	d \vol
			\\
		& \lesssim
			\varepsilon
			\int_{u'=0}^u
				\int_{\mathcal{C}_{u'}^t}
					(1 + t')^2
					\left|
						\left\lbrace
							\Lunit + \frac{1}{2} \mytr \upchi 
					  \right\rbrace
					  \mathscr{Z}^{\leq N} \Psi
					\right|^2
				 \, d \conevol
			\, du'	
				\notag \\
		&  \ \
			+ 
			\varepsilon
			\int_{t'=0}^t
					\int_{\Sigma_{t'}^u}
						\frac{1}{(1 + t')^2}
						\left|
							\rgeo \upmu \angdiff \mathscr{Z}^{N-1} \mytr \upchi^{(Small)}
						\right|^2
					\, d \tvol
				\, dt'
			\notag \\	
	& \lesssim
			\varepsilon
			\totonemax{\leq N}(t,u)	
			+ \varepsilon
			  \int_{\mathcal{M}_{t,u}}
				\frac{1}{(1 + t')^2}
				\left|
					\rgeo \upmu \angdiff \mathscr{Z}^{N-1} \mytr \upchi^{(Small)}
				\right|^2
				\, d \vol.
			\notag
	\end{align}
	The remainder of the proof now proceeds as in \eqref{E:REPTERMHARMLESSELLIPTICTOPORDERMULTERRORINTEGRAL}.
	
\end{proof}

In the next lemma, we derive estimates for the error integrals
corresponding to the deformation tensors 
of the multiplier vectorfields $\Mult$ and $\Mor.$
Thanks to the pointwise estimates we have previously derived,
the lemma is not difficult to prove.

\begin{lemma}[\textbf{The main estimates for the error integrals 
corresponding to the deformation tensors of $\Mult$ and $\Mor$
in terms of $\totzeromax{N},$ 
$\totonemax{N},$
and
$\totMormax{N}$}] \label{L:BASICENERGYERRORINTEGRALS}
	Assume that the small data and bootstrap assumptions 
	of Sects.~\ref{S:PSISOLVES}-\ref{S:C0BOUNDBOOTSTRAP},
	hold for $(t,u) \in [0,\Tboot) \times [0,U_0].$
	There exists a constant $C > 0$ such that 
	if $\varepsilon$ is sufficiently small, 
	then for $(t,u) \in [0,\Tboot) \times [0,U_0],$
	we have the following bounds for the non-$\waveinhom-$containing integrals appearing on
	the right-hand sides of the inequalities of Prop.~\ref{P:DIVTHMWITHCANCELLATIONS}
	(the terms $\basicenergyerror{\Mult}[\Psi]$ and $\basicenergyerror{\Mor}[\Psi]$ 
	are the energy error integrands from
	Def.~\ref{D:LK0K1ERRORINTEGRANDS} and Lemma~\ref{L:LK0K1ERRORINTEGRANDS}, 
	and the Morawetz spacetime integral with a good sign has been subtracted
	from the term $\basicenergyerror{\Mor}[\Psi]$ 
	on the left-hand side of \eqref{E:Q1BASICERRORINTEGRALESTIMATE}):
\begin{subequations}
\begin{align} \label{E:Q0BASICERRORINTEGRALESTIMATE}
		\int_{\mathcal{M}_{t,u}}
			\basicenergyerror{\Mult}[\mathscr{Z}^N \Psi]								
		\, d \vol
		& \leq
			C
			\int_{t'=0}^t 
				\frac{\ln^2(\myexp + t')}{(1 + t')^2} \totzeromax{N}(t',u) 
			\, dt'
				\\
		& \ \ + C
						\int_{t'=0}^t 
							\frac{\ln^2(\myexp + t')}{(1 + t')^2} \totonemax{N}(t',u) 
						\, dt' 
					\notag	\\
		& \ \ 
			+ C
				\varepsilon^{1/2}
				\int_{t'=0}^t 
					\frac{\ln(\myexp + t') + 1}{(1 + t')^2 \sqrt{\ln(\myexp + t) - \ln(\myexp + t')}}
					\totonemax{N}(t',u) 
				\, dt' 
				\notag \\
		& \ \ + C
						\varepsilon
						\int_{u'=0}^u
							\totzeromax{N}(t,u')		
						\, du'
					\notag \\
		& \ \ + C
						\frac{\ln^2(\myexp + t)}{(1 + t)^2} 
						\int_{u'=0}^u
							\totonemax{N}(t,u')		
						\, du'
					\notag \\
		& \ \ 
			+ C \varepsilon 
				\sup_{t' \in [0,t]}\frac{\totMormax{N}(t',u)}{(1 + t')^{1/2}},
			\notag
	\end{align}
	
	\begin{align} \label{E:Q1BASICERRORINTEGRALESTIMATE}
		\int_{\mathcal{M}_{t,u}}
			\left|
				\basicenergyerror{\Mor}[\mathscr{Z}^N \Psi]
				+ \frac{1}{2} \rgeo^2 [\Lunit \upmu]_- |\angdiff \mathscr{Z}^N \Psi|^2 
			\right|							
		& \leq
			C \int_{t'=0}^t 
				\frac{\ln^3(\myexp + t')}{(1 + t')} \totzeromax{N}(t',u) 
			\, dt'
				\\
		& \ \ + (1 + C \varepsilon)
						\int_{t'=0}^t 
							\frac{1}{\rgeo(t',u) \ln \left(\frac{\rgeo(t',u)}{\rgeo(0,u)} \right)} \totonemax{N}(t',u) 
						\, dt' 
					\notag	\\
		& \ \ + C 
						\frac{\ln^3(\myexp + t)}{(1 + t)^2} 
						\int_{u'=0}^u
								\totonemax{N}(t,u')		
						\, du'
					\notag \\
		& \ \ + C \varepsilon \totMormax{N}(t,u),
			\notag
	\end{align}
		
	\begin{align} \label{E:EASYERRORINTEGRANDTWOINTEGRALESTIMATE}
		& \int_{\Sigma_t^u} 
			\left|
				\left\lbrace
					\uLgood [\rgeo^2 \mytr \upchi]
					- \frac{1}{2} \rgeo^2  \upmu (\mytr \upchi)^2
					+ \rgeo^2 \mytr \upchi \angkuparg{(Trans-\Psi)}
					+ \rgeo^2 \upmu \mytr \upchi \angkuparg{(Tan-\Psi)}
				\right\rbrace
		\right|
		(\mathscr{Z}^N \Psi)^2
		\, d \tvol
			\\
		& \leq
			C
			\ln(\myexp + t)
			\totzeromax{N}(t,u),
			\notag
	\end{align}
	
	\begin{align} \label{E:EASYERRORINTEGRANDONEINTEGRALESTIMATE}
		\int_{\mathcal{M}_{t,u}} 
			\upmu 
			\left\lbrace 
				\square_{g(\Psi)} [\rgeo^2 \mytr \upchi] 
			\right\rbrace 
			(\mathscr{Z}^N \Psi)^2
		\, d \vol
		& \leq
			C
			\ln(\myexp + t)
			\totzeromax{N}(t,u).
	\end{align}
	\end{subequations}
\end{lemma}

\begin{remark}[\textbf{Isolating the Morawetz term $-\frac{1}{2} \rgeo^2 |\angdiff \Psi|^2 [\Lunit \upmu]_-$}]
	Recall that 
	we have isolated the spacetime integral of
	$-\frac{1}{2} \rgeo^2 |\angdiff \mathscr{Z}^N \Psi|^2 [\Lunit \upmu]_-$
	by ``removing'' it from the left-hand side of \eqref{E:Q1BASICERRORINTEGRALESTIMATE};
	this isolated term generates the coercive Morawetz spacetime integral
	\eqref{E:COERCIVEMORDEF}.
\end{remark}

\begin{proof}
To prove \eqref{E:Q0BASICERRORINTEGRALESTIMATE} and \eqref{E:Q1BASICERRORINTEGRALESTIMATE},
we integrate inequalities 
\eqref{E:MULTERRORINTEGRANDPOINTWISE}
and 
\eqref{E:MORAWETZERRORINTEGRANDPOINTWISE}
(with $\mathscr{Z}^N \Psi$ in the role of $\Psi$)
over the spacetime region $\mathcal{M}_{t,u}.$
The vast majority of the terms can be suitably bounded by
using Prop.~\ref{P:L2NORMSOFPSIINTERMSOFTHECOERCIVEQUANTITIES}
and the fact that
$\totzeromax{\leq N}$ and $\totonemax{\leq N}$ are increasing in their arguments.
There are only two terms that are slightly subtle to estimate. 
The first subtle term is the term
$\varepsilon \ln^2(\myexp + t') 
		\mathbf{1}_{\lbrace \upmu \leq 1/4 \rbrace} 
		\left|\angdiff \mathscr{Z}^N \Psi\right|^2$
from the right-hand side of \eqref{E:MULTERRORINTEGRANDPOINTWISE}.
To bound the corresponding spacetime integral,
we use the estimate \eqref{E:MORAWETZANGULARDERIVATIVESFRAMEDERIVATIVESINTERMSOFQ0ANDQ1}
to deduce that
\[
\varepsilon 
  \int_{\mathcal{M}_{t,u}}
  	\ln^2(\myexp + t') 
		\mathbf{1}_{\lbrace \upmu \leq 1/4 \rbrace} 
		\left|\angdiff \mathscr{Z}^N \Psi\right|^2
	\, d \vol
	\leq 
	C \varepsilon 
	\sup_{t' \in [0,t]}\frac{\totMormax{N}(t',u)}{(1 + t')^{1/2}}		
\]
as desired.
The second subtle integral is
$ \int_{\mathcal{M}_{t,u}}
		\varepsilon 
 		\frac{(1 + t')}{\ln(\myexp + t')} 
 		\mathbf{1}_{\lbrace \upmu \leq 1/4 \rbrace} 
 		|\angdiff \mathscr{Z}^N \Psi|^2
 	\, d \vol,
$
which is generated by the last term on the right-hand side of \eqref{E:MORAWETZERRORINTEGRANDPOINTWISE}.   
To bound the integral by $\lesssim \varepsilon \totMormax{N}(t,u),$
we simply use
Lemma~\ref{L:MORAWETZSPACETIMECOERCIVITY}.

The \emph{hypersurface integral} inequality \eqref{E:EASYERRORINTEGRANDTWOINTEGRALESTIMATE}
and the spacetime integral inequality \eqref{E:EASYERRORINTEGRANDONEINTEGRALESTIMATE} 
follow easily from the pointwise estimates
\eqref{E:EASYERRORINTEGRANDONEPOINTWISE} and \eqref{E:EASYERRORINTEGRANDTWOPOINTWISE},
Prop.~\ref{P:L2NORMSOFPSIINTERMSOFTHECOERCIVEQUANTITIES},
and the fact that $\totzeromax{\leq N}$ is increasing in its arguments.
	
\end{proof}

\section{Difficult top-order error integral estimates} \label{S:DIFFICULTTOPORDERINTEGRALESTIMATES}

We now derive suitable estimates for the difficult top-order error integrals,
which are generated when we commute the wave equation $\upmu \square_{g(\Psi)} \Psi = 0$ 
with a top-order
pure spatial commutator vectorfield operator
$\mathscr{S}^N$
(see Def.~\ref{D:DEFSETOFSPATIALCOMMUTATORVECTORFIELDS}).
These are the most important estimates in the monograph.
We state them in the next two ``primary'' lemmas. 
In order to break up the proofs into manageable pieces, 
we devote the remainder of the present section to proving a series auxiliary lemmas. 
We then give the proof of the two lemmas of primary interest,
along with Lemma~\ref{L:TOPORDERELLIPTICRECOVERY} above,
in Sects.~\ref{S:PROOFOFLEMMATOPORDERELLIPTICRECOVERY}-\ref{S:PROOFOFLEMMADANGEROUSTOPORDERMORERRORINTEGRAL}.

\subsection{Statement of the two primary lemmas}
We now state the two lemmas of primary interest.

\begin{lemma}[\textbf{The main estimates for the difficult top-order error spacetime integrals
corresponding to the multiplier $\Mult$}]
\label{L:DANGEROUSTOPORDERMULTERRORINTEGRAL}
Let $1 \leq N \leq 24$ be an integer, 
let $\mathscr{S}^{N-1}$ be an $(N-1)^{st}$ order pure spatial commutation vectorfield operator
(see Def.~\ref{D:DEFSETOFSPATIALCOMMUTATORVECTORFIELDS}),
and let $\varsigma > 0$ be a number.
Under the small-data and bootstrap assumptions 
of Sects.~\ref{S:PSISOLVES}-\ref{S:C0BOUNDBOOTSTRAP},
there exist a small constant $\Littleconone > 0$ and a large constant $C > 0$ such that
if $\varepsilon$ is sufficiently small, 
then the following estimates hold 
for $(t,u) \in [0,\Tboot) \times [0,U_0]$
(and the constants are \textbf{independent of} $\varsigma$):
\begin{align} \label{E:MULTMAINENERGYFLUXERRORINTEGRALESTIMATE}
	& \left|
			\int_{\mathcal{M}_{t,u}}
				\left\lbrace
					(1 + 2 \upmu) \Lunit \mathscr{S}^{N-1} \Rad \Psi
					+ 2 \Rad \mathscr{S}^{N-1} \Rad \Psi
				\right\rbrace
				(\Rad \Psi) \angLap \mathscr{S}^{N-1} \upmu
				\, d \vol
		 \right|,
		\\
	& \left|
			\int_{\mathcal{M}_{t,u}}
				\left\lbrace
					(1 + 2 \upmu) \Lunit \mathscr{S}^{N-1} \Rot \Psi
					+ 2 \Rad \mathscr{S}^{N-1} \Rot \Psi
				\right\rbrace
				(\Rad \Psi) \Rot \mathscr{S}^{N-1} \mytr \upchi^{(Small)}
				\, d \vol
		 \right|
		\notag \\
	& \leq 
			C \varepsilon \totzeromax{\leq N}(t,u)
			+ C \varepsilon^3 \frac{1}{\upmu_{\star}(t,u)}
			\notag \\
	& \ \ + C \varepsilon \ln^8(\myexp + t)  \frac{1}{\upmu_{\star}(t,u)} \totzeromax{\leq N-1}(t,u)
			+ C \varepsilon \frac{1}{\upmu_{\star}(t,u)} \totonemax{\leq N-1}(t,u)
			\notag \\
	& \ \  +	\boxed{9}
							\int_{t'=0}^t 
							\frac{\| [\Lunit \upmu]_- \|_{C^0(\Sigma_{t'}^u)}} 
									 {\upmu_{\star}(t',u)} 
							\totzeromax{\leq N}^{1/2}(t',u) 
						\int_{s=0}^{t'}
							\frac{\| [\Lunit \upmu]_- \|_{C^0(\Sigma_s^u)}} 
									{\upmu_{\star}(s,u)} 
							\totzeromax{\leq N}^{1/2}(s,u) 
						\, dt'
						\, ds	
						\notag \\		
		& \ \ + 
			C \varepsilon
				\int_{t'=0}^t
						\frac{1} 
								 {(1 + t')^{1 + \Littleconone} \upmu_{\star}(t',u)} 
						\totzeromax{\leq N}^{1/2}(t',u) 
						\int_{s=0}^{t'}
							\frac{1}{(1 + s)}
							\frac{1}{\upmu_{\star}(s,u)} 
							\totzeromax{\leq N}^{1/2}(s,u) 
						\, ds
				\, dt'
				\notag \\
		& \ \ + 
			C \varepsilon
				\int_{t'=0}^t
						\frac{1} 
								 {(1 + t')^{1 + \Littleconone} \upmu_{\star}(t',u)} 
						\totzeromax{\leq N}^{1/2}(t',u) 
						\int_{s=0}^{t'}
							\frac{1}{(1 + s)}
							\frac{1}{\upmu_{\star}(s,u)} 
							\totonemax{\leq N}^{1/2}(s,u) 
						\, ds
				\, dt'
				\notag \\
		& \ \ +   \boxed{9} 
							\int_{t'=0}^t 
								\frac{\| [\Lunit \upmu]_- \|_{C^0(\Sigma_{t'}^u)}} 
									   {\upmu_{\star}(t',u)} 
								\totzeromax{\leq N}(t',u)
							\, dt'	 
						\notag \\	
		& \ \ + 	9
							\int_{t'=0}^t 
							 \frac{1}{\rgeo(t',u) \left\lbrace 1 + \ln \left(\frac{\rgeo(t',u)}{\rgeo(0,u)} \right) \right\rbrace}
							 \totzeromax{\leq N}^{1/2}(t',u)
						\int_{s=0}^{t'} 
							\frac{\| [\Lunit \upmu]_- \|_{C^0(\Sigma_s^u)}} 
									{\upmu_{\star}(s,u)} 
							\totzeromax{\leq N}^{1/2}(s,u) 
						\, ds
						\, dt'	
						\notag \\	
		& \ \ +   9 
							\int_{t'=0}^t 
								\frac{1}{\rgeo(t',u) \left\lbrace 1 + \ln \left(\frac{\rgeo(t',u)}{\rgeo(0,u)} \right) \right\rbrace} 
								\totzeromax{\leq N}(t',u) 
							\, dt'	
							\notag \\
	& \ \ + C \varepsilon
				 \int_{t'=0}^t
					\frac{1} 
							 {(1 + t')^{1+\Littleconone} \upmu_{\star}(t',u)} 
				  \totzeromax{\leq N}(t',u)
				\, dt'
			\notag \\
	& \ \ + C \varepsilon
				\int_{t'=0}^t
					\frac{1} 
							 {(1 + t')^{1+\Littleconone} \upmu_{\star}(t',u)} 
				  \totonemax{\leq N}(t',u)
				\, dt'
		  \notag \\
	& \ \ + C 
				 (1 + \varsigma^{-1})
				 \int_{t'=0}^t
					\frac{1} 
							 {(1 + t')^{3/2}} 
				  \totzeromax{\leq N}(t',u)
				\, dt'
			\notag \\
	& \ \ + C 
				 (1 + \varsigma^{-1})
				 \int_{t'=0}^t
					\frac{1} 
							 {(1 + t')^{3/2}} 
				  \totonemax{\leq N}(t',u)
				\, dt'
		  \notag \\
	& \ \
		+ C (1 + \varsigma^{-1})
			\frac{1}{(1 + t)^{1/2}}
			\int_{u'=0}^u
				\totonemax{\leq N}(t,u') 
			\, du'
		+ C \varsigma 
			\sup_{t' \in [0,t]} 
			\frac{\totMormax{\leq N}(t',u)}{(1 + t')^{1/2}}.
			\notag
\end{align}	
\end{lemma}
We provide the proof of Lemma~\ref{L:DANGEROUSTOPORDERMULTERRORINTEGRAL}
in Sect.~\ref{S:PROOFOFLEMMADANGEROUSTOPORDERMULTERRORINTEGRAL}.

\begin{lemma}[\textbf{The main estimates for the difficult top-order error spacetime integrals
corresponding to the multiplier $\Mor$}]
\label{L:DANGEROUSTOPORDERMORERRORINTEGRAL}
Let $1 \leq N \leq 24$ be an integer, 
let $\mathscr{S}^{N-1}$ be an $(N-1)^{st}$ order pure spatial commutation vectorfield operator
(see Def.~\ref{D:DEFSETOFSPATIALCOMMUTATORVECTORFIELDS}),
and let $\widetilde{\varsigma} > 0$ be a number.
Under the small-data and bootstrap assumptions 
of Sects.~\ref{S:PSISOLVES}-\ref{S:C0BOUNDBOOTSTRAP},
there exist a small constant $a > 0$ and a large constant $C > 0$ such that
if $\varepsilon$ is sufficiently small, 
then the following estimates hold 
for  $(t,u) \in [0,\Tboot) \times [0,U_0]$
(and the constants are \textbf{independent of} $\widetilde{\varsigma}$):
\begin{align} \label{E:MORMAINENERGYFLUXERRORINTEGRALESTIMATE}
	& \left|
			\int_{\mathcal{M}_{t,u}}
				\rgeo^2
				\left\lbrace
					\Lunit \mathscr{S}^{N-1} \Rad \Psi 
					+ \frac{1}{2} \mytr \upchi \mathscr{S}^{N-1} \Rad \Psi
				\right\rbrace
				(\Rad \Psi) \angLap \mathscr{S}^{N-1} \upmu
				\, d \vol
		 \right|,
		\\
	& \left|
			\int_{\mathcal{M}_{t,u}}
				\rgeo^2
				\left\lbrace
					\Lunit \mathscr{S}^{N-1} \Rot \Psi 
					+ \frac{1}{2} \mytr \upchi \mathscr{S}^{N-1} \Rot \Psi
				\right\rbrace
				(\Rad \Psi) \Rot \mathscr{S}^{N-1} \mytr \upchi^{(Small)}
				\, d \vol
		 \right|
		\notag \\
	& \leq 
			C \varepsilon \totzeromax{\leq N}(t,u)
			+ C \varepsilon \totonemax{\leq N}(t,u)
			+ C \varepsilon^3 \frac{1}{\upmu_{\star}(t,u)}
			\notag \\
	& \ \ + C \varepsilon \frac{1}{\upmu_{\star}(t,u)} \ln^2(\myexp + t) \totzeromax{\leq N-1}(t,u)
			+ C \varepsilon \frac{1}{\upmu_{\star}(t,u)} \ln^2(\myexp + t) \totonemax{\leq N-1}(t,u)
			\notag \\
	& \ \				+ \boxed{5} 
							\int_{t'=0}^t 
								\frac{\| [\Lunit \upmu]_- \|_{C^0(\Sigma_{t'}^u)}} 
									   {\upmu_{\star}(t',u)} 
								\totonemax{\leq N}(t',u)
							\, dt'
			\notag \\		
	& \ \ 		+ 
							5
							\int_{t'=0}^t 
								\frac{1}{\rgeo(t',u) \left\lbrace 1 +  \ln \left(\frac{\rgeo(t',u)}{\rgeo(0,u)} \right) \right\rbrace}
								\totonemax{\leq N}(t',u)
							\, dt'
							\notag \\
	 & \ \ + C \varepsilon
	  					\int_{t'=0}^t
	  						\frac{1}{(1 + t')^{3/2} \upmu_{\star}^{1/2}(t',u)} 
								\totzeromax{\leq N}(t',u)
	  					\, dt'
			\notag \\				
	& \ \ +  C \varepsilon
	  					\int_{t'=0}^t
	  						\frac{1}{(1 + t')^{3/2} \upmu_{\star}(t',u)} 
	  						\totonemax{\leq N}(t',u)
	  					\, dt'
	  					\notag \\
	 & \ \ + C (1 + \widetilde{\varsigma}^{-1})
	  					\int_{t'=0}^t
	  						\frac{1}{(1 + t')^{3/2}} 
								\totzeromax{\leq N}(t',u)
	  					\, dt'
			\notag \\				
	& \ \ +  C (1 + \widetilde{\varsigma}^{-1})
	  					\int_{t'=0}^t
	  						\frac{1}{(1 + t')^{3/2}} 
	  						\totonemax{\leq N}(t',u)
	  					\, dt'
	  					\notag \\
		& \ \ + C \ln^4(\myexp + t) \totzeromax{N}(t,u)
					\notag \\
		& \ \ 
		+ \boxed{5}
		\frac{\| \Lunit \upmu \|_{C^0(\Sigmaminus{t}{t}{u})}} 
			   {\upmu_{\star}^{1/2}(t,u)}
		\totonemax{\leq N}^{1/2}(t,u)
		\int_{t' = 0}^t
			\frac{1}{\upmu_{\star}^{1/2}(t',u)}
			\totonemax{\leq N}^{1/2}(t',u)
		 \, dt'
		\notag \\
	& \ \
		+ 5
		\left\| 
			\frac{\Lunit \upmu}{\upmu} 
		\right \|_{C^0(\Sigmaplus{t}{t}{u})} 
		\totonemax{\leq N}^{1/2}(t,u)
		\int_{t' = 0}^t
			\left\|
					\sqrt{
					\frac{\upmu(t,\cdot)} 
					     {\upmu}
					     }
			\right\|_{C^0(\Sigmaplus{t'}{t}{u})}
			\totonemax{\leq N}^{1/2}(t',u)
		 \, dt'
		 \notag \\
	& \ \ + C \varepsilon \frac{1}{\upmu_{\star}^{1/2}(t,u)}
			\totonemax{\leq N}^{1/2}(t,u)	
			\int_{t'=0}^t
				\frac{1}{(1 + t')^{3/2}} \totzeromax{\leq N}^{1/2}(t',u)
			\, dt'
			\notag \\
	& \ \ + C \varepsilon \frac{1}{\upmu_{\star}^{1/2}(t,u)}
				\totonemax{\leq N}^{1/2}(t,u)
				\int_{t'=0}^t
					\frac{1}{(1 + t')^{3/2} \upmu_{\star}^{1/2}(t',u)} \totonemax{\leq N}^{1/2}(t',u)
				\, dt'
				\notag \\
		& \ \ + C (1 + \widetilde{\varsigma}^{-1})
						\int_{u'=0}^u
								\totonemax{\leq N}(t,u')		
						\, du' 
				+ C \widetilde{\varsigma} \totMormax{\leq N}(t,u).
					\notag 
		\end{align}	
\end{lemma}
We provide the proof of Lemma~\ref{L:DANGEROUSTOPORDERMORERRORINTEGRAL}
in Sect.~\ref{S:PROOFOFLEMMADANGEROUSTOPORDERMORERRORINTEGRAL}.

\subsection{Auxiliary lemmas}
We now state and prove a series of auxiliary lemmas 
in order to help us prove the previous two lemmas
as well as Lemma~\ref{L:TOPORDERELLIPTICRECOVERY} above.

\begin{lemma}[\textbf{Preliminary} $L^2$ \textbf{bounds for some top-order derivatives of} 
$\upchi^{(Small)}$ and $\Lunit_{(Small)}^i$]
	\label{L:L2COMMUTEDANGDLIECHIJUNKINTERMSOFTRACEFREEANDTRACEPART}
	Let $1 \leq N \leq 24$ be an integer.
	Under the small-data and bootstrap assumptions 
	of Sects.~\ref{S:PSISOLVES}-\ref{S:C0BOUNDBOOTSTRAP},
	if $\varepsilon$ is sufficiently small,  
	then the following estimates hold for $(t,u) \in [0,\Tboot) \times [0,U_0]:$
	\begin{align}  \label{E:L2COMMUTEDANGDLIECHIJUNKINTERMSOFTRACEFREEANDTRACEPART}
		\left\| \rgeo \upmu \angD \angLie_{\mathscr{Z}}^{N-1} \upchi^{(Small)} \right\|_{L^2(\Sigma_t^u)}
		& \lesssim
				\sum_{l=1}^3
				\left\| \upmu \Rot_{(l)} \mathscr{Z}^{N-1} \mytr \upchi^{(Small)} \right\|_{L^2(\Sigma_t^u)}
			+ \left\| \rgeo \upmu \angD \angfreeLietwoarg{\mathscr{Z}}{N-1} \hat{\upchi}^{(Small)} \right\|_{L^2(\Sigma_t^u)}
				\\
		& \ \ + \totzeromax{\leq N}^{1/2}(t,u)
			+ 
			\int_{t'=0}^t 
				\frac{1}{(1 + t')^{3/2} \upmu_{\star}^{1/2}(t',u)}
				\totonemax{\leq N}^{1/2}(t',u)
			\, d t'
			+ \varepsilon,
				\notag
	\end{align}
	\begin{align}  \label{E:L2COMMUTEDANGDIFFTRCHIJUNKINTERMSOFANGDCHIJUNK}
		& \left\| \rgeo \upmu \angdiff \mathscr{Z}^{N-1} \mytr \upchi^{(Small)} \right\|_{L^2(\Sigma_t^u)},
			\\ 
		& \left\| \upmu \Rot \mathscr{Z}^{N-1} \mytr \upchi^{(Small)} \right\|_{L^2(\Sigma_t^u)}
			\notag \\
		& \lesssim
			\left\| \rgeo \upmu \angD \angLie_{\mathscr{Z}}^{N-1} \upchi^{(Small)} \right\|_{L^2(\Sigma_t^u)}
			+ \totzeromax{\leq N}^{1/2}(t,u)
			+
			\int_{t'=0}^t 
				\frac{1}{(1 + t')^{3/2} \upmu_{\star}^{1/2}(t',u)}
				\totonemax{\leq N}^{1/2}(t',u)
			\, d t'
			+ \varepsilon,
				\notag
	\end{align}
	\begin{align}  \label{E:L2COMMUTEDANGDTRACEFREECHIJUNKINTERMSOFCHIJUNK}
		\left\| \rgeo \upmu \angD \angLie_{\mathscr{Z}}^{N-1} \hat{\upchi}^{(Small)} \right\|_{L^2(\Sigma_t^u)}
		& \lesssim
				\left\| \rgeo \upmu \angD \angLie_{\mathscr{Z}}^{N-1} \upchi^{(Small)} \right\|_{L^2(\Sigma_t^u)}
				\\
		& \ \ + \totzeromax{\leq N}^{1/2}(t,u)
			+ 
			\int_{t'=0}^t 
				\frac{1}{(1 + t')^{3/2} \upmu_{\star}^{1/2}(t',u)}
				\totonemax{\leq N}^{1/2}(t',u)
			\, d t'
			+ \varepsilon,
				\notag
	\end{align}
	\begin{align}  \label{E:L2COMMUTEDANGDTRACEFREECHIJUNKINTERMSOFTRACEPART} 
		\left\| \rgeo \upmu \angdiv \angfreeLietwoarg{\mathscr{Z}}{N-1} \hat{\upchi}^{(Small)} \right\|_{L^2(\Sigma_t^u)}
		& \lesssim
				\sum_{l=1}^3
				\left\| \upmu \Rot_{(l)} \mathscr{Z}^{N-1} \mytr \upchi^{(Small)} \right\|_{L^2(\Sigma_t^u)}
				\\
		& \ \ + \totzeromax{\leq N}^{1/2}(t,u)
			+ 
			\int_{t'=0}^t 
				\frac{1}{(1 + t')^{3/2} \upmu_{\star}^{1/2}(t',u)}
				\totonemax{\leq N}^{1/2}(t',u)
			\, d t'
			+ \varepsilon,
				\notag
\end{align}
\begin{align}  \label{E:L2LOWERORDERCOMMUTEDANGDTRACEFREECHIJUNKINTERMSOFTRACEPART} 
		\left\| \rgeo \angfreeLietwoarg{\mathscr{Z}}{N-1} \hat{\upchi}^{(Small)} \right\|_{L^2(\Sigma_t^u)}
		& \lesssim
			\ln(\myexp + t)
			\totzeromax{\leq N}^{1/2}(t,u)
			+ 
			\int_{t'=0}^t 
				\frac{1}{(1 + t') \upmu_{\star}^{1/2}(t',u)}
				\totonemax{\leq N}^{1/2}(t',u)
			\, d t'
			+ \varepsilon.
\end{align}

Furthermore, the following estimate holds for $(t,u) \in [0,\Tboot) \times [0,U_0]:$ 
\begin{align}  \label{E:L2COMMUTEDLIERADCHIJUNKINTERMSOFANGDSQUAREUPMUANDPSI}
		\left \|\rgeo \upmu \angLie_{\Rad} \angLie_{\mathscr{Z}}^{N-1} \upchi^{(Small)} \right \|_{L^2(\Sigma_t^u)}
		& \lesssim 
			\left\| \upmu \angD^2 \mathscr{Z}^{N-1} \upmu \right\|_{L^2(\Sigma_t^u)}
				\\
		& \ \ 
			+ \totzeromax{\leq N}^{1/2}(t,u)
			+ 
			\int_{t'=0}^t 
				\frac{1}{(1 + t')^{3/2} \upmu_{\star}^{1/2}(t',u)}
				\totonemax{\leq N}^{1/2}(t',u)
			\, d t'
			+ \varepsilon.
			\notag
\end{align}

Finally, the following estimate holds for $(t,u) \in [0,\Tboot) \times [0,U_0]$ (for $i=1,2,3$):
\begin{align}  \label{E:L2COMMUTEDANGLAPLJUNKIINTERMSOFANGDCHIANDPSI}
		\left \|\rgeo \upmu \angLap \mathscr{Z}^{\leq N-1} \Lunit_{(Small)}^i \right \|_{L^2(\Sigma_t^u)}
		& \lesssim 
			\sum_{l=1}^3 \left\| \upmu \Rot_{(l)} \mathscr{Z}^{N-1} \mytr \upchi^{(Small)} \right\|_{L^2(\Sigma_t^u)}
				\\
		& \ \ 
			+ \totzeromax{\leq N}^{1/2}(t,u)
			+ 
			\int_{t'=0}^t 
				\frac{1}{(1 + t')^{3/2} \upmu_{\star}^{1/2}(t',u)}
				\totonemax{\leq N}^{1/2}(t',u)
			\, d t'
			+ \varepsilon.
			\notag
\end{align}
\end{lemma}
\begin{proof}
We first prove \eqref{E:L2COMMUTEDANGDLIECHIJUNKINTERMSOFTRACEFREEANDTRACEPART}.
To this end, we multiply both sides of inequalities 
\eqref{E:POINTWISECOMMUTEDANGDCHIJUNKINTERMSOFTRACEFREEANDTRACEPART}
and 
\eqref{E:POINTWISECOMMUTEDANGDTRACEFREECHIJUNKINTERMSOFFREECOMMUTEDANGDTRACEFREECHI}
by $\rgeo \upmu$ and then take the $L^2(\Sigma_t^u)$ norm of both sides.
Using Prop.~\ref{P:L2NORMSOFPSIINTERMSOFTHECOERCIVEQUANTITIES},
inequality \eqref{E:EIKONALFUNCTIONQUANTITIESL2BOUNDSINTERMSOFQ0ANDQ1},
the estimate $\rgeo \upmu \lesssim (1 + t) \ln(\myexp + t)$
(that is, \eqref{E:C0BOUNDCRUCIALEIKONALFUNCTIONQUANTITIES}),
and the fact that $\totzeromax{\leq N}$ is increasing in its arguments,
we deduce that the products of $\rgeo \upmu$ and the right-hand sides of
\eqref{E:POINTWISECOMMUTEDANGDCHIJUNKINTERMSOFTRACEFREEANDTRACEPART}
and 
\eqref{E:POINTWISECOMMUTEDANGDTRACEFREECHIJUNKINTERMSOFFREECOMMUTEDANGDTRACEFREECHI}
are bounded in the norm $\| \cdot \|_{L^2(\Sigma_t^u)}$ by
\begin{align} \label{E:EASYTERMSL2COMMUTEDANGDLIECHIJUNKINTERMSOFTRACEFREEANDTRACEPART}
	& \lesssim 
		\totzeromax{\leq N}^{1/2}(t,u)
			+ 
			\int_{t'=0}^t 
				\frac{1}{(1 + t')^{3/2} \upmu_{\star}^{1/2}(t',u)}
				\totonemax{\leq N}^{1/2}(t',u)
			\, d t'
			+ \varepsilon.
\end{align}
We note that the right-hand side
of \eqref{E:EASYTERMSL2COMMUTEDANGDLIECHIJUNKINTERMSOFTRACEFREEANDTRACEPART}
is manifestly $\lesssim$ the right-hand side of 
\eqref{E:L2COMMUTEDANGDLIECHIJUNKINTERMSOFTRACEFREEANDTRACEPART}.
The desired estimate
\eqref{E:L2COMMUTEDANGDLIECHIJUNKINTERMSOFTRACEFREEANDTRACEPART}
now follows from 
inequalities 
\eqref{E:POINTWISECOMMUTEDANGDCHIJUNKINTERMSOFTRACEFREEANDTRACEPART}
and 
\eqref{E:POINTWISECOMMUTEDANGDTRACEFREECHIJUNKINTERMSOFFREECOMMUTEDANGDTRACEFREECHI},
the triangle inequality, and the estimate \eqref{E:FUNCTIONPOINTWISEANGDINTERMSOFANGLIEO},
which implies that 
$\left\| \upmu \rgeo \angdiff \mathscr{Z}^{N-1} \mytr \upchi^{(Small)}) \right\|_{L^2(\Sigma_t^u)}
\lesssim 
\sum_{l=1}^3
\left\| \upmu \Rot_{(l)} \mathscr{Z}^{N-1} \mytr \upchi^{(Small)} \right\|_{L^2(\Sigma_t^u)}.$

The estimate \eqref{E:L2COMMUTEDANGDIFFTRCHIJUNKINTERMSOFANGDCHIJUNK} for
$\left\| \upmu \Rot \mathscr{Z}^{N-1} \mytr \upchi^{(Small)} \right\|_{L^2(\Sigma_t^u)}$
can be proved in a similar fashion
with the help of the pointwise inequality \eqref{E:POINTWISEESTIMATEROTTRCHIJUNKINTERMSOFANGDCHIJUNKPLUSJUNK}.
The estimate \eqref{E:L2COMMUTEDANGDIFFTRCHIJUNKINTERMSOFANGDCHIJUNK} for
$\left\| \rgeo \upmu \angdiff \mathscr{Z}^{N-1} \mytr \upchi^{(Small)} \right\|_{L^2(\Sigma_t^u)}$
then follows from the estimate for
$\left\| \upmu \Rot \mathscr{Z}^{N-1} \mytr \upchi^{(Small)} \right\|_{L^2(\Sigma_t^u)}$
and the inequality
$\left\| \upmu \rgeo \angdiff \mathscr{Z}^{N-1} \mytr \upchi^{(Small)}) \right\|_{L^2(\Sigma_t^u)}
\lesssim 
\sum_{l=1}^3
\left\| \upmu \Rot_{(l)} \mathscr{Z}^{N-1} \mytr \upchi^{(Small)} \right\|_{L^2(\Sigma_t^u)}$
noted above.

To prove the estimate \eqref{E:L2COMMUTEDANGDTRACEFREECHIJUNKINTERMSOFCHIJUNK}, we again
base our argument on the pointwise inequality 
\eqref{E:POINTWISECOMMUTEDANGDCHIJUNKINTERMSOFTRACEFREEANDTRACEPART}, but this time we 
use the triangle inequality to bound
$\left\| \rgeo \upmu \angD \angLie_{\mathscr{Z}}^{N-1} \hat{\upchi}^{(Small)} \right\|_{L^2(\Sigma_t^u)}$
in terms of the remaining quantities in \eqref{E:POINTWISECOMMUTEDANGDCHIJUNKINTERMSOFTRACEFREEANDTRACEPART}, 
which we have already shown to be bounded
in the norm $\| \cdot \|_{L^2(\Sigma_t^u)}$
by the right-hand side of \eqref{E:L2COMMUTEDANGDTRACEFREECHIJUNKINTERMSOFCHIJUNK}. 

To prove
\eqref{E:L2COMMUTEDANGDTRACEFREECHIJUNKINTERMSOFTRACEPART},
we apply similar reasoning to the fourth pointwise inequality
in \eqref{E:DIVCHIJUNKANGDIFFTRCHIJUNKHIGHERORDERCOMMMUTOR}.

To prove \eqref{E:L2LOWERORDERCOMMUTEDANGDTRACEFREECHIJUNKINTERMSOFTRACEPART},
we apply similar reasoning to the pointwise inequality
\eqref{E:POINTWISELOWERORDERANGREELIETRACEFREECHIJUNKISJUNK}
with $N-1$ in the role of $N.$

To prove \eqref{E:L2COMMUTEDLIERADCHIJUNKINTERMSOFANGDSQUAREUPMUANDPSI},
we apply similar reasoning to the first pointwise inequality 
in \eqref{E:TOPORDERDERIVATIVESOFANGDSQUAREDUPMUINTERMSOFCONTROLLABLE}.

To prove \eqref{E:L2COMMUTEDANGLAPLJUNKIINTERMSOFANGDCHIANDPSI}, we apply similar
reasoning to the pointwise inequality \eqref{E:POINTWISEBOUNDFORANGLAPLJUNKICOMMUTED}.
As a preliminary step, we use the third pointwise inequality in \eqref{E:DIVCHIJUNKANGDIFFTRCHIJUNKHIGHERORDERCOMMMUTOR} to
bound the last term on the right-hand side of \eqref{E:POINTWISEBOUNDFORANGLAPLJUNKICOMMUTED}
in terms of $\angdiff \mathscr{Z}^{N-1} \mytr \upchi^{(Small)}$
and the error terms on the right-hand side of \eqref{E:DIVCHIJUNKANGDIFFTRCHIJUNKHIGHERORDERCOMMMUTOR}.
\end{proof}

\begin{lemma}[\textbf{Preliminary non-sharp} $L^2$ \textbf{bounds for}
$\|\upmu \Rot \mathscr{Z}^{N-1} \mytr \upchi^{(Small)} \|_{L^2(\Sigma_t^u)}$
\textbf{and} 
$\upmu \angLap \mathscr{Z}^{N-1} \upmu$ ]
	\label{L:SLIGHTLYLESSSHARPL2FORANLGAPUPMUANDROTTRCHI}
	Let $1 \leq N \leq 24$ be an integer.
	Under the small-data and bootstrap assumptions 
	of Sects.~\ref{S:PSISOLVES}-\ref{S:C0BOUNDBOOTSTRAP},
	if $\varepsilon$ is sufficiently small, 
	then the following estimates hold for $(t,u) \in [0,\Tboot) \times [0,U_0]:$
	\begin{align}  \label{E:SLIGHTLYLESSSHARPL2FORANLGAPUPMUANDROTTRCHI}
		& \|\upmu \Rot \mathscr{Z}^{N-1} \mytr \upchi^{(Small)} \|_{L^2(\Sigma_t^u)},
			\,
			\| \upmu \angLap \mathscr{Z}^{N-1} \upmu \|_{L^2(\Sigma_t^u)}
			\\
		& \lesssim 
			\int_{t'=0}^t 
					\frac{\| [\Lunit \upmu]_- \|_{C^0(\Sigma_{t'}^u)}} 
							 {\upmu_{\star}(t',u)} 
					\totzeromax{\leq N}^{1/2}(t',u) 
				\, dt'
			+ 
			\int_{t'=0}^t 
				\frac{1}{(1 + t')^{3/2} \upmu_{\star}(t',u)}
				\totzeromax{\leq N}^{1/2}(t',u)
			\, d t'		
				\notag \\
		& \ \ + \totzeromax{\leq N}^{1/2}(t,u)
			+ 
			\int_{t'=0}^t 
				\frac{1}{(1 + t')^{3/2} \upmu_{\star}(t',u)}
				\totonemax{\leq N}^{1/2}(t',u)
			\, d t'
			\notag \\
		& \ \ 
			+ \varepsilon
				\int_{t'=0}^t
					\frac{1}{(1 + t')^{3/2}}
					\myarray[\left\| \upmu \angfreeDsquared \mathscr{Z}^{\leq N-1} \upmu \right\|_{L^2(\Sigma_{t'}^u)}]
						{\left\| \rgeo \upmu \angD \angLie_{\mathscr{Z}}^{\leq N-1} \hat{\upchi}^{(Small)} \right\|_{L^2(\Sigma_{t'}^u)}}
				\, dt'
			+ \varepsilon
				\left\lbrace
					\ln \upmu_{\star}^{-1}(t,u) + 1
				\right\rbrace.
				\notag
	\end{align}
	
\end{lemma}

\begin{proof}
	We first assume that $\mathscr{Z}^{N-1} = \mathscr{S}^{N-1}$ is an $(N-1)^{st}$ order pure spatial commutation vectorfield operator.
	Then to prove \eqref{E:SLIGHTLYLESSSHARPL2FORANLGAPUPMUANDROTTRCHI},
	it suffices to bound the $L^2(\Sigma_t^u)$ norm of the right-hand side of 
	\eqref{E:SLIGHTLYLESSSHARPPOINTWISEANLGAPUPMUANDROTTRCHI}	
	by the right-hand side of \eqref{E:SLIGHTLYLESSSHARPL2FORANLGAPUPMUANDROTTRCHI}.
	To bound the norm $\| \cdot \|_{L^2(\Sigma_t^u)}$ of the first term on the right-hand side of 
	\eqref{E:SLIGHTLYLESSSHARPPOINTWISEANLGAPUPMUANDROTTRCHI} by
	$\lesssim \totzeromax{\leq N}^{1/2}(t,u),$ we use Prop.~\ref{P:L2NORMSOFPSIINTERMSOFTHECOERCIVEQUANTITIES}.
	
	To bound the norm $\| \cdot \|_{L^2(\Sigma_t^u)}$ of the second term on the right-hand side of 
	\eqref{E:SLIGHTLYLESSSHARPPOINTWISEANLGAPUPMUANDROTTRCHI}
	by $\lesssim \totzeromax{\leq N}^{1/2}(t,u) 
		+ \int_{t'=0}^t 
				\frac{1}{(1 + t')^{3/2} \upmu_{\star}^{1/2}(t',u)}
				\totonemax{\leq N}^{1/2}(t',u)
			\, d t'
			+ \varepsilon,$
	we use inequality \eqref{E:EIKONALFUNCTIONQUANTITIESL2BOUNDSINTERMSOFQ0ANDQ1}
	and the fact that $\totzeromax{\leq N}$ is increasing in both of its arguments.
	
	To bound the norm $\| \cdot \|_{L^2(\Sigma_t^u)}$ of the third term on the right-hand side of 
	\eqref{E:SLIGHTLYLESSSHARPPOINTWISEANLGAPUPMUANDROTTRCHI} 
	by the first term on the right-hand side of \eqref{E:SLIGHTLYLESSSHARPL2FORANLGAPUPMUANDROTTRCHI},
	we use Lemma~\ref{L:L2NORMSOFTIMEINTEGRATEDFUNCTIONS}
	and Prop.~\ref{P:L2NORMSOFPSIINTERMSOFTHECOERCIVEQUANTITIES}.
	
	To bound the norm $\| \cdot \|_{L^2(\Sigma_t^u)}$ of the fourth term on the right-hand side of 
	\eqref{E:SLIGHTLYLESSSHARPPOINTWISEANLGAPUPMUANDROTTRCHI},
	we use Lemma~\ref{L:L2NORMSOFTIMEINTEGRATEDFUNCTIONS},
	the estimate \eqref{E:ANNLOYINGSQRTMUOVERMUINTEGRATEDBOUND} with $\Littlecontwo=1/2,$
	Prop.~\ref{P:L2NORMSOFPSIINTERMSOFTHECOERCIVEQUANTITIES},
	and the fact that $\totzeromax{\leq N}$ is increasing in both of its arguments
	to bound it by
	\begin{align}
		\lesssim 
		\totzeromax{\leq N}^{1/2}(t,u) 
		\frac{1}{\rgeo(t,u)}
		\int_{t'=0}^t 
				\left\|
					\left(\frac{\upmu(t,\cdot)}{\upmu}\right)^2
				\right\|_{C^0(\Sigma_{t'}^u)}
				\, dt' \lesssim \totzeromax{\leq N}^{1/2}(t,u)
	\end{align}
	as desired.
	
	To bound the norm $\| \cdot \|_{L^2(\Sigma_t^u)}$ of the fifth term 
	$\frac{\ln^2(\myexp + t)}{(1 + t)^2}
	|\chifullmod^{[N]}|(0,u,\vartheta)$ on the right-hand side of 
	\eqref{E:SLIGHTLYLESSSHARPPOINTWISEANLGAPUPMUANDROTTRCHI}
	(which is defined just below \eqref{E:RENORMALIZEDTOPORDERTRCHIJUNKTRANSPORTINVERTED}),
	we use the definition of $\chifullmod^{[N]},$ 
	inequality \eqref{E:SPHEREVOLUMEFORMCOMPARISON},
	and Lemma~\ref{E:SMALLINITIALSOBOLEVNORMS}
	to deduce that
	\begin{align} \label{E:RENORMALZEDTRCHIJUNKDATATERML2NORMATTIMET}
	\left\| \chifullmod^{[N]}(0,\cdot) \right\|_{L^2(\Sigma_t^u)}^2
	& = \int_{u'=0}^u 
				\int_{S_{t,u}}
				\left|
					\chifullmod^{[N]}
				\right|^2
				(0,u',\vartheta)
				\, d \argspherevol{(t,u',\vartheta)}
			\, du'
				\\
		& \lesssim
			\rgeo^2(t,u)
			\int_{u'=0}^u 
				\int_{S_{t,u}}
				\left|
					\chifullmod^{[N]}
				\right|^2
				(0,u',\vartheta)
				\, d \argspherevol{(0,u',\vartheta)}
			\, du'
			\notag \\
		& = \rgeo^2(t,u) \left\| \chifullmod^{[N]}(0,\cdot) \right\|_{L^2(\Sigma_0^u)}^2
			\lesssim \varepsilon^2 \rgeo^2(t,u).
			\notag
\end{align}
From \eqref{E:RENORMALZEDTRCHIJUNKDATATERML2NORMATTIMET}, it follows that
$\frac{\ln^2(\myexp + t)}{(1 + t)^2} \| \chifullmod^{[N]}(0,\cdot)\|_{L^2(\Sigma_t^u)} 
\lesssim \varepsilon \frac{\ln^2(\myexp + t)}{1 + t}$ as desired.

	To bound the norm $\| \cdot \|_{L^2(\Sigma_t^u)}$ of the sixth term on the right-hand side of 
	\eqref{E:SLIGHTLYLESSSHARPPOINTWISEANLGAPUPMUANDROTTRCHI}, we use 
	Lemma~\ref{L:L2NORMSOFTIMEINTEGRATEDFUNCTIONS},
	Prop.~\ref{P:L2NORMSOFPSIINTERMSOFTHECOERCIVEQUANTITIES},
	the estimate \eqref{E:LOGLOSSLESSSINGULARTERMSMTHREEFOURTHSINTEGRALBOUND},
	and the fact that $\totzeromax{\leq N}$ and $\totonemax{\leq N}$
	are increasing in both of their arguments
	to bound it by
	\begin{align}
		& \lesssim 
		\frac{\ln^3(\myexp + t)}{1 + t}
		\int_{t'=0}^t 
			\frac{1}{(1 + t') \upmu_{\star}^{1/2}(t',u)} \totzeromax{\leq N}^{1/2}(t',u)
		\, dt'
		+
		\frac{\ln^3(\myexp + t)}{1 + t}
		\int_{t'=0}^t 
			\frac{1}{(1 + t') \upmu_{\star}^{1/2}(t',u)} \totonemax{\leq N}^{1/2}(t',u)
		\, dt'
			\\
		& \lesssim 
			\totzeromax{\leq N}^{1/2}(t,u)
			+
		\int_{t'=0}^t 
			\frac{1}{(1 + t')^{3/2} \upmu_{\star}^{1/2}(t',u)} \totonemax{\leq N}^{1/2}(t',u)
		\, dt'
		\notag
	\end{align}
	as desired.
	
	To bound the norm $\| \cdot \|_{L^2(\Sigma_t^u)}$ of the seventh term on the right-hand side of 
	\eqref{E:SLIGHTLYLESSSHARPPOINTWISEANLGAPUPMUANDROTTRCHI} 
	by the next-to-last term on the right-hand side of \eqref{E:SLIGHTLYLESSSHARPL2FORANLGAPUPMUANDROTTRCHI},
	we use Lemma~\ref{L:L2NORMSOFTIMEINTEGRATEDFUNCTIONS}.
	
	To bound the norm $\| \cdot \|_{L^2(\Sigma_t^u)}$ of the final term on the right-hand side of 
	\eqref{E:SLIGHTLYLESSSHARPPOINTWISEANLGAPUPMUANDROTTRCHI},
	we use Lemma~\ref{L:L2NORMSOFTIMEINTEGRATEDFUNCTIONS}, 
	the estimate \eqref{E:EIKONALFUNCTIONQUANTITIESL2BOUNDSINTERMSOFQ0ANDQ1}, 
	the estimates 
	\eqref{E:LOGLOSSMUINVERSEINTEGRALBOUND}
	and
	\eqref{E:LOGLOSSLESSSINGULARTERMSMTHREEFOURTHSINTEGRALBOUND},
	and the fact that $\totzeromax{\leq N}$ and $\totonemax{\leq N}$
	are increasing in both of their arguments
	to bound it by
	\begin{align}
	& \lesssim
		\varepsilon
		\frac{\ln^3(\myexp + t)}{1 + t}	
		\int_{t'=0}^t 
			\frac{1}{(1 + t')}
			\frac{1}{\upmu_{\star}(t',u)}
			\varepsilon
		\, dt'
		+
		\frac{\ln^3(\myexp + t)}{1 + t}	
		\int_{t'=0}^t 
			\frac{1}{(1 + t')}
			\frac{1}{\upmu_{\star}(t',u)}
			\int_{s=0}^{t'} 
				\frac{ \totzeromax{\leq N}^{1/2}(s,u)}{1 + s} 
			\, ds
		\, dt'
			\\
		& \ \
		+
		\frac{\ln^3(\myexp + t)}{1 + t}	
		\int_{t'=0}^t 
			\frac{1}{(1 + t')}
			\frac{1}{\upmu_{\star}(t',u)}
			\int_{s=0}^{t'} 
				\frac{\totonemax{\leq N}^{1/2}(s,u)}{(1 + s) \upmu_{\star}^{1/2}(s,u)} 
			\, ds
		\, dt'
		\notag
		\\
	& \lesssim
	\varepsilon
	\left\lbrace
		\ln \upmu_{\star}^{-1}(t,u) + 1
	\right\rbrace
	+ \int_{t'=0}^t 
			\frac{1}{(1 + t')^{3/2}}
			\frac{1}{\upmu_{\star}(t',u)}
			\totzeromax{\leq N}^{1/2}(t',u)
		\, dt'
			\notag \\
	& \ \ 
		+ \int_{t'=0}^t 
				\frac{1}{(1 + t')^{3/2}}
				\frac{1}{\upmu_{\star}(t',u)}
				\totonemax{\leq N}^{1/2}(t',u)
			\, dt'
		\notag
	\end{align}
	as desired. 
	We have thus proved the desired estimate \eqref{E:SLIGHTLYLESSSHARPL2FORANLGAPUPMUANDROTTRCHI}
	when $\mathscr{Z}^{N-1} = \mathscr{S}^{N-1}.$ 

	We now assume that $\mathscr{Z}^{N-1}$ contains a factor of $\rgeo \Lunit.$
	We prove the desired estimate \eqref{E:SLIGHTLYLESSSHARPL2FORANLGAPUPMUANDROTTRCHI}
	for $\| \upmu \angLap \mathscr{Z}^{N-1} \upmu \|_{L^2(\Sigma_t^u)}$ in detail.
	The estimate \eqref{E:SLIGHTLYLESSSHARPL2FORANLGAPUPMUANDROTTRCHI}
	for $\|\upmu \Rot \mathscr{Z}^{N-1} \mytr \upchi^{(Small)} \|_{L^2(\Sigma_t^u)}$
	can be proved in a similar fashion, and we omit those details.
	To proceed, we first use the 
	commutator estimates that we used to prove \eqref{E:COMMUTINGDONERADFIRSTWITHANRGEOLFACTORWANTSTOBEHARMLESS},
	inequality \eqref{E:LDERIVATIVECRUCICALTRANSPORTINTEQUALITIES},
	and the estimate $|\rgeo \mytr \upchi| \lesssim 1$
	(that is, \eqref{E:CRUDELOWERORDERC0BOUNDDERIVATIVESOFANGULARDEFORMATIONTENSORS})
	in order to deduce that
	\begin{align} \label{E:POINTWISEESTIMATESTOPORDERANGLAPUPMUTERMCONTAINSAGOODFACTOR}
		\left|
			\angLap \mathscr{Z}^{N-1} \upmu
		\right|
		& 
		\lesssim
		\frac{1}{1 + t}
			\left| 
				\fourmyarray[\rgeo \left\lbrace \Lunit + \frac{1}{2} \mytr \upchi \right\rbrace \mathscr{Z}^{\leq N-1} \Psi]
					{\Rad \mathscr{Z}^{\leq N-1} \Psi}
					{\rgeo \angdiff \mathscr{Z}^{\leq N-1} \Psi}
					{\mathscr{Z}^{\leq N-1} \Psi}
			\right|
			+ 
			\frac{\ln(\myexp + t)}{(1 + t)^3}
			\left|
				\myarray[\mathscr{Z}^{\leq N} (\upmu - 1)]
					{\sum_{a=1}^3 \rgeo |\mathscr{Z}^{\leq N} \Lunit_{(Small)}^a|} 
			\right|.
	\end{align}	
	We now multiply both sides of \eqref{E:POINTWISEESTIMATESTOPORDERANGLAPUPMUTERMCONTAINSAGOODFACTOR} by
	$\upmu,$ take the $L^2(\Sigma_t^u)$ norm of each side, 
	and use the estimates $\upmu \lesssim \ln(\myexp + t)$ (that is, \eqref{E:C0BOUNDCRUCIALEIKONALFUNCTIONQUANTITIES})
	and
	\eqref{E:EIKONALFUNCTIONQUANTITIESL2BOUNDSINTERMSOFQ0ANDQ1}
	and Prop.~\ref{P:L2NORMSOFPSIINTERMSOFTHECOERCIVEQUANTITIES} 
	to deduce that 
	\begin{align} \label{E:ALMOSTDONEL2ESTIMATETOPORDERANGLAPUPMUNOTPURESPATIAL}
		\| \upmu \angLap \mathscr{Z}^{N-1} \upmu \|_{L^2(\Sigma_t^u)}
		& \lesssim \totzeromax{\leq N}^{1/2}(t,u) 
			+ 
		\frac{\ln^2(\myexp + t)}{(1 + t)^2} 
		\int_{s=0}^t 
			\left\lbrace
				\frac{ \totzeromax{\leq N}^{1/2}(s,u)}{1 + s} 
				+ \frac{\totonemax{\leq N}^{1/2}(s,u)}{(1 + s)\upmu_{\star}^{1/2}(s,u)} 
			\right\rbrace
		\, ds
			\\
	& \ \ + \varepsilon
			\frac{\ln^2(\myexp + t)}{(1 + t)^2}.
			\notag
	\end{align}
	The desired estimate \eqref{E:SLIGHTLYLESSSHARPL2FORANLGAPUPMUANDROTTRCHI}
	for $\| \upmu \angLap \mathscr{Z}^{N-1} \upmu \|_{L^2(\Sigma_t^u)}$
	now follows easily from \eqref{E:ALMOSTDONEL2ESTIMATETOPORDERANGLAPUPMUNOTPURESPATIAL}
	and the fact that $\totzeromax{\leq N}$ is increasing in its arguments.

\end{proof}

\begin{lemma}[\textbf{Sharp $L^2$ estimates for a partially modified version of} $\angdiff \mathscr{S}^{N-1} \upmu$]
\label{L:ANGDIFFUPMURENORMALIZEDSHARPL2INTERMSOFQANDWIDETILDEQ}
Let $1 \leq N \leq 24$ and let $\mathscr{S}^{N-1}$ be an $(N-1)^{st}$ order pure spatial commutation vectorfield operator
(see Def.~\ref{D:DEFSETOFSPATIALCOMMUTATORVECTORFIELDS}).
Let $\mupartialmodarg{\mathscr{S}^{N-1}}$ be the partially modified 
$S_{t,u}$ one-form defined in \eqref{E:TRANSPORTPARTIALRENORMALIZEDUPMU}.
Recall the splitting $\Sigma_t^u = \Sigmaplus{t}{t}{u} \cup \Sigmaminus{t}{t}{u}.$
from Def.~\ref{D:REGIONSOFDISTINCTUPMUBEHAVIOR}.
Under the small-data and bootstrap assumptions 
of Sects.~\ref{S:PSISOLVES}-\ref{S:C0BOUNDBOOTSTRAP},
if $\varepsilon$ is sufficiently small,  
then the following estimates hold for $(t,u) \in [0,\Tboot) \times [0,U_0]:$
\begin{subequations}
\begin{align} \label{E:SIGMAMINUSANGDIFFUPMURENORMALIZEDSHARPL2INTERMSOFQANDWIDETILDEQ}
		& \left\| 
				\frac{1}{\sqrt{\upmu}}
				\rgeo
				(\Rad \Psi)  
				\mupartialmodarg{\mathscr{S}^{N-1}}
			\right\|_{L^2(\Sigmaminus{t}{t}{u})}
			\\
		& \leq 
			\sqrt{2}(1 + C \varepsilon)
			\| \rgeo \Lunit \upmu \|_{C^0(\Sigmaminus{t}{t}{u})}
			\frac{1}{\upmu_{\star}^{1/2}(t,u)}
			\int_{t'=0}^t
				\frac{1}{\rgeo(t',u) \upmu_{\star}^{1/2}(t',u)} \totonemax{\leq N}^{1/2}(t',u)
			\, dt'
				\notag \\
		& \ \ + C \varepsilon 
				\frac{1}{\upmu_{\star}^{1/2}(t,u)}
				\int_{t'=0}^t
					\frac{1}{(1 + t')^{3/2}} \totzeromax{\leq N}^{1/2}(t',u)
				\, dt'
			\notag \\
		& \ \ + C \varepsilon 
					\frac{1}{\upmu_{\star}^{1/2}(t,u)}
				\int_{t'=0}^t
					\frac{1}{(1 + t')^{3/2} \upmu_{\star}^{1/2}(t',u)} \totonemax{\leq N}^{1/2}(t',u)
				\, dt'
			\notag \\	
		& \ \ + C \varepsilon \ln(\myexp + t) \frac{1}{\upmu_{\star}^{1/2}(t,u)} \totzeromax{\leq N-1}^{1/2}(t,u)
				+ C \varepsilon \ln(\myexp + t) \frac{1}{\upmu_{\star}^{1/2}(t,u)} \totonemax{\leq N-1}^{1/2}(t,u)
				+ C \varepsilon^2 \frac{1}{\upmu_{\star}^{1/2}(t,u)},
				\notag
	\end{align}
	
	\begin{align} \label{E:SIGMAPLUSANGDIFFUPMURENORMALIZEDSHARPL2INTERMSOFQANDWIDETILDEQ}
		& \left\| 
				\frac{1}{\sqrt{\upmu}}
				\rgeo
				(\Rad \Psi)  
				\mupartialmodarg{\mathscr{S}^{N-1}}
			\right\|_{L^2(\Sigmaplus{t}{t}{u})}
			\\
		& \leq 
			\sqrt{2}(1 + C \varepsilon)
			\left\| \frac{\rgeo \Lunit \upmu}{\upmu} \right\|_{C^0(\Sigmaplus{t}{t}{u})}
			\int_{t'=0}^t
				\left\|
					\sqrt{
					\frac{\upmu(t,\cdot)} 
					     {\upmu}
					     }
				\right\|_{C^0(\Sigmaplus{t'}{t}{u})}
				\frac{1}{\rgeo(t',u)} \totonemax{\leq N}^{1/2}(t',u)
			\, dt'
				\notag \\
		& \ \ + C \varepsilon \frac{1}{\upmu_{\star}^{1/2}(t,u)}
				\int_{t'=0}^t
					\frac{1}{(1 + t')^{3/2}} \totzeromax{\leq N}^{1/2}(t',u)
				\, dt'
			\notag \\
		& \ \ + C \varepsilon \frac{1}{\upmu_{\star}^{1/2}(t,u)}
				\int_{t'=0}^t
					\frac{1}{(1 + t')^{3/2} \upmu_{\star}^{1/2}(t',u)} \totonemax{\leq N}^{1/2}(t',u)
				\, dt'
			\notag \\	
		& \ \ + C \varepsilon \ln(\myexp + t) \frac{1}{\upmu_{\star}^{1/2}(t,u)} \totzeromax{\leq N-1}^{1/2}(t,u)
				+ C \varepsilon \ln(\myexp + t) \frac{1}{\upmu_{\star}^{1/2}(t,u)} \totonemax{\leq N-1}^{1/2}(t,u)
				+ C \varepsilon^2 \frac{1}{\upmu_{\star}^{1/2}(t,u)}.
				\notag
	\end{align}
	\end{subequations}
	Furthermore, the following less sharp estimate also holds:
	\begin{align} \label{E:LESSSHARPSIGMAANGDIFFUPMURENORMALIZEDL2INTERMSOFQANDWIDETILDEQ}
		\left\| 
			  \rgeo
				\mupartialmodarg{\mathscr{S}^{N-1}}
		\right\|_{L^2(\Sigma_t^u)}
		& \le	  C
						(1 + t)
						\ln(\myexp + t)
						\totzeromax{\leq N}^{1/2}(t,u)
					+ C
						(1 + t)
						\ln(\myexp + t)
						\totonemax{\leq N}^{1/2}(t,u)
					\\
		& \  \ + C \varepsilon (1 + t).
					\notag 
	\end{align}
\end{lemma}

\begin{proof}
To prove \eqref{E:SIGMAMINUSANGDIFFUPMURENORMALIZEDSHARPL2INTERMSOFQANDWIDETILDEQ}, 
we first multiply inequality \eqref{E:ANGDIFFUPMUSHARPPOINTWISE} by 
$\frac{1}{\sqrt{\upmu}}(\Rad \Psi)$ and take the norm $\| \cdot \|_{L^2(\Sigmaminus{t}{t}{u})}$
of both sides.
It is only for the quantity arising from the first term on the right-hand side of \eqref{E:ANGDIFFUPMUSHARPPOINTWISE} 
that we use the subset norm $\| \cdot \|_{L^2(\Sigmaminus{t}{t}{u})}.$ We bound all other terms in the
(larger) norm $\| \cdot \|_{L^2(\Sigma_t^u)}.$ Next, we use Lemma~\ref{L:UPMUFIRSTTRANSPORT}, 
inequalities \eqref{E:LOWERORDERC0BOUNDLIEDERIVATIVESOFGRAME} and \eqref{E:C0BOUNDCRUCIALEIKONALFUNCTIONQUANTITIES},
and the bootstrap assumptions \eqref{E:PSIFUNDAMENTALC0BOUNDBOOTSTRAP} to deduce that
\begin{align} \label{E:GLLRADPSILUPMUINEQUALITY}
	\frac{1}{2} G_{\Lunit \Lunit} \Rad \Psi 
	& = \Lunit \upmu + \mathcal{O}\left(\varepsilon \frac{\ln(\myexp + t)}{(1 + t)^2} \right).
\end{align}
Using \eqref{E:GLLRADPSILUPMUINEQUALITY},
Lemma~\ref{L:L2NORMSOFTIMEINTEGRATEDFUNCTIONS},
and inequality \eqref{E:SQRTUPMUANGDIFFPSIL2INTERMSOFONE},
we derive the desired bound for the quantity arising from the first term on the right-hand side of 
\eqref{E:ANGDIFFUPMUSHARPPOINTWISE} as follows:
\begin{align} \label{E:KEYPROOFSTEPSIGMAMINUSANGDIFFUPMURENORMALIZEDSHARPL2INTERMSOFQANDWIDETILDEQ}
	& 
	\frac{1}{2} 
	\left\|
	\frac{1}{\sqrt{\upmu}}
	|(\Rad \Psi) G_{\Lunit \Lunit}|
		\int_{t'=0}^t
			\left|\rgeo \angdiff \mathscr{S}^N \Psi \right|(t',u,\vartheta)
		\, dt'
	\right\|_{L^2(\Sigmaminus{t}{t}{u})}
		\\
	& \leq
		\frac{1}{2}
		\left\|
			\frac{1}{\sqrt{\upmu}}
			(\Rad \Psi) G_{\Lunit \Lunit}
		\right \|_{C^0(\Sigmaminus{t}{t}{u})}
	\left\|
		\int_{t'=0}^t
			\left|\rgeo \angdiff \mathscr{S}^N \Psi \right|(t',u,\vartheta)
		\, dt'
	\right\|_{L^2(\Sigma_t^u)}
		\notag \\
	& \leq 
		(1 + C \varepsilon)
		\frac{1}{\upmu_{\star}^{1/2}(t,u)}
		\| \rgeo \Lunit \upmu \|_{C^0(\Sigmaminus{t}{t}{u})}
		\int_{t'=0}^t
			\left\| \angdiff \mathscr{S}^N \Psi \right\|_{L^2(\Sigma_t^u)}
		\, dt'
			\notag \\
		& \ \ + C \varepsilon \frac{\ln(\myexp + t)}{1 + t} \frac{1}{\upmu_{\star}^{1/2}(t,u)}
			\int_{t'=0}^t
				\left\| \angdiff \mathscr{S}^N \Psi \right\|_{L^2(\Sigma_t^u)}
			\, dt'
		\notag \\
		& \leq 
			\sqrt{2}(1 + C \varepsilon)
			\| \rgeo \Lunit \upmu \|_{C^0(\Sigmaminus{t}{t}{u})}
			\frac{1}{\upmu_{\star}^{1/2}(t,u)}
			\int_{t'=0}^t
				\frac{1}{\rgeo(t',u) \upmu_{\star}^{1/2}(t',u)} \totonemax{\leq N}^{1/2}(t',u)
			\, dt'
				\notag \\
			& \ \
				+ C \varepsilon \frac{1}{\upmu_{\star}^{1/2}(t,u)}
				\int_{t'=0}^t
					\frac{1}{(1 + t')^{3/2} \upmu_{\star}^{1/2}(t',u)} \totonemax{\leq N}^{1/2}(t',u)
				\, dt'
			\notag
\end{align}
as desired. 

Similarly, we use Lemma~\ref{L:L2NORMSOFTIMEINTEGRATEDFUNCTIONS},
Prop.~\ref{P:L2NORMSOFPSIINTERMSOFTHECOERCIVEQUANTITIES},
and the bootstrap assumption $\| \Rad \Psi \|_{C^0(\Sigma_t^u)} \leq \varepsilon (1 + t)^{-1}$
(that is, \eqref{E:PSIFUNDAMENTALC0BOUNDBOOTSTRAP}) 
to deduce the following
desired estimate for the quantity arising from the next-to-last term
on the right-hand side of \eqref{E:ANGDIFFUPMUSHARPPOINTWISE}:
\begin{align} \label{E:PSIBOUNDSPROOFSTEPSIGMAMINUSANGDIFFUPMURENORMALIZEDSHARPL2INTERMSOFQANDWIDETILDEQ}
	& \left\|
		\frac{1}{\sqrt{\upmu}}
		(\Rad \Psi)
		\int_{t'=0}^t
					\left| 
							\fourmyarray[\rgeo \left\lbrace \Lunit + \frac{1}{2} \mytr \upchi \right\rbrace \mathscr{Z}^{\leq N-1} \Psi]
								{\Rad \mathscr{Z}^{\leq N-1} \Psi}
								{\rgeo \angdiff \mathscr{Z}^{\leq N-1} \Psi}
								{\mathscr{Z}^{\leq N-1} \Psi}
						\right|(t',u,\vartheta)
		\, dt'
	\right\|_{L^2(\Sigma_t^u)}
	\\
	& \leq 
		C \varepsilon \frac{1}{\upmu_{\star}^{1/2}(t,u)}
		\int_{t'=0}^t
			\frac{1}{1 + t'} \totzeromax{\leq N-1}^{1/2}(t',u)
		\, dt'
		+ 
		C \varepsilon \frac{1}{\upmu_{\star}^{1/2}(t,u)}
		\int_{t'=0}^t
			\frac{1}{(1 + t') \upmu_{\star}^{1/2}(t',u)} \totonemax{\leq N-1}^{1/2}(t',u)
		\, dt'
		\notag
			\\
	& \leq 
		C \varepsilon \ln(\myexp + t) \frac{1}{\upmu_{\star}^{1/2}(t,u)} \totzeromax{\leq N-1}^{1/2}(t,u)
		+ C \varepsilon \ln(\myexp + t) \frac{1}{\upmu_{\star}^{1/2}(t,u)} \totonemax{\leq N-1}^{1/2}(t,u).
		\notag
\end{align}
In the last step of \eqref{E:PSIBOUNDSPROOFSTEPSIGMAMINUSANGDIFFUPMURENORMALIZEDSHARPL2INTERMSOFQANDWIDETILDEQ},
we used the fact that
$\totzeromax{\leq N-1}$ and $\totonemax{\leq N-1}$ are increasing in their arguments,
and we used the estimate \eqref{E:LOGLOSSLESSSINGULARTERMSMTHREEFOURTHSINTEGRALBOUND}
to annihilate the factor $\upmu_{\star}^{1/2}(t',u)$ in the denominator of the integrand.

Similarly, we bound the quantity arising from the last term
on the right-hand side of \eqref{E:ANGDIFFUPMUSHARPPOINTWISE}
with the help of the estimate
\eqref{E:EIKONALFUNCTIONQUANTITIESL2BOUNDSINTERMSOFQ0ANDQ1}
as follows:
\begin{align}
	& C 
	 \left\|
		\frac{1}{\sqrt{\upmu}(t,u)}
		(\Rad \Psi) 
					\int_{t'=0}^t
						\frac{\ln(\myexp + t')}{(1 + t')^2}
						\left| 
							\myarray
								[\mathscr{Z}^{\leq N} (\upmu - 1)]
								{\rgeo \sum_{a=1}^3 |\mathscr{Z}^{\leq N} \Lunit_{(Small)}^a|}
						 \right|(t',u,\vartheta)
					\, dt'	
	\right\|_{L^2(\Sigma_t^u)}
		\\
	& \lesssim 
		\varepsilon^2 
		\frac{1}{\upmu_{\star}^{1/2}(t,u)}
		\int_{t'=0}^t
			\frac{\ln(\myexp + t')}{(1 + t')^2}
		\, dt'
		+
		\varepsilon \frac{1}{\upmu_{\star}^{1/2}(t,u)}
		\int_{t'=0}^t
			\frac{\ln(\myexp + t')}{(1 + t')^2}
			\int_{s=0}^{t'}
				\frac{1}{1 + s} \totzeromax{\leq N}^{1/2}(s,u)
			\, ds
		\, dt'
			\notag \\
	& \ \	
		+
		\varepsilon \frac{1}{\upmu_{\star}^{1/2}(t,u)}
		\int_{t'=0}^t
			\frac{\ln(\myexp + t')}{(1 + t')^2}
			\int_{s=0}^{t'}
				\frac{1}{(1 + s)\upmu_{\star}^{1/2}(s,u)} \totonemax{\leq N}^{1/2}(s,u)
			\, ds
		\, dt'
			\notag \\	& \lesssim 
		\varepsilon^2 \frac{1}{\upmu_{\star}^{1/2}(t,u)}
		+
		\varepsilon \frac{1}{\upmu_{\star}^{1/2}(t,u)}
		\int_{t'=0}^t
			\frac{1}{(1 + t')^{3/2}}
			\totzeromax{\leq N}^{1/2}(t,u)
		\, dt'
		\notag \\
	& \ \ +
		\varepsilon \frac{1}{\upmu_{\star}^{1/2}(t,u)}
		\int_{t'=0}^t
			\frac{1}{(1 + t')^{3/2}}
			\totonemax{\leq N}^{1/2}(t,u)
		\, dt'.
		\notag
\end{align}

Next, using the bootstrap assumption $\| \Rad \Psi \|_{C^0(\Sigma_t^u)} \leq \varepsilon (1 + t)^{-1}$
and the estimate \eqref{E:ESTIMATEFORSPATIALL2NORMOF1},
we bound the quantity arising from the term
$C \varepsilon$ on the right-hand side of \eqref{E:ANGDIFFUPMUSHARPPOINTWISE} as follows:
\begin{align}
	\left\| 
		\frac{1}{\sqrt{\upmu}}
		(\Rad \Psi)  
			C \varepsilon
	\right\|_{L^2(\Sigma_t^u)}
	& \leq C \varepsilon^2 
				 \frac{1}{(1 + t)} 
				 \frac{1}{\upmu_{\star}^{1/2}(t,u)} 
				 \| 1 \|_{L^2(\Sigma_t^u)}
	\leq C \varepsilon^2 \frac{1}{\upmu_{\star}^{1/2}(t,u)}.
\end{align}

Finally, using the bootstrap assumption $\| \Rad \Psi \|_{C^0(\Sigma_t^u)} \leq \varepsilon (1 + t)^{-1},$
\eqref{E:LOWERORDERC0BOUNDLIEDERIVATIVESOFGRAME}, 
and Lemma~\ref{L:L2NORMSOFTIMEINTEGRATEDFUNCTIONS}, 
we deduce that the quantity arising from the term 
$C \varepsilon
					\int_{t'=0}^t
						\left|\angdiff \mathscr{S}^N \Psi \right|(t',u,\vartheta)
					\, dt'$
on the right-hand side of \eqref{E:ANGDIFFUPMUSHARPPOINTWISE} can be bounded as follows:
\begin{align} \label{E:LASTTERMPROOFSIGMAMINUSANGDIFFUPMURENORMALIZEDSHARPL2INTERMSOFQANDWIDETILDEQ}
	& 
	\left\|
		C \varepsilon
		\frac{1}{\sqrt{\upmu}}
		(\Rad \Psi)
		\int_{t'=0}^t
			\left|\angdiff \mathscr{S}^N \Psi \right|(t',u,\vartheta)
		\, dt'
		\right\|_{L^2(\Sigma_t^u)}
			\\
& \leq C \varepsilon^2 \frac{1}{\upmu_{\star}^{1/2}(t,u)}
	\int_{t'=0}^t
		\frac{1}{(1 + t')^{3/2} \upmu_{\star}^{1/2}(t',u)} \totonemax{\leq N}^{1/2}(t',u)
	\, dt'. \notag
\end{align}
We now observe that the right-hand side of \eqref{E:LASTTERMPROOFSIGMAMINUSANGDIFFUPMURENORMALIZEDSHARPL2INTERMSOFQANDWIDETILDEQ}
is $\lesssim$
the right-hand side of \eqref{E:SIGMAMINUSANGDIFFUPMURENORMALIZEDSHARPL2INTERMSOFQANDWIDETILDEQ} as desired.
We have thus proved the estimate \eqref{E:SIGMAMINUSANGDIFFUPMURENORMALIZEDSHARPL2INTERMSOFQANDWIDETILDEQ}.

The estimate \eqref{E:SIGMAPLUSANGDIFFUPMURENORMALIZEDSHARPL2INTERMSOFQANDWIDETILDEQ} can be proved by making
minor changes to our proof of \eqref{E:SIGMAMINUSANGDIFFUPMURENORMALIZEDSHARPL2INTERMSOFQANDWIDETILDEQ}.

The proof of \eqref{E:LESSSHARPSIGMAANGDIFFUPMURENORMALIZEDL2INTERMSOFQANDWIDETILDEQ} is 
much easier. We estimate the right-hand side of \eqref{E:ANGDIFFUPMUSHARPPOINTWISE} 
in the norm $\| \cdot \|_{L^2(\Sigma_t^u)}$ by using arguments similar to the ones
we used to deduce \eqref{E:SIGMAMINUSANGDIFFUPMURENORMALIZEDSHARPL2INTERMSOFQANDWIDETILDEQ},
but we treat the first term on the right-hand side rather crudely 
with the help of the estimate $\|G_{\Lunit \Lunit} \|_{C^0(\Sigma_t^u)} \lesssim 1$ (that is, \eqref{E:LOWERORDERC0BOUNDLIEDERIVATIVESOFGRAME}),
and without using an analog of \eqref{E:GLLRADPSILUPMUINEQUALITY}.
In doing so, we encounter the integrals 
$\int_{t'=0}^t
				\frac{1}{\rgeo(t',u) \upmu_{\star}^{1/2}(t',u)} \totonemax{\leq N}^{1/2}(t',u)
			\, dt'$
and $\int_{t'=0}^t
				\frac{1}{\rgeo(t',u)} \totzeromax{\leq N}^{1/2}(t',u)
			\, dt',$			
and we respectively bound them by 
$\lesssim \ln(\myexp + t) \totonemax{\leq N}^{1/2}(t,u)$ 
and 
$\lesssim \ln(\myexp + t) \totzeromax{\leq N}^{1/2}(t,u),$
where in the first estimate, we use inequality \eqref{E:LOGLOSSLESSSINGULARTERMSMTHREEFOURTHSINTEGRALBOUND}.

\end{proof}

\begin{lemma}[\textbf{Sharp} $L^2$ \textbf{estimates for a partially modified version of} $\mathscr{S}^{N-1} \mytr \upchi^{(Small)}$]
\label{L:TRCHIRENORMALIZEDSHARPL2INTERMSOFQANDWIDETILDEQ}
Let $1 \leq N \leq 24$ and let $\mathscr{S}^{N-1}$ be an $(N-1)^{st}$ order pure spatial commutation vectorfield operator
(see Def.~\ref{D:DEFSETOFSPATIALCOMMUTATORVECTORFIELDS}).
Let $\chipartialmodarg{\mathscr{S}^{N-1}}$
be the partially modified quantity defined in \eqref{E:TRANSPORTPARTIALRENORMALIZEDTRCHIJUNK}.
Recall the splitting $\Sigma_t^u = \Sigmaplus{t}{t}{u} \cup \Sigmaminus{t}{t}{u}.$
from Def.~\ref{D:REGIONSOFDISTINCTUPMUBEHAVIOR}.
Under the small-data and bootstrap assumptions 
of Sects.~\ref{S:PSISOLVES}-\ref{S:C0BOUNDBOOTSTRAP},
if $\varepsilon$ is sufficiently small,  
then the following estimates hold for $(t,u) \in [0,\Tboot) \times [0,U_0]:$
\begin{subequations}
\begin{align} \label{E:SIGMAMINUSTRCHIRENORMALIZEDSHARPL2INTERMSOFQANDWIDETILDEQ}
		& \left\| 
				\frac{1}{\sqrt{\upmu}}
				(\Rad \Psi)  
				\rgeo^2
				\chipartialmodarg{\mathscr{S}^{N-1}}
			\right\|_{L^2(\Sigmaminus{t}{t}{u})}
			\\
		& \leq 
			(\sqrt{12} + C \varepsilon)
			\| \rgeo \Lunit \upmu \|_{C^0(\Sigmaminus{t}{t}{u})}
			\frac{1}{\upmu_{\star}^{1/2}(t,u)}
			\int_{t'=0}^t
				\frac{1}{\rgeo(t',u) \upmu_{\star}^{1/2}(t',u)} \totonemax{\leq N}^{1/2}(t',u)
			\, dt'
				\notag \\
		& \ \
		 	+ C \varepsilon^2
			\frac{1}{\upmu_{\star}^{1/2}(t,u)}
			\int_{t'=0}^t
				\frac{1}{\rgeo(t',u) \upmu_{\star}^{1/2}(t',u)} \totonemax{\leq N}^{1/2}(t',u)
			\, dt'
				\notag \\
		& \ \ + C \varepsilon 
				\frac{1}{\upmu_{\star}^{1/2}(t,u)}
				\int_{t'=0}^t
					\frac{1}{(1 + t')^{3/2}} \totzeromax{\leq N}^{1/2}(t',u)
				\, dt'
			\notag \\
		& \ \ + C \varepsilon 
				\frac{1}{\upmu_{\star}^{1/2}(t,u)}
				\int_{t'=0}^t
					\frac{1}{(1 + t')^{3/2} \upmu_{\star}^{1/2}(t',u)} \totonemax{\leq N}^{1/2}(t',u)
				\, dt'
			\notag \\	
		& \ \ + C \varepsilon \ln(\myexp + t) \frac{1}{\upmu_{\star}^{1/2}(t,u)} \totzeromax{\leq N-1}^{1/2}(t,u)
				+ C \varepsilon \ln(\myexp + t) \frac{1}{\upmu_{\star}^{1/2}(t,u)} \totonemax{\leq N-1}^{1/2}(t,u)
				+ C \varepsilon^2 \frac{1}{\upmu_{\star}^{1/2}(t,u)},
				\notag
	\end{align}
	\begin{align} \label{E:SIGMAPLUSTRCHIRENORMALIZEDSHARPL2INTERMSOFQANDWIDETILDEQ}
		& \left \|
		 	\frac{1}{\sqrt{\upmu}} 
		 	(\Rad \Psi)
		 	\rgeo^2
			\chipartialmodarg{\mathscr{S}^{N-1}}
		\right\|_{L^2(\Sigmaplus{t}{t}{u})}
			\\
		& \leq 
			(\sqrt{12} + C \varepsilon)
			\left\| \frac{\rgeo \Lunit \upmu}{\upmu} \right\|_{C^0(\Sigmaplus{t}{t}{u})}
			\int_{t'=0}^t
				\left\|
					\sqrt{
					\frac{\upmu(t,\cdot)} 
					     {\upmu}
					     }
				\right\|_{C^0(\Sigmaplus{t'}{t}{u}}
				\frac{1}{\rgeo(t',u)} \totonemax{\leq N}^{1/2}(t',u)
			\, dt'
				\notag \\
		& \ \
			+ C \varepsilon^2
		  \frac{1}{\upmu_{\star}(t,u)}
			\int_{t'=0}^t
				\left\|
					\sqrt{
					\frac{\upmu(t,\cdot)} 
					     {\upmu}
					     }
				\right\|_{C^0(\Sigmaplus{t'}{t}{u})}
				\frac{1}{\rgeo(t',u)} \totonemax{\leq N}^{1/2}(t',u)
			\, dt'
				\notag \\
		& \ \ + C \varepsilon 
				\frac{1}{\upmu_{\star}^{1/2}(t,u)}
				\int_{t'=0}^t
					\frac{1}{(1 + t')^{3/2}} \totzeromax{\leq N}^{1/2}(t',u)
				\, dt'
			\notag \\
		& \ \ + C \varepsilon
				\frac{1}{\upmu_{\star}^{1/2}(t,u)}
				\int_{t'=0}^t
					\frac{1}{(1 + t')^{3/2} \upmu_{\star}^{1/2}(t',u)} \totonemax{\leq N}^{1/2}(t',u)
				\, dt'
			\notag \\	
		& \ \ + C \varepsilon \ln(\myexp + t) \frac{1}{\upmu_{\star}^{1/2}(t,u)}\totzeromax{\leq N-1}^{1/2}(t,u)
				+ C \varepsilon \ln(\myexp + t) \frac{1}{\upmu_{\star}^{1/2}(t,u)} \totonemax{\leq N-1}^{1/2}(t,u)
				+ C \varepsilon^2 \frac{1}{\upmu_{\star}^{1/2}(t,u)}.
				\notag
	\end{align}
	\end{subequations}
	Furthermore, the following less sharp estimate also holds:
	\begin{align} \label{E:LESSSHARPSIGMATRCHIRENORMALIZEDL2INTERMSOFQANDWIDETILDEQ}
	 \left\| 
			\rgeo^2 
			\chipartialmodarg{\mathscr{S}^{N-1}}
		\right\|_{L^2(\Sigma_t^u)}
		& \leq  C
						\ln(\myexp + t)
						(1 + t)
						\totzeromax{\leq N}^{1/2}(t,u)
					+ C
						\ln(\myexp + t)
						(1 + t)
						\totonemax{\leq N}^{1/2}(t,u)
					\\
		& \  \ + C 
						\varepsilon (1 + t).
					\notag 
	\end{align}
	
\end{lemma}

\begin{proof}
	We apply essentially the same reasoning used in the proof of
	Lemma~\ref{L:ANGDIFFUPMURENORMALIZEDSHARPL2INTERMSOFQANDWIDETILDEQ}
	to the pointwise estimate \eqref{E:TRCHIRENORMALIZEDSHARPPOINTWISE}.
	The factors of $\sqrt{12}$ arise because in order to estimate the first
	term on the right-hand side of \eqref{E:TRCHIRENORMALIZEDSHARPPOINTWISE},
	we use inequalities \eqref{E:ANGLAPFUNCTIONPOINTWISEINTERMSOFROTATIONS}
	and \eqref{E:SQRTUPMUANGDIFFPSIL2INTERMSOFONE}
	and Cor.~\ref{C:SQRTEPSILONREPLCEDWITHCEPSILON}.
\end{proof}

\begin{lemma}[\textbf{Sharp} $L^2$ \textbf{estimates for} $\angLie_{\Lunit} + \mytr \upchi$ 
\textbf{applied to a partially modified version of}  $\angdiff \mathscr{S}^{N-1} \upmu$]
\label{L:LDERIVANGDIFFUPMURENORMALIZEDSHARPL2INTERMSOFQANDWIDETILDEQ}
Let $1 \leq N \leq 24$ and let $\mathscr{S}^{N-1}$ be an $(N-1)^{st}$ order pure spatial commutation vectorfield operator
(see Def.~\ref{D:DEFSETOFSPATIALCOMMUTATORVECTORFIELDS}).
Let $\mupartialmodarg{\mathscr{S}^{N-1}}^{\#}$ be 
$S_{t,u}-$tangent vectorfield that is $\gsphere$ dual to
the partially modified 
$S_{t,u}$ one-form defined in \eqref{E:TRANSPORTPARTIALRENORMALIZEDUPMU}.
Recall the splitting $\Sigma_t^u = \Sigmaplus{t}{t}{u} \cup \Sigmaminus{t}{t}{u}.$
from Def.~\ref{D:REGIONSOFDISTINCTUPMUBEHAVIOR}.
Under the small-data and bootstrap assumptions 
of Sects.~\ref{S:PSISOLVES}-\ref{S:C0BOUNDBOOTSTRAP},
if $\varepsilon$ is sufficiently small,  
then the following estimate holds for $(t,u) \in [0,\Tboot) \times [0,U_0]:$
\begin{align} \label{E:LPLUSTRCHIAPPLIEDTORENORMALIZEDANGDIFFUPMUSHARPL2INTERMSOFQANDWIDETILDEQ}
	& \left\|
			\frac{1}{\sqrt{\upmu}}
			(\Rad \Psi)
			\rgeo
		\left\lbrace
			\angLie_{\Lunit}
			+ \mytr \upchi 
		\right\rbrace
		\mupartialmodarg{\mathscr{S}^{N-1}}^{\#}
		\right\|_{L^2(\Sigma_t^u)}
	  \\
	  & \leq 
	  	\sqrt{2}
	  	\frac
				{\left\| [\Lunit \upmu]_- \right\|_{C^0(\Sigma_t^u)}}
				{\upmu_{\star}(t,u)}
	  	\totonemax{\leq N}^{1/2}(t,u)
	  + 
	  	\sqrt{2}
	  	\left\| \frac{[\Lunit \upmu]_+}{\upmu} \right\|_{C^0(\Sigma_t^u)} 
	  	\totonemax{\leq N}^{1/2}(t,u)
	  		\notag \\
	  & \ \ + C \varepsilon 
	  					\frac{1}{(1 + t)^{3/2}} 
	  					\frac{1}{\upmu_{\star}^{1/2}(t,u)}
							\totzeromax{\leq N}^{1/2}(t,u)
	  			+ C \varepsilon
	  					\frac{1}{(1 + t)^{3/2}} 
	  					\frac{1}{\upmu_{\star}(t,u)}
	  					\totonemax{\leq N}^{1/2}(t,u)
	  					\notag \\
	  & \ \ + C \varepsilon 
	  				\frac{1}{(1 + t)} 
						\frac{1}{\upmu_{\star}^{1/2}(t,u)}
	  				\totzeromax{\leq N-1}^{1/2}(t,u)
	  			+ C \varepsilon 
	  				\frac{1}{(1 + t)} 
						\frac{1}{\upmu_{\star}(t,u)}
	  				\totonemax{\leq N-1}^{1/2}(t,u).
	  				\notag
\end{align}

Furthermore, the following less sharp estimate also holds:
\begin{align} \label{E:LESSSHARPLPLUSTRCHIAPPLIEDTORENORMALIZEDANGDIFFUPMUL2INTERMSOFQANDWIDETILDEQ}
		\left\|
			\rgeo
			\left\lbrace
				\angLie_{\Lunit}
				+ \mytr \upchi 
			\right\rbrace
			\mupartialmodarg{\mathscr{S}^{N-1}}^{\#}
		\right\|_{L^2(\Sigma_t^u)}
	  & \leq
	  	C
			\totzeromax{\leq N}^{1/2}(t,u)
	  	+	C
	 			\frac{1}{\upmu_{\star}^{1/2}(t,u)}
	  	 	\totonemax{\leq N}^{1/2}(t,u).
	 \end{align}

\end{lemma}

\begin{proof}
To prove \eqref{E:LPLUSTRCHIAPPLIEDTORENORMALIZEDANGDIFFUPMUSHARPL2INTERMSOFQANDWIDETILDEQ}, we start by multiplying
the second inequality in \eqref{E:UPMUANDTRCHIJUNKREFINEDPARTIALLYRENORMALIZEDTRANSPORTINEQUALITIES} by 
$\frac{1}{\sqrt{\upmu}} (\Rad \Psi)$ and taking the $L^2(\Sigma_t^u)$ norm.
To bound the norm $\| \cdot \|_{L^2(\Sigma_t^u)}$
of the quantity arising from the term 
$\frac{1}{2}\rgeo G_{\Lunit \Lunit} \angdiffuparg{\#} \mathscr{S}^{N-1} \Rad \Psi,$
we use the estimate \eqref{E:GLLRADPSILUPMUINEQUALITY}, 
the inequality 
\[
	\left\| \frac{\Lunit \upmu}{\upmu} \right\|_{C^0(\Sigma_t^u)}
	\leq 
	\frac{1}{\upmu_{\star}(t,u)}
	\left\| [\Lunit \upmu]_- \right\|_{C^0(\Sigma_t^u)}
	+ \left\| \frac{[\Lunit \upmu]_+}{\upmu} \right\|_{C^0(\Sigma_t^u)},
\]
and the inequality 
$\left\| \rgeo \sqrt{\upmu} \angdiffuparg{\#} \mathscr{S}^{N-1} \Rad \right\|_{L^2(\Sigma_t^u)} \leq \sqrt{2} \totonemax{\leq N}^{1/2}(t,u)$
(that is, \eqref{E:SQRTUPMUANGDIFFPSIL2INTERMSOFONE})
to deduce that
\begin{align} \label{E:IMPORTANTTERMSESTIMATEDLPLUSTRCHIAPPLIEDTORENORMALIZEDANGDIFFUPMUSHARPL2INTERMSOFQANDWIDETILDEQ}
	& \left\| 
			\frac{1}{\sqrt{\upmu}}
			(\Rad \Psi) 
			\rgeo 
			\frac{1}{2}
			G_{\Lunit \Lunit} 
			\angdiffuparg{\#} 
			\mathscr{S}^{N-1} \Rad \Psi 
		\right\|_{L^2(\Sigma_t^u)}
			\\
	& \leq 
	  	\sqrt{2}
	  	\frac{1}{\upmu_{\star}(t,u)}
	  	\left\| [\Lunit \upmu]_- \right\|_{C^0(\Sigma_t^u)} 
	  	\totonemax{\leq N}^{1/2}(t,u)
	  +  \sqrt{2}
	  	\left\| \frac{[\Lunit \upmu]_+}{\upmu} \right\|_{C^0(\Sigma_t^u)} 
	  	\totonemax{\leq N}^{1/2}(t,u)
	  		\notag \\
	 & \ \ + C \varepsilon \frac{\ln(\myexp + t)}{(1 + t)^2} \frac{1}{\upmu_{\star}(t,u)} \totonemax{\leq N}^{1/2}(t,u).
	 	\notag 
\end{align}
Clearly, the right-hand side of \eqref{E:IMPORTANTTERMSESTIMATEDLPLUSTRCHIAPPLIEDTORENORMALIZEDANGDIFFUPMUSHARPL2INTERMSOFQANDWIDETILDEQ}
is bounded by the right-hand side of \eqref{E:LPLUSTRCHIAPPLIEDTORENORMALIZEDANGDIFFUPMUSHARPL2INTERMSOFQANDWIDETILDEQ}
as desired.

The quantities arising from the remaining terms on the right-hand side of \eqref{E:UPMUANDTRCHIJUNKREFINEDPARTIALLYRENORMALIZEDTRANSPORTINEQUALITIES} 
can be bounded in the norm $\| \cdot \|_{L^2(\Sigma_t^u)}$
by the right-hand side of \eqref{E:LPLUSTRCHIAPPLIEDTORENORMALIZEDANGDIFFUPMUSHARPL2INTERMSOFQANDWIDETILDEQ}
with the help of 
Prop.~\ref{P:L2NORMSOFPSIINTERMSOFTHECOERCIVEQUANTITIES},
the bootstrap assumption $\| \Rad \Psi \|_{C^0(\Sigma_t^u)} \leq \varepsilon (1 + t)^{-1},$ 
the estimates 
\eqref{E:EIKONALFUNCTIONQUANTITIESL2BOUNDSINTERMSOFQ0ANDQ1}
and \eqref{E:C0BOUNDCRUCIALEIKONALFUNCTIONQUANTITIES}, 
and inequality \eqref{E:LOGLOSSLESSSINGULARTERMSMTHREEFOURTHSINTEGRALBOUND}.

The proof of \eqref{E:LESSSHARPLPLUSTRCHIAPPLIEDTORENORMALIZEDANGDIFFUPMUL2INTERMSOFQANDWIDETILDEQ} is simpler.
We bound the term $\frac{1}{2}\rgeo G_{\Lunit \Lunit} \angdiffuparg{\#} \mathscr{S}^{N-1} \Rad \Psi$ 
from inequality \eqref{E:UPMUANDTRCHIJUNKREFINEDPARTIALLYRENORMALIZEDTRANSPORTINEQUALITIES} 
by using Prop.~\ref{P:L2NORMSOFPSIINTERMSOFTHECOERCIVEQUANTITIES}
and the estimate $\|G_{\Lunit \Lunit} \|_{C^0(\Sigma_t^u)} \lesssim 1$ (that is, \eqref{E:LOWERORDERC0BOUNDLIEDERIVATIVESOFGRAME}) as follows:
\begin{align}
	\left\| 
			\frac{1}{2}
			\rgeo 
			G_{\Lunit \Lunit} 
			\angdiffuparg{\#} 
			\mathscr{S}^{N-1} \Rad \Psi 
		\right\|_{L^2(\Sigma_t^u)}
		& \lesssim \frac{1}{\upmu_{\star}^{1/2}(t,u)} \totonemax{\leq N}^{1/2}(t,u).
\end{align}
The remaining terms from inequality \eqref{E:UPMUANDTRCHIJUNKREFINEDPARTIALLYRENORMALIZEDTRANSPORTINEQUALITIES}
can be bounded by using arguments similar to the ones 
we used in our proof of \eqref{E:LPLUSTRCHIAPPLIEDTORENORMALIZEDANGDIFFUPMUSHARPL2INTERMSOFQANDWIDETILDEQ}.

\end{proof}

\begin{lemma}
[\textbf{Sharp} $L^2$ \textbf{estimates for} $\Lunit$ 
\textbf{applied to a partially modified version of} $\mathscr{S}^{N-1} \mytr \upchi^{(Small)}$] 
\label{L:LDERIVTRCHILRENORMALIZEDSHARPL2INTERMSOFQANDWIDETILDEQ}
Let $1 \leq N \leq 24$ and let $\mathscr{S}^{N-1}$ be an $(N-1)^{st}$ order pure spatial commutation vectorfield operator
(see Def.~\ref{D:DEFSETOFSPATIALCOMMUTATORVECTORFIELDS}).
Let $\chipartialmodarg{\mathscr{S}^{N-1}}$
be the partially modified quantity defined in \eqref{E:TRANSPORTPARTIALRENORMALIZEDTRCHIJUNK}.
Recall the splitting $\Sigma_t^u = \Sigmaplus{t}{t}{u} \cup \Sigmaminus{t}{t}{u}.$
from Def.~\ref{D:REGIONSOFDISTINCTUPMUBEHAVIOR}.
Under the small-data and bootstrap assumptions 
of Sects.~\ref{S:PSISOLVES}-\ref{S:C0BOUNDBOOTSTRAP},
if $\varepsilon$ is sufficiently small, 
then the following estimate holds for $(t,u) \in [0,\Tboot) \times [0,U_0]:$
\begin{align}  \label{E:TRCHILRENORMALIZEDSHARPL2INTERMSOFQANDWIDETILDEQ}
	& \left\|
			\frac{1}{\sqrt{\upmu}}
			(\Rad \Psi)
			\Lunit
			\left\lbrace
		 		\rgeo^2 
				\chipartialmodarg{\mathscr{S}^{N-1}}
			\right\rbrace
		\right\|_{L^2(\Sigma_t^u)}
	    \\
	  & \leq 
	  	\sqrt{12}(1 + C \varepsilon)
	  	\frac{\left\| [\Lunit \upmu]_- \right\|_{C^0(\Sigma_t^u)}}
			{\upmu(t,u)}
	  	\totonemax{\leq N}^{1/2}(t,u)
	  	\notag \\
	 & \ \ + 
	  	\sqrt{12}(1 + C \varepsilon)
	  	\left\| \frac{[\Lunit \upmu]_+}{\upmu} \right\|_{C^0(\Sigma_t^u)} 
	  	\totonemax{\leq N}^{1/2}(t,u)
	  		\notag \\
	  & \ \ + C \varepsilon 
	  					\frac{1}{(1 + t)^{3/2}} 
	  					\frac{1}{\upmu_{\star}^{1/2}(t,u)}
	  					\totzeromax{\leq N}^{1/2}(t,u)
	  			+ C \varepsilon
	  					\frac{1}{(1 + t)^{3/2}} 
	  					\frac{1}{\upmu_{\star}(t,u)}
	  					\totonemax{\leq N}^{1/2}(t,u)
	  					\notag \\
	  & \ \ + C \varepsilon 
	  				\frac{1}{(1 + t)} \frac{1}{\upmu_{\star}^{1/2}(t,u)} 
	  				\totzeromax{\leq N-1}^{1/2}(t,u)
	  			+ C \varepsilon 
	  				\frac{1}{(1 + t)} \frac{1}{\upmu_{\star}(t,u)} 
	  				\totonemax{\leq N-1}^{1/2}(t,u).
	  				\notag
\end{align}

Furthermore, the following less sharp estimates also hold:
\begin{align} \label{E:LESSSHARPLAPPLIEDTORENORMALIZEDTRCHIL2INTERMSOFQANDWIDETILDEQ}
	& \left\|
			\Lunit	
			\left\lbrace
		 		\rgeo^2 
				\chipartialmodarg{\mathscr{S}^{N-1}}
			\right\rbrace
		\right\|_{L^2(\Sigma_t^u)},
	   \\
	& \left\|
			\rgeo^2
			\left\lbrace
				\Lunit 
				+ \mytr \upchi
			\right\rbrace	
			\chipartialmodarg{\mathscr{S}^{N-1}}
		\right\|_{L^2(\Sigma_t^u)}
	  \notag \\
	  & \leq
			C
	  	\totzeromax{\leq N}^{1/2}(t,u)
				+
				C
	 			\frac{1}{\upmu_{\star}^{1/2}(t,u)}
	  	 	\totonemax{\leq N}^{1/2}(t,u).
	  	 	\notag 
\end{align}

\end{lemma}

\begin{proof}
We apply the same reasoning used in the proof of Lemma~\ref{L:LDERIVANGDIFFUPMURENORMALIZEDSHARPL2INTERMSOFQANDWIDETILDEQ}
to the third and fourth pointwise inequalities in \eqref{E:UPMUANDTRCHIJUNKREFINEDPARTIALLYRENORMALIZEDTRANSPORTINEQUALITIES}.
The factors of $\sqrt{12}$ arise because in order to bound the norm $\| \cdot \|_{L^2(\Sigma_t^u)}$ of the term 
$- \frac{1}{2} \rgeo^2 G_{\Lunit \Lunit} \angLap \mathscr{Z}^{N-1} \Psi$
from the pointwise inequalities, 
we use inequalities \eqref{E:ANGLAPFUNCTIONPOINTWISEINTERMSOFROTATIONS}
and \eqref{E:SQRTUPMUANGDIFFPSIL2INTERMSOFONE}
and Cor.~\ref{C:SQRTEPSILONREPLCEDWITHCEPSILON}.

\end{proof}

\begin{lemma}[\textbf{Estimates for some easy error integrals corresponding to the partially modified version of} 
$\Rot \mathscr{S}^{N-1} \mytr \upchi^{(Small)}$]
	\label{L:NONDAMAGINGTERMSRENORMALIZEDTRCHI}
	Let $1 \leq N \leq 24$ be an integer, and let
	$\mathscr{S}^{N-1}$ be an $(N-1)^{st}$ order pure spatial commutation vectorfield operator
	(see Def.~\ref{D:DEFSETOFSPATIALCOMMUTATORVECTORFIELDS}).
	Let $\chipartialmodarg{\mathscr{S}^{N-1}}$
	be the partially modified quantity defined in \eqref{E:TRANSPORTPARTIALRENORMALIZEDTRCHIJUNK}.
	Under the small-data and bootstrap assumptions 
	of Sects.~\ref{S:PSISOLVES}-\ref{S:C0BOUNDBOOTSTRAP},
	if $\varepsilon$ is sufficiently small, 	
	then we have the following spacetime integral estimate for
	the quantity $\IBPerror{1}[\mathscr{S}^{N-1} \Rot \Psi]$
	defined in \eqref{E:NONDAMAGINGSPACETIMERENORMALZEDTRCHI},
	with $\newsmoothfunction := \rgeo^2 \chipartialmodarg{\mathscr{S}^{N-1}},$
	for $(t,u) \in [0,\Tboot) \times [0,U_0]:$
	\begin{subequations}
	\begin{align} \label{E:NONDAMAGINGSPACETIMERENORMALZEDTRCHIESTIMATE}
			\int_{\mathcal{M}_{t,u}}
				\left|\IBPerror{1}[\mathscr{S}^{N-1} \Rot \Psi] \right|
			\, d \vol
			& \lesssim 
				\varepsilon \totzeromax{\leq N}(t,u)
				+
				\varepsilon
				\int_{t'=0}^t
					\frac{1}{(1 + t')^{3/2} \upmu_{\star}(t',u)}
					\totonemax{\leq N}(t',u)
				\, dt'
				\\
	& \ \ + \varepsilon \ln^2(\myexp + t) \totzeromax{\leq N-1}(t,u) 
				+ \varepsilon \ln^2(\myexp + t) \frac{1}{\upmu_{\star}^{1/2}(t,u)} \totonemax{\leq N-1}(t,u)
				\notag \\
	& \ \ + \varepsilon^3.
			\notag
	\end{align}
	
	Furthermore, we have the following hypersurface integral estimate for
	the quantity $\IBPerror{2}[\mathscr{S}^{N-1} \Rot \Psi]$
	defined in \eqref{E:NONDAMAGINGHYPERSURFACERENORMALZEDTRCHI}:
 \begin{align}  \label{E:NONDAMAGINGHYPERSURFACERENORMALZEDTRCHIESTIMATE}
			\int_{\Sigma_t^u}
				\left|\IBPerror{2}[\mathscr{S}^{N-1} \Rot \Psi] \right|
			\, d \tvol
			& \lesssim 
				\varepsilon \totzeromax{\leq N}(t,u) 
				+ \varepsilon \totonemax{\leq N}(t,u)
				\\
			& \ \ + \varepsilon \ln^2(\myexp + t) \frac{1}{\upmu_{\star}(t,u)} \totzeromax{\leq N-1}(t,u) 
				+ \varepsilon \ln^2(\myexp + t) \frac{1}{\upmu_{\star}(t,u)} \totonemax{\leq N-1}(t,u)
				\notag \\
			& \ \ + \varepsilon^3 \frac{1}{\upmu_{\star}(t,u)}.
				\notag
\end{align}

	Finally, we have the following hypersurface integral estimate for
	the quantity $\IBPerror{3}[\Rot \mathscr{S}^{N-1}]$
	defined in \eqref{E:NONDAMAGINGINITIALDATARENORMALZEDTRCHI}:
	\begin{align} \label{E:NONDAMAGINGINITIALDATARENORMALZEDTRCHIESTIMATE}
			\int_{\Sigma_0^u}
				\left|\IBPerror{3}[\mathscr{S}^{N-1} \Rot \Psi] \right|
			\, d \tvol
			& \lesssim \varepsilon^3.
	\end{align}
	\end{subequations}
\end{lemma}

\begin{proof}
 To prove \eqref{E:NONDAMAGINGSPACETIMERENORMALZEDTRCHIESTIMATE}, we first express 
	the spacetime integral as a time integral of
	integrals over $\Sigma_{t'}^u:$
	$\int_{\mathcal{M}_{t,u}} \cdots \, d\vol = 
	\int_{t'=0}^t \int_{\Sigma_{t'}^u} \cdots \, d \tvol \, d t'.$
	We bound all spatial integrals $\int_{\Sigma_{t'}^u} \cdots \, d \tvol $
	by using Cauchy-Schwarz in the form
	$\int_{\Sigma_{t'}^u} |v_1 v_2 v_3| \, d \tvol \, d t' 
	\lesssim \| v_1 \|_{C^0(\Sigma_{t'}^u)} \| v_2 \|_{L^2(\Sigma_{t'}^u)} \| v_3 \|_{L^2(\Sigma_{t'}^u)}.$
	The product $v_1 v_2 v_3$ represents any of the products
	from the right-hand side of \eqref{E:NONDAMAGINGSPACETIMERENORMALZEDTRCHI}
	with $\mathscr{Z}^N := \mathscr{S}^{N-1} \Rot,$
	where 
	$v_2$ is equal to either the factor 
	$\mathscr{S}^{N-1} \Rot \Psi,$ 
	$\Rot \mathscr{S}^{N-1} \Rot \Psi,$
	or $\angdiff \mathscr{S}^{N-1} \Rot \Psi,$
	$v_3$ is either $\Lunit F$ or $F,$ and $v_1$ is the product of the remaining factors.
	To bound $v_2,$ we use the estimates
	$\| \mathscr{S}^{N-1} \Rot \Psi \|_{L^2(\Sigma_{t'}^u)} \lesssim \totzeromax{\leq N-1}^{1/2}(t',u) 
	+ \upmu_{\star}^{-1/2}(t',u) \totonemax{\leq N-1}^{1/2}(t',u),$
	$\| \Rot \mathscr{S}^{N-1} \Rot \Psi \|_{L^2(\Sigma_{t'}^u)} \lesssim \upmu_{\star}^{-1/2}(t',u) \totonemax{\leq N}^{1/2}(t',u),$
	and
	$\| \angdiff \mathscr{S}^{N-1} \Rot \Psi \|_{L^2(\Sigma_{t'}^u)} \lesssim (1 + t')^{-1} \upmu_{\star}^{-1/2}(t',u) \totonemax{\leq N}^{1/2}(t',u),$
	which follow from 
	\eqref{E:FUNCTIONAVOIDINGCOMMUTING} with $f = \Psi,$
	the estimate $|\mytr \upchi| \lesssim (1 + t)^{-1}$ (that is, \eqref{E:CRUDELOWERORDERC0BOUNDDERIVATIVESOFANGULARDEFORMATIONTENSORS}),
	Prop.~\ref{P:L2NORMSOFPSIINTERMSOFTHECOERCIVEQUANTITIES}, 
	and
	\eqref{E:ROTATIONPOINTWISENORMESTIMATE}. To bound $v_3,$ we
	respectively bound $\| \Lunit F \|_{L^2(\Sigma_{t'}^u)}$ and $\| F \|_{L^2(\Sigma_{t'}^u)}$
	via the estimates \eqref{E:LESSSHARPLAPPLIEDTORENORMALIZEDTRCHIL2INTERMSOFQANDWIDETILDEQ}
	and \eqref{E:LESSSHARPSIGMATRCHIRENORMALIZEDL2INTERMSOFQANDWIDETILDEQ}.
	To bound $v_1,$ we use the 
	inequalities
	$\| \Rad \Psi \|_{C^0(\Sigma_{t'}^u)} \lesssim \varepsilon (1 + t')^{-1},$
	$\| \Rot \Rad \Psi \|_{C^0(\Sigma_{t'}^u)} \lesssim \varepsilon (1 + t')^{-1},$
	$\| \angdeformoneformupsharparg{\Rot}{\Lunit} \|_{C^0(\Sigma_{t'}^u)} \lesssim \varepsilon \ln(\myexp + t') (1 + t')^{-1},$
	$\| \angdiv \angdeformoneformupsharparg{\Rot}{\Lunit} \|_{C^0(\Sigma_{t'}^u)} \lesssim \varepsilon \ln(\myexp + t') (1 + t')^{-2},$
	$\| \mytr  \angdeform{\Rot} \|_{C^0(\Sigma_{t'}^u)} \lesssim \varepsilon \ln(\myexp + t') (1 + t')^{-1},$
	$\| \Lunit \Rad \Psi 
			 + \frac{1}{2} \mytr \upchi \Rad \Psi 
		\|_{C^0(\Sigma_{t'}^u)} \lesssim \varepsilon \ln(\myexp + t') (1 + t')^{-3},$
	and 
	$\left\| 
			\Rot
			\left\lbrace 
				\Lunit \Rad \Psi 
			 	+ \frac{1}{2} \mytr \upchi \Rad \Psi 
			\right\rbrace
		\right\|_{C^0(\Sigma_{t'}^u)}
		\lesssim \varepsilon \ln(\myexp + t') (1 + t')^{-3}.$
		The first two of these inequalities follow from
		the bootstrap assumptions \eqref{E:PSIFUNDAMENTALC0BOUNDBOOTSTRAP},
		the next two from 
		inequality \eqref{E:ONEFORMANGDINTERMSOFROTATIONALLIE},
		the estimate \eqref{E:LOWERORDERC0BOUNDROTDEFORMSPHEREL},
		and Cor.~\ref{C:SQRTEPSILONREPLCEDWITHCEPSILON},
		the fifth one from 
		\eqref{E:ROTDEFORMSPHERETRACE},
		\eqref{E:FIRSTESTIMATEFORANGULARDIFFERNTIALOFXI},
		\eqref{E:LOWERORDERC0BOUNDEUCLIDEANROTATIONRADCOMPONENT},
		\eqref{E:ROTATIONPOINTWISENORMESTIMATE},
		\eqref{E:LOWERORDERC0BOUNDLIEDERIVATIVESOFGRAME},
		\eqref{E:C0BOUNDCRUCIALEIKONALFUNCTIONQUANTITIES},
		Cor.~\ref{C:SQRTEPSILONREPLCEDWITHCEPSILON},
		and the bootstrap assumptions \eqref{E:PSIFUNDAMENTALC0BOUNDBOOTSTRAP},
		and the last two from \eqref{E:LPSIPLUSHALFTRCHIPSISTRONGC0BOUND}.
	
	Examining the products on the right-hand side of \eqref{E:NONDAMAGINGSPACETIMERENORMALZEDTRCHI},
	combining the above estimates for $v_1,$ $v_2,$ and $v_3,$ 
	using simple estimates of
	the form $ab \lesssim a^2 + b^2,$
	using the fact that 
	$\totzeromax{\leq N},$ 
	$\totonemax{\leq N},$
	$\totzeromax{\leq N-1},$ 
	and
	$\totonemax{\leq N-1}$
	are increasing in their arguments,
	and using
	the estimates \eqref{E:LESSSINGULARLOGUPMULOSSTERMSMUINTEGRALBOUND},
	\eqref{E:LESSSINGULARTERMSMUTHREEFOURTHSINTEGRALBOUND},
	and \eqref{E:LOGLOSSLESSSINGULARTERMSMTHREEFOURTHSINTEGRALBOUND},
	we deduce that all of the integrals 
	$\int_{\mathcal{M}_{t,u}} \cdots \, d \vol$
	are $\lesssim$ one of the following:
	\begin{align}
			& \varepsilon
				\int_{t'=0}^t
					\frac{1}{(1 + t')^{3/2} \upmu_{\star}^{1/2}(t',u)}
					\totzeromax{\leq N}(t',u)
				\, dt'
				\lesssim \varepsilon \totzeromax{\leq N}(t,u),
				\\
			& \varepsilon
				\int_{t'=0}^t
					\frac{1}{(1 + t')^{3/2} \upmu_{\star}(t',u)}
					\totonemax{\leq N}(t',u)
				\, dt',
				\\
		  & \varepsilon
				\int_{t'=0}^t
					\frac{1}{(1 + t')}
					\totzeromax{\leq N-1}^{1/2}(t',u)
					\totzeromax{\leq N}^{1/2}(t',u)
				\, dt'
				\lesssim 
				\varepsilon \totzeromax{\leq N}(t,u)
				+ \varepsilon \ln^2(\myexp + t) \totzeromax{\leq N-1}(t,u),
			\\
			 & \varepsilon
				\int_{t'=0}^t
					\frac{1}{(1 + t') \upmu_{\star}^{1/2}(t',u)}
					\totzeromax{\leq N-1}^{1/2}(t',u)
					\totonemax{\leq N}^{1/2}(t',u)
				\, dt'
				\lesssim 
				\varepsilon \totonemax{\leq N}(t,u)
				+ \varepsilon \ln^2(\myexp + t) \totzeromax{\leq N-1}(t,u),	
			\\
			 & \varepsilon
				\int_{t'=0}^t
					\frac{1}{(1 + t') \upmu_{\star}^{1/2}(t',u)}
					\totonemax{\leq N-1}^{1/2}(t',u)
					\totzeromax{\leq N}^{1/2}(t',u)
				\, dt'
				\lesssim 
				\varepsilon \totzeromax{\leq N}(t,u)
				+ \varepsilon \ln^2(\myexp + t) \totonemax{\leq N-1}(t,u),
				\\
			 & \varepsilon
				\int_{t'=0}^t
					\frac{1}{(1 + t') \upmu_{\star}(t',u)}
					\totonemax{\leq N-1}^{1/2}(t',u)
					\totonemax{\leq N}^{1/2}(t',u)
				\, dt'
				\lesssim 
				\varepsilon \totonemax{\leq N}(t,u)
				\\
				& \ \ + \varepsilon \ln^2(\myexp + t) \frac{1}{\upmu_{\star}^{1/2}(t,u)} \totonemax{\leq N-1}(t,u),
				\notag \\
			& \varepsilon^3
				\int_{t'=0}^t
					\frac{1}{(1 + t')^{3/2} \upmu_{\star}^{1/2}(t',u)}
				\, dt'
				\lesssim \varepsilon^3.
	\end{align}
	We have therefore proved the desired estimate \eqref{E:NONDAMAGINGSPACETIMERENORMALZEDTRCHIESTIMATE}.
	
	Applying similar reasoning to the terms on the right-hand side of \eqref{E:NONDAMAGINGHYPERSURFACERENORMALZEDTRCHI}
	(without having to integrate in time),
	we deduce \eqref{E:NONDAMAGINGHYPERSURFACERENORMALZEDTRCHIESTIMATE}.
	
	Finally, the estimate \eqref{E:NONDAMAGINGINITIALDATARENORMALZEDTRCHIESTIMATE} follows easily 
	from the estimates of Sect.~\ref{S:INITIALBEHAVIOROFQUANTITIES}
	since the integrands
	on the right-hand side of \eqref{E:NONDAMAGINGINITIALDATARENORMALZEDTRCHI} are cubic and depend on the data.
\end{proof}

\begin{lemma}[\textbf{Estimates for some easy error integrals corresponding to the partially modified version of}
$\angLap \mathscr{S}^{N-1} \upmu$]
	\label{L:NONDAMAGINGTERMSRENORMALIZEDANGLAPUPMU}
	Let $1 \leq N \leq 24$ be an integer, and let
	$\mathscr{S}^{N-1}$ be an $(N-1)^{st}$ order pure spatial commutation vectorfield operator
	(see Def.~\ref{D:DEFSETOFSPATIALCOMMUTATORVECTORFIELDS}).
	Let $\mupartialmodarg{\mathscr{S}^{N-1}}$ be the partially modified 
	$S_{t,u}$ one-form defined in \eqref{E:TRANSPORTPARTIALRENORMALIZEDUPMU}.
	Consider the $S_{t,u}-$tangent vectorfield 
	$Y:= \mupartialmodarg{\mathscr{S}^{N-1}}^{\#},$
	and define the weight function $w(t,u) := \rgeo^2(t,u).$
	Under the small-data and bootstrap assumptions 
	of Sects.~\ref{S:PSISOLVES}-\ref{S:C0BOUNDBOOTSTRAP},
	if $\varepsilon$ is sufficiently small, 
	then we have the following spacetime integral estimate for
	the quantity $\IBPseconderror{1}[\mathscr{S}^{N-1} \Rad \Psi;w]$
	defined in \eqref{E:NONDAMAGINGSPACETIMERENORMALZEDANGLAPUPMU}
	for $(t,u) \in [0,\Tboot) \times [0,U_0]:$
	\begin{subequations}
	\begin{align} \label{E:NONDAMAGINGSPACETIMERENORMALZEDANGLAPUPMUESTIMATE}
			\int_{\mathcal{M}_{t,u}}
				\left|\IBPseconderror{1}[\mathscr{S}^{N-1} \Rad \Psi;w]\right|
			\, d \vol
			& \lesssim 
				\varepsilon \totzeromax{\leq N}(t,u)
				+
				\varepsilon
				\int_{t'=0}^t
					\frac{1}{(1 + t')^{3/2} \upmu_{\star}(t',u)}
					\totonemax{\leq N}(t',u)
				\, dt'
				\\
	& \ \ + \varepsilon \ln^2(\myexp + t) \totzeromax{\leq N-1}(t,u) 
				+ \varepsilon \ln^2(\myexp + t) \frac{1}{\upmu_{\star}^{1/2}(t,u)} \totonemax{\leq N-1}(t,u)
				\notag \\
	& \ \ + \varepsilon^3.
			\notag
	\end{align}
	
	Furthermore, we have the following hypersurface integral estimate for
	the quantity $\IBPseconderror{2}[\Rad \mathscr{S}^{N-1};w]$
	defined in \eqref{E:NONDAMAGINGHYPERSURFACERENORMALZEDANGLAPUPMU}:
 \begin{align}  \label{E:NONDAMAGINGHYPERSURFACERENORMALZEDANGLAPUPMUESTIMATE}
			\int_{\Sigma_t^u}
				\left|\IBPseconderror{2}[\mathscr{S}^{N-1} \Rad \Psi;w] \right|
			\, d \tvol
			& \lesssim 
				\varepsilon \totzeromax{\leq N}(t,u) 
				+ \varepsilon \totonemax{\leq N}(t,u)
				\\
			& \ \ + \varepsilon \ln^2(\myexp + t) \frac{1}{\upmu_{\star}(t,u)} \totzeromax{\leq N-1}(t,u) 
				+ \varepsilon \ln^2(\myexp + t) \frac{1}{\upmu_{\star}(t,u)} \totonemax{\leq N-1}(t,u)
				\notag \\
			& \ \ + \varepsilon^3 \frac{1}{\upmu_{\star}(t,u)}.
			\notag
\end{align}

	Finally, we have the following hypersurface integral estimate for
	the quantity $\IBPseconderror{3}[\Rad \mathscr{S}^{N-1};w]$
	defined in \eqref{E:NONDAMAGINGINITIALDATARENORMALZEDANGLAPUPMU}:
	\begin{align}  \label{E:NONDAMAGINGINITIALDATARENORMALZEDANGLAPUPMUESTIMATE}
			\int_{\Sigma_0^u}
				\left|\IBPseconderror{3}[\mathscr{S}^{N-1} \Rad \Psi;w] \right|
			\, d \tvol
			& \lesssim \varepsilon^3.
	\end{align}
	\end{subequations}	
\end{lemma}

\begin{proof}
We apply reasoning similar to the reasoning that we used to prove Lemma~\ref{L:NONDAMAGINGTERMSRENORMALIZEDTRCHI}.
We use the estimates
\eqref{E:LESSSHARPSIGMAANGDIFFUPMURENORMALIZEDL2INTERMSOFQANDWIDETILDEQ}
and
\eqref{E:LESSSHARPLPLUSTRCHIAPPLIEDTORENORMALIZEDANGDIFFUPMUL2INTERMSOFQANDWIDETILDEQ}
in place of the estimates
\eqref{E:LESSSHARPSIGMATRCHIRENORMALIZEDL2INTERMSOFQANDWIDETILDEQ}
and \eqref{E:LESSSHARPLAPPLIEDTORENORMALIZEDTRCHIL2INTERMSOFQANDWIDETILDEQ}
used in the proof of Lemma~\ref{L:NONDAMAGINGTERMSRENORMALIZEDTRCHI}.
\end{proof}

In Sects.~\ref{S:PROOFOFLEMMATOPORDERELLIPTICRECOVERY}-\ref{S:PROOFOFLEMMADANGEROUSTOPORDERMORERRORINTEGRAL},
we use the auxiliary lemmas to prove the three lemmas of primary interest, 
namely Lemmas \ref{L:TOPORDERELLIPTICRECOVERY}, \ref{L:DANGEROUSTOPORDERMULTERRORINTEGRAL}, and \ref{L:DANGEROUSTOPORDERMORERRORINTEGRAL}.
We begin with the proof of Lemma~\ref{L:TOPORDERELLIPTICRECOVERY}.

\section{Proof of Lemma~\ref{L:TOPORDERELLIPTICRECOVERY}}
\label{S:PROOFOFLEMMATOPORDERELLIPTICRECOVERY}
We now prove Lemma~\ref{L:TOPORDERELLIPTICRECOVERY}.
To deduce \eqref{E:ELLIPTICRECOVERY}, we show that
\begin{align} \label{E:KEYGRONWALLABLEESTIMATEFORELLIPTICQUANTITIES}
	\fivemyarray
		[\| \upmu \angD^2 \mathscr{Z}^{\leq N-1} \upmu \|_{L^2(\Sigma_t^u)}]
		{\sum_{a=1}^3 \|\rgeo \upmu \angD^2 \mathscr{Z}^{\leq N-1} \Lunit_{(Small)}^a \|_{L^2(\Sigma_t^u)}}		
		{\| \rgeo \upmu \angD \angLie_{\mathscr{Z}}^{\leq N-1} \upchi^{(Small)} \|_{L^2(\Sigma_t^u)}}
		{\| \rgeo \upmu \angdiff \mathscr{Z}^{\leq N-1} \mytr \upchi^{(Small)} \|_{L^2(\Sigma_t^u)}}
		{\| \rgeo \upmu \angD \angLie_{\mathscr{Z}}^{\leq N-1} \hat{\upchi}^{(Small)} \|_{L^2(\Sigma_t^u)}}
	& \lesssim 
		\int_{t'=0}^t 
					\frac{\| [\Lunit \upmu]_- \|_{C^0(\Sigma_{t'}^u)}} 
							 {\upmu_{\star}(t',u)} 
					\totzeromax{\leq N}^{1/2}(t',u) 
				\, dt'	
				\\
		& \ \
			+ \int_{t'=0}^t 
					\frac{1}{(1 + t')^{3/2} \upmu_{\star}(t',u)} 
					\totzeromax{\leq N}^{1/2}(t',u)
				\, dt'
			+ \totzeromax{\leq N}^{1/2}(t,u)
			\notag \\
	& \ \ 
			+ \int_{t'=0}^t 
					\frac{1}{(1 + t')^{3/2} \upmu_{\star}(t',u)} 
					\totonemax{\leq N}^{1/2}(t',u)
				\, dt'
			+ 
			\varepsilon	
			\left\lbrace
				\upmu_{\star}^{-1}(t,u) + 1
			\right\rbrace
			\notag \\
	& \ \ 
		+ \varepsilon
			\int_{t'=0}^t
				\frac{1}{(1 + t')^{3/2}}
				\myarray[\| \upmu \angD^2 \mathscr{Z}^{\leq N-1} \upmu \|_{L^2(\Sigma_{t'}^u)}]
						{\| \rgeo \upmu \angD \angLie_{\mathscr{Z}}^{\leq N-1} \upchi^{(Small)} \|_{L^2(\Sigma_{t'}^u)}}
			\, dt'.
			\notag
\end{align}
The desired estimate \eqref{E:ELLIPTICRECOVERY} 
then follows from applying Gronwall's inequality to the quantity on the left-hand side of 
\eqref{E:KEYGRONWALLABLEESTIMATEFORELLIPTICQUANTITIES}.

It remains for us to prove \eqref{E:KEYGRONWALLABLEESTIMATEFORELLIPTICQUANTITIES}. We first prove the desired estimate for
$\myarray
			[\| \rgeo \upmu \angD \angLie_{\mathscr{Z}}^{\leq N-1} \upchi^{(Small)} \|_{L^2(\Sigma_t^u)}]
				{\| \rgeo \upmu \angD \angLie_{\mathscr{Z}}^{\leq N-1} \hat{\upchi}^{(Small)} \|_{L^2(\Sigma_t^u)}}.
$ 
To this end,
we first use \eqref{E:L2COMMUTEDANGDLIECHIJUNKINTERMSOFTRACEFREEANDTRACEPART} 
and \eqref{E:L2COMMUTEDANGDTRACEFREECHIJUNKINTERMSOFCHIJUNK}
to deduce that
\begin{align} \label{E:FIRSTSTEPTOPORDERANGDCHIJUNKDERIVATIVESELLIPTIC}
	\myarray
			[\| \rgeo \upmu \angD \angLie_{\mathscr{Z}}^{\leq N-1} \upchi^{(Small)} \|_{L^2(\Sigma_t^u)}]
				{\| \rgeo \upmu \angD \angLie_{\mathscr{Z}}^{\leq N-1} \hat{\upchi}^{(Small)} \|_{L^2(\Sigma_t^u)}}
	& \lesssim
				\sum_{l=1}^3
				\left\| \upmu \Rot_{(l)} \mathscr{Z}^{\leq N-1} \mytr \upchi^{(Small)} \right\|_{L^2(\Sigma_t^u)}
			+ \left\| \rgeo \upmu \angD \angfreeLietwoarg{\mathscr{Z}}{N-1} \hat{\upchi}^{(Small)} \right\|_{L^2(\Sigma_t^u)}
				\\
		& \ \ + \totzeromax{\leq N}^{1/2}(t,u)
			+ 
			\int_{t'=0}^t 
				\frac{1}{(1 + t')^{3/2} \upmu_{\star}^{1/2}(t',u)}
				\totonemax{\leq N}^{1/2}(t',u)
			\, d t'
			+ \varepsilon.
			\notag
\end{align}
Next, applying the elliptic estimate \eqref{E:TYPE02ELLIPTICESTIMATESPHERES}
to $\angD \angfreeLietwoarg{\mathscr{Z}}{N-1} \hat{\upchi}^{(Small)}$ and using
inequalities \eqref{E:L2COMMUTEDANGDTRACEFREECHIJUNKINTERMSOFTRACEPART} 
and
\eqref{E:L2LOWERORDERCOMMUTEDANGDTRACEFREECHIJUNKINTERMSOFTRACEPART},
we bound the second term on the right-hand side of \eqref{E:FIRSTSTEPTOPORDERANGDCHIJUNKDERIVATIVESELLIPTIC}
as follows:
\begin{align} \label{E:SECONDSTEPTOPORDERANGDCHIJUNKDERIVATIVESELLIPTIC}
 \left\| \rgeo \upmu \angD \angfreeLietwoarg{\mathscr{Z}}{N-1} \hat{\upchi}^{(Small)} \right\|_{L^2(\Sigma_t^u)}
 & \lesssim 
 	\left\| \rgeo \upmu \angdiv \angfreeLietwoarg{\mathscr{Z}}{N-1} \hat{\upchi}^{(Small)} \right\|_{L^2(\Sigma_t^u)}
 	+ \varepsilon \frac{\ln(\myexp + t)}{1 + t}
 		\left\| \rgeo \angfreeLietwoarg{\mathscr{Z}}{N-1} \hat{\upchi}^{(Small)} \right\|_{L^2(\Sigma_t^u)}
 		\\
 	& \lesssim
				\sum_{l=1}^3
				\left\| \upmu \Rot_{(l)} \mathscr{Z}^{\leq N-1} \mytr \upchi^{(Small)} \right\|_{L^2(\Sigma_t^u)}
				\notag \\
	& \ \ + \totzeromax{\leq N}^{1/2}(t,u)
			+ 
			\int_{t'=0}^t 
				\frac{1}{(1 + t')^{3/2} \upmu_{\star}^{1/2}(t',u)}
			\totonemax{\leq N}^{1/2}(t',u)
			\, d t'
		+ \varepsilon. 
		\notag
\end{align}
With the exception of
$\left\| \upmu \Rot \mathscr{Z}^{\leq N-1} \mytr \upchi^{(Small)} \right\|_{L^2(\Sigma_t^u)},$
all terms on the right-hand side of 
\eqref{E:FIRSTSTEPTOPORDERANGDCHIJUNKDERIVATIVESELLIPTIC}
and
\eqref{E:SECONDSTEPTOPORDERANGDCHIJUNKDERIVATIVESELLIPTIC}
are clearly $\lesssim$ the right-hand side of
\eqref{E:KEYGRONWALLABLEESTIMATEFORELLIPTICQUANTITIES} as desired.
To bound
$\left\| \upmu \Rot \mathscr{Z}^{\leq N-1} \mytr \upchi^{(Small)} \right\|_{L^2(\Sigma_t^u)},$
we use inequality
\eqref{E:SLIGHTLYLESSSHARPL2FORANLGAPUPMUANDROTTRCHI},
the simple inequality
$\left\| \upmu \angfreeDsquared \mathscr{Z}^{\leq N-1} \upmu \right\|_{L^2(\Sigma_t^u)}
\leq 
\left\| \upmu \angD^2 \mathscr{Z}^{\leq N-1} \upmu \right\|_{L^2(\Sigma_t^u)},$
inequality \eqref{E:L2COMMUTEDANGDTRACEFREECHIJUNKINTERMSOFCHIJUNK},
and the estimate
$\int_{s=0}^{t'} 
				\frac{1}{(1 + s)^{3/2} \upmu_{\star}^{1/2}(s,u)}
				\totonemax{\leq N}^{1/2}(s,u)
			\, ds
\lesssim \totonemax{\leq N}^{1/2}(t',u)$
(which follows from \eqref{E:LESSSINGULARTERMSMUTHREEFOURTHSINTEGRALBOUND}
and the fact that $\totonemax{\leq N}$ is increasing in its arguments)
to deduce that
\begin{align}   \label{E:THIRDSTEPTOPORDERANGDCHIJUNKDERIVATIVESELLIPTIC}
		& \left\| \upmu \Rot \mathscr{Z}^{\leq N-1} \mytr \upchi^{(Small)} \right\|_{L^2(\Sigma_t^u)}
			\\
		& \lesssim 
			\int_{t'=0}^t 
					\frac{\| [\Lunit \upmu]_- \|_{C^0(\Sigma_{t'}^u)}} 
							 {\upmu_{\star}(t',u)} 
					\totzeromax{\leq N}^{1/2}(t',u) 
				\, dt'
			+ 
			\int_{t'=0}^t 
				\frac{1}{(1 + t')^{3/2} \upmu_{\star}(t',u)}
				\totzeromax{\leq N}^{1/2}(t',u)
			\, d t'		
			\notag	\\
		& \ \ + \totzeromax{\leq N}^{1/2}(t,u)
			+ 
			\int_{t'=0}^t 
				\frac{1}{(1 + t')^{3/2} \upmu_{\star}(t',u)}
				\totonemax{\leq N}^{1/2}(t',u)
			\, d t'
			\notag \\
		& \ \ 
			+ \varepsilon
				\int_{t'=0}^t
					\frac{1}{(1 + t')^{3/2}}
					\myarray[\left\| \upmu \angD^2 \mathscr{Z}^{\leq N-1} \upmu \right\|_{L^2(\Sigma_{t'}^u)}]
						{\left\| \rgeo \upmu \angD \angLie_{\mathscr{Z}}^{\leq N-1} \upchi^{(Small)} \right\|_{L^2(\Sigma_{t'}^u)}}
				\, dt'
			+ \varepsilon
				\left\lbrace
					\ln \upmu_{\star}^{-1}(t,u) + 1
				\right\rbrace.
				\notag
	\end{align}
	Clearly, the right-hand side of \eqref{E:THIRDSTEPTOPORDERANGDCHIJUNKDERIVATIVESELLIPTIC}
	is $\lesssim$ the right-hand side of
	\eqref{E:KEYGRONWALLABLEESTIMATEFORELLIPTICQUANTITIES} as desired.
	We have thus bounded 
	$\| \rgeo \upmu \angD \angLie_{\mathscr{Z}}^{\leq N-1} \upchi^{(Small)} \|_{L^2(\Sigma_t^u)}$
	by the right-hand side of \eqref{E:KEYGRONWALLABLEESTIMATEFORELLIPTICQUANTITIES} as desired.
	Using in addition inequality \eqref{E:FUNCTIONPOINTWISEANGDINTERMSOFANGLIEO},
	we deduce from \eqref{E:THIRDSTEPTOPORDERANGDCHIJUNKDERIVATIVESELLIPTIC} that
	$\| \rgeo \upmu \angdiff \mathscr{Z}^{\leq N-1} \mytr \upchi^{(Small)} \|_{L^2(\Sigma_t^u)}$
	is also $\lesssim$ the right-hand side of
	\eqref{E:KEYGRONWALLABLEESTIMATEFORELLIPTICQUANTITIES} as desired.
	
	The proof that $\| \upmu \angD^2 \mathscr{Z}^{\leq N-1} \upmu \|_{L^2(\Sigma_t^u)}$
	is bounded by the right-hand side of \eqref{E:KEYGRONWALLABLEESTIMATEFORELLIPTICQUANTITIES} 
	is similar but simpler because the elliptic estimates involve a scalar quantity rather than a trace-free tensorial one. 
	Specifically, 
	we first use inequalities 
	\eqref{E:FUNCTIONPOINTWISEANGDINTERMSOFANGLIEO}, 
	\eqref{E:TYPE02TENSORANGDINTERMSOFROTATIONALLIE},
	\eqref{E:EIKONALFUNCTIONQUANTITIESL2BOUNDSINTERMSOFQ0ANDQ1},
	and the elliptic estimate \eqref{E:POISSONELLIPTIC} 
	to deduce that
	\begin{align}
	\left\| \upmu \angD^2 \mathscr{Z}^{\leq N-1}  \upmu \right\|_{L^2(\Sigma_t^u)}
	& \lesssim 
		\left\| \upmu \angLap \mathscr{Z}^{\leq N-1}  \upmu \right\|_{L^2(\Sigma_t^u)}
			\label{E:TOPORDERUPMUANGULARDERIVATIVESELLIPTIC}  \\
	& \ \ 
		+ \totzeromax{\leq N}^{1/2}(t,u)
		+ \int_{t'=0}^t
				\frac{1}{(1 + t')^{3/2} \upmu_{\star}^{1/2}(t',u)}
				\totonemax{\leq N}^{1/2}(t',u)
			\, dt'
		+ \varepsilon.
			\notag 
\end{align}
All terms on the right-hand side of \eqref{E:TOPORDERUPMUANGULARDERIVATIVESELLIPTIC}
except for the first one are manifestly bounded 
by the right-hand side of \eqref{E:KEYGRONWALLABLEESTIMATEFORELLIPTICQUANTITIES} as desired.
To bound the remaining term $\left\| \upmu \angLap \mathscr{Z}^{\leq N-1} \upmu \right\|_{L^2(\Sigma_t^u)}$
by the right-hand side of \eqref{E:KEYGRONWALLABLEESTIMATEFORELLIPTICQUANTITIES},
we use inequality \eqref{E:SLIGHTLYLESSSHARPL2FORANLGAPUPMUANDROTTRCHI}
and then argue as in the previous paragraph to bound this quantity
by the right-hand side of \eqref{E:THIRDSTEPTOPORDERANGDCHIJUNKDERIVATIVESELLIPTIC} as desired.

To complete the proof of \eqref{E:FIRSTSTEPTOPORDERANGDCHIJUNKDERIVATIVESELLIPTIC}, it only remains for
us to bound $\sum_{a=1}^3 \|\rgeo \upmu \angD^2 \mathscr{Z}^{\leq N-1} \Lunit_{(Small)}^a \|_{L^2(\Sigma_t^u)}$
by the right-hand side of \eqref{E:FIRSTSTEPTOPORDERANGDCHIJUNKDERIVATIVESELLIPTIC}.
The proof is based on the inequality 
\eqref{E:L2COMMUTEDANGLAPLJUNKIINTERMSOFANGDCHIANDPSI}
and the elliptic estimate \eqref{E:POISSONELLIPTIC} .
We omit the details because they are very similar to the bound for 
$\| \upmu \angD^2 \mathscr{Z}^{\leq N-1} \upmu \|_{L^2(\Sigma_t^u)}$ 
proved in the previous paragraph.

To prove \eqref{E:NONSHARPTOPORDERELLIPTICRECOVERY},
we first use inequalities
\eqref{E:KEYMUINVERSEINTEGRALBOUND}
and
\eqref{E:LESSSINGULARLOGUPMULOSSTERMSMUINTEGRALBOUND}
and the fact that $\totzeromax{\leq N}$
and $\totonemax{\leq N}$ are increasing in both of their arguments
to deduce that the right-hand side of 
\eqref{E:ELLIPTICRECOVERY} 
is 
\begin{align} \label{E:TIMEINTEGRATEDELLIPTICRECOVERYALMOSTPROVED}
\lesssim 
	\left\lbrace
		\ln \upmu_{\star}^{-1}(t,u) + 1
	\right\rbrace
	\totzeromax{\leq N}^{1/2}(t,u) 
	+ 
	\left\lbrace
		\ln \upmu_{\star}^{-1}(t,u) + 1
	\right\rbrace
	\totonemax{\leq N}^{1/2}(t,u) 
	+ \varepsilon
	\left\lbrace
			\ln \upmu_{\star}^{-1}(t,u) + 1	
	\right\rbrace.
\end{align}			
Next, we combine inequality \eqref{E:TIMEINTEGRATEDELLIPTICRECOVERYALMOSTPROVED} with 
Lemma~\ref{L:L2NORMSOFTIMEINTEGRATEDFUNCTIONS},
which allows us to bound the left-hand side of 
\eqref{E:NONSHARPTOPORDERELLIPTICRECOVERY} by 
\begin{align} \label{E:FIRSTESTIMATENONSHARPTOPORDERELLIPTICRECOVERY}
	&
	\lesssim
	(1 + t)
	\int_{t'=0}^t
		\frac{1}{1 + t'}
		(\ln \upmu_{\star}^{-1}(t,u) + \sqrt{\varepsilon})\totzeromax{\leq N}^{1/2}(t',u)
	\, dt'
	+ 
	(1 + t)
	\int_{t'=0}^t
		\frac{1}{1 + t'}
		\totonemax{\leq N}^{1/2}(t',u)
	\, dt'
		\\
	&
	\ \
	+ 
	\varepsilon
	(1 + t)
	\int_{t'=0}^t
		\frac{1}{1 + t'}
		\left\lbrace
			\ln \upmu_{\star}^{-1}(t,u) + 1
		\right\rbrace
	\, dt'.
	\notag
\end{align}
Again using 
the fact that $\totzeromax{\leq N}$
and $\totonemax{\leq N}$ are increasing in both of their arguments 
and also using \eqref{E:LOGLOSSLESSSINGULARTERMSMTHREEFOURTHSINTEGRALBOUND},
we deduce that
the right-hand side of
\eqref{E:FIRSTESTIMATENONSHARPTOPORDERELLIPTICRECOVERY}
is 
\[
\lesssim
	\ln(\myexp + t) (1 + t) \totzeromax{\leq N}^{1/2}(t,u)
	+ \ln(\myexp + t) (1 + t) \totonemax{\leq N}^{1/2}(t,u)
	+ \varepsilon \ln(\myexp + t) (1 + t)
\]
as desired.

$\hfill \qed$

\section{Proof of Lemma~\ref{L:DANGEROUSTOPORDERMULTERRORINTEGRAL}}
\label{S:PROOFOFLEMMADANGEROUSTOPORDERMULTERRORINTEGRAL}
We now use the auxiliary lemmas to prove the first lemma of primary interest, 
namely Lemma~\ref{L:DANGEROUSTOPORDERMULTERRORINTEGRAL}.
	We give the proof for the first integral on the left-hand side of \eqref{E:MULTMAINENERGYFLUXERRORINTEGRALESTIMATE}.
	The proof for the second one is nearly identical, and we describe the minor differences at the end of the proof. 
	To begin, we split the integrand into two pieces as follows:
	\begin{align} \label{E:DIFFICULTMULTINTEGRANDSPLITTING}
		& \left\lbrace
			(1 + 2 \upmu) \Lunit \mathscr{S}^{N-1} \Rad \Psi
				+ 2 \Rad \mathscr{S}^{N-1} \Rad  \Psi
		\right\rbrace
				(\Rad \Psi) \angLap \mathscr{S}^{N-1} \upmu
				\\
		& = (1 + 2 \upmu) (\Lunit \mathscr{S}^{N-1} \Rad \Psi) (\Rad \Psi) \angLap \mathscr{S}^{N-1} \upmu
				+ 2 (\Rad \mathscr{S}^{N-1} \Rad \Psi) (\Rad \Psi) \angLap \mathscr{S}^{N-1} \upmu.
				\notag
	\end{align}
				The integral of the first product on the right-hand side of 
				\eqref{E:DIFFICULTMULTINTEGRANDSPLITTING} is much easier to bound than that of the second; 
				we bound the integral of the easier product at the end of the proof using a separate argument.
	
	\noindent{\emph{Bound for} $\int_{\mathcal{M}_{t,u}}
				(\Rad \mathscr{S}^{N-1} \Rad \Psi)
				(2 \Rad \Psi) \angLap \mathscr{S}^{N-1} \upmu
				\, d \vol:$}
	Our goal is to bound the spacetime integral
	\begin{align} \label{E:DIFFICULTMULTINTEGRALBOUNDPROOF}
		\int_{\mathcal{M}_{t,u}}
				(2 \Rad \mathscr{S}^{N-1} \Rad \Psi)
				(\Rad \Psi) \angLap \mathscr{S}^{N-1} \upmu
				\, d \vol
	\end{align}
	in magnitude by the right-hand side of \eqref{E:MULTMAINENERGYFLUXERRORINTEGRALESTIMATE}. 
	To this end, we first use the final inequality 
	in \eqref{E:TOPORDERDERIVATIVESOFANGDSQUAREDUPMUINTERMSOFCONTROLLABLE},
	the bootstrap assumption $\| \Rad \Psi \|_{C^0(\Sigma_t^u)} \leq \varepsilon (1 + t)^{-1},$ 
	and Def.~\ref{D:HARMLESSTERMS} 
	to deduce that
	\begin{align} \label{E:ANGLAPUPMUISRADTRCHIJUNKPLUSHARMLESS}
		(\Rad \Psi) \angLap \mathscr{S}^{N-1} \upmu 
		= (\Rad \Psi) \Rad \mathscr{S}^{N-1} \mytr \upchi^{(Small)}
		+ Harmless^{\leq N}.
	\end{align}
	By Cor.~\ref{C:EASYHARMLESSMULTERRORINTEGRAL}, the part of the spacetime integral \eqref{E:DIFFICULTMULTINTEGRALBOUNDPROOF}
	involving the $Harmless^{\leq N}$ terms from \eqref{E:ANGLAPUPMUISRADTRCHIJUNKPLUSHARMLESS}
	has already been bounded by the right-hand side of \eqref{E:MULTMAINENERGYFLUXERRORINTEGRALESTIMATE}.
	It remains for us to bound the spacetime integral \eqref{E:DIFFICULTMULTINTEGRALBOUNDPROOF}
	with $\angLap \mathscr{S}^{N-1} \upmu$ replaced by 
	$\Rad \mathscr{S}^{N-1} \mytr \upchi^{(Small)}.$
	To this end, we first express it as a time integral of
	integrals over $\Sigma_{t'}^u:$
	$\int_{\mathcal{M}_{t,u}} \cdots \, d\vol = 
	\int_{t'=0}^t \int_{\Sigma_{t'}^u} \cdots \, d \tvol \, d t'.$
	Next, using Cauchy-Schwarz on $\Sigma_{t'}^u$
	and inequality \eqref{E:RADPSIL2INTERMSOFZERO},
	we deduce that
	\begin{align}  \label{E:DIFFICULTMULTINTEGRALOFCONSTANTTIMEINTEGRALBOUNDPROOF}
		&
		\left|
			\int_{\mathcal{M}_{t,u}}
				(2 \Rad \mathscr{S}^{N-1} \Rad \Psi) 
				(\Rad \Psi) \Rad \mathscr{S}^{N-1} \mytr \upchi^{(Small)}
			\, d \vol
		\right|
			\\
		& \leq
			2
			\int_{t'=0}^t
				\| \Rad \mathscr{S}^{N-1} \Rad \Psi \|_{L^2(\Sigma_{t'}^u)}
				\| (\Rad \Psi) \Rad \mathscr{S}^{N-1} \mytr \upchi^{(Small)} \|_{L^2(\Sigma_{t'}^u)}
			\, dt'
				\notag \\
		& \leq 
			2 \int_{t'=0}^t
					\totzeromax{\leq N}^{1/2}(t',u)	
					\| (\Rad \Psi) \Rad \mathscr{S}^{N-1} \mytr \upchi^{(Small)} \|_{L^2(\Sigma_{t'}^u)}
				\, dt'.
			 \notag
	\end{align}
	The main part of the argument consists of showing that there exists a small constant $\Littleconone > 0$ such that 
	the following key inequality holds for the integrand factor on the right-hand side of 
	\eqref{E:DIFFICULTMULTINTEGRALOFCONSTANTTIMEINTEGRALBOUNDPROOF}:
	\begin{align} \label{E:DIFFICULTMULTCONSTANTTIMEINTEGRALBOUNDKEYL2ESTIMATE}
		& \left\| 
			(\Rad \Psi) 
			\Rad \mathscr{S}^{N-1} \mytr \upchi^{(Small)} 
		\right\|_{L^2(\Sigma_t^u)}
			\\
		& \leq 
				\boxed{4.5}
							\frac{\| [\Lunit \upmu]_- \|_{C^0(\Sigma_t^u)}} 
									 {\upmu_{\star}(t,u)} 
						\int_{s=0}^t
							\frac{\| [\Lunit \upmu]_- \|_{C^0(\Sigma_s^u)}} 
									{\upmu_{\star}(s,u)} 
							\totzeromax{\leq N}^{1/2}(s,u) 
					\, ds	
				\notag		 \\			
	 & \ \  +	C \varepsilon 
								\frac{1}{(1 + t)^{1 + \Littleconone}} 
								\frac{1}{\upmu_{\star}(t,u)} 
						\int_{s=0}^t
							\frac{1}{(1 + s)}
							\frac{1}{\upmu_{\star}(s,u)} 
							\totzeromax{\leq N}^{1/2}(s,u) 
					\, ds	
						\notag \\	
	& \ \  +	C \varepsilon 
								\frac{1}{(1 + t)^{1 + \Littleconone}} 
								\frac{1}{\upmu_{\star}(t,u)} 
						\int_{s=0}^t
							\frac{1}{(1 + s)}
							\frac{1}{\upmu_{\star}(s,u)} 
							\totonemax{\leq N}^{1/2}(s,u) 
					\, ds	
						\notag \\	
		& \ \ +   \boxed{4.5} 
							\frac{\| [\Lunit \upmu]_- \|_{C^0(\Sigma_t^u)}} 
									   {\upmu_{\star}(t,u)} 
								\totzeromax{\leq N}^{1/2}(t,u)
						\notag \\	
		& \ \ + 	4.5
							\left\| 
								\frac{[\Lunit \upmu]_+}{\upmu} 
							\right\|_{C^0(\Sigma_t^u)}
							\int_{s=0}^t 
							\frac{\| [\Lunit \upmu]_- \|_{C^0(\Sigma_s^u)}} 
									{\upmu_{\star}(s,u)} 
							\totzeromax{\leq N}^{1/2}(s,u) 
						\, ds	
						\notag \\	
		& \ \ +   4.5 
							\left\| 
								\frac{[\Lunit \upmu]_+}{\upmu} 
							\right\|_{C^0(\Sigma_t^u)}
							\totzeromax{\leq N}^{1/2}(t,u) 
							\notag \\
		& \ \ + C \varepsilon  
							\frac{1}{(1 + t)^{1 + \Littleconone}} 
							\frac{1}{\upmu_{\star}(t,u)}
							\totzeromax{\leq N}^{1/2}(t,u)
			+ C \varepsilon 
					\frac{1}{(1 + t)^{1 + \Littleconone}} 
					\frac{1}{\upmu_{\star}(t,u)}
					\totonemax{\leq N}^{1/2}(t,u)
					\notag \\	
		& \ \ 	
					+ C \varepsilon
						\frac{\ln^3(\myexp + t)}{1 + t} \frac{1}{\upmu_{\star}^{3/2}(t,u)} \totzeromax{\leq N-1}(t,u)
					+ C \varepsilon
						\frac{1}{(1 + t)^{3/2}} \frac{1}{\upmu_{\star}^{3/2}(t,u)} \totonemax{\leq N-1}(t,u)
						\notag \\
		& \ \ + C \varepsilon^2 \frac{1}{(1 + t)^{3/2}} \frac{1}{\upmu_{\star}^{3/2}(t,u)}.
			\notag
\end{align}	
Once we have shown \eqref{E:DIFFICULTMULTCONSTANTTIMEINTEGRALBOUNDKEYL2ESTIMATE},
in order to bound the right-hand side of \eqref{E:DIFFICULTMULTINTEGRALOFCONSTANTTIMEINTEGRALBOUNDPROOF}
by the right-hand side of \eqref{E:MULTMAINENERGYFLUXERRORINTEGRALESTIMATE}, 
we insert inequality \eqref{E:DIFFICULTMULTCONSTANTTIMEINTEGRALBOUNDKEYL2ESTIMATE}
(with $t$ replaced by $t'$)
into the right-hand side of \eqref{E:DIFFICULTMULTINTEGRALOFCONSTANTTIMEINTEGRALBOUNDPROOF}
and then integrate the inequality from $t'=0$ to $t'=t.$ 
We use inequality \eqref{E:POSITIVEPARTOFLMUOVERMUISSMALL} to bound the factors 
$\left\| 
	 \frac{[\Lunit \upmu]_+}{\upmu} 
 \right\|_{C^0(\Sigma_{t'}^u)}
$
arising from the right-hand side of \eqref{E:DIFFICULTMULTCONSTANTTIMEINTEGRALBOUNDKEYL2ESTIMATE}.							
The factor of $2$ on the right-hand side of \eqref{E:DIFFICULTMULTINTEGRALOFCONSTANTTIMEINTEGRALBOUNDPROOF} 
results in the doubling of the ``boxed'' 
constants on the right-hand side of \eqref{E:DIFFICULTMULTCONSTANTTIMEINTEGRALBOUNDKEYL2ESTIMATE}, that is,
the boxed constants in \eqref{E:MULTMAINENERGYFLUXERRORINTEGRALESTIMATE}							
are twice as large as the ones in \eqref{E:DIFFICULTMULTCONSTANTTIMEINTEGRALBOUNDKEYL2ESTIMATE}.
These estimates allow us to bound the integrals corresponding to 
all but the last $5$ terms on the right-hand side of \eqref{E:DIFFICULTMULTCONSTANTTIMEINTEGRALBOUNDKEYL2ESTIMATE}
by the right-hand side of \eqref{E:MULTMAINENERGYFLUXERRORINTEGRALESTIMATE}.
We now explain how to suitably bound these last $5$ terms.
Specifically, to bound the integral arising from the term
$C \varepsilon \ln^3(\myexp + t)(1 + t)^{-1} \upmu_{\star}^{-3/2}(t,u) \totzeromax{\leq N-1}(t,u)$
on the right-hand side of \eqref{E:DIFFICULTMULTCONSTANTTIMEINTEGRALBOUNDKEYL2ESTIMATE} 
by the right-hand side of \eqref{E:MULTMAINENERGYFLUXERRORINTEGRALESTIMATE},
we use the fact that 
$\totzeromax{\leq N-1}$ 
and $\totzeromax{\leq N}$
are increasing in both of their arguments
and inequality \eqref{E:LOGLOSSKEYMUINTEGRALBOUND} to deduce that
\begin{align}
	&
	\varepsilon
	\int_{t'=0}^t 
			\frac{\ln^3(\myexp + t')}{(1 + t')} 
			\frac{1}{\upmu_{\star}^{3/2}(t',u)}
			\totzeromax{\leq N}^{1/2}(t',u)
			\totzeromax{\leq N-1}^{1/2}(t',u)
	\, dt'
		\\
	& \lesssim
		\varepsilon
		\ln^3(\myexp + t)
		\totzeromax{\leq N}^{1/2}(t,u)
		\totzeromax{\leq N-1}^{1/2}(t,u)
		\int_{t'=0}^t 
			\frac{1}{(1 + t')}
			\frac{1}{\upmu_{\star}^{3/2}(t',u)}
		\, dt'
		\notag	\\
	& \lesssim
		\varepsilon
		\ln^4(\myexp + t)
		\frac{1}{\upmu_{\star}^{1/2}(t',u)}
		\totzeromax{\leq N}^{1/2}(t,u)
		\totzeromax{\leq N-1}^{1/2}(t,u)
		\notag \\
	& \lesssim
		\varepsilon
		\ln^8(\myexp + t)
		\frac{1}{\upmu_{\star}(t',u)}
		\totzeromax{\leq N-1}(t,u)
		+
		\varepsilon
		\totzeromax{\leq N}(t,u),
		\notag
\end{align}
which is clearly $\lesssim$
the right-hand side of \eqref{E:MULTMAINENERGYFLUXERRORINTEGRALESTIMATE} as desired.
We can similarly bound the integrals generated by the second term on the next-to-last line 
of the right-hand side of \eqref{E:DIFFICULTMULTCONSTANTTIMEINTEGRALBOUNDKEYL2ESTIMATE}
and the term on the last line of the right-hand side of \eqref{E:DIFFICULTMULTCONSTANTTIMEINTEGRALBOUNDKEYL2ESTIMATE},
but we use inequality \eqref{E:LESSSINGULARTERMSMUINTEGRALBOUND}
in place of inequality \eqref{E:LOGLOSSKEYMUINTEGRALBOUND}.
Finally, it is straightforward to see that we can bound the integrals generated by the 
two terms on the third-to-last
line of the right-hand side of \eqref{E:DIFFICULTMULTCONSTANTTIMEINTEGRALBOUNDKEYL2ESTIMATE}
by the 
terms 
$
C \varepsilon
				 \int_{t'=0}^t
					\frac{1} 
							 {(1 + t')^{1+\Littleconone} \upmu_{\star}(t',u)} 
				  \totzeromax{\leq N}(t',u)
$
and
$
C \varepsilon
				 \int_{t'=0}^t
					\frac{1} 
							 {(1 + t')^{1+\Littleconone} \upmu_{\star}(t',u)} 
				  \totonemax{\leq N}(t',u)
$
on the right-hand side of \eqref{E:MULTMAINENERGYFLUXERRORINTEGRALESTIMATE}.
We have thus shown that the desired inequality \eqref{E:MULTMAINENERGYFLUXERRORINTEGRALESTIMATE}
follows from \eqref{E:DIFFICULTMULTCONSTANTTIMEINTEGRALBOUNDKEYL2ESTIMATE}.

It remains for us to prove \eqref{E:DIFFICULTMULTCONSTANTTIMEINTEGRALBOUNDKEYL2ESTIMATE}.	
By Prop.~\ref{P:MAINCOMMUTEDWAVEEQNINHOMOGENEOUSTERMPOINTWISEESTIMATES},
it suffices to bound the norm $\| \cdot \|_{L^2(\Sigma_t^u)}$
of the terms on the right-hand side of \eqref{E:TOPORDERTRCHIJUNKENERGYERRORTERMKEYPOINTWISEESTIMATE}
by the right-hand side of \eqref{E:DIFFICULTMULTCONSTANTTIMEINTEGRALBOUNDKEYL2ESTIMATE}.
We proceed by arguing one term at a time. We begin by 
using the estimate $\| \Rad \mathscr{S}^N \Psi \|_{L^2(\Sigma_t^u)} \leq \totzeromax{\leq N}^{1/2}(t,u),$
which follows from inequality \eqref{E:RADPSIL2INTERMSOFZERO},
to deduce that the norm $\| \cdot \|_{L^2(\Sigma_t^u)}$ of the first product on the
right-hand side of \eqref{E:TOPORDERTRCHIJUNKENERGYERRORTERMKEYPOINTWISEESTIMATE}
is
\begin{align} \label{E:KEYFIRSTTERMWHERETHECONSTANTMATTERS}
	= \boxed{2} \left \|
								\frac{\Lunit \upmu}{\upmu}
							\right \|_{C^0(\Sigma_t^u)}
						\| \Rad \mathscr{S}^N \Psi \|_{L^2(\Sigma_t^u)}
		\leq \boxed{2} 
				\frac{\| [\Lunit \upmu]_- \|_{C^0(\Sigma_t^u)}}
						 {\upmu_{\star}(t,u)}
				\totzeromax{\leq N}^{1/2}(t,u)
			+ 2 
				 \left \|
						\frac{[\Lunit \upmu]_+}{\upmu}
				 \right \|_{C^0(\Sigma_t^u)}
				 \totzeromax{\leq N}^{1/2}(t,u),
\end{align}
which in turn is manifestly bounded by 
the right-hand side of \eqref{E:DIFFICULTMULTCONSTANTTIMEINTEGRALBOUNDKEYL2ESTIMATE}.

To bound the norm $\| \cdot \|_{L^2(\Sigma_t^u)}$ of the second product on the
right-hand side of \eqref{E:TOPORDERTRCHIJUNKENERGYERRORTERMKEYPOINTWISEESTIMATE},
we first use Lemma~\ref{L:L2NORMSOFTIMEINTEGRATEDFUNCTIONS}
and inequality \eqref{E:RADPSIL2INTERMSOFZERO}
to bound the norm $\| \cdot \|_{L^2(\Sigma_t^u)}$ of the time integral term as follows:
\begin{align} \label{E:ANNOYINGTIMEINTEGRALL2ESTIMATE}
	& \left \|
		\int_{t'=0}^t 
					\frac{
						\left \|
							[\Lunit \upmu]_-
						\right \|_{C^0(\Sigma_{t'}^u)}}
							{\upmu_{\star}(t',u)}
					\rgeo(t',u)
					\left| \Rad \mathscr{S}^N \Psi \right|(t',u,\vartheta)
				\, dt'
		\right\|_{L^2(\Sigma_t^u)}
			\\
		& \leq
		(1 + C \varepsilon) 
		\rgeo(t,u)
		\int_{t'=0}^t
			\frac{
						\left \|
							[\Lunit \upmu]_-
						\right \|_{C^0(\Sigma_{t'}^u)}}
					 {\upmu_{\star}(t',u)}
			\totzeromax{\leq N}^{1/2}(t',u)
		\, dt'.
		\notag
\end{align}
The factor on the right-hand side of \eqref{E:TOPORDERTRCHIJUNKENERGYERRORTERMKEYPOINTWISEESTIMATE}
that multiplies the first time integral is bounded in the norm $\| \cdot \|_{C^0(\Sigma_t^u)}$ by
\begin{align} \label{E:ANNOYINGTIMEINTEGRALL2ESTIMATEPREFACTORBOUND}
	& \leq 4(1 + C \sqrt{\varepsilon}) 
		\frac{1}{\rgeo(t,u)}
		\frac{
					\left \|
						[\Lunit \upmu]_-
						\right \|_{C^0(\Sigma_t^u)}}
				{\upmu_{\star}(t,u)}
	+ 4 (1 + C \sqrt{\varepsilon}) 
		\frac{1}{\rgeo(t,u)}
		\left\| 
			\frac{[\Lunit \upmu]_+}
				   {\upmu}	
		\right\|_{C^0(\Sigma_t^u)}.
\end{align}
We now multiply \eqref{E:ANNOYINGTIMEINTEGRALL2ESTIMATEPREFACTORBOUND} by
the right-hand side of \eqref{E:ANNOYINGTIMEINTEGRALL2ESTIMATE} and note that the resulting
product is bounded by the right-hand side of \eqref{E:DIFFICULTMULTCONSTANTTIMEINTEGRALBOUNDKEYL2ESTIMATE}
as desired. Note that this argument does not exhaust the full amount of the
constant $\boxed{4.5}$ in front of the first term on the right-hand side of \eqref{E:DIFFICULTMULTCONSTANTTIMEINTEGRALBOUNDKEYL2ESTIMATE},
but rather only $4 (1 + C \sqrt{\varepsilon}).$

To bound the norm $\| \cdot \|_{L^2(\Sigma_t^u)}$ of the third product on the
right-hand side of \eqref{E:TOPORDERTRCHIJUNKENERGYERRORTERMKEYPOINTWISEESTIMATE},
we first use Lemma~\ref{L:L2NORMSOFTIMEINTEGRATEDFUNCTIONS},
inequality \eqref{E:ANNLOYINGSQRTMUOVERMUINTEGRATEDBOUND} 
with the parameter $\Littlecontwo > 0$ chosen to be small enough so that $1 + C \Littlecontwo \leq 1.1$
on the right-hand side of \eqref{E:ANNLOYINGSQRTMUOVERMUINTEGRATEDBOUND},
and inequality \eqref{E:RADPSIL2INTERMSOFZERO}
to bound the norm $\| \cdot \|_{L^2(\Sigma_t^u)}$ of the time integral 
(more precisely, the second time integral on the right-hand side of \eqref{E:TOPORDERTRCHIJUNKENERGYERRORTERMKEYPOINTWISEESTIMATE})
as follows:
\begin{align} \label{E:SECONDANNOYINGTIMEINTEGRALL2ESTIMATE}
	& \left \|
			\int_{t'=0}^t 
							\left\|
								\left(\frac{\upmu(t,\cdot)}{\upmu}\right)^2
							\right\|_{C^0(\Sigma_{t'}^u)}
							\rgeo(t',u)
							\left| \Rad \mathscr{S}^N \Psi \right|(t',u,\vartheta)
			\, dt'
		\right\|_{L^2(\Sigma_t^u)}
			\\
		& \leq 
		(1 + C \varepsilon)
		\rgeo(t,u)
		\int_{t'=0}^t 
			\left\|
				\left(\frac{\upmu(t,\cdot)}{\upmu}\right)^2
			\right\|_{C^0(\Sigma_{t'}^u)}
			\totzeromax{\leq N}^{1/2}(t',u)
			\, dt' 
				\notag \\
		& \leq 
			  (1.1)
			  (1 + C \varepsilon)
				\rgeo^2(t,u)
				\totzeromax{\leq N}^{1/2}(t,u)
			+ C \rgeo^2(t,u)
					\frac{\ln^2(\myexp + t)}{(1 + t)^{\Littlecontwo}}
					\totzeromax{\leq N}^{1/2}(t,u)
			\notag \\
		& \leq 
			  (1.2)
				\rgeo^2(t,u)
				\totzeromax{\leq N}^{1/2}(t,u)
			+ C \rgeo^2(t,u)
					\frac{1}{(1 + t)^{\Littleconone}}
					\totzeromax{\leq N}^{1/2}(t,u),
					\notag
\end{align}
where $0 < \Littleconone < \Littlecontwo$ is a constant.
The factor on the right-hand side of \eqref{E:TOPORDERTRCHIJUNKENERGYERRORTERMKEYPOINTWISEESTIMATE}
that multiplies the second time integral (which is under consideration) is bounded in the 
norm $\| \cdot \|_{C^0(\Sigma_t^u)}$ by
\eqref{E:ANNOYINGTIMEINTEGRALL2ESTIMATEPREFACTORBOUND}, except $\rgeo^{-1}$ is replaced by $\rgeo^{-2}$
and $4 (1 + C \sqrt{\varepsilon})$ is replaced with $2(1 + C \varepsilon).$ Multiplying this factor
by the first term on the right-hand side of \eqref{E:SECONDANNOYINGTIMEINTEGRALL2ESTIMATE}, we
see that the resulting product is bounded by the right-hand side of \eqref{E:DIFFICULTMULTCONSTANTTIMEINTEGRALBOUNDKEYL2ESTIMATE}
as desired. Note that in bounding this term and the first one, we did not exhaust the full amount of the
constant $\boxed{4.5}$ in front of the fourth term on the right-hand side of \eqref{E:DIFFICULTMULTCONSTANTTIMEINTEGRALBOUNDKEYL2ESTIMATE},
but rather only $2 + 2(1.2) + C \varepsilon = 4.4 + C \varepsilon.$
Next, we use the estimate \eqref{E:C0BOUNDLDERIVATIVECRUCICALEIKONALFUNCTIONQUANTITIES}
for $\Lunit \upmu$
to deduce the following easier estimate for the factor on the right-hand side of 
\eqref{E:TOPORDERTRCHIJUNKENERGYERRORTERMKEYPOINTWISEESTIMATE}
that multiplies the second time integral: it is bounded in the norm $\| \cdot \|_{C^0(\Sigma_t^u)}$ by
$\lesssim \varepsilon (1 + t)^{-3} \upmu_{\star}^{-1}(t,u).$ Hence, 
multiplying the factor by the second term on the right-hand side of \eqref{E:SECONDANNOYINGTIMEINTEGRALL2ESTIMATE}, 
we see that the resulting product is bounded by  
$\lesssim \varepsilon (1 + t)^{-(1 + \Littleconone)} \upmu_{\star}^{-1}(t,u) \totzeromax{\leq N}^{1/2}(t,u),$
which is in turn $\lesssim$ the right-hand side of \eqref{E:DIFFICULTMULTCONSTANTTIMEINTEGRALBOUNDKEYL2ESTIMATE} 
as desired.

The remaining terms on the right-hand of \eqref{E:TOPORDERTRCHIJUNKENERGYERRORTERMKEYPOINTWISEESTIMATE}
are relatively easy to bound in the norm $\| \cdot \|_{L^2(\Sigma_t^u)}$ and have only a tiny effect on the dynamics.
Many of the estimates we derive in this paragraph are non-optimal. Throughout this paragraph,
we silently use the coerciveness estimates provided by
Prop.~\ref{P:L2NORMSOFPSIINTERMSOFTHECOERCIVEQUANTITIES}
as well as the fact that
$\totzeromax{\leq N}$
and $\totonemax{\leq N}$
are increasing in their arguments.
To bound the norm $\| \cdot \|_{L^2(\Sigma_t^u)}$ of the product on the right-hand of \eqref{E:TOPORDERTRCHIJUNKENERGYERRORTERMKEYPOINTWISEESTIMATE}
involving the factor $|\chifullmod^{[N]}|(0,u,\vartheta)$ 
by $\lesssim \varepsilon^2 \ln^2(\myexp + t)(1 + t)^{-2}  \upmu_{\star}^{-1}(t,u) $ as desired,
we use inequality \eqref{E:RENORMALZEDTRCHIJUNKDATATERML2NORMATTIMET}.
Next, we see that the terms in the next line of \eqref{E:TOPORDERTRCHIJUNKENERGYERRORTERMKEYPOINTWISEESTIMATE}
(which involve an array of non-time-integrated top-order terms including
$\rgeo \sqrt{\upmu} \left\lbrace \Lunit + \frac{1}{2} \mytr \upchi \right\rbrace \mathscr{S}^{\leq N} \Psi$)
are bounded in the norm $\| \cdot \|_{L^2(\Sigma_t^u)}$
by $\lesssim \varepsilon (1 + t)^{-2} \upmu_{\star}^{-1/2}(t,u)  \totzeromax{\leq N}^{1/2}(t,u)
+ \varepsilon (1 + t)^{-2} \upmu_{\star}^{-1/2}(t,u)  \totonemax{\leq N}^{1/2}(t,u)$
as desired.
Similarly, we bound the terms on the next line of the 
right-hand side of \eqref{E:TOPORDERTRCHIJUNKENERGYERRORTERMKEYPOINTWISEESTIMATE}
(which involve two arrays of non-time-integrated below-top-order terms including $\Rad \mathscr{Z}^{\leq N-1} \Psi$)
in the norm $\| \cdot \|_{L^2(\Sigma_t^u)}$
by $\lesssim \varepsilon (1 + t)^{-1} \upmu_{\star}^{-1}(t,u) \totzeromax{\leq N-1}(t,u)
+ \varepsilon (1 + t)^{-3/2} \upmu_{\star}^{-3/2}(t,u) \totonemax{\leq N-1}(t,u).$
We remark in passing that the $\varepsilon (1 + t)^{-1} \upmu_{\star}^{-1}(t,u)  \totzeromax{\leq N-1}(t,u)$
term leads to some logarithmic-in-time growth 
in the top-order $L^2$ quantities $\totzeromax{\leq N}.$ 
To bound the terms in the next line
of the right-hand of \eqref{E:TOPORDERTRCHIJUNKENERGYERRORTERMKEYPOINTWISEESTIMATE}
(which involve an array of non-time integrated terms including $\mathscr{Z}^{\leq N} (\upmu - 1)$)
in the norm $\| \cdot \|_{L^2(\Sigma_t^u)}$ by
\[
	\lesssim \varepsilon \frac{1}{(1+t)^{3/2}  \upmu_{\star}(t,u)} \totzeromax{\leq N}^{1/2}(t,u)
	+ \varepsilon \frac{1}{(1+t)^{3/2} \upmu_{\star}(t,u)} \totonemax{\leq N}^{1/2}(t,u)
	+ \varepsilon^2 \frac{1}{(1+t)^{3/2} \upmu_{\star}(t,u)},
\]
we use \eqref{E:EIKONALFUNCTIONQUANTITIESL2BOUNDSINTERMSOFQ0ANDQ1}
and the estimate \eqref{E:LOGLOSSLESSSINGULARTERMSMTHREEFOURTHSINTEGRALBOUND}.
To bound the terms in the next line
of the right-hand of \eqref{E:TOPORDERTRCHIJUNKENERGYERRORTERMKEYPOINTWISEESTIMATE}
(which involve an array of time-integrated top-order terms that include the term 
$\rgeo \sqrt{\upmu} \left\lbrace \Lunit + \frac{1}{2} \mytr \upchi \right\rbrace \mathscr{S}^{\leq N} \Psi$)
in the norm $\| \cdot \|_{L^2(\Sigma_t^u)}$ by
the sum of the second and third terms on the right-hand side of 
\eqref{E:DIFFICULTMULTCONSTANTTIMEINTEGRALBOUNDKEYL2ESTIMATE}
(where $a > 0$ is a small constant),
we use Lemma~\ref{L:L2NORMSOFTIMEINTEGRATEDFUNCTIONS}.
To bound the terms in the next two lines
of the right-hand of \eqref{E:TOPORDERTRCHIJUNKENERGYERRORTERMKEYPOINTWISEESTIMATE}
(which involve two arrays of time-integrated below-top-order terms such as $\rgeo \Rad \mathscr{Z}^{\leq N-1} \Psi$)
in the norm $\| \cdot \|_{L^2(\Sigma_t^u)}$ by
\[
	\lesssim \varepsilon \frac{\ln^3(\myexp + t)}{(1 + t) \upmu_{\star}^{3/2}(t,u)} \totzeromax{\leq N-1}(t,u)
	+ \varepsilon \frac{1}{(1 + t)^{3/2}\upmu_{\star}^{3/2}(t,u)} \totonemax{\leq N-1}(t,u),
\]
we use Lemma~\ref{L:L2NORMSOFTIMEINTEGRATEDFUNCTIONS} and
the estimates  
\eqref{E:LOGLOSSKEYMUINTEGRALBOUND}
and
\eqref{E:LOGLOSSMUINVERSEINTEGRALBOUND}.
To bound the terms on the next-to-last line of \eqref{E:TOPORDERTRCHIJUNKENERGYERRORTERMKEYPOINTWISEESTIMATE} 
in the norm $\| \cdot \|_{L^2(\Sigma_t^u)}$ by the right-hand side of \eqref{E:DIFFICULTMULTCONSTANTTIMEINTEGRALBOUNDKEYL2ESTIMATE},
we use the estimate \eqref{E:NONSHARPTOPORDERELLIPTICRECOVERY}
and the simple inequality $|\hat{\angD}^2 \mathscr{S}^{\leq N-1} \upmu| \lesssim |\angD^2 \mathscr{S}^{\leq N-1} \upmu|.$
Finally, using Lemma~\ref{L:L2NORMSOFTIMEINTEGRATEDFUNCTIONS}, 
the estimate \eqref{E:EIKONALFUNCTIONQUANTITIESL2BOUNDSINTERMSOFQ0ANDQ1},
and the estimate \eqref{E:LOGLOSSMUINVERSEINTEGRALBOUND}, 
we see that the terms on the last line of \eqref{E:TOPORDERTRCHIJUNKENERGYERRORTERMKEYPOINTWISEESTIMATE} are
bounded in the norm $\| \cdot \|_{L^2(\Sigma_t^u)}$ by 
$\lesssim$ the sum of the second, third, and last terms on the right-hand side of
\eqref{E:DIFFICULTMULTCONSTANTTIMEINTEGRALBOUNDKEYL2ESTIMATE}.
We have thus proved the desired estimate \eqref{E:DIFFICULTMULTCONSTANTTIMEINTEGRALBOUNDKEYL2ESTIMATE}.

\ \\

\noindent{\emph{Bound for} $\int_{\mathcal{M}_{t,u}}
				(1 + 2 \upmu)
				(\Lunit \mathscr{S}^{N-1} \Rad \Psi)
				(\Rad \Psi) \angLap \mathscr{S}^{N-1} \upmu
				\, d \vol:$}
We now bound the integral corresponding to the easier piece in \eqref{E:DIFFICULTMULTINTEGRANDSPLITTING}.
More precisely, our goal is to bound the spacetime integral
	\begin{align} \label{E:PARTIIDIFFICULTMULTINTEGRALBOUNDPROOF}
		\int_{\mathcal{M}_{t,u}}
				(1 + 2 \upmu)(\Lunit \mathscr{S}^{N-1} \Rad \Psi)
				(\Rad \Psi) 
				\angLap \mathscr{S}^{N-1} \upmu
		\, d \vol
	\end{align}
	in magnitude by the right-hand side of \eqref{E:MULTMAINENERGYFLUXERRORINTEGRALESTIMATE}. 
	To this end, we first bound the integral of the integrand piece 
	$(\Lunit \mathscr{S}^{N-1} \Rad \Psi) (\Rad \Psi) \angLap \mathscr{S}^{N-1} \upmu;$
	the integral of the remaining piece 
	$2 \upmu (\Lunit \mathscr{S}^{N-1} \Rad \Psi) (\Rad \Psi) \angLap \mathscr{S}^{N-1} \upmu$ 
	can be handled similarly. To proceed, we further decompose the first factor as
	\begin{align} \label{E:FIRSTFACTORDECOMP}
		\Lunit \mathscr{S}^{N-1} \Rad \Psi 
		& = 
		\left\lbrace
			\Lunit \mathscr{S}^{N-1} \Rad \Psi
			+ \frac{1}{2} \mytr \upchi \mathscr{S}^{N-1} \Rad \Psi
		\right\rbrace
		- \frac{1}{2} \mytr \upchi \mathscr{S}^{N-1} \Rad \Psi.
	\end{align}
	We now bound the spacetime integral corresponding to the term 
	$\frac{1}{2} \mytr \upchi \mathscr{S}^{N-1} \Rad \Psi$
	in \eqref{E:FIRSTFACTORDECOMP}. Integrating by parts on the $S_{t,u}$
	in order to remove an angular derivative from $\angLap \mathscr{S}^{N-1} \upmu$
	and referring to definition \eqref{E:CHIJUNKDEF},
	we see that it suffices to bound the spacetime integrals
	\begin{align} \label{E:EASYTERMSIBP}
		& \int_{\mathcal{M}_{t,u}}
			\mytr \upchi 
			(\Rad \Psi) 
			(\angdiff \mathscr{S}^{N-1} \Rad \Psi)
			\angdiffuparg{\#} \mathscr{S}^{N-1} \upmu
		\, d \vol
		+ \int_{\mathcal{M}_{t,u}}
			(\angdiff \mytr \upchi^{(Small)})
			(\Rad \Psi) 
			(\mathscr{S}^{N-1} \Rad \Psi)
			\angdiffuparg{\#} \mathscr{S}^{N-1} \upmu
		\, d \vol
			\\
		& \ \ 
			+
			\int_{\mathcal{M}_{t,u}}
				\mytr \upchi 
				(\angdiff \Rad \Psi) 
				(\mathscr{S}^{N-1} \Rad \Psi)
				\angdiffuparg{\#} \mathscr{S}^{N-1} \upmu
			\, d \vol.
			\notag
	\end{align}
	To this end, we first express all spacetime integrals as integrals over $\Sigma_{t'}^u:$
	$\int_{\mathcal{M}_{t,u}} \cdots \, d\vol = \int_{t'=0}^t \int_{\Sigma_{t'}^u} \cdots \, d \tvol \, d t'.$
	None of the $\int_{\Sigma_{t'}^u}$ integrals are difficult to bound because of the favorable time decay that is available.
	To bound the integrals $\int_{\Sigma_{t'}^u} \cdots \, d \tvol,$ 
	we use Cauchy-Schwarz. The quantities
	$\angdiff \mathscr{S}^{N-1} \Rad \Psi,$
	$\mathscr{S}^{N-1} \Rad \Psi,$
	and
	$\angdiffuparg{\#} \mathscr{S}^{N-1} \upmu,$
	are bounded in the norm $\| \cdot \|_{L^2(\Sigma_{t'}^u)}$ as follows
	with the help of 
	Prop.~\ref{P:L2NORMSOFPSIINTERMSOFTHECOERCIVEQUANTITIES},
	inequalities \eqref{E:FUNCTIONPOINTWISEANGDINTERMSOFANGLIEO}
	and
	\eqref{E:EIKONALFUNCTIONQUANTITIESL2BOUNDSINTERMSOFQ0ANDQ1},
	and the estimate $\int_{s=0}^{t'} \frac{\totonemax{\leq N}^{1/2}(s,u)}{(1 + s) \upmu_{\star}^{1/2}(s,u)} \, ds
	\lesssim \ln(\myexp + t') \totonemax{\leq N}^{1/2}(t',u),$
	which follows from \eqref{E:LOGLOSSLESSSINGULARTERMSMTHREEFOURTHSINTEGRALBOUND}:
	\begin{align}
		\| \angdiff \mathscr{S}^{N-1} \Rad \Psi \|_{L^2(\Sigma_{t'}^u)}
		&  \lesssim \frac{1}{1 + t'} \upmu_{\star}^{-1/2}(t',u) \totonemax{\leq N}^{1/2}(t',u),
				\label{E:ANGDIFFRADZNMINUSONEPSIPROOFBOUND} \\
		\| \mathscr{S}^{N-1} \Rad \Psi \|_{L^2(\Sigma_{t'}^u)}
		& \lesssim \totzeromax{\leq N}^{1/2}(t',u),
			\label{E:RADZNMINUSONEPSIPROOFBOUND} \\
		\| \angdiffuparg{\#} \mathscr{S}^{N-1} \upmu \|_{L^2(\Sigma_{t'}^u)}
		& \lesssim \frac{1}{1 + t'} \|  \mathscr{Z}^N \upmu \|_{L^2(\Sigma_{t'}^u)}
			\lesssim 
			\varepsilon
			+ \ln(\myexp + t') \totzeromax{\leq N}^{1/2}(t',u)
			+ \ln(\myexp + t') \totonemax{\leq N}^{1/2}(t',u).
			\label{E:ANGDIFFZNMINUSONEUPMUPROOFBOUND}
	\end{align}	
	The remaining factors in $\int_{\Sigma_{t'}^u} \cdots \, d \tvol$
	are bounded in the norm $\| \cdot \|_{C^0(\Sigma_{t'}^u)}$ as follows:
	$\| \Rad \Psi \|_{C^0(\Sigma_{t'}^u)} \leq \varepsilon (1+t')^{-1},$
	$\| \angdiff \Rad \Psi \|_{C^0(\Sigma_{t'}^u)} \lesssim \varepsilon (1+t')^{-2},$
	$\| \mytr \upchi \|_{C^0(\Sigma_{t'}^u)} \lesssim (1+t')^{-1},$
	and
	$\| \angdiff \mytr \upchi^{(Small)} \|_{C^0(\Sigma_{t'}^u)} \lesssim \varepsilon \ln(\myexp + t') (1+t')^{-3}.$
	These estimates follow from
	the bootstrap assumptions \eqref{E:PSIFUNDAMENTALC0BOUNDBOOTSTRAP},
	\eqref{E:FUNCTIONPOINTWISEANGDINTERMSOFANGLIEO},
	\eqref{E:CRUDELOWERORDERC0BOUNDDERIVATIVESOFANGULARDEFORMATIONTENSORS},
	and
	\eqref{E:C0BOUNDCRUCIALEIKONALFUNCTIONQUANTITIES}.
	Also using simple estimates of the form $ab \lesssim a^2 + b^2,$
	we deduce the following bound for the spacetime integrals
	$\int_{\mathcal{M}_{t,u}} \cdots \, d\vol = \int_{t'=0}^t \int_{\Sigma_{t'}^u} \cdots \, d \tvol \, d t'$
	of interest:
	\begin{align}  \label{E:SPACETIMEINTEGRALBOUNDINVOLVING0ORDERTRCHISMALLTERM}
		& \lesssim
			\varepsilon
			\int_{t'=0}^t
				\frac{1}{(1 + t')^{3/2}}
				\totzeromax{\leq N}(t',u)
			\, dt'	
			+ 
			\varepsilon
			\int_{t'=0}^t
				\frac{1}{(1 + t')^{3/2} \upmu_{\star}^{1/2}(t',u)}
				\totonemax{\leq N}(t',u)
			\, dt'
			\\
		& \ \ 
			+ 
			\varepsilon^3
			\int_{t'=0}^t
				\frac{1}{(1 + t')^{3/2} \upmu_{\star}^{1/2}(t',u)}
			\, dt'.
			\notag
	\end{align}
	Using \eqref{E:LESSSINGULARTERMSMUTHREEFOURTHSINTEGRALBOUND}
	to destroy the factor $\upmu_{\star}^{-1/2}$ in the last integral in 
	\eqref{E:SPACETIMEINTEGRALBOUNDINVOLVING0ORDERTRCHISMALLTERM},
	we conclude that
	\eqref{E:SPACETIMEINTEGRALBOUNDINVOLVING0ORDERTRCHISMALLTERM}
	is $\lesssim$ the right-hand side of \eqref{E:MULTMAINENERGYFLUXERRORINTEGRALESTIMATE}
	as desired.
	
	We now bound the spacetime integral corresponding to the first term in 
	\eqref{E:FIRSTFACTORDECOMP}, that is, the integral
	\begin{align} \label{E:EASYTERMNEEDSIBP}
	&	- 
		\int_{\mathcal{M}_{t,u}}
			(\Rad \Psi) 
			\left\lbrace
				\Lunit \mathscr{S}^{N-1} \Rad \Psi
				+ \frac{1}{2} \mytr \upchi \mathscr{S}^{N-1} \Rad \Psi
			\right\rbrace
			\angLap \mathscr{S}^{N-1} \upmu
		\, d \vol.	
	\end{align}
	To this end, we first integrate by parts with Lemma~\ref{L:TOPORDERMORAWETZREORMALIZEDANGLAPUPMUIBP}, 
	where the weight function $w$ is equal to $1$ and $\angdiffuparg{\#} \mathscr{Z}^{N-1} \upmu$ plays the role of the vectorfield 
	$Y$ from the lemma.
	We therefore have to bound the integrals on the right-hand side of \eqref{E:TOPORDERMORAWETZREORMALIZEDANGLAPUPMUIBP}.
	To this end, we first express all spacetime integrals as integrals over $\Sigma_{t'}^u:$
	$\int_{\mathcal{M}_{t,u}} \cdots \, d\vol = \int_{t'=0}^t \int_{\Sigma_{t'}^u} \cdots \, d \tvol \, d t'.$
	None of the $\int_{\Sigma_{t'}^u}$ integrals are difficult to bound because of the favorable time decay that is available.
	Specifically, when bounding integrals over $\int_{\Sigma_{t'}^u},$ we use Cauchy-Schwarz and bound the 
	terms 
	$\angdiff \mathscr{S}^{N-1} \Rad \Psi,$ 
	$\mathscr{S}^{N-1} \Rad \Psi,$
	$\angdiffuparg{\#} \mathscr{S}^{N-1} \upmu,$
	and $\left\lbrace \angLie_{\Lunit} + \mytr \upchi \right\rbrace \angdiffuparg{\#} \mathscr{S}^{N-1} \upmu$
	in the norm $\| \cdot \|_{L^2(\Sigma_{t'}^u)}.$ 
	The first three of these terms are bounded with the estimates 
	\eqref{E:ANGDIFFRADZNMINUSONEPSIPROOFBOUND}-\eqref{E:ANGDIFFZNMINUSONEUPMUPROOFBOUND}.
	To bound the final term, we first use 
	Lemma~\ref{L:LANDRADCOMMUTEWITHANGDIFF},
	the identity
	$\left\lbrace
			\angLie_{\Lunit}
			+ \mytr \upchi 
	\right\rbrace \ginversesphere
	= - 2 \hat{\upchi}^{(Small)},$
	inequality \eqref{E:FUNCTIONPOINTWISEANGDINTERMSOFANGLIEO}, 
	and the bounds
	$\| \hat{\upchi}^{(Small)} \|_{C^0(\Sigma_{t'}^u)} \lesssim \varepsilon \frac{\ln(\myexp + t')}{(1 + t')^2}$
	and
	$\| [\Lunit, \Rot] \|_{C^0(\Sigma_{t'}^u)} \lesssim \varepsilon \ln(\myexp + t') (1 + t')^{-1},$
	which follow from
	\eqref{E:LCOMMUTETANGENTISTANGENT},
	\eqref{E:LOWERORDERC0BOUNDROTDEFORMSPHEREL},
	\eqref{E:C0BOUNDCRUCIALEIKONALFUNCTIONQUANTITIES},
	and Cor.~\ref{C:SQRTEPSILONREPLCEDWITHCEPSILON},
	to deduce the estimate
	\begin{align}	\label{E:LPLUSTRCHIDIFFERENTITATEDANGDIFFSHARPUPMUPOINTWISE}
		\left|
			\left\lbrace 
				\angLie_{\Lunit} + \mytr \upchi 
			\right\rbrace \angdiffuparg{\#} 
			\mathscr{S}^{N-1} \upmu 
		\right|
		& \lesssim 
			\left|
				\angdiff \Lunit \mathscr{S}^{N-1} \upmu 
			\right|
			+ \left|
					\hat{\upchi}^{(Small)}
				\right| 
				\left|
					\angdiff \mathscr{S}^{N-1} \upmu
				\right|
			\\
		& \lesssim 
			\frac{1}{1 + t'} \sum_{l = 1}^3 
			\left|
				\Lunit \Rot_{(l)} \mathscr{Z}^{N-1} \upmu
			\right| 
			+ \frac{1}{1 + t'} \sum_{l = 1}^3 
				\left|
					[\Lunit, \Rot_{(l)}]
				\right| 
				\left|
					\angdiff \mathscr{Z}^{N-1} \upmu
				\right|
			\notag \\
		& \ \ 
			+ \left|
					\hat{\upchi}^{(Small)}
				\right|
				\left|
					\angdiff \mathscr{Z}^{N-1} \upmu
				\right|
			\notag \\
		& \lesssim 
			\frac{1}{1 + t'} 
			\left|
				\Lunit \mathscr{Z}^{\leq N} \upmu
			\right|
			+ \varepsilon \frac{\ln(\myexp + t')}{(1 + t')^3} 
				\left|
					\mathscr{Z}^{\leq N} (\upmu - 1)
				\right|,
			\notag
	\end{align}
	where all terms above are evaluated at time $t'.$
	We then use inequalities
	\eqref{E:LOGLOSSLESSSINGULARTERMSMTHREEFOURTHSINTEGRALBOUND},
	\eqref{E:EIKONALFUNCTIONQUANTITIESL2BOUNDSINTERMSOFQ0ANDQ1},
	and \eqref{E:LDERIVATIVEEIKONALFUNCTIONQUANTITIESL2BOUNDSINTERMSOFQ0ANDQ1}
	to bound the norm $\| \cdot \|_{L^2(\Sigma_{t'}^u)}$ 
	of the right-hand side of \eqref{E:LPLUSTRCHIDIFFERENTITATEDANGDIFFSHARPUPMUPOINTWISE},
	which leads to the bound
	\begin{align}
		\left\|  \left\lbrace \angLie_{\Lunit} + \mytr \upchi \right\rbrace \angdiffuparg{\#} \mathscr{S}^{N-1} \upmu \right\|_{L^2(\Sigma_{t'}^u)}
		& \lesssim \varepsilon \frac{\ln(\myexp + t')}{(1 + t')^2}
			+  \frac{1}{1+t'} \totzeromax{\leq N}^{1/2}(t',u)
			+  \frac{1}{1+t'} \upmu_{\star}^{-1/2}(t',u) \totonemax{\leq N}^{1/2}(t',u).
	\end{align}	
	All of the remaining integrand factors on the right-hand side of \eqref{E:TOPORDERMORAWETZREORMALIZEDANGLAPUPMUIBP}
	are bounded in the norm $\| \cdot \|_{C^0(\Sigma_{t'}^u)}$
	via the estimates
	$\| \Rad \Psi \|_{C^0(\Sigma_{t'}^u)} \lesssim \varepsilon (1+t')^{-1},$
	$\| \angdiff \Rad \Psi \|_{C^0(\Sigma_{t'}^u)} \lesssim \varepsilon (1+t')^{-2},$
	$\| \mytr \upchi \|_{C^0(\Sigma_{t'}^u)} \lesssim (1+t')^{-1},$
	$\| \angdiff \mytr \upchi^{(Small)} \|_{C^0(\Sigma_{t'}^u)} \lesssim \varepsilon \ln(\myexp + t') (1+t')^{-3},$
	$\| \left\lbrace
				\Lunit 
				+ \frac{1}{2} \mytr \upchi
			 \right\rbrace
			 \Rad \Psi 
	\|_{C^0(\Sigma_{t'}^u)} 
		\lesssim \varepsilon \ln(\myexp + t') (1+t')^{-3},$
	$\| \angdiff
			\left\lbrace
				\Lunit 
				+ \frac{1}{2} \mytr \upchi
			 \right\rbrace
			 \Rad \Psi 
	\|_{C^0(\Sigma_{t'}^u)} 
	\lesssim \varepsilon \ln(\myexp + t') (1+t')^{-4}.$
	All of these $C^0$ estimates except for the last two were justified just above.
	The last two follow from \eqref{E:FUNCTIONPOINTWISEANGDINTERMSOFANGLIEO} 
	and \eqref{E:LPSIPLUSHALFTRCHIPSISTRONGC0BOUND}.
								
In total, this argument leads to the following bound for the \emph{spacetime integrals} on the right-hand side of \eqref{E:TOPORDERMORAWETZREORMALIZEDANGLAPUPMUIBP}:
\begin{align} \label{E:EASYTERMSBOUNDINTOTAL}
	& \lesssim 
	\varepsilon^2
	\int_{t'=0}^t
		\frac{1}{(1+t')^{3/2}} \totzeromax{\leq N}^{1/2}(t',u)
	\, dt'
	+
	\varepsilon
	\int_{t'=0}^t
		\frac{1}{(1+t')^{3/2} \upmu_{\star}^{1/2}} \totzeromax{\leq N}(t',u)
	\, dt'
		\\
	 & \ \  
	 	+ 
		\varepsilon^2
		\int_{t'=0}^t
			\frac{1}{(1+t')^{3/2} \upmu_{\star}^{1/2}(t',u)} \totonemax{\leq N}^{1/2}(t',u)
		\, dt'
		+ 
		\varepsilon
		\int_{t'=0}^t
			\frac{1}{(1+t')^{3/2} \upmu_{\star}(t',u)} \totonemax{\leq N}(t',u)
		\, dt'
		\notag \\
	 & \lesssim 
	\varepsilon
	\int_{t'=0}^t
		\frac{1}{(1+t')^{3/2} \upmu_{\star}^{1/2}} \totzeromax{\leq N}(t',u)
	+	
	\varepsilon
		\int_{t'=0}^t
			\frac{1}{(1+t')^{3/2} \upmu_{\star}(t',u)} \totonemax{\leq N}(t',u)
		\, dt'
	+ \varepsilon^3,
		\notag
\end{align}
where in the last step, we used the simple inequalities
$\varepsilon^2 \totzeromax{\leq N}^{1/2} \lesssim \varepsilon^3 + \varepsilon \totzeromax{\leq N}$
and
$\varepsilon^2 \totonemax{\leq N}^{1/2} \lesssim \varepsilon^3 + \varepsilon \totonemax{\leq N}$
as well as the estimate
$\varepsilon^3 
 \int_{t'=0}^t
		\frac{1}{(1+t')^{3/2} \upmu_{\star}^{1/2}(t',u)}
\, dt'
\lesssim \varepsilon^3,
$
which follows from \eqref{E:LESSSINGULARTERMSMUTHREEFOURTHSINTEGRALBOUND}.
We now note that the right-hand side of \eqref{E:EASYTERMSBOUNDINTOTAL} is manifestly $\lesssim$ the right-hand side of
	\eqref{E:MULTMAINENERGYFLUXERRORINTEGRALESTIMATE} as desired.
	
To bound the two $\Sigma_t^u$ integrals on the right-hand side of \eqref{E:TOPORDERMORAWETZREORMALIZEDANGLAPUPMUIBP}
(that is, the explicitly written one and the one corresponding to the integrand \eqref{E:NONDAMAGINGHYPERSURFACERENORMALZEDANGLAPUPMU}),
we use Cauchy-Schwarz,
\eqref{E:FUNCTIONPOINTWISEANGDINTERMSOFANGLIEO},
\eqref{E:EIKONALFUNCTIONQUANTITIESL2BOUNDSINTERMSOFQ0ANDQ1},
\eqref{E:ANGDIFFRADZNMINUSONEPSIPROOFBOUND}, \eqref{E:RADZNMINUSONEPSIPROOFBOUND}, 
the above $C^0(\Sigma_t^u)$ estimates for $\Rad \Psi$ and $\angdiff \Rad \Psi,$
and simple estimates of the form $ab \lesssim a^2 + b^2$
to bound them (in a non-optimal fashion) by
\begin{align} \label{E:LDERIVATIVETERMSPATIALINTEGRALSBOUND}
	& \lesssim 
	\varepsilon
	\frac{1}{\upmu_{\star}^{1/2}(t,u)}
	\totonemax{\leq N}^{1/2}(t,u)
	\int_{t'=0}^t
		\frac{1}{(1+t')^{3/2}} \totzeromax{\leq N}^{1/2}(t',u)
	\, dt'
		\\
	& 
	+ 
	\varepsilon
	\frac{1}{\upmu_{\star}^{1/2}(t,u)}
	\totonemax{\leq N}^{1/2}(t,u)
	\int_{t'=0}^t
		\frac{1}{(1+t')^{3/2} \upmu_{\star}^{1/2}(t',u)} \totonemax{\leq N}^{1/2}(t',u)
	\, dt'
		\notag \\
	& + \varepsilon \frac{1}{1 + t} \totzeromax{\leq N}(t,u)
	+ \varepsilon
		\frac{1}{1 + t} 
		\totonemax{\leq N}(t,u)
		+ \varepsilon^3 \frac{1}{(1 + t)} \frac{1}{\upmu_{\star}(t,u)},
		\notag
\end{align}
where the time integrals 
in \eqref{E:LDERIVATIVETERMSPATIALINTEGRALSBOUND}
arise from bounding $Y = \angdiffuparg{\#} \mathscr{Z}^{N-1} \upmu$ in
the norm $\| \cdot \|_{L^2(\Sigma_t^u)}$
with \eqref{E:EIKONALFUNCTIONQUANTITIESL2BOUNDSINTERMSOFQ0ANDQ1}.
Clearly, the right-hand side of 
\eqref{E:LDERIVATIVETERMSPATIALINTEGRALSBOUND} is 
$\lesssim$ the right-hand side of \eqref{E:MULTMAINENERGYFLUXERRORINTEGRALESTIMATE}
as desired.

The initial data hypersurface integrals on the right-hand side of \eqref{E:TOPORDERMORAWETZREORMALIZEDANGLAPUPMUIBP}
(corresponding to the integrand \eqref{E:NONDAMAGINGINITIALDATARENORMALZEDANGLAPUPMU})
are cubic and hence it follows easily 
from the estimates of Sect.~\ref{S:INITIALBEHAVIOROFQUANTITIES}
that they are bounded by $\lesssim \varepsilon^3$ as desired.
This completes our proof of the bound for the first integral on the left-hand side of \eqref{E:MULTMAINENERGYFLUXERRORINTEGRALESTIMATE}.

\ \\

\noindent{\emph{Minor changes needed to bound the second integral on the left-hand side of \eqref{E:MULTMAINENERGYFLUXERRORINTEGRALESTIMATE}.}}
The overall strategy for bounding the integral
				$
				\int_{\mathcal{M}_{t,u}}
				\left\lbrace
					(1 + 2 \upmu) \Lunit \mathscr{S}^{N-1} \Rot \Psi
					+ \Rad \mathscr{S}^{N-1} \Rot\Psi
				\right\rbrace
				(\Rad \Psi) \Rot \mathscr{S}^{N-1} \mytr \upchi^{(Small)}
				\, d \vol$ 
				is the same as the one we used for bounding the first integral.
We split the integrand as in \eqref{E:DIFFICULTMULTINTEGRANDSPLITTING}. The difficult part of the proof
is showing the analog of \eqref{E:DIFFICULTMULTCONSTANTTIMEINTEGRALBOUNDKEYL2ESTIMATE}, namely that 
$\| (\Rad \Psi) \Rot \mathscr{S}^{N-1} \mytr \upchi^{(Small)} \|_{L^2(\Sigma_t^u)}$ is bounded by the
right-hand side of \eqref{E:DIFFICULTMULTCONSTANTTIMEINTEGRALBOUNDKEYL2ESTIMATE}. However, this proof is 
exactly the same
as our proof of \eqref{E:DIFFICULTMULTCONSTANTTIMEINTEGRALBOUNDKEYL2ESTIMATE}. 
More precisely, again use Prop.~\ref{P:MAINCOMMUTEDWAVEEQNINHOMOGENEOUSTERMPOINTWISEESTIMATES}
 to reduce the proof of the bound for 
$\| (\Rad \Psi) \Rot \mathscr{S}^{N-1} \mytr \upchi^{(Small)} \|_{L^2(\Sigma_t^u)}$
to bounding the norm $\| \cdot \|_{L^2(\Sigma_t^u)}$ of the right-hand side of inequality 
\eqref{E:TOPORDERTRCHIJUNKENERGYERRORTERMKEYPOINTWISEESTIMATE},
which we already accomplished in the argument following \eqref{E:DIFFICULTMULTCONSTANTTIMEINTEGRALBOUNDKEYL2ESTIMATE}.
The rest of the proof of the bound for $\| (\Rad \Psi) \Rot \mathscr{S}^{N-1} \mytr \upchi^{(Small)} \|_{L^2(\Sigma_t^u)}$ 
is also identical to the proof of the bound for $\| (\Rad \Psi) \angLap \mathscr{S}^{N-1} \upmu \|_{L^2(\Sigma_t^u)}.$

Some minor changes are needed to bound the spacetime integral corresponding to the easier piece in \eqref{E:DIFFICULTMULTINTEGRANDSPLITTING}, namely the integral
$\int_{\mathcal{M}_{t,u}}
				(1 + 2 \upmu)
				(\Lunit \mathscr{S}^{N-1} \Rot \Psi)
				(\Rad \Psi) 
				\Rot \mathscr{S}^{N-1} \mytr \upchi^{(Small)} 
				\, d \vol.$
			We once again use the splitting 
$(1 + 2 \upmu)
	(\Lunit \mathscr{S}^{N-1} \Rot \Psi)
	= \Lunit \mathscr{S}^{N-1} \Rot \Psi + 2 \upmu \Lunit \mathscr{S}^{N-1} \Rot \Psi
$
and then further decompose $\Lunit \mathscr{S}^{N-1} \Rot \Psi$ as in \eqref{E:FIRSTFACTORDECOMP}.
To bound the spacetime integral corresponding to the term $\frac{1}{2} \mytr \upchi \mathscr{S}^{N-1} \Rot \Psi$
from the analog of \eqref{E:FIRSTFACTORDECOMP}, we use an analog of the integration by parts identity
\eqref{E:EASYTERMSIBP} in order to remove the vectorfield $\Rot$ from the factor
$\Rot \mathscr{S}^{N-1} \mytr \upchi^{(Small)}.$ 
This integration by parts also leads to the presence of 
some additional lower-order integrals containing
the factor $\mytr  \angdeform{\Rot}.$ 
This factor in fact enhances the decay of the corresponding integrals over $\Sigma_{t'}^u$
due to the bound $\| \mytr  \angdeform{\Rot} \|_{C^0(\Sigma_{t'}^u)} \lesssim \varepsilon \ln(\myexp + t')(1 + t')^{-1}$
noted in the proof of Lemma~\ref{L:NONDAMAGINGTERMSRENORMALIZEDTRCHI}.

In order to bound the spacetime integral corresponding to the term 
$\left\lbrace
			\Lunit \mathscr{S}^{N-1} \Rot \Psi
			+ \frac{1}{2} \mytr \upchi \mathscr{S}^{N-1} \Rot \Psi
		\right\rbrace
$
from the analog of \eqref{E:FIRSTFACTORDECOMP}, 
we integrate by parts to remove $\Rot$ from the factor
$\Rot \mathscr{S}^{N-1} \mytr \upchi^{(Small)}.$ This is the analog of
the integration by parts performed just after equation \eqref{E:EASYTERMNEEDSIBP} 
to remove an angular derivative off of $\angLap \mathscr{S}^{N-1} \upmu.$
To carry out this integration by parts, in place of Lemma~\ref{L:TOPORDERMORAWETZREORMALIZEDANGLAPUPMUIBP} used above, 
we now use Lemma~\ref{L:TOPORDERMORAWETZREORMALIZEDTRCHIIBP}
with the function $\newsmoothfunction$ from the lemma equal to $\mathscr{S}^{N-1} \mytr \upchi^{(Small)}$
and $\mathscr{S}^{N-1} \Rot \Psi$ in the role of $\mathscr{Z}^N \Psi.$
We then estimate the integrals on the right-hand side of
the integration by parts identity \eqref{E:TOPORDERMORAWETZREORMALIZEDTRCHIIBP}
by using arguments very similar to the ones 
we used in the discussion following equation \eqref{E:EASYTERMNEEDSIBP}.
In total, this line of reasoning allows us to deduce that
the right-hand side of \eqref{E:TOPORDERMORAWETZREORMALIZEDTRCHIIBP} 
is $\leq$ the right-hand side of \eqref{E:EASYTERMSBOUNDINTOTAL} as desired.

$\hfill \qed$

\section{Proof of Lemma~\ref{L:DANGEROUSTOPORDERMORERRORINTEGRAL}}
\label{S:PROOFOFLEMMADANGEROUSTOPORDERMORERRORINTEGRAL}

We now use the auxiliary lemmas to prove the second lemma of primary interest, 
namely Lemma~\ref{L:DANGEROUSTOPORDERMORERRORINTEGRAL}.
We provide complete details for the estimate of
\begin{align} \label{E:DIFFICULTMORTRCHIERRORINTEGRALRESTATED}
- \int_{\mathcal{M}_{t,u}}
							\left\lbrace
								\Lunit \mathscr{S}^{N-1} \Rot \Psi 
								+ \frac{1}{2} \mytr \upchi \mathscr{S}^{N-1} \Rot \Psi
							\right\rbrace
						(\Rad \Psi) 
						(\rgeo^2 \Rot \mathscr{S}^{N-1} \mytr \upchi^{(Small)})
						\,  d \vol.
\end{align}
At the end of the proof, we sketch the minor changes needed to estimate the other integral
on the left-hand side of \eqref{E:MORMAINENERGYFLUXERRORINTEGRALESTIMATE}. 
In order to avoid error integrals that lead to 
damaging top-order estimates,
we need to replace the factor $\rgeo^2 \mathscr{S}^{N-1} \mytr \upchi^{(Small)}$ in \eqref{E:DIFFICULTMORTRCHIERRORINTEGRALRESTATED}
with the partially modified quantity $\rgeo^2 \chipartialmodarg{\mathscr{S}^{N-1}}.$ 
We recall that $\chipartialmodarg{\mathscr{S}^{N-1}}$ is
defined in \eqref{E:TRANSPORTPARTIALRENORMALIZEDTRCHIJUNK}; 
for convenience, we state:
\begin{align} \label{E:PROOFTRCHIRENORMALIZED}
	\rgeo^2 \chipartialmodarg{\mathscr{S}^{N-1}}
	:= \rgeo^2 \mathscr{S}^{N-1} \mytr \upchi^{(Small)} 
			- \frac{1}{2} \rgeo^2 \angGmixedarg{A}{A} \Lunit \mathscr{S}^{N-1} \Psi
			- \frac{1}{2} \rgeo^2 G_{\Lunit \Lunit} \Lunit \mathscr{S}^{N-1} \Psi
			+ \rgeo^2 \angGmixedarg{\Lunit}{A} \angdiffarg{A} \mathscr{S}^{N-1} \Psi.
\end{align}
The reason that we must estimate the quantity
\eqref{E:PROOFTRCHIRENORMALIZED}
rather than $\rgeo^2 \mathscr{S}^{N-1} \mytr \upchi^{(Small)}$
is that we need to use the sharp $L^2$ estimates 
\eqref{E:SIGMAMINUSTRCHIRENORMALIZEDSHARPL2INTERMSOFQANDWIDETILDEQ},
\eqref{E:SIGMAPLUSTRCHIRENORMALIZEDSHARPL2INTERMSOFQANDWIDETILDEQ},
and \eqref{E:TRCHILRENORMALIZEDSHARPL2INTERMSOFQANDWIDETILDEQ} for 
$\rgeo^2 \chipartialmodarg{\mathscr{S}^{N-1}};$
if we did not exploit these sharp $L^2$ estimates, 
then we would encounter error integrals that could grow in time at 
a rate that is just damaging enough to spoil our top-order 
a priori $L^2$ estimates.

Upon making this replacement, we generate an additional error
integral equal to the vectorfield $\Rot$ applied to
difference between 
$\mathscr{S}^{N-1} \mytr \upchi^{(Small)}$
and $\rgeo^2 \chipartialmodarg{\mathscr{S}^{N-1}}.$
Equivalently, in view of definitions 
\eqref{E:TRANSPORTPARTIALRENORMALIZEDTRCHIJUNK}  
and \eqref{E:TRANSPORTPARTIALRENORMALIZEDTRCHIJUNKDISCREPANCY},
we generate the following additional error integrals:
\begin{align} \label{E:TRCHIRENORMALIZEDADDITIONALERRORINTEGRALS}
	 & \int_{\mathcal{M}_{t,u}}
	 	 \rgeo^2
		 \left\lbrace
			 \Lunit \mathscr{S}^{N-1} \Rot \Psi 
			 + \frac{1}{2} \mytr \upchi \mathscr{S}^{N-1} \Rot \Psi
		 \right\rbrace
		 (\Rad \Psi) 
		 \Rot \chipartialmodinhomarg{\mathscr{S}^{N-1}}
		\,  d \vol.
\end{align}
From the estimate \eqref{E:TRCHIJUNKMORAWETZREPLACEMENTTERMSAREHARMLESS}
and Cor.~\ref{C:EASYHARMLESSMORERRORINTEGRAL},
it follows that the magnitude of the integral \eqref{E:TRCHIRENORMALIZEDADDITIONALERRORINTEGRALS} can be bounded 
by the non-boxed-constant-multiplied integrals on the right-hand side of \eqref{E:MORMAINENERGYFLUXERRORINTEGRALESTIMATE}.

Thus, to complete the proof, we have to estimate the main error integral
\begin{align} \label{E:REMAININGDIFFICULTERRORINTEGRAL}
- \int_{\mathcal{M}_{t,u}} 
							\left\lbrace
								\Lunit \mathscr{S}^{N-1} \Rot  \Psi 
								+ \frac{1}{2} \mytr \upchi \mathscr{S}^{N-1} \Rot  \Psi
							\right\rbrace
						(\Rad \Psi) 
						\Rot (\rgeo^2 \chipartialmodarg{\mathscr{S}^{N-1}})
						\, d \vol.
\end{align}
We first integrate by parts with Lemma~\ref{L:TOPORDERMORAWETZREORMALIZEDTRCHIIBP} 
(with $\rgeo^2 \chipartialmodarg{\mathscr{S}^{N-1}}$ in the role of $\newsmoothfunction$ from the lemma)
and therefore have to estimate both the spacetime and the hypersurface integrals on the right-hand side of
\eqref{E:TOPORDERMORAWETZREORMALIZEDTRCHIIBP}. All of these integrals except the 
first two (the most difficult ones) were bounded 
by the non-boxed-constant-multiplied integrals on the right-hand side of \eqref{E:MORMAINENERGYFLUXERRORINTEGRALESTIMATE}
in Lemma~\ref{L:NONDAMAGINGTERMSRENORMALIZEDTRCHI}.
 
We now bound the first difficult integral, which is the 
first spatial integral $\int_{\Sigma_t^u} \cdots$ on the right-hand side of \eqref{E:TOPORDERMORAWETZREORMALIZEDTRCHIIBP}.

\noindent \emph{Estimate of}
$\int_{\Sigma_t^u}
									(\Rot \mathscr{S}^{N-1} \Rot \Psi)
									(\Rad \Psi) 
									\left(
	 									\rgeo^2 \chipartialmodarg{\mathscr{S}^{N-1}}
	 								\right)
\, d \tvol:$
To bound this integral, we first use  
inequality \eqref{E:SQRTUPMUANGDIFFPSIL2INTERMSOFONE},
the estimate \eqref{E:ROTATIONPOINTWISENORMESTIMATE},
and Cor.~\ref{C:SQRTEPSILONREPLCEDWITHCEPSILON}
to deduce that
\begin{align} \label{E:SIMPLEENERGYNORMCOMPARISON}
\left \|
	\sqrt{\upmu} \Rot \mathscr{S}^{N-1} \Rot \Psi
\right\|_{L^2(\Sigma_t^u)}
& \leq (1 + C \varepsilon) 
	\left \|
		\rgeo \sqrt{\upmu} \angdiff \mathscr{S}^{N-1} \Rot \Psi
	\right\|_{L^2(\Sigma_t^u)}
	\leq \sqrt{2} (1 + C \varepsilon) \totonemax{\leq N}^{1/2}(t,u).
\end{align}
Then from Cauchy-Schwarz and \eqref{E:SIMPLEENERGYNORMCOMPARISON},
we deduce that
\begin{align} \label{E:DIFFICULTRREENORMALIZEDTRCHIHYPERSURFACEINTEGRALFIRSTESTIMATE}
	\left|
		\int_{\Sigma_t^u}
			(\Rot \mathscr{S}^{N-1} \Rot \Psi)
			(\Rad \Psi)
			\left(
	 			\rgeo^2 \chipartialmodarg{\mathscr{S}^{N-1}}
	 		\right)
		\, d \tvol
	\right|
& \leq
\sqrt{2}
(1 + C \varepsilon)
\left \|
	 \frac{1}{\sqrt{\upmu}}
	 (\Rad \Psi) 
	 \left(
	 	\rgeo^2 \chipartialmodarg{\mathscr{S}^{N-1}}
	 \right)
\right\|_{L^2(\Sigma_t^u)}
\totonemax{\leq N}^{1/2}(t,u).
\end{align}
Hence, splitting $\Sigma_t^u = \Sigmaminus{t}{t}{u} \cup \Sigmaplus{t}{t}{u}$
(see Def.~\ref{D:REGIONSOFDISTINCTUPMUBEHAVIOR}),
applying Lemma~\ref{L:TRCHIRENORMALIZEDSHARPL2INTERMSOFQANDWIDETILDEQ}, 
and using simple inequalities of the form $ab \lesssim a^2 + b^2,$
we bound the right-hand 
of \eqref{E:DIFFICULTRREENORMALIZEDTRCHIHYPERSURFACEINTEGRALFIRSTESTIMATE} side by
\begin{align} \label{E:DIFFICULTRREENORMALIZEDTRCHIHYPERSURFACEINTEGRALSECONDESTIMATE}
	& \leq 
		(\sqrt{24} + C \varepsilon)
			\totonemax{\leq N}^{1/2}(t,u)
			\| \rgeo \Lunit \upmu \|_{C^0(\Sigmaminus{t}{t}{u})}
			\frac{1}{\upmu_{\star}^{1/2}(t,u)}
			\int_{t'=0}^t
				\frac{1}{\rgeo(t',u) \upmu_{\star}^{1/2}(t',u)} \totonemax{\leq N}^{1/2}(t',u)
			\, dt'
			\\
	& \ \ + (\sqrt{24} + C \varepsilon)
			\totonemax{\leq N}^{1/2}(t,u)
			\left\| \frac{\rgeo \Lunit \upmu}{\upmu} \right\|_{C^0(\Sigmaplus{t}{t}{u})}
			\int_{t'=0}^t
				\left\|
					\sqrt{
					\frac{\upmu(t,\cdot)} 
					     {\upmu}
					     }
				\right\|_{C^0(\Sigmaplus{t'}{t}{u})}
				\frac{1}{\rgeo(t',u)} \totonemax{\leq N}^{1/2}(t',u)
		\, dt'
				\notag \\
	& \ \ + C \varepsilon \frac{1}{\upmu_{\star}^{1/2}(t,u)}
			\totonemax{\leq N}^{1/2}(t,u)	
			\int_{t'=0}^t
				\frac{1}{(1 + t')^{3/2}} \totzeromax{\leq N}^{1/2}(t',u)
			\, dt'
			\notag \\
	& \ \ + C \varepsilon \frac{1}{\upmu_{\star}^{1/2}(t,u)}
				\totonemax{\leq N}^{1/2}(t,u)
				\int_{t'=0}^t
					\frac{1}{(1 + t')^{3/2} \upmu_{\star}^{1/2}(t',u)} \totonemax{\leq N}^{1/2}(t',u)
				\, dt'
				\notag \\
		& \ \ + C \varepsilon \totonemax{\leq N}(t,u)
					\notag \\
		& \ \ + C \varepsilon \frac{1}{\upmu_{\star}(t,u)} \ln^2(\myexp + t) \totzeromax{\leq N-1}(t,u)
					+ C \varepsilon \frac{1}{\upmu_{\star}(t,u)} \ln^2(\myexp + t) \totonemax{\leq N-1}(t,u)
					+ C \varepsilon^3 \frac{1}{\upmu_{\star}(t,u)}.
				\notag
\end{align}
Since $\sqrt{24} < \boxed{5},$ we observe that all terms on the right-hand side of 
\eqref{E:DIFFICULTRREENORMALIZEDTRCHIHYPERSURFACEINTEGRALSECONDESTIMATE} 
are manifestly bounded by the right-hand side of \eqref{E:MORMAINENERGYFLUXERRORINTEGRALESTIMATE}
as desired. 

To complete the proof, it remains only for us to estimate the second difficult integral, which is the
first spacetime integral $\int_{\mathcal{M}_{t,u}} \cdots$ on the right-hand side of \eqref{E:TOPORDERMORAWETZREORMALIZEDTRCHIIBP}
(with $\rgeo^2 \chipartialmodarg{\mathscr{S}^{N-1}}$ in the role of $\newsmoothfunction$ from the lemma).

\noindent \emph{Estimate of}
$- \int_{\mathcal{M}_{t,u}}
								(\Rot \mathscr{S}^{N-1} \Rot \Psi)
								(\Rad \Psi)
								\Lunit
								\left(
	 								\rgeo^2 \chipartialmodarg{\mathscr{S}^{N-1}}
	 							\right)
\, d \vol:$
To bound this spacetime integral, we first express it as a time integral of
integrals over $\Sigma_{t'}^u:$
$\int_{\mathcal{M}_{t,u}} \cdots \, d\vol = 
\int_{t'=0}^t \int_{\Sigma_{t'}^u} \cdots \, d \tvol \, d t'.$
We then use Cauchy-Schwarz and the estimate \eqref{E:SIMPLEENERGYNORMCOMPARISON}
to deduce that
\begin{align} \label{E:DIFFICULTRREENORMALIZEDTRCHISPACETIMEINTEGRALFIRSTESTIMATE}
	& \left|
		\int_{\Sigma_{t'}^u}
			(\Rot \mathscr{S}^{N-1} \Rot \Psi)
			(\Rad \Psi)
			\Lunit
			\left(
	 			\rgeo^2 \chipartialmodarg{\mathscr{S}^{N-1}}
	 		\right)
		\, d \vol
	\right|
		\\
	& \leq \sqrt{2} (1 + C \varepsilon)
				\totonemax{\leq N}^{1/2}(t',u)
				\left\| 
					\frac{1}{\sqrt{\upmu}}
					(\Rad \Psi)
					\Lunit
					\left(
	 					\rgeo^2 \chipartialmodarg{\mathscr{S}^{N-1}}
	 				\right)
				\right \|_{L^2(\Sigma_{t'}^u)}.
				\notag
\end{align}
We now use the key $L^2$ estimate \eqref{E:TRCHILRENORMALIZEDSHARPL2INTERMSOFQANDWIDETILDEQ},
simple estimates of the form $ab \lesssim a^2 + b^2,$
and inequality \eqref{E:POSITIVEPARTOFLMUOVERMUISSMALL}
to bound the last factor on the right-hand side of \eqref{E:DIFFICULTRREENORMALIZEDTRCHISPACETIMEINTEGRALFIRSTESTIMATE},
thereby concluding that the right-hand side of \eqref{E:DIFFICULTRREENORMALIZEDTRCHISPACETIMEINTEGRALFIRSTESTIMATE} 
is bounded by
\begin{align} \label{E:DIFFICULTRREENORMALIZEDTRCHISPACETIMEINTEGRALSECONDESTIMATE}
	& \leq \sqrt{24}(1 + C \varepsilon)
	  	\frac{\left\| [\Lunit \upmu]_- \right\|_{C^0(\Sigma_{t'}^u)}}{\upmu_{\star}(t',u)} 
	  	\totonemax{\leq N}(t',u)
	  	\\
	 & \ \ + 
	  	\sqrt{24}(1 + C \varepsilon)
	  	\frac{1}{\rgeo(t',u) \left\lbrace 1 + \ln \left(\frac{\rgeo(t',u)}{\rgeo(0,u)} \right) \right\rbrace} 
	  	\totonemax{\leq N}(t',u)
	  	\notag \\
	 & \ \ + C \varepsilon
	  					\frac{1}{(1 + t')^{3/2}} 
	  					\totzeromax{\leq N}(t',u)
	  			+ C \varepsilon
	  					\frac{1}{(1 + t')^{3/2}} 
	  					\frac{1}{\upmu_{\star}(t',u)}
	  					\totonemax{\leq N}(t',u)
	  					\notag
	  	\\
	 & \ \ + C \varepsilon 
	  				\frac{1}{(1 + t')} 
						\frac{1}{\upmu_{\star}^{1/2}(t',u)} 
	  				\totzeromax{\leq N-1}^{1/2}(t',u) \totonemax{\leq N}^{1/2}(t',u)
	  			+ C \varepsilon 
	  				\frac{1}{(1 + t')} \frac{1}{\upmu_{\star}(t',u)} 
	  				\totonemax{\leq N-1}^{1/2}(t',u) \totonemax{\leq N}^{1/2}(t',u).
	  				\notag
\end{align}
We now integrate inequality \eqref{E:DIFFICULTRREENORMALIZEDTRCHISPACETIMEINTEGRALSECONDESTIMATE} 
$dt'.$ Noting that $\sqrt{24} < \boxed{5},$ we observe that the time integrals of
the first four products on
the right-hand side of \eqref{E:DIFFICULTRREENORMALIZEDTRCHISPACETIMEINTEGRALSECONDESTIMATE}
are bounded by the right-hand side of \eqref{E:MORMAINENERGYFLUXERRORINTEGRALESTIMATE}
as desired. To bound the time integrals of the terms on the last line of
\eqref{E:DIFFICULTRREENORMALIZEDTRCHISPACETIMEINTEGRALSECONDESTIMATE}, 
we first note that we can bound all factors 
$\totzeromax{\leq N}^{1/2}(t',u),$ 
$\totonemax{\leq N}^{1/2}(t',u),$ 
$\totzeromax{\leq N-1}^{1/2}(t',u),$ 
and
$\totonemax{\leq N-1}^{1/2}(t',u)$
by their values at time $t$ and then pull these factors out of the time integrals
(since these factors are increasing in $t'$).
We then use the estimates 
\eqref{E:LOGLOSSMUINVERSEINTEGRALBOUND}
and \eqref{E:LOGLOSSLESSSINGULARTERMSMTHREEFOURTHSINTEGRALBOUND} to deduce that
the remaining time integrals are bounded as follows:
$\int_{t'=0}^t 
	\frac{1}{(1 + t')}
	\upmu_{\star}^{-1/2}(t',u)
\, dt' \leq C \ln(\myexp + t)
$
and
$\int_{t'=0}^t 
	\frac{1} 
			 {1 + t'}
	\upmu_{\star}^{-1}(t',u)
\, dt' \leq C \ln(\myexp + t) \left\lbrace \ln \upmu_{\star}^{-1}(t,u) + 1 \right\rbrace.$
Also using simple inequalities of the form $ab \lesssim a^2 + b^2,$
we see that in total, the time integrals of the terms on the last line of
\eqref{E:DIFFICULTRREENORMALIZEDTRCHISPACETIMEINTEGRALSECONDESTIMATE}
are bounded by the terms on the first two lines on right-hand side of \eqref{E:MORMAINENERGYFLUXERRORINTEGRALESTIMATE}. 
We have thus bounded the integral $\int_{\mathcal{M}_{t,u}} \cdots \, d\vol$ by the right-hand side of \eqref{E:MORMAINENERGYFLUXERRORINTEGRALESTIMATE} as desired.
This completes our proof of the desired bound for the second integral on the left-hand side of \eqref{E:MORMAINENERGYFLUXERRORINTEGRALESTIMATE}.

\ \\

\noindent{\emph{Minor changes needed to bound}
					$- \int_{\mathcal{M}_{t,u}}
							\rgeo^2
							\left\lbrace
								\Lunit \mathscr{S}^{N-1} \Rad \Psi 
								+ \frac{1}{2} \mytr \upchi \mathscr{S}^{N-1} \Rad \Psi
							\right\rbrace
							(\Rad \Psi) 
							\angLap \mathscr{S}^{N-1} \upmu
				\, d \vol.$}
To bound the first spacetime integral on the left-hand side of \eqref{E:MORMAINENERGYFLUXERRORINTEGRALESTIMATE}, we use the same
overall strategy that we used in bounding the second one. We again need to use a modified quantity in analogy with the quantity 
\eqref{E:PROOFTRCHIRENORMALIZED}. More precisely, we view 
$\angLap \mathscr{S}^{N-1} \upmu = \angdiv \angdiffuparg{\#} \mathscr{S}^{N-1} \upmu$
and replace 
$\angdiffuparg{\#} \mathscr{S}^{N-1} \upmu$
with the partially modified $S_{t,u}-$tangent vectorfield 
$\mupartialmodargsharp{\mathscr{S}^{N-1}},$
which is $\gsphere-$dual to the one-form $\mupartialmodarg{\mathscr{S}^{N-1}}$
defined in \eqref{E:TRANSPORTPARTIALRENORMALIZEDUPMU}.
For convenience, we state the following identity:
\begin{align} \label{E:PROOFANGLAPRENORMALIZED}
	\mupartialmodargsharp{\mathscr{S}^{N-1}}
	& 	=	\angdiffuparg{\#} \mathscr{S}^{N-1} \upmu
					+ \frac{1}{2} \upmu  G_{\Lunit \Lunit} \angdiffuparg{\#} \mathscr{S}^{N-1} \Psi
					+ \upmu G_{\Lunit \Radunit} \angdiffuparg{\#} \mathscr{S}^{N-1} \Psi.
\end{align}
The reason that we must estimate the quantity
\eqref{E:PROOFANGLAPRENORMALIZED}
rather than $\angdiffuparg{\#} \mathscr{S}^{N-1} \upmu$
is that we need to use the sharp $L^2$ estimates 
of Lemma~\ref{L:ANGDIFFUPMURENORMALIZEDSHARPL2INTERMSOFQANDWIDETILDEQ}
and Lemma~\ref{L:LDERIVANGDIFFUPMURENORMALIZEDSHARPL2INTERMSOFQANDWIDETILDEQ}.
This replacement leads to the generation of an additional error integral
that is of the form \eqref{E:TRCHIRENORMALIZEDADDITIONALERRORINTEGRALS}
but with $\angdiv \mupartialmodinhomargsharp{\mathscr{S}^{N-1}}$
(see definition \eqref{E:TRANSPORTPARTIALRENORMALIZEDUPMUDISCREPANCY})
in place of the factor $\Rot \chipartialmodinhomarg{\mathscr{S}^{N-1}}.$
This additional error integral can be suitably bounded 
with the help of the estimate \eqref{E:UPMUMORAWETZREPLACEMENTTERMSAREHARMLESS}
and Cor.~\ref{C:EASYHARMLESSMORERRORINTEGRAL}, in analogy 
with the way we bounded the integral \eqref{E:TRCHIRENORMALIZEDADDITIONALERRORINTEGRALS}.

The difficult part of the analysis is bounding the 
remaining spacetime integral involving $\mupartialmodargsharp{\mathscr{S}^{N-1}},$
which is the analog of the integral \eqref{E:REMAININGDIFFICULTERRORINTEGRAL}.
To estimate this integral, we integrate by parts 
to remove the $\angdiv$ operator from
$\mupartialmodargsharp{\mathscr{S}^{N-1}}.$
Rather than integrating by parts with Lemma~\ref{L:TOPORDERMORAWETZREORMALIZEDTRCHIIBP} as above, this time
we use Lemma~\ref{L:TOPORDERMORAWETZREORMALIZEDANGLAPUPMUIBP}, where
the role of $Y$ from the lemma is played by $\mupartialmodargsharp{\mathscr{S}^{N-1}}$
and the weight function $w$ is equal to $\rgeo^2$
because of the presence of this weight in the integrand.
After integrating by parts, we have to estimate both the spacetime and the hypersurface integrals on the right-hand side of
\eqref{E:TOPORDERMORAWETZREORMALIZEDANGLAPUPMUIBP}. 
All of these integrals except the 
first two (the most difficult ones) were bounded 
by the non-boxed-constant-multiplied integrals on the right-hand side of \eqref{E:MORMAINENERGYFLUXERRORINTEGRALESTIMATE}
in Lemma~\ref{L:NONDAMAGINGTERMSRENORMALIZEDANGLAPUPMU}. This leaves the two difficult error integrals
to estimate: \textbf{a)} the first hypersurface error integral on the right-hand side of \eqref{E:TOPORDERMORAWETZREORMALIZEDANGLAPUPMUIBP},
which we first bound via Cauchy-Schwarz
and \eqref{E:SQRTUPMUANGDIFFPSIL2INTERMSOFONE} as follows: 
\begin{align} \label{E:OTHERPRELIMINARYHYPERSURFACEINTEGRALBOUND}
\left|
\int_{\Sigma_t^u}
	\rgeo^2
	(\Rad \Psi)
	(\angdiff \mathscr{S}^{N-1} \Rad \Psi)
	\mupartialmodargsharp{\mathscr{S}^{N-1}}
\, d \tvol
\right|
& \leq 
\sqrt{2}(1 + C \varepsilon) 
	\totonemax{\leq N}^{1/2}(t,u) 
	\left \| 
		\frac{1}{\sqrt{\upmu}} 
			\rgeo (\Rad \Psi) \mupartialmodarg{\mathscr{S}^{N-1}} 
	\right\|_{L^2(\Sigma_t^u)},
\end{align}
and \textbf{b)} the first spacetime integral on the right-hand side of \eqref{E:TOPORDERMORAWETZREORMALIZEDANGLAPUPMUIBP},
which we first bound via Cauchy-Schwarz and \eqref{E:SQRTUPMUANGDIFFPSIL2INTERMSOFONE} as follows:
\begin{align} \label{E:OTHERPRELIMINARYSPACETIMEINTEGRALBOUND}
&
\left|
\int_{\mathcal{M}_{t,u}}
	\rgeo^2
	(\Rad \Psi)
	(\angdiff \mathscr{S}^{N-1} \Rad \Psi)
	\left\lbrace
		\angLie_{\Lunit} 
			+ \mytr \upchi
	\right\rbrace
	\mupartialmodarg{\mathscr{S}^{N-1}}
\, d \vol
\right|
	\\
& \leq
\sqrt{2}(1 + C \varepsilon)
\int_{t'=0}^t
	\totonemax{\leq N}^{1/2}(t',u)
	\left\|	
		\frac{1}{\sqrt{\upmu}}
		\rgeo 
		(\Rad \Psi)
		\left\lbrace
			\angLie_{\Lunit} 
			+ \mytr \upchi
		\right\rbrace
		\mupartialmodarg{\mathscr{S}^{N-1}}
	\right\|_{L^2(\Sigma_{t'}^u)}
\, dt'.
\notag
\end{align}

To derive a suitable bound for the right-hand side of \eqref{E:OTHERPRELIMINARYHYPERSURFACEINTEGRALBOUND}, 
we argue as in our proof of \eqref{E:DIFFICULTRREENORMALIZEDTRCHIHYPERSURFACEINTEGRALSECONDESTIMATE}, 
using Lemma~\ref{L:ANGDIFFUPMURENORMALIZEDSHARPL2INTERMSOFQANDWIDETILDEQ}
in place of Lemma~\ref{L:TRCHIRENORMALIZEDSHARPL2INTERMSOFQANDWIDETILDEQ}
to estimate the last factor 
$\| \cdots \|_{L^2(\Sigma_t^u)}$
on the right-hand side of \eqref{E:OTHERPRELIMINARYHYPERSURFACEINTEGRALBOUND}.
This line of reasoning allows us to bound 
the right-hand side of \eqref{E:OTHERPRELIMINARYHYPERSURFACEINTEGRALBOUND}
by $\leq$ the right-hand side of \eqref{E:MORMAINENERGYFLUXERRORINTEGRALESTIMATE}
as desired. 

To derive a suitable bound for the integrand on the right-hand side of \eqref{E:OTHERPRELIMINARYSPACETIMEINTEGRALBOUND}, 
we argue as in our proof of \eqref{E:DIFFICULTRREENORMALIZEDTRCHISPACETIMEINTEGRALSECONDESTIMATE}, using 
Lemma~\ref{L:LDERIVANGDIFFUPMURENORMALIZEDSHARPL2INTERMSOFQANDWIDETILDEQ}
in place of Lemma~\ref{L:LDERIVTRCHILRENORMALIZEDSHARPL2INTERMSOFQANDWIDETILDEQ}
to bound the last factor 
$\| \cdots \|_{L^2(\Sigma_{t'}^u)}$
in the integrand on the right-hand side of \eqref{E:OTHERPRELIMINARYSPACETIMEINTEGRALBOUND}.
This line of reasoning allows us to bound  
the right-hand side of \eqref{E:OTHERPRELIMINARYSPACETIMEINTEGRALBOUND}
by $\leq$ the right-hand side of 
\eqref{E:MORMAINENERGYFLUXERRORINTEGRALESTIMATE}
as desired. 

$\hfill \qed$

\section{Proof of Prop.~\ref{P:MAINTOPORDERENERGYANDFLUXINTEGRALINEQUALITIES}}
\label{S:PROOFOFPROPOSITIONMAINTOPORDERENERGYANDFLUXINTEGRALINEQUALITIES}
We now use the previously derived estimates estimates to prove Prop.~\ref{P:MAINTOPORDERENERGYANDFLUXINTEGRALINEQUALITIES}.
Let $0 \leq N \leq 24$ be an integer.
Let $Z \in \mathscr{Z} = \lbrace \rgeo \Lunit, \Rad, \Rot_{(1)}, \Rot_{(2)}, \Rot_{(3)} \rbrace$
and for $N \geq 1,$ let $\mathscr{Z}^N$ be an $N^{th}$ order commutation vectorfield differential operator
of the form $\mathscr{Z}^N = \mathscr{Z}^{N-1} Z.$ We recall that by 
\eqref{E:WAVENTIMESCOMMUTEDBASICSTRUCTURE}-\eqref{E:LOWERORDERINHOMOGENEOUSTERMSFIRSTPOINTWISE}
and Cor.~\ref{C:POINTWISEESTIMATESFOREASYCOMMUTATORTERMS}, we have
\begin{align} 
	\upmu \square_{g(\Psi)} \mathscr{Z}^N \Psi
	& = \inhomarg{\mathscr{Z}^N},
		\\	
	\inhomarg{\mathscr{Z}^N}
	& = \mathscr{Z}^{N-1}(\upmu \D_{\alpha} \Jcurrent{Z}^{\alpha})
			+ Harmless^{\leq N}.
			\label{E:FINALTOPORDERPROOFFNDECOMP}
\end{align}
	We now state and separately consider four cases depending on the structure of $\mathscr{Z}^N.$
	It is easy to see that the four cases exhaust all possibilities.
\begin{center}
	\underline{\large{The four cases for the structure of $\mathscr{Z}^N$}}
\end{center}
		Using
		Lemma~\ref{L:WAVENTIMESCOMMUTEDBASICSTRUCTURE},
		Prop.~\ref{P:IDOFKEYDIFFICULTENREGYERRORTERMS}, 
		Cor.~\ref{C:POINTWISEESTIMATESFOREASYCOMMUTATORTERMS},
		and Cor.~\ref{C:REDUCTIONOFPROOFTOPURESPATIALCOMMUTATORS}, 
		we deduce the following estimates in the four cases.
\begin{enumerate}
	\item $\mathscr{Z}^N$ contains more than one factor of $\rgeo \Lunit.$
		Then 	
		\begin{align*}
			\inhomarg{\mathscr{Z}^N} = Harmless^{\leq N}.
		\end{align*}
	\item $\mathscr{Z}^N$ contains precisely one factor of $\rgeo \Lunit.$
		Then one of the following two possibilities must occur: 
			\begin{align*}	
				\inhomarg{\mathscr{Z}^N} 
				& = Harmless^{\leq N},
					\\
				\inhomarg{\mathscr{Z}^N} 
				& = \rgeo (\angdiff \Psi^{\#}) \cdot (\upmu \angdiff \mathscr{S}^{N-1} \mytr \upchi^{(Small)})
			+ Harmless^{\leq N},
			\end{align*}	
		where $\mathscr{S}^{N-1}$ is an $(N-1)^{st}$ order pure spatial commutation vectorfield operator
		(see definition \eqref{E:DEFSETOFSPATIALCOMMUTATORVECTORFIELDS}).
	\item $\mathscr{Z}^N = \mathscr{S}^{N-1} \Rad,$
		where $\mathscr{S}^{N-1}$ is an $(N-1)^{st}$ order pure spatial commutation vectorfield operator.
		Then 
			\begin{align*}
				\inhomarg{\mathscr{Z}^N} 
				& = (\Rad \Psi) \angLap \mathscr{S}^{N-1} \upmu
					+ (\upmu \angdiffuparg{\#} \Psi) \cdot (\upmu \angdiff \mathscr{S}^{N-1} \mytr \upchi^{(Small)})
					+ Harmless^{\leq N}.
			\end{align*}
	\item $\mathscr{Z}^N = \mathscr{S}^{N-1} \Rot_{(l)},$
		where $\mathscr{S}^{N-1}$ is an $(N-1)^{st}$ order pure spatial commutation vectorfield operator.
		Then 
		\begin{align*}
			\inhomarg{\mathscr{Z}^N} 
			& = (\Rad \Psi) \Rot_{(l)} \mathscr{Z}^{N-1} \mytr \upchi^{(Small)}
				+ \RotRadcomponent{l}(\angdiffuparg{\#} \Psi) \cdot (\upmu \angdiff \mathscr{S}^{N-1} \mytr \upchi)
				+ Harmless^{\leq N}.
		\end{align*}
		
\end{enumerate}

We now give an overview of the proof of the proposition.
We must derive suitable estimates for the error integrals on the right-hand sides of
the energy-flux inequalities \eqref{E:E0DIVID} and \eqref{E:E1DIVID},
where $\mathscr{Z}^N \Psi$ is in the role of $\Psi,$ and
the inhomogeneous term $\waveinhom$ appearing in \eqref{E:E0DIVID} and \eqref{E:E1DIVID} 
is equal to the term $\inhomarg{\mathscr{Z}^N}$ given in equation \eqref{E:FINALTOPORDERPROOFFNDECOMP}.
This difficult analysis has already been carried
out in the previous lemmas.
That analysis allows us to bound
$\enzero[\mathscr{Z}^N \Psi]$
and 
$\flzero[\mathscr{Z}^N \Psi]$
by the right-hand side of
\eqref{E:Q0TOPORDERGRONWALLREADYINEQUALITY},
and
$\enone[\mathscr{Z}^N \Psi],$
$\flone[\mathscr{Z}^N \Psi],$
and 
$\Morint[\mathscr{Z}^N \Psi]$
by the right-hand side of
\eqref{E:Q1TOPORDERGRONWALLREADYINEQUALITY}.
We will then take the max of these estimates
over all such operators of the form $\mathscr{Z}^N,$
$0 \leq N \leq 24,$ and then the sup over $t$ and $u.$ 
This will immediately imply the desired 
top-order inequalities
\eqref{E:Q0TOPORDERGRONWALLREADYINEQUALITY}
and
\eqref{E:Q1TOPORDERGRONWALLREADYINEQUALITY}
for 
$\totzeromax{\leq N},$ 
$\totonemax{\leq N},$
and 
$\totMormax{\leq N}.$

\ \\

\noindent{\emph{Estimate for} $\totzeromax{\leq N}:$}
We now derive the desired inequality \eqref{E:Q0TOPORDERGRONWALLREADYINEQUALITY}.
In our analysis, we implicitly use the quantities defined in
Defs.~\ref{D:TOPORDERMULTBOUNDINGQUANTITIES}
and \ref{D:BELOWORDERMULTBOUNDINGQUANTITIES}.
We now carry out the main step, which is deriving suitable bounds for the quantities
$\enzero[\mathscr{Z}^N \Psi](t,u) + \flzero[\mathscr{Z}^N \Psi](t,u)$
on the left-hand side of \eqref{E:E0DIVID}
(where $\mathscr{Z}^N \Psi$ is in the role of $\Psi$). 
To obtain the desired bounds, we
bound the right-hand side of \eqref{E:E0DIVID} by 
the right-hand side of \eqref{E:Q0TOPORDERGRONWALLREADYINEQUALITY}.
To this end, we first analyze the difficult error integral
\[- \int_{\mathcal{M}_{t,u}}
		\left\lbrace
			(1 + 2 \upmu) (\Lunit \mathscr{Z}^N \Psi) 
			+ 2 \Rad \mathscr{Z}^N \Psi
		\right\rbrace 
		\inhomarg{\mathscr{Z}^N} 
\, d \vol
\]
on the right-hand side of \eqref{E:E0DIVID},
which is absent when $N=0.$
We separately consider the exhaustive cases $(1)-(4)$ 
(depending on the structure of $\mathscr{Z}^N$)
stated above. We bound all error integrals of the form
\[ \int_{\mathcal{M}_{t,u}}
		\left|
			(1 + 2 \upmu) (\Lunit \mathscr{Z}^N \Psi) 
			+ 2 \Rad \mathscr{Z}^N \Psi
		\right|
		Harmless^{\leq N} 
\, d \vol
\]
in magnitude by the right-hand side of \eqref{E:Q0TOPORDERGRONWALLREADYINEQUALITY}
with Cor.~\ref{C:EASYHARMLESSMULTERRORINTEGRAL}.

In case $(2),$ we bound the error integral
\[ \int_{\mathcal{M}_{t,u}}
		\left\lbrace
			(1 + 2 \upmu) (\Lunit \mathscr{Z}^N \Psi) 
			+ 2 \Rad \mathscr{Z}^N \Psi
		\right\rbrace 
		\rgeo (\angdiff \Psi^{\#})
		\cdot 
		(\upmu \angdiff \mathscr{Z}^{N-1} \mytr \upchi^{(Small)}) 
\, d \vol
\]
in magnitude by the right-hand side of \eqref{E:Q0TOPORDERGRONWALLREADYINEQUALITY}
by Lemma~\ref{L:HARMLESSELLIPTICTOPORDERMULTERRORINTEGRAL}.

In case $(3),$ we bound the dangerous error integral
\[ \int_{\mathcal{M}_{t,u}}
		\left\lbrace
			(1 + 2 \upmu) (\Lunit \mathscr{S}^{N-1} \Rad \Psi) 
			+ 2 \Rad \mathscr{S}^{N-1} \Rad \Psi
		\right\rbrace 
		(\Rad \Psi) \angLap \mathscr{S}^{N-1} \upmu 
\, d \vol
\]
in magnitude by the right-hand side of \eqref{E:Q0TOPORDERGRONWALLREADYINEQUALITY}
by Lemma~\ref{L:DANGEROUSTOPORDERMULTERRORINTEGRAL}.
We stress that the error terms on the right-hand side 
of Lemma~\ref{L:DANGEROUSTOPORDERMULTERRORINTEGRAL} 
are the only ones that result
in the difficult capital Roman numeral error terms 
$\topboxedmulterrorone{\leq N},$
$\topmulterrortwo{\leq N},$
$\cdots$ 
on the right-hand side of \eqref{E:Q0TOPORDERGRONWALLREADYINEQUALITY}.
Furthermore, we bound the error integral
\[ \int_{\mathcal{M}_{t,u}}
		\left\lbrace
			(1 + 2 \upmu) (\Lunit \mathscr{S}^{N-1} \Rad \Psi) 
			+ 2 \Rad \mathscr{S}^{N-1} \Rad \Psi
		\right\rbrace 
		(\upmu \angdiffuparg{\#} \Psi) 
		\cdot
		(\upmu \angdiff \mathscr{S}^{N-1} \mytr \upchi^{(Small)}) 
\, d \vol
\]
in magnitude by the right-hand side of \eqref{E:Q0TOPORDERGRONWALLREADYINEQUALITY}
by Lemma~\ref{L:HARMLESSELLIPTICTOPORDERMULTERRORINTEGRAL}.

In case $(4),$ we bound the dangerous error integral
\[ \int_{\mathcal{M}_{t,u}}
		\left\lbrace
			(1 + 2 \upmu) (\Lunit \mathscr{S}^{N-1} \Rot \Psi) 
			+ 2 \Rad \mathscr{S}^{N-1} \Rot \Psi
		\right\rbrace 
		(\Rad \Psi) \Rot_{(l)} \mathscr{S}^{N-1} \mytr \upchi^{(Small)} 
\, d \vol
\]
in magnitude by the right-hand side of \eqref{E:Q0TOPORDERGRONWALLREADYINEQUALITY}
by Lemma~\ref{L:DANGEROUSTOPORDERMULTERRORINTEGRAL}.
Furthermore, we bound the error integral
\[
 \int_{\mathcal{M}_{t,u}}
		\left\lbrace
			(1 + 2 \upmu) (\Lunit \mathscr{S}^{N-1} \Rot \Psi) 
			+ 2 \Rad \mathscr{S}^{N-1} \Rot \Psi
		\right\rbrace 
		\rgeo 
		(\angdiff \Psi^{\#}) 
		\cdot
		(\upmu \angdiff \mathscr{S}^{N-1} \mytr \upchi^{(Small)}) 
\, d \vol
\]
in magnitude by the right-hand side of \eqref{E:Q0TOPORDERGRONWALLREADYINEQUALITY}
by Lemma~\ref{L:HARMLESSELLIPTICTOPORDERMULTERRORINTEGRAL}.

We now bound the error integral on the right-hand side of \eqref{E:E0DIVID} that 
is not related to the inhomogeneous term $\inhomarg{\mathscr{Z}^N},$ that is, the integral
\[
- \frac{1}{2} 
						\int_{\mathcal{M}_{t,u}}
							\upmu \enmomtensor^{\alpha \beta} \deformarg{\Mult}{\alpha}{\beta}
						\, d \vol.
\]
This integral is bounded
by the right-hand side of \eqref{E:Q0TOPORDERGRONWALLREADYINEQUALITY} by \eqref{E:Q0BASICERRORINTEGRALESTIMATE}
(see definition \eqref{E:MULTERRORINT}).

Finally, we note that the terms $\enzero[\mathscr{Z}^N \Psi](0,u)$
corresponding to the first term on the right-hand side of \eqref{E:E0DIVID}
are trivially bounded by $\leq \totzeromax{\leq N}(0,u).$

We have thus bounded the sum $\enzero[\mathscr{Z}^N \Psi](t,u) + \flzero[\mathscr{Z}^N \Psi](t,u)$
by the right-hand side of \eqref{E:Q0TOPORDERGRONWALLREADYINEQUALITY} as desired.
Taking the max of these estimates over all such operators of the form $\mathscr{Z}^N,$
$0 \leq N \leq 24,$ taking the sup over $t$ and $u,$
and recalling Def.~\ref{D:MAINCOERCIVEQUANT},
we conclude the desired inequality \eqref{E:Q0TOPORDERGRONWALLREADYINEQUALITY}.

\ \\

\noindent{\emph{Estimate for} $\totonemax{\leq N}:$}
We now derive the desired inequality \eqref{E:Q1TOPORDERGRONWALLREADYINEQUALITY}.
In our analysis, we implicitly use the quantities defined in
Defs.~\ref{D:TOPORDERMORBOUNDINGQUANTITIES} and \ref{D:BELOWORDERMORBOUNDINGQUANTITIES}.
We now carry out the main step, which is
deriving suitable bounds for the quantities
$\enone[\mathscr{Z}^N \Psi](t,u) + \flone[\mathscr{Z}^N \Psi](t,u)$ on the left-hand side of \eqref{E:E1DIVID}
(where $\mathscr{Z}^N \Psi$ is in the role of $\Psi$). 
To obtain the desired bounds, we
\textbf{a)} show that on the right-hand side of \eqref{E:E1DIVID},
the non-positive integral $-\Morint[\mathscr{Z}^N \Psi](t,u)$ 
(see definition \eqref{E:COERCIVEMORDEF}) is present
and hence we can bring it over to the left-hand side to generate a coercive spacetime integral
and \textbf{b)} bound the remaining integrals on
the right-hand side of \eqref{E:E1DIVID} by 
the right-hand side of \eqref{E:Q1TOPORDERGRONWALLREADYINEQUALITY}.
To accomplish $\textbf{a)},$ we simply appeal to definition \eqref{E:COERCIVEMORDEF}
(with $\mathscr{Z}^N \Psi$ in the role of $\Psi$),
which shows that $\Morint[\mathscr{Z}^N \Psi]$ is equal to 
the first term on the right-hand side of the error term \eqref{E:MORAWETZENERGYERRORINTEGRANDS}.
We emphasize that by \eqref{E:MORERRORINT}, $\Morint[\mathscr{Z}^N \Psi]$ is in fact part of the integral
$- \frac{1}{2} \int_{\mathcal{M}_{t,u}}
								\upmu \enmomtensor^{\alpha \beta}[\mathscr{Z}^N \Psi] 
									\left\lbrace
										\deformarg{\Mor}{\alpha}{\beta} 
										- \rgeo^2 \mytr \upchi g_{\alpha \beta}
									\right\rbrace
								\, d \vol$
and hence it would have appeared on the right-hand side of \eqref{E:Q1TOPORDERGRONWALLREADYINEQUALITY}
if we did not bring it to the left.

We now accomplish $\textbf{b)};$ the proof is similar to the bound we derived for $\totzeromax{\leq N}$ above.
Our goal is to bound
the right-hand side of \eqref{E:E1DIVID} (except for the coercive term $\Morint[\mathscr{Z}^N \Psi]$ that was addressed in \textbf{a)})
by the right-hand side of \eqref{E:Q1TOPORDERGRONWALLREADYINEQUALITY}.
To this end, we first analyze the difficult error integral
\[
	- \int_{\mathcal{M}_{t,u}}
			\rgeo^2
				\left\lbrace
					\Lunit \mathscr{Z}^N \Psi 
					+ \frac{1}{2} \mytr \upchi \mathscr{Z}^N \Psi
				\right\rbrace 
			\inhomarg{\mathscr{Z}^N} 
		\, d \vol
\]
on the right-hand side of \eqref{E:E1DIVID},
which is absent when $N=0.$
We consider the same four cases
(based on the structure of $\mathscr{Z}^N$)
that we did in our estimates for $\totzeromax{\leq N}.$
We bound all error integrals of the form
\[ \int_{\mathcal{M}_{t,u}}
		\rgeo^2
		\left|
			\Lunit \mathscr{Z}^N \Psi
				+ \frac{1}{2} \mytr \upchi \mathscr{Z}^N \Psi 
		\right|
		Harmless^{\leq N} 
\, d \vol
\]
in magnitude by the right-hand side of \eqref{E:Q1TOPORDERGRONWALLREADYINEQUALITY}
with the help of Cor.~\ref{C:EASYHARMLESSMORERRORINTEGRAL}.

In case $(2),$ we bound the error integral
\[ \int_{\mathcal{M}_{t,u}}
		\rgeo^2
		\left\lbrace
			\Lunit \mathscr{Z}^N \Psi 
			+ \frac{1}{2} \mytr \upchi \mathscr{Z}^N \Psi
		\right\rbrace
		\rgeo 
		(\angdiff \Psi^{\#}) 
		\cdot
		(\upmu \angdiff \mathscr{Z}^{N-1} \mytr \upchi^{(Small)}) 
\, d \vol
\]
in magnitude by the right-hand side of \eqref{E:Q1TOPORDERGRONWALLREADYINEQUALITY}
by Lemma~\ref{L:HARMLESSELLIPTICTOPORDERMORERRORINTEGRAL}.

In case $(3),$ we bound the dangerous error integral
\[ \int_{\mathcal{M}_{t,u}}
		\left\lbrace
			\Lunit \mathscr{S}^{N-1} \Rad \Psi 
			+ \frac{1}{2} \mytr \upchi \mathscr{S}^{N-1} \Rad \Psi
		\right\rbrace
		(\Rad \Psi) \angLap \mathscr{S}^{N-1} \upmu 
\, d \vol
\]
in magnitude by the right-hand side of \eqref{E:Q1TOPORDERGRONWALLREADYINEQUALITY}
by Lemma~\ref{L:DANGEROUSTOPORDERMORERRORINTEGRAL}.
We stress that the error terms on the right-hand side of 
Lemma~\ref{L:DANGEROUSTOPORDERMORERRORINTEGRAL}
are the only ones
that result in the difficult capital Roman numeral error terms
$\topboxedmorerrorone{\leq N},$
$\mathbf{\widetilde{II}})_{\leq N},$
$\cdots$
on the right-hand side of \eqref{E:Q1TOPORDERGRONWALLREADYINEQUALITY}.
Furthermore, 
we bound the error integral
\[ \int_{\mathcal{M}_{t,u}}
		\left\lbrace
			\Lunit \mathscr{S}^{N-1} \Rad \Psi 
			+ \frac{1}{2} \mytr \upchi \mathscr{S}^{N-1} \Rad \Psi
		\right\rbrace
		(\upmu \angdiffuparg{\#} \Psi) 
		\cdot
		(\upmu \angdiff \mathscr{S}^{N-1} \mytr \upchi^{(Small)}) 
\, d \vol
\]
in magnitude by the right-hand side of \eqref{E:Q1TOPORDERGRONWALLREADYINEQUALITY}
by Lemma~\ref{L:HARMLESSELLIPTICTOPORDERMORERRORINTEGRAL}.

In case $(4),$ we bound the dangerous error integral
\[ \int_{\mathcal{M}_{t,u}}
		\left\lbrace
			\Lunit \mathscr{S}^{N-1} \Rot \Psi 
			+ \frac{1}{2} \mytr \upchi \mathscr{S}^{N-1} \Rot \Psi
		\right\rbrace 
		(\Rad \Psi) \Rot_{(l)} \mathscr{S}^{N-1} \mytr \upchi^{(Small)} 
\, d \vol
\]
in magnitude by the right-hand side of \eqref{E:Q1TOPORDERGRONWALLREADYINEQUALITY}
by Lemma~\ref{L:DANGEROUSTOPORDERMORERRORINTEGRAL}.
Furthermore, we bound the error integral
\[ \int_{\mathcal{M}_{t,u}}
		\left\lbrace
			\Lunit \mathscr{S}^{N-1} \Rot \Psi 
			+ \frac{1}{2} \mytr \upchi \mathscr{S}^{N-1} \Rot \Psi
		\right\rbrace 
		\rgeo 
		(\angdiff \Psi^{\#}) 
		\cdot
		(\upmu \angdiff \mathscr{S}^{N-1} \mytr \upchi^{(Small)}) 
\, d \vol
\]
in magnitude by the right-hand side of \eqref{E:Q1TOPORDERGRONWALLREADYINEQUALITY}
by Lemma~\ref{L:HARMLESSELLIPTICTOPORDERMORERRORINTEGRAL}.

We now bound the error integrals on the right-hand side of \eqref{E:E1DIVID} that 
are not related to the inhomogeneous term $\inhomarg{\mathscr{Z}^N},$
that is, the first, second, fourth, and fifth (final)
integrals on the right-hand side of \eqref{E:E1DIVID}.
The coercive Morawetz part of the third integral was handled in \textbf{a)},
while the remaining part (see definitions \eqref{E:MORERRORINT} and \eqref{E:MORAWETZENERGYERRORINTEGRANDS})
is bounded by the right-hand side of \eqref{E:Q1TOPORDERGRONWALLREADYINEQUALITY} by 
\eqref{E:Q1BASICERRORINTEGRALESTIMATE}.
The remaining two error integrals are
bounded by
\eqref{E:EASYERRORINTEGRANDTWOINTEGRALESTIMATE}
and
\eqref{E:EASYERRORINTEGRANDONEINTEGRALESTIMATE},
where in bounding the 
hypersurface integral $- \frac{1}{4} \int_{\Sigma_0^u} \cdots,$
we use the fact that by Lemma~\ref{L:CONTROLLINGQUANTITIESAREINTIALLYSMALL}, we have
$\totonemax{\leq N}(0,u) \leq C \totzeromax{\leq N}(0,u).$

Finally, we note that 
by Lemma~\ref{L:CONTROLLINGQUANTITIESAREINTIALLYSMALL},
the terms $\enone[\mathscr{Z}^N \Psi](0,u)$
corresponding to the first term on the right-hand side of \eqref{E:E1DIVID}
are bounded by $\leq \totonemax{\leq N}(0,u) \leq C \totzeromax{\leq N}(0,u).$

We have thus bounded the sum 
$\enone[\mathscr{Z}^N \Psi](t,u) + \flone[\mathscr{Z}^N \Psi](t,u) + \frac{1}{2} \Morint[\mathscr{Z}^N \Psi](t,u)$
by the right-hand side of \eqref{E:Q1TOPORDERGRONWALLREADYINEQUALITY} as desired.
Taking the max of these estimates over all such operators of the form $\mathscr{Z}^N,$
$0 \leq N \leq 24,$
taking the sup over $t$ and $u,$
and recalling Def.~\ref{D:MAINCOERCIVEQUANT},
we conclude the desired inequality \eqref{E:Q1TOPORDERGRONWALLREADYINEQUALITY}.
\hfill $\qed$

\section{Proof of Prop.~\ref{P:MAINBELOWTOPORDERENERGYANDFLUXINTEGRALINEQUALITIES}}
\label{S:PROOFOFPROPOSITIONMAINBELOWTOPORDERENERGYANDFLUXINTEGRALINEQUALITIES}
We now use the previously derived estimates to prove Prop.~\ref{P:MAINBELOWTOPORDERENERGYANDFLUXINTEGRALINEQUALITIES}.
The proof is very similar to the proof of Prop.~\ref{P:MAINTOPORDERENERGYANDFLUXINTEGRALINEQUALITIES}
with only one key change: we estimate all of the error integrals that cause top-order
$L^2$ degeneracy with respect to $\upmu_{\star}^{-1}$
in a \emph{different way}. Specifically, we estimate these integrals in terms of higher-order
$L^2$ quantities (that is, these estimates lose one derivative), 
which will result in the presence of the terms
$\multerrorzero{\leq N+1}$ and $\morerrorzero{\leq N+1}$
(see definitions \eqref{E:MULTLOSSOFONEDERIVERRORINTEGRALBOUND} and \eqref{E:MORLOSSOFONEDERIVERRORINTEGRALTERM})
on the right-hand sides of \eqref{E:Q0BELOWTOPORDERGRONWALLREADYINEQUALITY} and \eqref{E:Q1BELOWTOPORDERGRONWALLREADYINEQUALITY}.
The major gain is the following: because of this alternate strategy,
the top-order error terms $\topboxedmulterrorone{\leq N}-\topmulterroreight{\leq N}$
and $\topboxedmorerrorone{\leq N}-\topmorerroreight{\leq N}$
which are present on the right-hand sides of the top-order estimates
\eqref{E:Q0TOPORDERGRONWALLREADYINEQUALITY}-\eqref{E:Q1TOPORDERGRONWALLREADYINEQUALITY},
are \emph{not present on the right-hand sides of the below-top-order estimates
\eqref{E:Q0BELOWTOPORDERGRONWALLREADYINEQUALITY} and \eqref{E:Q1BELOWTOPORDERGRONWALLREADYINEQUALITY}}.
For the same reason, many terms 
on the first two lines of the 
right-hand sides of \eqref{E:Q0BELOWTOPORDERGRONWALLREADYINEQUALITY} and \eqref{E:Q1BELOWTOPORDERGRONWALLREADYINEQUALITY}
are absent from the below-top-order estimates
\eqref{E:Q0BELOWTOPORDERGRONWALLREADYINEQUALITY} and \eqref{E:Q1BELOWTOPORDERGRONWALLREADYINEQUALITY}.

We now provide the details for the proof of \eqref{E:Q0BELOWTOPORDERGRONWALLREADYINEQUALITY}.
We repeat the proof of \eqref{E:Q0TOPORDERGRONWALLREADYINEQUALITY} and note that 
the only reason that the worst terms $\topboxedmulterrorone{\leq N}-\topmulterroreight{\leq N}$
are present on the right-hand side of \eqref{E:Q0TOPORDERGRONWALLREADYINEQUALITY}
is because we had to use them as an upper bound for the dangerous error integrals
\begin{align} \label{E:MULTERRORINTEGRALALLOWONEDERIVATIVELOSS}
- \int_{\mathcal{M}_{t,u}}
		\left\lbrace
			(1 + 2 \upmu) (\Lunit \mathscr{Z}^N \Psi) 
			+ (\Rad \mathscr{Z}^N \Psi)
		\right\rbrace 
		(\Rad \Psi)
		\myarray
			[\angLap \mathscr{S}^{N-1} \upmu]
			{\Rot_{(l)} \mathscr{S}^{N-1} \mytr \upchi^{(Small)}} 
\, d \vol
\end{align}
from cases $(3)$ and $(4)$ above.
Similarly,
the error integrals
\eqref{E:MULTERRORINTEGRALALLOWONEDERIVATIVELOSS}
are the only reason that the terms on the first two lines
of the right-hand side of \eqref{E:Q0TOPORDERGRONWALLREADYINEQUALITY}
are present, aside from the initial data term 
$\totzeromax{\leq N}(0,u)$
(although the term $C \varepsilon^3 \upmu_{\star}^{-1}(t,u)$ from the first line appears without the factor
of $\upmu_{\star}^{-1}(t,u)$ in many other error integrals).
Since we are no longer bounding a top-order quantity, 
to bound the integral \eqref{E:MULTERRORINTEGRALALLOWONEDERIVATIVELOSS},
rather than using the argument given in the proof of \eqref{E:Q0TOPORDERGRONWALLREADYINEQUALITY}, 
we instead allow a permissible loss of one derivative by
using the bootstrap assumption $\| \Rad \Psi \|_{C^0(\Sigma_t^u)} \leq \varepsilon (1 + t)^{-1}$ 
(that is, \eqref{E:PSIFUNDAMENTALC0BOUNDBOOTSTRAP})
and the derivative-losing Lemma~\ref{L:MULTERRORINTEGRALESTIMATELOSSOFONEDERIVATIVE}
to deduce that \eqref{E:MULTERRORINTEGRALALLOWONEDERIVATIVELOSS} is bounded in magnitude by
\begin{align} 
	& 
	\lesssim
	\varepsilon
	\int_{t'=0}^t
	\frac{1}{(1 + t')^{3/2} \upmu_{\star}^{1/2}(t',u)}
		\totzeromax{\leq N}^{1/2}(t',u)	
		\int_{s=0}^{t'}
			\frac{1}{1 + s}
			\totzeromax{\leq N+1}^{1/2}(s,u)
		\, ds
\, dt'	
	\\
& \ \ 
+ \varepsilon
 	\int_{t'=0}^t
	\frac{1}{(1 + t')^{3/2} \upmu_{\star}^{1/2}(t',u)}
		\totzeromax{\leq N}^{1/2}(t',u)	
		\int_{s=0}^{t'}
			\frac{1}{(1 + s)\upmu_{\star}^{1/2}(s,u)}
			\totonemax{\leq N+1}^{1/2}(s,u)
		\, ds
\, dt'
	\notag \\
& \ \
+ \varepsilon
\int_{t'=0}^t
	\frac{1}{(1 + t')^{3/2} \upmu_{\star}^{1/2}(t',u)}
		\totzeromax{\leq N}(t',u)	
\, dt'
+ \varepsilon^3
	\notag 
\end{align}
as desired (see definition \eqref{E:MULTLOSSOFONEDERIVERRORINTEGRALBOUND}).
We have thus proved inequality \eqref{E:Q0BELOWTOPORDERGRONWALLREADYINEQUALITY}.

The proof of \eqref{E:Q1BELOWTOPORDERGRONWALLREADYINEQUALITY} is similar. 
More precisely, we repeat the proof of \eqref{E:Q1TOPORDERGRONWALLREADYINEQUALITY} and make the following key
change, which completely eliminates the dangerous terms 
$\topboxedmorerrorone{\leq N}-\topmorerroreight{\leq N}$ that are present on the right-hand side
of the top-order estimate \eqref{E:Q1TOPORDERGRONWALLREADYINEQUALITY}
as well as the terms on the first two lines
of the right-hand side of \eqref{E:Q1TOPORDERGRONWALLREADYINEQUALITY}
(aside from the initial data term $C \totzeromax{\leq N}(0,u)$
and the term $C \varepsilon^3 \upmu_{\star}^{-1}(t,u),$ which appears in
the ameliorated form $C \varepsilon^3$):
we bound the integrals
\begin{align} \label{E:MORERRORINTEGRALALLOWONEDERIVATIVELOSS}
- \int_{\mathcal{M}_{t,u}}
		\rgeo^2
		\left\lbrace
			\Lunit \mathscr{Z}^N \Psi 
			+ \frac{1}{2} \mytr \upchi \mathscr{Z}^N \Psi
		\right\rbrace 
		(\Rad \Psi)
		\myarray
			[\angLap \mathscr{S}^{N-1} \upmu]
			{\Rot_{(l)} \mathscr{S}^{N-1} \mytr \upchi^{(Small)}} 
\, d \vol
\end{align}
by using the bootstrap assumption $\| \Rad \Psi \|_{C^0(\Sigma_t^u)} \lesssim \frac{\varepsilon}{1 + t}$
and Lemma~\ref{L:MORERRORINTEGRALESTIMATELOSSOFONEDERIVATIVE}
to deduce that \eqref{E:MORERRORINTEGRALALLOWONEDERIVATIVELOSS} is bounded in magnitude by
\begin{align} 
	& 
	\lesssim
	\varepsilon
	\int_{t'=0}^t
		\frac{1}{(1 + t')^2}
		\left(
			\int_{s=0}^{t'}
				\frac{ \totzeromax{\leq N+1}^{1/2}(s,u)}{1 + s} 
			\, ds
		\right)^2
	\, dt'
		+ 
		\varepsilon
		\int_{t'=0}^t
			\frac{1}{(1 + t')^2}
				\left(
					\int_{s=0}^{t'}
						\frac{\totonemax{\leq N+1}^{1/2}(s,u)}{(1 + s)\upmu_{\star}^{1/2}(s,u)} 
					\, ds
				\right)^2
		\, dt'
			\\
		& \ \ 
		+ 
		\varepsilon
		\int_{u'=0}^u
			\totonemax{\leq N}(t,u') 
		\, du'
		+	
		\varepsilon^3
		\notag
\end{align}
as desired (see definition \eqref{E:MORLOSSOFONEDERIVERRORINTEGRALTERM}).
We have thus proved inequality \eqref{E:Q1BELOWTOPORDERGRONWALLREADYINEQUALITY}.
\hfill $\qed$

\section{Proof of Lemma~\ref{L:FUNDAMENTALGRONWALL}} 
\label{S:PROOFOFLEMMAFUNDAMENTALGRONWALL}

Before proving Lemma~\ref{L:FUNDAMENTALGRONWALL}, we first provide the following simple lemma, which 
will be used to estimate the right-hand side of \eqref{E:STRANGEINTEGRATINGFACTORINTEGRALBOUND}.

\begin{lemma} [\textbf{Integral estimate for an unusual Gronwall factor}] \label{L:UNUSUALINTEGRATINGFACTORESTIMATE} 
Let $t \geq 0.$ There exists a constant $C > 0$ \textbf{independent of} $t$ such that
\begin{align} \label{E:UNUSUALINTEGRATINGFACTORESTIMATE}
	\int_{t'=0}^t 
		\frac{\ln^5(\myexp + t')}{(1 + t')^2 \sqrt{\ln(\myexp + t) - \ln(\myexp + t')}}
	\, dt'
	& \leq C.
\end{align}

\end{lemma}

\begin{proof}[\textbf{Proof of Lemma} \ref{L:UNUSUALINTEGRATINGFACTORESTIMATE}]
We make the change of variables
$\uptau' := \ln(\myexp + t'),$ 
$\uptau := \ln(\myexp + t),$
$d \uptau' = \frac{dt'}{\myexp + t'}.$
We then split the integration domain into the
intervals $[1, \uptau/2]$ and $[\uptau/2,\uptau]$
and bound each piece as follows:
\begin{align}
& \int_{t'=0}^t 
		\frac{\ln^5(\myexp + t')}{(1 + t')^2 \sqrt{\ln(\myexp + t) - \ln(\myexp + t')}}
	\, dt'
	\\
& \leq 
	C
	\int_{t'=0}^t 
		\frac{\ln^5(\myexp + t')}{(\myexp + t')^2 \sqrt{\ln(\myexp + t) - \ln(\myexp + t')}}
	\, dt'
		\notag \\
& \leq C 
	\int_{\uptau'= 1}^{\uptau/2} 
	\frac{\exp(-\uptau') \uptau'^5}{\sqrt{\uptau - \uptau'}}
\, d \uptau'
+ \int_{\uptau'=\uptau/2}^{\uptau} 
	\frac{\exp(-\uptau') \uptau'^5}{\sqrt{\uptau - \uptau'}}
\, d \uptau'
	\notag \\
& \leq 
	\frac{C}{\sqrt{\uptau}}
	\int_{\uptau'= 1}^{\infty} 
		\exp(-\uptau') \uptau'^5
\, d \uptau'
	+ C
	\exp(-\uptau/2) \uptau^5
	\int_{\uptau'=\uptau/2}^{\uptau} 
		\frac{1}{\sqrt{\uptau - \uptau'}}
	\, d \uptau'
	\notag 
		\\
& \leq \frac{C}{\sqrt{\uptau}} 
	+ C \exp(-\uptau/2) \uptau^{11/2}
	\leq C.
	\notag
\end{align}

\end{proof}

We now prove Lemma~\ref{L:FUNDAMENTALGRONWALL}.

\ \\

\noindent{\emph{Proof of the estimates for the top-order quantities.}}
Throughout the proof, we often use the 
estimate $\rgeo(t,u) \approx 1 + t$ 
(on $\mathcal{M}_{\Tboot,U_0}$)
and the fact that 
$\totzeromax{\leq N}(t,u)$
and $\totonemax{\leq N}(t,u)$
are increasing in their arguments
without explicitly mentioning them each time.
We first prove the desired (difficult) estimates for the top-order quantities
$\totzeromax{\leq 24}(t,u),$
$\totonemax{\leq 24}(t,u),$
and
$\totMormax{\leq 24}(t,u).$ Our argument involves the following functions:
\begin{align}
	\iota_1(t,u)
	& := \exp(u), 
		\\
	\iota_2(t,u)
	& := \exp\left(\int_{s=0}^{t} \frac{1}{(1 + s)^{1 + \Littleconone}} \, ds \right),
		\\
	\iota_3(t,u)
	& := 1 + \ln \left(\frac{\rgeo(t,u)}{\rgeo(0,u)} \right),
		\\
	\iota(t,u)
	& :=
		\iota_1^{2 \Conthree}(t,u)
		\iota_2^{2 \Conthree}(t,u)
		\iota_3^{2 \Cononestar}(t,u)
		\upmu_{\star}^{-17.5}(t,u),
		 \label{E:TOPORDERQ0INTEGRATINGFACTORPRODUCT} \\
	\widetilde{\iota}(t,u)
	& := 
		\iota_1^{2 \Conthree}(t,u)
		\iota_2^{2 \Conthree}(t,u)
		\iota_3^{2 \Cononestar + 4}(t,u)
		\upmu_{\star}^{-17.5}(t,u),
		\label{E:TOPORDERQ1INTEGRATINGFACTORPRODUCT} 
\end{align}
where $\Conthree > 0$ and $\Cononestar > 4$ are constants to be determined below, and 
$0 < \Littleconone < 1/2$ is the small positive constant on the right-hand side of 
\eqref{E:TOPORDERMULTTERMSMALLA},
\eqref{E:TOPORDERQ0LARGETIMEDECAYCOUNTERSUPMULOSSTERM},
and \eqref{E:TOPORDERQ0INVOLVINGQ1LARGETIMEDECAYCOUNTERSUPMULOSSTERM}.
Note the important discrepancy factor $\iota_3^4(t,u)$
between $\iota(t,u)$ and $\widetilde{\iota}(t,u).$
The functions $\iota^{-1}(t,u)$ and $\widetilde{\iota}^{-1}(t,u)$ can be thought of as
products of approximate integrating factors for the inequalities of Prop.~\ref{P:MAINTOPORDERENERGYANDFLUXINTEGRALINEQUALITIES}
and hence 
$\iota(t,u)$ and $\widetilde{\iota}(t,u)$ 
are connected to the expected behavior of
$\totzeromax{\leq 24}(t,u),$
$\totonemax{\leq 24}(t,u),$
and
$\totMormax{\leq 24}(t,u).$
Note that 
$\iota_1(t,u), \iota_2(t,u), \iota_3(t,u), \upmu_{\star}^{-1}(t,u), \iota(t,u), \widetilde{\iota}(t,u)$ 
are \emph{increasing} in $u,$
that $\iota_1(t,u), \iota_2(t,u), \iota_3(t,u)$ are also increasing in $t,$
and that $\upmu_{\star}^{-1}(t,u), \iota(t,u)$ and $\widetilde{\iota}(t,u)$
are \emph{approximately increasing in} $t$
in the sense that
\begin{align} \label{E:INTEGRATINGFACTORAPPROXIMATELYINCREASING}
	\upmu_{\star}^{-1}(t_1,u)
	& \leq (1 + C \sqrt{\varepsilon}) \upmu_{\star}^{-1}(t_2,u), & & \mbox{if \ } t_1 \leq t_2,
		\\
	\iota(t_1,u) 
	& \leq (1 + C \sqrt{\varepsilon}) \iota_1(t_2,u), & & \mbox{if \ } t_1 \leq t_2,
		\label{E:SECONDINTEGRATINGFACTORAPPROXIMATELYINCREASING} \\
	\widetilde{\iota}(t_1,u) 
	& \leq (1 + C \sqrt{\varepsilon}) \widetilde{\iota}(t_2,u), & & \mbox{if \ } t_1 \leq t_2.
	\label{E:THIRDINTEGRATINGFACTORAPPROXIMATELYINCREASING}
\end{align} 
To deduce \eqref{E:INTEGRATINGFACTORAPPROXIMATELYINCREASING}, 
we have used the approximate monotonicity inequality \eqref{E:MUSTARINVERSETOPOWERAMUSTGROWUPTOACONSTANT}. 
Note also that for a fixed $\Conthree,$ 
$\iota_1^{\Conthree}(t,u)$ and $\iota_2^{\Conthree}(t,u)$ are uniformly bounded from above by a positive constant
for $(t,u) \in [0,\Tboot) \times [0,U_0].$
Throughout this proof, we often use these increasing/approximately increasing properties without explicitly mentioning
them every time.

In order to generate sufficient smallness that will allow us to absorb error integrals
during this Gronwall-type argument, we use the crucially important estimates of Prop.~\ref{P:MUINVERSEINTEGRALESTIMATES}
together with the following simple inequalities:
\begin{align}
	\int_{u'=0}^u
		\iota_1^{\Conthree}(t,u')	
	\, du'
	& \leq \frac{1}{\Conthree} \iota_1^{\Conthree}(t,u),
		\label{E:MAINGRONWALLPROOFKEYEASYUINTEGRALBOUND} \\
	\int_{t'=0}^t
		\frac{1}{(1 + t')^{1 + \Littleconone}} 
		\iota_2^{\Conthree}(t',u)	
	\, dt'
	& \leq \frac{1}{\Conthree} \iota_2^{\Conthree}(t,u),
		\label{E:MAINGRONWALLPROOFKEYLARGETDECAYINTEGRALBOUND} \\
	\int_{t'=0}^t
		\frac{1}{\rgeo(t',u) \left\lbrace 1 + \ln \left(\frac{\rgeo(t',u)}{\rgeo(0,u)} \right) \right\rbrace}
		\iota_3^{\Cononestar}(t',u) 	
	\, dt'
	& \leq \frac{1}{\Cononestar} \iota_3^{\Cononestar}(t,u).
		\label{E:MAINGRONWALLPROOFKEYLOGGROWTHGENERATINGINTEGRALBOUND} 
\end{align}
The above estimates are straightforward to verify by explicit computation. The ``smallness factors''
that we exploit are the factors $\frac{1}{\Conthree}$ and $\frac{1}{\Cononestar}$
as well as the bootstrap parameter $\varepsilon.$

To proceed, we define the following rescaled functions $q(t,u)$ and $\widetilde{q}(t,u),$
which we would like to show are small on the domain of interest:
\begin{align}
	q(t,u) 
	& : = 
		\sup_{(\hat{t},\hat{u}) \in [0,t] \times [0,u]}
			\iota^{-1}(\hat{t},\hat{u})
			\totzeromax{\leq 24}(\hat{t},\hat{u}),
		\\
\widetilde{q}(t,u) 
&: = 	\sup_{(\hat{t},\hat{u}) \in [0,t] \times [0,u]}
				\widetilde{\iota}^{-1}(\hat{t},\hat{u})
				\left\lbrace 
					\totonemax{\leq 24}(\hat{t},\hat{u}) 
					+ \Morint(\hat{t},\hat{u})
				\right\rbrace.
\end{align}
In order to prove \eqref{E:Q0TOPORDERKEYBOUND} and \eqref{E:Q1TOPORDERKEYBOUND}, it suffices to prove that
the following estimates hold for $(t,u) \in [0,\Tboot) \times [0,U_0]:$
\begin{align}
	q(t,u) \label{E:KEYTOPORDERLITTLEFBOUND}
	& \leq C 
				\left\lbrace
					\mathring{\upepsilon}^2
					+ \varepsilon^3
				\right\rbrace,
				\\
	\widetilde{q}(t,u) 
	& \leq C 
				\left\lbrace
					\mathring{\upepsilon}^2
					+ \varepsilon^3
				\right\rbrace.
				\label{E:KEYTOPORDERWIDETILDELITTLEFBOUND}
\end{align}
The reason that the bounds \eqref{E:KEYTOPORDERLITTLEFBOUND} and \eqref{E:KEYTOPORDERWIDETILDELITTLEFBOUND} are sufficient 
for proving \eqref{E:Q0TOPORDERKEYBOUND} and \eqref{E:Q1TOPORDERKEYBOUND}
is that $\iota_1(t,u)$ and $\iota_2(t,u)$
are each bounded from above by a uniform constant $C$ on the domain $(t,u) \in [0,\infty) \times [0,U_0],$
while $\iota_3(t,u) \approx \ln(\myexp + t).$

Our first goal is to derive a suitable bound for $\widetilde{q}(t,u)$ in terms of 
$q(t,u).$ More precisely, we will show that if $\Conthree > 1$ and $\Cononestar > 4$ are large enough, 
then the following bound holds for $(t,u) \in [0,\Tboot) \times [0,U_0]:$
\begin{align} \label{E:DESIREDDIFFICULTQ1BOUNDINTERMSOFQ0}
	\widetilde{q}(t,u)
	& \leq 
				C 
				\left\lbrace
					\mathring{\upepsilon}^2
					+ \varepsilon^3
					+ q(t,u)
				\right\rbrace.
\end{align}
To this end, we will prove that
\begin{align} \label{E:PREDESIREDDIFFICULTQ1BOUNDINTERMSOFQ0}
	\widetilde{q}(t,u) & \leq 
				C 
				\left\lbrace
					\mathring{\upepsilon}^2
					+ \varepsilon^3
					+ q(t,u)
				\right\rbrace
			 + \upalpha 
			   \widetilde{q}(t,u),
\end{align}
where $0 < \upalpha < 1$ is a constant. The desired estimate 
\eqref{E:DESIREDDIFFICULTQ1BOUNDINTERMSOFQ0} easily follows from
\eqref{E:PREDESIREDDIFFICULTQ1BOUNDINTERMSOFQ0} by absorbing the 
product $\upalpha \widetilde{q}(t,u)$ 
on the right-hand side of \eqref{E:PREDESIREDDIFFICULTQ1BOUNDINTERMSOFQ0} back into the left.

In order to prove \eqref{E:PREDESIREDDIFFICULTQ1BOUNDINTERMSOFQ0}, 
we evaluate both sides of \eqref{E:Q1TOPORDERGRONWALLREADYINEQUALITY}
at $(\hat{t},\hat{u}),$
multiply both sides of \eqref{E:Q1TOPORDERGRONWALLREADYINEQUALITY}
by $\widetilde{\iota}^{-1}(\hat{t},\hat{u}),$
and then take $\sup_{(\hat{t},\hat{u}) \in [0,t] \times [0,u]}.$ 
The left-hand side of the resulting inequality is 
precisely the term $\widetilde{q}(t,u)$ on the left-hand side of \eqref{E:PREDESIREDDIFFICULTQ1BOUNDINTERMSOFQ0}.
The most difficult terms on the right-hand side 
of the resulting inequality are the ones involving ``boxed constants terms,'' 
that is,
terms arising from 
$\topboxedmorerrorone{\leq 24}$ 
and 
$\topboxedmorerrorfive{\leq 24}$ 
(see Def.~\ref{D:TOPORDERMORBOUNDINGQUANTITIES}).
We now derive a suitable bound for these difficult terms.
We give complete details for bounding the term corresponding to $\topboxedmorerrorone{\leq 24}.$
We then provide abbreviated proofs for the remaining terms.
To handle the term corresponding to $\topboxedmorerrorone{\leq 24}$ 
(see definition \eqref{E:FIRSTHARDMORTERM}), we have to bound
\begin{align} \label{E:FIRSTHARDMORTERMGRONWALL}
	\sup_{(\hat{t},\hat{u}) \in [0,t] \times [0,u]}
	\left\lbrace
		\widetilde{\iota}^{-1}(\hat{t},\hat{u})	
		\topboxedmorerrorone{\leq 24}(\hat{t},\hat{u})
	\right\rbrace.
\end{align}
We now multiply and divide by $\upmu_{\star}^{17.5}(t',\hat{u})$ 
in the integral on the right-hand side of the bound \eqref{E:FIRSTHARDMORTERM} for
$\topboxedmorerrorone{\leq 24}(\hat{t},\hat{u})$
and use the fact that $\iota_1,$ $\iota_2,$ $\iota_3,$ and $\upmu_{\star}^{-1}$ 
are increasing in both of their arguments,
thereby bounding the terms in braces in \eqref{E:FIRSTHARDMORTERMGRONWALL} by
\begin{align} \label{E:ABSORBABLEBOUNDFIRSTHARDMORTERMGRONWALL}
	& \leq 
		\boxed{5}
		\widetilde{\iota}^{-1}(\hat{t},\hat{u})
		\int_{t'=0}^{\hat{t}}
			\frac{\| [\Lunit \upmu]_- \|_{C^0(\Sigma_{t'}^{\hat{u}})}} 
				  {\upmu_{\star}(t',\hat{u})} 
			\totonemax{\leq 24}(t',\hat{u})
		\, dt'
		\\
	& \leq 
		\boxed{5}
				\widetilde{\iota}^{-1}(\hat{t},\hat{u})
				\left\lbrace
					\sup_{(t',u') \in [0,\hat{t}] \times [0,\hat{u}]}
					\upmu_{\star}^{17.5}(t',u')
					\totonemax{\leq 24}(t',u')
				\right\rbrace
				\int_{t'=0}^{\hat{t}}
					\frac{\| [\Lunit \upmu]_- \|_{C^0(\Sigma_{t'}^{\hat{u}})}} 
							 {\upmu_{\star}(t',\hat{u})} 
					\upmu_{\star}^{-17.5}(t',\hat{u})
				\, dt'
			\notag \\
	& \leq 
		\boxed{5}
				\widetilde{\iota}^{-1}(\hat{t},\hat{u})
				\iota_1^{2 \Conthree}(\hat{t},\hat{u})
				\iota_2^{2 \Conthree}(\hat{t},\hat{u})
				\iota_3^{2 \Cononestar}(\hat{t},\hat{u})
				\widetilde{q}(\hat{t},\hat{u})
				\int_{t'=0}^{\hat{t}}
					\frac{\| [\Lunit \upmu]_- \|_{C^0(\Sigma_{t'}^{\hat{u}})}} 
							 {\upmu_{\star}(t',\hat{u})} 
					\upmu_{\star}^{-17.5}(t',\hat{u})
				\, dt'
			\notag \\
	& =
		\boxed{5}
		\widetilde{q}(\hat{t},\hat{u})
		\upmu_{\star}^{17.5}(\hat{t},\hat{u})
		 	\int_{t'=0}^{\hat{t}}
					\frac{\| [\Lunit \upmu]_- \|_{C^0(\Sigma_{t'}^{\hat{u}})}} 
							 {\upmu_{\star}(t',\hat{u})} 
							\upmu_{\star}^{-17.5}(t',\hat{u})
				\, dt'
				\notag \\
	& \leq (1 + C \sqrt{\varepsilon})
		\boxed{5} \times \frac{1}{17.5}
		\widetilde{q}(t,u).
		\notag
\end{align}
In the last step of \eqref{E:ABSORBABLEBOUNDFIRSTHARDMORTERMGRONWALL},
we used the crucially important integral estimate \eqref{E:KEYMUTOAPOWERINTEGRALBOUND}
and the fact that $\widetilde{q}$ is increasing in both of its arguments.
We have therefore bounded \eqref{E:FIRSTHARDMORTERMGRONWALL}
by the right-hand side of \eqref{E:ABSORBABLEBOUNDFIRSTHARDMORTERMGRONWALL}.

A similar argument based on the crucially important estimate \eqref{E:KEYMUDECAYHYPERSURFACEMUTOAPOWERINTEGRALBOUND}
leads to the following bound (see definition \eqref{E:SECONDHARDMORTERM}):
\begin{align} \label{E:ABSORBABLEBOUNDDECONDHARDMORTERMGRONWALL}
	\sup_{(\hat{t},\hat{u}) \in [0,t] \times [0,u]}
	\left\lbrace
	 	\widetilde{\iota}^{-1}(\hat{t},\hat{u})	
		\topboxedmorerrorfive{\leq 24}(\hat{t},\hat{u})
	\right\rbrace
	& \leq 
		(1 + C \sqrt{\varepsilon})
		\boxed{5} \times \frac{1}{8.25}
		\widetilde{q}(t,u).
\end{align}
The important point is that the constants $\frac{5}{17.5}$ and $\frac{5}{8.25}$ from 
the estimates \eqref{E:ABSORBABLEBOUNDFIRSTHARDMORTERMGRONWALL} and \eqref{E:ABSORBABLEBOUNDDECONDHARDMORTERMGRONWALL} 
verify $\frac{5}{17.5} + \frac{5}{8.25} < .9 < 1.$ This sum makes up the bulk of the constant
$\upalpha$ on the right-hand side of \eqref{E:PREDESIREDDIFFICULTQ1BOUNDINTERMSOFQ0}.
The estimates \eqref{E:ABSORBABLEBOUNDFIRSTHARDMORTERMGRONWALL} and \eqref{E:ABSORBABLEBOUNDDECONDHARDMORTERMGRONWALL}
are the only two that 
force us to prove bounds for 
$\totonemax{\leq 24}(t,u)$ that involve large degeneracy with respect to powers of $\upmu_{\star}^{-1}.$ 
As we will see, the remaining terms on the right-hand side of \eqref{E:Q1TOPORDERGRONWALLREADYINEQUALITY}
can be suitably bounded by choosing $\Conthree$ and $\Cononestar$ to be sufficiently large.

To bound the term corresponding to $\topmorerrortwo{\leq 24}$
(see definition \eqref{E:TOPORDERMORTERCREATESLOGGROWTH}),
we use a similar argument based on the estimate
\eqref{E:MAINGRONWALLPROOFKEYLOGGROWTHGENERATINGINTEGRALBOUND} 
and the increasing/approximately increasing properties of the integrating factors
to deduce that
\begin{align}  \label{E:TOPORDERQ1LOGGROWTHTERMBOUNDED}
	\sup_{(\hat{t},\hat{u}) \in [0,t] \times [0,u]}
	\left\lbrace
		\widetilde{\iota}^{-1}(\hat{t},\hat{u})
		\topmorerrortwo{\leq 24}(\hat{t},\hat{u})
	\right\rbrace
	& \leq 
		\frac{C}{2 \Cononestar + 4}
		\widetilde{q}(t,u).
\end{align}
By choosing $\Cononestar$ to be large in \eqref{E:TOPORDERQ1LOGGROWTHTERMBOUNDED}, we can ensure 
that $\frac{C}{2 \Cononestar + 4}$ is as small as we need it to be.

Next, using a similar argument based on the estimate 
\eqref{E:LESSSINGULARTERMSMUINTEGRALBOUND}
(with $\Conone = 0$ and $\Littleconone = 1/2$),
we deduce that
(see definitions 
\eqref{E:TOPORDERMORDANGEROUSUPMUBUTLOTSOFTIMEDECAYINTERMSOFQ0},
\eqref{E:TOPORDERMORDANGEROUSUPMUBUTLOTSOFTIMEDECAYINTERMSOFQ1},
\eqref{E:TOPORDERMORHYPERSURFACEDANGEROUSUPMUBUTLOTSOFTIMEDECAYINTERMSOFQ0},
and
\eqref{E:TOPORDERMORHYPERSURFACEDANGEROUSUPMUBUTLOTSOFTIMEDECAYINTERMSOFQ1})
\begin{align}
\sup_{(\hat{t},\hat{u}) \in [0,t] \times [0,u]}
	\left\lbrace
		\widetilde{\iota}^{-1}(\hat{t},\hat{u})
	\topmorerrorthree{\leq 24}(\hat{t},\hat{u})
	\right\rbrace
	& \leq C \varepsilon q,
			\label{E:TOPORDERMORDANGEROUSUPMUBUTLOTSOFTIMEDECAYINTERMSOFQ0ESTIMATED} \\
		\sup_{(\hat{t},\hat{u}) \in [0,t] \times [0,u]}
	\left\lbrace
		\widetilde{\iota}^{-1}(\hat{t},\hat{u})
	\topmorerrorfour{\leq 24}(\hat{t},\hat{u})
	\right\rbrace
	& \leq C \varepsilon \widetilde{q}, 
			\label{E:TOPORDERMORDANGEROUSUPMUBUTLOTSOFTIMEDECAYINTERMSOFQ1ESTIMATED} \\
		\sup_{(\hat{t},\hat{u}) \in [0,t] \times [0,u]}
	\left\lbrace
		\widetilde{\iota}^{-1}(\hat{t},\hat{u})
	\topmorerrorseven{\leq 24}(\hat{t},\hat{u})
	\right\rbrace
	& \leq 
		C \varepsilon q^{1/2} \widetilde{q}^{1/2}
		\leq 
			C \varepsilon q
			+ C \varepsilon\widetilde{q},
			\label{E:TOPORDERMORHYPERSURFACEDANGEROUSUPMUBUTLOTSOFTIMEDECAYINTERMSOFQ0ESTIMATED} \\
	\sup_{(\hat{t},\hat{u}) \in [0,t] \times [0,u]}
	\left\lbrace
		\widetilde{\iota}^{-1}(\hat{t},\hat{u})
	\topmorerroreight{\leq 24}(\hat{t},\hat{u})
	\right\rbrace
	& \leq 
		C \varepsilon \widetilde{q}.
		\label{E:TOPORDERMORHYPERSURFACEDANGEROUSUPMUBUTLOTSOFTIMEDECAYINTERMSOFQ1ESTIMATED}
\end{align}

To bound the term corresponding to $\topmorerrorsix{\leq 24}$
(see definition \eqref{E:ANNOYINGTOPORDERMORTERCREATESLOGGROWTH}),
we use a similar argument based on the estimate \eqref{E:KEYMUNODECAYHYPERSURFACEMUTOAPOWERINTEGRALBOUND}
to deduce that
\begin{align} 
	\sup_{(\hat{t},\hat{u}) \in [0,t] \times [0,u]}
	\left\lbrace
		\widetilde{\iota}^{-1}(\hat{t},\hat{u})
		\topmorerrorsix{\leq 24}(\hat{t},\hat{u})
	\right\rbrace
	& \leq 
		\frac{C}{\Cononestar + 1/2}
		\widetilde{q}(t,u).
\end{align}
By choosing $\Cononestar$ to be large in the previous estimate, we can ensure 
that $\frac{C}{\Cononestar + 1/2}$ is as small as we need it to be.
We have thus accounted for all of the capital roman numeral terms on
the right-hand side of \eqref{E:Q1TOPORDERGRONWALLREADYINEQUALITY}.

To bound the terms corresponding to the easy terms 
$\morerrorone{\leq 24}-\morerrorfive{\leq 24}$
(see Def.~\ref{D:BELOWORDERMORBOUNDINGQUANTITIES})
on the right-hand side of \eqref{E:Q1TOPORDERGRONWALLREADYINEQUALITY}, 
we use similar arguments to deduce that
\begin{align} 
	(1 + \widetilde{\varsigma}^{-1})
	\sup_{(\hat{t},\hat{u}) \in [0,t] \times [0,u]}
	\left\lbrace
	 	\widetilde{\iota}^{-1}(\hat{t},\hat{u})	
		\morerrorone{\leq 24}(\hat{t},\hat{u})
	\right\rbrace
	& \leq C (1 + \widetilde{\varsigma}^{-1}) q(t,u),
		\label{E:EASYBELOWTOPORDERQ1THATINVOLVESQ0TERMESTIMATE} \\
	\sup_{(\hat{t},\hat{u}) \in [0,t] \times [0,u]}
	\left\lbrace
	 	\widetilde{\iota}^{-1}(\hat{t},\hat{u})
		\morerrortwo{\leq 24}(\hat{t},\hat{u})
	\right\rbrace
	& \leq \frac{C}{2\Cononestar + 4} \widetilde{q}(t,u),
		\label{E:EASYBELOWTOPORDERQ1LOGLOSSTERMESTIMATE} \\
	(1 + \widetilde{\varsigma}^{-1})
	\sup_{(\hat{t},\hat{u}) \in [0,t] \times [0,u]}
		\left\lbrace
			\widetilde{\iota}^{-1}(\hat{t},\hat{u})
			\morerrorthree{\leq 24}(\hat{t},\hat{u})
		\right\rbrace
	& \leq \frac{C}{\Conthree} (1 + \widetilde{\varsigma}^{-1}) \widetilde{q}(t,u),
		\label{E:EASYBELOWTOPORDERLARGETIMEDECAYTERMESTIMATE} \\
	(1 + \widetilde{\varsigma}^{-1})
	\sup_{(\hat{t},\hat{u}) \in [0,t] \times [0,u]}
	\left\lbrace
		\widetilde{\iota}^{-1}(\hat{t},\hat{u})
		\morerrorfour{\leq 24}(\hat{t},\hat{u})
	\right\rbrace
	& \leq \frac{C}{\Conthree} (1 + \widetilde{\varsigma}^{-1}) \widetilde{q}(t,u),
		\label{E:EASYBELOWTOPORDERQ1UTERMESTIMATE}
		\\
	(\widetilde{\varsigma} + \varepsilon)
	\sup_{(\hat{t},\hat{u}) \in [0,t] \times [0,u]}
	\left\lbrace
		\widetilde{\iota}^{-1}(\hat{t},\hat{u})
		\morerrorfive{\leq 24}(\hat{t},\hat{u})
	\right\rbrace
	& \leq C (\widetilde{\varsigma} + \varepsilon) \widetilde{q}(t,u).
		\label{E:EASYBELOWTOPORDERQ1MORAWETZTERMESTIMATE}
\end{align}
We make the following clarifying remarks concerning the above estimates.
The estimate \eqref{E:EASYBELOWTOPORDERQ1THATINVOLVESQ0TERMESTIMATE}
relies on the inequality $\widetilde{\iota}^{-1}(t,u) \ln^4(\myexp + t) \leq C \iota^{-1}(t,u).$
The estimate \eqref{E:EASYBELOWTOPORDERQ1LOGLOSSTERMESTIMATE}
was already proved in \eqref{E:TOPORDERQ1LOGGROWTHTERMBOUNDED}.
The estimate \eqref{E:EASYBELOWTOPORDERLARGETIMEDECAYTERMESTIMATE} relies on inequality
\eqref{E:MAINGRONWALLPROOFKEYLARGETDECAYINTEGRALBOUND}
in analogy with the way that
\eqref{E:KEYMUTOAPOWERINTEGRALBOUND} 
was used to deduce
\eqref{E:ABSORBABLEBOUNDFIRSTHARDMORTERMGRONWALL},
while the estimate \eqref{E:EASYBELOWTOPORDERQ1UTERMESTIMATE} relies on inequality
\eqref{E:MAINGRONWALLPROOFKEYEASYUINTEGRALBOUND}.
The estimate involving \eqref{E:EASYBELOWTOPORDERQ1MORAWETZTERMESTIMATE} follows easily from the definitions.

In order to complete the proof of \eqref{E:PREDESIREDDIFFICULTQ1BOUNDINTERMSOFQ0},
it remains for us to bound the terms arising from the terms on the first two lines of the right-hand side of 
\eqref{E:Q1TOPORDERGRONWALLREADYINEQUALITY}.
We begin by bounding the terms generated by the terms on the first line.
By Lemma~\ref{L:CONTROLLINGQUANTITIESAREINTIALLYSMALL},
the term corresponding to $C \totzeromax{\leq 24}(0,u)$ is $\leq C \mathring{\upepsilon}^2.$
The term corresponding to $C \varepsilon^3 \upmu_{\star}^{-1}(t,u)$ can easily be bounded by $\leq C \varepsilon^3.$
The term corresponding to $C \varepsilon \totzeromax{\leq N}(t,u)$
can easily be bounded by $\leq C \varepsilon q(t,u),$
while the term corresponding to $C \varepsilon \totonemax{\leq N}(t,u)$
can easily bounded by $\leq C \varepsilon \widetilde{q}(t,u).$ 
To bound the term arising from the first term on the second line of the right-hand side of \eqref{E:Q1TOPORDERGRONWALLREADYINEQUALITY},
we use the bootstrap assumption \eqref{E:Q0MIDBOOT} for $\totzeromax{\leq 23}(t,u)$ 
to deduce that
\begin{align} \label{E:TERMSBOUNDEDSWITHBOOTSTRAPASSUMPTIONS}
	\sup_{(\hat{t},\hat{u}) \in [0,t] \times [0,u]}
	& 
	\left\lbrace
		C \varepsilon
		\widetilde{\iota}^{-1}(\hat{t},\hat{u})
		\upmu_{\star}^{-1}(\hat{t},\hat{u}) 
		\ln^2(\myexp + \hat{t})
		\totzeromax{\leq 23}(\hat{t},\hat{u})
	\right\rbrace
		\\
	& \leq C \varepsilon^3 \sup_{(\hat{t},\hat{u}) \in [0,t] \times [0,u]} \upmu_{\star}(\hat{t},\hat{u})
		\leq C \varepsilon^3.
		\notag 
\end{align}
Using a similar argument based on the bootstrap assumption \eqref{E:Q1MIDBOOT} for $\totonemax{\leq 23}(t,u),$
we also bound the term arising from the second term on the second line of the right-hand side of \eqref{E:Q1TOPORDERGRONWALLREADYINEQUALITY}
by $\leq C \varepsilon^3.$

Combining all of the above estimates, we deduce that
\begin{align} \label{E:MOREEXPLICITPREDESIREDDIFFICULTQ1BOUNDINTERMSOFQ0}
	\widetilde{q}(t,u) & \leq 
				C 
				\left\lbrace
					\mathring{\upepsilon}^2
					+ \varepsilon^3
					+ \widetilde{\varsigma}^{-1}
				\right\rbrace
				q(t,u)
					\\
			& \ \ 
					+ \left\lbrace
						(1 + C \sqrt{\epsilon})
						\frac{5}{17.5} 
						+ (1 + C \sqrt{\epsilon})
							\frac{5}{8.25}
						+ \frac{C}{\Cononestar}
						+ \frac{C}{\Conthree} (1 + \widetilde{\varsigma}^{-1})
						+ C \widetilde{\varsigma}
						+ C \varepsilon
					\right\rbrace
					\widetilde{q}(t,u).
				\notag
\end{align}
The bound \eqref{E:PREDESIREDDIFFICULTQ1BOUNDINTERMSOFQ0} now follows from \eqref{E:MOREEXPLICITPREDESIREDDIFFICULTQ1BOUNDINTERMSOFQ0}
if we first choose $\varepsilon$ and $\widetilde{\varsigma}$ to be sufficiently small
and we then choose $\Cononestar$ and $\Conthree$ to be sufficiently large.
We remark that later in the proof, just below inequality 
\eqref{E:MOREEXPLICITPREDESIREDDIFFICULTQ0BOUND}, 
we may need to further enlarge $\Cononestar$ and $\Conthree.$

We now use the estimate \eqref{E:DESIREDDIFFICULTQ1BOUNDINTERMSOFQ0} 
to help us derive the desired estimate \eqref{E:KEYTOPORDERLITTLEFBOUND} for $q(t,u).$
The desired estimate \eqref{E:KEYTOPORDERWIDETILDELITTLEFBOUND} for $\widetilde{q}$
then follows easily from \eqref{E:DESIREDDIFFICULTQ1BOUNDINTERMSOFQ0}
and \eqref{E:KEYTOPORDERLITTLEFBOUND}. 
To prove \eqref{E:KEYTOPORDERLITTLEFBOUND},
we will argue as in our proof of \eqref{E:PREDESIREDDIFFICULTQ1BOUNDINTERMSOFQ0}
to show that there exists a constant $0 < \upbeta < 1$ such that
\begin{align} \label{E:PREDESIREDDIFFICULTQ0BOUND}
	q(t,u) & \leq 
				C 
				\left\lbrace
					\mathring{\upepsilon}^2
					+ \varepsilon^3
				\right\rbrace
				+ \upbeta q(t,u).
\end{align}
Clearly, the desired bound \eqref{E:KEYTOPORDERLITTLEFBOUND} follows once we have shown \eqref{E:PREDESIREDDIFFICULTQ0BOUND}.
Our proof of \eqref{E:PREDESIREDDIFFICULTQ0BOUND} is very similar to our proof of \eqref{E:PREDESIREDDIFFICULTQ1BOUNDINTERMSOFQ0}.
More precisely, in order to prove \eqref{E:PREDESIREDDIFFICULTQ0BOUND}, 
for each term on the right-hand side of \eqref{E:Q0TOPORDERGRONWALLREADYINEQUALITY} involving an integral, we find an
effective approximate integrating factor and then multiply both sides of \eqref{E:Q0TOPORDERGRONWALLREADYINEQUALITY}
by the product of all the integrating factors. Specifically, we evaluate both sides of \eqref{E:Q0TOPORDERGRONWALLREADYINEQUALITY}
at $(\hat{t},\hat{u}),$
multiply both sides of \eqref{E:Q0TOPORDERGRONWALLREADYINEQUALITY}
by $\iota^{-1}(\hat{t},\hat{u})$
(see \eqref{E:TOPORDERQ0INTEGRATINGFACTORPRODUCT})
and then take $\sup_{(\hat{t},\hat{u}) \in [0,t] \times [0,u]}.$ 
The left-hand side of the resulting inequality is 
precisely the term $q(t,u)$ on the left-hand side of \eqref{E:PREDESIREDDIFFICULTQ0BOUND}.
The most difficult terms on the right-hand side 
of the resulting inequality are the ones involving ``boxed constants terms''
that is, terms arising from 
$\topboxedmulterrorone{\leq 24}$ 
and 
$\topboxedmulterrorfive{\leq 24}$
(see Def.~\ref{D:TOPORDERMULTBOUNDINGQUANTITIES}).
We now derive a suitable bound for these difficult terms.
To handle the term corresponding to $\topboxedmulterrorone{\leq 24}$
(see definition \eqref{E:FIRSTHARDMULTTERM}),
we have to bound
\begin{align} \label{E:FIRSTHARDMULTTERMGRONWALL}
	\sup_{(\hat{t},\hat{u}) \in [0,t] \times [0,u]}
	\left\lbrace
		\iota(\hat{t},\hat{u})	
		\topboxedmulterrorone{\leq 24}(\hat{t},\hat{u})
	\right\rbrace.
\end{align}
Using an argument similar to the one we used to prove inequality \eqref{E:ABSORBABLEBOUNDFIRSTHARDMORTERMGRONWALL},
and in particular using the crucially important integral estimate \eqref{E:KEYMUTOAPOWERINTEGRALBOUND} \emph{twice}
since the right-hand side of \eqref{E:FIRSTHARDMULTTERM} involves two time integrations,
we deduce that the right-hand side of \eqref{E:FIRSTHARDMULTTERMGRONWALL} is
\begin{align}  \label{E:FIRSTHARDMULTTERMGRONWALLCARRIEDOUT}
	\leq (1 + C \sqrt{\varepsilon}) 
				\boxed{9} \times \frac{1}{8.75} \times \frac{1}{8.75} q(t,u).
\end{align}

A similar argument based on the same crucially important estimate \eqref{E:KEYMUTOAPOWERINTEGRALBOUND}
leads to the following bound for the term involving $\topboxedmulterrorfive{\leq 24}$ (see definition \eqref{E:SECONDHARDMULTTERM}):
\begin{align} \label{E:SECONDHARDMULTTERMGRONWALL}
	\sup_{(\hat{t},\hat{u}) \in [0,t] \times [0,u]}
	\left\lbrace
		\iota(\hat{t},\hat{u})	\topboxedmulterrorfive{\leq 24}(\hat{t},\hat{u})
	\right\rbrace
	& \leq \boxed{9} \times \frac{1}{17.5}.
\end{align}
As before, the important point is that the constants 
$\boxed{9} \times \frac{9}{8.75} \times \frac{1}{8.75}$ 
and $\boxed{9} \times \frac{1}{17.5}$ from 
the estimates 
\eqref{E:FIRSTHARDMULTTERMGRONWALLCARRIEDOUT} and \eqref{E:SECONDHARDMULTTERMGRONWALL} 
verify $\frac{9}{(8.75)^2} + \frac{9}{17.5} < .65 < 1.$ 
As we will see, the remaining terms on the right-hand side of \eqref{E:Q0TOPORDERGRONWALLREADYINEQUALITY}
can be suitably bounded by choosing $\Conthree$ and $\Cononestar$ to be sufficiently large. 

In our analysis of the remaining terms,
in order to connect the quantity 
$\totonemax{\leq 24}(t,u)$ to the quantity $q(t,u)$
and to account for the 
$\iota_3^4(t,u)$
discrepancy
between $\iota(t,u)$ and $\widetilde{\iota}(t,u),$
we use the following estimate, which follows from
the definitions of the quantities involved, 
the approximate monotonicity of the integrating factors,
and \eqref{E:DESIREDDIFFICULTQ1BOUNDINTERMSOFQ0}:
\begin{align} \label{E:Q1INTERMSOFLITTLEF}
		\iota^{-1}(t,u)
		& \sup_{(t',u') \in [0,t] \times [0,u]}
			\left\lbrace
				\iota_1^{-2 \Conthree}(t',u')
				\iota_3^{-4}(t',u')
				\totonemax{\leq 24}(t',u') 
			\right\rbrace
			\\
	& \leq 
		C
		\iota_1^{-2 \Conthree}(t,u)
		\sup_{(t',u') \in [0,t] \times [0,u]}
		\left\lbrace
				\iota_1^{-2 \Conthree}(t',u')
				\iota_2^{-2 \Conthree}(t',u')
				\iota_3^{-(2 \Cononestar + 4)}(t',u')
				\upmu_{\star}^{17.5}(t',u')
				\totonemax{\leq 24}(t',u')
		\right\rbrace
			\notag \\
	& \leq 
		C
		\iota_1^{-2 \Conthree}(t,u)
		\left\lbrace
				\mathring{\upepsilon}^2
				+ \varepsilon^3
				+ q(t,u)
		\right\rbrace.
		\notag
\end{align}
Furthermore,
\begin{align}
	& \mbox{\eqref{E:Q1INTERMSOFLITTLEF} also holds with any of \ }
	1, \iota_2^{-2 \Conthree}(t,u),
	\iota_3^{-2 \Cononestar}(t,u),
	\mbox{\ or \ }
	\upmu_{\star}^{17.5}(t,u)
	\mbox{\ in the role of \ }
	\iota_1^{-2 \Conthree}(t,u),
		\label{E:MOREQ1INTERMSOFLITTLEF} \\
	& \mbox{and the same estimates hold with \ }
		\totMormax{\leq 24} 
		\mbox{\ in place of \ }
		\totonemax{\leq 24}.
		\label{E:MORAWETZTERMINTERMSOFLITTLEF}
\end{align}

We now claim that the following estimates hold for the remaining terms on the last three lines of the right-hand side of 
\eqref{E:Q0TOPORDERGRONWALLREADYINEQUALITY}
(see Def.~\ref{D:TOPORDERMULTBOUNDINGQUANTITIES} and Def.~\ref{D:BELOWORDERMULTBOUNDINGQUANTITIES}):
\begin{align} 
	\sup_{(\hat{t},\hat{u}) \in [0,t] \times [0,u]}
	\left\lbrace
	 	\iota^{-1}(\hat{t},\hat{u})
		\topmulterrortwo{\leq 24}(\hat{t},\hat{u})
	\right\rbrace
	& \leq C \varepsilon q(t,u),
		\label{E:TOPORDERQ0SMALLASMIXEDTERMESTIMATE} \\
	\sup_{(\hat{t},\hat{u}) \in [0,t] \times [0,u]}
	\left\lbrace
	 	\iota^{-1}(\hat{t},\hat{u})
		\topmulterrorthree{\leq 24}(\hat{t},\hat{u})
	\right\rbrace
	& \leq C
			\left\lbrace
					\mathring{\upepsilon}^2
					+ \varepsilon^3
					+ \varepsilon q(t,u)
				\right\rbrace,
		\label{E:TOPORDERSECONDQ0SMALLASMIXEDTERMESTIMATE} \\
	\sup_{(\hat{t},\hat{u}) \in [0,t] \times [0,u]}
		\left\lbrace
			\iota^{-1}(\hat{t},\hat{u})
			\topmulterrorfour{\leq 24}(\hat{t},\hat{u})
		\right\rbrace
	& \leq \frac{C}{\Cononestar} q(t,u),
		\label{E:TOPORDERQ0DANGEROUSMIXEDTERMESTIMATE} \\
	\sup_{(\hat{t},\hat{u}) \in [0,t] \times [0,u]}
	\left\lbrace
		\iota^{-1}(\hat{t},\hat{u})
		\topmulterrorsix{\leq 24}(\hat{t},\hat{u})
	\right\rbrace
	& \leq \frac{C}{\Cononestar} q(t,u),
		\label{E:TOPORDERQ0LOWGROWTHTERMESTIMATE} \\
	\sup_{(\hat{t},\hat{u}) \in [0,t] \times [0,u]}
	\left\lbrace
	 	\iota^{-1}(\hat{t},\hat{u})
		\topmulterrorseven{\leq 24}(\hat{t},\hat{u})
	\right\rbrace
	& \leq C \varepsilon q(t,u),
		\label{E:TOPORDERQ0LARGETIMEDECAYCOUNTERSUPMULOSSTERMESTIMATE} \\
	\sup_{(\hat{t},\hat{u}) \in [0,t] \times [0,u]}
		\left\lbrace
			\iota^{-1}(\hat{t},\hat{u})
			\topmulterroreight{\leq 24}(\hat{t},\hat{u})
		\right\rbrace
	& \leq    C
						\left\lbrace
							\mathring{\upepsilon}^2
							+ \varepsilon^3
							+ \varepsilon q(t,u)
						\right\rbrace,
		\label{E:TOPORDERQ0INVOLVINGQ1LARGETIMEDECAYCOUNTERSUPMULOSSTERMESTIMATE} 
\end{align}
and
\begin{align} 
	(1 + \varsigma^{-1})
	\sup_{(\hat{t},\hat{u}) \in [0,t] \times [0,u]}
	\left\lbrace
	 	\iota^{-1}(\hat{t},\hat{u})	
		\multerrorone{\leq 24}(\hat{t},\hat{u})
	\right\rbrace
	& \leq \frac{C}{\Conthree} (1 + \varsigma^{-1}) q(t,u),
		\label{E:BELOWTOPORDEREASYQ0LARGETIMEDECAYTERMESTIMATE} \\
	\sup_{(\hat{t},\hat{u}) \in [0,t] \times [0,u]}
	\left\lbrace
	 	\iota^{-1}(\hat{t},\hat{u})
		\multerrortwo{\leq 24}(\hat{t},\hat{u})
	\right\rbrace
	& \leq  C 
			 	 	\left\lbrace
						\mathring{\upepsilon}^2
						+ \varepsilon^3
					\right\rbrace
				+ C \varepsilon^{1/2} q(t,u),
			\label{E:UNUSUALINTEGRATINGFACTORTERMGRONWALLED} \\
	(1 + \varsigma^{-1})
	\sup_{(\hat{t},\hat{u}) \in [0,t] \times [0,u]}
		\left\lbrace
			\iota^{-1}(\hat{t},\hat{u})
			\multerrorthree{\leq 24}(\hat{t},\hat{u})
		\right\rbrace
	& \leq \frac{C}{\Conthree} 
	        (1 + \varsigma^{-1})
	         \left\lbrace
							\mathring{\upepsilon}^2
							+ \varepsilon^3
							+ q(t,u)
					 \right\rbrace,
		\label{E:EASYBELOWTOPORDERQ0LARGETIMETERMTHATINVOLVESQ1ESTIMATE} \\
	\sup_{(\hat{t},\hat{u}) \in [0,t] \times [0,u]}
	\left\lbrace
		\iota^{-1}(\hat{t},\hat{u})
		\multerrorfour{\leq 24}(\hat{t},\hat{u})
	\right\rbrace
	& \leq \frac{C}{\Conthree} q(t,u),
		\label{E:EASYBELOWTOPORDERQ0UTERMESTIMATE}
		\\
	(1 + \varsigma^{-1})
	\sup_{(\hat{t},\hat{u}) \in [0,t] \times [0,u]}
	\left\lbrace
	 	\iota^{-1}(\hat{t},\hat{u})	
		\multerrorfive{\leq 24}(\hat{t},\hat{u})
	\right\rbrace
	& \leq \frac{C}{\Conthree} 
	        (1 + \varsigma^{-1})
	         \left\lbrace
							\mathring{\upepsilon}^2
							+ \varepsilon^3
							+ q(t,u)
					 \right\rbrace,
		\label{E:EASYBELOWTOPORDERQ0UTERMTHATINVOLVESQ1ESTIMATE}  \\
	(\varsigma + \varepsilon)
	\sup_{(\hat{t},\hat{u}) \in [0,t] \times [0,u]}
	\left\lbrace
	 	\iota^{-1}(\hat{t},\hat{u})
		\multerrorsix{\leq 24}(\hat{t},\hat{u})
	\right\rbrace
	& \leq 
			C 	
			\left\lbrace
				\mathring{\upepsilon}^2
				+ \varepsilon^3
				+ (\varsigma + \varepsilon) q(t,u)
			\right\rbrace.
		\label{E:BELOWTOPORDEREASYMORAWETZINTEGRALTERMESTIMATE}
\end{align}

We now explain how to derive the bounds 
\eqref{E:TOPORDERQ0SMALLASMIXEDTERMESTIMATE}-
\eqref{E:TOPORDERQ0INVOLVINGQ1LARGETIMEDECAYCOUNTERSUPMULOSSTERMESTIMATE}.

The bound \eqref{E:TOPORDERQ0DANGEROUSMIXEDTERMESTIMATE} 
(see definition \eqref{E:MIXEDTYPEHARDMULTTERM})
can be proved by using an argument similar to the one used to prove 
\eqref{E:ABSORBABLEBOUNDFIRSTHARDMORTERMGRONWALL}.
More precisely, the inner time integral on the right-hand side of \eqref{E:MIXEDTYPEHARDMULTTERM}
is handled with the key integral estimate \eqref{E:KEYMUTOAPOWERINTEGRALBOUND},
and then the outer time integral on the right-hand side of \eqref{E:MIXEDTYPEHARDMULTTERM}
is bounded with the help of inequality \eqref{E:MAINGRONWALLPROOFKEYLOGGROWTHGENERATINGINTEGRALBOUND}, 
which provides the smallness factor $\frac{1}{\Cononestar}$ on the right-hand side of \eqref{E:TOPORDERQ0DANGEROUSMIXEDTERMESTIMATE}.

To obtain the bound \eqref{E:TOPORDERSECONDQ0SMALLASMIXEDTERMESTIMATE} 
(see definition \eqref{E:TOPORDERMULTTERMSECONDSMALLA}), 
we first use the approximate monotonicity of the integrating factors
to bound the term in braces on the left-hand side of \eqref{E:TOPORDERSECONDQ0SMALLASMIXEDTERMESTIMATE}  
as follows:
\begin{align} \label{E:FIRSTESTIMATETOPORDERSECONDQ0SMALLASMIXEDTERMESTIMATE}
	& \leq
	C \varepsilon
	\iota^{-1}(\hat{t},\hat{u})
				\int_{t'=0}^{\hat{t}}
						\frac{1} 
								 {(1 + t')^{1 + a} \upmu_{\star}(t',\hat{u})} 
						\totzeromax{\leq N}^{1/2}(t',u) 
						\int_{s=0}^{t'}
							\frac{1}{(1 + s)} 
							\frac{1}{\upmu_{\star}(s,\hat{u})} 
							\totonemax{\leq N}^{1/2}(s,\hat{u}) 
						\, ds
				\, dt'
				\\
	& \leq
		C \varepsilon
		\iota^{-1}(\hat{t},\hat{u})
		\left\lbrace
			\sup_{(t',u') \in [0,\hat{t}] \times [0,\hat{u}]}
			\upmu_{\star}^{8.75}(t',u')
			\totzeromax{\leq 24}^{1/2}(t',u')
		\right\rbrace
		\notag \\
	& \ \ \ \ \
		\times
		\left\lbrace
			\sup_{(t',u') \in [0,\hat{t}] \times [0,\hat{u}]}
			\iota_3^{-2}(t',u')
			\upmu_{\star}^{8.75}(t',u')
			\totonemax{\leq 24}^{1/2}(t',u')
		\right\rbrace
			\notag \\
	& \ \ \ \ \
		\times
			\int_{t'=0}^{\hat{t}}
						\frac{1} 
								 {(1 + t')^{1 + a} \upmu_{\star}^{9.75}(t',\hat{u})} 
								\iota_3^2(t',\hat{u})
						\int_{s=0}^{t'}
							\frac{1}{(1 + s)} 
							\frac{1}{\upmu_{\star}^{9.75}(s,\hat{u})} 
						\, ds
				\, dt'.
		\notag
\end{align}
We bound the inner time integral on the right-hand side of \eqref{E:FIRSTESTIMATETOPORDERSECONDQ0SMALLASMIXEDTERMESTIMATE}
with the estimate \eqref{E:LOGLOSSKEYMUINTEGRALBOUND}
and then the outer time integral with the estimate
\eqref{E:LESSSINGULARTERMSMUINTEGRALBOUND},
which yields that the right-hand side of \eqref{E:FIRSTESTIMATETOPORDERSECONDQ0SMALLASMIXEDTERMESTIMATE}
is
\begin{align} \label{E:SECONDESTIMATETOPORDERSECONDQ0SMALLASMIXEDTERMESTIMATE}
	& \leq
		C \varepsilon
		\upmu_{\star}^{-17.5}(\hat{t},\hat{u})
		\iota^{-1}(\hat{t},\hat{u})
		\left\lbrace
			\sup_{(t',u') \in [0,\hat{t}] \times [0,\hat{u}]}
			\upmu_{\star}^{8.75}(t',u')
			\totzeromax{\leq 24}^{1/2}(t',u')
		\right\rbrace
			\\
 & \ \ \ \ \ 
	\times
	\left\lbrace
			\sup_{(t',u') \in [0,\hat{t}] \times [0,\hat{u}]}
			\iota_3^{-2}(t',u')
			\upmu_{\star}^{8.75}(t',u')
			\totonemax{\leq 24}^{1/2}(t',u')
		\right\rbrace.
		\notag
\end{align}
Using the monotonicity of $\iota_1,$ $\iota_2,$ and $\iota_3$ in both of their arguments
and the bound \eqref{E:DESIREDDIFFICULTQ1BOUNDINTERMSOFQ0}, we 
deduce that the right-hand side of
\eqref{E:SECONDESTIMATETOPORDERSECONDQ0SMALLASMIXEDTERMESTIMATE}
is
\begin{align} \label{E:THIRDESTIMATETOPORDERSECONDQ0SMALLASMIXEDTERMESTIMATE}
& \leq C \varepsilon q^{1/2}(\hat{t},\hat{u}) \widetilde{q}^{1/2}(\hat{t},\hat{u})
	 \leq C \left\lbrace
					\mathring{\upepsilon}^2
					+ \varepsilon^3
					+ q(t,u)
				\right\rbrace.
\end{align}
The desired estimate \eqref{E:TOPORDERSECONDQ0SMALLASMIXEDTERMESTIMATE} now 
follows from taking 
$\sup_{(\hat{t},\hat{u}) \in [0,t] \times [0,u]}$
in inequality \eqref{E:THIRDESTIMATETOPORDERSECONDQ0SMALLASMIXEDTERMESTIMATE}.

The bound \eqref{E:TOPORDERQ0SMALLASMIXEDTERMESTIMATE}
(see definition \eqref{E:TOPORDERMULTTERMSMALLA})
can be proved by using arguments similar to the ones we used
to prove \eqref{E:TOPORDERSECONDQ0SMALLASMIXEDTERMESTIMATE}.


The bound \eqref{E:TOPORDERQ0LOWGROWTHTERMESTIMATE} 
(see definition \eqref{E:TOPORDERMULTTERCREATESLOGGROWTH})
can be proved by using the same argument used to prove \eqref{E:TOPORDERQ1LOGGROWTHTERMBOUNDED}.

The bound \eqref{E:TOPORDERQ0LARGETIMEDECAYCOUNTERSUPMULOSSTERMESTIMATE} 
(see definition \eqref{E:TOPORDERQ0LARGETIMEDECAYCOUNTERSUPMULOSSTERM}) 
can similarly be proved 
by using inequality \eqref{E:LESSSINGULARTERMSMUINTEGRALBOUND} with $\Conone = 0$
and $\Contwo = 17.5+1 = 18.5.$ The bound \eqref{E:TOPORDERQ0INVOLVINGQ1LARGETIMEDECAYCOUNTERSUPMULOSSTERMESTIMATE}
(see definition \eqref{E:TOPORDERQ0INVOLVINGQ1LARGETIMEDECAYCOUNTERSUPMULOSSTERM})
can similarly be proved 
with the help of the estimate \eqref{E:DESIREDDIFFICULTQ1BOUNDINTERMSOFQ0}.

We now explain how to derive the easier estimates \eqref{E:BELOWTOPORDEREASYQ0LARGETIMEDECAYTERMESTIMATE}-\eqref{E:BELOWTOPORDEREASYMORAWETZINTEGRALTERMESTIMATE}.
The bound \eqref{E:BELOWTOPORDEREASYQ0LARGETIMEDECAYTERMESTIMATE} (see definition 
\eqref{E:BELOWTOPORDEREASYQ0LARGETIMEDECAYTERM})
can be proved by using an argument similar to the one we used to prove 
\eqref{E:ABSORBABLEBOUNDFIRSTHARDMORTERMGRONWALL},
but with $\iota_2^{-2 \Conthree}$ in place of $\upmu_{\star}^{-17.5}$ and the estimate 
\eqref{E:MAINGRONWALLPROOFKEYLARGETDECAYINTEGRALBOUND} in place of
the estimate \eqref{E:KEYMUTOAPOWERINTEGRALBOUND}.

To obtain the bound \eqref{E:UNUSUALINTEGRATINGFACTORTERMGRONWALLED} 
(see definition \eqref{E:STRANGEINTEGRATINGFACTORINTEGRALBOUND}),
we use inequality \eqref{E:MOREQ1INTERMSOFLITTLEF} in the case of the constant function $1$
to bound the term in braces on the left-hand side of \eqref{E:UNUSUALINTEGRATINGFACTORTERMGRONWALLED} as follows:
\begin{align}
	& C	\varepsilon^{1/2}		
			 \iota^{-1}(\hat{t},\hat{u})
			 \int_{t'=0}^{\hat{t}} 
					\frac{\ln(\myexp + t')}{(1 + t')^2 \sqrt{\ln(\myexp + \hat{t}) - \ln(\myexp + t')}}
					\totonemax{\leq N}(t',\hat{u}) 
			 \, dt'
			\label{E:UNUSUALINTEGRATINGFACTORTERMALMOSTGRONWALLED}	\\
	& \leq C \varepsilon^{1/2}	
			 \left\lbrace
					\mathring{\upepsilon}^2
					+ \varepsilon^3
					+ q(t,u)
				\right\rbrace
			\int_{t'=0}^{\hat{t}} 
					\frac{\ln^5(\myexp + t')}{(1 + t')^2 \sqrt{\ln(\myexp + t) - \ln(\myexp + t')}}
			\, dt'.
				\notag 
\end{align}
The desired estimate \eqref{E:UNUSUALINTEGRATINGFACTORTERMGRONWALLED} now follows from
\eqref{E:UNUSUALINTEGRATINGFACTORTERMALMOSTGRONWALLED} and Lemma~\ref{L:UNUSUALINTEGRATINGFACTORESTIMATE}.

To obtain the bound \eqref{E:EASYBELOWTOPORDERQ0LARGETIMETERMTHATINVOLVESQ1ESTIMATE}
(see definition \eqref{E:EASYBELOWTOPORDERQ0LARGETIMETERMTHATINVOLVESQ1}),
we use inequality \eqref{E:MOREQ1INTERMSOFLITTLEF} in the case of the function $\iota_2^{-2 \Conthree}$
to bound the term in braces on the left-hand side of \eqref{E:EASYBELOWTOPORDERQ0LARGETIMETERMTHATINVOLVESQ1ESTIMATE} as follows:
\begin{align}
	& C	 (1 + \varsigma^{-1})	
			 \iota^{-1}(\hat{t},\hat{u})
			 \int_{t'=0}^{\hat{t}} 
					\frac{1}{(1 + t')^{3/2}}
					\totonemax{\leq N}(t',\hat{u}) 
			 \, dt'
			\label{E:EASYBELOWTOPORDERQ0LARGETIMETERMTHATINVOLVESQ1ALMOSTESTIMATE}	\\
	& \leq C (1 + \varsigma^{-1})	
			 \left\lbrace
					\mathring{\upepsilon}^2
					+ \varepsilon^3
					+ q(t,u)
				\right\rbrace
			\iota_2^{-2 \Conthree}(\hat{t},\hat{u})
			\int_{t'=0}^{\hat{t}} 
				 \frac{\iota_2^4(t',u)}{(1 + t')^{3/2}}
				 \iota_2^{2 \Conthree}(t',u)
			\, dt'.
				\notag 
\end{align}
The desired bound \eqref{E:EASYBELOWTOPORDERQ0LARGETIMETERMTHATINVOLVESQ1ESTIMATE} now follows from
\eqref{E:EASYBELOWTOPORDERQ0LARGETIMETERMTHATINVOLVESQ1ALMOSTESTIMATE} and inequality \eqref{E:MAINGRONWALLPROOFKEYLARGETDECAYINTEGRALBOUND}.

The bound \eqref{E:EASYBELOWTOPORDERQ0UTERMESTIMATE} (see definition \eqref{E:EASYBELOWTOPORDERQ0UTERM})
can be proved by using an argument similar to the one we used to prove \eqref{E:BELOWTOPORDEREASYQ0LARGETIMEDECAYTERMESTIMATE},
but with $\iota_1^{-2 \Conthree}$ in place of $\iota_2^{-2 \Conthree}$ and the estimate \eqref{E:MAINGRONWALLPROOFKEYEASYUINTEGRALBOUND} in place of
the estimate \eqref{E:MAINGRONWALLPROOFKEYLARGETDECAYINTEGRALBOUND}.

To obtain the bound \eqref{E:BELOWTOPORDEREASYMORAWETZINTEGRALTERMESTIMATE} (see definition \eqref{E:BELOWTOPORDEREASYMORAWETZINTEGRALTERM}),
we first use inequality \eqref{E:MORAWETZTERMINTERMSOFLITTLEF} in the case of the constant function $1$
to deduce that
\begin{align}
	\iota^{-1}(t',\hat{u}) \frac{\Morint_{(\leq N)}(t',\hat{u})}{(1 + t')^{1/2}}
	& \leq C \frac{\iota_3^4(t',\hat{u})}{(1 + t')^{1/2}} 
					\left(
						\mathring{\upepsilon}^2
						+ \varepsilon^3
						+ q(t',\hat{u})
					\right).
\end{align}
Hence, the term in braces on the left-hand side of \eqref{E:BELOWTOPORDEREASYMORAWETZINTEGRALTERMESTIMATE}
can be bounded as follows:
\begin{align}
	\leq
			C	
			\iota^{-1}(\hat{t},\hat{u})
			\sup_{t' \in [0,\hat{t}]}\frac{\Morint_{(\leq N)}(t',\hat{u})}{(1 + t')^{1/2}}
	& \leq 
			C	
			\sup_{t' \in [0,\hat{t}]} 
			\left\lbrace
				\iota^{-1}(t',\hat{u}) 
				\frac{\Morint_{(\leq N)}(t',\hat{u})}{(1 + t')^{1/2}}
			\right\rbrace
			\label{E:BELOWTOPORDEREASYMORAWETZINTEGRALTERMALMOSTESTIMATE}	\\
	& \leq C 
			 \sup_{t' \in [0,\hat{t}]}
			 \left\lbrace
			 		\frac{\iota_3^4(t',u)}{(1 + t')^{1/2}}
			 		\left(
						\mathring{\upepsilon}^2
						+ \varepsilon^3
						+ q(t',u)
					\right)
				\right\rbrace
				\notag 
				\\
	& \leq C 	
			 \left\lbrace
					\mathring{\upepsilon}^2
					+ \varepsilon^3
					+ q(t,u)
				\right\rbrace.
				\notag
\end{align}
The desired bound \eqref{E:BELOWTOPORDEREASYMORAWETZINTEGRALTERMESTIMATE} now easily follows from inequality 
\eqref{E:BELOWTOPORDEREASYMORAWETZINTEGRALTERMALMOSTESTIMATE}.

To obtain the bound \eqref{E:EASYBELOWTOPORDERQ0UTERMTHATINVOLVESQ1ESTIMATE} 
(see definition \eqref{E:EASYBELOWTOPORDERQ0UTERMTHATINVOLVESQ1}),
we use an argument similar to the one used to prove 
\eqref{E:BELOWTOPORDEREASYMORAWETZINTEGRALTERMESTIMATE},
but we use inequality \eqref{E:MOREQ1INTERMSOFLITTLEF} in the case of the constant function $1$
in place of inequality \eqref{E:MORAWETZTERMINTERMSOFLITTLEF}.

In order to complete the proof of \eqref{E:PREDESIREDDIFFICULTQ0BOUND},
it remains for us to bound the terms arising from the terms on the first two lines of the right-hand side of 
\eqref{E:Q0TOPORDERGRONWALLREADYINEQUALITY}.
The five corresponding terms are respectively bounded from above by $\leq$ 
$C \mathring{\upepsilon}^2,$ 
$C \varepsilon^3,$ 
$C \varepsilon q(t,u),$
$C \varepsilon^3,$ 
and 
$C \varepsilon^3.$ 
The first of these bounds follows from Lemma~\ref{L:CONTROLLINGQUANTITIESAREINTIALLYSMALL}, 
while the next two are easy to derive. 
To bound the two terms on the second line of the right-hand side of 
\eqref{E:Q0TOPORDERGRONWALLREADYINEQUALITY}
by $\leq C \varepsilon^3,$
we argue as in our proof of \eqref{E:TERMSBOUNDEDSWITHBOOTSTRAPASSUMPTIONS}.

Combining all of the above bounds, we arrive at the following analog of 
\eqref{E:MOREEXPLICITPREDESIREDDIFFICULTQ1BOUNDINTERMSOFQ0}:
\begin{align} \label{E:MOREEXPLICITPREDESIREDDIFFICULTQ0BOUND}
	q(t,u) & \leq 
				C (1 + \varsigma^{-1})
				\left\lbrace
					\mathring{\upepsilon}^2
					+ \varepsilon^3
				\right\rbrace
					\\
			& \ \ 
					+ \left\lbrace
						(1 + C \sqrt{\epsilon})
						\frac{6}{(8.75)^2} 
						+ (1 + C \sqrt{\epsilon})
							\frac{2}{17.5}
						+ \frac{C}{\Cononestar}
						+ \frac{C}{\Conthree}
							(1 + \varsigma^{-1})
						+ C \varsigma
						+ C \varepsilon^{1/2}
					\right\rbrace
					q(t,u).
				\notag
\end{align}
The desired bound \eqref{E:PREDESIREDDIFFICULTQ0BOUND} thus follows from \eqref{E:MOREEXPLICITPREDESIREDDIFFICULTQ0BOUND}
if we first choose $\varsigma$ and $\varepsilon$ to be sufficiently small
and we then choose $\Cononestar$ and $\Conthree$ to be sufficiently large
(at least as large as they were chosen to be in the 
part of the proof following inequality \eqref{E:MOREEXPLICITPREDESIREDDIFFICULTQ1BOUNDINTERMSOFQ0}).

\medskip

\noindent{\emph{Proof of the estimates for the just-below-top-order quantities.}}
We now prove the desired estimates for the below top-order quantities
$\totzeromax{\leq N}(t,u),$
$\totonemax{\leq N}(t,u),$
$\totMormax{\leq N}(t,u),$
$0 \leq N \leq 23.$ These estimates are much easier to prove than the top-order estimates.
We give complete details for the estimates of
$\totzeromax{\leq 23}(t,u),$
$\totonemax{\leq 23}(t,u),$
and $\totMormax{\leq 23}(t,u).$
We then indicate the minor changes in the proof needed to derive the desired estimates
for $\totzeromax{\leq N}(t,u),$
$\totonemax{\leq N}(t,u),$
and $\totMormax{\leq N}(t,u)$ when $0 \leq N \leq 22.$

To begin, in place of \eqref{E:TOPORDERQ0INTEGRATINGFACTORPRODUCT} and \eqref{E:TOPORDERQ1INTEGRATINGFACTORPRODUCT},
we define
\begin{align}
	\iota(t,u)
	& :=
		\iota_1^{2 \Conthree}(t,u)
		\iota_2^{2 \Conthree}(t,u)
		\upmu_{\star}^{-15.5}(t,u),
		 \label{E:LOWERORDERQ0INTEGRATINGFACTORPRODUCT} \\
	\widetilde{\iota}(t,u)
	& := 
		\iota_1^{2 \Conthree}(t,u)
		\iota_2^{2 \Conthree}(t,u)
		\iota_3^4(t,u)
		\upmu_{\star}^{-15.5}(t,u).
		\label{E:LOWERORDERQ1INTEGRATINGFACTORPRODUCT} 
\end{align}
Notice in particular that the power of $\upmu_{\star}^{-1}$ in 
\eqref{E:LOWERORDERQ0INTEGRATINGFACTORPRODUCT} and \eqref{E:LOWERORDERQ1INTEGRATINGFACTORPRODUCT}
has been reduced by $2$ and that there is no logarithmic factor 
$\iota_3(t,u)$ in \eqref{E:LOWERORDERQ0INTEGRATINGFACTORPRODUCT}. 
Next, in place of \eqref{E:TOPORDERQ0INTEGRATINGFACTORPRODUCT} and \eqref{E:TOPORDERQ1INTEGRATINGFACTORPRODUCT}, we define
\begin{align}
	q(t,u) 
	& : = 
		\sup_{(\hat{t},\hat{u}) \in [0,t] \times [0,u]}
			\iota^{-1}(\hat{t},\hat{u})
			\totzeromax{\leq 23}(\hat{t},\hat{u}),
		\label{E:LOWERORDERRENORMALIZEDQ0} \\
\widetilde{q}(t,u) 
&: = 	\sup_{(\hat{t},\hat{u}) \in [0,t] \times [0,u]}
				\widetilde{\iota}^{-1}(\hat{t},\hat{u})
				\left\lbrace 
					\totonemax{\leq 23}(\hat{t},\hat{u}) 
					+ \totMormax{\leq 23}(\hat{t},\hat{u})
				\right\rbrace.
				\label{E:LOWERORDERRENORMALIZEDQ1}
\end{align}
Our goal is to prove the following analogs of \eqref{E:KEYTOPORDERLITTLEFBOUND} and \eqref{E:KEYTOPORDERWIDETILDELITTLEFBOUND}:
\begin{align}
	q(t,u) \label{E:KEYLOWERORDERLITTLEFBOUND}
	& \leq C 
				\left\lbrace
					\mathring{\upepsilon}^2
					+ \varepsilon^3
				\right\rbrace,
				\\
	\widetilde{q}(t,u) 
	& \leq C 
				\left\lbrace
					\mathring{\upepsilon}^2
					+ \varepsilon^3
				\right\rbrace.
				\label{E:KEYLOWERORDERWIDETILDELITTLEFBOUND}
\end{align}
The desired estimates \eqref{E:Q0MIDORDERKEYBOUND} and \eqref{E:Q1MIDORDERKEYBOUND} for
$\totzeromax{\leq 23},$ $\totonemax{\leq 23},$ and $\totMormax{\leq 23}$ easily follow from
\eqref{E:KEYLOWERORDERLITTLEFBOUND}-\eqref{E:KEYLOWERORDERWIDETILDELITTLEFBOUND} and the definitions
of the quantities involved.

In order to prove \eqref{E:KEYLOWERORDERLITTLEFBOUND} and \eqref{E:KEYLOWERORDERWIDETILDELITTLEFBOUND},
we first prove the following analog of \eqref{E:DESIREDDIFFICULTQ1BOUNDINTERMSOFQ0}:
\begin{align} \label{E:LOWERORDERDESIREDDIFFICULTQ1BOUNDINTERMSOFQ0}
	\widetilde{q}(t,u)
	& \leq 
				C 
				\left\lbrace
					\mathring{\upepsilon}^2
					+ \varepsilon^3
					+ q(t,u)
				\right\rbrace.
\end{align}
In order to prove \eqref{E:LOWERORDERDESIREDDIFFICULTQ1BOUNDINTERMSOFQ0}, 
we will prove the following analog of \eqref{E:PREDESIREDDIFFICULTQ1BOUNDINTERMSOFQ0}:
\begin{align} \label{E:LOWERORDERPREDESIREDDIFFICULTQ1BOUNDINTERMSOFQ0}
	\widetilde{q}(t,u) & \leq 
				C 
				\left\lbrace
					\mathring{\upepsilon}^2
					+ \varepsilon^3
					+ q(t,u)
				\right\rbrace
				+ \upalpha
					\widetilde{q}(t,u),
\end{align}
where $0 < \upalpha < 1$ is a constant. 

We now claim that the following estimates hold for the terms on
the second through fourth lines of the right-hand side of \eqref{E:Q1BELOWTOPORDERGRONWALLREADYINEQUALITY}
in the case $N=23$ (see Def.~\ref{D:BELOWORDERMORBOUNDINGQUANTITIES}):
\begin{align}
	\sup_{(\hat{t},\hat{u}) \in [0,t] \times [0,u]}
	\left\lbrace
	 	\widetilde{\iota}^{-1}(\hat{t},\hat{u})	
		\morerrorzero{\leq 24}(\hat{t},\hat{u})
	\right\rbrace
	& \leq C \varepsilon^3,
		 \label{E:MORLOSSOFONEDERIVERRORINTEGRALTERMESTIMATE} \\
	(1 + \widetilde{\varsigma}^{-1})
	\sup_{(\hat{t},\hat{u}) \in [0,t] \times [0,u]}
	\left\lbrace
	 	\widetilde{\iota}^{-1}(\hat{t},\hat{u})	
		\morerrorone{\leq 23}(\hat{t},\hat{u})
	\right\rbrace
	& \leq C (1 + \widetilde{\varsigma}^{-1}) q(t,u),
		 \label{E:DESCENTEASYBELOWTOPORDERQ1THATINVOLVESQ0TERMESTIMATE} \\
	\sup_{(\hat{t},\hat{u}) \in [0,t] \times [0,u]}
	\left\lbrace
	 	\widetilde{\iota}^{-1}(\hat{t},\hat{u})
		\morerrortwo{\leq 23}(\hat{t},\hat{u})
	\right\rbrace
	& \leq \frac{1}{2} \widetilde{q}(t,u),
				\label{E:DESCENTEASYBELOWTOPORDERQ1LOGLOSSTERMESTIMATE} \\
	(1 + \widetilde{\varsigma}^{-1})
	\sup_{(\hat{t},\hat{u}) \in [0,t] \times [0,u]}
		\left\lbrace
			\widetilde{\iota}^{-1}(\hat{t},\hat{u})
			\morerrorthree{\leq 23}(\hat{t},\hat{u})
		\right\rbrace
	& \leq \frac{C}{\Conthree} (1 + \widetilde{\varsigma}^{-1}) \widetilde{q}(t,u),
		\label{E:DESCENTEASYBELOWTOPORDERLARGETIMEDECAYTERM} \\
	(1 + \widetilde{\varsigma}^{-1})
	\sup_{(\hat{t},\hat{u}) \in [0,t] \times [0,u]}
	\left\lbrace
		\widetilde{\iota}^{-1}(\hat{t},\hat{u})
		\morerrorfour{\leq 23}(\hat{t},\hat{u})
	\right\rbrace
	& \leq \frac{C}{\Conthree} (1 + \widetilde{\varsigma}^{-1}) \widetilde{q}(t,u),
		\label{E:DESCENTEASYBELOWTOPORDERQ1UTERMESTIMATE}
		\\
	(\widetilde{\varsigma} + \varepsilon)
	\sup_{(\hat{t},\hat{u}) \in [0,t] \times [0,u]}
	\left\lbrace
		\widetilde{\iota}^{-1}(\hat{t},\hat{u})
		\morerrorfive{\leq 23}(\hat{t},\hat{u})
	\right\rbrace
	& \leq C (\widetilde{\varsigma} + \varepsilon) \widetilde{q}(t,u).
		\label{E:DESCENTEASYBELOWTOPORDERQ1MORAWETZTERMESTIMATE}
\end{align}
The bounds \eqref{E:DESCENTEASYBELOWTOPORDERQ1THATINVOLVESQ0TERMESTIMATE}-\eqref{E:DESCENTEASYBELOWTOPORDERQ1UTERMESTIMATE}
can be proved by using arguments similar to the ones we used to prove 
\eqref{E:EASYBELOWTOPORDERQ1THATINVOLVESQ0TERMESTIMATE}-\eqref{E:EASYBELOWTOPORDERQ1MORAWETZTERMESTIMATE}.
We note that the factor $\frac{1}{2}$ on the right-hand side of \eqref{E:DESCENTEASYBELOWTOPORDERQ1LOGLOSSTERMESTIMATE}
arises from the fact that the constant on the right-hand side of \eqref{E:EASYBELOWTOPORDERQ1LOGLOSSTERM}
is precisely $2.$

To obtain \eqref{E:MORLOSSOFONEDERIVERRORINTEGRALTERMESTIMATE}, we set
$N = 23$ in \eqref{E:MORLOSSOFONEDERIVERRORINTEGRALTERM},
insert the estimates
$\totzeromax{\leq 24}^{1/2} \leq C 
				\left\lbrace
					\mathring{\upepsilon}
					+ \varepsilon^{3/2}
				\right\rbrace
				\ln^{\Cononestar} (\myexp + t) \upmu_{\star}^{-8.75}(t,u)$ 
and $\totonemax{\leq 24}^{1/2} \leq C 
				\left\lbrace
					\mathring{\upepsilon}
					+ \varepsilon^{3/2}
				\right\rbrace
				\ln^{\Cononestar + 2} (\myexp + t) \upmu_{\star}^{-8.75}(t,u),$
which follow from 
\eqref{E:KEYTOPORDERLITTLEFBOUND}-\eqref{E:KEYTOPORDERWIDETILDELITTLEFBOUND}
(and the constant $\Cononestar$ has already been chosen above),
and use the integral estimates \eqref{E:LESSSINGULARTERMSMUINTEGRALBOUND} and \eqref{E:LOGLOSSKEYMUINTEGRALBOUND}
to deduce that
\begin{align}
\morerrorzero{\leq 24}(\hat{t},\hat{u})
& \leq 
		C \varepsilon^3
		\int_{t'=0}^{\hat{t}} 
			\frac{\ln^{2 \Cononestar + 4}(\myexp + t')}{(1 + t')^2}
			\left(
			\int_{s=0}^{t'}
				\frac{1}{(1 + s)\upmu_{\star}^{8.75 + .5}(s,u)} 
			\, ds
			\right)^2
		\, dt' 
	 	 \label{E:BELOWTOPORDERQ1LOSEADERIVATIVEGAININUPMUESTIMATE} \\
	& \leq 	C \varepsilon^3
		\int_{t'=0}^{\hat{t}} 
			\frac{1}{(1 + t')^{3/2}}
			\left(
				\frac{1}{\upmu_{\star}^{8.25}(t',u)}
			\right)^2
		\, dt'
		\leq C \varepsilon^3 \upmu_{\star}^{-15.5}(\hat{t},\hat{u}).
		\notag
\end{align}
The desired estimate \eqref{E:MORLOSSOFONEDERIVERRORINTEGRALTERMESTIMATE} now follows from \eqref{E:BELOWTOPORDERQ1LOSEADERIVATIVEGAININUPMUESTIMATE}
and the definitions of the quantities involved.

In order to complete the proof of \eqref{E:LOWERORDERPREDESIREDDIFFICULTQ1BOUNDINTERMSOFQ0},
it remains for us to bound the terms arising from the terms on the first line of the right-hand side of 
\eqref{E:Q1BELOWTOPORDERGRONWALLREADYINEQUALITY}.
Arguing as in our proof of the bounds for the terms on the first line of the
right-hand side of \eqref{E:Q1TOPORDERGRONWALLREADYINEQUALITY},
we deduce that the two corresponding terms are respectively bounded by $\leq$
$C \mathring{\upepsilon}^2$ 
and 
$C \varepsilon^3.$

Combining all of the above bounds, we arrive at the following analog of 
\eqref{E:MOREEXPLICITPREDESIREDDIFFICULTQ1BOUNDINTERMSOFQ0}:
\begin{align} \label{E:DESCENTMOREEXPLICITPREDESIREDDIFFICULTQ1BOUNDINTERMSOFQ0}
	\widetilde{q}(t,u) & \leq 
				C 
				\left\lbrace
					\mathring{\upepsilon}^2
					+ \varepsilon^3
					+ q(t,u)
				\right\rbrace
				+ \left\lbrace
						\frac{1}{2}
						+ \frac{C}{\Conthree}
							(1 + \widetilde{\varsigma}^{-1})	
						+ C \widetilde{\varsigma}
						+ C \varepsilon
					\right\rbrace
					\widetilde{q}(t,u),
\end{align}
and the desired bound \eqref{E:LOWERORDERPREDESIREDDIFFICULTQ1BOUNDINTERMSOFQ0} thus follows
if first $\widetilde{\varsigma}$ and $\varepsilon$ are chosen to be sufficiently small
and then $\Conthree$ is chosen to be sufficiently large.

We now use inequality \eqref{E:LOWERORDERDESIREDDIFFICULTQ1BOUNDINTERMSOFQ0} to help us derive 
the estimate \eqref{E:KEYLOWERORDERLITTLEFBOUND}. To prove \eqref{E:KEYLOWERORDERLITTLEFBOUND},
we will argue as in our proof of \eqref{E:LOWERORDERDESIREDDIFFICULTQ1BOUNDINTERMSOFQ0}
to show that there exists a constant $0 < \upbeta < 1$ such that
\begin{align} \label{E:PREDESIREDKEYLOWERORDERLITTLEFBOUND}
	q(t,u) & \leq 
				C 
				\left\lbrace
					\mathring{\upepsilon}^2
					+ \varepsilon^3
				\right\rbrace
				+ \upbeta q(t,u).
\end{align}
The desired bound \eqref{E:KEYLOWERORDERLITTLEFBOUND} follows easily once we have shown \eqref{E:PREDESIREDKEYLOWERORDERLITTLEFBOUND}.

We claim that the following estimates hold for the terms on
the second through fourth lines of the right-hand side of \eqref{E:Q0BELOWTOPORDERGRONWALLREADYINEQUALITY}
in the case $N=23$ (see Def.~\ref{D:BELOWORDERMULTBOUNDINGQUANTITIES}):
\begin{align} 
	\sup_{(\hat{t},\hat{u}) \in [0,t] \times [0,u]}
	\left\lbrace
	 	\iota^{-1}(\hat{t},\hat{u})	
		\multerrorzero{\leq 24}(\hat{t},\hat{u})
	\right\rbrace
	& \leq C \varepsilon^3,
		 \label{E:MULTLOSSOFONEDERIVERRORINTEGRALTERMESTIMATE} \\
	(1 + \varsigma^{-1})
	\sup_{(\hat{t},\hat{u}) \in [0,t] \times [0,u]}
	\left\lbrace
	 	\iota^{-1}(\hat{t},\hat{u})	
		\multerrorone{\leq 23}(\hat{t},\hat{u})
	\right\rbrace
	& \leq \frac{C}{\Conthree} (1 + \varsigma^{-1}) q(t,u),
		\label{E:LOWERORDEREASYQ0LARGETIMEDECAYTERMESTIMATE} \\
	\sup_{(\hat{t},\hat{u}) \in [0,t] \times [0,u]}
	\left\lbrace
	 	\iota^{-1}(\hat{t},\hat{u})
		\multerrortwo{\leq 23}(\hat{t},\hat{u})
	\right\rbrace
	& \leq  C 
			 	 	\left\lbrace
						\mathring{\upepsilon}^2
						+ \varepsilon^3
						+ \varepsilon^{1/2} q(t,u)
					\right\rbrace,
			\label{E:LOWERORDERUNUSUALINTEGRATINGFACTORTERMGRONWALLED} \\
	(1 + \varsigma^{-1})
	\sup_{(\hat{t},\hat{u}) \in [0,t] \times [0,u]}
		\left\lbrace
			\iota^{-1}(\hat{t},\hat{u})
			\multerrorthree{\leq 23}(\hat{t},\hat{u})
		\right\rbrace
	& \leq  \frac{C}{\Conthree} 
					(1 + \varsigma^{-1})
					\left\lbrace
							\mathring{\upepsilon}^2
							+ \varepsilon^3
							+ q(t,u)
					 \right\rbrace,
		\label{E:EASYLOWERORDERQ0LARGETIMETERMTHATINVOLVESQ1ESTIMATE} \\
	\sup_{(\hat{t},\hat{u}) \in [0,t] \times [0,u]}
	\left\lbrace
		\iota^{-1}(\hat{t},\hat{u})
		\multerrorfour{\leq 23}(\hat{t},\hat{u})
	\right\rbrace
	& \leq \frac{C}{\Conthree} q(t,u),
		\label{E:EASYLOWERORDERQ0UTERMESTIMATE}
		\\
	(1 + \varsigma^{-1})
	\sup_{(\hat{t},\hat{u}) \in [0,t] \times [0,u]}
	\left\lbrace
	 	\iota^{-1}(\hat{t},\hat{u})	
		\multerrorfive{\leq 23}(\hat{t},\hat{u})
	\right\rbrace
	& \leq 
			 \frac{C}{\Conthree}  
			(1 + \varsigma^{-1})	
			\left\lbrace
				\mathring{\upepsilon}^2
				+ \varepsilon^3
				+ q(t,u)
			\right\rbrace,
		\label{E:EASYLOWERORDERQ0UTERMTHATINVOLVESQ1ESTIMATE}  \\
	\sup_{(\hat{t},\hat{u}) \in [0,t] \times [0,u]}
	(\varsigma + \varepsilon)
	\left\lbrace
	 	\iota^{-1}(\hat{t},\hat{u})
		\multerrorsix{\leq 23}(\hat{t},\hat{u})
	\right\rbrace
	& \leq 
			C 	
			\left\lbrace
				\mathring{\upepsilon}^2
				+ \varepsilon^3
				+ (\varsigma + \varepsilon) q(t,u)
			\right\rbrace.
		\label{E:LOWERORDEREASYMORAWETZINTEGRALTERMESTIMATE}
\end{align}
The estimate \eqref{E:MULTLOSSOFONEDERIVERRORINTEGRALTERMESTIMATE} can be proved by using an argument similar to the
one used to prove \eqref{E:MORLOSSOFONEDERIVERRORINTEGRALTERMESTIMATE}.
The estimates \eqref{E:LOWERORDEREASYQ0LARGETIMEDECAYTERMESTIMATE}-\eqref{E:LOWERORDEREASYMORAWETZINTEGRALTERMESTIMATE}
can be proved by using arguments similar to the ones we used to prove
\eqref{E:BELOWTOPORDEREASYQ0LARGETIMEDECAYTERMESTIMATE}-\eqref{E:BELOWTOPORDEREASYMORAWETZINTEGRALTERMESTIMATE}.

In order to complete the proof of \eqref{E:LOWERORDERPREDESIREDDIFFICULTQ1BOUNDINTERMSOFQ0},
it remains for us to bound the terms arising from the terms on the first line of the right-hand side of 
\eqref{E:Q0BELOWTOPORDERGRONWALLREADYINEQUALITY}.
Arguing as in our proof of the bounds for the terms on the first line of the
right-hand side of \eqref{E:Q0TOPORDERGRONWALLREADYINEQUALITY},
we deduce that the two corresponding terms are respectively bounded by
$\leq$ $C\mathring{\upepsilon}^2$ and $C \varepsilon^3.$

Combining all of the above bounds, we arrive at the following analog of 
\eqref{E:MOREEXPLICITPREDESIREDDIFFICULTQ0BOUND}:
\begin{align} \label{E:DESCENTMOREEXPLICITPREDESIREDDIFFICULTQ0BOUND}
	q(t,u) & \leq 
				C (1 + \varsigma^{-1})
				\left\lbrace
					\mathring{\upepsilon}^2
					+ \varepsilon^3
				\right\rbrace
				+ \left\lbrace
						\frac{C}{\Conthree}
						(1 + \varsigma^{-1})
						+ C \varsigma
						+ C \varepsilon
					\right\rbrace
					q(t,u),
\end{align}
and the desired bound \eqref{E:PREDESIREDKEYLOWERORDERLITTLEFBOUND} thus follows
if first $\varsigma$ and $\varepsilon$ are chosen to be sufficiently small
and then $\Conthree$ is chosen to be sufficiently large.

\medskip

\noindent{\emph{Further descent.}}
We now explain how to inductively derive the desired estimates for the remaining lower-order quantities
$\totzeromax{\leq N}(t,u),$
$\totonemax{\leq N}(t,u),$
$\totMormax{\leq N}(t,u),$
$0 \leq N \leq 22$
by descending. At each step in the descent, we define the approximating integrating factors 
\eqref{E:LOWERORDERQ0INTEGRATINGFACTORPRODUCT}-\eqref{E:LOWERORDERQ1INTEGRATINGFACTORPRODUCT},
except we reduce the power of $\upmu_{\star}^{-1}$ by two at each step. Starting at 
the case $N = 15,$ the factor involving $\upmu_{\star}^{-1}$ is absent.
At each step, we prove the estimates 
\eqref{E:MORLOSSOFONEDERIVERRORINTEGRALTERMESTIMATE}-\eqref{E:DESCENTEASYBELOWTOPORDERQ1MORAWETZTERMESTIMATE}
and \eqref{E:MULTLOSSOFONEDERIVERRORINTEGRALTERMESTIMATE}-\eqref{E:LOWERORDEREASYMORAWETZINTEGRALTERMESTIMATE}
(with $24$ replaced by $N+1$ and $23$ replaced by $N$)
using essentially
the same arguments used in the case $N=23.$ One small change is needed starting at $N=15;$
$N=15$ is the first instance in which the $\upmu_{\star}^{-1}$ degeneracy is completely absent. 
Specifically, when proving the desired estimates for $N \leq 15,$
we use the inequalities \eqref{E:LESSSINGULARTERMSMUTHREEFOURTHSINTEGRALBOUND} and \eqref{E:LOGLOSSLESSSINGULARTERMSMTHREEFOURTHSINTEGRALBOUND}
to help prove the corresponding analogs of the
estimates \eqref{E:MORLOSSOFONEDERIVERRORINTEGRALTERMESTIMATE} and \eqref{E:MULTLOSSOFONEDERIVERRORINTEGRALTERMESTIMATE};
inequalities \eqref{E:LESSSINGULARTERMSMUTHREEFOURTHSINTEGRALBOUND} and \eqref{E:LOGLOSSLESSSINGULARTERMSMTHREEFOURTHSINTEGRALBOUND}
are the ones that break the $\upmu_{\star}^{-1}$ degeneracy. 

\hfill $\qed$


\chapter{Local Well-Posedness and Continuation Criteria}
\label{C:LOCALWELLPOSEDNESS}
\thispagestyle{fancy}
In Chapter~\ref{C:LOCALWELLPOSEDNESS},
we sketch the proof of a proposition that
provides local well-posedness and related continuation
criteria for the covariant wave equation $\square_{g(\Psi)} \Psi = 0.$
We state some results in terms of the rectangular coordinates
$(t,x^1,x^2,x^3)$ and others in terms of the geometric coordinates 
$(t,u,\vartheta^1,\vartheta^2).$
The rectangular coordinates are a natural coordinate system for showing that the solution exists 
on some region and thus for initiating the bootstrap argument
that we use in the proof of the sharp classical lifespan theorem
(Theorem~\ref{T:LONGTIMEPLUSESTIMATES}). 
On the other hand, the geometric coordinates are the ones we have used throughout 
the monograph to derive sharp estimates.

\section{Local well-posedness and continuation criteria}
\label{S:LOCALWP}
We now provide the proposition. For convenience, 
we assume the amount of regularity on the data that we use in proving
our sharp classical lifespan theorem; this assumption is highly non-optimal. 

\begin{proposition}[\textbf{Local well-posedness and continuation criteria}] 
\label{P:CLASSICALLOCAL}
Let $N =24,$ let $0 < U_0 < 1$ be a constant,
and let $\Sigma_0^{U_0} = \cup_{u \in [0,U_0]} S_{0,u} \subset \mathbb{R}^3$ be 
the annular region foliated by the level sets $S_{0,u}$ of the function $u = 1 - r$ defined on $\Sigma_0^{U_0}.$ 
Let $(\mathring{\Psi} := \Psi|_{\Sigma_0^{U_0}}, \mathring{\Psi}_0 := \partial_t \Psi|_{\Sigma_0^{U_0}})$ 
be initial data for the covariant wave equation $\square_{g(\Psi)} \Psi = 0$
(that is, equation \eqref{E:WAVEGEO})
that are defined on $\Sigma_0^{U_0}$ and that have vanishing trace on the outer sphere $S_{0,0}.$
Assume that the metric $g(\Psi)$ verifies \eqref{E:LITTLEGDECOMPOSED} and $(g^{-1})^{00} = -1.$
Let $\mathring{\upepsilon} := \| \mathring{\Psi} \|_{H_{\Euct}^{N+1}(\Sigma_0^1)} + \| \mathring{\Psi}_0 \|_{H_{\Euct}^N(\Sigma_0^1)}$
denote the size of the data\footnote{Because we are only studying the influence of the 
nontrivial portion of the data belonging to the subset $\Sigma_0^{U_0}$ of $\Sigma_0^1,$
we could replace the data norms
$\| \cdot \|_{H_{\Euct}^{\cdot}(\Sigma_0^1)}$ with $\| \cdot \|_{H_{\Euct}^{\cdot}(\Sigma_0^{U_0})}$
without altering any of the conclusions of the proposition.} 
as defined in Def.~\ref{D:SMALLDATA},
and let $\mathcal{H}$ be the set of real numbers\footnote{
$\mathcal{H}$ can be viewed as the set of $\Psi$ for which
the metric $g(\Psi)$ is Lorentzian and 
for which the hypersurfaces $\Sigma_t$ are spacelike.
} $b$ such that the following conditions hold:
\begin{itemize}
	\item The rectangular components $g_{\mu \nu}(\cdot),$ $(\mu, \nu = 0,1,2,3),$ are smooth on a neighborhood of $b.$
	\item $g_{00}(b) < 0.$
	\item The eigenvalues of the $3 \times 3$ matrix $\gt_{ij}(b)$ 
	(see Def.~\ref{D:FIRSTFUND}),
	$(i,j = 1,2,3),$ are positive.
\end{itemize}

\medskip

\noindent \textbf{Part I): Statements relative to the rectangular coordinates.}

\medskip

\noindent \underline{\textbf{Local well-posedness.}}
Assume that there is a compact subset $\mathfrak{K} \subset \mbox{interior}(\mathcal{H})$
such that $\mathring{\Psi}(\Sigma_0^{U_0}) \subset \mathfrak{K}.$
If $\mathring{\upepsilon} < \infty,$ 
then these data launch a unique classical solution
$\Psi$ to the equation $\square_g \Psi = 0$
and a unique outgoing eikonal function $u$ that is a classical solution to
$(g^{-1})^{\alpha \beta}(\Psi) \partial_{\alpha} u \partial_{\beta} u = 0,$
that takes on the initial conditions $1 - r$ along $\Sigma_0^{U_0},$
and that verifies $(g^{-1})^{\alpha \beta}(\Psi) \partial_{\alpha} t \partial_{\beta} u < 0$
along $\Sigma_0^{U_0}.$
The solution exists on a nontrivial spacetime region of the form 
$\mathcal{M}_{\Tlocal,U_0}$ (see definition \eqref{E:MTUDEF})
for some $\Tlocal > 0.$
There exists a compact subset 
$\mathfrak{K}'$ such that
$\mathfrak{K} \subset \mathfrak{K}' \subset \mbox{interior}(\mathcal{H})$
and such that on $\mathcal{M}_{\Tlocal,U_0},$
we have $\Psi \in \mathfrak{K}'.$
Furthermore, on $\mathcal{M}_{\Tlocal,U_0},$
we have
$\sum_{a=1}^3 |\partial_a u| > 0,$
the one-form with rectangular components 
$(\partial_1 u, \partial_2 u, \partial_3 u)$ on $\Sigma_0^{U_0}$
is inward-pointing relative to $S_{t,u},$
$(g^{-1})^{\alpha \beta}(\Psi) \partial_{\alpha} t \partial_{\beta} u < 0,$
$0 < \upmu < \infty,$ and each $S_{t,u}$ is an embedded two-dimensional sphere.
In addition, on $\mathcal{M}_{\Tlocal,U_0},$ the scalar-valued functions
	$\upmu$
	and
	$\Lunit_{(Small)}^i,$ $(i=1,2,3),$
	which are defined by 
	\eqref{E:GEOMETRICRADIAL},
	\eqref{E:UPMUDEF}, \eqref{E:LUNIT}, \eqref{E:LUNITJUNK},
	and the geometric angular coordinates 
	$(\vartheta^1,\vartheta^2)$ constructed in Chapter~\ref{C:BASICGEOMETRICCONSTRUCTIONS}
	are $C^{N-2}$ functions of the rectangular coordinates.
	A similar statement holds 
	(with, in some cases, a different degree of differentiability)
	for the rectangular components
	$\Xi^{\mu},$
	$\upchi_{\mu \nu}^{(Small)},$
	$\Lunit^{\mu},$
	$\Radunit_{(Small)}^i,$
	$\Radunit^{\mu},$
	$\Rad^{\mu},$
	$\uLgood^{\mu},$
	and
	$\Rot^{\mu},$
	$(\mu, \nu = 0,1,2,3),$
	and for the rectangular components of 
	all of the other geometric quantities defined throughout the monograph.

The (open-at-the-top) region 
$\mathcal{M}_{\Tlocal,U_0} = 
\lbrace	
	\cup_{u \in [0,U_0]} \mathcal{C}_u^{\Tlocal}
\rbrace 
\cap 
\lbrace
	\cup_{t \in [0, \Tlocal)} \Sigma_t
\rbrace$ 
is foliated by level sets $\mathcal{C}_u^{\Tlocal}$ of the eikonal function $u,$ 
where each $\mathcal{C}_u^{\Tlocal}$ is a truncated null hypersurface of the metric $g(\Psi),$
which itself is foliated by spheres $S_{t,u} = \mathcal{C}_u^{\Tlocal} \cap \Sigma_t^{U_0}.$ 
That is, we have $\mathcal{C}_u^{\Tlocal} = \cup_{t \in [0,\Tlocal]} S_{t,u}.$

The solution has the following regularity properties relative to the rectangular coordinates:
\begin{subequations}
\begin{align} 
	\Psi, \ u & \in C^{N-1}(\mathcal{M}_{\Tlocal,U_0}),
		\\
	\partial^{\vec{I}} \Psi, \ \partial^{\vec{I}} u & \in C([0,\Tlocal),H_{\Euct}^{N+1-|\vec{I}|}(\Sigma_t^{U_0})),
		\qquad |\vec{I}| \leq N+1, 
		\label{E:PSIANDEIKONALREGULARITYRELATIVETORECTANGULAR}
\end{align}
\end{subequations}
where $\vec{I}$ denotes a multi-index corresponding to repeated differentiation with respect to the rectangular 
spacetime coordinate vectorfields $\partial_{\nu},$ 
$(\nu = 0,1,2,3),$
$\Sigma_t^{U_0} = \cup_{u \in [0,U_0]} S_{t,u},$
and $H_{\Euct}^N(\Sigma_t^{U_0})$ is the standard Euclidean Sobolev space involving 
order $\leq N$ rectangular spatial derivatives along $\Sigma_t^{U_0}.$

In addition, the solution depends continuously on the data. Furthermore,
if $\mathring{\upepsilon}$ is sufficiently small, then the
existence time $\Tlocal$ from above can be bounded from below by $f(\mathring{\upepsilon}^{-1}),$
where $f$ is a continuous increasing function such that $f \uparrow \infty$ as $\mathring{\upepsilon} \downarrow 0.$

\medskip

\noindent \underline{\textbf{Regularity of the change of variables map.}}
	The change of variables map 
	$\Upsilon: [0,\Tlocal) \times [0,U_0] \times \mathbb{S}^2 \rightarrow \mathcal{M}_{\Tlocal,U_0}$
	from geometric to rectangular coordinates
	(see Lemma~\ref{L:JACOBIAN}) 
	is a $C^{22}$ diffeomorphism 
	with an everywhere positive Jacobian determinant. 
	Hence, on $\mathcal{M}_{\Tlocal,U_0},$ 
	we have that the scalar-valued functions
	$\Psi,$ 
	$u,$ 
	$\upmu,$
	$\Lunit_{(Small)}^i,$ etc.
	are many times continuously differentiable 
	with respect to the geometric coordinates $(t,u,\vartheta).$
	The same statement holds for the 
	rectangular
	vectorfield components
	$\Lunit^{\mu},$
	$\Radunit^{\mu},$
	$\Rad^{\mu},$
	$\uLgood^{\mu},$
	$\Rot^{\mu},$
	and the rectangular components of all of the other geometric quantities defined throughout the monograph.
	
	\medskip
	
	\noindent \underline{\textbf{Regularity of the geometric norms and energies.}}
	The geometric norm and energy quantities appearing on the left-hand sides of 
	\eqref{E:PSIMAINTHEROEMLINFINIFTYESTIMATES}-\eqref{E:TOPORDERCHIJUNKANGULARMAINTHEROEML2ESTIMATES}
	are well-defined, continuous functions of $(t,u)$ on $[0,\Tlocal) \times [0,U_0].$
	Furthermore, if $\mathring{\upepsilon}$ is sufficiently small, then
	at $(t,u) = (0,U_0),$ these quantities are all bounded by 
	$\lesssim \mathring{\upepsilon}.$

	\medskip

\noindent \textbf{Part II): Continuation criteria.}
Assume that none of the following $4$ breakdown scenarios occur on $\mathcal{M}_{\Tlocal,U_0}:$
\begin{enumerate}
	\item $\inf_{\mathcal{M}_{\Tlocal,U_0}} \upmu = 0.$ 
	\item $\sup_{\mathcal{M}_{\Tlocal,U_0}} \upmu = \infty.$
	\item There exists a sequence $p_n \in \mathcal{M}_{\Tlocal,U_0}$
		such that $\Psi(p_n)$ escapes every compact subset of $\mathcal{H}$ as $n \to \infty.$
	\item $\sup_{\mathcal{M}_{\Tlocal,U_0}} 
				\max_{\kappa=0,1,2,3}
				\left|
					\partial_{\kappa} \Psi
				\right|
		= \infty.$
\end{enumerate}

In addition, assume that the following condition is verified:
\begin{enumerate}
	\setcounter{enumi}{4}
	\item The change of variables map $\Upsilon$
		extends to the compact set $[0,\Tlocal] \times [0,U_0] \times \mathbb{S}^2$
		as a (global) $C^1$ diffeomorphism onto its image.
\end{enumerate}

Then there exists a $\Delta > 0$ such that 
$\Psi,$ 
$u,$ 
$\upmu,$
$\Upsilon,$
$\vartheta^1,$
$\vartheta^2,$
$\Lunit^{\mu},$
$\Lunit_{(Small)}^i,$
$\Xi^{\mu},$
and all of the
other quantities can be extended 
(where $\Psi$ and $u$ are solutions)
to a strictly larger region of the form
$\mathcal{M}_{\Tlocal + \Delta,U_0}$ on which they have all of the properties
stated in \textbf{Part I)}.

\end{proposition}

\begin{remark}[\textbf{Controlling top-order derivatives of $u$}]
\label{R:EIKONALGAINOFONEDERIATIVE}
	The proof of the finiteness of the quantities 
	\eqref{E:TOPORDERLUPMUMAINTHEROEML2ESTIMATES}-\eqref{E:TOPORDERCHIJUNKANGULARMAINTHEROEML2ESTIMATES},
	which play an essential role in our proof of the sharp classical lifespan theorem,
	is highly nontrivial. These quantities involve some $26^{th}$ geometric derivatives of $u,$
	whereas \eqref{E:PSIANDEIKONALREGULARITYRELATIVETORECTANGULAR} 
	yields only the finiteness of the 
	$L^2$ norms of the $25^{th}$
	rectangular derivatives of $u.$ The gain of one derivative 
	is made possible by the following key ingredients:
	\begin{itemize}
		\item The special structure of the right-hand side of 
		the divergence identity \eqref{E:DIVCOMMUTATIONCURRENTDECOMPOSITION}
		and of the deformation tensors 
		for the commutation vectorfields $\mathscr{Z}$
		(see Prop.~\ref{P:DEFORMATIONTENSORFRAMECOMPONENTS}).
		\item The special structure of the
		equations verified by the modified quantities from 
		Ch.~\ref{C:RENORMALIZEDEIKONALFUNCTIONQUANTITIES}.
		\item The availability of the elliptic estimates of Sect.~\ref{S:ELLIPTICESTIMATES}.
	\end{itemize}
\end{remark}

\begin{proof}[Sketch of a proof] \

\noindent \textbf{Discussion of aspects of the proof of Part I) involving rectangular coordinates.}

\medskip

These aspects of the proposition can be proved with rather standard techniques.
For the main ideas behind the proof of local well-posedness 
for the wave equation \eqref{E:WAVEGEO}
relative to the rectangular 
coordinates, readers may consult, for example, \cite{lH1997}*{Ch. VI}. 
After $\Psi$ has been solved for relative to rectangular coordinates, we can
then solve for the eikonal function $u$ and deduce its regularity properties relative 
to the rectangular coordinates. To achieve this, we first use the eikonal equation
$(g^{-1})^{\alpha \beta}(\Psi) \partial_{\alpha} u \partial_{\beta} u = 0$
to solve for $\partial_t u = f(\Psi,\partial_1 u, \partial_2 u, \partial_3 u),$
where $f$ is smooth. We then apply the standard $L^2-$type energy methods to this 
scalar equation, which leads to the existence of $u$ and its properties.
The remaining quantities 
$\upmu,$
$\Lunit_{(Small)}^{\mu},$ 
etc. are constructed out of $\Psi$ and $u.$
Hence, their existence and regularity
properties are easy consequences
of the properties of $\Psi$ and $u.$

\medskip

\noindent \textbf{Discussion of aspects of the proof of Part I) involving geometric coordinates.}

\medskip

We first note that we have already sketched a proof that
$\Psi,$ 
$u,$ 
$\upmu,$
$\Lunit_{(Small)}^{\mu},$ 
etc.
are many times differentiable
with respect to the rectangular coordinates
on $\mathcal{M}_{\Tlocal,U_0}.$
	
	We now show that (shrinking $\Tlocal$ if necessary)
	the change of variables map $\Upsilon$ 
	from geometric to rectangular coordinates
	(see Lemma~\ref{L:JACOBIAN})
	is a bijection from $[0,\Tlocal) \times [0,U_0] \times \mathbb{S}^2$ to $\mathcal{M}_{\Tlocal,U_0}.$
	Since $t$ is the same function in both coordinate systems, 
	we have to show only that for any $t \in [0,\Tlocal),$
	$\Upsilon(t,\cdot): [0,U_0] \times \mathbb{S}^2 \rightarrow \Sigma_t^{U_0}$ 
 	is bijective. The main point is that
 	the estimates proved relative to rectangular coordinates guarantee
	that the rectangular components $\Lunit^{\nu} = \Lunit x^{\nu}$
	are $C^{21}$ functions of the rectangular coordinates on 
	$\mathcal{M}_{t,U_0}.$ 
	Hence, since $\Lunit t = 1,$ it follows from the standard theory of ODEs
	that the restriction $\varphi_t|_{\Sigma_0^{U_0}}$ to $\Sigma_0^{U_0}$
	of the flow map $\varphi_t$ of $\Lunit$
	(see the proof of Lemma~\ref{L:STUDERIVATIVES})
	is, relative to the rectangular spatial coordinates, a bijection from $\Sigma_0^{U_0}$ to $\Sigma_t^{U_0}.$
	In view of the manner in which we constructed
	the geometric coordinates in
	Chapter~\ref{C:BASICGEOMETRICCONSTRUCTIONS}, 
	we see that
	$\Upsilon(t,u,\vartheta)$ is equal to $\varphi_t|_{\Sigma_0^{U_0}}$
	pre-composed with the smooth diffeomorphism $[0,U_0] \times \mathbb{S}^2 \rightarrow \Sigma_0^{U_0}$
	that maps the geometric coordinates $(u,\vartheta) \in [0,U_0] \times \mathbb{S}^2$
	along $\Sigma_0^{U_0}$ to the rectangular spatial coordinates $(x^1,x^2,x^3)$ along $\Sigma_0^{U_0}.$
	We have thus shown the bijectivity of $\Upsilon.$
	
	Furthermore, examining the proof of Lemma~\ref{L:JACOBIAN},
	we see that at $t = 0,$ the Jacobian of $\Upsilon^{-1}$ 
	is a $C^{21}$ matrix
	$\frac{\partial (t,u,\vartheta^1, \vartheta^2)}{\partial (x^0,x^1,x^2,x^3)}$
	of the form 
	$\begin{pmatrix}
  1 & 0 \\
  * & M
 	\end{pmatrix}^{-1},$
 	where $M$ is an invertible $3 \times 3$ matrix.
	Hence, at least for short times, the Jacobian 
	of $\Upsilon^{-1}$ remains invertible, and by
	the inverse function theorem, 
	$\Upsilon$ has the same regularity as $\Upsilon^{-1}.$
	We remark that the 
	$\Xi^i$ are the terms in the Jacobian (see equation \eqref{E:CHOVMATRIXCOMPUTATIONS}) 
	with the least regularity (that is, $C^{21}$);
	their regularity can be derived
	with the help of the transport equation \eqref{E:XIEVOLUTION},
	expressed relative to the rectangular coordinates.
	
	The statements concerning the finiteness and the $(t,u)-$continuity of the geometric norms and energies
	\eqref{E:PSIMAINTHEROEMLINFINIFTYESTIMATES}-\eqref{E:TOPORDERCHIJUNKANGULARMAINTHEROEML2ESTIMATES}	
	are difficult to derive. However, the only difficult step is obtaining a priori estimates for these
	quantities and in particular, avoiding derivative loss in some of the top-order eikonal function quantities 
	(see Remark~\ref{R:EIKONALGAINOFONEDERIATIVE}).
	In our proof of Theorem~\ref{T:LONGTIMEPLUSESTIMATES}, we show how to derive such
	a priori estimates and hence we do not repeat the lengthy argument here.
	To bound the geometric norms and energies at $(t,u) = (0,U_0)$ by
	$\lesssim \mathring{\upepsilon},$ we use the small-data estimates
	derived in Sect.~\ref{S:INITIALBEHAVIOROFQUANTITIES}.	

	\medskip
	
	\noindent \textbf{Discussion of the proof of Part II).}
	
	\medskip
	
The continuation principle of \textbf{Part II)} can also be proved
by using mostly standard arguments; see, for example, \cite{jS2008b} for the main ideas behind a proof.
The main idea is that once we rule out the breakdown in the hyperbolic character of the equations
and the possible blow-up of first rectangular derivatives of various quantities, 
we can then ensure that the solutions $\Psi$ and $u$ can be continued relative to rectangular coordinates.
Furthermore, assuming that the condition $(5)$ holds
and using the extendibility of the solution relative to the rectangular coordinates,
it is straightforward to show that we can extend the change of variables map $\Upsilon$ 
to a strictly larger set $[0,\Tlocal + \Delta) \times [0,U_0] \times \mathbb{S}^2$
(for some $\Delta > 0$)
as a diffeomorphism onto its image. We can then derive
the desired properties, on the domain
$[0,\Tlocal + \Delta) \times [0,U_0] \times \mathbb{S}^2,$
for all of the quantities of interest relative to the geometric
coordinates by using the same arguments as in the proof of 
local well-posedness.

We now provide the proof of the one somewhat subtle aspect, 
which is showing, under the assumption that none of the $4$ breakdown scenarios occur,
that $u$ remains regular in the following sense:
\begin{align} \label{E:EIKONALFUNCTIONRECTANGULARDERIVATIVESUPPERANDLOWERBOUND}
	0
	& 
	< \inf_{\mathcal{M}_{\Tlocal,U_0}} \sum_{a=1}^3 |\partial_a u| 
	\leq \sup_{\mathcal{M}_{\Tlocal,U_0}} \sum_{a=1}^3 |\partial_a u|
	< \infty.
\end{align}
Assuming that none of the $4$ breakdown scenarios occur,
we have in particular that the $3 \times 3$ 
matrices $\gt_{ij}(\Psi)$ 
and $(\gtinverse)^{ij}(\Psi)$
are uniformly positive definite  
on $\mathcal{M}_{\Tlocal,U_0}$
(their eigenvalues are bounded from above and uniformly from below away from $0$). 
Hence, the inequalities in \eqref{E:EIKONALFUNCTIONRECTANGULARDERIVATIVESUPPERANDLOWERBOUND}
follow from the identity $(\gt^{-1})^{ab} \partial_a u \partial_b u = \upmu^{-2}$
(see \eqref{E:ALTERNATEEIKONAL}).

	
\end{proof}

\chapter{The Sharp Classical Lifespan Theorem} \label{C:LONGTIMEEXISTENCE}
\thispagestyle{fancy}
In Chapter~\ref{C:LONGTIMEEXISTENCE}, 
we state and prove our sharp classical lifespan theorem, which is the main theorem
of the monograph. It guarantees that the
solution persists unless $\upmu_{\star}$ becomes $0$ in finite time, 
in which case some rectangular derivatives of $\Psi$ blow-up and a shock singularity 
has started to form. We also prove Cor.~\ref{C:ANGULARDDERIVATIVEESTIMATES},
which shows that if the data have ``very small'' angular derivatives, then 
this property is propagated by the solution. We use the theorem and the corollary
in Chapter~\ref{C:PROOFOFSHOCKFORMATION}, when we prove finite-time shock
formation for an open set of nearly spherically symmetric small data.

\section{The sharp classical lifespan theorem}
\label{S:SHARPCLASSICALLIFESPAN}

We now state and prove the main theorem of the monograph.

\begin{theorem}[\textbf{The sharp classical lifespan theorem together with estimates}]
\label{T:LONGTIMEPLUSESTIMATES} 
Let $(\mathring{\Psi} := \Psi|_{\Sigma_0}, \mathring{\Psi}_0 := \partial_t \Psi|_{\Sigma_0})$ 
be initial data for the covariant wave equation \eqref{E:WAVEGEO} 
under the assumption\footnote{As we described in Chapter~\ref{C:INTRO}, 
this assumption is easy to eliminate.} \eqref{E:GINVERSE00ISMINUSONE}.
Assume that the data are supported in the 
Euclidean unit ball $\Sigma_0^1.$
Let $\mathring{\upepsilon} = \| \mathring{\Psi} \|_{H_{\Euct}^{25}(\Sigma_0^1)} + \| \mathring{\Psi}_0 \|_{H_{\Euct}^{24}(\Sigma_0^1)}$ 
be the size of the data as defined in Def.~\ref{D:SMALLDATA}.
Assume that the data verify the hypotheses of Prop.~\ref{P:CLASSICALLOCAL}
(the local well-posedness proposition), and let
$0 < U_0 < 1$ be a fixed constant.
Let $\Psi$ denote the solution corresponding to the data
existing on a nontrivial region of the form $\mathcal{M}_{\Tlocal,U_0}$
(see definition \eqref{E:MTUDEF})
Recall that $\upmu_{\star}(t,u) := \min\lbrace 1, \min_{\Sigma_t^u} \upmu \rbrace,$
$t$ denotes the Minkowski time coordinate, $u$ is the eikonal function
(with initial data $u|_{\Sigma_0} = 1 - r,$ where $r = \sqrt{\sum_{a=1}^3 (x^a)^2}$),
and $\rgeo(t,u) := 1 - u + t$ is the geometric radial coordinate. 
There exist large constants $C > 0,$ $C_{(Lower-Bound)} > 0,$ and $\Cononestar > 4$ 
and a small constant $\upepsilon_0 > 0$ such that if 
$\mathring{\upepsilon} < \upepsilon_0,$
then the following statements hold true.
In the statements, the constants can depend on $U_0$
and the nonlinearities in \eqref{E:WAVEGEO} 
and in particular, $C$ can blow-up as $U_0 \uparrow 1.$ 

\medskip

\noindent \underline{\textbf{Existence as long as $\upmu_{\star} > 0$ and a classical lifespan lower bound.}}
Let  
\begin{align} \label{E:TLIFESPAN}
	T_{(Lifespan);U_0}:= \sup \lbrace t \ | \ \inf_{s \in [0,t)} \upmu_{\star}(s,U_0) > 0 \rbrace.
\end{align}
Then $T_{(Lifespan);U_0}$ is the classical lifespan of the solution
in the region determined by the portion of the data in $\Sigma_0$
belonging to the exterior of the Euclidean sphere $S_{0,U_0}$ of radius $1-U_0$
(see \eqref{E:EIKONALFUNCTIONINITIALCONDITIONS} and \eqref{E:STU}).
That is, $\Psi$ can be extended 
as a classical solution (relative to both the geometric and the rectangular coordinates)
to the region $\mathcal{M}_{T_{(Lifespan);U_0},U_0}$
on which it has all of the properties stated in 
Prop.~\ref{P:CLASSICALLOCAL}.
Furthermore, if $T_{(Lifespan);U_0}< \infty,$ then
\begin{align} \label{E:RECTANGULARDERIVATIVESBLOWUP}
	\sup_{\mathcal{M}_{T_{(Lifespan);U_0},U_0}} \max_{\nu=0,1,2,3} |\partial_{\nu} \Psi| = \infty.
\end{align}

In addition, with $\mathring{\upmu}(u,\vartheta) = \upmu(0,u,\vartheta),$
$G_{\Lunit \Lunit} = \frac{d}{d \Psi} g_{\alpha \beta}(\Psi) \Lunit^{\alpha} \Lunit^{\beta},$
and $\InitialFutFailFac(\vartheta)$ as in Def.~\ref{D:FAILUREFACTOR},
we have the following estimates for $0 \leq s \leq t < T_{(Lifespan);U_0}:$
\begin{subequations}
\begin{align} \label{E:MAINTHEOREMKEYUPMUBEHAVIOR}
		\left|
		\upmu(s,u,\vartheta)
		-
		\left\lbrace	 
			\mathring{\upmu}(u,\vartheta)
			+ \frac{1}{2} 
				\ln \left(\frac{\rgeo(s,u)}{\rgeo(0,u)} \right)
				[\rgeo G_{\Lunit \Lunit} \Rad \Psi](t,u,\vartheta)
		\right\rbrace
		\right|	
		& \leq C \mathring{\upepsilon},	
			\\
	\left|
		G_{\Lunit \Lunit}(t,u,\vartheta)
		- \InitialFutFailFac(\vartheta)
	\right|
	& \leq C \mathring{\upepsilon}
		\label{E:MAINTHEOREMGLLNEARNULLCONDFAILFACT},
		\\
	\left|
		\mathring{\upmu}(u,\vartheta) 
		- 1
	\right|
	& \leq C \mathring{\upepsilon},
		\label{E:MUISNEARONEINITIALLYMAINTHEOREMSTATEMENT} \\
	\left|
		\rgeo \Rad \Psi
	\right|(t,u,\vartheta)
	& \leq C \mathring{\upepsilon}.
	\label{E:RGEOTIMESRADPSIMAINTHEOREMSTATEMENT}
\end{align}
\end{subequations}

Furthermore, 
$\upmu_{\star}(t,U_0) > 0$ whenever $t < \exp\left( \left\lbrace C_{(Lower-Bound)} \mathring{\upepsilon} \right\rbrace^{-1} \right),$ 
and hence
\begin{align}  \label{E:TLIFESPANLOWERBOUND}
	T_{(Lifespan);U_0}> \exp\left(\frac{1}{C_{(Lower-Bound)}\mathring{\upepsilon}} \right).
\end{align}

\medskip

\noindent \underline{\textbf{$\upmu$ remains large along outgoing null curves with $G_{\Lunit \Lunit}(t,u,\vartheta) = 0$.}}
For $(t,u) \in [0,T_{(Lifespan);U_0}) \times [0,U_0],$ the
following pointwise estimate holds along the portions of the integral curves of the 
outgoing null vectorfield $\Lunit$
that terminate at the (possibly empty) set of points in $\Sigma_t^{U_0}$
with geometric coordinates $(t,u,\vartheta)$ verifying $G_{\Lunit \Lunit}(t,u,\vartheta)=0:$
\begin{align} \label{E:UPMUREMAINSLARGEWHENNULLCONDITIONFAILUREFACTORVANISHES}
	\upmu(t,u,\vartheta) \geq 1 - C \mathring{\upepsilon}.
\end{align}

\medskip

\noindent \underline{\textbf{The lower-order $\mathscr{Z}$ derivatives of $\Psi$ remain regular relative to $(t,u,\vartheta)$.}}
The following $C^0$ estimates\footnote{See Remark~\ref{R:C0NORMCONVENTION} concerning our use of the norm
$\| \cdot \|_{C^0(\Sigma_t^u)}.$} for 
the lower-order derivatives of $\Psi$ hold on the domain 
$(t,u) \in [0,T_{(Lifespan);U_0}) \times [0,U_0]:$
	\begin{align} \label{E:PSIMAINTHEROEMLINFINIFTYESTIMATES}
		\| \mathscr{Z}^N \Psi \|_{C^0(\Sigma_t^u)}
		& \leq C \mathring{\upepsilon} \frac{1}{1 + t},
		&& (N \leq 13).
	\end{align}
In \eqref{E:PSIMAINTHEROEMLINFINIFTYESTIMATES}, $\mathscr{Z}^N$ denotes an arbitrary 
$N^{th}$ order differential operator corresponding to repeated differentiation with respect to
commutation vectorfields belonging to the set $\mathscr{Z}$ defined in
\eqref{E:DEFSETOFCOMMUTATORVECTORFIELDS}.
In particular, 
if $T_{(Lifespan);U_0}< \infty,$
then the quantities $\mathscr{Z}^{\leq 12} \Psi$ extend as continuous functions
of $(t,u,\vartheta)$ to $\Sigma_{T_{(Lifespan);U_0}}^{U_0}.$

\medskip

\noindent \underline{\textbf{The lower-order $\mathscr{Z}$ derivatives of the rectangular metric components 
remain regular relative to $(t,u,\vartheta)$.}}
The following $C^0$ estimates for 
the lower-order derivatives of the rectangular components 
$g_{\mu \nu},$ $(\mu,\nu=0,1,2,3),$ of the spacetime metric
and the rectangular spatial components $\gsphere_{ij},$
$(i,j=1,2,3),$ of the metric induced by $g$ on $S_{t,u}$ 
(see Def.~\ref{D:FIRSTFUND})
hold on the domain $(t,u) \in [0,T_{(Lifespan);U_0}) \times [0,U_0]:$
	\begin{subequations}
	\begin{align} \label{E:RECTANGULARSPACETIMEMETRICCOMPONENTSMAINTHEROEMLINFINIFTYESTIMATES}
		\left \| \mathscr{Z}^N \left\lbrace g_{\mu \nu} - m_{\mu \nu} \right\rbrace \right\|_{C^0(\Sigma_t^u)}
		& \leq C \mathring{\upepsilon} \frac{1}{1 + t},
		&& (N \leq 13),
			\\
		\left\| 
			\mathscr{Z}^N 
					\left\lbrace
						\gsphere_{ij} 
						- \left(\delta_{ij} - \frac{x^i x^j}{\rgeo^2} \right) 
					\right\rbrace
		\right\|_{C^0(\Sigma_t^u)}
		& \leq C \mathring{\upepsilon} \frac{\ln(\myexp + t)}{1 + t},
		&& (N \leq 12).
			\label{E:RECTANGULARSPHEREMETRICCOMPONENTSMAINTHEROEMLINFINIFTYESTIMATES}
	\end{align}
	\end{subequations}
In \eqref{E:RECTANGULARSPACETIMEMETRICCOMPONENTSMAINTHEROEMLINFINIFTYESTIMATES},
$m_{\mu \nu} = \mbox{diag}(-1,1,1,1)$ denotes the Minkowski metric.
In particular, 
if $T_{(Lifespan);U_0}< \infty,$
then the quantities $\mathscr{Z}^{\leq 12} \left\lbrace g_{\mu \nu} - m_{\mu \nu} \right\rbrace$ 
and 
$\mathscr{Z}^{\leq 11} \left\lbrace
						\gsphere_{ij} 
						- \left(\delta_{ij} - \frac{x^i x^j}{\rgeo^2} \right) 
					\right\rbrace$
extend as continuous functions
of $(t,u,\vartheta)$ to $\Sigma_{T_{(Lifespan);U_0}}^{U_0}.$

\medskip

\noindent \underline{\textbf{The lower-order $\mathscr{Z}$ derivatives of the eikonal function quantities 
remain regular relative to $(t,u,\vartheta)$.}}
For $(t,u) \in [0,T_{(Lifespan);U_0}) \times [0,U_0],$ 
the following $C^0$ estimates hold
for the inverse foliation density $\upmu,$ the rectangular components 
$\Lunit_{(Small)}^i = \Lunit^i - \frac{x^i}{\rgeo}$
and
$\Radunit_{(Small)}^i = \Radunit^i + \frac{x^i}{\rgeo},$
and the $S_{t,u}$ tensorfield $\upchi^{(Small)} = \chi - \frac{\gsphere}{\rgeo}:$
\begin{subequations}
\begin{align}
	\| \mathscr{Z}^N (\upmu - 1) \|_{C^0(\Sigma_t^u)}
		& \leq C \mathring{\upepsilon} \ln(\myexp + t),  
			&& (N \leq 12),
		\label{E:UPMUMAINTHEROEMLINFINIFTYESTIMATES} \\	
	\| \mathscr{Z}^N \Lunit_{(Small)}^i \|_{C^0(\Sigma_t^u)},
		\, \| \mathscr{Z}^N \Radunit_{(Small)}^i \|_{C^0(\Sigma_t^u)} 	
		& \leq C \mathring{\upepsilon} \frac{\ln(\myexp + t)}{1 + t},  
			&& (N \leq 12), 
		\label{E:LUNITJUNKMAINTHEROEMLINFINIFTYESTIMATES} \\
	\| \angLie_{\mathscr{Z}}^N \upchi^{(Small)} \|_{C^0(\Sigma_t^u)}
		& \leq C \mathring{\upepsilon} \frac{\ln(\myexp + t)}{(1 + t)^2},  
			&& (N \leq 11).
			\label{E:TRCHIJUNKMAINTHEROEMLINFINIFTYESTIMATES}
\end{align}
\end{subequations}
In particular, 
if $T_{(Lifespan);U_0}< \infty,$ then
the quantities $\mathscr{Z}^{\leq 11} \upmu,$ 
$\mathscr{Z}^{\leq 11} \Lunit_{(Small)}^i,$
$\mathscr{Z}^{\leq 11} \Radunit_{(Small)}^i,$
and 
$\angLie_{\mathscr{Z}}^{\leq 10} \upchi^{(Small)}$
extend as continuous functions
of $(t,u,\vartheta)$ to $\Sigma_{T_{(Lifespan);U_0}}^{U_0}.$

\medskip

\noindent \underline{\textbf{Behavior of the change of variables map $\Upsilon:$}}
If $T_{(Lifespan);U_0} < \infty,$
then the change of variables map (from geometric to rectangular coordinates)
$\Upsilon: [0,T_{(Lifespan);U_0}) \times [0,U_0] \times \mathbb{S}^2 \rightarrow \mathcal{M}_{T_{(Lifespan);U_0},U_0},$
$\Upsilon(t,u,\vartheta) = (t,x^1,x^2,x^3),$
extends as a $C^{10}$ function
to $[0,T_{(Lifespan);U_0}] \times [0,U_0] \times \mathbb{S}^2.$
In addition, 
$\Upsilon$ is a bijection
from\footnote{Recall that $\mathcal{M}_{\Tlocal,U_0}$ is, by definition, ``open at the top.''} $[0,\Tlocal) \times [0,U_0] \times \mathbb{S}^2$ to $\mathcal{M}_{\Tlocal,U_0}$
with a positive Jacobian determinant.
Furthermore, if $T_{(Lifespan);U_0}< \infty,$
then on $[0,T_{(Lifespan);U_0}] \times [0,U_0] \times \mathbb{S}^2,$
its Jacobian determinant vanishes precisely on the subset 
\begin{align} \label{E:MUVANISHESSET}
	\lbrace (T_{(Lifespan);U_0},u,\vartheta) \ | \ \upmu(T_{(Lifespan);U_0},u,\vartheta) = 0 \rbrace.
\end{align}

\medskip

\noindent \underline{\textbf{What happens when $\upmu \to 0$}:}
In the subset of $\mathcal{M}_{T_{(Lifespan);U_0},U_0}$ 
such that the following two conditions hold:
\textbf{i)} $\upmu \leq 1/4$ and 
\textbf{ii)} $G_{\Lunit \Lunit} \neq 0$ 
we have
\begin{align} 
	\Lunit \upmu(t,u,\vartheta) 
	& \leq
		- \frac{2}{\rgeo(t,u) \left\lbrace 1 + \ln \left(\frac{\rgeo(t,u)}{\rgeo(0,u)} \right) \right \rbrace},
	  \label{E:MAINTHEOREMSMALLMUIMPLIESLMUISNEGATIVE} \\
	|\Radunit \Psi|(t,u,\vartheta) 
	& \geq \frac{1}{\upmu(t,u,\vartheta)} \frac{1}{\rgeo(t,u) \left\lbrace 1 + \ln \left(\frac{\rgeo(t,u)}{\rgeo(0,u)} \right) \right\rbrace} 
		\frac{1}{\left| G_{\Lunit \Lunit}(t,u,\vartheta) \right|},
	\label{E:RADUNITPSIBLOWSUP}
\end{align}
where the vectorfield $\Radunit$ verifies the Euclidean estimate 
$\left| \Radunit - (- \partial_r) \right|_{\Euct} \lesssim \mathring{\upepsilon} \ln(\myexp + t) (1 + t)^{-1}.$
Here, $|V|_{\Euct}^2 := \delta_{ab} V^a V^b$ 
and $\partial_r$ is the standard Euclidean radial derivative.

If $T_{(Lifespan);U_0} < \infty,$ then let
\begin{align} \label{E:SIGMATLIFESPANBLOWUP}
	\Sigma_{T_{(Lifespan);U_0};(Blow-up)}^{U_0}
	& := \lbrace
				(T_{(Lifespan);U_0}, u, \vartheta)
				\ | \
				\upmu(T_{(Lifespan);U_0}, u, \vartheta)
				= 0
			\rbrace.
\end{align}
In particular, since \eqref{E:UPMUREMAINSLARGEWHENNULLCONDITIONFAILUREFACTORVANISHES} shows that
no point $p$ with $G_{\Lunit \Lunit}(p) = 0$ can belong to $\Sigma_{T_{(Lifespan);U_0};(Blow-up)}^{U_0},$
it follows from \eqref{E:RADUNITPSIBLOWSUP} that 
at any point in $\Sigma_{T_{(Lifespan);U_0};(Blow-up)}^{U_0},$
\textbf{the near-Euclidean-unit-length derivative $\Radunit \Psi$ blows up.}
\medskip

\noindent \underline{\textbf{A hierarchy of $L^2$ estimates for $\Psi$ with no $\upmu_{\star}^{-1}$ degeneracy at the lower orders.}}
The following estimates hold for the $L^2-$based quantities 
$\totzeromax{N}$ and $\totonemax{N}$ defined in Def.~\ref{D:MAINCOERCIVEQUANT} 
and the spacetime Morawetz integral $\totMormax{N}$
defined in Def.~\ref{D:COERCIVEMORAWETZINTEGRAL}
on the domain $(t,u) \in [0,T_{(Lifespan);U_0}) \times [0,U_0]:$
\begin{subequations}
\begin{align}
	\totzeromax{N}^{1/2}(t,u)
	& \leq C \mathring{\upepsilon}, 
		&& (0 \leq N \leq 15), 
		\label{E:Q0LOWESTORDERIMPROVED} \\  
	\totonemax{N}^{1/2}(t,u)
		+ \totMormax{N}^{1/2}(t,u)
	& \leq C \mathring{\upepsilon}
	 	\ln^2(\myexp + t), 
		&& (0 \leq N \leq 15),
		\label{E:Q1LOWESTORDERIMPROVED} \\
	\totzeromax{16 + M}^{1/2}(t,u)
	& \leq C \mathring{\upepsilon}
				\upmu_{\star}^{- .75 - M}(t,u), 
		&& (0 \leq M \leq 7), 
		\label{E:Q0MIDORDERIMPROVED} \\  
	\totonemax{16 + M}^{1/2}(t,u)
		+ \totMormax{16 + M}^{1/2}(t,u)
	& \leq C \mathring{\upepsilon}
			\ln^2 (\myexp + t) 
			\upmu_{\star}^{-.75 - M}(t,u), 
		&& (0 \leq M \leq 7),
		\label{E:Q1MIDORDERIMPROVED} \\
		\totzeromax{24}^{1/2}(t,u)
	& \leq C \mathring{\upepsilon}
				\ln^{\Cononestar} (\myexp + t) 
				\upmu_{\star}^{-8.75}(t,u), 
		&&  
		\label{E:Q0TOPORDERIMPROVED} \\  
	\totonemax{24}^{1/2}(t,u)
		+ \totMormax{24}^{1/2}(t,u)
	& \leq C \mathring{\upepsilon}
				\ln^{\Cononestar + 2} (\myexp + t) \upmu_{\star}^{-8.75}(t,u).
				\label{E:Q1TOPORDERIMPROVED}
\end{align}
\end{subequations}

\medskip

\noindent \underline{\textbf{A hierarchy of $L^2$ estimates for the eikonal function quantities with no $\upmu_{\star}^{-1}$ 
degeneracy at the lower orders.}}
The following $L^2$ estimates hold for the eikonal function quantities on the domain 
$(t,u) \in [0,T_{(Lifespan);U_0}) \times [0,U_0]:$
\begin{subequations}
	\begin{align}
		\| \mathscr{Z}^N(\upmu - 1) \|_{L^2(\Sigma_t^u)}
		& \leq C \mathring{\upepsilon} (1 + t) \ln(\myexp + t), 
			&& (0 \leq N \leq 15),
			 \label{E:NODEGENERACYCUPMUMAINTHEROEML2ESTIMATES} 
			 	\\
		\| \mathscr{Z}^N \Lunit_{(Small)}^i \|_{L^2(\Sigma_t^u)}
		& \leq C \mathring{\upepsilon} \ln(\myexp + t), 
			&& (0 \leq N \leq 15),
			 \label{E:NODEGENERACYLJUNKMAINTHEROEML2ESTIMATES} \\
		\| \angLie_{\mathscr{Z}}^N \upchi^{(Small)} \|_{L^2(\Sigma_t^u)}
		& \leq C \mathring{\upepsilon} \frac{\ln(\myexp + t)}{1 + t}, 
			&& (0 \leq N \leq 14),
			\label{E:NODEGENERACYCHIJUNKMAINTHEROEML2ESTIMATES} 
	\end{align}
	
	\begin{align}
	\| \mathscr{Z}^{16 + M} \upmu \|_{L^2(\Sigma_t^u)}
		& \leq C \mathring{\upepsilon} (1 + t) \ln^3(\myexp + t)\upmu_{\star}^{-.25 - M}(t,u), 
			&& (0 \leq M \leq 7),
		 \label{E:MIDORDERUPMUMAINTHEROEML2ESTIMATES} \\
	\| \mathscr{Z}^{16 + M} \Lunit_{(Small)}^i \|_{L^2(\Sigma_t^u)}
		& \leq C \mathring{\upepsilon} \ln^3(\myexp + t)\upmu_{\star}^{-.25 - M}(t,u), 
			&& (0 \leq M \leq 7),
			\label{E:MIDORDERLJUNKMAINTHEROEML2ESTIMATES} \\
	\| \angLie_{\mathscr{Z}}^{15 + M} \upchi^{(Small)} \|_{L^2(\Sigma_t^u)}
		& \leq C \mathring{\upepsilon} \frac{\ln^3(\myexp + t)}{1 + t} \upmu_{\star}^{-.25 - M}(t,u) , 
			\label{E:MIDORDERCHIJUNKMAINTHEROEML2ESTIMATES}  && (0 \leq M \leq 7),
	\end{align}
	
	\begin{align}
	\| \mathscr{Z}^{24} \upmu \|_{L^2(\Sigma_t^u)}
		& \leq C \mathring{\upepsilon} (1 + t) \ln^{\Cononestar + 3}(\myexp + t)\upmu_{\star}^{-8.25}(t,u) , 
			&& 
			\label{E:JUSTBELOWTOPORDERUPMUMAINTHEROEML2ESTIMATES} \\
	\| \mathscr{Z}^{24} \Lunit_{(Small)}^i \|_{L^2(\Sigma_t^u)}
		& \leq C \mathring{\upepsilon} \ln^{\Cononestar + 3}(\myexp + t)\upmu_{\star}^{-8.25}(t,u), 
			&& 
				\label{E:JUSTBELOWTOPORDERLJUNKMAINTHEROEML2ESTIMATES} \\
		\| \angLie_{\mathscr{Z}}^{23} \upchi^{(Small)} \|_{L^2(\Sigma_t^u)}
		& \leq C \mathring{\upepsilon} \frac{\ln^{\Cononestar + 3}(\myexp + t)}{1 + t} \upmu_{\star}^{-8.25}(t,u), 
			&& 
			\label{E:JUSTBELOWTOPORDERCHIJUNKMAINTHEROEML2ESTIMATES}  \\
		\| \Lunit \mathscr{Z}^{24} \upmu \|_{L^2(\Sigma_t^u)}
		& \leq C \mathring{\upepsilon} \ln^{\Cononestar + 2} (\myexp + t) \upmu_{\star}^{-9.25}(t,u) , 
			&& 
			\label{E:TOPORDERLUPMUMAINTHEROEML2ESTIMATES} \\
	\| \Lunit\left[ \rgeo \mathscr{Z}^{24} \Lunit_{(Small)}^i \right] \|_{L^2(\Sigma_t^u)}
		& \leq C \mathring{\upepsilon} \ln^{\Cononestar + 2} (\myexp + t) \upmu_{\star}^{-9.25}(t,u), 
			&& 
			\label{E:TOPORDERLLJUNKMAINTHEROEML2ESTIMATES} \\
		\| \angLie_{\Lunit} \left[\rgeo^2 \angLie_{\mathscr{Z}}^{23} \upchi^{(Small)} \right] \|_{L^2(\Sigma_t^u)}
		& \leq C \mathring{\upepsilon} \ln^{\Cononestar + 2} (\myexp + t) \upmu_{\star}^{-9.25}(t,u), 
			&& 
		\label{E:TOPORDERLCHIJUNKMAINTHEROEML2ESTIMATES} \\
		\| \upmu \angLie_{\Rad} \angLie_{\mathscr{Z}}^{23} \upchi^{(Small)} \|_{L^2(\Sigma_t^u)}
		& \leq C \mathring{\upepsilon} \ln^{\Cononestar + 2} (\myexp + t) \upmu_{\star}^{-8.75}(t,u), 
			&& 
			\label{E:TOPORDERRADCHIJUNKMAINTHEROEML2ESTIMATES} \\
		\| \upmu \angD^2 \mathscr{Z}^{23} \upmu \|_{L^2(\Sigma_t^u)}
		& \leq C \mathring{\upepsilon} \ln^{\Cononestar}(\myexp + t) \upmu_{\star}^{-8.75}(t,u), 
			\label{E:TOPORDERANGLAPUPMUMAINTHEROEML2ESTIMATES} \\
		\| \upmu \angD^2 \mathscr{Z}^{23} \Lunit_{(Small)}^i \|_{L^2(\Sigma_t^u)}
		& \leq C \mathring{\upepsilon} \frac{\ln^{\Cononestar}(\myexp + t)}{1 + t} \upmu_{\star}^{-8.75}(t,u), 
			 \label{E:TOPORDERANGLAPLJUNKMAINTHEROEML2ESTIMATES} \\
		\| \upmu \angD \angLie_{\mathscr{Z}}^{23} \upchi^{(Small)} \|_{L^2(\Sigma_t^u)}
		& \leq C \mathring{\upepsilon} \frac{\ln^{\Cononestar}(\myexp + t)}{1 + t} \upmu_{\star}^{-8.75}(t,u).
			\label{E:TOPORDERCHIJUNKANGULARMAINTHEROEML2ESTIMATES}
	\end{align}
	\end{subequations}				
\end{theorem}

\begin{remark}[\textbf{Global existence when the classic null condition holds}]
	When Klainerman's classic null condition \cite{sk1984} is verified,
	that is, when $\FutFailFac \equiv 0,$ it is straightforward to show that 
	the estimate \eqref{E:MAINTHEOREMGLLNEARNULLCONDFAILFACT} can be improved to
	$\left| G_{\Lunit \Lunit} \right|(t,u,\vartheta) \leq C \mathring{\upepsilon} \ln(\myexp + t) (1 + t)^{-1}.$
	Consequently, Theorem~\ref{T:LONGTIMEPLUSESTIMATES} implies, 
	in particular via the estimates \eqref{E:MAINTHEOREMKEYUPMUBEHAVIOR}
	and \eqref{E:RGEOTIMESRADPSIMAINTHEOREMSTATEMENT},
	that $\upmu_{\star}$ never vanishes and hence the solution exists globally in the spacetime region
	bounded by the outgoing null cone $\mathcal{C}_{U_0}$ and the flat outgoing null cone
	$\mathcal{C}_0.$
\end{remark}

\begin{remark}[\textbf{Extensions to data with additional regularity}]
	Theorem~\ref{T:LONGTIMEPLUSESTIMATES} can of course be extended to apply to data of higher
	Sobolev regularity, that is, with the smallness of
	$\| \mathring{\Psi} \|_{H_{\Euct}^{25}(\Sigma_0^1)} + \| \mathring{\Psi}_0 \|_{H_{\Euct}^{24}(\Sigma_0^1)}$
	replaced by the smallness of
	$\| \mathring{\Psi} \|_{H^{N_{(Top)}+1}(\Sigma_0^1)} + \| \mathring{\Psi}_0 \|_{H^{N_{(Top)}}(\Sigma_0^1)},$
	where $N_{(Top)} \geq 25$ is an integer.
	In this case, the estimates 
	\eqref{E:Q0TOPORDERIMPROVED}-\eqref{E:Q1TOPORDERIMPROVED}
	hold with $24$ replaced by $N_{(Top)},$
	the estimates \eqref{E:Q0MIDORDERIMPROVED}-\eqref{E:Q1MIDORDERIMPROVED}
	hold with $16+M$ replaced by $N_{(Top)}-8+M,$
	the estimates \eqref{E:Q0LOWESTORDERIMPROVED}-\eqref{E:Q0LOWESTORDERIMPROVED}
	hold with $0 \leq N \leq 15$ replaced by $0 \leq N \leq N_{(Top)}-9,$
	and similarly for the other estimates in the theorem.
	
	Alternatively, we could assume the smallness of 
	$\| \mathring{\Psi} \|_{H_{\Euct}^{25}(\Sigma_0^1)} + \| \mathring{\Psi}_0 \|_{H_{\Euct}^{24}(\Sigma_0^1)}$
	and that
	$\| \mathring{\Psi} \|_{\dot{H}_{\Euct}^{M+1}(\Sigma_0^1)} + \| \mathring{\Psi}_0 \|_{\dot{H}_{\Euct}^M(\Sigma_0^1)} < \infty$
	(without smallness)
	for $M=25,26,\cdots,N_{(Top)}.$
	Here, $\| \cdot \|_{\dot{H}_{\Euct}^M(\Sigma_0^1)}$ denotes the standard homogeneous Sobolev
	norm corresponding to $M^{th}$ order rectangular derivatives along $\Sigma_0^1.$
	The results of Theorem~\ref{T:LONGTIMEPLUSESTIMATES} would of course hold verbatim,
	and we could also prove estimates of the form
	$\totzeromax{25+M}^{1/2}(t,u) \leq 
	f(\| \mathring{\Psi} \|_{H_{\Euct}^{25+M}(\Sigma_0^1)} + \| \mathring{\Psi}_0 \|_{H_{\Euct}^{24+M}(\Sigma_0^1)})
	\ln^{\Conone} (\myexp + t) \upmu_{\star}^{-\Contwo}(t,u),$
	where $f$ is a smooth increasing function 
	and $\Conone,$ $\Contwo$ are positive constants. 
	We could also prove similar estimates for $\totonemax{25+M}^{1/2}(t,u).$
	We could prove these higher-order estimates by
	using arguments similar to the ones we used to prove 
	\eqref{E:Q0TOPORDERKEYBOUND}-\eqref{E:Q1TOPORDERKEYBOUND}.
	In particular, in order to avoid derivative loss,
	the proofs of these estimates must rely on the modified quantities
	of Ch.~\ref{C:RENORMALIZEDEIKONALFUNCTIONQUANTITIES}, at least at the top order.
	
\end{remark}

\begin{proof}[Proof of Theorem~\ref{T:LONGTIMEPLUSESTIMATES}]
	We set
	\begin{align}
		T_{(Max);U_0} &:= \mbox{\ The supremum of the set of times \ } \Tboot \geq 0 \mbox{ \ such that:} 
			\label{E:LIFESPANPROOF}	\\
			& \bullet \mbox{$\Psi,$ $u,$ $\upmu,$ $\Lunit_{(Small)}^i,$ $\Upsilon,$
				and all of the other quantities} 
					\notag \\
			& \ \ \mbox{defined throughout the monograph
				  exist classically on \ } \mathcal{M}_{\Tboot,U_0}
					\notag \\
			& \mbox{\ and all of the solution properties from \textbf{Part I)} of
					Prop.~\ref{P:CLASSICALLOCAL} 
					hold on \ } \mathcal{M}_{\Tboot,U_0}.
					\notag \\
		& \bullet \inf \left\lbrace \upmu_{\star}(t,U_0) \ | \ t \in [0,\Tboot) \right\rbrace > 0.
				\notag \\
		& \bullet \mbox{The fundamental \ } C^0 \mbox{ \ bootstrap assumptions \ } 
						\eqref{E:PSIFUNDAMENTALC0BOUNDBOOTSTRAP}
						\notag \\
		& \ \ \mbox{\ hold with \ } \varepsilon := C_* \mathring{\upepsilon}
								\mbox{\ for \ } (t,u) \in \times [0,\Tboot) \times [0,U_0].
					\notag \\
			& \bullet \mbox{The $L^2-$type bootstrap assumptions \ } \eqref{E:Q0LOWBOOT}-\eqref{E:Q1TOPBOOT}
								\notag \\
			& \ \ \mbox{\ hold with \ } \varepsilon := C_* \mathring{\upepsilon}
						\mbox{\ for \ } (t,u) \in \times [0,\Tboot) \times [0,U_0].
						\notag
		\end{align}
	Actually, the property involving $\upmu_{\star}(t,U_0)$
	is redundant in the sense that one of the 
	solution properties stated in Prop.~\ref{P:CLASSICALLOCAL} 
	is that $\upmu > 0$ on $\mathcal{M}_{\Tboot,U_0}.$
	However, we have nonetheless included the property in the definition of
	$T_{(Max);U_0}$ in order to clarify our discussion below.
	From Proposition
	\ref{P:CLASSICALLOCAL},
	we deduce that if $\mathring{\upepsilon}$ is sufficiently small and $C_* > 1$ 
	is sufficiently large, then $T_{(Max);U_0}> 0.$ 
	
	In the first part of the proof, we focus mostly on deriving the quantitative
	estimates stated in the theorem, but with the 
	time $T_{(Lifespan);U_0}$ defined by \eqref{E:TLIFESPAN}
	replaced by the time $T_{(Max);U_0}$ defined by \eqref{E:LIFESPANPROOF}.
	In the second part of the proof, 
	we focus mostly on proving the qualitative statements,
	showing that $T_{(Max);U_0} = T_{(Lifespan);U_0},$
	and showing that $T_{(Lifespan);U_0}$ is the classical lifespan of the solution.
	Note that the definitions imply 
	that $T_{(Max);U_0} \leq T_{(Lifespan);U_0}$
	and hence to show that $T_{(Max);U_0} = T_{(Lifespan);U_0},$
	we have to show only that $T_{(Max);U_0} \geq T_{(Lifespan);U_0}.$
	
	To proceed, we use Lemma~\ref{L:FUNDAMENTALGRONWALL} and enlarge $C_*$ if necessary
	to deduce that if $C_* \varepsilon^{1/2} = C_*^{3/2} \mathring{\upepsilon}^{1/2} < 1,$
	then the $L^2$ bootstrap assumptions \eqref{E:Q0LOWBOOT}-\eqref{E:Q1TOPBOOT} 
	hold for $(t,u) \in \times [0,T_{(Max);U_0}) \times [0,U_0]$ 
	with the factor $\varepsilon = C_* \mathring{\upepsilon}$ 
	on the right-hand side replaced by the following \emph{smaller} quantity:
	\begin{align} \label{E:KEYGRONWALLBOUNDRECALLED}
			 	\frac{1}{4}
			 	C_*
				\left\lbrace
					\mathring{\upepsilon}
					+ \varepsilon^{3/2}
				\right\rbrace
				= \frac{1}{4}
			 	\left\lbrace
					1
					+ C_* \varepsilon^{1/2}
				\right\rbrace
				\varepsilon
				\leq \frac{1}{2} \varepsilon
				= \frac{1}{2} C_* \mathring{\upepsilon}.
	\end{align}
	In particular,
	the estimates \eqref{E:Q0LOWESTORDERIMPROVED}-\eqref{E:Q1TOPORDERIMPROVED}
	follow from \eqref{E:KEYGRONWALLBOUNDRECALLED}.
	Furthermore,
	using Cor.~\ref{E:C0BOUNDSOBOLEVINTERMSOFENERGIES},
	Def.~\ref{D:MAINCOERCIVEQUANT}, 
	and inequality \eqref{E:FUNCTIONPOINTWISEANGDINTERMSOFANGLIEO},
	and enlarging the constant $C_*$ from \eqref{E:KEYGRONWALLBOUNDRECALLED} if necessary,
	we deduce that the $C^0$ bootstrap assumptions 
	\eqref{E:PSIFUNDAMENTALC0BOUNDBOOTSTRAP}
	for $\Psi$ hold for $(t,u) \in \times [0,T_{(Max);U_0}) \times [0,U_0]$ 
	with $\varepsilon$
	on the right-hand side of \eqref{E:PSIFUNDAMENTALC0BOUNDBOOTSTRAP} replaced by 
	the right-hand side of \eqref{E:KEYGRONWALLBOUNDRECALLED}
	and the norm $C^0(\Sigma_t^u)$ on the left-hand side replaced by $C^0(S_{t,u}).$
	To recover the fact that for $N \leq 13,$ 
	$\mathscr{Z}^N \Psi$ is an element of $C^0(\Sigma_t^u)$
	and not just $C^0(S_{t,u}),$ we simply need non-quantitative
	estimates showing that $\mathscr{Z}^N \Psi$ is jointly continuous in 
	the eikonal and angular coordinates. 
	To this end, we note that at each point in $\mathcal{M}_{T_{(Max);U_0},U_0},$ 
	the vectorfields belonging to the commutation set $\mathscr{Z}$ 
	have span equal to that of the geometric coordinate derivative frame
	$\lbrace 
		\frac{\partial}{\partial t}, 
		\frac{\partial}{\partial u}, 
		\frac{\partial}{\partial \vartheta^1}, 
		\frac{\partial}{\partial \vartheta^2} 
	\rbrace$
	and have components relative to this frame
	that are many times differentiable with respect to 
	$(t,u,\vartheta^1,\vartheta^2).$
	Hence, the desired non-quantitative fact follows from
	the $L^2$ bootstrap assumption \eqref{E:Q0LOWBOOT}, 
	the assumption that $\upmu_{\star} > 0,$
	and standard Sobolev embedding 
	$H^2(\Sigma_t^u) \hookrightarrow C^0(\Sigma_t^u).$
	 
	Throughout the remainder of the proof, we silently use the 
	estimate $\rgeo(t,u) \approx 1 + t$ 
	(on $\mathcal{M}_{T_{(Max);U_0},U_0}$)
	and the fact that 
	$\totzeromax{\leq N}$
	and $\totonemax{\leq N}$
	are increasing in their arguments.

	The estimate \eqref{E:PSIMAINTHEROEMLINFINIFTYESTIMATES} now follows from the above reasoning and 	
	\eqref{E:KEYGRONWALLBOUNDRECALLED}. Furthermore,
	the estimate \eqref{E:PSIMAINTHEROEMLINFINIFTYESTIMATES} implies that
	$\| \Lunit \mathscr{Z}^{\leq 12} \Psi \|_{C^0(\Sigma_t^{U_0})}$
	is uniformly bounded for $0 \leq t < T_{(Max);U_0}.$
	Hence, recalling that $\Lunit = \frac{\partial}{\partial t},$ 
	we conclude that if $T_{(Max);U_0} < \infty,$
	then the functions $\mathscr{Z}^{\leq 12} \Psi$ extend as continuous functions
	of $(t,u,\vartheta)$ to $\Sigma_{T_{(Max);U_0}}^{U_0}.$
	
	The estimates \eqref{E:RECTANGULARSPACETIMEMETRICCOMPONENTSMAINTHEROEMLINFINIFTYESTIMATES}
	and \eqref{E:RECTANGULARSPHEREMETRICCOMPONENTSMAINTHEROEMLINFINIFTYESTIMATES} 
	now follow from 
	inequalities
	\eqref{E:GRECTCOMPONENTC0BOUND}
	and \eqref{E:GSPHERERECTCOMPONENTC0BOUND},
	Cor.~\ref{C:SQRTEPSILONREPLCEDWITHCEPSILON},
	and \eqref{E:KEYGRONWALLBOUNDRECALLED}.
	As in the previous paragraph, the estimates
	\eqref{E:RECTANGULARSPACETIMEMETRICCOMPONENTSMAINTHEROEMLINFINIFTYESTIMATES}
	and \eqref{E:RECTANGULARSPHEREMETRICCOMPONENTSMAINTHEROEMLINFINIFTYESTIMATES}
	imply that
	if $T_{(Max);U_0} < \infty,$
	then the quantities $\mathscr{Z}^{\leq 12} \left\lbrace g_{\mu \nu} - m_{\mu \nu} \right\rbrace,$ etc.
	extend as continuous functions
	of $(t,u,\vartheta)$ to $\Sigma_{T_{(Max);U_0}}^{U_0}.$
	
	The estimates \eqref{E:UPMUMAINTHEROEMLINFINIFTYESTIMATES}-\eqref{E:TRCHIJUNKMAINTHEROEMLINFINIFTYESTIMATES}
	now follow from \eqref{E:C0BOUNDCRUCIALEIKONALFUNCTIONQUANTITIES}, \eqref{E:PSIMAINTHEROEMLINFINIFTYESTIMATES},
	and the fact that $\Radunit_{(Small)}^i = - \Lunit_{(Small)}^i$ plus a smooth function of $\Psi$ that vanishes at $\Psi = 0$
	(that is, \eqref{E:RADUNITJUNKLIKELMINUSUNITJUNK}).
	As in the previous two paragraphs, the estimates
	\eqref{E:UPMUMAINTHEROEMLINFINIFTYESTIMATES}-\eqref{E:TRCHIJUNKMAINTHEROEMLINFINIFTYESTIMATES}
	imply that
	if $T_{(Max);U_0}< \infty,$
	then the quantities $\mathscr{Z}^{\leq 11} \upmu,$ etc.
	extend as continuous functions
	of $(t,u,\vartheta)$ to $\Sigma_{T_{(Max);U_0}}^{U_0}.$
	
	The estimates \eqref{E:TOPORDERANGLAPUPMUMAINTHEROEML2ESTIMATES}-\eqref{E:TOPORDERCHIJUNKANGULARMAINTHEROEML2ESTIMATES}
	now follow from inserting the already proven estimates for $\totzeromax{\leq 24}$ and $\totonemax{\leq 24}$
	into inequality \eqref{E:ELLIPTICRECOVERY} and using the inequalities 
	\eqref{E:KEYMUTOAPOWERINTEGRALBOUND} and \eqref{E:LESSSINGULARTERMSMUINTEGRALBOUND}.
	
	The estimate \eqref{E:TOPORDERRADCHIJUNKMAINTHEROEML2ESTIMATES} 
	now follows from inserting the already proven estimates for 
	$\totzeromax{\leq 24}$ and $\totonemax{\leq 24}$
	and \eqref{E:TOPORDERANGLAPUPMUMAINTHEROEML2ESTIMATES}
	into inequality \eqref{E:L2COMMUTEDLIERADCHIJUNKINTERMSOFANGDSQUAREUPMUANDPSI}
	and using inequality \eqref{E:LESSSINGULARTERMSMUINTEGRALBOUND}.
	
	The estimates \eqref{E:TOPORDERLUPMUMAINTHEROEML2ESTIMATES}-\eqref{E:TOPORDERLCHIJUNKMAINTHEROEML2ESTIMATES}
	now follow from inserting the already proven estimates for $\totzeromax{\leq 24}$ and $\totonemax{\leq 24}$
	into inequality \eqref{E:LDERIVATIVEEIKONALFUNCTIONQUANTITIESL2BOUNDSINTERMSOFQ0ANDQ1}.
	
	The estimates \eqref{E:JUSTBELOWTOPORDERUPMUMAINTHEROEML2ESTIMATES}-\eqref{E:JUSTBELOWTOPORDERCHIJUNKMAINTHEROEML2ESTIMATES}
	now follow from inserting the already proven estimates for $\totzeromax{\leq 24}$ and $\totonemax{\leq 24}$
	into inequality \eqref{E:EIKONALFUNCTIONQUANTITIESL2BOUNDSINTERMSOFQ0ANDQ1} and using 
	inequality \eqref{E:LOGLOSSKEYMUINTEGRALBOUND}.
	
	The estimates \eqref{E:MIDORDERUPMUMAINTHEROEML2ESTIMATES}-\eqref{E:MIDORDERCHIJUNKMAINTHEROEML2ESTIMATES}	
	now follow from inserting the already proven estimates for $\totzeromax{\leq 16 + M}$ and $\totonemax{\leq 16 + M}$
	into inequality \eqref{E:EIKONALFUNCTIONQUANTITIESL2BOUNDSINTERMSOFQ0ANDQ1} and using 
	inequality \eqref{E:LOGLOSSKEYMUINTEGRALBOUND}.
	
	The estimates \eqref{E:NODEGENERACYCUPMUMAINTHEROEML2ESTIMATES}-\eqref{E:NODEGENERACYCHIJUNKMAINTHEROEML2ESTIMATES}
	now follow from inserting the already proven estimates for $\totzeromax{\leq 16}$ and $\totonemax{\leq 16}$
	into inequality \eqref{E:EXTRADERIVATIVELOSSEIKONALFUNCTIONQUANTITIESL2BOUNDSINTERMSOFQ0ANDQ1} and,
	in the case of  
	$\| \mathscr{Z}^{15}(\upmu - 1) \|_{L^2(\Sigma_t^u)},$
	$\| \mathscr{Z}^{15} \Lunit_{(Small)}^i \|_{L^2(\Sigma_t^u)},$
	and $\| \angLie_{\mathscr{Z}}^{14} \upchi^{(Small)} \|_{L^2(\Sigma_t^u)},$
	using inequality \eqref{E:LOGLOSSLESSSINGULARTERMSMTHREEFOURTHSINTEGRALBOUND}.
	
	The estimate \eqref{E:MAINTHEOREMKEYUPMUBEHAVIOR}
	now follows from
	\eqref{E:UPMUALMOSTSHARPBOUND}
	and \eqref{E:KEYGRONWALLBOUNDRECALLED}.
	The estimate \eqref{E:MAINTHEOREMGLLNEARNULLCONDFAILFACT}
	now follows from
	\eqref{E:NULLCONDFACTANNOYINGSTIMATE},
	\eqref{E:IMPORTANTGLLDIFFERENCEESTIMATE},
	and \eqref{E:KEYGRONWALLBOUNDRECALLED}.
	The estimate \eqref{E:MUISNEARONEINITIALLYMAINTHEOREMSTATEMENT} was proved as \eqref{E:UPMUMINUSONESMALLDATA}.
	The estimate \eqref{E:RGEOTIMESRADPSIMAINTHEOREMSTATEMENT} is a special case of \eqref{E:PSIMAINTHEROEMLINFINIFTYESTIMATES},
	which we have already proved.
	
	The estimate \eqref{E:MAINTHEOREMSMALLMUIMPLIESLMUISNEGATIVE}
	was proved as \eqref{E:SMALLMUIMPLIESLMUISNEGATIVE}.
	Inequality \eqref{E:RADUNITPSIBLOWSUP}
	then follows from
	the identity $\Radunit = \upmu^{-1} \Rad,$
	definition \eqref{E:BIGOMEGADEF},
	and
	the estimates \eqref{E:BIGOMEGAISWELLAPPROXIMATEDBYNULLCONDITIONFAILUREFACTORTERM}
	and \eqref{E:MAINTHEOREMSMALLMUIMPLIESLMUISNEGATIVE}.
	To estimate the Euclidean length of the difference of $\Radunit$ 
	and $- \partial_r,$ we use the identity
	$\partial_r = \frac{x^a}{r} \partial_a,$
	the identity \eqref{E:RADUNITJUNK}, 
	and the estimates
	\eqref{E:EASYEUCLIDEANRADIALVARIABLEBOUND}
	(by Cor.~\ref{C:SQRTEPSILONREPLCEDWITHCEPSILON} and \eqref{E:KEYGRONWALLBOUNDRECALLED},
	this estimate is valid with $\varepsilon^{1/2}$ replaced by $C \mathring{\upepsilon}$)
	and \eqref{E:LUNITJUNKMAINTHEROEMLINFINIFTYESTIMATES} 
	to deduce that 
		\begin{align} \label{E:RADUNITNEAREUCLIDEANRADIAL}
			\left| \Radunit - (- \partial_r) \right|_{\Euct}
			& \lesssim
				\sum_{a=1}^3 \left| \frac{x^a}{r} \right|
				\left| 1 - \frac{r}{\rgeo} \right|
				+ \sum_{a=1}^3 \left| \Radunit_{(Small)}^{a} \right|
			\lesssim  
			\mathring{\upepsilon} 
			\frac{\ln(\myexp + t)}{1 + t}.
		\end{align}
	
		We now prove the statements concerning the behavior of the change of variables map $\Upsilon$
		in the case $T_{(Max);U_0} < \infty.$
		Using equation \eqref{E:CHOVMATRIXCOMPUTATIONS}
		and the argument that we used in the proof of Cor.~\ref{C:VECTORFIEDLDIFFERENTIABILITYVSGEOCOORDINATEDIFFERENTIABILITY},
		it is straightforward to see that
		the extendibility of $\Upsilon$ 
		to $\Sigma_{T_{(Max);U_0}}^{U_0}$ as a
		$C^{10}$ function
		of $(t,u,\vartheta^1,\vartheta^2)$
		follows after we prove that
		\textbf{i)} 
		$\mathscr{Z}^{\leq 10} \Xi^i,$
		$\mathscr{Z}^{\leq 11} \Lunit^i,$
		$\mathscr{Z}^{\leq 11} \Rad^i,$
		$\mathscr{Z}^{\leq 11} \Xi^i,$
		and $\mathscr{Z}^{\leq 11} X_A^i$
		extend as continuous functions
		of $(t,u,\vartheta^1,\vartheta^2)$ to $\Sigma_{T_{(Max);U_0}}^{U_0}$
		and 
		\textbf{ii)} for each $Z \in \mathscr{Z},$ 
		the components of $Z$ relative to the geometric coordinates
		extend as $C^{10}$ functions
		of $(t,u,\vartheta^1,\vartheta^2)$ to $\Sigma_{T_{(Max);U_0}}^{U_0}.$
		The result \textbf{ii)} follows from
		the estimates of Lemma~\ref{L:GEOMETRICCOMPONENTS},
		the argument that we used in the proof of
		Cor.~\ref{C:VECTORFIEDLDIFFERENTIABILITYVSGEOCOORDINATEDIFFERENTIABILITY},
		and the argument that we used to prove
		the continuous extendibility of $\mathscr{Z}^{\leq 12} \Psi.$
	 	We note for later use that \textbf{ii)} 
	 	and Lemma~\ref{L:VECTORFIELDCOMMUTATORS}
	 	imply that
		the vectorfields $\angLie_{\mathscr{Z}}^{\leq 10} X_A = \Lie_{\mathscr{Z}}^{\leq 10} X_A$
		also continuously extend
		as functions of $(t,u,\vartheta^1,\vartheta^2)$ to $\Sigma_{T_{(Max);U_0}}^{U_0}.$
		To prove \textbf{i)}, we first note that
		the previous arguments have shown that 
		$\mathscr{Z}^{\leq 11} \Lunit^i$ and $\mathscr{Z}^{\leq 11} \Rad^i = \mathscr{Z}^{\leq 11} (\upmu \Radunit^i)$				
		extend as continuous functions
		of $(t,u,\vartheta^1,\vartheta^2)$ to $\Sigma_{T_{(Max);U_0}}^{U_0}.$
		Using similar reasoning and the estimate \eqref{E:C0BOUNDSXIRECTANGULARCOMPONENTS},
		we conclude that $\mathscr{Z}^{\leq 10} \Xi^i$ 
		extends as a continuous function
		of $(t,u,\vartheta^1,\vartheta^2)$ to $\Sigma_{T_{(Max);U_0}}^{U_0}.$
		To handle the rectangular components $X_A^i,$
		we first use similar reasoning and
		the estimate \eqref{E:LOWERORDERPOINTWISEBOUNDPROJECTEDLIEDERIVATIVESANGDIFFCOORDINATEX}
		to deduce that the $S_{t,u}$ one-form $\angdiff \mathscr{Z}^{\leq 11} x^i$
		extends as a continuous function
		of $(t,u,\vartheta^1,\vartheta^2)$ to $\Sigma_{T_{(Max);U_0}}^{U_0}.$
		Using the aforementioned continuous extendibility of the vectorfields $\angLie_{\mathscr{Z}}^{\leq 10} X_A,$
		and the identity $X_A^i = X_A \cdot \angdiff x^i,$
		we conclude that $\mathscr{Z}^{\leq 10} X_A^i$
		extends as a continuous function
		of $(t,u,\vartheta^1,\vartheta^2)$ to $\Sigma_{T_{(Max);U_0}}^{U_0}.$ 
		We have thus proved the extendibility of $\Upsilon$
		as a $C^{10}$ function of $(t,u,\vartheta^1,\vartheta^2).$
		
		
	With the help of the above results,
	we now show that either 
	\textbf{a)} $T_{(Max);U_0} = \infty$ or
	\textbf{b)} $T_{(Max);U_0} < \infty$ and
		$\inf \left\lbrace \upmu_{\star}(t,U_0) \ | \ t \in [0,T_{(Max);U_0}) \right\rbrace = 0.$
	We note that after we have shown that either \textbf{a)} or \textbf{b)} must hold, 
	it follows from the comments made
	in the second paragraph of the proof that $T_{(Max);U_0} = T_{(Lifespan);U_0}$
	(where $T_{(Lifespan);U_0}$ is defined in \eqref{E:TLIFESPAN}).
	Furthermore, in the case $T_{(Lifespan);U_0} < \infty,$
	it follows from the estimates
	\eqref{E:RADUNITPSIBLOWSUP} 
	and \eqref{E:RADUNITNEAREUCLIDEANRADIAL}
	that the vanishing of $\upmu_{\star}(t,U_0)$ at time $T_{(Lifespan);U_0}$
	leads to the blow-up result \eqref{E:RECTANGULARDERIVATIVESBLOWUP}.
	
	It remains for us show that \textbf{a)} or \textbf{b)} must hold.
	It suffices to show that it is impossible to have
	$T_{(Max);U_0} < \infty$ unless 
	$\inf \left\lbrace \upmu_{\star}(t,U_0) \ | \ t \in [0,T_{(Max);U_0}) \right\rbrace = 0.$
	To proceed, we assume for the sake of deriving a contradiction that: 
	$T_{(Max);U_0} < \infty$ 
	and $\inf \left\lbrace \upmu_{\star}(t,U_0) \ | \ t \in [0,T_{(Max);U_0}) \right\rbrace > 0.$
	We now rule out, 
	on the spacetime domain $\mathcal{M}_{T_{(Max);U_0},U_0},$
	the $4$ breakdown scenarios from \textbf{Part II)} of Prop.~\ref{P:CLASSICALLOCAL},
	and we also show that the condition $(5)$ for $\Upsilon$ is verified
	(when $\varepsilon$ is sufficiently small).
	Scenario $(1)$ is ruled out by assumption.
	Scenario $(2)$ is ruled out by the estimate \eqref{E:C0BOUNDCRUCIALEIKONALFUNCTIONQUANTITIES}.
	Scenario $(3)$ is ruled out by the bootstrap assumptions \eqref{E:PSIFUNDAMENTALC0BOUNDBOOTSTRAP}
	and the fact that $g_{\mu \nu}^{(Small)}(\cdot)$ is a smooth function of $\Psi$
	that vanishes at $\Psi = 0$ (see \eqref{E:LITTLEGDECOMPOSED} and \eqref{E:METRICNONLINEARITYVANISHESWHENPSIIS0}). 
	
	In the next paragraph, we show that $\Upsilon$ is a 
	$C^{10}$ diffeomorphism from $[0,T_{(Max);U_0}] \times [0,U_0] \times \mathbb{S}^2$
	onto its image and hence $\Upsilon$ verifies the condition $(5).$
	Given this fact, we can easily rule out the breakdown scenario $(4).$
	To this end, we first
	use the previously shown fact that
	the functions $\mathscr{Z}^{\leq 12} \Psi$ extend as continuous functions
	of $(t,u,\vartheta)$ to $\Sigma_{T_{(Max);U_0}}^{U_0}$
	and argue as in our proof that $\Upsilon$ extends
	as a $C^{10}$ function of $(t,u,\vartheta)$ 
	in order to conclude that $\Psi$ extends to $\Sigma_{T_{(Max);U_0}}^{U_0}$ 
	as a $C^{12}$ function of $(t,u,\vartheta).$ 
	Hence, we deduce from the regular nature of $\Upsilon^{-1}$ 
	that $\partial_{\kappa} \Psi,$ 
	$(\kappa = 0,1,2,3),$ 
	extends to $\mathcal{M}_{T_{(Max);U_0},{U_0}} \cup \Sigma_{T_{(Max);U_0}}^{U_0}$
	as a continuous function
	of the \emph{rectangular coordinates} and that
	$\sup_{\mathcal{M}_{T_{(Max);U_0},U_0} \cup \Sigma_{T_{(Max);U_0}}^{U_0}}
	\max_{\kappa =0,1,2,3} |\partial_{\kappa} \Psi| < \infty.$
	We have thus ruled out the breakdown scenario $(4).$ 
	
	It remains for us to show that condition 
	$(5)$ from \textbf{Part II)} of Prop.~\ref{P:CLASSICALLOCAL} is verified.
	We first recall that we have already shown that $\Upsilon$ 
	is a bijection on $[0,T_{(Max);U_0}) \times [0,U_0] \times \mathbb{S}^2$ 
	that extends as a $C^{10}$ function of $(t,u,\vartheta)$
	to the compact set $[0,T_{(Max);U_0}] \times [0,U_0] \times \mathbb{S}^2.$
	Furthermore, from equation \eqref{E:JACOBIAN},
	the estimates 
	\eqref{E:C0BOUNDCRUCIALEIKONALFUNCTIONQUANTITIES},
	\eqref{E:SQRTDETANGSPHEREMETRICPOINTWISEESTIMATE},
	and
	\eqref{E:PSIMAINTHEROEMLINFINIFTYESTIMATES},
	and the assumption $\inf \left\lbrace \upmu_{\star}(t,U_0) \ | \ t \in [0,T_{(Max);U_0}) \right\rbrace > 0,$
	we conclude that the Jacobian determinant of $\Upsilon$
	is strictly positive on $[0,T_{(Max);U_0}] \times [0,U_0] \times \mathbb{S}^2.$
	From the inverse function theorem, we conclude that 
	$\Upsilon$ extends to $[0,T_{(Max);U_0}] \times [0,U_0] \times \mathbb{S}^2$
	as a locally invertible $C^{10}$ function.  
	We now show that $\Upsilon$ is a $C^{10}$ (global) diffeomorphism from
 	$[0,T_{(Max);U_0}] \times [0,U_0] \times \mathbb{S}^2$ onto its image.
 	To this end, we first note that inequality \eqref{E:EIKONALFUNCTIONRECTANGULARDERIVATIVES},
 	the assumption $\inf \left\lbrace \upmu_{\star}(t,U_0) \ | \ t \in [0,T_{(Max);U_0}) \right\rbrace > 0,$
 	and the estimate $r \geq (1 - \mathcal{O}(\varepsilon)) \rgeo \geq (1 - \mathcal{O}(\varepsilon))(1 - U_0)$
	(which follows from \eqref{E:EASYEUCLIDEANRADIALVARIABLEBOUND} and Cor.~\ref{C:SQRTEPSILONREPLCEDWITHCEPSILON})
 	imply that
 	$\sum_{a=1}^3 
		\left|
			\partial_a u
		\right|
	$
	is uniformly bounded from above and below strictly away from $0.$
	It follows that the outgoing null cone portions $\mathcal{C}_u^{T_{(Max);U_0}}$
	corresponding to two distinct values $u \in [0,U_0]$ cannot intersect.
	Hence, to prove the global invertibility of $\Upsilon,$ it remains only for us
	to show that\footnote{Recall that the angular coordinates $\vartheta^A$ are, by construction,
	constant along the integral curves of $\Lunit.$} for each $u \in [0,U_0],$
	the integral curves of $\Lunit,$  
	which rule $\mathcal{C}_u^{T_{(Max);U_0}},$
	do not intersect\footnote{Equivalently, we must show that
	if $q \in \mathcal{C}_u^{T_{(Max);U_0}},$ then
	no point in $\mathcal{C}_u^{T_{(Max);U_0}}$ can belong to the null cut locus of $q.$} 
 	at time $T_{(Max);U_0}.$ 
 	Equivalently, we show that for each fixed $(t,u) \in [0,T_{(Max);U_0}] \times [0,U_0],$
 	the map $\Upsilon(t,u,\cdot),$ defined on the domain $\mathbb{S}^2,$ is injective.
	To this end, we define, relative to the rectangular coordinates, 
	the map $F(t,x^1,x^2,x^3): \Sigma_t \setminus \lbrace (t,0,0,0) \rbrace \rightarrow \mathbb{S}^2 \subset \mathbb{R}^3$
	by $F(t,x^1,x^2,x^3)=(x^1/r,x^2/r,x^3/r).$
	Consider the composition $F \circ \Upsilon(t,u,\cdot): \mathbb{S}^2 \rightarrow \mathbb{S}^2,$ which is $C^{10}$
	in view of the estimate $r \geq (1 - \mathcal{O}(\varepsilon))(1 - U_0)$
	mentioned above. 
	Just below we show that $F \circ \Upsilon(t,u,\cdot)$
	is a $C^{10}$ immersion.\footnote{Recall that $F \circ \Upsilon(t,u,\cdot)$ is said to be an immersion
	if the Jacobian $\frac{\partial (F \circ \Upsilon)}{\partial (\vartheta^1,\vartheta^2)}$
	is rank $2$ at each point in its domain $\mathbb{S}^2.$}
	Given this fact,
	it follows from the compactness of $\mathbb{S}^2$ that
	$F \circ \Upsilon(t,u,\cdot)$ is a covering map.\footnote{
	A continuous surjective function 
	$f: X \rightarrow \mathbb{S}^2$ is said to be a covering map of $\mathbb{S}^2$
	if for every $p \in \mathbb{S}^2,$ there exists a neighborhood $\Omega$ of $p$ such 
	that $f^{-1}(\Omega)$ is a disjoint union of open sets in $X,$
	each of which is mapped homeomorphically onto $\Omega$ by $f.$} 
	Then, since $\mathbb{S}^2$ is simply connected, 
	it is a basic result
	of algebraic topology that $F \circ \Upsilon(t,u,\cdot)$ must be a $C^{10}$ 
	diffeomorphism from $\mathbb{S}^2$ to $\mathbb{S}^2,$
	which yields the desired injectivity of $\Upsilon(t,u,\cdot).$
	It remains for us to show that for each fixed $(t,u) \in [0,T_{(Max);U_0}] \times [0,U_0],$
	$F \circ \Upsilon(t,u,\cdot)$ is an immersion. 
	That is, we must show that
	the $3 \times 2$ matrix
	$\frac{\partial (F \circ \Upsilon)}{\partial (\vartheta^1,\vartheta^2)}$
	is rank $2.$
	To proceed, we first compute that the $3 \times 4$ matrix
	$\frac{\partial F}{\partial (t,x^1,x^2,x^3)}$ is 
	rank $2$ with kernel $K$ given by
	$K = \mbox{\upshape span} \lbrace (1,0,0,0)^{\intercal}, (0,x^1/r,x^2/r,x^3/r)^{\intercal} \rbrace,$
	where $\intercal$ denotes the transpose.
	Next, we note that inequality \eqref{E:EIKONALFUNCTIONRECTANGULARDERIVATIVES}
	implies that the Euclidean outward normal to $S_{t,u}$ in $\Sigma_t$
	at the point $(t,x^1,x^2,x^3) = \Upsilon(t,u,\vartheta)$
	has rectangular components equal to $(0,x^1/r,x^2/r,x^3/r)^{\intercal}$ 
	up to an error vector of Euclidean length 
	$\mathcal{O}\left(\varepsilon \ln(\myexp + t)(1 + t)^{-1} \right).$
	Hence, since the two columns of the $4 \times 2$ matrix
	$\frac{\partial \Upsilon}{\partial (\vartheta^1,\vartheta^2)}$
	are tangent to $S_{t,u},$ they do not belong to $K.$
	Since the two columns are also linearly dependent, it follows that 
	their union with $K$ forms a basis of $\mathbb{R}^4.$
	It thus follows from the chain rule that the $4 \times 2$ matrix
	$\frac{\partial (F \circ \Upsilon)}{\partial (\vartheta^1,\vartheta^2)}$
	is rank $2$ as desired.

	In total, we conclude from \textbf{Part II)} of Prop.~\ref{P:CLASSICALLOCAL}
	that there exists a $\Delta > 0$ such that $\Psi,$ $u,$ $\upmu,$ $\Lunit_{(Small)}^i,$
	etc.  
	can be extended as solutions to a strictly larger region of the form
	$\mathcal{M}_{T_{(Max);U_0}+ \Delta,U_0}$ for some $\Delta > 0.$ 
	Furthermore, the proposition also yields
	that the quantities appearing on the left-hand sides
	of \eqref{E:Q0LOWBOOT}-\eqref{E:Q1TOPBOOT} are continuous in $t$ and $u.$
	In view of this fact and the bound \eqref{E:KEYGRONWALLBOUNDRECALLED}, 
	we deduce that
	if $\Delta$ is small enough, then 
	\eqref{E:Q0LOWBOOT}-\eqref{E:Q1TOPBOOT} 
	hold for $(t,u) \in \times [0,T_{(Max);U_0}+ \Delta) \times [0,U_0]$ with $\varepsilon$
	replaced by $\frac{3}{4} \varepsilon.$ We have thus arrived at the
	desired contradiction of the definition of $T_{(Max);U_0}.$
	
	To prove \eqref{E:TLIFESPANLOWERBOUND}, that is, 
	that there exists a constant $C_{(Lower-Bound)} > 0$ such that
	$\upmu_{\star}(t,U_0) > 0$ whenever 
	$t < \exp\left( \left\lbrace C_{(Lower-Bound)} \mathring{\upepsilon} \right\rbrace^{-1} \right),$
	we only have to combine 
	\eqref{E:MAINTHEOREMKEYUPMUBEHAVIOR},
	\eqref{E:MUISNEARONEINITIALLYMAINTHEOREMSTATEMENT},
	and \eqref{E:RGEOTIMESRADPSIMAINTHEOREMSTATEMENT}.
	Similar reasoning also yields \eqref{E:UPMUREMAINSLARGEWHENNULLCONDITIONFAILUREFACTORVANISHES}.
	
	Finally, to conclude,
	in the case $T_{(Lifespan);U_0} < \infty,$
	that the Jacobian determinant of $\Upsilon$ vanishes precisely on the subset 
	\eqref{E:MUVANISHESSET},
	we use equation \eqref{E:JACOBIAN},
	the estimate \eqref{E:PSIMAINTHEROEMLINFINIFTYESTIMATES} for $\Psi,$
	and inequality 
	\eqref{E:SQRTDETANGSPHEREMETRICPOINTWISEESTIMATE}
	(which, by \eqref{E:KEYGRONWALLBOUNDRECALLED}, 
	is valid with $\varepsilon$ replaced by $C \mathring{\upepsilon}$).

\end{proof}

\section{More precise control over angular derivatives.}
\label{S:PRECISECONTROLOVERANGULARDERIVATIVES}
We now state and prove a corollary of Theorem~\ref{T:LONGTIMEPLUSESTIMATES}  
that provides additional information about some of the 
angular derivatives of $\Psi.$ Roughly, the corollary allows us to infer that nearly spherically symmetric data 
launch nearly spherically solutions, even though the wave equation 
\eqref{E:WAVEGEO} is not necessarily invariant under Euclidean rotations.
The estimates from the corollary play an important role in our proof of 
small data shock formation, which is based on arguments that are valid 
for nearly spherically symmetric solutions.

\begin{corollary}[\textbf{More precise control over the lower-order angular derivatives}] 
\label{C:ANGULARDDERIVATIVEESTIMATES}
Assume the hypotheses and conclusions of Theorem~\ref{T:LONGTIMEPLUSESTIMATES}.
In particular, let $\mathring{\upepsilon}$ be the size of the data, 
and let $T_{(Lifespan);U_0}$ and $U_0$ be the numbers appearing in the statement of the theorem.
Let 
$\enzero[\cdots],$
$\flzero[\cdots],$
$\enone[\cdots],$
$\flone[\cdots],$
and $\Morint[\cdots]$
be the energies, fluxes, and Morawetz spacetime integral from
Defs.~\ref{D:ENERGIESANDFLUXES} and \ref{D:COERCIVEMORAWETZINTEGRAL}.
Let
\begin{align} \label{E:SMALLANGULARDERIVATIVES}
	\mathring{\updelta} := \max_{2 \leq |\vec{I}| \leq 4} \enzero^{1/2}[\mathscr{O}^{\vec{I}} \Psi](0,U_0).
\end{align}
There exists a constant $C > 0$ such that the following estimates hold for $(t,u) \in [0,T_{(Lifespan);U_0}) \times [0,U_0]:$ 
\begin{subequations}
\begin{align}
	\max_{2 \leq |\vec{I}| \leq 4}
		\enzero^{1/2}[\mathscr{O}^{\vec{I}} \Psi](t,u)
		+ \flzero^{1/2}[\mathscr{O}^{\vec{I}} \Psi](t,u)
	& \leq C (\mathring{\updelta} + \mathring{\upepsilon}^{3/2}),
			\label{E:L2ANGULARDERIVATIESEVENSMALLER} \\
	\max_{2 \leq |\vec{I}| \leq 4}
		\enone^{1/2}[\mathscr{O}^{\vec{I}} \Psi](t,u)
		+ \flone^{1/2}[\mathscr{O}^{\vec{I}} \Psi](t,u)
		+ \Morint^{1/2}[\mathscr{O}^{\vec{I}} \Psi](t,u)
	& \leq C (\mathring{\updelta} + \mathring{\upepsilon}^{3/2}) \ln^2(\myexp + t),
			\label{E:MORL2ANGULARDERIVATIESEVENSMALLER} \\
	\| \angLap \Psi \|_{C^0(\Sigma_t^u)}
	& \leq C (\mathring{\updelta} + \mathring{\upepsilon}^{3/2}) \frac{1}{(1 + t)^3}.
		\label{E:ANGULARLAPLACIANEVENSMALLER}
\end{align}
\end{subequations}

\end{corollary}

\begin{proof}
The estimate \eqref{E:ANGULARLAPLACIANEVENSMALLER} follows from
\eqref{E:L2ANGULARDERIVATIESEVENSMALLER}
and inequalities
\eqref{E:FUNCTIONPOINTWISEANGDINTERMSOFANGLIEO},
\eqref{E:ANGLAPFUNCTIONPOINTWISEINTERMSOFROTATIONS},
and \eqref{E:C0BOUNDSOBOLEVINTERMSOFENERGIES}.
Hence, we need to prove only inequalities \eqref{E:L2ANGULARDERIVATIESEVENSMALLER}
and \eqref{E:MORL2ANGULARDERIVATIESEVENSMALLER}.

To this end, we commute the wave equation \eqref{E:WAVEGEO} with 
iterated rotational differential operators $\mathscr{O}^{\vec{I}},$
where $2 \leq |\vec{I}| \leq 4,$ to derive an equation of the form
\begin{align}
	\upmu \square_{g(\Psi)} \mathscr{O}^{\vec{I}} \Psi 
	& = \inhomarg{\mathscr{O}^{\vec{I}}}.
\end{align}
Using 
\eqref{E:NTIMESCOMMUTEDWAVEINHOMOGENEOUSTERMFIRSTSPLITTING},
\eqref{E:LOWERORDERINHOMOGENEOUSTERMSFIRSTPOINTWISE},
\eqref{E:ROTDEFORMLL}-\eqref{E:ROTDEFORMSPHERETRACE},
\eqref{E:DIVCOMMUTATIONCURRENTDECOMPOSITION},
and the $C^0$ estimates of Theorem~\ref{T:LONGTIMEPLUSESTIMATES},
we see that
$\inhomarg{\mathscr{O}^{\vec{I}}}$ is \emph{quadratically small}
in the sense that its magnitude is bounded by a decaying function of time 
multiplied by $\mathring{\upepsilon}^2.$
This quadratic smallness is 
in particular based on the fact that 
the deformation tensors
$\deform{\Rot}$ of the rotations
completely vanish for the background solution $\Psi \equiv 0$
(see \eqref{E:ROTDEFORMLL}-\eqref{E:ROTDEFORMSPHERETRACE}).
The reason is that the rotations $\Rot$ are equal to the Euclidean rotations when 
$\Psi \equiv 0,$ and 
the Euclidean rotations are Killing fields\footnote{That is, we have $\Lie_{\Roteuc} m = 0.$}  of the background Minkowski metric. 
We remark that this stands in contrast to $\deform{\Rad}$ and $\deform{\rgeo \Lunit},$
which do not vanish even when $\Psi \equiv 0.$

We now revisit the energy-flux identities \eqref{E:E0DIVID}-\eqref{E:E1DIVID},
where $\mathscr{O}^{\vec{I}} \Psi$ plays the role of $\Psi$
and $\inhomarg{\mathscr{O}^{\vec{I}}}$ plays the role of $\waveinhom.$
To simplify the notation, we set
\begin{align}
	\mathbb{Q}_{(\mathscr{\Rot})}(t,u)
	 & := \max_{2 \leq |\vec{I}| \leq 4} 
		\enzero[\mathscr{O}^{\vec{I}} \Psi](t,u)
		+ \flzero[\mathscr{O}^{\vec{I}} \Psi](t,u),
		\\
	\widetilde{\mathbb{Q}}_{(\mathscr{\Rot})}(t,u)
	& := \max_{2 \leq |\vec{I}| \leq 4} 
		\enone[\mathscr{O}^{\vec{I}} \Psi](t,u)
		+ \flone[\mathscr{O}^{\vec{I}} \Psi](t,u),
			\\
	\widetilde{\mathbb{K}}_{(\mathscr{\Rot})}(t,u)
	& := \max_{2 \leq |\vec{I}| \leq 4}
			\Morint[\mathscr{O}^{\vec{I}} \Psi](t,u),
\end{align}
where $\Morint[\cdot](t,u)$ is the Morawetz integral from \eqref{E:COERCIVEMORDEF}.
Our goal is to use the already proven non-$\upmu_{\star}^{-1}-$degenerate estimates 
for the lower-order derivatives of various quantities implied by Theorem~\ref{T:LONGTIMEPLUSESTIMATES}
in order to prove that the following estimates hold for $(t,u) \in [0,T_{(Lifespan);U_0}) \times [0,U_0]:$
\begin{align}
	\mathbb{Q}_{(\mathscr{\Rot})}(t,u)
	& \leq
		C \mathring{\updelta}^2
		+
		C \mathring{\upepsilon}^3
		+ C
			\int_{t'=0}^t 
				\frac{\ln^2(\myexp + t')}{(1 + t')^2} \mathbb{Q}_{(\mathscr{\Rot})}(t',u) 
			\, dt'
			 \label{E:MULTLOWERORDERROTATIONALGRONWALLREADY} \\
	& \ \ 
		+ C
			\int_{t'=0}^t 
				\frac{\ln^2(\myexp + t')}{(1 + t')^2} \widetilde{\mathbb{Q}}_{(\mathscr{\Rot})}(t',u) 
			\, dt'
		\notag \\
	& \ \ 
			+ C
				\mathring{\upepsilon}^{1/2}
				\int_{t'=0}^t 
					\frac{\ln(\myexp + t')}{(\myexp + t')^2 \sqrt{\ln(\myexp + t) - \ln(\myexp + t')}}
					\widetilde{\mathbb{Q}}_{(\mathscr{\Rot})}(t',u)   
				\, dt' 
				\notag \\
	& \ \ +  	C
						\frac{\ln^2(\myexp + t)}{(1 + t)^2} 
						\int_{u'=0}^u
							\widetilde{\mathbb{Q}}_{(\mathscr{\Rot})}(t,u')
						\, du',
				\notag \\
	\widetilde{\mathbb{Q}}_{(\mathscr{\Rot})}(t,u)
		+ \widetilde{\mathbb{K}}_{(\mathscr{\Rot})}(t,u)
	& \leq
			C \mathring{\updelta}^2
			+
			C \mathring{\upepsilon}^3
			+
			C \ln^4(\myexp + t) \mathbb{Q}_{(\mathscr{\Rot})}(t,u)
				\label{E:MORLOWERORDERROTATIONALGRONWALLREADY} \\
		& \ \ + (1 + C \varepsilon)
						\int_{t'=0}^t 
							\frac{1}{\rgeo(t',u) \ln \left(\frac{\rgeo(t',u)}{\rgeo(0,u)} \right)} 
							\widetilde{\mathbb{Q}}_{(\mathscr{\Rot})}(t',u) 
						\, dt' 
					\notag	\\
		& \ \ + C 
						\frac{\ln^3(\myexp + t)}{(1 + t)^2} 
						\int_{u'=0}^u
							\widetilde{\mathbb{Q}}_{(\mathscr{\Rot})}(t,u')		
						\, du'.
					\notag
\end{align}
Once we have shown \eqref{E:MULTLOWERORDERROTATIONALGRONWALLREADY}-\eqref{E:MORLOWERORDERROTATIONALGRONWALLREADY},
the desired estimates \eqref{E:L2ANGULARDERIVATIESEVENSMALLER}-\eqref{E:MORL2ANGULARDERIVATIESEVENSMALLER}
follow from applying a Gronwall-type argument to the system of inequalities \eqref{E:MULTLOWERORDERROTATIONALGRONWALLREADY}-\eqref{E:MORLOWERORDERROTATIONALGRONWALLREADY}
in the variables 
$\mathbb{Q}_{(\mathscr{\Rot})},$
$\widetilde{\mathbb{Q}}_{(\mathscr{\Rot})},$
and $\widetilde{\mathbb{K}}_{(\mathscr{\Rot})}.$
Because the argument is much simpler than the
Gronwall estimates proved for the below-top-order quantities 
$\totzeromax{\leq 23},$
$\totonemax{\leq 23},$
and $\totMormax{\leq 23}(t,u)$
(the proof of those estimates starts just above equation \eqref{E:LOWERORDERQ0INTEGRATINGFACTORPRODUCT}),
we omit the details. However, we remark that we again use
Lemma~\ref{L:UNUSUALINTEGRATINGFACTORESTIMATE} to handle the next-to-last term on the right-hand side of 
\eqref{E:MULTLOWERORDERROTATIONALGRONWALLREADY}.

In order to prove \eqref{E:MULTLOWERORDERROTATIONALGRONWALLREADY}-\eqref{E:MORLOWERORDERROTATIONALGRONWALLREADY},
we have to estimate the integrals on the right-hand sides of
\eqref{E:E0DIVID}-\eqref{E:E1DIVID},
where $\mathscr{O}^{\vec{I}} \Psi$ plays the role of $\Psi$
and $\inhomarg{\mathscr{O}^{\vec{I}}}$ plays the role of $\waveinhom.$
For each fixed $\vec{I},$ we show how to deduce an inequality of the form 
\eqref{E:MULTLOWERORDERROTATIONALGRONWALLREADY}-\eqref{E:MORLOWERORDERROTATIONALGRONWALLREADY},
but with $\enzero[\mathscr{O}^{\vec{I}} \Psi](t,u) + \flzero[\mathscr{O}^{\vec{I}} \Psi](t,u)$ 
on the left-hand side of \eqref{E:MULTLOWERORDERROTATIONALGRONWALLREADY}
in place of $\mathbb{Q}_{(\mathscr{\Rot})}(t,u)$
and $\enone[\mathscr{O}^{\vec{I}} \Psi](t,u) + \flone[\mathscr{O}^{\vec{I}} \Psi](t,u) + \Morint[\mathscr{O}^{\vec{I}} \Psi](t,u)$
on the left-hand side of \eqref{E:MORLOWERORDERROTATIONALGRONWALLREADY}
in place of $\widetilde{\mathbb{Q}}_{(\mathscr{\Rot})}(t,u) + \widetilde{\mathbb{K}}_{(\mathscr{\Rot})}(t,u).$
We then take the max of both sides of the inequalities
over all relevant quantities with
$2 \leq |\vec{I}| \leq 4$ 
in order to deduce \eqref{E:MULTLOWERORDERROTATIONALGRONWALLREADY}-\eqref{E:MORLOWERORDERROTATIONALGRONWALLREADY}.

To begin the detailed analysis, 
we remark that we will handle the integrals involving $\inhomarg{\mathscr{O}^{\vec{I}}}$
at the end of the proof. 
From the assumption \eqref{E:SMALLANGULARDERIVATIVES}, 
we see that the terms $\enzero[\mathscr{O}^{\vec{I}} \Psi](0,u)$
arising from the first term on the right-hand side of
\eqref{E:E0DIVID} are $\leq \mathring{\updelta}^2$
and hence are $\leq$ the right-hand side of \eqref{E:MULTLOWERORDERROTATIONALGRONWALLREADY}
as desired.
We now explain how we bound the integral 
$\int_{\mathcal{M}_{t,u}} \upmu \enmomtensor^{\alpha \beta}[\mathscr{O}^{\vec{I}} \Psi] 
\deformarg{\Mult}{\alpha}{\beta} \, d \vol$
from the right-hand side of \eqref{E:E0DIVID} 
by $\leq$ the right-hand side of \eqref{E:MULTLOWERORDERROTATIONALGRONWALLREADY}.
To derive the desired bound, we examine the proof of inequality
\eqref{E:Q0BASICERRORINTEGRALESTIMATE}
and conclude that the same inequality holds true
with $\mathscr{O}^{\vec{I}} \Psi$ in place of $\mathscr{Z}^N \Psi$ on the left-hand side and
$\mathbb{Q}_{(\mathscr{\Rot})},$
$\widetilde{\mathbb{Q}}_{(\mathscr{\Rot})},$
and $\widetilde{\mathbb{K}}_{(\mathscr{\Rot})}$
respectively in place of 
$\totzeromax{N},$
$\totonemax{N},$
and $\totMormax{N}$
on the right-hand side.
Theorem~\ref{T:LONGTIMEPLUSESTIMATES}
implies that
Lemma~\ref{L:BASICENERGYERRORINTEGRALS} holds with $\varepsilon$ replaced by $C \mathring{\upepsilon}$
and furthermore, by \eqref{E:Q0LOWESTORDERIMPROVED} and \eqref{E:Q1LOWESTORDERIMPROVED}, 
the fourth and last integrals on the right-hand side of \eqref{E:Q0BASICERRORINTEGRALESTIMATE}
(with $\mathbb{Q}_{(\mathscr{\Rot})}$ in the role of $\totzeromax{N}$ and $\widetilde{\mathbb{K}}_{(\mathscr{\Rot})}$ in the role of $\totMormax{N}$)
are $\lesssim \mathring{\upepsilon}^3$ and hence are bounded by the right-hand side
of \eqref{E:MULTLOWERORDERROTATIONALGRONWALLREADY} as desired.


We now address the terms on the right-hand side of \eqref{E:E1DIVID}
that do not involve the integrand factor $\inhomarg{\mathscr{O}^{\vec{I}}},$
where $\mathscr{O}^{\vec{I}} \Psi$ plays the role of $\Psi.$
The integral over the hypersurface $\Sigma_0^u$ is quadratic in 
$\mathscr{O}^{\vec{I}} \Psi$ 
and is controlled by the data.
Hence, by \eqref{E:MULTENERGYCOERCIVITY}
and the assumption \eqref{E:SMALLANGULARDERIVATIVES},
the magnitude of this integral is $\lesssim \mathring{\updelta}^2$ 
and thus $\leq$ the right-hand side of \eqref{E:MORLOWERORDERROTATIONALGRONWALLREADY}
as desired. 
Clearly, the same bound holds for the first term 
$\enone[\cdots](0,u)$
on the right-hand side of \eqref{E:E1DIVID}.
To bound the integral over the hypersurface $\Sigma_t^u$ 
and the last spacetime integral on the right-hand side of \eqref{E:E1DIVID}
by $\leq$ the right-hand side of \eqref{E:MORLOWERORDERROTATIONALGRONWALLREADY},
we examine the proofs of inequalities
\eqref{E:EASYERRORINTEGRANDTWOINTEGRALESTIMATE} 
and
\eqref{E:EASYERRORINTEGRANDONEINTEGRALESTIMATE}
and conclude that the same inequalities hold true
with $\mathscr{O}^{\vec{I}} \Psi$ in place of $\mathscr{Z}^N \Psi$ on the left-hand side and
$\mathbb{Q}_{(\mathscr{\Rot})}$
in place of 
$\totzeromax{N}$
on the right-hand side.
Hence, the integrals under consideration
are in magnitude
$\lesssim \ln(\myexp + t) \mathbb{Q}_{(\mathscr{\Rot})}(t,u)$ as desired.
To handle the spacetime integral
$- \frac{1}{2} \int_{\mathcal{M}_{t,u}}
								\upmu \enmomtensor^{\alpha \beta}[\mathscr{O}^{\vec{I}} \Psi] 
									\left\lbrace
										\deformarg{\Mor}{\alpha}{\beta}
										- \rgeo^2 \mytr \upchi g_{\alpha \beta}
									\right\rbrace
								\, d \vol$
from the right-hand side of \eqref{E:E1DIVID},
we appeal to definition \eqref{E:COERCIVEMORDEF}
(with $\mathscr{O}^{\vec{I}} \Psi$ in the role of $\Psi$)
and \eqref{E:MORAWETZENERGYERRORINTEGRANDS},
which show that $-\Morint[\mathscr{O}^{\vec{I}}\Psi]$ is equal to 
the one of the terms that make up this integral.
We bring the term $\Morint[\mathscr{O}^{\vec{I}}\Psi]$
over to the left-hand side of \eqref{E:E1DIVID}
as a positive integral.
We then bound the remaining part of this spacetime integral 
in magnitude by 
$\leq$
the right-hand side of \eqref{E:MORLOWERORDERROTATIONALGRONWALLREADY}
with the help of \eqref{E:Q1BASICERRORINTEGRALESTIMATE}.
More precisely, examining the proof of \eqref{E:Q1BASICERRORINTEGRALESTIMATE},
we see that the
inequality holds 
true with
with $\mathscr{O}^{\vec{I}} \Psi$ in place of $\mathscr{Z}^N \Psi$ on the left-hand side and
$\mathbb{Q}_{(\mathscr{\Rot})},$
$\widetilde{\mathbb{Q}}_{(\mathscr{\Rot})},$
and $\widetilde{\mathbb{K}}_{(\mathscr{\Rot})}$
respectively in place of 
$\totzeromax{N},$
$\totonemax{N},$
and $\totMormax{N}$
on the right-hand side,
which yields the desired bounds.
We remark that using the argument given at the end of the previous paragraph,
we have bounded the last term on the right-hand side of \eqref{E:Q1BASICERRORINTEGRALESTIMATE} 
(with $\widetilde{\mathbb{K}}_{(\mathscr{\Rot})}$ in the role of $\totMormax{N}$)
by $\lesssim \mathring{\upepsilon}^3$ as desired.

To complete the proof of 
\eqref{E:MULTLOWERORDERROTATIONALGRONWALLREADY}-\eqref{E:MORLOWERORDERROTATIONALGRONWALLREADY},
it remains for us to show that
the second term on
the right-hand side of \eqref{E:E0DIVID} 
and the fourth term on the right-hand side of
\eqref{E:E1DIVID}
are bounded as follows:
\begin{align} \label{E:ROTATIONCOMMUTEDINHOMOGENEOUSERRORINTEGRALCUBICALLYSMALL}
 \left|
		\int_{\mathcal{M}_{t,u}}
			(1 + 2 \upmu) (\Lunit \mathscr{O}^{\vec{I}}\Psi) \inhomarg{\mathscr{O}^{\vec{I}}}
				+ 2 (\Rad \mathscr{O}^{\vec{I}}\Psi) \inhomarg{\mathscr{O}^{\vec{I}}} 
		\, d \vol
	\right|
	& \lesssim \mathring{\upepsilon}^3,
		\\
	\left|
	\int_{\mathcal{M}_{t,u}}
		\rgeo^2 
			\left\lbrace
				\Lunit \mathscr{O}^{\vec{I}} \Psi 
				+ \frac{1}{2} \mytr \upchi \mathscr{O}^{\vec{I}} \Psi
			\right\rbrace
			\inhomarg{\mathscr{O}^{\vec{I}}}
		\,  d \vol
	\right|
	& \lesssim \mathring{\upepsilon}^3.
		\label{E:MORAWETZROTATIONCOMMUTEDINHOMOGENEOUSERRORINTEGRALCUBICALLYSMALL}
\end{align}
The important point is that all integrands 
on the left-hand sides of
\eqref{E:ROTATIONCOMMUTEDINHOMOGENEOUSERRORINTEGRALCUBICALLYSMALL}
and \eqref{E:MORAWETZROTATIONCOMMUTEDINHOMOGENEOUSERRORINTEGRALCUBICALLYSMALL}
are in fact \emph{cubically} small. 
We stress again that we have cubic smallness only because all commutation vectorfields 
under consideration are rotations.
More precisely, to prove 
\eqref{E:ROTATIONCOMMUTEDINHOMOGENEOUSERRORINTEGRALCUBICALLYSMALL}
and \eqref{E:MORAWETZROTATIONCOMMUTEDINHOMOGENEOUSERRORINTEGRALCUBICALLYSMALL},
we first note that the estimates
\eqref{E:NTIMESCOMMUTEDWAVEINHOMOGENEOUSTERMFIRSTSPLITTING}
and
\eqref{E:LOWERORDERINHOMOGENEOUSTERMSFIRSTPOINTWISE}
imply that
for $|\vec{I}| \leq 4,$
the terms $\inhomarg{\mathscr{O}^{\vec{I}}}$
are $\lesssim$ the right-hand side of
\eqref{E:LOWERORDERINHOMOGENEOUSTERMSFIRSTPOINTWISE}
with $N=5.$
Hence, by \eqref{E:LOWERORDERDERIVATIVESOFCOMMUTATIONCURRENTSAREHARMLESS},
we deduce that for 
$|\vec{I}| \leq 4,$
the terms
$\inhomarg{\mathscr{O}^{\vec{I}}}$
are $Harmless^{\leq 5};$
see definition \eqref{E:HARMLESSTERMS} and note
that by Theorem~\ref{T:LONGTIMEPLUSESTIMATES}, the estimate 
$|\inhomarg{\mathscr{O}^{\vec{I}}}| = Harmless^{\leq 5}$
holds
with $C \mathring{\upepsilon}$ in place of $\varepsilon$
in \eqref{E:HARMLESSTERMS}.
However, in the present context, there is additional structure present
in the terms beyond the structure indicated on
the right-hand side of \eqref{E:HARMLESSTERMS}.
The new feature in the present context is
that $|\inhomarg{\mathscr{O}^{\vec{I}}}|$ 
is $\lesssim$ a modified version of the right-hand side of
\eqref{E:HARMLESSTERMS} (with $N=5$)
in which the terms on the first line of the right-hand side are
also multiplied by a $\mathring{\upepsilon}$ factor.
This factor is simply a reflection of the quadratic smallness of 
$\inhomarg{\mathscr{O}^{\vec{I}}}.$
It is straightforward to verify that thanks 
to this $\mathring{\upepsilon}$ factor,
the left-hand sides of
\eqref{E:ROTATIONCOMMUTEDINHOMOGENEOUSERRORINTEGRALCUBICALLYSMALL}
and \eqref{E:MORAWETZROTATIONCOMMUTEDINHOMOGENEOUSERRORINTEGRALCUBICALLYSMALL}
are $\lesssim \mathring{\upepsilon}$
times the terms on the right-hand sides of
\eqref{E:EASYHARMLESSMULTERRORINTEGRAL}
and \eqref{E:EASYHARMLESSMORERRORINTEGRAL}
with $N=5$
and with $C \mathring{\upepsilon}$ in place of $\varepsilon.$
Also using the estimates
\eqref{E:Q0LOWESTORDERIMPROVED} and \eqref{E:Q1LOWESTORDERIMPROVED}
to bound the terms by $\lesssim \mathring{\upepsilon}^2,$
we conclude the
desired estimates
\eqref{E:ROTATIONCOMMUTEDINHOMOGENEOUSERRORINTEGRALCUBICALLYSMALL}
and \eqref{E:MORAWETZROTATIONCOMMUTEDINHOMOGENEOUSERRORINTEGRALCUBICALLYSMALL}.

\end{proof}


\chapter{Proof of Shock Formation for Nearly Spherically Symmetric Data}
\label{C:PROOFOFSHOCKFORMATION}
\setcounter{equation}{0}
\thispagestyle{fancy}
In Chapter~\ref{C:PROOFOFSHOCKFORMATION}, 
we use the estimates of Theorem~\ref{T:LONGTIMEPLUSESTIMATES}
and Cor.~\ref{C:ANGULARDDERIVATIVEESTIMATES}
to show that there exists an open set $\shockset$
of nearly spherically symmetric small data 
that launch shock-forming solutions. 
For technical reasons to be explained,
the set $\shockset$ consists of data
on $\Sigma_{-1/2}$ that are compactly supported
in the Euclidean ball of radius $1/2$ centered at the origin.
We state the main result as Theorem~\ref{T:STABLESHOCKFORMATION}.
As we described in Sect.~\ref{SS:INTROSHOCKFORMINGDATACOMPARISON},
the theorem can be extended
to show shock formation for a significantly larger set of
compactly supported small data.

By Theorem~\ref{T:LONGTIMEPLUSESTIMATES},
in order to show that a shock forms,
we must show that $\upmu$ vanishes in finite time.
The main idea of the proof is to use the estimates 
\eqref{E:MAINTHEOREMKEYUPMUBEHAVIOR}-\eqref{E:RGEOTIMESRADPSIMAINTHEOREMSTATEMENT}
of Theorem~\ref{T:LONGTIMEPLUSESTIMATES} 
to effectively reduce the problem to proving that the term
$\rgeo G_{\Lunit \Lunit} \Rad \Psi$
from inequality \eqref{E:MAINTHEOREMKEYUPMUBEHAVIOR}
becomes negative along some integral curve of $\Lunit;$
see Sects.~\ref{S:LOWERBOUNDLEADINGTOBLOWUP} and \ref{S:UPMUHEURISTICS} for a sketch of the argument.
The main term that we have yet to suitably estimate is the product $\rgeo \Rad \Psi.$
To estimate it, we consider the wave equation in the 
form \eqref{E:WAVEEQUATIONSHOCKFORMATIONVERSION}
(with $f= \Psi$),
which, thanks to the estimates of Theorem~\ref{T:LONGTIMEPLUSESTIMATES}
and Cor.~\ref{C:ANGULARDDERIVATIVEESTIMATES},
can be viewed as a transport equation 
along the integral curves of $\Lunit = \frac{\partial}{\partial t}|_{u,\vartheta}$
with small error sources. 
The role of Cor.~\ref{C:ANGULARDDERIVATIVEESTIMATES} 
is that it allows us to prove that 
for nearly spherically symmetric data,
the linear term $\rgeo \angLap \Psi$ on the right-hand side of
\eqref{E:WAVEEQUATIONSHOCKFORMATIONVERSION} has a sufficiently small amplitude and
can be treated as a negligible error term, even
near $t=0;$ all of the remaining terms on the right-hand side are quadratic and 
are therefore much more easily seen to be negligible error terms. 
In total, for data in $\shockset,$  
this line of reasoning allows us to prove that for sufficiently large times,
the following estimate holds
along some integral curve of $\Lunit$ 
(that is, a curve with fixed geometric coordinates $(u_*,\vartheta_*)$):
\begin{align} \label{E:EXPLANATORYBLOWUPODE}
	\upmu(t,u_*,\vartheta_*)
	& \leq 
	 1 
	 + C \mathring{\upepsilon} 
	 	- \shockparameter \ln \left( \frac{\rgeo(t,u_*)}{\rgeo(0,u_*)} \right),
\end{align}
where $\shockparameter > 0$ is a small constant depending on 
a neighborhood of the data.
It easily follows from \eqref{E:EXPLANATORYBLOWUPODE} that 
$\upmu_{\star}(t,U_0) := \min\lbrace \inf_{\Sigma_t^{U_0}} \upmu, 1 \rbrace$
vanishes by the time $t \sim \exp \left( \frac{1 + C \mathring{\upepsilon}}{\shockparameter} \right).$

\section{Preliminary pointwise estimates based on transport equations}
\label{S:PRELIMINARYPOINTWISE}
We begin our detailed analysis with the following simple lemma, which shows that
for the solutions $\Psi$ under consideration, the equation
$\upmu \square_{g(\Psi)} \Psi = 0$ can be treated as a transport equation with small error sources.

\begin{lemma} [\textbf{Wave equation approximately reduces to a transport equation}]
\label{L:ALMOSTSPHERICALLYSYMMETRICWAVEEQUATION}
Assume the hypotheses and conclusions of Theorem~\ref{T:LONGTIMEPLUSESTIMATES}
and Cor.~\ref{C:ANGULARDDERIVATIVEESTIMATES}.
In particular, let $\mathring{\upepsilon}$ be the size of the data
as defined in Def.~\ref{D:SMALLDATA}, 
let $\mathring{\updelta}$ be the initial size of the geometric rotational derivatives of $\Psi$ defined in
\eqref{E:SMALLANGULARDERIVATIVES}, and let $T_{(Lifespan);U_0}$ and $U_0$ be the numbers appearing in the statement of the theorem and corollary. Then the following estimates hold on the spacetime domain $\mathcal{M}_{T_{(Lifespan);U_0},U_0}:$
		\begin{subequations}
		\begin{align} \label{E:WAVEEQUATIONAPPROXIMATETRANSPORT}
		|\Lunit \uLgood [\rgeo \Psi]|(t,u,\vartheta)
			& \leq C (\mathring{\updelta} + \mathring{\upepsilon}^{3/2})
					\frac{\ln(\myexp + t)}{(1 + t)^2},
				\\
		\left|
			[\rgeo \Rad \Psi](t,u,\vartheta)
			- \frac{1}{2} \uLgood [\rgeo \Psi](0,u,\vartheta)
		\right|
		& \leq C (\mathring{\updelta} + \mathring{\upepsilon}^{3/2})
			+ C \mathring{\upepsilon} \frac{\ln(\myexp + t)}{1 + t},
			\label{E:KEYTERMISNEARITSDATAVALUE}
				\\
		\left|
			[\rgeo \Rad \Psi](t,u,\vartheta)
			- 
			\frac{1}{2}
			\left\lbrace 
				r \mathring{\Psi}_0 
				- r \partial_r \mathring{\Psi}  
				- \mathring{\Psi} 
			\right\rbrace (u,\vartheta)
		\right|
		& \leq C (\mathring{\updelta} + \mathring{\upepsilon}^{3/2})
			+ C \mathring{\upepsilon} \frac{\ln(\myexp + t)}{1 + t},
			\label{E:KEYTERMINTERMSOFWAVEEQUATIONDATAISNEARITSDATAVALUE}
	\end{align}
	\end{subequations}
	where $\partial_r = \frac{x^a}{r}\partial_a$ is the standard Euclidean radial derivative.
\end{lemma}

\begin{proof}
	To deduce \eqref{E:WAVEEQUATIONAPPROXIMATETRANSPORT}, we use the wave equation in the form
	\eqref{E:WAVEEQUATIONSHOCKFORMATIONVERSION} with $f = \Psi.$ 
	By Cor.~\ref{C:ANGULARDDERIVATIVEESTIMATES}, 
	the linear term $\rgeo \angLap \Psi$ on the right-hand side of
	\eqref{E:WAVEEQUATIONSHOCKFORMATIONVERSION} verifies the pointwise bound
	\begin{align}	
		|\rgeo \angLap \Psi|
		& \lesssim
			(\mathring{\updelta} + \mathring{\upepsilon}^{3/2})
			\frac{\ln(\myexp + t)}{(1 + t)^2},		
	\end{align}
	which is $\lesssim$ the right-hand side of \eqref{E:WAVEEQUATIONAPPROXIMATETRANSPORT}.
	The remaining terms on the right-hand side of
	\eqref{E:WAVEEQUATIONSHOCKFORMATIONVERSION} are of the schematic form
	\begin{align}   \label{E:SMALLWAVEEQNQUADRATICTERMSSCHEMATIC}
			&
			\rgeo (\upmu - 1) \angLap \Psi
			+ (\upmu - 1) \Lunit \Psi
			+ G_{(Frame)} 
			\myarray
				[\upmu \Lunit \Psi]
				{\Rad \Psi}
			\Psi
			\\
		& \ \
			+
			\rgeo
			G_{(Frame)}
			\ginversesphere
			\threemyarray
				[\upmu \Lunit \Psi]
				{\Rad \Psi}	
				{\upmu \angdiff \Psi}	
			\myarray
				[\Lunit \Psi]
				{\angdiff \Psi}
		+ \rgeo \mytr \upchi^{(Small)} \Rad \Psi.
		\notag
\end{align}
	In particular, the terms in \eqref{E:SMALLWAVEEQNQUADRATICTERMSSCHEMATIC}
	are quadratically small, and
	from inequality
\eqref{E:ANGLAPFUNCTIONPOINTWISEINTERMSOFROTATIONS},
the estimates \eqref{E:LOWERORDERC0BOUNDLIEDERIVATIVESOFGRAME},
\eqref{E:CRUDELOWERORDERC0BOUNDDERIVATIVESOFANGULARDEFORMATIONTENSORS},
and
\eqref{E:C0BOUNDCRUCIALEIKONALFUNCTIONQUANTITIES},
the bootstrap assumptions \eqref{E:PSIFUNDAMENTALC0BOUNDBOOTSTRAP}, 
and the fact that Theorem~\ref{T:LONGTIMEPLUSESTIMATES}
implies that these estimates hold with
$\varepsilon$ replaced by $C \mathring{\upepsilon},$
we deduce that these terms are bounded in magnitude by
$\lesssim \mathring{\upepsilon}^2 \ln(\myexp + t)(1 + t)^{-2}$		
as desired. We have thus proved \eqref{E:WAVEEQUATIONAPPROXIMATETRANSPORT}.
	
	To prove \eqref{E:KEYTERMISNEARITSDATAVALUE}, we first integrate \eqref{E:WAVEEQUATIONAPPROXIMATETRANSPORT}
	along the integral curves of $\Lunit$ to deduce that
	\begin{align} \label{E:ULRADPSIKEYSHOCKESTIMATE}
		\left|
			\uLgood [\rgeo \Psi](t,u,\vartheta)
			- \uLgood [\rgeo \Psi](0,u,\vartheta)
		\right|
		& \lesssim \mathring{\updelta} + \mathring{\upepsilon}^{3/2}.
	\end{align}
	Next, noting that the identities
	$\uLgood \Psi = 2 \Rad \Psi + \upmu \Lunit \Psi$
	and
	$\uLgood \rgeo = \upmu - 2$
	imply that
	$\uLgood [\rgeo \Psi] 
	=  2 \rgeo \Rad \Psi 
		- (\upmu - 2) \Psi
		+ \rgeo \upmu \Lunit \Psi,
	$
	and using the estimates 
	$|(\upmu - 1) \Psi|, \, |\rgeo \upmu \Lunit \Psi| 
	\lesssim \mathring{\upepsilon} \ln(\myexp + t)(1 + t)^{-1},$
	which follow from the estimates stated just after \eqref{E:SMALLWAVEEQNQUADRATICTERMSSCHEMATIC},
	we deduce that 
	\begin{align} \label{E:RGEORADPSINEARHALFULRADPSI}
		\left|
			\rgeo \Rad \Psi(t,u,\vartheta)
			- 
			\frac{1}{2} \uLgood [\rgeo \Psi](t,u,\vartheta)
		\right|
		& \lesssim \mathring{\upepsilon} \frac{\ln(\myexp + t)}{1 + t}.
	\end{align}
	The desired estimate \eqref{E:KEYTERMISNEARITSDATAVALUE} now follows from
	\eqref{E:ULRADPSIKEYSHOCKESTIMATE} and \eqref{E:RGEORADPSINEARHALFULRADPSI}.
	
	Next, once we have shown that
	\begin{align} \label{E:KEYSHOCKTERMDATAEXPLICITLYEXPANDED}
		\uLgood [\rgeo \Psi](0,u,\vartheta)
		& =  r \mathring{\Psi}_0 
				- r \partial_r \mathring{\Psi}  
				- \mathring{\Psi} 
			+ \mathcal{O}(\mathring{\upepsilon}^2),
	\end{align}
	the estimate \eqref{E:KEYTERMINTERMSOFWAVEEQUATIONDATAISNEARITSDATAVALUE}
	follows from 
	inequalities
	\eqref{E:KEYTERMISNEARITSDATAVALUE} 
	and 
	\eqref{E:KEYSHOCKTERMDATAEXPLICITLYEXPANDED}.
	To prove \eqref{E:KEYSHOCKTERMDATAEXPLICITLYEXPANDED}, we first further expand
	the term $\uLgood [\rgeo \Psi]$ as follows:
	\begin{align} \label{E:EXPANSIONOFULRGEOPSI}
		\uLgood [\rgeo \Psi] 
		& = 2 \rgeo \Rad \Psi 
			+ \rgeo \Lunit \Psi
			- \Psi
			+ (\upmu - 1) \Psi
			+ \rgeo (\upmu - 1) \Lunit \Psi. 
	\end{align}
	The last two terms on the right-hand side of \eqref{E:EXPANSIONOFULRGEOPSI} 
	are quadratically small and hence from the estimates stated just after \eqref{E:SMALLWAVEEQNQUADRATICTERMSSCHEMATIC},
	they are $\lesssim \mathring{\upepsilon}^2.$
	By virtue of the estimates
	\begin{align}
		\rgeo(0,u,\vartheta) 
		& = r(0,u,\vartheta) = 1 - u,
			\\
		\Rad^0 & = 0,
			\\
		\Rad^i(0,u,\vartheta)
		& = - \frac{x^i}{r} + \mathcal{O}(\mathring{\upepsilon}),
			\label{E:RADIDATAEXPANSION} \\
		\Lunit^0 & = 1,
			\\
		\Lunit^i(0,u,\vartheta) 
		& = \frac{x^i}{r} + \mathcal{O}(\mathring{\upepsilon})
			\label{E:LUNITIDATAEXPANSION}
	\end{align}
	(where \eqref{E:RADIDATAEXPANSION} and \eqref{E:LUNITIDATAEXPANSION} follow from 
	Def.~\ref{D:RENORMALIZEDVARIABLES} and Lemma~\ref{L:INITIALEXPRESSIONSFORUPMUANDLJUNKI}),
	we can further expand the first two terms on the right-hand side of \eqref{E:EXPANSIONOFULRGEOPSI} as
	\begin{align}
		[2 \rgeo \Rad \Psi](0,u,\vartheta) 
		& = - 2 r \partial_r \Psi(0,u,\vartheta)
			+ \mathcal{O}(\mathring{\upepsilon}),
				\label{E:RRADPSIDATAEXPANDED} \\
		[\rgeo \Lunit \Psi](0,u,\vartheta) 
		& = r \partial_t \Psi(0,u,\vartheta) 
			+ r \partial_r \Psi(0,u,\vartheta)
			+ \mathcal{O}(\mathring{\upepsilon}).
			\label{E:RLUNITPSIDATAEXPANDED}
	\end{align}
	The estimate \eqref{E:KEYSHOCKTERMDATAEXPLICITLYEXPANDED} now follows from
	\eqref{E:EXPANSIONOFULRGEOPSI}, the sentence after \eqref{E:EXPANSIONOFULRGEOPSI}, 
	\eqref{E:RRADPSIDATAEXPANDED}, 
	and
	\eqref{E:RLUNITPSIDATAEXPANDED}.
\end{proof}

As we will see in Lemma~\ref{L:KEYESTIMATEFORTHETERMTHATDRIVESTHESHOCKFORMATION},
the following data-dependent function 
appearing in Lemma~\ref{L:ALMOSTSPHERICALLYSYMMETRICWAVEEQUATION}
plays a central role in our proof of shock formation.

\begin{definition} [\textbf{The data-dependent function that drives shock formation}]
	\label{D:SHOCKFORMINGFUNCTION}
Let 
$\uLunit_{(Flat)} := \partial_t - \partial_r$ denote the standard Minkowskian ingoing
null vectorfield.
Let $(\mathring{\Psi} := \Psi|_{\Sigma_0^{U_0}}, \mathring{\Psi}_0 := \partial_t \Psi|_{\Sigma_0^{U_0}})$ 
be initial data for the wave equation \eqref{E:WAVEGEO}.
We define the following data-dependent function: 
\begin{align} \label{E:INITIALDATASHOCKFORMATIONFUNCTION}
	\shockfunction[(\mathring{\Psi}, \mathring{\Psi}_0)] 
	& :=
		\frac{1}{4} 
		\uLunit_{(Flat)} (r \Psi)|_{(\Psi,\partial_t \Psi) = (\mathring{\Psi},\mathring{\Psi}_0)}
		=
		 \frac{1}{4}
		 \left\lbrace 
		 	r \mathring{\Psi}_0 - r \partial_r \mathring{\Psi}  - \mathring{\Psi} 
		\right\rbrace.
\end{align}
\end{definition}

In the next lemma, we insert the estimates of 
Lemma~\ref{L:ALMOSTSPHERICALLYSYMMETRICWAVEEQUATION}
into some estimates proved in Theorem~\ref{T:LONGTIMEPLUSESTIMATES}.
We arrive at the key pointwise estimate \eqref{E:KEYUPMUVERYSHARPPOINTWISEESTIMATE} 
verified by $\upmu,$ which is the main ingredient
in our proof of finite-time shock formation.  

\begin{lemma}[\textbf{The key estimate that drives the shock formation}] \label{L:KEYESTIMATEFORTHETERMTHATDRIVESTHESHOCKFORMATION}
Let $(\mathring{\Psi} := \Psi|_{\Sigma_0^{U_0}}, \mathring{\Psi}_0 := \partial_t \Psi|_{\Sigma_0^{U_0}})$ 
be initial data for the wave equation \eqref{E:WAVEGEO}.
Assume the hypotheses and conclusions of Lemma~\ref{L:ALMOSTSPHERICALLYSYMMETRICWAVEEQUATION}
and Lemma~\ref{L:KEYESTIMATEFORTHETERMTHATDRIVESTHESHOCKFORMATION}.
In particular, let $\mathring{\upepsilon}$
be the size of the data as defined in Def.~\ref{D:SMALLDATA},  
and let $\mathring{\updelta}$ be the initial size of the geometric rotational derivatives of $\Psi$ 
defined by \eqref{E:SMALLANGULARDERIVATIVES}.
Let $\InitialFutFailFac$ be as in Def.~\ref{D:FAILUREFACTOR} 
and let $\shockfunction[(\mathring{\Psi}, \mathring{\Psi}_0)](\cdots),$
viewed as a function of the geometric coordinates $(u,\vartheta),$
be as in Def.~\ref{E:INITIALDATASHOCKFORMATIONFUNCTION}.
Then the following pointwise estimate holds
on the spacetime domain $\mathcal{M}_{T_{(Lifespan);U_0},U_0}:$
\begin{align} \label{E:KEYUPMUVERYSHARPPOINTWISEESTIMATE}
	\left|
		\upmu(t,u,\vartheta)
	 	- 1 
	 	-
	 		\InitialFutFailFac(\vartheta)
	 		\shockfunction[(\mathring{\Psi}, \mathring{\Psi}_0)](u,\vartheta)
			\ln \left( \frac{\rgeo(t,u)}{\rgeo(0,u)} \right)
	\right|
	& \leq
			C \left(\mathring{\updelta} + \mathring{\upepsilon}^{3/2} \right)
				\ln \left( \frac{\rgeo(t,u)}{\rgeo(0,u)} \right)
		+ C \mathring{\upepsilon}.
\end{align}
\end{lemma}

\begin{proof}
	The estimate \eqref{E:KEYUPMUVERYSHARPPOINTWISEESTIMATE}
	follows from 
	the estimates
	\eqref{E:MAINTHEOREMKEYUPMUBEHAVIOR}-\eqref{E:RGEOTIMESRADPSIMAINTHEOREMSTATEMENT}
	and
	\eqref{E:KEYTERMINTERMSOFWAVEEQUATIONDATAISNEARITSDATAVALUE}.
	
\end{proof}

\section{Existence of small, stable, shock-forming data}
\label{S:SHOCKFORMINGDATA}
In Lemma~\ref{L:EXISTENCEOFSHOCKFORMINGDATA}, 
we construct an open set of small initial data with properties
that lead to finite-time shock formation in their solutions. 
We provide the proof of shock formation 
in the subsequent Theorem~\ref{T:STABLESHOCKFORMATION},
which is located in Sect.~\ref{S:STABLESHOCKFORMATION}.

As a preliminary step, we prove the following lemma, 
which is based on an argument of John given in \cite{fJ1985}.
The lemma sets the stage for proving our shock formation
result, roughly by showing that compactly supported spherically symmetric data
always induce the ``shock-driving'' sign in the quantity
that drives shock formation.
It is a spherically symmetric
analog of Prop.~\ref{P:JOHNSCRITERIONISALWAYSSATISFIEDFORCOMPACTLYSUPPORTEDDATA}
for wave equations of the form \eqref{E:WAVEGEO}.

\begin{lemma}[\textbf{Nontrivial spherically symmetric data force} $\timeminushalfshockfunction$ \textbf{to take on both signs}]
 	\label{L:SPHERICALLYSYMMETRICDATAGIVESHOCKFUNCTIONBOTHSIGNS}
 	Let $(\halfdata, \halftimedata) \in C^1(\mathbb{R}^3) \times C^0(\mathbb{R}^3)$ 
 	be data for the wave equation \eqref{E:WAVEGEO}
 	\textbf{given along} $\Sigma_{-1/2}$
 	with support contained in the Euclidean ball of radius $1/2$ centered at the origin. 
 	Assume that the data are spherically symmetric in the sense 
	that they depend on only\footnote{Since equation \eqref{E:WAVEGEO} is not generally invariant under Euclidean rotations,
	spherically symmetric data do not generally launch spherically symmetric solutions.} 
	the Euclidean radial variable $r$ and not on the Euclidean angular variable $\theta,$
	that is, $\halfdata = \halfdata(r)$ and $\halftimedata = \halftimedata(r).$ 
	For $p \in [-1/2,1/2],$ we define the data-dependent function
	$\timeminushalfshockfunction[(\halfdata, \halftimedata)]$ by
	\begin{align} \label{E:SEXTENDED}
	\timeminushalfshockfunction[(\halfdata, \halftimedata)](p)
	: =
	\frac{1}{4}
		 \left\lbrace 
		 	p \halftimedata(|p|) - |p| \halfdata'(|p|) - \halfdata(|p|) 
		\right\rbrace.
	\end{align}
	Then if the data are nontrivial, we have
	\begin{align} \label{E:SPHERICALLYSYMMETRICDATAGIVESHOCKFUNCTIONBOTHSIGNS}
		\inf_{p \in [-1/2,1/2]} \timeminushalfshockfunction[(\halfdata, \halftimedata)](p) 
		< 0 
		< \sup_{p \in [-1/2,1/2]} \timeminushalfshockfunction[(\halfdata, \halftimedata)](p).
	\end{align}
\end{lemma}

\begin{remark}[\textbf{Relationship between $\timeminushalfshockfunction$ and $\shockfunction$}]
	\label{R:RELATIONSHIPBETWEENTWOSHOCFUNCTIONS}
	Note that for $p > 0,$
	$\timeminushalfshockfunction[(\halfdata, \halftimedata)](p)
	= \shockfunction[(\halfdata, \halftimedata)](p),$
	where the latter function is defined by \eqref{E:INITIALDATASHOCKFORMATIONFUNCTION}.
	For $p < 0,$
	$\timeminushalfshockfunction[(\halfdata, \halftimedata)](p)$
	can be viewed as an extension of 
	$\shockfunction[(\halfdata, \halftimedata)](\cdot)$
	constructed by extending the data 
	$(\halfdata, \halftimedata)$
	to be even functions of their argument.
\end{remark}

\begin{proof}[Proof of Lemma \ref{L:SPHERICALLYSYMMETRICDATAGIVESHOCKFUNCTIONBOTHSIGNS}]
Throughout, we suppress the dependence of $\timeminushalfshockfunction$ on $(\halfdata, \halftimedata).$
We view the data as \emph{even} functions of 
the real variable $p \in [-1/2,1/2]$ that verify $(\halfdata, \halftimedata) \in C^1([-1/2,1/2]) \times C^0([-1/2,1/2]).$
Note that 
$4 \timeminushalfshockfunction(p) =
	p \halftimedata(p) - \partial_p (p \halfdata(p)).
$
It easily follows from this identity and the previous discussion that
\begin{align} \label{E:SINTEGRATESTO0}
	\int_{p = -1/2}^{1/2} \timeminushalfshockfunction(p) \, dp = 0.
\end{align}
We now show that if
$\min_{p \in [-1/2,1/2]} \timeminushalfshockfunction(p) > 0,$
then the data must completely vanish.
We therefore assume that $\min_{p \in [-1/2,1/2]} \timeminushalfshockfunction(p) > 0.$
It then follows from \eqref{E:SINTEGRATESTO0} that $\timeminushalfshockfunction(p) \equiv 0$ for $p \in [-1/2,1/2].$
But since $\halfdata(p)$ and $\halftimedata(p)$ are even functions,
we have the following identities:
\begin{align}
	0 & = 4 \left\lbrace \timeminushalfshockfunction(p) + \timeminushalfshockfunction(-p) \right\rbrace 
		= - 2 \partial_p (p \halfdata),
		\label{E:DERIVATIVEOFRTIMESDATAIS0} \\
	0 & = 4 \left\lbrace \timeminushalfshockfunction(p) - \timeminushalfshockfunction(-p) \right\rbrace
		= 2 p \halftimedata(p).
			\label{E:RTIMESSECONDDATA}
\end{align}
Since the data are compactly supported, we easily conclude from 
\eqref{E:DERIVATIVEOFRTIMESDATAIS0} and \eqref{E:RTIMESSECONDDATA} that they
are trivial. We have thus proved that the first inequality in \eqref{E:SPHERICALLYSYMMETRICDATAGIVESHOCKFUNCTIONBOTHSIGNS}
must hold when the data are nontrivial. The second inequality can be proved in a similar fashion.
\end{proof}

In the next lemma, we show that
for small spherically symmetric data given on $\Sigma_{-1/2},$
the function $\timeminushalfshockfunction[(\halfdata, \halftimedata)]$
effectively determines the function
$\shockfunction[(\mathring{\Psi}, \mathring{\Psi}_0)]$
corresponding to the ``data'' induced by the solution on $\Sigma_0^1.$

\begin{lemma}[\textbf{Quantitative relationship between $\shockfunction[(\mathring{\Psi}, \mathring{\Psi}_0)]$
	and $\timeminushalfshockfunction[(\halfdata, \halftimedata)]$}]
	\label{L:PROPAGATINGSHOCKFUNCTION}
	Let $(\halfdata, \halftimedata)$ be spherically symmetric data
	for the wave equation \eqref{E:WAVEGEO}
	\textbf{given along} $\Sigma_{-1/2}$
 	with support contained in the Euclidean ball of radius $1/2$
 	centered at the origin, as in Lemma \ref{L:SPHERICALLYSYMMETRICDATAGIVESHOCKFUNCTIONBOTHSIGNS}.
	Let 
	$\epsilon 
	= \epsilon[(\halfdata, \halftimedata)] 
	:= 
	\| \halfdata \|_{H_{\Euct}^{25}(\Sigma_{-1/2})} + \| \halftimedata \|_{H_{\Euct}^{24}(\Sigma_{-1/2})} < \infty$
	denote their size,\footnote{As always in this monograph, $H_{\Euct}^N$ is a standard Euclidean Sobolev space involving 
	rectangular spatial derivatives.}  
	and let $\Psi$ denote the corresponding solution.
	Let $(\mathring{\Psi} := \Psi|_{\Sigma_0^1}, \mathring{\Psi}_0 := \partial_t \Psi|_{\Sigma_0^1})$
	be the ``data'' along $\Sigma_0$ induced by $\Psi.$ Then if $\epsilon$ is sufficiently small,
	$\Psi$ persists beyond time $0$ and along $\Sigma_0^1,$ we have
	\begin{align} \label{E:NORMDOESNOTGROWTOOMUCH}
		\mathring{\upepsilon}[(\mathring{\Psi},\mathring{\Psi}_0)]
		:=
		\| \mathring{\Psi} \|_{H_{\Euct}^{25}(\Sigma_0^1)} + \| \mathring{\Psi}_0 \|_{H_{\Euct}^{24}(\Sigma_0^1)} 
		\lesssim \epsilon.
	\end{align}
	
	Furthermore, 
	let $\shockfunction[(\mathring{\Psi}, \mathring{\Psi}_0)]$ 
	be the data-dependent function from Def.~\ref{D:SHOCKFORMINGFUNCTION}
	viewed as a function of
	standard spherical coordinates $(r,\theta)$ on $\Sigma_0^1,$
	and let $\timeminushalfshockfunction[(\halfdata, \halftimedata)]$
	be the function of the real variable $p \in [-1/2,1/2]$ as defined by
	\eqref{E:SEXTENDED}.
	If $\epsilon$ is sufficiently small, then the following estimates hold for 
	$(r,\theta) \in [0,1] \times \mathbb{S}^2:$ 
	\begin{align} \label{E:SFUNCTIONSARECLOSE}
		\shockfunction[(\mathring{\Psi}, \mathring{\Psi}_0)](r,\theta)
		& = \timeminushalfshockfunction[(\halfdata, \halftimedata)](r - 1/2)
		+ \mathcal{O}(\epsilon^2),
	\end{align}
	where the implicit constants in $\mathcal{O}$ can be chosen to be independent of 
	$(r,\theta).$
\end{lemma}

\begin{proof}
	The persistence result
	and inequality \eqref{E:NORMDOESNOTGROWTOOMUCH} are standard results
	that can be proved using the same techniques 
	that we described in sketch of the proof of the local well-posedness proposition
	(that is, Prop.~\ref{P:CLASSICALLOCAL}).
	More generally, for $t \in [-1/2,0],$ we have the following standard estimate:
	\begin{align} \label{E:PROOFNORMDOESNOTGROWTOOMUCH}
		\sum_{M=0}^{25} \| \partial_t^M \Psi \|_{H_{\Euct}^{25-M}(\Sigma_t)}  
		\lesssim \epsilon.
	\end{align}
	
	Next, we reduce the analysis to that of an effectively 
	spherically symmetric problem over the time interval $t \in [-1/2,0].$
	To this end, we first write the nonlinear wave equation \eqref{E:WAVEGEO}
	as a perturbation of the spherically symmetric linear wave equation as follows:
	\begin{align} \label{E:NONLINEARISLINEARPLUSJUNK}
		 \partial_t^2 (r \Psi)
		- \partial_r^2 (r \Psi)
		& = r \angFlatLap \Psi
			+ \mathcal{N}(\Psi, \partial^2 \Psi)
			+ \mathcal{N}(\Psi)(\partial \Psi, \partial \Psi),
	\end{align}
	where $\angFlatLap$ is the angular Laplacian corresponding
	to the metric $\msphere$ induced 
	by the Minkowski metric $m$
	on the Euclidean spheres of constant $t$ and $r$ 
	and $\mathcal{N}(\cdot,\cdot),\mathcal{N}(\Psi)(\cdot,\cdot)$
	are nonlinear terms that are quadratic in their arguments.
	We now show that for $t \in [-1/2,0],$ we have
	\begin{align} \label{E:ALMOSTRADIAL}
		\partial_t^2 (r \Psi)
		- \partial_r^2 (r \Psi)
		& = \mathcal{O}(\epsilon^2).
	\end{align}
	To prove \eqref{E:ALMOSTRADIAL}, we first
	deduce from \eqref{E:PROOFNORMDOESNOTGROWTOOMUCH} 
	that for $t \in [-1/2,0],$ 
	the two quadratic terms $\mathcal{N}$ on the right-hand side of \eqref{E:NONLINEARISLINEARPLUSJUNK}
	are $\mathcal{O}(\epsilon^2).$
	To show that the term $r \angFlatLap \Psi$ is $\mathcal{O}(\epsilon^2),$
	we can argue as in the proof of Cor.~\ref{C:ANGULARDDERIVATIVEESTIMATES}, but
	using the Euclidean rotations $\Roteuc$ 
	(see Def.~\ref{D:ROTATION})
	as commutators instead of the geometric ones.
	More precisely, since the data are spherically symmetric,
	we have $\| \mathscr{O}_{(Flat)}^{\leq 25} \Psi \|_{L_{\Euct}^2(\Sigma_t)} = 0,$
	where $\mathscr{O}_{(Flat)}$
	denotes the set of (three) Euclidean rotations.
	Furthermore, for $t \in [-1/2,0],$ we have the following estimate:
	\begin{align} \label{E:EUCLIDEANROTSMALLFORSHORTNEGATIVETIME}
		\sum_{\alpha = 0}^3 \| \partial_{\alpha} \mathscr{O}_{(Flat)}^{\leq 24} \Psi \|_{L_{\Euct}^2(\Sigma_t)} \lesssim \epsilon^2,
	\end{align}
	To prove \eqref{E:EUCLIDEANROTSMALLFORSHORTNEGATIVETIME}, 
	we commute equation \eqref{E:WAVEGEO}
	with $\leq 24$ Euclidean rotations, derive standard Minkowskian energy estimates, 
	and use the fact that all error terms in the commuted equations are quadratic
	(in particular, they vanish in the case of the linear wave equation).
	The desired estimate $r \angFlatLap \Psi = \mathcal{O}(\epsilon^2)$
	then follows from 
	\eqref{E:EUCLIDEANROTSMALLFORSHORTNEGATIVETIME},
	the Minkowskian analog of Lemma~\ref{L:POINTWISEANGDINTERMSOFANGLIEO},
	and Sobolev embedding. We have thus shown \eqref{E:ALMOSTRADIAL}.
	
	We now recall (see Def.~\ref{D:SHOCKFORMINGFUNCTION}) that 
	$
	\shockfunction[(\mathring{\Psi},\mathring{\Psi}_0)] 
	= 
	\frac{1}{4}
	\uLunit_{(Flat)} (r \Psi)|_{(\Psi,\partial_t \Psi) = (\mathring{\Psi}, \mathring{\Psi}_0)}.
	$
	In view of this identity and \eqref{E:SEXTENDED},
	we see that the desired estimate \eqref{E:SFUNCTIONSARECLOSE} will follow once
	we establish the following estimate for $t \in [-1/2,0]$ 
	(and set $t=0$ to conclude \eqref{E:SFUNCTIONSARECLOSE}):
	\begin{align} \label{E:ULUNITFLATOFRPSIALMOSTRADIALFORMULA}
		[\uLunit_{(Flat)} (r \Psi)](t,r,\theta)
		& =  
		 (r - t - 1/2) \halftimedata(|1/2 + t - r|)	
			- |(1/2 + t - r)| \halfdata'(|1/2 + t - r|)
			\\
	& \ \ 
			- \halfdata(|1/2 + t - r|)
			+ \mathcal{O}(\epsilon^2).
	\notag
	\end{align}
	To prove \eqref{E:ULUNITFLATOFRPSIALMOSTRADIALFORMULA},
	we view the solution $\Psi$ of \eqref{E:ALMOSTRADIAL} 
	as the sum of two solutions: \textbf{i)} a truly spherically symmetric solution
	to the homogeneous linear wave equation with data $(\halfdata,\halftimedata)$
	and \textbf{ii)} a solution to the inhomogeneous equation \eqref{E:ALMOSTRADIAL}
	(whose source term is size $\epsilon^2$) but with trivial data.
	The desired expression \eqref{E:ULUNITFLATOFRPSIALMOSTRADIALFORMULA}
	then follows from the standard formula for
	solutions to the spherically symmetric linear wave equation in 
	$3$ spatial dimensions in terms of its data, 
	where in the formula \eqref{E:ULUNITFLATOFRPSIALMOSTRADIALFORMULA}, we have accounted for
	the fact that the data $(\halfdata,\halftimedata)$ are given at $t = - 1/2.$
	We stress that the inhomogeneous solution $\textbf{ii)}$ is completely responsible for the
	term $\mathcal{O}(\epsilon^2)$ in \eqref{E:ULUNITFLATOFRPSIALMOSTRADIALFORMULA}.
	
\end{proof}

We now construct an open set $\shockset$ of small data 
along $\Sigma_{-1/2}$ that lead to finite-time shock formation.

\begin{lemma}[\textbf{Small, nontrivial spherically symmetric data are stable, shock-forming data}] \label{L:EXISTENCEOFSHOCKFORMINGDATA}
	We assume that the null condition fails 
	in the covariant wave equation \eqref{E:WAVEGEO}, and we define 
	\begin{align} \label{E:NULLCONDITIONFAILREFACTORMAX}
		\FutFailFac_* := \max_{\vartheta \in \mathbb{S}^2} |\InitialFutFailFac(\vartheta)| > 0,
	\end{align}
	where $\InitialFutFailFac$ is as in Def.~\ref{D:FAILUREFACTOR}.
	Let $\vartheta_* \in \mathbb{S}^2$ be any angular direction such that 
	$|\InitialFutFailFac(\vartheta_*)| = \FutFailFac_*.$
	Then there exists a non-empty \textbf{open}\footnote{
	By ``open,'' we mean open relative to the topology corresponding to the function space 
	$H_{\Euct}^{25}(\Sigma_{-1/2}) \times H_{\Euct}^{24}(\Sigma_{-1/2}).$} 
	set 
	$\shockset \subset H_{\Euct}^{25}(\Sigma_{-1/2}) \times H_{\Euct}^{24}(\Sigma_{-1/2})$
	of $\vartheta_*-$dependent data
	\textbf{given along} 
	$\Sigma_{-1/2}$
 	with support contained in the Euclidean ball of radius $1/2$
	and featuring the following properties.
	In the statements below,
	$(\halfdata, \halftimedata)$ denotes data on $\Sigma_{-1/2},$
	$(\mathring{\Psi} := \Psi|_{\Sigma_0}, \mathring{\Psi}_0 := \partial_t \Psi|_{\Sigma_0})$ 
	denotes the ``data'' induced on $\Sigma_0^1$ by the solution $\Psi,$
	and $\mathring{\upepsilon} = \mathring{\upepsilon}[(\mathring{\Psi},\mathring{\Psi}_0)]
	:= \| \mathring{\Psi} \|_{H_{\Euct}^{25}(\Sigma_0^1)}
	+ \| \mathring{\Psi}_0 \|_{H_{\Euct}^{24}(\Sigma_0^1)}$
	denotes the size of the latter, as in \eqref{E:SMALLDATA}.
\begin{enumerate}
	\item $\shockset$ is by definition the union of
		open neighborhoods $\Omega[(\halfdata^{(Radial)}, \halftimedata^{(Radial)})]$
		of spherically symmetric data $(\halfdata^{(Radial)}, \halftimedata^{(Radial)}),$
		as described in item \textbf{3}.
	\item If $(\halfdata^{(Radial)}, \halftimedata^{(Radial)}) \in H_{\Euct}^{25}(\Sigma_{-1/2}) \times H_{\Euct}^{24}(\Sigma_{-1/2})$ 
		are nontrivial and spherically symmetric in the sense that they depend on only		 
		the Euclidean radial variable $r,$ 
		and $\uplambda > 0$ is sufficiently small 
		(where the required smallness depends on $(\halfdata^{(Radial)}, \halftimedata^{(Radial)})$),
		then $(\uplambda \halfdata^{(Radial)}, \uplambda \halftimedata^{(Radial)}) \in \shockset.$
	\item 
		For any spherically symmetric data $(\halfdata^{(Radial)}, \halftimedata^{(Radial)}) \in \shockset,$
		there exist 
		a number $U_0 \in (0,1),$ 
		an open (in $H_{\Euct}^{25}(\Sigma_{-1/2}) \times H_{\Euct}^{24}(\Sigma_{-1/2})$)
		subset $\Omega[(\halfdata^{(Radial)}, \halftimedata^{(Radial)})] \subset \shockset$ containing it,
		a data-dependent constant 
		$\shockparameter := \shockparameter[(\halfdata^{(Radial)}, \halftimedata^{(Radial)})] > 0,$ 
		and a number $u_* \in (0, U_0)$ 
		such that if 
		$(\halfdata, \halftimedata) \in \Omega[(\mathring{\Psi}^{(Radial)}, \mathring{\Psi}_0^{(Radial)})],$
		then the following inequalities hold for 
		the data $(\mathring{\Psi}, \mathring{\Psi}_0)$ induced on $\Sigma_0^{U_0}$ by the solution $\Psi:$
		\begin{align} \label{E:KEYDATABLOWUPTERMHASBLOWUPSIGNANDSIZE}
			\InitialFutFailFac(\vartheta_*)
			\shockfunction[(\mathring{\Psi}, \mathring{\Psi}_0)](u_*,\vartheta_*) + C \left(\mathring{\updelta} + \mathring{\upepsilon}^{3/2} 		
				\right) & <
				- \shockparameter
				< 0.
			\end{align}
		In \eqref{E:KEYDATABLOWUPTERMHASBLOWUPSIGNANDSIZE},
		$\shockfunction[(\mathring{\Psi}, \mathring{\Psi}_0)]$ is the function defined in
		\eqref{E:INITIALDATASHOCKFORMATIONFUNCTION}
		and $C \left(\mathring{\updelta} + \mathring{\upepsilon}^{3/2} \right)$ 
		is the factor in the first product on the right-hand side of \eqref{E:KEYUPMUVERYSHARPPOINTWISEESTIMATE}.
		Furthermore, $\shockparameter > 0$ 
		(see Lemma \ref{L:SPHERICALLYSYMMETRICDATAGIVESHOCKFUNCTIONBOTHSIGNS})
		is defined by
		\begin{align} \label{E:CONSTANTDEPENDSONDATANEIGHBORHOOD}
				\shockparameter
				:=
				\begin{cases}
					\frac{1}{2}
					\FutFailFac_*
					\sup_{p \in [-1/2,1/2]}
					\timeminushalfshockfunction[(\halfdata^{(Radial)}, \halftimedata^{(Radial)})](p)
					& 
					\mbox{if} \ \InitialFutFailFac(\vartheta_*) < 0,
						\\
					-
					\frac{1}{2}
					\FutFailFac_*
					\inf_{p \in [-1/2,1/2]}
					\timeminushalfshockfunction[(\halfdata^{(Radial)}, \halftimedata^{(Radial)})](p)
					& 
					\mbox{if} \ \InitialFutFailFac(\vartheta_*) > 0.
				\end{cases}
		\end{align}
		There exists a constant $C > 0$ such that
		\begin{align} \label{E:SHOCKPARAMETERBOUNDEDBYDATAATTIMEZERO}
			\shockparameter
			& \leq C \mathring{\upepsilon}.
		\end{align}
	\item For all data in $\shockset,$ 
	the ``data'' $(\mathring{\Psi} := \Psi|_{\Sigma_0^1}, \mathring{\Psi}_0 :=  \partial_t \Psi|_{\Sigma_0^1})$ induced
	by the solution $\Psi$ are small and nearly spherically symmetric in the sense that
	the hypotheses and conclusions of Theorem~\ref{T:LONGTIMEPLUSESTIMATES}
	and Cor.~\ref{C:ANGULARDDERIVATIVEESTIMATES} are verified.
	\item All constants can depend on $U_0.$
\end{enumerate}
\end{lemma}

\begin{proof}
	The main idea of the proof is to first use an approximate scaling argument to show
	that given any spherically symmetric data $(\halfdata^{(Radial)}, \halftimedata^{(Radial)})$ on $\Sigma_{-1/2},$ 
	we can rescale it so that the desired inequalities are verified by the rescaled data.
	We then show that the inequalities are stable under small general perturbations.
	The set $\shockset$ will therefore be the union of 
	open neighborhoods of rescaled spherically symmetric data.
	
	We assume throughout the proof that $\InitialFutFailFac(\vartheta_*) < 0;$
	the case $\InitialFutFailFac(\vartheta_*) > 0$ can be treated using similar arguments.
	To begin, we let $(\halfdata^{(Radial)}, \halftimedata^{(Radial)}) \in H_{\Euct}^{25}(\Sigma_{-1/2}) \times 
	H_{\Euct}^{24}(\Sigma_{-1/2})$ 
	be spherically symmetric data on $\Sigma_{-1/2}$ supported in the Euclidean ball of radius 
	$1/2$ centered at the origin.
	For $\uplambda > 0,$ let $\scaledPsi{\uplambda}$ denote the solution
	launched by the rescaled data $(\uplambda \halfdata^{(Radial)}, \uplambda \halftimedata^{(Radial)}),$
	and let
	$(\scaleddata{\uplambda} := \scaledPsi{\uplambda}|_{\Sigma_0^1}, 
	\scaledtimedata{\uplambda} := \partial_t \scaledPsi{\uplambda}|_{\Sigma_0^1})$
	denote the ``data'' induced on $\Sigma_0^1$ by $\scaledPsi{\uplambda}.$
	For the remainder of the proof, we will assume that $\uplambda$ is sufficiently small,
	where the required smallness is allowed to depend on $(\halfdata^{(Radial)}, \halftimedata^{(Radial)}).$
	From
	\eqref{E:SPHERICALLYSYMMETRICDATAGIVESHOCKFUNCTIONBOTHSIGNS},
	\eqref{E:SFUNCTIONSARECLOSE},	
	and the fact that 
	$\timeminushalfshockfunction[\cdots]$ scales linearly in its functional argument $\cdots,$
	we deduce that there exist a parameter $U_0 \in (0,1)$ and a number $u_* \in (0,U_0)$ such that
	relative to the geometric coordinates $(u,\vartheta)$ along $\Sigma_0^1,$ we have
	\begin{align} \label{E:LAMBDASQUAREDERROR}
		\InitialFutFailFac(\vartheta_*)
		\shockfunction[(\scaleddata{\uplambda}, \scaledtimedata{\uplambda})](u_*,\vartheta_*)
		& = - \uplambda 
				\FutFailFac_*
				\sup_{p \in [-1/2,1/2]}
				\timeminushalfshockfunction[(\halfdata^{(Radial)}, \halftimedata^{(Radial)})](p)
				+ \mathcal{O}(\uplambda^2)
				< 0.
	\end{align}
	We now claim that the 
	factor $C \left(\mathring{\updelta} + \mathring{\upepsilon}^{3/2} \right)$ 
	in the first product on the right-hand side of \eqref{E:KEYUPMUVERYSHARPPOINTWISEESTIMATE}
	is $\mathcal{O}(\uplambda^{3/2}).$
	It then follows from this estimate
	and \eqref{E:LAMBDASQUAREDERROR}
	that inequality \eqref{E:KEYDATABLOWUPTERMHASBLOWUPSIGNANDSIZE} 
	holds for the rescaled spherically symmetric data 
	$(\halfdata^{(Radial)}, \halftimedata^{(Radial)}),$
	where $\shockparameter$ is defined in \eqref{E:CONSTANTDEPENDSONDATANEIGHBORHOOD}
	(note the factor of $1/2$ in the definition).
	It remains for us to establish that $C \left(\mathring{\updelta} + \mathring{\upepsilon}^{3/2} \right)$
	is $\lesssim \uplambda^{3/2}.$
	To bound the factor $C \mathring{\upepsilon}^{3/2}$
	by $\lesssim \uplambda^{3/2},$ we simply use the estimate
	\eqref{E:NORMDOESNOTGROWTOOMUCH}.
	
	We now show that 
	the quantity $\mathring{\updelta},$ which is defined by \eqref{E:SMALLANGULARDERIVATIVES}
	and appears in the estimate \eqref{E:KEYUPMUVERYSHARPPOINTWISEESTIMATE}, 
	verifies the desired bound
	\begin{align} \label{E:SPHERICALLYSYMMETRICDATAHAVESMALLDELTA}
		\mathring{\updelta} \lesssim \uplambda^2.
	\end{align}
	To prove \eqref{E:SPHERICALLYSYMMETRICDATAHAVESMALLDELTA}, we first use
	\eqref{E:EUCLIDEANROTSMALLFORSHORTNEGATIVETIME} to deduce that
	\begin{align} \label{E:EUCLIDEANROTSMALLSHOCFORMINGPROOF}
		\sum_{\alpha = 0}^3 \| \partial_{\alpha} \mathscr{O}_{(Euc)}^{\leq 24} \scaledPsi{\uplambda} \|_{L_{\Euct}^2(\Sigma_0^1)} 
		\lesssim \uplambda^2.
	\end{align}
	Furthermore, from \eqref{E:NORMDOESNOTGROWTOOMUCH}, we deduce that
	\begin{align} \label{E:FULLDATAARESMALL}
		\| \scaledPsi{\uplambda} \|_{H_{\Euct}^{25}(\Sigma_0^1)} 
		+ \| \partial_t \scaledPsi{\uplambda} \|_{H_{\Euct}^{24}(\Sigma_0^1)}
		& \lesssim \uplambda,
	\end{align}
	and thus the solution is small along $\Sigma_0.$
	Hence, it follows from the analysis of Sect.~\ref{S:INITIALBEHAVIOROFQUANTITIES} that along $\Sigma_0^1,$
	all of the tensorfields defined throughout the monograph
	have rectangular components equal to their
	Minkowskian counterparts plus an error bounded by $\mathcal{O}(\uplambda),$
	that is, $\upmu = 1 + \mathcal{O}(\uplambda),$
	$\Lunit^{\nu} = \Lunit_{(Flat)}^{\nu} + \mathcal{O}(\uplambda),$
	$\Rad^{\nu} = \Radunit_{(Flat)}^{\nu} + \mathcal{O}(\uplambda),$
	$\Rot^i = \Roteuc^i + \mathcal{O}(\uplambda),$ etc.
	Furthermore, similar estimates hold for the rectangular derivatives of these components.
	The desired bound \eqref{E:SPHERICALLYSYMMETRICDATAHAVESMALLDELTA}
	easily follows from \eqref{E:EUCLIDEANROTSMALLSHOCFORMINGPROOF} and these observations.
	
	To prove \eqref{E:SHOCKPARAMETERBOUNDEDBYDATAATTIMEZERO},
	we solve, relative to the rectangular coordinates, 
	the wave equation \eqref{E:WAVEGEO}
	with data of small size $\mathring{\upepsilon}$ 
	on $\Sigma_0^1$ backwards in time and deduce that
	the ``data'' $(\halfdata^{(Radial)}, \halftimedata^{(Radial)})$
	appearing in definition \eqref{E:CONSTANTDEPENDSONDATANEIGHBORHOOD},
	which we view to be induced by the solution on $\Sigma_{-1/2},$ 
	have size $\lesssim \mathring{\upepsilon};$
	this is a basic bound that can be proved using the ideas mentioned in
	our sketch of a proof of Prop.~\ref{P:CLASSICALLOCAL}.
	Inequality \eqref{E:SHOCKPARAMETERBOUNDEDBYDATAATTIMEZERO} follows easily 
	from this bound.

	To finish the proof of the lemma, it remains only for us to show that the inequalities proved
	above for the small spherically symmetric data 
	on $\Sigma_{-1/2}$
	are stable with respect to small
	non-symmetric perturbations 
	belonging to $H_{\Euct}^{25}(\Sigma_{-1/2}) \times H_{\Euct}^{24}(\Sigma_{-1/2})$
	and supported in the Euclidean ball of radius $1/2$
	centered at the origin.
	The stability follows easily from Sobolev embedding
	without changing the parameter $U_0$ associated to the spherically symmetric data
	since inequality \eqref{E:KEYDATABLOWUPTERMHASBLOWUPSIGNANDSIZE} is strict,
	and since Prop.~\ref{P:CLASSICALLOCAL}
	implies that on the interval $t \in [-1/2,0],$
	the solution varies continuously with respect to the data
	given at $\Sigma_{-1/2}.$
\end{proof}

\section{Proof of shock formation for small, nearly spherically symmetric data}
\label{S:STABLESHOCKFORMATION}
In this section, we prove our main shock-formation theorem. 
Specifically, we prove that the open set $\shockset$ 
of data from Lemma~\ref{L:EXISTENCEOFSHOCKFORMINGDATA}
launches solutions that develop a shock singularity in finite time.

\begin{theorem}[\textbf{Finite-time shock formation for data in} $\shockset$]
\label{T:STABLESHOCKFORMATION}
Let 
$(\halfdata := \Psi|_{\Sigma_{-1/2}}, \halftimedata := \partial_t \Psi|_{\Sigma_{-1/2}})$
be data for the covariant wave equation \eqref{E:WAVEGEO}
under the assumption\footnote{As we described in Chapter~\ref{C:INTRO}, this assumption is easy to eliminate.}
\eqref{E:GINVERSE00ISMINUSONE}.
Assume that $(\halfdata,\halftimedata)$
belong to the open set $\shockset$ (of small, nearly spherically symmetric data)
from Lemma~\ref{L:EXISTENCEOFSHOCKFORMINGDATA}. 
In particular, the data
are supported in the Euclidean ball of radius $1/2$ centered at the origin.
Let $\Psi$ denote the corresponding solution.
Let $\mathring{\upepsilon} = \mathring{\upepsilon}[(\mathring{\Psi},\mathring{\Psi}_0)]
	:= \| \mathring{\Psi} \|_{H_{\Euct}^{25}(\Sigma_0^1)}
	+ \| \mathring{\Psi}_0 \|_{H_{\Euct}^{24}(\Sigma_0^1)}$
denote the size of the ``data'' 
$(\mathring{\Psi} := \Psi|_{\Sigma_0^1}, \mathring{\Psi}_0 := \partial_t \Psi|_{\Sigma_0^1})$
(which are supported in the Euclidean unit ball $\Sigma_0^1$)
induced by the solution on $\Sigma_0^1,$
and let $\mathring{\updelta}$ be the size of the geometric rotational derivatives of $\Psi$ 
on $\Sigma_0^1$ as defined in \eqref{E:SMALLANGULARDERIVATIVES}. 
Note that by the definition of $\shockset,$ 
the hypotheses and conclusions of 
Theorem~\ref{T:LONGTIMEPLUSESTIMATES}
and Cor.~\ref{C:ANGULARDDERIVATIVEESTIMATES} hold for data in $\shockset.$ 
In particular, we have $\mathring{\upepsilon} < \upepsilon_0,$ where
$\upepsilon_0$ is the small positive constant from Theorem~\ref{T:LONGTIMEPLUSESTIMATES}.
Let $U_0 \in (0,1)$ and $\shockparameter > 0$ be the 
\underline{data-neighborhood-dependent} parameters appearing in the statement of Lemma~\ref{L:EXISTENCEOFSHOCKFORMINGDATA}
and let $T_{(Lifespan);U_0}$ be the classical lifespan of $\Psi$ appearing in the statement of 
Theorem~\ref{T:LONGTIMEPLUSESTIMATES}
and Cor.~\ref{C:ANGULARDDERIVATIVEESTIMATES}.
Then $T_{(Lifespan);U_0} < \infty$ 
precisely because $\upmu$ vanishes for the first time at one more points in $\Sigma_{T_{(Lifespan);U_0}},$ 
thus, according to Theorem~\ref{T:LONGTIMEPLUSESTIMATES}, 
signifying the onset of shock formation.
Furthermore, there exists a constant $C > 0$ depending on $U_0$ such that
\begin{align} \label{E:SHOCKTIMEUPPERBOUND}
	T_{(Lifespan);U_0}
	& \leq \exp \left(\frac{1 + C \mathring{\upepsilon}}{\shockparameter} \right).
\end{align}
\end{theorem}

\begin{remark}[\textbf{Beyond nearly spherically symmetric data.}]
	\label{R:SPHERICALLYSYMMETRICDATADONOTGENERATEASPHERICALLYSYMMETRICSOLUTION}
	Although the data in $\shockset$ are nearly spherically symmetric,
	equation \eqref{E:WAVEGEO} is not generally invariant under Euclidean rotations.
	Thus, even when the data are exactly spherically symmetric, 
	the solution does not generally inherit this property.
	It might therefore seem a bit unnatural to prove shock formation for the data
	in $\shockset.$ It is only for technical convenience 
	that we have treated these data in detail;
	in Sect.~\ref{SS:INTROSHOCKFORMINGDATACOMPARISON},
	we outlined how to extend our results to prove shock formation 
	for a much larger set of compactly supported small data.
\end{remark}

\begin{proof}
Recall that (see Theorem~\ref{T:LONGTIMEPLUSESTIMATES})
$T_{(Lifespan);U_0}= \sup \lbrace t \ | \ \inf_{[0,t)} \upmu_{\star}(t,U_0) > 0 \rbrace,$
where $\upmu_{\star}$ is defined in \eqref{E:MUSTARDEF}.
Let $u_*$ and $\vartheta_*$ be as in Lemma~\ref{L:EXISTENCEOFSHOCKFORMINGDATA}.
From inequality \eqref{E:KEYUPMUVERYSHARPPOINTWISEESTIMATE}, the properties of the data
in $\shockset$ that are stated in Lemma~\ref{L:EXISTENCEOFSHOCKFORMINGDATA},
and the definition of $\upmu_{\star},$ it follows that
\begin{align} \label{E:UPMUSTARKEYPOINTWISEBOUNDPROOFOFSHOCKFORMATION}
			\upmu_{\star}(t,U_0)
			\leq
			\upmu(t,u_*,\vartheta_*)
	 		& \leq 
	 			1 
	 			+ C \mathring{\upepsilon} 
	 			- \shockparameter \ln \left( \frac{\rgeo(t,u_*)}{\rgeo(0,u_*)} \right).
	\end{align}
		It now easily follows from 
		\eqref{E:SHOCKPARAMETERBOUNDEDBYDATAATTIMEZERO},
		\eqref{E:UPMUSTARKEYPOINTWISEBOUNDPROOFOFSHOCKFORMATION},
		and the definition $\rgeo(t,u) = 1 - u + t$
		that $\upmu_{\star}(t,U_0)$
		must vanish before the time on the right-hand side of \eqref{E:SHOCKTIMEUPPERBOUND}.

\end{proof}



\appendix

\chapter{Extension of the Results to a Class of Non-Covariant Wave Equations} 
\label{A:EQUIVALENTPROBLEM}
\thispagestyle{fancy}
In Appendix \ref{A:EQUIVALENTPROBLEM}, we sketch how to extend the results of
Theorems~\ref{T:LONGTIMEPLUSESTIMATES} and~\ref{T:STABLESHOCKFORMATION}
(Remark~\ref{R:SPHERICALLYSYMMETRICDATADONOTGENERATEASPHERICALLYSYMMETRICSOLUTION} also applies in the present context)
to non-covariant equations of the form
\begin{align} \label{E:HOMOGENEOUSPHIDERIVATIVEQUASILINEARWAVE}
	(h^{-1})^{\alpha \beta}(\partial \Phi) \partial_{\alpha} \partial_{\beta} \Phi
	& = 0,
\end{align}
where
\begin{align}
	h_{\mu \nu}
	& = m_{\mu \nu} + h_{\mu \nu}^{(Small)}(\partial \Phi)
\end{align}	
and $h_{\mu \nu}^{(Small)}(0) = 0.$ As in 
\eqref{E:GINVERSE00ISMINUSONE}, in order to simplify the calculations, 
we make the non-essential assumption $(h^{-1})^{00} = -1.$ 

\begin{remark}[\textbf{A big difference between equation \eqref{E:WAVEGEO} and equation \eqref{E:HOMOGENEOUSPHIDERIVATIVEQUASILINEARWAVE}}]
	\label{R:BIGDIFFERENCE}
	Even though equation \eqref{E:HOMOGENEOUSPHIDERIVATIVEQUASILINEARWAVE}
	does not contain any semilinear terms, it nonetheless exhibits 
	small-data future shock formation\footnote{That is, an analog of Theorem~\ref{T:STABLESHOCKFORMATION},
	showing that nearly spherically symmetric small data lead to finite future time shock formation,
	holds for equation  \eqref{E:HOMOGENEOUSPHIDERIVATIVEQUASILINEARWAVE}.} 
	whenever $\FutFailFac \not \equiv 0,$
	where $\FutFailFac$ is defined in \eqref{E:NEWNULLCONDITIONFAILUREFACTOR}.
	We recall that in contrast, small-data global existence holds for equation \eqref{E:WAVEGEO}
	if we delete the semilinear terms 
	that are present when the equation is expanded relative to the rectangular coordinates;
	see Remark~\ref{R:QUASILINEARTERMSAREHARMLESSBYTHEMSELVES}.
\end{remark}

\section{From the scalar quasilinear wave equation to the equivalent system of wave equations}
The main idea is that the overwhelming majority of the analysis 
of solutions to \eqref{E:HOMOGENEOUSPHIDERIVATIVEQUASILINEARWAVE}
is essentially the same as the analysis 
we used on the covariant scalar equation $\square_{g(\Psi)} \Psi = 0.$
As we will see, the reason is that we can
differentiate equation \eqref{E:HOMOGENEOUSPHIDERIVATIVEQUASILINEARWAVE} 
to transform it into a coupled system of equations that are closely 
related to equation \eqref{E:WAVEGEO}.
Specifically, we will commute the wave equation 
\eqref{E:HOMOGENEOUSPHIDERIVATIVEQUASILINEARWAVE} with the rectangular
derivatives $\frac{\partial}{\partial x^{\nu}} = \partial_{\nu}.$ 
This motivates the following definition.
\begin{definition}[\textbf{The quantities $\Psi_{\nu}$}]
	We define $\Psi_{\nu},$ $(\nu = 0,1,2,3),$ and $\vec{\Psi}$ as follows:
	\begin{align}
		\Psi_{\nu} 
		& := \partial_{\nu} \Phi,
			\\
		\vec{\Psi} 
		&:= (\Psi_0,\Psi_1,\Psi_2,\Psi_3).
	\end{align}
\end{definition}
We view the $\Psi_{\nu}$ as the new unknowns.
The following functions of $\vec{\Psi},$ 
which are analogs of the functions $G_{\mu \nu}$ from \eqref{E:BIGGDEF},
play a fundamental role in our analysis.
\begin{definition}[\textbf{Metric component derivative functions}]
	We define
	\begin{align} \label{E:BIGHDEF}
		H_{\mu \nu}^{\lambda}
		= H_{\mu \nu}^{\lambda}(\vec{\Psi})
		& := \frac{\partial h_{\mu \nu}^{(Small)}(\vec{\Psi})}{\partial \Psi_{\lambda}}.
	\end{align}
\end{definition}
It is straightforward to calculate that the
analog of equation \eqref{E:RECTANGULARCHRISTOFFEL} for the Christoffel symbols 
relative to the rectangular coordinates is
\begin{align} \label{E:NEWCHRIST}
	\Gamma_{\alpha \kappa \beta}
	= \frac{1}{2} 
			\left\lbrace
				\partial_{\alpha} h_{\kappa \beta}
				+ \partial_{\beta} h_{\alpha \kappa}
				- \partial_{\kappa} h_{\alpha \beta}
			\right\rbrace
	& = \frac{1}{2}
			\left\lbrace
				H_{\kappa \beta}^{\lambda} \partial_{\alpha} \Psi_{\lambda}
				+ H_{\alpha \kappa}^{\lambda} \partial_{\beta} \Psi_{\lambda}
				- H_{\alpha \beta}^{\lambda} \partial_{\kappa} \Psi_{\lambda}
			\right\rbrace.
\end{align}

Recall that for equation \eqref{E:WAVEGEO},
the future null condition failure factor is provided in equation \eqref{E:FAILUREFACTOR}.
We now provide,
for equation \eqref{E:HOMOGENEOUSPHIDERIVATIVEQUASILINEARWAVE},
the appropriate analog in the region 
$\lbrace t \geq 0 \rbrace.$

\begin{definition}[\textbf{Future null condition failure factor for the equation} 
	$(h^{-1})^{\alpha \beta}(\partial \Phi) \partial_{\alpha} \partial_{\beta} \Phi = 0$]
	\label{D:FUTUREFAILUREFACTORFORNONGEOMETRIC}
We define $\FutFailFac,$ the future null condition failure factor corresponding to equation \eqref{E:HOMOGENEOUSPHIDERIVATIVEQUASILINEARWAVE},
by
\begin{align} \label{E:NEWNULLCONDITIONFAILUREFACTOR}
		\FutFailFac
		&:=
		m_{\kappa \lambda} H_{\mu \nu}^{\kappa}(\vec{\Psi} = 0) 
		\Lunit_{(Flat)}^{\lambda} \Lunit_{(Flat)}^{\mu} \Lunit_{(Flat)}^{\nu},
	\end{align}
where $m_{\mu \nu} = \mbox{diag}(-1,1,1,1)$ is the Minkowski metric,
$\Lunit_{(Flat)} = \partial_t + \partial_r,$
and $\partial_r = x^a/r \partial_a$ is the standard Euclidean radial derivative.
\end{definition}

We now explain the relevance of Def.~\ref{D:FUTUREFAILUREFACTORFORNONGEOMETRIC}.
Our discussion here closely parallels the discussion in Sect.~\ref{S:WHENNULLCONDITIONFAILS}.
We first note that the nonlinear terms in equation \eqref{E:HOMOGENEOUSPHIDERIVATIVEQUASILINEARWAVE}
verify Klainerman's classic null condition \cite{sk1984}
(see the discussion in Sect.~\ref{S:STRUCTUREOFNONLINEARITIES})
if and only if
$m_{\kappa \lambda} H_{\mu \nu}^{\kappa}(\vec{\Psi} = 0) 
		\ell^{\lambda} \ell^{\mu} \ell^{\nu} = 0$
for all Minkowski-null vectors $\ell.$ 
Furthermore, the previous equation holds if and only if $\FutFailFac \equiv 0;$
see Footnote \ref{FN:ALLPOSSIBLENULLVECTORS} on pg. \pageref{FN:ALLPOSSIBLENULLVECTORS}.
In particular, when $\FutFailFac \equiv 0,$
the methods of \cite{sK1986} and \cite{dC1986a} yield small-data global existence.
Moreover, upon Taylor expanding equation \eqref{E:HOMOGENEOUSPHIDERIVATIVEQUASILINEARWAVE}
to quadratic order around the $\textbf{0}$ solution and then decomposing the
quadratic terms relative to the flat frame
\eqref{E:MINKOWSKIFRAME}, we find that $\FutFailFac$ is the coefficient of
the term $(\Radunit_{(Flat)} \Phi) \Radunit_{(Flat)} (\Radunit_{(Flat)} \Phi).$
As we mentioned in Remark~\ref{R:BIGDIFFERENCE},
when $\FutFailFac$ is nontrivial, 
equation \eqref{E:HOMOGENEOUSPHIDERIVATIVEQUASILINEARWAVE}
exhibits small-data future shock formation
caused precisely by the term 
$(\Radunit_{(Flat)} \Phi) \Radunit_{(Flat)} (\Radunit_{(Flat)} \Phi);$
the methods of \cite{sK1986} and \cite{dC1986a} imply
that in the region $\lbrace t \geq 0 \rbrace,$
the other quadratic terms are compatible
with small-data future-global existence.

\begin{remark}[\textbf{Past null condition failure factor for equation} \eqref{E:HOMOGENEOUSPHIDERIVATIVEQUASILINEARWAVE}]
	We could also study shock formation in the region $\lbrace t \leq 0 \rbrace.$
	As we described in Remark~\ref{R:PASTFAILUREFACTOR}, 
	in this region, the relevant past null condition failure factor 
	$\PastFailFac$ is defined by replacing
	$\Lunit_{(Flat)}$ with $- \partial_t + \partial_r$ in equation
	\eqref{E:NEWNULLCONDITIONFAILUREFACTOR}.
\end{remark}

In the next lemma, we provide the wave equations verified by the
$\Psi_{\nu} = \partial_{\nu} \Phi,$ $(\nu = 0,1,2,3).$
For reasons to be explained, we also further commute the equations
with the \emph{Minkowski scaling vectorfield} $S := x^{\alpha} \partial_{\alpha},$ 
which is a conformal Killing field\footnote{That is, $\Lie_S m$ is equal to a scalar function multiple of $m.$} of the Minkowski metric.

\begin{lemma}[\textbf{The structure of the commuted equations}] \label{L:COMMUTED}
Assume that $\Phi$ verifies equation \eqref{E:HOMOGENEOUSPHIDERIVATIVEQUASILINEARWAVE}. Let
$\Psi \in \lbrace \partial_{\nu} \Phi \rbrace_{\nu=0,1,2,3}.$ Then the scalar-valued 
function $\Psi$ verifies the following equation
relative to the \textbf{rectangular Minkowski} coordinates:
\begin{align} \label{E:COMMUTED}
	\partial_{\alpha} 
	\left(
		(h^{-1})^{\alpha \beta} \partial_{\beta} \Psi 
	\right) 
	& = H^{\mu \alpha \beta}(\vec{\Psi})
			\left\lbrace
				\partial_{\beta} \Psi_{\alpha} \partial_{\mu} \Psi
				- \partial_{\mu} \Psi_{\alpha} \partial_{\beta} \Psi
			\right\rbrace,
\end{align}	

\begin{align} \label{E:ALMOSTNULLFORM}
	H^{\mu \alpha \beta}(\vec{\Psi})	
	& := (h^{-1})^{\alpha \alpha'} (h^{-1})^{\beta \beta'} H_{\alpha' \beta'}^{\mu}.
\end{align}

Similarly, let $S$ be the Minkowski scaling vectorfield
\begin{align} \label{E:MINKOWSKISCALING}
	S 
	&:= x^{\alpha} \partial_{\alpha},
\end{align}
and let 
\begin{align} \label{E:SCALINGRENORMALIZEDPSI}
	\Psi 
	&:= S \Phi - \Phi
		= x^{\alpha} \Psi_{\alpha} - \Phi.
\end{align}
Then $\Psi$ verifies the following equation
relative to the \textbf{rectangular Minkowski} coordinates:
\begin{align} \label{E:SCALINGCOMMUTED}
	\partial_{\alpha} 
	\left(
		(h^{-1})^{\alpha \beta} \partial_{\beta} \Psi 
	\right) 
	& = H^{\mu \alpha \beta}(\vec{\Psi})
			\left\lbrace
				\partial_{\beta} \Psi_{\alpha} \partial_{\mu} \Psi
				- \partial_{\mu} \Psi_{\alpha} \partial_{\beta} \Psi
			\right\rbrace.
\end{align}	



\end{lemma}

\begin{proof}
Equation \eqref{E:COMMUTED} follows from commuting \eqref{E:HOMOGENEOUSPHIDERIVATIVEQUASILINEARWAVE} with $\partial_{\nu}$ 
and performing straightforward calculations. Similarly, 
\eqref{E:SCALINGCOMMUTED}
follows from 
commuting \eqref{E:HOMOGENEOUSPHIDERIVATIVEQUASILINEARWAVE} 
with 
$S$ 
and performing straightforward calculations.
\end{proof}

Just below, we verify that \eqref{E:COMMUTED} and \eqref{E:SCALINGCOMMUTED}
can be rewritten as covariant wave equations whose 
\emph{right-hand} sides have a good null structure.
The following functions play a role in our identification of this geometric structure.

\begin{definition}
	\label{D:GOODNULLSTRUCTUREWAVEEQUATIONTERM}
	We define the rectangular components $\Omega^{\nu},$ 
	$(\nu = 0,1,2,3),$
	to be the following 
	scalar-valued functions of $\vec{\Psi}:$
	\begin{align} \label{E:RECTANGULARVOLUMEFORMFACTOR}
		\Omega^{\nu} 
		& = \Omega^{\nu}(\vec{\Psi})
		:= \frac{1}{\sqrt{|\mbox{\upshape{det}} h|(\vec{\Psi})}} 
			\frac{\partial \sqrt{|\mbox{\upshape{det}} h|(\vec{\Psi})}}{\partial \Psi_{\nu}}
		= \frac{1}{2} 
			\frac{\partial \ln \left(|\mbox{\upshape{det}} h|(\vec{\Psi})\right)}{\partial \Psi_{\nu}},
	\end{align}
	where the determinant in \eqref{E:RECTANGULARVOLUMEFORMFACTOR} is taken
	relative to rectangular coordinates.

\end{definition}

In the next lemma, we rewrite the equations from
Lemma~\ref{L:COMMUTED} as covariant wave equations
with inhomogeneous terms.

\begin{lemma} [\textbf{Rewriting the equations as covariant wave equations}] 
\label{L:ECOMMUTEDINTERMSOFHCON}
Relative to the rectangular coordinates, the system of wave equations \eqref{E:COMMUTED} 
can be equivalently expressed as the following system of covariant wave equations in the unknowns 
$\vec{\Psi} := (\Psi_0, \Psi_1, \Psi_2, \Psi_3),$ where $\Psi_{\nu} := \partial_{\nu} \Phi$ 
is viewed to be a scalar-valued function:
\begin{align} \label{E:PSIGEOMETRIC}
	\square_{h(\vec{\Psi})} \Psi_{\nu} 
	& = \mathscr{Q}(\partial \vec{\Psi}, \partial \Psi_{\nu}), & (\nu & = 0,1,2,3),
\end{align}
\begin{align} \label{E:GEOMETRICWAVEOPERATORACTION}
	\square_{h(\vec{\Psi})} \Psi := \frac{1}{\sqrt{|\mbox{\upshape{det}} h|}} 
		\partial_{\alpha} \left(\sqrt{|\mbox{\upshape{det}} h|} (h^{-1})^{\alpha \beta} \partial_{\beta} \Psi \right)
\end{align}
is the covariant wave operator of $h$ applied to $\Psi,$ 
and 
\begin{align}
	\mathscr{Q}(\partial \vec{\Psi}, \partial \Psi) 
	& :=	H^{\mu \alpha \beta}
				\left\lbrace
					\partial_{\beta} \Psi_{\alpha} \partial_{\mu} \Psi
					- \partial_{\mu} \Psi_{\alpha} \partial_{\beta} \Psi
				\right\rbrace
			+ (h^{-1})^{\alpha \beta} \Omega^{\lambda} \partial_{\alpha} \Psi_{\lambda} 
			\partial_{\beta} \Psi.
\end{align}
The quantities $H^{\alpha \beta \mu}$ and $\Omega^{\lambda}$
are defined in \eqref{E:ALMOSTNULLFORM} and \eqref{E:RECTANGULARVOLUMEFORMFACTOR}.

Similarly, the wave equation \eqref{E:SCALINGCOMMUTED} for the 
scalar-valued function $\Psi := S \Phi - \Phi$ can be expressed as
\begin{align} \label{E:SCALINGPSIGEOMETRIC}
	\square_{h(\vec{\Psi})} \Psi
	& = \mathscr{Q}(\partial \vec{\Psi}, \partial \Psi).
\end{align}


\end{lemma}


\begin{proof}
	Lemma~\ref{L:ECOMMUTEDINTERMSOFHCON} follows from Lemma~\ref{L:COMMUTED},
	Def.~\ref{D:GOODNULLSTRUCTUREWAVEEQUATIONTERM},
	and straightforward computations.
\end{proof}

The next lemma is at the heart of the connection between
equations \eqref{E:HOMOGENEOUSPHIDERIVATIVEQUASILINEARWAVE} and
\eqref{E:WAVEGEO}. It shows shows that the quadratic
terms $\mathscr{Q}(\cdot,\cdot)$ from Lemma~\ref{L:ECOMMUTEDINTERMSOFHCON}
have a special null structure and hence have a negligible effect
on the dynamics of small-data solutions.
Before we state the lemma, we first state the 
definitions of the standard null forms relative to $h.$

\begin{definition}[\textbf{Standard null forms relative to $h$}]
	\label{D:STANDARDNULLFORMS}
	We define the following bilinear forms
	$\mathscr{Q}_{(0)}$ and $\mathscr{Q}_{(\alpha \beta)},$
	which act on pairs of functions $\phi$ and $\widetilde{\phi}$
	as follows:
	\begin{subequations}
	\begin{align}
		\mathscr{Q}_{(0)}(\partial \phi, \partial \widetilde{\phi})
		& := (h^{-1})^{\alpha \beta} \partial_{\alpha} \phi \partial_{\beta} \widetilde{\phi},
			\label{E:STANDARD0NULLFORM} \\
		\mathscr{Q}_{(\alpha \beta)}(\partial \phi, \partial \widetilde{\phi})
		& := \partial_{\alpha} \phi \partial_{\beta} \widetilde{\phi}
			- \partial_{\alpha} \widetilde{\phi} \partial_{\beta} \phi.
			\label{E:STANDARDALPHABETANULLFORM}
	\end{align}
	\end{subequations}
\end{definition}

We now recall the important property of the standard null forms.
By expanding $h^{-1}(\vec{\Psi})$ relative to the frame 
$\lbrace \Lunit, \Radunit, X_1, X_2 \rbrace$
as
$(h^{-1})^{\mu \nu} = 
- \Lunit^{\mu} \Lunit^{\nu}
- (\Lunit^{\mu} \Radunit^{\nu} + \Radunit^{\mu} \Lunit^{\nu})
	+ (\hinversesphere)^{\mu \nu}$
(as in \eqref{E:GINVERSEFRAMEWITHRECTCOORDINATESFORGSPHEREINVERSE})
and noting that
$(\hinversesphere)^{\mu \nu}$
is a smooth function of $\vec{\Psi}$ and
the rectangular components $\Lunit^{\mu},$ $\Radunit^{\mu},$
we see that
\begin{align} \label{E:Q0NULLFORMEXPANDED}
	\mathscr{Q}_{(0)}(\partial \phi, \partial \widetilde{\phi})
	& = f_1 \Lunit \phi \Lunit \widetilde{\phi}
		+ f_2 \Lunit \phi \Radunit \widetilde{\phi}
		+ f_3 \Radunit \phi \Lunit \widetilde{\phi}
		+ f_4^{AB} \angdiffarg{A} \phi \angdiffarg{B} \widetilde{\phi},
\end{align}
where $f_1,$ $f_2,$ and $f_3$ are smooth scalar-valued functions of $\vec{\Psi}$
and the rectangular components $\Lunit^{\mu},$ $\Radunit^{\mu},$
and the same statement holds for the rectangular components $f_4^{\alpha \beta}$
of the symmetric type $\binom{2}{0}$ $S_{t,u}$ tensor $f_4.$
The important point is that there are no quadratic terms on the right-hand side of
\eqref{E:Q0NULLFORMEXPANDED} involving two $\mathcal{C}_u^t-$transversal derivatives,
that is, no terms proportional to $\Radunit \phi \Radunit \widetilde{\phi}.$
Similarly, upon expanding 
$\mathscr{Q}_{(\alpha \beta)}(\partial \phi, \partial \widetilde{\phi})$
relative to the non-rescaled frame $\lbrace \Lunit, \Radunit, X_1, X_2 \rbrace,$
we find that the coefficients depend smoothly on $\vec{\Psi}$
and the rectangular components $\Lunit^{\mu},$ $\Radunit^{\mu},$
and that there are no terms proportional to 
$\Radunit \phi \Radunit \widetilde{\phi}.$ 

We capture the absence of the component $\Radunit \phi \Radunit \widetilde{\phi}$
with the following definition, which is relevant in the region 
$\lbrace t \geq 0 \rbrace.$

\begin{definition}[\textbf{Future strong null condition}]
	\label{D:STRONGNULL}
	We say that $\mathscr{Q}(\partial \phi, \partial \widetilde{\phi})$ verifies the 
	future strong null condition if,
	relative to the non-rescaled frame $\lbrace \Lunit, \Radunit, X_1, X_2 \rbrace,$
	we have
	\begin{align} \label{E:STRONGNULL}
		\mathscr{Q}(\partial \phi, \partial \widetilde{\phi})
		& = f_1 \Lunit \phi \Lunit \widetilde{\phi}
			+ f_2 \Lunit \phi \Radunit \widetilde{\phi}
			+ f_3 \Radunit \phi \Lunit \widetilde{\phi}
				 \\
		& \ \ + f_4^A \angdiffarg{A} \phi \Lunit \widetilde{\phi}
			+ f_5^A \angdiffarg{A} \phi \Radunit \widetilde{\phi}
			+ f_6^A \Lunit \phi \angdiffarg{A} \widetilde{\phi} 
			+ f_7^A \Radunit \phi \angdiffarg{A} \widetilde{\phi} 
			+ f_8^{AB} \angdiffarg{A} \phi \angdiffarg{B} \widetilde{\phi},
			\notag
	\end{align}
	where the scalars $f_1,$ $f_2$ and $f_3$
	and the rectangular components 
	$f_4^{\alpha},$ $f_5^{\alpha},$ $f_6^{\alpha},$ $f_7^{\alpha},$
	and $f_8^{\alpha \beta}$
	are smooth scalar-valued functions of $\vec{\Psi}$
	and the rectangular components of the vectorfields $\Lunit$ and $\Radunit.$
\end{definition}	

\begin{remark}[\textbf{The future strong null condition is nonlinear in nature}]
	Note that the future strong null condition is a nonlinear condition tied to the 
	dynamic metric $h.$ In particular, it is not necessarily 
	verified by cubic terms, which could contain dangerous factors 
	involving two or more derivatives in the $\Radunit$ direction.
\end{remark}

\begin{remark}[\textbf{Past strong null condition}]
	We could also formulate a past strong null condition,
	which would be the correct condition in the region
	$\lbrace t \leq 0 \rbrace.$
	To define such a condition, one would first need 
	to replace the frame
	$\lbrace \Lunit, \Radunit, X_1, X_2 \rbrace$
	with an analogous frame in which the null vectorfield
	is outgoing as we head towards the past.
	We would then define the past strong null condition 
	by replacing, in definition \eqref{E:STRONGNULL},
	$\Lunit$ with the new outgoing null vectorfield.
\end{remark}

We now provide the aforementioned lemma.
 
\begin{lemma}[\textbf{Special null structure of the inhomogeneous terms}]
	\label{L:SPECIALNULLSTRUCTUREINHOMOGENEOUS}
	The quadratic term $\mathscr{Q}(\partial \vec{\Psi}, \partial \Psi)$ on the right-hand sides of 
	\eqref{E:PSIGEOMETRIC}
	and \eqref{E:SCALINGPSIGEOMETRIC} 
	verifies the future strong null condition.
\end{lemma}

\begin{proof}
	Lemma~\ref{L:SPECIALNULLSTRUCTUREINHOMOGENEOUS} follows from
	comparing the right-hand sides of 
	\eqref{E:PSIGEOMETRIC}
	and \eqref{E:SCALINGPSIGEOMETRIC} 
	to the null forms
	$\mathscr{Q}_{(0)}(\cdot,\cdot)$
	and $\mathscr{Q}_{(\alpha \beta)}(\cdot,\cdot)$
	from
	Def.~\ref{D:STANDARDNULLFORMS}.
\end{proof}

The important consequence of Lemma~\ref{L:SPECIALNULLSTRUCTUREINHOMOGENEOUS}
is that the quadratic terms 
$\mathscr{Q}(\partial \vec{\Psi}, \partial \Psi)$ 
on the right-hand sides of \eqref{E:PSIGEOMETRIC}
and \eqref{E:SCALINGPSIGEOMETRIC} 
have a negligible influence on the
future dynamics of small-data solutions, even those that form shocks.
Specifically, because of their special null structure,
the $\mathscr{Z}^N$ derivatives of 
$\upmu$ times\footnote{Recall that in our analysis, we always work with the 
$\upmu-$weighted wave operator $\upmu \square_{h(\vec{\Psi})} \Psi_{\nu}.$} 
$\mathscr{Q}(\cdot,\cdot)$
are $Harmless^{\leq N}$ in the sense of Def.~\ref{D:HARMLESSTERMS}
and are therefore easy to treat.

\section{The main new estimate needed at the top order}
\label{S:MAINNEWESTIMATE}
Lemma~\ref{L:SPECIALNULLSTRUCTUREINHOMOGENEOUS} implies that
$\vec{\Psi}$ verifies the coupled system \eqref{E:PSIGEOMETRIC} 
of covariant wave equations with quadratic inhomogeneous terms that
have a special null structure.
We now explain how to prove an analog of Theorem~\ref{T:LONGTIMEPLUSESTIMATES}
for the system. The main point is that
we can derive $L^2$ estimates for each component
$\Psi_{\nu},$ which verifies a scalar wave equation,
by using essentially the same arguments that we used to  
treat the scalar equation \eqref{E:WAVEGEO}. 
To control $\vec{\Psi},$ 
we replace the controlling quantity $\totzeromax{N}$
from Def.~\ref{D:MAINCOERCIVEQUANT}
with
\begin{align}
	\totzeromax{N}(t,u)
	& := \max_{|\vec{I}| = N} \max_{\nu=0,1,2,3} \sup_{(t',u') \in [0,t] \times [0,u]} 
		\left\lbrace
			\enzero[\mathscr{Z}^{\vec{I}} \Psi_{\nu}](t',u')	
			+ \flzero[\mathscr{Z}^{\vec{I}} \Psi_{\nu}](t',u')
		\right\rbrace,
\end{align}
and similarly for the other controlling quantities.
As we explained in the previous paragraph, 
the quadratic terms on the right-hand side of \eqref{E:PSIGEOMETRIC} 
are easy to treat without invoking any new ideas.
However, we do need one new ingredient to close the estimates,
one that is needed for the
top-order $L^2$ estimates for the $\Psi_{\nu}.$
To clarify matters, we recall the analogous situation
that arose in our analysis of the scalar equation \eqref{E:WAVEGEO}.
At various points in our argument for deriving top-order $L^2$ estimates,
it was essential that we could use equation \eqref{E:UPMUFIRSTTRANSPORT} 
to make the replacement 
\begin{align} \label{E:KEYREPLACE}
	G_{\Lunit \Lunit} \Rad \Psi = 2 \Lunit \upmu + \mbox{Error},
\end{align}
where the error terms are small and decay at an integrable-in-time rate.
One example of the importance of this identity
was explained in detail just before
equation \eqref{E:GRONWALLREADYINTROTOPORDERENERGYIDCARICATURE}.
The main point is that we can treat the coupled system $\vec{\Psi}$
by separately analyzing each $\Psi_{\nu},$
but we need an analog of \eqref{E:KEYREPLACE} 
\emph{for each of the four} $\Psi_{\nu}.$ 
An analog of the factor $G_{\Lunit \Lunit} \Rad \Psi$ from \eqref{E:KEYREPLACE}
arises when we derive top-order $L^2$ estimates for 
the scalar component $\Psi_{\nu}$ in equation \eqref{E:PSIGEOMETRIC}.
Specifically, the analog is $-H_{\Lunit \Lunit}^{\Lunit} \Rad \Psi_{\nu},$
where $H_{\Lunit \Lunit}^{\Lunit}:= H_{\alpha \beta}^{\kappa} \Lunit^{\alpha} \Lunit^{\beta} \Lunit_{\kappa}.$
On the other hand, as we will see
in equation \eqref{E:LITTLEOMEGATRANSVERSALPART} below, 
$\upmu$ verifies a transport equation of the form
$\Lunit \upmu = - \frac{1}{2} H_{\Lunit \Lunit}^{\Lunit} \Radunit^a \Rad \Psi_a
+ \mbox{Error}.$ 
Hence, it is not entirely obvious that for solutions to the system \eqref{E:PSIGEOMETRIC},
the following analog of \eqref{E:KEYREPLACE},
which is sufficient for closing the top-order $L^2$ estimates,
should hold for $\nu = 0,1,2,3:$
\begin{align} \label{E:FIRSTSTATEMENTHARDSHARPRADPSIPOINTWISEESTIMATE}
	\left|H_{\Lunit \Lunit}^{\Lunit} \Rad \Psi_{\nu} \right| 
	& \leq 2 \left|\Lunit \upmu \right| + \mbox{Error}.
\end{align}
We stress that inequality \eqref{E:FIRSTSTATEMENTHARDSHARPRADPSIPOINTWISEESTIMATE} is exactly
the new ingredient mentioned above that we need to close the top-order $L^2$ estimates.
We dedicate the remainder of this section to deriving the
estimate \eqref{E:FIRSTSTATEMENTHARDSHARPRADPSIPOINTWISEESTIMATE}.
We prove the final estimate of interest 
in Prop.~\ref{P:EACHRADPSINUFACTORISESSENTIALLYTWICELUPMU}.

In the next lemma, we provide the analog of equation \eqref{E:UPMUFIRSTTRANSPORT},
that is, the transport equation equation verified by $\upmu.$
Our main goal is to identify the dangerous slow-decaying component
$\upomega^{(Trans-\vec{\Psi})},$ which arises from the failure of the classic null condition 
and which can drive $\upmu$ to $0.$

\begin{lemma} [\textbf{The transport equation verified by $\upmu$}] \label{L:NEWUPMUTRANSPORT}
In the case of equation \eqref{E:HOMOGENEOUSPHIDERIVATIVEQUASILINEARWAVE},
the inverse foliation density $\upmu$ defined in \eqref{E:UPMUDEF} verifies the following transport equation:
\begin{align} \label{E:NEWUPMUTRANSPORT}
	\Lunit \upmu 
	& = \upomega^{(Trans-\vec{\Psi})} 
		+ \upmu \upomega^{(Tan-\vec{\Psi})}
		:= \upomega,
\end{align}	
where
\begin{subequations}
\begin{align}
	\upomega^{(Trans-\vec{\Psi})} 
	& := - \frac{1}{2} H_{\Lunit \Lunit}^{\Lunit} \Radunit^a \Rad \Psi_a,
		\label{E:LITTLEOMEGATRANSVERSALPART} \\
	\upomega^{(Tan-\vec{\Psi})} 
	& := 
		- \frac{1}{2} H_{\Lunit \Lunit}^{\Lunit} \Radunit^a \Lunit \Psi_a
		- \frac{1}{2} H_{\Lunit \Lunit}^{\Radunit} \Radunit^a \Lunit \Psi_a
		+ \frac{1}{2} \angH_{\Lunit \Lunit}^A \Radunit^a \angdiffarg{A} \Psi_a
		- \frac{1}{2} H_{\Lunit \Lunit}^{\lambda} \Lunit \Psi_{\lambda}
		- H_{\Lunit \Radunit}^{\lambda} \Lunit \Psi_{\lambda}.
		\label{E:LITTLEOMEGATANGENTIALPART}
\end{align}
\end{subequations}
The scalar-valued function $H_{\mu \nu}^{\lambda}$ above is defined in \eqref{E:BIGHDEF},
$H_{\Lunit \Lunit}^{\Lunit}:= H_{\alpha \beta}^{\kappa} \Lunit^{\alpha} \Lunit^{\beta} \Lunit_{\kappa},$
$\angH_{\Lunit \Lunit}^{\lambda} 
:= H_{\alpha \beta}^{\kappa} \Lunit^{\alpha} \Lunit^{\beta} \sphereproject_{\ \kappa}^{\lambda},$
$\angH_{\Lunit \Lunit}^A 
:= H_{\alpha \beta}^{\kappa} \Lunit^{\alpha} \Lunit^{\beta} \sphereproject_{\ \kappa}^{\lambda}
\gsphere_{\lambda \mu} X_A^{\mu},$
etc.
\end{lemma}

\begin{proof}
	The proof is very similar to the proof of Lemma~\ref{L:UPMUFIRSTTRANSPORT}.
	The main difference is that we use the identity
	$\partial_{\alpha} \Psi_{\beta} = \partial_{\beta} \Psi_{\alpha}$
	to rewrite 
	\begin{align}
		H_{\Lunit \Lunit}^{\lambda} \Rad \Psi_{\lambda}
		& = \upmu H_{\Lunit \Lunit}^{\lambda} \Radunit^a \partial_{\lambda} \Psi_a
			\\
		& = - H_{\Lunit \Lunit}^{\Lunit} \Radunit^a \Rad \Psi_a
			- \upmu H_{\Lunit \Lunit}^{\Lunit} \Radunit^a \Lunit \Psi_a
			- \upmu H_{\Lunit \Lunit}^{\Radunit} \Radunit^a \Lunit \Psi_a
			+ \upmu \angH_{\Lunit \Lunit}^A \Radunit^a \angdiffarg{A} \Psi_a.
			\notag
	\end{align}
\end{proof}

The next lemma provides the key technical estimates that
allow us to conclude \eqref{E:FIRSTSTATEMENTHARDSHARPRADPSIPOINTWISEESTIMATE}.

\begin{lemma}[\textbf{Approximate identities verified by $\Rad \Psi_{\nu}$}]
	Under the small-data and bootstrap assumptions 
	of Sects.~\ref{S:PSISOLVES}-\ref{S:C0BOUNDBOOTSTRAP}
	(with $\vec{\Psi}$ in place of $\Psi$),
	if $\varepsilon$ is sufficiently small, 
	then the following estimates hold for $(t,u) \in [0,\Tboot) \times [0,U_0]$
	($i = 1,2,3$):
	\begin{subequations}
	\begin{align}  \label{E:RADPSI0ISNEARRADIALCONTRACTIONRADPSI}
		\left \| 
			\Rad \Psi_0
			- \Radunit^a \Rad \Psi_a
		\right \|_{C^0(\Sigma_t^u)}
		& \leq C \varepsilon \frac{\ln(\myexp + t)}{(1 + t)^2},
			\\
		\left \| 
			\Rad \Psi_i
			+ \frac{t x^i}{r^2} \Radunit^a \Rad \Psi_a
		\right \|_{C^0(\Sigma_t^u)}
		& \leq C \varepsilon \frac{\ln(\myexp + t)}{(1 + t)^2}.
		\label{E:RADPSIIISNEARWEIGHTEDRADIALCONTRACTIONRADSPI}
	\end{align}
	\end{subequations}
\end{lemma}

\begin{proof}
	We first sketch the proof of \eqref{E:RADPSI0ISNEARRADIALCONTRACTIONRADPSI}.
	From
	\eqref{E:RADUNITJUNK},
	\eqref{E:EASYCOORDINATEBOUND},	
	and
	\eqref{E:FRAMECOMPONENTSBOOT}
	(which together imply that $|\Radunit^a| \lesssim 1$)
	and \eqref{E:PSIFUNDAMENTALC0BOUNDBOOTSTRAP},
	we see that
	\eqref{E:RADPSI0ISNEARRADIALCONTRACTIONRADPSI} easily follows for $t \leq 1.$ Hence, we 
	have to address only the case $t > 1.$
	To this end, 
	we use
	direct computation and 
	\eqref{E:RADUNITJUNK}
	to compute that
	\begin{align} \label{E:RADSCALINGPHIMINUSPHIID}
	\Rad (S \Phi - \Phi)
		& = t \Rad \Psi_0
			+ x^a \Rad \Psi_a
				\\
		& = t \Rad \Psi_0
			- \rgeo \Radunit^a \Rad \Psi_a
			+ \rgeo \Radunit_{(Small)}^a \Rad \Psi_a,
			\notag
	\end{align}
	where $S = x^{\alpha} \partial_{\alpha}$ is the Minkowskian scaling vectorfield.
	By 
	\eqref{E:RADUNITJUNKLIKELMINUSUNITJUNK} (with $h$ in place of $g$),
	\eqref{E:PSIFUNDAMENTALC0BOUNDBOOTSTRAP},
	\eqref{E:RGEOMINKOWSKITIMECOMPARISON},
	and
	\eqref{E:C0BOUNDCRUCIALEIKONALFUNCTIONQUANTITIES},
	we see that the last
	term on the right-hand side of \eqref{E:RADSCALINGPHIMINUSPHIID}
	is bounded in magnitude by $\lesssim \varepsilon \ln(\myexp + t) (1 + t)^{-1}.$
	Also using the fact that $\rgeo = 1 - u + t,$
	we see that \eqref{E:RADPSI0ISNEARRADIALCONTRACTIONRADPSI} for $t \geq 1$
	will follow once we show that
	\begin{align} \label{E:SCALINGPHIMINUSPHIJUNKESTIMATE}
		\|
		\Rad 
		(S \Phi - \Phi)
		\|_{C^0(\Sigma_t^u)}
		& \leq C \varepsilon \frac{1}{1 + t}.
	\end{align}
	To this end, we first multiply both sides of \eqref{E:SCALINGPSIGEOMETRIC} by $\upmu.$
	By Lemma~\ref{L:SPECIALNULLSTRUCTUREINHOMOGENEOUS},
	we have 
	\begin{align} \label{E:SCALINGPHIMINUSPHI}
		\upmu \square_{h(\vec{\Psi})} 
		(S \Phi - \Phi)
		& = \waveinhom,
	\end{align}
	where $\waveinhom$
	can be expressed as $f(\vec{\psi},\upmu,\Lunit^1,\Lunit^2,\Lunit^3)$ times
	terms that are quadratic in the derivatives of $S \Phi - \Phi$ 
	relative to the rescaled frame $\Lunit, \Rad, X_1, X_2,$
	with at least one
	factor involving a derivative that is tangential to $\mathcal{C}_u^t,$ where $f$ is a smooth function.
	Hence, we can apply essentially the same reasoning that we used to derive
	the estimates \eqref{E:Q0LOWESTORDERIMPROVED}-\eqref{E:Q1LOWESTORDERIMPROVED}
	in order to deduce that
	\begin{align} 
		\enzero^{1/2}[\mathscr{Z}^{\leq 2} \Rad (S \Phi - \Phi)]
		+ 
		\flzero^{1/2}[\mathscr{Z}^{\leq 2} \Rad (S \Phi - \Phi)]
		& \leq C \mathring{\upepsilon},
			\label{E:LOWERORDERENERGIESOFSCALINGPHIMINUSPHIARENONDEGENERATE} \\
		\enone^{1/2}[\mathscr{Z}^{\leq 2} \Rad (S \Phi - \Phi)]
		+ 
		\flone^{1/2}[\mathscr{Z}^{\leq 2} \Rad (S \Phi - \Phi)]
		+ \Morint^{1/2}[\mathscr{Z}^{\leq 2} \Rad (S \Phi - \Phi)]
		& \leq C \mathring{\upepsilon} \ln^2(\myexp + t).
	\end{align}
	Then from 
	\eqref{E:C0BOUNDSOBOLEVINTERMSOFENERGIES}
	and
	\eqref{E:LOWERORDERENERGIESOFSCALINGPHIMINUSPHIARENONDEGENERATE},
	we deduce the desired bound \eqref{E:SCALINGPHIMINUSPHIJUNKESTIMATE}.

	
	
	
	To prove \eqref{E:RADPSIIISNEARWEIGHTEDRADIALCONTRACTIONRADSPI}, we
	first use the identity $\partial_j \Psi_i = \partial_i \Psi_j$
	to compute that
	\begin{align} \label{E:RADPSIIRADIALPLUSJUNK}
		\Rad \Psi_i 
		& = \frac{x^i}{r} x^a \Rad \Psi_a 
			+ \upmu \Radunit^a V \Psi_a,
	\end{align}
	where $V$ is the
	$\Sigma_t^{U_0}-$tangent vectorfield
	(depending on $i$)
	with rectangular components
	\begin{align}
		V^j 
		& = \delta_i^j
				- \frac{x^i}{r} \frac{x^j}{r}.
	\end{align}
	From the bound $|\Radunit^a| \lesssim 1$ noted above,
	\eqref{E:RADUNITJUNK},
	\eqref{E:RADUNITJUNKLIKELMINUSUNITJUNK},
	and \eqref{E:RADPSIIRADIALPLUSJUNK},
	the estimates
	\eqref{E:EASYEUCLIDEANRADIALVARIABLEBOUND}
	and \eqref{E:C0BOUNDCRUCIALEIKONALFUNCTIONQUANTITIES},
	and the bootstrap assumptions
	\eqref{E:PSIFUNDAMENTALC0BOUNDBOOTSTRAP},
	we see that the desired bound \eqref{E:RADPSIIISNEARWEIGHTEDRADIALCONTRACTIONRADSPI}
	will therefore follow once we show that for $a=1,2,3,$ we have
	\begin{align} \label{E:UPMUVPSIKEYPOINTWISEBOUND}
		\upmu |V \Psi_a|
		& \lesssim \varepsilon \frac{\ln(\myexp + t)}{(1 + t)^2}.
	\end{align}
	To prove \eqref{E:UPMUVPSIKEYPOINTWISEBOUND}, we first note that
	since the Euclidean length of $V$ is $\lesssim 1$ and since $V$ is tangent
	to the spheres in $\Sigma_t^{U_0}$ of constant Euclidean radius,
	it is a standard fact
	(that can be proved by arguing as in our proof of \eqref{E:FUNCTIONPOINTWISEANGDINTERMSOFANGLIEO})
	that
	\begin{align} \label{E:VPSIBOUNDEDBYEUCRINVERSETIMESEUCROTATIONS}
		|V \Psi_a|
		\lesssim  
		\frac{1}{r} \sum_{l=1}^3 |\Roteucarg{l} \Psi_a|,
	\end{align}
	where $\Roteucarg{l}$ are the Euclidean rotations from Def.~\ref{D:ROTATION}.
	The desired bound \eqref{E:UPMUVPSIKEYPOINTWISEBOUND} now follows from
	\eqref{E:ROTATIONDECOMPOSITIONINTOEUCLIDEANPLUSRADCOMPONENT},
	\eqref{E:PSIFUNDAMENTALC0BOUNDBOOTSTRAP},
	\eqref{E:EASYEUCLIDEANRADIALVARIABLEBOUND},
	\eqref{E:LOWERORDERC0BOUNDEUCLIDEANROTATIONRADCOMPONENT},
	\eqref{E:C0BOUNDCRUCIALEIKONALFUNCTIONQUANTITIES},
	Cor.~\ref{C:SQRTEPSILONREPLCEDWITHCEPSILON}, 
	and \eqref{E:VPSIBOUNDEDBYEUCRINVERSETIMESEUCROTATIONS}.

	
	

\end{proof}

We now use the above estimates to derive the main estimate of interest,
namely \eqref{E:FIRSTSTATEMENTHARDSHARPRADPSIPOINTWISEESTIMATE}.

\begin{proposition}[\textbf{Sharp pointwise estimates for $H_{\Lunit \Lunit}^{\Lunit} \Rad \Psi_{\nu}$ in terms of $\Lunit \upmu$}]
	\label{P:EACHRADPSINUFACTORISESSENTIALLYTWICELUPMU}
	Under the small-data and bootstrap assumptions 
	of Sects.~\ref{S:PSISOLVES}-\ref{S:C0BOUNDBOOTSTRAP}
	(with $\vec{\Psi}$ in place of $\Psi$),
	if $\varepsilon$ is sufficiently small, 
	then the following estimates hold on the spacetime region $\mathcal{M}_{\Tboot,U_0}$
	($\nu=0,1,2,3$):
	\begin{align} \label{E:HARDSHARPRADPSIPOINTWISEESTIMATE}
		\left|
			H_{\Lunit \Lunit}^{\Lunit} \Rad \Psi_{\nu}
		\right|
		& \leq
			2
			\left|
				\Lunit \upmu
			\right|
			+ C \varepsilon \frac{\ln(\myexp + t)}{(1 + t)^2}.
	\end{align}
	
	\end{proposition}

\begin{proof}
	We consider equation \eqref{E:NEWUPMUTRANSPORT}. Using an argument similar 
	to the one we used to prove \eqref{E:KEYGLLRADPSIOVERUPMUTIMESRADZNPSIBOUND}, we deduce that
	$|\upomega^{(Tan-\vec{\Psi})}| \lesssim \varepsilon (1 + t)^{-2}$
	and $|H_{\Lunit \Lunit}^{\Lunit}| \lesssim 1.$
	Using these two inequalities and inequality \eqref{E:RADPSI0ISNEARRADIALCONTRACTIONRADPSI}
	and examining the expression \eqref{E:LITTLEOMEGATRANSVERSALPART} for $\upomega^{(Trans-\vec{\Psi})},$
	we immediately deduce the desired estimate \eqref{E:HARDSHARPRADPSIPOINTWISEESTIMATE} in the case $\nu = 0.$	
	
	To prove the desired estimate \eqref{E:HARDSHARPRADPSIPOINTWISEESTIMATE} in the cases $\nu = 1,2,3,$
	we use the following inequalities:
	\begin{align} \label{E:MINKOWSKITCOORDINATEXIOVEREUCLIDEANRSQUAREDPOINTWISE}
	\left|
			\frac{t x^i}{r^2}
		\right|
		& \leq
			\left|
				\frac{t}{\rgeo}
			\right|
			\left|
				\frac{\rgeo}{r}
			\right|
		\leq 1 + \frac{\ln(\myexp + t)}{1 + t}.
	\end{align}
	The first inequality in \eqref{E:MINKOWSKITCOORDINATEXIOVEREUCLIDEANRSQUAREDPOINTWISE} is trivial,
	while the second follows from \eqref{E:EASYEUCLIDEANRADIALVARIABLEBOUND},
	the definition $\rgeo = 1 - u + t,$ and the bound $0 \leq u \leq U_0 < 1.$
	Using inequalities 
	\eqref{E:RADPSI0ISNEARRADIALCONTRACTIONRADPSI},
	\eqref{E:RADPSIIISNEARWEIGHTEDRADIALCONTRACTIONRADSPI},
	and \eqref{E:MINKOWSKITCOORDINATEXIOVEREUCLIDEANRSQUAREDPOINTWISE},
	the estimate $|H_{\Lunit \Lunit}^{\Lunit}| \lesssim 1$ mentioned above,
	and the bound $|\Radunit^a \Rad \Psi_a| \lesssim \varepsilon (1 + t)^{-1}$
	(which follows from the bound $|\Radunit^a| \lesssim 1$ proved above
	and the bootstrap assumption $\|\Rad \Psi_i \|_{C^0(\Sigma_t^u)} \leq \varepsilon (1 + t)^{-1}$),
	we deduce that
	\begin{align} \label{E:HLLLRADPSIIINTERMSOFHLLRADPSI0PLUSJUNK}
		\left|
			H_{\Lunit \Lunit}^{\Lunit} \Rad \Psi_i
		\right|
		\leq 
		\left|
			H_{\Lunit \Lunit}^{\Lunit} \Rad \Psi_0
		\right|
		+ C \varepsilon \frac{\ln(\myexp + t)}{(1 + t)^2}.
	\end{align}
	The desired estimate now follows from \eqref{E:HLLLRADPSIIINTERMSOFHLLRADPSI0PLUSJUNK}
	and the estimate \eqref{E:HARDSHARPRADPSIPOINTWISEESTIMATE} in the case $\nu = 0.$
	
\end{proof}


\chapter{Summary of Notation} \label{A:NOTATION}
\thispagestyle{fancy}
In Appendix \ref{A:NOTATION}, we collect some important notation 
and conventions that we use throughout the monograph
so that the reader can refer to it as needed.

\section{Coordinates}
\begin{itemize}
	\item $x^0,x^1,x^2,x^3$ denote the rectangular spacetime coordinates.
	\item $x^1,x^2,x^3$ denote the rectangular spatial coordinates.
	\item We often use the alternate notation $t = x^0.$
	\item $(t,u,\vartheta^1,\vartheta^2)$ are the geometric coordinates
		(where $t$ is the rectangular time coordinate and $u$ is the eikonal function).
	\item $r = \sqrt{\sum_{a=1}^3 (x^a)^2}$ is the Euclidean radial variable.
	\item $\theta = (\theta^1,\theta^2)$ are local angular coordinates 
		corresponding to standard spherical coordinates on Minkowski spacetime.
	\item $\rgeo = 1 - u + t$ is the geometric radial variable.
\end{itemize}

\section{Indices}
Lowercase Greek indices $\mu, \nu,$ etc. correspond to components with respect
to the \emph{rectangular spacetime} coordinates $x^0,x^1,x^2,x^3,$ and lowercase
Latin indices $i,j,$ etc. correspond to components with respect
to the \emph{rectangular spatial} coordinates $x^1,x^2,x^3.$ 
That is, lowercase Greek indices vary over $0,1,2,3$ and lowercase
Latin indices vary over $1,2,3.$
All lowercase Greek indices are lowered and raised with the spacetime metric
$g$ and its inverse $g^{-1},$ and \emph{not with the Minkowski metric}.
Capital Latin letters $A,B$ etc. correspond to components with respect to the 
geometric coordinate derivative frame 
$\lbrace X_1 = \frac{\partial}{\partial \vartheta^1} ,X_2 = \frac{\partial}{\partial \vartheta^2} \rbrace$ 
on the spheres $S_{t,u}.$ 
These indices are lowered and raised with $\gsphere$ and $\ginversesphere.$
We use Einstein's summation convention in that repeated indices are summed 
over their respective ranges.

\section{Constants}
\begin{itemize}
\item $C$ denotes a uniform constant that is free to vary from line to line.
\item If we want to emphasize that the constant $C$ depends on an a quantity
	$Q,$ then we use notation such as ``$C_Q.$''
\item We use the notation 
\[
	f_1 \lesssim f_2
\]
to indicate that there exists a uniform constant $C > 0$ such that
$f_1 \leq C f_2.$ 
We sometimes use the alternate notation $\mathcal{O}(f_2)$ to denote a quantity $f_1$
that verifies $f_1 \lesssim f_2.$
\item If we want to emphasize that the implicit constant $C$ depends on an a quantity
$Q,$ then we use the alternate notation
\[
	f_1 \overset{Q}{\lesssim} f_2.
\]
\end{itemize}

\section{Spacetime subsets}
\begin{itemize}
	\item $\Sigma_t = \lbrace (x^0,x^1,x^2,x^3) \ | \ x^0 = t \rbrace;$ pg. \pageref{E:SIGMAT}.
	\item $\mathcal{C}_u =$ the outgoing null hypersurface equal to the corresponding level set of the eikonal function;
		pg. \pageref{E:CUT}.
	\item $S_{t,u}=$ a topological two-dimensional sphere equal to $\Sigma_t \cap \mathcal{C}_u;$
		pg. \pageref{E:STU}.
	\item $\Sigma_t^u = $ the portion of $\Sigma_t$ in between $\mathcal{C}_0$ and $\mathcal{C}_u;$
		pg. \pageref{E:SIGMATU}.
	\item $\mathcal{C}_u^t = $ the portion of $\mathcal{C}_u$ in between $\Sigma_0$ and $\Sigma_t;$
		pg. \pageref{E:CUT}.
	\item $\mathcal{M}_{t,u} =$ the spacetime region trapped in between
		$\Sigma_0^u,$ $\Sigma_t^u,$ $\mathcal{C}_u^t,$ and $\mathcal{C}_0^t;$
		pg. \pageref{E:MTUDEF}.
\end{itemize}

\section{Metrics and musical notation}

\begin{itemize}
	\item $g = g(\Psi)$ denotes the spacetime metric.
	\item Relative to rectangular coordinates,
		$g_{\mu \nu} = m_{\mu \nu} + g_{\mu \nu}^{(Small)},$
		where $m_{\mu \nu} = \mbox{diag}(-1,1,1,1);$
		pg. \pageref{E:LITTLEGDECOMPOSED}.
	\item Relative to rectangular coordinates,
		$G_{\mu \nu}(\Psi) = \frac{d}{d \Psi} g_{\mu \nu}(\Psi)$
		and
		$G_{\mu \nu}'(\Psi) = \frac{d}{d \Psi} G_{\mu \nu}(\Psi);$
 		pg. \pageref{E:BIGGDEF}.
	\item $\gt$ denotes the first fundamental form of $\Sigma_t,$ that is, $\gt_{ij} = g_{ij};$
		pg. \pageref{D:FIRSTFUND}.
	\item $\gt^{-1}$ denotes the inverse first fundamental form of $\Sigma_t;$
		pg. \pageref{D:FIRSTFUND}. 
	\item $\gsphere$ denotes the first fundamental form of $S_{t,u};$
		pg. \pageref{D:FIRSTFUND}.
	\item $\gsphere^{-1}$ denotes the inverse first fundamental form of $S_{t,u};$
		pg. \pageref{D:FIRSTFUND}.
	\item We denote the $\gsphere-$dual of a $S_{t,u}$ one-form $\xi$
		by $\xi^{\#}.$ 
		In particular, $(\xi^{\#})^A = (\gsphere^{-1})^{AB} \xi_B.$
		We often use the abbreviated notation
		$\xi^A := (\gsphere^{-1})^{AB} \xi_B,$
		which is the standard convention of differential geometry. 
		Similarly, $Y$ is an $S_{t,u}-$tangent vector, then
		$Y_{\flat}$ denotes the $\gsphere-$dual of $Y,$
		which is an $S_{t,u}$ covector. That is,
		$(Y_{\flat})_A = (\ginversesphere)_{AB} Y^B.$
		We often abbreviate $Y_A := (Y_{\flat})_A.$
		Similarly, we denote the $\gsphere-$dual of a symmetric type $\binom{0}{2}$ tensor $\xi$
		by $\xi^{\#}$ 
		and we denote its double $\gsphere-$dual
		by $\xi^{\# \#}.$
		In particular, 
		$(\xi^{\#})_{\ B}^A = (\gsphere^{-1})^{AC} \xi_{CB}$
		and 
		$(\xi^{\# \#})^{AB} = (\gsphere^{-1})^{AC} (\gsphere^{-1})^{BD} \xi_{CD}.$
		We use similar notation to denote the $\gsphere-$duals of general type 
		$\binom{m}{0}$ and type $\binom{0}{n}$ $S_{t,u}$ tensors, and we use
		abbreviations similar to the ones mentioned above for vectors and covectors;
		pg. \pageref{D:SHARPNOTATION}.
\end{itemize}

\section{Eikonal function quantities}
	\begin{itemize}
		\item The eikonal function $u$ verifies the eikonal equation $(g^{-1})^{\alpha \beta}\partial_{\alpha} u \partial_{\beta} u = 0$
			and has initial conditions $u|_{t=0} 1 - r;$ pg. \pageref{E:OUTGOINGEIKONAL}.
		\item $\upmu = - \left\lbrace (g^{-1})^{\alpha \beta}\partial_{\alpha} u \partial_{\beta} t \right\rbrace^{-1}$ 
			denotes the inverse foliation density; pg. \pageref{E:UPMUDEF}.
		\item $\Lgeo^{\nu} = - (g^{-1})^{\nu \alpha} \partial_{\alpha} u$ denotes the rectangular components of the outgoing
			null geodesic vectorfield; pg. \pageref{E:LGEOEQUATION}.
		\item $\Lunit^{\nu} = \upmu \Lgeo^{\nu}$ denotes a rescaled outgoing null vectorfield; pg. \pageref{E:LUNIT}.
		\item $\Lunit_{(Small)}^i = \Lunit^i - \rgeo^{-1} x^i$ is a re-centered version of the rectangular component $\Lunit^i;$
			pg. \pageref{E:LUNITJUNK}.
		\item $\Radunit_{(Small)}^i = \Radunit^i + \rgeo^{-1} x^i$ is a re-centered version of the rectangular component $\Radunit^i;$
			pg. \pageref{E:RADUNITJUNK}.
		\item $\upchi_{AB} = \frac{1}{2} \angLie_{\Lunit} \gsphere_{AB}$
			is the null second fundamental form of $\mathcal{C}_u$
			relative to $g;$
			pg. \pageref{E:CHIDEF}.
		\item $\upchi_{AB}^{(Small)} = \upchi_{AB} - \rgeo^{-1} \gsphere_{AB}$ is a re-centered version of $\upchi_{AB};$
			pg. \pageref{E:CHIJUNKDEF}.
	\end{itemize}

\section{Additional \texorpdfstring{$S_{t,u}$}{spherical} tensorfields related to the frame connection coefficients}
\begin{itemize}
	\item $k= \frac{1}{2} \SigmatLie_{\Timenormal} \gt$ 
		is the second fundamental form of $\Sigma_t$ relative to $g;$ pg. \pageref{E:SECONDFUNDSIGMATDEF}.
	\item $\upzeta_A = g(\D_A \Lunit, \Radunit) = \upmu^{-1} g(\D_A \Lunit, \Rad);$
		pg. \pageref{E:ZETADEF}.
	\item $\upzeta_A = \upmu^{-1} \upzeta_A^{(Trans-\Psi)} + \upzeta_A^{(Tan-\Psi)};$
		pg. \pageref{E:ZETADECOMPOSED}.
	\item $\angkdoublearg{A}{B} = \upmu^{-1} \angktriplearg{A}{B}{(Trans-\Psi)}
		+ \angktriplearg{A}{B}{(Tan-\Psi)};$
		pg. \pageref{E:ANGKDECOMPOSED}.
	\item $\upeta_A = \upzeta_A + \upmu^{-1} \angdiff_A \upmu;$
		pg. \pageref{E:UPETATENSOR}.
	\item $\underline{\upchi}_{AB} = 2 \upmu \angktriplearg{A}{B} - \upmu \upchi_{AB};$ pg. \pageref{E:UCHI}.
\end{itemize}

\section{Vectorfields}
	\begin{itemize}
		\item $\Lunit$ is $\mathcal{C}_u^t-$tangent, outward pointing, 
			and verifies $g(\Lunit,\Lunit) = 0$ and $\Lunit t = 1.$
		\item $\Lunit = \frac{\partial}{\partial t}$ relative to the geometric coordinates;
			pg. \pageref{E:LISDDT}.
		\item $\Rad$ is $\Sigma_t^u-$tangent, $S_{t,u}-$orthogonal,
			inward pointing, and verifies $g(\Rad,\Rad) = \upmu^2,$ $\Rad u = 1;$
			pg. \pageref{E:RADINTERMSOFLANDULGOOD}.			
		\item $\Rad = \frac{\partial}{\partial u} - \Xi$ relative to the geometric coordinates,
			where $\Xi$ is $S_{t,u}-$tangent;
			pg. \pageref{E:RADINTERMSOFGEOMETRICCOORDINATEPARTIALDERIVATIVES}.
		\item $\Radunit = \upmu^{-1} \Rad;$
			pg. \pageref{E:RADUNIT}. 
		\item $\uLgood = \upmu \Lunit + 2 \Rad$ is an inward pointing null
			vectorfield that verifies $g(\Lunit, \uLgood) = - 2 \upmu;$
			pg. \pageref{D:NULLVECTORFIELDS}.
		\item For $A=1,2,$ $X_A = \frac{\partial}{\partial \vartheta^A}$ 
			is the geometric angular coordinate derivative vectorfield;
			pg. \pageref{E:STUFRAME}. 
		\item $\Timenormal = \Lunit + \Radunit$ is the 
			future-directed normal to $\Sigma_t;$ pg. \pageref{E:UNITSIGMATNORMAL}.
		\item $\lbrace \Lunit, \Radunit, X_1, X_2 \rbrace$
			denotes the non-rescaled frame.
		\item $\lbrace \Lunit, \Rad, X_1, X_2 \rbrace$
			denotes the rescaled frame.	
		\item $\lbrace \Lunit, \uLgood, X_1, X_2 \rbrace$
			denotes the rescaled null frame.
	\end{itemize}

\section{Contraction and component notation}
	\begin{itemize}
		\item If $\xi$ is a symmetric type $\binom{0}{2}$ spacetime tensor and 
			$V,$ $W$ are
			vectors, then 
		$\xi_{VW} = \xi_{\alpha \beta}V^{\alpha} W^{\beta}$
			and $\xi_V$ is the one-form with rectangular components
			$(\xi_V)_{\nu} = \xi_{\alpha \nu}V^{\alpha};$
			pg. \pageref{S:CONTRACTIONNOTATION}.
	\item If $\xi$ is a symmetric type $\binom{2}{0}$ spacetime tensor,
		then $\xi_{VW} = \xi^{\alpha \beta} V_{\alpha} W_{\beta}.$
	\item We use similar contraction notation for tensors $\xi$ of any type.
	\item We use abbreviations such as $\xi_A = \xi_{X_A}$ when contracting against the $S_{t,u}$ frame
		vectors $X_1$ and $X_2.$
	\item We sometimes use a ``$\cdot$'' to emphasize contraction between two tensors.
		For example, if $T$ is a symmetric type $\binom{2}{0}$ $S_{t,u}$ tensorfield and $\xi$ is a $S_{t,u}$ one-form, 
		then $(\angdiv T) \cdot \xi = (\angDarg{A} T^{AB}) \xi_B.$
	\item If $V$ is a spacetime vector, then
			$V = V^{\Lunit} \Lunit + V^{\Rad} \Rad + V^A X_A$
		denotes its decomposition relative to the rescaled frame
		$\lbrace \Lunit, \Rad, X_1, X_2 \rbrace.$
	\item We have
		$V^{\Lunit} = - V_{\Lunit} - \upmu^{-1} V_{\Rad},$
		$V^{\Rad} = - \upmu^{-1} V_{\Lunit},$
		and $V^A = (\ginversesphere)^{AB}V_B.$
\end{itemize}

\section{Projection operators and frame components}
	\begin{itemize}
		\item $\Sigmatproject$ denotes the type $\binom{1}{1}$ tensorfield that projects onto $\Sigma_t^u;$
			pg. \pageref{E:SIGMATPROJECTION}.
		\item $\sphereproject$ denotes the type $\binom{1}{1}$ tensorfield that projects onto $S_{t,u};$
			pg. \pageref{E:SPHEREPROJECTION}.
		\item If $\xi$ is a spacetime tensor, then $\angxi = \sphereproject \xi$
			is the projection of $\xi$ onto $S_{t,u};$
			pg. \pageref{D:PROJECTIONOFATENSOR}.
		\item If $\xi$ is a type $\binom{0}{2}$ spacetime tensor, then 
			$\angxiarg{V} = \sphereproject (\xi_V);$
			pg. \pageref{D:PROJECTIONOFATENSOR}.
		\item $G_{(Frame)} = 
			\left(
				G_{\Lunit \Lunit},
				G_{\Lunit \Radunit},
				G_{\Radunit \Radunit},
				\angGarg{\Lunit},
				\angGarg{\Radunit},
				\angG
			\right)$
			is the array of components of $G_{\mu \nu}$
			relative to the non-rescaled frame $\lbrace \Lunit, \Radunit, X_1, X_2 \rbrace;$
			pg. \pageref{D:GFRAMECOMPONENTS}.
		\item $G_{(Frame)}' = 
			\left(
				G_{\Lunit \Lunit}',
				G_{\Lunit \Radunit}',
				G_{\Radunit \Radunit}',
				\angGprimearg{\Lunit},
				\angGprimearg{\Radunit},
				\angGprime
			\right)$
			is the array of components of $G_{\mu \nu}'$
			relative to the non-rescaled frame $\lbrace \Lunit, \Radunit, X_1, X_2 \rbrace;$
			pg. \pageref{D:GFRAMECOMPONENTS}.
	\end{itemize}

\section{Tensor products and the trace and trace-free parts of tensors}
\begin{itemize}
	\item $\myspacetimetr \pi = (g^{-1})^{\alpha \beta} \pi_{\alpha \beta}$
		denotes the $g-$trace of the type $\binom{0}{2}$ spacetime tensor $\pi_{\mu \nu}.$
	\item $\mytr \xi = (\ginversesphere)^{AB} \xi_{AB}$ denotes the $\gsphere-$trace 
		of the type $\binom{0}{2}$ $S_{t,u}$ tensor $\xi_{AB};$
		pg. \pageref{E:TRACEDEF}. 
	\item $\hat{\xi}_{AB} = \xi_{AB} - \frac{1}{2} \mytr \xi \gsphere_{AB}$ denotes the trace-free part 
		of the type $\binom{0}{2}$ $S_{t,u}$ tensor $\xi_{AB};$
		pg. \pageref{E:TRACEFREE}.
	\item $(\xi \otimes \omega)_{AB} = \xi_A \omega_B$ denotes the tensor product of the $S_{t,u}$ one-forms $\xi_A$ and $\omega_A.$
	\item $(\xi \hat{\otimes} \omega)_{AB} = \frac{1}{2}(\xi_A \omega_B + \xi_B \omega_A) - \frac{1}{2} (\ginversesphere)^{CD} \xi_C \omega_D \gsphere_{AB}$
		denotes the symmetrized trace-free part of the type $\binom{0}{2}$ $S_{t,u}$ tensor 
		$(\xi \otimes \omega)_{AB};$ pg. \pageref{E:SYMMETRIZEDTRACEFREETENSORPRODUCT}.
\end{itemize}

\section{Energy-momentum tensor, multiplier vectorfields, and compatible currents}
\begin{itemize}
	\item $\enmomtensor_{\mu \nu}[\Psi] = \D_{\mu} \Psi \D_{\nu} \Psi
	- \frac{1}{2} g_{\mu \nu} (g^{-1})^{\alpha \beta} \D_{\alpha} \Psi \D_{\beta} \Psi$
		denotes the energy-momentum tensor associated to $\Psi;$
		pg. \pageref{E:ENERGYMOMENTUMTENSOR}.
	\item $\Mult = (1 + 2 \upmu) \Lunit + 2 \Rad$ denotes the timelike multiplier vectorfield;
		pg. \pageref{E:DEFINITIONMULT}.
	\item $\Mor = \rgeo^2 \Lunit$ denotes the Morawetz multiplier vectorfield;
		pg. \pageref{E:DEFINITIONMOR}.
	\item $\Jenergycurrent{\Mult}^{\nu}[\Psi] = \enmomtensor_{\ \alpha}^{\nu}[\Psi] \Mult^{\alpha}$
			denotes the compatible current associated to $\Mult;$
			pg. \pageref{E:MULTCURRENTDEF}.
	\item $\Jenergycurrent{\Mor}^{\nu}[\Psi]
					= \enmomtensor_{\ \alpha}^{\nu}[\Psi]\Mor^{\alpha}$
					denotes the compatible current associated to $\Mor;$
				pg. \pageref{E:MORCURRENTDEF}.
	\item $\Jenergycurrent{Correction}^{\nu}[\Psi]
					= \frac{1}{2} 
						\left\lbrace
							\rgeo^2 \mytr \upchi \Psi \D^{\nu} \Psi
							- \frac{1}{2} \Psi^2 \D^{\nu}[\rgeo^2 \mytr \upchi]
					  \right\rbrace$
				denotes a lower-order correction current;
				pg. \pageref{E:MODCURRENTDEF}.
	\item $\Jenergycurrent{\Mor+Correction}^{\nu}[\Psi] = {\Jenergycurrent{\Mor}^{\nu}[\Psi]}	
				+ {\Jenergycurrent{Correction}^{\nu}[\Psi]};$
				pg. \pageref{E:MORPLUSMODCURRENT}.
\end{itemize}

\section{Commutation vectorfields}
\begin{itemize}
	\item For $l=1,2,3,$ 
		$\Roteucarg{l}$ denotes the Euclidean rotation vectorfield
		with rectangular spatial components 
		$\Roteucarg{l}^j = \epsilon_{laj} x^a;$
		pg. \pageref{E:EULCIDEANROTATION}.
	\item $\Rot_{(l)} = \sphereproject \Roteuc$ denotes an $S_{t,u}-$tangent rotation vectorfield;	
		pg. \pageref{E:GEOMETRICROTATION}.
	\item $\mathscr{Z} = \lbrace \rgeo \Lunit, \Rad, \Rot_{(1)}, \Rot_{(2)}, \Rot_{(3)} \rbrace$
		denotes the full set of commutation vectorfields;
		pg. \pageref{E:DEFSETOFCOMMUTATORVECTORFIELDS}.
	\item $\mathscr{S} = \lbrace \Rad, \Rot_{(1)}, \Rot_{(2)}, \Rot_{(3)} \rbrace$
		are the spatial commutators;
		pg. \pageref{E:DEFSETOFSPATIALCOMMUTATORVECTORFIELDS}.
	\item $\mathscr{O} = \lbrace \Rot_{(1)}, \Rot_{(2)}, \Rot_{(3)} \rbrace$
		are the rotation commutators;
		pg. \pageref{E:DEFSETOFROTATIONCOMMUTATORVECTORFIELDS}.
\end{itemize}

\section{Differential operators and commutator notation}
\begin{itemize}
	\item $\partial_{\mu}$ denotes the rectangular coordinate derivative $\frac{\partial}{\partial x^{\mu}}.$
	\item $\frac{\partial}{\partial t},\frac{\partial}{\partial u},\frac{\partial}{\vartheta^1},\frac{\partial}{\vartheta^2}$ 
		denote the geometric coordinate derivatives.
	\item $V f = V^{\alpha} \partial_{\alpha} f$ 
		denotes the $V-$directional derivative of a function $f.$
	\item $d f =$ standard differential of a function $f$ on spacetime.
	\item $\angdiff f = \sphereproject d f,$ where $f$ is a function on spacetime and $\sphereproject$ denotes projection onto $S_{t,u};$
		pg. \pageref{E:ANGDIFFDEF}. Alternatively, $\angdiff f$ can be viewed as the inherent angular differential of a function $f$ on the $S_{t,u}.$
	\item $\D =$ Levi-Civita connection of $g.$
	\item $\nabla = $ Levi-Civita connection of the Minkowski metric $m.$
	\item $\angD =$ Levi-Civita connection of $\gsphere.$
	\item $\angdiffarg{A} = X_A \cdot \angdiff,$ 
	$\D_A = \D_{X_A},$ 
		$\angDarg{A} = \angDarg{X_A},$
		and similarly for other differential operators.
	\item $\angD^2 \xi$ denotes the second $S_{t,u}$ covariant derivative of $\xi.$
	\item $\angLap f =  \mytr \angD^2 f$ denotes the angular Laplacian of $f.$
	\item $\angfreeDsquared f$ denotes the trace-free part of the second $S_{t,u}$ covariant derivative of the function $f;$
		pg. \pageref{E:ANGFREEDDEF}.
	\item If $Y$ is an $S_{t,u}-$tangent vectorfield,
		then $\angdiv Y = \angDarg{A} Y^A$ denotes its angular divergence;
		pg. \pageref{E:ANGDIVVECTORFIELD}.
		Similarly, if $\xi$ is a type $\binom{1}{1}$ $S_{t,u}$ tensorfield, then 
		$\angdiv \xi$ is the $S_{t,u}$ one-form with components
		$(\angdiv \xi)_A
		= \angDarg{B} \xi_{\ A}^B;$
		pg. \pageref{E:ANGDIVYPE11}.
		We similarly define the angular 
		divergence of symmetric type $\binom{2}{0}$
		or symmetric type $\binom{0}{2}$ $S_{t,u}$ tensorfields.
	\item If $\xi$ is a symmetric type $\binom{0}{2}$ $S_{t,u}$ tensorfield, then
			$\angcheckDarg{A} \xi_{BC} = \frac{1}{2} \left\lbrace
				\angDarg{A} \xi_{BC}
			+ \angDarg{B} \xi_{AC}
			- \angDarg{C} \xi_{AB}
			\right\rbrace;$
			pg. \pageref{E:ANGDCHECKDEF}.
	\item $\Lie_V \xi$ denotes the Lie derivative of $\xi$ with respect to $V;$
		pg. \pageref{E:LIEDERIVATIVE}.
	\item $[V,W] = \Lie_V W$ when $V$ and $W$ are vectorfields.
	\item More generally, if $P$ and $Q$ are two operators,
		then $[P,Q]$ denotes their commutator.
	\item $\SigmatLie_V \xi = \Sigmatproject \Lie_V \xi$ is the $\Sigma_t-$projected 
		Lie derivative of $\xi$ with respect to $V;$
		pg. \pageref{E:SIGMATTPROJECTEDLIE}.
	\item $\angLie_V \xi = \sphereproject \Lie_V \xi$ is the $S_{t,u}-$projected 
		Lie derivative of $\xi$ with respect to $V;$
		pg. \pageref{E:STUPROJECTEDLIE}.
	\item If $\xi$ is type $\binom{0}{2}$ tensorfield, then
		$\angfreeLiearg{V} \xi$ denotes the trace-free $S_{t,u}-$projected Lie derivative operator;
		pg. \pageref{E:ANGFREELIEDEF}.
\end{itemize}

\section{Floor and ceiling functions and repeated differentiation}
\begin{itemize}
\item If $M$ is a non-negative integer,
	then $\lfloor M/2 \rfloor = M/2$ for $M$ even and 
	$\lfloor M/2 \rfloor = (M-1)/2$ for $M$ odd,
	while
	$\lceil M/2 \rceil = M/2$ for $M$ even and 
	$\lceil M/2 \rceil = (M+1)/2$ for $M$ odd;
	pg. \pageref{R:FLOORFUNCTION}.
\item We fix a labeling $\mathscr{Z} = \lbrace Z_{(i)} \rbrace_{i=1}^5$ of the commutation vectorfields.
		If $\vec{I} = (\iota_{(1)}, \iota_{(2)}, \cdots, \iota_{(N)})$ is a multi-index
		of length $N$ with $\iota_{(1)}, \iota_{(2)}, \cdots, \iota_{(N)} \in \lbrace 1,2,3,4,5 \rbrace,$
		then $\mathscr{Z}^{\vec{I}} = Z_{\iota_{(1)}} Z_{\iota_{(2)}} \cdots Z_{\iota_{(N)}}$ 
		denotes the corresponding $N^{th}$ order differential operator.
	\item When we are not concerned with the precise structure of 
		the multi-index $\vec{I},$ we abbreviate $\mathscr{Z}^N = Z_{\iota_{(1)}} Z_{\iota_{(2)}} \cdots Z_{\iota_{(N)}}.$ 
	\item Similarly, $\angLie_{\mathscr{Z}}^N = \angLie_{Z_{\iota_{(1)}}} \angLie_{Z_{\iota_{(2)}}} \cdots \angLie_{Z_{\iota_{(N)}}}$	
		denotes an $N^{th}$ order $S_{t,u}-$projected Lie derivative operator.
		\item Similarly, $\angfreeLietwoarg{\mathscr{Z}}{N} = \angfreeLiearg{Z_{\iota_{(1)}}} \angfreeLiearg{Z_{\iota_{(2)}}} \cdots \angfreeLiearg
		{Z_{\iota_{(N)}}}$	
		denotes an $N^{th}$ order trace-free $S_{t,u}-$projected Lie derivative operator.
	\item If $f$ is a function, then
			$\left|\mathscr{Z}^{\leq N} f \right|$ denotes any term that is
			$\leq \sum_{|\vec{I}| \leq N} c_{\vec{I}} \left| \mathscr{Z}^{\vec{I}} f \right|,$
			where the $c_{\vec{I}}$ are non-negative constants that verify
			$\sum_{|\vec{I}| \leq N} c_{\vec{I}} \leq 1.$
			Similarly, $\left|\mathscr{Z}^N f \right|$ denotes any term that is
			$\leq \sum_{|\vec{I}| = N} c_{\vec{I}} \left| \mathscr{Z}^{\vec{I}} f \right|,$
			where the $c_{\vec{I}}$ are non-negative constants that verify
			$\sum_{|\vec{I}| = N} c_{\vec{I}} \leq 1.$
			Similarly, if $V$ is a vectorfield, then 
			$\left|V \mathscr{Z}^{\leq N} f \right|$ 
			denotes any term that is
			$\leq \sum_{|\vec{I}| \leq N} c_{\vec{I}} \left| V \mathscr{Z}^{\vec{I}} f \right|,$
			where the $c_{\vec{I}}$ are non-negative constants that verify
			$\sum_{|\vec{I}| \leq N} c_{\vec{I}} \leq 1.$
			Similarly,  
			$\left|V \mathscr{Z}^N f \right|$ 
			denotes any term that is
			$\leq \sum_{|\vec{I}| = N} c_{\vec{I}} \left| V \mathscr{Z}^{\vec{I}} f \right|,$
			where the $c_{\vec{I}}$ are non-negative constants that verify
			$\sum_{\vec{I}} c_{\vec{I}} \leq 1.$
			Similarly,
			$\left|\angdiff \mathscr{Z}^{\leq N} f \right|$ 
			denotes any term that is
			$\leq \sum_{|\vec{I}| \leq N} c_{\vec{I}} \left| \angdiff \mathscr{Z}^{\vec{I}} f \right|,$
			where the $c_{\vec{I}}$ are non-negative constants that verify
			$\sum_{\vec{I}} c_{\vec{I}} \leq 1.$
			We use similar notation for other expressions.
		\item If $\xi$ is an $S_{t,u}$ tensor, then 
			$\left| \angLie_{\mathscr{Z}}^{\leq N} \xi \right|$ 
			denotes any term that is
			$\leq \sum_{|\vec{I}| \leq N} c_{\vec{I}} \left| \angLie_{\mathscr{Z}}^{\vec{I}} \xi \right|,$
			where the $c_{\vec{I}}$ are non-negative constants that verify
			$\sum_{|\vec{I}| \leq N} c_{\vec{I}} \leq 1.$
			Similarly, $\left| \angLie_{\mathscr{Z}}^N \xi \right|$ 
			denotes any term that is
			$\leq \sum_{|\vec{I}| = N} c_{\vec{I}} \left| \angLie_{\mathscr{Z}}^{\vec{I}} \xi \right|,$
			where the $c_{\vec{I}}$ are non-negative constants that verify
			$\sum_|{\vec{I}| = N} c_{\vec{I}} \leq 1.$
			Similarly, if $V$ is a vectorfield, then
			$\left| \angLie_V \angLie_{\mathscr{Z}}^{\leq N} \xi \right|$ 
			denotes any term that is
			$\leq \sum_{|\vec{I}| \leq N} c_{\vec{I}} \left| \angLie_V \angLie_{\mathscr{Z}}^{\vec{I}} \xi \right|,$
			where the $c_{\vec{I}}$ are non-negative constants that verify
			$\sum_{|\vec{I}| \leq N} c_{\vec{I}} \leq 1.$
			Similarly,
			$\left| \angLie_V \angLie_{\mathscr{Z}}^N \xi \right|$ 
			denotes any term that is
			$\leq \sum_{|\vec{I}| = N} c_{\vec{I}} \left| \angLie_V \angLie_{\mathscr{Z}}^{\vec{I}} \xi \right|,$
			where the $c_{\vec{I}}$ are non-negative constants that verify
			$\sum_{|\vec{I}| = N} c_{\vec{I}} \leq 1.$
			We use similar notation for other expressions,
			including the case in which trace-free $S_{t,u}-$projected 
			Lie derivatives $\angfreeLie$ are present instead
			of ordinary $S_{t,u}-$projected 
			Lie derivatives $\angLie.$
		\item We use similar notation with $\mathscr{S}$ or $\mathscr{O}$ in place of $\mathscr{Z}$ 
			when all derivatives are spatial derivatives or rotation derivatives.
\end{itemize}

\section{Area and volume forms}
	\begin{itemize}
		\item $d \argspherevol{(t',u',\vartheta)} = \sqrt{\mbox{\upshape{det}} \gsphere(t',u',\vartheta)} d \vartheta$
			denotes the area form on $S_{t,u}$ induced by $g;$
			pg. \pageref{E:SPHEREAREAFORM}.
		\item $d \tvol = d \argspherevol{(t',u',\vartheta)} \, du'$
			denotes a volume form on $\Sigma_t^{U_0};$
			$\upmu \, d \tvol$ is the volume form on $\Sigma_t^u$
				induced by $g;$
        pg. \pageref{E:RESCALEDSIGMATVOLUMEFORM}.			
		\item $d \conevol = d \argspherevol{(t',u',\vartheta)} \, dt'$
			denotes a volume form on $\mathcal{C}_u^t;$
			pg. \pageref{E:RESCALEDCONEVOLUMEFORM}.
		\item $d \vol = d \argspherevol{(t',u',\vartheta)} \, du' \, dt'$
			denotes a volume form on $\mathcal{M}_{t,u};$
			$\upmu \, d \vol$ is the volume form on $\mathcal{M}_{t,u}$
			induced by $g;$
			pg. \pageref{E:RESCALEDSPACETIMEVOLUMEFORM}.
	\end{itemize}

\section{Norms}
	\begin{itemize}
		\item 
			$|\xi|^2 
				= 
				(\ginversesphere)^{B_1 \widetilde{B}_1} \cdots (\ginversesphere)^{B_n \widetilde{B}_n}
				\gsphere_{A_1 \widetilde{A}_1} \cdots \gsphere_{A_m \widetilde{A}_m}
				\xi_{B_1 \cdots B_n}^{A_1 \cdots A_m}
				\xi_{\widetilde{B}_1 \cdots \widetilde{B}_n}^{\widetilde{A}_1 \cdots \widetilde{A}_m}
			$ 
			denotes the square of the norm of the type $\binom{m}{n}$ $S_{t,u}$ tensor $\xi;$
			pg. \pageref{E:POINTWISENORMS}.
		\item $\| f \|_{C^0(\Omega)} = \sup_{p \in \Omega} |f(p)|;$
			pg. \pageref{E:C0NORM}.
		\item Throughout the monograph,
			whenever we state an estimate
			implying that $\| f \|_{C^0(\Omega)} < \infty,$
			we are also implicitly indicating that $f \in C^0(\Omega),$
			where $C^0(\Omega)$ is the set of functions 
			of the geometric coordinates $(t,u,\vartheta)$
			that are continuous on $\Omega.$
		\item $\| f \|_{L^2(S_{t,u})}^2 = \int_{\mathbb{S}^2} f^2(t,u,\vartheta) d \spherevol
			= \int_{S_{t,u}} f^2 d \spherevol;$
			pg. \pageref{E:STUL2NORMDEF}.
		\item $\| f \|_{L^2(\Sigma_t^u)}^2
			= \int_{u'=0}^u \int_{\mathbb{S}^2} f^2(t,u',\vartheta) d \spherevol \, du'
			= \int_{\Sigma_t^u} f^2 \, d \tvol;$
			pg. \pageref{E:SIGMATL2NORMDEF}.
		\item $\| f \|_{L^2(\mathcal{C}_u^t)}^2
			= \int_{t'=0}^t \int_{\mathbb{S}^2} f^2(t',u,\vartheta) d \spherevol \, dt'
			= \int_{\mathcal{C}_u^t} f^2 \, d \conevol;$
			pg. \pageref{E:CONEL2NORMDEF}.
		\item $\| f \|_{\mathcal{M}_{t,u}}^2
			= \int_{t'=0}^t \int_{u'=0}^u \int_{\mathbb{S}^2} f^2(t',u',\vartheta) d \spherevol \, du' \, dt'
			= \int_{\mathcal{M}_{t,u}} f^2 \, d \vol;$
			pg. \pageref{E:SPACETIMEL2NORMDEF}.
		\item We use similar notation for the norms $\| \cdot \|_{L^2(\Omega)}$ of functions defined on subsets $\Omega$
			of $S_{t,u},$ $\Sigma_t^u,$ $\mathcal{C}_u^t,$ or $\mathcal{M}_{t,u}.$
	\end{itemize}

\section{\texorpdfstring{$L^2-$}{square integral-}controlling quantities}
\begin{itemize}
	\item $\enzero[\Psi](t,u)$ denotes the energy of $\Psi$ along $\Sigma_t^u$ corresponding the multiplier $\Mult;$
		pg. \pageref{E:E0DEF}.
	\item $\flzero[\Psi](t,u)$ denotes the flux of $\Psi$ along $\mathcal{C}_u^t$ corresponding the multiplier $\Mult;$
		pg. \pageref{E:F0DEF}.
	\item $\enone[\Psi](t,u)$ denotes the energy of $\Psi$ along $\Sigma_t^u$ corresponding the multiplier $\Mor,$
		including a lower-order correction term; pg. \pageref{E:E1DEF}.
	\item $\flone[\Psi](t,u)$ denotes the flux of $\Psi$ along $\mathcal{C}_u^t$ corresponding the multiplier $\Mor,$
		including a lower-order correction term; pg. \pageref{E:F1DEF}.
	\item $\totzeromax{N}(t,u) = \max_{|\vec{I}| = N} \sup_{(t',u') \in [0,t] \times [0,u]} 
		\left\lbrace
			\enzero[\mathscr{Z}^{\vec{I}} \Psi](t',u')	
			+ \flzero[\mathscr{Z}^{\vec{I}} \Psi](t',u')
		\right\rbrace;$
		pg. \pageref{E:Q0LEQNDEF}.
	\item $\totonemax{N}(t,u) = \max_{|\vec{I}| = N} \sup_{(t',u') \in [0,t] \times [0,u]} 
		\left\lbrace
			\enone[\mathscr{Z}^{\vec{I}} \Psi](t',u')	
			+ \flone[\mathscr{Z}^{\vec{I}} \Psi](t',u')
		\right\rbrace;$
		pg. \pageref{E:Q1LEQNDEF}.
	\item $\totzeromax{\leq N}(t,u) = \max_{0 \leq M \leq N} \totzeromax{M}(t,u);$
		pg. \pageref{E:MAXEDQ0LEQNDEF}.
	\item $\totonemax{\leq N}(t,u) = \max_{0 \leq M \leq N} \totonemax{M}(t,u);$
	pg. \pageref{E:MAXEDQ1LEQNDEF}.
	\item $\Morint[\Psi](t,u) =
	 	\frac{1}{2}
	 	\int_{\mathcal{M}_{t,u}}
			\rgeo^2
			[\Lunit \upmu]_{-}
			|\angdiff \Psi|^2
		\, d \vol$
		denotes the coercive Morawetz spacetime integral associated to $\Psi;$
		pg. \pageref{E:COERCIVEMORDEF}.
	\item $\totMormax{N}(t,u) = \max_{|\vec{I}|=N} \Morint[\mathscr{Z}^{\vec{I}} \Psi](t,u);$
		pg. \pageref{E:MAXEDQ0LEQNDEF}.
	\item $\totMormax{\leq N}(t,u) = \max_{0 \leq M \leq N} \totMormax{M}(t,u);$
		pg. \pageref{E:MAXEDQ0LEQNDEF}.
\end{itemize}

\section{Modified quantities}
\begin{itemize}
	\item $\chifullmod = \upmu \mytr \upchi^{(Small)} + \mathfrak{X}$
		is the lowest-order fully modified version of $\mytr \upchi^{(Small)};$ 
		pg. \pageref{E:LOWESTORDERTRANSPORTRENORMALIZEDTRCHIJUNK}.
	\item $\chifullmodarg{\mathscr{S}^N}
		= \upmu \mathscr{S}^N \mytr \upchi^{(Small)} + \mathscr{S}^N \mathfrak{X}$
		is a higher-order fully modified version of 
		$\mytr \upchi^{(Small)};$ 
		pg. \pageref{E:TRANSPORTRENORMALIZEDTRCHIJUNK}.
	\item $\chipartialmod
			= \mytr \upchi^{(Small)} + \chipartialmodinhom$
		is the lowest-order partially modified version of $\mytr \upchi^{(Small)};$
		\pageref{E:LOWESTORDERTRANSPORTPARTIALRENORMALIZEDTRCHIJUNK}.
	\item $\chipartialmodarg{\mathscr{Z}^N}
		= \mathscr{Z}^N \mytr \upchi^{(Small)} 
		+ \chipartialmodinhomarg{\mathscr{Z}^N}$
		is a higher-order partially modified version of $\mytr \upchi^{(Small)};$
		pg. \pageref{E:TRANSPORTPARTIALRENORMALIZEDTRCHIJUNK}.
	\item $\mupartialmod = \angdiff \upmu + \mupartialmodinhom$
		is the lowest-order partially modified version of $\angdiff \upmu;$
		pg. \pageref{E:LOWESTORDERTRANSPORTPARTIALRENORMALIZEDUPMU}.
	\item $\mupartialmodarg{\mathscr{Z}^N}
		= \angdiff \mathscr{Z}^N \upmu 
				+ \mupartialmodinhomarg{\mathscr{Z}^N}$
		is a higher-order partially modified version of $\angdiff \upmu;$
		pg. \pageref{E:TRANSPORTPARTIALRENORMALIZEDUPMU}.
\end{itemize}

\section{Curvatures}
\begin{itemize}
	\item $\Cur_{\mu \nu \alpha \beta} $ is the Riemann curvature tensor of $g;$ 
		pg. \pageref{D:SPACETIMERIEMANN}.
	\item $\Rsphere_{ABCD}$ is the Riemann curvature tensor of $\gsphere;$ 
		pg. \pageref{L:RICCIISGAUSSTIMEMETRIC}.
	\item $\Rsphere_{AB}$ is the Ricci curvature tensor of $\gsphere;$
		pg. \pageref{L:FIRSTVARIATIONOFSCALARCURVATURE}.
	\item $\Gauss $ is the Gaussian curvature of $\gsphere;$
		pg. \pageref{L:RICCIISGAUSSTIMEMETRIC}.
\end{itemize}

\section{Omission of the independent variables in some expressions}
\begin{itemize}
	\item
	Many of our pointwise estimates are stated in the form
	\[
		|f_1| \lesssim h(t,u)|f_2|
	\]
	for some function $h.$
	Unless we otherwise indicate, it is understood that both $f_1$ and $f_2$
	are evaluated at the point with geometric coordinates
	$(t,u,\vartheta).$
	\item Unless we otherwise indicate,
		in integrals $\int_{S_{t,u}} f \, d \spherevol,$
		the integrand $f$ and the area form $d \spherevol$ are viewed
		as functions of $(t,u,\vartheta)$ and $\vartheta$ is the integration variable.
	\item Unless we otherwise indicate,
		in integrals $\int_{\Sigma_t^u} f \, d \tvol,$
		the integrand $f$ and the volume form $d \tvol$ are viewed
		as functions of $(t,u',\vartheta)$ and $(u',\vartheta)$ are the integration variables.
	\item Unless we otherwise indicate,
		in integrals $\int_{\mathcal{C}_u^t} f \, d \conevol,$
		the integrand $f$ and the volume form $d \conevol$ are viewed
		as functions of $(t',u,\vartheta)$ and $(t',\vartheta)$ are the integration variables.
	\item Unless we otherwise indicate,
		in integrals $\int_{\mathcal{M}_{t,u}} f \, d \vol,$
		the integrand  $f$ and the volume form $d \vol$ are viewed
		as functions of $(t',u',\vartheta)$ and
		$(t',u',\vartheta)$ are the integration variables.
\end{itemize}

\section{Functions that drive shock formation}
\label{S:FUNCTOINSTHATDRIVESHOCKFORMATION}
\begin{itemize}
	\item $(\mathring{\Psi}, \mathring{\Psi}_0)$ are data
		on $\Sigma_0^1$ that are compactly supported in the Euclidean unit
		ball centered at the origin;
		pg. \pageref{S:DATA}.
	\item $(\halfdata, \halftimedata)$
		are spherically symmetric data
		on $\Sigma_{-1/2}$ that are compactly supported in the Euclidean
		ball of radius $1/2$ centered at the origin;
		pg. \pageref{L:SPHERICALLYSYMMETRICDATAGIVESHOCKFUNCTIONBOTHSIGNS}.
	\item $\shockfunction[(\mathring{\Psi}, \mathring{\Psi}_0)]$
		is the data-dependent function that drives 
		nearly spherically symmetric 
		small-data shock formation
		pg. \pageref{E:INITIALDATASHOCKFORMATIONFUNCTION}.
	\item $\timeminushalfshockfunction[(\halfdata, \halftimedata)]$
		is the data-dependent function
		that effectively determines 
		the behavior of $\shockfunction[(\mathring{\Psi}, \mathring{\Psi}_0)];$
		pg. \pageref{E:SEXTENDED}.
\end{itemize}


\bibliographystyle{amsalpha}
\bibliography{JBib}

\end{document}